%% file: thesis.tex
\documentclass{SPhdThesis}
\usepackage{fullpage}

\SgSetTitle{Distributed and Stochastic Optimization Methods\\ with Gradient Compression and Local Steps}
\SgSetAuthor{Eduard Gorbunov}
\SgSetYear{2021}
\SgSetDegree{Doctor of Philosophy}
\SgSetDepartment{Phystech School of Applied Mathematics and Informatics}
\SgSetUniversity{Moscow Institute of Physics and Technology}
\SgSetDeclarationDate{XXXX 2021}

\usepackage{float}
\usepackage{amssymb, amsmath, amsthm, latexsym}
\usepackage{fullpage, color}
\usepackage{url}
\usepackage{algorithm}
\usepackage{algpseudocode}
\usepackage{tabularx}
\usepackage{paralist}
\usepackage{mathtools}
\usepackage{nicefrac}

\newcommand{\mytextstyle}{} 

 \usepackage{adjustbox}

\usepackage[utf8]{inputenc} 

\usepackage[T1]{fontenc}    
\usepackage{lmodern}

\usepackage{hyperref}       
\usepackage{url}            
\usepackage{booktabs}       
\usepackage{amsfonts}       
\usepackage{nicefrac}       
\usepackage{microtype}      

\usepackage{graphicx}
\usepackage{grffile}

\usepackage{tabularx}

\usepackage{amssymb, amsmath, amsthm, latexsym}
\usepackage{color}
\usepackage{url}
\usepackage{algorithm}
\usepackage{algpseudocode}
\usepackage{tabularx}
\usepackage{paralist}
\usepackage{mathtools}

\usepackage{colortbl}
\definecolor{bgcolor}{rgb}{0.8,1,1}
\definecolor{bgcolor2}{rgb}{0.8,1,0.8}

\usepackage{makecell}
\newcommand{\myred}[1]{{\color{red}#1}}
\newcommand{\myblue}[1]{{\color{blue}#1}}

\usepackage[numbers]{natbib}
\usepackage{sidecap}

\usepackage[flushleft]{threeparttable} 

\usepackage{caption}
\usepackage{multirow}

\usepackage[dvipsnames]{xcolor}
\usepackage{pifont}
\newcommand{\cmark}{{\color{PineGreen}\ding{51}}}%
\newcommand{\xmark}{{\color{BrickRed}\ding{55}}}%

\newcommand{\Exp}{\mathbb{E}}
\newcommand{\Prob}{\mathbb{P}}
\newcommand{\R}{\mathbb{R}}
\newcommand{\eqdef}{\stackrel{\text{def}}{=}}

\def\<#1,#2>{\left\langle #1,#2\right\rangle}

\usepackage{mdframed} 
\usepackage{thmtools}

\definecolor{shadecolor}{gray}{0.9}
\declaretheoremstyle[
headfont=\normalfont\bfseries,
notefont=\mdseries, notebraces={(}{)},
bodyfont=\normalfont,
postheadspace=0.5em,
spaceabove=1pt,
mdframed={
  skipabove=8pt,
  skipbelow=8pt,
  hidealllines=true,
  backgroundcolor={shadecolor},
  innerleftmargin=4pt,
  innerrightmargin=4pt}
]{shaded}

\declaretheorem[style=shaded,within=section]{definition}
\declaretheorem[style=shaded,sibling=definition]{theorem}
\declaretheorem[style=shaded,sibling=definition]{proposition}
\declaretheorem[style=shaded,sibling=definition]{assumption}
\declaretheorem[style=shaded,sibling=definition]{corollary}

\declaretheorem[style=shaded,sibling=definition]{lemma}
\declaretheorem[style=shaded,sibling=definition]{remark}
\declaretheorem[style=shaded,sibling=definition]{example}

\newcommand{\argmin}{\mathop{\arg\!\min}}

\usepackage[colorinlistoftodos,bordercolor=orange,backgroundcolor=orange!20,linecolor=orange,textsize=scriptsize]{todonotes}

\newcommand{\cC}{{\cal C}}
\newcommand{\cD}{{\cal D}}
\newcommand{\cE}{{\cal E}}

\newcommand{\cL}{{\cal L}}
\newcommand{\cM}{{\cal M}}
\newcommand{\cN}{{\cal N}}
\newcommand{\cO}{{\cal O}}

\newcommand{\cQ}{{\cal Q}}

\newcommand{\Var}{\mathrm{Var}}

\newcommand{\mA}{{\bf A}}

\newcommand{\mI}{{\bf I}}
\newcommand{\mJ}{{\bf J}}

\newcommand{\mM}{{\bf M}}

\newcommand{\mS}{{\bf S}}

\newcommand{\mW}{{\bf W}}
\newcommand{\mX}{{\bf X}}

\newcommand{\EE}{\mathbb{E}}
\newcommand{\PP}{\mathbb{P}}

\newcommand{\prox}{\mathop{\mathrm{prox}}\nolimits}
\newcommand{\proxR}{\prox_{\gamma R}}

\newcommand{\tx}{\tilde{x}}

\newcommand{\psvrg}{q}

\usepackage{xspace}

\newcommand{\algname}[1]{{\tt #1}\xspace}
\newcommand{\algnamex}[1]{{\tt #1}\xspace}

\newcommand{\tA}{\widetilde{A}}
\newcommand{\hA}{\widehat{A}}
\newcommand{\tB}{\widetilde{B}}
\newcommand{\hB}{\widehat{B}}
\newcommand{\tF}{\widetilde{F}}
\newcommand{\hF}{\widehat{F}}
\newcommand{\tD}{\widetilde{D}}
\newcommand{\hD}{\widehat{D}}

\newcommand{\txi}{\widetilde{\xi}}
\newcommand{\oxi}{\overline{\xi}}

\newcommand{\Jac}{{ \bf \nabla F}}
\newcommand{\Proj}{{\bf \Pi}}
\newcommand{\norm}[1]{\left \| #1 \right\|}
\newcommand{\ED}[1]{\mathbb{E}_{\cD}\left[#1\right] }
\newcommand{\Tr}[1]{\mbox{Tr}\left( #1\right)}
\newcommand{\ones}{e}
\newcommand{\one}{{\bf 1}}

\newcommand{\x}{x}

\usepackage{hyperref}

\newcommand{\Mod}[1]{\ \mathrm{mod}\ #1}

\usepackage{accents}
\newlength{\dhatheight}

\graphicspath{{plots/}}

\begin{document}
	\begin{frontmatter}
		\SgAddTitle%
        \input{to_my_parents}
		\input{acknowledgments}%
		\SgAddToc
		\SgAddLof
		\SgAddLot
	\end{frontmatter}
   
	\input{ch1_Intro}
	\input{ch2_sigma_k}

	\input{ch3_ef_sigma_k}

	\input{ch4_local_sigma_k}

	\input{ch5_marina}
	\input{ch6_moshpit}
	
	\SgIncludeBib{biblio}
	
	\appendix
	
	\input{Appendix_Basic_ineqs}
	\input{Appendix_ef_sigma_k}
	\input{Appendix_local_sigma_k}
	\input{Appendix_marina}
	\input{Appendix_moshpit}
\end{document}

%% file: to_my_parents.tex
\thispagestyle{plain}
	\begin{center}
		\vspace*{5.5cm}
		{%
			To my parents
		}
	\end{center}
	\vspace{0.5cm}
\SgIntClearDoublePage

%% file: acknowledgments.tex
\begin{acknowledgments}
I express my deepest gratitude to my supervisors Alexander Gasnikov and\newline Peter Richt\'arik. I have learned a lot from both of you about various aspects of being a researcher. Thank you a lot for your guidance, encouragement, and opportunities that you provided. This all allowed me to realize my potential.

Next, I am grateful to Pavel Dvurechensky for all his help and guidance, especially during the work on our first joint papers.

I thank all my co-authors for their work, fruitful discussions and great impact on my research (in the order of appearance of joint works): Evgeniya~Vorontsova, Dmitry~Kovalev, Elnur~Gasanov, Ahmed~Mohammed, Elena~Chernousova,\newline {Konstantin}~Mishchenko, Martin~Tak\'{a}\v{c}, El~Houcine~Bergou, Darina~Dvinskikh, C\'esar~A.~Uribe, Filip~Hanzely, Adel~Bibi, Ozan~Sener, Sergey~Guz,\newline Maksim~Shirobokov, Egor~Shulgin, Aleksandr~Beznosikov, Marina~Danilova,\newline Dmitry~Makarenko, Alexander~Rogozin, Sergey~Guminov,
Dmitry~Kamzolov,\newline Innokentiy~Shibaev, Konstantin~Burlachenko, Zhize~Li, Max~Ryabinin,\newline Vsevolod~Plokhotnyuk, Gennady~Pekhimenko, {Alexander}~Borzunov, Michael~Diskin, Ilyas Fatkhullin, Igor Sokolov, Gauthier~Gidel, Nicolas Loizou, and Hugo Berard.

I also express my gratitude to Artem Babenko, Francis Bach, Aymeric Dieuleveut, Samuel Horv\'ath, Praneeth Karimireddy, Eric Moulines, Anton Osokin,\newline {Alexander}~Panin, Liudmila Prokhorenkova, Adrien Taylor, and Aleksei~Ustimenko for fruitful discussions. Further, I owe a great thanks to my internship advisor\newline Gauthier Gidel and to Nicolas Loizou, who I actively collaborated with during my internship. I learned a lot from you!

Finally, it is hard to express how much I appreciate all the support I received from MIPT and, in particular, from Andrei M. Raigorodskii and the Phystech School of Applied Mathematics and Informatics.

\end{acknowledgments}

%% file: ch1_Intro.tex
\chapter{Introduction}
In\footnote{The work on this thesis was partially supported by RFBR 19-31-51001 and was partially supported by the Ministry of Science and Higher Education of the Russian Federation (Goszadaniye) 075-00337-20-03, project no. 0714-2020-0005.} this chapter, we give a general introduction with an overview of the developed results in this thesis. All subsequent chapters also contain their own detailed introductions.
\section{Stochastic First-Order Methods}
Stochastic optimization \cite{shapiro2014lectures, lan2020first} is a young but rapidly developing branch of optimization. Stochastic optimization methods are at the heart of various applications of statistics \cite{spokoiny2012parametric} and machine learning \cite{Goodfellow-et-al-2016,shalev2014understanding}. Sometimes the use of stochasticity is dictated by the nature of the optimization problem, in other situations, people artificially introduce stochasticity to solve the problem faster, e.g., in randomized coordinate-wise methods \cite{RCDM, richtarik2020stochastic, kovalev2018stochastic} and stochastic derivative-free approaches \cite{nesterov2017random, gorbunov2018accelerated, dvurechensky2021accelerated, bergou2020stochastic, gorbunov2020stochastic}.

Due to their practical efficiency and simplicity in implementation, stochastic first-order methods are the most popular stochastic optimization methods. The simplest and brightest example of such a method is Stochastic Gradient Descent ({\tt SGD}) \cite{RobbinsMonro:1951}. In its basic form, {\tt SGD} applied to the unconstrained minimization problem
\begin{equation}
    \min\limits_{x\in \R^d} f(x) \label{eq:problem_intro}
\end{equation}
has the update rule
\begin{equation}
    x^{k+1} = x^k - \gamma_k g^k, \label{eq:SGD_intro}
\end{equation}
where $\{x^k\}_{k\ge 0}$ is the sequence of optimization variables, $\{\gamma_k\}_{k\ge 0}$ is the sequence of stepsizes, and $\{g^k\}_{k\ge 0}$ are {\it stochastic gradients} -- the key ingredient in {\tt SGD}. In a nutshell, stochastic gradient $g^k$ is a random vector that, in some sense, approximates the true gradient $\nabla f(x^k)$ of the objective function $f$ at the point $x^k$. Of course, in each particular situation, it should be clarified in what sense $g^k$ approximates $\nabla f(x^k)$. Typically, this means that $g^k$ is an unbiased estimate of $\nabla f(x^k)$ for fixed $x^k$:
\begin{equation}
    \EE\left[g^k\mid x^k\right] = \nabla f(x^k). \label{eq:unbiasedness_intro}
\end{equation}
Although this assumption is natural, it is not enough to ensure the convergence of {\tt SGD} to some solution of the problem \eqref{eq:problem_intro}. Therefore, it is necessary to introduce additional assumptions on the stochastic gradient. Moreover, before that, it needs to be clarified what we mean by the ``convergence'' of a stochastic method.

As in the majority of papers on stochastic optimization, in this thesis, we focus on the convergence in expectation, i.e., we study the convergence rates of the considered methods to achieve a desired accuracy of the solution (in terms of functional suboptimality/squared distance to the solution/squared norm of the gradient) in expectation. In many real-world problems, ``in-expectation'' convergence guarantees are in good correspondence with behavior of the method during a particular run and they are often easier to derive than their high-probability counterparts. However, we emphasize that for a deeper understanding of the stochastic methods, it is also highly important to analyze their high-probability convergence rates \cite{Nemirovski-Juditsky-Lan-Shapiro-2009, ghadimi2012optimal, ghadimi2013optimal, nazin2019algorithms, davis2021low, gorbunov2020clipped_sstm, gorbunov2021near, cutkosky2021high}, as well as limit distributions \cite{polyak1992acceleration, hodgkinson2021multiplicative, gurbuzbalaban2021heavy, wang2021convergence} and almost-surely convergence guarantees \cite{bottou2003stochastic, zhou2017stochastic, nguyen2018sgd, mertikopoulos2020almost, sebbouh2021almost, patel2020stopping}.

The classical convergence guarantees for {\tt SGD} \cite{RobbinsMonro:1951, Nemirovski-Juditsky-Lan-Shapiro-2009} rely on the bounded second moment assumption:
\begin{equation*}
    \EE\left[\|g^k\|^2\mid x^k\right] \le G^2
\end{equation*}
for some constant $G > 0$. Although this assumption is reasonable for convex non-smooth objectives, it does not hold for strongly convex problems and for several smooth convex problems. To resolve this issue for smooth problems, one can analyze \cite{ghadimi2013stochastic} {\tt SGD} assuming only boundedness of the variance:
\begin{equation*}
    \EE\left[\|g^k - \nabla f(x^k)\|^2\mid x^k\right] \le \sigma^2
\end{equation*}
for some $\sigma \ge 0$. Next, if stochastic realizations of the objective function $f$ are smooth, then one can relax this assumption even further \cite{nguyen2018sgd, gower2019sgd}. Moreover, taking into account some structural properties of the problem one can construct $g^k$ in such a way that it will satisfy certain inequalities needed to derive the convergence of the resulting method. For example, in finite-sum optimization, one can consider variance reduced methods \cite{SAG, SVRG, SAGA}, in distributed optimization, one might be interested in designing parallel stochastic methods with communication compression \cite{alistarh2017qsgd, seide20141}, and, when the dimension of the problem is an issue, one can use coordinate-wise randomization \cite{RCDM}.

As a result, a lot of different stochastic methods appeared in the literature and were analyzed under various assumptions. However, a large group of {\tt SGD} methods have update rules of the form \eqref{eq:SGD_intro} with gradient estimates satisfying \eqref{eq:unbiasedness_intro}. Therefore, it is important to have a clean systematic way to analyze all of them, i.e., have a general theoretical framework that provides tight analysis for all of these methods.

\subsubsection{The First Contribution: Unified Theory of SGD}
Our first contribution is a general analysis of {\tt SGD} in the strongly convex case with proximable regularization. That is, we propose a unified assumption on the stochastic gradients and the problem that covers various existing methods in different settings. Whenever we recover a known method, our general theorem provides the tightest know rate for this method. Moreover, inspired by the proposed theoretical framework, we generalize several existing methods and develop new stochastic methods.

Chapter~\ref{ch:sigma_k} is devoted to the first contribution of this thesis and is based on the following paper:
\begin{itemize}
    \item[\cite{gorbunov2019unified}] Eduard Gorbunov, Filip Hanzely, and Peter Richt\'arik. \textit{A Unified Theory of SGD: Variance Reduction, Sampling, Quantization and Coordinate Descent}. In Silvia Chiappa and Roberto Calandra, editors, Proceedings of the Twenty Third International Conference on Artificial Intelligence and Statistics, volume 108 of Proceedings of Machine Learning Research, pages 680–690. PMLR, 26–28 Aug 2020.
\end{itemize}

\section{Centralized Distributed Stochastic Optimization}
As we mentioned earlier, stochastic optimization methods are widely used in machine learning applications. With the growth of data and complexity of models it became inevitable to consider ways of solving the problems in a parallel/distributed way. Indeed, training modern deep neural networks would take a prohibitively long time (e.g., days or even years of computations) if executed on a single machine, even if this machine is a top-of-the-line GPU server \cite{gpt3costlambda}. Therefore, \textit{distributed} stochastic methods are usually applied in such problems \cite{goyal2017accurate, You2020Large}, where parallel computations help to reduce the training time significantly. Moreover, distributed methods are the natural choice when the data is private and/or distributed across multiple devices, e.g., in federated learning \cite{FEDLEARN, FL2017-AISTATS}.

In its general form, distributed unconstrained optimization problem can be defined in the following way: $n$ devices/peers/workers/nodes/machines solve the minimization problem
\begin{equation}
    \min\limits_{x\in\R^d}\left\{f(x) := \frac{1}{n}\sum\limits_{i=1}^n f_i(x)\right\}, \label{eq:distr_problem_intro}
\end{equation}
where function $f_i$ is known for worker $i$ only but not necessarily to other workers, meaning that worker $i$ can compute some specific quantities such as functional value or (stochastic) gradient of $f_i$ but other workers do not necessarily have an access to this information. For example, in federated learning, the information about function $f_i$ is privately stored on device $i$, and $f_1,\ldots,f_n$ are naturally heterogeneous. In large-batch training of deep neural networks, all functions $f_i$ can be equal to $f$.

Perhaps the simplest {\tt SGD} variant for solving \eqref{eq:distr_problem_intro} is Parallel {\tt SGD} \cite{localsgd_first}:
\begin{equation}
    x^{k+1} = x^k - \gamma_k g^k = x^k - \frac{\gamma_k}{n}\sum\limits_{i=1}^n g_i^k, \label{eq:Parallel_SGD_intro}
\end{equation}
where $g_i^k$ is a stochastic gradient of function $f_i$ at point $x^k$. That is, at each iteration of Parallel {\tt SGD}, workers first compute stochastic gradients $g_i^k$, and, after that, vectors $g_i^k$ for $i=1,\ldots,n$ are aggregated and new point $x^{k+1}$ is computed. Here the following natural question arises: \textit{how are the stochastic gradients aggregated?}

The classical and historically first way of gradients aggregation is to use the Parameter Server architecture \cite{parameter_server_first}. In this approach, workers cannot communicate between each other directly, and instead are only allowed to communicate with a dedicated machine: a server or master. Therefore, to update $x^{k+1}$ via \eqref{eq:Parallel_SGD_intro}, workers need to send the gradients $g_i^k$ to the server. After that, the server averages the received vectors, computes $x^{k+1}$, and broadcasts the result back to the workers.

Despite its simplicity, this idea works quite well in practice. However, Parallel {\tt SGD} has a significant issue that rapidly becomes evident with the growth of the number of workers $n$ and/or growth of the dimension of the problem $d$. This issue is called \textit{communication bottleneck}. It means that for large enough $n$ or $d$, communication may take much more time than computation. This happens because of several reasons: 1) stochastic gradients $g_i^k$ can be dense and huge-dimensional, 2) workers communicate at each iteration of the method, and 3) a single machine (server) is responsible for aggregating a large amount of information at each iteration. In this thesis, we address all these three problems separately.

\subsection{Communication Compression}
The natural way of addressing the communication bottleneck is to use \textit{communication compression} \cite{seide20141, suresh2017distributed}, which is based on applying compression to the gradient vectors or tensors that workers need to send to the master. For example, one can modify Parallel {\tt SGD} in the following way \cite{alistarh2017qsgd}:
\begin{equation}
    x^{k+1} = x^k - \frac{\gamma_k}{n}\sum\limits_{i=1}^n \cC(g_i^k), \label{eq:QSGD_intro}
\end{equation}
where $\cC:\R^d \to \R^d$ is some (possibly randomized) operator called \textit{compression operator}. This method is usually called Compressed or Quantized {\tt SGD} ({\tt QSGD}). In this scheme, instead of sending $g_i^k$, the workers send the compressed message $\cC(g_i^k)$ to the server. Therefore, the operator $\cC$ is designed in such a way that transmitting $\cC(g_i^k)$ requires much less time than transmitting $g_i^k$. For example, one can use the so-called Rand$K$ operator that picks $K$ components of the input uniformly at random and scales the result to ensure unbiasedness:
\begin{equation}
    \text{Rand}K(x) = \frac{d}{K}\sum\limits_{i \in S}x_ie_i. \notag
\end{equation}
Here $(e_1,e_2,\ldots, e_d)$ is a standard basis in $\R^d$, $x = (x_1,\ldots,x_d)^\top \in \R^d$, and $S$ is a random set uniformly distributed on the family of $K$-element subsets of $\{1,2,\ldots,d\}$. When $K \ll d$, the per-iteration communication cost of {\tt QSGD} is significantly smaller than for Parallel {\tt SGD}.

Moreover, in the (strongly) convex case, one can prove that {\tt QSGD} converges to the solution with any predefined accuracy if the operator $\cC$ satisfies
\begin{equation}
    \EE\left[\cC(x)\right] = x,\quad \EE\left[\|\cC(x) - x\|^2\right] \le \omega\|x\|^2 \label{eq:unbiased_compr_intro}
\end{equation}
with some $\omega \ge 0$ for all $x\in\R^d$. Compression operators satisfying \eqref{eq:unbiased_compr_intro} are usually called  \textit{unbiased compressors}. Although inequality \eqref{eq:unbiased_compr_intro} is satisfied for a wide range of compression operators, it does not cover several practically important \textit{biased} compression operators such as the Top$K$ compression operator that picks $K$ components of the input with the largest absolute values. Usually, when the compression operator $\cC$ is biased, it is assumed that
\begin{equation}
    \EE\left[\|\cC(x) - x\|^2\right] \le (1-\delta)\|x\|^2 \label{eq:biased_compr_intro}
\end{equation}
with some $\delta \in (0,1]$ for all $x\in \R^d$. Interestingly, Compressed {\tt SGD} \eqref{eq:QSGD_intro} with biased compression $\cC$ may diverge exponentially fast even for strongly convex problems \cite{beznosikov2020biased}. To circumvent this issue, one can use the so-called \textit{error compensation} mechanism \cite{seide20141}. The resulting method is usually called Error Compensated {\tt SGD} ({\tt EC-SGD}), and has the following update rule:
\begin{equation}
    x^{k+1} = x^k - \frac{1}{n}\sum\limits_{i=1}^nv_i^k,\quad v_i^k = \cC(\gamma_k g_i^k + e_i^k),\quad e_i^{k+1} = \gamma_k g_i^k + e_i^k - v_i^k. \label{eq:EC_SGD_intro}
\end{equation}
Here, each worker $i$ ``memorizes'' the unsent information $e_i^{k+1} = \gamma_k g_i^k + e_i^k - v_i^k$ in order to use it during the next iterations.

{\tt EC-SGD} was analyzed in many papers under different assumptions \cite{stich2018sparsified, stich2020error, beznosikov2020biased}. However, before this thesis, there were several important gaps in the theory of stochastic methods with error compensation in the (strongly) convex case. In particular, there were no \textit{full-gradient} methods ($g_i^k = \nabla f_i(x^k)$) with error compensation that have \textit{linear} convergence in the strongly convex case. Moreover, there were no variance reduced variants of {\tt EC-SGD} and variants with arbitrary sampling was never analyzed.

\subsubsection{The Second Contribution: Unified Theory of Error Compensated Methods}
Our second contribution in this thesis can be seen as an extension of the first contribution to the class of methods with error compensation. That is, we propose a new unified theoretical framework for the analysis of stochastic first-order methods supporting \textit{error compensation}. Using this framework, we develop new efficient error-compensated methods. In particular, we develop the first full-gradient methods with error compensation that have linear convergence in the strongly convex case and the first variance reduced method with error compensation that also enjoys linear convergence on strongly convex problems. Moreover, our framework covers methods with \textit{delayed} updates. Overall, using this new framework we develop 16 new optimization methods.

Chapter~\ref{ch:ef_sigma_k} is devoted to the second contribution of this thesis, and is based on the following paper:
\begin{itemize}
    \item[\cite{gorbunov2020linearly}] Eduard Gorbunov, Dmitry Kovalev, Dmitry Makarenko, and Peter Richt\'arik. \textit{Linearly Converging Error Compensated SGD}. In H. Larochelle, M. Ranzato, R. Hadsell, M. F. Balcan, and H. Lin, editors, Advances in Neural Information Processing Systems, volume 33, pages 20889–20900. Curran Associates, Inc., 2020.
\end{itemize}

\subsection{Local Updates}
Another popular way of addressing communication bottleneck is to use more computations locally on workers between two sequential communication rounds. For example, workers can perform several ($\tau \ge 1$) {\tt SGD} steps between two neighboring communications rounds rather than a single ($\tau = 1$) step. Formally, the update of resulting method can be written in the form:
\begin{equation}
    x_i^{k+1} = \begin{cases}x_i^k - \gamma_k g_i^k,& \text{if } k+1 \mod \tau \neq 0,\\ \frac{1}{n}\sum\limits_{i=1}^n(x_i^k - \gamma_k g_i^k),& \text{if } k+1 \mod \tau = 0, \end{cases} \label{eq:Local_SGD_intro}
\end{equation}
where $x_i^k$ denotes the local iterate stored on node $i\in \{1,\ldots,n\}$ at iteration $k$. This method is known as {\tt Local-SGD}/Federated Averaging ({\tt FedAvg}) \cite{FEDLEARN, FL2017-AISTATS, stich2018local}. {\tt Local-SGD} and its different variants gained a lot of attention and were studied in a number of papers \cite{localsgd_first, mcmahan2016federated, stich2018local, LinSPJ2018local,liang2019variance, wu2019federated, karimireddy2019scaffold, khaled2020tighter, woodworth2020local}. However, several promising directions, such as better understanding of so-called \textit{local shifts}, more sophisticated local gradient estimators allowing importance sampling, variance reduction or coordinate descent, variable number of local steps, and general theory supporting different data similarity types, were unexplored in the previous works.

\subsubsection{The Third Contribution: Unified Theory of Methods with Local Updates}
Motivated by the first two contributions, we propose yet another unified theoretical framework, this time for the analysis of {\tt Local-SGD}-type methods, in the regime when the objective function is (strongly) convex. We recover multiple known local optimizers as a special case of our general framework, along with their convergence rates (up to small constant factors). To demonstrate the strengths of our approach we develop a new method called {\tt S-Local-SVRG} fitting our general framework. Moreover, using our general theorem we prove that {\tt S-Local-SVRG} converges \textit{linearly} even when the local loss functions are \textit{arbitrarily heterogeneous}. This is the first variance reduced linearly converging {\tt Local-SGD} method. Moreover, to obtain this result, we did not need to rely on any restrictive assumptions such as gradient boundedness or gradients similarity.

Chapter~\ref{ch:local_sigma_k} is devoted to the third contribution of this thesis and based on the following paper:
\begin{itemize}
    \item[\cite{gorbunov2020local}] Eduard Gorbunov, Filip Hanzely, and Peter Richt\'arik. \textit{Local SGD: Unified Theory and New Efficient Methods}. In Arindam Banerjee and Kenji Fukumizu, editors, Proceedings of The 24th International Conference on Artificial Intelligence and Statistics, volume 130 of Proceedings of Machine Learning Research, pages 3556–3564. PMLR, 13–15 Apr 2021.
\end{itemize}

\subsection{Non-Convex Distributed Optimization with Compression}
In the previous sections, we focus on (quasi-strongly) convex problems. However, there are many practically important problems that are non-convex, including training deep neural networks \cite{Goodfellow-et-al-2016}, and matrix completion and recovery \cite{ma2018implicit, bhojanapalli2016dropping}. Clearly, it is important to design efficient {\tt SGD}-type methods for solving non-convex problems \cite{danilova2020recent}.

Nowadays one of the most popular example of non-convex optimization problems is training of deep neural networks. As we mentioned before, some of these tasks are so computationally hard that even top-of-the-line GPU servers \cite{gpt3costlambda} may require years of computations to solve them. Therefore, such problems are necessarily solved in a distributed manner.

As is the case in the convex regime, communication bottleneck appears in non-convex distributed optimization too, and one can handle this issue using \textit{communication compression}. The optimization and machine learning communities have exerted considerable effort in recent years to design  distributed methods  supporting compressed communication. From the many methods proposed, we emphasize here \algname{VR-DIANA} \cite{horvath2019stochastic},  \algname{FedCOMGATE} \cite{haddadpour2020federated}, and \algname{FedSTEPH} \cite{das2020improved} because they are supported by the state-of-the-art theoretical complexity results in the setup when the local loss functions are allowed to be arbitrarily heterogeneous.

\subsubsection{The Fourth Contribution: Faster Methods for Non-Convex Distributed Optimization with Compression}
We develop and analyze \algname{MARINA}: a new communication efficient method for non-convex distributed learning over heterogeneous datasets.   \algname{MARINA} employs a novel communication compression strategy based on the compression of gradient differences that is reminiscent of but different from the strategy employed in the \algname{DIANA} method \cite{mishchenko2019distributed}. Unlike  virtually all competing distributed first-order methods,  including \algname{DIANA}, ours is based on a carefully designed {\em biased} gradient estimator, which is the key to its superior theoretical and practical performance. The communication complexity bounds we prove for \algname{MARINA} are evidently better than those of all previous first-order methods. Further, we develop and analyze two variants of \algname{MARINA}: \algname{VR-MARINA} and \algname{PP-MARINA}. The first method  is designed for the case when the local loss functions owned by clients are either of a finite sum or of an expectation form, and the second method allows for a partial participation of clients -- a feature important in federated learning. All our methods are superior to previous state-of-the-art methods in terms of oracle/communication complexity. Finally, we provide a convergence analysis of all methods  for problems satisfying the Polyak-{\L}ojasiewicz condition.

Chapter~\ref{ch:marina} is devoted to the fourth contribution of this thesis, and is based on the following paper:
\begin{itemize}
    \item[\cite{pmlr-v139-gorbunov21a}] Eduard Gorbunov, Konstantin P. Burlachenko, Zhize Li, and Peter Richt\'arik. \textit{MARINA: Faster Non-Convex Distributed Learning with Compression}. In Marina Meila and Tong Zhang, editors, Proceedings of the 38th International Conference on Machine Learning, volume 139 of Proceedings of Machine Learning Research, pages 3788–3798. PMLR, 18–24 Jul 2021.
\end{itemize}

\section{Distributed Optimization Without a Central Server}
In situations when it is possible to engineer the \textit{network} that defines communication links among the machines, one can handle communication bottleneck even without compressed communications and local updates. As we explained in the previous sections, in the parameter-server architecture the communication bottleneck arises mainly because of the existence of a machine (server) that aggregates a lot of data at each iteration. To alleviate this issue, one can change the communication protocol in such a way that no machine is required to aggregate too much data at any iteration.

One of the most popular decentralized communication protocols is \textit{gossip} \cite{boyd2006randomized,tsitsiklis1984problems,lian2017can}. For any given network structure, and initial vectors $x_1^0, x_2^0, \ldots, x_n^0 \in \R^d$, gossip generates the sequence of points $\{x_i^{k}\}_{k\ge 0}$ on each worker $i = 1,\ldots, n$ such that
\begin{equation}
    x_i^{k+1} = \sum\limits_{j=1}^n \mM_{ij} x_j^k, \label{eq:gossip_intro}
\end{equation}
where $\mM_{ij}$ is the $i,j$-th element of a \textit{mixing} matrix $\mM$. The key property of a mixing matrix is that $\mM_{i,j} = 0$ iff $i\neq j$ and $(i,j)\not\in \cE$, where $\cE$ denotes the set of edges in the communication network. Further,  for $(i,j) \in \cE$ it satisfies $\mM_{ij} > 0$ and $\mM_{ii} > 0$ for all $i = 1,\ldots,n$. Moreover, it is usually assumed that $\mM$ is symmetric $\mM = \mM^\top$, $\mM \one = \one$, where $\one = (1,\ldots,1)^\top \in \R^n$, and $\lambda_2(\mM) < 1$, where $\lambda_2(\mM)$ is the absolute value of the second largest (in absolute value) eigenvalue of $\mM$ \cite{gorbunov2020recent}. Under these assumptions gossip converges linearly to the exact average of $x_1^0, x_2^0, \ldots, x_n^0$ as follows:
\begin{equation}
    \|\mX^k - \overline{\mX}\|_2 \le (\lambda_2(\mM))^k\|\mX^0 - \overline{\mX}\|_2, \notag
\end{equation}
where $\mX^k = [x_1^k, x_2^k, \ldots , x_n^k]\in \R^{d\times n}$ and $\overline{\mX} = [\overline{x}, \overline{x}, \ldots , \overline{x}]\in \R^{d\times n}$, $\overline{x} = \frac{1}{n}\sum_{i=1}^n x_i^0$. That is, gossip finds approximate average on nodes with accuracy $\|\mX^k - \overline{\mX}\|_2 \le \varepsilon$ after $\cO\left((1-\lambda_2(\mM))^{-1}\log(\varepsilon^{-1})\right)$ iterations. The quantity $\eta = 1-\lambda_2(\mM)$ is called the spectral gap of the mixing matrix $\mM$, and $\eta^{-1}$ is typically a polynomial of the total number of nodes $n$ when the maximal degree of the node is $\cO(1)$. For example, for uniformly averaging $\mM$ one can show that $\eta^{-1} = \cO(n^2)$ for the ring topology (node degree $2$), $\eta^{-1} = \cO(n)$ for the two-dimensional torus topology (node degree  $2$), and $\eta^{-1} = \cO(1)$ for the fully connected graph (node degree $n-1$) \cite{aldous2002reversible}.

One or several steps of gossip can be used in distributed optimization algorithms as an alternative to aggregation through the central server, e.g., in Parallel {\tt SGD}. Choosing the communication graph in such way that there are no ``overloaded'' nodes, i.e., each node has a degree $\cO(1)$, one can significantly reduce the cost of one communication round in comparison to the parameter-server architecture. However, the communication complexity of gossip-based decentralized optimization methods often has multiplicative dependence on either $\cO(\eta^{-1})$ (see \cite{xu2020distributed} and references therein) or $\cO(\eta^{-\nicefrac{1}{2}})$ \cite{scaman2019optimal,uribe2020dual,fallah2019robust,kovalev2020optimal}, which is not improvable for gossip-based methods~\cite{arjevani2015communication,scaman2017optimal}. Since in the practically interesting cases we have $\eta = \Omega(n)$, it means that the overall number of communication rounds needed to achieve the desired accuracy of the solution grows with the number of workers $n$ as $\Omega(n)$ or $\Omega(\sqrt{n})$.

As an alternative to gossip, many practical distributed training systems perform averaging with All-Reduce ~\cite{goyal2017accurate,mikami2019massively,shoeybi2019megatron,You2020Large}. This name refers to a collection of protocols originally developed for HPC applications. Workers can follow these protocols to collectively compute the average gradient more efficiently than with a central server. The simplest variant of All-Reduce is known as Butterfly All-Reduce~\cite{bandwidth_optimal_allreduce}. Each of $n$ participants splits its local vector into $n$ chunks. Then, the $i$-th worker aggregates the $i$-th chunk of data from all peers and sends back the averaged chunk. As long as the vector size $s$ is greater than $n$, this protocol uses $\cO\left(s \times \frac{n - 1}{n}\right)$ total bandwidth on each worker. However, it requires all-to-all communication, which is not always practical for the HPC infrastructure. Real-world systems typically use Ring or Tree All-Reduce, where each worker only communicates with a small subset of its peers. These protocols enable highly efficient and scalable averaging with $\cO(1)$ or $\cO(\log N)$ total communication per worker.

As a result, All-Reduce Parallel {\tt SGD} enjoys the benefits of two worlds: the number of communication rounds does not grow with $n$, and each worker handles $\cO(s)$ amount data only, where $s$ is the size of one vector. However, All-Reduce protocols share a common drawback: they cannot tolerate node failures or network instability. If any single participant fails to execute its part or takes long to respond, this paralyzes all other workers. In contrast, gossip-based algorithms are more robust to such changes, which makes them applicable to time-varying networks~\cite{nedic2014distributed,nedic2016stochastic,nedic2018network,rogozin2019projected} and federated learning~\cite{ram2009asynchronous,yan2012distributed,yuan2016convergence}.

\subsubsection{The Fifth Contribution: Fault-Tolerant and Communication-Efficient Decentralized Optimization Method}
In this thesis, we lift the above restrictions by proposing {\tt Moshpit All-Reduce} --- an iterative averaging protocol that exponentially converges to the global average even with unreliable communication-constrained devices. According to our analysis, this method has exponential convergence independent of the network topology. Armed with this averaging protocol, we develop {\tt Moshpit SGD} for distributed optimization. We derive convergence rates for this algorithm and establish its equivalence to Centralized (Local) {\tt SGD} for (strongly) convex and non-convex problems.

Chapter~\ref{ch:moshpit} is devoted to the fifth contribution of this thesis, and is based on the following paper:
\begin{itemize}
    \item[\cite{ryabinin2021moshpit}] Max Ryabinin*, Eduard Gorbunov*, Vsevolod Plokhotnyuk, and Gennady Pekhimenko (*equal contribution). \textit{Moshpit SGD: Communication-Efficient Decentralized Training on Heterogeneous Unreliable Devices}. Advances in Neural Information Processing Systems, volume 34 (accepted), 2021.
\end{itemize}


\section{Scientific Novelty}
All results are new. They are summarized as follows:
\begin{itemize}
    \item We propose new general analysis of {\tt SGD} in the strongly convex case with proximable regularization. Our approach covers various existing methods in different settings. Whenever we recover a known method, our general theorem provides the tightest know rate for this method. Moreover, inspired by the proposed theoretical framework, we develop new stochastic methods ({\tt SGD-MB, SGD-star, N-SEGA, N-SAGA, Q-SGD-SR}).
    \item We propose a new unified theoretical framework for the analysis of stochastic first-order methods with error compensation and delayed updates. Using this framework, we develop 16 new methods. In particular, we develop the first full-gradient methods with error compensation that have linear convergence in the strongly convex case ({\tt EC-SGD-DIANA}) and the first variance reduced method with error compensation that also enjoys linear convergence on strongly convex problems ({\tt EC-LSVRG-DIANA}).
    \item We develop a new unified theoretical framework for the analysis of {\tt Local-SGD}-type methods when the objective function is (strongly) convex. We recover multiple known local optimizers as a special case of our general framework, along with their convergence rates (up to small constant factors). To demonstrate the strengths of our approach, we develop a new method called {\tt S-Local-SVRG} fitting our general framework. Moreover, using our general theorem we prove that {\tt S-Local-SVRG} converges linearly even when the local loss functions are arbitrarily heterogeneous. That is, we propose the first variance reduced linearly converging method without any restrictive assumptions.
    \item We develop and analyze \algname{MARINA}: a new communication efficient method for non-convex distributed learning over heterogeneous datasets. \algname{MARINA} employs a novel communication compression strategy based on the compression of gradient differences. Unlike  virtually all competing distributed first-order methods, ours is based on a carefully designed {\em biased} gradient estimator. Further, we develop and analyze two variants of \algname{MARINA}: \algname{VR-MARINA} and \algname{PP-MARINA}. The first method  is designed for the case when the local loss functions owned by clients are either of a finite sum or of an expectation form, and the second method allows for partial participation of clients. The proposed methods are superior to previous state-of-the-art methods in terms of oracle/communication complexity. Finally, we provide a convergence analysis of all methods  for problems satisfying the Polyak-{\L}ojasiewicz condition.
    \item We develop {\tt Moshpit All-Reduce} --- an iterative averaging protocol that exponentially converges to the global average even with unreliable communication-constrained devices. According to our analysis, this method has exponential convergence independent of the network topology. Armed with this averaging protocol, we develop {\tt Moshpit SGD} for distributed optimization. We derive convergence rates for this algorithm and establish its equivalence to Centralized (Local) {\tt SGD} for (strongly) convex and non-convex problems.
\end{itemize}

\section{Presentations and Validation of Research Results}
The results of this thesis were presented at the following conferences and seminars.

\begin{itemize}
    \item Neural Information Processing Systems 34 (NeurIPS 2021), ``Moshpit SGD: Communication-Efficient Decentralized Training on Heterogeneous Unreliable Devices'', online, 10 December, 2021.
    \item 38th International Conference on Machine Learning (ICML 2021), ``MARINA: Faster Non-Convex Distributed Learning with Compression'', online, 21 July, 2021.
    \item 24th International Conference on Artificial Intelligence and Statistics (AISTATS 2021), ``Local SGD: Unified Theory and New Efficient Methods'', online, 14 April, 2021.
    \item Federated Learning One-World Seminar, ``MARINA: Faster Non-Convex Distributed Learning with Compression'', online, 10 March, 2021.
    \item Neural Information Processing Systems 33 (NeurIPS 2020), ``Linearly Converging Error Compensated SGD'', online, 9 December, 2020.
    \item Federated Learning One-World Seminar and All-Russian Optimization Seminar, ``Linearly Converging Error Compensated SGD'', online, 7 October, 2020.
    \item 23rd International Conference on Artificial Intelligence and Statistics (AISTATS 2020), ``A Unified Theory of SGD: Variance Reduction, Sampling, Quantization and Coordinate Descent'', online, 26--28 August, 2020.
\end{itemize}

\section{Publications}
Chapters \ref{ch:sigma_k}-\ref{ch:moshpit} are based on the following papers, respectively:

\textbf{Published papers:}
\begin{itemize}
    \item[\cite{gorbunov2019unified}] Eduard Gorbunov, Filip Hanzely, and Peter Richt\'arik. \textit{A Unified Theory of SGD: Variance Reduction, Sampling, Quantization and Coordinate Descent}. In Silvia Chiappa and Roberto Calandra, editors, Proceedings of the Twenty Third International Conference on Artificial Intelligence and Statistics, volume 108 of Proceedings of Machine Learning Research, pages 680–690. PMLR, 26–28 Aug 2020.
    \item[\cite{gorbunov2020linearly}] Eduard Gorbunov, Dmitry Kovalev, Dmitry Makarenko, and Peter Richt\'arik. \textit{Linearly Converging Error Compensated SGD}. In H. Larochelle, M. Ranzato, R. Hadsell, M. F. Balcan, and H. Lin, editors, Advances in Neural Information Processing Systems, volume 33, pages 20889–20900. Curran Associates, Inc., 2020.
    \item[\cite{gorbunov2020local}] Eduard Gorbunov, Filip Hanzely, and Peter Richt\'arik. \textit{Local SGD: Unified Theory and New Efficient Methods}. In Arindam Banerjee and Kenji Fukumizu, editors, Proceedings of The 24th International Conference on Artificial Intelligence and Statistics, volume 130 of Proceedings of Machine Learning Research, pages 3556–3564. PMLR, 13–15 Apr 2021.
    \item[\cite{pmlr-v139-gorbunov21a}] Eduard Gorbunov, Konstantin P. Burlachenko, Zhize Li, and Peter Richt\'arik. \textit{MARINA: Faster Non-Convex Distributed Learning with Compression}. In Marina Meila and Tong Zhang, editors, Proceedings of the 38th International Conference on Machine Learning, volume 139 of Proceedings of Machine Learning Research, pages 3788–3798. PMLR, 18–24 Jul 2021.
\end{itemize}

\textbf{In print:}
\begin{itemize}
    \item[\cite{ryabinin2021moshpit}] Max Ryabinin*, Eduard Gorbunov*, Vsevolod Plokhotnyuk, and Gennady Pekhimenko (*equal contribution). \textit{Moshpit SGD: Communication-Efficient Decentralized Training on Heterogeneous Unreliable Devices}. Advances in Neural Information Processing Systems, volume 34, 2021.
\end{itemize}

Appendix~\ref{app:ef_sigma_k} contains extra plots, some missing proofs, and the results for the methods with delayed updates from \cite{gorbunov2020linearly} (Chapter~\ref{ch:ef_sigma_k}). Extra experiments and missing proofs of the general results from Chapter~\ref{ch:local_sigma_k} are deferred to Appendix~\ref{app:local_sigma_k}. Finally, missing proofs and additional technical details from Chapters~\ref{ch:marina}~and~\ref{ch:moshpit} are given in Appendices~\ref{app:marina}~and~\ref{app:moshpit}, respectively.

\subsection{Excluded Papers}
During my PhD studies, I was also fortunate to co-author two papers on stochastic optimization with heavy-tailed noise in stochastic gradients \cite{gorbunov2020clipped_sstm,gorbunov2021near}, a paper on Byzantine-tolerant distributed optimization without parameter server \cite{gorbunov2021secure}, two review-papers on non-convex optimization \cite{danilova2020recent} and decentralized distributed optimization \cite{gorbunov2020recent}, a paper on extensions of modern error feedback \cite{fatkhullin2021ef21}, a paper on the last-iterate convergence analysis of Extragradient method \cite{gorbunov2021extragradient}, and a paper on new analysis of its stochastic versions \cite{gorbunov2021stochastic}.

\section{Thesis Structure}
The thesis consists of an introduction, 5 main chapters, list of 247 references, and 5 chapters in the Appendix with technical details, some proofs, and auxiliary results.

%% file: ch2_sigma_k.tex
\chapter{A Unified Theory of SGD: Variance Reduction, Sampling, Quantization  and Coordinate Descent}\label{ch:sigma_k}

\section{Introduction}

In this chapter, we are interested in  the optimization  problem
    \begin{equation}\label{eq:problem_gen}
        \min_{x\in\R^d} f(x) + R(x),
    \end{equation}
    where $f$ is convex, differentiable with Lipschitz gradient, and $R:\R^d\to \R\cup \{+\infty\}$ is a proximable (proper closed convex) regularizer. In particular, we  focus on situations when  it is prohibitively expensive to  compute the gradient of $f$, while an unbiased estimator of the gradient can be computed efficiently. This is typically the case for stochastic optimization problems, i.e., when \begin{equation} \label{eq:f_exp}
    f(x)=\EE_{\xi \sim \cD} \left[ f_\xi(x)\right] ,\end{equation} where $\xi$ is a random variable, and  $f_\xi:\R^d\to \R$  is smooth for all $\xi$. Stochastic optimization problems are of key importance in statistical supervised learning theory. In this setup, $x$ represents a machine learning model described by $d$ parameters (e.g., logistic regression or a deep neural network), $\cD$  is an unknown distribution of labelled examples,  $f_\xi(x)$ represents the loss of model $x$ on datapoint $\xi$, and $f$ is the generalization error. Problem \eqref{eq:problem_gen} seeks to find the model $x$ minimizing the generalization error. In statistical learning theory one assumes that while $\cD$ is not known, samples $\xi\sim \cD$ are available. In such a case, $\nabla f(x)$ is not computable, while $\nabla f_\xi(x)$, which is  an unbiased estimator of the gradient of $f$ at $x$, is easily computable. 

Another prominent example, one of special interest in this chapter, are functions $f$ which arise as averages of a very large number of smooth functions: \begin{equation}\label{eq:f_sum} 
 f(x)=\frac{1}{n}\sum \limits_{i=1}^n f_i(x).\end{equation} This problem often arises by approximation of the stochastic optimization loss function \eqref{eq:f_exp} via Monte Carlo integration, and is in this context known as the empirical risk minimization (ERM) problem.  ERM is currently the dominant paradigm for solving supervised learning problems \citep{shalev2014understanding}.  If index $i$ is chosen uniformly at random from $[n]\eqdef \{1,2,\dots,n\}$, $\nabla f_i(x)$ is an unbiased estimator of $\nabla f(x)$. Typically, $\nabla f(x)$ is about $n$ times more expensive to compute than $\nabla f_i(x)$. 

Lastly, in some applications, especially in distributed training of supervised models, one considers problem \eqref{eq:f_sum}, with $n$ being the number of machines, and each $f_i$ also having a finite sum structure, i.e.,
\begin{equation}\label{eq:f_i_sum}
    f_i(x) = \frac{1}{m}\sum\limits_{j=1}^m  f_{ij}(x),
\end{equation}
where $m$ corresponds to the number of training examples stored on machine $i$.

\section{The Many Faces of Stochastic Gradient Descent} \label{sec:many_faces_of_SGD}

Stochastic gradient descent ({\tt SGD}) \citep{RobbinsMonro:1951, Nemirovski-Juditsky-Lan-Shapiro-2009, Vaswani2019-overparam} is a state-of-the-art algorithmic paradigm for solving optimization problems \eqref{eq:problem_gen} in situations when $f$ is either of structure 
 \eqref{eq:f_exp} or \eqref{eq:f_sum}. In its generic form, (proximal) {\tt SGD} defines the new iterate by subtracting a multiple of a stochastic gradient from the current iterate, and subsequently applying the proximal operator of $R$:
\begin{equation} \label{eq:SGD} x^{k+1} = \proxR(x^k - \gamma g^k).\end{equation} 
Here, $g^k$ is an unbiased estimator of the gradient (i.e., a stochastic gradient), \begin{equation}\label{eq:stoch_grad}\EE \left[g^k \;|\; x^k\right] = \nabla f(x^k),\end{equation}
and $\proxR(x)\eqdef \argmin_u \{\gamma R(x) + \frac{1}{2}\norm{u-x}^2\}$. However, and this is the starting point of our journey in this paper, there are {\em infinitely many} ways of obtaining a random vector $g^k$ satisfying \eqref{eq:stoch_grad}. On the one hand, this gives algorithm designers the flexibility to {\em construct}  stochastic gradients in various ways in order to target desirable properties such as convergence speed, iteration cost, parallelizability and generalization.  On the other hand, this poses considerable challenges in terms of convergence analysis. Indeed, if one aims to, as one should, obtain the sharpest bounds possible, dedicated analyses are needed to handle each of the particular variants of {\tt SGD}.

\textbf{Vanilla\footnote{In this thesis, by {\em vanilla} {\tt SGD} we refer to {\tt SGD} variants with or without importance sampling and mini-batching, but {\em excluding} variance-reduced variants, such as {\tt SAGA} \citep{SAGA} and {\tt SVRG} \citep{SVRG}.} {\tt SGD}.} The flexibility in the design of efficient strategies for constructing $g^k$ has led to a creative renaissance in the optimization and machine learning communities, yielding a large number of immensely powerful new variants of {\tt SGD}, such as those employing  {\em importance sampling} \citep{IProx-SDCA, NeedellWard2015}, and {\em mini-batching} \citep{mS2GD}. These efforts are subsumed by the recently developed and remarkably sharp analysis of {\tt SGD} under  {\em arbitrary sampling} paradigm \citep{gower2019sgd}, first introduced in the study of randomized coordinate descent methods by \citep{NSync}. The arbitrary sampling paradigm covers virtually all stationary mini-batch and importance sampling strategies in a unified way, thus making headway towards theoretical unification of two separate strategies for constructing stochastic gradients. For strongly convex $f$, the {\tt SGD} methods analyzed in \citep{gower2019sgd} converge linearly to a neighbourhood of the solution $x^* = \arg \min_x f(x)$ for a fixed stepsize $\gamma^k=\gamma$. The size of the neighbourhood is proportional to the second moment of the stochastic gradient at the optimum ($\sigma^2 \eqdef \tfrac{1}{n}\sum_{i=1}^n \norm{\nabla f_i(x^*)}^2$), to the stepsize ($\gamma$), and inversely proportional to the modulus of strong convexity. The effect of various sampling strategies, such as importance sampling and mini-batching, is twofold: i) improvement of the linear convergence rate by enabling larger stepsizes, and ii) modification of $\sigma^2$. However, none of these strategies\footnote{Except for the full batch strategy, which is prohibitively expensive.} is able to completely eliminate the adverse effect of $\sigma^2$. That is,  {\tt SGD} with a fixed stepsize does not reach the optimum, unless one happens to be in the overparameterized case characterized by the identity $\sigma^2=0$.

\textbf{Variance reduced {\tt SGD}.} While sampling strategies such as importance sampling and mini-batching reduce the variance of the stochastic gradient, in the finite-sum case \eqref{eq:f_sum} a new type of {\em variance reduction} strategies has been developed over the last few years \citep{SAG, SAGA, SVRG, SDCA, QUARTZ,nguyen2017sarah, kovalev2019don}. These variance-reduced {\tt SGD} methods differ from the sampling strategies discussed before in a significant way: they can iteratively {\em learn} the stochastic gradients at the optimum,  and in so doing are able to eliminate the adverse effect of the gradient noise $\sigma^2>0$ which, as mentioned above, prevents the iterates of vanilla {\tt SGD} from converging to the optimum. As a result, for strongly convex $f$, these new variance-reduced {\tt SGD} methods converge linearly to $x^*$, with a fixed stepsize. At the moment, these variance-reduced variants require a markedly different convergence theory from the vanilla variants of {\tt SGD}. An exception to this is the situation when $\sigma^2=0$ as then variance reduction is not needed; indeed, vanilla {\tt SGD} already converges to the optimum, and with a fixed stepsize. We end the discussion here by remarking that this {\em hints} at a possible existence of a more unified theory, one that would include both vanilla and variance-reduced {\tt SGD}.

\textbf{Distributed {\tt SGD}, quantization and variance reduction.} When {\tt SGD} is implemented in a distributed fashion,  the problem is often expressed in the form \eqref{eq:f_sum}, where $n$ is the number of workers/nodes, and  $f_i$ corresponds to the loss based on data stored on node $i$. Depending on the number of data points stored on each node, it may or may not be efficient to compute the gradient of $f_i$ in each iteration. In general, {\tt SGD} is implemented in this way: each  node $i$ first computes a stochastic gradient $g_i^k$ of $f_i$ at the current point $x^k$ (maintained individually by each node). These gradients are then aggregated by a master node \citep{DANE, RDME}, in-network by a switch~\citep{switchML}, or a different technique best suited to the architecture used. To alleviate the communication bottleneck, various lossy update compression strategies such as quantization~\citep{seide20141, Gupta:2015limited, zipml}, sparsification~\citep{RDME, alistarh2018convergence, wangni2018gradient} and dithering~\citep{alistarh2017qsgd} were proposed. The basic idea is for each worker to apply a randomized transformation $Q:\R^d\to \R^d$ to $g_i^k$, resulting in a vector which is still an unbiased estimator of the gradient, but one that can be communicated with fewer bits. Mathematically, this amounts to injecting additional noise into the already noisy stochastic gradient $g_i^k$. The field of quantized {\tt SGD} is still young, and even some basic questions remained open until recently. For instance, there was no distributed quantized {\tt SGD} capable of provably solving  \eqref{eq:problem_gen} until the {\tt DIANA} algorithm \citep{mishchenko2019distributed} was introduced. {\tt DIANA} applies quantization to {\em gradient differences}, and in so doing is able to learn the gradients at the optimum, which makes it able to work for any regularizer $R$. {\tt DIANA} has some structural similarities with {\tt SEGA}~\citep{hanzely2018sega}---the first coordinate descent type method which works for non-separable regularizers---but a more precise relationship remains elusive.  When the functions of $f_i$ are of a finite-sum structure as in \eqref{eq:f_i_sum}, one can apply variance reduction to reduce the variance of the stochastic gradients $g_i^k$ together with quantization, resulting in the {\tt VR-DIANA} method~\citep{horvath2019stochastic}. This is the first distributed quantized {\tt SGD} method which provably converges to the solution of \eqref{eq:problem_gen}+\eqref{eq:f_i_sum} with a fixed stepsize.


\textbf{Randomized coordinate descent ({\tt RCD}).} Lastly, in a distinctly separate strain, there are {\tt SGD} methods for the coordinate/subspace descent variety~\citep{RCDM}. While it is possible to see  {\em some} {\tt RCD} methods as special cases of \eqref{eq:SGD}+\eqref{eq:stoch_grad}, most of them do not follow this algorithmic template. First, standard {\tt RCD} methods use different stepsizes for updating different coordinates~\citep{ALPHA}, and this seems to be crucial to their success. Second,  until the recent discovery of the {\tt SEGA} method, {\tt RCD} methods were not able to converge with non-separable regularizers. Third, {\tt RCD} methods are naturally variance-reduced in the $R=0$ case as partial derivatives at the optimum are all zero. As a consequence, attempts at creating variance-reduced {\tt RCD} methods seem to be futile. Lastly, {\tt RCD} methods are typically analyzed using different techniques. While there are deep links between standard {\tt SGD} and {\tt RCD} methods, these are often indirect and rely on duality \citep{SDCA, FACE-OFF, SDA}. 


\section{Contributions}

As outlined in the previous section, the world of {\tt SGD} is vast and beautiful. It is formed by many largely disconnected islands populated by elegant and efficient methods, with their own applications, intuitions, and convergence analysis techniques. While some links already exist (e.g., the unification of importance sampling and mini-batching variants under the arbitrary sampling umbrella), there is no comprehensive general theory. 
It is becoming increasingly difficult for the community to understand the relationships between these variants, both in theory and practice. New variants are yet to be discovered, but it is not clear what tangible principles one should adopt beyond intuition to aid the discovery. This situation is exacerbated by the fact that a number of different assumptions on the stochastic gradient, of various levels of strength, is being used in the literature.


The main contributions of this work include: 

$\bullet$ {\bf  Unified analysis.} In this work we propose a {\em unifying theoretical framework}
 which covers all of the  variants of {\tt SGD} outlined in Section~\ref{sec:many_faces_of_SGD}. As a by-product, we obtain the {\em first unified analysis} of vanilla and variance-reduced {\tt SGD} methods.  For instance, our analysis covers as special cases vanilla {\tt SGD} methods from~\citep{nguyen2018sgd} and~\citep{gower2019sgd}, variance-reduced {\tt SGD} methods such as {\tt SAGA}~\citep{SAGA}, {\tt L-SVRG}~\citep{hofmann2015variance, kovalev2019don} and {\tt JacSketch}~\citep{gower2018stochastic}. Another by-product  is {\em the unified analysis of {\tt SGD} methods which include {\tt RCD}.} For instance, our theory covers the subspace descent method {\tt SEGA}~\citep{hanzely2018sega} as a special case. Lastly, our framework is general enough to capture the phenomenon of {\em quantization}. For instance, we obtain the {\tt DIANA} and {\tt VR-DIANA} methods in special cases. 


$\bullet$ {\bf Generalization of existing methods.} An important yet {\em relatively} minor contribution of our work is that it enables  {\em generalization} of knowns methods.  For instance, some particular methods we consider, such as {\tt L-SVRG} (Alg~\ref{alg:L-SVRG})~\citep{kovalev2019don}, were not analyzed in the proximal ($R\neq 0$) case before. To illustrate how this can be done within our framework, we do it here for {\tt L-SVRG}. Further, most\footnote{Our analysis allows for arbitrary sampling of all methods except of those using partial derivatives such as {\tt SEGA} or {\tt N-SEGA}. We shall note that arbitrary sampling for {\tt SEGA} was developed concurrently in~\citep{hanzely2019one}. Note that~\citep{hanzely2019one} proposes many novel variance reduced algorithms, for some of which we can obtain best rates. A detailed discussion and comparison to~\citep{hanzely2019one} is provided in Remark~\ref{rem:gjs} in the Appendix} of the methods we analyze can be extended to the {\em arbitrary sampling} paradigm. 

$\bullet$ {\bf Sharp rates.} In all known special cases, the rates obtained from our general theorem (Theorem~\ref{thm:main_gsgm}) are the {\em best known rates} for these methods.

$\bullet$ {\bf New methods.}   Our general analysis provides estimates for a possibly infinite array of new and yet-to-be-developed variants of {\tt SGD}. One only needs to verify that Assumption~\ref{as:general_stoch_gradient} holds, and a complexity estimate is readily furnished by Theorem~\ref{thm:main_gsgm}. Selected existing and new methods that fit our framework are summarized in Table~\ref{tbl:special_cases2}. This list is for illustration only, we believe that future work by us and others will lead to its rapid expansion.

$\bullet$ {\bf Experiments.} We show through extensive experimentation that some of the {\em new} and {\em generalized} methods proposed here and analyzed via our framework have some intriguing practical properties when compared against appropriately selected existing methods.

\section{Main Result \label{sec:main_res}}

We first introduce the key assumption on the stochastic gradients $g^k$ enabling our general analysis (Assumption~\ref{as:general_stoch_gradient}), then state our assumptions on $f$ (Assumption~\ref{as:mu_strongly_quasi_convex}), and finally  state and comment on our  unified convergence result (Theorem~\ref{thm:main_gsgm}).

\subsection{Key Assumption}
Our first assumption is of key importance. It is mainly an assumption on the sequence of stochastic gradients $\{g^k\}$ generated by an arbitrary randomized algorithm. Besides unbiasedness (see \eqref{eq:general_stoch_grad_unbias}), we require two recursions to hold for the iterates $x^k$ and the stochastic gradients $g^k$ of a randomized method. We allow for flexibility by casting these inequalities in a parametric manner.

\begin{assumption}\label{as:general_stoch_gradient} Let $\{x^k\}$ be the random iterates produced by proximal {\tt SGD} (Algorithm in Eq~\eqref{eq:SGD}).
We first assume that the stochastic gradients $g^k$ are unbiased
    \begin{equation}\label{eq:general_stoch_grad_unbias}
        \EE\left[ g^k\mid x^k\right] = \nabla f(x^k),
    \end{equation}
    for all $k\geq 0$. Further, we assume that there exist
non-negative constants $A, B, C, D_1, D_2, \rho$ and a (possibly) random sequence $\{\sigma_k^2\}_{k\ge 0}$ such that the following two relations hold\footnote{For convex and $L$-smooth $f$, one can show that $    \norm{\nabla f(x) - \nabla f(y)}^2 \le 2LD_{f}(x,y).$ Hence, $D_f$ can be used as a measure of proximity for the gradients.}
    \begin{equation}\label{eq:general_stoch_grad_second_moment}
        \EE\left[\norm{g^k -\nabla f(x^*)}^2\mid x^k\right] \le 2AD_f(x^k,x^*) + B\sigma_k^2 + D_1,
    \end{equation}
    \begin{equation}\label{eq:gsg_sigma}
        \EE\left[\sigma_{k+1}^2 \, \mid \, \sigma^2_k\right] \le (1-\rho) \sigma_k^2 + 2CD_f(x^k,x^*)  + D_2,
    \end{equation} 
The expectation above is with respect to the randomness of the algorithm. 

\end{assumption}

The unbiasedness assumption \eqref{eq:general_stoch_grad_unbias} is standard.  The key innovation we bring is inequality \eqref{eq:general_stoch_grad_second_moment} coupled with \eqref{eq:gsg_sigma}. We argue, and justify this statement by furnishing many examples in Section~\ref{sec:main-paper-special-cases}, that these inequalities capture the essence of a wide array of existing and some new {\tt SGD} methods, including vanilla, variance reduced, arbitrary sampling, quantized and coordinate descent variants. Note that in the case when $\nabla f(x^*) = 0$ (e.g., when $R=0$),  the inequalities in Assumption~\ref{as:general_stoch_gradient}  reduce to
\begin{equation}\label{eq:general_stoch_grad_second_moment_special}
        \EE\left[\norm{g^k}^2\mid x^k\right] \le 2A(f(x^k) - f(x^*)) + B\sigma_k^2 + D_1,
\end{equation}
\begin{equation}\label{eq:gsg_sigma_special}
    \EE\left[\sigma_{k+1}^2 \, \mid \, \sigma^2_k\right] \le (1-\rho) \sigma_k^2 + 2C(f(x^k) - f(x^*)) + D_2.
\end{equation}
Similar inequalities can be found in the analysis of stochastic first-order methods. However, this is the first time that such inequalities are generalized, equipped with parameters, and elevated to the status of an assumption that can be used on its own, independently from any other details defining the underlying method that generated them.

To give a further intuition about inequalities~\eqref{eq:general_stoch_grad_second_moment} and~\eqref{eq:gsg_sigma},    we shall note that sequence $\sigma_k$ usually represents the portion of noise that can gradually decrease over the course of optimization while constants $D_1, D_2$ represent a static noise. On the other hand, constants $A, C$ are usually related to some measure of smoothness of the objective. For instance, the parameters for (deterministic) gradient descent can be chosen as $A = L, B = C = D_1 =D_2 = \sigma_k^2 = \rho =0$. For an overview of parameter choices for specific instances of~\eqref{eq:SGD}, see Table~\ref{tbl:special_cases-parameters}. Note also that the choice of parameters of~\eqref{eq:general_stoch_grad_second_moment} and~\eqref{eq:gsg_sigma} is not unique, however this has no impact on convergence rates we provide.

\subsection{Main Theorem}

 For simplicity, we shall assume throughout that $f$ is $(\mu,x^*)$-strongly quasi-convex, which is a generalization of $\mu$-strong convexity. 

\begin{assumption}[$(\mu,x^*)$-strong quasi-convexity]\label{as:mu_strongly_quasi_convex}
There exists $\mu>0$ such that $f:\R^d\to \R$ satisfies the following inequality for all $ x\in\R^d$:
    \begin{equation}\label{eq:mu_strongly_quasi_convex}
        f(x^*) \ge f(x) + \langle \nabla f(x), x^* - x\rangle + \frac{\mu}{2}\norm{x^* - x}^2.    
    \end{equation}
\end{assumption}

We are now ready to present the key lemma of this paper which states per iteration recurrence to analyze~\eqref{eq:SGD}.
 
\begin{lemma}  \label{lem:iter_dec}
Let Assumptions~\ref{as:general_stoch_gradient}~and~\ref{as:mu_strongly_quasi_convex}  be satisfied. Then the following inequality holds for all $k\geq 0$:
\begin{eqnarray*}
    \EE\left[\norm{x^{k+1}-x^*}^2\right] + M\gamma^2\EE\left[\sigma_{k+1}^2\right]  &\le& (1-\gamma\mu)\EE \left[\norm{x^k - x^*}^2 \right] + \left(1 - \rho + \frac{B}{M}\right)M\gamma^2\EE\left[\sigma_k^2\right] \\
    && \quad - 2\gamma\left(1-\gamma(A+CM)\right)\EE\left[D_f(x^k,x^*)\right]\\
    &&\quad + (D_1+MD_2)\gamma^2.
\end{eqnarray*}
\end{lemma}
\begin{proof}
    We start with estimating the first term of the Lyapunov function. Let $r^k = x^k - x^*$. Then
    \begin{eqnarray*}
        \norm{r^{k+1}}^2 &=& \norm{ \prox_{\gamma R} (x^k- \gamma  g^k) - \prox_{\gamma R} ( x^* - \gamma\nabla f(x^*)) }^2 \\
         &\leq & 
          \norm{x^k- x^* - \gamma (  g^k - \nabla f(x^*)) }^2
          \\
          &= & 
           \norm{r^k}^2 - 2\gamma\langle r^k,g^k  - \nabla f(x^*)\rangle + \gamma^2\norm{ g^k  - \nabla f(x^*)}^2.
    \end{eqnarray*}
    Taking expectation conditioned on $x^k$ we get
    \begin{eqnarray*}
        \EE\left[\norm{r^{k+1}}^2\mid x^k\right] &=& \norm{r^k}^2 - 2\gamma\langle r^k,\nabla f(x^k)-\nabla f(x^*)\rangle + \gamma^2\EE\left[\norm{ g^k - \nabla f(x^*)}^2\mid x^k\right]
        \\
        &\overset{\eqref{eq:mu_strongly_quasi_convex}}{\le}& 
        (1-\gamma\mu)\norm{r^k}^2 - 2\gamma D_f(x^k,x^*) + \gamma^2\EE\left[\norm{g^k - \nabla f(x^*)}^2\mid x^k\right]
        \\
        &\overset{\eqref{eq:general_stoch_grad_unbias}+\eqref{eq:general_stoch_grad_second_moment}}{\le}& 
        (1-\gamma\mu)\norm{r^k}^2 + 2\gamma\left(A\gamma - 1\right)D_f(x^k,x^*) + B\gamma^2\sigma_k^2 + \gamma^2 D_1.
    \end{eqnarray*}
    Using this we estimate the full expectation of $V^{k+1}$ in the following way:
    \begin{eqnarray}
    \EE\norm{x^{k+1}-x^*}^2 + M\gamma^2\EE\sigma_{k+1}^2 &\overset{\eqref{eq:gsg_sigma}}{\le}& (
        1-\gamma\mu)\EE\norm{x^k-x^*}^2 + 2\gamma\left(A\gamma - 1\right)\EE\left[D_f(x^k,x^*)\right] \notag
        \\
        &&\quad 
        + (1-\rho)M\gamma^2\EE\sigma_k^2  + 2CM\gamma^2\EE\left[D_f(x^k,x^*)\right] \notag\\
        &&\quad + B\gamma^2\EE\sigma_k^2 + (D_1+MD_2)\gamma^2 \notag\\
        &=& (1-\gamma\mu)\EE\norm{x^k - x^*}^2 + \left(1 + \frac{B}{M} - \rho\right)M\gamma^2\EE\sigma_k^2\notag
        \\
        &&\quad
         + 2\gamma\left(\gamma(A+CM)-1\right)\EE\left[D_f(x^k,x^*)\right] \notag\\
        &&\quad + (D_1+MD_2)\gamma^2.\notag
        \label{eq:gsgm_recurrence}
    \end{eqnarray}
It remains to rearrange the terms.
\end{proof}

Using recursively Lemma~\ref{lem:iter_dec}, we obtain the convergence rate of proximal SGD, which we state as Theorem~\ref{thm:main_gsgm}.

\begin{theorem}\label{thm:main_gsgm}
Let Assumptions~\ref{as:general_stoch_gradient}~and~\ref{as:mu_strongly_quasi_convex} be satisfied. Choose constant $M$ such that $M > \frac{B}{\rho}$. Choose a stepsize  satisfying
    \begin{equation}\label{eq:gamma_condition_gsg}0 < \gamma \le \min\left\{\frac{1}{\mu}, \frac{1}{A+CM}\right\}.
    \end{equation}
    Then the iterates $\{x^k\}_{k\geq 0}$ of proximal {\tt SGD} (Algorithm~\eqref{eq:SGD}) satisfy
    \begin{equation}
        \EE\left[V^k\right] \le \max\left\{(1-\gamma\mu)^k,\left(1+\frac{B}{M}-\rho\right)^k\right\} V^0 + \frac{(D_1+MD_2)\gamma^2 }{\min\left\{\gamma\mu, \rho - \frac{B}{M}\right\}},\label{eq:main_result_gsgm}
    \end{equation}
    where the Lyapunov function $V^k$ is defined by $V^k \eqdef \norm{x^k - x^*}^2 + M\gamma^2\sigma_k^2$.
\end{theorem}
\begin{proof}
    Note first that due to~\eqref{eq:gamma_condition_gsg} we have $2\gamma\left(1-\gamma(A+CM)\right)\EE D_f(x^k,x^*)>0$, thus we can omit the term.

    Unrolling the recurrence from Lemma~\ref{lem:iter_dec} and using the Lyapunov function notation gives us
    \begin{eqnarray*}
        \EE V^{k} &\le& \max\left\{(1-\gamma\mu)^k,\left(1+\frac{B}{M}-\rho\right)^k\right\}V^0\\
        &&\quad + (D_1+MD_2)\gamma^2\sum\limits_{l=0}^{k-1}\max\left\{(1-\gamma\mu)^l,\left(1+\frac{B}{M}-\rho\right)^l\right\}\\
        &\le& \max\left\{(1-\gamma\mu)^k,\left(1+\frac{B}{M}-\rho\right)^k\right\}V^0\\
        &&\quad + (D_1+MD_2)\gamma^2\sum\limits_{l=0}^{\infty}\max\left\{(1-\gamma\mu)^l,\left(1+\frac{B}{M}-\rho\right)^l\right\}\\
        &\le& \max\left\{(1-\gamma\mu)^k,\left(1+\frac{B}{M}-\rho\right)^k\right\} V^0 + \frac{(D_1+MD_2)\gamma^2}{\min\left\{\gamma\mu, \rho - \frac{B}{M}\right\}}.
    \end{eqnarray*}
\end{proof}

This theorem establishes a linear rate for a wide range of proximal {\tt SGD} methods up to a certain oscillation radius, controlled by the additive term in \eqref{eq:main_result_gsgm}, and namely, by parameters $D_1$ and $D_2$. As we shall see  in Section~\ref{sec:special_cases} (refer to Table~\ref{tbl:special_cases-parameters}), the main difference between the vanilla and variance-reduced {\tt SGD} methods is that while the former  satisfy inequality \eqref{eq:gsg_sigma} with $D_1>0$ or $D_2>0$, which in view of \eqref{eq:main_result_gsgm} prevents them from reaching the optimum $x^*$ (using a fixed stepsize), the latter methods satisfy inequality \eqref{eq:gsg_sigma} with $D_1=D_2=0$, which in view of \eqref{eq:main_result_gsgm} enables them to reach the optimum.

\section{The Classic, The Recent and The Brand New} \label{sec:main-paper-special-cases}

In this section we deliver on the promise from the introduction and show how many existing and some new variants of {\tt SGD} fit our general framework (see Table~\ref{tbl:special_cases2}). 

{\bf An overview.} As claimed, our framework is powerful enough to include vanilla methods (\xmark\; in the ``VR'' column) as well as variance-reduced methods (\cmark\; in the ``VR'' column), methods which generalize to arbitrary sampling (\cmark\; in the ``AS'' column), methods supporting gradient quantization (\cmark\; in the ``Quant'' column) and finally, also {\tt RCD} type methods (\cmark\; in the ``RCD'' column). 

\begin{table}[H]
\caption{List of specific existing (in some cases generalized) and new methods which fit our general analysis framework. VR = variance reduced method, AS = arbitrary sampling, Quant = supports gradient quantization, RCD = randomized coordinate descent type method. ${}^{a}$ Special case of {\tt SVRG} with  1 outer loop only;  ${}^{b}$ Special case of {\tt DIANA} with $1$ node and quantization of exact gradient. }
\label{tbl:special_cases2}
\begin{center}
\footnotesize
\begin{tabular}{|c|c|c|c|c|c|c|c|c|c|}
\hline
Problem & Method &  Alg \# & Citation &   VR?  & AS? & Quant? &  RCD? & Section  & Result \\
\hline
\eqref{eq:problem_gen}+\eqref{eq:f_exp} & {\tt SGD}  & Alg \ref{alg:sgd_prox} & {\tiny\cite{nguyen2018sgd}}  & \xmark &  \xmark & \xmark &  \xmark & \ref{sec:SGD} & 
Cor~\ref{cor:recover_sgd_rate} \\
\eqref{eq:problem_gen}+\eqref{eq:f_sum}  & {\tt SGD-SR} & Alg \ref{alg:sgdas} & {\tiny\cite{gower2019sgd}} & \xmark &  \cmark & \xmark & \xmark & \ref{SGD-AS} &   Cor~\ref{cor:recover_sgd-as_rate} \\
\eqref{eq:problem_gen}+\eqref{eq:f_sum} &  {\tt SGD-MB} & Alg~\ref{alg:SGD-MB} & {\bf NEW} & \xmark & \xmark & \xmark & \xmark & \ref{sec:SGD-MB} & Cor~\ref{cor:mb}   \\
\eqref{eq:problem_gen}+\eqref{eq:f_sum} &  {\tt SGD-star} & Alg \ref{alg:SGD-star} & {\bf NEW} & \cmark & \cmark  & \xmark & \xmark & \ref{sec:SGD-star} & Cor~\ref{cor:SGD-star} \\
\eqref{eq:problem_gen}+\eqref{eq:f_sum}  & {\tt SAGA} & Alg~\ref{alg:SAGA} & {\tiny\cite{SAGA}} & \cmark & \xmark  & \xmark & \xmark & \ref{sec:saga} & Cor~\ref{thm:recover_saga_rate} \\
\eqref{eq:problem_gen}+\eqref{eq:f_sum}  & {\tt N-SAGA} & Alg~\ref{alg:N-SAGA} & {\bf NEW} & \xmark &  \xmark & \xmark & \xmark & \ref{N-SAGA} & Cor~\ref{cor:N-SAGA} \\
\eqref{eq:problem_gen} & {\tt SEGA}  & Alg~\ref{alg:SEGA} &  {\tiny\cite{hanzely2018sega}}  & \cmark &   \xmark & \xmark & \cmark & \ref{sec:sega} & Cor~\ref{cor:sega} \\
\eqref{eq:problem_gen} & {\tt N-SEGA}  & Alg~\ref{alg:N-SEGA} &  {\bf NEW}  & \xmark & \xmark  & \xmark  & \cmark & \ref{N-SEGA} &
Cor~\ref{cor:N-SEGA}\\
\eqref{eq:problem_gen}+\eqref{eq:f_sum}  & {\tt SVRG}${}^{a}$  & Alg~\ref{alg:SVRG} & {\tiny\cite{SVRG}} & \cmark & \xmark & \xmark & \xmark & \ref{sec:svrg} & Cor~\ref{cor:svrg} \\
\eqref{eq:problem_gen}+\eqref{eq:f_sum}  & {\tt L-SVRG} & Alg~\ref{alg:L-SVRG} & {\tiny\cite{hofmann2015variance}} & \cmark &  \xmark & \xmark & \xmark & \ref{sec:L-SVRG} & Cor~\ref{cor:recover_l-svrg_rate}  \\
\eqref{eq:problem_gen}+\eqref{eq:f_sum}  & {\tt DIANA} & Alg~\ref{alg:diana} & 
{\tiny\cite{mishchenko2019distributed}} & \xmark & \xmark & \cmark &  \xmark & \ref{sec:diana} & Cor~\ref{cor:main_diana}  \\
\eqref{eq:problem_gen}+\eqref{eq:f_sum}  & {\tt DIANA}${}^b$  & Alg~\ref{alg:diana_case}  & {\tiny\cite{mishchenko2019distributed}} & \cmark &  \xmark & \cmark  &  \xmark &\ref{sec:diana}  & Cor~\ref{cor:main_diana_special_case}\\
\eqref{eq:problem_gen}+\eqref{eq:f_sum}  & {\tt Q-SGD-SR} & Alg \ref{alg:qsgdas} & {\bf NEW} & \xmark &  \cmark & \cmark & \xmark & \ref{Q-SGD-AS} &   Cor~\ref{cor:recover_q-sgd-as_rate} \\ 
\eqref{eq:problem_gen}+\eqref{eq:f_sum}+\eqref{eq:f_i_sum}  & {\tt VR-DIANA} & Alg~\ref{alg:vr-diana_sigma_k} & {\tiny\cite{horvath2019stochastic}}& \cmark & \xmark & \cmark & \xmark & \ref{sec:VR-DIANA} & Cor~\ref{cor:main_vr_diana} \\
\eqref{eq:problem_gen}+\eqref{eq:f_sum}& {\tt JacSketch} & Alg~\ref{alg:jacsketch} & {\tiny\cite{gower2018stochastic}} & \cmark &  \cmark \xmark & \xmark & \xmark & \ref{sec:JacSketch} & Cor~\ref{thm:main_jacsketch}\\
\hline
\end{tabular}
\end{center}
\end{table}

For existing methods we provide a citation; new methods developed in this paper are marked accordingly.  We provide a link to the appropriate section for easy navigation. While these details are important, the main message of this chapter, i.e., the generality of our approach, is captured by Table~\ref{tbl:special_cases2}. The ``Result'' column of Table~\ref{tbl:special_cases2} points to a corollary of Theorem~\ref{thm:main_gsgm}; these corollaries state in detail the convergence  statements for the various methods. In all cases where known methods are recovered, these corollaries of Theorem~\ref{thm:main_gsgm} recover the best known rates.

{\bf Parameters.} From the point of view of Assumption~\ref{as:general_stoch_gradient}, the methods listed in Table~\ref{tbl:special_cases2} exhibit certain patterns.  To shed some light on this, in Table~\ref{tbl:special_cases-parameters} we summarize the values of these parameters. 

\begin{table}[H]
\caption{The parameters for which the methods from Table~\ref{tbl:special_cases2} (special cases of~\eqref{eq:SGD}) satisfy Assumption~\ref{as:general_stoch_gradient}. The meaning of the expressions appearing in the table, as well as their justification is defined in detail in Section~\ref{sec:special_cases}. }
\label{tbl:special_cases-parameters}
\begin{center}
\begin{tabular}{|c|c|c|c|c|c|c|}
\hline
 Method &   $A$ & $B$ & $\rho$ & $C$ & $D_1$ & $D_2$\\
\hline
 {\tt SGD}   &  $2L $ & $0$ & $1$ & $0$ & $2\sigma^2$ & $0$ \\
 {\tt SGD-SR} &  $2 \cL$ & $0$ & $1$ & $0$ & $2\sigma^2$ & $0$ \\
  {\tt SGD-MB}  &  $\frac{A' + L(\tau-1)}{\tau}$ & 0 & $1$ & $0$ & $ \frac{D'}{\tau}$ & $0$ \\
  {\tt SGD-star}  &  $2 \cL$ & $0$ & $1$ & $0$ & $0$ & $0$ \\
 {\tt SAGA}   &  $2L$ & $2$ & $1/n$ & $L/n$ & $0$ & $0$\\
{\tt N-SAGA}   &  $2L$ & $2$ & $1/n$ & $L/n$ & $2\sigma^2$ & $\frac{\sigma^2}{n}$\\
{\tt SEGA}   &   $2dL$ & $2d$ & $1/d$ & $L/d$ & $0$ & $0$\\
{\tt N-SEGA}  &   $2dL$ & $2d$ & $1/d$ & $L/d$ & $2d\sigma^2$ & $\frac{\sigma^2}{d}$\\
 {\tt SVRG}${}^{a}$   & $2L$  & $2$ & $0$ & $0$& $0$ & $0$ \\
 {\tt L-SVRG}   &  $2 L$ & $2$ & $p$ & $Lp$ & $0$ & $0$ \\
 {\tt DIANA}   & $\left(1+\frac{2\omega}{n}\right)L$ & $\frac{2\omega}{n}$ & $\alpha$ & $L\alpha$ & $\frac{(1+\omega)\sigma^2}{n}$ & $\alpha\sigma^2$  \\
 {\tt DIANA}${}^b$  &  $\left(1+2\omega\right)L$ & $2\omega$ & $\alpha$ & $L\alpha$ & $0$ & $0$ \\
{\tt Q-SGD-SR} &  $2(1+\omega)\cL$ & $0$ & $1$ & $0$ & $2(1+\omega)\sigma^2$ & $0$ \\
 {\tt VR-DIANA} &  $\left(1+\frac{4\omega + 2}{n}\right)L$ & $\frac{2(\omega+1)}{n}$ & $\alpha$ & $\left(\frac{1}{m}+4\alpha\right)L$ & $0$ & $0$ \\
 {\tt JacSketch}   &  $2\cL_1$ & $\frac{2\lambda_{\max}}{n}$ & $\lambda_{\min}$ & $\frac{\cL_2}{n}$ & $0$ & $0$ \\
\hline
\end{tabular}
\end{center}
\end{table}

Note, for example, that for all methods the parameter $A$ is non-zero. Typically, this a multiple of an appropriately defined smoothness parameter (e.g., $L$ is the Lipschitz constant of the gradient of $f$, $\cL$ and $\cL_1$ in {\tt SGD-SR}\footnote{{\tt SGD-SR} is first {\tt SGD} method analyzed in the {\em arbitrary sampling} paradigm. It was developed using the  {\em stochastic reformulation} approach (whence the ``SR'') pioneered in \citep{richtarik2020stochastic} in a numerical linear algebra setting, and later extended to develop the {\tt JacSketch} variance-reduction technique for finite-sum optimization  \citep{gower2018stochastic}.}, {\tt SGD-star} and {\tt JacSketch} are {\em expected smoothness} parameters). In the three variants of the {\tt DIANA} method, $\omega$ captures the variance of the  quantization operator $Q$. That is, one assumes that $\Exp{Q(x)}=x$ and $\Exp{\norm{Q(x)-x}^2} \leq \omega \norm{x}^2$ for all $x\in \R^d$. In view of \eqref{eq:gamma_condition_gsg}, large $A$ means a smaller stepsize, which slows down the rate. Likewise, the variance $\omega$ also affects the parameter $B$, which in view of \eqref{eq:main_result_gsgm} also has an adverse effect on the rate.
Further, as predicted by Theorem~\ref{thm:main_gsgm}, whenever either $D_1>0$ or $D_2>0$, the corresponding method converges to an oscillation region only. These methods are not variance-reduced. All symbols used in Table~\ref{tbl:special_cases-parameters} are defined in the appendix, in the same place where the methods are described and analyzed.

{\bf Five new methods.} To illustrate the usefulness of our general framework, we develop  {\em 5 new variants} of {\tt SGD} never explicitly considered in the literature before (see Table~\ref{tbl:special_cases2}). Here we briefly motivate them; details can be found in the corresponding sections. 

$\bullet$ {\tt SGD-MB} (Algorithm~\ref{alg:SGD-MB}).
This method is specifically designed for functions of the finite-sum structure  \eqref{eq:f_i_sum}. As we show through experiments, this is a powerful mini-batch {\tt SGD} method, with mini-batches formed with replacement as follows: in each iteration, we repeatedly ($\tau$ times) and independently pick $i\in [n]$ with probability $p_i>0$. Stochastic gradient $g^k$ is then formed by averaging the stochastic gradients $\nabla f_i(x^k)$ for all selected indices $i$ (including each $i$ as many times as this index was selected). This allows for a more practical importance mini-batch sampling implementation than what was until now possible (see Remark~\ref{rem:sgdmb} for more details and experiment in Figure~\ref{fig:SGDMB}).

$\bullet$ {\tt SGD-star} (Algorithm~\ref{alg:SGD-star}). 
This new method forms a bridge between vanilla and variance-reduced {\tt SGD} methods. While  not practical, it sheds light on the role of variance reduction. Again, we consider functions of the finite-sum form~\eqref{eq:f_i_sum}. This methods answers the following question:  assuming that the gradients $\nabla f_i(x^*)$, $i\in [n]$ are {\em known}, can they be used to design a more powerful {\tt SGD} variant? The answer is {\em yes}, and {\tt SGD-star} is the method. In its most basic form, {\tt SGD-star}  constructs the stochastic gradient via $g^k=\nabla f_i(x^k) - \nabla f_i(x^*)+\nabla f(x^*)$, where $i\in [n]$ is chosen uniformly at random. 
Inferring from Table~\ref{tbl:special_cases-parameters}, where $D_1=D_2=0$, this method converges to $x^*$, and not merely to some oscillation region.  Variance-reduced methods essentially work by iteratively constructing increasingly more accurate {\em estimates} of $\nabla f_i(x^*)$. Typically, the term $\sigma_k^2$ in the Lyapunov function of variance reduced methods will contain a term of the form $\sum_i \norm{h_i^k - \nabla f_i(x^*)}^2$, with $h_i^k$ being the estimators maintained by the method. Remarkably, {\tt SGD-star} was never explicitly considered in the literature before.

$\bullet$ {\tt N-SAGA} (Algorithm~\ref{alg:N-SAGA}). This is a novel variant of {\tt SAGA}~\citep{SAGA}, one in which one does not have access to the gradients of $f_i$, but instead only has access to {\em noisy} stochastic estimators thereof (with noise $\sigma^2$). Like {\tt SAGA}, {\tt N-SAGA} is able to reduce the variance inherent in the finite sum structure~\eqref{eq:f_i_sum} of the problem. However, it necessarily pays the price of noisy estimates of $\nabla f_i $, and hence, just like vanilla {\tt SGD} methods, is ultimately unable to converge to $x^*$. The oscillation region is governed by the noise level $\sigma^2$ (refer to $D_1$ and $D_2$ in Table~\ref{tbl:special_cases-parameters}). This method will be of practical importance for problems where each $f_i$ is of the form~\eqref{eq:f_exp}, i.e., for problems of the ``average of expectations'' structure. Batch versions of {\tt N-SAGA} would be well suited for distributed optimization, where each $f_i$ is owned by a different worker, as in such a case one wants the workers to work in parallel.

$\bullet$ {\tt N-SEGA} (Algorithm~\ref{alg:N-SEGA}). This is a {\em noisy} extension of the {\tt RCD}-type method {\tt SEGA}, in complete analogy with the relationship between {\tt SAGA} and {\tt N-SAGA}. Here we assume that we only have noisy estimates of partial derivatives (with noise $\sigma^2$). This situation is common in derivative-free optimization, where such a noisy estimate can be obtained by taking (a random) finite difference approximation~\citep{nesterov2017randomDFO}. Unlike {\tt SEGA}, {\tt N-SEGA} only converges to an oscillation region the size of which is governed by $\sigma^2$.

$\bullet$ {\tt Q-SGD-SR} (Algorithm~\ref{alg:qsgdas}). This is a quantized version of {\tt SGD-SR}, which is the first {\tt SGD} method analyzed in the  arbitrary sampling paradigm. As such, {\tt Q-SGD-SR}  is a vast generalization of the celebrated {\tt QSGD} method \citep{alistarh2017qsgd}.

\section{Special Cases}\label{sec:special_cases}

\subsection{Proximal {\tt SGD} for Stochastic Optimization}\label{sec:SGD}

\begin{algorithm}[H]
    \caption{{\tt SGD}}
    \label{alg:sgd_prox}
    \begin{algorithmic}
        \Require learning rate $\gamma>0$, starting point $x^0\in\R^d$, distribution $\cD$ over $\xi $
        \For{ $k=0,1,2,\ldots$ }
        \State{Sample $\xi \sim \cD$}
        \State{$g^k = \nabla f_\xi (x^k)$}
        \State{$x^{k+1} = \proxR(x^k - \gamma g^k)$}
        \EndFor
    \end{algorithmic}
\end{algorithm}

We start with stating the problem, the assumptions on the objective and on the stochastic gradients for {\tt SGD} \citep{nguyen2018sgd}. Consider the expectation minimization problem
\begin{equation}\label{eq:main_sgd}
    \min_{x\in\R^d} f(x) + R(x),\quad f(x) \eqdef \ED{f_\xi(x)}
\end{equation}
where $\xi \sim \cD$, $f_\xi(x)$ is differentiable and $L$-smooth almost surely in $\xi$.

 Lemma~\ref{lem:lemma1_sgd} shows that the stochastic gradient $g^k = \nabla f_\xi(x^k)$ satisfies Assumption~\ref{as:general_stoch_gradient}. The corresponding choice of parameters can be found in Table~\ref{tbl:special_cases-parameters}. 

\begin{lemma}[Generalization of Lemmas 1,2 from~\citep{nguyen2018sgd}]\label{lem:lemma1_sgd}
    Assume that $f_\xi(x)$ is convex in $x$ for every $\xi$. Then for every $x\in\R^d$
    \begin{equation}\label{eq:lemma1_sgd}
        \ED{\norm{\nabla f_\xi(x)- \nabla f(x^*)}^2} \le 4L(D_f(x,x^*)) + 2\sigma^2,
    \end{equation}
    where $\sigma^2\eqdef \EE_\xi\left[\norm{\nabla f_\xi(x^*)}^2\right]$. If further $f(x)$ is $\mu$-strongly convex with possibly non-convex $f_\xi$, then for every $x\in\R^d$
    \begin{equation}\label{eq:lemma2_sgd}
        \ED{\norm{\nabla f_\xi(x) - \nabla f(x^*)}^2} \le 4L\kappa(D_f(x,x^*)) + 2\sigma^2,
    \end{equation}
    where $\kappa = \frac{L}{\mu}$.
\end{lemma}

\begin{corollary}\label{cor:recover_sgd_rate}
    Assume that $f_\xi(x)$ is convex in $x$ for every $\xi$ and $f$ is $\mu$-strongly quasi-convex. Then {\tt SGD} with $\gamma \le \frac{1}{2L}$ satisfies
    \begin{equation}\label{eq:recover_sgd_rate}
        \EE\left[\norm{x^k - x^*}^2\right] \le (1-\gamma\mu)^k\norm{x^0-x^*}^2 + \frac{2\gamma \sigma^2}{\mu}.
    \end{equation}
    If we further assume that $f(x)$ is $\mu$-strongly convex with possibly non-convex $f_\xi(x)$,   {\tt SGD} with $\gamma \le \frac{1}{2L\kappa}$ satisfies~\eqref{eq:recover_sgd_rate} as well. 
\end{corollary}
\begin{proof}
    It suffices to plug parameters from Table~\ref{tbl:special_cases-parameters} into Theorem~\ref{thm:main_gsgm}. 
\end{proof}

\subsubsection*{Proof of Lemma~\ref{lem:lemma1_sgd}}
The proof is a direct generalization to the one from~\citep{nguyen2018sgd}.
Note that
\begin{eqnarray*}
&& \frac12 \ED{\norm{ \nabla f_\xi (x) - \nabla f(x^*)}^2} - \ED{\norm{\nabla f_\xi (x^*) - \nabla f(x^*)}^2} \\
&& \qquad \qquad=
\frac12 \ED{\norm{ \nabla f_\xi (x) - \nabla f(x^*)}^2 - \norm{\nabla f_\xi (x^*) - \nabla f(x^*)}^2} \\
&& \qquad \qquad\overset{\eqref{eq:1/2a_minus_b}}{\le}
 \ED{\norm{ \nabla f_\xi (x) - \nabla f_{\xi}(x^*)}^2}  \\
&& \qquad \qquad\leq
2LD_f(x,x^*).
\end{eqnarray*}
It remains to rearrange the above to get~\eqref{eq:lemma1_sgd}. To obtain~\eqref{eq:lemma2_sgd}, we shall proceed similarly:
\begin{eqnarray*}
&& \frac12 \ED{\norm{ \nabla f_\xi (x) - \nabla f(x^*)}^2} - \ED{\norm{ \nabla f_\xi (x^*) - \nabla f(x^*)}^2} \\
&& \qquad \qquad=
\frac12 \ED{\norm{ \nabla f_\xi (x) - \nabla f(x^*)}^2 - \norm{ \nabla f_\xi (x^*) - \nabla f(x^*)}^2} \\
&& \qquad \qquad\overset{\eqref{eq:1/2a_minus_b}}{\le}
 \ED{\norm{\nabla f_\xi (x) - \nabla f_{\xi}(x^*)}^2}  \\
&& \qquad \qquad\leq
L^2\norm{ x - x^*}^2 \\ 
&& \qquad \qquad\leq
2\frac{L^2}{\mu}D_f(x,x^*).
\end{eqnarray*}
Again, it remains to rearrange the terms.

\subsection{{\tt SGD-SR}}\label{SGD-AS}

In this section, we recover convergence result of {\tt SGD} under expected smoothness property from~\citep{gower2019sgd}. This setup allows obtaining tight convergence rates of {\tt SGD} under arbitrary stochastic reformulation of finite sum minimization\footnote{For technical details on how to exploit expected smoothness for specific reformulations, see~\citep{gower2019sgd}}. 

The stochastic reformulation is a special instance of~\eqref{eq:main_sgd}:
\begin{equation}\label{eq:problem_sgd-as}
    \min\limits_{x\in\R^d}f(x) + R(x),\quad f(x) =\ED{f_\xi(x)}, \quad  f_\xi(x) \eqdef \frac1n \sum_{i=1}^n \xi_i f_i(x) 
\end{equation}
where $\xi$ is a random vector from distribution $\cD$ such that for all $i$: $\ED{\xi_i}=1$ and $f_i$ (for all $i$) is smooth, possibly non-convex function. We next state the expextes smoothness assumption. A specific instances of this assumption allows to get tight convergence rates of SGD, which we recover in this section.

\begin{algorithm}[H]
    \caption{{\tt SGD-SR}}
    \label{alg:sgdas}
    \begin{algorithmic}
        \Require learning rate $\gamma>0$, starting point $x^0\in\R^d$, distribution $\cD$ over $\xi \in\R^n$ such that $\ED{\xi}$ is vector of ones
        \For{ $k=0,1,2,\ldots$ }
        \State{Sample $\xi \sim \cD$}
        \State{$g^k = \nabla f_\xi (x^k)$}
        \State{$x^{k+1} = \proxR(x^k - \gamma g^k)$}
        \EndFor
    \end{algorithmic}
\end{algorithm}

\begin{assumption}[Expected smoothness]\label{as:exp_smoothness_sgd-as}
    We say that $f$ is $\cL$-smooth in expectation with respect to distribution $\cD$ if there exists $\cL = \cL(f,\cD) > 0$ such that
    \begin{equation}\label{eq:exp_smoothness_sgd-as}
        \ED{\norm{\nabla f_\xi(x) - \nabla f_\xi(x^*)}^2} \le 2\cL D_f(x,x^*),
    \end{equation}
    for all $x\in\R^d$. For simplicity, we will write $(f,\cD) \sim ES(\cL)$ to say that \eqref{eq:exp_smoothness_sgd-as} holds.
\end{assumption}

Next, we present Lemma~\ref{lem:exp_smoothness_grad_up_bound_sgd-as} which shows that choice of constants for Assumption~\ref{as:general_stoch_gradient} from Table~\ref{tbl:special_cases-parameters} is valid. 

\begin{lemma}[Generalization of Lemma~2.4, \citep{gower2019sgd}]\label{lem:exp_smoothness_grad_up_bound_sgd-as}
    If $(f,\cD)\sim ES(\cL)$, then
    \begin{equation}\label{eq:exp_smoothness_grad_up_bound_sgd-as}
        \ED{\norm{\nabla f_\xi(x) - \nabla f(x^*)}^2} \le 4\cL D_f(x,x^*) + 2\sigma^2.
    \end{equation}
    where $\sigma^2 \eqdef \ED{\norm{\nabla f_\xi(x^*) - \nabla f(x^*	)}^2}$.
\end{lemma}

A direct consequence of Theorem~\ref{thm:main_gsgm} in this setup is Corollary~\ref{cor:recover_sgd-as_rate}.

\begin{corollary}\label{cor:recover_sgd-as_rate}
    Assume that $f(x)$ is $\mu$-strongly quasi-convex and $(f,\cD)\sim ES(\cL)$. Then {\tt SGD-SR} with $\gamma^k \equiv\gamma \le \frac{1}{2\cL}$ satisfies
    \begin{equation}\label{eq:recover_sgd-as_rate}
        \EE\left[\norm{x^k - x^*}^2\right] \le (1-\gamma\mu)^k\norm{x^0-x^*}^2 + \frac{2\gamma\sigma^2}{\mu}.
    \end{equation}
\end{corollary}


\subsubsection*{Proof of Lemma~\ref{lem:exp_smoothness_grad_up_bound_sgd-as}}
Here we present the generalization of the proof of Lemma~2.4 from \citep{gower2019sgd} for the case when $\nabla f(x^*) \neq 0$. In this proof all expectations are conditioned on $x^k$.
\begin{eqnarray*}
    \EE\left[\norm{\nabla f_{\xi} (x) - \nabla f(x^*)}^2\right] &=& \EE\left[ \norm{ \nabla f_\xi(x) - \nabla f_{\xi}(x^*) + \nabla f_{\xi}(x^*) - \nabla f(x^*) }^2\right] \\
    &\overset{\eqref{eq:a_b_norm_squared}}{\le}&  2 \EE\left[\norm{ \nabla f_{\xi}(x) - \nabla f_{\xi}(x^*)}^2\right] 
 + 2  \EE\left[\norm{\nabla f_{\xi}(x^*) - \nabla f(x^*)}^2\right] \notag\\ 
    &\overset{\eqref{eq:exp_smoothness_sgd-as}}{\le}& 4 \cL D_f(x,x^*) + 2 \sigma^2.
\end{eqnarray*}

\subsection{{\tt SGD-MB}}\label{sec:SGD-MB}
In this section, we present a specific practical formulation of~\eqref{eq:problem_sgd-as} which was not considered in~\citep{gower2019sgd}. The resulting algorithm (Algorithm~\ref{alg:SGD-MB}) is novel; it was not considered in~\citep{gower2019sgd} as a specific instance of {\tt SGD-SR}. The key idea behind {\tt SGD-MB} is constructing unbiased gradient estimate via with-replacement sampling.

Consider random variable $\nu \sim \cD$ such that 
\begin{equation}\label{eq:nb87fg87f}
 \Prob(\nu = i) =p_i;  \qquad \sum_{i=1}^np_i=1.
\end{equation} Notice that if we define \begin{equation}
\label{eq:reform_function}\psi_i(x)\eqdef \frac{1}{n p_i} f_i(x), \qquad i=1,2,\dots,n,\end{equation}
then 
\begin{equation} \label{eq:reform_problem}f(x) = \frac{1}{n} \sum_{i=1}^n f_i(x)  \overset{\eqref{eq:reform_function}}{=} \sum_{i=1}^n p_i \psi_i(x) \overset{\eqref{eq:nb87fg87f}}{=} \ED{\psi_\nu (x)}.\end{equation}
So, we have rewritten the finite sum problem \eqref{eq:f_sum} into the {\em equivalent stochastic optimization problem}  
\begin{equation}\label{eq:reform_opt} \min_{x\in \R^d} \ED{\psi_\nu (x)}.\end{equation}

We are now ready to describe our method. At each iteration $k$ we sample $\nu^k_i,\dots,\nu^k_{\tau}\sim \cD$ independently ($1\leq \tau \leq n$), and define $g^k\eqdef  \frac{1}{\tau}\sum_{i=1}^\tau \nabla \psi_{\nu^k_i}(x^k) $. Further, we use $g^k$ as a stochastic gradient, resulting in Algorithm~\ref{alg:SGD-MB}.

\begin{algorithm}[H]
    \caption{{\tt SGD-MB}}
    \label{alg:SGD-MB}
    \begin{algorithmic}
        \Require learning rate $\gamma>0$, starting point $x^0\in\R^d$, distribution $\cD$ over $\nu$ such that~\eqref{eq:nb87fg87f} holds.
        \For{ $k=0,1,2,\ldots$ }
        \State{Sample $\nu^k_i,\dots,\nu^k_{\tau}\sim \cD$ independently}
        \State{$g^k =  \frac{1}{\tau}\sum_{i=1}^\tau \nabla \psi_{\nu^k_i}(x^k) $}
        \State{$x^{k+1} = x^k - \gamma g^k$}
        \EndFor
    \end{algorithmic}
\end{algorithm}

To remain in full generality, consider the following Assumption.

\begin{assumption}\label{ass:main_mb} There exists constants $A'>0$ and $D'\geq 0$ such that \begin{equation}\label{eq:ABassumpt} \ED{ \norm{\nabla \psi_{\nu}(x)}^2 }\leq 2 A' (f(x) - f(x^*)) + D'\end{equation} for all $x\in \R^d$.  
\end{assumption}

Note that it is sufficient to have convex and smooth $f_i$ in order to satisfy Assumption~\ref{ass:main_mb}, as Lemma~\ref{lem:b98h0f} states. 

\begin{lemma}\label{lem:b98h0f}  Let $\sigma^2\eqdef \ED{\norm{\nabla \psi_\nu (x^*)}^2}$. If $f_i$ are convex and $L_i$-smooth, then Assumption~\ref{ass:main_mb} holds for $A'=2\cL$ and $D'=2\sigma^2$, where \begin{equation}\label{eq:ineq90f9f}\cL \leq \max_i \frac{L_i}{n p_i}.\end{equation} 
If moreover $\nabla f_i(x^*)=0$ for all $i$,  then Assumption~\ref{ass:main_mb} holds for $A'=\cL$ and $D'=0$. 
\end{lemma}

Next, Lemma~\ref{lem:mb} states that Algorithm~\ref{alg:SGD-MB} indeed satisfies Assumption~\ref{as:general_stoch_gradient}. 

\begin{lemma} \label{lem:mb}
Suppose that Assumption~\ref{ass:main_mb} holds. Then $g^k$ is unbiased; i.e. $\ED{g^k} = \nabla f(x^k)$. Further, 
\begin{eqnarray*}
\ED{\norm{g^k}^2} \leq  \frac{2A' + 2L(\tau-1)}{\tau}(f(x^k) - f(x^*)) +  \frac{  D'}{\tau}.
\end{eqnarray*}
\end{lemma}
Thus, parameters from Table~\ref{tbl:special_cases-parameters} are validated. As a direct consequence of Theorem~\ref{thm:main_gsgm} we get Corollary~\ref{cor:mb}.

\begin{corollary}\label{cor:mb}
    As long as $0< \gamma \leq \frac{\tau}{A' + L(\tau-1)}$, we have
    \begin{equation}\label{eq:mb_rate}
\EE \norm{x^k-x^*}^2 \leq (1-\gamma \mu)^k \norm{x^0-x^*}^2 + \frac{\gamma D'}{\mu \tau}.    \end{equation}
\end{corollary}

\begin{remark} \label{rem:sgdmb}
For $\tau=1$, {\tt SGD-MB} is a special of the method from~\citep{gower2019sgd}, Section~3.2. However, for $\tau>1$, this is a different method; the difference lies in the with-replacement sampling. Note that with-replacement trick allows for efficient and implementation of independent importance sampling~\footnote{Distribution of random sets $S$ for which random variables $i\in S$ and $j\in S$ are independent for $j\neq i$.} with complexity $\cO(\tau \log(n))$. In contrast, implementation of without-replacement importance sampling has complexity $\cO(n)$, which can be significantly more expensive to the cost of evaluating $\sum_{i\in S} \nabla f_i(x)$.
 \end{remark}

\subsubsection*{Proof of Lemma~\ref{lem:mb}}
Notice first that
\begin{eqnarray*}
\ED{g^k}
&\overset{\eqref{eq:reform_function}}{=}&
  \frac{1}{\tau}\sum_{i=1}^\tau\ED{\frac{1}{n p_{\nu^k_{i}}} \nabla f_{\nu^k_i}(x^k)}
  \\
&=& 
\ED{\frac{1}{n p_{\nu}} \nabla f_{\nu}(x^k)}
  \\
&\overset{\eqref{eq:nb87fg87f}}{=}& 
\sum_{i=1}^n p_i \frac{1}{n p_i} \nabla f_i(x^k) 
\\
&=& \nabla f(x_k).
\end{eqnarray*}

So, $g^k$ is an unbiased estimator of the gradient $\nabla f(x^k)$. Next, 
\begin{eqnarray*}
\ED{\norm{g^k}^2}&=& \ED{\norm{\frac{1 }{\tau}\sum_{i=1}^\tau \nabla \psi_{\nu^k_i}(x^k)}^2 }\\ 
&=& \frac{1}{\tau^2}\ED{\sum_{i=1}^\tau \norm{\nabla \psi_{\nu^k_i}(x^k)}^2  +  2    \sum_{i< j} \< \nabla \psi_{\nu^k_i}(x^k),  \nabla \psi_{\nu^k_j}(x^k)>  } \\
&=& \frac{1}{\tau}\ED{ \norm{\nabla \psi_{\nu}(x^k)}^2}+  \frac{2}{\tau^2}    \sum_{i< j} \< \ED{\nabla \psi_{\nu^k_i}(x^k) },  \ED{\nabla \psi_{\nu^k_j}(x^k) }>     \\
&=& \frac{1}{\tau} \ED{\norm{\nabla \psi_{\nu}(x^k)}^2} + \frac{\tau - 1}{\tau} \norm{\nabla f(x^k)}^2 \\
&\overset{\eqref{eq:ABassumpt} }{\leq}& \frac{2A' (f(x^k) - f(x^*)) + D' + 2L(\tau-1)(f(x^k) - f(x^*))}{\tau}.
\end{eqnarray*}

\subsubsection*{Proof of Lemma~\ref{lem:b98h0f}}
Let $\cL = \cL(f,\cD)>0$ be any constant for which
\begin{equation}\label{eq:ES}\EE_{\xi\sim \cD} \norm{\nabla \phi_{\xi}(x) - \nabla \phi_{\xi}(x^*)}^2 \leq 2\cL (f(x) - f(x^*))\end{equation}
 holds for all $x\in \R^d$. This is the expected smoothness property (for a single item sampling) from \citep{gower2019sgd}. It was shown in \cite[Proposition 3.7]{gower2019sgd}  that \eqref{eq:ES} holds, and that $ \cL$ satisfies \eqref{eq:ineq90f9f}. The claim now follows by applying \cite[Lemma~2.4]{gower2019sgd}. 

\subsection{{\tt SGD-star}}\label{sec:SGD-star}
Consider problem~\eqref{eq:problem_sgd-as}. Suppose that $\nabla f_i(x^*)$ is known for all $i$. In this section we present a novel algorithm~---~{\tt SGD-star}~---~which is {\tt SGD-SR} shifted by the stochastic gradient in the optimum. The method is presented under Expected Smoothness Assumption~\eqref{eq:exp_smoothness_sgd-as}, obtaining general rates under arbitrary sampling.  The algorithm is presented as Algorithm~\ref{alg:SGD-star}. 

\begin{algorithm}[H]
    \caption{{\tt SGD-star}}
    \label{alg:SGD-star}
    \begin{algorithmic}
        \Require learning rate $\gamma>0$, starting point $x^0\in\R^d$, distribution $\cD$ over $\xi \in \R^n$ such that $\ED{\xi}$ is vector of ones
        \For{ $k=0,1,2,\ldots$ }
        \State{Sample $\xi \sim \cD$}
        \State{$g^k = \nabla f_\xi(x^k) - \nabla f_\xi(x^*) + \nabla f(x^*)$}
        \State{$x^{k+1} = \proxR(x^k - \gamma g^k)$}
        \EndFor
    \end{algorithmic}
\end{algorithm}
Suppose that $(f,\cD) \sim ES(\cL)$. 
Note next that {\tt SGD-star} is just  {\tt SGD-SR} applied on objective $D_f(x,x^*)$ instead of $f(x)$ when $\nabla f(x^*) = 0$. This careful design of the objective yields $(D_f(\cdot, x^*),\cD) \sim ES(\cL)$ and $\ED{\norm{\nabla_x D_{f_\xi}(x,x^*)}^2 \, \mid x=x^*} = 0 $, and thus Lemma~\eqref{lem:exp_smoothness_grad_up_bound_sgd-as} becomes
\begin{lemma}[Lemma~2.4, \citep{gower2019sgd}]\label{lem:SGD-star}
    If $(f,\cD)\sim ES(\cL)$, then
    \begin{equation}\label{eq:sgd-star_lemma}
        \ED{\norm{g^k -\nabla f(x^*)}^2} \le 4\cL D_f(x^k,x^*).
    \end{equation}
\end{lemma}

A direct consequence of Corollary (thus also a direct consequence of Theorem~\ref{thm:main_gsgm}) in this setup is Corollary~\ref{cor:SGD-star}.

\begin{corollary}\label{cor:SGD-star}
    Suppose that $(f,\cD)\sim ES(\cL)$. Then {\tt SGD-star} with $\gamma = \frac{1}{2\cL}$ satisfies
    \begin{equation}\label{eq:SGD-shift-xx}
        \EE\left[\norm{x^k - x^*}^2\right] \le \left(1-\frac{\mu}{2\cL} \right)^k\norm{x^0-x^*}^2.
    \end{equation}
\end{corollary}


\begin{remark}
Note that results from this section are obtained by applying results from~\ref{SGD-AS}. Since Section~\ref{sec:SGD-MB} presets a specific sampling algorithm for {\tt SGD-SR}, the results can be thus extended to {\tt SGD-star} as well.
\end{remark}

\subsubsection*{Proof of Lemma~\ref{lem:SGD-star}}
In this proof all expectations are conditioned on $x^k$.
\begin{eqnarray*}
    \ED{\norm{g^k -\nabla f(x^*)}^2}  &=& \ED{\norm{\nabla f_\xi(x^k) -\nabla f_\xi(x^*)}^2} \\ 
    &\overset{\eqref{eq:exp_smoothness_sgd-as}}{\le}& 4 \cL D_f(x^k,x^*).
\end{eqnarray*}

\subsection{{\tt SAGA}}\label{sec:saga}
In this section we show that our approach is suitable for {\tt SAGA} \citep{SAGA} (see Algorithm~\ref{alg:SAGA}). Consider the finite-sum minimization problem
\begin{equation}\label{eq:main_l-svrg}
    f(x) = \frac{1}{n}\sum\limits_{i=1}^n f_i(x) + R(x),
\end{equation}
where $f_i$ is convex, $L$-smooth for each $i$ and $f$ is $\mu$-strongly convex. 

\begin{algorithm}[H]
    \caption{{\tt SAGA} \citep{SAGA}}
    \label{alg:SAGA}
    \begin{algorithmic}
        \Require learning rate $\gamma>0$, starting point $x^0\in\R^d$
        \State Set $\phi_j^0 = x^0$ for each $j\in[n]$
        \For{ $k=0,1,2,\ldots$ }
        \State{Sample  $j \in [n]$ uniformly at random}
        \State{Set $\phi_j^{k+1} = x^k$ and $\phi_i^{k+1} = \phi_i^{k}$ for $i\neq j$}
        \State{$g^k = \nabla f_j(\phi_j^{k+1}) - \nabla f_j(\phi_j^k) + \frac{1}{n}\sum\limits_{i=1}^n\nabla f_i(\phi_i^k)$}
        \State{$x^{k+1} = \proxR\left(x^k - \gamma g^k\right)$}
        \EndFor
    \end{algorithmic}
\end{algorithm}

\begin{lemma}\label{lem:stoch_grad_second_moment_saga}
    We have
    \begin{equation}\label{eq:stoch_grad_second_moment_saga1}
        \EE\left[\norm{g^k- \nabla f(x^*)}^2\mid x^k\right] \le 4LD_f(x^k,x^*)+ 2\sigma_k^2
    \end{equation}
    and 
        \begin{equation}\label{eq:stoch_grad_second_moment_saga2}
        \EE\left[\sigma_{k+1}^2\mid x^k\right] \le \left(1 - \frac{1}{n}\right)\sigma_k^2 + \frac{2L}{n}D_f(x^k,x^*),
    \end{equation}
    where $\sigma_k^2 = \frac{1}{n}\sum\limits_{i=1}^n\norm{\nabla f_i(\phi_i^k) - \nabla f_i(x^*)}^2$.
\end{lemma}

Clearly, Lemma~\ref{lem:stoch_grad_second_moment_saga} shows that Algorithm~\ref{alg:SAGA} satisfies Assumption~\ref{as:general_stoch_gradient}; the corresponding parameter choice can be found in Table~\ref{tbl:special_cases-parameters}. Thus, as a direct consequence of Theorem~\ref{thm:main_gsgm} with $M=4n$ we obtain the next corollary.

\begin{corollary}\label{thm:recover_saga_rate}
    {\tt SAGA} with $\gamma = \frac{1}{6L}$ satisfies
    \begin{equation}\label{eq:recover_saga_rate_coro}
        \EE V^k \le \left(1-\min\left\{\frac{\mu}{6L},\frac{1}{2n}\right\}\right)^kV^0.
    \end{equation}
\end{corollary}

\subsubsection*{Proof of Lemma~\ref{lem:stoch_grad_second_moment_saga}}
Note that  Lemma~\ref{lem:stoch_grad_second_moment_saga} is a special case of Lemmas 3,4 from~\citep{mishchenko202099} without prox term. We reprove it with prox for completeness.

    Let all expectations be conditioned on $x^k$ in this proof. Note that $L$-smoothness and convexity of $f_i$ implies
\begin{equation}\label{eq:norm_diff_grads}
    \frac{1}{2L}\norm{\nabla f_i(x) - \nabla f_i(y)}^2 \le f_i(x) - f_i(y) - \<\nabla f_i(y),x-y>,\quad \forall x,y\in\R^d, i\in[n].
\end{equation}

 By definition of $g^k$ we have
    \begin{eqnarray*}
        \EE\left[\norm{g^k - \nabla f(x^*)}^2\right] 
        &=& 
        \EE\left[ \norm{\nabla f_j(\phi_j^{k+1}) - \nabla f_j(\phi_j^k) + \frac{1}{n}\sum\limits_{i=1}^n\nabla f_i(\phi_i^k) - \nabla f(x^*)}^2\right]
        \\
        &\overset{\eqref{eq:a_b_norm_squared}}{\le}& 
        2\EE\left[\norm{\nabla f_j(x^k) - \nabla f_j(x^*)}^2\mid x^k\right]
        \\
        &&\quad + 2\EE\left[\norm{\nabla f_j(x^*) - \nabla f_j(\phi_j^k) - \EE\left[\nabla f_j(x^*) - \nabla f_j(\phi_j^k)\right] }^2\right]
        \\
        &\overset{\eqref{eq:variance_decomposition}+\eqref{eq:norm_diff_grads}}{\le}&
         \frac{4L}{n}\sum\limits_{i=1}^nD_{f_i}(x^k,x^*)+ 2\EE\left[\norm{\nabla f_j(x^*) - \nabla f_j(\phi_j^k)}^2\mid x^k\right]\\
        &=& 4LD_f(x^k,x^*) + 2\underbrace{\frac{1}{n}\sum\limits_{i=1}^n\norm{\nabla f_i(\phi_i^k) - \nabla f_i(x^*)}^2}_{\sigma_k^2}.
    \end{eqnarray*}
    
To proceed with~\eqref{eq:stoch_grad_second_moment_saga2}, we have
    \begin{eqnarray*}
        \EE\left[\sigma_{k+1}^2\right] &=& \frac{1}{n}\sum\limits_{i=1}^n\EE\left[\norm{\nabla f_i(\phi_i^{k+1}) - \nabla f_i(x^*)}^2\right]\\
        &=& \frac{1}{n}\sum\limits_{i=1}^n\left(\frac{n-1}{n}\norm{\nabla f_i(\phi_i^k) - \nabla f_i(x^*)}^2 + \frac{1}{n}\norm{\nabla f_i(x^k) - \nabla f_i(x^*)}^2\right)\\
        &\overset{\eqref{eq:norm_diff_grads}}{\le}& \left(1 - \frac{1}{n}\right)\frac{1}{n}\sum\limits_{i=1}^n\norm{\nabla f_i(\phi_i^k) - \nabla f_i(x^*)}^2\\
        &&\quad + \frac{2L}{n^2}\sum\limits_{i=1}^n D_{f_i}(x^k,x^*)\\
        &=& \left(1 - \frac{1}{n}\right)\sigma_k^2 + \frac{2L}{n}D_f(x^k,x^*).
    \end{eqnarray*}

\subsection{{\tt N-SAGA}} \label{N-SAGA}

\begin{algorithm}[H]
    \caption{Noisy {\tt SAGA} ({\tt N-SAGA})}
    \label{alg:N-SAGA}
    \begin{algorithmic}
        \Require learning rate $\gamma>0$, starting point $x^0\in\R^d$
        \State Set $\psi_j^0 = x^0$ for each $j\in[0]$
        \For{ $k=0,1,2,\ldots$ }
        \State{Sample  $j \in [n]$ uniformly at random and $\zeta$}
        \State{Set $g_j^{k+1} = g_j(x^k,\xi)$ and $g_i^{k+1} = g_i^{k} $ for $i\neq j$}
        \State{$g^k =g_j(x^k,\xi) - g_j^k + \frac{1}{n}\sum\limits_{i=1}^ng_i^k$}
        \State{$x^{k+1} = \prox_{\gamma R}(x^k - \gamma g^k)$}
        \EndFor
    \end{algorithmic}
\end{algorithm}

Note that it can in practice happen that instead of $\nabla f_i(x)$ one can query $g_i(x,\zeta)$ such that $\EE_\xi g_i(\cdot,\xi)=  \nabla f_i(\cdot) $ and $\EE_\xi \norm{g_i(\cdot,\xi)}^2 \leq  \sigma^2$. This leads to a variant of {\tt SAGA} which only uses noisy  estimates of the stochastic gradients $\nabla_i(\cdot)$. We call this variant {\tt N-SAGA} (see Algorithm~\ref{alg:N-SAGA}).

\begin{lemma}\label{lem:saga_inexact1}
    We have
    \begin{equation}\label{eq:stoch_grad_second_moment_saga_noise1}
        \EE\left[\norm{g^k - \nabla f(x^*)}^2\mid x^k\right] \le 4LD_f(x^k,x^*) + 2\sigma_k^2 + 2\sigma^2,
    \end{equation}
    and
        \begin{equation}\label{eq:stoch_grad_second_moment_saga_noise1}
        \EE\left[\sigma_{k+1}^2\mid x^k\right] \le \left(1 - \frac{1}{n}\right)\sigma_k^2 + \frac{2L}{n}D_f(x^k,x^*) + \frac{\sigma^2}{n},
    \end{equation}
    where $\sigma_k^2 \eqdef \frac{1}{n}\sum\limits_{i=1}^n\norm{ g_i^k - \nabla f_i(x^*)}^2$.
\end{lemma}

\begin{corollary}\label{cor:N-SAGA}
Let $\gamma = \frac{1}{6L}$. Then, iterates of Algorithm~\ref{alg:N-SAGA} satisfy
\[
\EE{V^k}\leq \left( 1- \min\left( \frac{\mu}{6L}, \frac{1}{2n} \right)\right)^k V^0 + \frac{\sigma^2 }{ L \min(\mu, \frac{3 L }{n})}.
\]
\end{corollary}

Analogous results can be obtained for {\tt L-SVRG}.

\subsubsection*{Proof of Lemma~\ref{lem:saga_inexact1}}

    Let all expectations be conditioned on $x^k$. By definition of $g^k$ we have
    \begin{eqnarray*}
        \EE\left[\norm{g^k - \nabla f(x^*)}^2\right] 
        &\le& 
        \EE\left[\norm{g_j(x^k,\zeta)- g_j^k + \frac{1}{n}\sum\limits_{i=1}^n g_i^k - \nabla f(x^*)}^2\right]
        \\
        &=&
         \EE\left[\norm{g_j(x^k,\zeta) - \nabla f_j(x^*) + \nabla f_j(x^*) - g_j^k + \frac{1}{n}\sum\limits_{i=1}^n g_i^k  - \nabla f(x^*)}^2\right]
         \\
        &\overset{\eqref{eq:a_b_norm_squared}}{\le}& 
        2\EE\left[\norm{g_j(x^k,\zeta) - \nabla f_j(x^*)}^2\right]
        \\
        &&\quad 
        + 2\EE\left[\norm{\nabla f_j(x^*) -g_j^k - \EE\left[\nabla f_j(x^*) - g_j^k\right] }^2\right]
        \\
        &\overset{\eqref{eq:variance_decomposition}}{\le}&
         2\EE\left[\norm{g_j(x^k,\zeta) - \nabla f_j(x^*)}^2\right] + 2\EE\left[\norm{\nabla f_j(x^*) - g_j^k }^2\right]
         \\
        &=& 
        2\EE\left[\norm{g_j(x^k,\zeta) - \nabla f_j(x^*)}^2\right] + 2\underbrace{\frac{1}{n}\sum\limits_{i=1}^n\norm{g_i^k - \nabla f_i(x^*)}^2}_{\sigma_k^2} 
        \\
        &\overset{\eqref{eq:variance_decomposition}}{\le}& 
                 2\EE\left[ \norm{ \nabla f_j(x^k) - \nabla f_j(x^*)}^2\right] + 2\sigma^2 + 2\sigma_k^2
                \\
                &\overset{\eqref{eq:norm_diff_grads}}{\le}&
                4LD_f(x^k,x^*) + 2\sigma_k^2 + 2\sigma^2
    \end{eqnarray*}

    For the second inequality, we have
    \begin{eqnarray*}
        \EE\left[\sigma_{k+1}^2\right] &=& \frac{1}{n}\sum\limits_{i=1}^n\EE\left[\norm{g_i^{k+1} - \nabla f_i(x^*)}^2\right]
        \\
        &=& 
        \frac{1}{n}\sum\limits_{i=1}^n\left(\frac{n-1}{n}\norm{g_i^k  - \nabla f_i(x^*)}^2 + \frac{1}{n}\EE\left[ \norm{ g_i(x^k, \zeta) - \nabla f_i(x^*)}^2\right]\right)
        \\
        &\leq&
                \frac{1}{n}\sum\limits_{i=1}^n\left(\frac{n-1}{n}\norm{g_i^k  - \nabla f_i(x^*)}^2 + \frac{1}{n} \norm{ \nabla f_i(x^k) - \nabla f_i(x^*)}^2 + \frac{\sigma^2}{n}\right)
        \\
        &\overset{\eqref{eq:norm_diff_grads}}{\le}&
        \left(1 - \frac{1}{n}\right)\sigma_k^2 + \frac{2L}{n}D_f(x^k,x^*)+ \frac{\sigma^2}{n}.
    \end{eqnarray*}

\subsection{{\tt SEGA}}\label{sec:sega}

\begin{algorithm}[H]
    \caption{{\tt SEGA} \citep{hanzely2018sega}}
    \label{alg:SEGA}
    \begin{algorithmic}
        \Require learning rate $\gamma>0$, starting point $x^0\in\R^d$
        \State Set $h^0 = 0$
        \For{ $k=0,1,2,\ldots$ }
        \State{Sample  $j \in [d]$ uniformly at random}
        \State{Set $h^{k+1} = h^{k} +e_i( \nabla_i f(x^k) - h^{k}_i)$}
        \State{$g^k = de_i (\nabla_i f(x^k) - h_i^k) + h^k$}
        \State{$x^{k+1} = \prox_{\gamma R}(x^k - \gamma g^k)$}
        \EndFor
    \end{algorithmic}
\end{algorithm}

We show that the framework recovers the simplest version of {\tt SEGA} (i.e., setup from Theorem D1 from~\citep{hanzely2018sega}) in the proximal setting\footnote{General version for arbitrary gradient sketches instead of partial derivatives can be recovered as well, however, we omit it for simplicity}. 

\begin{lemma} (Consequence of Lemmas A.3., A.4. from~\citep{hanzely2018sega})
We have
\[
\EE\left[\norm{g^{k}-\nabla f(x^*) \mid x^k}^{2}\right] \leq 2d\norm{\nabla f\left(x^{k}\right)-\nabla f(x^*)}^{2}+2d\sigma_k^2
\]
and
\[
\EE
\left[\sigma_{k+1}^2 \mid x^k \right]= \left(1-\frac1d\right)\sigma_k^2+\frac1d \norm{\nabla f\left(x^{k}\right)-\nabla f(x^*)}^{2},
\]
where $\sigma_k^2 \eqdef \norm{h^{k}-\nabla f(x^*)}^2$.
\end{lemma}

Given that we have from convexity and smoothness $\norm{\nabla f(x^{k})-\nabla f(x^*)}^{2} \leq 2L D_f(x^k,x^*)$, Assumption~\ref{as:general_stoch_gradient} holds the parameter choice as per Table~\ref{tbl:special_cases-parameters}. Setting further $M = 4d^2$, we get the next corollary.
\begin{corollary}\label{cor:sega}
{\tt SEGA} with $\gamma =\frac{1}{6dL} $ satisfies
\[
\EE V^k \leq \left( 1- \frac{\mu}{6dL}\right)^kV^0.
\] 
\end{corollary}

\subsection{{\tt N-SEGA}} \label{N-SEGA}

\begin{algorithm}[H]
    \caption{Noisy {\tt SEGA} ({\tt N-SEGA})}
    \label{alg:N-SEGA}
    \begin{algorithmic}
        \Require learning rate $\gamma>0$, starting point $x^0\in\R^d$
        \State Set $h^0 = 0$
        \For{ $k=0,1,2,\ldots$ }
        \State{Sample  $i \in [d]$ uniformly at random and sample $\xi$}
        \State{Set $h^{k+1} = h^{k} +e_i( g_i(x,\xi)- h^{k}_i)$}
        \State{$g^k = de_i (g_i(x,\xi) - h_i^k) +  h^k$}
        \State{$x^{k+1} = x^k - \gamma g^k$}
        \EndFor
    \end{algorithmic}
\end{algorithm}

Here we assume that $g_i(x,\zeta)$ is a noisy estimate of the partial derivative $\nabla_i f(x)$ such that $\EE_\zeta g_i(x,\zeta) = \nabla_i f(x)$ and $\EE_\zeta | g_i(x,\zeta) - \nabla_i f(x)|^2 \leq \frac{\sigma^2}{d}$. 

\begin{lemma} \label{lem:sega_noise}
The following inequalities hold:
\[
\EE\left[\norm{g^{k}-\nabla f(x^*)}^{2}\right]
 \leq 
 4dLD_f(x^k,x^*)+2d\sigma_k^2 +2d\sigma^2,
\]
\[
\EE
\left[\sigma_{k+1}^2\right] \leq \left(1-\frac1d\right)\sigma_k^2+\frac{2L}{d}D_f(x^k,x^*) + \frac{\sigma^2}{d}, 
\]
where $\sigma_k^2 = \norm{h^{k}-\nabla f(x^*)}^2$.
\end{lemma}

\begin{corollary}\label{cor:N-SEGA}
Let $\gamma = \frac{1}{6Ld}$. Applying Theorem~\ref{thm:main_gsgm} with $M = 4d^2$, iterates of Algorithm~\ref{alg:N-SEGA} satisfy
\[
\EE{V^k}\leq \left( 1- \frac{\mu}{6dL}\right)^k V^0 + \frac{\sigma^2 }{ L \mu}.
\]
\end{corollary}

\subsubsection*{Proof of Lemma~\ref{lem:sega_noise}}
 
Let all expectations be conditioned on $x^k$. 
For the first bound, we write \[g^k  - \nabla f(x^*)= \underbrace{h^k - \nabla f(x^*)- d  h^k_i e_{i}+d\nabla_i f(x^*)e_i}_{a} +   \underbrace{ d g_i(x^k,\xi)e_i- d\nabla_i f(x^*)e_i}_{b}.\]
 Let us bound the expectation of each term individually. The first term can be bounded as
\begin{eqnarray*}
\EE{\norm{a}^2} &=& \EE{\norm{ \left(\mI - d e_{i} e_{i}^\top \right) (h^k - \nabla f(x^*))  }_{2}^2}\\
&=& (d-1)\norm{h^k- \nabla f(x^*)}^2\\
&\leq& d\norm{h^k- \nabla f(x^*)}^2.
\end{eqnarray*}
The second term can be bounded as
\begin{eqnarray*}
\EE{\norm{b}^2} &=&  \EE_i \EE_{\xi}{\norm{ d g_i(x,\xi)e_i- d\nabla f_i(x^*)e_i}^2}\\
&=& 
 \EE_{i} \EE_{\xi} \norm{ d g_i(x^k,\xi)e_i- d\nabla_i f(x^k)e_i }^2 + \EE_i\norm{ d\nabla_i f(x^k)e_i-d\nabla f_i(x^*)e_i}^2 
\\
&\leq & d\sigma^2 + d  \norm{ \nabla f(x^k)- \nabla f(x^*)}^2 \\
&\leq& d\sigma^2 + 2Ld D_f(x^k,x^*),
\end{eqnarray*}
where in the last step we used $L$--smoothness of $f$. It remains to combine the two bounds.

For the second bound, we have
\begin{eqnarray*}
\EE{\norm{ h^{k+1}  - \nabla f(x^*)}^2} &=& \EE{ \norm{h^k + g_i(x^k,\xi)e_i - h^k_i - \nabla f(x^*) }^2 }\\
 &=&  \EE{ \norm{\left(\mI - e_{i} e_{i}^\top \right)h^k + g_i(x^k,\xi)e_i -\nabla f(x^*)  }^2 }\\
 &=& \EE{\norm{ \left(\mI - e_{i} e_{i}^\top \right) ( h^k - \nabla f(x^*)) }^2} + \EE{\norm{g_i(x^k,\xi)e_i - \nabla_i f(x^*)e_i }^2 }\\
  &=& \left(1-\frac{1}{d}\right) \norm{h^k - \nabla f(x^*)}^2 + \EE{\norm{g_i(x^k,\xi)e_i - \nabla_i f(x^k)e_i }^2 }   \\
  && \qquad + \EE{\norm{\nabla_i f(x^k)e_i  - \nabla_i f(x^*)e_i}^2 } \\
 &=& \left(1-\frac{1}{d}\right) \norm{h^k - \nabla f(x^*)}^2 + \frac{\sigma^2}{d} + \frac{1}{d} \norm{\nabla f(x^k) - \nabla f(x^*)}^2 \\
  &\leq&
   \left(1-\frac{1}{d}\right) \norm{h^k - \nabla f(x^*)}^2 + \frac{\sigma^2}{d} + \frac{2L}{d}D_f(x^k,x^*).
\end{eqnarray*}

\subsection{{\tt SVRG}} \label{sec:svrg}

\begin{algorithm}[H]
    \caption{{\tt SVRG} \citep{SVRG}}
    \label{alg:SVRG}
    \begin{algorithmic}
        \Require learning rate $\gamma>0$, epoch length $m$, starting point $x^0\in\R^d$
        \State  $\phi = x^0$
        \For{ $s=0,1,2,\ldots$ }
        \For{ $k=0,1,2,\ldots, m-1$ }
        \State{Sample  $i \in \{1,\ldots, n\}$ uniformly at random}
        \State{$g^k = \nabla f_i(x^k) - \nabla f_i(\phi) + \nabla f(\phi)$}
        \State{$x^{k+1} = \proxR(x^k - \gamma g^k)$}
        \EndFor
        \State  $\phi = x^0 = \frac1m \sum_{k=1}^m x^k$
        \EndFor
    \end{algorithmic}
\end{algorithm}
Let $\sigma_k^2 \eqdef \frac{1}{n}\sum\limits_{i=1}^n\norm{\nabla f_i(\phi) - \nabla f_i(x^*)}^2$. We will show that Lemma~\ref{lem:iter_dec} recovers per-epoch analysis of {\tt SVRG} in a special case.

\begin{lemma}\label{lem:svrg_lemma_1}
For $k \mod m \neq 0$ we have    
        \begin{equation}\label{eq:svrg1}
        \EE\left[\norm{g^k -\nabla f(x^*)}^2\mid x^k\right] \le 4LD_f(x^k,x^*) + 2\sigma_k^2
    \end{equation}
    and
    \begin{equation}\label{eq:svrg3}
        \EE\left[\sigma_{k+1}^2\mid x^k\right] = \sigma_{k+1}^2 = \sigma_k^2.
    \end{equation}
\end{lemma}
\begin{proof}
The proof of~\eqref{eq:svrg1} is identical to the proof of~\eqref{eq:stoch_grad_second_moment_saga1}. Next,~\eqref{eq:svrg3} holds since $\sigma_k$ does not depend on $k$. 
\end{proof}

Thus, Assumption~\ref{as:general_stoch_gradient} holds with parameter choice as per Table~\ref{tbl:special_cases-parameters} and Lemma~\ref{lem:iter_dec} implies the next corollary.
\begin{corollary}\label{cor:svrg}
\begin{equation}\label{eq:svrg_iter_dec}
 \EE\norm{x^{k+1}-x^*}^2 + \gamma (1-2\gamma L)\EE D_{f}(x^k,x^*) \leq 
        (1-\gamma\mu)\EE\norm{x^k - x^*}^2 +2\gamma^2\EE\sigma_k^2.
\end{equation}
\end{corollary}

\subsubsection*{Recovering SVRG rate}
Summing~\eqref{eq:svrg_iter_dec} for $k=0, \dots, m-1$ using $\sigma_k = \sigma_0$ we arrive at
\begin{eqnarray*}
\EE\norm{x^{m}-x^*}^2+ \sum_{k=1}^m \gamma (1-2\gamma L)\EE D_{f}(x^k,x^*)
&\leq& (1-\gamma\mu)\EE\norm{x^0 - x^*}^2 + 2m\gamma^2\EE\sigma_0^2 \\
&\leq &2 \left( \mu^{-1} +  2m\gamma^2   L\right) D_f(x^0,x^*) \;.
\end{eqnarray*}

Since $D_f$ is convex in the first argument, we have 
\[
m \gamma (1-2\gamma L)D_{f}\left( \frac1m \sum_{k=1}^m x^k,x^*\right) \leq \norm{x^{m}-x^*}^2+ \sum_{k=1}^m \gamma (1-2\gamma L)D_{f}(x^k,x^*)  
\]
and thus
\[
D_{f}\left( \frac1m \sum_{k=1}^m x^k,x^*\right) \leq \frac{ 2 \left( \mu^{-1} +  2m\gamma^2   L\right)}
{m \gamma (1-2\gamma L)} D_f(x^0,x^*) ,
 \]
which recovers rate from Theorem 1 in~\citep{SVRG}.

\subsection{{\tt L-SVRG}}\label{sec:L-SVRG}

In this section we show that our approach also covers {\tt L-SVRG} analysis from \citep{hofmann2015variance, kovalev2019don} (see Algorithm~\ref{alg:L-SVRG}) with a minor extension -- it allows for proximable regularizer $R$. Consider the finite-sum minimization problem
\begin{equation}\label{eq:main_l-svrg}
    f(x) = \frac{1}{n}\sum\limits_{i=1}^n f_i(x) + R(x),
\end{equation}
where each $f_i$ convex and  $L$-smooth for each $i$ and $f$ is $\mu$-strongly convex.

\begin{algorithm}[H]
    \caption{{\tt L-SVRG} (\cite{hofmann2015variance, kovalev2019don})}
    \label{alg:L-SVRG}
    \begin{algorithmic}
        \Require learning rate $\gamma>0$, probability $p\in (0,1]$, starting point $x^0\in\R^d$
        \State $w^0 = x^0$
        \For{ $k=0,1,2,\ldots$ }
        \State{Sample  $i \in \{1,\ldots, n\}$ uniformly at random}
        \State{$g^k = \nabla f_i(x^k) - \nabla f_i(w^k) + \nabla f(w^k)$}
        \State{$x^{k+1} = x^k - \gamma g^k$}
        \State{$w^{k+1} = \begin{cases}
            x^{k}& \text{with probability } p\\
            w^k& \text{with probability } 1-p
            \end{cases}$
        }
        \EndFor
    \end{algorithmic}
\end{algorithm}
Note that the gradient estimator is again unbiased, i.e. $\EE\left[g^k\mid x^k\right] = \nabla f(x^k)$. Next, Lemma~\ref{lem:l-svrg} provides with the remaining constants for Assumption~\ref{as:general_stoch_gradient}. The corresponding choice is stated in Table~\ref{tbl:special_cases-parameters}.

\begin{lemma}[Lemma~4.2 and Lemma~4.3 from~\citep{kovalev2019don} extended to prox setup]\label{lem:l-svrg}
    We have
    \begin{equation}\label{eq:lemma4.2_l-svrg}
        \EE\left[\norm{g^k - \nabla f(x^*)}^2\mid x^k\right] \le 4LD_f(x^k,x^*) + 2\sigma_k^2
    \end{equation}
    and 
        \begin{equation}\label{eq:lemma4.3_l-svrg}
        \EE\left[\sigma_{k+1}^2\mid x^k\right] \le (1-p)\sigma_k^2 + 2LpD_f(x^k, x^*),
    \end{equation}
    where $\sigma_k^2 \eqdef \frac{1}{n}\sum\limits_{i=1}^n\norm{\nabla f_i(w^k) - \nabla f_i(x^*)}^2$.
\end{lemma}

Next, applying Theorem~\ref{thm:main_gsgm} on Algorithm~\ref{alg:L-SVRG} with $M=\frac{4}{p}$ we get Corollary~\ref{cor:recover_l-svrg_rate}. 

\begin{corollary}\label{cor:recover_l-svrg_rate}
    {\tt L-SVRG} with $\gamma = \frac{1}{6L}$ satisfies
    \begin{equation}\label{eq:recover_l-svrg_rate}
        \EE V^k \le \left(1-\min\left\{\frac{\mu}{6L}, \frac{p}{2}\right\}\right)^kV^0.
    \end{equation}
\end{corollary}

\subsubsection*{Proof of Lemma~\ref{lem:l-svrg}}\label{sec:proofs_l-svrg}

    Let all expectations be conditioned on $x^k$. Using definition of $g^k$
    \begin{eqnarray*}
        \EE\left[\norm{g^k - \nabla f(x^*)}^2\right] &\overset{\text{Alg.}~\ref{alg:L-SVRG}}{=}&
        \EE\left[\norm{
            \nabla f_i(x^k)  - \nabla f_i(w^k) + \nabla f(w^k)
        - \nabla f(x^*)}^2\right]
        \\
        &\overset{\eqref{eq:a_b_norm_squared}}{\leq}&
        2\EE\left[\norm{\nabla f_i(x^k) - \nabla f_i(x^*)}^2\right]\\
        &&\quad+ 
        2\EE\left[\norm{\nabla f_i(x^*) - \nabla f_i(w^k) - \EE\left[\nabla f_i(x^*) - \nabla f_i(w^k)\mid x^k\right]}^2\right]\\
        &\overset{\eqref{eq:norm_diff_grads},\eqref{eq:variance_decomposition}}{\leq}&
        4LD_f(x^k,x^*) + 2 \EE\left[\norm{\nabla f_i(w^k) - \nabla f_i(x^*)}^2\right] \\
        &=&
        4LD_f(x^k,x^*) + 2 \sigma_k^2.
    \end{eqnarray*}
    For the second bound, we shall have
        \begin{eqnarray*}
        \EE\left[\sigma_{k+1}^2\right] &\overset{\text{Alg.}~\ref{alg:L-SVRG}}{=}& (1-p) \sigma_k^2 +  \frac{p}{n} \sum\limits_{i=1}^{n} \norm{\nabla f_i(x^k) - \nabla f_i(x^*)}^2\\
        &\overset{\eqref{eq:norm_diff_grads}}{\leq}& (1-p) \sigma_k^2 + 2Lp D_f(x^k,x^*).
    \end{eqnarray*}

\subsection{{\tt DIANA}}\label{sec:diana}
In this section we consider a distributed setup where each function $f_i$ from~\eqref{eq:f_sum} is owned by $i$-th machine (thus, we have all together $n$ machines). 

We show that our approach covers the analysis of {\tt DIANA} from \citep{mishchenko2019distributed, horvath2019stochastic}. {\tt DIANA} is a specific algorithm for distributed optimization with { \em quantization} -- lossy compression of gradient updates, which reduces the communication between the server and workers\footnote{It is a well-known problem in distributed optimization that the communication between machines often takes more time than actual computation.}. 

In particular, {\tt DIANA} quantizes gradient differences instead of the actual gradients. This trick allows for the linear convergence to the optimum once the full gradients are evaluated on each machine, unlike other popular quantization methods such as {\tt QSGD} \citep{alistarh2017qsgd} or {\tt TernGrad} \citep{wen2017terngrad}. In this case, {\tt DIANA} behaves as variance reduced method -- it reduces a variance that was injected due to the quantization. However, {\tt DIANA} also allows for evaluation of stochastic gradients on each machine, as we shall further see.

First of all, we introduce the notion of quantization operator.

\begin{definition}[Quantization]\label{def:quantization}
    We say that $\hat \Delta$ is a \textit{quantization} of vector $\Delta\in\R^d$ and write $\hat \Delta \sim {\rm Q}(\Delta)$ if
    \begin{equation}\label{eq:quantization}
        \EE\hat\Delta = \Delta, \qquad \EE\norm{\hat \Delta - \Delta}^2 \le \omega \norm{\Delta}^2
    \end{equation}
    for some $\omega > 0$.
\end{definition}

\begin{algorithm}[H]
   \caption{{\tt DIANA} \citep{mishchenko2019distributed, horvath2019stochastic}}
   \label{alg:diana}
\begin{algorithmic}[1]
   \Require learning rates $\alpha>0$ and $\gamma>0$, initial vectors $x^0, h_1^0,\dotsc, h_n^0 \in \R^d$ and $h^0 = \frac{1}{n}\sum_{i=1}^n h_i^0$
   \For{$k=0,1,\dotsc$}
       \State Broadcast $x^{k}$ to all workers
        \For{$i=1,\dotsc,n$ in parallel}
            \State Sample $g^{k}_i$ such that $\EE [g^k_i \;|\; x^k]  =\nabla f_i(x^k)$ 
            \State $\Delta^k_i = g^k_i - h^k_i$
            \State Sample $\hat \Delta^k_i \sim {\rm Q}(\Delta^k_i)$
            \State $h_i^{k+1} = h_i^k + \alpha \hat \Delta_i^k$
            \State $\hat g_i^k = h_i^k + \hat \Delta_i^k$
        \EndFor
       \State $\hat \Delta^k = \frac{1}{n}\sum_{i=1}^n \hat \Delta_i^k$
       \State $ g^k = \frac{1}{n}\sum_{i=1}^n \hat g_i^k = h^k + \hat \Delta^k $
        \State $x^{k+1} = \proxR\left(x^k - \gamma  g^k \right)$
        \State $h^{k+1}  = \tfrac{1}{n}\sum_{i=1}^n h_i^{k+1} = h^k + \alpha \hat \Delta^k$
   \EndFor
\end{algorithmic}
\end{algorithm}

The aforementioned method is applied to solve problem \eqref{eq:problem_gen}+\eqref{eq:f_sum} where each $f_i$ is convex and $L$-smooth and $f$ is $\mu$-strongly convex. 

\begin{lemma}[Lemma 1 and consequence of Lemma 2 from \citep{horvath2019stochastic}]\label{lem:lemma1_diana}
    Suppose that $\alpha \le \frac{1}{1+\omega}$. For all iterations $k\ge 0$ of Algorithm~\ref{alg:diana} it holds
    \begin{eqnarray}
        \EE\left[ g^k\mid x^k\right] &=& \nabla f(x^k), \label{eq:unbiased_diana}\\
        \EE\left[\norm{g^k - \nabla f(x^*)}^2\mid x^k\right] &\le& \left(1+\frac{2\omega}{n}\right)\frac{1}{n}\sum\limits_{i=1}^n\norm{\nabla f_i(x^k) - \nabla f_i(x^*)}^2\notag\\
        &&\quad + \frac{2\omega\sigma_k^2}{n} + \frac{(1+\omega)\sigma^2}{n},\label{eq:second_moment_diana}
        \\
    \EE\left[\sigma_{k+1}^2 \mid x^k\right] &\le& (1-\alpha)\sigma_k^2 + \frac{\alpha}{n}\sum\limits_{i=1}^n\norm{\nabla f_i(x^k) - \nabla f_i(x^*)}^2 + \alpha \sigma^2.\label{eq:h_i_k_sec_moment_diana}
    \end{eqnarray}
    where $\sigma_k^2  = \frac{1}{n}\sum\limits_{i=1}^n\norm{h_i^k -  \nabla f_i(x^*)}^2$ and $\sigma^2$ is such that $\frac{1}{n}\sum\limits_{i=1}^n\EE\left[\norm{g_i^k - \nabla f_i(x^k)}^2\mid x^k\right]\le \sigma^2$.
\end{lemma}
Bounding further $\frac1n \sum_{i=1}^n\norm{\nabla f_i(x^k) - \nabla f_i(x^*)}^2 \leq 2L D_{f}(x^k,x^*)$ in the above Lemma, we see that Assumption~\ref{as:general_stoch_gradient} as per Table~\ref{tbl:special_cases-parameters} is valid. Thus, as a special case of Theorem~\ref{thm:main_gsgm}, we obtain the following corollary.

\begin{corollary}\label{cor:main_diana}
    Assume that $f_i$ is convex and $L$-smooth for all $i\in[n]$ and $f$ is $\mu$ strongly convex, $\alpha \le \frac{1}{\omega+1}$, $\gamma \le \frac{1}{\left(1+\frac{2\omega}{n}\right)L + ML\alpha}$ where $M > \frac{2\omega}{n\alpha}$. Then the iterates of {\tt DIANA} satisfy
    \begin{equation}\label{eq:convergence_diana}
        \EE\left[V^k\right] \le \max\left\{(1-\gamma\mu)^k, \left(1 + \frac{2\omega}{nM} - \alpha\right)^k\right\}V^0 + \frac{\left(\frac{1+\omega}{n} + M\alpha\right)\sigma^2\gamma^2}{\min\left\{\gamma\mu, \alpha - \frac{2\omega}{nM}\right\}},
    \end{equation}
    where the Lyapunov function $V^k$ is defined by $V^k \eqdef \norm{x^k - x^*}^2 + M\gamma^2\sigma_k^2$. For the particular choice $\alpha = \frac{1}{\omega+1}$, $M = \frac{4\omega(\omega+1)}{n}$, $\gamma = \frac{1}{\left(1 + \frac{6\omega}{n}\right)L}$, then {\tt DIANA} converges to a solution neighborhood and the leading iteration complexity term is
    \begin{equation}\label{eq:diana_leading_term}
        \max\left\{\frac{1}{\gamma\mu}, \frac{1}{\alpha - \frac{2\omega}{nM}}\right\} = \max\left\{\kappa + \kappa \frac{6\omega}{n}, 2(\omega+1)\right\},
    \end{equation}
    where $\kappa = \frac{L}{\mu}$.
\end{corollary}

As mentioned, once the full (deterministic) gradients are evaluated on each machine, {\tt DIANA} converges linearly to the exact optimum. In particular, in such case we have $\sigma^2 = 0$. Corollary~\ref{cor:main_diana_special_case} states the result in the case when $n=1$, i.e. there is only a single node~\footnote{node = machine}. For completeness, we present the mentioned simple case of {\tt DIANA} as Algorithm~\ref{alg:diana_case}.

\begin{algorithm}[H]
   \caption{{\tt DIANA}: 1 node $\&$ exact gradients \citep{mishchenko2019distributed, horvath2019stochastic}}
   \label{alg:diana_case}
\begin{algorithmic}[1]
   \Require learning rates $\alpha>0$ and $\gamma>0$, initial vectors $x^0, h^0 \in \R^d$
   \For{$k=0,1,\dotsc$}
       \State $\Delta^k = \nabla f(x^k) - h^k$
       \State Sample $\hat \Delta^k \sim {\rm Q}(\Delta^k)$
       \State $h^{k+1} = h^k + \alpha \hat \Delta^k$
       \State $g^k = h^k + \hat \Delta^k$
       \State $x^{k+1} = \proxR\left(x^k - \gamma  g^k \right)$
   \EndFor
\end{algorithmic}
\end{algorithm}

\begin{corollary}\label{cor:main_diana_special_case}
    Assume that $f_i$ is $\mu$-strongly convex and $L$-smooth for all $i\in[n]$, $\alpha \le \frac{1}{\omega+1}$, $\gamma \le \frac{1}{\left(1+2\omega\right)L + ML\alpha}$ where $M > \frac{2\omega}{\alpha}$. Then the stochastic gradient $\hat g^k$ and the objective function $f$ satisfy Assumption~\ref{as:general_stoch_gradient} with $A = \left(1+2\omega\right)L, B = 2\omega, \sigma_k^2 = \norm{h^k - h^*}^2, \rho = \alpha, C = L\alpha, D_1 = 0, D_2 = 0$ and 
    \begin{equation}\label{eq:convergence_diana_special_case}
        \EE\left[V^k\right] \le \max\left\{(1-\gamma\mu)^k, \left(1 + \frac{2\omega}{M} - \alpha\right)^k\right\}V^0,
    \end{equation}
    where the Lyapunov function $V^k$ is defined by $V^k \eqdef \norm{x^k - x^*}^2 + M\gamma^2\sigma_k^2$. For the particular choice $\alpha = \frac{1}{\omega+1}$, $M = 4\omega(\omega+1)$, $\gamma = \frac{1}{\left(1 + 6\omega\right)L}$ the leading term in the iteration complexity bound is
    \begin{equation}\label{eq:diana_leading_term_special_case}
        \max\left\{\frac{1}{\gamma\mu}, \frac{1}{\alpha - \frac{2\omega}{M}}\right\} = \max\left\{\kappa + 6\kappa\omega, 2(\omega+1)\right\},
    \end{equation}
    where $\kappa = \frac{L}{\mu}$.
\end{corollary}

\subsection{{\tt Q-SGD-SR}}\label{Q-SGD-AS}

In this section, we consider a quantized version of {\tt SGD-SR}. 

\begin{algorithm}[H]
    \caption{{\tt Q-SGD-SR}}
    \label{alg:qsgdas}
    \begin{algorithmic}
        \Require learning rate $\gamma>0$, starting point $x^0\in\R^d$, distribution $\cD$ over $\xi \in\R^n$ such that $\ED{\xi}$ is vector of ones
        \For{ $k=0,1,2,\ldots$ }
        \State{Sample $\xi \sim \cD$}
        \State{$g^k \sim {\rm Q}(\nabla f_\xi (x^k))$}
        \State{$x^{k+1} = \proxR(x^k - \gamma g^k)$}
        \EndFor
    \end{algorithmic}
\end{algorithm}

\begin{lemma}[Generalization of Lemma~2.4, \citep{gower2019sgd}]\label{lem:exp_smoothness_grad_up_bound_q-sgd-as}
    If $(f,\cD)\sim ES(\cL)$, then
    \begin{equation}\label{eq:exp_smoothness_grad_up_bound_sgd-as}
        \ED{\norm{g^k - \nabla f(x^*)}^2} \le 4\cL(1+\omega)D_f(x^k,x^*) + 2\sigma^2(1+\omega).
    \end{equation}
    where $\sigma^2 \eqdef \ED{\norm{\nabla f_\xi(x^*)}^2}$.
\end{lemma}

A direct consequence of Theorem~\ref{thm:main_gsgm} in this setup is Corollary~\ref{cor:recover_q-sgd-as_rate}. 

\begin{corollary}\label{cor:recover_q-sgd-as_rate}
    Assume that $f(x)$ is $\mu$-strongly quasi-convex and $(f,\cD)\sim ES(\cL)$. Then {\tt Q-SGD-SR} with $\gamma^k \equiv\gamma \le \frac{1}{2(1+\omega)\cL}$ satisfies
    \begin{equation}\label{eq:recover_q-sgd-as_rate}
        \EE\left[\norm{x^k - x^*}^2\right] \le (1-\gamma\mu)^k\norm{x^0-x^*}^2 + \frac{2\gamma(1+\omega)\sigma^2}{\mu}.
    \end{equation}
\end{corollary}

\subsubsection*{Proof of Lemma~\ref{lem:exp_smoothness_grad_up_bound_q-sgd-as}}
In this proof all expectations are conditioned on $x^k$. First of all, from Lemma~\ref{lem:exp_smoothness_grad_up_bound_sgd-as} we have
\begin{eqnarray*}
    \ED{\norm{\nabla f_\xi(x^k) - \nabla f(x^*)}^2} \le 4\cL D_f(x^k,x^*) + 2\sigma^2.
\end{eqnarray*}
The remaining step is to understand how quantization of $\nabla f_\xi(x^k)$ changes the above inequality if we put $g^k\sim {\rm Q}(\nabla f_\xi(x^k))$ instead of $\nabla f_\xi(x^k)$. Let us denote mathematical expectation with respect randomness coming from quantization by $\EE_Q\left[\cdot\right]$. Using tower property of mathematical expectation we get
\begin{eqnarray*}
    \EE\left[\|g^k - \nabla f(x^*)\|^2\right] &=& \EE_{\cD}\left[\EE_Q\|g^k - \nabla f(x^*)\|^2\right]\\
    &\overset{\eqref{eq:variance_decomposition}}{=}& \EE\left[\|g^k - \nabla f_\xi(x^k)\|^2\right] + \EE\left[\|\nabla f_\xi(x^k) - \nabla f(x^*)\|^2\right]\\
    &\overset{\eqref{eq:exp_smoothness_grad_up_bound_sgd-as}}{\le}&  \EE\left[\|g^k - \nabla f_\xi(x^k)\|^2\right] + 4\cL D_f(x^k,x^*) + 2\sigma^2.
\end{eqnarray*} 
Next, we estimate the first term in the last row of the previous inequality
\begin{eqnarray*}
    \EE\left[\|g^k - \nabla f_\xi(x^k)\|^2\right] &\overset{\eqref{eq:quantization}}{\le}& \omega\EE\left[\|\nabla f_\xi(x^k)\|^2\right]\\
    &\overset{\eqref{eq:a_b_norm_squared}}{\le}& 2\omega\EE\left[\|\nabla f_\xi(x^k) - \nabla f_\xi(x^*)\|^2\right] + 2\omega\EE\left[\|\nabla f_\xi(x^*)\|^2\right]\\
    &\le& 4\omega\cL D_f(x^k,x^*) + 2\omega\sigma^2.
\end{eqnarray*}
Putting all together we get the result.

\subsection{{\tt VR-DIANA}} \label{sec:VR-DIANA}
Corollary~\ref{cor:main_diana} shows that once each machine evaluates a stochastic gradient instead of the full gradient, {\tt DIANA} converges linearly only to a certain neighborhood. In contrast, {\tt VR-DIANA}~\citep{horvath2019stochastic} uses a variance reduction trick within each machine, which enables linear convergence to the exact solution. In this section, we show that our approach recovers {\tt VR-DIANA} as well. 

\begin{algorithm}[H]
   \caption{{\tt VR-DIANA} based on L-SVRG (Variant 1), SAGA (Variant 2), \citep{horvath2019stochastic}}
   \label{alg:vr-diana_sigma_k}
\begin{algorithmic}[1]
        \Require{learning rates $\alpha > 0$ and $\gamma > 0$, initial vectors $x^0, h_{1}^0, \dots, h_{n}^0$, $h^0 = \frac{1}{n}\sum_{i=1}^n h_i^0$}
        \For{$k = 0,1,\ldots$}
        \State Sample random 
            $
                u^k = \begin{cases}
                    1,& \text{with probability } \frac{1}{m}\\
                    0,& \text{with probability } 1 - \frac{1}{m}\\
                \end{cases}
            $ \Comment{only for Variant 1}
        \State Broadcast $x^k$, $u^k$ to all workers\;
            \For{$i = 1, \ldots, n$ in parallel} \Comment{Worker side}
            \State Pick random $j_i^k \sim_{\rm u.a.r.} [m]$\;
            \State $\mu_i^k = \frac{1}{m} \sum\limits_{j=1}^{m} \nabla f_{ij}(w_{ij}^k)$\label{ln:mu} \;
            \State $g_i^k = \nabla f_{ij_i^k}(x^k) - \nabla f_{ij_i^k}(w_{ij_i^k}^k) + \mu_i^k$\;
            \State $\hat{\Delta}_i^k = Q(g_i^k - h_i^k)$\;
            \State $h_i^{k+1} = h_i^k + \alpha \hat{\Delta}_i^k$\;
                \For{$j = 1, \ldots, m$}
                    \State
                    $
                    w_{ij}^{k+1} =
                    \begin{cases}
                        x^k, & \text{if } u^k = 1 \\
                        w_{ij}^k, &\text{if } u^k = 0\\
                    \end{cases}
                    $ \Comment{Variant 1 (L-SVRG): update epoch gradient if $u^k = 1$}
                    \State
                    $
                    w_{ij}^{k+1} =
                    \begin{cases}
                    x^k, & j = j_i^k\\
                    w_{ij}^k, & j \neq j_i^k\\
                    \end{cases}
                    $ \Comment{Variant 2 (SAGA): update gradient table}
                \EndFor
            \EndFor
            \State $h^{k+1} \! = \! h^k \!+\! \frac{\alpha}{n} \displaystyle \sum_{i=1}^n \hat{\Delta}_i^k$ \Comment{Gather quantized updates} 
            \State $g^k = \frac{1}{n}\sum\limits_{i=1}^{n} (\hat{\Delta}_i^k + h_i^k)$\;
            \State $x^{k+1} = x^k - \gamma g^k$\;
        \EndFor

\end{algorithmic}  
\end{algorithm}

The aforementioned method is applied to solve problem \eqref{eq:problem_gen}+\eqref{eq:f_sum} where each $f_i$ is also of a finite sum structure, as in \eqref{eq:f_i_sum}, with  each $f_{ij}(x)$ being convex and $L$-smooth, and $f_i(x)$ being $\mu$-strongly convex. Note that $\nabla f(x^*) = 0$ and, in particular, $D_f(x,x^*) = f(x) - f(x^*)$ since the problem is considered without regularization.

\begin{lemma}[Lemmas 3, 5, 6 and 7 from \citep{horvath2019stochastic}]\label{lemmas_vr_diana}
    Let $\alpha \le \frac{1}{\omega+1}$. Then for all iterates $k\ge 0$ of Algorithm~\ref{alg:vr-diana_sigma_k} the following inequalities hold:
    \begin{eqnarray}
        \EE\left[g^k\mid x^k\right] &=& \nabla f(x^k),\label{eq:unbiased_g_k_vr_diana}\\
        \EE\left[H^{k+1}\mid x^k\right] &\le& \left(1-\alpha\right)H^k + \frac{2\alpha}{m}D^k + 8\alpha Ln\left(f(x^k) - f(x^*)\right),\label{eq:H_k+1_bound_vr_diana}\\
        \EE\left[D^{k+1}\mid x^k\right] &\le& \left(1 - \frac{1}{m}\right)D^k + 2Ln\left(f(x^k) - f(x^*)\right),\label{eq:D_k+1_bound_vr_diana}\\
        \EE\left[\norm{g^k}^2\mid x^k\right] &\le& 2L\left(1+\frac{4\omega + 2}{n}\right)\left(f(x^k)-f(x^*)\right) + \frac{2\omega}{n^2}\frac{D^k}{m} + \frac{2(\omega+1)}{n^2}H^k,\label{eq:second_moment_g_k_vr_diana}
    \end{eqnarray}
    where $H^k = \sum\limits_{i=1}^n\norm{h_i^k - \nabla f_i(x^*)}^2$ and $D^k = \sum\limits_{i=1}^n\sum\limits_{j=1}^m\norm{\nabla f_{ij}(w_{ij}^k) - \nabla f_{ij}(x^*)}^2$.
\end{lemma}

\begin{corollary}\label{cor:vr_diana_meets_assumption}
    Let $\alpha \le \min\left\{\frac{1}{3m},\frac{1}{\omega+1}\right\}$. Then stochastic gradient $\hat g^k$ (Algorithm~\ref{alg:vr-diana}) and the objective function $f$ satisfy Assumption~\ref{as:general_stoch_gradient} with $A = \left(1+\frac{4\omega + 2}{n}\right)L, B = \frac{2(\omega+1)}{n}, \rho = \alpha, C = L\left(\frac{1}{m}+4\alpha\right), D_1 = 0, D_2 = 0$ and
    \[
        \sigma_k^2 = \frac{H^k}{n} + \frac{D^k}{nm} = \frac{1}{n}\sum\limits_{i=1}^n\norm{h_i^k - \nabla f_i(x^*)}^2 + \frac{1}{nm}\sum\limits_{i=1}^n\sum\limits_{j=1}^m\norm{\nabla f_{ij}(w_{ij}^k) - \nabla f_{ij}(x^*)}^2.
    \]
\end{corollary}
\begin{proof}
    Indeed, \eqref{eq:general_stoch_grad_unbias} holds due to \eqref{eq:unbiased_g_k_vr_diana}. Inequality \eqref{eq:general_stoch_grad_second_moment} follows from \eqref{eq:second_moment_g_k_vr_diana} with $A = \left(1+\frac{4\omega + 2}{n}\right)L, B = \frac{2(\omega+1)}{n}, D_1 = 0, \sigma_k^2 = \frac{H^k}{n} + \frac{D^k}{nm}$ if we take into account that $\frac{2\omega}{n^2}\frac{D^k}{m} + \frac{2(\omega+1)}{n^2}H^k \le \frac{2(\omega+1)}{n}\left(\frac{D^k}{nm} + \frac{H^k}{n}\right)$. Finally, summing inequalities \eqref{eq:H_k+1_bound_vr_diana} and \eqref{eq:D_k+1_bound_vr_diana} and using $\alpha\le\frac{1}{3m}$
    \begin{eqnarray*}
        \EE\left[\sigma_k^2\mid x^k\right] &=& \frac{1}{n}\EE\left[H^{k+1}\mid x^k\right] + \frac{1}{nm}\EE\left[D^{k+1}\mid x^k\right]\\
        &\overset{\eqref{eq:H_k+1_bound_vr_diana}+\eqref{eq:D_k+1_bound_vr_diana}}{\le}& \left(1-\alpha\right)\frac{H^k}{n} + \left(1+2\alpha-\frac{1}{m}\right)\frac{D^k}{nm} + 2L\left(\frac{1}{m}+4\alpha\right)\left(f(x^k)-f(x^*)\right)\\
        &\le& \left(1-\alpha\right)\sigma_k^2 + 2L\left(\frac{1}{m}+4\alpha\right)\left(f(x^k)-f(x^*)\right)
    \end{eqnarray*}
    we get \eqref{eq:gsg_sigma} with $\rho = \alpha, C = L\left(\frac{1}{m}+4\alpha\right), D_2 = 0$.
\end{proof}

\begin{corollary}\label{cor:main_vr_diana}
    Assume that $f_i$ is $\mu$-strongly convex and $f_{ij}$ is convex and $L$-smooth for all $i\in[n], j\in[m]$, $\alpha \le \min\left\{\frac{1}{3m},\frac{1}{\omega+1}\right\}$, $\gamma \le \frac{1}{\left(1+\frac{4\omega + 2}{n}\right)L + ML\left(\frac{1}{m}+4\alpha\right)}$ where $M > \frac{2(\omega+1)}{n\alpha}$. Then the iterates of {\tt VR-DIANA} satisfy
    \begin{equation}\label{eq:convergence_vr_diana}
        \EE\left[V^k\right] \le \max\left\{(1-\gamma\mu)^k, \left(1 + \frac{2(\omega+1)}{nM} - \alpha\right)^k\right\}V^0,
    \end{equation}
    where the Lyapunov function $V^k$ is defined by $V^k \eqdef \norm{x^k - x^*}^2 + M\gamma^2\sigma_k^2$. Further, if we set  $\alpha = \min\left\{\frac{1}{3m},\frac{1}{\omega+1}\right\}$, $M = \frac{4(\omega+1)}{n\alpha}$, $\gamma = \frac{1}{\left(1 + \frac{20\omega+18}{n} + \frac{4\omega+4}{n\alpha m}\right)L}$, then to achieve precision $\EE\left[\norm{x^k-x^*}^2\right] \le \varepsilon V^0$ {\tt VR-DIANA} needs $\cO\left(\max\left\{\kappa+\kappa\frac{\omega+1}{n}+\kappa\frac{(\omega+1)\max\left\{m,\omega+1\right\}}{nm},m,\omega+1\right\}\log\frac{1}{\varepsilon}\right)$ iterations, where $\kappa = \frac{L}{\mu}$.
\end{corollary}
\begin{proof}
    Using Corollary~\ref{cor:vr_diana_meets_assumption} we apply Theorem~\ref{thm:main_gsgm} and get the result.
\end{proof}

\begin{remark}
{\tt VR-DIANA} can be easily extended to the proximal setup in our framework.
\end{remark}

\subsection{{\tt JacSketch}} \label{sec:JacSketch}

In this section, we show that our approach covers the analysis of {\tt JacSketch} from \citep{gower2018stochastic}. {\tt JacSketch} is a generalization of {\tt SAGA} in the following manner. {\tt SAGA} observes every iteration $\nabla f_i(x)$ for random index $i$ and uses it to build both stochastic gradient as well as the control variates on the stochastic gradient in order to progressively decrease variance. In contrast, {\tt JacSketch}  observes every iteration the random sketch of the Jacobian, which is again used to build both stochastic gradient as well as the control variates on the stochastic gradient.

For simplicity, we do not consider proximal setup, since~\citep{gower2018stochastic} does not either.

We first introduce the necessary notation (same as in \citep{gower2018stochastic}). Denote first the Jacobian the objective \begin{equation}\label{eq:jac_def}\Jac(x) \eqdef [\nabla f_1(x), \ldots, \nabla f_n(x)] \in \R^{d\times n}.\end{equation} 
Every iteration of the method, a random sketch of Jacobian $\nabla F(x^k)\mS$ (where $\mS\sim \cD$) is observed. Then, the method builds a variable $\mJ^k$, which is the current Jacobian estimate, updated using so-called sketch and project iteration~\citep{gower2015randomized}:
\[
\mJ^{k+1}  = \mJ^k(\mI - \Proj_{\mS_k}) + \Jac(x^k)\Proj_{\mS_k},
\]

where $\Proj_\mS$ is a projection under $\mW$ norm\footnote{Weighted Frobenius norm of matrix $\mX\in\R^{n\times n}$ with a positive definite weight matrix $\mW\in \R^{n\times n}$ is defined as 
$\norm{\mX}_{\mW^{-1}} \eqdef \sqrt{\Tr{ \mX \mW^{-1} \mX^\top}}.$
} ($\mW\in \R^{n\times n}$ is some positive definite weight matrix) defined as
$\Proj_\mS \eqdef  \mS (\mS^\top \mW \mS)^{\dagger} \mS^\top \mW$\footnote{Symbol $\dagger$ stands for Moore-Penrose pseudoinverse.}.

Further, in order to construct unbiased stochastic gradient, an access to the random scalar $\theta_{\mS}$ such that
\begin{equation}\label{eq:unbiased}
\ED{\theta_{\mS} \Proj_\mS} \ones  =  \ones,
\end{equation}
where $e$ is the vector of all ones.

Next, the simplest option for the choice of the stochastic gradient is $\nabla f_{\mS}(x)$ -- an unbiased estimate of $\nabla f$ directly constructed using $\mS,\theta_{\mS} $:
\begin{equation}
\label{eq:stochgradplain}
\nabla f_{\mS}(x)   = \frac{\theta_{\mS}}{n}\Jac(x) \Proj_\mS \ones.
\end{equation}

However, one can build a smarter estimate $ \nabla f_{\mS,\mJ}(x) $ via control variates constructed from $\mJ$:
\begin{equation}\label{eq:controlgradJ} 
 \nabla f_{\mS,\mJ}(x) = \frac{\theta_{\mS}}{n} (\Jac(x)-\mJ)   \Proj_\mS\ones  + \frac{1}{n} \mJ \ones.
\end{equation}
The resulting algorithm is stated as Algorithm~\ref{alg:jacsketch}.

\begin{algorithm}
    \begin{algorithmic}[1]
        \Require $\left(\cD, \mW, \theta_{\mS} \right)$, $x^0\in \R^d$, Jacobian estimate $\mJ^0 \in \R^{d \times n}$, stepsize $\gamma>0$         
        
        \For {$k =  0, 1, 2, \dots$}
        \State Sample a fresh copy $\mS_k\sim \cD$

        \State $ \mJ^{k+1}  = \mJ^k(\mI - \Proj_{\mS_k}) + \Jac(x^k)\Proj_{\mS_k}$
         \label{ln:jacupdate}        
        \State $g^{k} =  \nabla f_{\mS_k, \mJ^k}(x^k)$       \label{ln:gradupdate}    
                    
        \State $x^{k+1} = x^k - \gamma g^{k}$ \label{ln:xupdate}    
         
        \EndFor
    \end{algorithmic}
    \caption{{\tt JacSketch} \citep{gower2018stochastic}}
    \label{alg:jacsketch}
\end{algorithm}

Next we present Lemma~\ref{lem:lemmas39_310_jacsketch} which directly justifies the parameter choice from Table~\ref{tbl:special_cases2}. 

\begin{lemma}[Lemmas 2.5, 3.9 and 3.10 from \citep{gower2018stochastic}]\label{lem:lemmas39_310_jacsketch}
    Suppose that there are constants $\cL_1, \cL_2>0$ such that 
    \begin{eqnarray*}
      \ED{ \norm{ \nabla f_{\mS}(x) - \nabla f_{\mS}(x^*)}_2^2 } &\leq& 2  \cL_1 (f(x)-f(x^*)), \qquad \forall x\in \R^d \\
      \ED{\norm{(\Jac(x)-\Jac(x^*)) \Proj_{\mS}  }_{\mW^{-1}}^2 }& \leq & 2\cL_2 (f(x) -f(x^*)), \qquad \forall x\in \R^d, \label{eq:ES2}
\end{eqnarray*}    
    
     Then
    \begin{equation}\label{eq:jacs_contraction}
     \ED{\norm{\mJ^{k+1} -\Jac(x^*)}_{\mW^{-1}}^2} \le (1-\lambda_{\min}) \norm{\mJ^{k}-\Jac(x^*)}_{\mW^{-1}}^2 +  2\cL_2(f(x^k) -f(x^*)),
     \end{equation}
     \begin{equation}\label{eq:gradbndsubdeltaXX}
    \ED{\norm{g^k}_2^2 } \le 4 \cL_1  (f(x^k) -f(x^*)) +  2 \frac{\lambda_{\max}}{n^2} \norm{\mJ^{k} -\Jac(x^*)}_{\mW^{-1}}^2,
    \end{equation}
    where $\lambda_{\min} = \lambda_{\min}\left(\ED{\Proj_{\mS}}\right)$ and $\lambda_{\max} = \lambda_{\max}\left( \mW^{1/2}\left( \ED{\theta_{\mS}^2 \Proj_{\mS} \ones \ones^\top \Proj_{\mS}} -\ones \ones^\top\right) \mW^{1/2}\right)$. Further,     $
    \ED{ \nabla f_{\mS,\mJ}(x)} = \nabla f(x)$.
\end{lemma}

Thus, as a direct consequence of Theorem~\ref{thm:main_gsgm}, we obtain the next corollary.

\begin{corollary}\label{thm:main_jacsketch}
Consider the setup from Lemma~\ref{lem:lemmas39_310_jacsketch}. Suppose that $f$ is $\mu$-strongly convex and choose $\gamma \le \min\left\{\frac{1}{\mu},\frac{1}{2\cL_1 + M\frac{\cL_2}{n}}\right\}$ where $M > \frac{2\lambda_{\max}}{n\lambda_{\min}}$. Then the iterates of {\tt JacSketch} satisfy
    \begin{equation}\label{eq:convergence_jacsketch}
        \EE\left[V^k\right] \le \max\left\{(1-\gamma\mu)^k, \left(1 + \frac{2\lambda_{\max}}{nM} - \lambda_{\min}\right)^k\right\}V^0.
    \end{equation}

\end{corollary}

\begin{remark}\label{rem:gjs}
We shall note that concurrently with this work, a more general version of JacSketch with refined analysis was proposed in~\citep{hanzely2019one}, obtaining many new methods in special case (such as {\tt LSVRG}, {\tt SEGA} and several new ones), with best known rate in each special case. As mentioned in the main body of the paper, the rates from~\citep{hanzely2019one} for methods that have randomness in partial derivatives and non-uniform smoothness are better to what can Theorem~\ref{thm:main_gsgm} achieve. On the other hand, \citep{hanzely2019one} only focuses on variance reduced methods, while this paper analyzes also methods with extra noise. 
\end{remark}

\subsection{Interpolation Between Methods ~\label{sec:interpol}}

Given that a set of stochastic gradients satisfy Assumption~\ref{as:general_stoch_gradient}, we show that an any convex combination of the mentioned stochastic gradients satisfy Assumption~\ref{as:general_stoch_gradient} as well.

\begin{lemma}\label{lem:convex_comb}
    Assume that sequences of stochastic gradients $\{g_1^k\}_{k\ge 0}, \ldots, \{g_m^k\}_{k\ge 0}$ at the common iterates $\{x^k\}_{k\ge 0}$ satisfy the Assumption~\ref{as:general_stoch_gradient} with parameters $A(j),B(j),\{\sigma_k^2(j)\}_{k\ge 0},$ $ C(j),\rho(j),D_1(j),D_2(j)$, $j\in[m]$ respectively. Then for any vector $\tau = (\tau_1,\ldots,\tau_m)^\top$ such as $\sum\limits_{j=1}^m\tau_j = 1$ and $\tau_j \ge 0, j\in[m]$ stochastic gradient $g_\tau^k = \sum\limits_{j=1}^m\tau_j g_j^k$ satisfies the Assumption~\ref{as:general_stoch_gradient} with parameters:
    \begin{eqnarray}
        A_\tau = \sum\limits_{j=1}^m\tau_j A(j),\quad B_\tau = 1,\quad \sigma_{\tau,k}^2 = \sum\limits_{j=1}^m B(j)\tau_j \sigma_k^2(j),\quad \rho_\tau = \min\limits_{j\in[m]}\rho(j),\notag\\
        C_\tau = \sum\limits_{j=1}^m\tau_j C(j) B(j),\quad D_{\tau,1} = \sum\limits_{j=1}^m\tau_j D_1(j),\quad D_{\tau,2} = \sum\limits_{j=1}^m\tau_j D_2(j) B(j).\label{eq:conv_comb_params}
    \end{eqnarray}
    Furthermore, if stochastic gradients $g_1^k, \dots, g_m^k$ are independent for all $k$, Assumption~\ref{as:general_stoch_gradient} is satisfied with parameters
        \begin{eqnarray}
        A_\tau = L+\sum\limits_{j=1}^m\tau_j^2 A(j),\quad B_\tau = 1,\quad \sigma_{\tau,k}^2 = \sum\limits_{j=1}^m B(j)\tau_j^2 \sigma_k^2(j),\quad \rho_\tau = \min\limits_{j\in[m]}\rho(j),\notag\\
        C_\tau = \sum\limits_{j=1}^m\tau_j^2 C(j) B(j),\quad D_{\tau,1} = \sum\limits_{j=1}^m\tau_j^2 D_1(j),\quad D_{\tau,2} = \sum\limits_{j=1}^m\tau_j^2 D_2(j) B(j).\label{eq:conv_comb_params_indep}
    \end{eqnarray}
\end{lemma}

What is more, instead of taking convex combination one can choose stochastic gradient at random. Lemma~\ref{lem:flipping_a_coin} provides the result. 
\begin{lemma}\label{lem:flipping_a_coin}
    Assume that sequences of stochastic gradients $\{g_1^k\}_{k\ge 0}, \ldots, \{g_m^k\}_{k\ge 0}$ at the common iterates $\{x^k\}_{k\ge 0}$ satisfy the Assumption~\ref{as:general_stoch_gradient} with parameters $A(j),B(j),\{\sigma_k^2(j)\}_{k\ge 0},$ $ C(j),\rho(j),D_1(j),D_2(j)$, $j\in[m]$ respectively. Then for any vector $\tau = (\tau_1,\ldots,\tau_m)^\top$ such as $\sum\limits_{j=1}^m\tau_j = 1$ and $\tau_j \ge 0, j\in[m]$ stochastic gradient $g_\tau^k$ which equals $g_j^k$ with probability $\tau_j$ satisfies the Assumption~\ref{as:general_stoch_gradient} with parameters:
    \begin{eqnarray}
        A_\tau = \sum\limits_{j=1}^m\tau_j A(j),\quad B_\tau =1,\quad \sigma_{\tau,k}^2 = \sum\limits_{j=1}^m\tau_j B(j) \sigma_k^2(j),\quad \rho_\tau = \min\limits_{j\in[m]}\rho(j),\notag\\
        C_\tau = \sum\limits_{j=1}^m\tau_j  B(j)C(j),\quad D_{\tau,1} = \sum\limits_{j=1}^m\tau_j D_1(j),\quad D_{\tau,2} = \sum\limits_{j=1}^mB(j)\tau_j D_2(j).\label{eq:flipping_a_coin}
    \end{eqnarray}
    Furthermore, if stochastic gradients $g_1^k, \dots, g_m^k$ are independent for all $k$, Assumption~\ref{as:general_stoch_gradient} is satisfied with parameters
        \begin{eqnarray}
        A_\tau = L+\sum\limits_{j=1}^m\tau_j^2 A(j),\quad B_\tau = 1,\quad \sigma_{\tau,k}^2 = \sum\limits_{j=1}^m B(j)\tau_j^2 \sigma_k^2(j),\quad \rho_\tau = \min\limits_{j\in[m]}\rho(j),\notag\\
        C_\tau = \sum\limits_{j=1}^m\tau_j^2 C(j) B(j),\quad D_{\tau,1} = \sum\limits_{j=1}^m\tau_j^2 D_1(j),\quad D_{\tau,2} = \sum\limits_{j=1}^m\tau_j^2 D_2(j) B(j).\label{eq:flipping_a_coin_indep}
    \end{eqnarray}
\end{lemma}

\begin{example}[{\tt $\tau$-L-SVRG}]
    Consider the following method~--- {\tt $\tau$-L-SVRG}~--- which interpolates between vanilla {\tt SGD} and {\tt L-SVRG}.
    \begin{algorithm}[H]
    \caption{{\tt $\tau$-L-SVRG}}
    \label{alg:tau-L-SVRG}
    \begin{algorithmic}
        \Require learning rate $\gamma>0$, probability $p\in (0,1]$, starting point $x^0\in\R^d$, convex combination parameter $\tau\in[0,1]$
        \State $w^0 = x^0$
        \For{ $k=0,1,2,\ldots$ }
        \State{Sample  $i \in \{1,\ldots, n\}$ uniformly at random}
        \State{$g^k_{L-SVRG} = \nabla f_i(x^k) - \nabla f_i(w^k) + \nabla f(w^k)$}
        \State{Sample  $j \in \{1,\ldots, n\}$ uniformly at random}
        \State{$g^k_{SGD} = \nabla f_j(x^k)$}
        \State{$g^k = \tau g^k_{SGD} + (1-\tau)g^k_{L-SVRG}$}        
        \State{$x^{k+1} = x^k - \gamma g^k$}
        \State{$w^{k+1} = \begin{cases}
            x^{k}& \text{with probability } p\\
            w^k& \text{with probability } 1-p
            \end{cases}$
        }
        \EndFor
    \end{algorithmic}
    \end{algorithm}
    When $\tau = 0$ the Algorithm~\ref{alg:tau-L-SVRG} becomes {\tt L-SVRG} and when $\tau = 1$ it is just {\tt SGD} with uniform sampling.    Notice that Lemmas~\ref{lem:l-svrg} and~\ref{lem:exp_smoothness_grad_up_bound_sgd-as} still hold as they does not depend on the update rule for $x^{k+1}$.
    
    Thus, sequences $\{g_{SGD}^k\}_{k\ge 0}$ and $\{g_{L-SVRG}^k\}_{k\ge 0}$ satisfy the conditions of Lemma~\ref{lem:convex_comb} and, as a consequence, stochastic gradient $g^k$ from {\tt $\tau$-L-SVRG} meets the Assumption~\ref{as:general_stoch_gradient} with the following parameters:
    \begin{eqnarray*}
        A_\tau =  L + 2\tau^2\cL + 2(1-\tau)^2L,\quad B_\tau = 1,\quad \sigma_{\tau,k}^2 = 2\frac{(1-\tau)^2}{n}\sum\limits_{i=1}^n\norm{\nabla f_i(w^k) - \nabla f_i(x^*)}^2,\notag\\
         \rho_\tau = p,\quad C_\tau = 2(1-\tau)^2Lp,\quad D_{\tau,1} = 2\tau^2\sigma^2,\quad D_{\tau,2} = 0.
    \end{eqnarray*}
\end{example}

\begin{remark}
Similar interpolation with the analogous analysis can be considered between {\tt SGD} and {\tt SAGA}, or {\tt SGD} and {\tt SVRG}. 
\end{remark}

\subsubsection*{Proof of Lemma~\ref{lem:convex_comb}}

    Indeed, \eqref{eq:general_stoch_grad_unbias} holds due to linearity of mathematical expectation. Next, summing inequalities \eqref{eq:general_stoch_grad_second_moment} for $g_1^k,\ldots,g_m^k$ and using convexity of $\norm{\cdot}^2$ we get
    \begin{eqnarray*}
        \EE\left[\norm{g_\tau^k -\nabla f(x^*)}^2\mid x^k\right] &\le& \sum\limits_{j=1}^m\tau_j\EE\left[\norm{g_j^k-\nabla f(x^*)}^2\mid x^k\right]
        \\
        & \overset{\eqref{eq:general_stoch_grad_second_moment}}{\le}& 2\sum\limits_{j=1}^m\tau_j A(j) D_f(x^k,x^*) + \sum\limits_{j=1}^mB(j)\tau_j\sigma_k^2(j) +     \sum\limits_{j=1}^m\tau_j D_1(j),
    \end{eqnarray*}    
    which implies \eqref{eq:general_stoch_grad_second_moment} for $g_\tau^k$ with $A_\tau = \sum\limits_{j=1}^m\tau_j A(j), B_\tau =1, \sigma_{\tau,k}^2 = \sum\limits_{j=1}^m\tau_j B(j)\sigma_k^2(j), D_{\tau,1} = \sum\limits_{j=1}^m\tau_j D_1(j)$.
    Finally, summing \eqref{eq:gsg_sigma} for $g_1^k,\ldots,g_m^k$ gives us
    \begin{equation*}
        \EE\left[\sigma_{\tau,k+1}^2\mid\sigma_{\tau,k}^2\right] \overset{\eqref{eq:gsg_sigma}}{\le} \left(1-\min\limits_{j\in[m]}\rho(j)\right)\sigma_{\tau,k}^2 + 2\sum\limits_{j=1}^m\tau_j B(j)C(j)D_f(x^k,x^*) + \sum\limits_{j=1}^m\tau_j B(j)D_2(j),
    \end{equation*}
    which is exactly \eqref{eq:gsg_sigma} for $\sigma_{\tau,k}^2$ with $\rho =\min\limits_{j\in[m]}\rho(j), C_\tau = \sum\limits_{j=1}^m\tau_j C(j), D_{\tau,2} = \sum\limits_{j=1}^m\tau_j D_2(j)$.

Next, for independent gradients we have
\begin{eqnarray}
\EE\left[\norm{g_\tau^k -\nabla f(x^*)}^2\mid x^k\right]   &=&
 \sum_{j=1}^m \tau_j^2 \EE \left[\norm{g_j^k -\nabla f(x^*)}^2\mid x^k\right]\notag\\
 &&\quad+ 2\sum_{i< j} \tau_i \tau_j\EE\<g_j^k -\nabla f(x^*),g_i^k -\nabla f(x^*) >
 \nonumber
 \\
  &=&
 \sum_{j=1}^m \tau_j^2 \EE \left[\norm{g_j^k -\nabla f(x^*)}^2\mid x^k\right] + 2\sum_{i< j}\tau_i \tau_j \norm{\nabla f(x^k) -\nabla f(x^*)}^2
  \nonumber
 \\
  &\leq&
 \sum_{j=1}^m \tau_j^2 \EE \left[\norm{g_j^k -\nabla f(x^*)}^2\mid x^k\right] + \left(\sum_{ j=1}^m \tau_j\right)^2 \norm{\nabla f(x^k) -\nabla f(x^*)}^2
  \nonumber
  \\
  &=&
 \sum_{j=1}^m \tau_j^2 \EE \left[\norm{g_j^k -\nabla f(x^*)}^2\mid x^k\right] +  \norm{\nabla f(x^k) -\nabla f(x^*)}^2
  \nonumber
   \\
  &\leq&
 \sum_{j=1}^m \tau_j^2 \EE \left[\norm{g_j^k -\nabla f(x^*)}^2\mid x^k\right] +  2LD_f(x^k,x^*).
 \label{eq:indp_bounding}
\end{eqnarray}
and further the bounds follow. 

\subsubsection*{Proof of Lemma~\ref{lem:flipping_a_coin}}

    Indeed, \eqref{eq:general_stoch_grad_unbias} holds due to linearity and tower property of mathematical expectation. Next, using tower property of mathematical expectation and inequalities \eqref{eq:general_stoch_grad_second_moment} for $g_1^k,\ldots,g_m^k$ we get
    \begin{eqnarray*}
        \EE\left[\norm{g_\tau^k-\nabla f(x^*)}^2\mid x^k\right] &=& \EE\left[\EE_\tau\left[\norm{g_\tau^k-\nabla f(x^*)}^2\right]\mid x^k\right] = \sum\limits_{j=1}^m\tau_j\EE\left[\norm{g_j^k-\nabla f(x^*)}^2\mid x^k\right]\\ &\overset{\eqref{eq:general_stoch_grad_second_moment}}{\le}& 2\sum\limits_{j=1}^m\tau_j A(j) D_f(x^k,x^*) + \sum\limits_{j=1}^m B(j)\tau_j\sigma_k^2(j) +     \sum\limits_{j=1}^m\tau_j D_1(j),
    \end{eqnarray*}
    which implies \eqref{eq:general_stoch_grad_second_moment} for $g_\tau^k$ with $A_\tau = \sum\limits_{j=1}^m\tau_j A(j), B_\tau = 1, \sigma_{\tau,k}^2 = \sum\limits_{j=1}^m\tau_j B(j) \sigma_k^2(j), D_{\tau,1} = \sum\limits_{j=1}^m\tau_j D_1(j)$.
    Finally, summing \eqref{eq:gsg_sigma} for $g_1^k,\ldots,g_m^k$ gives us
    \begin{equation*}
        \EE\left[\sigma_{\tau,k+1}^2\mid\sigma_{\tau,k}^2\right] \overset{\eqref{eq:gsg_sigma}}{\le} \left(1-\min\limits_{j\in[m]}\rho(j)\right)\sigma_{\tau,k}^2 + 2\sum\limits_{j=1}^m\tau_j B(j) C(j)D_f(x^k,x^*) + \sum\limits_{j=1}^m\tau_jB(j)D_2(j),
    \end{equation*}
    which is exactly \eqref{eq:gsg_sigma} for $\sigma_{\tau,k}^2$ with $\rho =\min\limits_{j\in[m]}\rho(j), C_\tau = \sum\limits_{j=1}^m\tau_j B(j)C(j), D_{\tau,2} = \sum\limits_{j=1}^m\tau_j B(j)D_2(j)$.
To show~\eqref{eq:flipping_a_coin_indep}, it suffices to combine above bounds with the trick~\eqref{eq:indp_bounding}.

\begin{remark}
Recently, \citep{tran2019hybrid} demonstrated in that the convex combination of {\tt SGD} and {\tt SARAH}~\citep{nguyen2017sarah} performs very well on non-convex problems. 
\end{remark}

\section{Experiments \label{sec:exp}}

\subsection{Experiments on {\tt SGD-MB}}

\begin{figure}[H]
\centering
\begin{minipage}{0.24\textwidth}
  \centering
\includegraphics[width =  \textwidth ]{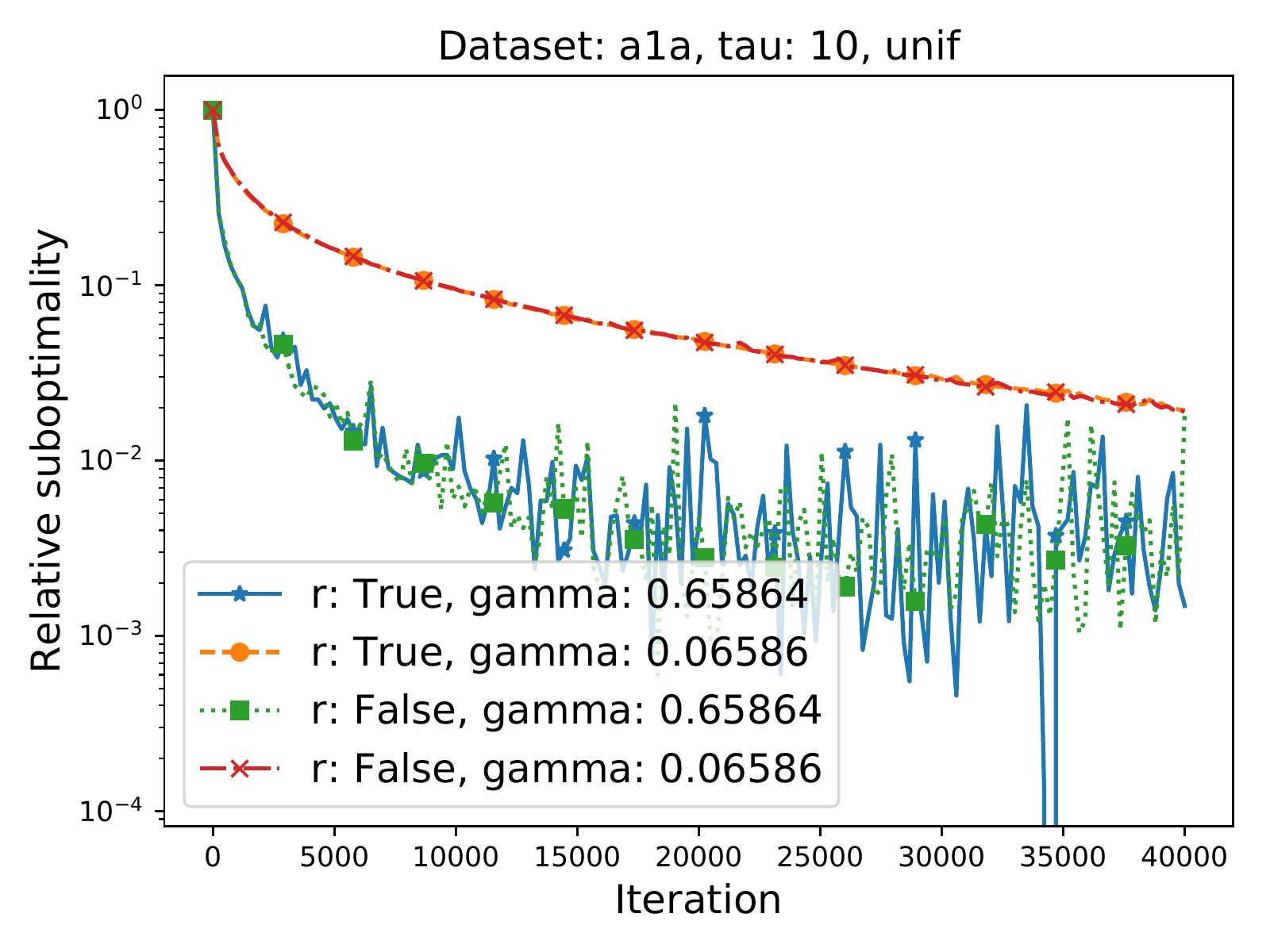}
\end{minipage}%
\begin{minipage}{0.24\textwidth}
  \centering
\includegraphics[width =  \textwidth ]{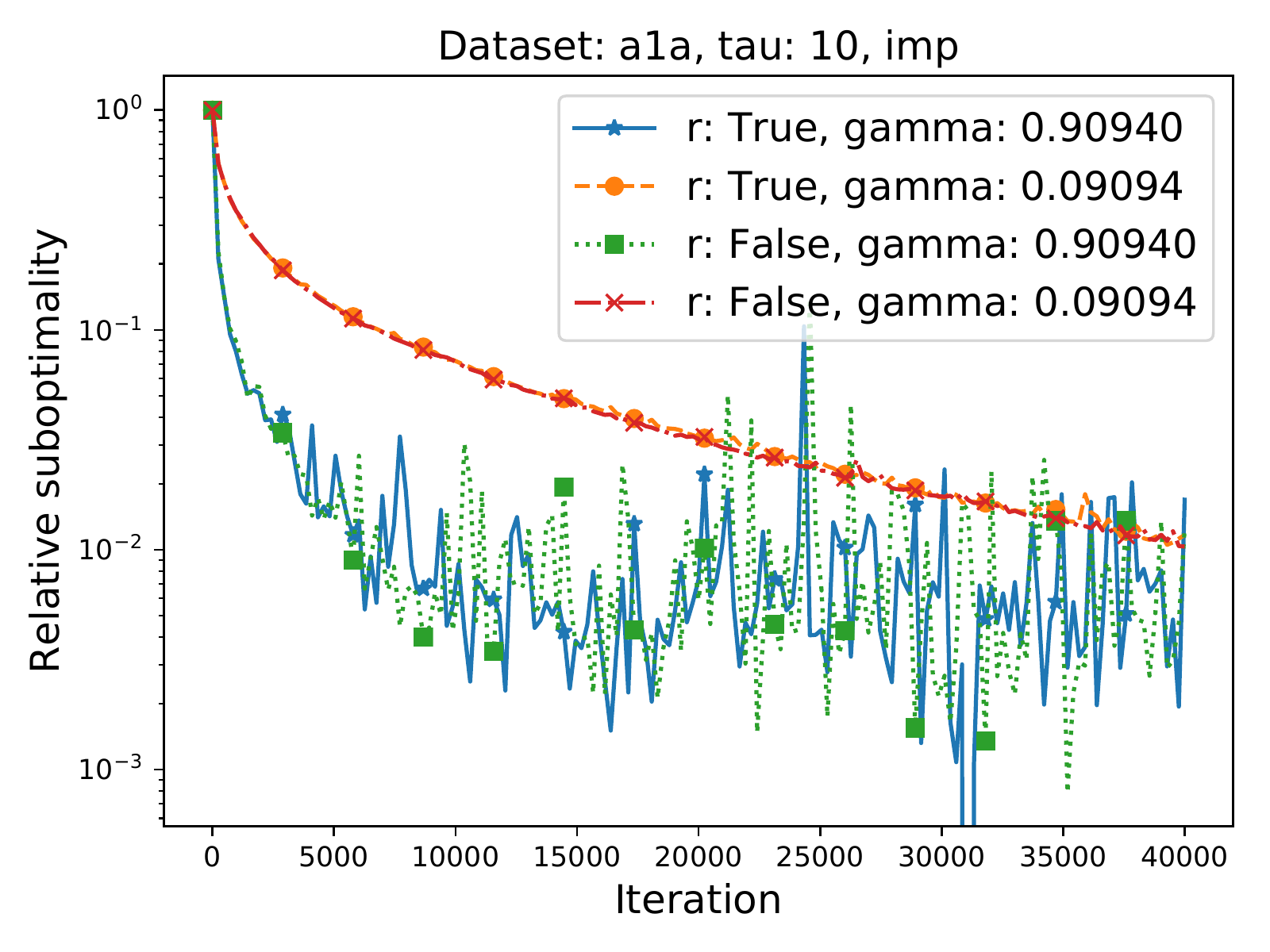}
\end{minipage}%
\begin{minipage}{0.24\textwidth}
  \centering
\includegraphics[width =  \textwidth ]{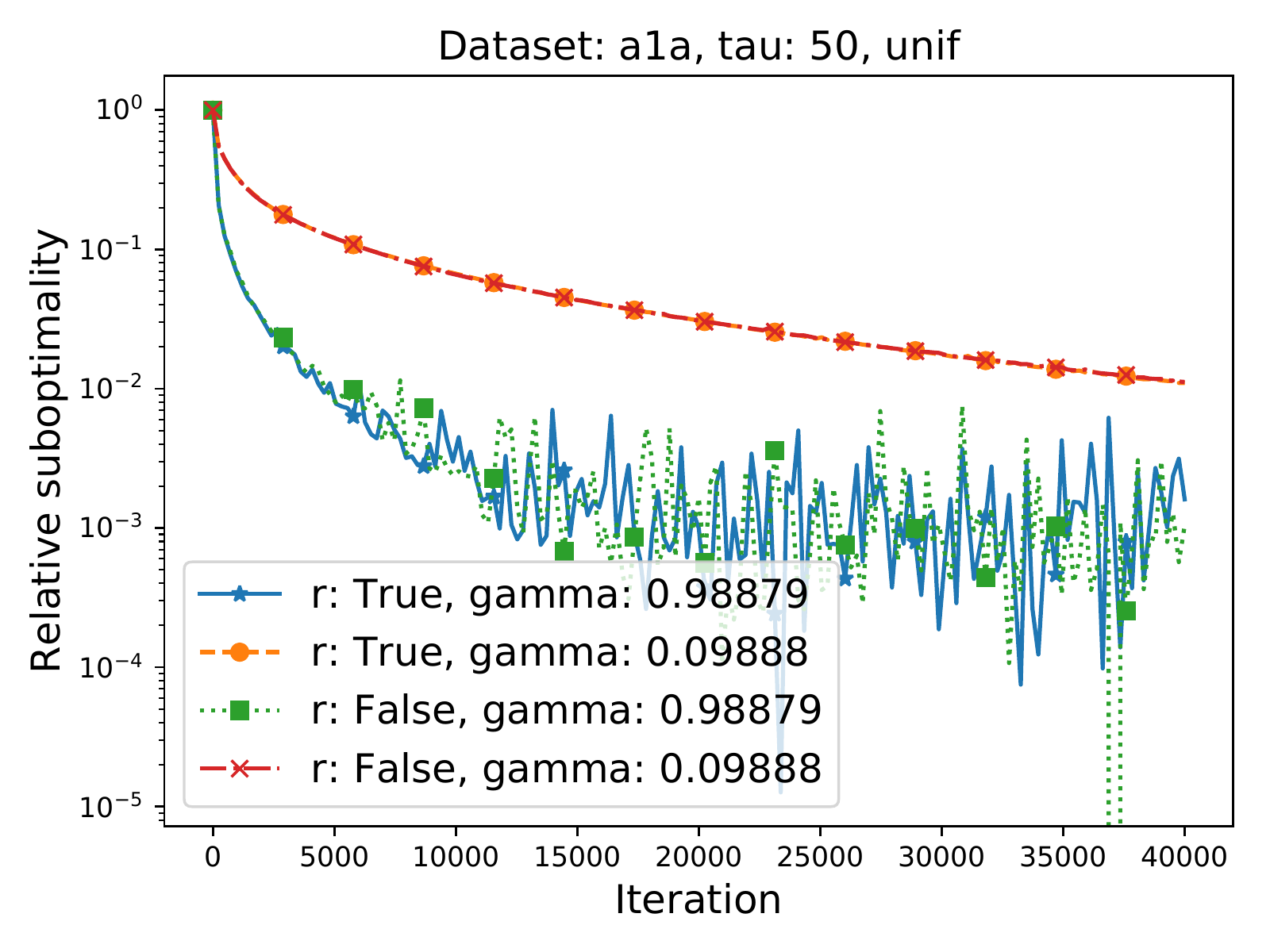}
\end{minipage}%
\begin{minipage}{0.24\textwidth}
  \centering
\includegraphics[width =  \textwidth ]{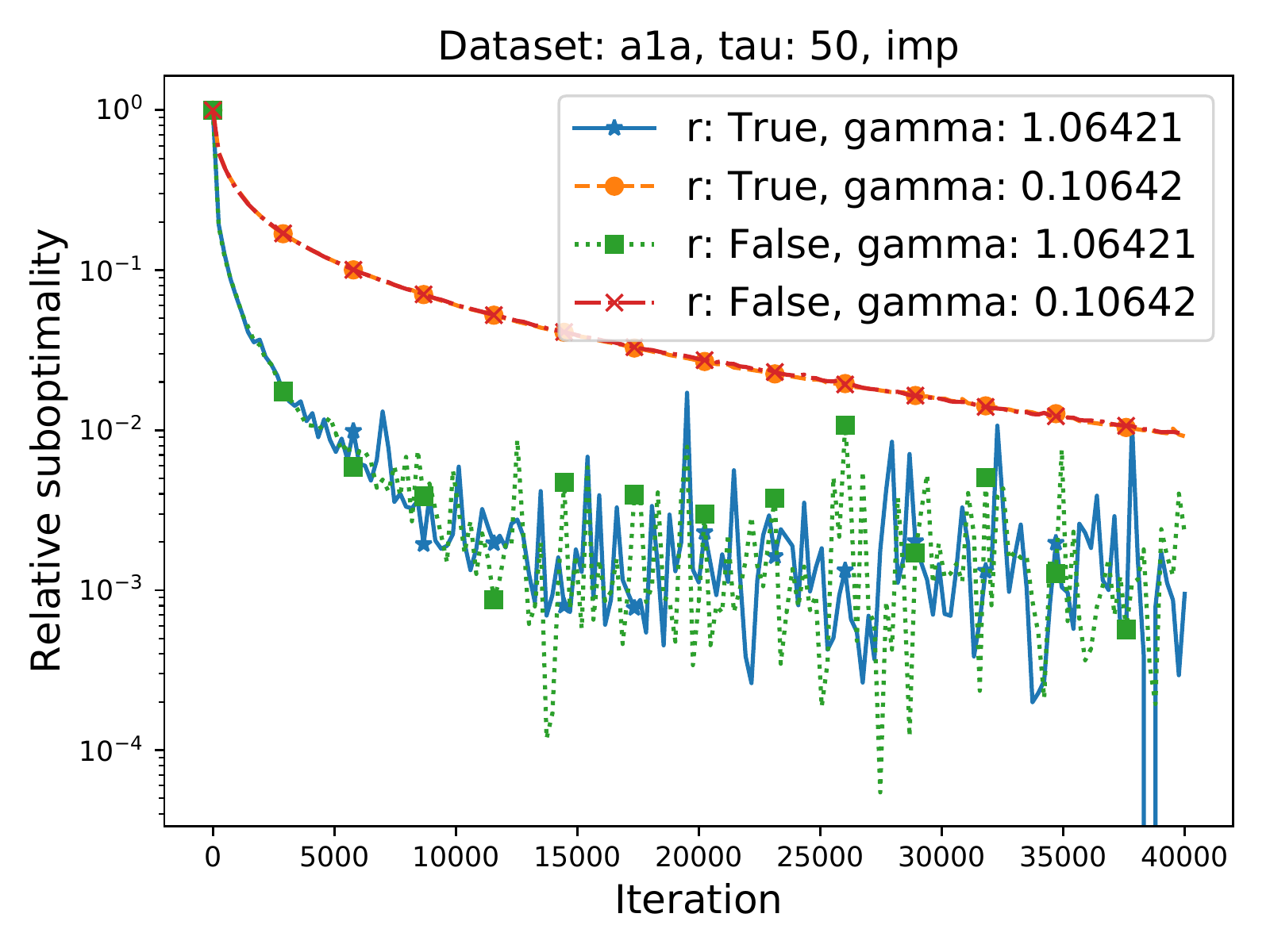}
\end{minipage}%
\\
\begin{minipage}{0.24\textwidth}
  \centering
\includegraphics[width =  \textwidth ]{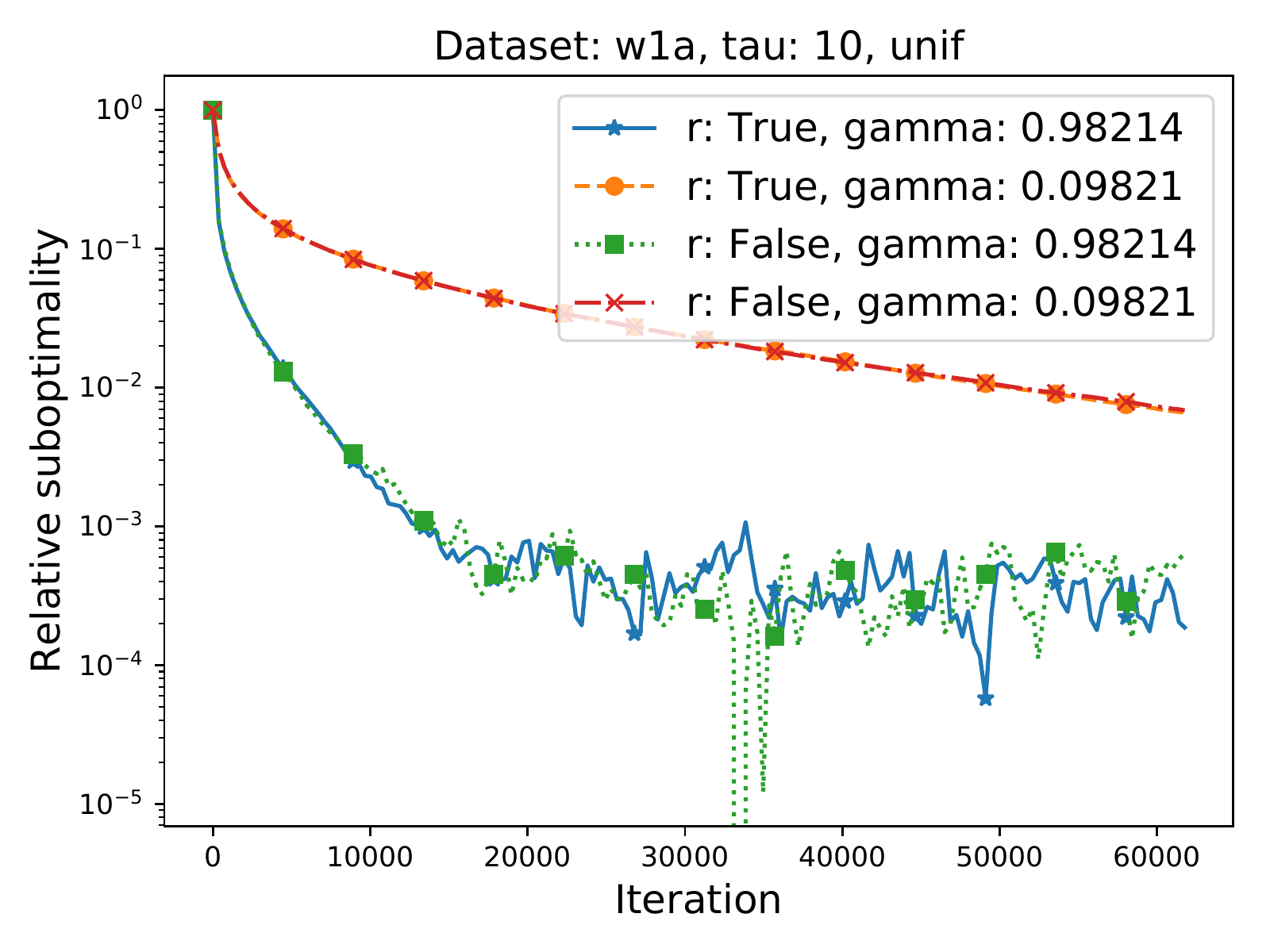}
\end{minipage}%
\begin{minipage}{0.24\textwidth}
  \centering
\includegraphics[width =  \textwidth ]{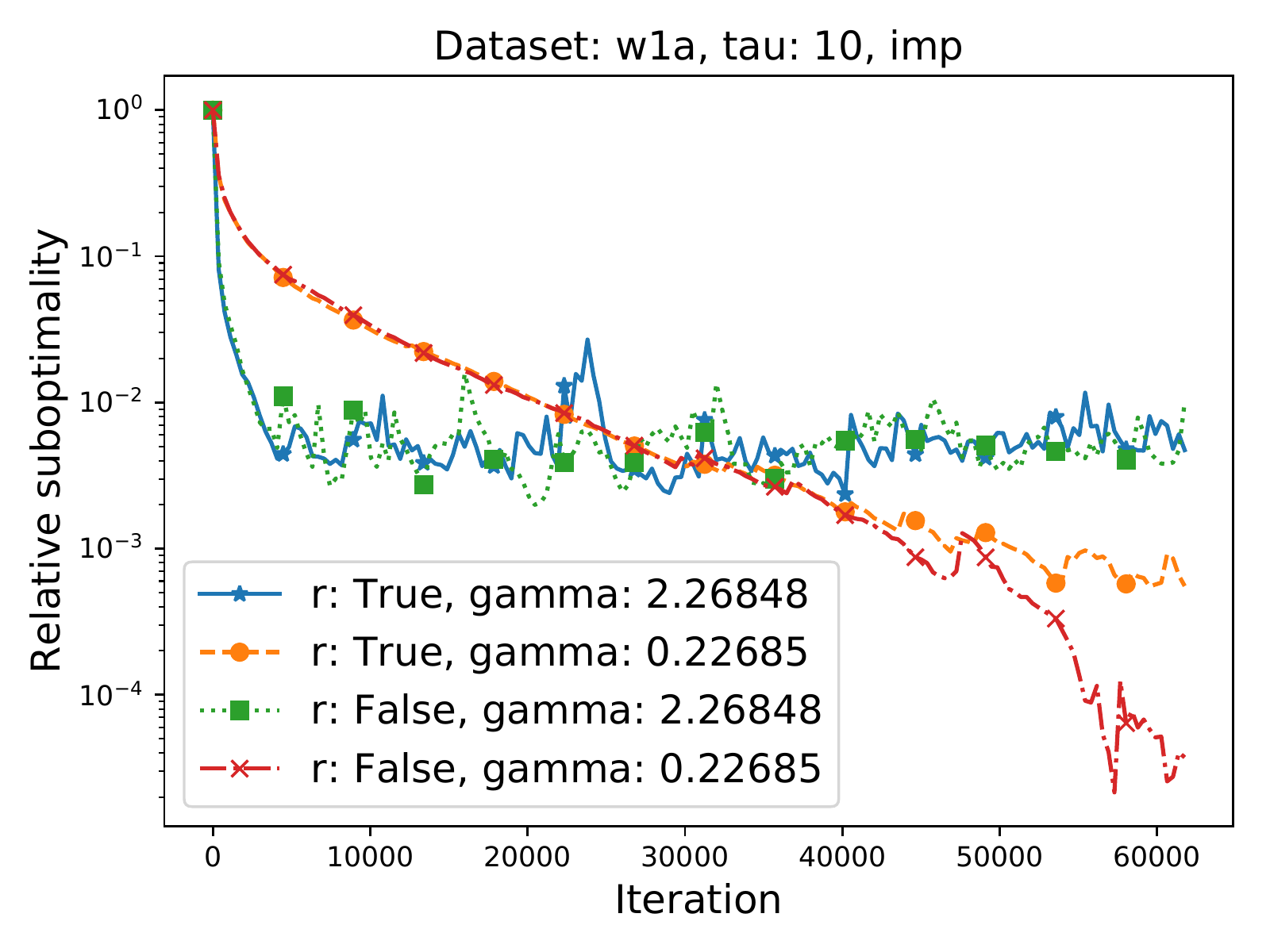}
\end{minipage}%
\begin{minipage}{0.24\textwidth}
  \centering
\includegraphics[width =  \textwidth ]{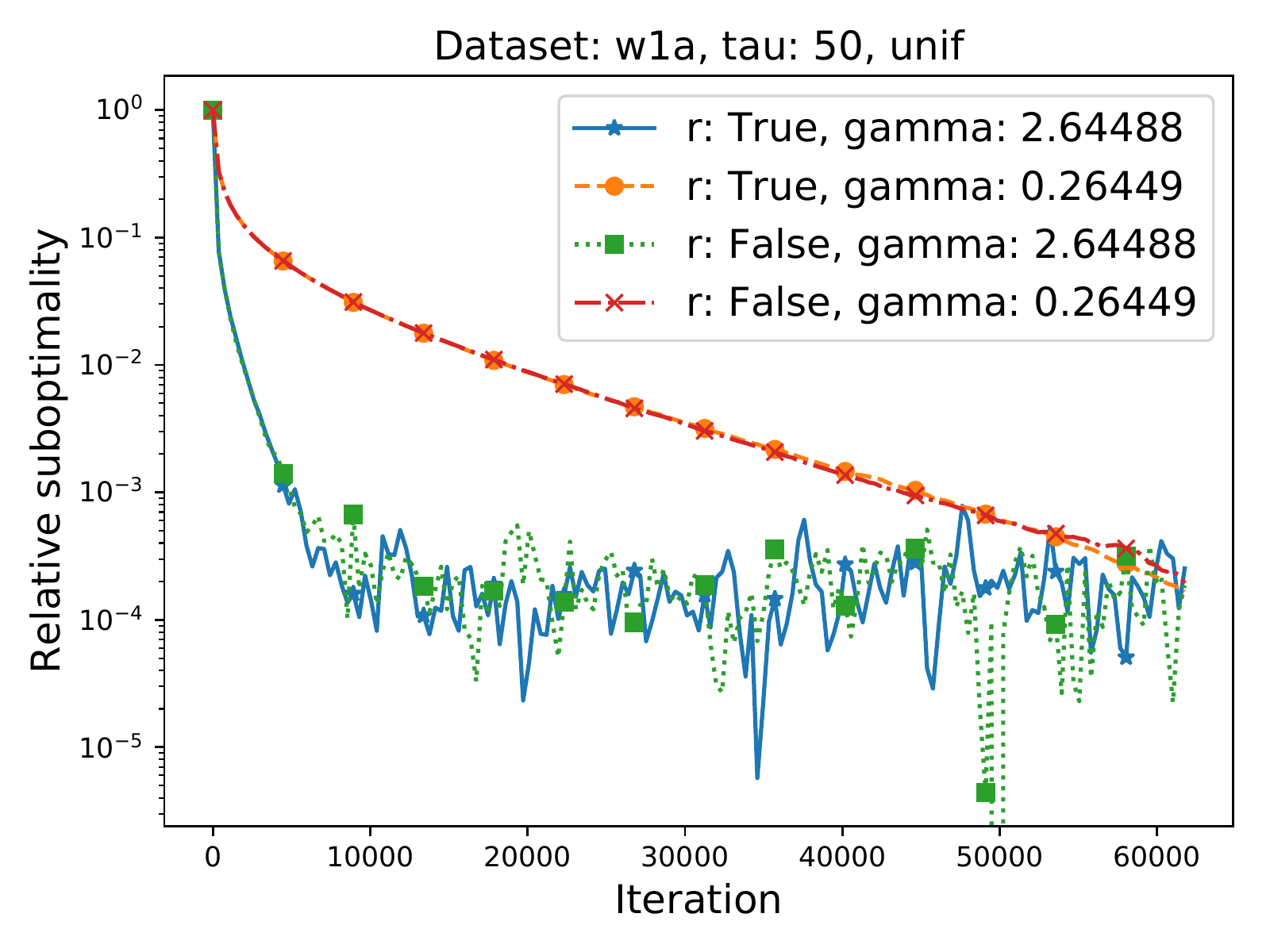}
\end{minipage}%
\begin{minipage}{0.24\textwidth}
  \centering
\includegraphics[width =  \textwidth ]{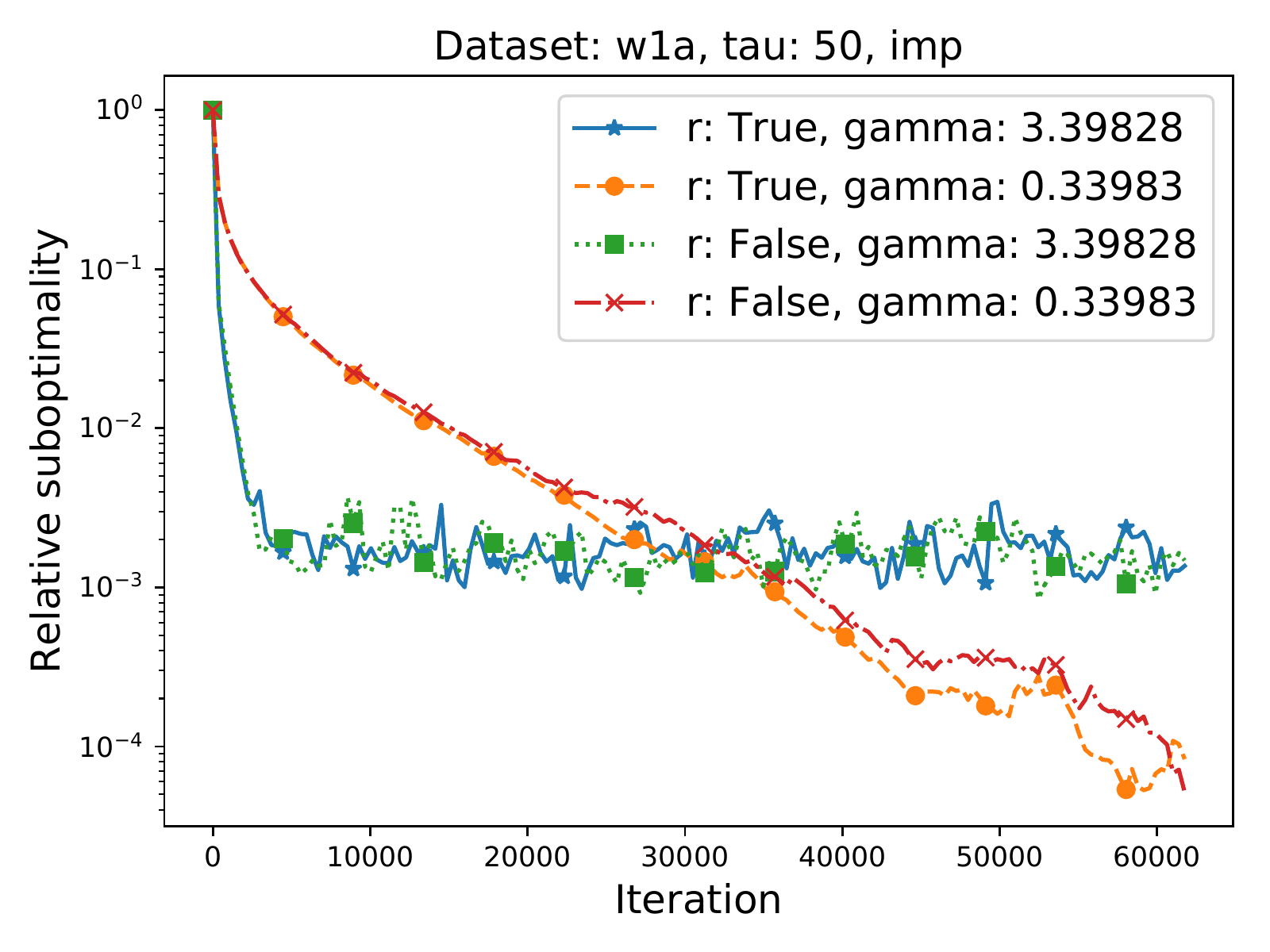}
\end{minipage}%
\\
\begin{minipage}{0.24\textwidth}
  \centering
\includegraphics[width =  \textwidth ]{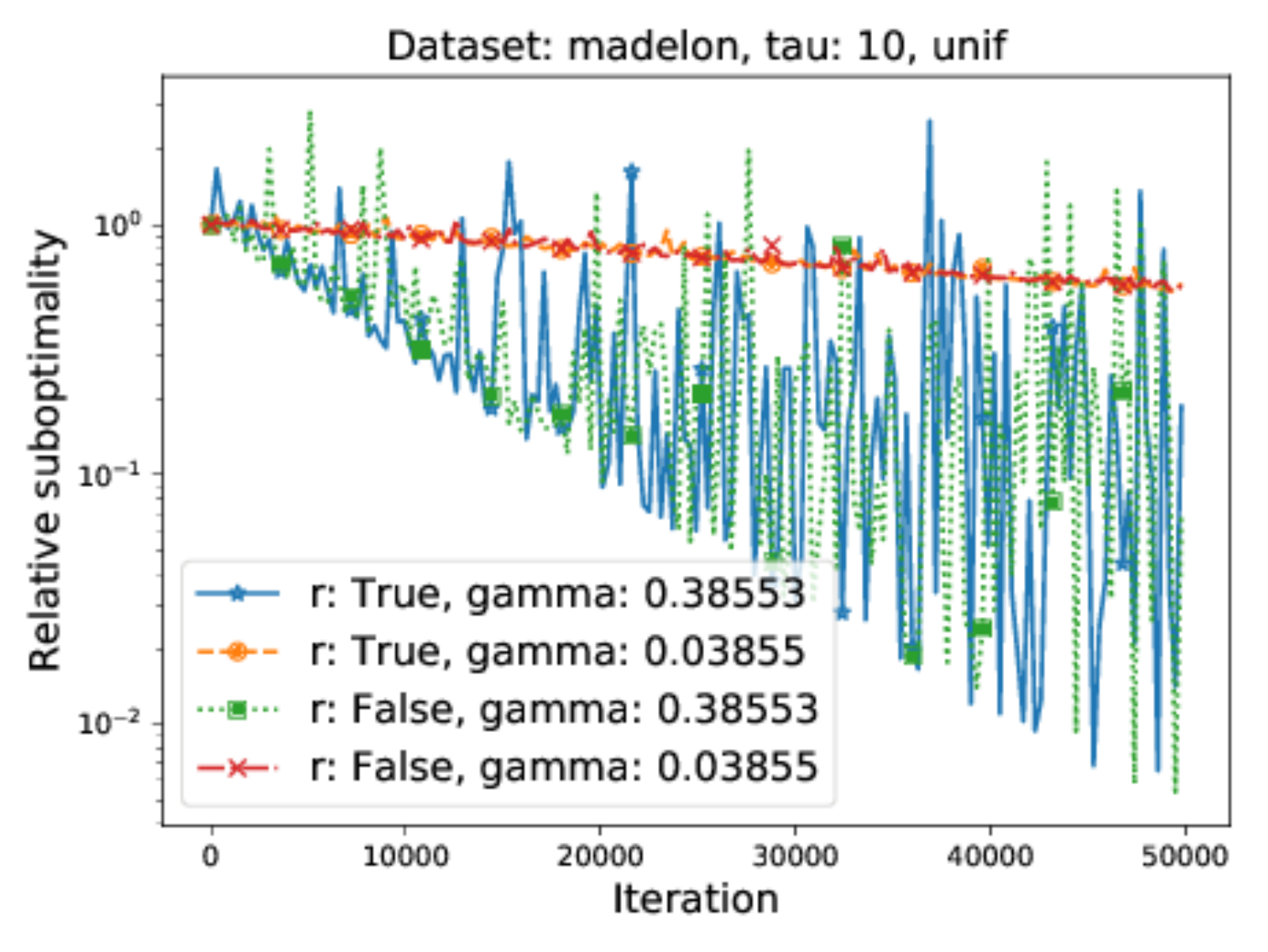}
\end{minipage}%
\begin{minipage}{0.24\textwidth}
  \centering
\includegraphics[width =  \textwidth ]{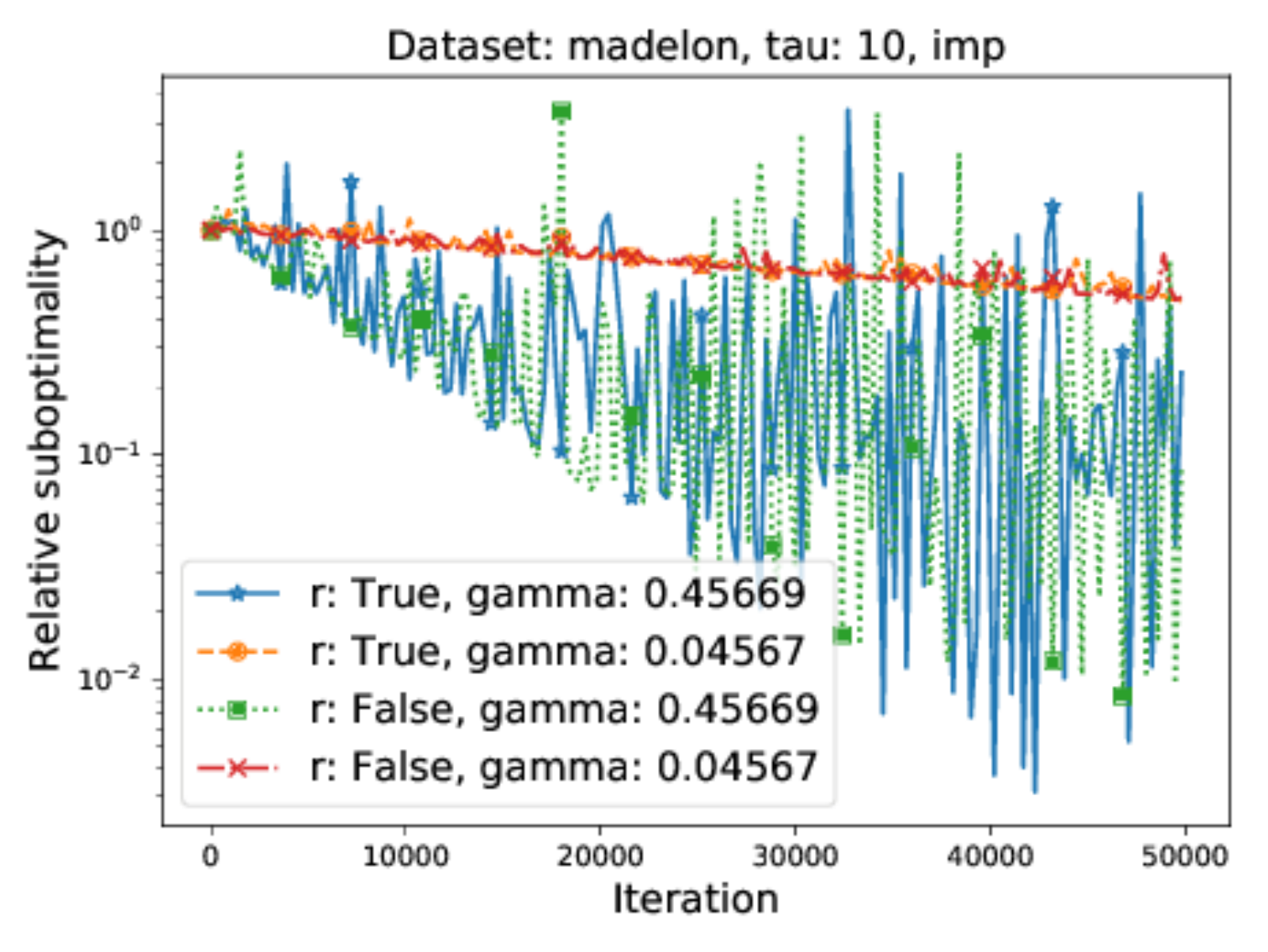}
\end{minipage}%
\begin{minipage}{0.24\textwidth}
  \centering
\includegraphics[width =  \textwidth ]{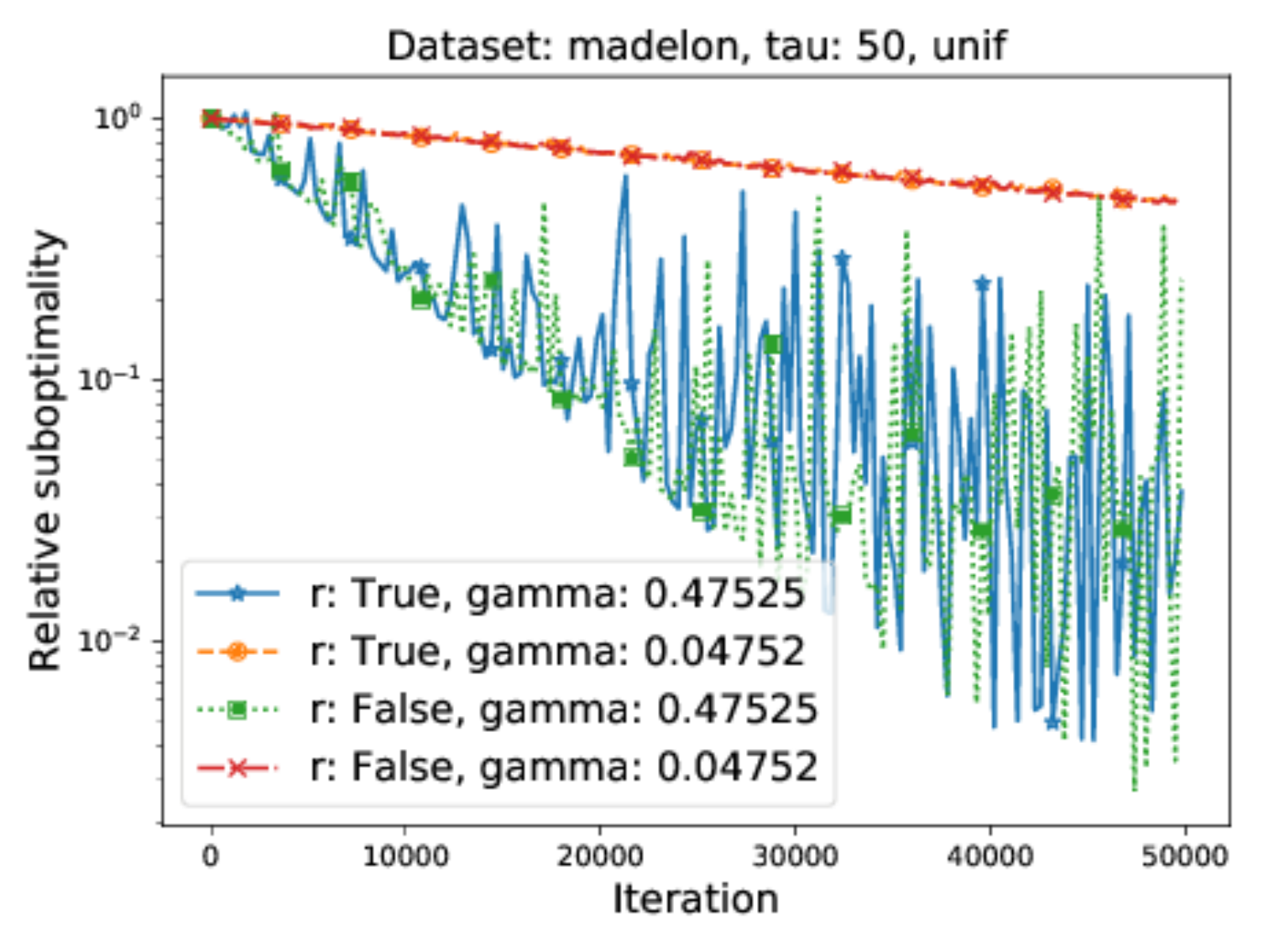}
\end{minipage}%
\begin{minipage}{0.24\textwidth}
  \centering
\includegraphics[width =  \textwidth ]{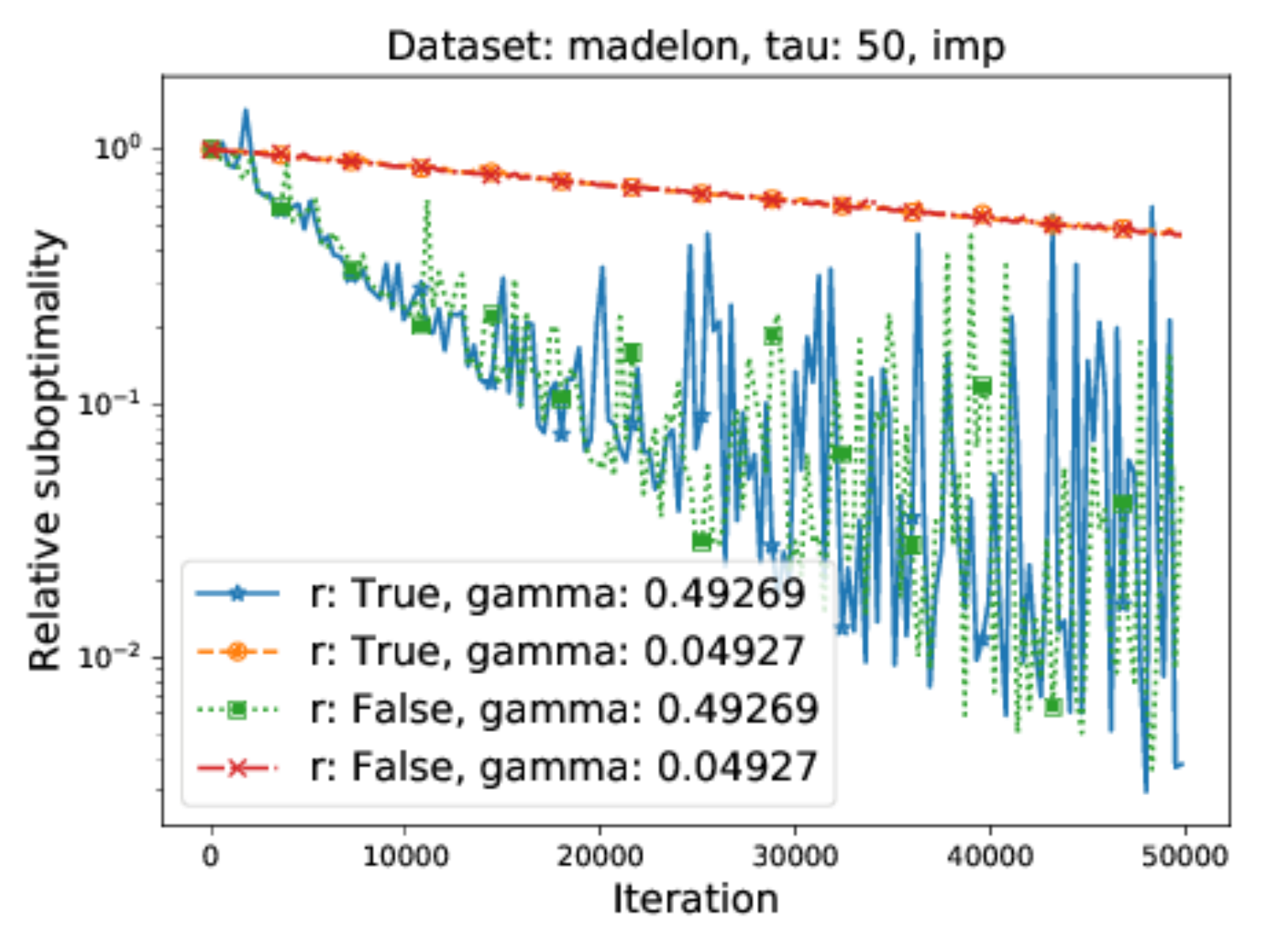}
\end{minipage}%
\\
\begin{minipage}{0.24\textwidth}
  \centering
\includegraphics[width =  \textwidth ]{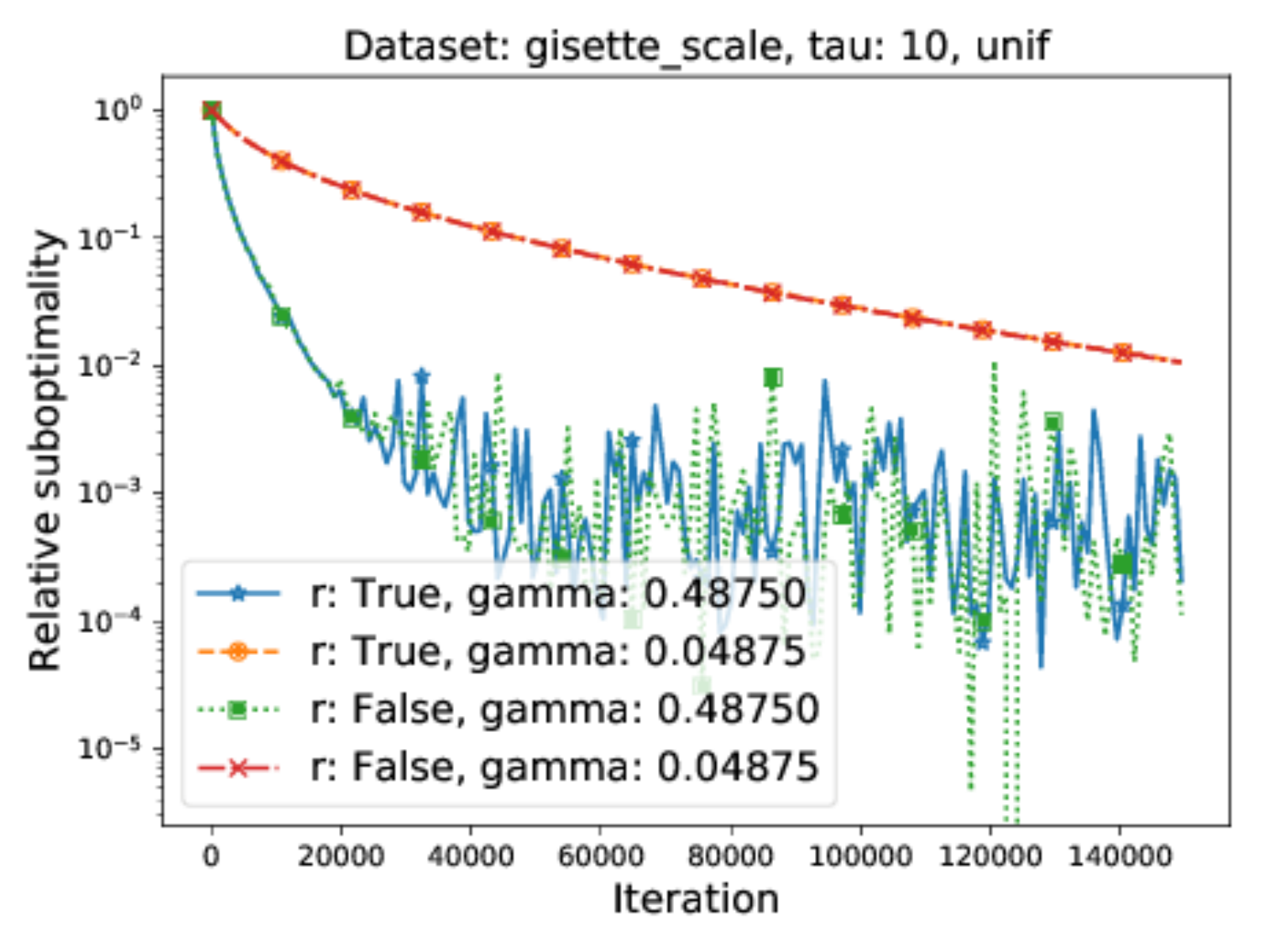}
\end{minipage}%
\begin{minipage}{0.24\textwidth}
  \centering
\includegraphics[width =  \textwidth ]{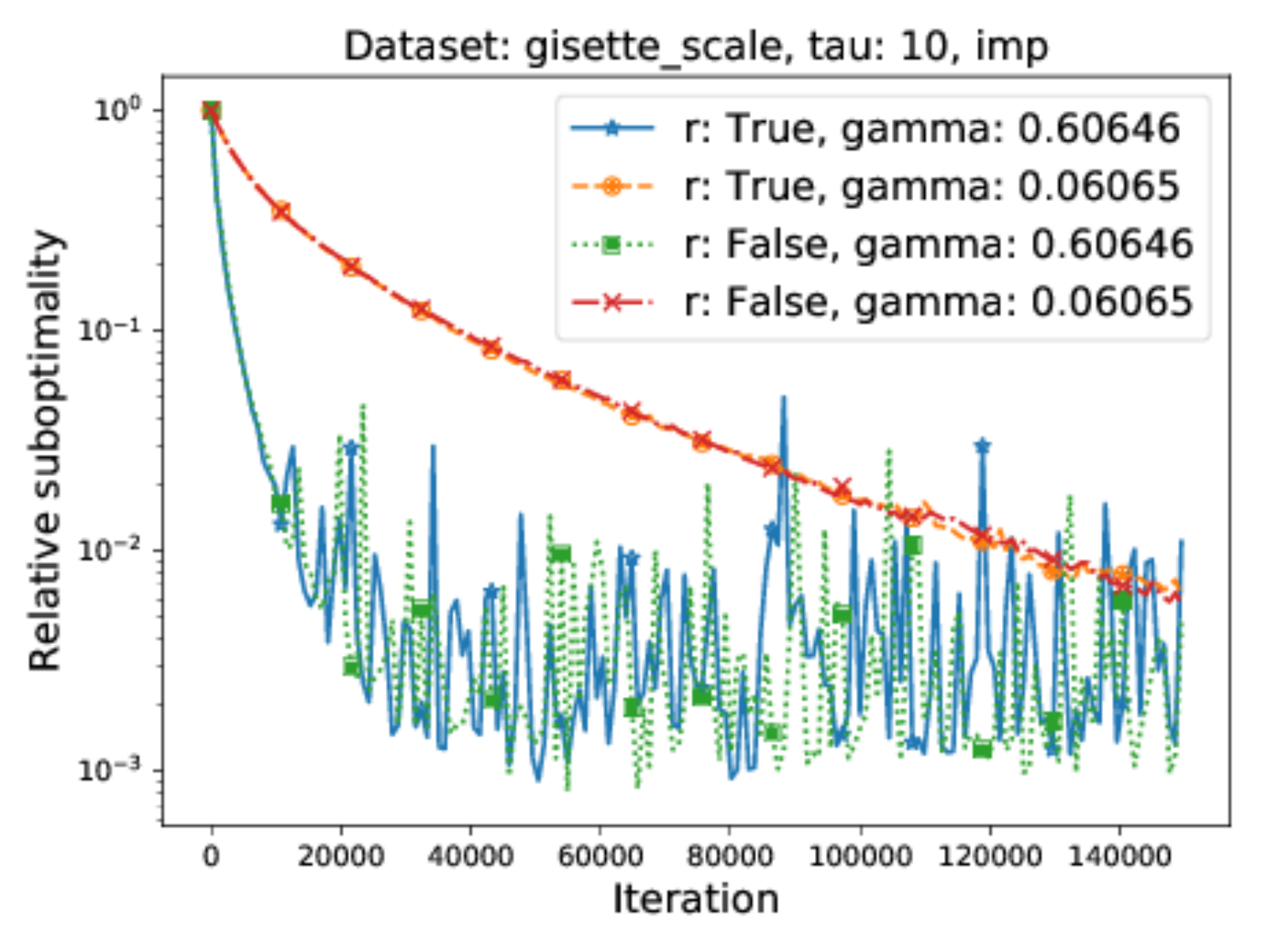}
\end{minipage}%
\begin{minipage}{0.24\textwidth}
  \centering
\includegraphics[width =  \textwidth ]{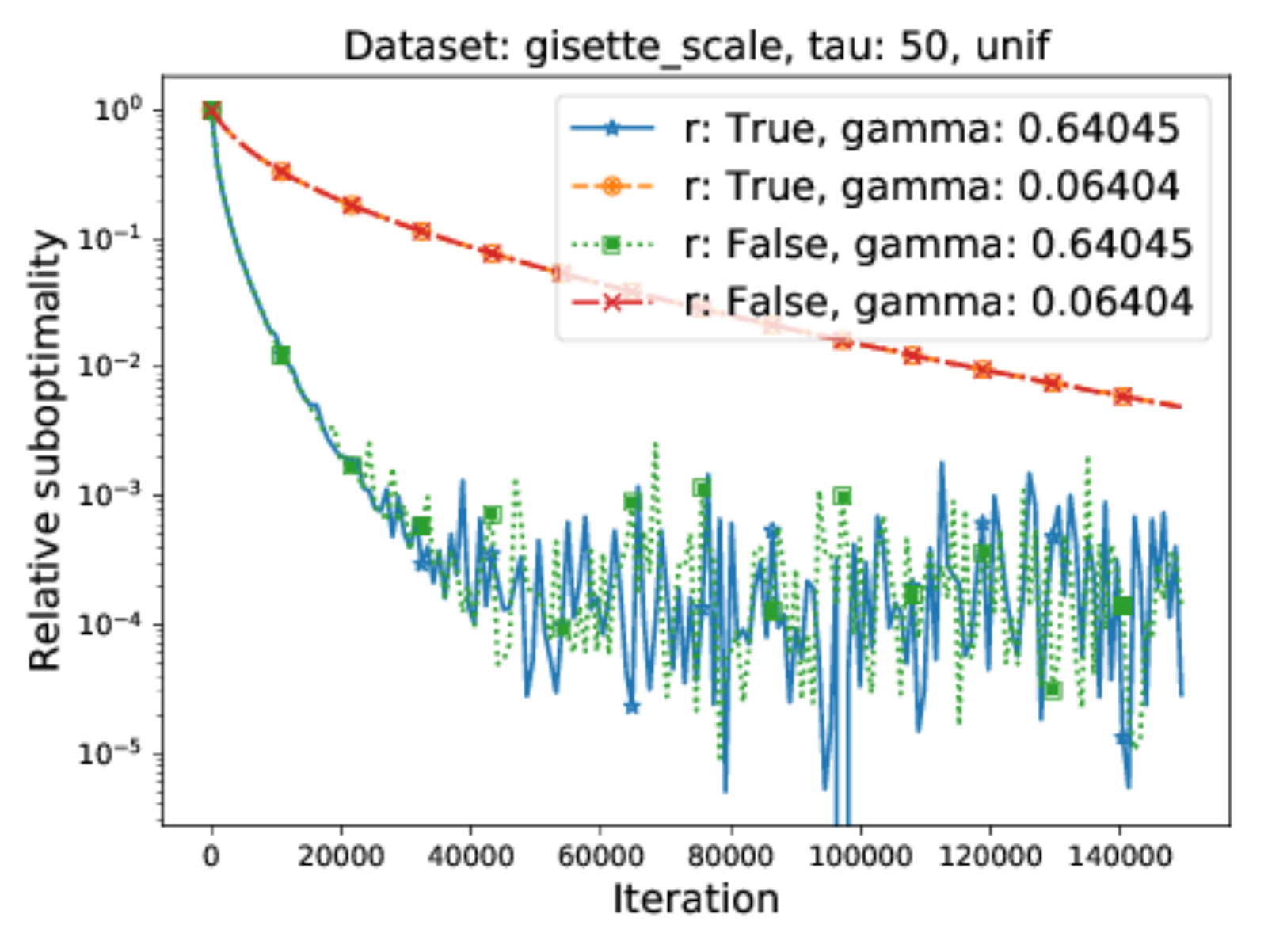}
\end{minipage}%
\begin{minipage}{0.24\textwidth}
  \centering
\includegraphics[width =  \textwidth ]{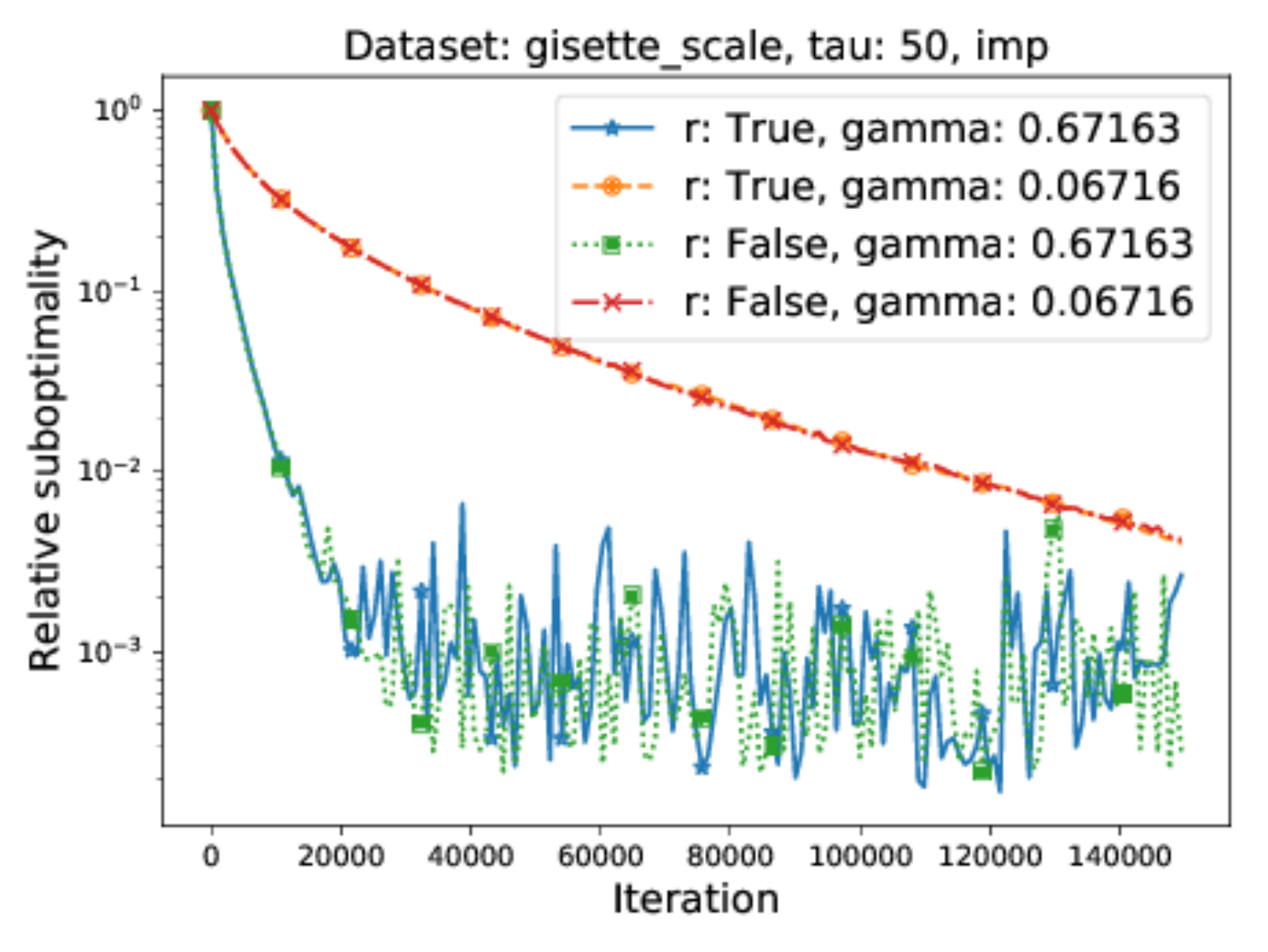}
\end{minipage}%
\\
\begin{minipage}{0.24\textwidth}
  \centering
\includegraphics[width =  \textwidth ]{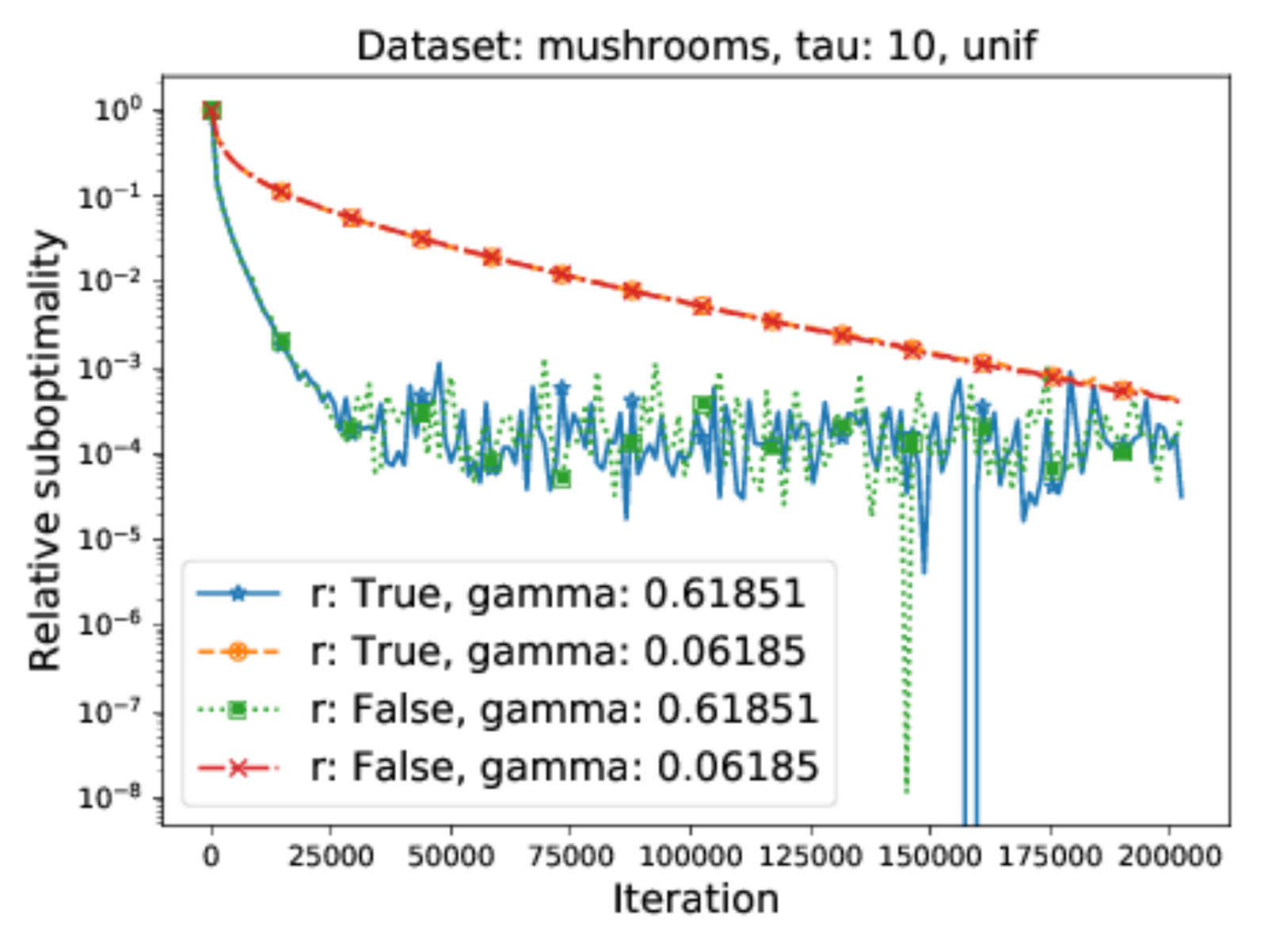}
\end{minipage}%
\begin{minipage}{0.24\textwidth}
  \centering
\includegraphics[width =  \textwidth ]{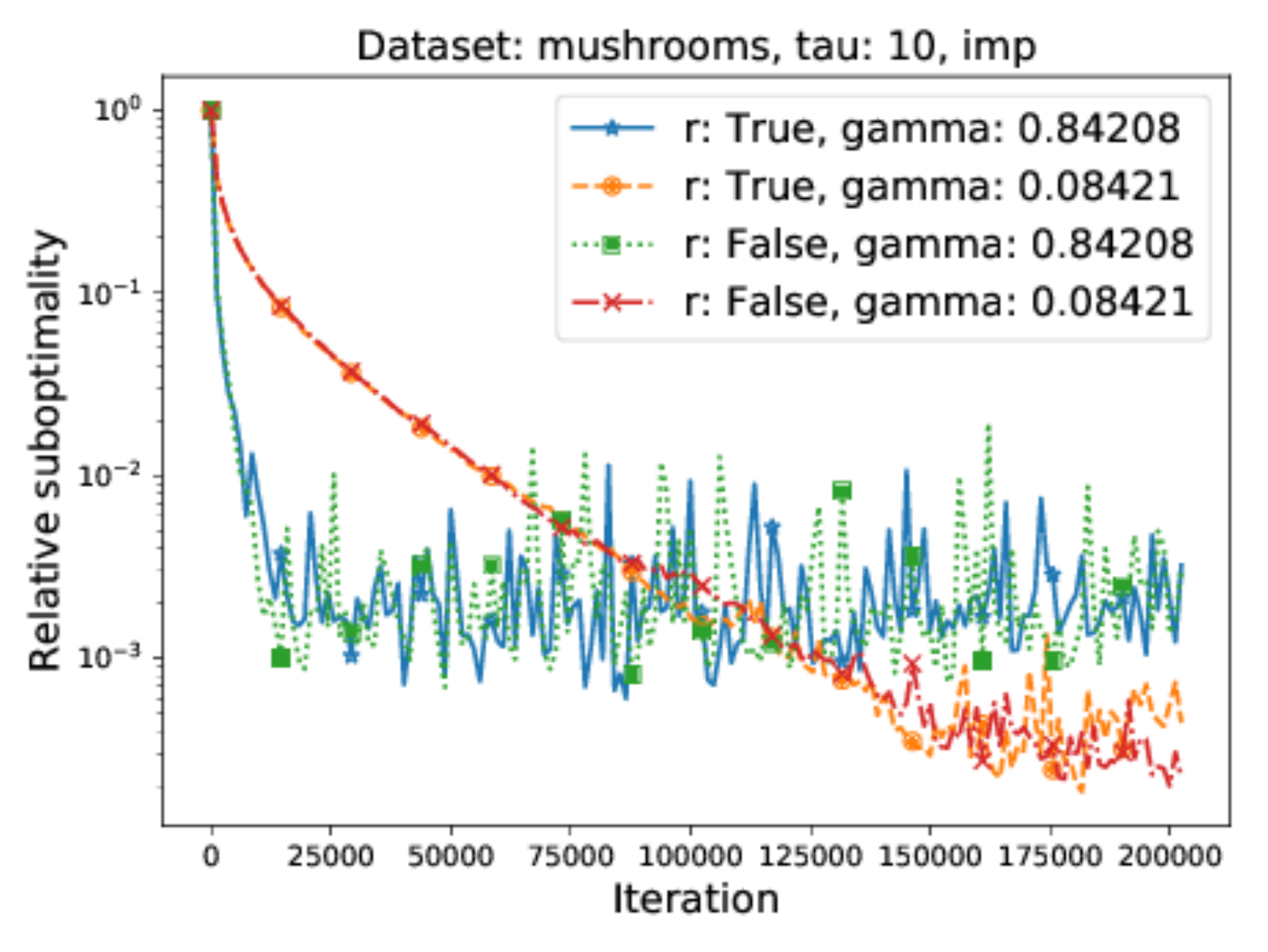}
\end{minipage}%
\begin{minipage}{0.24\textwidth}
  \centering
\includegraphics[width =  \textwidth ]{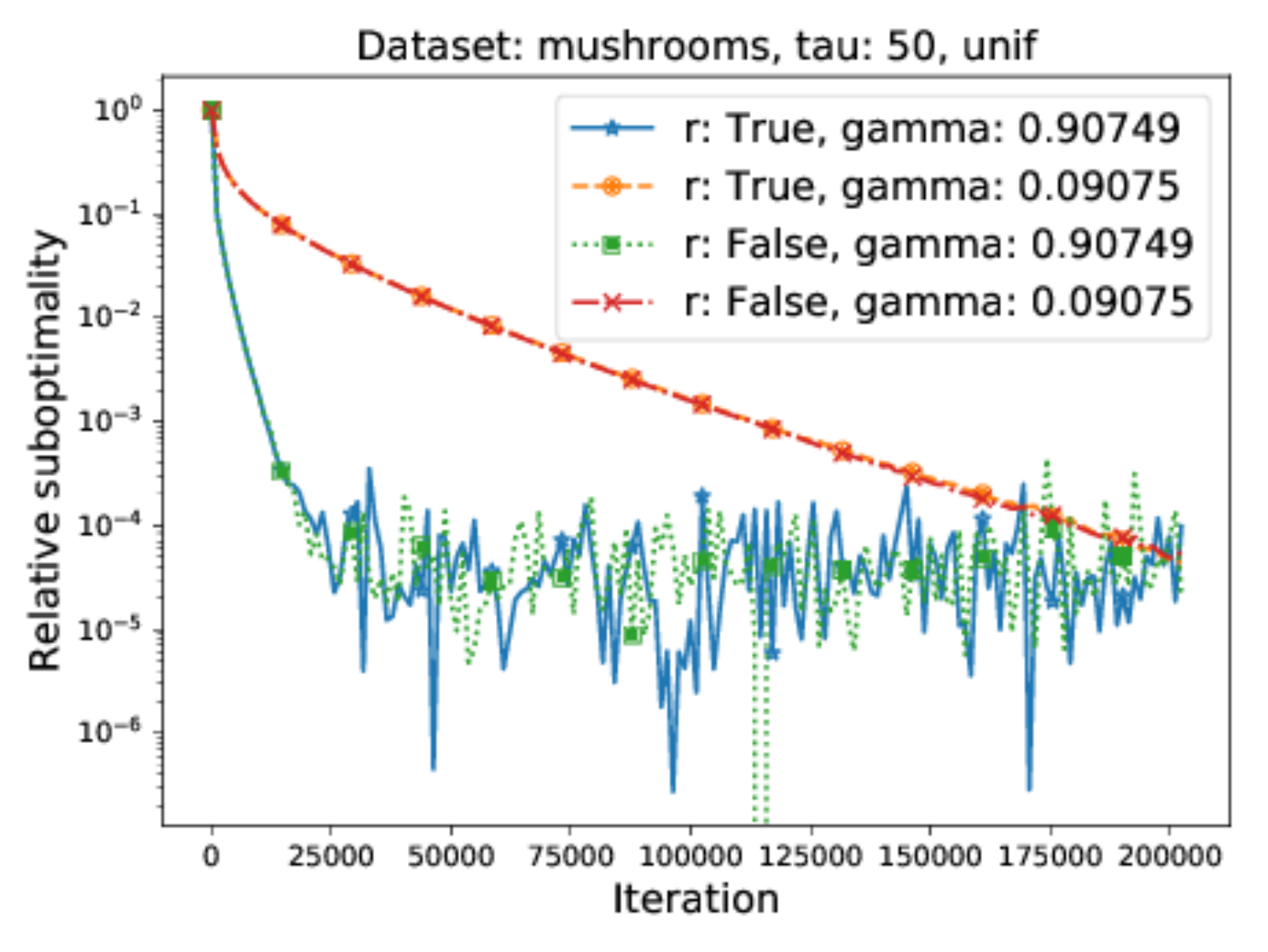}
\end{minipage}%
\begin{minipage}{0.24\textwidth}
  \centering
\includegraphics[width =  \textwidth ]{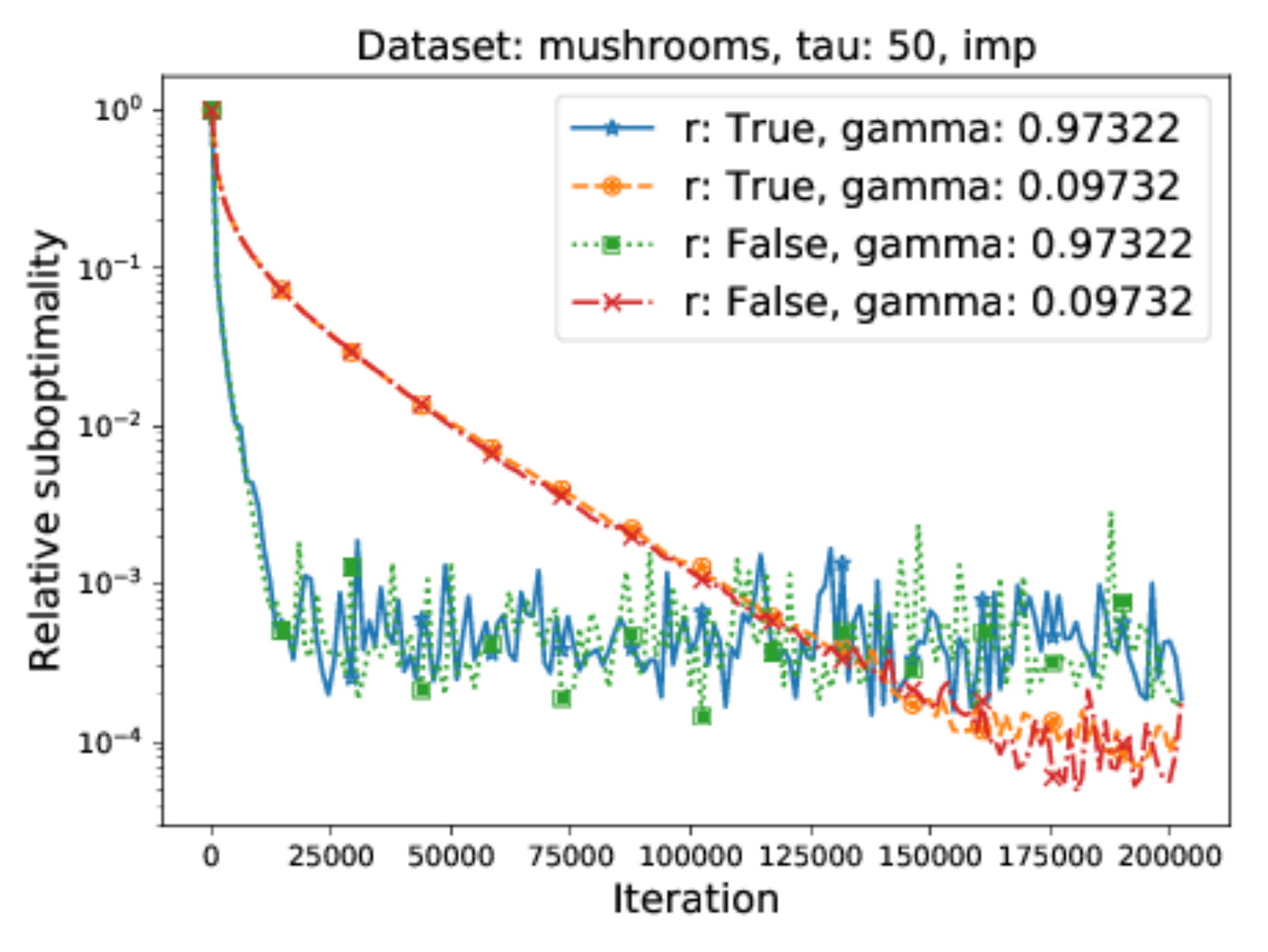}
\end{minipage}%
\caption{{\tt SGD-MB} and independent {\tt SGD} applied on LIBSVM~\citep{chang2011libsvm} datasets with regularization parameter $\lambda = 10^{-5}$. Axis $y$ stands for relative suboptimality, i.e. $\frac{f(x^k)-f(x^*)}{f(x^k)-f(x^0)}$. Title label ``unif'' corresponds to probabilities chosen by~\ref{item:unif} while label ``imp'' corresponds to probabilities chosen by~\ref{item:imp}. Lastly, legend label ``r'' corresponds to ``replacement'' with value ``True'' for {\tt SGD-MB} and value ``False'' for independent {\tt SGD}.}
\label{fig:SGDMB}
\end{figure}

In Section~\ref{sec:SGD-MB}, we describe in detail the {\tt SGD-MB} method already outlined before. The main advantage of {\tt SGD-MB} is that the sampling procedure it employs can be implemented in just $\cO(\tau \log n)$ time. In contrast, even the simplest without-replacement sampling which selects each function into the minibatch with a prescribed probability independently (we will refer to it as independent {\tt SGD}) requires $n$ calls of a uniform random generator.  We demonstrate numerically that {\tt SGD-MB} has essentially identical iteration complexity to  independent {\tt SGD} in practice. We consider logistic regression with Tikhonov regularization:
\begin{equation}\label{eq:logreg}
\frac1n \sum_{i=1}^n   \log \left(1+\exp\left(a_i^\top x\cdot  b_i\right) \right)+\frac{\lambda}{2} \norm{ x}^2,
\end{equation}
For a fixed  expected sampling size $\tau$, consider two options for the probability of sampling the $i$-th function:
\begin{enumerate}
\item \label{item:unif} $\frac{\tau}{n}$, or
\item \label{item:imp} $\frac{\norm{a_i}^2+\lambda}{\delta+\norm{a_i}^2+\lambda}$, where $\delta$ is such that\footnote{An {\tt RCD} version of this sampling was proposed in~\citep{AccMbCd}; it was shown to be superior to uniform sampling both in theory and practice.} $\sum_{i=1}^n\frac{\norm{a_i}^2+\lambda}{\delta+\norm{a_i}^2+\lambda}=1$. 
\end{enumerate}
The results can be found in Figure~\ref{fig:SGDMB}, where we also report the choice of stepsize $\gamma$ and the choice of $\tau$ in the legend and title of the plot, respectively.

Indeed, iteration complexity of {\tt SGD-MB} and independent {\tt SGD} is almost identical. Since the cost of each iteration of {\tt SGD-MB} is cheaper\footnote{The relative difference between iteration costs of {\tt SGD-MB} and independent {\tt SGD} can be arbitrary, especially for the case when cost of evaluating $\nabla f_i(x)$ is cheap, $n$ is huge and $n\gg \tau$. In such case, cost of one iteration of {\tt SGD-MB} is $\tau \text{Cost}(\nabla f_i) +\tau\log(n)$ while the cost of one iteration of independent {\tt SGD} is $\tau \text{Cost}(\nabla f_i) + n$.}, we conclude superiority of {\tt SGD-MB} to independent {\tt SGD}.


\subsection{Experiments on {\tt SGD-star} \label{sec:exp_star}}
In this section, we study {\tt SGD-star} and numerically verify claims from Section~\ref{sec:SGD-star}. In particular, Corollary~\ref{cor:SGD-star} shows that {\tt SGD-star} enjoys linear convergence rate which is constant times better to the rate of {\tt SAGA} (given that problem condition number is high enough). We compare 3 methods -- {\tt SGD-star}, {\tt SGD} and {\tt SAGA}. We consider simple and well-understood least squares problem $\min_x \frac12 \| \mA x-b\|^2$ where elements of $\mA,b$ were generated (independently) from standard normal distribution. Further, rows of $\mA$ were normalized so that $\|\mA_{i:}\|=1$. Thus, denoting $f_i(x) = \frac12 (\mA_{i:}^\top x-b_i )^2$, $f_i$ is 1-smooth. For simplicity, we consider {\tt SGD-star} with uniform serial sampling, i.e. $\cL=1$.

\begin{figure}[H]
\centering
\begin{minipage}{0.3\textwidth}
  \centering
\includegraphics[width =  \textwidth ]{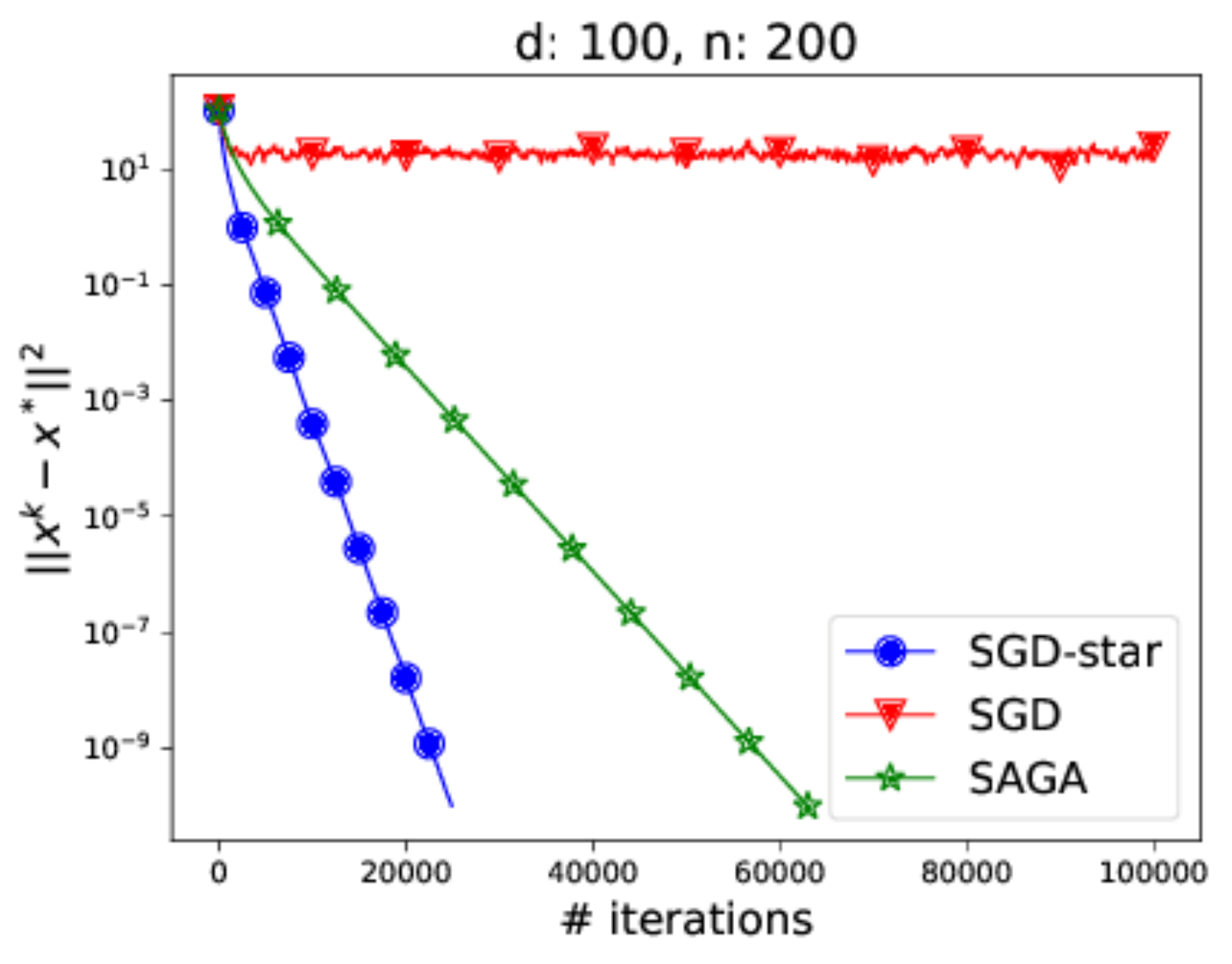}
\end{minipage}%
\begin{minipage}{0.3\textwidth}
  \centering
\includegraphics[width =  \textwidth ]{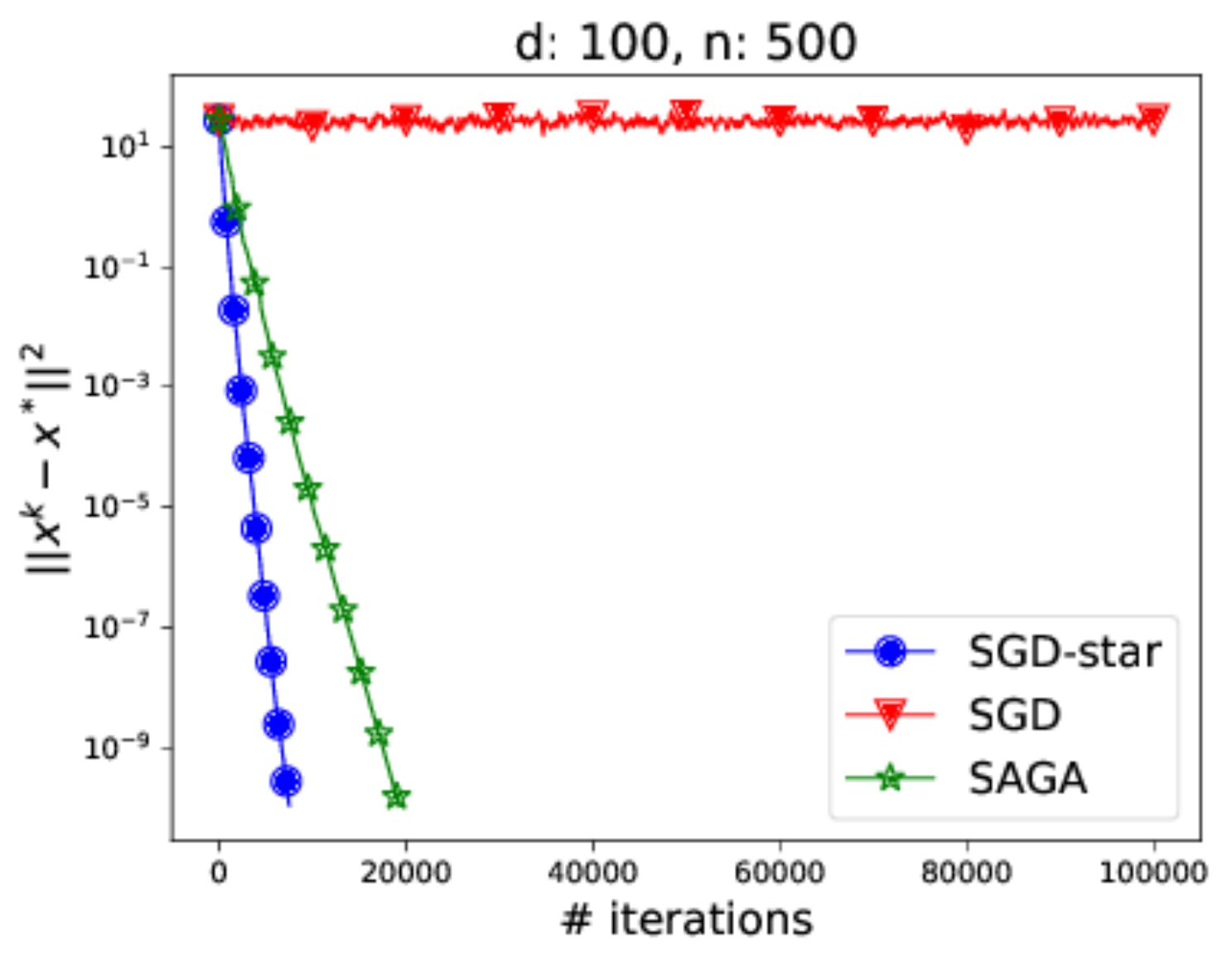}
\end{minipage}%
\begin{minipage}{0.3\textwidth}
  \centering
\includegraphics[width =  \textwidth ]{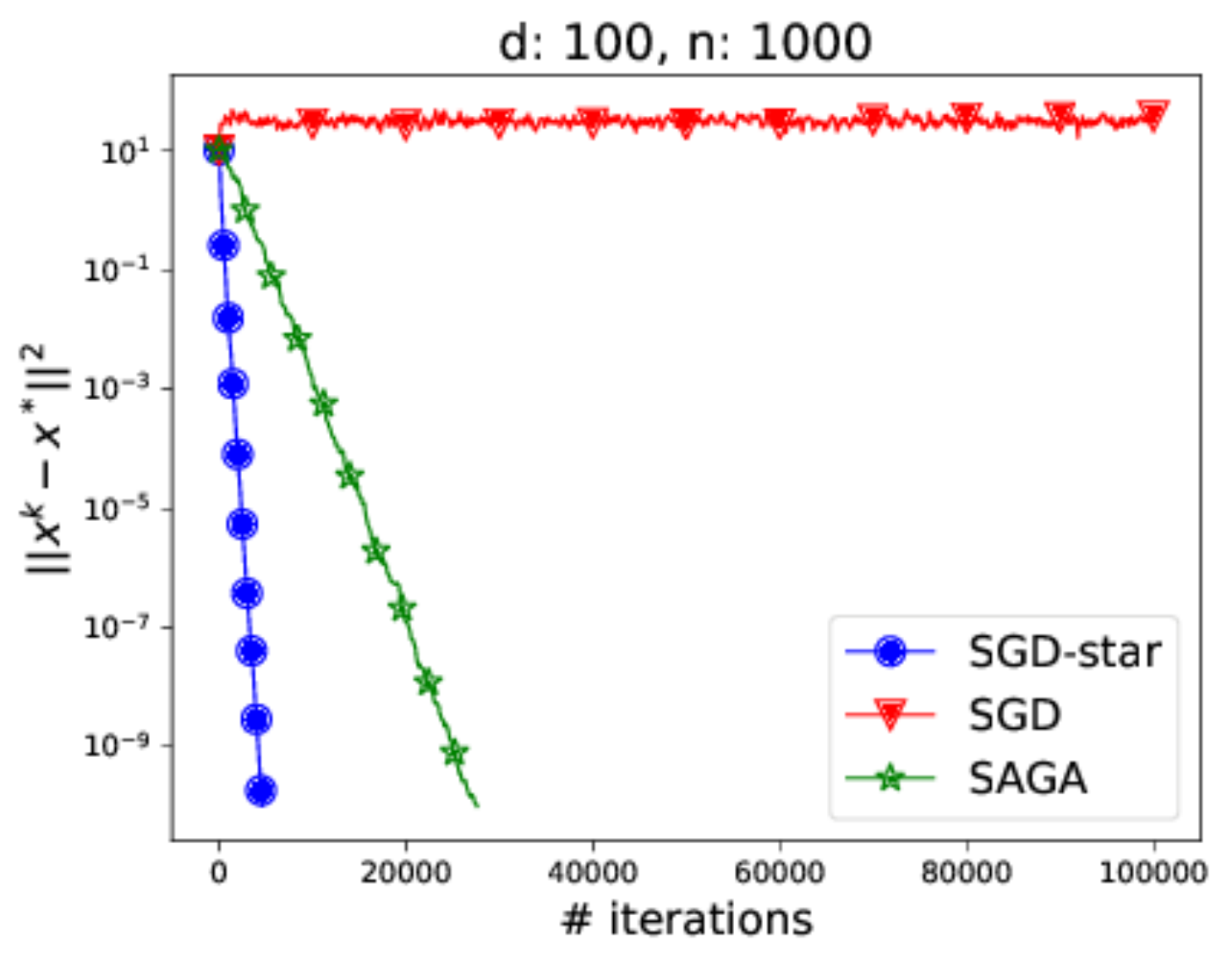}
\end{minipage}%
\\
\begin{minipage}{0.3\textwidth}
  \centering
\includegraphics[width =  \textwidth ]{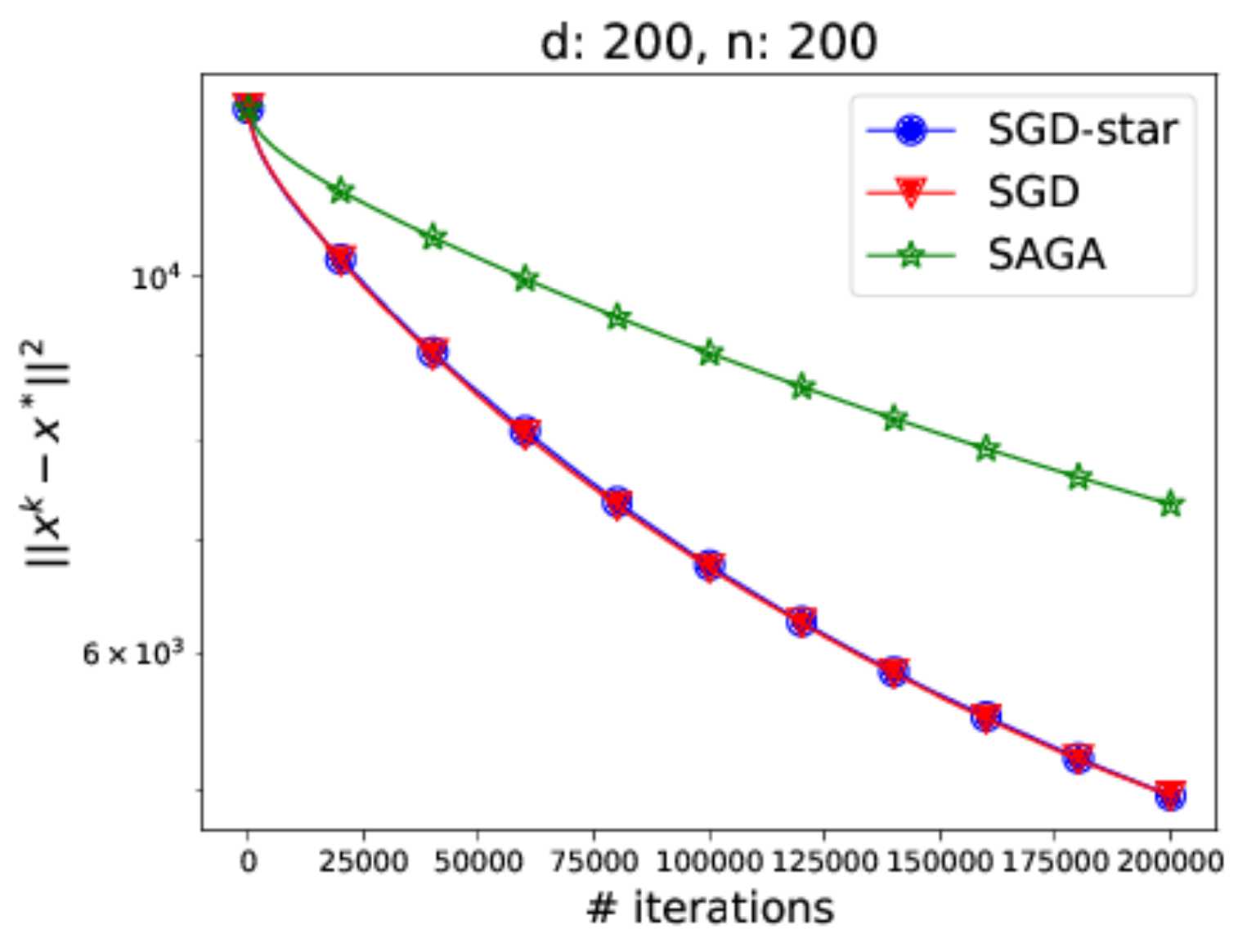}
\end{minipage}%
\begin{minipage}{0.3\textwidth}
  \centering
\includegraphics[width =  \textwidth ]{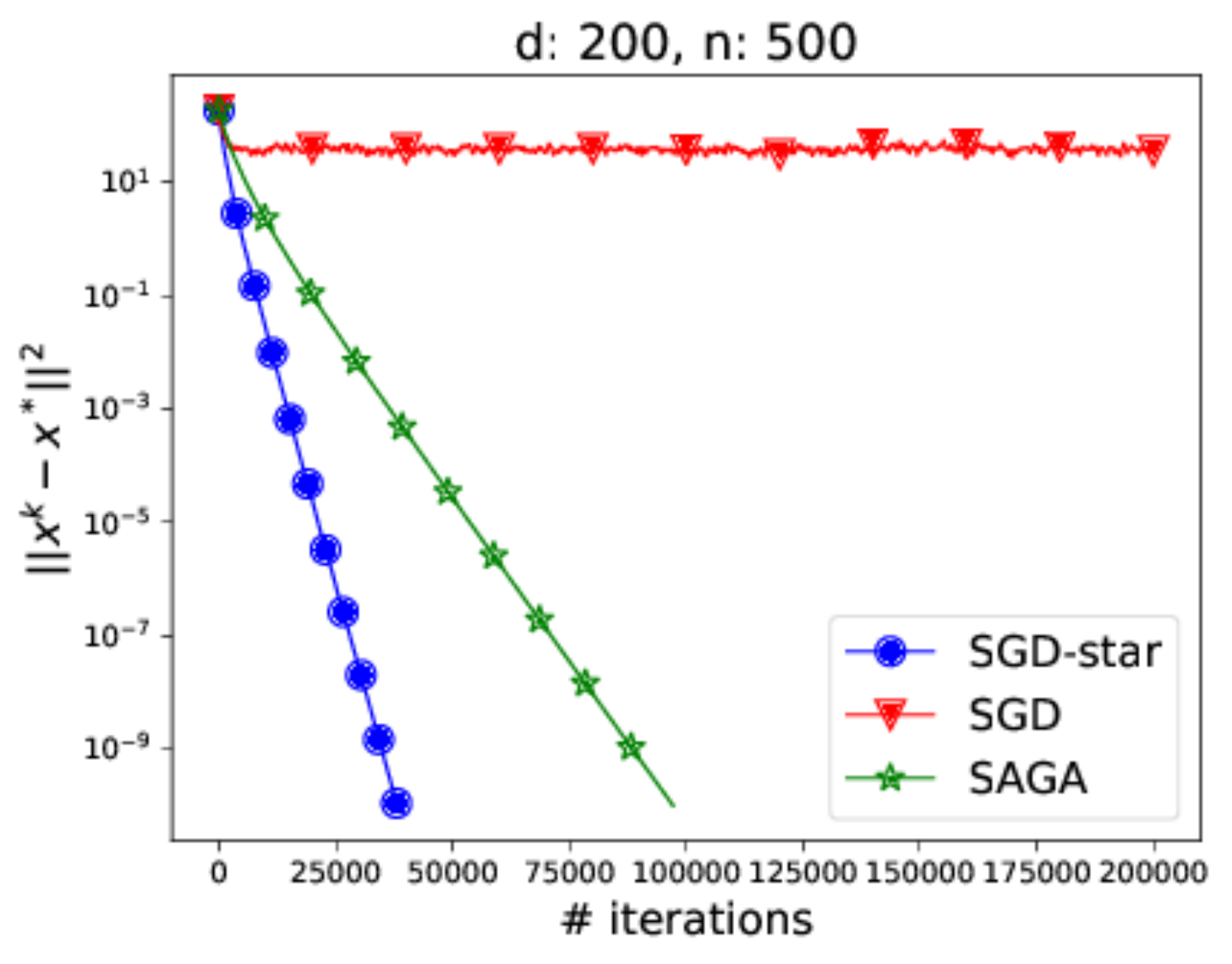}
\end{minipage}%
\begin{minipage}{0.3\textwidth}
  \centering
\includegraphics[width =  \textwidth ]{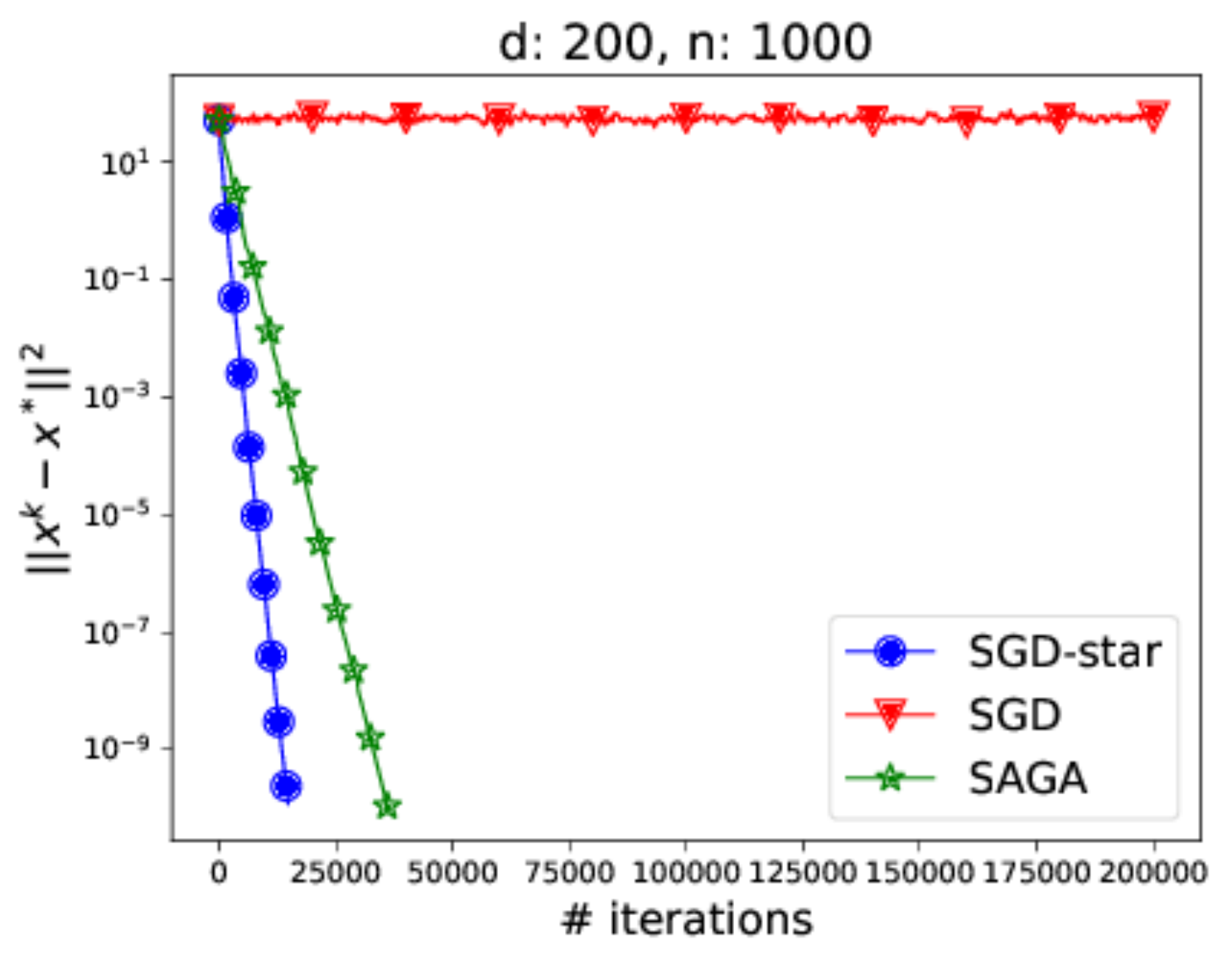}
\end{minipage}%
\\
\begin{minipage}{0.3\textwidth}
  \centering
\includegraphics[width =  \textwidth ]{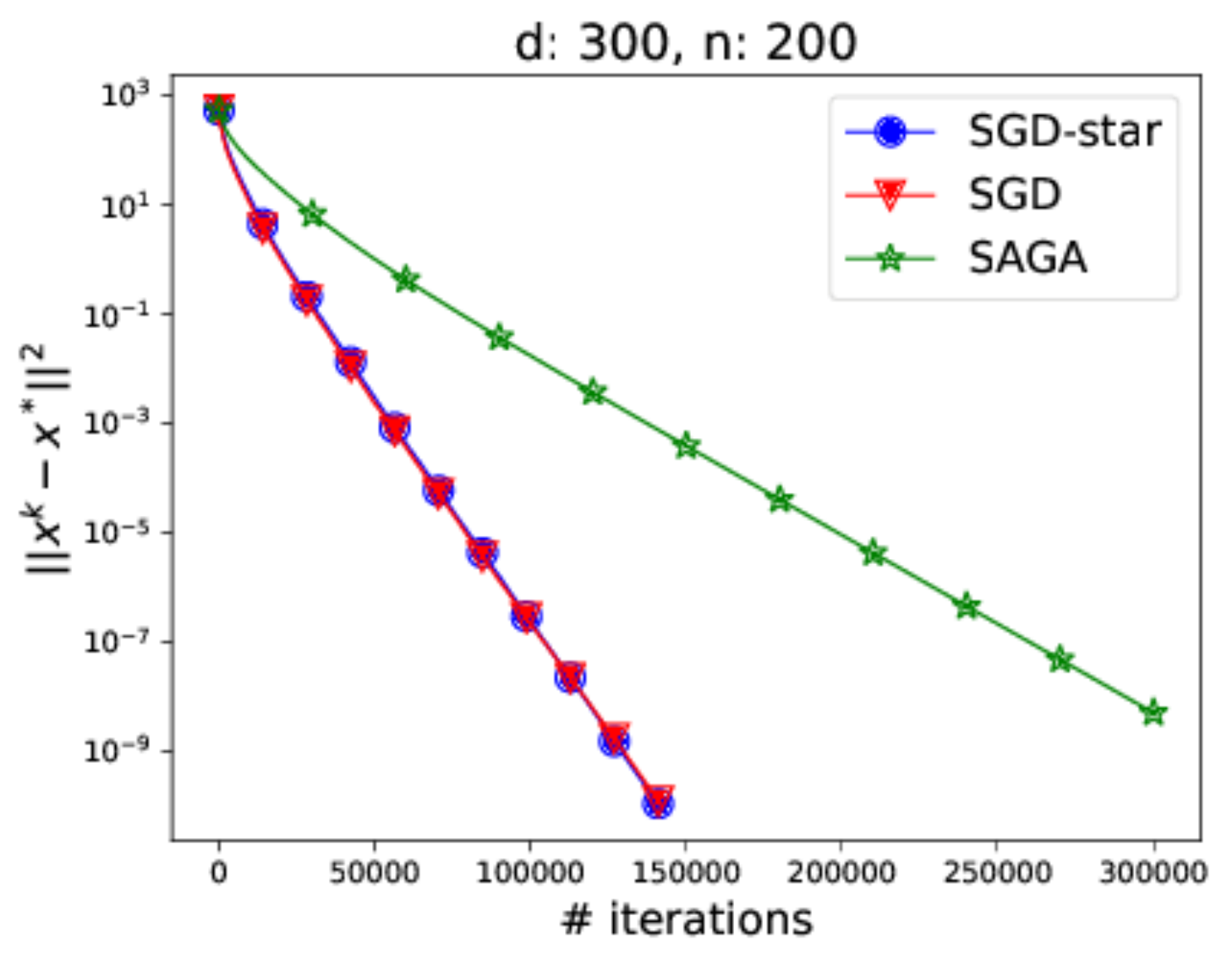}
\end{minipage}%
\begin{minipage}{0.3\textwidth}
  \centering
\includegraphics[width =  \textwidth ]{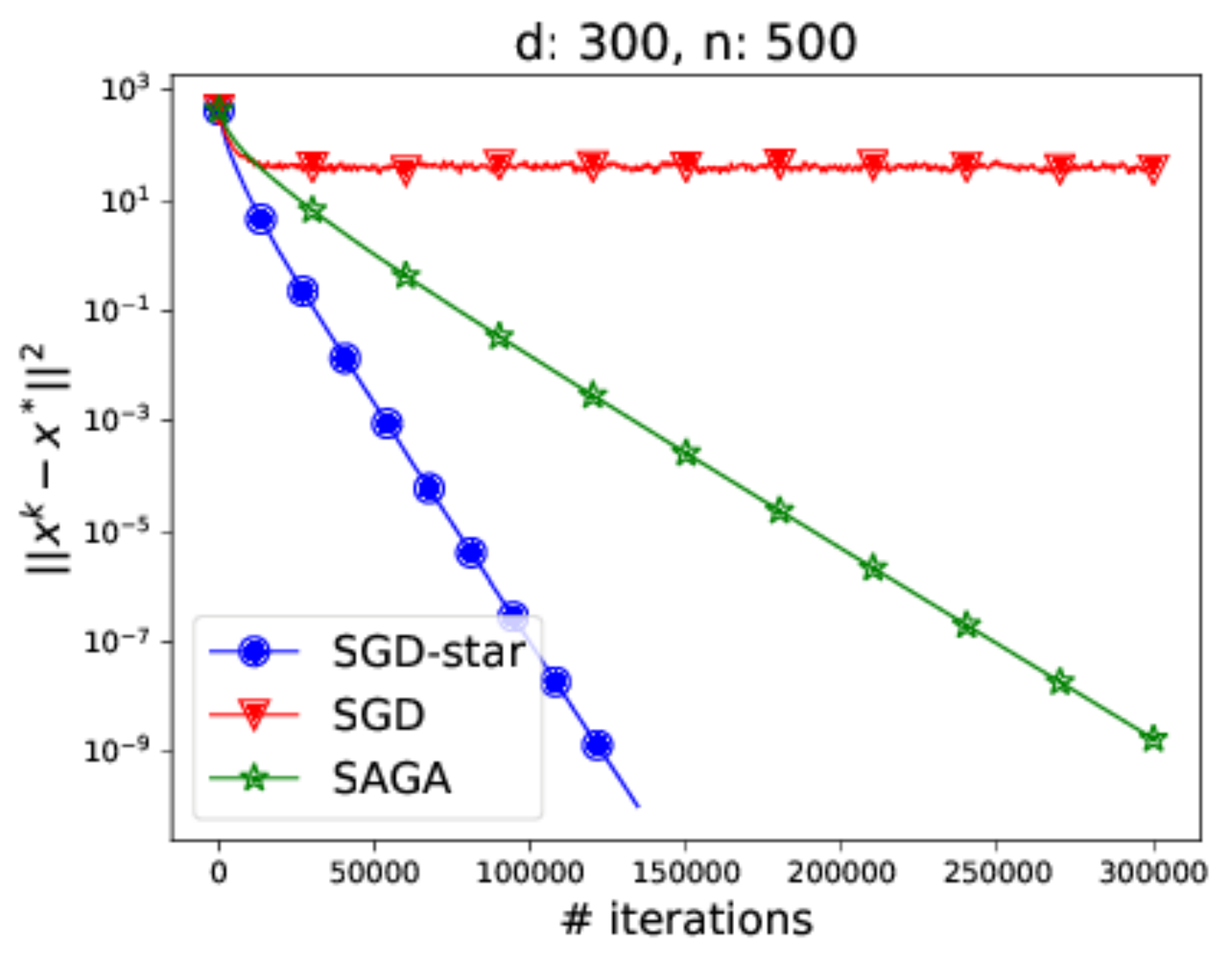}
\end{minipage}%
\begin{minipage}{0.3\textwidth}
  \centering
\includegraphics[width =  \textwidth ]{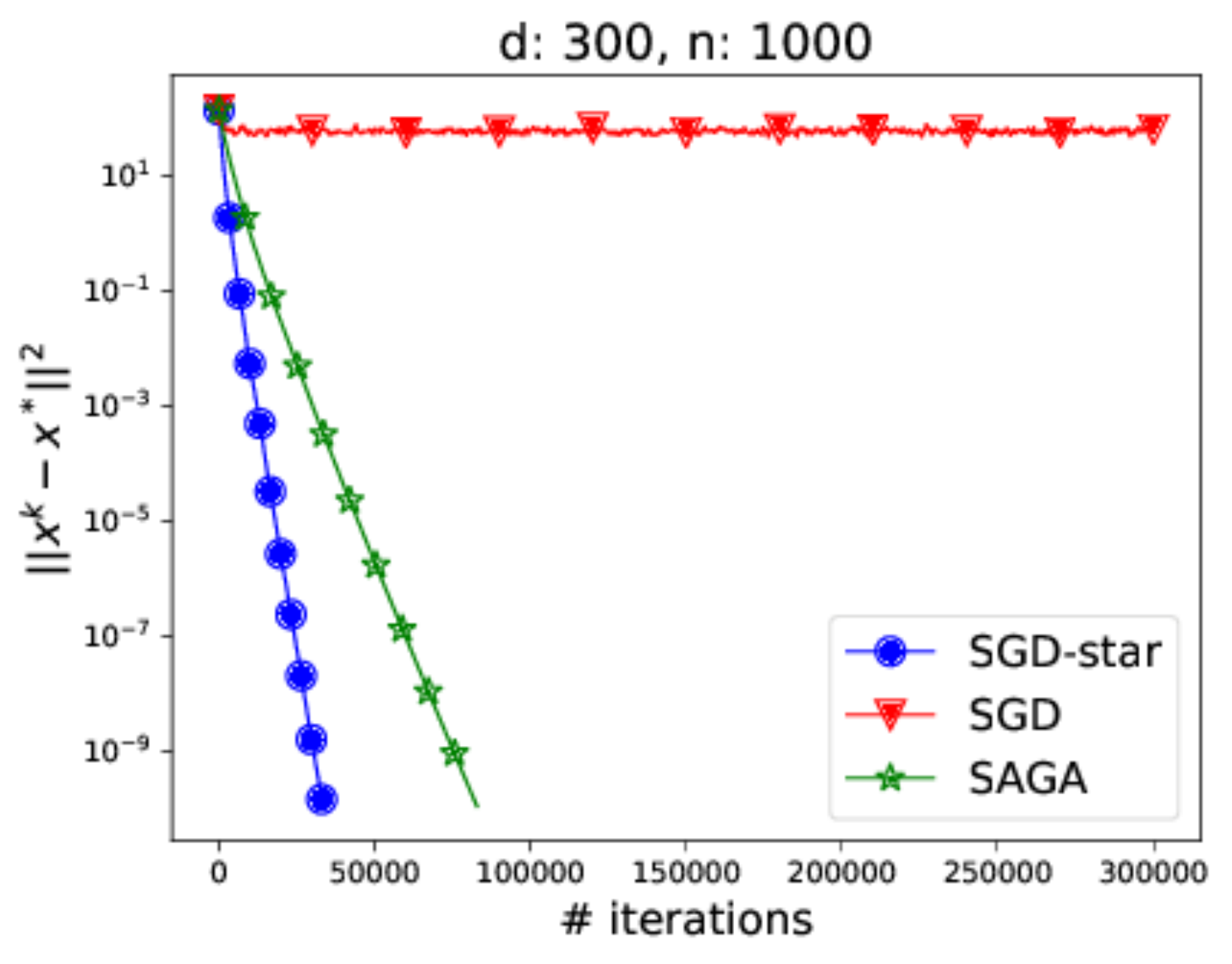}
\end{minipage}%
\caption{Comparison of {\tt SGD-star}, {\tt SGD} and {\tt SAGA} on least squares problem.}
\label{fig:star}
\end{figure}

 Next, for both {\tt SGD-star} and {\tt SGD} we use stepsize $\gamma = \frac{1}{2}$ (theory supported stepsize for {\tt SGD-star}), while for {\tt SAGA} we set $\gamma = \frac{1}{5}$ (almost theory supported stepsize). Figure~\ref{fig:star} shows the results.

Note that, as theory predicts, {\tt SGD-star} is always faster to {\tt SAGA}, although only constant times. Further, in the cases where $d\geq n$, performance of {\tt SGD} seems identical to the performance of {\tt SGD-shift}. This is due to a simple reason: if $d\geq n$, we must have $\nabla f_i(x^*) = 0$ for all $i$, and thus {\tt SGD} and {\tt SGD-shift} are in fact identical algorithms. 

\subsection{Experiments on {\tt N-SEGA} \label{sec:exp_nsega}}
In this experiment we study the effect of noise on {\tt N-SEGA}. We consider unit ball constrained least squares problem: $\min_{\|x\|\leq 1} f(x)$ where $f(x)=\|\mA x-b\|^2$. and we suppose that there is an oracle providing us with 
noised partial derivative $g_i(x,\zeta) = \nabla_i f(x) +\zeta$, where $\zeta \sim N(0,\sigma^2)$. For each problem instance (i.e. pair $\mA, b$), we compare performance of  {\tt N-SEGA} under various noise magnitudes $\sigma^2$.

\begin{figure}[H]
\centering
\begin{minipage}{0.24\textwidth}
  \centering
\includegraphics[width =  \textwidth ]{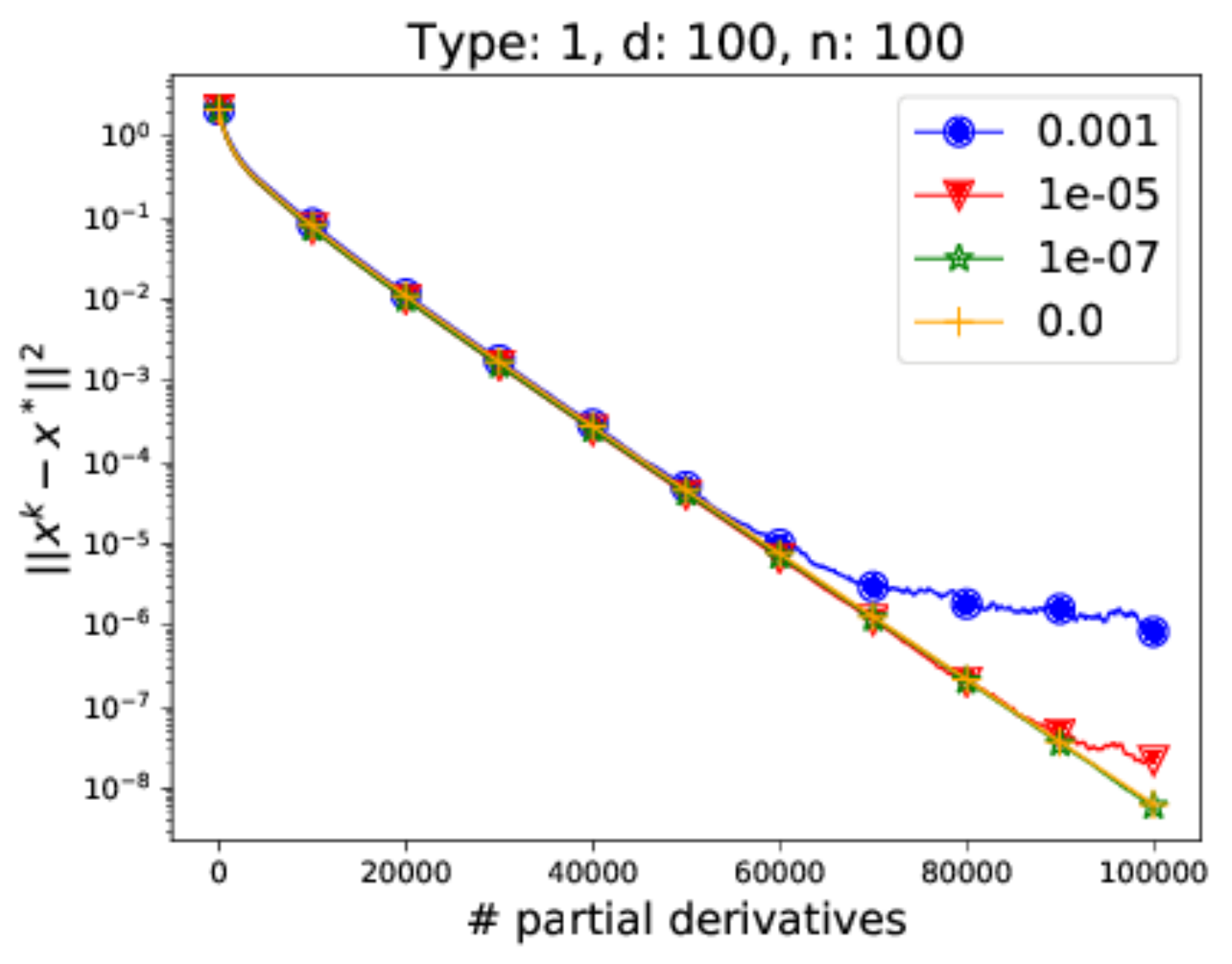}
\end{minipage}%
\begin{minipage}{0.24\textwidth}
  \centering
\includegraphics[width =  \textwidth ]{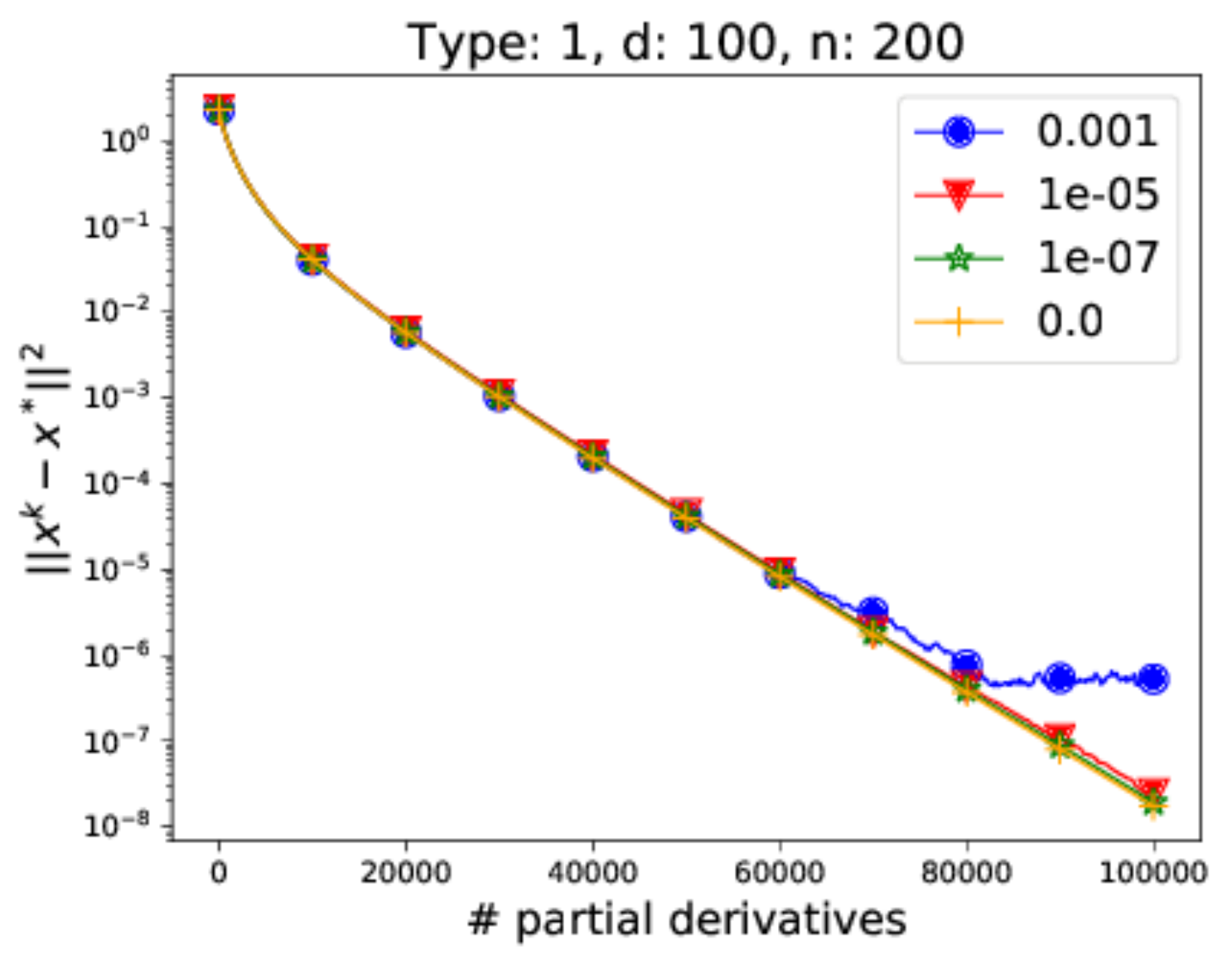}
\end{minipage}%
\begin{minipage}{0.24\textwidth}
  \centering
\includegraphics[width =  \textwidth ]{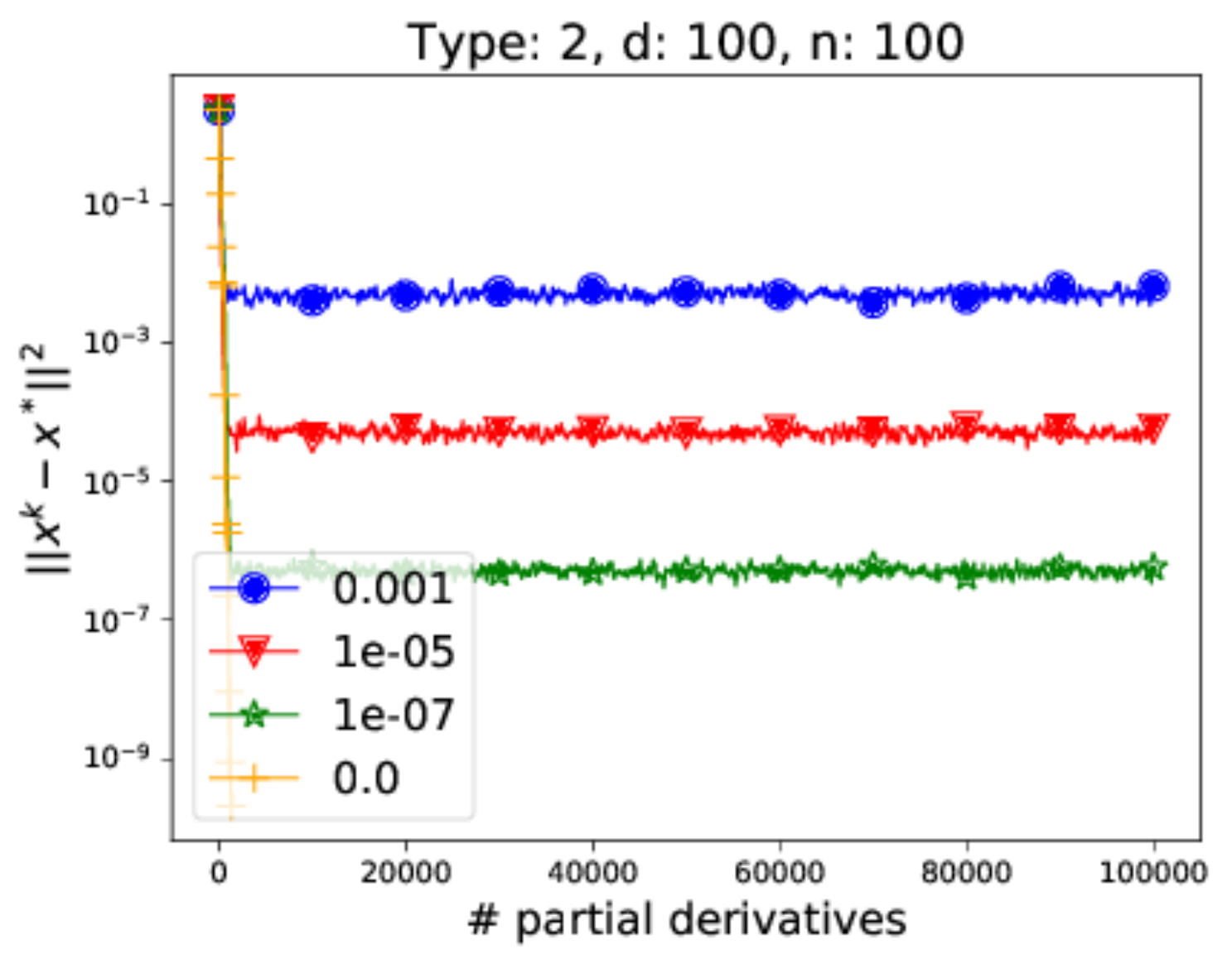}
\end{minipage}%
\begin{minipage}{0.24\textwidth}
  \centering
\includegraphics[width =  \textwidth ]{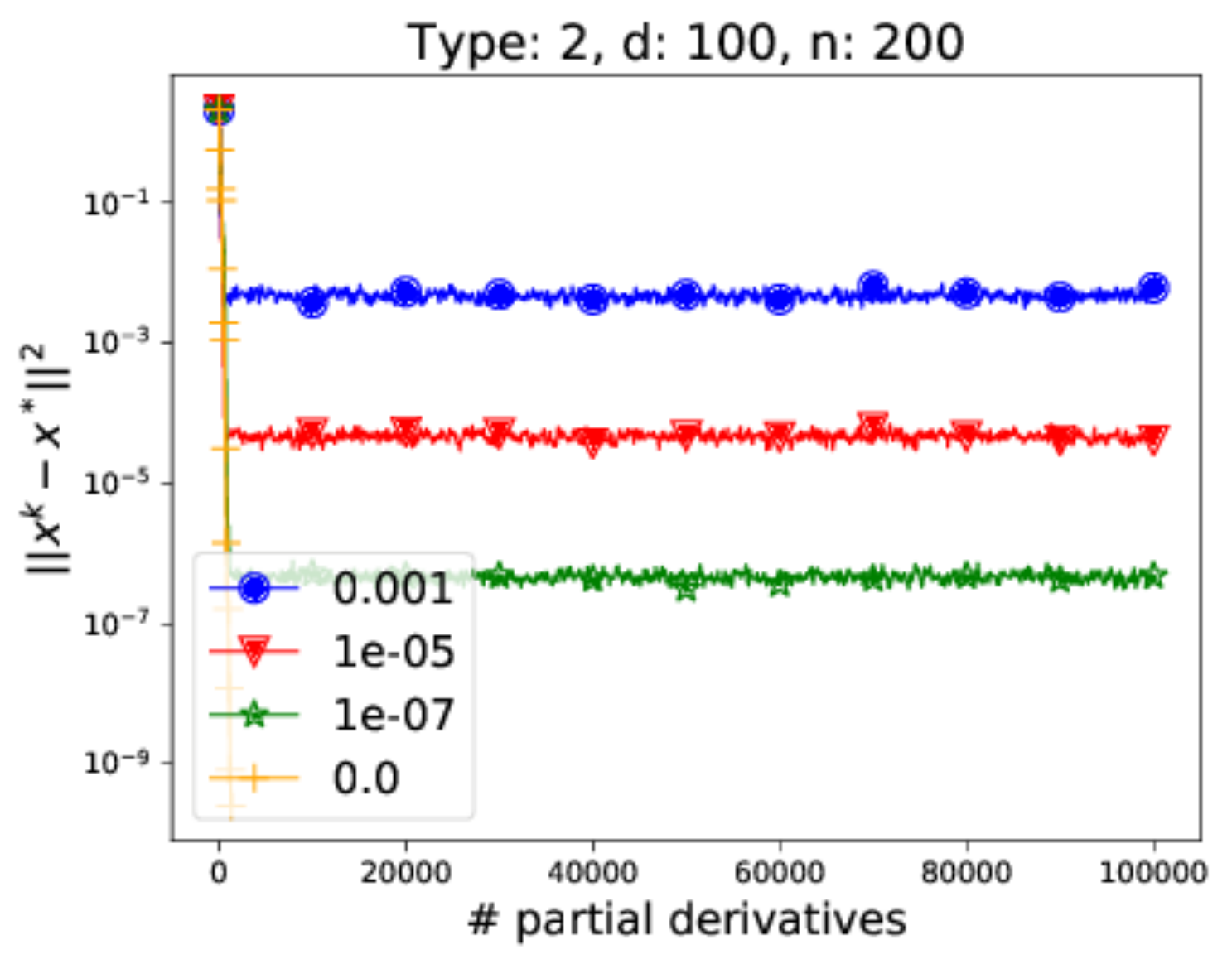}
\end{minipage}%
\\
\begin{minipage}{0.24\textwidth}
  \centering
\includegraphics[width =  \textwidth ]{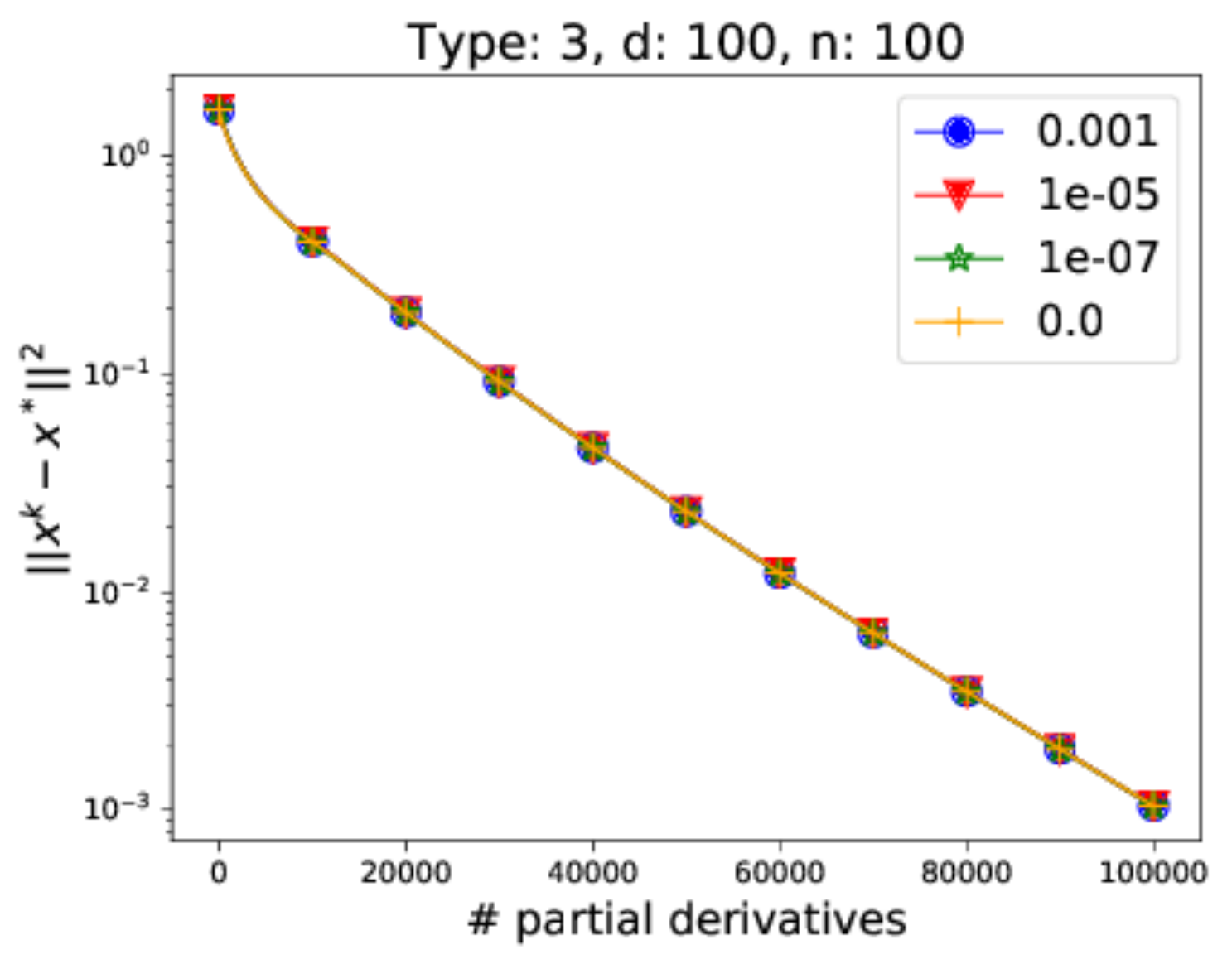}
\end{minipage}%
\begin{minipage}{0.24\textwidth}
  \centering
\includegraphics[width =  \textwidth ]{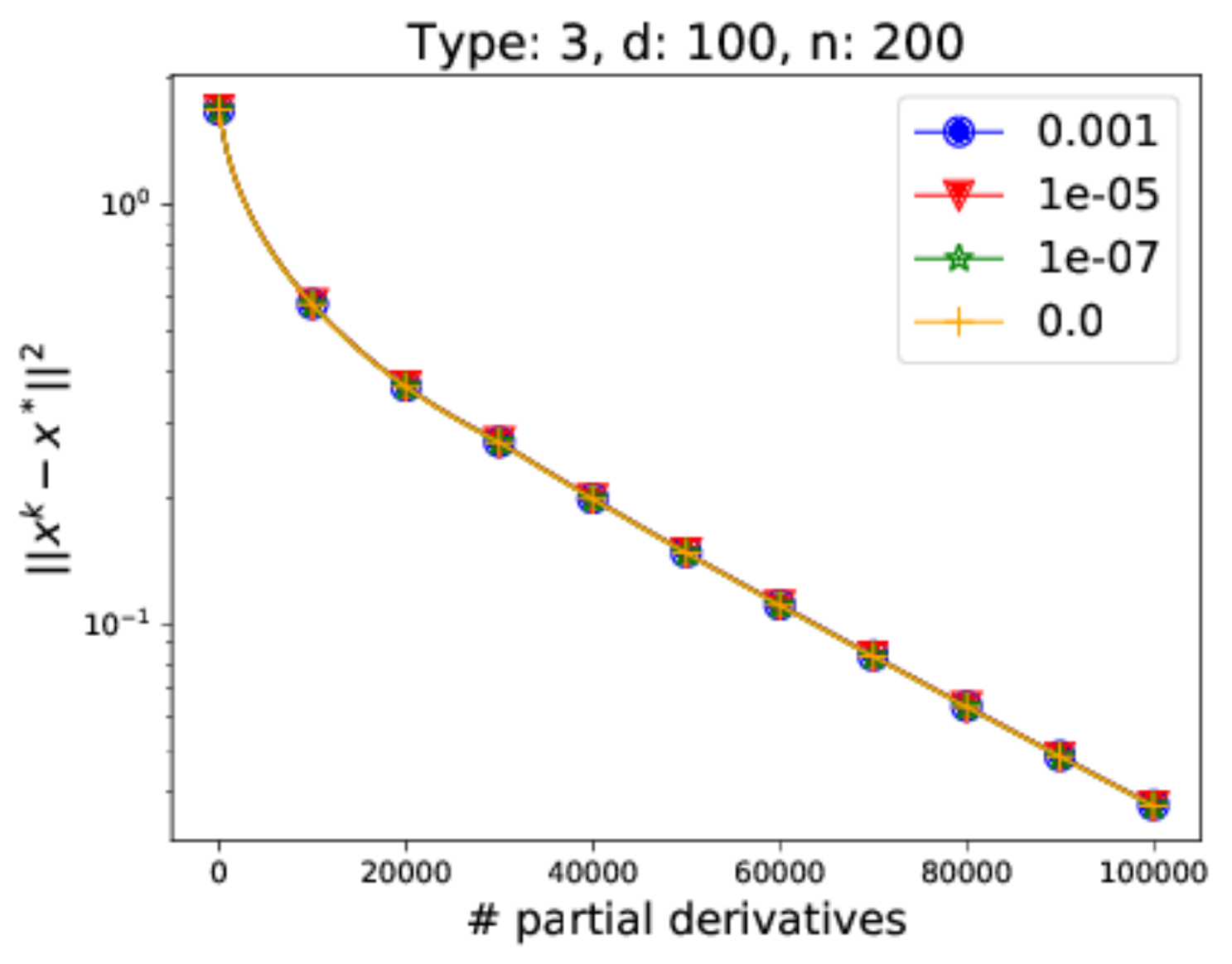}
\end{minipage}%
\begin{minipage}{0.24\textwidth}
  \centering
\includegraphics[width =  \textwidth ]{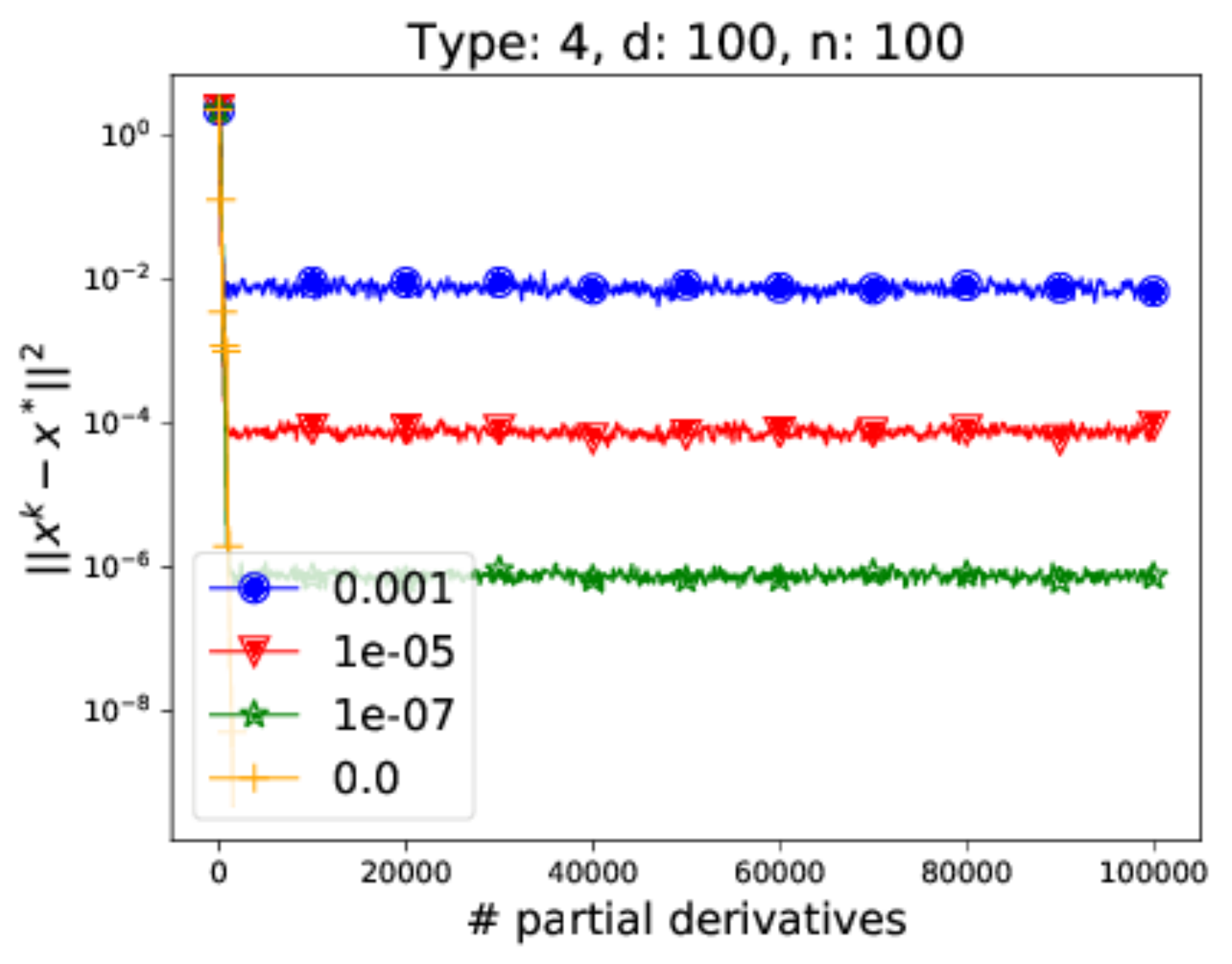}
\end{minipage}%
\begin{minipage}{0.24\textwidth}
  \centering
\includegraphics[width =  \textwidth ]{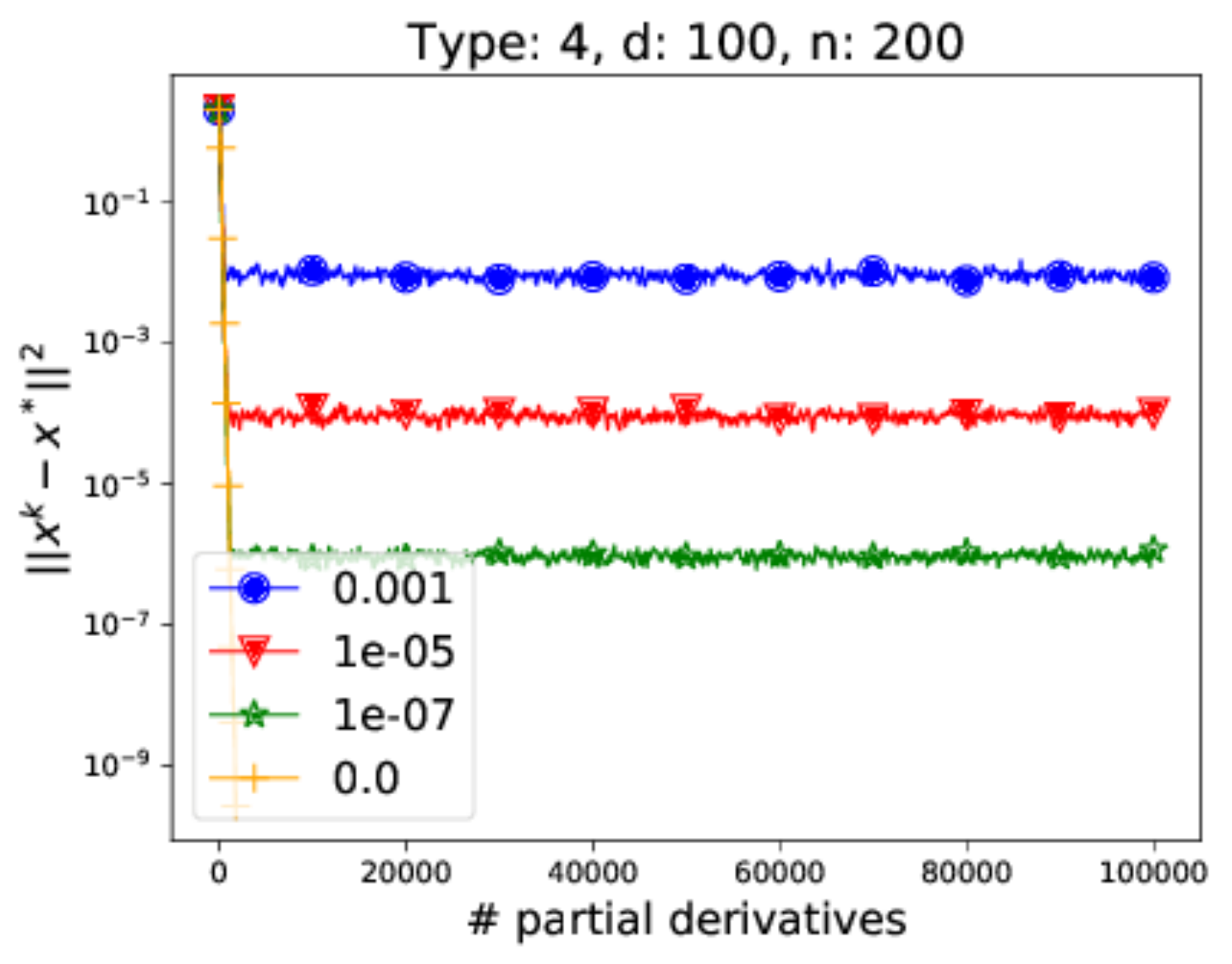}
\end{minipage}%
\caption{ {\tt N-SEGA} applied on constrained least squares problem with noised partial derivative oracle. Legend labels stand for the magnitude $\sigma^2$ of the oracle noise. }
\label{fig:nsega}
\end{figure}

The specific problem instances are presented in Table~\ref{tbl:leastsquares}. Figure~\ref{fig:nsega} shows the results. 

\begin{table}[H]
\begin{center}
\begin{tabular}{|c|c|c|}
\hline
Type & $\mA $ & $b$ \\
 \hline
 \hline
 1   & $\mA_{ij}\sim N(0,1)$ (independently)  & vector of ones  \\
 \hline
  2   & Same as 1, but scaled so that $\lambda_{\max}(A^\top A)=1$ & vector of ones \\
\hline
 3   & $\mA_{ij} = \varrho_{ij}\varpi_{j}$ $\forall i,j:\varrho_{ij},\varpi_{j} \sim N(0,1)$ (independently)  & vector of ones  \\
 \hline
  4   & Same as 3, but scaled so that $\lambda_{\max}(A^\top A)=1$ & vector of ones \\
\hline
\end{tabular}
\end{center}
\caption{Four types of least squares. }
\label{tbl:leastsquares}
\end{table}

We shall mention that this experiment serves to support and give a better intuition about the results from Section~\ref{N-SEGA} and is by no means practical. The results show, as predicted by theory, linear convergence to a specific neighborhood of the objective. The effect of the noise varies, however, as a general rule, the larger strong convexity $\mu$ is (i.e. problems 1,3 where scaling was not applied), the smaller the effect of noise is.

\section{Discussion}

Although our approach is rather general, we still see several possible directions for  future extensions, including: 

$\bullet$ Recently, our results were extended to generally convex \cite{khaled2020unified} and non-convex functions \cite{khaled2020tighter, li2020unified}. 

$\bullet$ It would be further interesting to unify our theory with {\em biased} gradient estimators. If this was possible, one could recover methods as {\tt SAG}~\citep{SAG} in special cases, or obtain rates for the zero-order optimization. We have some preliminary results in this direction already.

$\bullet$ Although our theory allows for non-uniform stochasticity, it does not recover the best known rates for {\tt RCD} type methods with {\em importance sampling}. It would be thus interesting to provide a more refined  analysis capable of capturing importance sampling phenomena more accurately.

$\bullet$ An extension of Assumption~\ref{as:general_stoch_gradient} to {\em iteration dependent} parameters $A,B,C,D_1, D_2, \rho$ would enable an array of new methods, such as {\tt SGD} with decreasing stepsizes. Such an extension is rather very straightforward. 

$\bullet$ It would be interesting to provide a unified analysis of stochastic methods with {\em acceleration} and {\em momentum}. In fact,~\citep{kulunchakov2019estimate} provide (separately) a unification of some methods with and without variance reduction. Hence, an attempt to combine our insights with their approach seems to be a promising starting point in these efforts.

%% file: ch3_ef_sigma_k.tex
\chapter{Linearly Converging Error Compensated SGD}\label{ch:ef_sigma_k}
\section{Introduction}\label{sec:intro}
We\footnote{Part of the work was done while I was a research intern at KAUST.} consider distributed optimization problems of the form
\begin{equation}
\mytextstyle
\min\limits_{x\in\R^d} \left\{ f(x) = \frac{1}{n}\sum\limits_{i=1}^n f_i(x) \right\}, \label{eq:main_problem_ef}
\end{equation}
where  $n$ is the number of  workers/devices/clients/nodes. The  information about function $f_i$ is stored on the $i$-th worker only.  Problems of this form appear in the distributed or federated training of supervised machine learning models \cite{DANE, FEDLEARN}.  In such applications, $x\in \R^d$ describes the parameters identifying a statistical model we wish to train, and $f_i$ is the (generalization or empirical) loss of model $x$ on the data accessible by worker $i$. If worker $i$ has access to data with distribution $\cD_i$, then $f_i$ is assumed to have the structure
\begin{equation}
	f_i(x) = \EE_{\xi_i\sim \cD_i}\left[f_{\xi_i}(x)\right]. \label{eq:f_i_expectation_ef}
\end{equation}
Dataset $\cD_i$ may or may not be available to worker $i$ in its entirety. Typically, we assume that worker $i$ has  only access to samples from $\cD_i$. If the dataset is fully available, it is typically finite, in which case we can assume that $f_i$ has the finite-sum form\footnote{The implicit assumption that each worker contains exactly $m$ data points is for simplicity only; all our results have direct analogues in the general setting with $m_i$ data points on worker $i$.}:
\begin{equation}
	\mytextstyle f_i(x) = \frac{1}{m}\sum\limits_{j=1}^m f_{ij}(x). \label{eq:f_i_sum_ef}
\end{equation}



\textbf{Communication bottleneck.}
  The key bottleneck in practical distributed \cite{goyal2017accurate} and federated \cite{FEDLEARN,kairouz2019advances} systems comes from the high cost of communication of messages among the clients needed to find a solution of sufficient qualities. Several approaches to addressing this communication bottleneck have been proposed in the literature.
  
In the very rare situation when it is possible to adjust the network architecture connecting the clients, one may consider a fully decentralized setup \cite{bertsekas1989parallel}, and allow each client in each iteration to  communicate to their neighbors only. One can argue that in some circumstances and in a certain sense,  decentralized architecture may be preferable to centralized architectures \cite{lian2017can}.  Another natural way to address the communication bottleneck is to do more meaningful (which typically means more expensive) work on each client before each communication round. This is done in the hope that such extra work will produce more valuable messages to be communicated, which hopefully results in the need for fewer communication rounds. A popular technique of this type which is particularly relevant to Federated Learning is based in applying multiple {\em local updates} instead of a single update only.  This is the main idea behind  {\tt Local-SGD} \cite{Stich18local}; see also \cite{basu2019qsparse, haddadpour2019convergence, karimireddy2020scaffold, khaled2020tighter, koloskova2020unified, stich2020error, woodworth2020local}.  However, in this chapter, we contribute to the line work which aims to  resolve the communication bottleneck issue via {\em communication compression}. That is, the  information that is normally exchanged---be it iterates, gradients or some more sophisticated vectors/tensors---is compressed in a lossy manner before communication. By applying compression,  fewer bits are transmitted  in each communication round, and one hopes that the increase in the number of communication rounds necessary to solve the problem, if any, is compensated by the savings, leading to a more efficient method overall.

\textbf{Error-feedback framework.} In order to address these issues, in this chapter we study a broad class of distributed stochastic first order methods for solving problem \eqref{eq:main_problem_ef} described by the  iterative framework
\begin{eqnarray}
	x^{k+1} &=& \mytextstyle x^k - \frac{1}{n} \sum \limits_{i=1}^n v_i^k, \label{eq:x^k+1_update}\\
	e_i^{k+1} &=&  e_i^k +  \gamma g_i^k -  v_i^k , \qquad i=1,2,\dots , n.\label{eq:error_update}
\end{eqnarray}
In this scheme,  $x^k$ represents the key  iterate, $v_i^k$ is the contribution of worker $i$ towards the update in iteration $k$, $g_i^k$ is an unbiased estimator of $\nabla f_i(x^k)$ computed by worker $i$, $\gamma>0$ is a fixed stepsize and $e_i^k$ is the error accumulated at node $i$ prior to iteration $k$ (we set to $e_i^0
= 0$ for all $i$). In order to better understand the role of the vectors $v_i^k$ and $e_i^k$, first consider the simple special case with $v_i^k \equiv \gamma g_i^k$. In this case, $e_i^k=0$ for all $i$ and $k$, and method \eqref{eq:x^k+1_update}--\eqref{eq:error_update} reduces to distributed {\tt SGD}:
\begin{equation}\label{eq:SGD-ss}\mytextstyle x^{k+1} = x^k - \frac{ \gamma}{n}\sum \limits_{i=1}^n g_i^k. \end{equation}
However, by allowing to chose the vectors $v_i^k$ in a different manner,  we obtain a more general update rule  than what the {\tt SGD} update \eqref{eq:SGD-ss} can offer. \cite{stich2020error},  who studied  \eqref{eq:x^k+1_update}--\eqref{eq:error_update} in the $n=1$ regime, show that this  flexibility allows to capture several types of methods,  including those employing  i) compressed communication, ii)  delayed gradients, and iii) minibatch gradient updates. While our general results apply to all these special cases and more, in order to not dilute the focus of the chapter,  in the main body of this chapter we concentrate on the first use case---compressed communication---which we now describe.

\textbf{Error-compensated compressed gradient methods.}
Note that in distributed {\tt SGD} \eqref{eq:SGD-ss}, each worker needs to know the aggregate gradient $g^k = \frac{1}{n}\sum_{i=1}^n g_i^k$ to form $x^{k+1}$, which is needed before the next iteration can start. This can be achieved, for example, by each worker $i$ communicating their gradient $g_i^k$ to all other workers. Alternatively, in a parameter server setup, a dedicated master node  collects the gradients from all workers, and broadcasts their average $g^k$  to all workers.  Instead of communicating the gradient vectors $g_i^k$, which is  expensive in distributed learning in general and in federated learning in particular, and especially if $d$ is large, we wish to communicate other but  closely related vectors which can be represented with fewer bits. To this effect, each worker $i$  sends  the vector 
\begin{equation}v_i^k = \cC(e_i^k +  \gamma g_i^k), \qquad \forall i\in [n]\label{eq:bu98gf}\end{equation} instead, where $\cC:\R^d\to \R^d$ is a (possibly randomized, and in such a case, drawn independently of all else in iteration $k$)  compression operator used  to reduce communication. We assume throughout that there exists $\delta \in (0,1]$ such  that  the following inequality holds for all $x\in\R^d$
	\begin{equation}
		\EE\left[\|\cC(x) - x\|^2\right] \le (1 - \delta)\|x\|^2. \label{eq:compression_def}
	\end{equation}

For any $k \geq 0$, the vector $e_i^{k+1} = \sum_{t=0}^{k} \gamma g_i^t - v_i^t $ captures the {\em error} accumulated by worker $i$ up to iteration $k$. This is the difference between the ideal {\tt SGD} update $\sum_{t=0}^{k} \gamma g_i^t$ and the applied update $ \sum_{t=0}^{k} v_i^t $. As we see in \eqref{eq:bu98gf},  at iteration $k$ the current error $e_i^k$ is added to the gradient update $\gamma g_i^k$---this is referred to as {\em error feedback}---and subsequently compressed, which defines the update vector $v_i^k$. Compression introduces additional error, which is added to $e_i^k$, and the process is repeated.


\textbf{Compression operators.} For a rich collection of specific  operators satisfying \eqref{eq:compression_def}, we refer the reader to \cite{stich2020error} and \cite{beznosikov2020biased}. These include various unbiased or contractive sparsification operators such as RandK and TopK, respectively,  and quantization operators  such as natural compression  and natural dithering \cite{Cnat}. Several additional comments related to compression operators are included in Section~\ref{sec:compressions}.


%

\section{Summary of Contributions}\label{sec:contrib}

We now summarize the key contributions of this chapter.

\textbf{$\diamond$ General theoretical framework.} 
In this work we propose a  {\em general theoretical framework} for analyzing a wide class of methods that can be written in the the error-feedback form \eqref{eq:x^k+1_update}-\eqref{eq:error_update}. We perform {\em complexity analysis}  under $\mu$-strong quasi convexity (Assumption~\ref{ass:quasi_strong_convexity})  and $L$-smoothness (Assumption~\ref{ass:L_smoothness})  assumptions on the functions $f$ and $\{f_i\}$, respectively. Our analysis is based on an additional {\em parametric assumption}  (using parameters $A$, $A'$, $B_1$, $B_1'$, $B_2$, $B_2'$, $C_1$, $C_2$, $D_1$, $D_1'$, $D_2$, $D_3$, $\eta$, $\rho_1$, $\rho_2$, $F_1$, $F_2$, $G$) on the relationship between the iterates $x^k$,  stochastic gradients $g^k$, errors $e^k$ and a few other quantities (see Assumption~\ref{ass:key_assumption_new}, and the stronger Assumption~\ref{ass:key_assumption_finite_sums_new}). We prove a single theorem (Theorem~\ref{thm:main_result_new}) from which all our complexity results follow as special cases. That is, for each existing or new specific method, we {\em prove} that one (or both) of our parametric assumptions holds, and specify the parameters for which it holds. These parameters have direct impact on the theoretical rate of the method. A summary of the values of the parameters for all methods developed in this chapter is provided in Table~\ref{tbl:special_cases-parameters} in the appendix. We remark that the values of the parameters $A, A', B_1, B_1', B_2, B_2', C_1, C_2$ and $\rho_1, \rho_2$ influence the theoretical stepsize. 

\textbf{$\diamond$ Sharp rates.} For  existing methods covered by our general framework, our main convergence result (Theorem~\ref{thm:main_result_new}) recovers the best known rates for these methods up to constant factors. 

\textbf{$\diamond$ Eight new error-compensated (EC) methods.} We study several specific EC methods for solving problem \eqref{eq:main_problem_ef}. First, we recover  the {\tt EC-SGD} method first analyzed in the $n=1$ case by \cite{stich2020error} and later in the general $n\geq 1$ case by \cite{beznosikov2020biased}. More importantly, we develop {\em eight new methods}: {\tt EC-SGDsr}, {\tt EC-GDstar}, {\tt EC-SGD-DIANA}\footnote{Inspired by personal communication with D.~Kovalev in November 2019 who shared a key algorithm and preliminary results of our work, \cite{stich2020communication} studied almost the same algorithm and also other related methods and independently derived convergence rates. Our work was finalized and submitted to NeurIPS 2020 in June 2020, while the results in \cite{stich2020communication} were obtained in Summer 2020 and appeared on arXiv in September 2020.  Moreover, in our work, we obtain tighter rates (see Table~\ref{tbl:special_cases2_ef} for the details). }, {\tt EC-SGDsr-DIANA}, {\tt EC-GD-DIANA}, {\tt EC-LSVRG}, {\tt EC-LSVRGstar} and {\tt EC-LSVRG-DIANA}. 
Some of these methods are designed to work with the expectation  structure of the local functions $f_i$ given in \eqref{eq:f_i_expectation_ef}, and others are specifically designed
 to exploit the  finite-sum structure \eqref{eq:f_i_sum_ef}. All these methods follow the error-feedback framework \eqref{eq:x^k+1_update}--\eqref{eq:error_update}, with $v_i^k$ chosen as in \eqref{eq:bu98gf}. They differ in how the gradient estimator $g_i^k$ is {\em constructed} (see Table~\ref{tbl:EC_methods_summary} for a compact description of all these methods;  formal descriptions can be found in the appendix). As we shall see, the existing {\tt EC-SGD} method uses a rather  naive gradient estimator, which renders it less efficient in theory and practice when compared to the best of our new methods. A key property of our parametric assumption described above is that it allows for the construction and modeling of more elaborate gradient estimators, which leads to new EC methods with superior theoretical and practical properties when compared to prior state of the art.

\begin{table}[H]
\caption{Complexity of Error-Compensated {\tt SGD} methods established in this chapter. Symbols: $\varepsilon = $ error tolerance; $\delta = $ contraction factor of compressor $\cC$; $\omega = $ variance parameter of compressor $\cQ$; $\kappa = \nicefrac{L}{\mu}$; $\cL =$ expected smoothness constant; $\sigma_*^2 = $ variance of the stochastic gradients in the solution; $\zeta_*^2 =$ average of $\|\nabla f_i(x^*)\|^2$; $\sigma^2 =$ average of the uniform bounds for the variances of stochastic gradients of workers. {\tt EC-GDstar}, {\tt EC-LSVRGstar} and {\tt EC-LSVRG-DIANA} are the first EC methods with a linear convergence rate without assuming that $\nabla f_i(x^*)=0$ for all $i$. {\tt EC-LSVRGstar} and {\tt EC-LSVRG-DIANA} are the first EC methods with a linear convergence rate which do not require the computation of the full gradient $\nabla f_i(x^k)$ by all workers in each iteration. Out of these three methods, only {\tt EC-LSVRG-DIANA} is practical. $^\dagger${\tt EC-GD-DIANA} is a special case of {\tt EC-SGD-DIANA} where each worker $i$ computes the full gradient  $\nabla f_i(x^k)$.
}
\label{tbl:special_cases2_ef}
\begin{center}
\scriptsize
\begin{tabular}{|c|l|c|c|c|c|}
\hline
\bf Problem & \bf Method &   \bf Alg \# &  \bf Citation &  \bf  Sec \#  
& \bf Rate (constants ignored)\\
\hline
\eqref{eq:main_problem_ef}+\eqref{eq:f_i_sum_ef} & {\tt EC-SGDsr}  & Alg \ref{alg:ec-SGDsr} & {\color{red}\bf new} & \ref{sec:ec_SGDsr} 
& {\color{red}$\widetilde{\cO}\left(\frac{\cL}{\mu} + \frac{L+\sqrt{\delta L\cL}}{\delta\mu} + \frac{\sigma_*^2}{n\mu\varepsilon} + \frac{\sqrt{L(\sigma_*^2 + \nicefrac{\zeta_*^2}{\delta})}}{\mu\sqrt{\delta\varepsilon}}\right)$}\\
\hline
\eqref{eq:main_problem_ef}+\eqref{eq:f_i_expectation_ef} & {\tt EC-SGD}  & Alg \ref{alg:ec-sgd} & {\cite{stich2020error}}  & \ref{sec:ec_sgd_pure} 
& $\widetilde{\cO}\left(\frac{\kappa}{\delta} + \frac{\sigma_*^2}{n\mu\varepsilon} + \frac{\sqrt{L(\sigma_*^2 + \nicefrac{\zeta_*^2}{\delta})}}{\delta\mu\sqrt{\varepsilon}}\right)$\\
\hline
\eqref{eq:main_problem_ef}+\eqref{eq:f_i_sum_ef} & {\tt EC-GDstar}  & Alg \ref{alg:EC-GDstar} & {\color{red}\bf new}
& \ref{sec:ec_SGDstar} 
& {\color{red} $\cO\left( \frac{\kappa}{\delta}  \log \frac{1}{\varepsilon} \right)$ } \\
\hline
\eqref{eq:main_problem_ef}+\eqref{eq:f_i_expectation_ef} & {\tt EC-SGD-DIANA}  & Alg \ref{alg:EC-SGD-DIANA} & {\color{red}\bf new}
& \ref{sec:ec_diana} 
& \begin{tabular}{c}
	{\color{red}Opt.\ I: $\widetilde{\cO}\left(\omega + \frac{\kappa}{\delta} + \frac{\sigma^2}{n\mu\varepsilon} + \frac{\sqrt{L\sigma^2}}{\delta\mu\sqrt{\varepsilon}}\right)$}\\
	{\color{red}Opt.\ II: $\widetilde{\cO}\left(\frac{1+\omega}{\delta} + \frac{\kappa}{\delta} + \frac{\sigma^2}{n\mu\varepsilon} + \frac{\sqrt{L\sigma^2}}{\mu\sqrt{\delta\varepsilon}}\right)$}
\end{tabular}\\
\hline
\eqref{eq:main_problem_ef}+\eqref{eq:f_i_sum_ef} & {\tt EC-SGDsr-DIANA}  & Alg \ref{alg:EC-SGDsr-DIANA} & {\color{red}\bf new} & \ref{sec:ec_sgdsr_DIANA} 
& \begin{tabular}{c}
	{\color{red} Opt.\ I: $\widetilde{\cO}\left(\omega + \frac{\cL}{\mu} + \frac{\sqrt{L\cL}}{\delta\mu} + \frac{\sigma_*^2}{n\mu\varepsilon} + \frac{\sqrt{L\sigma_*^2}}{\delta\mu\sqrt{\varepsilon}}\right)$}\\
	{\color{red} Opt.\ II: $\widetilde{\cO}\left(\frac{1+\omega}{\delta} + \frac{\cL}{\mu} + \frac{\sqrt{L\cL}}{\delta\mu} + \frac{\sigma_*^2}{n\mu\varepsilon} + \frac{\sqrt{L\sigma_*^2}}{\mu\sqrt{\delta\varepsilon}}\right)$}
\end{tabular}\\
\hline
\eqref{eq:main_problem_ef}+\eqref{eq:f_i_expectation_ef} & {\tt EC-GD-DIANA}$^\dagger$  & Alg \ref{alg:EC-SGD-DIANA} & {\color{red}\bf new}
& \ref{sec:ec_diana} 
& {\color{red}$\cO\left(\left(\omega + \frac{\kappa}{\delta}\right) \log \frac{1}{\varepsilon}\right)$} \\
\hline
\eqref{eq:main_problem_ef}+\eqref{eq:f_i_sum_ef} & {\tt EC-LSVRG}  & Alg \ref{alg:ec-LSVRG} & {\color{red}\bf new}
& \ref{sec:ec_LSVRG} 
& {\color{red}$\widetilde{\cO}\left(m + \frac{\kappa}{\delta} + \frac{\sqrt{L\zeta_*^2}}{\delta\mu\sqrt{\varepsilon}}\right)$}\\
\hline
\eqref{eq:main_problem_ef}+\eqref{eq:f_i_sum_ef} & {\tt EC-LSVRGstar}  & Alg \ref{alg:ec-LSVRGstar} & {\color{red}\bf new}
 & \ref{sec:ec_LSVRGstar} 
 & {\color{red}$\cO\left(\left(m + \frac{\kappa}{\delta}\right)  \log \frac{1}{\varepsilon}\right)$} \\
\hline
\eqref{eq:main_problem_ef}+\eqref{eq:f_i_sum_ef} & {\tt EC-LSVRG-DIANA}  & Alg \ref{alg:ec-LSVRG-diana} & {\color{red}\bf new}
& \ref{sec:ec_LSVRG-diana} 
& {\color{red}$\cO \left(\left( \omega + m + \frac{\kappa }{\delta} \right) \log \frac{1}{\varepsilon}\right)$} \\
\hline
\end{tabular}
\end{center}
\end{table}

\begin{table}[H]\caption{Error compensated methods developed in this work. In all cases, $v_i^k = \cC(e_i^k + \gamma g_i^k)$. The full descriptions of the algorithms are included in Section~\ref{sec:special_cases_ef}.} \label{tbl:EC_methods_summary}
\begin{center}
\footnotesize
\begin{tabular}{|c|c|c|c|c|} 
\hline
Problem & Method & $g_i^k $ & Comment \\
\hline
\eqref{eq:main_problem_ef} + \eqref{eq:f_i_sum_ef}       &   {\tt EC-SGDsr}           & $\frac{1}{m} 
\sum \limits_{j=1}^m \xi_{ij}\nabla f_{ij}(x^k)$ & \begin{tabular}{c}$\Exp\left[ \xi_{ij} \right] = 1 $ \\  $\EE_{\cD_i}\left[\|\nabla f_{\xi_i}(x) - \nabla f_{\xi_i}(x^*)\|^2\right]$ \\
	$~~~~~\le 2\cL D_{f_i}(x,x^*)$
\end{tabular}  \\  
\hline
\eqref{eq:main_problem_ef} + \eqref{eq:f_i_expectation_ef}       &   {\tt EC-SGD}           & $\nabla f_{\xi_i}(x^k)$ &  \\  
\hline

\eqref{eq:main_problem_ef}             & {\tt EC-GDstar}           & $\nabla f_i(x^k) - \nabla f_i(x^*)$ & known $\nabla f_i(x^*) \; \forall i$  \\  
\hline
\eqref{eq:main_problem_ef} + \eqref{eq:f_i_expectation_ef}              & {\tt EC-SGD-DIANA}   & $\hat{g}_i^k - h_i^k + h^k$ & \begin{tabular}{c}$\EE\left[\hat{g}_i^k\right] = \nabla f_i(x^k)$ \\
 			$\Exp_k \left[ \|\hat{g}_i^k - \nabla f_i(x^k)\|^2\right] \leq D_{1,i}$\\
             $h_i^{k+1} = h_i^k + \alpha \cQ(\hat{g}_i^k - h_i^k)$\\
             $h^k = \frac{1}{n} \sum \limits_{i=1}^n h_i^k$ \end{tabular}  \\
\hline
\eqref{eq:main_problem_ef} +  \eqref{eq:f_i_sum_ef}               & {\tt EC-SGDsr-DIANA}   & $\nabla f_{\xi_i^k}(x^k) - h_i^k + h^k$ & \begin{tabular}{c}$\EE\left[\nabla f_{\xi_i^k}(x^k)\right] = \nabla f_i(x^k)$ \\
 			$\EE_{\cD_i}\left[\|\nabla f_{\xi_i}(x) - \nabla f_{\xi_i}(x^*)\|^2\right]$ \\
	$~~~~~\le 2\cL D_{f_i}(x,x^*)$\\
             $h_i^{k+1} = h_i^k + \alpha \cQ(\nabla f_{\xi_i^k}(x^k) - h_i^k)$\\
             $h^k = \frac{1}{n} \sum \limits_{i=1}^n h_i^k$ \end{tabular}  \\
\hline
\eqref{eq:main_problem_ef} +  \eqref{eq:f_i_sum_ef}         &   {\tt EC-LSVRG}           & {$ \begin{aligned} \nabla & f_{il}(x^k) - \nabla f_{il}(w_i^k) \\ &+ \nabla f_i(w_i^k) \end{aligned} $}  & \begin{tabular}{c}$l$ chosen uniformly from $[m]$\\ $w_i^{k+1} = \begin{cases}x^k,& \text{with prob. } p,\\ w_i^k,& \text{with prob. } 1-p\end{cases}$\end{tabular}  \\  
\hline
\eqref{eq:main_problem_ef} +  \eqref{eq:f_i_sum_ef}         &   {\tt EC-LSVRGstar}           & { $ \begin{aligned}\nabla & f_{il}(x^k) - \nabla f_{il}(w_i^k) \\ & + \nabla f_i(w_i^k) - \nabla f_i(x^*) \end{aligned}$ }  & \begin{tabular}{c}$l$ chosen uniformly from $[m]$\\ $w_i^{k+1} = \begin{cases}x^k,& \text{with prob. } p,\\ w_i^k,& \text{with prob. } 1-p\end{cases}$\end{tabular}  \\  
\hline
\eqref{eq:main_problem_ef} +  \eqref{eq:f_i_sum_ef}         &   {\tt EC-LSVRG-DIANA}           & \begin{tabular}{c} $\hat{g}_i^k - h_i^k + h^k$ \\ where \\ { $ \begin{aligned}& \hat{g}_i^k =  \nabla  f_{il}(x^k) \\ &- \nabla f_{il}(w_i^k)  + \nabla f_i(w_i^k)  \end{aligned}$ } \end{tabular} & \begin{tabular}{c}  $h_i^{k+1} = h_i^k + \alpha \cQ(\hat{g}_i^k - h_i^k)$\\
             $h^k = \frac{1}{n} \sum \limits_{i=1}^n h_i^k$ \\ $l$ chosen uniformly from $[m]$\\ $w_i^{k+1} = \begin{cases}x^k,& \text{with prob. } p,\\ w_i^k,& \text{with prob. } 1-p\end{cases}$\end{tabular} \\  
\hline
\end{tabular}
\end{center}
\end{table}

\textbf{$\diamond$ First linearly converging EC methods.}  The key theoretical consequence of our general framework is the development of the {\em first linearly converging} error-compensated {\tt SGD}-type methods for distributed training with biased communication compression. In particular, we design four such methods: two simple  but impractical methods, {\tt EC-GDstar} and {\tt EC-LSVRGstar}, with rates $\cO\left(\frac{\kappa}{\delta}\ln\frac{1}{\varepsilon}\right)$ and $\cO\left(\left(m+\frac{\kappa}{\delta}\right)\ln\frac{1}{\varepsilon}\right)$, respectively,
and two  practical but more elaborate methods, {\tt EC-GD-DIANA}, with rate  $\cO\left(\left(\omega + \frac{\kappa }{\delta}\right) \ln \frac{1}{\varepsilon}\right)$, and  {\tt EC-LSVRG-DIANA}, with rate $\cO\left(\left(\omega + m + \frac{\kappa }{\delta}\right) \ln \frac{1}{\varepsilon}\right).$
In these rates, $\kappa=\nicefrac{L}{\mu}$ is the condition number,  $0<\delta \leq 1$ is the contraction parameter associated with the compressor $\cC$ used in \eqref{eq:bu98gf}, and $\omega$ is the variance parameter associated with a {\em secondary unbiased compressor\footnote{We assume that $\Exp{\cQ(x)}=x$ and $\Exp{\|\cQ(x)-x\|^2} \leq \omega \|x\|^2$ for all $x\in \R^d$.} $\cQ$} which plays a key role in the construction of the gradient estimator $g_i^k$. The complexity of the  first and third methods does not depend on $m$ as they require the computation of the full gradient $\nabla f_i(x^k)$ for each $i$. The remaining two methods only need to compute $\cO(1)$ stochastic gradients $\nabla f_{ij}(x^k)$ on each worker $i$.  

The first two methods, while impractical, provided us with the intuition which enabled us to develop the practical variant. We include them in this chapter due to their simplicity, because of the added insights they offer, and to showcase the flexibility of our general theoretical framework, which is able to describe them. {\tt EC-GDstar} and {\tt EC-LSVRGstar} are impractical since they require the knowledge of the gradients $\{\nabla f_i(x^*)\}$, where $x^*$ is an optimal solution of \eqref{eq:main_problem_ef}, which are obviously not known since $x^*$ is not known. 

The only known linear convergence result for an error compensated {\tt SGD} method is due to \cite{beznosikov2020biased}, who require the computation of the full gradient of $f_i$ by each  machine $i$ (i.e., $m$ stochastic gradients), and the additional assumption that $\nabla f_i(x^*) = 0$ for all $i$. We do not need such assumptions, thereby resolving a major theoretical issue with EC methods.

\textbf{$\diamond$ Results in the convex case.} 
Our theoretical analysis goes beyond distributed optimization and recovers the results from \cite{gorbunov2019unified,khaled2020unified} (without regularization) in the  special case when $v_i^k \equiv \gamma g_i^k$. As we have seen, in this case $e_i^k\equiv 0$ for all $i$ and $k$, and  the error-feedback framework \eqref{eq:x^k+1_update}--\eqref{eq:error_update}  reduces to distributed {\tt SGD} \eqref{eq:SGD-ss}. In this regime, the relation \eqref{eq:sum_of_errors_bound_new} in Assumption~\ref{ass:key_assumption_new} becomes void, while relations \eqref{eq:second_moment_bound_new} and \eqref{eq:sigma_k+1_bound_1} with $\sigma_{2,k}^2\equiv 0$ are precisely those used by \cite{gorbunov2019unified} to analyze a wide array of {\tt SGD} methods, including  vanilla {\tt SGD} \cite{RobbinsMonro:1951}, {\tt SGD} with arbitrary sampling \cite{gower2019sgd}, as well as variance reduced methods such as {\tt SAGA} \cite{SAGA}, {\tt SVRG} \cite{SVRG}, {\tt LSVRG} \cite{hofmann2015variance, kovalev2019don}, {\tt JacSketch} \cite{gower2018stochastic},  {\tt SEGA} \cite{hanzely2018sega} and {\tt DIANA} \cite{mishchenko2019distributed, horvath2019stochastic}. Our theorem recovers the rates of all the methods just listed in both the convex case $\mu = 0$ \cite{khaled2020unified} and the strongly-convex case $\mu > 0$ \cite{gorbunov2019unified} under the more general Assumption~\ref{ass:key_assumption_new}. 


\textbf{$\diamond$ {\tt DIANA} with bi-directional quantization.}
To illustrate how our framework can be used even in the case when $v_i^k \equiv \gamma g_i^k$, $e_i^k \equiv 0$, we develop analyze a new version of {\tt DIANA} called {\tt DIANAsr-DQ} that uses arbitrary sampling on every node and double quantization\footnote{In the concurrent work (which appeared on arXiv after we have submitted our paper to NeurIPS) a similar method was independently proposed under the name of {\tt Artemis} \cite{philippenko2020artemis}. However, our analysis is  more general, see all the details on this method in the appendix. This footnote was added to the paper during the preparation of the camera-ready version of our paper.}, i.e., unbiased compression not only on the workers' side but also on the master's one.

\textbf{$\diamond$ Methods with delayed updates.}
Following \cite{stich2019unified}, we also show that our approach covers {\tt SGD} with delayed updates \cite{agarwal2011distributed, arjevani2020tight, feyzmahdavian2016asynchronous} ({\tt D-SGD}), and our analysis shows the best-known rate for this method. Due to the flexibility of our framework, we are able develop several new variants of {\tt D-SGD} with and without quantization, variance reduction, and arbitrary sampling. Again, due to space limitations, we put these methods together with their convergence analyses in the appendix.

\section{Main Result}\label{sec:main_res}
In this section we present the main theoretical result of our chapter.  First,  we introduce our assumption on $f$, which is a relaxation of $\mu$-strong convexity (see also Assumption~\ref{as:mu_strongly_quasi_convex}).
\begin{assumption}[$\mu$-strong quasi-convexity]\label{ass:quasi_strong_convexity}
	Assume that function $f$ has a unique minimizer $x^*$. We say that function $f$ is strongly quasi-convex with parameter $\mu\ge 0$ if for all $x\in\R^d$
	\begin{equation}
	\mytextstyle	f(x^*) \ge f(x) + \langle\nabla f(x), x^* - x\rangle + \frac{\mu}{2}\|x - x^*\|^2. \label{eq:str_quasi_cvx}
	\end{equation}
\end{assumption}
We allow $\mu$ to be zero, in which case  $f$ is sometimes called {\em weakly quasi-convex} (see \cite{stich2019unified} and references therein).

We now introduce our key parametric assumption on  the stochastic gradient $g^k$. This is a generalization of the assumption introduced by \cite{gorbunov2019unified} for the particular class of methods described covered by the EF framework \eqref{eq:x^k+1_update}--\eqref{eq:error_update}.

\begin{assumption}\label{ass:key_assumption_finite_sums_new}
	For all $k\ge 0$, the stochastic gradient $g^k$ is an average of stochastic gradients $g_i^k$ such that
	\begin{equation}
\mytextstyle		g^k = \frac{1}{n}\sum\limits_{i=1}^ng_i^k,\qquad \EE\left[g^k\mid x^k\right] = \nabla f(x^k). \label{eq:unbiasedness_g_i^k_new}
	\end{equation}
	Moreover,  there exist  constants $A,\widetilde{A}, A', B_1, B_2, \widetilde{B}_1, \widetilde{B}_2, B_1', B_2', C_1, C_2, G, D_1, \widetilde{D}_1, D_1', D_2, D_3 \ge 0$, and $\rho_1,\rho_2 \in [0,1]$ and two sequences of (probably random) variables $\{\sigma_{1,k}\}_{k\ge 0}$ and $\{\sigma_{2,k}\}_{k\ge 0}$, such  that the following recursions hold:
	\begin{eqnarray}
	\mytextstyle	\frac{1}{n}\sum\limits_{i=1}^n\left\|\bar{g}_i^k\right\|^2 &\le& 2A(f(x^k) - f(x^*)) + B_1\sigma_{1,k}^2 + B_2\sigma_{2,k}^2 + D_1, \label{eq:second_moment_bound_g_i^k_new}\\
\mytextstyle		\frac{1}{n}\sum\limits_{i=1}^n\EE\left[\left\|g_i^k-\bar{g}_i^k\right\|^2\mid x^k\right] &\le& 2\widetilde{A}(f(x^k) - f(x^*)) + \widetilde{B}_1\sigma_{1,k}^2 + \widetilde{B}_2\sigma_{2,k}^2 + \widetilde{D}_1, \label{eq:variance_bound_g_i^k_new}\\
		\EE\left[\|g^k\|^2\mid x^k\right] &\le& 2A'(f(x^k) - f(x^*)) + B_1'\sigma_{1,k}^2 + B_2'\sigma_{2,k}^2 + D_1', \label{eq:second_moment_bound_new}\\
		\EE\left[\sigma_{1,k+1}^2\mid \sigma_{1,k}^2, \sigma_{2,k}^2\right] &\le& (1-\rho_1)\sigma_{1,k}^2 + 2C_1\left(f(x^k) - f(x^*)\right) + G\rho_1 \sigma_{2,k}^2 + D_2,\label{eq:sigma_k+1_bound_1}\\
		\EE\left[\sigma_{2,k+1}^2\mid \sigma_{2,k}^2\right] &\le& (1-\rho_2)\sigma_{2,k}^2 + 2C_2\left(f(x^k) - f(x^*)\right),\label{eq:sigma_k+1_bound_2}
	\end{eqnarray}
	where $\bar{g}_i^k = \EE\left[g_i^k\mid x^k\right]$.
\end{assumption}

Let us briefly explain the intuition behind the assumption and the meaning of the introduced parameters. First of all, we assume that the stochastic gradient at iteration $k$ is conditionally unbiased estimator of $\nabla f(x^k)$, which is a natural and commonly used assumption on the stochastic gradient in the literature. However, we explicitly do {\em not} require unbiasedness of $g_i^k$, which is very useful in some special cases. Secondly, let us consider the simplest special case when $g^k \equiv \nabla f(x^k)$ and $f_1 = \ldots = f_n = f$, i.e., there is no stochasticity/randomness in the method and the workers have the same functions. Then due to $\nabla f(x^*) = 0$, we have that
\begin{equation*}
	\|\nabla f(x^k)\|^2 \overset{\eqref{eq:L_smoothness_cor_3}}{\le} 2L(f(x^k) - f(x^*)),
\end{equation*}
which implies that Assumption~\ref{ass:key_assumption_finite_sums_new} holds in this case with $A = A' = L$, $\widetilde{A}=0$ and $B_1 = B_2 = \widetilde{B}_1 = \widetilde{B}_2 = B_1' = B_2' = C_1 = C_2 = D_1 = \widetilde{D}_1 = D_1' = D_2 = 0$, $\rho = 1$, $\sigma_{1,k}^2 \equiv \sigma_{2,k}^2  \equiv 0$.

In general, if $g^k$ satisfies Assumption~\ref{ass:key_assumption_new}, then  parameters $A$, $\widetilde{A}$ and $A'$ are usually connected with the smoothness properties of $f$ and typically they are just multiples of $L$, whereas terms $B_1\sigma_{1,k}^2$, $B_2\sigma_{2,k}^2$, $\widetilde{B}_1\sigma_{1,k}^2$, $\widetilde{B}_2\sigma_{2,k}^2$, $B_1'\sigma_{1,k}^2$, $B_2'\sigma_{2,k}^2$ and $D_1$, $\widetilde{D}_1$, $D_1'$ appear due to the stochastic nature of $g_i^k$. Moreover, $\{\sigma_{1,k}^2\}_{k\ge 0}$ and $\{\sigma_{2,k}^2\}_{k\ge 0}$ are sequences connected with variance reduction processes and for the methods; without any kind of variance reduction these sequences contains only zeros. Parameters $B_1$ and $B_2$ are often $0$ or small positive constants, e.g., $B_1 = B_2 = 2$, and $D_1$ characterizes the remaining variance in the estimator $g^k$ that is not included in the first two terms. 

Inequalities \eqref{eq:sigma_k+1_bound_1} and \eqref{eq:sigma_k+1_bound_2} describe the variance reduction processes: one can interpret $\rho_1$ and $\rho_2$ as the {\em rates} of the variance reduction processes, $2C_1(f(x^k) - f(x^*))$ and $2C_2(f(x^k) - f(x^*))$ are ``optimization'' terms and, similarly to $D_1$, $D_2$ represents the remaining variance that is not included in the first two terms. Typically, $\sigma_{1,k}^2$ controls the variance coming from compression and $\sigma_{2,k}^2$ controls the variance taking its origin in finite-sum type randomization (i.e., subsampling) by each worker. In the case $\rho_1 = 1$ we assume that $B_1 = B_1' = C_1 = G = 0, D_2 = 0$ (for $\rho_2 = 1$ analogously), since inequality \eqref{eq:sigma_k+1_bound_1} becomes superfluous.

However, in our main result we need a slightly different assumption.
\begin{assumption}\label{ass:key_assumption_new}
	For all $k\ge 0$, the stochastic gradient $g^k$ is an unbiased estimator of $\nabla f(x^k)$:
	\begin{equation}
		\EE\left[g^k\mid x^k\right] = \nabla f(x^k). \label{eq:unbiasedness_new}
	\end{equation}
	Moreover, there exist non-negative constants $A',B_1', B_2',C_1, C_2,F_1, F_2, G, D_1',D_2, D_3 \ge 0, \rho_1, \rho_2 \in [0,1]$ and two sequences of (probably random) variables $\{\sigma_{1,k}\}_{k\ge 0}$ and $\{\sigma_{2,k}\}_{k\ge 0}$ such that inequalities \eqref{eq:second_moment_bound_new}, \eqref{eq:sigma_k+1_bound_1} and \eqref{eq:sigma_k+1_bound_2} hold and
	\begin{eqnarray}
	\mytextstyle	3L\sum\limits_{k=0}^K w_k\EE\|e^k\|^2 &\le& \mytextstyle \frac{1}{4}\sum\limits_{k=0}^K w_k\EE\left[f(x^k) - f(x^*)\right] + F_1\sigma_{1,0}^2 + F_2\sigma_{2,0}^2 + \gamma D_3 W_K \label{eq:sum_of_errors_bound_new}
	\end{eqnarray}
for all $k,	K\ge 0$, where $e^k = \frac{1}{n}\sum_{i=1}^n e_i^k$ and $\{W_K\}_{K\ge 0}$ and $\{w_k\}_{k\ge 0}$ are defined as 
	\begin{equation}
	\mytextstyle			W_K = \sum\limits_{k=0}^K w_k,\quad w_k = (1 - \eta)^{-(k+1)},\quad \eta = \min\left\{\frac{\gamma\mu}{2}, \frac{\rho_1}{4}, \frac{\rho_2}{4}\right\}. \label{eq:w_k_definition_new}
	\end{equation}		
\end{assumption}
This assumption is more flexible than Assumption~\ref{ass:key_assumption_finite_sums_new} and helps us to obtain a unified analysis of all methods falling in the error-feedback framework. We emphasize that in this assumption we do not assume that \eqref{eq:second_moment_bound_g_i^k_new} and \eqref{eq:variance_bound_g_i^k_new} hold \textit{explicitly}. Instead of this, we introduce inequality \eqref{eq:sum_of_errors_bound_new}, which is the key tool that helps us to analyze the effect of error-feedback and comes from the analysis from \cite{stich2020error} with needed adaptations connected with the first three inequalities. As we show in the appendix, this inequality can be derived for {\tt SGD} with error compensation and delayed updates under Assumption~\ref{ass:key_assumption_finite_sums_new} and, in particular, using \eqref{eq:second_moment_bound_g_i^k_new} and \eqref{eq:variance_bound_g_i^k_new}. As before, $D_3$ hides a variance that is not handled by variance reduction processes and $F_1$ and $F_2$ are some constants that typically depend on $L, B_1, B_2, \rho_1, \rho_2$ and $\gamma$.

We now proceed to stating our main theorem.

\begin{theorem}\label{thm:main_result_new}
	Let Assumptions~\ref{ass:quasi_strong_convexity},~\ref{ass:L_smoothness} and ~\ref{ass:key_assumption_new} be satisfied and~$\gamma \le \nicefrac{1}{4(A'+C_1M_1+C_2M_2)}$. Then for all $K\ge 0$ we have
	\begin{equation}
	\mytextstyle	\EE\left[f(\bar x^K) - f(x^*)\right] \le \left(1 - \eta\right)^K\frac{4(T^0 + \gamma F_1 \sigma_{1,0}^2+ \gamma F_2 \sigma_{2,0}^2)}{\gamma} + 4\gamma\left(D_1' + M_1D_2 + D_3\right) \label{eq:main_result_new}
	\end{equation}	
	when $\mu > 0$ and
	\begin{equation}
	\mytextstyle	\EE\left[f(\bar x^K) - f(x^*)\right] \le \frac{4(T^0 + \gamma F_1 \sigma_{1,0}^2+ \gamma F_2 \sigma_{2,0}^2)}{\gamma K} + 4\gamma\left(D_1' + M_1D_2 + D_3\right) \label{eq:main_result_new_cvx}
	\end{equation}
	when $\mu = 0$, where $\eta = \min\left\{\nicefrac{\gamma\mu}{2},\nicefrac{\rho_1}{4},\nicefrac{\rho_2}{4}\right\}$, $T^k \eqdef \|\tx^k - x^*\|^2 + M_1\gamma^2 \sigma_{1,k}^2 + M_2\gamma^2 \sigma_{2,k}^2$ and $M_1 = \frac{4B_1'}{3\rho_1}$, $M_2 = \frac{4\left(B_2' + \frac{4}{3}G\right)}{3\rho_2}$.
\end{theorem}

All the complexity results summarized in Table~\ref{tbl:special_cases2_ef} follow from this theorem; the detailed proofs of the main results are included in the appendix. Furthermore, in the appendix we include similar results but for methods employing {\em delayed} updates.

\section{Further Notation}

In what follows it will be useful to denote
$$\mytextstyle v^k \eqdef \frac{1}{n}\sum_i  v_i^k, \quad  g^k \eqdef \frac{1}{n}\sum_i g_i^k, \quad e^k \eqdef \frac{1}{n}\sum_i e_i^k.$$ By aggregating identities  \eqref{eq:error_update} across all $i$, we get $e^{k+1} = e^k + \gamma g^k - v^k.$ In our proofs we also use the perturbed iterates technique \cite{leblond2018improved,mania2017perturbed} based on the analysis of the following sequence
\begin{equation}
	\tx^k = x^k - e^k. \label{eq:perturbed_uterate}
\end{equation}
This sequence satisfies very useful for the analysis relation:
\begin{equation}
	\tx^{k+1} \overset{\eqref{eq:perturbed_uterate}}{=} x^{k+1} - e^{k+1} \overset{\eqref{eq:x^k+1_update},\eqref{eq:error_update}}{=} x^k - v^k - (e^k + \gamma g^k - v^k) = x^k - e^k - \gamma g^k \overset{\eqref{eq:perturbed_uterate}}{=} \tx^k - \gamma g^k. \label{eq:perturbed_iterates_key_relation}
\end{equation}

\section{{\tt SGD} as a Special Case}\label{sec:sgd}
In this section we want to show that our approach is general enough to cover many existing methods of {\tt SGD} type. Consider the following situation:
\begin{equation}
	v^k = \gamma g^k,\quad e^0 = 0.\label{eq:vk_ek_plain_sgd}
\end{equation}
It implies that $e^k = 0$ for all $k\ge 0$ and the updates rules \eqref{eq:x^k+1_update}-\eqref{eq:error_update} gives us a simple {\tt SGD}:
\begin{equation}
	x^{k+1} = x^k - \gamma g^k.\label{eq:plain_sgd}
\end{equation}
The following lemma formally shows that {\tt SGD} under general enough assumptions satisfies Assumption~\ref{ass:key_assumption_new}.
\begin{lemma}\label{lem:sgd_as_a_special_case}
	Let Assumptions~\ref{ass:quasi_strong_convexity}~and~\ref{ass:L_smoothness} be satisfies and inequalities \eqref{eq:unbiasedness_new}, \eqref{eq:second_moment_bound_new}, \eqref{eq:sigma_k+1_bound_1} and \eqref{eq:sigma_k+1_bound_2} hold. Then for the method \eqref{eq:plain_sgd} inequality \eqref{eq:sum_of_errors_bound_new} holds with $F_1 = F_2 = 0$ and $D_3 = 0$ for all $k\ge 0$.
\end{lemma}
\begin{proof}
	Since $e^k = 0$ and $f(x^k)\ge f(x^*)$ for all $k\ge 0$ we get
	\begin{equation*}
		3L\sum\limits_{k=0}^K w_k\EE\|e^k\|^2 = 0 \le \frac{1}{4}\sum\limits_{k=0}^K w_k\EE\left[f(x^k) - f(x^*)\right]
	\end{equation*}
	which concludes the proof.	 
\end{proof}

It implies that all methods considered in Chapter~\ref{ch:sigma_k} fit our framework. Moreover, using Theorem~\ref{thm:main_result_new} we derive the following result.
\begin{theorem}\label{thm:main_result_sgd}
	Let Assumptions~\ref{ass:quasi_strong_convexity}~and~\ref{ass:L_smoothness} be satisfied, inequalities \eqref{eq:unbiasedness_new}, \eqref{eq:second_moment_bound_new}, \eqref{eq:sigma_k+1_bound_1}, \eqref{eq:sigma_k+1_bound_2} hold and $\gamma \le \nicefrac{1}{4(A'+C_1M_1 + C_2M_2)}$. Then for the method \eqref{eq:plain_sgd} for all $K\ge 0$ we have
	\begin{equation*}
		\EE\left[f(\bar x^K) - f(x^*)\right] \le \left(1 - \min\left\{\frac{\gamma\mu}{2},\frac{\rho_1}{4},\frac{\rho_2}{4}\right\}\right)^K\frac{4T^0}{\gamma} + 4\gamma\left(D_1' + M_1D_2\right),
	\end{equation*}	
	when $\mu > 0$ and
	\begin{equation}
		\EE\left[f(\bar x^K) - f(x^*)\right] \le \frac{4T^0}{\gamma K} + 4\gamma\left(D_1' + M_1D_2\right) \notag
	\end{equation}
	when $\mu = 0$, where $T^k \eqdef \|x^k - x^*\|^2 + M_1\gamma^2 \sigma_{1,k}^2 + M_2\gamma^2 \sigma_{2,k}^2$ and $M_1 = \frac{4B_1'}{3\rho_1}$, $M_2 = \frac{4\left(B_2' + \frac{4}{3}G\right)}{3\rho_2}$.
\end{theorem}
In particular, if $\sigma_{2,k}^2 \equiv 0$, then our assumption coincides with the key assumption from \cite{gorbunov2019unified} and our theorem recovers the same rates as in \cite{gorbunov2019unified} when $\mu > 0$. The case when $\mu = 0$ was not considered in \cite{gorbunov2019unified}, while in our analysis we get it for free.

\section{Special Cases: {\tt SGD}}\label{sec:special_cases_sgd}
To illustrate the generality of our approach, we develop and analyse a new special case of {\tt SGD} without error-feedback and show that in some cases, our framework recovers tighter rates than the framework from \cite{gorbunov2019unified}.

\subsection{{\tt DIANA} with Arbitrary Sampling and Double Quantization}\label{sec:diana_arbitrary_sampling}
In this section we consider problem \eqref{eq:main_problem_ef} with $f(x)$ being $\mu$-quasi strongly convex and $f_i(x)$ satisfying \eqref{eq:f_i_sum_ef} where functions $f_{ij}(x)$ are differentiable, but not necessary convex. Following \cite{gower2019sgd} we construct a stochastic reformulation of this problem:
\begin{equation}
	f(x) = \EE_{\cD}\left[f_\xi(x)\right],\quad f_\xi(x) = \frac{1}{n}\sum\limits_{i=1}^n f_{\xi_i}(x),\quad f_{\xi_i}(x) = \frac{1}{m}\sum\limits_{j=1}^m \xi_{ij}f_{ij}(x), \label{eq:sr_def}
\end{equation}
where $\xi = (\xi_1^\top,\ldots, \xi_n^\top), \xi_i = (\xi_{i1},\ldots, \xi_{im})^\top$ is a random vector with distribution $\cD_i$ such that $\EE_{\cD_i}[\xi_{ij}] = 1$ for all $i\in[n], j\in[m]$ and the following assumption holds.
\begin{assumption}[Expected smoothness]\label{ass:exp_smoothness}
	We assume that functions $f_1,\ldots, f_n$ are $\cL$-smooth in expectation w.r.t.\ distributions $\cD_1,\ldots,\cD_n$, i.e., there exists constant $\cL = \cL(f,\cD_1,\ldots,\cD_n)$ such that
	\begin{equation}
		\EE_{\cD_i}\left[\|\nabla f_{\xi_i}(x) - \nabla f_{\xi_i}(x^*)\|^2\right] \le 2\cL D_{f_i}(x,x^*)\label{eq:exp_smoothness}
	\end{equation}
	for all $i\in [n]$ and $x\in\R^d$.
\end{assumption}

To solve this problem, we consider {\tt DIANA} \cite{mishchenko2019distributed, horvath2019stochastic}~--- a distributed stochastic method using unbiased compressions or \textit{quantizations} for communication between workers and master. We start with the formal definition of quantization. In \cite{mishchenko2019distributed, horvath2019stochastic} {\tt DIANA} was analyzed under the assumption that stochastic gradients $g_i^k$ have uniformly bounded variances which is not very practical.

Therefore, we consider a slightly different method called {\tt DIANAsr-DQ} which works with the stochastic reformulation \eqref{eq:sr_def} of problem \eqref{eq:main_problem_ef}+\eqref{eq:f_i_sum_ef}, see Algorithm~\ref{alg:DIANAsr-DQ}.
\begin{algorithm}[t]
   \caption{{\tt DIANAsr} with Double Compression ({\tt DIANAsr-DQ})}\label{alg:DIANAsr-DQ}
\begin{algorithmic}[1]
   \Require learning rates $\gamma>0$, $\alpha \in (0,1]$, initial vectors $x^0, h_1^0,\ldots, h_n^0 \in \R^d$
   \State Set $h^0 = \frac{1}{n}\sum_{i=1}^n h_i^0$   
   \For{$k=0,1,\dotsc$}
       \State Broadcast $g^{k-1}$ to all workers \Comment{If $k=0$, then broadcast $x^0$}
        \For{$i=1,\dotsc,n$ in parallel}
        	\State $x^{k} = x^{k-1} - \gamma g^{k-1}$ \Comment{Ignore this line if $k=0$}
			\State Sample $g_i^{k,1} = \nabla f_{\xi_i^k}(x^k)$ satisfying Assumption~\ref{ass:exp_smoothness} independtently from other workers
            \State $\hat\Delta_i^k = g_i^{k,1} - h_i^k$
            \State Sample $\Delta_i^k \sim Q_1(\hat\Delta_i^k)$ indepently from other workers
            \State $g_i^{k,2} = h_i^k + \Delta_i^k$
            \State $h_i^{k+1} = h_i^k + \alpha \Delta_i^k$
        \EndFor
        \State $g^{k,2} = \frac{1}{n}\sum_{i=1}^ng_i^{k,2} = h^k + \frac{1}{n}\sum_{i=1}^n\Delta_i^k$
        \State $h^{k+1} = \frac{1}{n}\sum\limits_{i=1}^n h_i^{k+1} = h^k + \alpha\frac{1}{n}\sum\limits_{i=1}^n\Delta_i^k$
       \State Sample $g^k \sim Q_2(g^{k,2})$
       \State $x^{k+1} = x^{k} - \gamma g^{k-1}$
   \EndFor
\end{algorithmic}
\end{algorithm}
Moreover, to illustrate the flexibility of our approach, we consider compression not only on the workers' side but also on the master side. To perform an update of {\tt DIANAsr-DQ} master needs to gather quantized gradient differences $\Delta_i^k$ and the to broadcast quantized stochastic gradient $g^k$ to all workers. Clearly, in this case, only compressed vectors participate in communication.  

In the concurrent work \cite{philippenko2020artemis} the same method was independently proposed under the name of {\tt Artemis}. However, our analysis is slightly more general: it is based on Assumption~\ref{ass:exp_smoothness} while in \cite{philippenko2020artemis} authors assume $L$-cocoercivity of stochastic gradients almost surely. Next, a very similar approach was considered in \cite{tang2019doublesqueeze}, where authors present a method with error compensation on master and worker sides. Moreover, recently another method called {\tt DORE} was developed in \cite{liu2019double}, which uses {\tt DIANA}-trick on the worker side and error compensation on the master side. However, in these methods, compression operators are the same on both sides, despite the fact that gathering the information often costs much more than broadcasting. Therefore, the natural idea is in using different quantization for gathering and broadcasting, and it is what {\tt DIANAsr-DQ} does. Moreover, we do not assume uniform boundedness of the second moment of the stochastic gradient like in \cite{tang2019doublesqueeze}, and we also do not assume uniform boundedness of the variance of the stochastic gradient like in \cite{liu2019double}. Assumption~\ref{ass:exp_smoothness} is more natural and always holds for the problems \eqref{eq:main_problem_ef}+\eqref{eq:f_i_sum_ef} when $f_{ij}$ are convex and $L$-smooth for each $i\in[n]$, $j\in[m]$. In contrast, in the same setup, there exist such problems that the variance of the stochastic gradients is not uniformly upper bounded by any finite constant.

We assume that $Q_1$ and $Q_2$ satisfy \eqref{eq:quantization_def} with parameters $\omega_1$ and $\omega_2$ respectively.
\begin{lemma}\label{lem:diana_second_moment_bound}
	Let Assumption~\ref{ass:exp_smoothness} be satisfied. Then, for all $k\ge 0$ we have
	\begin{eqnarray}
		\EE\left[g^k\mid x^k\right] &=& \nabla f(x^k), \label{eq:diana_unbiasedness}\\
		\EE\left[\|g^k\|^2\mid x^k\right] &\le& 2\cL(1+\omega_2)\left(2+\frac{3\omega_1}{n}\right)\left(f(x^k) - f(x^*)\right) + \frac{3\omega_1(1+\omega_2)}{n}\sigma_k^2 + D_1', \label{eq:diana_second_moment_bound}
	\end{eqnarray}
	where $\sigma_k^2 = \frac{1}{n}\sum_{i=1}^n\|h_i^k - \nabla f(x^*)\|^2$ and $D_1' = \frac{(2+3\omega_1)(1+\omega_2)}{n^2}\sum\limits_{i=1}^n\EE_{\cD_i}\left[\|\nabla f_{\xi_i}(x^*) - \nabla f_i(x^*)\|^2\right]$.
\end{lemma}
\begin{proof}
	First of all, we show inbiasedness of $g^k$:
	\begin{eqnarray*}
		\EE\left[g^k\mid x^k\right] &\overset{\eqref{eq:tower_property},\eqref{eq:quantization_def}}{=}& \EE\left[g^{k,2}\mid x^k\right] = h^k + \frac{1}{n}\sum\limits_{i=1}^n\EE\left[\Delta_i^k\mid x^k\right]\\
		&\overset{\eqref{eq:tower_property},\eqref{eq:quantization_def}}{=}& h^k + \frac{1}{n}\sum\limits_{i=1}^n\EE\left[\hat \Delta_i^k\mid x^k\right]\\
		&=& h^k + \frac{1}{n}\sum\limits_{i=1}^n\left(\nabla f_i(x^k) - h_i^k\right) = \nabla f(x^k).
	\end{eqnarray*}
	Next, to denote mathematical expectation w.r.t.\ the randomness coming from quantizations $Q_1$ and $Q_2$ at iteration $k$ we use $\EE_{Q_1^k}[\cdot]$ and $\EE_{Q_2^k}[\cdot]$ respectively. Using these notations and the definition of quantization we derive
	\begin{eqnarray*}
		\EE_{Q_2^k}[\|g^k\|^2] &\overset{\eqref{eq:variance_decomposition},\eqref{eq:quantization_def}}{=}& \|g^{k,1}\|^2 + \EE_{Q_2^k}\left[\|g^{k,2}-g^{k,1}\|^2\right]\\
		&\overset{\eqref{eq:quantization_def}}{\le}& (1+\omega_2)\|g^{k,1}\|^2. 
	\end{eqnarray*}
	Taking the conditopnal mathematical expectation $\EE_{Q_1^k}[\cdot]$ from the both sides of previous inequality and using the  independence of $\Delta_i^1,\ldots,\Delta_i^n$ we get
	\begin{eqnarray*}
		\EE_{Q_1^k,Q_2^k}\left[\|g^k\|^2\right] &\overset{\eqref{eq:tower_property}}{=}& (1+\omega_2)\EE_{Q_1^k}\left[\|g^{k,1}\|^2\right]  = (1+\omega_2)\EE_{Q_1^k}\left[\left\|\frac{1}{n}\sum\limits_{i=1}^n (h_i^k + \Delta_i^k)\right\|^2\right]\\
		&\overset{\eqref{eq:variance_decomposition}}{=}& (1+\omega_2)\left\|\frac{1}{n}\sum\limits_{i=1}^n\left(h_i^k + \hat \Delta_i^k\right)\right\|^2 + (1+\omega_2)\EE_{Q_1^k}\left[\left\|\frac{1}{n}\sum\limits_{i=1}^n(\Delta_i^k - \hat\Delta_i^k)\right\|^2\right]\\
		&=& (1+\omega_2)\left\|\frac{1}{n}\sum\limits_{i=1}^n\left(\nabla f_{\xi_i^k}(x^k) - \nabla f_{\xi_i^k}(x^*) + \nabla f_{\xi_i^k}(x^*) - \nabla f_i(x^*)\right)\right\|^2\\
		&&\quad + \frac{(1+\omega_2)}{n^2}\sum\limits_{i=1}^n\EE_{Q_1^k}\left[\|\Delta_i^k - \hat\Delta_i^k\|^2\right]\\
		&\overset{\eqref{eq:a_b_norm_squared},\eqref{eq:quantization_def}}{\le}& \frac{2(1+\omega_2)}{n}\sum\limits_{i=1}^n\|\nabla f_{\xi_i^k}(x^k) - \nabla f_{\xi_i^k}(x^*)\|^2\\
		&&\quad + 2(1+\omega_2)\left\|\frac{1}{n}\sum\limits_{i=1}^n\left(\nabla f_{\xi_i^k}(x^*) - \nabla f_i(x^*)\right)\right\|^2\\
		&&\quad + \frac{\omega_1(1+\omega_2)}{n^2}\sum\limits_{i=1}^n\|\nabla f_{\xi_i^k}(x^k) - h_i^k\|^2\\
		&\overset{\eqref{eq:a_b_norm_squared}}{\le}& \frac{2(1+\omega_2)}{n}\sum\limits_{i=1}^n\|\nabla f_{\xi_i^k}(x^k) - \nabla f_{\xi_i^k}(x^*)\|^2\\
		&&\quad  + 2(1+\omega_2)\left\|\frac{1}{n}\sum\limits_{i=1}^n\left(\nabla f_{\xi_i^k}(x^*) - \nabla f_i(x^*)\right)\right\|^2\\
		&&\quad + \frac{3\omega_1(1+\omega_2)}{n^2}\sum\limits_{i=1}^n\|\nabla f_{\xi_i^k}(x^k) - \nabla f_{\xi_i^k}(x^*)\|^2\\
		&&\quad + \frac{3\omega_1(1+\omega_2)}{n^2}\sum\limits_{i=1}^n\|\nabla f_{\xi_i^k}(x^*) - \nabla f_i(x^*)\|^2\\
		&&\quad + \frac{3\omega_1(1+\omega_2)}{n^2}\sum\limits_{i=1}^n\|h_i^k - \nabla f_i(x^*)\|^2.
	\end{eqnarray*}
	Finally, we take conditional mathematical expectation $\EE[\cdot\mid x^k]$ from the both sides of the inequality above and use the independece of $\xi_1^k,\ldots,\xi_n^k$:
	\begin{eqnarray*}
		\EE\left[\|g^k\|^2\mid x^k\right] &\overset{\eqref{eq:exp_smoothness}}{\le}& 2\cL(1+\omega_2)\left(2+\frac{3\omega_1}{n}\right)(f(x^k) - f(x^*))+\frac{3\omega_1(1+\omega_2)}{n}\sigma_{k}^2\\
		&&\quad + 2(1+\omega_2)\EE\left[\left\|\frac{1}{n}\sum\limits_{i=1}^n\left(\nabla f_{\xi_i^k}(x^*) - \nabla f_i(x^*)\right)\right\|^2\mid x^k\right]\\
		&&\quad + \frac{3\omega_1(1+\omega_2)}{n^2}\sum\limits_{i=1}^n\EE_{\cD_i}\left[\|\nabla f_{\xi_i}(x^*) - \nabla f_i(x^*)\|^2\right]\\
		&=& 2\cL(1+\omega_2)\left(2+\frac{3\omega_1}{n}\right)(f(x^k) - f(x^*))+\frac{3\omega_1(1+\omega_2)}{n}\sigma_{k}^2\\
		&&\quad + \frac{(1+\omega_2)(2+3\omega_1)}{n^2}\sum\limits_{i=1}^n\EE_{\cD_i}\left[\|\nabla f_{\xi_i}(x^*) - \nabla f_i(x^*)\|^2\right].
	\end{eqnarray*}
\end{proof}

\begin{lemma}\label{lem:diana_sigma_k+1_bound}
	Let $f_i$ be convex and $L$-smooth, Assumption~\ref{ass:exp_smoothness} holds and $\alpha \le \nicefrac{1}{(\omega_1+1)}$. Then, for all $k\ge 0$ we have
	\begin{equation}
		\EE\left[\sigma_{k+1}^2\mid x^k\right] \le (1 - \alpha)\sigma_k^2 + 2\alpha(3\cL+4L)(f(x^k) - f(x^*)) + D_2, \label{eq:diana_sigma_k+1_bound}
	\end{equation}
	where $\sigma_k^2 = \frac{1}{n}\sum_{i=1}^n\|h_i^k - \nabla f_i(x^*)\|^2$ and $D_2 = \frac{3\alpha}{n}\sum_{i=1}^n \EE_{\cD_i}\left[\|\nabla f_{\xi_i}(x^*) - \nabla f_i(x^*)\|^2\right]$.
\end{lemma}
\begin{proof}
	For simplicity, we introduce new notation: $h_i^* \eqdef \nabla f_i(x^*)$. Using this we derive an upper bound for the second moment of $h_i^{k+1} - h_i^*$:
	\begin{eqnarray*}
		\EE\left[\|h_i^{k+1} - h_i^*\|^2\mid x^k\right] &=& \EE\left[\left\|h_i^k - h_i^* + \alpha \Delta_i^k \right\|^2\mid x^k\right]\\
		&\overset{\eqref{eq:quantization_def}}{=}& \|h_i^k - h_i^*\|^2 +2\alpha\langle h_i^k - h_i^*, \nabla f_i(x^k) - h_i^k \rangle + \alpha^2\EE\left[\|\Delta_i^k\|^2\mid x^k\right]\\
		&\overset{\eqref{eq:quantization_def},\eqref{eq:tower_property}}{\le}& \|h_i^k - h_i^*\|^2 +2\alpha\langle h_i^k - h_i^*, \nabla f_i(x^k) - h_i^k \rangle\\
		&&\quad + \alpha^2(\omega_1+1)\EE\left[\|\nabla f_{\xi_i^k}(x^k) - h_i^k\|^2\mid x^k\right].
	\end{eqnarray*}
	Using variance decomposition \eqref{eq:variance_decomposition} and $\alpha \le \nicefrac{1}{(\omega_1+1)}$ we get
	\begin{eqnarray*}
		\alpha^2(\omega_1+1)\EE_{\cD_i}\left[\|\nabla f_{\xi_i^k}(x^k) - h_i^k\|^2\right] &\overset{\eqref{eq:variance_decomposition}}{=}& \alpha^2(\omega_1+1)\EE_{\cD_i}\left[\|\nabla f_{\xi_i^k}(x^k) - \nabla f_i(x^k)\|^2\right]\\
		&&\quad + \alpha^2(\omega_1+1)\|\nabla f_i(x^k) - h_i^k\|^2\\
		&\overset{\eqref{eq:a_b_norm_squared}}{\le}& 3\alpha\EE_{\cD_i}\left[\|\nabla f_{\xi_i^k}(x^k) - \nabla f_{\xi_i^k}(x^*)\|^2\right]\\
		&&\quad +3\alpha\EE_{\cD_i}\left[\|\nabla f_{\xi_i^k}(x^*) - \nabla f_i(x^*)\|^2\right]\\
		&&\quad +3\alpha\|\nabla f_i(x^k) - \nabla f_i(x^*)\|^2\\
		&&\quad + \alpha\|\nabla f_i(x^k) - h_i^k\|^2\\
		&\overset{\eqref{eq:L_smoothness_cor_3},\eqref{eq:exp_smoothness}}{\le}& 6\alpha(\cL + L)D_{f_i}(x^k,x^*) + \alpha\|\nabla f_i(x^k) - h_i^k\|^2\\
		&&\quad +3\alpha\EE_{\cD_i}\left[\|\nabla f_{\xi_i^k}(x^*) - \nabla f_i(x^*)\|^2\right]
	\end{eqnarray*}
	Putting all together we obtain
	\begin{eqnarray*}
		\EE\left[\|h_i^{k+1} - h_i^*\|^2\mid x^k\right] &\le& \|h_i^k - h_i^*\|^2 + \alpha\left\langle \nabla f_i(x^k) - h_i^k, f_i(x^k) + h_i^k - 2h_i^* \right\rangle\\
		&&\quad + 6\alpha(\cL + L)D_{f_i}(x^k,x^*) +3\alpha\EE_{\cD_i}\left[\|\nabla f_{\xi_i^k}(x^*) - \nabla f_i(x^*)\|^2\right]\\
		&\overset{\eqref{eq:a-b_a+b}}{=}& \|h_i^k - h_i^*\|^2 + \alpha\|\nabla f_i(x^k) - h_i^*\|^2 - \alpha\|h_i^k - h_i^*\|^2\\
		&&\quad + 6\alpha(\cL + L)D_{f_i}(x^k,x^*) +3\alpha\EE_{\cD_i}\left[\|\nabla f_{\xi_i^k}(x^*) - \nabla f_i(x^*)\|^2\right]\\
		&\overset{\eqref{eq:L_smoothness_cor_3}}{\le}& (1-\alpha)\|h_i^k - h_i^*\|^2 + \alpha(6\cL + 8L)D_{f_i}(x^k,x^*)\\
		&&\quad +3\alpha\EE_{\cD_i}\left[\|\nabla f_{\xi_i^k}(x^*) - \nabla f_i(x^*)\|^2\right].
	\end{eqnarray*}
	Summing up the above inequality for $i=1,\ldots, n$ we derive
	\begin{eqnarray*}
		\frac{1}{n}\sum\limits_{i=1}^n\EE\left[\|h_i^{k+1} - h_i^*\|^2\mid x^k\right] &\le& \frac{1-\alpha}{n}\sum\limits_{i=1}^n\|h_i^k - h_i^*\|^2 + \alpha(6\cL + 8L)(f(x^k) - f(x^*))\\
		&&\quad + \frac{3\alpha}{n}\sum\limits_{i=1}^n \EE_{\cD_i}\left[\|\nabla f_{\xi_i^k}(x^*) - \nabla f_i(x^*)\|^2\right].
	\end{eqnarray*}
\end{proof}

\begin{theorem}\label{thm:diana}
	Assume that $f_i(x)$ is convex and $L$-smooth for all $i=1,\ldots, n$, $f(x)$ is $\mu$-quasi strongly convex and Assumption~\ref{ass:exp_smoothness} holds. Then {\tt DIANAsr-DQ} satisfies Assumption~\ref{ass:key_assumption_new} with
	\begin{gather*}
		A' = \cL(1+\omega_2)\left(2+\frac{3\omega_1}{n}\right),\quad B_1' = \frac{3\omega_1(1+\omega_2)}{n},\\
		D_1' = \frac{(2+3\omega_1)(1+\omega_2)}{n^2}\sum\limits_{i=1}^n\EE_{\cD_i}\left[\|\nabla f_{\xi_i}(x^*) - \nabla f_i(x^*)\|^2\right],\\
		\sigma_{1,k}^2 = \sigma_k^2 = \frac{1}{n}\sum\limits_{i=1}^n\|h_i^k - \nabla f_i(x^*)\|^2,\quad B_2' = 0,\quad \sigma_{2,k}^2\equiv 0,\quad \rho_1 = \alpha,\quad \rho_2 = 1,\\
		C_1 = \alpha(3\cL+4L),\quad C_2 = 0,\quad D_2 = \frac{3\alpha}{n}\sum_{i=1}^n \EE_{\cD_i}\left[\|\nabla f_{\xi_i}(x^*) - \nabla f_i(x^*)\|^2\right],\\
		 G = 0,\quad F_1 = F_2 = 0,\quad D_3 = 0,
	\end{gather*}
	with $\gamma$ and $\alpha$ satisfying
	\begin{equation*}
		\gamma \le \frac{1}{4(1+\omega_2)\left(\cL\left(2+\frac{15\omega_1}{n}\right)+\frac{16L\omega_1}{n}\right)},\quad \alpha \le \frac{1}{\omega+1},\quad M_1 = \frac{4\omega_1(1+\omega_2)}{n\alpha},\quad M_2 = 0
	\end{equation*}
	and for all $K \ge 0$
	\begin{equation*}
		\EE\left[f(\bar x^K) - f(x^*)\right] \le \left(1 - \min\left\{\frac{\gamma\mu}{2},\frac{\alpha}{4}\right\}\right)^K\frac{4T^0}{\gamma} + 4\gamma\left(D_1' + M_1D_2\right),
	\end{equation*}	
	when $\mu > 0$ and
	\begin{equation*}
		\EE\left[f(\bar{x}^K) - f(x^*)\right] \le \frac{4T^0}{\gamma K} + 4\gamma\left(D_1' + M_1D_2\right)
	\end{equation*}
	when $\mu=0$, where $T^k \eqdef \|x^k - x^*\|^2 + M_1\gamma^2 \sigma_{1,k}^2$.
\end{theorem}
In other words, if 
\begin{equation*}
		\gamma = \frac{1}{4(1+\omega_2)\left(\cL\left(2+\frac{15\omega_1}{n}\right)+\frac{16L\omega_1}{n}\right)},\quad \alpha = \frac{1}{\omega+1}
\end{equation*}
and $D_1 = 0$, i.e., $\nabla f_{\xi_i^k}(x^k) = \nabla f_i(x^k)$ almost surely, {\tt DIANAsr-DQ} converges with the linear rate
\begin{equation*}
	\cO\left(\left(\omega_1 + \frac{\cL}{\mu}(1+\omega_2)\left(1+\frac{\omega_1}{n}\right)\right)\ln\frac{1}{\varepsilon}\right)
\end{equation*}
to the exact solution. Applying Lemma~\ref{lem:lemma2_stich} we establish the rate of convergence to $\varepsilon$-solution.
\begin{corollary}\label{cor:diana_str_cvx_cor}
	Let the assumptions of Theorem~\ref{thm:diana} hold and $\mu > 0$. Then after $K$ iterations of {\tt DIANAsq-DQ} with the stepsize
	\begin{eqnarray*}
		\gamma_0 &=& \frac{1}{4(1+\omega_2)\left(\cL\left(2+\frac{15\omega_1}{n}\right)+\frac{16L\omega_1}{n}\right)}\\
		\gamma &=& \min\left\{\gamma_0, \frac{\ln\left(\max\left\{2,\frac{\mu^2K^2(\|x^0-x^*\|^2+M_1\gamma_0^2\sigma_{1,0}^2)}{D_1'+M_1D_2}\right\}\right)}{\mu K}\right\},\quad M_1 = \frac{4\omega_1(1+\omega_2)}{n\alpha}
	\end{eqnarray*}		
	and $\alpha = \frac{1}{\omega+1}$ we have
	\begin{equation*}
		\EE\left[f(\bar{x}^K) - f(x^*)\right] = \widetilde\cO\left(A'\|x^0 - x^*\|^2\exp\left(-\min\left\{\frac{\mu}{A'},\frac{1}{\omega_1}\right\}K\right) + \frac{D_1'+M_1D_2}{\mu K}\right).
	\end{equation*}
	That is, to achive $\EE\left[f(\bar{x}^K) - f(x^*)\right] \le \varepsilon$ {\tt DIANAsq-DQ} requires
	\begin{equation*}
		\widetilde{\cO}\left(\omega_1 + \frac{\cL\left(1+\frac{\omega_1}{n}\right)(1+\omega_2)}{\mu} + \frac{(1+\omega_1)(1+\omega_2)}{n^2\mu\varepsilon}\sum\limits_{i=1}^n\EE_{\cD_i}\|\nabla f_{\xi_i}(x^*)-\nabla f_i(x^*)\|^2\right) \text{ iterations.}
	\end{equation*}
\end{corollary}

Applying Lemma~\ref{lem:lemma_technical_cvx} we get the complexity result in the case when $\mu = 0$.
\begin{corollary}\label{cor:diana_cvx_cor}
	Let the assumptions of Theorem~\ref{thm:diana} hold and $\mu = 0$. Then after $K$ iterations of {\tt DIANAsq-DQ} with the stepsize
	\begin{eqnarray*}
		\gamma_0 &=& \frac{1}{4(1+\omega_2)\left(\cL\left(2+\frac{15\omega_1}{n}\right)+\frac{16L\omega_1}{n}\right)}\\
		\gamma &=& \min\left\{\gamma_0, \sqrt{\frac{\|x^0 - x^*\|^2}{M_1\sigma_{1,0}^2}}, \sqrt{\frac{\|x^0 - x^*\|^2}{(D_1'+M_1D_2) K}}\right\},\quad M_1 = \frac{4\omega_1(1+\omega_2)}{n\alpha}
	\end{eqnarray*}		
	and $\alpha = \frac{1}{\omega+1}$ we have $\EE\left[f(\bar{x}^K) - f(x^*)\right]$ of order
	\begin{equation*}
		\cO\left(\frac{\cL R_0^2(1+\omega_2)\left(1+\frac{\omega_1}{n}\right)}{K} + \frac{R_0\sigma_{1,0}(1+\omega_1)\sqrt{1+\omega_2}}{\sqrt{n}K} +  \frac{R_0\sqrt{(1+\omega_1)(1+\omega_2)D_{\text{opt}}}}{\sqrt{nK}}\right)
	\end{equation*}
	where $R_0 = \|x^0 - x^*\|^2, D_{\text{opt}} = \frac{1}{n}\sum\limits_{i=1}^n\EE_{\cD_i}\|\nabla f_{\xi_i}(x^*)-\nabla f_i(x^*)\|^2$. That is, to achive \newline $\EE\left[f(\bar{x}^K) - f(x^*)\right] \le \varepsilon$ {\tt DIANAsq-DQ} requires
	\begin{equation*}
		\cO\left(\frac{\cL R_0^2(1+\omega_2)\left(1+\frac{\omega_1}{n}\right)}{\varepsilon} + \frac{R_0\sigma_{1,0}(1+\omega_1)\sqrt{1+\omega_2}}{\sqrt{n}\varepsilon} +  \frac{R_0^2(1+\omega_1)(1+\omega_2)D_{\text{opt}}}{n\varepsilon^2}\right)
	\end{equation*}
	iterations.
\end{corollary}

\subsection{Recovering Tight Complexity Bounds for {\tt VR-DIANA}}
In this section we consider the same problem \eqref{eq:main_problem_ef}+\eqref{eq:f_i_sum_ef} and variance reduced version of {\tt DIANA} called {\tt VR-DIANA} \cite{horvath2019stochastic}, see Algorithm~\ref{alg:vr-diana}.
\begin{algorithm}[t]
   \caption{{\tt VR-DIANA} based on {\tt LSVRG} (Variant 1), {\tt SAGA} (Variant 2), \cite{horvath2019stochastic}}
   \label{alg:vr-diana}
\begin{algorithmic}[1]
        \Require{learning rates $\alpha > 0$ and $\gamma > 0$, initial vectors $x^0, h_{1}^0, \dots, h_{n}^0$, $h^0 = \frac{1}{n}\sum_{i=1}^n h_i^0$}
        \For{$k = 0,1,\ldots$}
        \State Sample random 
            $
                u^k = \begin{cases}
                    1,& \text{with probability } \frac{1}{m}\\
                    0,& \text{with probability } 1 - \frac{1}{m}\\
                \end{cases}
            $ \Comment{only for Variant 1}
        \State Broadcast $x^k$, $u^k$ to all workers\;
            \For{$i = 1, \ldots, n$ in parallel} \Comment{Worker side}
            \State Pick $j_i^k$ uniformly at random from $[m]$\;
            \State $\mu_i^k = \frac{1}{m} \sum\limits_{j=1}^{m} \nabla f_{ij}(w_{ij}^k)$\label{ln:mu} \;
            \State $g_i^k = \nabla f_{ij_i^k}(x^k) - \nabla f_{ij_i^k}(w_{ij_i^k}^k) + \mu_i^k$\;
            \State $\hat{\Delta}_i^k = Q(g_i^k - h_i^k)$\;
            \State $h_i^{k+1} = h_i^k + \alpha \hat{\Delta}_i^k$\;
                \For{$j = 1, \ldots, m$}
                    \State
                    $
                    w_{ij}^{k+1} =
                    \begin{cases}
                        x^k, & \text{if } u^k = 1 \\
                        w_{ij}^k, &\text{if } u^k = 0\\
                    \end{cases}
                    $ \Comment{Variant 1 (L-SVRG): update epoch gradient if $u^k = 1$}
                    \State
                    $
                    w_{ij}^{k+1} =
                    \begin{cases}
                    x^k, & j = j_i^k\\
                    w_{ij}^k, & j \neq j_i^k\\
                    \end{cases}
                    $ \Comment{Variant 2 (SAGA): update gradient table}
                \EndFor
            \EndFor
            \State $h^{k+1} \! = \! h^k \!+\! \frac{\alpha}{n} \displaystyle \sum_{i=1}^n \hat{\Delta}_i^k$ \Comment{Gather quantized updates} 
            \State $g^k = \frac{1}{n}\sum\limits_{i=1}^{n} (\hat{\Delta}_i^k + h_i^k)$\;
            \State $x^{k+1} = x^k - \gamma g^k$\;
        \EndFor
\end{algorithmic}  
\end{algorithm}
For simplicity we assume that each $f_{ij}$ is convex and $L$-smooth and $f_i$ is additionally $\mu$-strongly convex.
\begin{lemma}[Lemmas 3, 5, 6 and 7 from \cite{horvath2019stochastic}]\label{lemmas_vr_diana}
    Let $\alpha \le \frac{1}{\omega+1}$. Then for all iterates $k\ge 0$ of Algorithm~\ref{alg:vr-diana} the following inequalities hold:
    \begin{eqnarray}
        \EE\left[g^k\mid x^k\right] &=& \nabla f(x^k),\label{eq:unbiased_g_k_vr_diana}\\
        \EE\left[H^{k+1}\mid x^k\right] &\le& \left(1-\alpha\right)H^k + \frac{2\alpha}{m}D^k + 8\alpha Ln\left(f(x^k) - f(x^*)\right),\label{eq:H_k+1_bound_vr_diana}\\
        \EE\left[D^{k+1}\mid x^k\right] &\le& \left(1 - \frac{1}{m}\right)D^k + 2Ln\left(f(x^k) - f(x^*)\right),\label{eq:D_k+1_bound_vr_diana}\\
        \EE\left[\|g^k\|^2\mid x^k\right] &\le& 2L\left(1+\frac{4\omega + 2}{n}\right)\left(f(x^k)-f(x^*)\right) + \frac{2\omega}{n^2}\frac{D^k}{m} + \frac{2(\omega+1)}{n^2}H^k,\label{eq:second_moment_g_k_vr_diana}
    \end{eqnarray}
    where $H^k = \sum\limits_{i=1}^n\|h_i^k - \nabla f_i(x^*)\|^2$ and $D^k = \sum\limits_{i=1}^n\sum\limits_{j=1}^m\|\nabla f_{ij}(w_{ij}^k) - \nabla f_{ij}(x^*)\|^2$.
\end{lemma}
This lemma shows that {\tt VR-DIANA} satisfies \eqref{eq:second_moment_bound_new}, \eqref{eq:sigma_k+1_bound_1} and \eqref{eq:sigma_k+1_bound_2}. Applying Theorem~\ref{thm:main_result_sgd} we get the following result.
\begin{theorem}\label{thm:vr-diana}
	Assume that $f_{ij}(x)$ is convex and $L$-smooth for all $i=1,\ldots, n$ and $f_i(x)$ is $\mu$-strongly convex for all $i=1,\ldots,n$. Then {\tt VR-DIANA} satisfies Assumption~\ref{ass:key_assumption_new} with
	\begin{gather*}
		A' = L\left(1+\frac{4\omega+2}{n}\right),\quad B_1' = \frac{2(\omega+1)}{n},\quad D_1' = 0,\\
		\sigma_{1,k}^2 = H^{k} = \frac{1}{n}\sum\limits_{i=1}^n\|h_i^k - \nabla f_i(x^*)\|^2,\quad B_2' = \frac{2\omega}{n},\\
		\sigma_{2,k}^2 =  D^k = \frac{1}{nm}\sum\limits_{i=1}^n\sum\limits_{j=1}^m\|\nabla f_{ij}(w_{ij}^k) - \nabla f_{ij}(x^*)\|^2,\quad \rho_1 = \alpha,\quad \rho_2 = \frac{1}{m},\\
		C_1 = 4\alpha L,\quad C_2 = \frac{L}{m},\quad D_2 = 0,\quad G = 2,\quad F_1 = F_2 = 0,\quad D_3 = 0,
	\end{gather*}
	with $\gamma$ and $\alpha$ satisfying
	\begin{equation*}
		\gamma \le \frac{3}{L\left(\frac{41}{3}+\frac{52\omega+35}{n}\right)},\quad \alpha \le \frac{1}{\omega+1},\quad M_1 = \frac{8(\omega+1)}{3n\alpha},\quad M_2 = \frac{8\omega m}{3n} + \frac{32m}{9}
	\end{equation*}
	and for all $K \ge 0$
	\begin{equation*}
		\EE\left[f(\bar x^K) - f(x^*)\right] \le \left(1 - \min\left\{\frac{\gamma\mu}{2},\frac{\alpha}{4},\frac{1}{4m}\right\}\right)^K\frac{4T^0}{\gamma},
	\end{equation*}	
	when $\mu > 0$ and
	\begin{equation*}
		\EE\left[f(\bar{x}^K) - f(x^*)\right] \le \frac{4T^0}{\gamma K}
	\end{equation*}
	when $\mu=0$, where $T^k \eqdef \|x^k - x^*\|^2 + M_1\gamma^2 \sigma_{1,k}^2 + M_2\gamma^2\sigma_{2,k}^2$.
\end{theorem}
In other words, if $\mu > 0$ and
\begin{equation*}
		\gamma =\frac{3}{L\left(\frac{41}{3}+\frac{52\omega+35}{n}\right)},\quad \alpha = \frac{1}{\omega+1},
\end{equation*}
then {\tt VR-DIANA} converges with the linear rate
\begin{equation*}
	\cO\left(\left(\omega + m + \kappa\left(1+\frac{\omega}{n}\right)\right)\ln\frac{1}{\varepsilon}\right)
\end{equation*}
to the exact solution which coincides with the rate obtained in \cite{horvath2019stochastic}. We notice that the framework from \cite{gorbunov2019unified} establishes slightly worse guarantee:
\begin{equation*}
	\cO\left(\left(\omega + m + \kappa\left(1+\frac{\omega}{n}\right)\frac{\max\{m,\omega+1\}}{m}\right)\ln\frac{1}{\varepsilon}\right)
\end{equation*}
This guarantee is strictly worse than our bound when $m \le 1+\omega$. The key tool that helps us to improve the rate is two sequences of $\{\sigma_{1,k}^2\}_{k\ge 0}$, $\{\sigma_{2,k}^2\}_{k\ge 0}$ instead of one sequence $\{\sigma_{k}^2\}_{k\ge 0}$ as in \cite{gorbunov2019unified}.

Applying Lemma~\ref{lem:lemma_technical_cvx} we get the complexity result in the case when $\mu = 0$.
\begin{corollary}\label{cor:vr-diana_cvx_cor}
	Let the assumptions of Theorem~\ref{thm:vr-diana} hold and $\mu = 0$. Then after $K$ iterations of {\tt VR-DIANA} with the stepsize
	\begin{eqnarray*}
		\gamma_0 &=&\frac{3}{L\left(\frac{41}{3}+\frac{52\omega+35}{n}\right)}\\
		\gamma &=& \min\left\{\gamma_0, \sqrt{\frac{\|x^0 - x^*\|^2}{M_1\sigma_{1,0}^2+ M_2\sigma_{2,0}^2}}\right\},\quad  M_1 = \frac{8(\omega+1)}{3n\alpha},\quad M_2 = \frac{8\omega m}{3n} + \frac{32m}{9}
	\end{eqnarray*}		
	and $\alpha = \frac{1}{\omega+1}$ we have $\EE\left[f(\bar{x}^K) - f(x^*)\right]$ of order
	\begin{equation*}
		\cO\left(\frac{L R_0^2\left(1+\frac{\omega}{n}\right)}{K} + \frac{R_0\sqrt{\frac{(1+\omega)^2}{n}\sigma_{1,0}^2 + \left(1+\frac{\omega}{n}\right)m\sigma_{2,0}^2}}{K}\right)
	\end{equation*}
	where $R_0 = \|x^0 - x^*\|^2$. That is, to achive $\EE\left[f(\bar{x}^K) - f(x^*)\right] \le \varepsilon$ {\tt VR-DIANA} requires
	\begin{equation*}
		\cO\left(\frac{L R_0^2\left(1+\frac{\omega}{n}\right)}{\varepsilon} + \frac{R_0\sqrt{\frac{(1+\omega)^2}{n}\sigma_{1,0}^2 + \left(1+\frac{\omega}{n}\right)m\sigma_{2,0}^2}}{\varepsilon}\right)
	\end{equation*}
	iterations.
\end{corollary}

\section{{Distributed \tt SGD} with Compression and Error Compensation}\label{sec:ec_sgd}
In this section we consider the scenario when compression and error-feedback is applied in order to reduce the communication cost of the method, i.e., we consider {\tt SGD} with error compensation and compression ({\tt EC-SGD}) which has updates of the form \eqref{eq:x^k+1_update}-\eqref{eq:error_update} with
\begin{eqnarray}
	g^k &=& \frac{1}{n}\sum\limits_{i=1}^ng_i^k\notag\\
	v^k &=& \frac{1}{n}\sum\limits_{i=1}^n v_i^k,\quad v_i^k = C(e_i^k + \gamma g_i^k)\label{eq:v^k_def_ec_sgd}\\
	e^k &=& \frac{1}{n}\sum\limits_{i=1}^n e_i^k,\quad e_i^{k+1} = e_i^k + \gamma g_i^k - v_i^k = e_i^k + \gamma g_i^k - C(e_i^k + \gamma g_i^k).\label{eq:e^k_def_ec_sgd}
\end{eqnarray}
Moreover, we assume that $e_i^0 = 0$ for $i = 1,\ldots,n$.

\begin{lemma}\label{lem:ec_sgd_key_lemma_new}
	Let Assumptions~\ref{ass:quasi_strong_convexity}~and~\ref{ass:L_smoothness} be satisfied,  Assumption~\ref{ass:key_assumption_finite_sums_new} holds and\footnote{When $\rho_1 = 1$ and $\rho_2=1$ one can always set the parameters in such a way that $B_1 = \widetilde{B}_1 = B_2 = \widetilde{B}_2 = C_1 = C_2 = 0$, $D_2 = 0$. In this case we assume that $\frac{2}{1-\rho_1}\left(\frac{C_1}{\rho_1}+\frac{2GC_2}{\rho_2(1-\rho_2)}\right)\left(\frac{2B_1}{\delta}+\widetilde{B}_1\right) + \frac{2C_2\left(\frac{2B_2}{\delta}+\widetilde{B}_2\right)}{\rho_2(1-\rho_2)} = 0$.}
	\begin{equation}
		\gamma \le \min\left\{\frac{\delta}{4\mu}, \sqrt{\frac{\delta}{96L\left(\frac{2A}{\delta} + \widetilde{A} + \frac{2}{1-\rho_1}\left(\frac{C_1}{\rho_1}+\frac{2GC_2}{\rho_2(1-\rho_2)}\right)\left(\frac{2B_1}{\delta}+\widetilde{B}_1\right) + \frac{2C_2\left(\frac{2B_2}{\delta}+\widetilde{B}_2\right)}{\rho_2(1-\rho_2)}\right)}}\right\},\label{eq:gamma_condition_ec_sgd_new}
	\end{equation}
	where $ M_1 = \frac{4B_1'}{3\rho_1}$ and $M_2 = \frac{4\left(B_2' + \frac{4}{3}G\right)}{3\rho_2}$. Then {\tt EC-SGD} satisfies Assumption~\ref{ass:key_assumption_new}, i.e., inequality \eqref{eq:sum_of_errors_bound_new} holds with the following parameters:
	\begin{equation}
		F_1 = \frac{24L\gamma^2}{\delta\rho_1(1-\eta)}\left(\frac{2B_1}{\delta}+\widetilde{B}_1\right),\quad F_2 = \frac{24L\gamma^2}{\delta\rho_2(1-\eta)}\left(\frac{2G}{1-\rho_1}\left(\frac{2B_1}{\delta}+\tilde{B}_1\right)+\frac{2B_2}{\delta} + \widetilde{B}_2\right), \label{eq:ec_sgd_parameters_new}
	\end{equation}
	\begin{equation}
		D_3 = \frac{6L\gamma}{\delta}\left(\frac{D_2}{\rho_1}\left(\frac{2B_1}{\delta}+\widetilde{B}_1\right) + \frac{2D_1}{\delta} + \widetilde{D}_1\right).\label{eq:ec_sgd_parameters_new_2}
	\end{equation}
\end{lemma}

That is, Assumption~\ref{ass:key_assumption_finite_sums_new} implies Assumption~\ref{ass:key_assumption_new} in the case of error compensation. As a direct application of Lemma~\ref{lem:ec_sgd_key_lemma_new} and Theorem~\ref{thm:main_result_new} we get the following result.
\begin{theorem}\label{thm:ec_sgd_main_result_new}
	Let Assumptions~\ref{ass:quasi_strong_convexity}~and~\ref{ass:L_smoothness} be satisfied, Assumption~\ref{ass:key_assumption_finite_sums_new} holds and
	\begin{eqnarray*}
		\gamma &\le& \frac{1}{4(A'+C_1M_1+C_2M_2)},\\
		\gamma &\le& \min\left\{\frac{\delta}{4\mu}, \sqrt{\frac{\delta}{96L\left(\frac{2A}{\delta} + \widetilde{A} + \frac{2}{1-\rho_1}\left(\frac{C_1}{\rho_1}+\frac{2GC_2}{\rho_2(1-\rho_2)}\right)\left(\frac{2B_1}{\delta}+\widetilde{B}_1\right) + \frac{2C_2\left(\frac{2B_2}{\delta}+\widetilde{B}_2\right)}{\rho_2(1-\rho_2)}\right)}}\right\},
	\end{eqnarray*}
	where $ M_1 = \frac{4B_1'}{3\rho_1}$ and $M_2 = \frac{4\left(B_2' + \frac{4}{3}G\right)}{3\rho_2}$. Then for all $K\ge 0$ we have
	\begin{equation*}
		\EE\left[f(\bar x^K) - f(x^*)\right] \le \left(1 - \eta\right)^K\frac{4(T^0 + \gamma F_1 \sigma_{1,0}^2+ \gamma F_2 \sigma_{2,0}^2)}{\gamma} + 4\gamma\left(D_1' + M_1D_2 + D_3\right), 
	\end{equation*}
	when $\mu > 0$ and
	\begin{equation*}
		\EE\left[f(\bar x^K) - f(x^*)\right] \le \frac{4(T^0 + \gamma F_1 \sigma_{1,0}^2+ \gamma F_2 \sigma_{2,0}^2)}{\gamma K} + 4\gamma\left(D_1' + M_1D_2 + D_3\right) 
	\end{equation*} 
	when $\mu = 0$, where $\eta = \min\left\{\nicefrac{\gamma\mu}{2},\nicefrac{\rho_1}{4},\nicefrac{\rho_2}{4}\right\}$, $T^k \eqdef \|\tx^k - x^*\|^2 + M_1\gamma^2 \sigma_{1,k}^2 + M_2\gamma^2 \sigma_{2,k}^2$ and 
	\begin{equation*}
		F_1 = \frac{24L\gamma^2}{\delta\rho_1(1-\eta)}\left(\frac{2B_1}{\delta}+\widetilde{B}_1\right),\quad F_2 = \frac{24L\gamma^2}{\delta\rho_2(1-\eta)}\left(\frac{2G}{1-\rho_1}\left(\frac{2B_1}{\delta}+\tilde{B}_1\right)+\frac{2B_2}{\delta} + \widetilde{B}_2\right),
	\end{equation*}
	\begin{equation*}
		D_3 = \frac{6L\gamma}{\delta}\left(\frac{D_2}{\rho_1}\left(\frac{2B_1}{\delta}+\widetilde{B}_1\right) + \frac{2D_1}{\delta} + \widetilde{D}_1\right).
	\end{equation*}
\end{theorem}

\section{Special Cases: Error Compensated Methods}\label{sec:special_cases_ef}

\subsection{{\tt EC-SGDsr}}\label{sec:ec_SGDsr}
In this section we consider the same setup as in Section~\ref{sec:diana_arbitrary_sampling} and assume additionally that $f_1,\ldots,f_n$ are $L$-smooth.
\begin{algorithm}[t]
   \caption{{\tt EC-SGDsr}}\label{alg:ec-SGDsr}
\begin{algorithmic}[1]
   \Require learning rate $\gamma>0$, initial vector $x^0 \in \R^d$
	\State Set $e_i^0 = 0$ for all $i=1,\ldots, n$   
   \For{$k=0,1,\dotsc$}
       \State Broadcast $x^{k}$ to all workers
        \For{$i=1,\dotsc,n$ in parallel}
            \State Sample $g^{k}_i = \nabla f_{\xi_i}(x^k)$
            \State $v_i^k = C(e_i^k + \gamma g_i^k)$
            \State $e_i^{k+1} = e_i^k + \gamma g_i^k - v_i^k$
        \EndFor
        \State $e^k = \frac{1}{n}\sum_{i=1}^ne_i^k$, $g^k = \frac{1}{n}\sum_{i=1}^ng_i^k$, $v^k = \frac{1}{n}\sum_{i=1}^nv_i^k$
       \State $x^{k+1} = x^k - v^k$
   \EndFor
\end{algorithmic}
\end{algorithm}

\begin{lemma}\label{lem:key_lemma_ec-SGDsr}
	For all $k\ge 0$ we have
	\begin{eqnarray*}
		\frac{1}{n}\sum\limits_{i=1}^n\EE\left[\|g_i^k\|^2\mid x^k\right] &\le& 4L\left(f(x^k) - f(x^*)\right) + \frac{2}{n}\sum\limits_{i=1}^n\|\nabla f_{i}(x^*)\|^2, \notag\\
		\frac{1}{n}\sum\limits_{i=1}^n\EE\left[\|g_i^k-\bar{g}_i^k\|^2\mid x^k\right] &\le& 6(\cL + L)\left(f(x^k)-f(x^*)\right) + \frac{3}{n}\sum\limits_{i=1}^n\EE_{\cD}\left[\|\nabla f_{\xi_i}(x^*)-\nabla f_{i}(x^*)\|^2\right],\\
		\EE\left[\|g^k\|^2\mid x^k\right] &\le& 4\cL\left(f(x^k) - f(x^*)\right) + \frac{2}{n^2}\sum\limits_{i=1}^n\EE_{\cD}\left[\|\nabla f_{\xi_i}(x^*) - \nabla f_i(x^*)\|^2\right]. \notag
	\end{eqnarray*}
\end{lemma}
\begin{proof}
	Applying straightforward inequality $\|a+b\|^2 \le 2\|a\|^2 + 2\|b\|^2$ for $a,b\in\R^d$ we get
	\begin{eqnarray}
		\frac{1}{n}\sum\limits_{i=1}^n\|\bar{g}_i^k\|^2 &=& \frac{1}{n}\sum\limits_{i=1}^n\|\nabla f_i(x^k)-\nabla f_i(x^*)+\nabla f_i(x^*)\|^2\notag\\
		&\overset{\eqref{eq:a_b_norm_squared}}{\le}& \frac{1}{n}\sum\limits_{i=1}^n\|\nabla f_{i}(x^k) - \nabla f_{i}(x^*)\|^2 + \frac{2}{n}\sum\limits_{i=1}^n\|\nabla f_{i}(x^*)\|^2\notag\\
		&\overset{\eqref{eq:L_smoothness_cor_3}}{\le}& 4L\left(f(x^k) - f(x^*)\right) + \frac{2}{n}\sum\limits_{i=1}^n\|\nabla f_{i}(x^*)\|^2.\label{eq:ec_useful_technical_stuff}
	\end{eqnarray}
	Similarly we obtain
	\begin{eqnarray*}
		\frac{1}{n}\sum\limits_{i=1}^n\EE\left[\|g_i^k-\bar{g}_i^k\|^2\mid x^k\right] &=& \frac{1}{n}\sum\limits_{i=1}^n\EE_{\cD}\left[\|\nabla f_{\xi_i}(x^k)  - \nabla f_i(x^k)\|^2\right]\\
		&\overset{\eqref{eq:a_b_norm_squared}}{\le}& \frac{3}{n}\sum\limits_{i=1}^n\EE_{\cD}\left[\|\nabla f_{\xi_i}(x^k)-\nabla f_{\xi_i}(x^*)\|^2\right]\\
		&&\quad + \frac{3}{n}\sum\limits_{i=1}^n\EE_{\cD}\left[\|\nabla f_{\xi_i}(x^*)-\nabla f_{i}(x^*)\|^2\right] \\
		&&\quad + \frac{3}{n}\sum\limits_{i=1}^n\|\nabla f_i(x^*) - \nabla f_i(x^k)\|^2\\
		&\overset{\eqref{eq:L_smoothness_cor_3},\eqref{eq:exp_smoothness}}{\le}& 6(\cL + L)\left(f(x^k)-f(x^*)\right)\\
		&&\quad + \frac{3}{n}\sum\limits_{i=1}^n\EE_{\cD}\left[\|\nabla f_{\xi_i}(x^*)-\nabla f_{i}(x^*)\|^2\right].
	\end{eqnarray*}
	Next, using the independence of $\xi_1^k,\ldots, \xi_n^k$ we derive
	\begin{eqnarray*}
		\EE\left[\left\|g^k\right\|^2\mid x^k\right] &=& \EE\left[\left\|\frac{1}{n}\sum\limits_{i=1}^n\left(\nabla f_{\xi_i^k}(x^k) - \nabla f_{\xi_i^k}(x^*) + \nabla f_{\xi_i^k}(x^*) - \nabla f_i(x^*)\right)\right\|^2\mid x^k\right]\\
		&\overset{\eqref{eq:a_b_norm_squared}}{\le}& \frac{2}{n}\sum\limits_{i=1}^n\EE\left[\left\|\nabla f_{\xi_i^k}(x^k) - \nabla f_{\xi_i^k}(x^*)\right\|^2\mid x^k\right]\\
		&&\quad + 2\EE\left[\left\|\frac{1}{n}\sum\limits_{i=1}^n\left(\nabla f_{\xi_i^k}(x^*) - \nabla f_i(x^*)\right)\right\|^2\mid x^k\right]\\
		&\overset{\eqref{eq:exp_smoothness}}{\le}& 4\cL\left(f(x^k) - f(x^*)\right) + \frac{2}{n^2}\sum\limits_{i=1}^n\EE_{\cD_i}\left[\left\|\nabla f_{\xi_i}(x^*) - \nabla f_i(x^*)\right\|^2\right].
	\end{eqnarray*}
\end{proof}

Applying Theorem~\ref{thm:ec_sgd_main_result_new} we get the following result.
\begin{theorem}\label{thm:ec_SGDsr}
	Assume that $f(x)$ is $\mu$-quasi strongly convex, $f_1,\ldots,f_n$ are $L$-smooth and Assumption~\ref{ass:exp_smoothness} holds. Then {\tt EC-SGDsr} satisfies Assumption~\ref{ass:key_assumption_finite_sums_new} with
	\begin{gather*}
		A = 2L,\quad \widetilde{A} = 3(\cL+L),\quad A' = 2\cL,\quad B_1 = \widetilde{B}_1 = B_1' = B_2 = \widetilde{B}_2 = B_2' = 0,\\
		D_1 = \frac{2}{n}\sum\limits_{i=1}^n\|\nabla f_{i}(x^*)\|^2,\quad \widetilde{D}_1 = \frac{3}{n}\sum\limits_{i=1}^n\EE_{\cD}\left[\|\nabla f_{\xi_i}(x^*) - \nabla f_i(x^*)\|^2\right],\quad \sigma_{1,k}^2 \equiv \sigma_{2,k}^2 \equiv 0,\\
		D_1' = \frac{2}{n^2}\sum\limits_{i=1}^n\EE_{\cD}\left[\|\nabla f_{\xi_i}(x^*) - \nabla f_i(x^*)\|^2\right],\quad \rho_1 = \rho_2 = 1,\quad C_1 = C_2 = 0,\quad G = 0,\quad D_2 = 0,\\
		F_1 = F_2 = 0,\quad D_3 = \frac{6L\gamma}{\delta}\left(\frac{2D_1}{\delta}+\widetilde{D}_1\right),
	\end{gather*}
	with $\gamma$ satisfying
	\begin{equation*}
		\gamma \le \min\left\{\frac{1}{8\cL},\frac{\delta}{4\sqrt{6L\left(4L+3\delta(\cL+L)\right)}}\right\}
	\end{equation*}
	and for all $K \ge 0$
	\begin{equation*}
		\EE\left[f(\bar{x}^K) - f(x^*)\right] \le \left(1 - \frac{\gamma\mu}{2}\right)^K\frac{4\|x^0 - x^*\|^2}{\gamma} + 4\gamma\left(D_1' + \frac{12L\gamma}{\delta^2}D_1 + \frac{6L\gamma}{\delta}\widetilde{D}_1\right)
	\end{equation*}
	when $\mu > 0$ and
	\begin{equation*}
		\EE\left[f(\bar{x}^K) - f(x^*)\right] \le \frac{4\|x^0 - x^*\|^2}{K\gamma} + 4\gamma\left(D_1' + \frac{12L\gamma}{\delta^2}D_1 + \frac{6L\gamma}{\delta}\widetilde{D}_1\right)
	\end{equation*}
	when $\mu=0$.
\end{theorem}
In other words, {\tt EC-SGDsr} converges with linear rate $\cO\left(\left(\frac{\cL}{\mu} + \frac{L + \sqrt{\delta L\cL}}{\mu\delta}\right)\ln\frac{1}{\varepsilon}\right)$ to the neighbourhood of the solution when $\mu > 0$. Applying Lemma~\ref{lem:lemma2_stich} we establish the rate of convergence to $\varepsilon$-solution.
\begin{corollary}\label{cor:ec_SGDsr_str_cvx_cor}
	Let the assumptions of Theorem~\ref{thm:ec_SGDsr} hold and $\mu > 0$. Then after $K$ iterations of {\tt EC-SGDsr} with the stepsize
	\begin{equation*}
		\gamma = \min\left\{\frac{1}{8\cL},\frac{\delta}{4\sqrt{6L\left(4L+3\delta(\cL+L)\right)}}, \frac{\ln\left(\max\left\{2,\min\left\{\frac{\|x^0-x^*\|^2\mu^2K^2}{D_1'}, \frac{\delta\|x^0-x^*\|^2\mu^3K^3}{6L(\nicefrac{2D_1}{\delta}+\widetilde{D}_1)}\right\}\right\}\right)}{\mu K}\right\}
	\end{equation*}	 
	we have $\EE\left[f(\bar{x}^K) - f(x^*)\right]$ of order
	\begin{equation*}
		\widetilde\cO\left(\left(\cL + \frac{L+\sqrt{\delta L\cL}}{\delta}\right)\|x^0 - x^*\|^2\exp\left(-\frac{\mu}{\cL+\frac{L+\sqrt{\delta L\cL}}{\delta}}K\right) + \frac{D_1'}{\mu K} + \frac{L(\widetilde{D}_1 + \nicefrac{D_1}{\delta})}{\delta\mu^2 K^2}\right).
	\end{equation*}
	That is, to achive $\EE\left[f(\bar{x}^K) - f(x^*)\right] \le \varepsilon$ {\tt EC-SGDsr} requires
	\begin{equation*}
		\widetilde{\cO}\left(\frac{\cL}{\mu} + \frac{L+\sqrt{\delta L\cL}}{\delta\mu} + \frac{D_1'}{\mu\varepsilon} + \frac{\sqrt{L(\widetilde{D}_1 + \nicefrac{D_1}{\delta})}}{\mu\sqrt{\delta\varepsilon}}\right) \text{ iterations.}
	\end{equation*}
\end{corollary}

Applying Lemma~\ref{lem:lemma_technical_cvx} we get the complexity result in the case when $\mu = 0$.
\begin{corollary}\label{cor:ec_sgdsr_cvx_cor}
	Let the assumptions of Theorem~\ref{thm:ec_SGDsr} hold and $\mu = 0$. Then after $K$ iterations of {\tt EC-SGDsr} with the stepsize
	\begin{eqnarray*}
		\gamma_0 &=&\min\left\{\frac{1}{8\cL},\frac{\delta}{4\sqrt{6L\left(4L+3\delta(\cL+L)\right)}}\right\}\\
		\gamma &=& \min\left\{\gamma_0, \sqrt{\frac{\|x^0 - x^*\|^2}{D_1' K}}, \sqrt[3]{\frac{\|x^0 - x^*\|^2\delta}{6L(\nicefrac{2D_1}{\delta}+\widetilde{D}_1) K}}\right\}	\end{eqnarray*}		
	we have $\EE\left[f(\bar{x}^K) - f(x^*)\right]$ of order
	\begin{equation*}
		\cO\left(\frac{R_0^2\left(\cL+\frac{L+\sqrt{\delta L\cL}}{\delta}\right)}{K} + \sqrt{\frac{R_0^2 D_1'}{K}} + \frac{\sqrt[3]{LR_0^4(\nicefrac{2D_1}{\delta}+\widetilde{D}_1)}}{\left(\delta K^2\right)^{\nicefrac{1}{3}}}\right)
	\end{equation*}
	where $R_0 = \|x^0 - x^*\|^2$. That is, to achive $\EE\left[f(\bar{x}^K) - f(x^*)\right] \le \varepsilon$ {\tt EC-SGDsr} requires
	\begin{equation*}
		\cO\left(\frac{R_0^2\left(\cL+\frac{L+\sqrt{\delta L\cL}}{\delta}\right)}{\varepsilon} + \frac{R_0^2 D_1'}{\varepsilon^2} + \frac{R_0^2\sqrt{L(\nicefrac{2D_1}{\delta}+\widetilde{D}_1)}}{\sqrt{\delta \varepsilon^3}}\right)
	\end{equation*}
	iterations.
\end{corollary}

\subsection{{\tt EC-SGD}}\label{sec:ec_sgd_pure}
In this section we consider problem \eqref{eq:main_problem_ef} with $f_i(x)$ satisfying \eqref{eq:f_i_expectation_ef} where functions $f_{\xi_i}(x)$ are differentiable and $L$-smooth almost surely in $\xi_i$, $i=1,\ldots,n$.
\begin{algorithm}[t]
   \caption{{\tt EC-SGD}}\label{alg:ec-sgd}
\begin{algorithmic}[1]
   \Require learning rate $\gamma>0$, initial vector $x^0 \in \R^d$
	\State Set $e_i^0 = 0$ for all $i=1,\ldots, n$   
   \For{$k=0,1,\dotsc$}
       \State Broadcast $x^{k}$ to all workers
        \For{$i=1,\dotsc,n$ in parallel}
            \State Sample $g^{k}_i = \nabla f_{\xi_i}(x^k)$ independently from other workers
            \State $v_i^k = C(e_i^k + \gamma g_i^k)$
            \State $e_i^{k+1} = e_i^k + \gamma g_i^k - v_i^k$
        \EndFor
        \State $e^k = \frac{1}{n}\sum_{i=1}^ne_i^k$, $g^k = \frac{1}{n}\sum_{i=1}^ng_i^k$, $v^k = \frac{1}{n}\sum_{i=1}^nv_i^k$
       \State $x^{k+1} = x^k - v^k$
   \EndFor
\end{algorithmic}
\end{algorithm}

\begin{lemma}[See also Lemmas 1,2 from~\cite{nguyen2018sgd}]\label{lem:lemma_sgd}
    Assume that $f_{\xi_i}(x)$ are convex in $x$ for every $\xi_i$, $i=1,\ldots,n$. Then for every $x\in\R^d$ and $i=1,\ldots, n$
	\begin{eqnarray*}
		\frac{1}{n}\sum\limits_{i=1}^n\|\nabla f_{i}(x)\|^2 &\le& 4L\left(f(x) - f(x^*)\right) + \frac{2}{n}\sum\limits_{i=1}^n\|\nabla f_{i}(x^*)\|^2,\\
		\frac{1}{n}\sum\limits_{i=1}^n\EE_{\xi_i\sim\cD_i}{\|\nabla f_{\xi_i}(x)-\nabla f_i(x)\|^2} &\le& 12L\left(f(x) - f(x^*)\right) + \frac{3}{n}\sum\limits_{i=1}^n\EE\left[\|\nabla f_{\xi_i}(x^*)-\nabla f_{i}(x^*)\|^2\right],\\
		\EE_{\xi_1,\ldots,\xi_n}{\left\|\frac{1}{n}\sum\limits_{i=1}^n\nabla f_{\xi_i}(x)\right\|^2} &\le& 4L\left(f(x) - f(x^*)\right) + \frac{2}{n^2}\sum\EE\left[\|\nabla f_{\xi_i}(x^*)-\nabla f_i(x^*)\|^2\right].
	\end{eqnarray*}	    
    If further $f(x)$ is $\mu$-strongly convex with $\mu > 0$ and possibly non-convex $f_i,f_{\xi_i}$, then for every $x\in\R^d$ and $i = 1,\ldots, n$
    \begin{eqnarray*}
    	\frac{1}{n}\sum\limits_{i=1}^n\|\nabla f_{i}(x)\|^2 &\le& 4L\kappa\left(f(x) - f(x^*)\right) + \frac{2}{n}\sum\limits_{i=1}^n\|\nabla f_{i}(x^*)\|^2,\\
		\frac{1}{n}\sum\limits_{i=1}^n\EE_{\xi_i\sim\cD_i}{\|\nabla f_{\xi_i}(x)-\nabla f_i(x)\|^2} &\le& 12L\kappa\left(f(x) - f(x^*)\right) \\
		&&\quad + \frac{3}{n}\sum\limits_{i=1}^n\EE\left[\|\nabla f_{\xi_i}(x^*)-\nabla f_{i}(x^*)\|^2\right],\\
    	\EE_{\xi_1,\ldots,\xi_n}{\left\|\frac{1}{n}\sum\limits_{i=1}^n\nabla f_{\xi_i}(x)\right\|^2} &\le& 4L\kappa\left(f(x) - f(x^*)\right)\\
    	&&\quad + \frac{2}{n^2}\sum\limits_{i=1}^n\EE\left[\|\nabla f_{\xi_i}(x^*)-\nabla f_i(x^*)\|^2\right].
    \end{eqnarray*}
    where $\kappa = \frac{L}{\mu}$.
\end{lemma}
\begin{proof}
	We start with the case when functions $f_{\xi_i}(x)$ are convex in $x$ for every $\xi_i$. The first inequality follows from \eqref{eq:ec_useful_technical_stuff}. Next, we derive
	\begin{eqnarray*}
		\frac{1}{n}\sum\limits_{i=1}^n\EE_{\xi_i\sim\cD_i}{\|\nabla f_{\xi_i}(x)-\nabla f_i(x)\|^2} &\overset{\eqref{eq:a_b_norm_squared}}{\le}& \frac{3}{n} \sum\limits_{i=1}^n\EE_{\xi_i\sim\cD_i}{\|\nabla f_{\xi_i}(x) - \nabla f_{\xi_i}(x^*)\|^2} \\
		&&\quad + \frac{3}{n}\sum\limits_{i=1}^n\EE_{\xi_i\sim\cD_i}{\|\nabla f_{\xi_i}(x^*) - \nabla f_i(x^*)\|^2}\\
		&&\quad + \frac{3}{n}\sum\limits_{i=1}^n{\|\nabla f_{i}(x^*) - \nabla f_i(x)\|^2}\\
		&\overset{\eqref{eq:L_smoothness_cor_3}}{\le}& 12L\left(f(x)-f(x^*)\right) + \frac{3}{n}\sum\limits_{i=1}^n\EE{\|\nabla f_{\xi_i}(x^*)\|^2}.
	\end{eqnarray*}
	Due to independence of $\xi_1^k,\ldots,\xi_n^k$ we get
	\begin{eqnarray*}
		\EE_{\xi_1,\ldots,\xi_n}{\left\|\frac{1}{n}\sum\limits_{i=1}^n\nabla f_{\xi_i}(x)\right\|^2} &=& \EE_{\xi_1,\ldots,\xi_n}{\left\|\frac{1}{n}\sum\limits_{i=1}^n\left(\nabla f_{\xi_i}(x)-\nabla f_{\xi_i}(x^*) + \nabla f_{\xi_i}(x^*) - \nabla f_i(x^*)\right)\right\|^2}\\
		&\overset{\eqref{eq:a_b_norm_squared}}{\le}& \frac{2}{n}\sum\limits_{i=1}^n\EE_{\xi_i\sim\cD_i}\left[\|\nabla f_{\xi_i}(x) - \nabla f_{\xi_i}(x^*)\|^2\right]\\
		&&\quad + 2\EE_{\xi_1,\ldots,\xi_n}{\left\|\frac{1}{n}\sum\limits_{i=1}^n\left(\nabla f_{\xi_i}(x^*) - \nabla f_i(x^*)\right)\right\|^2}\\
		&\overset{\eqref{eq:L_smoothness_cor_3}}{\le}& 4L\left(f(x)-f(x^*)\right) + \frac{2}{n^2}\sum\EE\left[\|\nabla f_{\xi_i}(x^*)-\nabla f_i(x^*)\|^2\right].
	\end{eqnarray*}
	
	Next, we consider the second case: $f(x)$ is $\mu$-strongly convex with possibly non-convex $f_i,f_{\xi_i}$. In this case
	\begin{eqnarray*}
		\frac{1}{n}\sum\limits_{i=1}^n{\|\nabla f_{i}(x)\|^2} &\overset{\eqref{eq:a_b_norm_squared}}{\le}& \frac{2}{n} \sum\limits_{i=1}^n{\|\nabla f_{i}(x) - \nabla f_{i}(x^*)\|^2} + \frac{2}{n}\sum\limits_{i=1}^n{\|\nabla f_{i}(x^*)\|^2}\\
		&\overset{\eqref{eq:L_smoothness_def}}{\le}& 2L^2\|x - x^*\|^2 + \frac{2}{n}\sum\limits_{i=1}^n{\|\nabla f_{i}(x^*)\|^2}\\
		&\le& \frac{4L^2}{\mu}\left(f(x)-f(x^*)\right) + \frac{2}{n}\sum\limits_{i=1}^n{\|\nabla f_{i}(x^*)\|^2}
	\end{eqnarray*}
	where the last inequality follows from $\mu$-strong convexity of $f$. Similarly, we get
	\begin{eqnarray*}
		\frac{1}{n}\sum\limits_{i=1}^n\EE_{\xi_i\sim\cD_i}\left[{\|\nabla f_{\xi_i}(x)-\nabla f_{i}(x)\|^2}\right] &\overset{\eqref{eq:a_b_norm_squared}}{\le}& \frac{3}{n} \sum\limits_{i=1}^n\EE_{\xi_i\sim\cD_i}\left[{\|\nabla f_{\xi_i}(x)-\nabla f_{\xi_i}(x^*)\|^2}\right]\\
		&&\quad + \frac{3}{n} \sum\limits_{i=1}^n\EE_{\xi_i\sim\cD_i}\left[{\|\nabla f_{\xi_i}(x^*)-\nabla f_{i}(x^*)\|^2}\right]\\
		&&\quad + \frac{3}{n} \sum\limits_{i=1}^n{\|\nabla f_{i}(x^*)-\nabla f_{i}(x)\|^2}\\
		&\overset{\eqref{eq:L_smoothness_def}}{\le}& 6L^2\|x - x^*\|^2 \\
		&&\quad + \frac{3}{n} \sum\limits_{i=1}^n\EE_{\xi_i\sim\cD_i}\left[{\|\nabla f_{\xi_i}(x^*)-\nabla f_{i}(x^*)\|^2}\right]\\
		&\le& \frac{12L^2}{\mu}\left(f(x)-f(x^*)\right)\\
		&&\quad + \frac{3}{n} \sum\limits_{i=1}^n\EE_{\xi_i\sim\cD_i}\left[{\|\nabla f_{\xi_i}(x^*)-\nabla f_{i}(x^*)\|^2}\right].
	\end{eqnarray*}
	Finally, using independence of $\xi_1^k,\ldots,\xi_n^k$ we derive
	\begin{eqnarray*}
		\EE_{\xi_1,\ldots,\xi_n}{\left\|\frac{1}{n}\sum\limits_{i=1}^n\nabla f_{\xi_i}(x)\right\|^2} &=& \EE_{\xi_1,\ldots,\xi_n}{\left\|\frac{1}{n}\sum\limits_{i=1}^n\left(\nabla f_{\xi_i}(x)-\nabla f_{\xi_i}(x^*) + \nabla f_{\xi_i}(x^*) - \nabla f_i(x^*)\right)\right\|^2}\\
		&\overset{\eqref{eq:a_b_norm_squared}}{\le}& \frac{2}{n}\sum\limits_{i=1}^n\EE_{\xi_i\sim\cD_i}\left[\|\nabla f_{\xi_i}(x) - \nabla f_{\xi_i}(x^*)\|^2\right]\\
		&&\quad + 2\EE_{\xi_1,\ldots,\xi_n}{\left\|\frac{1}{n}\sum\limits_{i=1}^n\left(\nabla f_{\xi_i}(x^*) - \nabla f_i(x^*)\right)\right\|^2}\\
		&\overset{\eqref{eq:L_smoothness_def}}{\le}& 2L^2\|x-x^*\|^2 + \frac{2}{n^2}\sum\EE\left[\|\nabla f_{\xi_i}(x^*)-\nabla f_i(x^*)\|^2\right]\\
		&\le& \frac{4L^2}{\mu}\left(f(x)-f(x^*)\right) + \frac{2}{n^2}\sum\EE\left[\|\nabla f_{\xi_i}(x^*)-\nabla f_i(x^*)\|^2\right].
	\end{eqnarray*}
\end{proof}

Applying Theorem~\ref{thm:ec_sgd_main_result_new} we get the following result.
\begin{theorem}\label{thm:ec_sgd_pure}
	Assume that $f_\xi(x)$ is convex and $L$-smooth in $x$ for every $\xi$ and $f(x)$ is $\mu$-quasi strongly convex. Then {\tt EC-SGD} satisfies Assumption~\ref{ass:key_assumption_finite_sums_new} with
	\begin{gather*}
		A = A' = 2L,\quad \widetilde{A} = 6L,\quad B_1 = \widetilde{B}_1 = B_1' = B_2 = \widetilde{B}_2 = B_2' = 0,\\
		D_1 = \frac{2}{n}\sum\limits_{i=1}^n\|\nabla f_{i}(x^*)\|^2,\quad \widetilde{D}_1 = \frac{2}{n}\sum\limits_{i=1}^n\EE\left[\|\nabla f_{\xi_i}(x^*) - \nabla f_i(x^*)\|^2\right],\quad \sigma_{1,k}^2 \equiv \sigma_{2,k}^2 \equiv 0,\\
		D_1' = \frac{2}{n^2}\sum\limits_{i=1}^n\EE\left[\|\nabla f_{\xi_i}(x^*) - \nabla f_i(x^*)\|^2\right],\quad \rho_1 = \rho_2 = 1,\quad C_1 = C_2 = 0,\quad G = 0,\quad D_2 = 0,\\
		F_1 = F_2 = 0,\quad D_3 = \frac{6L\gamma}{\delta}\left(\frac{2D_1}{\delta}+\widetilde{D}_1\right),
	\end{gather*}
	with $\gamma$ satisfying
	\begin{equation*}
		\gamma \le \frac{\delta}{8L\sqrt{6+9\delta}}
	\end{equation*}
	and for all $K \ge 0$
	\begin{equation*}
		\EE\left[f(\bar{x}^K) - f(x^*)\right] \le \left(1 - \frac{\gamma\mu}{2}\right)^K\frac{4\|x^0 - x^*\|^2}{\gamma} + 4\gamma\left(D_1' + \frac{12L\gamma}{\delta^2}D_1 + \frac{6L\gamma}{\delta}\widetilde{D}_1\right)
	\end{equation*}
	when $\mu > 0$ and
	\begin{equation*}
		\EE\left[f(\bar{x}^K) - f(x^*)\right] \le \frac{4\|x^0 - x^*\|^2}{K\gamma} + 4\gamma\left(D_1' + \frac{12L\gamma}{\delta^2}D_1 + \frac{6L\gamma}{\delta}\widetilde{D}_1\right)
	\end{equation*}
	when $\mu=0$. If further $f(x)$ is $\mu$-strongly convex with $\mu > 0$ and possibly non-convex $f_i,f_{\xi_i}$, then
	{\tt EC-SGD} satisfies Assumption~\ref{ass:key_assumption_finite_sums_new} with
	\begin{gather*}
		A = A' = 2L\kappa,\quad \widetilde{A} = 6L\kappa,\quad B_1 = \widetilde{B}_1 = B_1' = B_2 = \widetilde{B}_2 = B_2' = 0,\\
		D_1 = \frac{2}{n}\sum\limits_{i=1}^n\|\nabla f_{i}(x^*)\|^2,\quad \widetilde{D}_1 = \frac{2}{n}\sum\limits_{i=1}^n\EE\left[\|\nabla f_{\xi_i}(x^*) - \nabla f_i(x^*)\|^2\right],\quad \sigma_{1,k}^2 \equiv \sigma_{2,k}^2 \equiv 0,\\
		D_1' = \frac{2}{n^2}\sum\limits_{i=1}^n\EE\left[\|\nabla f_{\xi_i}(x^*) - \nabla f_i(x^*)\|^2\right],\quad \rho_1 = \rho_2 = 1,\quad C_1 = C_2 = 0,\quad G = 0,\quad D_2 = 0,\\
		F_1 = F_2 = 0,\quad D_3 = \frac{6L\gamma}{\delta}\left(\frac{2D_1}{\delta}+\widetilde{D}_1\right),
	\end{gather*}
	with $\gamma$ satisfying
	\begin{equation*}
		\gamma \le \min\left\{\frac{1}{8\kappa L},\frac{\delta}{8L\sqrt{3\kappa(2+3\delta)}}\right\}
	\end{equation*}
	and for all $K \ge 0$
	\begin{equation*}
		\EE\left[f(\bar{x}^K) - f(x^*)\right] \le \left(1 - \frac{\gamma\mu}{2}\right)^K\frac{4\|x^0 - x^*\|^2}{\gamma} + 4\gamma\left(D_1' + \frac{12L\gamma}{\delta^2}D_1 + \frac{6L\gamma}{\delta}\widetilde{D}_1\right).
	\end{equation*}
\end{theorem}

In other words, {\tt EC-SGD} converges with linear rate $\cO\left(\frac{\kappa}{\delta}\ln\frac{1}{\varepsilon}\right)$ to the neighbourhood of the solution when $f_{\xi}(x)$ are convex for each $\xi$ and $\mu > 0$. Applying Lemma~\ref{lem:lemma2_stich} we establish the rate of convergence to $\varepsilon$-solution.
\begin{corollary}\label{cor:ec_SGD_pure_str_cvx_cor}
	Let the assumptions of Theorem~\ref{thm:ec_sgd_pure} hold, $f_{\xi}(x)$ are convex for each $\xi$ and $\mu > 0$. Then after $K$ iterations of {\tt EC-SGD} with the stepsize
	\begin{equation*}
		\gamma = \min\left\{\frac{\delta}{8L\sqrt{6+9\delta}}, \frac{\ln\left(\max\left\{2,\min\left\{\frac{\|x^0-x^*\|^2\mu^2K^2}{D_1'}, \frac{\delta\|x^0-x^*\|^2\mu^3K^3}{6L(\nicefrac{2D_1}{\delta}+\widetilde{D}_1)}\right\}\right\}\right)}{\mu K}\right\}
	\end{equation*}	 
	we have
	\begin{equation*}
		\EE\left[f(\bar{x}^K) - f(x^*)\right] = \widetilde\cO\left(\frac{L}{\delta}\|x^0 - x^*\|^2\exp\left(-\frac{\delta\mu}{L}K\right) + \frac{D_1'}{\mu K} + \frac{L(\widetilde{D}_1 + \nicefrac{D_1}{\delta})}{\delta\mu^2 K^2}\right).
	\end{equation*}
	That is, to achive $\EE\left[f(\bar{x}^K) - f(x^*)\right] \le \varepsilon$ {\tt EC-SGD} requires
	\begin{equation*}
		\widetilde{\cO}\left(\frac{L}{\delta\mu} + \frac{D_1'}{\mu\varepsilon} + \frac{\sqrt{L(\widetilde{D}_1 + \nicefrac{D_1}{\delta})}}{\mu\sqrt{\delta\varepsilon}}\right) \text{ iterations.}
	\end{equation*}
\end{corollary}
\begin{corollary}\label{cor:ec_SGD_pure_str_cvx_cor_2}
	Let the assumptions of Theorem~\ref{thm:ec_sgd_pure} hold and $f(x)$ is $\mu$-strongly convex with $\mu > 0$ and possibly non-convex $f_i,f_{\xi_i}$. Then after $K$ iterations of {\tt EC-SGD} with the stepsize
	\begin{equation*}
		\gamma = \min\left\{\frac{1}{8\kappa L}, \frac{\delta}{8L\sqrt{3\kappa(2+3\delta)}}, \frac{\ln\left(\max\left\{2,\min\left\{\frac{\|x^0-x^*\|^2\mu^2K^2}{D_1'}, \frac{\delta\|x^0-x^*\|^2\mu^3K^3}{6L(\nicefrac{2D_1}{\delta}+\widetilde{D}_1)}\right\}\right\}\right)}{\mu K}\right\}
	\end{equation*}	 
	we have $\EE\left[f(\bar{x}^K) - f(x^*)\right]$ of order
	\begin{equation*}
		\widetilde\cO\left(\left(L\kappa+\frac{L\sqrt{\kappa}}{\delta}\right)\|x^0 - x^*\|^2\exp\left(-\min\left\{\frac{\delta\mu}{L\sqrt{\kappa}}, \frac{1}{\kappa^2}\right\}K\right) + \frac{D_1'}{\mu K} + \frac{L(\widetilde{D}_1 + \nicefrac{D_1}{\delta})}{\delta\mu^2 K^2}\right).
	\end{equation*}
	That is, to achive $\EE\left[f(\bar{x}^K) - f(x^*)\right] \le \varepsilon$ {\tt EC-SGD} requires
	\begin{equation*}
		\widetilde{\cO}\left(\kappa^2 + \frac{\kappa^{\nicefrac{3}{2}}}{\delta} + \frac{D_1'}{\mu\varepsilon} + \frac{\sqrt{L(\widetilde{D}_1 + \nicefrac{D_1}{\delta})}}{\mu\sqrt{\delta\varepsilon}}\right) \text{ iterations.}
	\end{equation*}
\end{corollary}

Applying Lemma~\ref{lem:lemma_technical_cvx} we get the complexity result in the case when $\mu = 0$.
\begin{corollary}\label{cor:ec_sgd_cvx_cor}
	Let the assumptions of Theorem~\ref{thm:ec_sgd_pure} hold, $f_{\xi}(x)$ are convex for each $\xi$ and $\mu = 0$. Then after $K$ iterations of {\tt EC-SGD} with the stepsize
	\begin{eqnarray*}
		\gamma &=& \min\left\{\frac{\delta}{8L\sqrt{6+9\delta}}, \sqrt{\frac{\|x^0 - x^*\|^2}{D_1' K}}, \sqrt[3]{\frac{\|x^0 - x^*\|^2\delta}{6L(\nicefrac{2D_1}{\delta}+\widetilde{D}_1) K}}\right\}	\end{eqnarray*}		
	we have $\EE\left[f(\bar{x}^K) - f(x^*)\right]$ of order
	\begin{equation*}
		\cO\left(\frac{LR_0^2}{\delta K} + \sqrt{\frac{R_0^2 D_1'}{K}} + \frac{\sqrt[3]{LR_0^4(\nicefrac{2D_1}{\delta}+\widetilde{D}_1)}}{\left(\delta K^2\right)^{\nicefrac{1}{3}}}\right)
	\end{equation*}
	where $R_0 = \|x^0 - x^*\|^2$. That is, to achive $\EE\left[f(\bar{x}^K) - f(x^*)\right] \le \varepsilon$ {\tt EC-SGD} requires
	\begin{equation*}
		\cO\left(\frac{LR_0^2}{\delta\varepsilon} + \frac{R_0^2 D_1'}{\varepsilon^2} + \frac{R_0^2\sqrt{L(\nicefrac{2D_1}{\delta}+\widetilde{D}_1)}}{\sqrt{\delta \varepsilon^3}}\right)
	\end{equation*}
	iterations.
\end{corollary}

\subsection{{\tt EC-GDstar}}\label{sec:ec_SGDstar}
We assume that $i$-th node has access to the gradient of $f_i$ at the optimality, i.e., to the $\nabla f_i(x^*)$. It is unrealistic scenario but it gives some insights that we will use next in order to design the method that converges asymptotically to \textit{the exact solution}.
\begin{algorithm}[t]
   \caption{{\tt EC-GDstar} (see also \cite{gorbunov2019unified})}\label{alg:EC-GDstar}
\begin{algorithmic}[1]
   \Require learning rate $\gamma>0$, initial vector $x^0 \in \R^d$
	\State Set $e_i^0 = 0$ for all $i=1,\ldots, n$   
   \For{$k=0,1,\dotsc$}
       \State Broadcast $x^{k}$ to all workers
        \For{$i=1,\dotsc,n$ in parallel}
            \State $g^{k}_i = \nabla f_{i}(x^k) - \nabla f_i(x^*)$
            \State $v_i^k = C(e_i^k + \gamma g_i^k)$
            \State $e_i^{k+1} = e_i^k + \gamma g_i^k - v_i^k$
        \EndFor
        \State $e^k = \frac{1}{n}\sum_{i=1}^ne_i^k$, $g^k = \frac{1}{n}\sum_{i=1}^ng_i^k$, $v^k = \frac{1}{n}\sum_{i=1}^nv_i^k$
       \State $x^{k+1} = x^k - v^k$
   \EndFor
\end{algorithmic}
\end{algorithm}

Assume that $f(x)$ is $\mu$-quasi strongly convex and each $f_i$ is convex and $L$-smooth. By definition of $g_i^k$ it trivially follows that
\begin{equation*}
	g^k = \frac{1}{n}\sum\limits_{i=1}^n g_i^k = \frac{1}{n}\sum\limits_{i=1}^n \left(\nabla f_i(x^k) - \nabla f_i(x^*)\right) = \nabla f(x^k) - \nabla f(x^*) = \nabla f(x^k),
\end{equation*}
$g_i^k = \bar{g}_i^k$, and
\begin{eqnarray*}
	\frac{1}{n}\sum\limits_{i=1}^n\|g_i^k\|^2 &=& \frac{1}{n}\sum\limits_{i=1}^n\|\nabla f_i(x^k) - \nabla f_i(x^*)\|^2\\
	&\overset{\eqref{eq:L_smoothness_cor_3}}{\le}& \frac{2L}{n}\sum\limits_{i=1}^n\left(f_i(x^k) - f_i(x^*) - \langle\nabla f_i(x^*), x^k - x^*\rangle\right) = 2L\left(f(x^k) - f(x^*)\right),\\
	\|g^k\|^2 &=& \|\nabla f(x^k)\|^2 \overset{\eqref{eq:L_smoothness_cor_3}}{\le} 2L\left(f(x^k) - f(x^*)\right).
\end{eqnarray*}

Applying Theorem~\ref{thm:ec_sgd_main_result_new} we get the following result.
\begin{theorem}\label{thm:ec_sgd_star}
	Assume that $f_i(x)$ is convex and $L$-smooth for all $i=1,\ldots, n$ and $f(x)$ is $\mu$-quasi strongly convex. Then {\tt EC-GDstar} satisfies Assumption~\ref{ass:key_assumption_finite_sums_new} with
	\begin{gather*}
		A = A' = L,\quad \widetilde{A} = 0,\quad B_1 = B_2 = \widetilde{B}_1 = \widetilde{B}_2 = B_1' = B_2' = 0,\\
		 D_1 = \widetilde{D}_1 = D_1' = 0,\quad \sigma_{1,k}^2 \equiv \sigma_{2,k}^2 \equiv 0,\\
		\rho_1 = \rho_2 = 1,\quad C_1 = C_2 = 0,\quad G = 0,\quad D_2 = 0,\quad F_1 = F_2 = 0,\quad D_3 = 0,
	\end{gather*}
	with $\gamma$ satisfying
	\begin{equation*}
		\gamma \le \frac{\delta}{8L\sqrt{3}}
	\end{equation*}
	and for all $K \ge 0$
	\begin{equation*}
		\EE\left[f(\bar{x}^K) - f(x^*)\right] \le \left(1 - \frac{\gamma\mu}{2}\right)^K\frac{4\|x^0 - x^*\|^2}{\gamma},
	\end{equation*}
	when $\mu > 0$ and
	\begin{equation*}
		\EE\left[f(\bar{x}^K) - f(x^*)\right] \le \frac{4\|x^0 - x^*\|^2}{K\gamma}
	\end{equation*}
	when $\mu=0$.
\end{theorem}
In other words, {\tt EC-GDstar} converges with linear rate $\cO\left(\frac{\kappa}{\delta}\ln\frac{1}{\varepsilon}\right)$ to the exact solution when $\mu > 0$ removing the drawback of {\tt EC-SGD} and {\tt EC-GD}. If $\mu = 0$ then the rate of convergence is $\cO\left(\frac{L\|x^0-x^*\|^2}{\delta\varepsilon}\right)$. However, {\tt EC-GDstar} relies on the fact that $i$-th node knows $\nabla f_i(x^*)$ which is not realistic.

\subsection{{\tt EC-SGD-DIANA}}\label{sec:ec_diana}
In this section we present a new method that converges to the exact optimum asymptotically but does not need to know $\nabla f_i(x^*)$ and instead of this it learns the gradients at the optimum. This method is inspired by another method called {\tt DIANA} (see \cite{mishchenko2019distributed, horvath2019stochastic}).
\begin{algorithm}[t]
   \caption{{\tt EC-SGD-DIANA}}\label{alg:EC-SGD-DIANA}
\begin{algorithmic}[1]
   \Require learning rates $\gamma>0$, $\alpha \in (0,1]$, initial vectors $x^0, h_1^0,\ldots, h_n^0 \in \R^d$
	\State Set $e_i^0 = 0$ for all $i=1,\ldots, n$   
	\State Set $h^0 = \frac{1}{n}\sum_{i=1}^n h_i^0$   
   \For{$k=0,1,\dotsc$}
       \State Broadcast $x^{k}, h^k$ to all workers
        \For{$i=1,\dotsc,n$ in parallel}
			\State Sample $\hat g_i^k$ such that $\EE[\hat g_i^k\mid x^k] = \nabla f_i(x^k)$ and $\EE\left[\|\hat g_i^k - \nabla f_i(x^k)\|^2\mid x^k\right] \le \widetilde{D}_{1,i}$ independently from other workers        
            \State $g^{k}_i = \hat g_i^k - h_i^k + h^k$
            \State $v_i^k = C(e_i^k + \gamma g_i^k)$
            \State $e_i^{k+1} = e_i^k + \gamma g_i^k - v_i^k$
            \State $h_i^{k+1} = h_i^k + \alpha Q(\hat g_i^k - h_i^k)$
        \EndFor
        \State $e^k = \frac{1}{n}\sum_{i=1}^ne_i^k$, $g^k = \frac{1}{n}\sum_{i=1}^ng_i^k$, $v^k = \frac{1}{n}\sum_{i=1}^nv_i^k$, $h^{k+1} = \frac{1}{n}\sum\limits_{i=1}^n h_i^{k+1} = h^k + \alpha\frac{1}{n}\sum\limits_{i=1}^n Q(\hat g_i^k - h_i^k)$
       \State $x^{k+1} = x^k - v^k$
   \EndFor
\end{algorithmic}
\end{algorithm}

We notice that master needs to gather only $C(e_i^k + \gamma g_i^k)$ and $Q(\hat g_i^k - h_i^k)$ from all nodes in order to perform an update.

\begin{lemma}\label{lem:ec_diana_second_moment_bound}
	Assume that $f_i(x)$ is convex and $L$-smooth for all $i=1,\ldots,n$. Then, for all $k\ge 0$ we have
	\begin{eqnarray}
		\EE\left[g^k\mid x^k\right] &=& \nabla f(x^k), \label{eq:ec_diana_unbiasedness}\\
		\frac{1}{n}\sum\limits_{i=1}^n\|\bar{g}_i^k\|^2 &\le& 4L\left(f(x^k) - f(x^*)\right) + 2\sigma_k^2, \label{eq:ec_diana_second_moment_bound}\\
		\frac{1}{n}\sum\limits_{i=1}^n\EE\left[\|g_i^k-\bar{g}_i^k\|^2\mid x^k\right] &\le& \widetilde{D}_1, \label{eq:ec_diana_variance_bound}\\
		\EE\left[\|g^k\|^2\mid x^k\right] &\le& 2L\left(f(x^k) - f(x^*)\right) + \frac{\widetilde{D}_1}{n} \label{eq:ec_diana_second_moment_bound_2}
	\end{eqnarray}
	where $\widetilde{D}_1 = \frac{1}{n}\sum_{i=1}^n \widetilde{D}_{1,i}$ and $\sigma_k^2 = \frac{1}{n}\sum_{i=1}^n\|h_i^k - \nabla f(x^*)\|^2$.
\end{lemma}
\begin{proof}
	First of all, we show unbiasedness of $g^k$:
	\begin{equation*}
		\EE\left[g^k\mid x^k\right] = \frac{1}{n}\sum\limits_{i=1}^n\EE\left[g_i^k\mid x^k\right] = \frac{1}{n}\sum\limits_{i=1}^n\left(\nabla f_i(x^k) - h_i^k + h^k\right) = \nabla f(x^k).
	\end{equation*}
	Next, we derive the upper bound for $\|\bar{g}_i^k\|^2$:
	\begin{eqnarray*}
		\|\bar{g}_i^k\|^2 &=& \|\nabla f_i(x^k) - h_i^k - h^k\|^2\\
		&\overset{\eqref{eq:a_b_norm_squared}}{\le}& 2\|\nabla f_i(x^k) - \nabla f_i(x^*)\|^2 + 2\left\|h_i^k - \nabla f_i(x^*) -\left(h^k + \nabla f(x^*)\right)\right\|^2	\\
		&\overset{\eqref{eq:L_smoothness_cor_3}}{\le}& 4L\left(f_i(x^k) - \nabla f_i(x^*) - \langle\nabla f_i(x^*), x^k - x^*\rangle\right)\\
		&&\quad +2\left\|h_i^k - \nabla f_i(x^*) -\left(h^k + \nabla f(x^*)\right)\right\|^2.
	\end{eqnarray*}
	Summing up previous inequality for $i = 1,\ldots, n$ we get
	\begin{eqnarray}
		\frac{1}{n}\sum\limits_{i=1}^n\|\bar{g}_i^k\|^2 &\le&4L(f(x^k) - f(x^*)) + \frac{2}{n}\sum\limits_{i=1}^n\left\|h_i^k - \nabla f_i(x^*) - \left(\frac{1}{n}\sum\limits_{i=1}^n(h_i^k - \nabla f_i(x^*))\right)\right\|^2\notag\\
		&\overset{\eqref{eq:variance_decomposition}}{\le}& 4L\left(f(x^k) - f(x^*)\right) + \frac{2}{n}\sum\limits_{i=1}^n\|h_i^k - \nabla f(x^*)\|^2.\label{eq:ec_sgd_diana_first_ineq}
	\end{eqnarray}
	Using the unbiasedness of $\hat g_i^k$ we derive
	\begin{equation*}
		\frac{1}{n}\sum\limits_{i=1}^n\EE\left[\|g_i^k-\bar{g}_i^k\|^2\mid x^k\right] = \frac{1}{n}\sum\limits_{i=1}^n\EE\left[\|\hat g_i^k-\nabla f_i(x^k)\|^2\mid x^k\right] \le \frac{1}{n}\sum\limits_{i=1}^n\widetilde{D}_{1,i} = \widetilde{D}_1.
	\end{equation*}
	Finally, we obtain the upper bound for the second moment of $g^k$ using the independence of $\hat g_1^k,\ldots, \hat g_n^k$:
	\begin{eqnarray*}
		\EE\left[\|g^k\|^2\mid x^k\right] &\overset{\eqref{eq:variance_decomposition}}{=}& \|\nabla f(x^k)\|^2 + \EE\left[\|g^k - \nabla f(x^k)\|^2\right]\\
		&\overset{\eqref{eq:L_smoothness_cor_3}}{\le}&  2L(f(x^k)-f(x^*)) + \EE\left[\left\|\frac{1}{n}\sum\limits_{i=1}^n(\hat g_i^k - \nabla f_i(x^k))\right\|^2\mid x^k\right]\\
		&=&2L(f(x^k)-f(x^*)) + \frac{1}{n^2}\sum\limits_{i=1}^n\EE\left[\left\|\hat g_i^k - \nabla f_i(x^k)\right\|^2\mid x^k\right]\\
		&\le& 2L(f(x^k)-f(x^*)) + \frac{1}{n^2}\sum\limits_{i=1}^n\widetilde{D}_{1,i}.
	\end{eqnarray*}
\end{proof}

\begin{lemma}\label{lem:ec_diana_sigma_k+1_bound}
	Let assumptions of Lemma~\ref{lem:ec_diana_second_moment_bound} hold and $\alpha\le\nicefrac{1}{(\omega+1)}$. Then, for all $k\ge 0$ we have
	\begin{equation}
		\EE\left[\sigma_{k+1}^2\mid x^k\right] \le (1 - \alpha)\sigma_k^2 + 2L\alpha(f(x^k) - f(x^*)) + \alpha^2(\omega+1) \widetilde{D}_1, \label{eq:ec_diana_sigma_k+1_bound}
	\end{equation}
	where $\sigma_k^2 = \frac{1}{n}\sum_{i=1}^n\|h_i^k - \nabla f_i(x^*)\|^2$ and $\widetilde{D}_1 = \frac{1}{n}\sum_{i=1}^n \widetilde{D}_{1,i}$.
\end{lemma}
\begin{proof}
	For simplicity, we introduce new notation: $h_i^* \eqdef \nabla f_i(x^*)$. Using this we derive an upper bound for the second moment of $h_i^{k+1} - h_i^*$:
	\begin{eqnarray*}
		\EE\left[\|h_i^{k+1} - h_i^*\|^2\mid x^k\right] &=& \EE\left[\left\|h_i^k - h_i^* + \alpha Q(\hat g_i^k - h_i^k) \right\|^2\mid x^k\right]\\
		&\overset{\eqref{eq:quantization_def}}{=}& \|h_i^k - h_i^*\|^2 +2\alpha\langle h_i^k - h_i^*, \nabla f_i(x^k) - h_i^k \rangle\\
		&&\quad + \alpha^2\EE\left[\|Q(\hat g_i^k - h_i^k)\|^2\mid x^k\right]\\
		&\overset{\eqref{eq:quantization_def},\eqref{eq:tower_property}}{\le}& \|h_i^k - h_i^*\|^2 +2\alpha\langle h_i^k - h_i^*, \nabla f_i(x^k) - h_i^k \rangle\\
		&&\quad + \alpha^2(\omega+1)\EE\left[\|\hat g_i^k - h_i^k\|^2\mid x^k\right].
	\end{eqnarray*}
	Using variance decomposition \eqref{eq:variance_decomposition} and $\alpha\le\nicefrac{1}{(\omega+1)}$ we get
	\begin{eqnarray*}
		\alpha^2(\omega+1)\EE\left[\|\hat g_i^k - h_i^k\|^2\mid x^k\right] &\overset{\eqref{eq:variance_decomposition}}{=}& \alpha^2(\omega+1)\EE\left[\|\hat g_i^k - \nabla f_i(x^k)\|^2\mid x^k\right]\\
		&&\quad + \alpha^2(\omega+1)\|\nabla f_i(x^k) - h_i^k\|^2\\
		&\le& \alpha^2(\omega+1) \widetilde{D}_{1,i} + \alpha\|\nabla f_i(x^k) - h_i^k\|^2.
	\end{eqnarray*}
	Putting all together we obtain
	\begin{eqnarray*}
		\EE\left[\|h_i^{k+1} - h_i^*\|^2\mid x^k\right] &\le& \|h_i^k - h_i^*\|^2 + \alpha\left\langle \nabla f_i(x^k) - h_i^k, f_i(x^k) + h_i^k - 2h_i^* \right\rangle + \alpha^2(\omega+1) \widetilde{D}_{1,i}\\
		&\overset{\eqref{eq:a-b_a+b}}{=}& \|h_i^k - h_i^*\|^2 + \alpha\|\nabla f_i(x^k) - h_i^*\|^2 - \alpha\|h_i^k - h_i^*\|^2 + \alpha^2(\omega+1) \widetilde{D}_{1,i}\\
		&\overset{\eqref{eq:L_smoothness_cor_3}}{\le}& (1-\alpha)\|h_i^k - h_i^*\|^2 + 2L\alpha\left(f_i(x^k) - f_i(x^*) - \langle\nabla f_i(x^*), x^k - x^* \rangle\right)\\
		&&\quad +\alpha^2(\omega+1) \widetilde{D}_{1,i}.
	\end{eqnarray*}
	Summing up the above inequality for $i=1,\ldots, n$ we derive
	\begin{equation*}
		\frac{1}{n}\sum\limits_{i=1}^n\EE\left[\|h_i^{k+1} - h_i^*\|^2\mid x^k\right] \le \frac{1-\alpha}{n}\sum\limits_{i=1}^n\|h_i^k - h_i^*\|^2 + 2L\alpha(f(x^k) - f(x^*)) + \frac{\alpha^2(\omega+1)}{n}\sum\limits_{i=1}^n\widetilde{D}_{1,i}.
	\end{equation*}
\end{proof}

Applying Theorem~\ref{thm:ec_sgd_main_result_new} we get the following result.
\begin{theorem}\label{thm:ec_diana}
	Assume that $f_i(x)$ is convex and $L$-smooth for all $i=1,\ldots, n$ and $f(x)$ is $\mu$-quasi strongly convex. Then {\tt EC-SGD-DIANA} satisfies Assumption~\ref{ass:key_assumption_finite_sums_new} with
	\begin{gather*}
		A = 2L,\quad \widetilde{A} = 0,\quad A' = L,\quad B_1 = 2,\quad \widetilde{D}_1 = \frac{1}{n}\sum\limits_{i=1}^n \widetilde{D}_{1,i},\quad \sigma_{1,k}^2 = \sigma_k^2 = \frac{1}{n}\sum\limits_{i=1}^n\|h_i^k - \nabla f_i(x^*)\|^2,\\
		B_1' = B_2' = B_2 = \widetilde{B}_1 = \widetilde{B}_2 = 0,\quad \sigma_{2,k}^2\equiv 0,\quad \rho_1 = \alpha,\quad \rho_2 = 1,\quad C_1 = L\alpha,\quad C_2 = 0,\quad D_1 = 0,\\
		D_2 = \alpha^2(\omega+1) \widetilde{D}_1,\quad D_1' = \frac{D_1}{n},\quad G = 0,\\
		F_1 = \frac{96L\gamma^2}{\delta^2\alpha\left(1-\min\left\{\frac{\gamma\mu}{2},\frac{\alpha}{4}\right\}\right)},\quad F_2 = 0,\quad D_3 = \frac{6L\gamma}{\delta}\left(\frac{4\alpha(\omega+1)}{\delta}+1\right)\widetilde{D}_1,
	\end{gather*}
	with $\gamma$ and $\alpha$ satisfying
	\begin{equation*}
		\gamma \le \min\left\{\frac{1}{4L}, \frac{\delta\sqrt{1-\alpha}}{8L\sqrt{6(3-\alpha)}}\right\},\quad \alpha \le \frac{1}{\omega+1},\quad M_1 = M_2 = 0
	\end{equation*}
	and for all $K \ge 0$
	\begin{equation*}
		\EE\left[f(\bar x^K) - f(x^*)\right] \le \left(1 - \min\left\{\frac{\gamma\mu}{2},\frac{\alpha}{4}\right\}\right)^K\frac{4(\|x^0 - x^*\|^2 + \gamma F_1 \sigma_0^2)}{\gamma} + 4\gamma\left(D_1' + D_3\right),
	\end{equation*}	
	when $\mu > 0$ and
	\begin{equation*}
		\EE\left[f(\bar{x}^K) - f(x^*)\right] \le \frac{4(\|x^0 - x^*\|^2 + \gamma F_1 \sigma_0^2)}{\gamma K} + 4\gamma\left(D_1' + D_3\right)
	\end{equation*}
	when $\mu=0$.
\end{theorem}
In other words, if 
\begin{equation*}
		\gamma = \min\left\{\frac{1}{4L}, \frac{\delta\sqrt{1-\alpha}}{8L\sqrt{6(3-\alpha)}}\right\},\quad \alpha = \min\left\{\frac{1}{\omega+1},\frac{1}{2}\right\}
\end{equation*}
and $\widetilde{D}_1 = 0$, i.e., $\hat g_i^k = \nabla f_i(x^k)$ almost surely (this is the setup of {\tt EC-GD-DIANA}), {\tt EC-SGD-DIANA} converges with the linear rate
\begin{equation*}
	\cO\left(\left(\omega + \frac{\kappa}{\delta}\right)\ln\frac{1}{\varepsilon}\right)
\end{equation*}
to the exact solution.  Applying Lemma~\ref{lem:lemma2_stich} we establish the rate of convergence to $\varepsilon$-solution in the case when $\mu > 0$.
\begin{corollary}\label{cor:ec_diana_str_cvx_cor}
	Let the assumptions of Theorem~\ref{thm:ec_diana} hold and $\mu > 0$. Then after $K$ iterations of {\tt EC-SGD-DIANA} with the stepsize
	\begin{eqnarray*}
		\gamma_0 &=& \min\left\{\frac{1}{4L}, \frac{\delta\sqrt{1-\alpha}}{8L\sqrt{6(3-\alpha)}}\right\},\quad R_0 = \|x^0-x^*\|,\quad  \tilde{F_1} = \frac{784L\gamma^2}{7\delta^2\alpha},\\
		\gamma &=& \min\left\{\gamma_0, \frac{\ln\left(\max\left\{2,\min\left\{\frac{n\left(R_0^2 + \tilde{F}_1\gamma_0\sigma_{1,0}^2\right)\mu^2K^2}{\widetilde{D}_1}, \frac{\delta\left(R_0^2 + \tilde{F}_1\gamma_0\sigma_{1,0}^2\right)\mu^3K^3}{6L\widetilde{D}_1(\nicefrac{4\alpha(\omega+1)}{\delta}+1)}\right\}\right\}\right)}{\mu K}\right\},
	\end{eqnarray*}
	and $\alpha \le \frac{1}{\omega+1}$ we have
	\begin{equation*}
		\EE\left[f(\bar{x}^K) - f(x^*)\right] = \widetilde\cO\left(\frac{L}{\delta}R_0^2\exp\left(-\min\left\{\frac{\delta\mu}{L},\alpha\right\}K\right) + \frac{\widetilde{D}_1}{n\mu K} + \frac{L\widetilde{D}_1\left(\nicefrac{\alpha(\omega+1)}{\delta}+1\right)}{\delta\mu^2 K^2}\right).
	\end{equation*}
	That is, to achive $\EE\left[f(\bar{x}^K) - f(x^*)\right] \le \varepsilon$ {\tt EC-SGD-DIANA} requires
	\begin{equation*}
		\widetilde{\cO}\left(\frac{1}{\alpha} + \frac{L}{\delta\mu} + \frac{D_1}{n\mu\varepsilon} + \frac{\sqrt{L\widetilde{D}_1\left(\nicefrac{\alpha(\omega+1)}{\delta}+1\right)}}{\mu\sqrt{\delta\varepsilon}}\right) \text{ iterations.}
	\end{equation*}
	In particular, if $\alpha = \frac{1}{\omega+1}$, then to achive $\EE\left[f(\bar{x}^K) - f(x^*)\right] \le \varepsilon$ {\tt EC-SGD-DIANA} requires
	\begin{equation*}
		\widetilde{\cO}\left(\omega + \frac{L}{\delta\mu} + \frac{\widetilde{D}_1}{n\mu\varepsilon} + \frac{\sqrt{L\widetilde{D}_1}}{\delta\mu\sqrt{\varepsilon}}\right) \text{ iterations,}
	\end{equation*}
	and if $\alpha = \frac{\delta}{\omega+1}$, then to achive $\EE\left[f(\bar{x}^K) - f(x^*)\right] \le \varepsilon$ {\tt EC-SGD-DIANA} requires
	\begin{equation*}
		\widetilde{\cO}\left(\frac{\omega+1}{\delta} + \frac{L}{\delta\mu} + \frac{\widetilde{D}_1}{n\mu\varepsilon} + \frac{\sqrt{L\widetilde{D}_1}}{\mu\sqrt{\delta\varepsilon}}\right) \text{ iterations.}
	\end{equation*}
\end{corollary}

Applying Lemma~\ref{lem:lemma_technical_cvx} we get the complexity result in the case when $\mu = 0$.
\begin{corollary}\label{cor:ec_diana_cvx_cor}
	Let the assumptions of Theorem~\ref{thm:ec_diana} hold and $\mu = 0$. Then after $K$ iterations of {\tt EC-SGD-DIANA} with the stepsize
	\begin{eqnarray*}
		\gamma_0 &=& \min\left\{\frac{1}{4L}, \frac{\delta\sqrt{1-\alpha}}{8L\sqrt{6(3-\alpha)}}\right\},\quad R_0 = \|x^0-x^*\|,\\
		\gamma &=& \min\left\{\gamma_0, \sqrt[3]{\frac{R_0^2\delta^2\alpha\left(1-\min\left\{\frac{\gamma_0\mu}{2},\frac{\alpha}{4}\right\}\right)}{96L\sigma_0^2}}, \sqrt{\frac{nR_0^2}{\widetilde{D}_1 K}}, \sqrt[3]{\frac{\delta R_0^2}{6L\widetilde{D}_1\left(\frac{4\alpha(\omega+1)}{\delta}+1\right) K}}\right\},
	\end{eqnarray*}	
	and $\alpha \le \frac{1}{\omega+1}$ we have $\EE\left[f(\bar{x}^K) - f(x^*)\right]$ of order
	\begin{equation*}
		\cO\left(\frac{LR_0^2}{\delta K} + \frac{\sqrt[3]{LR_0^4\sigma_0^2}}{K\sqrt[3]{\delta^2\alpha}} + \sqrt{\frac{R_0^2\widetilde{D}_1}{nK}}+ \sqrt[3]{\frac{LR_0^4\widetilde{D}_1\left(\frac{\alpha(\omega+1)}{\delta}+1\right)}{\delta K^2}}\right).
	\end{equation*}
	That is, to achive $\EE\left[f(\bar{x}^K) - f(x^*)\right] \le \varepsilon$ {\tt EC-SGD-DIANA} requires
	\begin{equation*}
		\cO\left(\frac{LR_0^2}{\delta \varepsilon} + \frac{\sqrt[3]{LR_0^4\sigma_0^2}}{\varepsilon\sqrt[3]{\delta^2\alpha}} + \frac{R_0^2\widetilde{D}_1}{n\varepsilon^2}+ \frac{R_0^2\sqrt{L\widetilde{D}_1\left(\frac{\alpha(\omega+1)}{\delta}+1\right)}}{\sqrt{\delta\varepsilon^3}}\right)
	\end{equation*}
	iterations. In particular, if $\alpha = \frac{1}{\omega+1}$, then to achive $\EE\left[f(\bar{x}^K) - f(x^*)\right] \le \varepsilon$ {\tt EC-SGD-DIANA} requires
	\begin{equation*}
		\cO\left(\frac{LR_0^2}{\delta \varepsilon} + \frac{\sqrt[3]{LR_0^4(\omega+1)\sigma_0^2}}{\varepsilon\sqrt[3]{\delta^2}} + \frac{R_0^2\widetilde{D}_1}{n\varepsilon^2}+ \frac{R_0^2\sqrt{L\widetilde{D}_1}}{\delta\sqrt{\varepsilon^3}}\right) \text{ iterations,}
	\end{equation*}
	and if $\alpha = \frac{\delta}{\omega+1}$, then to achive $\EE\left[f(\bar{x}^K) - f(x^*)\right] \le \varepsilon$ {\tt EC-SGD-DIANA} requires
	\begin{equation*}
		\cO\left(\frac{LR_0^2}{\delta \varepsilon} + \frac{\sqrt[3]{LR_0^4(\omega+1)\sigma_0^2}}{\delta\varepsilon} + \frac{R_0^2\widetilde{D}_1}{n\varepsilon^2}+ \frac{R_0^2\sqrt{L\widetilde{D}_1}}{\sqrt{\delta\varepsilon^3}}\right) \text{ iterations.}
	\end{equation*}
\end{corollary}

\subsection{{\tt EC-SGDsr-DIANA}}\label{sec:ec_sgdsr_DIANA}
In this section we consider the same setup as in Section~\ref{sec:diana_arbitrary_sampling} and consider {\tt EC-SGD-DIANA} adjusted to this setup. The resulting algorithm is called {\tt EC-SGDsr-DIANA}, see
\begin{algorithm}[t]
   \caption{{\tt EC-SGDsr-DIANA}}\label{alg:EC-SGDsr-DIANA}
\begin{algorithmic}[1]
   \Require learning rates $\gamma>0$, $\alpha \in (0,1]$, initial vectors $x^0, h_1^0,\ldots, h_n^0 \in \R^d$
	\State Set $e_i^0 = 0$ for all $i=1,\ldots, n$   
	\State Set $h^0 = \frac{1}{n}\sum_{i=1}^n h_i^0$   
   \For{$k=0,1,\dotsc$}
       \State Broadcast $x^{k}, h^k$ to all workers
        \For{$i=1,\dotsc,n$ in parallel}
			\State Sample $\hat g_i^{k} = \nabla f_{\xi_i^k}(x^k)$ satisfying Assumption~\ref{ass:exp_smoothness} independtently from other workers            
            \State $g^{k}_i = \hat g_i^k - h_i^k + h^k$
            \State $v_i^k = C(e_i^k + \gamma g_i^k)$
            \State $e_i^{k+1} = e_i^k + \gamma g_i^k - v_i^k$
            \State $h_i^{k+1} = h_i^k + \alpha Q(\hat g_i^k - h_i^k)$ \Comment{$Q(\cdot)$ is calculated independtly from other workers}
        \EndFor
        \State $e^k = \frac{1}{n}\sum_{i=1}^ne_i^k$, $g^k = \frac{1}{n}\sum_{i=1}^ng_i^k$, $v^k = \frac{1}{n}\sum_{i=1}^nv_i^k$, $h^{k+1} = \frac{1}{n}\sum\limits_{i=1}^n h_i^{k+1} = h^k + \alpha\frac{1}{n}\sum\limits_{i=1}^n Q(\hat g_i^k - h_i^k)$
       \State $x^{k+1} = x^k - v^k$
   \EndFor
\end{algorithmic}
\end{algorithm}
\begin{lemma}\label{lem:ec_sgdsr_diana_second_moment_bound}
	Let Assumption~\ref{ass:exp_smoothness} be satisfied and $f_i$ be convex and $L$-smooth for all $i\in[n]$. Then, for all $k\ge 0$ we have
	\begin{eqnarray}
		\EE\left[g^k\mid x^k\right] &=& \nabla f(x^k), \label{eq:ec_sgdsr_diana_unbiasedness}\\
		\frac{1}{n}\sum\limits_{i=1}^n\|\bar{g}_i^k\|^2 &\le& 4L\left(f(x^k) - f(x^*)\right) + 2\sigma_k^2, \label{eq:ec_sgdsr_diana_second_moment_bound}\\
		\frac{1}{n}\sum\limits_{i=1}^n\EE\left[\|g_i^k-\bar{g}_i^k\|^2\mid x^k\right] &\le& 6(\cL+L)\left(f(x^k) - f(x^*)\right) + \widetilde{D}_1, \label{eq:ec_sgdsr_diana_variance_bound}\\
		\EE\left[\|g^k\|^2\mid x^k\right] &\le& 4\cL\left(f(x^k) - f(x^*)\right) + D_1' \label{eq:ec_sgssr_diana_second_moment_bound_2}
	\end{eqnarray}
	where $\sigma_k^2 = \frac{1}{n}\sum_{i=1}^n\|h_i^k - \nabla f(x^*)\|^2$, $\widetilde{D}_1 = \frac{3}{n}\sum_{i=1}^n \EE_{\cD_i}\left[\|\nabla f_{\xi_i}(x^*) - \nabla f_i(x^*)\|^2\right]$ and\newline $D_1' = \frac{2}{n^2}\sum\limits_{i=1}^n\EE_{\cD_i}\left[\|\nabla f_{\xi_i}(x^*) - \nabla f_{i}(x^*)\|^2\right]$.
\end{lemma}
\begin{proof}
	First of all, we show unbiasedness of $g^k$:
	\begin{equation*}
		\EE\left[g^k\mid x^k\right] = \frac{1}{n}\sum\limits_{i=1}^n\EE\left[g_i^k\mid x^k\right] = \frac{1}{n}\sum\limits_{i=1}^n\left(\nabla f_i(x^k) - h_i^k + h^k\right) = \nabla f(x^k).
	\end{equation*}
	Following the same steps as in the proof of \eqref{eq:ec_sgd_diana_first_ineq} we derive \eqref{eq:ec_sgdsr_diana_second_moment_bound}. Next, we establish \eqref{eq:ec_sgdsr_diana_variance_bound}:
	\begin{eqnarray*}
		\frac{1}{n}\sum\limits_{i=1}^n\EE\left[\|g_i^k-\bar{g}_i^k\|^2\mid x^k\right] &=& \frac{1}{n}\sum\limits_{i=1}^n\EE_{\cD_i}\left[\|\nabla f_{\xi_i^k}(x^k)-\nabla f_i(x^k)\|^2\right]\\
		&\overset{\eqref{eq:a_b_norm_squared}}{\le}& \frac{3}{n}\sum\limits_{i=1}^n\EE_{\cD_i}\left[\|\nabla f_{\xi_i^k}(x^k)-\nabla f_{\xi_i^k}(x^*)\|^2\right]\\
		&&\quad + \frac{3}{n}\sum\limits_{i=1}^n\EE_{\cD_i}\left[\|\nabla f_{\xi_i^k}(x^*)-\nabla f_i(x^*)\|^2\right]\\
		&&\quad + \frac{3}{n}\sum\limits_{i=1}^n\|\nabla f_{i}(x^*)-\nabla f_i(x^k)\|^2\\
		&\overset{\eqref{eq:L_smoothness_cor_3},\eqref{eq:exp_smoothness}}{\le}& 6(\cL + L)\left(f(x^k)-f(x^*)\right) \\
		&&\quad + \frac{3}{n}\sum\limits_{i=1}^n\EE_{\cD_i}\left[\|\nabla f_{\xi_i}(x^*)-\nabla f_i(x^*)\|^2\right].
	\end{eqnarray*}
	Finally, we obtain the upper bound for the second moment of $g^k$ using the independence of $\xi_1^k,\ldots,\xi_n^k$:
	\begin{eqnarray*}
		\EE\left[\|g^k\|^2\mid x^k\right] &=& \EE\left[\left\|\frac{1}{n}\sum\limits_{i=1}^n(\nabla f_{\xi_i^k}(x^k) - \nabla f_{\xi_i^k}(x^*) + \nabla f_{\xi_i^k}(x^*) - \nabla f_i(x^*))\right\|^2\mid x^k\right]\\
		&\overset{\eqref{eq:a_b_norm_squared}}{\le}& \frac{2}{n}\sum\limits_{i=1}^n\EE\left[\|\nabla f_{\xi_i^k}(x^k) - \nabla f_{\xi_i^k}(x^*)\|^2\mid x^k\right]\\
		&&\quad + 2\EE\left[\left\|\frac{1}{n}\sum\limits_{i=1}^n(\nabla f_{\xi_i^k}(x^*) - \nabla f_i(x^*))\right\|^2\mid x^k\right]\\
		&\overset{\eqref{eq:exp_smoothness}}{\le}& 4\cL\left(f(x^k) - f(x^*)\right) + \frac{2}{n^2}\sum\limits_{i=1}^n\EE_{\cD_i}\left[\|\nabla f_{\xi_i}(x^*) - \nabla f_{i}(x^*)\|^2\right].
	\end{eqnarray*}
\end{proof}

\begin{lemma}\label{lem:ec_sgdsr_diana_sigma_k+1_bound}
	Let $f_i$ be convex and $L$-smooth, Assumption~\ref{ass:exp_smoothness} holds and $\alpha \le \nicefrac{1}{(\omega+1)}$. Then, for all $k\ge 0$ we have
	\begin{equation}
		\EE\left[\sigma_{k+1}^2\mid x^k\right] \le (1 - \alpha)\sigma_k^2 + 2\alpha(3\cL+4L)(f(x^k) - f(x^*)) + D_2, \label{eq:ec_sgdsr_diana_sigma_k+1_bound}
	\end{equation}
	where $\sigma_k^2 = \frac{1}{n}\sum_{i=1}^n\|h_i^k - \nabla f_i(x^*)\|^2$ and $D_2 = \alpha^2(\omega+1)\widetilde{D}_1$.
\end{lemma}
\begin{proof}
The proof is identical to the proof of Lemma~\ref{lem:diana_sigma_k+1_bound} up to the following changes in the notation: $\omega_1 = \omega$, $\Delta_i^k = Q(\hat g_i^k - h_i^k)$ and $\hat \Delta_i^k = \hat g_i^k - h_i^k$.
\end{proof}

Applying Theorem~\ref{thm:ec_sgd_main_result_new} we get the following result.
\begin{theorem}\label{thm:ec_sgdsr_diana}
	Assume that $f_i(x)$ is convex and $L$-smooth for all $i=1,\ldots, n$, $f(x)$ is $\mu$-quasi strongly convex and Assumption~\ref{ass:exp_smoothness} holds. Then {\tt EC-SGDsr-DIANA} satisfies Assumption~\ref{ass:key_assumption_finite_sums_new} with
	\begin{gather*}
		A = 2L,\quad \widetilde{A}=3(\cL+L),\quad A' = 2\cL,\quad B_1 = 2,\quad \widetilde{D}_1 = \frac{3}{n}\sum_{i=1}^n \EE_{\cD_i}\left[\|\nabla f_{\xi_i}(x^*) - \nabla f_i(x^*)\|^2\right],\\
		\sigma_{1,k}^2 = \sigma_k^2 = \frac{1}{n}\sum\limits_{i=1}^n\|h_i^k - \nabla f_i(x^*)\|^2,\quad D_1 = 0,\quad D_1' = \frac{2}{3n}\widetilde{D}_1,\quad D_2 = \alpha^2(\omega+1) \widetilde{D}_1\\
		\widetilde{B}_1 = B_1' = B_2' = B_2 = \widetilde{B}_2 = 0,\quad \sigma_{2,k}^2\equiv 0,\quad \rho_1 = \alpha,\quad \rho_2 = 1,\quad C_1 = 2\alpha(3\cL+4L),\quad C_2 = 0,\\
		G = 0,\quad F_1 = \frac{96L\gamma^2}{\delta^2\alpha\left(1-\min\left\{\frac{\gamma\mu}{2},\frac{\alpha}{4}\right\}\right)},\quad F_2 = 0,\quad D_3 = \frac{6L\gamma}{\delta}\left(\frac{4\alpha(\omega+1)}{\delta}+1\right)\widetilde{D}_1,
	\end{gather*}
	with $\gamma$ and $\alpha$ satisfying
	\begin{equation*}
		\gamma \le \min\left\{\frac{1}{4\cL}, \frac{\delta}{4\sqrt{6L\left(4L + 3\delta(\cL + L) + \frac{16(3\cL + 4L)}{1-\alpha}\right)}}\right\},\quad \alpha \le \frac{1}{\omega+1},\quad M_1 = M_2 = 0.
	\end{equation*}
	and for all $K \ge 0$
	\begin{equation*}
		\EE\left[f(\bar x^K) - f(x^*)\right] \le \left(1 - \min\left\{\frac{\gamma\mu}{2},\frac{\alpha}{4}\right\}\right)^K\frac{4(\|x^0 - x^*\|^2 + \gamma F_1 \sigma_0^2)}{\gamma} + 4\gamma\left(D_1' + D_3\right),
	\end{equation*}	
	when $\mu > 0$ and
	\begin{equation*}
		\EE\left[f(\bar{x}^K) - f(x^*)\right] \le \frac{4(\|x^0 - x^*\|^2 + \gamma F_1 \sigma_0^2)}{\gamma K} + 4\gamma\left(D_1' + D_3\right)
	\end{equation*}
	when $\mu=0$.
\end{theorem}
Applying Lemma~\ref{lem:lemma2_stich} we establish the rate of convergence to $\varepsilon$-solution in the case when $\mu > 0$.
\begin{corollary}\label{cor:ec_sgdsr_diana_str_cvx_cor}
	Let the assumptions of Theorem~\ref{thm:ec_sgdsr_diana} hold and $\mu > 0$. Then after $K$ iterations of {\tt EC-SGDsr-DIANA} with the stepsize
	\begin{eqnarray*}
		\gamma_0 &=& \min\left\{\frac{1}{4\cL}, \frac{\delta}{4\sqrt{6L\left(4L + 3\delta(\cL + L) + \frac{16(3\cL + 4L)}{1-\alpha}\right)}}\right\},\\
		R_0 &=& \|x^0-x^*\|,\quad  \tilde{F_1} = \frac{96L\gamma_0^2}{\delta^2\alpha\left(1-\min\left\{\frac{\gamma_0\mu}{2},\frac{\alpha}{4}\right\}\right)},\\
		\gamma &=& \min\left\{\gamma_0, \frac{\ln\left(\max\left\{2,\min\left\{\frac{3n\left(R_0^2 + \tilde{F}_1\gamma_0\sigma_{1,0}^2\right)\mu^2K^2}{2\widetilde{D}_1}, \frac{\delta\left(R_0^2 + \tilde{F}_1\gamma_0\sigma_{1,0}^2\right)\mu^3K^3}{6L\widetilde{D}_1\left(\frac{4\alpha(\omega+1)}{\delta}+1\right)}\right\}\right\}\right)}{\mu K}\right\},
	\end{eqnarray*} 
	and $\alpha \le \frac{1}{\omega+1}$ we have $\EE\left[f(\bar{x}^K) - f(x^*)\right]$ of order
	\begin{equation*}
		\widetilde\cO\left(\left(\cL +\frac{\sqrt{L\cL}}{\delta}\right)R_0^2\exp\left(-\min\left\{\frac{\mu}{\cL + \frac{\sqrt{L\cL}}{\delta}},\alpha\right\}K\right) + \frac{\widetilde{D}_1}{n\mu K} + \frac{L\widetilde{D}_1\left(\frac{\alpha(\omega+1)}{\delta}+1\right)}{\delta\mu^2 K^2}\right)
	\end{equation*}
	That is, to achive $\EE\left[f(\bar{x}^K) - f(x^*)\right] \le \varepsilon$ {\tt EC-SGDsr-DIANA} requires
	\begin{equation*}
		\widetilde{\cO}\left(\frac{1}{\alpha} + \frac{\cL}{\mu} + \frac{\sqrt{L\cL}}{\delta\mu} + \frac{\widetilde{D}_1}{n\mu\varepsilon} + \frac{\sqrt{L\widetilde{D}_1\left(\frac{\alpha(\omega+1)}{\delta}+1\right)}}{\mu\sqrt{\delta\varepsilon}}\right) \text{ iterations.}
	\end{equation*}
	In particular, if $\alpha = \frac{1}{\omega+1}$, then to achive $\EE\left[f(\bar{x}^K) - f(x^*)\right] \le \varepsilon$ {\tt EC-SGDsr-DIANA} requires
	\begin{equation*}
		\widetilde{\cO}\left(\omega + \frac{\cL}{\mu} + \frac{\sqrt{L\cL}}{\delta\mu} + \frac{\widetilde{D}_1}{n\mu\varepsilon} + \frac{\sqrt{L\widetilde{D}_1}}{\delta\mu\sqrt{\varepsilon}}\right) \text{ iterations,}
	\end{equation*}
	and if $\alpha = \frac{\delta}{\omega+1}$, then to achive $\EE\left[f(\bar{x}^K) - f(x^*)\right] \le \varepsilon$ {\tt EC-SGDsr-DIANA} requires
	\begin{equation*}
		\widetilde{\cO}\left(\frac{\omega+1}{\delta} + \frac{\cL}{\mu} + \frac{\sqrt{L\cL}}{\delta\mu} + \frac{\widetilde{D}_1}{n\mu\varepsilon} + \frac{\sqrt{L\widetilde{D}_1}}{\mu\sqrt{\delta\varepsilon}}\right) \text{ iterations.}
	\end{equation*}
\end{corollary}

Applying Lemma~\ref{lem:lemma_technical_cvx} we get the complexity result in the case when $\mu = 0$.
\begin{corollary}\label{cor:ec_sgdsr_diana_cvx_cor}
	Let the assumptions of Theorem~\ref{thm:ec_sgdsr_diana} hold and $\mu = 0$. Then after $K$ iterations of {\tt EC-SGDsr-DIANA} with the stepsize
	\begin{eqnarray*}
		\gamma_0 &=& \min\left\{\frac{1}{4\cL}, \frac{\delta}{4\sqrt{6L\left(4L + 3\delta(\cL + L) + \frac{16(3\cL + 4L)}{1-\alpha}\right)}}\right\},\quad R_0 = \|x^0-x^*\|,\\
		\gamma &=& \min\left\{\gamma_0, \sqrt[3]{\frac{R_0^2\delta^2\alpha\left(1-\min\left\{\frac{\gamma_0\mu}{2},\frac{\alpha}{4}\right\}\right)}{96L\sigma_0^2}}, \sqrt{\frac{3nR_0^2}{2\widetilde{D}_1 K}}, \sqrt[3]{\frac{\delta R_0^2}{6L\widetilde{D}_1\left(\frac{4\alpha(\omega+1)}{\delta}+1\right) K}}\right\},
	\end{eqnarray*}	
	and $\alpha \le \frac{1}{\omega+1}$ we have $\EE\left[f(\bar{x}^K) - f(x^*)\right]$ of order
	\begin{equation*}
		\cO\left(\frac{\cL R_0^2}{K} + \frac{\sqrt{\cL L} R_0^2}{\delta K} + \frac{\sqrt[3]{LR_0^4\sigma_0^2}}{K\sqrt[3]{\delta^2\alpha}} + \sqrt{\frac{R_0^2\widetilde{D}_1}{nK}}+ \sqrt[3]{\frac{LR_0^4\widetilde{D}_1\left(\frac{\alpha(\omega+1)}{\delta}+1\right)}{\delta K^2}}\right).
	\end{equation*}
	That is, to achive $\EE\left[f(\bar{x}^K) - f(x^*)\right] \le \varepsilon$ {\tt EC-SGDsr-DIANA} requires
	\begin{equation*}
		\cO\left(\frac{\cL R_0^2}{\varepsilon} + \frac{\sqrt{\cL L} R_0^2}{\delta \varepsilon} + \frac{\sqrt[3]{LR_0^4\sigma_0^2}}{\varepsilon\sqrt[3]{\delta^2\alpha}} + \frac{R_0^2\widetilde{D}_1}{n\varepsilon^2}+ \frac{R_0^2\sqrt{L\widetilde{D}_1\left(\frac{\alpha(\omega+1)}{\delta}+1\right)}}{\sqrt{\delta\varepsilon^3}}\right)
	\end{equation*}
	iterations. In particular, if $\alpha = \frac{1}{\omega+1}$, then to achive $\EE\left[f(\bar{x}^K) - f(x^*)\right] \le \varepsilon$ {\tt EC-SGDsr-DIANA} requires
	\begin{equation*}
		\cO\left(\frac{\cL R_0^2}{\varepsilon} + \frac{\sqrt{\cL L} R_0^2}{\delta \varepsilon} + \frac{\sqrt[3]{LR_0^4(\omega+1)\sigma_0^2}}{\varepsilon\sqrt[3]{\delta^2}} + \frac{R_0^2\widetilde{D}_1}{n\varepsilon^2}+ \frac{R_0^2\sqrt{L\widetilde{D}_1}}{\delta\sqrt{\varepsilon^3}}\right) \text{ iterations,}
	\end{equation*}
	and if $\alpha = \frac{\delta}{\omega+1}$, then to achive $\EE\left[f(\bar{x}^K) - f(x^*)\right] \le \varepsilon$ {\tt EC-SGDsr-DIANA} requires
	\begin{equation*}
		\cO\left(\frac{\cL R_0^2}{\varepsilon} + \frac{\sqrt{\cL L} R_0^2}{\delta \varepsilon} + \frac{\sqrt[3]{LR_0^4(\omega+1)\sigma_0^2}}{\delta\varepsilon} + \frac{R_0^2\widetilde{D}_1}{n\varepsilon^2}+ \frac{R_0^2\sqrt{L\widetilde{D}_1}}{\sqrt{\delta\varepsilon^3}}\right) \text{ iterations.}
	\end{equation*}
\end{corollary}

\subsection{{\tt EC-LSVRG}}\label{sec:ec_LSVRG}
In this section we consider problem \eqref{eq:main_problem_ef} with $f(x)$ being $\mu$-quasi strongly convex and $f_i(x)$ satisfying \eqref{eq:f_i_sum_ef} where functions $f_{ij}(x)$ are convex and $L$-smooth. For this problem we propose a new method called {\tt EC-LSVRG} which takes for the origin another method called {\tt LSVRG} (see \cite{hofmann2015variance,kovalev2019don}). 
\begin{algorithm}[t]
   \caption{{\tt EC-LSVRG}}\label{alg:ec-LSVRG}
\begin{algorithmic}[1]
   \Require learning rate $\gamma>0$, initial vector $x^0 \in \R^d$
   \State Set $e_i^0 = 0$ for all $i=1,\ldots, n$   
   \For{$k=0,1,\dotsc$}
       \State Broadcast $x^{k}$ to all workers
        \For{$i=1,\dotsc,n$ in parallel}
        	\State Pick $l$ uniformly at random from $[m]$
            \State Set $g^{k}_i = \nabla f_{il}(x^k) - \nabla f_{il}(w_i^k) + \nabla f_i(w_i^k)$
            \State $v_i^k = C(e_i^k + \gamma g_i^k)$
            \State $e_i^{k+1} = e_i^k + \gamma g_i^k - v_i^k$
            \State $w_i^{k+1} = \begin{cases}x^k,& \text{with probability } p,\\ w_i^k,& \text{with probability } 1-p\end{cases}$
        \EndFor
        \State $e^k = \frac{1}{n}\sum_{i=1}^ne_i^k$, $g^k = \frac{1}{n}\sum_{i=1}^ng_i^k$, $v^k = \frac{1}{n}\sum_{i=1}^nv_i^k$
       \State $x^{k+1} = x^k - v^k$
   \EndFor
\end{algorithmic}
\end{algorithm}

\begin{lemma}\label{lem:second_moment_bound_ec-LSVRG}
	For all $k\ge 0$, $i\in [n]$ we have
	\begin{equation}
		\bar{g}_i^k = \EE\left[g_i^k\mid x^k\right] = \nabla f_i(x^k) \label{eq:unbiasedness_g_i^k_ec-LSVRG}
	\end{equation}		
	and
	\begin{eqnarray}
		\frac{1}{n}\sum\limits_{i=1}^n\|\bar{g}_i^k\|^2 &\le& 4L\left(f(x^k) - f(x^*)\right) + D_1, \label{eq:second_moment_bound_ec-LSVRG}\\
		\frac{1}{n}\sum\limits_{i=1}^n\EE\left[\|g_i^k - \bar{g}_i^k\|^2\mid x^k\right] &\le& 12L\left(f(x^k) - f(x^*)\right) + 3\sigma_k^2, \label{eq:variance_bound_ec-LSVRG}\\
		\EE\left[\|g^k\|^2\mid x^k\right] &\le& 4L\left(f(x^k) - f(x^*)\right) + 2\sigma_k^2 \label{eq:second_moment_bound_ec-LSVRG_2}
	\end{eqnarray}
	where $\sigma_k^2 = \frac{1}{nm}\sum_{i=1}^n\sum_{j=1}^n\|\nabla f_{ij}(w_i^k) - \nabla f_{ij}(x^*)\|^2$ and $D_1 = \frac{2}{n}\sum_{i=1}^n\|\nabla f_{i}(x^*)\|^2$.
\end{lemma}
\begin{proof}
	First of all, we derive unbiasedness of $g_i^k$:
	\begin{equation*}
		\EE\left[g_i^k\mid x^k\right] = \frac{1}{m}\sum\limits_{j=1}^m\left(\nabla f_{ij}(x^k) - \nabla f_{ij}(w_i^k) + \nabla f_i(w_i^k)\right) = \nabla f_i(x^k).
	\end{equation*}
	Next, we get an upper bound for $\frac{1}{n}\sum\limits_{i=1}^n\|\bar{g}_i^k\|^2$:
	\begin{eqnarray*}
		\frac{1}{n}\sum\limits_{i=1}^n\|\bar{g}_i^k\|^2 &=& \frac{1}{n}\sum\limits_{i=1}^n\|\nabla f_i(x^k)\|^2\\
		&\overset{\eqref{eq:a_b_norm_squared}}{\le}& \frac{2}{n}\sum\limits_{i=1}^n\|\nabla f_i(x^k) - \nabla f_i(x^*)\|^2 + \frac{2}{n}\sum\limits_{i=1}^n\|\nabla f_i(x^*)\|^2\\
		&\overset{\eqref{eq:L_smoothness_cor_3}}{\le}& 4L\left(f(x^k)-f(x^*)\right) + \frac{2}{n}\sum\limits_{i=1}^n\|\nabla f_i(x^*)\|^2.
	\end{eqnarray*}
	Using \eqref{eq:unbiasedness_g_i^k_ec-LSVRG} we establish the following inequality:
	\begin{eqnarray*}
		\frac{1}{n}\sum\limits_{i=1}^n\EE\left[\|g_i^k - \bar{g}_i^k\|^2\mid x^k\right] &\overset{\eqref{eq:a_b_norm_squared}}{\le}& \frac{3}{n}\sum\limits_{i=1}^n\EE\left[\left\|\nabla f_{il}(w_i^k)- \nabla f_{il}(x^*) - \left(\nabla f_i(w_i^k) - \nabla f_i(x^*)\right)\right\|^2\mid x^k\right] \\
		&&\quad +  \frac{3}{n}\sum\limits_{i=1}^n\EE\left[\|\nabla f_{il}(x^k)- \nabla f_{il}(x^*)\|^2\mid x^k\right]\\
		&&\quad + \frac{3}{n}\sum\limits_{i=1}^n\|\nabla f_{i}(x^*)-\nabla f_i(x^k)\|^2\\
		&\overset{\eqref{eq:L_smoothness_cor_3},\eqref{eq:variance_decomposition}}{\le}& 12L\left(f(x^k)-f(x^*)\right) + \frac{3}{nm}\sum\limits_{i=1}^n\sum\limits_{j=1}^m\|\nabla f_{ij}(w_i^k) - \nabla f_{ij}(x^*)\|^2.
	\end{eqnarray*}
	Finally, we derive \eqref{eq:second_moment_bound_ec-LSVRG_2}:
	\begin{eqnarray*}
		\EE\left[\|g^k\|^2\mid x^k\right] &=& \EE\left[\left\|\frac{1}{n}\sum\limits_{i=1}^n\left(\nabla f_{il}(x^k)-\nabla f_{il}(w_i^k) + \nabla f_i(w_i^k) - \nabla f_i(x^*)\right)\right\|^2\mid x^k\right]\\
		&\overset{\eqref{eq:a_b_norm_squared}}{\le}& \frac{2}{n}\sum\limits_{i=1}^n\EE\left[\|\nabla f_{il}(x^k)-\nabla f_{il}(x^*)\|^2\mid x^k\right]\\
		&&\quad + \frac{2}{n}\sum\limits_{i=1}^n\EE\left[\left\|\nabla f_{il}(w_i^k)- \nabla f_{il}(x^*) - \left(\nabla f_i(w_i^k) - \nabla f_i(x^*)\right)\right\|^2\mid x^k\right]\\
		&=& \frac{2}{nm}\sum\limits_{i=1}^n\sum\limits_{j=1}^m\left\|\nabla f_{ij}(w_i^k) - \nabla f_{ij}(x^*) - \frac{1}{m}\sum\limits_{j=1}^m\left(\nabla f_{ij}(w_i^k) - \nabla f_{ij}(x^*)\right)\right\|^2\\
		&&\quad + \frac{2}{nm}\sum\limits_{i=1}^n\sum\limits_{j=1}^m\|\nabla f_{ij}(x^k) - \nabla f_{ij}(x^*)\|^2\\
		&\overset{\eqref{eq:L_smoothness_cor_3},\eqref{eq:variance_decomposition}}{\le}& 4L\left(f(x^k)-f(x^*)\right) + \frac{2}{nm}\sum\limits_{i=1}^n\sum\limits_{j=1}^m\left\|\nabla f_{ij}(w_i^k) - \nabla f_{ij}(x^*)\right\|^2.
	\end{eqnarray*}
\end{proof}

\begin{lemma}\label{lem:sigma_k+1_bound_ec-LSVRG}
	For all $k\ge 0$, $i\in [n]$ we have
	\begin{equation}
		\EE\left[\sigma_{k+1}^2\mid x^k\right] \le (1-p)\sigma_k^2 + 2Lp\left(f(x^k) - f(x^*)\right), \label{eq:sigma_k+1_ec-LSVRG} 
	\end{equation}
	where $\sigma_k^2 = \frac{1}{nm}\sum_{i=1}^n\sum_{j=1}^n\|\nabla f_{ij}(w_i^k) - \nabla f_{ij}(x^*)\|^2$.
\end{lemma}
\begin{proof}
	By definition of $w_i^{k+1}$ we get
	\begin{eqnarray*}
		\EE\left[\sigma_{k+1}^2\mid x^k\right] &=& \frac{1}{nm}\sum\limits_{i=1}^n\sum\limits_{j=1}^m\EE\left[\|\nabla f_{ij}(w_i^{k+1}) - \nabla f_{ij}(x^*)\|^2\mid x^k\right]\\
		&=& \frac{1-p}{nm}\sum\limits_{i=1}^n\sum\limits_{j=1}^m\|\nabla f_{ij}(w_i^{k}) - \nabla f_{ij}(x^*)\|^2 + \frac{p}{nm}\sum\limits_{i=1}^n\sum\limits_{j=1}^m\|\nabla f_{ij}(x^{k}) - \nabla f_{ij}(x^*)\|^2\\
		&\overset{\eqref{eq:L_smoothness_cor_3}}{\le}& (1-p)\sigma_k^2 + \frac{2Lp}{nm}\sum\limits_{i=1}^n\sum\limits_{j=1}^m D_{f_{ij}}(x^k,x^*)\\
		&=& (1-p)\sigma_k^2 + 2Lp\left(f(x^k) - f(x^*)\right).
	\end{eqnarray*}
\end{proof}

Applying Theorem~\ref{thm:ec_sgd_main_result_new} we get the following result.
\begin{theorem}\label{thm:ec_LSVRG}
	Assume that $f(x)$ is $\mu$-quasi strongly convex and functions $f_{ij}$ are convex and $L$-smooth for all $i\in[n],j\in[m]$. Then {\tt EC-LSVRG} satisfies Assumption~\ref{ass:key_assumption_finite_sums_new} with
	\begin{gather*}
		A = 2L,\quad \widetilde{A} = 12L,\quad A' = 2L,\quad B_1 = \widetilde{B}_1 = B_1' = B_2 = 0,\quad D_1 = \frac{2}{n}\sum\limits_{i=1}^n\|\nabla f_{i}(x^*)\|^2,\\
		D_1' = \widetilde{D}_1 = 0,\quad \widetilde{B}_2 = 3,\quad B_2' = 2,\quad \sigma_{1,k}^2 \equiv 0,\quad C_1 = 0,\\
		\sigma_{2,k}^2 = \sigma_{k}^2 = \frac{1}{nm}\sum\limits_{i=1}^n\sum\limits_{j=1}^m\|\nabla f_{ij}(w_i^{k}) - \nabla f_{ij}(x^*)\|^2,\quad \rho_1 = 1,\quad \rho_2 = p,\quad C_2 = Lp,\quad D_2 = 0,\\
		G = 0,\quad F_1 = 0,\quad F_2 = \frac{72L\gamma^2}{\delta p\left(1-\min\left\{\frac{\gamma\mu}{2},\frac{p}{4}\right\}\right)},\quad D_3 = \frac{12L\gamma}{\delta^2}D_1,
	\end{gather*}
	with $\gamma$ satisfying
	\begin{equation*}
		\gamma \le \min\left\{\frac{1}{24L}, \frac{\delta}{8L\sqrt{3\left(2 + 3\delta\left(2+\frac{1}{1-p}\right)\right)}}\right\}, \quad M_2 = \frac{4}{p}.
	\end{equation*}
	and for all $K \ge 0$
	\begin{equation*}
		\EE\left[f(\bar x^K) - f(x^*)\right] \le \left(1 - \min\left\{\frac{\gamma\mu}{2},\frac{p}{4}\right\}\right)^K\frac{4(T^0 + \gamma F_2 \sigma_0^2)}{\gamma} + \frac{48L\gamma^2}{\delta^2}D_1
	\end{equation*}
	when $\mu > 0$ and
	\begin{equation*}
		\EE\left[f(\bar x^K) - f(x^*)\right] \le \frac{4(T^0 + \gamma F_2 \sigma_0^2)}{\gamma K} + \frac{48L\gamma^2}{\delta^2}D_1
	\end{equation*}
	when $\mu = 0$, where $T^k \eqdef \|x^k - x^*\|^2 + M_2\gamma^2 \sigma_k^2$.
\end{theorem}

In other words, {\tt EC-LSVRG} converges with linear rate $\cO\left(\left(\frac{1}{p} + \frac{\kappa}{\delta\sqrt{1-p}}\right)\ln\frac{1}{\varepsilon}\right)$ to the neighbourhood of the solution. If $m\ge 2$ then taking $p = \frac{1}{m}$ we get that in expectation the sample complexity of one iteration of {\tt EC-LSVRG} is $\cO(1)$ gradients calculations per node as for {\tt EC-SGDsr} with standard sampling and the rate of convergence to the neighbourhood becomes $\cO\left(\left(m + \frac{\kappa}{\delta}\right)\ln\frac{1}{\varepsilon}\right)$. We notice that the size of this neighbourhood is typically smaller than for {\tt EC-SGDsr}, but still the method fails to converge to the exact solution with linear rate. Applying Lemma~\ref{lem:lemma2_stich} we establish the rate of convergence to $\varepsilon$-solution in the case when $\mu > 0$.
\begin{corollary}\label{cor:ec_lsvrg_str_cvx_cor}
	Let the assumptions of Theorem~\ref{thm:ec_LSVRG} hold and $\mu > 0$. Then after $K$ iterations of {\tt EC-LSVRG} with the stepsize
	\begin{eqnarray*}
		\gamma_0 &=&  \min\left\{\frac{1}{24L}, \frac{\delta}{8L\sqrt{3\left(2 + 3\delta\left(2+\frac{1}{1-p}\right)\right)}}\right\},\\
		\tilde{T^0} &=& \|x^0-x^*\|^2 + M_2\gamma_0^2\sigma_0^2,\quad  \tilde{F_2} = \frac{72L\gamma_0^2}{\delta p\left(1-\min\left\{\frac{\gamma_0\mu}{2},\frac{p}{4}\right\}\right)},\\
		\gamma &=& \min\left\{\gamma_0, \frac{\ln\left(\max\left\{2,\frac{\delta^2\left(\tilde{T^0} + \tilde{F}_2\gamma_0\sigma_{0}^2\right)\mu^3K^3}{48LD_1}\right\}\right)}{\mu K}\right\},
	\end{eqnarray*}	 
	and $p = \frac{1}{m}$, $m\ge 2$ we have
	\begin{equation*}
		\EE\left[f(\bar{x}^K) - f(x^*)\right] = \widetilde\cO\left(\frac{L}{\delta}\left(\tilde{T^0} + \tilde{F}_2\gamma_0\sigma_{0}^2\right)\exp\left(-\min\left\{\frac{\delta\mu}{L},\frac{1}{m}\right\}K\right) + \frac{LD_1}{\delta^2\mu^2 K^2}\right).
	\end{equation*}
	That is, to achive $\EE\left[f(\bar{x}^K) - f(x^*)\right] \le \varepsilon$ {\tt EC-LSVRG} requires
	\begin{equation*}
		\widetilde{\cO}\left(m + \frac{L}{\delta\mu} + \frac{\sqrt{LD_1}}{\delta\mu\sqrt{\varepsilon}}\right) \text{ iterations.}
	\end{equation*}
\end{corollary}

Applying Lemma~\ref{lem:lemma_technical_cvx} we get the complexity result in the case when $\mu = 0$.
\begin{corollary}\label{cor:ec_lsvrg_cvx_cor}
	Let the assumptions of Theorem~\ref{thm:ec_LSVRG} hold and $\mu = 0$. Then after $K$ iterations of {\tt EC-LSVRG} with the stepsize
	\begin{eqnarray*}
		\gamma_0 &=& \min\left\{\frac{1}{24L}, \frac{\delta}{8L\sqrt{3\left(2 + 3\delta\left(2+\frac{1}{1-p}\right)\right)}}\right\},\quad R_0 = \|x^0-x^*\|,\\
		\gamma &=& \min\left\{\gamma_0, \sqrt{\frac{R_0^2p}{4\sigma_0^2}}, \sqrt[3]{\frac{R_0^2\delta p\left(1-\min\left\{\frac{\gamma_0\mu}{2},\frac{p}{4}\right\}\right)}{72L\sigma_0^2}}, \sqrt[3]{\frac{\delta^2 R_0^2}{12LD_1 K}}\right\},
	\end{eqnarray*}	
	and $p = \frac{1}{m}$, $m\ge 2$ we have $\EE\left[f(\bar{x}^K) - f(x^*)\right]$ of order
	\begin{equation*}
		\cO\left(\frac{L R_0^2}{\delta K} + \frac{\sqrt{m R_0^2\sigma_0^2}}{K} + \frac{\sqrt[3]{LR_0^4m\sigma_0^2}}{\sqrt[3]{\delta} K} + \frac{\sqrt[3]{LR_0^4}}{(\delta K)^{\nicefrac{2}{3}}}\sqrt[3]{\frac{1}{n}\sum\limits_{i=1}^n\|\nabla f_i(x^*)\|^2}\right).
	\end{equation*}
	That is, to achive $\EE\left[f(\bar{x}^K) - f(x^*)\right] \le \varepsilon$ {\tt EC-LSVRG} requires
	\begin{equation*}
		\cO\left(\frac{L R_0^2}{\delta\varepsilon} + \frac{\sqrt{m R_0^2\sigma_0^2}}{\varepsilon} + \frac{\sqrt[3]{LR_0^4m\sigma_0^2}}{\sqrt[3]{\delta} \varepsilon}+ \frac{R_0^2}{\delta\varepsilon^{\nicefrac{3}{2}}}\sqrt{\frac{L}{n}\sum\limits_{i=1}^n\|\nabla f_i(x^*)\|^2}\right)
	\end{equation*}
	iterations.
\end{corollary}

\subsection{{\tt EC-LSVRGstar}}\label{sec:ec_LSVRGstar}
In the setup of Section~\ref{sec:ec_LSVRG} we now assume that $i$-th node has an access to the $\nabla f_i(x^*)$. Under this unrealistic assumption we construct the method called {\tt EC-LSVRGstar} that asymptotically converges to the exact solution.

\begin{algorithm}[t]
   \caption{{\tt EC-LSVRGstar}}\label{alg:ec-LSVRGstar}
\begin{algorithmic}[1]
   \Require learning rate $\gamma>0$, initial vector $x^0 \in \R^d$
   \State Set $e_i^0 = 0$ for all $i=1,\ldots, n$   
   \For{$k=0,1,\dotsc$}
       \State Broadcast $x^{k}$ to all workers
        \For{$i=1,\dotsc,n$ in parallel}
        	\State Pick $l$ uniformly at random from $[m]$
            \State Set $g^{k}_i = \nabla f_{il}(x^k) - \nabla f_{il}(w_i^k) + \nabla f_i(w_i^k) - \nabla f_i(x^*)$
            \State $v_i^k = C(e_i^k + \gamma g_i^k)$
            \State $e_i^{k+1} = e_i^k + \gamma g_i^k - v_i^k$
            \State $w_i^{k+1} = \begin{cases}x^k,& \text{with probability } p,\\ w_i^k,& \text{with probability } 1-p\end{cases}$
        \EndFor
        \State $e^k = \frac{1}{n}\sum_{i=1}^ne_i^k$, $g^k = \frac{1}{n}\sum_{i=1}^ng_i^k$, $v^k = \frac{1}{n}\sum_{i=1}^nv_i^k$
       \State $x^{k+1} = x^k - v^k$
   \EndFor
\end{algorithmic}
\end{algorithm}

\begin{lemma}\label{lem:second_moment_bound_ec-LSVRGstar}
	For all $k\ge 0$, $i\in [n]$ we have
	\begin{equation}
		\EE\left[g^k\mid x^k\right] = \nabla f(x^k) \label{eq:unbiasedness_g^k_ec-LSVRGstar}
	\end{equation}		
	and
	\begin{equation}
		\frac{1}{n}\sum\limits_{i=1}^n\|\bar{g}_i^k\|^2 \le 2L\left(f(x^k) - f(x^*)\right), \label{eq:second_moment_bound_ec-LSVRGstar} 
	\end{equation}
	\begin{equation}
		\frac{1}{n}\sum\limits_{i=1}^n\EE\left[\|g_i^k-\bar{g}_i^k\|^2\mid x^k\right] \le 4L\left(f(x^k) - f(x^*)\right) + 2\sigma_k^2, \label{eq:variance_bound_ec-LSVRGstar} 
	\end{equation}
	\begin{equation}
		\EE\left[\|g^k\|^2\mid x^k\right] \le 4L\left(f(x^k) - f(x^*)\right) + 2\sigma_k^2, \label{eq:second_moment_bound_ec-LSVRGstar_2} 
	\end{equation}
	where $\sigma_k^2 = \frac{1}{nm}\sum_{i=1}^n\sum_{j=1}^n\|\nabla f_{ij}(w_i^k) - \nabla f_{ij}(x^*)\|^2$.
\end{lemma}
\begin{proof}
	First of all, we derive unbiasedness of $g^k$:
	\begin{eqnarray*}
		\EE\left[g^k\mid x^k\right] &=& \frac{1}{n}\sum\limits_{i=1}^n\EE\left[\nabla f_{il}(x^k) - \nabla f_{il}(w_i^k) + \nabla f_i(w_i^k) - \nabla f_i(x^*)\mid x^k\right]\\
		&=& \frac{1}{nm}\sum\limits_{i=1}^n\sum\limits_{j=1}^m\left(\nabla f_{ij}(x^k) - \nabla f_{ij}(w_i^k) + \nabla f_i(w_i^k) - \nabla f_i(x^*)\right)\\
		&=& \nabla f(x^k) + \frac{1}{n}\sum\limits_{i=1}^n\left(-\nabla f_i(w_i^k) + \nabla f_i(w_i^k)\right) - \nabla f(x^*) = \nabla f(x^k).
	\end{eqnarray*}
	Next, we get an upper bound for $\frac{1}{n}\sum\limits_{i=1}^n\|\bar{g}_i^k\|^2$:
	\begin{equation*}
		\frac{1}{n}\sum\limits_{i=1}^n\|\bar{g}_i^k\|^2 = \frac{1}{n}\sum\limits_{i=1}^n\|\nabla f_i(x^k)-\nabla f_i(x^*)\|^2 \overset{\eqref{eq:L_smoothness_cor_3}}{\le} 2L\left(f(x^k)-f(x^*)\right).
	\end{equation*}		
	Since the variance of random vector is not greater than its second moment we obtain:
	\begin{eqnarray}
		\frac{1}{n}\sum\limits_{i=1}^n\EE\left[\|g_i^k-\bar{g}_i^k\|^2\mid x^k\right] &\overset{\eqref{eq:variance_decomposition}}{\le}& \frac{1}{n}\sum\limits_{i=1}^n\EE\left[\|g_i^k\|^2\mid x^k\right]\notag\\
		&\overset{\eqref{eq:a_b_norm_squared}}{\le}& \frac{2}{n}\sum\limits_{i=1}^n\EE\left[\left\|\nabla f_{il}(w_i^k)- \nabla f_{il}(x^*) - \left(\nabla f_i(w_i^k) - \nabla f_i(x^*)\right)\right\|^2\mid x^k\right] \notag\\
		&&\quad + \frac{2}{n}\sum\limits_{i=1}^n\EE\left[\|\nabla f_{il}(x^k)- \nabla f_{il}(x^*)\|^2\mid x^k\right]\notag\\
		&\overset{\eqref{eq:L_smoothness_cor_3},\eqref{eq:variance_decomposition}}{\le}& 4L\left(f(x^k)-f(x^*)\right) + \frac{2}{nm}\sum\limits_{i=1}^n\sum\limits_{j=1}^m\|\nabla f_{ij}(w_i^k) - \nabla f_{ij}(x^*)\|^2.\notag
	\end{eqnarray}
	Inequality \eqref{eq:second_moment_bound_ec-LSVRGstar_2} trivially follows from the inequality above by Jensen's inequality and convexity of $\|\cdot\|^2$. 
\end{proof}

\begin{lemma}\label{lem:sigma_k+1_bound_ec-LSVRGstar}
	For all $k\ge 0$, $i\in [n]$ we have
	\begin{equation}
		\EE\left[\sigma_{k+1}^2\mid x^k\right] \le (1-p)\sigma_k^2 + 2Lp\left(f(x^k) - f(x^*)\right), \label{eq:sigma_k+1_ec-LSVRGstar} 
	\end{equation}
	where $\sigma_k^2 = \frac{1}{nm}\sum_{i=1}^n\sum_{j=1}^n\|\nabla f_{ij}(w_i^k) - \nabla f_{ij}(x^*)\|^2$.
\end{lemma}
\begin{proof} The proof of this lemma is identical to the proof of Lemma~\ref{lem:sigma_k+1_bound_ec-LSVRG}.
\end{proof}

Applying Theorem~\ref{thm:ec_sgd_main_result_new} we get the following result.
\begin{theorem}\label{thm:ec_LSVRGstar}
	Assume that $f(x)$ is $\mu$-quasi strongly convex and functions $f_{ij}$ are convex and $L$-smooth for all $i\in[n],j\in[m]$. Then {\tt EC-LSVRGstar} satisfies Assumption~\ref{ass:key_assumption_finite_sums_new} with
	\begin{gather*}
		A = L,\quad \widetilde{A} = A' = 2L,\quad B_1 = \widetilde{B}_1 = B_1' = B_2 = 0,\quad \widetilde{B}_2 = B_2' = 2,\quad D_1 = D_1' = 0,\\
		\sigma_{1,k}^2 \equiv 0,,\quad C_1 = 0,\quad\sigma_{2,k}^2 = \sigma_{k}^2 = \frac{1}{nm}\sum\limits_{i=1}^n\sum\limits_{j=1}^m\|\nabla f_{ij}(w_i^{k}) - \nabla f_{ij}(x^*)\|^2,\quad \rho_1 = 1,\\
		\rho_2 = p,\quad C_2 = Lp,\quad D_2 = 0,\quad G = 0,\quad F_1 = 0,\quad F_2 = \frac{48L\gamma^2(2+p)}{\delta p},\quad D_3 = 0,
	\end{gather*}
	with $\gamma$ satisfying
	\begin{equation*}
		\gamma \le \min\left\{\frac{3}{56L}, \frac{\delta}{8L\sqrt{3\left(1+\delta\left(1 + \frac{2}{1-p}\right)\right)}}\right\}, \quad M_2 = \frac{8}{3p}.
	\end{equation*}
	and for all $K \ge 0$
	\begin{equation*}
		\EE\left[f(\bar x^K) - f(x^*)\right] \le \left(1 - \min\left\{\frac{\gamma\mu}{2},\frac{p}{4}\right\}\right)^K\frac{4(T^0 + \gamma F_2 \sigma_0^2)}{\gamma}
	\end{equation*}	
	when $\mu > 0$ and
	\begin{equation*}
		\EE\left[f(\bar x^K) - f(x^*)\right] \le \frac{4(T^0 + \gamma F_2 \sigma_0^2)}{\gamma K}
	\end{equation*}
	when $\mu = 0$, where $T^k \eqdef \|x^k - x^*\|^2 + M_2\gamma^2 \sigma_k^2$.
\end{theorem}

In other words, {\tt EC-LSVRGstar} converges with linear rate $\cO\left(\left(\frac{1}{p} + \frac{\kappa}{\delta\sqrt{1-p}}\right)\ln\frac{1}{\varepsilon}\right)$ exactly to the solution when $\mu > 0$. If $m\ge 2$ then taking $p = \frac{1}{m}$ we get that in expectation the sample complexity of one iteration of {\tt EC-LSVRGstar} is $\cO(1)$ gradients calculations per node as for {\tt EC-SGDsr} with standard sampling and the rate of convergence becomes $\cO\left(\left(m + \frac{\kappa}{\delta}\right)\ln\frac{1}{\varepsilon}\right)$. 

Applying Lemma~\ref{lem:lemma_technical_cvx} we get the complexity result in the case when $\mu = 0$.
\begin{corollary}\label{cor:ec_lsvrg_star_cvx_cor}
	Let the assumptions of Theorem~\ref{thm:ec_LSVRGstar} hold and $\mu = 0$. Then after $K$ iterations of {\tt EC-LSVRGstar} with the stepsize
	\begin{eqnarray*}
		\gamma_0 &=& \min\left\{\frac{3}{56L}, \frac{\delta}{8L\sqrt{3\left(1+\delta\left(1 + \frac{2}{1-p}\right)\right)}}\right\},\quad R_0 = \|x^0-x^*\|,\\
		\gamma &=& \min\left\{\gamma_0, \sqrt{\frac{3pR_0^2}{8\sigma_0^2}}, \sqrt[3]{\frac{R_0^2\delta p\left(1-\min\left\{\frac{\gamma_0\mu}{2},\frac{p}{4}\right\}\right)}{72L\sigma_0^2}}\right\},
	\end{eqnarray*}	
	and $p = \frac{1}{m}$, $m\ge 2$ we have $\EE\left[f(\bar{x}^K) - f(x^*)\right]$ of order
	\begin{equation*}
		\cO\left(\frac{L R_0^2}{\delta K} + \frac{\sqrt{R_0^2m\sigma_0^2}}{K} + \frac{\sqrt[3]{LR_0^4m\sigma_0^2}}{\sqrt[3]{\delta}K}\right).
	\end{equation*}
	That is, to achive $\EE\left[f(\bar{x}^K) - f(x^*)\right] \le \varepsilon$ {\tt EC-LSVRGstar} requires
	\begin{equation*}
		\cO\left(\frac{L R_0^2}{\delta \varepsilon} + \frac{\sqrt{R_0^2m\sigma_0^2}}{\varepsilon} + \frac{\sqrt[3]{LR_0^4m\sigma_0^2}}{\sqrt[3]{\delta}\varepsilon}\right)
	\end{equation*}
	iterations.
\end{corollary}
However, such convergence guarantees are obtained under very restrictive assumption: the method requires to know vectors $\nabla f_i(x^*)$.

\subsection{{\tt EC-LSVRG-DIANA}}\label{sec:ec_LSVRG-diana}
In the setup of Section~\ref{sec:ec_LSVRG} we construct a new method called {\tt EC-LSVRG-DIANA} which does not require to know $\nabla f_i(x^*)$ and has linear convergence to the exact solution.
\begin{algorithm}[t]
   \caption{{\tt EC-LSVRG-DIANA}}\label{alg:ec-LSVRG-diana}
\begin{algorithmic}[1]
   \Require learning rates $\gamma>0$, $\alpha \in (0,1]$, initial vectors $x^0, h_1^0,\ldots, h_n^0 \in \R^d$
	\State Set $e_i^0 = 0$ for all $i=1,\ldots, n$   
	\State Set $h^0 = \frac{1}{n}\sum_{i=1}^n h_i^0$   
   \For{$k=0,1,\dotsc$}
       \State Broadcast $x^{k}, h^k$ to all workers
        \For{$i=1,\dotsc,n$ in parallel}
			\State Pick $l$ uniformly at random from $[m]$
            \State Set $\hat g^{k}_i = \nabla f_{il}(x^k) - \nabla f_{il}(w_i^k) + \nabla f_i(w_i^k)$           
            \State $g^{k}_i = \hat g_i^k - h_i^k + h^k$
            \State $v_i^k = C(e_i^k + \gamma g_i^k)$
            \State $e_i^{k+1} = e_i^k + \gamma g_i^k - v_i^k$
            \State $h_i^{k+1} = h_i^k + \alpha Q(\hat g_i^k - h_i^k)$
            \State $w_i^{k+1} = \begin{cases}x^k,& \text{with probability } p,\\ w_i^k,& \text{with probability } 1-p\end{cases}$
        \EndFor
        \State $e^k = \frac{1}{n}\sum_{i=1}^ne_i^k$, $g^k = \frac{1}{n}\sum_{i=1}^ng_i^k$, $v^k = \frac{1}{n}\sum_{i=1}^nv_i^k$, $h^{k+1} = \frac{1}{n}\sum\limits_{i=1}^n h_i^{k+1} = h^k + \alpha\frac{1}{n}\sum\limits_{i=1}^n Q(\hat g_i^k - h_i^k)$
       \State $x^{k+1} = x^k - v^k$
   \EndFor
\end{algorithmic}
\end{algorithm}
As in {\tt EC-SGD-DIANA} the master needs to gather only $C(e_i^k + \gamma g_i^k)$ and $Q(\hat g_i^k - h_i^k)$ from all nodes in order to perform an update.

\begin{lemma}\label{lem:ec_LSVRG-diana_second_moment_bound}
	Assume that $f_{ij}(x)$ is convex and $L$-smooth for all $i=1,\ldots,n$, $j=1,\ldots,m$. Then, for all $k\ge 0$ we have
	\begin{eqnarray}
		\EE\left[g^k\mid x^k\right] &=& \nabla f(x^k), \label{eq:ec_LSVRG-diana_unbiasedness}\\
		\frac{1}{n}\sum\limits_{i=1}^n\|\bar{g}_i^k\|^2 &\le& 4L\left(f(x^k) - f(x^*)\right) + 2\sigma_{1,k}^2, \label{eq:ec_LSVRG-diana_second_moment_bound}\\
		\frac{1}{n}\sum\limits_{i=1}^n\EE\left[\|g_i^k-\bar{g}_i^k\|^2\mid x^k\right] &\le& 6L\left(f(x^k) - f(x^*)\right) + 3\sigma_{1,k}^2 + 3\sigma_{2,k}^2, \label{eq:ec_LSVRG-diana_variance_bound}\\
		\EE\left[\|g^k\|^2\mid x^k\right] &\le& 4L\left(f(x^k) - f(x^*)\right) + 2\sigma_{2,k}^2 \label{eq:ec_LSVRG-diana_second_moment_bound_2}
	\end{eqnarray}
	where $$\sigma_{1,k}^2 = \frac{1}{n}\sum_{i=1}^n\|h_i^k - \nabla f(x^*)\|^2,\quad \sigma_{2,k}^2 =  \frac{1}{nm}\sum_{i=1}^n\sum_{j=1}^m\|\nabla f_{ij}(w_i^k) - \nabla f_{ij}(x^*)\|^2.$$
\end{lemma}
\begin{proof}
	First of all, we show unbiasedness of $g^k$:
	\begin{eqnarray*}
		\EE\left[g^k\mid x^k\right] &=& \frac{1}{n}\sum\limits_{i=1}^n\EE\left[\hat g_i^k - h_i^k + h^k\mid x^k\right]\\
		&=& \frac{1}{nm}\sum\limits_{i=1}^n\sum\limits_{j=1}^m\left(\nabla f_{ij}(x^k) - \nabla f_{ij}(w_i^k) + \nabla f_i(w_i^k) - h_i^k + h^k\right) = \nabla f(x^k).	
	\end{eqnarray*}
	Next, we derive the upper bound for $\frac{1}{n}\sum\limits_{i=1}^n\|\bar{g}_i^k\|^2$:
	\begin{eqnarray*}
		\frac{1}{n}\sum\limits_{i=1}^n\|\bar{g}_i^k\|^2 &=& \frac{1}{n}\sum\limits_{i=1}^n\|\nabla f_i(x^k)-h_i^k + h^k\|^2\\
		&\overset{\eqref{eq:a_b_norm_squared}}{\le}& \frac{2}{n}\sum\limits_{i=1}^n\|\nabla f_i(x^k)-\nabla f_i(x^*)\|^2 + \frac{2}{n}\sum\limits_{i=1}^n\left\|h_i^k - \nabla f_i(x^*)-\left(h^k-\nabla f(x^*)\right)\right\|^2\\
		&\overset{\eqref{eq:L_smoothness_cor_3},\eqref{eq:variance_decomposition}}{\le}& 4L\left(f(x^k)-f(x^*)\right) + \frac{2}{n}\sum\limits_{i=1}^n\|h_i^k-\nabla f_i(x^*)\|^2.
	\end{eqnarray*}
	Since the variance of random vector is not greater than its second moment we obtain:
	\begin{eqnarray*}
		\frac{1}{n}\sum\limits_{i=1}^n\EE\left[\|g_i^k-\bar{g}_i^k\|^2\mid x^k\right] &\overset{\eqref{eq:variance_decomposition}}{\le}& \frac{1}{n}\sum\limits_{i=1}^n\EE\left[\|g_i^k\|^2\mid x^k\right]\\
		&=& \frac{1}{n}\sum\limits_{i=1}^n\EE\left[\|\nabla f_{il}(x^k) - \nabla f_{il}(w_i^k) + \nabla f_i(w_i^k) - h_i^k + h^k\|^2\mid x^k\right]\\
		&\overset{\eqref{eq:a_b_norm_squared}}{\le}& \frac{3}{n}\sum\limits_{i=1}^n\EE\left[\left\|\nabla f_{il}(w_i^k) - \nabla f_{il}(x^*) -\left(\nabla f_i(w_i^k) - \nabla f_i(x^*)\right)\right\|^2\mid x^k\right]\\
		&&\quad + \frac{3}{n}\sum\limits_{i=1}^n\EE\left[\left\|\nabla f_{il}(x^k) - \nabla f_{il}(x^*)\right\|^2\mid x^k\right]\\
		&&\quad + \frac{3}{n}\sum\limits_{i=1}^n\left\|h_i^k - \nabla f_{i}(x^*) -\left(h^k - \nabla f(x^*)\right)\right\|^2\\
		&\overset{\eqref{eq:L_smoothness_cor_3},\eqref{eq:variance_decomposition}}{\le}& 6L\left(f(x^k)-f(x^*)\right) + \frac{3}{nm}\sum\limits_{i=1}^n\sum\limits_{j=1}^m\|\nabla f_{ij}(w_i^k) - \nabla f_{ij}(x^*)\|^2\\
		&&\quad + \frac{3}{n}\sum\limits_{i=1}^n\left\|h_i^k - \nabla f_{i}(x^*)\right\|^2.
	\end{eqnarray*}
	Finally, we obtain an upper boud for the second moment of $g^k$:
	\begin{eqnarray*}
		\EE\left[\|g^k\|^2\mid x^k\right] &=& \EE\left[\left\|\frac{1}{n}\sum\limits_{i=1}^n\left(\nabla f_{il}(x^k)-\nabla f_{il}(w_i^k) + \nabla f_i(w_i^k) - \nabla f_i(x^*)\right)\right\|^2\mid x^k\right]\\
		&\overset{\eqref{eq:a_b_norm_squared}}{\le}& \frac{2}{n}\sum\limits_{i=1}^n\EE\left[\|\nabla f_{il}(x^k)-\nabla f_{il}(x^*)\|^2\mid x^k\right]\\
		&&\quad + \frac{2}{n}\sum\limits_{i=1}^n\EE\left[\left\|\nabla f_{il}(w_i^k)- \nabla f_{il}(x^*) - \left(\nabla f_i(w_i^k) - \nabla f_i(x^*)\right)\right\|^2\mid x^k\right]\\
		&=& \frac{2}{nm}\sum\limits_{i=1}^n\sum\limits_{j=1}^m\left\|\nabla f_{ij}(w_i^k) - \nabla f_{ij}(x^*) - \frac{1}{m}\sum\limits_{j=1}^m\left(\nabla f_{ij}(w_i^k) - \nabla f_{ij}(x^*)\right)\right\|^2\\
		&&\quad + \frac{2}{nm}\sum\limits_{i=1}^n\sum\limits_{j=1}^m\|\nabla f_{ij}(x^k) - \nabla f_{ij}(x^*)\|^2\\
		&\overset{\eqref{eq:L_smoothness_cor_3},\eqref{eq:variance_decomposition}}{\le}& 4L\left(f(x^k)-f(x^*)\right) + \frac{2}{nm}\sum\limits_{i=1}^n\sum\limits_{j=1}^m\left\|\nabla f_{ij}(w_i^k) - \nabla f_{ij}(x^*)\right\|^2.
	\end{eqnarray*}
\end{proof}

\begin{lemma}\label{lem:ec_LSVRG-diana_sigma_k+1_bound}
	Assume that $\alpha \le \nicefrac{1}{(\omega+1)}$. Then, for all $k\ge 0$ we have
	\begin{equation}
		\EE\left[\sigma_{1,k+1}^2\mid x^k\right] \le (1 - \alpha)\sigma_{1,k}^2 + 6L\alpha(f(x^k) - f(x^*)) + 2\alpha\sigma_{2,k}^2, \label{eq:ec_LSVRG-diana_sigma_k+1_bound}
	\end{equation}
	\begin{equation}
		\EE\left[\sigma_{2,{k+1}}^2\mid x^k\right] \le (1 - p)\sigma_{k,2}^2 + 2Lp\left(f(x^k)-f(x^*)\right)
	\end{equation}
	where $\sigma_{1,k}^2 = \frac{1}{n}\sum_{i=1}^n\|h_i^k - \nabla f_i(x^*)\|^2$ and $\sigma_{2,k}^2= \frac{1}{nm}\sum_{i=1}^n\sum_{j=1}^m\|\nabla f_{ij}(w_i^k) - \nabla f_{ij}(x^*)\|^2$.
\end{lemma}
\begin{proof}
	First of all, we derive an upper bound for the second moment of $h_i^{k+1} - h_i^*$:
	\begin{eqnarray*}
		\EE\left[\|h_i^{k+1} - h_i^*\|^2\mid x^k\right] &=& \EE\left[\left\|h_i^k - h_i^* + \alpha Q(\hat g_i^k - h_i^k) \right\|^2\mid x^k\right]\\
		&\overset{\eqref{eq:quantization_def}}{=}& \|h_i^k - h_i^*\|^2 +2\alpha\langle h_i^k - h_i^*, \nabla f_i(x^k) - h_i^k \rangle\\
		&&\quad + \alpha^2\EE\left[\|Q(\hat g_i^k - h_i^k)\|^2\mid x^k\right]\\
		&\overset{\eqref{eq:quantization_def},\eqref{eq:tower_property}}{\le}& \|h_i^k - h_i^*\|^2 +2\alpha\langle h_i^k - h_i^*, \nabla f_i(x^k) - h_i^k \rangle\\
		&&\quad + \alpha^2(\omega+1)\EE\left[\|\hat g_i^k - h_i^k\|^2\mid x^k\right].
	\end{eqnarray*}
	Using variance decomposition \eqref{eq:variance_decomposition} and $\alpha \le \nicefrac{1}{(\omega+1)}$ we get
	\begin{eqnarray*}
		\alpha^2(\omega+1)\EE\left[\|\hat g_i^k - h_i^k\|^2\mid x^k\right] &\overset{\eqref{eq:variance_decomposition}}{=}& \alpha^2(\omega+1)\EE\left[\|\hat g_i^k - \nabla f_i(x^k)\|^2\mid x^k\right] + \alpha^2(\omega+1)\|\nabla f_i(x^k) - h_i^k\|^2\\
		&\le& \alpha\EE\left[\|\hat g_i^k - \nabla f_i(x^k)\|^2\mid x^k\right] + \alpha\|\nabla f_i(x^k) - h_i^k\|^2\\
		&\overset{\eqref{eq:a_b_norm_squared}}{\le}& 2\alpha\EE\left[\left\|\nabla f_{il}(x^k) - \nabla f_{il}(x^*) -\left(\nabla f_i(x^k) - \nabla f_i(x^*)\right)\right\|^2\mid x^k\right]\\
		&&\quad + 2\alpha\EE\left[\left\|\nabla f_{il}(w_i^k) - \nabla f_{il}(x^*) -\left(\nabla f_i(w_i^k) - \nabla f_i(x^*)\right)\right\|^2\mid x^k\right]\\
		&&\quad+ \alpha\|\nabla f_i(x^k) - h_i^k\|^2\\
		&\overset{\eqref{eq:variance_decomposition}}{\le}&	2\alpha\EE\left[\left\|\nabla f_{il}(x^k) - \nabla f_{il}(x^*)\right\|^2\mid x^k\right]\\
		&&\quad + 2\alpha\EE\left[\left\|\nabla f_{il}(w_i^k) - \nabla f_{il}(x^*)\right\|^2\mid x^k\right] + \alpha\|\nabla f_i(x^k) - h_i^k\|^2\\
		&\overset{\eqref{eq:L_smoothness_cor_3}}{\le}& 4L\alpha D_{f_i}(x^k,x^*) + \frac{2\alpha}{m}\sum\limits_{j=1}^m\|\nabla f_{ij}(w_i^k) - \nabla f_{ij}(x^*)\|^2\\
		&&\quad + \alpha\|\nabla f_i(x^k) - h_i^k\|^2
	\end{eqnarray*}
	Putting all together we obtain
	\begin{eqnarray*}
		\EE\left[\|h_i^{k+1} - h_i^*\|^2\mid x^k\right] &\le& \|h_i^k - h_i^*\|^2 + \alpha\left\langle \nabla f_i(x^k) - h_i^k, f_i(x^k) + h_i^k - 2h_i^* \right\rangle\\
		&&\quad + 4L\alpha D_{f_i}(x^k,x^*) + \frac{2\alpha}{m}\sum\limits_{j=1}^m\|\nabla f_{ij}(w_i^k) - \nabla f_{ij}(x^*)\|^2\\
		&\overset{\eqref{eq:a-b_a+b}}{=}& \|h_i^k - h_i^*\|^2 + \alpha\|\nabla f_i(x^k) - h_i^*\|^2 - \alpha\|h_i^k - h_i^*\|^2\\
		&&\quad + 4L\alpha D_{f_i}(x^k,x^*) + \frac{2\alpha}{m}\sum\limits_{j=1}^m\|\nabla f_{ij}(w_i^k) - \nabla f_{ij}(x^*)\|^2\\
		&\overset{\eqref{eq:L_smoothness_cor_3}}{\le}& (1-\alpha)\|h_i^k - h_i^*\|^2 + 6L\alpha D_{f_i}(x^k,x^*)\\
		&&\quad + \frac{2\alpha}{m}\sum\limits_{j=1}^m\|\nabla f_{ij}(w_i^k) - \nabla f_{ij}(x^*)\|^2.
	\end{eqnarray*}
	Summing up the above inequality for $i=1,\ldots, n$ we derive
	\begin{eqnarray}
		\EE\left[\sigma_{1,k+1}^2\mid x^k\right] &\le& (1-\alpha)\sigma_{1,k}^2 + 6L\alpha(f(x^k) - f(x^*)) + 2\alpha\sigma_{2,k}^2.\notag
	\end{eqnarray}
	Similarly to the proof of Lemma~\ref{lem:sigma_k+1_bound_ec-LSVRG} we get
	\begin{eqnarray}
		\EE\left[\sigma_{2,k+1}^2\mid x^k\right] &=& \frac{1}{nm}\sum\limits_{i=1}^n\sum\limits_{j=1}^m\EE\left[\|\nabla f_{ij}(w_i^{k+1}) - \nabla f_{ij}(x^*)\|^2\mid x^k\right]\notag\\
		&=& \frac{1-p}{nm}\sum\limits_{i=1}^n\sum\limits_{j=1}^m\|\nabla f_{ij}(w_i^{k}) - \nabla f_{ij}(x^*)\|^2\notag\\
		&&\quad + \frac{p}{nm}\sum\limits_{i=1}^n\sum\limits_{j=1}^m\|\nabla f_{ij}(x^{k}) - \nabla f_{ij}(x^*)\|^2\notag\\
		&\overset{\eqref{eq:L_smoothness_cor_3}}{\le}& (1-p)\sigma_{2,k}^2 + \frac{2Lp}{nm}\sum\limits_{i=1}^n\sum\limits_{j=1}^m D_{f_{ij}}(x^k,x^*)\notag\\
		&=& (1-p)\sigma_{2,k}^2 + 2Lp\left(f(x^k) - f(x^*)\right).\notag
	\end{eqnarray}
\end{proof}

Applying Theorem~\ref{thm:ec_sgd_main_result_new} we get the following result.
\begin{theorem}\label{thm:ec_LSVRG-diana}
	Assume that $f_{ij}(x)$ is convex and $L$-smooth for all $i=1,\ldots, n$, $j=1,\ldots,m$ and $f(x)$ is $\mu$-quasi strongly convex. Then {\tt EC-LSVRG-DIANA} satisfies Assumption~\ref{ass:key_assumption_finite_sums_new} with
	\begin{gather*}
		A = A' = 2L,\quad B_1' = B_2 = 0,\quad B_1 = B_2' = 2,\quad D_1 = \widetilde{D}_1 = D_1' = D_2 = D_3 = 0,\\
		\widetilde{A} = 3L,\quad \widetilde{B}_1 = \widetilde{B}_2 = 3,\quad \sigma_{1,k}^2 = \frac{1}{n}\sum\limits_{i=1}^n\|h_i^k - \nabla f_i(x^*)\|^2,\quad \rho_1 = \alpha,\\
		\sigma_{2,k}^2 = \frac{1}{nm}\sum\limits_{i=1}^n\sum_{j=1}^m\|\nabla f_{ij}(w_i^k) - \nabla f_{ij}(x^*)\|^2,\quad \rho_2 = p,\quad C_1 = 3L\alpha,\quad C_2 = Lp,\\
		G = 2,\quad F_1 = \frac{24L\gamma^2\left(\frac{4}{\delta}+3\right)}{\delta\alpha\left(1-\min\left\{\frac{\gamma\mu}{2},\frac{\alpha}{4},\frac{p}{4}\right\}\right)},\quad F_2 = \frac{24L\gamma^2\left(\frac{4}{1-\alpha}\left(\frac{4}{\delta}+3\right) + 3\right)}{\delta p\left(1-\min\left\{\frac{\gamma\mu}{2},\frac{\alpha}{4},\frac{p}{4}\right\}\right)},
	\end{gather*}
	with $\gamma$ and $\alpha$ satisfying
	\begin{equation*}
		\gamma \le \min\left\{\frac{9}{296L}, \frac{\delta}{4L\sqrt{6\left(4+3\delta+\frac{2}{1-\alpha}\left(3+\frac{4}{1-p}\right)(4+3\delta)+\frac{6\delta}{1-p}\right)}}\right\},\quad \alpha \le \frac{1}{\omega+1}
	\end{equation*}
	with $M_1 = 0$ and $M_2 = \frac{8}{3p} + \frac{32}{9p}$ and for all $K \ge 0$
	\begin{equation*}
		\EE\left[f(\bar x^K) - f(x^*)\right] \le \left(1 - \min\left\{\frac{\gamma\mu}{2},\frac{\alpha}{4},\frac{p}{4}\right\}\right)^K\frac{4(T^0 + \gamma F_1 \sigma_{1,0}^2 + \gamma F_2 \sigma_{2,0}^2)}{\gamma},
	\end{equation*}	
	when $\mu > 0$ and
	\begin{equation*}
		\EE\left[f(\bar x^K) - f(x^*)\right] \le \frac{4(T^0 + \gamma F_1 \sigma_{1,0}^2 + \gamma F_2 \sigma_{2,0}^2)}{K\gamma}
	\end{equation*}
	when $\mu = 0$, where $T^k \eqdef \|x^k - x^*\|^2+ M_2\gamma^2 \sigma_{2,k}^2$.
\end{theorem}
In other words, if $p = \nicefrac{1}{m}$, $m \ge 2$ and
\begin{equation*}
		\gamma = \min\left\{\frac{9}{296L}, \frac{\delta}{4L\sqrt{6\left(4+3\delta+\frac{2}{1-\alpha}\left(3+\frac{4}{1-p}\right)(4+3\delta)+\frac{6\delta}{1-p}\right)}}\right\},\quad \alpha = \min\left\{\frac{1}{\omega+1},\frac{1}{2}\right\},
\end{equation*}
then {\tt EC-LSVRG-DIANA} converges with the linear rate
\begin{equation*}
	\cO\left(\left(\omega + m + \frac{\kappa}{\delta}\right)\ln\frac{1}{\varepsilon}\right)
\end{equation*}
to the exact solution when $\mu > 0$.

Applying Lemma~\ref{lem:lemma_technical_cvx} we get the complexity result in the case when $\mu = 0$.
\begin{corollary}\label{cor:ec_lsvrg_diana_cvx_cor}
	Let the assumptions of Theorem~\ref{thm:ec_LSVRG-diana} hold and $\mu = 0$. Then after $K$ iterations of {\tt EC-LSVRG-DIANA} with the stepsize
	\begin{eqnarray*}
		\gamma_0 &=& \min\left\{\frac{9}{296L}, \frac{\delta}{4L\sqrt{6\left(4+3\delta+\frac{2}{1-\alpha}\left(3+\frac{4}{1-p}\right)(4+3\delta)+\frac{6\delta}{1-p}\right)}}\right\},\quad R_0 = \|x^0-x^*\|,\\
		\gamma &=& \min\left\{\gamma_0, \sqrt{\frac{9pR_0^2}{56\sigma_{2,0}^2}}, \sqrt[3]{\frac{R_0^2}{\frac{24L\left(\frac{4}{\delta}+3\right)}{\delta\alpha\left(1-\min\left\{\frac{\gamma_0\mu}{2},\frac{\alpha}{4},\frac{p}{4}\right\}\right)}\sigma_{1,0}^2 + \frac{24L\left(\frac{4}{1-\alpha}\left(\frac{4}{\delta}+3\right) + 3\right)}{\delta p\left(1-\min\left\{\frac{\gamma_0\mu}{2},\frac{\alpha}{4},\frac{p}{4}\right\}\right)}\sigma_{2,0}^2}}\right\},
	\end{eqnarray*}	
	and $p = \frac{1}{m}$, $m\ge 2$, $\alpha = \min\left\{\frac{1}{\omega+1},\frac{1}{2}\right\}$ we have $\EE\left[f(\bar{x}^K) - f(x^*)\right]$ of order
	\begin{equation*}
		\cO\left(\frac{L R_0^2}{\delta K} + \frac{\sqrt{R_0^2 m\sigma_{2,0}^2}}{K} + \frac{\sqrt[3]{LR_0^4((\omega+1)\sigma_{1,0}^2 + m\sigma_{2,0}^2)}}{\delta^{\nicefrac{2}{3}}K}\right).
	\end{equation*}
	That is, to achive $\EE\left[f(\bar{x}^K) - f(x^*)\right] \le \varepsilon$ {\tt EC-LSVRG-DIANA} requires
	\begin{equation*}
		\cO\left(\frac{L R_0^2}{\delta \varepsilon} + \frac{\sqrt{R_0^2 m\sigma_{2,0}^2}}{\varepsilon} + \frac{\sqrt[3]{LR_0^4((\omega+1)\sigma_{1,0}^2 + m\sigma_{2,0}^2)}}{\delta^{\nicefrac{2}{3}}\varepsilon}\right)
	\end{equation*}
	iterations.
\end{corollary}


\section{Numerical Experiments}\label{sec:numerical_exp}
To justify our theory, we conduct several numerical experimentson logistic regression problem with $\ell_2$-regularization:
\begin{equation}
\mytextstyle	\min\limits_{x\in\R^d} \left\{ f(x) = \frac{1}{N}\sum\limits_{i=1}^N \log\left(1 + \exp\left(-y_i\cdot (Ax)_i\right)\right) + \frac{\mu}{2}\|x\|^2 \right\}, \label{eq:log_loss}
\end{equation}
where $N$ is a number of features, $x\in\R^d$ represents the weights of the model, $A\in\R^{N\times d}$ is a feature matrix, vector $y\in\{-1,1\}^N$ is a vector of labels and $(Ax)_i$ denotes the $i$-th component of vector $Ax$. Clearly, this problem is $L$-smooth and $\mu$-strongly convex with $L = \mu + \nicefrac{\lambda_{\max}(A^\top A)}{4N}$, where $\lambda_{\max}(A^\top A)$ is a largest eigenvalue of $A^\top A$. The datasets were taken from LIBSVM library \cite{chang2011libsvm}, and the code was written in Python 3.7 using standard libraries. Our code is available at  \url{https://github.com/eduardgorbunov/ef_sigma_k}.

We simulate parameter-server architecture using one machine with Intel(R) Core(TM) i7-9750 CPU \@ 2.60 GHz in the following way. First of all, we always use such $N$ that $N = n\cdot m$ and consider  $n = 20$ and $n=100$ workers. The choice of $N$ for each dataset that we consider is stated in Table~\ref{tab:logreg_datasets}.
\begin{table}[H]
    \centering \footnotesize
    \caption{Summary of datasets: $N =$  total \# of data samples; $d =$ \# of features.}
    \label{tab:logreg_datasets}
    \begin{tabular}{|c|c|c|c|c|c|c|}
        \hline
         & {\tt a9a} & {\tt w8a} & {\tt gisette} & {\tt mushrooms} & {\tt madelon} & {\tt phishing} \\
        \hline
        $N$ & $32,000$ & $49,700$ & $6,000$ & $8,000$ & $2,000$ & $11,000$ \\
        \hline
        $d$ & $123$ & $300$ & $5,000$ & $112$ & $500$ & $68$ \\
        \hline
    \end{tabular}
\end{table}
Next, we shuffle the data and split in $n$ groups of size $m$. To emulate the work of workers, we use a single machine and run the methods with the parallel loop in series. Since in these experiments we study sample complexity and number of bits used for communication, this setup is identical to the real parameter-server setup in this sense.

In all experiments we use the stepsize $\gamma = \nicefrac{1}{L}$ and $\ell_2$-regularization parameter $\mu = \nicefrac{10^{-4}\lambda_{\max}(A^\top A)}{4N}$. The starting point $x^0$ for each dataset was chosen so that $f(x^0) - f(x^*) \sim 10$. In experiments with stochastic methods we used batches of size $1$ and uniform sampling for simplicity. For {\tt LSVRG}-type methods we choose $p = \nicefrac{1}{m}$. 

\begin{figure}[H]
    \centering
    \includegraphics[width=0.32\textwidth]{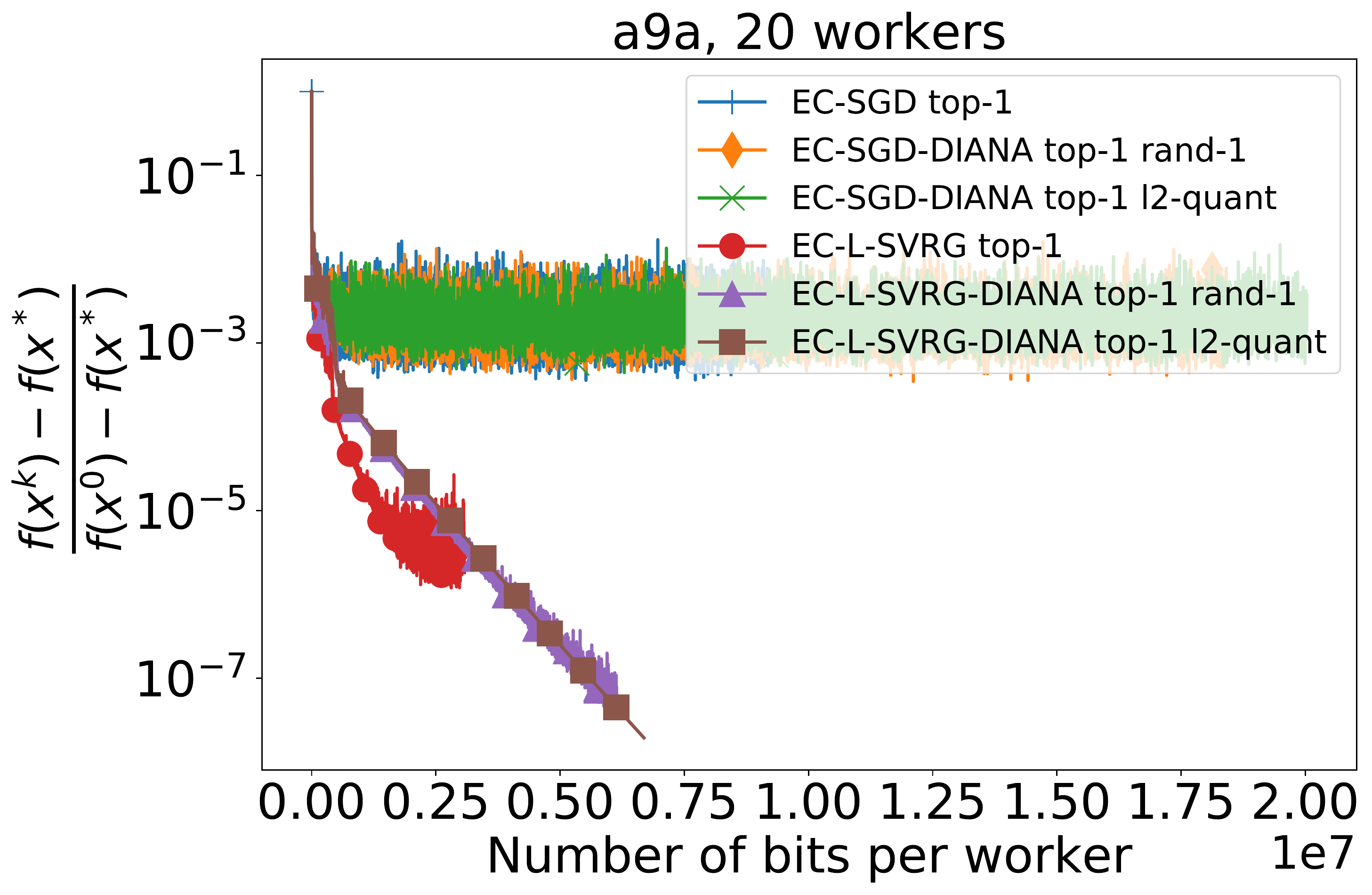}
	\includegraphics[width=0.32\textwidth]{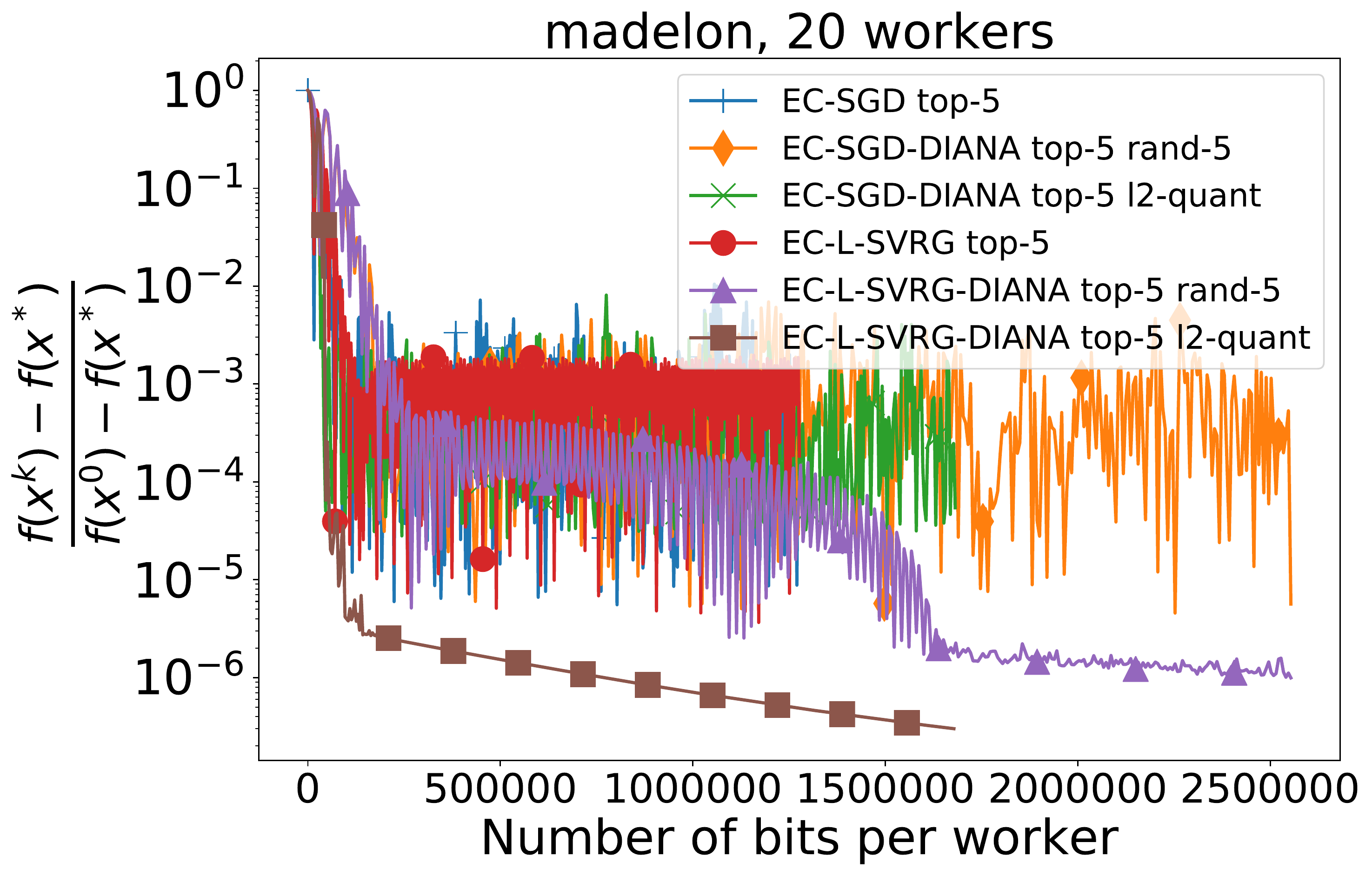}    
	\includegraphics[width=0.32\textwidth]{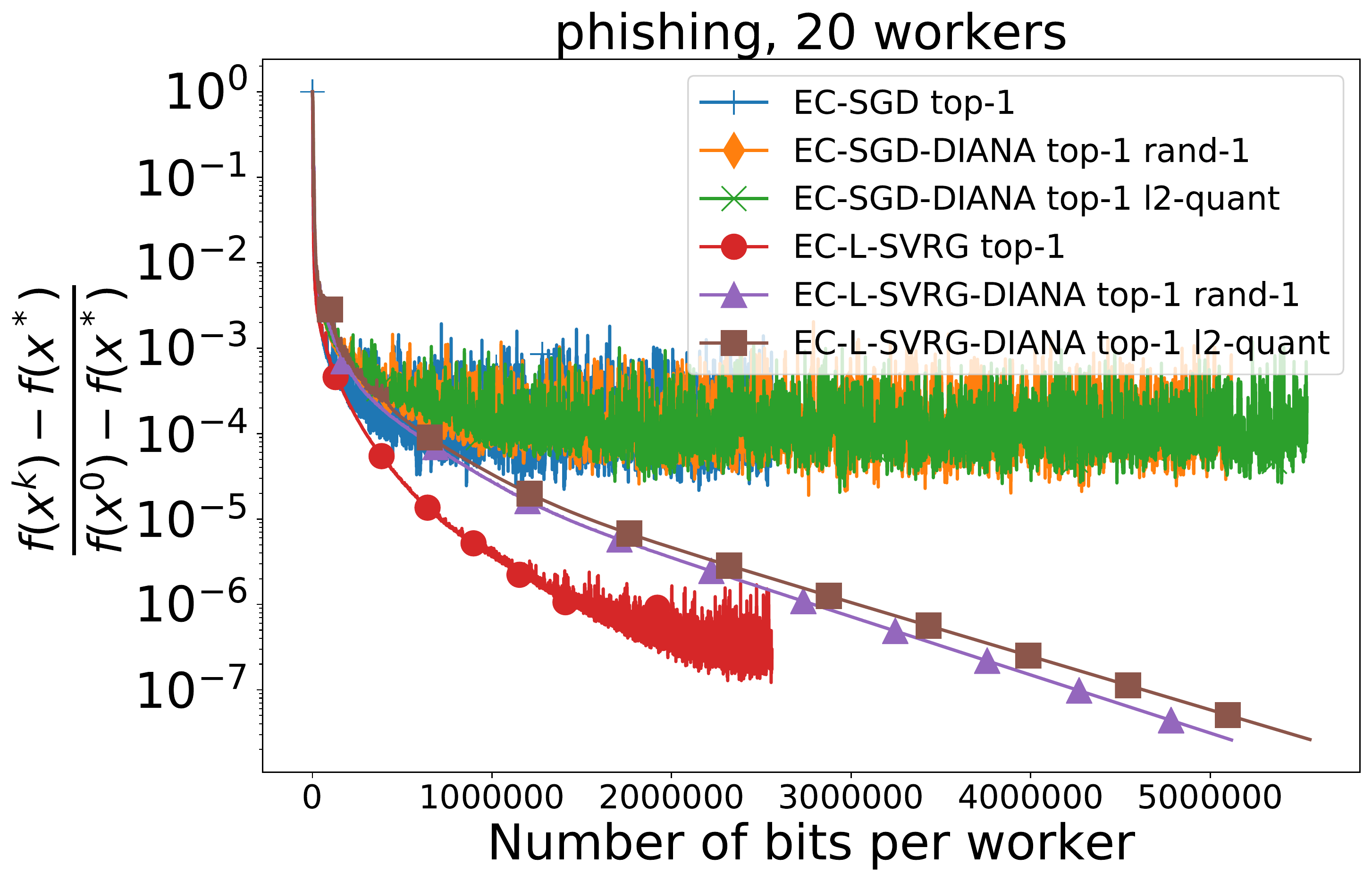}
    \\
    \includegraphics[width=0.32\textwidth]{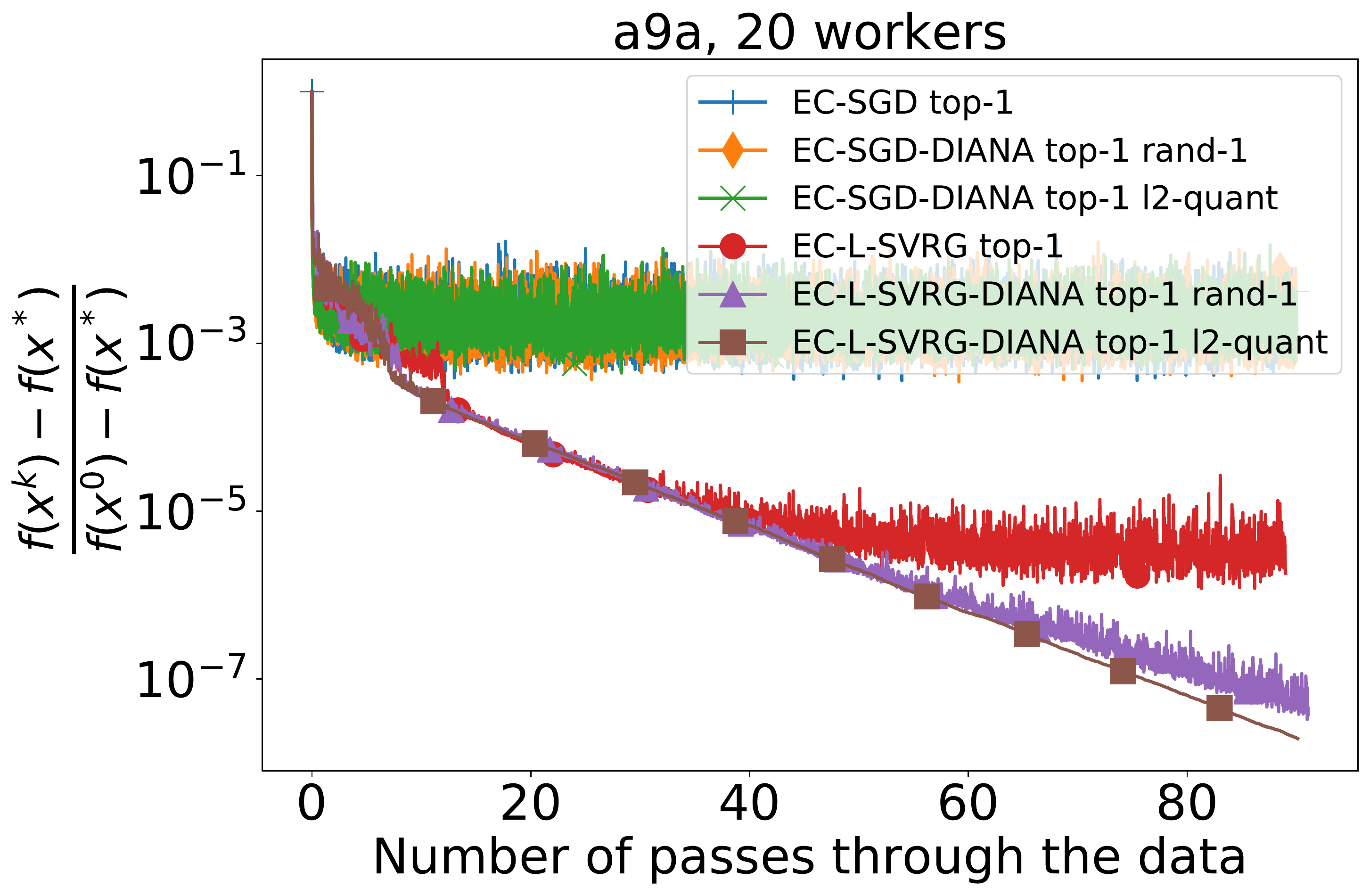}
    \includegraphics[width=0.32\textwidth]{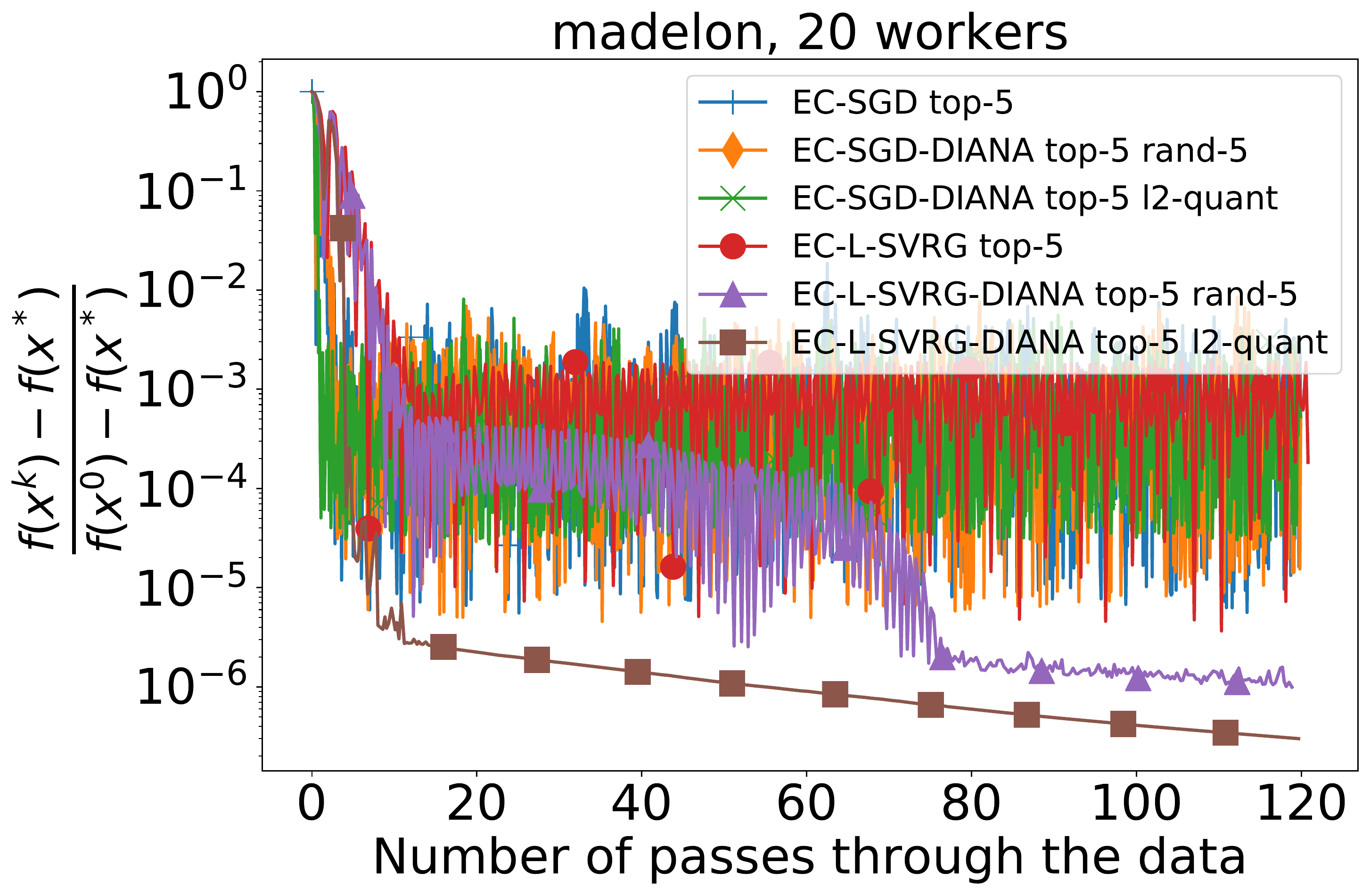}    
	\includegraphics[width=0.32\textwidth]{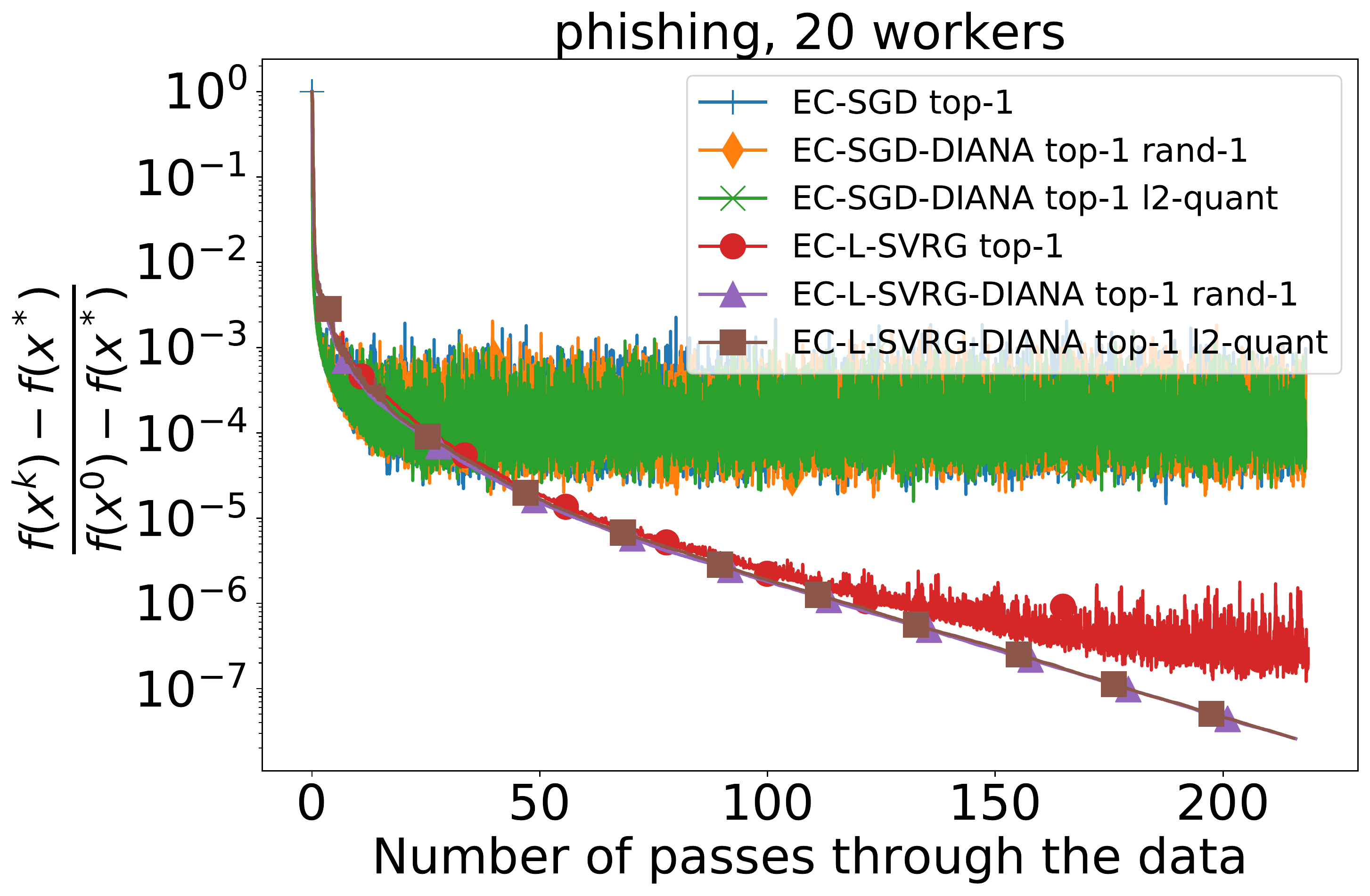} 
     
    \caption{Trajectories of {\tt EC-SGD}, {\tt EC-SGD-DIANA}, {\tt EC-LSVRG} and {\tt EC-LSVRG-DIANA} applied to solve logistic regression problem with $20$ workers.}
    \label{fig:sgd_logreg}
\end{figure}

\textbf{Compressing stochastic gradients.}
The results for {\tt a9a}, {\tt madelon} and {\tt phishing} can be found in Figure~\ref{fig:sgd_logreg} (included here) and for {\tt w8a}, {\tt mushrooms} and {\tt gisette} in Figure~\ref{fig:sgd_logreg_extra} (in the Appendix). We choose number of components for TopK operator of the order $\max\{1,\nicefrac{d}{100}\}$. Clearly, in these experiments we see two levels of noise. For some datasets, like {\tt a9a}, {\tt phishing} or {\tt mushrooms}, the noise that comes from the stochasticity of the gradients dominates the noise coming from compression. Therefore,  methods such as {\tt EC-SGD} and {\tt EC-SGD-DIANA} start to oscillate around a larger value of the loss function than other methods we consider. {\tt EC-LSVRG} reduces the largest source of noise and, as a result, finds a better approximation of the solution. However, at some point, it reaches another level of the loss function and starts to oscillate there due to the noise coming from compression. Finally, {\tt EC-LSVRG-DIANA} reduces the variance of both types, and as a result, finds an even better approximation of the solution. In contrast, for  the {\tt madelon} dataset, both noises are of the same order, and therefore, {\tt EC-LSVRG} and {\tt EC-SGD-DIANA} behave similarly to {\tt EC-SGD}. However, {\tt EC-LSVRG-DIANA} again reduces both types of noise effectively and finds a better approximation of the solution after a given number of epochs.   In the experiments with {\tt w8a} and {\tt gisette} datasets, the noise produced by compression is dominated by the noise coming from the stochastic gradients. As a result, we see that the {\tt DIANA}-trick is not needed here.

\textbf{Compressing full gradients.}
\begin{figure}[H]
    \centering
    \includegraphics[width=0.32\textwidth]{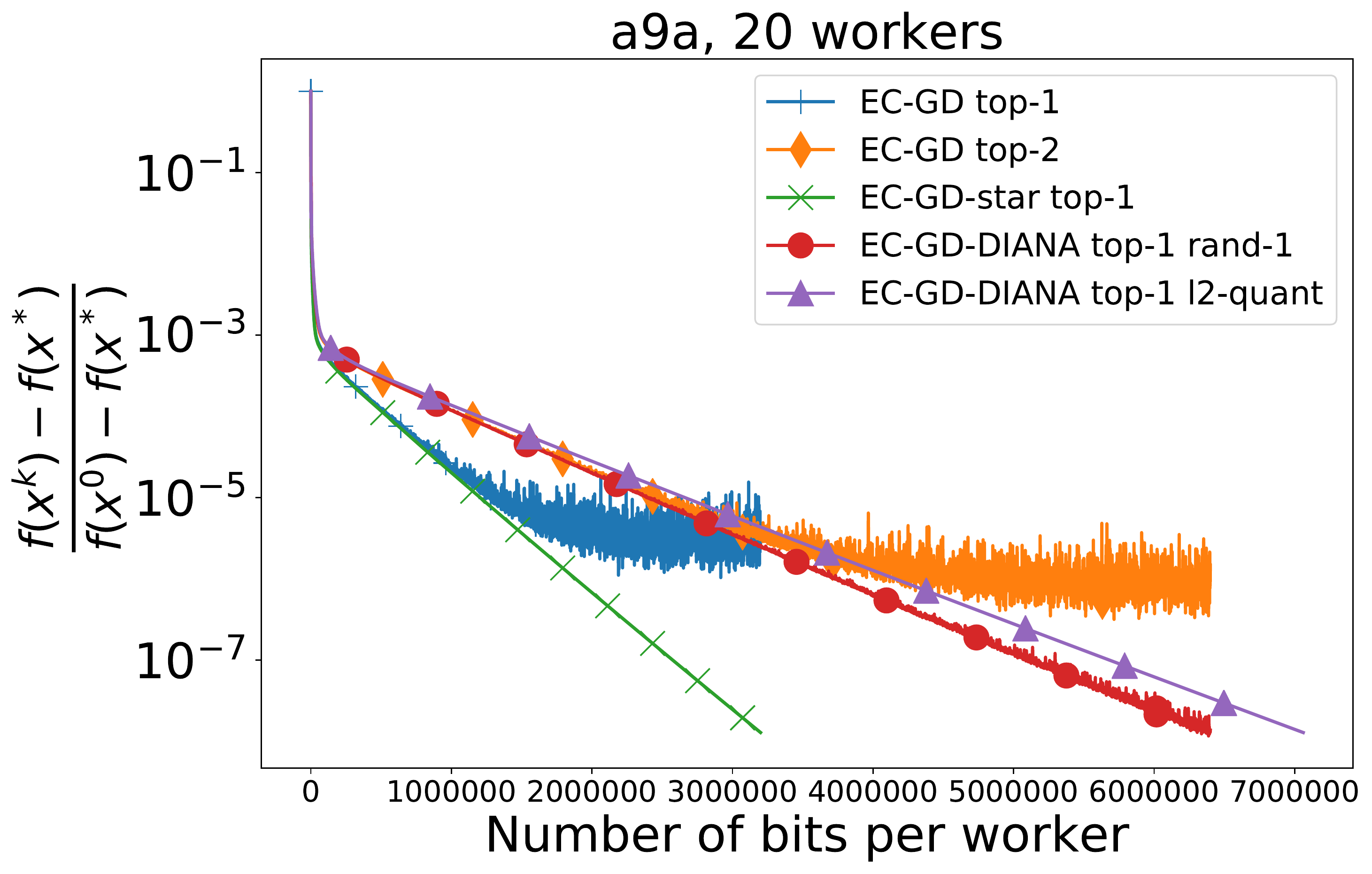}
	\includegraphics[width=0.32\textwidth]{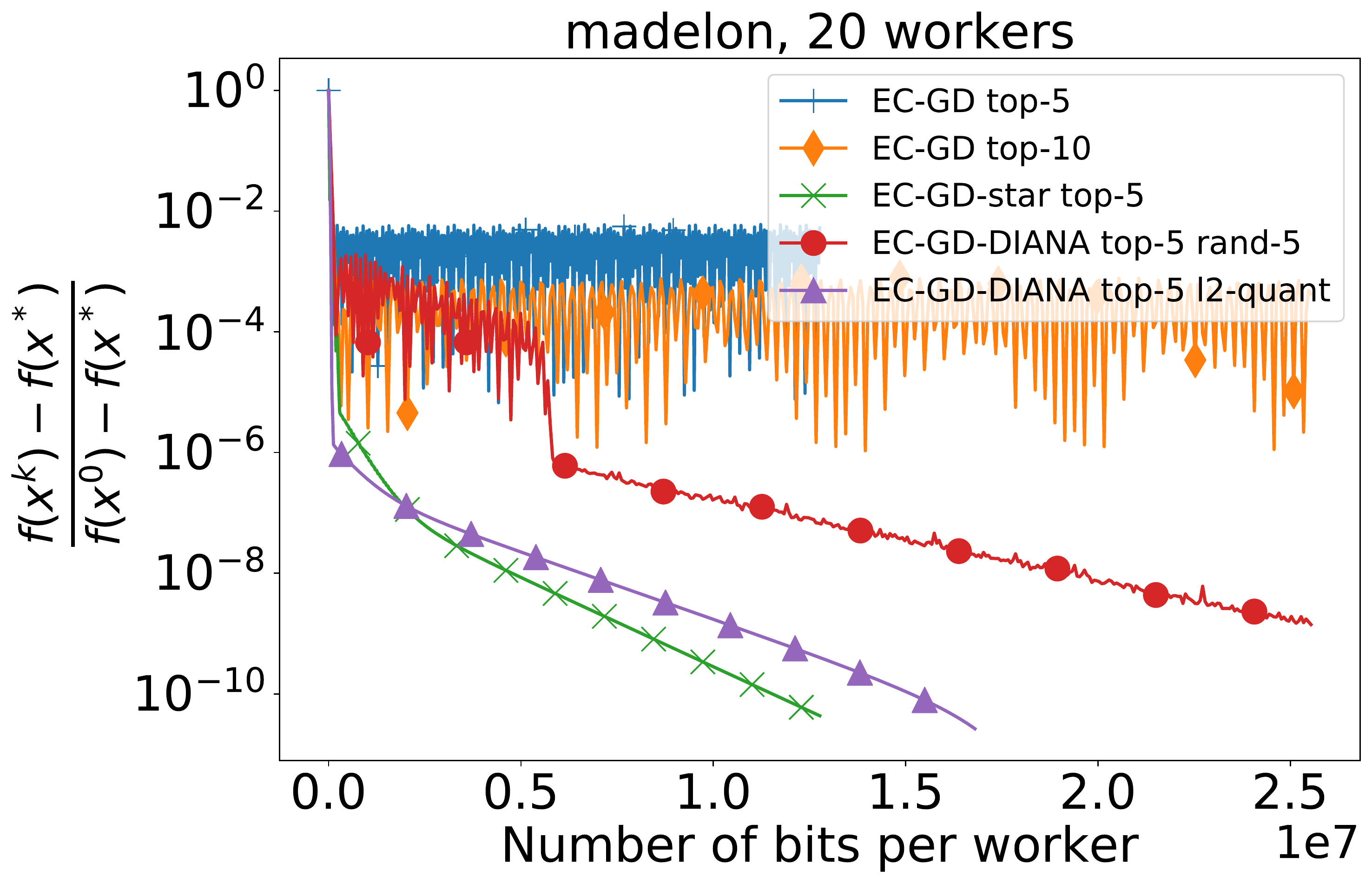}    
	\includegraphics[width=0.32\textwidth]{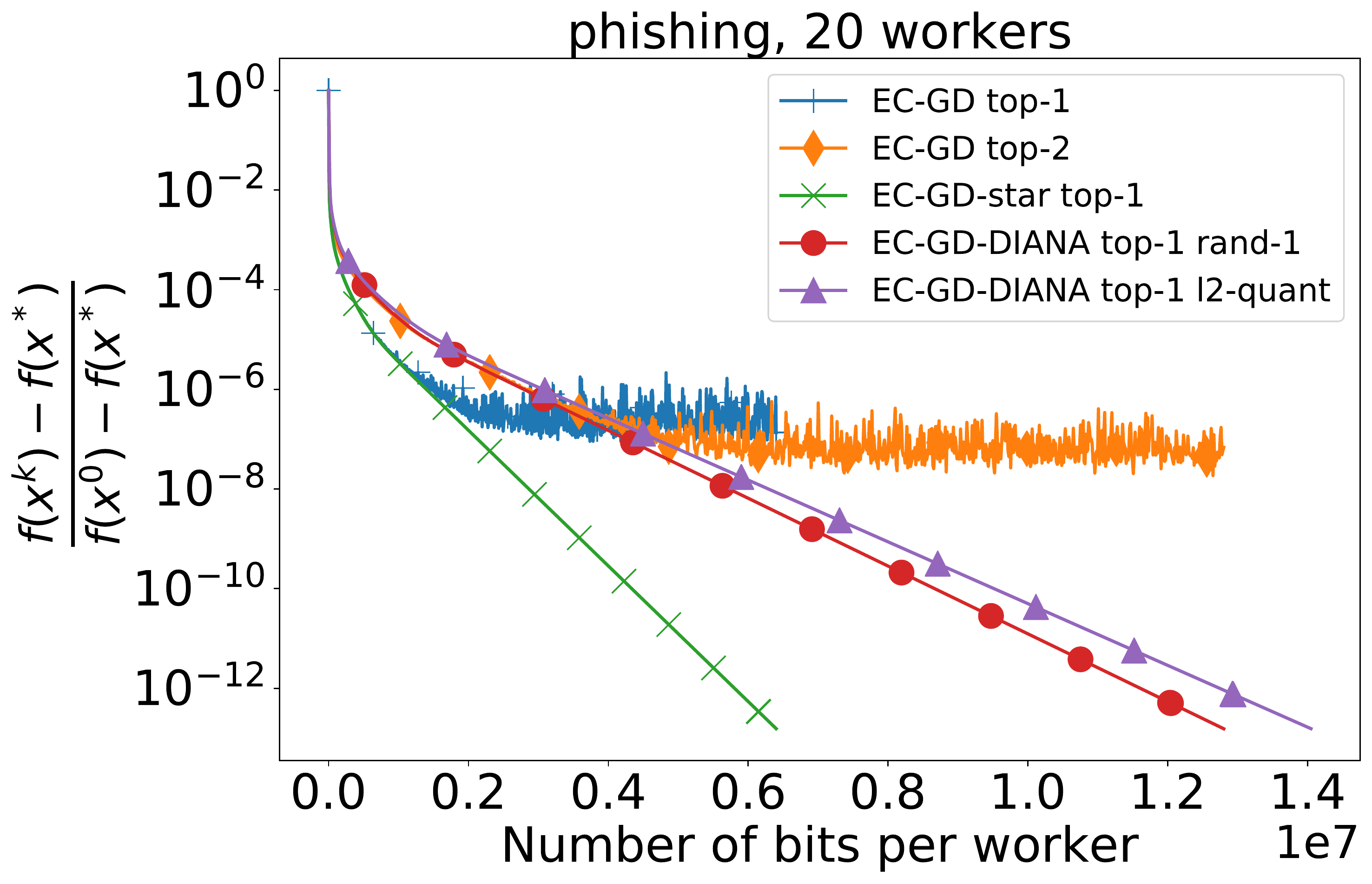}    
    \\
    \includegraphics[width=0.32\textwidth]{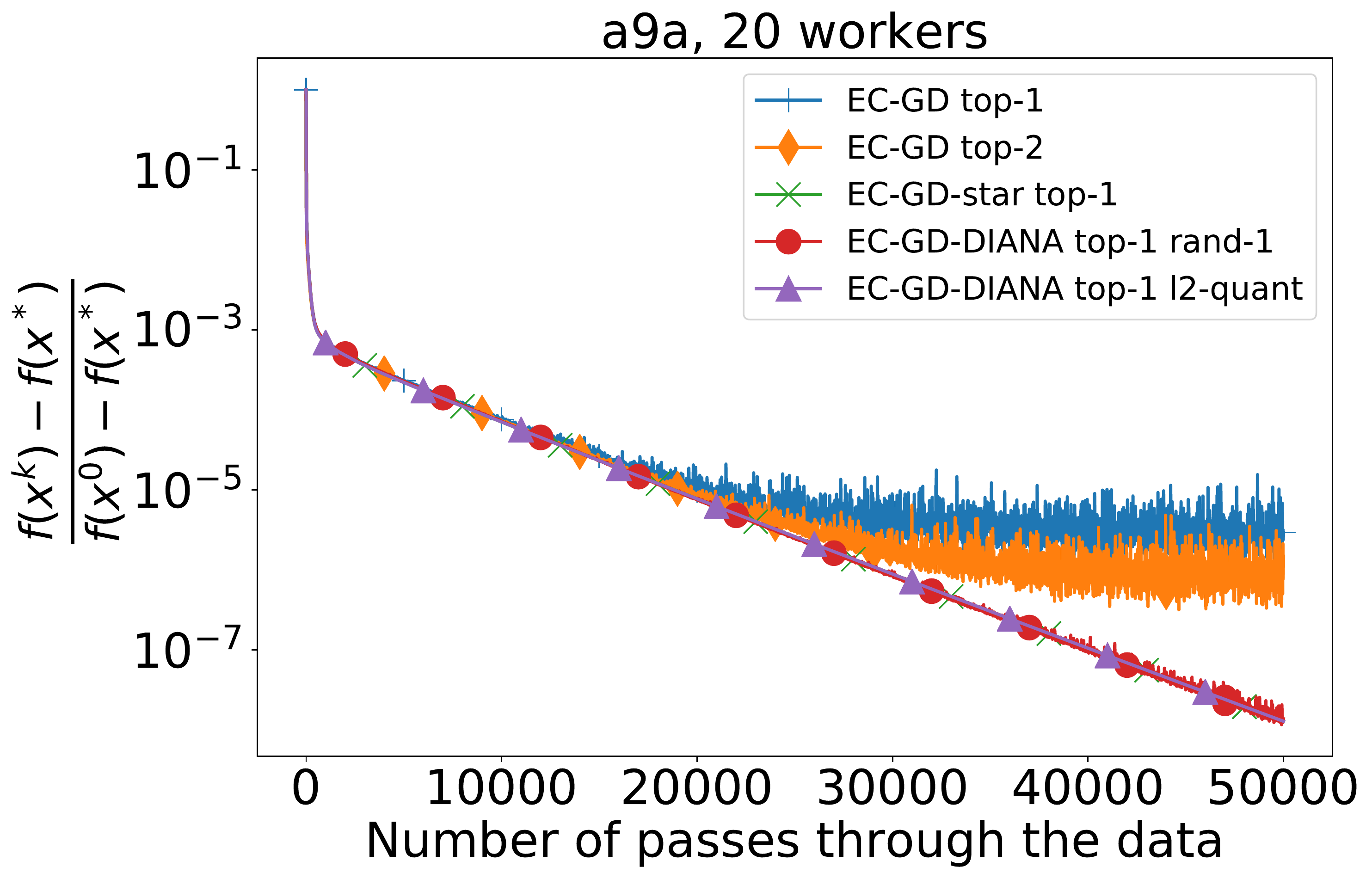}    
	\includegraphics[width=0.32\textwidth]{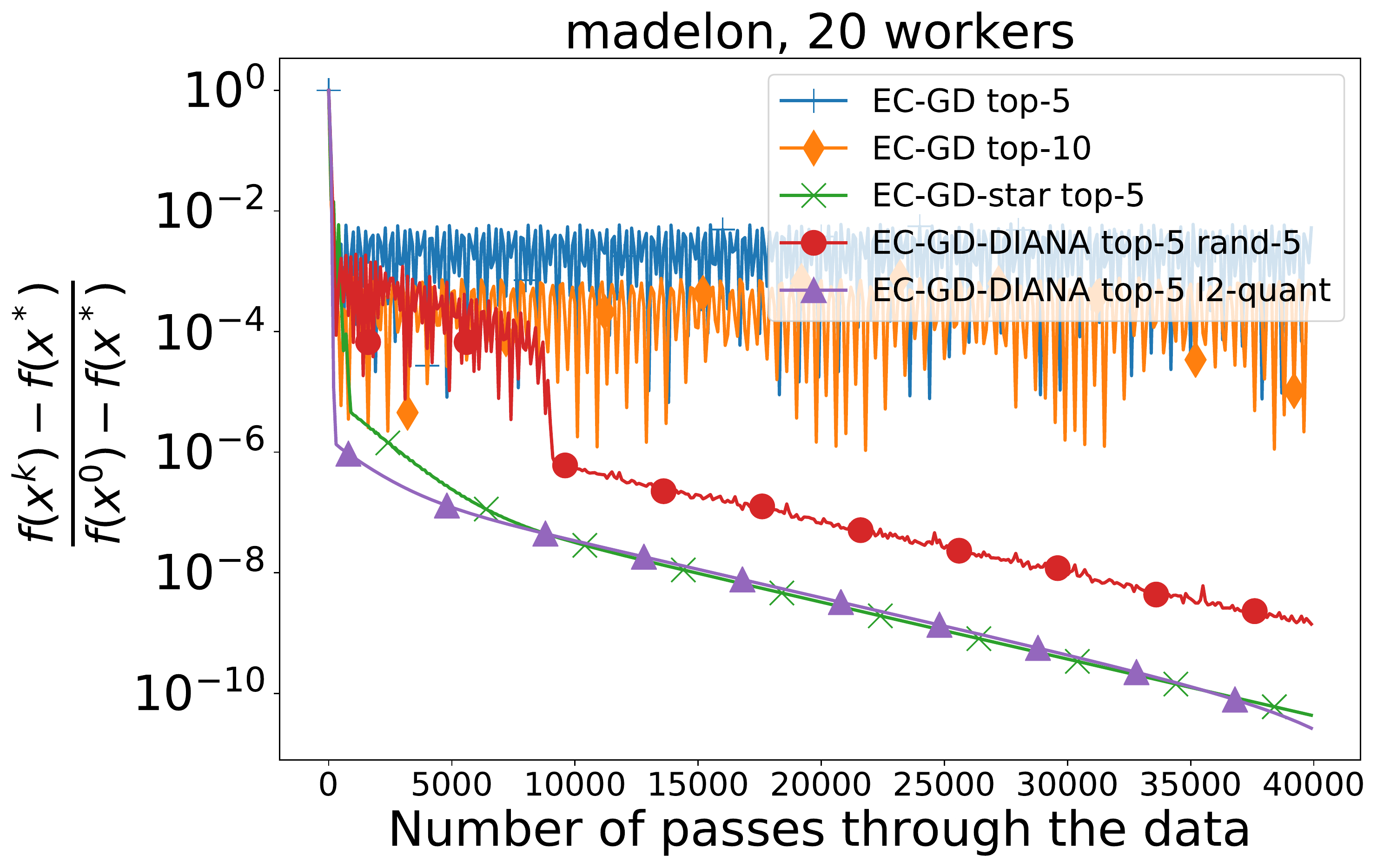}    
	\includegraphics[width=0.32\textwidth]{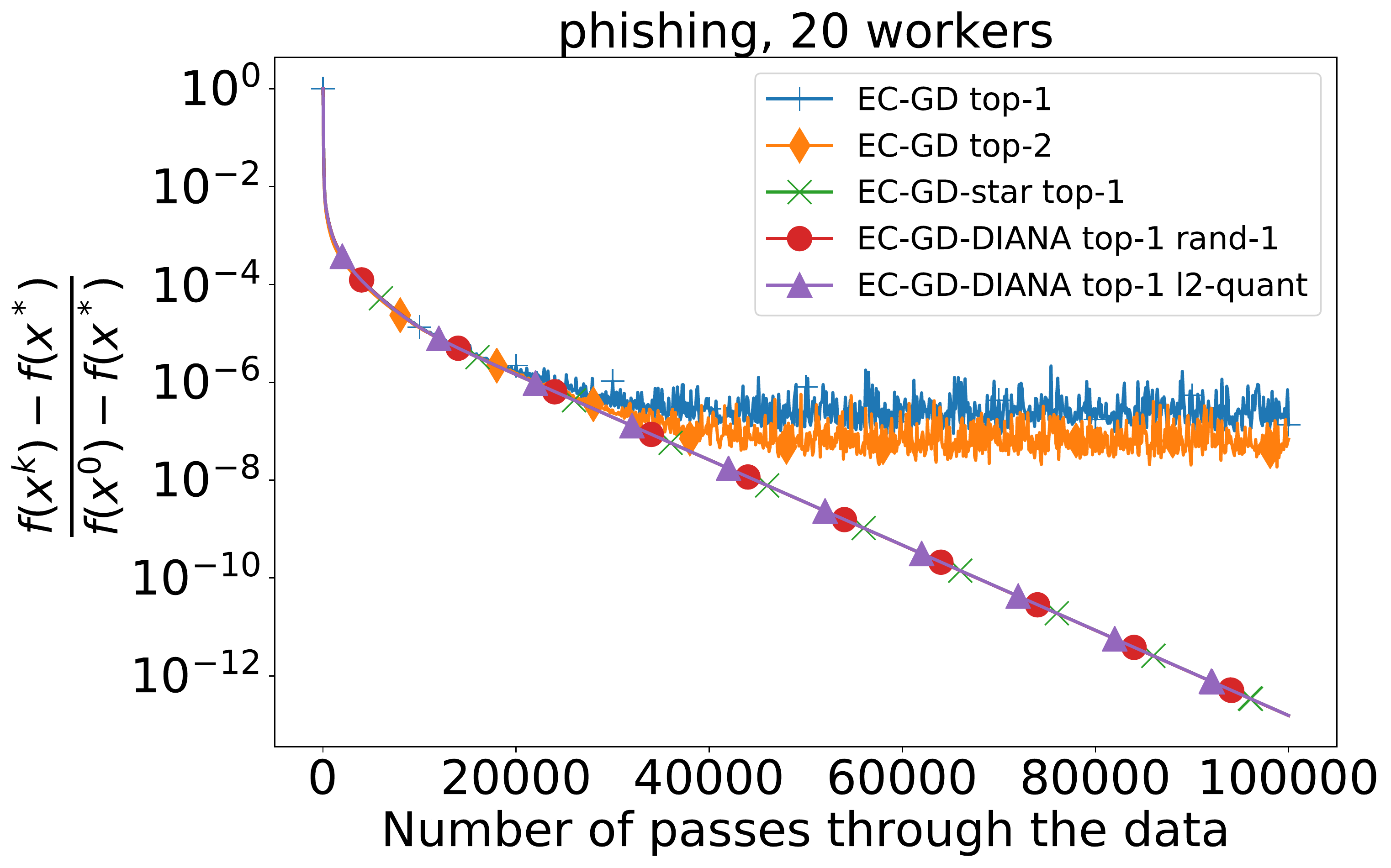}    	
	\caption{Trajectories of {\tt EC-GD}, {\tt EC-GD-star} and {\tt EC-DIANA} applied to solve logistic regression problem with $20$ workers.}
    \label{fig:gd_logreg_20_workers}
\end{figure}
In order to show the effect of {\tt DIANA}-type variance reduction itself, we consider the case when all workers compute the full gradients of their functions, see Figure~\ref{fig:gd_logreg_20_workers} (included here) and Figures~\ref{fig:gd_logreg_20_workers_appendix}--\ref{fig:gd_logreg_100_workers_id} (in the Appendix). Clearly, for all datasets except {\tt mushrooms}, {\tt EC-GD} with constant stepsize converges  to a neighborhood of the solution only, while {\tt EC-GDstar} and {\tt EC-GD-DIANA} converge with linear rate asymptotically to the exact solution. {\tt EC-GDstar} always show the best performance, however, it is impractical: we used a very good approximation of the solution to apply this method. In contrast, {\tt EC-DIANA} converges slightly slower and requires more bits for communication; but it is practical and shows better performance than {\tt EC-GD}. On the {\tt mushrooms} datasets, {\tt EC-GD} does not reach the oscillation region after the given number of epochs, therefore, it is preferable there.

%% file: ch4_local_sigma_k.tex
\chapter{Local SGD: Unified Theory and New Efficient Methods}\label{ch:local_sigma_k}

\section{Introduction}\label{sec:intro}
In this chapter\footnote{Part of this work was done while I was a research intern at KAUST.}, we are interested in a centralized distributed optimization problem of the form
\begin{equation}
\mytextstyle
	\min\limits_{x\in\R^d} f(x) = \frac{1}{n}\sum\limits_{i=1}^n f_i(x), \label{eq:main_problem}
\end{equation}
where $n$ is the number of devices/clients/nodes/workers. We assume that $f_i$ can be represented either as a) an expectation, i.e.,
\begin{equation}
\mytextstyle
	f_i(x) = \EE_{\xi_i\sim \cD_i}\left[f_{\xi_i}(x)\right], \label{eq:f_i_expectation}
\end{equation}
where $\cD_i$ describes the distribution of data on device $i$,  or b) as a finite sum, i.e.,
\begin{equation}
\mytextstyle
	f_i(x) = \frac{1}{m}\sum\limits_{j=1}^m f_{ij}(x). \label{eq:f_i_sum}
\end{equation} 
While our theory allows the number of functions $m$ to vary across the devices, for simplicity of exposition, we restrict the narrative to this simpler case.

Federated learning (FL)---an emerging subfield of machine learning \citep{mcmahan2016federated, FEDLEARN, FL2017-AISTATS}---is traditionally cast as an instance of problem~\eqref{eq:main_problem} with several idiosyncrasies. First, the number of devices $n$ is very large: tens of thousands to millions. Second,  the devices (e.g., mobile phones) are often very heterogeneous in their compute, connectivity, and storage capabilities. The data defining each function $f_i$ reflects the usage patterns of the device owner, and as such, it is either unrelated or at best related only weakly. Moreover, device owners desire to protect their local private data, and for that reason, training needs to take place with the data remaining on the devices.  Finally, and this is of key importance for the development in this work, communication among the workers, typically conducted via a trusted aggregation server, is very expensive. 


\paragraph{Communication bottleneck.} There are two main directions in the literature for tackling the communication cost issue in FL. The first approach consists of algorithms that aim to reduce the number of transmitted bits by applying a  carefully chosen gradient compression scheme, such as quantization~\citep{alistarh2017qsgd, pmlr-v80-bernstein18a, mishchenko2019distributed, horvath2019stochastic,ramezani2019nuqsgd, reisizadeh2020fedpaq}, sparsification~\citep{aji2017sparse, deep, alistarh2018convergence, wangni2018gradient, wang2018atomo, mishchenko202099}, or other more sophisticated strategies~\citep{karimireddy2019error, stich2020error, pmlr-v80-wu18d, vogels2019powersgd, beznosikov2020biased,gorbunov2020linearly}. The second approach---one that we investigate in this chapter---instead focuses on increasing the total amount of local computation in between the communication rounds in the hope that this will reduce the total number of communication rounds needed to build a model of sufficient quality~\citep{DANE,zhang2015disco,  reddi2016aide, li2018federated, pathak2020fedsplit}. These two approaches, {\em communication compression} and {\em local computation}, can be combined for a better practical performance~\citep{basu2019qsparse}.

\paragraph{Local first-order algorithms.} Motivated by recent development in the field~\citep{localsgd_first, mcmahan2016federated, stich2018local, LinSPJ2018local,liang2019variance, wu2019federated, karimireddy2020scaffold, khaled2020tighter, woodworth2020local}, in this chapter we perform an in-depth and general study of {\em local first-order algorithms}. Contrasted with zero or higher order local methods, local first order methods perform several gradient-type steps in between the communication rounds.  In particular, we consider the following family of methods:
\begin{equation}
	x_i^{k+1} = \begin{cases} x_i^k - \gamma g_i^k,& \text{if } c_{k+1} = 0,\\ \frac{1}{n}\sum\limits_{i=1}^n \left(x_i^k - \gamma g_i^k\right),& \text{if } c_{k+1} = 1, \end{cases} \label{eq:local_sgd_def}
\end{equation}
where $x_i^k$ represents the local variable maintained by the $i$-th device, $g_i^k$ represents local first order direction\footnote{Vector $g_i^k$ can be a simple unbiased estimator of $\nabla f_i(x_i^k)$, but  can also involve a local ``shift'' designed to correct the (inherently wrong) fixed point of local methods. We elaborate on this point later.} and (possibly random) sequence $\{c_{k}\}_{k\ge 1}$ with $c_k \in\{0,1\}$ encoding the times when communication takes place. 

Both the classical {\tt Local-SGD/FedAvg}~\citep{mcmahan2016federated, stich2018local, khaled2020tighter, woodworth2020local} and shifted local {\tt SGD}~\citep{liang2019variance, karimireddy2020scaffold} methods fall into this category of algorithms. However, most of the existing methods have been analyzed with limited flexibility only, leaving many potentially fruitful directions unexplored. The most important unexplored questions include i) better understanding of the local shift that aims to correct the fixed point of local methods, ii) support for more sophisticated local gradient estimators that allow for importance sampling, variance reduction, or coordinate descent, iii) variable number of local steps, and iv) general theory supporting multiple data similarity types, including identical, heterogeneous and partially heterogeneous ($\zeta$-heterogeneous - defined later). 

Consequently, there is a need for a single framework unifying the theory of local stochastic first order methods, ideally one capable of pointing to new and more efficient variants. This is what we do in this work.

\paragraph{Unification of stochastic algorithms.} There have been multiple recent papers aiming to unify the theory of first-order optimization algorithms. The closest to our work is the unification of (non-local) stochastic algorithms in~\citep{gorbunov2019unified} that proposes a relatively simple yet powerful framework for analyzing variants of {\tt SGD} that allow for minibatching, arbitrary sampling,\footnote{A tight convergence rate given any sampling strategy and any smoothness structure of the objective.} variance reduction, subspace gradient oracle, and quantization. We recover this framework as a special case in a non-local regime. Next, a framework for analyzing error compensated or delayed SGD methods was recently proposed in~\citep{gorbunov2020linearly}. Another relevant approach covers the unification of decentralized {\tt SGD} algorithms~\citep{koloskova2020unified}, which is able to recover the basic variant of {\tt Local-SGD} as well. While our framework matches their rate for basic {\tt Local-SGD}, we cover a broader range of local methods in this work as we focus on the centralized setting.

\subsection{Our Contributions}

In this chapter, we propose a general framework for analyzing a broad family of local stochastic gradient methods of the form~\eqref{eq:local_sgd_def}. Given that a particular local algorithm satisfies a specific parametric assumption (Assumption~\ref{ass:key_assumption}) in a certain scenario, we provide a tight convergence rate of such a method. 

Let us give a glimpse of our results and their generality. A local algorithm of the  form~\eqref{eq:local_sgd_def} is allowed to consist of an {\em arbitrary} local stochastic gradient estimator (see Section~\ref{sec:local_solver} for details), a possible {\em drift/shift} to correct for the non-stationarity of local methods\footnote{Basic local algorithms such as {\tt FedAvg}/{\tt Local-SGD} or {\tt FedProx}~\citep{li2018federated} have incorrect fixed points~\citep{pathak2020fedsplit}. To eliminate this issue, a strategy of adding an extra ``drift'' or ``shift'' to the local gradient has been proposed recently~\citep{liang2019variance, karimireddy2020scaffold}.} and a fixed or random local loop size. Further, we provide a tight convergence rate in both the identical and heterogeneous data regimes for strongly (quasi) convex and convex objectives. Consequently, our framework is capable of:

$\bullet$  {\bf Recovering known optimizers along with their tight rates.} We recover multiple known local optimizers as a special case of our general framework, along with their convergence rates (up to small constant factors). This includes {\tt FedAvg/Local-SGD}~\citep{mcmahan2016federated, stich2018local} with currently the best-known convergence rate~\citep{khaled2020tighter, woodworth2020local, koloskova2020unified, woodworth2020minibatch} and {\tt SCAFFOLD}~\citep{karimireddy2020scaffold}. Moreover, in a special case we recover a general framework for analyzing non-local {\tt SGD} method developed in ~\citep{gorbunov2019unified}, and consequently we recover multiple variants of {\tt SGD} with and without variance reduction, including {\tt SAGA}~\citep{SAGA}, {\tt L-SVRG}~\citep{kovalev2019don},  {\tt SEGA}~\citep{hanzely2018sega}, gradient compression methods~\citep{mishchenko2019distributed, horvath2019stochastic} and many more. 

$\bullet$  {\bf Filling missing gaps for known methods.} Many of the recovered optimizers have only been analyzed under specific and often limiting circumstances and regimes. Our framework allows us to extend known methods into multiple hitherto unexplored settings. For instance, for each (local) method our framework encodes, we allow for a random/fixed local loop size, identical/heterogeneous/$\zeta$-heterogeneous data (introduced soon), and convex/strongly convex objective.

$\bullet$  {\bf Extending the established optimizers.} To the best of our knowledge, none of the known local methods have been analyzed under arbitrary smoothness structure of the local objectives\footnote{By this we mean that function $f_{i,j}$ from~\eqref{eq:f_i_sum} is $\mM_{i,j}$-smooth with $\mM_{i,j}\in \R^{d\times d}, \mM_{i,j}\succeq 0$, i.e., for all $x,y\in  \R^d$ we have $f_{i,j}(x)\leq f_{i,j}(y) + \langle\nabla  f_{i,j}(y),x-y \rangle + \frac{1}{2} (x-y)^\top \mM_{i,j} (x-y)$. As an example, logistic regression possesses naturally such a structure with matrices $\mM_{i,j}$ of rank 1.} and consequently, our framework is the first to allow for the local stochastic gradient to be constructed via importance (possibly minibatch) sampling. Next, we allow for a local loop with a random length, which is a new development contrasting with the classical fixed-length regime. We discuss advantages of of the random loop in Section~\ref{sec:data_and_loop}.

$\bullet$  {\bf New efficient algorithms.} Perhaps most importantly, our framework is powerful enough to point to a range of novel methods. A notable example is {\tt S-Local-SVRG}, which is a local variance reduced {\tt SGD} method able to learn the optimal drift. This is the first time that local variance reduction is successfully combined with an on-the-fly learning of the local drift. Consequently, this is the first method which enjoys a linear convergence rate to the exact optimum (as opposed to a neighborhood of the solution only) without any restrictive assumptions and is thus superior in theory to the convergence of all existing local first order methods. We also develop another linearly converging method: {\tt S*-Local-SGD*}. Albeit not of practical significance as it depends on the a-priori knowledge of the optimal solution $x^*$, it is of theoretical interest as it enabled us to discover {\tt S-Local-SVRG}. See Table~\ref{tbl:special_cases} which summarizes all our complexity results.

\paragraph{Notation.} Due to its generality, our chapter is heavy in notation. For the reader's convenience, we present a notation table in Section~\ref{sec:notation_table} of the appendix.

\section{Our Framework}\label{sec:main_res}
In this section we present the main result of the chapter. Let us first introduce the key assumptions that we impose on our objective~\eqref{eq:main_problem}. We start with a relaxation of $\mu$-strong convexity (see also Assumptions~\ref{as:mu_strongly_quasi_convex} and \ref{ass:quasi_strong_convexity}).
\begin{assumption}[$(\mu,x^*)$-strong quasi-convexity]\label{ass:quasi_strong_convexity_local}
Let $x^*$ be a minimizer of $f$. We assume that  $f_i$ is $(\mu,x^*)$-strongly quasi-convex for all $i\in[n]$ with $\mu\geq 0$, i.e.\ for all $x\in\R^d$:
	\begin{equation}
	\mytextstyle	f_i(x^*) \ge f_i(x) + \langle\nabla f_i(x), x^* - x\rangle + \frac{\mu}{2}\|x - x^*\|^2. \label{eq:str_quasi_cvx}
	\end{equation}
\end{assumption}

Next, we require classical $L$-smoothness\footnote{While we require $L$-smoothness of $f_i$ to establish the main convergence theorem, some of the parameters of Assumption~\ref{ass:key_assumption} can be tightened considering a more complex smoothness structure of the local objective.} of local objectives, or equivalently, $L$-Lipschitzness of their gradients.
\begin{assumption}[$L$-smoothness]\label{ass:L_smoothness}
	Functions $f_i$ are $L$-smooth for all $i\in[n]$ with $L\geq 0$, i.e., 
	\begin{equation}
		\|\nabla f_i(x) - \nabla f_i(y)\| \le L\|x-y\|, \quad \forall x,y\in\R^d. \label{eq:L_smoothness}
	\end{equation}
\end{assumption}

In order to simplify our notation, it will be convenient to introduce the notion of virtual iterates $x^k$ defined as a mean of the local iterates~\citep{stich2020error}:
$
	\mytextstyle x^k \eqdef \frac{1}{n}\sum_{i=1}^n x_i^k. 
$
Despite the fact that $x^k$ is being physically computed only for $k$ for which  $c_k = 1$, virtual iterates are a very useful tool facilitating the convergence analysis. Next, we shall measure the discrepancy between the local and virtual iterates via the quantity $V_k$ defined as $
	\mytextstyle 	V_k \eqdef \frac{1}{n}\sum_{i=1}^n\|x_i^k - x^k\|^2. $

We are now ready to introduce the parametric assumption on both stochastic gradients $g_i^k$ and function $f$. This is a non-trivial generalization of the assumption from~\citep{gorbunov2019unified} to the class of local stochastic methods of the form~\eqref{eq:local_sgd_def}, and forms the heart of this work.\footnote{Recently, the assumption from~\citep{gorbunov2019unified} was generalized in a different way to cover  the class of the methods with error compensation and delayed updates~\citep{gorbunov2020linearly}.}

\begin{assumption}[Key parametric assumption]\label{ass:key_assumption}
	Assume that for all $k\ge 0$ and $i\in[n]$, local stochastic directions $g_i^k$ satisfy
	\begin{equation}
		\frac{1}{n}\sum\limits_{i=1}^n\EE_k\left[g_i^k\right] = \frac{1}{n}\sum\limits_{i=1}^n\nabla f_i(x_i^k), \label{eq:unbiasedness}
	\end{equation}
	where $\EE_k[\cdot]$ defines the expectation w.r.t.\ randomness coming from the $k$-th iteration only. Further, assume that there exist  non-negative constants $A, A', B, B', C, C', F, F',G, H, D_1, D_1',D_2,$ $D_3 \ge 0, \rho \in (0,1]$ and a sequence of (possibly random) variables $\{\sigma_k^2\}_{k\ge 0}$ such that
	\begin{align}
		\frac{1}{n}\sum\limits_{i=1}^n\EE\left[\|g_i^k\|^2\right] \le & 2A\EE\left[f(x^k) - f(x^*)\right] + B\EE\left[\sigma_k^2\right]  + F\EE\left[V_k\right] + D_1, \label{eq:second_moment_bound}\\
		\EE\left[\left\|\frac{1}{n}\sum\limits_{i=1}^ng_i^k\right\|^2\right] \le& 2A'\EE\left[f(x^k) - f(x^*)\right] + B'\EE\left[\sigma_k^2\right]  + F'\EE\left[V_k\right] + D_1', \label{eq:second_moment_bound_2}\\
		\EE\left[\sigma_{k+1}^2\right] \le& (1-\rho)\EE\left[\sigma_k^2\right] + 2C\EE\left[f(x^k) - f(x^*)\right]  + G\EE\left[V_k\right] + D_2,\label{eq:sigma_k+1_bound}\\
		2L\sum\limits_{k=0}^K w_k\EE[V_k] \le& \frac{1}{2}\sum\limits_{k=0}^Kw_k\EE\left[f(x^k) - f(x^*)\right]  + 2LH\EE\sigma_0^2+ 2LD_3\gamma^2 W_K ,\label{eq:sum_V_k_bounds}
	\end{align}
	where sequences $\{W_K\}_{K\ge 0}$, $\{w_k\}_{k\ge 0}$ are defined as
	\begin{equation}
		W_K \eqdef \sum\limits_{k=0}^K w_k,\quad w_k \eqdef \frac{1}{\left( 1 - \min\left\{\gamma\mu,\frac{\rho}{4}\right\}  \right)^{k+1}}, \label{eq:w_k_definition}
	\end{equation}

\end{assumption}

Admittedly, with its many parameters (whose meaning will become clear from the rest of the chapter), Assumption~\ref{ass:key_assumption} is not easy to parse on first reading. Several comments are due at this point. First, while the complexity of this assumption may be misunderstood as being problematic, the opposite is true. This assumption enables us to prove a single theorem (Thm.~\ref{thm:main_result}) capturing the convergence behavior, in a tight manner, of all local first-order methods described by our framework \eqref{eq:local_sgd_def}. So, the parametric and structural complexity of this assumption is paid for by the unification aspect it provides. Second, for each specific method we consider in this work, we {\em prove} that Assumption~\ref{ass:key_assumption} is satisfied, and each such proof is based on much simpler and generally accepted assumptions. So, Assumption~\ref{ass:key_assumption} should be seen as a ``meta-assumption'' forming an intermediary and abstract step in the analysis, one revealing the structure of the inequalities needed to obtain a general and tight convergence result for local first-order methods. We dedicate the rest of the chapter to explaining these parameters and to describing the algorithms and the associate rates their combination encodes. We are now ready to present our main convergence result.



\begin{theorem}\label{thm:main_result}
	Let Assumption~\ref{ass:quasi_strong_convexity_local},~\ref{ass:L_smoothness} and~\ref{ass:key_assumption} be satisfied and assume the stepsize satisfies $0<
		\gamma \le \min\left\{\frac{1}{2(A'+ \frac{4CB'}{3\rho})}, \frac{L}{F'+\frac{4GB'}{3\rho}}\right\}$. Define $\overline{x}^K \eqdef \frac{1}{W_K}\sum_{k=0}^K w_k x^k$,\newline $\Phi^0 \eqdef \frac{2\|x^0 - x^*\|^2 +   \frac{8B'}{3\rho}\gamma^2 \EE\sigma_0^2 + 4LH\gamma\EE\sigma_0^2}{\gamma}$ and $\Psi^0 \eqdef 2\left(D_1' + \frac{4B'}{3\rho}D_2 + 2L\gamma D_3\right)$. Let $\theta \eqdef 1 - \min\left\{\gamma\mu,\frac{\rho}{4}\right\}$. Then if $\mu>0$, we have
	\begin{align}
		\EE\left[f(\overline{x}^K) \right] - f(x^*) \le& \theta^K\Phi^0 + \gamma \Psi^0, \label{eq:main_result_1}
	\end{align}
	and in the case when $\mu = 0$, we have
	\begin{align}
	\mytextstyle
		\EE\left[f(\overline{x}^K) \right] - f(x^*)\le& \mytextstyle\frac{\Phi^0}{ K}  + \gamma \Psi^0. \label{eq:main_result_2} 
	\end{align}
\end{theorem}

As already mentioned, Thm.~\ref{thm:main_result} serves as a general, unified theory for local stochastic gradient algorithms.  The strongly convex case provides a linear convergence rate up to a specific neighborhood of the optimum. On the other hand, the weakly convex case yields an $\cO(K^{-1})$ convergence rate up to a particular neighborhood. One might easily derive  $\cO(K^{-1})$ and $\cO(K^{-\nicefrac{1}{2}})$ convergence rates to the exact optimum in the strongly and weakly convex case, respectively, by using a particular decreasing stepsize rule. The next corollary gives an example of such a result in the strongly convex scenario, where the estimate of $D_3$ does not depend on the stepsize $\gamma$. A detailed result that covers all cases is provided in Section~\ref{sec:corollaries} of the appendix.

\begin{corollary}\label{cor:main_complexity_cor_str_cvx}
Consider the setup from Thm.~\ref{thm:main_result} and by  $\frac1\nu$ denote  the resulting upper bound on $\gamma$.\footnote{In order to get tight estimate of $D_3$ and $H$, we will impose further bounds on $\gamma$ (see Tbl.~\ref{tbl:data_loop}). Assume that these extra bounds are included in parameter $h$.} Suppose that $\mu> 0$ and $D_3$ does not depend on $\gamma$. Let
\[
\mytextstyle
		\gamma = \min\left\{\frac{1}{\nu},  \frac{\ln\left(\max\left\{2,\min\left\{\frac{\Upsilon_1\mu^2K^2}{\Upsilon_2},\frac{\Upsilon_1\mu^3K^3}{\Upsilon_3}\right\}\right\}\right)}{\mu K}\right\},
\]
where $\Upsilon_1 = 2\|x^0 - x^*\|^2 +   \frac{8B'\EE\sigma_0^2}{3\nu^2\rho} + \frac{4LH\EE\sigma_0^2}{\nu}$, $\Upsilon_2 = 2D_1' + \frac{4B'D_2}{3\rho}$, $\Upsilon_3 = 4LD_3$. Then, the procedure~\eqref{eq:local_sgd_def} achieves $$\EE\left[f(\overline{x}^K) \right] - f(x^*)\le \varepsilon$$ as long as 
	\begin{equation*}
	\mytextstyle
		K\geq \widetilde\cO\left(\left(\frac{1}{\rho} + \frac{\nu}{\mu}\right)\log\left(\frac{\nu \Upsilon_1}{\varepsilon}\right) + \frac{\Upsilon_2}{\mu \varepsilon} + \sqrt{\frac{\Upsilon_3}{\mu^2 \varepsilon}}\right).
	\end{equation*}

\end{corollary}

\begin{remark}\label{rem:2nd_moment_decomposition}
Admittedly, Thm.~\ref{thm:main_result} does not yield the tightest known convergence rate in the heterogeneous setup under Assumption~\ref{thm:main_result}. Specifically, the neighborhood to which {\tt Local-SGD} converges can be slightly smaller~\citep{koloskova2020unified}. While we provide a tighter theory that matches the best-known results, we have deferred it to the appendix for the sake of clarity. In particular, to get the tightest rate, one shall replace the bound on the second moment of the stochastic direction~\eqref{eq:second_moment_bound} with two analogous bounds -- first one for the variance and the second one for the squared expectation. See Assumption~\ref{ass:hetero_second_moment} for details. Fortunately, Thm.~\ref{thm:main_result} does not need to change as it does not require parameters from~\eqref{eq:second_moment_bound}; these are only used later to derive $D_3, H, \gamma$ based on the data type. Therefore, only a few extra parameters should be determined in the specific scenario to get the tightest rate. 
\end{remark}

\begin{remark}
As we show in the appendix when looking at particular special cases,  local gradient methods are only as good as their non-local counterparts (i.e., when $\tau=1$) in terms of the communication complexity in the fully heterogeneous setup. Furthermore, the non-local methods outperform local ones in terms of  computation complexity. While one might think that this observation is a byproduct of our analysis, our observations are supported by findings in recent literature on this topic~\citep{karimireddy2020scaffold, khaled2020tighter}. To rise to the defense of local methods, we remark that they might be preferable to their non-local cousins in the homogeneous data setup~\citep{woodworth2020local} or for personalized federated learning~\citep{hanzely2020federated}.
\end{remark}

The parameters that drive both the convergence speed and the neighborhood size are determined by Assumption~\ref{ass:key_assumption}. In order to see through the provided rates, we shall discuss the value of these parameters in various scenarios. In general, we would like to have $\rho \in (0,1]$ as large as possible, while all other parameters are desired to be small so as to make the inequalities as tight as possible.

Let us start with studying data similarity and inner loop type as these can be decoupled from the type of the local direction that the method~\eqref{eq:local_sgd_def} takes.

\section{Data Similarity and Local Loop~\label{sec:data_and_loop}}
We now explain how our framework supports fixed and random local loop, and several data similarity regimes. 

\paragraph{Local loop.} Our framework supports \emph{local loop of a fixed length} $\tau \geq 1$ (i.e., we support local methods performing $\tau$ local iterations in between communications). This option, which is the de facto standard for local methods in theory and practice~\citep{mcmahan2016federated}, is recovered by setting $c_{a\tau } = 1$ for all non-negative integers $a$ and $c_k = 0$ for $k$ that are not divisible by $\tau$ in~\eqref{eq:local_sgd_def}. However, our framework also captures the very rarely considered \emph{local loop with a random length}. We recover this when $c_k$ are random samples from the Bernoulli distribution $\text{Be}(p)$ with parameter $p\in(0,1]$.

\paragraph{Data similarity.} We look at various possible data similarity regimes. The first option we consider is the fully heterogeneous setting where we do not assume any similarity between the local objectives whatsoever. Secondly, we consider the identical data regime with $f_1=\ldots=f_n$. Lastly, we consider the $\zeta$-heterogeneous data setting, which bounds the dissimilarity between the full and the local gradients~\citep{woodworth2020minibatch} (see Def.~\ref{def:zeta_hetero}).

\begin{definition}[$\zeta$-heterogeneous functions]\label{def:zeta_hetero}
	We say that functions $f_1,\ldots,f_n$ are\newline $\zeta$-heterogeneous for some $\zeta\geq 0$ if the following inequality holds for all $x\in\R^d$:
	\begin{equation}
\mytextstyle
		\frac{1}{n}\sum\limits_{i=1}^n\|\nabla f_i(x) - \nabla f(x)\|^2 \le \zeta^2. \label{eq:bounded_data_dissimilarity}
	\end{equation}
\end{definition}

The $\zeta$-heterogeneous data regime recovers the heterogeneous data for $\zeta =\infty$ and identical data for $\zeta = 0$.

In Section~\ref{sec:a_data_and_loop} of the appendix, we show that  the local loop type and the data similarity type affect parameters $H$ and $D_3$ from Assumption~\ref{ass:key_assumption} only. However, in order to obtain an efficient bound on these parameters, we impose additional constraints on the stepsize $\gamma$. While we do not have space to formally state our results in the main body, we provide a comprehensive summary in Tbl.~\ref{tbl:data_loop}.

\begin{table}[!t]
\caption{The effect of data similarity and local loop on Assumption~\ref{ass:key_assumption}. Constant factors are ignored. Homogeneous data are recovered as a special case of  $\zeta$-heterogeneous data with $\zeta=0$. Heterogeneous case is slightly loose in light of Remark~\ref{rem:2nd_moment_decomposition}. If one replaces the bound on the second moments~\eqref{eq:second_moment_bound} with a analogous bound on variance squared expectation (see Assumption~\ref{ass:hetero_second_moment}), the bounds on $\gamma$, $D_3$ and $H$ will have $(\tau-1)$ times better dependence on the variance parameters (or $\frac{1-p}{p}$ times for the random loop). See Section~\ref{sec:const_loop_hetero} and~\ref{sec:random_local_loop_hetero} of appendix for more details.}
\label{tbl:data_loop}
\begin{center}
\footnotesize
\begin{tabular}{|c|c|c|c|c|}
\hline
 Data &   Loop & Extra upper bounds on $\gamma$ & $D_3$ & $H$   \\
\hline
\hline
het  & fixed & {   $\frac{1}{\tau\mu}, \frac{1}{\tau\sqrt{\left(F+\frac{BG}{\rho(1-\rho)}\right)}}, \frac{1}{\tau\sqrt{2L\left(A + \frac{BC}{\rho(1-\rho)}\right)}}$} & {  $ (\tau-1)^2\left(D_1 + \frac{BD_2}{\rho}\right)$} & {  $\frac{B(\tau-1)^2\gamma^2}{\rho}$}  \\
 \hline
 $\zeta$-het  & fixed & {   $\frac{1}{\tau\mu}, \frac{1}{\sqrt{\tau\left(F+\frac{BG}{\rho(1-\rho)}\right)}}, \frac{1}{\sqrt{L\tau\left(A + \frac{BC}{\rho(1-\rho)}\right)}}$} & {  $  (\tau-1)\left(D_1 + \frac{\zeta^2}{\gamma\mu} + \frac{BD_2}{\rho}\right)$} & {  $\frac{B(\tau-1)\gamma^2}{\rho}$}  \\
 \hline
 het  & random & { $\frac{p}{\mu}$,   $\frac{p}{\sqrt{(1-p)F}}, \frac{p\sqrt{\rho(1-\rho)}}{\sqrt{BG(1-p)}}, \frac{p}{\sqrt{L(1-p)\left(A + \frac{BC}{\rho(1-\rho)}\right)}}$} & {  $ \frac{(1-p)\left(D_1 + \frac{BD_2}{\rho}\right) }{p^2}$} & {  $\frac{B(1-p)\gamma^2}{p^2\rho}$}  \\
 \hline
 $\zeta$-het  & radnom & {  $\frac{p}{\mu}$,  $\sqrt{\frac{p}{F(1-p)}}, \sqrt{\frac{p\rho(1-\rho)}{BG(1-p)}}, \sqrt{\frac{p}{L(1-p)\left(A + \frac{BC}{\rho(1-\rho)}\right)}}$} & {  $ \frac{(1-p)}{p}\left(D_1 + \frac{\zeta^2}{\gamma\mu} + \frac{BD_2}{\rho}\right)$} & {  $\frac{B(1-p)\gamma^2}{p\rho}$}  \\
 \hline
\end{tabular}
\end{center}
\end{table}

Methods with a random loop communicate once per $p^{-1}$ iterations on average, while the fixed loop variant communicates once every $\tau$ iterations. Consequently, we shall compare the two loop types for $\tau= p^{-1}$. In such a case, parameters $D_3$ and $H$ and the extra conditions on stepsize $\gamma$ match exactly, meaning that the loop type does not influence the convergence rate. 
Having said that, random loop choice provides more flexibility compared to the fixed loop. Indeed, one might want the local direction $g_i^k$ to be synchronized with the communication time-stamps in some special cases. However, our framework does not allow such synchronization for a fixed loop since we assume that the local direction $g_i^k$ follows some stationary distribution over stochastic gradients. The random local loop comes in handy here; the random variable that determines the communication follows a stationary distribution, thus possibly synchronized with the local computations.

\section{Local Stochastic Direction \label{sec:local_solver}}

This section discusses how the choice of $g_i^k$ allows us to obtain the remaining parameters from Assumption~\ref{ass:key_assumption} that were not covered in the previous section. To cover the most practical scenarios, we set $g_i^k$ to be a difference of two components $a_i^k,b_i^k\in \R^d$, which we explain next. We stress that the construction of $g_i^k$ is very general: we recover various state-of-the-art methods along with their rates while covering many new interesting algorithms. We will discuss this in more detail in Section~\ref{sec:special_cases}.

\subsection{Unbiased Local Gradient Estimator $a_i^k$}
The first component of the local direction that the method~\eqref{eq:local_sgd_def} takes is $a_i^k$ -- an unbiased, possibly variance reduced, estimator of the local gradient, i.e., $\EE_k[a_i^k] = \nabla f_i(x^k_i)$. Besides the unbiasedness, $a_i^k$ is allowed to be anything that satisfies the parametric recursive relation from~\citep{gorbunov2019unified}, which tightly covers many variants of {\tt SGD} including non-uniform, minibatch, and variance reduced stochastic gradient. The parameters of such a relation are capable of encoding both the general smoothness structure of the objective and the gradient estimator's properties that include a diminishing variance, for example. We state the adapted version of this recursive relation as Assumption~\ref{ass:sigma_k_original}. 

\begin{assumption} \label{ass:sigma_k_original}
Let the unbiased local gradient estimator $a_i^k$ be such that
\begin{align*}
&		\EE_k\left[\|a_i^k- \nabla f_i(x^*)\|^2 \right] \le 2A_iD_{f_i} (x^k_i, x^*)+ B_i\sigma_{i,k}^2 + D_{1,i}, 
		\\
&		\EE_k\left[\sigma_{i,k+1}^2\right] \le (1-\rho_i)\sigma_{ik}^2 + 2C_iD_{f_i} (x^k_i, x^*)  +  D_{2,i}
	\end{align*}
for $A_i\geq0, B_i\geq0, D_{1,i} \geq0, 0\leq \rho_i \leq 1,C_i \geq0, D_{2,i} \geq0$ and a non-negative sequence $ \{\sigma^2_{i,k}\}_{k=0}^\infty$.\footnote{By $D_{f_i}(x_i^k, x^k)$ we mean Bregman distance between $x_i^k, x^k$ defined as $D_{f_i}(x_i^k, x^k) \eqdef f_i(x_i^k) - f_i(x^k) - \langle \nabla f_i(x^k), x_i^k - x^k\rangle$.}
\end{assumption}

Note that the parameters of Assumption~\ref{ass:sigma_k_original} can be taken directly from~\citep{gorbunov2019unified} and offer a broad range of unbiased local gradient estimators $a_i^k$ in different scenarios. The most interesting setups covered include minibatching, importance sampling, variance reduction, all either under the classical smoothness assumption or under a uniform bound on the stochastic gradient variance. 

Our next goal is to derive the parameters of Assumption~\ref{ass:key_assumption} from the parameters of Assumption~\ref{ass:sigma_k_original}. However, let us first discuss the second component of the local direction -- the local shift $b_i^k$.

\subsection{Local Shift $b_i^k$} 

The local update rule~\eqref{eq:local_sgd_def} can include the local shift/drift $b_i^k$ allowing us to eliminate the infamous non-stationarity of the local methods. The general requirement for the choice of $b_i^k$ is so that it sums up to zero ($\sum_{i=1}^n b_i^k = 0$) to avoid unnecessary extra bias. For the sake of simplicity (while maintaining generality), we will consider three choices of $b_i^k$ -- zero, ideal shift ($=\nabla f_i(x^*)$) and on-the-fly shift via a possibly outdated local stochastic non-variance reduced gradient estimator that satisfies a similar bound as Assumption~\ref{ass:sigma_k_original}.

\begin{assumption} \label{ass:bik}
Consider the following choices:\\
Case I: $b_i^k = 0 $, \\
Case II: $b_i^k = \nabla f_i(x^*)  $, \\
Case III: $b_i^k = h_i^k  - \frac1n \sum_{i=1}^n h_i^k$ where $h_i^k\in \R^d$ is a delayed local gradient estimator defined recursively as
$$
h_i^{k+1} = \begin{cases}
h_i^k & \text{with probability } 1-\rho_i' \\
l_i^k &  \text{with probability } \rho_i' \
\end{cases},
$$
 where $0\leq \rho'_i \leq 1$ and $l_i^k\in \R^d$ is an unbiased non-variance reduced possibly stochastic gradient estimator of $\nabla f_i(x^k)$ such that for some $A'_i, D_{3,i}\geq 0$ we have
\begin{equation}\label{eq:bdef}
\EE_k\left[\|l_i^k- \nabla f_i(x^*)\|^2 \right] \le 2A'_iD_{f_i} (x^k_i, x^*)+ D_{3,i}.
\end{equation}

\end{assumption}

Let us look closer at Case III as this one is the most interesting. Note that what we assume about $l_i^k$ (i.e., ~\eqref{eq:bdef}) is essentially a variant of Assumption~\ref{ass:bik} with $\sigma^2_{i,k}$ parameters set to zero. This is achievable for a broad range of non-variance reduced gradient estimators that includes minibatching and importance sampling~\citep{gower2019sgd}. An intuitive choice of $l_i^k$ is to set it to $a_i^k$ given that $a_i^k$ is not variance reduced. In such a case, the scheme~\eqref{eq:local_sgd_def} reduces to {\tt SCAFFOLD}~\citep{karimireddy2020scaffold} along with its rate.  

However, our framework can do much more beyond this example. First, we cover the local variance reduced gradient $a_i^k$ with $l_i^k$ constructed as its non-variance reduced part. In such a case, the neighborhood of the optimum from Thm.~\ref{thm:main_result} to which the method~\eqref{eq:local_sgd_def} converges shrinks. There is a way to get rid of this neighborhood, noticing that $l_i^k$ is used only once in a while. Indeed, the combination of the full local gradient $l_i^k$ together with the variance reduced $a_i^k$ leads to a linear rate in the strongly (quasi) convex case or $\cO(K^{-1})$ rate in the weakly convex case. We shall remark that the variance reduced gradient might require a sporadic computation of the full local gradient -- it makes  sense to synchronize it with the update rule for $h_i^k$. In such a case, the computation of $l_i^k$ is for free. We have just described the {\tt S-Local-SVRG} method (Algorithm~\ref{alg:l_local_svrg_fs}). 

\subsection{Parameters of Assumption~\ref{ass:key_assumption}}

We proceed with a key lemma that provides us with the remaining parameters of Assumption~\ref{ass:key_assumption} that were not covered in Section~\ref{sec:data_and_loop}. These parameters will be chosen purely based on the selection of $a_i^k$ and $b_i^k$ discussed earlier.

\begin{lemma}\label{lem:local_solver}
For all $i \in [n]$ suppose that $a_i^k$ satisfies Assumption~\ref{ass:sigma_k_original}, while $b_i^k$ was chosen as per Assumption~\ref{ass:bik}. Then,~\eqref{eq:second_moment_bound}, \eqref{eq:second_moment_bound_2} and \eqref{eq:sigma_k+1_bound} hold with
	\begin{align*}
	A &= 4\max_i A_{i}, B = 2, F = 4 L \max_i A_{i},\\ 
	D_1 &= \begin{cases}
	\frac{2}{n} \sum_{i=1}^n \left( D_{1,i}  + \| \nabla f_i(x^*)\|^2\right) & \text{Case I}, \\
	\frac{2}{n} \sum_{i=1}^n  D_{1,i} & \text{Case II, III},
	\end{cases}
	\\
 B' &\mytextstyle= \frac1n, F' = \frac{2 L \max_i A_{i}}{n}+2L^2, D_1' = \frac{1}{n^2} \sum_{i=1}^n D_{1,i} \\
  A' &\mytextstyle = \frac{2\max_i A_{i}}{n} + L , G = \nicefrac{CL}{2}, \\
	\rho &=\begin{cases} \min_i \rho_i  & \text{Case I, II}, \\
	  \min_i \min \left\{ \rho_i, \rho'_i\right\} & \text{Case III,}
	 \end{cases}    \\
	 D_2  &=\begin{cases}\frac2n\sum\limits_{i=1}^n B_iD_{2,i} , & \text{Case I, II}, \\
  \frac1n\sum\limits_{i=1}^n\left( 2B_iD_{2,i} + \rho_i' D_{3,i} \right)& \text{Case III,}
	 \end{cases}    \\
	 C &=\begin{cases}
	4 \max_{i}\{B_iC_i\}   & \text{Case I, II}, \\
		4 \max_{i}\{B_iC_i \} +4\max_i\{\rho_i' A'_i\}& \text{Case III}.
	\end{cases}
	\end{align*}

	\end{lemma}

We have just broken down the parameters of Assumption~\ref{ass:key_assumption} based on the optimization objective and the particular instance of~\eqref{eq:local_sgd_def}. However, it might still be hard to understand particular rates based on these choices. In the appendix, we state a range of methods and decouple their convergence rates. A summary of the key parameters from Assumption~\ref{ass:key_assumption} is provided in Tbl.~\ref{tbl:special_cases-parameters}.

\section{Special Cases}\label{sec:special_cases}

Our theory covers a broad range of local stochastic gradient algorithms. While we are able to recover multiple known methods along with their rates, we also introduce several new methods along with extending the analysis of known algorithms. As already mentioned, our theory covers convex and strongly convex cases, identical and heterogeneous data regimes. From the algorithmic point of view, we cover the fixed and random loop, various shift types, and arbitrary local stochastic gradient estimator. We stress that our framework gives a tight convergence rate under any circumstances.

While we might not cover all of these combinations in a deserved detail, we thoroughly study a subset of them in the following subsections. An overview of these methods is presented in Tbl.~\ref{tbl:special_cases} together with their convergence rates in the strongly convex case (see Tbl.~\ref{tbl:special_cases_weakly_convex} for the rates in the weakly convex setting). Next, we describe a selected number of special cases of our framework.

$\bullet$  {\bf Non-local stochastic methods.} Our theory recovers a broad range of non-local stochastic methods. In particular, if $n=1$, we have $V_k = 0$, and consequently we can choose $A=A', B=B', D_1=  D_1', F=F'=G= H=D_3=0$. With such a choice, our theory matches\footnote{Up to the non-smooth regularization/proximal steps and small constant factors.} the general analysis of stochastic gradient methods from~\citep{gorbunov2019unified} for $\tau=1$. Consequently, we recover a broad range of algorithms as a special case along with their convergence guarantees, namely {\tt SGD}~\citep{RobbinsMonro:1951} with its best-known rate on smooth objectives~\citep{nguyen2018sgd, gower2019sgd}, variance reduced finite sum algorithms such as {\tt SAGA}~\citep{SAGA}, {\tt SVRG}~\citep{SVRG},  {\tt L-SVRG}~\citep{hofmann2015variance, kovalev2019don}, variance reduced subspace descent methods such as {\tt SEGA/SVRCD}~\citep{hanzely2018sega, hanzely2019one}, quantized methods~\citep{mishchenko2019distributed, horvath2019stochastic} and others. 

$\bullet$  {\bf ``Star''-shifted local methods.} As already mentioned, local methods have inherently incorrect fixed points~\citep{pathak2020fedsplit}; and one can fix these by shifting the local gradients. Star-shifted local methods employ the ideal stationary shift using the local gradients at the optimum $b_i^k = \nabla f_i(x^*)$ (i.e.,  Case II from Assumption~\ref{ass:bik}) and serve as a transition from the plain local methods (Case I from Assumption~\ref{ass:bik}) to the local methods that shift using past gradients such as {\tt SCAFFOLD} (Case III from Assumption~\ref{ass:bik}). In the appendix, we present two such methods: {\tt S*-Local-SGD} (Algorithm~\ref{alg:local_sgd_star}) and  {\tt S*-Local-SGD*} (Algorithm~\ref{alg:local_sgd_star_star}).  While being impractical in most cases since $\nabla f_i(x^*)$ is not known, star-shifted local methods give new insights into the role and effect of the shift for local algorithms. Specifically, these methods enjoy superior convergence rate when compared to methods without local shift (Case I) and methods with a shift constructed from observed gradients (Case III), while their rate serves as an aspiring goal for local methods in general. Fortunately, in several practical scenarios, one can match the rate of star methods using an approach from Case III, as we shall see in the next point. 

$\bullet$  {\bf Shifted Local {\tt SVRG} ({\tt S-Local-SVRG}).} As already mentioned, local {\tt SGD} suffers from convergence to a neighborhood of the optimum only, which is credited to i) inherent variance of the local stochastic gradient, and ii) incorrect fixed point of local {\tt GD}. We propose a way to correct both issues. To the best of our knowledge, this is the first time that on-device variance reduction was combined with the trick for reducing the non-stationarity of local methods. Specifically, the latter is achieved by selecting $b_i^k$ as a particular instance of Case III from Assumption~\ref{ass:bik} such that $l_i^k$ is the full local gradient, which in turns yields $D'_{1,i}=0, A'_i = L$. In order to not waste local computation, we synchronize the evaluation of $l_i^k$ with the computation of the full local gradient for the {\tt L-SVRG}~\citep{hofmann2015variance, kovalev2019don} estimator, which we use to construct $a_i^k$. Consequently, some terms cancel out, and we obtain a simple, fast, linearly converging local {\tt SGD} method, which we present as Algorithm~\ref{alg:l_local_svrg_fs}. We believe that this is remarkable since only a very few local methods converge linearly to the exact optimum.\footnote{A linearly converging local {\tt SGD} variant can be recovered from stochastic decoupling~\citep{mishchenko2019stochastic}, although this was not considered therein. Besides that, FedSplit~\citep{pathak2020fedsplit} achieves a linear rate too, however, with a much stronger local oracle.}

\begin{table}[H]
\caption{A selection of methods that can be analyzed using our framework, which we detail in the appendix. A choice of $a_i^k, b_i^k$ and $l_i^k$ is presented along with the established complexity bounds (= number of iterations to find such $\hat x$ that $\EE[f(\hat{x}) - f(x^*)]\le \varepsilon$) and a specific setup under which the methods are analyzed. For Algorithms 1-4 we suppress constants and $\log \tfrac{1}{\varepsilon}$ factors. Since Algorithms 5 and 6 converge linearly, we suppress constants only while keeping $\log \tfrac{1}{\varepsilon}$ factors. All rates are provided in the \textbf{strongly convex} setting. UBV stands for the ``Uniform Bound on the Variance'' of local stochastic gradient, which is often assumed when $f_i$ is of the form~\eqref{eq:f_i_expectation}. ES stands for the ``Expected Smoothness''~\citep{gower2019sgd}, which does not impose any extra assumption on the objective/noise, but rather can be derived given the sampling strategy and the smoothness structure of $f_i$. Consequently, such a setup allows us to obtain local methods with importance sampling. Next, the simple setting is a special case of ES when we uniformly sample a single index on each node each iteration. $^\clubsuit$: {\tt Local-SGD} methods have never been analyzed under ES assumption. Notation: $\sigma^2$ -- averaged (within nodes) uniform upper bound for the variance of local stochastic gradient, $\sigma_*^2$ -- averaged variance of local stochastic gradients at the solution, $\zeta_*^2 \eqdef \frac1n\sum_{i=1}^n \| \nabla f_i(x^*)\|^2$, $\max L_{ij}$ -- the worst smoothness of $f_{i,j}, i\in[n],j\in[m]$, $\cL$ -- the worst ES constant for all nodes.
}
\label{tbl:special_cases}
\begin{center}
\scriptsize
\begin{tabular}{|c|c|c|c|c|c|}
\hline
{\bf Method}  &   $\myred{a_i^k}, \myblue{b_i^k}, l_i^k$   & {\bf Complexity}  &   {\bf Setting}  & {\bf Sec}  \\
\hline
\hline
 \begin{tabular}{c}
	{\tt Local-SGD} \\ Alg.~\ref{alg:local_sgd}, {\tiny\citep{woodworth2020minibatch}} 
\end{tabular}   & $\myred{f_{\xi_i}(x_i^k)}, \myblue{0}, - $& $\frac{L}{\mu}+\frac{\sigma^2}{n\mu\varepsilon}+\sqrt{\frac{L\tau(\sigma^2 + \tau\zeta^2)}{\mu^2\varepsilon}}$ &  \begin{tabular}{c}
	UBV,\\
	$\zeta$-Het
\end{tabular} & \ref{sec:sgd_bounded_var}  \\
\hline
 \begin{tabular}{c}
 {\tt Local-SGD} \\ Alg.~\ref{alg:local_sgd}, {\tiny\citep{koloskova2020unified}}
\end{tabular}   & $\myred{f_{\xi_i}(x_i^k)},\myblue{0}, - $& $\frac{\tau L}{\mu}+\frac{\sigma^2}{n\mu\varepsilon}+\sqrt{\frac{L(\tau-1)(\sigma^2 + (\tau-1)\zeta_*^2)}{\mu^2\varepsilon}}$ &  \begin{tabular}{c}
	UBV,\\
	Het
\end{tabular} & \ref{sec:sgd_bounded_var}  \\
\hline
 \begin{tabular}{c}
 {\tt Local-SGD}\\ Alg.~\ref{alg:local_sgd},   {\tiny\citep{khaled2020tighter}}$^{\clubsuit}$
\end{tabular}   & $\myred{f_{\xi_i}(x_i^k)}, \myblue{0}, - $  & \begin{tabular}{c}
 $\frac{L+\nicefrac{\cL}{n}+\sqrt{(\tau-1) L\cL}}{\mu}+\frac{\sigma_*^2}{n\mu\varepsilon}\quad\quad$\\
 $+ \frac{L\zeta^2(\tau-1)}{\mu^2\varepsilon}+\sqrt{\frac{L(\tau-1)(\sigma_*^2 + \zeta_*^2)}{\mu^2\varepsilon}}$
\end{tabular}  &   \begin{tabular}{c}
	ES,\\
	$\zeta$-Het
\end{tabular} & \ref{sec:sgd_es}  \\
\hline
 \begin{tabular}{c}
 {\tt Local-SGD}\\ Alg.~\ref{alg:local_sgd}, {\tiny\citep{khaled2020tighter}}$^{\clubsuit}$
\end{tabular}  & $\myred{f_{\xi_i}(x_i^k)}, \myblue{0}, - $  & \begin{tabular}{c}
	$\frac{L\tau+\nicefrac{\cL}{n}+\sqrt{(\tau-1)L\cL}}{\mu}+\frac{\sigma_*^2}{n\mu\varepsilon}\quad\quad$\\
	$\quad\quad+\sqrt{\frac{L(\tau-1)(\sigma_*^2 + (\tau-1)\zeta_*^2)}{\mu^2\varepsilon}}$ 
\end{tabular}  &   \begin{tabular}{c}
	ES,\\
	Het
\end{tabular} & \ref{sec:sgd_es}  \\
\hline
\rowcolor{bgcolor} \begin{tabular}{c}
	{\tt Local-SVRG}\\ Alg.~\ref{alg:local_svrg},  {\tiny(NEW)}
\end{tabular}  & \begin{tabular}{c}
	\myred{$\nabla f_{i,j_i}(x^k_i) -  \nabla f_{i,j_i}(y_i^k)$}\\\myred{$ +  \nabla f_{i}(y_i^k)$},\\ \myblue{$0$}, $- $ 
\end{tabular}  & \begin{tabular}{c}
	 $m+\frac{L+\nicefrac{\max L_{ij}}{n}+\sqrt{(\tau-1) L\max L_{ij}}}{\mu}$\\
	 $\quad+ \frac{L\zeta^2(\tau-1)}{\mu^2\varepsilon}+\sqrt{\frac{L(\tau-1)\zeta_*^2}{\mu^2\varepsilon}}$
\end{tabular} & \begin{tabular}{c}
	simple,\\
	$\zeta$-Het
\end{tabular} & \ref{sec:llsvrg}   \\
\hline
\rowcolor{bgcolor} \begin{tabular}{c}
	{\tt Local-SVRG}\\ Alg.~\ref{alg:local_svrg},  {\tiny(NEW)}
\end{tabular}  & \begin{tabular}{c}
	\myred{$\nabla f_{i,j_i}(x^k_i) -  \nabla f_{i,j_i}(y_i^k)$}\\\myred{$ +  \nabla f_{i}(y_i^k)$},\\ \myblue{$0$}, $- $ 
\end{tabular}  & \begin{tabular}{c}
	 $m+\frac{L\tau+\nicefrac{\max L_{ij}}{n}+\sqrt{(\tau-1)L\max L_{ij}}}{\mu}$\\
	 $+\sqrt{\frac{L(\tau-1)^2\zeta_*^2}{\mu^2\varepsilon}}$
\end{tabular} & \begin{tabular}{c}
	simple,\\
	Het
\end{tabular} & \ref{sec:llsvrg}   \\
\hline
\rowcolor{bgcolor} \begin{tabular}{c}
	{\tt S*-Local-SGD}\\ Alg.~\ref{alg:local_sgd_star},  {\tiny(NEW)}
\end{tabular}  & $\myred{f_{\xi_i}(x_i^k)}, \myblue{\nabla f_i(x^*)}, - $ & $\frac{\tau L}{\mu}+\frac{\sigma^2}{n\mu\varepsilon}+\sqrt{\frac{L(\tau-1)\sigma^2}{\mu^2\varepsilon}}$  &  \begin{tabular}{c}
	UBV,\\
	Het
\end{tabular} & \ref{sec:sgd_star_bounded_var}  \\
\hline
 \begin{tabular}{c}
	{\tt SS-Local-SGD}\\ Alg.~\ref{alg:l_local_svrg}, {\tiny\citep{karimireddy2020scaffold}} 
\end{tabular}   & \begin{tabular}{c}
$\myred{f_{\xi_i}(x_i^k)}, \myblue{h_i^k  - \frac1n \sum_{i=1}^n h_i^k},$\\
$ \nabla f_{\tilde{\xi}_i^k} (y_i^k)$ 
\end{tabular}  & $\frac{L}{p\mu}+\frac{\sigma^2}{n\mu\varepsilon}+\sqrt{\frac{L(1-p)\sigma^2}{p\mu^2\varepsilon}}$ &  \begin{tabular}{c}
	UBV,\\
	Het
\end{tabular} & \ref{sec:loopless_local_svrg}  \\
\hline
\rowcolor{bgcolor} \begin{tabular}{c}
	{\tt SS-Local-SGD}\\ Alg.~\ref{alg:l_local_svrg}, {\tiny(NEW)}
\end{tabular}  & \begin{tabular}{c}
$\myred{f_{\xi_i}(x_i^k)}, \myblue{h_i^k  - \frac1n \sum_{i=1}^n h_i^k},$\\
$ \nabla f_{\tilde{\xi}_i^k} (y_i^k)$ 
\end{tabular}  &  \begin{tabular}{c}
$\frac{L}{p\mu} + \frac{\cL}{n\mu} + \frac{\sqrt{L\cL(1-p)}}{p\mu}$\\ $+\frac{\sigma_*^2}{n\mu\varepsilon}+\sqrt{\frac{L(1-p)\sigma_*^2}{p\mu^2\varepsilon}}$
\end{tabular} &  \begin{tabular}{c}
	ES,\\
	Het
\end{tabular} & \ref{sec:loopless_local_svrg_es}  \\
\hline
\rowcolor{bgcolor2} \begin{tabular}{c}
	{\tt S*-Local-SGD*}\\ Alg.~\ref{alg:local_sgd_star_star}, {\tiny(NEW)}
\end{tabular}  & \begin{tabular}{c}
 	\myred{$\nabla f_{i,j_i}(x^k_i) -  \nabla f_{i,j_i}(x^*)$}\\
 	\myred{$ +  \nabla f_{i}(x^*)$}, $\myblue{\nabla f_i(x^*)}, - $\\ 	
 \end{tabular} &  \begin{tabular}{c}
$\Big(\frac{\tau L}{\mu}+\frac{\max L_{ij}}{n\mu}\quad\quad\quad\quad\quad\quad\quad\quad$\\ $\quad\quad\quad\quad+ \frac{\sqrt{(\tau-1)L\max L_{ij}}}{\mu}\Big)\log\frac{1}{\varepsilon}$ 
\end{tabular}   &  \begin{tabular}{c}
	simple,\\
	Het
\end{tabular} & \ref{sec:S*-Local-SGD*}  \\
\hline
\rowcolor{bgcolor2} \begin{tabular}{c}
	{\tt S-Local-SVRG}\\ Alg.~\ref{alg:l_local_svrg_fs},  {\tiny(NEW)}
\end{tabular}  &\begin{tabular}{c}
	$\myred{ \nabla f_{i,j_i}(x^k_i) -  \nabla f_{i,j_i}(y^k)}$\\
	\myred{$ +  \nabla f_{i}(y^k) $},\\
	$\myblue{h_i^k  - \frac1n \sum_{i=1}^n h_i^k}, \nabla f_{i} (y^k)$ 
\end{tabular} & \begin{tabular}{c}
$\Big(m + \frac{L}{p\mu}+\frac{\max L_{ij}}{n\mu}\quad\quad\quad\quad\quad\quad$\\ $\quad\quad\quad\quad+ \frac{\sqrt{L\max L_{ij}(1-p)}}{p\mu}\Big)\log\frac{1}{\varepsilon}$ 
\end{tabular} & \begin{tabular}{c}
	simple,\\
	Het
\end{tabular} &   \ref{sec:loopless_local_svrg_fs} \\
\hline
\end{tabular}
\end{center}
\vskip -0.2cm
\end{table}

\begin{table}[H]
\caption{A selection of methods that can be analyzed using our framework. A choice of $a_i^k, b_i^k$ and $l_i^k$ is presented along with the established complexity bounds (= number of iterations to find such $\hat x$ that $\EE[f(\hat{x}) - f(x^*)]\le \varepsilon$) and a specific setup under which the methods are analyzed. For all algorithms we suppress constants factors. All rates are provided in the \textbf{weakly convex} setting. UBV stands for the ``Uniform Bound on the Variance'' of local stochastic gradient, which is often assumed when $f_i$ is of the form~\eqref{eq:f_i_expectation}. ES stands for the ``Expected Smoothness''~\citep{gower2019sgd}, which does not impose any extra assumption on the objective/noise, but rather can be derived given the sampling strategy and the smoothness structure of $f_i$. Consequently, such a setup allows us to obtain local methods with importance sampling. Next, the simple setting is a special case of ES when we uniformly sample a single index on each node each iteration. $^\clubsuit$: {\tt Local-SGD} methods have never been analyzed under ES assumption. Notation: $\sigma^2$ -- averaged (within nodes) uniform upper bound for the variance of local stochastic gradient, $\sigma_*^2$ -- averaged variance of local stochastic gradients at the solution, $\zeta_*^2 \eqdef \frac1n\sum_{i=1}^n \| \nabla f_i(x^*)\|^2$, $\max L_{ij}$ -- the worst smoothness of $f_{i,j}, i\in[n],j\in[m]$, $\cL$ -- the worst ES constant for all nodes, $R_0\eqdef \|x^0 - x^*\|$ -- distance of the starting point $x^0$ from the closest solution $x^*$, $\Delta_0 \eqdef f(x^0)-f(x^*)$.
}
\label{tbl:special_cases_weakly_convex}
\begin{center}
\vskip -0.2cm
\scriptsize
\begin{tabular}{|c|c|c|c|c|c|}
\hline
{\bf Method} &   $\myred{a_i^k}, \myblue{b_i^k}, l_i^k$   & {\bf Complexity}  &   {\bf Setting}  & {\bf Sec}  \\
\hline
\hline
 \begin{tabular}{c}
	{\tt Local-SGD}\\  Alg.~\ref{alg:local_sgd}, {\tiny\citep{woodworth2020minibatch}} 
\end{tabular} & $\myred{f_{\xi_i}(x_i^k)}, \myblue{0}, - $& \begin{tabular}{c}
$\frac{L R_0^2}{\varepsilon}+\frac{\sigma^2R_0^2}{n\varepsilon^2}\quad\quad\quad\quad\quad\quad\quad\quad$\\
$\quad\quad\quad+\frac{R_0^2\sqrt{L\tau(\sigma^2 + \tau\zeta^2)}}{\varepsilon^{\nicefrac{3}{2}}}$ 
\end{tabular} &  \begin{tabular}{c}
	UBV,\\
	$\zeta$-Het
\end{tabular} & \ref{sec:sgd_bounded_var}  \\
\hline
 \begin{tabular}{c}
 {\tt Local-SGD}\\ Alg.~\ref{alg:local_sgd},  {\tiny\citep{koloskova2020unified}}
\end{tabular}  & $\myred{f_{\xi_i}(x_i^k)},\myblue{0}, - $& \begin{tabular}{c}
$\frac{\tau L R_0^2}{\varepsilon}+\frac{\sigma^2R_0^2}{n\varepsilon^2}\quad\quad\quad\quad\quad\quad\quad\quad$\\
$\quad\quad\quad+\frac{R_0^2\sqrt{L(\tau-1)(\sigma^2 + (\tau-1)\zeta_*^2)}}{\varepsilon^{\nicefrac{3}{2}}}$ 
\end{tabular}   &  \begin{tabular}{c}
	UBV,\\
	Het
\end{tabular} & \ref{sec:sgd_bounded_var}  \\
\hline
 \begin{tabular}{c}
 {\tt Local-SGD}\\ Alg.~\ref{alg:local_sgd},  {\tiny\citep{khaled2020tighter}}$^{\clubsuit}$
\end{tabular}  & $\myred{f_{\xi_i}(x_i^k)}, \myblue{0}, - $  & \begin{tabular}{c}
	$\frac{\left(L+\nicefrac{\cL}{n}+\sqrt{(\tau-1)L\cL}\right)R_0^2}{\varepsilon}+\frac{\sigma_*^2R_0^2}{n\varepsilon^2}\quad$\\
	$+\frac{L\zeta^2(\tau-1)R_0^2}{\mu\varepsilon^2}+\frac{R_0^2\sqrt{L(\tau-1)(\sigma_*^2 + \zeta_*^2)}}{\varepsilon^{\nicefrac{3}{2}}}$ 
\end{tabular}  &   \begin{tabular}{c}
	ES,\\
	$\zeta$-Het
\end{tabular} & \ref{sec:sgd_es}  \\
\hline
 \begin{tabular}{c}
 {\tt Local-SGD}\\ Alg.~\ref{alg:local_sgd},  {\tiny\citep{khaled2020tighter}}$^{\clubsuit}$
\end{tabular} & $\myred{f_{\xi_i}(x_i^k)}, \myblue{0}, - $  & \begin{tabular}{c}
	$\frac{\left(L\tau+\nicefrac{\cL}{n}+\sqrt{(\tau-1)L\cL}\right)R_0^2}{\varepsilon}+\frac{\sigma_*^2R_0^2}{n\varepsilon^2}\quad$\\
	$\quad\quad+\frac{R_0^2\sqrt{L(\tau-1)(\sigma_*^2 + (\tau-1)\zeta_*^2)}}{\varepsilon^{\nicefrac{3}{2}}}$ 
\end{tabular}  &   \begin{tabular}{c}
	ES,\\
	Het
\end{tabular} & \ref{sec:sgd_es}  \\
\hline
\rowcolor{bgcolor} \begin{tabular}{c}
	{\tt Local-SVRG}\\ Alg.~\ref{alg:local_svrg},  {\tiny(NEW)}
\end{tabular} & \begin{tabular}{c}
	\myred{$\nabla f_{i,j_i}(x^k_i) -  \nabla f_{i,j_i}(y_i^k)$}\\\myred{$ +  \nabla f_{i}(y_i^k)$},\\ \myblue{$0$}, $- $ 
\end{tabular}  & \begin{tabular}{c}
	 $\frac{\left(L + \max L_{ij}\sqrt{\nicefrac{m}{n}} + \sqrt{(\tau-1)L\max L_{ij}}\right)R_0^2}{\varepsilon}$\\
	 $\frac{\sqrt[3]{(\tau-1)mL\max L_{ij}}R_0^2}{\varepsilon}+ \frac{L\zeta^2(\tau-1)R_0^2}{\mu\varepsilon^2}$\\
	 $+\frac{R_0^2\sqrt{L(\tau-1)\zeta_*^2}}{\varepsilon^{\nicefrac{3}{2}}}$
\end{tabular} & \begin{tabular}{c}
	simple,\\
	$\zeta$-Het
\end{tabular} & \ref{sec:llsvrg}   \\
\hline
\rowcolor{bgcolor} \begin{tabular}{c}
	{\tt Local-SVRG}\\ Alg.~\ref{alg:local_svrg},  {\tiny(NEW)}
\end{tabular} & \begin{tabular}{c}
	\myred{$\nabla f_{i,j_i}(x^k_i) -  \nabla f_{i,j_i}(y_i^k)$}\\\myred{$ +  \nabla f_{i}(y_i^k)$},\\ \myblue{$0$}, $- $ 
\end{tabular}  & \begin{tabular}{c}
	 $\frac{\left(L\tau + \max L_{ij}\sqrt{\nicefrac{m}{n}} + \sqrt{(\tau-1)L\max L_{ij}}\right)R_0^2}{\varepsilon}$\\
	 $\frac{\sqrt[3]{(\tau-1)mL\max L_{ij}}R_0^2}{\varepsilon}+\frac{R_0^2\sqrt{L(\tau-1)^2\zeta_*^2}}{\varepsilon^{\nicefrac{3}{2}}}$
\end{tabular} & \begin{tabular}{c}
	simple,\\
	Het
\end{tabular} & \ref{sec:llsvrg}   \\
\hline
\rowcolor{bgcolor} \begin{tabular}{c}
	{\tt S*-Local-SGD}\\ Alg.~\ref{alg:local_sgd_star},  {\tiny(NEW)}
\end{tabular} & $\myred{f_{\xi_i}(x_i^k)}, \myblue{\nabla f_i(x^*)}, - $ & $\frac{\tau LR_0^2}{\varepsilon}+\frac{\sigma^2R_0^2}{n\varepsilon^2}+\frac{R_0^2\sqrt{L(\tau-1)\sigma^2}}{\varepsilon^{\nicefrac{3}{2}}}$  &  \begin{tabular}{c}
	UBV,\\
	Het
\end{tabular} & \ref{sec:sgd_star_bounded_var}  \\
\hline
 \begin{tabular}{c}
	{\tt SS-Local-SGD}\\ Alg.~\ref{alg:l_local_svrg}, {\tiny\citep{karimireddy2020scaffold}} 
\end{tabular}  & \begin{tabular}{c}
$\myred{f_{\xi_i}(x_i^k)}, \myblue{h_i^k  - \frac1n \sum_{i=1}^n h_i^k},$\\
$ \nabla f_{\tilde{\xi}_i^k} (y_i^k)$ 
\end{tabular}  & $\frac{LR_0^2}{p\varepsilon}+\frac{\sigma^2R_0^2}{n\varepsilon^2}+\frac{R_0^2\sqrt{L(1-p)\sigma^2}}{p^{\nicefrac{1}{2}}\varepsilon^{\nicefrac{3}{2}}}$ &  \begin{tabular}{c}
	UBV,\\
	Het
\end{tabular} & \ref{sec:loopless_local_svrg}  \\
\hline
\rowcolor{bgcolor} \begin{tabular}{c}
	{\tt SS-Local-SGD}\\ Alg.~\ref{alg:l_local_svrg}, {\tiny(NEW)}
\end{tabular}  & \begin{tabular}{c}
$\myred{f_{\xi_i}(x_i^k)}, \myblue{h_i^k  - \frac1n \sum_{i=1}^n h_i^k},$\\
$ \nabla f_{\tilde{\xi}_i^k} (y_i^k)$ 
\end{tabular}  &  \begin{tabular}{c}
$\frac{\left(L + \nicefrac{p\cL}{n} + \sqrt{p(1-p)L\cL}\right)R_0^2}{p\varepsilon}$\\$+ \frac{\sqrt[3]{(1-p)L(L+p\cL)R_0^4\Delta_0}}{p\varepsilon}$\\ $+\frac{\sqrt[3]{(1-p)L\sigma_*^2R_0^4}}{p^{\nicefrac{2}{3}}\varepsilon}+\frac{\sigma_*^2R_0^2}{n\varepsilon^2}$\\$+\frac{R_0^2\sqrt{L(1-p)\sigma_*^2}}{p^{\nicefrac{1}{2}}\varepsilon^{\nicefrac{3}{2}}}$
\end{tabular} &  \begin{tabular}{c}
	ES,\\
	Het
\end{tabular} & \ref{sec:loopless_local_svrg_es}  \\
\hline
\rowcolor{bgcolor2} \begin{tabular}{c}
	{\tt S*-Local-SGD*}\\ Alg.~\ref{alg:local_sgd_star_star},  {\tiny(NEW)}
\end{tabular} & \begin{tabular}{c}
 	\myred{$\nabla f_{i,j_i}(x^k_i) -  \nabla f_{i,j_i}(x^*)$}\\
 	\myred{$ +  \nabla f_{i}(x^*)$}, $\myblue{\nabla f_i(x^*)}, - $\\ 	
 \end{tabular} &  $\frac{\left(L\tau + \nicefrac{\max L_{ij}}{n} + \sqrt{(\tau-1)L\max L_{ij}}\right)R_0^2}{\varepsilon}$  &  \begin{tabular}{c}
	simple,\\
	Het
\end{tabular} & \ref{sec:S*-Local-SGD*}  \\
\hline
\rowcolor{bgcolor2} \begin{tabular}{c}
	{\tt S-Local-SVRG}\\ Alg.~\ref{alg:l_local_svrg_fs},  {\tiny(NEW)}
\end{tabular}  &\begin{tabular}{c}
	$\myred{ \nabla f_{i,j_i}(x^k_i) -  \nabla f_{i,j_i}(y^k)}$\\
	\myred{$ +  \nabla f_{i}(y^k) $},\\
	$\myblue{h_i^k  - \frac1n \sum_{i=1}^n h_i^k}, \nabla f_{i} (y^k)$ 
\end{tabular} & \begin{tabular}{c}
$\frac{\left(L + pL\sqrt{\nicefrac{m}{n}} + \sqrt{(1-p)L\max L_{ij}}\right)R_0^2}{p\varepsilon}$\\ $+ \frac{R_0^2\sqrt[3]{L\max L_{ij}^2}}{p^{\nicefrac{2}{3}}\varepsilon}$ 
\end{tabular} & \begin{tabular}{c}
	simple,\\
	Het
\end{tabular} &   \ref{sec:loopless_local_svrg_fs} \\
\hline
\end{tabular}
\end{center}
\vskip -0.2cm
\end{table}

\subsection{{\tt Local-SGD}} \label{sec:local_sgd_app}
We start with the analysis of {\tt Local-SGD} (see Algorithm~\ref{alg:local_sgd}) under different assumptions of stochastic gradients and data similarity.
\begin{algorithm}[h]
   \caption{{\tt Local-SGD}}\label{alg:local_sgd}
\begin{algorithmic}[1]
   \Require learning rate $\gamma>0$, initial vector $x^0 \in \R^d$, communication period $\tau \ge 1$
	\For{$k=0,1,\dotsc$}
       	\For{$i=1,\dotsc,n$ in parallel}
            \State Sample $g^{k}_i = \nabla f_{\xi_i^k}(x_i^k)$ independently from other nodes
            \If{$k+1 \mod \tau = 0$}
            \State $x_i^{k+1} = x^{k+1} = \frac{1}{n}\sum\limits_{i=1}^n\left(x_i^k - \gamma g_i^k\right)$ \Comment{averaging}
            \Else
            \State $x_i^{k+1} = x_i^k - \gamma g_i^k$ \Comment{local update}
            \EndIf
        \EndFor
   \EndFor
\end{algorithmic}
\end{algorithm}

\subsubsection{Uniformly Bounded Variance}\label{sec:sgd_bounded_var}

In this section we assume that $f_i$ has a form of expectation (see \eqref{eq:f_i_expectation}) and stochastic gradients $\nabla f_{\xi_i}(x)$ satisfy
\begin{equation}
	\EE_{\xi_i}\left[\|\nabla f_{\xi_i}(x) - \nabla f_i(x)\|^2\right] \le D_{1,i},\quad \forall\; x\in\R^d,\;\forall\; i\in [n].\label{eq:bounded_variance}
\end{equation}
We also introduce the average variance $\sigma^2$ and the parameter of heterogeneity at the solution $\zeta_*^2$ in the following way:
\begin{equation*}
	\sigma^2 = \frac{1}{n}\sum\limits_{i=1}^nD_{1,i},\quad \zeta_*^2 = \frac{1}{n}\sum\limits_{i=1}^n\|\nabla f_i(x^*)\|^2.
\end{equation*}

\begin{lemma}\label{lem:local_sgd_interesting_labels}
Assume that functions $f_i$ are convex and $L$-smooth for all $i\in[n]$. Then
\begin{equation}\label{eq:dnaossniadd}
\frac{1}{n}\sum\limits_{i=1}^n\|\nabla f_i(x_i^k)\|^2\leq
		6L\left(f(x^k) - f(x^*)\right) + 3L^2 V_k + 3\zeta_*^2
\end{equation}
and
\begin{equation}\label{eq:vdgasvgda}
	 \left\|\frac{1}{n}\sum\limits_{i=1}^n\nabla f_i(x_i^k)\right\|^2
		\leq  4L\left(f(x^k)-f(x^*)\right) + 2L^2 V_k .
\end{equation}
\end{lemma}
\begin{proof}
First, to show~\eqref{eq:dnaossniadd} we shall have
	\begin{eqnarray*}
 \frac{1}{n}\sum\limits_{i=1}^n\|\nabla f_i(x_i^k)\|^2
		&\overset{\eqref{eq:a_b_norm_squared}}{\le}&
		\frac{3}{n}\sum\limits_{i=1}^n\|\nabla f_i(x_i^k) - \nabla f_i(x^k)\|^2 + \frac{3}{n}\sum\limits_{i=1}^n\|\nabla f_i(x^k) - \nabla f_i(x^*)\|^2 \\
		&& \qquad
		+ \frac{3}{n}\sum\limits_{i=1}^n\|\nabla f_i(x^*)\|^2
		\\
		&\overset{\eqref{eq:L_smoothness},\eqref{eq:L_smoothness_cor}}{\le}& \frac{3L^2}{n}\sum\limits_{i=1}^n\|x_i^k - x^k\|^2+ \frac{6L}{n}\sum\limits_{i=1}^nD_{f_i}(x^k,x^*) + 3\zeta_*^2
		\\
		&=&
		6L\left(f(x^k) - f(x^*)\right) + 3L^2 V_k + 3\zeta_*^2.
	\end{eqnarray*}
Next, to establish~\eqref{eq:vdgasvgda}, we have
	\begin{eqnarray*}
	 \left\|\frac{1}{n}\sum\limits_{i=1}^n\nabla f_i(x_i^k)\right\|^2
		&=&  \left\|\frac{1}{n}\sum\limits_{i=1}^n\left(\nabla f_i(x_i^k)-\nabla f_i(x^*)\right)\right\|^2
		\\
		&\overset{\eqref{eq:a_b_norm_squared}}{\le}&  \frac{2}{n}\sum\limits_{i=1}^n\|\nabla f_i(x_i^k) - \nabla f(x^k)\|^2 + \frac{2}{n}\sum\limits_{i=1}^n\|\nabla f_i(x^k) - \nabla f(x^*)\|^2
		\\
		&\overset{\eqref{eq:L_smoothness},\eqref{eq:L_smoothness_cor}}{\le}& \frac{2L^2}{n}\sum\limits_{i=1}^n\|x_i^k - x^k\|^2 + \frac{4L}{n}\sum\limits_{i=1}^nD_{f_i}(x^k,x^*)\\
		&=& 4L\left(f(x^k)-f(x^*)\right) + 2L^2 V_k .
	\end{eqnarray*}
\end{proof}

\begin{lemma}\label{lem:local_sgd_second_moment}
	Let $f_i$ be convex and $L$-smooth for all $i\in[n]$. Then for all $k\ge 0$
	\begin{eqnarray}
		\frac{1}{n}\sum\limits_{i=1}^n \EE\left[\|g_i^k\|^2\mid x^k\right] &\le& 6L\left(f(x^k) - f(x^*)\right) + 3L^2 V_k + \sigma^2 + 3\zeta_*^2,\label{eq:second_moment_local_sgd}\\
		\frac{1}{n}\sum\limits_{i=1}^n \EE\left[\|g_i^k-\bar{g}_i^k\|^2\mid x^k\right] &\le& \sigma^2,\label{eq:variance_local_sgd}\\
		\EE\left[\left\|\frac{1}{n}\sum\limits_{i=1}^ng_i^k\right\|^2\mid x^k\right] &\le& 4L\left(f(x^k) - f(x^*)\right) + 2L^2 V_k + \frac{\sigma^2}{n},\label{eq:second_moment_local_sgd_2}
	\end{eqnarray}
	where $\EE[\cdot\mid x^k]\eqdef \EE[\cdot\mid x_1^k,\ldots,x_n^k]$.
\end{lemma}
\begin{proof}
	First of all, we notice that $\bar{g}_i^k = \EE\left[g_i^k\mid x^k\right] = \nabla f_i(x_i^k)$. Using this we get
	\begin{eqnarray*}
		\frac{1}{n}\sum\limits_{i=1}^n \EE\left[\|g_i^k-\bar g_i^k\|^2\mid x_i^k\right] &=& \frac{1}{n}\sum\limits_{i=1}^n \EE_{\xi_i^k}\|\nabla f_{\xi_i^k}(x_i^k) - \nabla f_i(x_i^k)\|^2 \overset{\eqref{eq:bounded_variance}}{\le} \frac{1}{n}\sum\limits_{i=1}^n D_{1,i},\\
		\frac{1}{n}\sum\limits_{i=1}^n \EE\left[\|g_i^k\|^2\mid x_i^k\right] &\overset{\eqref{eq:variance_decomposition}}{=}& \frac{1}{n}\sum\limits_{i=1}^n \EE_{\xi_i^k}\|\nabla f_{\xi_i^k}(x_i^k) - \nabla f_i(x_i^k)\|^2 + \frac{1}{n}\sum\limits_{i=1}^n\|\nabla f_i(x_i^k)\|^2\\
		&\overset{\eqref{eq:bounded_variance},\eqref{eq:dnaossniadd}}{\le}&
		 6L\left(f(x^k) - f(x^*)\right) + 3L^2 V_k + \frac{1}{n}\sum\limits_{i=1}^n\left(D_{1,i} + 3\|\nabla f_i(x^*)\|^2\right).
	\end{eqnarray*}
	Finally, using independence of $g_1^k,g_2^k,\ldots,g_n^k$ we obtain
	\begin{eqnarray*}
		\EE\left[\left\|\frac{1}{n}\sum\limits_{i=1}^ng_i^k\right\|^2\mid x^k\right] &\overset{\eqref{eq:variance_decomposition}}{\le}& \EE\left[\left\|\frac{1}{n}\sum\limits_{i=1}^n\left(g_i^k - \nabla f_i(x_i^k)\right)\right\|^2\mid x^k\right] + \left\|\frac{1}{n}\sum\limits_{i=1}^n\nabla f_i(x_i^k)\right\|^2\\
		&=& \frac{1}{n^2}\sum\limits_{i=1}^n\EE\left[\|g_i^k - \nabla f_i(x_i^k)\|^2\mid x_i^k\right] +  \left\|\frac{1}{n}\sum\limits_{i=1}^n\nabla f_i(x_i^k)\right\|^2\\
		&\overset{\eqref{eq:bounded_variance},\eqref{eq:vdgasvgda}}{\le}&
		4L\left(f(x^k)-f(x^*)\right) + 2L^2 V_k + \frac{1}{n^2}\sum\limits_{i=1}^nD_{1,i}.
	\end{eqnarray*}
\end{proof}

\subsubsection*{Heterogeneous Data}
Applying Corollary~\ref{cor:const_loop} and Lemmas~\ref{lem:local_sgd_interesting_labels}~and~\ref{lem:local_sgd_second_moment} we get the following result.
\begin{theorem}\label{thm:local_sgd}
	Assume that $f_i(x)$ is $\mu$-strongly convex and $L$-smooth for every $i\in[n]$. Then {\tt Local-SGD} satisfies Assumption~\ref{ass:hetero_second_moment} with
	\begin{gather*}
		\tA = 3L,\quad \hA= 0,\quad \tB = \hB = 0,\quad \tF = 3L^2,\quad \hF = 0, \quad \tD_1 = 3\zeta_*^2,\quad \hD = \sigma^2,\\
		A' = 2L,\quad B' = 0,\quad F' = 2L^2, \quad D_1' = \frac{\sigma^2}{n},\quad \sigma_k^2 \equiv 0,\quad \rho = 1,\quad C = 0,\quad G = 0,\quad D_2 = 0,\\
		H = 0,\quad D_3 = 2e(\tau-1)\left(3(\tau-1)\zeta_*^2 + \sigma^2\right)
	\end{gather*}
	with $\gamma$ satisfying
	\begin{eqnarray*}
		\gamma &\le& \min\left\{\frac{1}{4L}, \frac{1}{4\sqrt{6e}(\tau-1)L}\right\}.
	\end{eqnarray*}
	and for all $K \ge 0$
	\begin{eqnarray}
		\EE\left[f(\overline{x}^K) - f(x^*)\right] &\le& \frac{2\|x^0-x^*\|^2}{\gamma W_K} + 2\gamma\left(\nicefrac{\sigma^2}{n} + 4eL(\tau-1)\gamma \left(\sigma^2 + 3(\tau-1)\zeta_*^2\right)\right). \notag
	\end{eqnarray}
	In particular, if $\mu > 0$ then
	\begin{eqnarray}
		\EE\left[f(\overline{x}^K) - f(x^*)\right] &\le& \left(1 - \gamma\mu\right)^K\frac{2\|x^0-x^*\|^2}{\gamma}\notag\\
		&&\quad + 2\gamma\left(\nicefrac{\sigma^2}{n} + 4eL(\tau-1)\gamma \left(\sigma^2 + 3(\tau-1)\zeta_*^2\right)\right) \label{eq:local_sgd_str_cvx}
	\end{eqnarray}
	and when $\mu = 0$ we have
	\begin{eqnarray}
		\EE\left[f(\overline{x}^K) - f(x^*)\right] &\le& \frac{2\|x^0-x^*\|^2}{\gamma K} + 2\gamma\left(\nicefrac{\sigma^2}{n} + 4eL(\tau-1)\gamma \left(\sigma^2 + 3(\tau-1)\zeta_*^2\right)\right). \label{eq:local_sgd_cvx}
	\end{eqnarray}
\end{theorem}

The theorem above together with Lemma~\ref{lem:lemma2_stich} implies the following result.
\begin{corollary}\label{cor:local_sgd_str_cvx}
	Let assumptions of Theorem~\ref{thm:local_sgd} hold with $\mu > 0$. Then for
	\begin{eqnarray*}
	    \gamma &=& \min\left\{\frac{1}{4L}, \frac{1}{4\sqrt{6e}(\tau-1)L},  \gamma_K\right\},\\
		\gamma_K &=& \frac{\ln\left(\max\left\{2,\min\left\{\nicefrac{\|x^0 - x^*\|^2n\mu^2K^2}{\sigma^2},\nicefrac{\|x^0 - x^*\|^2\mu^3K^3}{4eL(\tau-1)\left(\sigma^2 + 3(\tau-1)\zeta_*^2\right)}\right\}\right\}\right)}{\mu K}
	\end{eqnarray*}
	for all $K$ such that 
	\begin{eqnarray*}
		\text{either}&&\mu\gamma_K \le 1\\
		\text{or}&&\min\left\{\frac{1}{4L}, \frac{1}{4\sqrt{6e}(\tau-1)L}\right\} \le \gamma_K
	\end{eqnarray*}
	we have that $\EE\left[f(\overline{x}^K)-f(x^*)\right]$ equals
	\begin{equation}
		 \widetilde\cO\left(\tau L\|x^0 - x^*\|^2\exp\left(- \frac{\mu}{\tau L} K\right) + \frac{\sigma^2}{n\mu K} + \frac{L(\tau-1)\left(\sigma^2+(\tau-1)\zeta_*^2\right)}{\mu^2K^2}\right).\label{eq:local_sgd_str_cvx_1}
	\end{equation}
	That is, to achieve $\EE\left[f(\overline{x}^K)-f(x^*)\right] \le \varepsilon$ in this case {\tt Local-SGD} requires
	\begin{equation*}
		\widetilde\cO\left(\frac{\tau L}{\mu} + \frac{\sigma^2}{n\mu\varepsilon} + \sqrt{\frac{L(\tau-1)\left(\sigma^2+(\tau-1)\zeta_*^2\right)}{\mu^2\varepsilon}}\right)
	\end{equation*}
	iterations/oracle calls per node and $\tau$ times less communication rounds.
\end{corollary}
Now we consider some special cases. First of all, if $D_{1,i} = 0$ for all $i\in [n]$, i.e.\ $g_i^k = \nabla f_i(x_i^k)$ almost surely, then our result implies that for {\tt Local-SGD} it is enough to perform
\begin{equation*}
	\widetilde\cO\left(\frac{\tau L}{\mu} + \sqrt{\frac{L(\tau-1)^2\zeta_*^2}{\mu^2\varepsilon}}\right)
\end{equation*}
iterations in order to achieve $\EE\left[f(\overline{x}^K)-f(x^*)\right] \le \varepsilon$. It is clear that for this scenario the optimal choice for $\tau$ is $\tau = 1$ which recovers\footnote{We notice that for this particular case our analysis doesn't give extra logarithmical factors if we apply \eqref{eq:local_sgd_str_cvx} instead of \eqref{eq:local_sgd_str_cvx_1}.} the rate of Gradient Descent.

Secondly, if $\tau = 1$ then we recover the rate of parallel {\tt SGD}:
\begin{eqnarray*}
	\widetilde\cO\left(\frac{L}{\mu} + \frac{\sigma^2}{n\mu\varepsilon}\right)&& \text{communication rounds/oracle calls per node}
\end{eqnarray*}
in order to achieve $\EE\left[f(\overline{x}^K)-f(x^*)\right] \le \varepsilon$.

Finally, our result gives a negative answer to the following question: is {\tt Local-SGD} always worse then Parallel Minibatch {\tt SGD} ({\tt PMSGD}) for heterogeneous data? To achieve $\EE\left[f(\overline{x}^K)-f(x^*)\right] \le \varepsilon$ {\tt Local-SGD} requires
\begin{equation*}
	\widetilde\cO\left(\frac{\tau L}{\mu}+ \frac{\sigma^2}{n\mu\varepsilon}+ \sqrt{\frac{L(\tau-1)\left(\sigma^2+(\tau-1)\zeta_*^2\right)}{\mu^2\varepsilon}}\right)\quad \text{oracle calls per node.}
\end{equation*}
It means that if $\frac{\sigma^2}{n\sqrt{L(\tau-1)\left(\sigma^2+(\tau-1)\zeta_*^2\right)\varepsilon}} \ge 1$ for given $\tau > 1$ and $\varepsilon$ and $\sigma^2$ are such that the first term in the complexity bound is dominated by other terms, then the second term corresponding to the complexity of {\tt PMSGD} dominates the third term. Informally speaking, if the variance is large or $\varepsilon$ is small then {\tt Local-SGD} with $\tau > 1$ has the same complexity bounds as {\tt PMSGD}.

Combining Theorem~\ref{thm:local_sgd} and Lemma~\ref{lem:lemma_technical_cvx} we derive the following result for the convergence of {\tt Local-SGD} in the case when $\mu = 0$.
\begin{corollary}
	\label{cor:local_sgd_cvx}
	Let assumptions of Theorem~\ref{thm:local_sgd} hold with $\mu = 0$. Then for
	\begin{equation*}
		\gamma = \min\left\{\frac{1}{4L}, \frac{1}{4\sqrt{6e}(\tau-1)L}, \sqrt{\frac{nR_0^2}{\sigma^2 K}}, \sqrt[3]{\frac{R_0^2}{4eL(\tau-1)\left(\sigma^2+(\tau-1)\zeta_*^2\right)K}}\right\},
	\end{equation*}
	where $R_0 = \|x^0 - x^*\|$, we have that
	\begin{equation}
		\EE\left[f(\overline{x}^K)-f(x^*)\right] = \cO\left(\frac{\tau LR_0^2}{K} + \sqrt{\frac{R_0^2\sigma^2}{nK}} + \frac{\sqrt[3]{LR_0^4(\tau-1)\left(\sigma^2+(\tau-1)\zeta_*^2\right)}}{K^{\nicefrac{2}{3}}} \right).\label{eq:local_sgd_cvx_1}
	\end{equation}
	That is, to achieve $\EE\left[f(\overline{x}^K)-f(x^*)\right] \le \varepsilon$ in this case {\tt Local-SGD} requires
	\begin{equation*}
		\cO\left(\frac{\tau LR_0^2}{\varepsilon} + \frac{R_0^2\sigma^2}{n\varepsilon^2} + \frac{R_0^2\sqrt{L(\tau-1)\left(\sigma^2+(\tau-1)\zeta_*^2\right)}}{\varepsilon^{\nicefrac{3}{2}}}\right)
	\end{equation*}
	iterations/oracle calls per node and $\tau$ times less communication rounds.
\end{corollary}

\subsubsection*{Homogeneous Data}
In this case we modify the approach a little bit and apply the following result.
\begin{lemma}[Lemma~1 from \citep{khaled2020tighter}]
	Under the homogeneous data assumption for {\tt Local-SGD} we have
	\begin{equation}
		\EE\left[V_k\right] \le (\tau-1)\gamma^2\sigma^2 \label{eq:homocase_V_k_ubv}
	\end{equation}
	for all $k \ge 0$.
\end{lemma}
Using this we derive the following inequality for the weighted sum of $V_k$:
\begin{equation*}
	2L\sum\limits_{k=0}^Kw_k\EE[V_k] \le 2L(\tau-1)\gamma^2\sigma^2\sum\limits_{k=0}^Kw_k = 2L(\tau-1)\gamma^2\sigma^2 W_K.
\end{equation*}
Together with Lemmas~\ref{lem:local_sgd_interesting_labels}~and~\ref{lem:local_sgd_second_moment} and Theorem~\ref{thm:main_result} it gives the following result.
\begin{theorem}\label{thm:local_sgd_homo}
	Assume that $f(x)$ is $\mu$-strongly convex and $L$-smooth and $f_1 = \ldots = f_n = f$. Then {\tt Local-SGD} satisfies Assumption~\ref{ass:key_assumption} with
	\begin{gather*}
		A = 3L,\quad B = 0,\quad F = 3L^2, \quad D_1 = \sigma^2,\quad A' = 2L,\quad B' = 0,\quad F' = 2L^2, \quad D_1' = \frac{\sigma^2}{n},\\
		\sigma_k^2 \equiv 0,\quad \rho = 1,\quad C = 0,\quad G = 0,\quad D_2 = 0,\quad H = 0,\quad D_3 = (\tau-1)\sigma^2
	\end{gather*}
	with $\gamma$ satisfying
	\begin{eqnarray*}
		\gamma &\le& \frac{1}{4L}.
	\end{eqnarray*}
	and for all $K \ge 0$
	\begin{eqnarray}
		\EE\left[f(\overline{x}^K) - f(x^*)\right] &\le& \frac{2\|x^0-x^*\|^2}{\gamma W_K} + 2\gamma\left(\nicefrac{\sigma^2}{n} + 2L(\tau-1)\gamma \sigma^2\right). \notag
	\end{eqnarray}
	In particular, if $\mu > 0$ then
	\begin{eqnarray}
		\EE\left[f(\overline{x}^K) - f(x^*)\right] &\le& \left(1 - \gamma\mu\right)^K\frac{2\|x^0-x^*\|^2}{\gamma} + 2\gamma\left(\nicefrac{\sigma^2}{n} + 2L(\tau-1)\gamma \sigma^2\right) \label{eq:local_sgd_str_cvx_homo}
	\end{eqnarray}
	and when $\mu = 0$ we have
	\begin{eqnarray}
		\EE\left[f(\overline{x}^K) - f(x^*)\right] &\le& \frac{2\|x^0-x^*\|^2}{\gamma K} + 2\gamma\left(\nicefrac{\sigma^2}{n} + 2L(\tau-1)\gamma \sigma^2\right). \label{eq:local_sgd_cvx_homo}
	\end{eqnarray}
\end{theorem}

The theorem above together with Lemma~\ref{lem:lemma2_stich} implies the following result.
\begin{corollary}\label{cor:local_sgd_str_cvx_homo}
	Let assumptions of Theorem~\ref{thm:local_sgd_homo} hold with $\mu > 0$. Then for
	\begin{equation*}
		\gamma = \min\left\{\frac{1}{4L}, \frac{\ln\left(\max\left\{2,\min\left\{\nicefrac{\|x^0 - x^*\|^2n\mu^2K^2}{\sigma^2},\nicefrac{\|x^0 - x^*\|^2\mu^3K^3}{2L(\tau-1)\sigma^2}\right\}\right\}\right)}{\mu K}\right\}
	\end{equation*}
	for all $K$ such that 
	\begin{eqnarray*}
		\text{either} && \frac{\ln\left(\max\left\{2,\min\left\{\nicefrac{\|x^0 - x^*\|^2n\mu^2K^2}{\sigma^2},\nicefrac{\|x^0 - x^*\|^2\mu^3K^3}{2L(\tau-1)\sigma^2}\right\}\right\}\right)}{ K} \le 1\\
		\text{or} && \frac{1}{4L} \le \frac{\ln\left(\max\left\{2,\min\left\{\nicefrac{\|x^0 - x^*\|^2n\mu^2K^2}{\sigma^2},\nicefrac{\|x^0 - x^*\|^2\mu^3K^3}{2L(\tau-1)\sigma^2}\right\}\right\}\right)}{\mu K}
	\end{eqnarray*}
	we have that
	\begin{equation}
		\EE\left[f(\overline{x}^K)-f(x^*)\right] = \widetilde\cO\left(L\|x^0 - x^*\|^2\exp\left(- \frac{\mu}{L} K\right) + \frac{\sigma^2}{n\mu K} + \frac{L(\tau-1)\sigma^2}{\mu^2K^2}\right).\label{eq:local_sgd_str_cvx_1_homo}
	\end{equation}
	That is, to achieve $\EE\left[f(\overline{x}^K)-f(x^*)\right] \le \varepsilon$ in this case {\tt Local-SGD} requires
	\begin{equation*}
		\widetilde\cO\left(\frac{L}{\mu}\ln\left(\frac{L\|x^0 - x^*\|^2}{\varepsilon}\right) + \frac{\sigma^2}{n\mu\varepsilon} + \sqrt{\frac{L(\tau-1)\sigma^2}{\mu^2\varepsilon}}\right)
	\end{equation*}
	iterations/oracle calls per node and $\tau$ times less communication rounds.
\end{corollary}
It means that if $\frac{\sigma^2}{n^2L\varepsilon} \ge 1$, $\tau \le 1 + \frac{\sigma^2}{n^2L\varepsilon}$ and $\varepsilon$ and $\sigma^2$ are such that the first term in the complexity bound is dominated by other terms, then the second term corresponding to the complexity of {\tt PMSGD} dominates the third term. Informally speaking, if the variance is large or $\varepsilon$ is small then {\tt Local-SGD} with $\tau > 1$ has the same complexity bounds as {\tt PMSGD}.

Combining Theorem~\ref{thm:local_sgd_homo} and Lemma~\ref{lem:lemma_technical_cvx} we derive the following result for the convergence of {\tt Local-SGD} in the case when $\mu = 0$.
\begin{corollary}
	\label{cor:local_sgd_cvx_homo}
	Let assumptions of Theorem~\ref{thm:local_sgd_homo} hold with $\mu = 0$. Then for
	\begin{equation*}
		\gamma = \min\left\{\frac{1}{4L}, \sqrt{\frac{nR_0^2}{\sigma^2 K}}, \sqrt[3]{\frac{R_0^2}{2L(\tau-1)\sigma^2 K}}\right\},
	\end{equation*}
	where $R_0 = \|x^0 - x^*\|$, we have that
	\begin{equation}
		\EE\left[f(\overline{x}^K)-f(x^*)\right] = \cO\left(\frac{LR_0^2}{K} + \sqrt{\frac{R_0^2\sigma^2}{nK}} + \frac{\sqrt[3]{LR_0^4(\tau-1)\sigma^2}}{K^{\nicefrac{2}{3}}} \right).\label{eq:local_sgd_cvx_1_homo}
	\end{equation}
	That is, to achieve $\EE\left[f(\overline{x}^K)-f(x^*)\right] \le \varepsilon$ in this case {\tt Local-SGD} requires
	\begin{equation*}
		\cO\left(\frac{LR_0^2}{\varepsilon} + \frac{R_0^2\sigma^2}{n\varepsilon^2} + \frac{R_0^2\sqrt{L(\tau-1)\sigma^2}}{\varepsilon^{\nicefrac{3}{2}}}\right)
	\end{equation*}
	iterations/oracle calls per node and $\tau$ times less communication rounds.
\end{corollary}

\subsubsection*{$\zeta$-Heterogeneous Data}
In this setup we also use an external result to bound $\EE[V_k]$.
\begin{lemma}[Lemma~8 from \citep{woodworth2020minibatch}]
	If $f_1,f_2,\ldots,f_n$ are $\zeta$-heterogeneous then for {\tt Local-SGD} we have
	\begin{equation}
		\EE\left[V_k\right] \le 3\tau\gamma^2\sigma^2 + 6\tau^2\gamma^2\zeta^2 \label{eq:zeta_heterocase_V_k_ubv}
	\end{equation}
	for all $k \ge 0$.
\end{lemma}
Using this we derive the following inequality for the weighted sum of $V_k$:
\begin{equation*}
	2L\sum\limits_{k=0}^Kw_k\EE[V_k] \le 6\tau L\gamma^2\left(\sigma^2 + 2\tau\zeta^2\right)\sum\limits_{k=0}^Kw_k = 6\tau L\gamma^2\left(\sigma^2 + 2\tau\zeta^2\right) W_K.
\end{equation*}
Together with Lemmas~\ref{lem:local_sgd_interesting_labels}~and~\ref{lem:local_sgd_second_moment} and Theorem~\ref{thm:main_result} it gives the following result.
\begin{theorem}\label{thm:local_sgd_zeta_hetero}
	Assume that $f_1,\ldots,f_n$ are $\zeta$-heterogeneous, $\mu$-strongly convex and $L$-smooth functions. Then {\tt Local-SGD} satisfies Assumption~\ref{ass:key_assumption} with
	\begin{gather*}
		A = 3L,\quad B = 0,\quad F = 3L^2, \quad D_1 = \sigma^2 + 3\zeta_*^2,\quad A' = 2L,\quad B' = 0,\quad F' = 2L^2, \quad D_1' = \frac{\sigma^2}{n},\\
		\sigma_k^2 \equiv 0,\quad \rho = 1,\quad C = 0,\quad G = 0,\quad D_2 = 0,\quad H = 0,\quad D_3 = 3\tau\left(\sigma^2 + 2\tau\zeta^2\right)
	\end{gather*}
	with $\gamma$ satisfying
	\begin{eqnarray*}
		\gamma &\le& \frac{1}{4L}.
	\end{eqnarray*}
	and for all $K \ge 0$
	\begin{eqnarray}
		\EE\left[f(\overline{x}^K) - f(x^*)\right] &\le& \frac{2\|x^0-x^*\|^2}{\gamma W_K} + 2\gamma\left(\nicefrac{\sigma^2}{n} + 6L\tau\gamma \left(\sigma^2 + 2\tau\zeta^2\right)\right). \notag
	\end{eqnarray}
	In particular, if $\mu > 0$ then
	\begin{eqnarray}
		\EE\left[f(\overline{x}^K) - f(x^*)\right] &\le& \left(1 - \gamma\mu\right)^K\frac{2\|x^0-x^*\|^2}{\gamma} + 2\gamma\left(\nicefrac{\sigma^2}{n} + 6L\tau\gamma \left(\sigma^2 + 2\tau\zeta^2\right)\right) \label{eq:local_sgd_str_cvx_zeta_hetero}
	\end{eqnarray}
	and when $\mu = 0$ we have
	\begin{eqnarray}
		\EE\left[f(\overline{x}^K) - f(x^*)\right] &\le& \frac{2\|x^0-x^*\|^2}{\gamma K} + 2\gamma\left(\nicefrac{\sigma^2}{n} + 6L\tau\gamma \left(\sigma^2 + 2\tau\zeta^2\right)\right). \label{eq:local_sgd_cvx_zeta_hetero}
	\end{eqnarray}
\end{theorem}
The theorem above together with Lemma~\ref{lem:lemma2_stich} implies the following result.
\begin{corollary}\label{cor:local_sgd_str_cvx_zeta_hetero}
	Let assumptions of Theorem~\ref{thm:local_sgd_zeta_hetero} hold with $\mu > 0$. Then for
	\begin{equation*}
		\gamma = \min\left\{\frac{1}{4L}, \frac{\ln\left(\max\left\{2,\min\left\{\nicefrac{\|x^0 - x^*\|^2n\mu^2K^2}{\sigma^2},\nicefrac{\|x^0 - x^*\|^2\mu^3K^3}{6L\tau(\sigma^2+2\tau\zeta^2)}\right\}\right\}\right)}{\mu K}\right\}
	\end{equation*}
	for all $K$ such that 
	\begin{eqnarray*}
		\text{either} && \frac{\ln\left(\max\left\{2,\min\left\{\nicefrac{\|x^0 - x^*\|^2n\mu^2K^2}{\sigma^2},\nicefrac{\|x^0 - x^*\|^2\mu^3K^3}{6L\tau(\sigma^2+2\tau\zeta^2)}\right\}\right\}\right)}{K} \le 1\\
		\text{or} && \frac{1}{4L} \le \frac{\ln\left(\max\left\{2,\min\left\{\nicefrac{\|x^0 - x^*\|^2n\mu^2K^2}{\sigma^2},\nicefrac{\|x^0 - x^*\|^2\mu^3K^3}{6L\tau(\sigma^2+2\tau\zeta^2)}\right\}\right\}\right)}{\mu K}
	\end{eqnarray*}
	we have that
	\begin{equation}
		\EE\left[f(\overline{x}^K)-f(x^*)\right] = \widetilde\cO\left(L\|x^0 - x^*\|^2\exp\left(- \frac{\mu}{L} K\right) + \frac{\sigma^2}{n\mu K} + \frac{L\tau(\sigma^2+\tau\zeta^2)}{\mu^2K^2}\right).\label{eq:local_sgd_str_cvx_1_zeta_hetero}
	\end{equation}
	That is, to achieve $\EE\left[f(\overline{x}^K)-f(x^*)\right] \le \varepsilon$ in this case {\tt Local-SGD} requires
	\begin{equation*}
		\widetilde\cO\left(\frac{L}{\mu}\ln\left(\frac{L\|x^0 - x^*\|^2}{\varepsilon}\right) + \frac{\sigma^2}{n\mu\varepsilon} + \sqrt{\frac{L\tau (\sigma^2+\tau\zeta^2)}{\mu^2\varepsilon}}\right)
	\end{equation*}
	iterations/oracle calls per node and $\tau$ times less communication rounds.
\end{corollary}

Combining Theorem~\ref{thm:local_sgd_zeta_hetero} and Lemma~\ref{lem:lemma_technical_cvx} we derive the following result for the convergence of {\tt Local-SGD} in the case when $\mu = 0$.
\begin{corollary}
	\label{cor:local_sgd_cvx_zeta_hetero}
	Let assumptions of Theorem~\ref{thm:local_sgd_zeta_hetero} hold with $\mu = 0$. Then for
	\begin{equation*}
		\gamma = \min\left\{\frac{1}{4L}, \sqrt{\frac{nR_0^2}{\sigma^2 K}}, \sqrt[3]{\frac{R_0^2}{6L\tau(\sigma^2+2\tau\zeta^2) K}}\right\},
	\end{equation*}
	where $R_0 = \|x^0 - x^*\|$, we have that
	\begin{equation}
		\EE\left[f(\overline{x}^K)-f(x^*)\right] = \cO\left(\frac{LR_0^2}{K} + \sqrt{\frac{R_0^2\sigma^2}{nK}} + \frac{\sqrt[3]{LR_0^4\tau(\sigma^2+\tau\zeta^2)}}{K^{\nicefrac{2}{3}}} \right).\label{eq:local_sgd_cvx_1_zeta_hetero}
	\end{equation}
	That is, to achieve $\EE\left[f(\overline{x}^K)-f(x^*)\right] \le \varepsilon$ in this case {\tt Local-SGD} requires
	\begin{equation*}
		\cO\left(\frac{LR_0^2}{\varepsilon} + \frac{R_0^2\sigma^2}{n\varepsilon^2} + \frac{R_0^2\sqrt{L\tau(\sigma^2+\tau\zeta^2)}}{\varepsilon^{\nicefrac{3}{2}}}\right)
	\end{equation*}
	iterations/oracle calls per node and $\tau$ times less communication rounds.
\end{corollary}

\subsubsection{Expected Smoothness and Arbitrary Sampling}\label{sec:sgd_es}
In this section we continue our consideration of {\tt Local-SGD} but now we make another assumption on stochastic gradients $\nabla f_{\xi_i}(x)$.
\begin{assumption}[Expected Smoothness]\label{ass:expected_smoothness}
	We assume that for all $i\in[n]$ stochastic gradients $\nabla f_{\xi_i}(x)$ are unbiased estimators of $\nabla f_i(x)$ and there exists such constant $\cL > 0$ that $\forall x,y\in\R^d$
	\begin{equation}
		\EE_{\xi_i\sim \cD_i}\left[\left\|\nabla f_{\xi_i}(x) - \nabla f_{\xi_i}(x^*)\right\|^2\right] \le 2\cL D_{f_i}(x,x^*) \label{eq:expected_smoothness_1}
	\end{equation}
	where $D_{f_i}(x,y) \eqdef f_i(x) - f_i(y) - \langle\nabla f_i(y), x-y\rangle$.
\end{assumption}

In particular, let us consider the following special case. Assume that $f_i(x)$ has a form of finite sum (see \eqref{eq:f_i_sum}) and consider the following stochastic reformulation:
\begin{equation}
	f_i(x) = \EE_{\xi_i}\left[f_{\xi_i}(x)\right],\quad f_{\xi_i}(x) = \frac{1}{m}\sum\limits_{j=1}^m \xi_{i,j}f_{i,j}(x), \label{eq:stochastic_reformulation}
\end{equation}
where $\EE[\xi_{i,j}] = 1$ and $\EE[\xi_{i,j}^2] < \infty$. In this case, $\EE_{\xi_i}[\nabla f_{\xi_i}] = \nabla f_i(x)$. If each $f_{i,j}(x)$ is $L_{i,j}$-smooth then there exists such $\cL \le \max_{j\in[m]}L_{i,j}$ that Assumption~\ref{ass:expected_smoothness} holds. Clearly, $\cL$ depends on the sampling strategy and in some cases one can make $\cL$ much smaller than $\max_{j\in[m]}L_{i,j}$ via good choice of this strategy. Our analysis works for an arbitrary sampling strategy that satisfies Assumption~\ref{ass:expected_smoothness}.

\begin{lemma}\label{lem:local_sgd_es_second_moment}
	Let $f_i$ be convex and $L$-smooth for all $i\in[n]$. Then for all $k\ge 0$
	\begin{eqnarray}
		\frac{1}{n}\sum\limits_{i=1}^n \EE\left[\|g_i^k\|^2\mid x^k\right] &\le& 8\cL\left(f(x^k) - f(x^*)\right) + 4\cL L V_k + 2\sigma_*^2 + 2\zeta_*^2,\label{eq:second_moment_local_sgd_es}\\
		\frac{1}{n}\sum\limits_{i=1}^n \EE\left[\|g_i^k-\bar{g}_i^k\|^2\mid x^k\right] &\le& 8\cL\left(f(x^k) - f(x^*)\right) + 4\cL L V_k + 2\sigma_*^2,\label{eq:var_local_sgd_es}\\
		\EE\left[\left\|\frac{1}{n}\sum\limits_{i=1}^ng_i^k\right\|^2\mid x^k\right] &\le& 4\left(\nicefrac{2\cL}{n} + L\right)(f(x^k) - f(x^*)) + 2L\left(\nicefrac{2\cL}{n} + L\right) V_k\notag\\
		&&\quad + \frac{2\sigma_*^2}{n},\label{eq:second_moment_local_sgd_es_2}
	\end{eqnarray}
	where $\sigma_*^2 = \frac{1}{n}\sum_{i=1}^n\EE\|\nabla f_{\xi_i}(x^*) - \nabla f_i(x^*)\|^2$, $\zeta_*^2 = \frac{1}{n}\sum_{i=1}^n\|\nabla f_i(x^*)\|^2$ and $\EE[\cdot\mid x^k]\eqdef \EE[\cdot\mid x_1^k,\ldots,x_n^k]$.
\end{lemma}
\begin{proof}
	First of all, we notice that $\bar{g}_i^k = \EE\left[g_i^k\mid x^k\right] = \nabla f_i(x_i^k)$. Using this we get
	\begin{eqnarray*}
		\frac{1}{n}\sum\limits_{i=1}^n \EE\left[\|g_i^k\|^2\mid x^k\right] &\overset{\eqref{eq:a_b_norm_squared}}{\le}& \frac{2}{n}\sum\limits_{i=1}^n \EE_{\xi_i^k}\|\nabla f_{\xi_i^k}(x_i^k) - \nabla f_{\xi_i^k}(x^*)\|^2 + \frac{2}{n}\sum\limits_{i=1}^n \EE_{\xi_i^k}\|\nabla f_{\xi_i^k}(x^*)\|^2\\
		&\overset{\eqref{eq:expected_smoothness_1},\eqref{eq:variance_decomposition}}{\le}&
		\frac{4\cL}{n}\sum\limits_{i=1}^nD_{f_i}(x_i^k,x^*) + \frac{2}{n}\sum\limits_{i=1}^n \EE_{\xi_i}\left[\|\nabla f_{\xi_i}(x^*)-\nabla f_i(x^*)\|^2\right]\\
		&&\quad + \frac{2}{n}\sum\limits_{i=1}^n \|\nabla f_i(x^*)\|^2\\
		&\overset{\eqref{eq:poiouhnkj}}{\le}& 8\cL\left(f(x^k) - f(x^*)\right) + 4\cL L V_k + 2\sigma_*^2 + 2\zeta_*^2
	\end{eqnarray*}
	and
	\begin{eqnarray}
		\frac{1}{n}\sum\limits_{i=1}^n \EE\left[\|g_i^k-\bar{g}_i^k\|^2\mid x^k\right] &=& \frac{1}{n}\sum\limits_{i=1}^n\EE_{\xi_i^k}\|\nabla f_{\xi_i^k}(x_i^k)-\nabla f_i(x_i^k)\|^2\notag\\
		&\overset{\eqref{eq:variance_decomposition}}{\le}& \frac{1}{n}\sum\limits_{i=1}^n\EE_{\xi_i^k}\|\nabla f_{\xi_i^k}(x_i^k)-\nabla f_i(x^*)\|^2 \notag \\
		&\overset{\eqref{eq:a_b_norm_squared}}{\le}& \frac{2}{n}\sum\limits_{i=1}^n\EE_{\xi_i^k}\|\nabla f_{\xi_i^k}(x_i^k)-\nabla f_{\xi_i^k}(x^*)\|^2\notag\\
		&&\quad + \frac{2}{n}\sum\limits_{i=1}^n\EE_{\xi_i^k}\|\nabla f_{\xi_i^k}(x^*)-\nabla f_{i}(x^*)\|^2\notag\\
		&\overset{\eqref{eq:expected_smoothness_1}}{\le}& \frac{4\cL}{n}\sum\limits_{i=1}^nD_{f_i}(x_i^k,x^*) + 2\sigma_*^2\notag\\
		&\overset{\eqref{eq:poiouhnkj}}{\le}& 8\cL\left(f(x^k) - f(x^*)\right) + 4\cL L V_k + 2\sigma_*^2. \label{eq:hbdsfbhdbfvfdh}
	\end{eqnarray}		
	Finally, using independence of $\xi_1^k,\xi_2^k,\ldots,\xi_n^k$ we obtain
	\begin{eqnarray*}
		\EE\left[\left\|\frac{1}{n}\sum\limits_{i=1}^ng_i^k\right\|^2\mid x^k\right] &\overset{\eqref{eq:variance_decomposition}}{=}& \EE_{\xi_i^k}\left[\left\|\frac{1}{n}\sum\limits_{i=1}^n(\nabla f_{\xi_i^k}(x_i^k) - \nabla f_{i}(x_i^k))\right\|^2\right] + \left\|\frac{1}{n}\sum\limits_{i=1}^n\nabla f_{i}(x_i^k)\right\|^2\\
		&=& \frac{1}{n^2}\sum\limits_{i=1}^n\EE_{\xi_i^k}\left[\|\nabla f_{\xi_i^k}(x_i^k) - \nabla f_{i}(x_i^k)\|^2\right] + \left\|\frac{1}{n}\sum\limits_{i=1}^n\nabla f_{i}(x_i^k)\right\|^2\\
		&\overset{\eqref{eq:hbdsfbhdbfvfdh},\eqref{eq:vdgasvgda}}{\le}& 4\left(\nicefrac{2\cL}{n} + L\right)(f(x^k) - f(x^*)) + 2L\left(\nicefrac{2\cL}{n} + L\right) V_k + \frac{2\sigma_*^2}{n}.
	\end{eqnarray*}
\end{proof}

\subsubsection*{Heterogeneous Data}
Applying Corollary~\ref{cor:const_loop} and Lemmas~\ref{lem:local_sgd_interesting_labels}~and~\ref{lem:local_sgd_es_second_moment} we get the following result.
\begin{theorem}\label{thm:local_sgd_es}
	Assume that $f_i(x)$ is $\mu$-strongly convex and $L$-smooth for $i\in[n]$. Let Assumption~\ref{ass:expected_smoothness} holds. Then {\tt Local-SGD} satisfies Assumption~\ref{ass:hetero_second_moment} with
	\begin{gather*}
		\tA = 3L,\quad \hA = 4\cL,\quad \tB = \hB = 0,\quad \tF = 3L^2,\quad \hF = 4\cL L, \quad \tD_1 = 3\zeta_*^2,\quad \hD_1 = 2\sigma_*^2\\
		A' = \frac{4\cL}{n} + 2L,\quad B' = 0,\quad F' = \frac{4\cL L}{n} + 2L^2, \quad D_1' = \frac{2\sigma_*^2}{n},\\
		\sigma_k^2 \equiv 0,\quad \rho = 1,\quad C = 0,\quad G = 0,\quad D_2 = 0,\\
		H = 0,\quad D_3 = 2e(\tau-1)\left(2\sigma_*^2+3(\tau-1)\zeta_*^2\right)
	\end{gather*}
	with $\gamma$ satisfying
	\begin{eqnarray*}
		\gamma &\le& \min\left\{\frac{1}{\nicefrac{8\cL}{n} + 4L}, \frac{1}{4\sqrt{2eL(\tau-1)\left(3L(\tau-1)+4\cL \right)}}\right\}.
	\end{eqnarray*}
	and for all $K \ge 0$
	\begin{eqnarray}
		\EE\left[f(\overline{x}^K) - f(x^*)\right] &\le& \frac{2\|x^0-x^*\|^2}{\gamma W_K} + 2\gamma\left(\nicefrac{2\sigma_*^2}{n} + 4eL(\tau-1)\gamma\left(2\sigma_*^2+3(\tau-1)\zeta_*^2\right)\right). \notag
	\end{eqnarray}
	In particular, if $\mu > 0$ then
	\begin{eqnarray}
		\EE\left[f(\overline{x}^K) - f(x^*)\right] &\le& \left(1 - \gamma\mu\right)^K\frac{2\|x^0-x^*\|^2}{\gamma}\notag\\
		&&\quad + 2\gamma\left(\nicefrac{2\sigma_*^2}{n} + 4eL(\tau-1)\gamma\left(2\sigma_*^2+3(\tau-1)\zeta_*^2\right)\right) \label{eq:local_sgd_es_str_cvx}
	\end{eqnarray}
	and when $\mu = 0$ we have
	\begin{eqnarray}
		\EE\left[f(\overline{x}^K) - f(x^*)\right] &\le& \frac{2\|x^0-x^*\|^2}{\gamma K}\notag \\
		&&\quad + 2\gamma\left(\nicefrac{2\sigma_*^2}{n} + 4eL(\tau-1)\gamma\left(2\sigma_*^2+3(\tau-1)\zeta_*^2\right)\right). \label{eq:local_sgd_es_cvx}
	\end{eqnarray}
\end{theorem}

The theorem above together with Lemma~\ref{lem:lemma2_stich} implies the following result.
\begin{corollary}\label{cor:local_sgd_es_str_cvx}
	Let assumptions of Theorem~\ref{thm:local_sgd_es} hold with $\mu > 0$. Then for 
	\begin{eqnarray*}
		\gamma_0 &=& \min\left\{\frac{1}{\nicefrac{8\cL}{n} + 4L}, \frac{1}{4\sqrt{2eL(\tau-1)\left(3L(\tau-1)+4\cL \right)}}\right\},\\
		\gamma &=& \min\left\{\gamma_0,\frac{\ln\left(\max\left\{2, \min\left\{\nicefrac{n\|x^0 - x^*\|^2\mu^2K^2}{2\sigma_*^2}, \nicefrac{\|x^0 - x^*\|^2\mu^3K^3}{4eL(\tau-1)\gamma\left(2\sigma_*^2+3(\tau-1)\zeta_*^2\right)}\right\}\right\}\right)}{\mu K}\right\}
	\end{eqnarray*}
	for all $K$ such that 
	\begin{eqnarray*}
		\text{either} && \frac{\ln\left(\max\left\{2, \min\left\{\nicefrac{n\|x^0 - x^*\|^2\mu^2K^2}{2\sigma_*^2}, \nicefrac{\|x^0 - x^*\|^2\mu^3K^3}{4eL(\tau-1)\gamma\left(2\sigma_*^2+3(\tau-1)\zeta_*^2\right)}\right\}\right\}\right)}{ K} \le 1\\
		\text{or} && \gamma_0 \le \frac{\ln\left(\max\left\{2, \min\left\{\nicefrac{n\|x^0 - x^*\|^2\mu^2K^2}{2\sigma_*^2}, \nicefrac{\|x^0 - x^*\|^2\mu^3K^3}{4eL(\tau-1)\gamma\left(2\sigma_*^2+3(\tau-1)\zeta_*^2\right)}\right\}\right\}\right)}{\mu K}
	\end{eqnarray*}
	we have that $\EE\left[f(\overline{x}^K)-f(x^*)\right]$ is of the order
	\begin{eqnarray}
	    \widetilde\cO\Bigg(\left(L\tau + \nicefrac{\cL}{n}+\sqrt{(\tau-1)\cL L}\right)R_0^2\exp\left(- \frac{\mu}{L\tau + \nicefrac{\cL}{n}+\sqrt{(\tau-1)\cL L}} K\right) &&\notag\\
		 &&\hspace{-5cm}+ \frac{\sigma_*^2}{n\mu K} + \frac{L(\tau-1)\left(\sigma_*^2+(\tau-1)\zeta_*^2\right)}{\mu^2 K^2}\Bigg),\notag
	\end{eqnarray}
	where $R_0 = \|x^0-x^*\|$. That is, to achieve $\EE\left[f(\overline{x}^K)-f(x^*)\right] \le \varepsilon$ in this case {\tt Local-SGD} requires
	\begin{equation*}
		\widetilde{\cO}\left(\frac{L\tau}{\mu}+\frac{\cL}{n\mu}+\frac{\sqrt{(\tau-1)\cL L}}{\mu} + \frac{\sigma_*^2}{n\mu\varepsilon} + \sqrt{\frac{L(\tau-1)\left(\sigma_*^2 + (\tau-1)\zeta_*^2\right)}{\mu^2\varepsilon}}\right)
	\end{equation*}
	iterations/oracle calls per node and $\tau$ times less communication rounds.	
\end{corollary}

Combining Theorem~\ref{thm:local_sgd_es} and Lemma~\ref{lem:lemma_technical_cvx} we derive the following result for the convergence of {\tt Local-SGD} in the case when $\mu = 0$.
\begin{corollary}
	\label{cor:local_sgd_es_cvx}
	Let assumptions of Theorem~\ref{thm:local_sgd_es} hold with $\mu = 0$. Then for
	\begin{eqnarray*}
		\gamma_0 &=& \min\left\{\frac{1}{\nicefrac{8\cL}{n} + 4L}, \frac{1}{4\sqrt{2eL(\tau-1)\left(3L(\tau-1)+4\cL \right)}}\right\},\\
		\gamma &=& \min\left\{\gamma_0, \sqrt{\frac{nR_0^2}{2\sigma_*^2 K}}, \sqrt[3]{\frac{R_0^2}{4eL(\tau-1)\left(2\sigma_*^2+3(\tau-1)\zeta_*^2\right) K}}\right\},
	\end{eqnarray*}
	where $R_0 = \|x^0 - x^*\|$, we have that $\EE\left[f(\overline{x}^K)-f(x^*)\right]$ equals
	\begin{equation}
		\cO\left(\frac{\left(L\tau + \nicefrac{\cL}{n}+\sqrt{(\tau-1)\cL L}\right)R_0^2}{K} + \sqrt{\frac{R_0^2\sigma_*^2}{nK}} + \frac{\sqrt[3]{LR_0^4(\tau-1)\left(\sigma_*^2+(\tau-1)\zeta_*^2\right)}}{K^{\nicefrac{2}{3}}} \right).\label{eq:local_sgd_es_cvx_1}
	\end{equation}
	That is, to achieve $\EE\left[f(\overline{x}^K)-f(x^*)\right] \le \varepsilon$ in this case {\tt Local-SGD} requires
	\begin{equation*}
		\cO\left(\frac{\left(L\tau + \nicefrac{\cL}{n}+\sqrt{(\tau-1)\cL L}\right)R_0^2}{\varepsilon} + \frac{R_0^2\sigma_*^2}{n\varepsilon^2} + \frac{R_0^2\sqrt{L(\tau-1)\left(\sigma_*^2+(\tau-1)\zeta_*^2\right)}}{\varepsilon^{\nicefrac{3}{2}}}\right)
	\end{equation*}
	iterations/oracle calls per node and $\tau$ times less communication rounds.
\end{corollary}

\subsubsection*{$\zeta$-Heterogeneous Data}
Applying Corollary~\ref{cor:const_loop_homo} and Lemma~\ref{lem:local_sgd_es_second_moment} we get the following result.
\begin{theorem}\label{thm:local_sgd_es_homo}
	Assume that $f_i(x)$ is $L$-smooth for $i\in[n]$ and $f_1, \ldots, f_n$ are $\zeta$-heterogeneous and $\mu$-strongly convex. Let Assumption~\ref{ass:expected_smoothness} holds. Then {\tt Local-SGD} satisfies Assumption~\ref{ass:key_assumption} with
	\begin{gather*}
		A = 4\cL,\quad B = 0,\quad F = 4\cL L, \quad D_1 = 2\sigma_*^2 + 2\zeta_*^2,\\
		A' = \frac{4\cL}{n} + 2L,\quad B' = 0,\quad F' = \frac{4\cL L}{n} + 2L^2, \quad D_1' = \frac{2\sigma_*^2}{n},\\
		\sigma_k^2 \equiv 0,\quad \rho = 1,\quad C = 0,\quad G = 0,\quad D_2 = 0,\\
		H = 0,\quad D_3 = 2(\tau-1)\left(2\sigma_*^2 + 2\zeta_*^2 + \frac{\zeta^2}{\gamma\mu}\right)
	\end{gather*}
	with $\gamma$ satisfying
	\begin{eqnarray*}
		\gamma &\le& \min\left\{\frac{1}{\nicefrac{8\cL}{n} + 4L}, \frac{1}{8\sqrt{2L\cL(\tau-1)}}\right\}.
	\end{eqnarray*}
	and for all $K \ge 0$
	\begin{eqnarray}
		\EE\left[f(\overline{x}^K) - f(x^*)\right] &\le& \frac{2\|x^0-x^*\|^2}{\gamma W_K} + 2\gamma\left(\frac{2\sigma_*^2}{n} + \frac{4L\zeta^2(\tau-1)}{\mu} + 8L(\tau-1)\gamma\left(\sigma_*^2+\zeta_*^2\right) \right). \notag
	\end{eqnarray}
	In particular, if $\mu > 0$ then
	\begin{eqnarray}
		\EE\left[f(\overline{x}^K) - f(x^*)\right] &\le& \left(1 - \gamma\mu\right)^K\frac{2\|x^0-x^*\|^2}{\gamma}\notag\\
		&&\quad + 2\gamma\left(\frac{2\sigma_*^2}{n} + \frac{4L\zeta^2(\tau-1)}{\mu} + 8L(\tau-1)\gamma\left(\sigma_*^2+\zeta_*^2\right) \right) \label{eq:local_sgd_es_str_cvx_homo}
	\end{eqnarray}
	and when $\mu = 0$ we have
	\begin{eqnarray}
		\EE\left[f(\overline{x}^K) - f(x^*)\right] &\le& \frac{2\|x^0-x^*\|^2}{\gamma K}\notag\\
		&&\quad + 2\gamma\left(\frac{2\sigma_*^2}{n} + \frac{4L\zeta^2(\tau-1)}{\mu} + 8L(\tau-1)\gamma\left(\sigma_*^2+\zeta_*^2\right) \right). \label{eq:local_sgd_es_cvx_homo}
	\end{eqnarray}
\end{theorem}

The theorem above together with Lemma~\ref{lem:lemma2_stich} implies the following result.
\begin{corollary}\label{cor:local_sgd_es_str_cvx_homo}
	Let assumptions of Theorem~\ref{thm:local_sgd_es_homo} hold with $\mu > 0$. Then for 
	\begin{eqnarray*}
		\gamma_0 &=& \min\left\{\frac{1}{\nicefrac{8\cL}{n} + 4L}, \frac{1}{8\sqrt{2L\cL(\tau-1)}}\right\},\\
		\gamma_K &=& \frac{\ln\left(\max\left\{2, \min\left\{\nicefrac{\|x^0 - x^*\|^2\mu^2K^2}{\left(\nicefrac{2\sigma_*^2}{n}+\nicefrac{4L\zeta^2(\tau-1)}{\mu}\right)}, \nicefrac{\|x^0 - x^*\|^2\mu^3K^3}{8L(\tau-1)\left(\sigma_*^2+\zeta_*^2\right) }\right\}\right\}\right)}{\mu K},\\
		\gamma &=& \min\left\{\gamma_0,\gamma_K\right\}
	\end{eqnarray*}
	for all $K$ such that either $\mu\gamma_K\le 1$ or $\gamma_0 \le \gamma_K$ we have that $\EE\left[f(\overline{x}^K)-f(x^*)\right]$ is of the order
	\begin{eqnarray}
	    \widetilde\cO\Bigg(\left(L + \nicefrac{\cL}{n}+ \sqrt{(\tau-1)\cL L}\right)R_0^2\exp\left(- \frac{\mu}{L + \nicefrac{\cL}{n}+ \sqrt{(\tau-1)\cL L}} K\right)&&\notag\\
		&&\hspace{-6cm}+ \frac{\sigma_*^2}{n\mu K} + \frac{L\zeta^2(\tau-1)}{\mu^2 K} + \frac{L(\tau-1)\left(\sigma_*^2+\zeta_*^2\right)}{\mu^2 K^2}\Bigg),\notag
	\end{eqnarray}
	where $R_0 = \|x^0 - x^*\|$. That is, to achieve $\EE\left[f(\overline{x}^K)-f(x^*)\right] \le \varepsilon$ in this case {\tt Local-SGD} requires
	\begin{equation*}
		\widetilde{\cO}\left(\frac{L}{\mu}+\frac{\cL}{n\mu}+\frac{\sqrt{(\tau-1)\cL L}}{\mu} + \frac{\sigma_*^2}{n\mu\varepsilon} + \frac{L\zeta^2(\tau-1)}{\mu^2 \varepsilon} + \sqrt{\frac{L(\tau-1)\left(\sigma_*^2 + \zeta_*^2\right)}{\mu^2\varepsilon}}\right)
	\end{equation*}
	iterations/oracle calls per node and $\tau$ times less communication rounds.	
\end{corollary}

Combining Theorem~\ref{thm:local_sgd_es_homo} and Lemma~\ref{lem:lemma_technical_cvx} we derive the following result for the convergence of {\tt Local-SGD} in the case when $\mu = 0$.
\begin{corollary}
	\label{cor:local_sgd_es_cvx_homo}
	Let assumptions of Theorem~\ref{thm:local_sgd_es_homo} hold with $\mu = 0$. Then for
	\begin{eqnarray*}
	\gamma_0 &=& \min\left\{\frac{1}{\nicefrac{8\cL}{n} + 4L}, \frac{1}{8\sqrt{2L\cL(\tau-1)}}\right\},\\
		\gamma &=& \min\left\{\gamma_0, \sqrt{\frac{R_0^2}{\left(\nicefrac{2\sigma_*^2}{n} + \nicefrac{4L\zeta^2(\tau-1)}{\mu}\right) K}}, \sqrt[3]{\frac{R_0^2}{8L(\tau-1)\left(\sigma_*^2+\zeta_*^2\right) K}}\right\},
	\end{eqnarray*}
	where $R_0 = \|x^0 - x^*\|$, we have that $\EE\left[f(\overline{x}^K)-f(x^*)\right]$ equals
	\begin{equation}
		 \cO\left(\frac{\left(L + \nicefrac{\cL}{n}+ \sqrt{(\tau-1)\cL L}\right)R_0^2}{K} + \sqrt{\frac{R_0^2\left(\nicefrac{\sigma_*^2}{n} + \nicefrac{L\zeta^2(\tau-1)}{\mu}\right)}{K}} + \frac{\sqrt[3]{LR_0^4(\tau-1)\left(\sigma_*^2+\zeta_*^2\right)}}{K^{\nicefrac{2}{3}}} \right).\notag
	\end{equation}
	That is, to achieve $\EE\left[f(\overline{x}^K)-f(x^*)\right] \le \varepsilon$ in this case {\tt Local-SGD} requires
	\begin{equation*}
		\cO\left(\frac{\left(L + \nicefrac{\cL}{n}+ \sqrt{(\tau-1)\cL L}\right)R_0^2}{\varepsilon} + \frac{\left(\nicefrac{\sigma_*^2}{n} + \nicefrac{L\zeta^2(\tau-1)}{\mu}\right)R_0^2}{\varepsilon^2} + \frac{R_0^2\sqrt{L(\tau-1)\left(\sigma_*^2+\zeta_*^2\right)}}{\varepsilon^{\nicefrac{3}{2}}}\right)
	\end{equation*}
	iterations/oracle calls per node and $\tau$ times less communication rounds.
\end{corollary}

\subsection{{\tt Local-SVRG} \label{sec:llsvrg}}

As an alternative to {\tt Local-SGD} when the local objective is of a finite sum structure~\eqref{eq:f_i_sum}, we propose {\tt L-SVRG}~\citep{hofmann2015variance, kovalev2019don} stochastic gradient as a local direction instead of the plain stochastic gradient. Specifically, we consider 
\[
a_i^k \eqdef  \nabla f_{i,j_i}(x_i^k) -\nabla f_{i,j_i}(w_i^k) + \nabla f_{i}(w_i^k), \qquad  b_i^k = 0,
\]
where index $1\leq j_i \leq m$ is selected uniformly at random and $w_i^k$ is a particular iterate from the local history updated as follows: 
\[
w_i^{k+1} =  \begin{cases}
            x^{k}_i & \text{w.p. } \psvrg \\
            w_i^k & \text{w.p. } 1- \psvrg.
            \end{cases}
\]

Next, we will assume that the local functions $f_{i,j}$ are $\max L_{ij}$-smooth.\footnote{It is easy to see that we must have $\max L_{ij}\geq L \geq \frac1m\max L_{ij}$.} Lastly, we will equip the mentioned method with the fixed local loop. The formal statement of the described instance of~\eqref{eq:local_sgd_def} is given as Algorithm~\ref{alg:local_svrg}.

\begin{algorithm}[h]
   \caption{{\tt Local-SVRG}}\label{alg:local_svrg}
\begin{algorithmic}[1]
   \Require learning rate $\gamma>0$, initial vector $x^0 \in \R^d$, communication period $\tau \ge 1$
	\For{$k=0,1,\dotsc$}
       	\For{$i=1,\dotsc,n$ in parallel}
            \State Choose $j_i$ uniformly at random, independently across nodes
            \State $g_i^k = \nabla f_{i,j_i}(x_i^k) -\nabla f_{i,j_i}(w_i^k) + \nabla f_{i}(w_i^k)$
            \State $w_i^{k+1} =  \begin{cases}
            x^{k}_i & \text{w.p. } \psvrg \\
            w_i^k & \text{w.p. } 1- \psvrg
            \end{cases}$
            \If{$k+1 \mod \tau = 0$}
            \State $x_i^{k+1} = x^{k+1} = \frac{1}{n}\sum\limits_{i=1}^n\left(x_i^k - \gamma g_i^k\right)$ \Comment{averaging}
            \Else
            \State $x_i^{k+1} = x_i^k - \gamma g_i^k$ \Comment{local update}
            \EndIf
        \EndFor
   \EndFor
\end{algorithmic}
\end{algorithm}

Let us next provide the details on the convergence rate. In order to do so, let us identify the parameters of Assumption~\ref{ass:sigma_k_original}.

\begin{proposition}[see \citep{gorbunov2019unified}]\label{prop:local_svrg} 
Gradient estimator $a_i^k$ satisfies Assumption~\ref{ass:sigma_k_original} with parameters $A_i = 2\max L_{ij}, B_i=2, D_{1,i} = 0, \rho_i = \psvrg, C_i = \max L_{ij} \psvrg, D_{2,i}= 0$, and $\sigma_{i,k}^2 = \frac{1}{m}\sum\limits_{j=1}^m\|\nabla f_{ij}(w_i^k)-\nabla f_{ij}(x^*)\|^2$.
\end{proposition}

\subsubsection{$\zeta$-Heterogeneous Data}
It remains to use Lemma~\ref{lem:local_solver} along with Corollary~\ref{cor:const_loop_homo} to recover all parameters of Assumption~\ref{ass:key_assumption} and obtain a convergence rate of Algorithm~\ref{alg:local_svrg} in $\zeta$-heterogeneous case.
\begin{theorem}\label{thm:local_svrg_homo}
	Assume that $f_i(x)$ is $\mu$-strongly convex and $L$-smooth for $i\in[n]$ and $f_1, \ldots, f_n$ are $\zeta$-heterogeneous, convex and $\max L_{ij}$-smooth. Then {\tt Local-SVRG} satisfies Assumption~\ref{ass:key_assumption} with
	\begin{gather*}
		A = 8\max L_{ij},\quad B = 2,\quad F = 8L\max L_{ij}, \quad D_1 = 2\zeta_*^2,\\
		A' = \frac{4 \max L_{ij}}{n} + L,\quad B' = \frac{1}{n},\quad F' = \frac{4L \max L_{ij}}{n} + 2L^2, \quad D_1' = 0,\\
		\sigma_k^2 = \frac{4}{nm}\sum\limits_{i=1}^n\sum\limits_{j=1}^m\|\nabla f_{ij}(w_i^k) - \nabla f_{ij}(x^*)\|^2,\quad \rho = q,\quad C = 8q\max L_{ij},\quad G = 4qL\max L_{ij},\\
		D_2 = 0,\quad H = \frac{8(\tau-1)(2+q)\gamma^2}{q},\quad D_3 = 2(\tau-1)\left(2\zeta_*^2 + \frac{\zeta^2}{\gamma\mu}\right)
	\end{gather*}
	with $\gamma$ satisfying
	\begin{eqnarray*}
		\gamma &\le& \min\left\{\frac{1}{2\left(\nicefrac{44\max L_{ij}}{n}+L\right)}, \frac{1}{16\sqrt{L\max L_{ij}(\tau-1)\left(1+\nicefrac{4}{(1-q)}\right)}}\right\}.
	\end{eqnarray*}
	and for all $K \ge 0$
	\begin{eqnarray}
		\EE\left[f(\overline{x}^K) - f(x^*)\right] &\le& \frac{\Phi^0}{\gamma W_K} + 8L(\tau-1)\gamma\left(\frac{\zeta^2}{\mu}+2\gamma\zeta_*^2\right), \notag
	\end{eqnarray}
	where $\Phi^0 = 2\|x^0 - x^*\|^2 +   \frac{8}{3nq}\gamma^2 \sigma_0^2 + \frac{32L(\tau-1)(2+q)\gamma^3}{q}\sigma_0^2$. In particular, if $\mu > 0$ then
	\begin{eqnarray}
		\EE\left[f(\overline{x}^K) - f(x^*)\right] &\le& \left(1 - \min\left\{\gamma\mu,\frac{q}{4}\right\}\right)^K\frac{\Phi^0}{\gamma} + 8L(\tau-1)\gamma\left(\frac{\zeta^2}{\mu}+2\gamma\zeta_*^2\right) \label{eq:local_svrg_str_cvx_homo}
	\end{eqnarray}
	and when $\mu = 0$ we have
	\begin{eqnarray}
		\EE\left[f(\overline{x}^K) - f(x^*)\right] &\le& \frac{\Phi^0}{\gamma K} + 8L(\tau-1)\gamma\left(\frac{\zeta^2}{\mu}+2\gamma\zeta_*^2\right). \label{eq:local_svrg_cvx_homo}
	\end{eqnarray}
\end{theorem}

The theorem above together with Lemma~\ref{lem:lemma2_stich} implies the following result.
\begin{corollary}\label{cor:local_svrg_str_cvx_homo}
	Let assumptions of Theorem~\ref{thm:local_svrg_homo} hold with $\mu > 0$. Then for 
	\begin{eqnarray*}
		\gamma_0 &=& \min\left\{\frac{1}{2\left(\nicefrac{44\max L_{ij}}{n}+L\right)}, \frac{1}{16\sqrt{L\max L_{ij}(\tau-1)\left(1+\nicefrac{4}{(1-q)}\right)}}\right\},\quad q = \frac{1}{m},\quad m > 1,\\
		\widetilde{\Phi}^0 &=& 2\|x^0 - x^*\|^2 +   \frac{8}{3nq}\gamma_0^2 \sigma_0^2 + \frac{32L(\tau-1)(2+q)\gamma_0^3}{q}\sigma_0^2,\\
		\gamma &=& \min\left\{\gamma_0,\frac{\ln\left(\max\left\{2, \min\left\{\nicefrac{\widetilde{\Phi}^0\mu^3K^2}{8L\zeta^2(\tau-1)}, \nicefrac{\widetilde{\Phi}^0\mu^3K^3}{16L(\tau-1)\zeta_*^2 }\right\}\right\}\right)}{\mu K}\right\},
	\end{eqnarray*}
	for all $K$ such that 
	\begin{eqnarray*}
		\text{either} && \frac{\ln\left(\max\left\{2, \min\left\{\nicefrac{\widetilde{\Phi}^0\mu^3K^2}{8L\zeta^2(\tau-1)}, \nicefrac{\widetilde{\Phi}^0\mu^3K^3}{16L(\tau-1)\zeta_*^2 }\right\}\right\}\right)}{ K} \le \frac{1}{m}\\
		\text{or} && \gamma_0 \le \frac{\ln\left(\max\left\{2, \min\left\{\nicefrac{\widetilde{\Phi}^0\mu^3K^2}{8L\zeta^2(\tau-1)}, \nicefrac{\widetilde{\Phi}^0\mu^3K^3}{16L(\tau-1)\zeta_*^2 }\right\}\right\}\right)}{\mu K}
	\end{eqnarray*}
	we have that $\EE\left[f(\overline{x}^K)-f(x^*)\right]$ is of the order
	\begin{equation}
		 \widetilde\cO\left(\frac{\widetilde{\Phi}^0}{\gamma_0}\exp\left(- \min\left\{m^{-1}, \gamma_0\mu\right\} K\right) + \frac{\zeta^2L(\tau-1)}{\mu^2 K} + \frac{L(\tau-1)\zeta_*^2}{\mu^2 K^2}\right).\notag
	\end{equation}
	That is, to achieve $\EE\left[f(\overline{x}^K)-f(x^*)\right] \le \varepsilon$ in this case {\tt Local-SVRG} requires
	\begin{equation*}
		\widetilde{\cO}\left(m + \frac{L}{\mu}+\frac{\max L_{ij}}{n\mu}+\frac{\sqrt{(\tau-1)L\max L_{ij}}}{\mu} + \frac{L\zeta^2(\tau-1)}{\mu^2 \varepsilon} + \sqrt{\frac{L(\tau-1)\zeta_*^2}{\mu^2\varepsilon}}\right)
	\end{equation*}
	iterations/oracle calls per node and $\tau$ times less communication rounds.	
\end{corollary}

Combining Theorem~\ref{thm:local_svrg_homo} and Lemma~\ref{lem:lemma_technical_cvx} we derive the following result for the convergence of {\tt Local-SVRG} in the case when $\mu = 0$.
\begin{corollary}
	\label{cor:local_svrg_cvx_homo}
	Let assumptions of Theorem~\ref{thm:local_svrg_homo} hold with $\mu = 0$. Then for
	\begin{eqnarray*}
		\gamma_0 &=& \min\left\{\frac{1}{2\left(\nicefrac{44\max L_{ij}}{n}+L\right)}, \frac{1}{16\sqrt{L\max L_{ij}(\tau-1)\left(1+\nicefrac{4}{(1-q)}\right)}}\right\},\quad q = \frac{1}{m},\quad m > 1,\\	
		\gamma &=& \min\left\{\gamma_0, \sqrt{\frac{3nR_0^2}{4m\sigma_0^2}}, \sqrt[3]{\frac{R_0^2}{16Lm(\tau-1)(2+\nicefrac{1}{m})\sigma_0^2}}, \sqrt{\frac{\mu R_0^2}{4L\zeta^2(\tau-1) K}}, \sqrt[3]{\frac{R_0^2}{8L(\tau-1)\zeta_*^2 K}}\right\},
	\end{eqnarray*}
	where $R_0 = \|x^0 - x^*\|$, we have that $\EE\left[f(\overline{x}^K)-f(x^*)\right]$ is of the order
	\begin{eqnarray*}
		\cO\Bigg(\frac{(L +\nicefrac{\max L_{ij}}{n} +\sqrt{(\tau-1)L\max L_{ij}})R_0^2 + \sqrt{\nicefrac{m\sigma_0^2R_0^2}{n}} + \sqrt[3]{Lm(\tau-1)\sigma_0^2 R_0^4}}{K}& \\
		&\hspace{-3cm}+ \sqrt{\frac{LR_0^2\zeta^2(\tau-1)}{\mu K}} + \frac{\sqrt[3]{LR_0^4(\tau-1)\zeta_*^2}}{K^{\nicefrac{2}{3}}} \Bigg).
	\end{eqnarray*}
	That is, to achieve $\EE\left[f(\overline{x}^K)-f(x^*)\right] \le \varepsilon$ in this case {\tt Local-SVRG} requires
	\begin{eqnarray*}
		\cO\Bigg(\frac{(L +\nicefrac{\max L_{ij}}{n} +\sqrt{(\tau-1)L\max L_{ij}})R_0^2 + \sqrt{\nicefrac{m\sigma_0^2R_0^2}{n}} + \sqrt[3]{Lm(\tau-1)\sigma_0^2 R_0^4}}{\varepsilon}&\\
		&\hspace{-3cm} + \frac{L\zeta^2(\tau-1)R_0^2}{\mu\varepsilon^2} + \frac{R_0^2\sqrt{L(\tau-1)\zeta_*^2}}{\varepsilon^{\nicefrac{3}{2}}}\Bigg)
	\end{eqnarray*}
	iterations/oracle calls per node and $\tau$ times less communication rounds.
\end{corollary}

\begin{remark}
	To get the rate from Tbl.~\ref{tbl:special_cases_weakly_convex} it remains to apply the following inequality:
	\begin{eqnarray*}
		\sigma_0^2 = \frac{4}{nm}\sum\limits_{i=1}^n\sum\limits_{j=1}^m\|\nabla f_{ij}(x^0)-\nabla f_{ij}(x^*)\|^2 \overset{\eqref{eq:L_smoothness}}{\le} 4\max L_{ij}^2 \|x^0-x^*\|^2.
	\end{eqnarray*}
\end{remark}

\subsubsection{Heterogeneous Data}
First of all, we need the following lemma.
\begin{lemma}\label{lem:small_lemma_local_svrg}
	Assume that $f_i(x)$ is $L$-smooth for $i\in[n]$ and $f_{ij}$ is convex and $\max L_{ij}$-smooth for $i\in[n], j\in [m]$. Then for {\tt Local-SVRG} we have
	\begin{eqnarray}
		\frac{1}{n}\sum\limits_{i=1}^n\EE\left[\|\bar g_i^k\|^2\right] &\le& 6L\EE\left[f(x^k) - f(x^*)\right] + 3L^2 \EE[V_k] + 3\zeta_*^2,\label{eq:hetero_ineq_1_local_svrg}\\
		\frac{1}{n}\sum\limits_{i=1}^n\EE\left[\left\|g_i^k-\bar g_i^k\right\|^2\right] &\le& 8\max L_{ij}\EE\left[f(x^k)-f(x^*)\right] + \frac{1}{2}\EE[\sigma_k^2] + 4L\max L_{ij}\EE[V_k], \label{eq:hetero_ineq_2_local_svrg}
	\end{eqnarray}
	where $\sigma_k^2 = \frac{4}{nm}\sum\limits_{i=1}^n\sum\limits_{j=1}^m\|\nabla f_{ij}(w_i^k) - \nabla f_{ij}(x^*)\|^2$.
\end{lemma}
\begin{proof}
	Inequality \eqref{eq:hetero_ineq_1_local_svrg} follows from $\bar g_i^k = \EE\left[g_i^k\mid x^k\right] = \nabla f_i(x_i^k)$ and inequality \eqref{eq:dnaossniadd}. Next, using Young's inequality we derive
	\begin{eqnarray*}
		\frac{1}{n}\sum\limits_{i=1}^n\EE\left[\left\|g_i^k-\bar g_i^k\right\|^2\right] &\overset{\eqref{eq:variance_decomposition}}{\le}& \frac{1}{n}\sum\limits_{i=1}^n\EE\left[\left\|g_i^k- \nabla f_i(x^*)\right\|^2\right]\\
		&\overset{\eqref{eq:a_b_norm_squared}}{\le}& \frac{2}{n}\sum\limits_{i=1}^n\EE\left[\|\nabla f_{ij_i}(w_i^k)-\nabla f_{ij_i}(x^*) - (\nabla f_i(w_i^k)-\nabla f_i(x^*))\|^2\right]\\
		&&\quad + \frac{2}{n}\sum\limits_{i=1}^n\EE\left[\|\nabla f_{ij_i}(x_i^k)-\nabla f_{ij_i}(x^*)\|^2\right]\\
		&\overset{\eqref{eq:tower_property}}{=}& \frac{2}{nm}\sum\limits_{i=1}^n\sum\limits_{j=1}^m\EE\left[\|\nabla f_{ij}(w_i^k)-\nabla f_{ij}(x^*) - (\nabla f_i(w_i^k)-\nabla f_i(x^*))\|^2\right]\\
		&&\quad + \frac{2}{nm}\sum\limits_{i=1}^n\sum\limits_{j=1}^m\EE\left[\|\nabla f_{ij}(x_i^k)-\nabla f_{ij}(x^*)\|^2\right]\\
		&\overset{\eqref{eq:L_smoothness_cor},\eqref{eq:variance_decomposition}}{\le}& \frac{4\max L_{ij}}{n}\sum\limits_{i=1}^n\EE\left[D_{f_i}(x_i^k,x^*)\right]\\
		&&\quad + \frac{2}{nm}\sum\limits_{i=1}^n\sum\limits_{j=1}^m\EE\left[\|\nabla f_{ij}(w_i^k)-\nabla f_{ij}(x^*)\|^2\right]\\
		&\overset{\eqref{eq:poiouhnkj}}{\le}& 8\max L_{ij}\EE\left[f(x^k)-f(x^*)\right] + \frac{1}{2}\EE[\sigma_k^2] + 4L\max L_{ij}\EE[V_k].	
	\end{eqnarray*}
\end{proof}

Applying Corollary~\ref{cor:const_loop}, Lemma~\ref{lem:small_lemma_local_svrg}, Proposition~\ref{prop:local_svrg} and Lemma~\ref{lem:local_solver} we get the following result.
\begin{theorem}\label{thm:local_svrg}
	Assume that $f_i(x)$ is $\mu$-strongly convex and $L$-smooth for $i\in[n]$ and $f_{ij}$ is convex and $\max L_{ij}$-smooth for $i\in[n], j\in [m]$. Then {\tt Local-SVRG} satisfies Assumption~\ref{ass:hetero_second_moment} with
	\begin{gather*}
		\tA = 3L,\quad \hA = 4\max L_{ij},\quad \tB = 0,\quad \hB = \frac{1}{2},\quad \tF = 3L^2,\quad \hF = 4L\max L_{ij}, \quad \tD_1 = 3\zeta_*^2,\\
		\hD_1 = 0,\quad A' = \frac{4 \max L_{ij}}{n} + L,\quad B' = \frac{1}{n},\quad F' = \frac{4L \max L_{ij}}{n} + 2L^2, \quad D_1' = 0,\\
		\sigma_k^2 = \frac{4}{nm}\sum\limits_{i=1}^n\sum\limits_{j=1}^m\|\nabla f_{ij}(w_i^k) - \nabla f_{ij}(x^*)\|^2,\quad \rho = q,\quad C = 8q\max L_{ij},\quad G = 4qL\max L_{ij},\\
		D_2 = 0,\quad H = \frac{2e(\tau-1)(2+q)\gamma^2}{q},\quad D_3 = 6e(\tau-1)^2\zeta_*^2
	\end{gather*}
	with $\gamma$ satisfying
	\begin{eqnarray*}
		\gamma &\le& \min\left\{\frac{1}{2\left(\nicefrac{44\max L_{ij}}{n}+L\right)}, \frac{1}{4\sqrt{2eL(\tau-1)\left(3L(\tau-1)+4\max L_{ij} + \nicefrac{8\max L_{ij}}{(1-q)}\right)}}\right\}.
	\end{eqnarray*}
	and for all $K \ge 0$
	\begin{eqnarray}
		\EE\left[f(\overline{x}^K) - f(x^*)\right] &\le& \frac{\Phi^0}{\gamma W_K} + 24e L(\tau-1)^2\zeta_*^2\gamma^2, \notag
	\end{eqnarray}
	where $\Phi^0 = 2\|x^0 - x^*\|^2 +   \frac{8}{3nq}\gamma^2 \sigma_0^2 + \frac{8eL(\tau-1)(2+q)\gamma^3}{q}\sigma_0^2$
	In particular, if $\mu > 0$ then
	\begin{eqnarray}
		\EE\left[f(\overline{x}^K) - f(x^*)\right] &\le& \left(1 - \min\left\{\gamma\mu,\frac{q}{4}\right\}\right)^K\frac{\Phi^0}{\gamma} + 24e L(\tau-1)^2\zeta_*^2\gamma^2 \label{eq:local_svrg_str_cvx}
	\end{eqnarray}
	and when $\mu = 0$ we have
	\begin{eqnarray}
		\EE\left[f(\overline{x}^K) - f(x^*)\right] &\le& \frac{\Phi^0}{\gamma K} + 24e L(\tau-1)^2\zeta_*^2\gamma^2. \label{eq:local_svrg_cvx}
	\end{eqnarray}
\end{theorem}

The theorem above together with Lemma~\ref{lem:lemma2_stich} implies the following result.
\begin{corollary}\label{cor:local_svrg_str_cvx}
	Let assumptions of Theorem~\ref{thm:local_svrg} hold with $\mu > 0$. Then for 
	\begin{eqnarray*}
		\gamma_0 &=& \min\left\{\frac{1}{2\left(\nicefrac{44\max L_{ij}}{n}+L\right)}, \frac{1}{4\sqrt{2eL(\tau-1)\left(3L(\tau-1)+4\max L_{ij} + \nicefrac{8\max L_{ij}}{(1-q)}\right)}}\right\},\\
		\widetilde{\Phi}^0 &=& 2\|x^0 - x^*\|^2 +   \frac{8}{3nq}\gamma_0^2 \sigma_0^2 + \frac{8eL(\tau-1)(2+q)\gamma_0^3}{q}\sigma_0^2,\quad q = \frac{1}{m},\quad m > 1,\\
		\gamma &=& \min\left\{\gamma_0,\frac{\ln\left(\max\left\{2, \nicefrac{\widetilde{\Phi}^0\mu^3K^3}{24eL(\tau-1)^2\zeta_*^2 }\right\}\right)}{\mu K}\right\},
	\end{eqnarray*}
	for all $K$ such that 
	\begin{eqnarray*}
		\text{either } \frac{\ln\left(\max\left\{2, \nicefrac{\widetilde{\Phi}^0\mu^3K^3}{24eL(\tau-1)^2\zeta_*^2 }\right\}\right)}{ K} \le \frac{1}{m} \text{ or }  \gamma_0 \le \frac{\ln\left(\max\left\{2, \nicefrac{\widetilde{\Phi}^0\mu^3K^3}{24eL(\tau-1)^2\zeta_*^2 }\right\}\right)}{\mu K}
	\end{eqnarray*}
	we have that $\EE\left[f(\overline{x}^K)-f(x^*)\right]$ is of the order
	\begin{equation}
		 \widetilde\cO\left(\frac{\widetilde{\Phi}^0}{\gamma_0}\exp\left(- \min\left\{m^{-1}, \gamma_0\mu\right\} K\right)  + \frac{L(\tau-1)^2\zeta_*^2}{\mu^2 K^2}\right).\notag
	\end{equation}
	That is, to achieve $\EE\left[f(\overline{x}^K)-f(x^*)\right] \le \varepsilon$ in this case {\tt Local-SVRG} requires
	\begin{equation*}
		\widetilde{\cO}\left(m + \frac{L\tau}{\mu}+\frac{\max L_{ij}}{n\mu}+\frac{\sqrt{(\tau-1)L\max L_{ij}}}{\mu} + \sqrt{\frac{L(\tau-1)^2\zeta_*^2}{\mu^2\varepsilon}}\right)
	\end{equation*}
	iterations/oracle calls per node and $\tau$ times less communication rounds.	
\end{corollary}

Combining Theorem~\ref{thm:local_svrg} and Lemma~\ref{lem:lemma_technical_cvx} we derive the following result for the convergence of {\tt Local-SVRG} in the case when $\mu = 0$.
\begin{corollary}
	\label{cor:local_svrg_cvx}
	Let assumptions of Theorem~\ref{thm:local_svrg} hold with $\mu = 0$. Then for $q = \frac{1}{m},$ $m > 1$ and
	\begin{eqnarray*}
		\gamma_0 &=& \min\left\{\frac{1}{2\left(\nicefrac{44\max L_{ij}}{n}+L\right)}, \frac{1}{4\sqrt{2eL(\tau-1)\left(3L(\tau-1)+4\max L_{ij} + \nicefrac{8\max L_{ij}}{(1-q)}\right)}}\right\},\\	
		\gamma &=& \min\left\{\gamma_0, \sqrt{\frac{3nR_0^2}{4m\sigma_0^2}}, \sqrt[3]{\frac{R_0^2}{4eLm(\tau-1)(2+\nicefrac{1}{m})\sigma_0^2}}, \sqrt[3]{\frac{R_0^2}{12eL(\tau-1)^2\zeta_*^2 K}}\right\},
	\end{eqnarray*}
	where $R_0 = \|x^0 - x^*\|$, we have that $\EE\left[f(\overline{x}^K)-f(x^*)\right]$ is of the order
	\begin{eqnarray*}
		\cO\Bigg(\frac{(L\tau +\nicefrac{\max L_{ij}}{n} +\sqrt{(\tau-1)L\max L_{ij}})R_0^2 + \sqrt{\nicefrac{m\sigma_0^2R_0^2}{n}} + \sqrt[3]{Lm(\tau-1)\sigma_0^2 R_0^4}}{K} &&\\
		&&\hspace{-2cm} + \frac{\sqrt[3]{LR_0^4(\tau-1)^2\zeta_*^2}}{K^{\nicefrac{2}{3}}} \Bigg).
	\end{eqnarray*}
	That is, to achieve $\EE\left[f(\overline{x}^K)-f(x^*)\right] \le \varepsilon$ in this case {\tt Local-SVRG} requires
	\begin{eqnarray*}
		\cO\Bigg(\frac{(L\tau +\nicefrac{\max L_{ij}}{n} +\sqrt{(\tau-1)L\max L_{ij}})R_0^2 + \sqrt{\nicefrac{m\sigma_0^2R_0^2}{n}} + \sqrt[3]{Lm(\tau-1)\sigma_0^2 R_0^4}}{\varepsilon}&&\\
		&&\hspace{-2cm} + \frac{R_0^2\sqrt{L(\tau-1)^2\zeta_*^2}}{\varepsilon^{\nicefrac{3}{2}}}\Bigg)
	\end{eqnarray*}
	iterations/oracle calls per node and $\tau$ times less communication rounds.
\end{corollary}

\begin{remark}
	To get the rate from Tbl.~\ref{tbl:special_cases_weakly_convex} it remains to apply the following inequality:
	\begin{eqnarray*}
		\sigma_0^2 = \frac{4}{nm}\sum\limits_{i=1}^n\sum\limits_{j=1}^m\|\nabla f_{ij}(x^0)-\nabla f_{ij}(x^*)\|^2 \overset{\eqref{eq:L_smoothness}}{\le} 4\max L_{ij}^2 \|x^0-x^*\|^2.
	\end{eqnarray*}
\end{remark}

\subsection{{\tt S*-Local-SGD}}\label{sec:sgd_star_bounded_var}
In this section we consider the same settings as in Section~\ref{sec:sgd_bounded_var} and our goal is to remove one of the main drawbacks of {\tt Local-SGD} in heterogeneous case which in the case of $\mu$-strongly convex $f_i$ with $\mu > 0$ converges with linear rate only to the neighbourhood of the solution even in the full-gradients case, i.e.\ when $D_{1,i} = 0$ for all $i\in[n]$. However, we start with unrealistic assumption that $i$-th node has an access to $\nabla f_i(x^*)$ for all $i\in[n]$. Under this assumption we present a new method called Star-Shifted Local-SGD ({\tt S*-Local-SGD}, see Algorithm~\ref{alg:local_sgd_star}).

\begin{algorithm}[h]
   \caption{{\tt S*-Local-SGD}}\label{alg:local_sgd_star}
\begin{algorithmic}[1]
   \Require learning rate $\gamma>0$, initial vector $x^0 \in \R^d$, communication period $\tau \ge 1$
	\For{$k=0,1,\dotsc$}
       	\For{$i=1,\dotsc,n$ in parallel}
            \State Sample $\hat g^{k}_i = \nabla f_{\xi_i^k}(x_i^k)$ independently from other nodes
            \State $g_i^k = \hat g_i^k - \nabla f_i(x^*)$
            \If{$k+1 \mod \tau = 0$}
            \State $x_i^{k+1} = x^{k+1} = \frac{1}{n}\sum\limits_{i=1}^n\left(x_i^k - \gamma g_i^k\right)$ \Comment{averaging}
            \Else
            \State $x_i^{k+1} = x_i^k - \gamma g_i^k$ \Comment{local update}
            \EndIf
        \EndFor
   \EndFor
\end{algorithmic}
\end{algorithm}

\begin{lemma}\label{lem:local_sgd_star_second_moment}
	Let $f_i$ be convex and $L$-smooth for all $i\in[n]$. Then for all $k\ge 0$
	\begin{eqnarray}
		\frac{1}{n}\sum\limits_{i=1}^n \EE\left[g_i^k\mid x_i^k\right] &=& \frac{1}{n}\sum\limits_{i=1}^n\nabla f_i(x_i^k), \label{eq:unbiasedness_local_sgd_star}\\
		\frac{1}{n}\sum\limits_{i=1}^n \|\bar{g}_i^k\|^2 &\le& 4L\left(f(x^k) - f(x^*)\right) + 2L^2 V_k,\label{eq:second_moment_2_local_sgd_star}\\
		\frac{1}{n}\sum\limits_{i=1}^n \EE\left[\|g_i^k - \bar{g}_i^k\|^2\mid x_i^k\right] &\le& \sigma^2,\label{eq:variance_local_sgd_star}\\
		\EE\left[\left\|\frac{1}{n}\sum\limits_{i=1}^ng_i^k\right\|^2\mid x^k\right] &\le& 4L\left(f(x^k) - f(x^*)\right) + 2L^2 V_k + \frac{\sigma^2}{n},\label{eq:second_moment_local_sgd_star_2}
	\end{eqnarray}
	where $\sigma^2 \eqdef \frac{1}{n}\sum_{i=1}^nD_{1,i}$ and $\EE[\cdot\mid x^k]\eqdef \EE[\cdot\mid x_1^k,\ldots,x_n^k]$.
\end{lemma}
\begin{proof}
	First of all, we notice that $\EE\left[g_i^k\mid x_i^k\right] = \nabla f_i(x_i^k)-\nabla f_i(x^*)$ and
	\begin{equation*}
		\frac{1}{n}\sum\limits_{i=1}^n\EE\left[g_i^k\mid x_i^k\right] = \frac{1}{n}\sum\limits_{i=1}^n\left(\nabla f_i(x_i^k)-\nabla f_i(x^*)\right) = \frac{1}{n}\sum\limits_{i=1}^n\nabla f_i(x_i^k).
	\end{equation*}
	Using this we get
	\begin{eqnarray*}
		\frac{1}{n}\sum\limits_{i=1}^n \|\bar{g}_i^k\|^2 &=& \frac{1}{n}\sum\limits_{i=1}^n \|\nabla f_i(x_i^k) - \nabla f_i(x^*)\|^2 \overset{\eqref{eq:L_smoothness_cor}}{\le} \frac{2L}{n}\sum\limits_{i=1}^nD_{f_i}(x_i^k,x^*)\\
		&\overset{\eqref{eq:poiouhnkj}}{\le}& 4L\left(f(x^k) - f(x^*)\right) + 2L^2 V_k
	\end{eqnarray*}
	and
	\begin{equation*}
		\frac{1}{n}\sum\limits_{i=1}^n \EE\left[\|g_i^k - \bar{g}_i^k\|^2\mid x_i^k\right] = \frac{1}{n}\sum\limits_{i=1}^n \EE\left[\|\nabla f_{\xi_i^k}(x_i^k) - \nabla f_i(x_i^k)\|^2\right] \overset{\eqref{eq:bounded_variance}}{\le} \frac{1}{n}\sum\limits_{i=1}^nD_{1,i} =: \sigma^2.
	\end{equation*}
	Finally, using independence of $g_1^k,g_2^k,\ldots,g_n^k$ and $\frac{1}{n}\sum_{i=1}^n\nabla f_i(x^*) = \nabla f(x^*) = 0$ we obtain
	\begin{eqnarray*}
		\EE\left[\left\|\frac{1}{n}\sum\limits_{i=1}^ng_i^k\right\|^2\mid x^k\right] &\overset{\eqref{eq:variance_decomposition},\eqref{eq:unbiasedness_local_sgd_star}}{=}& \EE\left[\left\|\frac{1}{n}\sum\limits_{i=1}^n\left(g_i^k - \nabla f_i(x_i^k)\right)\right\|^2\mid x^k\right] + \left\|\frac{1}{n}\sum\limits_{i=1}^n\nabla f_i(x_i^k)\right\|^2\\
		&=& \EE\left[\left\|\frac{1}{n}\sum\limits_{i=1}^n\left(\nabla f_{\xi_i^k}(x_i^k) - \nabla f_i(x_i^k)\right)\right\|^2\mid x^k\right] + \left\|\frac{1}{n}\sum\limits_{i=1}^n\nabla f_i(x_i^k)\right\|^2\\
		&=& \frac{1}{n^2}\sum\limits_{i=1}^n\EE_{\xi_i^k}\left[\|\nabla f_{\xi_i^k}(x_i^k) - \nabla f_i(x_i^k)\|^2\right] + \left\|\frac{1}{n}\sum\limits_{i=1}^n\nabla f_i(x_i^k)\right\|^2\\
		&\overset{\eqref{eq:bounded_variance},\eqref{eq:vdgasvgda}}{\le}& 4L\left(f(x^k) - f(x^*)\right) + 2L^2 V_k + \frac{\sigma^2}{n}.
	\end{eqnarray*}
\end{proof}

Applying Corollary~\ref{cor:const_loop} and Lemma~\ref{lem:local_sgd_star_second_moment} we get the following result.
\begin{theorem}\label{thm:local_sgd_star}
	Assume that $f_i(x)$ is $\mu$-strongly convex and $L$-smooth for every $i\in[n]$. Then {\tt S*-Local-SGD} satisfies Assumption~\ref{ass:hetero_second_moment} with
	\begin{gather*}
		\tA = 2L,\quad \hA = 0,\quad \tB = \hB = 0,\quad \tF = 2L^2,\quad \hF = 0,\quad \tD_1 = 0,\quad \hD_1 = \sigma^2 := \frac{1}{n}\sum\limits_{i=1}^nD_{1,i}\\
		A' = 2L,\quad B' = 0,\quad F' = 2L^2, \quad D_1' = \frac{\sigma^2}{n},\quad \sigma_k^2 \equiv 0,\quad \rho = 1,\quad C = 0,\quad G = 0,\quad D_2 = 0,\\
		H = 0,\quad D_3 = 2e(\tau-1)\sigma^2.
	\end{gather*}
	Consequently, if
	\begin{eqnarray*}
		\gamma &\le& \min\left\{\frac{1}{4L}, \frac{1}{8\sqrt{e}(\tau-1)L}\right\}.
	\end{eqnarray*}
we have for $\mu > 0$
	\begin{eqnarray}
		\EE\left[f(\overline{x}^K) - f(x^*)\right] &\le& \left(1 - \gamma\mu\right)^K\frac{2\|x^0-x^*\|^2}{\gamma} + 2\gamma\left(\frac{\sigma^2}{n} + 4eL(\tau-1)\gamma \sigma^2\right) \notag
	\end{eqnarray}
	and when $\mu = 0$ we have
	\begin{eqnarray}
		\EE\left[f(\overline{x}^K) - f(x^*)\right] &\le& \frac{2\|x^0-x^*\|^2}{\gamma K} + 2\gamma\left(\frac{\sigma^2}{n} + 4eL(\tau-1)\gamma \sigma^2\right). \notag
	\end{eqnarray}
\end{theorem}

In the special case when $\nabla f_{\xi_i^k}(x_i^k) = \nabla f_i(x_i^k)$ for all $i\in[n]$ and $k\ge 0$ we obtain {\tt S*-Local-GD} which converges with $\cO\left(\tau\kappa\ln \frac{1}{\varepsilon}\right)$ rate when $\mu > 0$ and with $\cO\left(\frac{L\tau\|x^0 - x^*\|^2}{\varepsilon}\right)$ rate when $\mu = 0$ to the exact solution asymptotically.

The theorem above together with Lemma~\ref{lem:lemma2_stich} implies the following result.
\begin{corollary}\label{cor:sstarsgd}
	Let assumptions of Theorem~\ref{thm:local_sgd_star} hold with $\mu > 0$. Then for
	\begin{eqnarray*}
	    \gamma &=& \min\left\{\frac{1}{4L}, \frac{1}{8\sqrt{e}(\tau-1)L},  \gamma_K\right\},\\
		\gamma_K &=& \frac{\ln\left(\max\left\{2,\min\left\{\nicefrac{\|x^0 - x^*\|^2n\mu^2K^2}{\sigma^2},\nicefrac{\|x^0 - x^*\|^2\mu^3K^3}{4eL(\tau-1)\sigma^2}\right\}\right\}\right)}{\mu K}
	\end{eqnarray*}
	for all $K$ such that 
	\begin{eqnarray*}
		\text{either}\quad \mu\gamma_K \le 1\quad
		\text{or}\quad \min\left\{\frac{1}{4L}, \frac{1}{8\sqrt{e}(\tau-1)L}\right\} \le \gamma_K
	\end{eqnarray*}
	we have that
	\begin{equation}
		\EE\left[f(\overline{x}^K)-f(x^*)\right] = \widetilde\cO\left(\tau L\|x^0 - x^*\|^2\exp\left(- \frac{\mu}{\tau L} K\right) + \frac{\sigma^2}{n\mu K} + \frac{L(\tau-1)\sigma^2}{\mu^2K^2}\right).\notag
	\end{equation}
	That is, to achieve $\EE\left[f(\overline{x}^K)-f(x^*)\right] \le \varepsilon$ in this case {\tt S*-Local-SGD} requires
	\begin{equation*}
		\widetilde\cO\left(\frac{\tau L}{\mu} + \frac{\sigma^2}{n\mu\varepsilon} + \sqrt{\frac{L(\tau-1)\sigma^2}{\mu^2\varepsilon}}\right)
	\end{equation*}
	iterations/oracle calls per node and $\tau$ times less communication rounds.
\end{corollary}

Combining Theorem~\ref{thm:local_sgd_star} and Lemma~\ref{lem:lemma_technical_cvx} we derive the following result for the convergence of {\tt S*-Local-SGD} in the case when $\mu = 0$.
\begin{corollary}
	\label{cor:local_sgd_star_cvx}
	Let assumptions of Theorem~\ref{thm:local_sgd_star} hold with $\mu = 0$. Then for
	\begin{equation*}
		\gamma = \min\left\{\frac{1}{4L}, \frac{1}{8\sqrt{e}(\tau-1)L}, \sqrt{\frac{nR_0^2}{\sigma^2 K}}, \sqrt[3]{\frac{R_0^2}{4eL(\tau-1)\sigma^2K}}\right\},
	\end{equation*}
	where $R_0 = \|x^0 - x^*\|$, we have that
	\begin{equation}
		\EE\left[f(\overline{x}^K)-f(x^*)\right] = \cO\left(\frac{\tau LR_0^2}{K} + \sqrt{\frac{R_0^2\sigma^2}{nK}} + \frac{\sqrt[3]{LR_0^4(\tau-1)\sigma^2}}{K^{\nicefrac{2}{3}}} \right).\notag
	\end{equation}
	That is, to achieve $\EE\left[f(\overline{x}^K)-f(x^*)\right] \le \varepsilon$ in this case {\tt S*-Local-SGD} requires
	\begin{equation*}
		\cO\left(\frac{\tau LR_0^2}{\varepsilon} + \frac{R_0^2\sigma^2}{n\varepsilon^2} + \frac{R_0^2\sqrt{L(\tau-1)\sigma^2}}{\varepsilon^{\nicefrac{3}{2}}}\right)
	\end{equation*}
	iterations/oracle calls per node and $\tau$ times less communication rounds.
\end{corollary}

\subsection{{\tt SS-Local-SGD}}

\subsubsection{Uniformly Bounded Variance}\label{sec:loopless_local_svrg}
In this section we consider the same settings as in Section~\ref{sec:sgd_bounded_var}

\begin{algorithm}[h]
   \caption{Stochastically Shifted {\tt Local-SGD} ({\tt {\tt SS-Local-SGD}})}\label{alg:l_local_svrg}
\begin{algorithmic}[1]
   \Require learning rate $\gamma>0$, initial vector $x^0 \in \R^d$, probability of communication $p\in(0,1]$, probability of the shift's update $q\in(0,1]$, batchsize $r$ for computing shifts
   \State $y^0 = x^0$
   \State For $i\in[n]$ compute $r$ independent samples $\nabla f_{\oxi_{i,1}^0}(y^0), \nabla f_{\oxi_{i,2}^0}(y^0), \ldots, \nabla f_{\oxi_{i,r}^0}(y^0)$, set $\nabla f_{\oxi_i^0}(y^0) = \frac{1}{r}\sum_{j=1}^r\nabla f_{\oxi_{i,j}^0}(y^0)$ and $\nabla f_{\oxi^0}(y^0) = \frac{1}{n}\sum_{i=1}^n\nabla f_{\oxi_i^0}(y^0)$
	\For{$k=0,1,\dotsc$}
       	\For{$i=1,\dotsc,n$ in parallel}
       		\State Sample $\nabla f_{\xi_i^k}(x_i^k)$ independently from other nodes
            \State $g_i^k = \nabla f_{\xi_i^k}(x_i^k) - \nabla f_{\txi_i^k}(y^k) + \nabla f_{\txi^k}(y^k)$, where $\nabla f_{\oxi_i^k}(y^k) = \frac{1}{r}\sum_{j=1}^r\nabla f_{\oxi_{i,j}^k}(y^k)$ and $\nabla f_{\oxi^k}(y^k) = \frac{1}{n}\sum_{i=1}^n\nabla f_{\oxi_i^k}(y^k)$
            \State $x_i^{k+1} = \begin{cases}x^{k+1},&\text{w.p. } p,\\
            x_i^k - \gamma g_i^k,& \text{w.p. } 1 - p, \end{cases}$ where $x^{k+1} = \frac{1}{n}\sum\limits_{i=1}^n(x_i^k - \gamma g_i^k)$
            \State $y^{k+1} = \begin{cases}x^k,&\text{w.p. } q,\\
            y^k,& \text{w.p. } 1 - q, \end{cases}$ and for all $i\in[n]$, $j\in[r]\;$ $\oxi_{i,j}^{k+1}$ is $\begin{cases}\text{a fresh sample},&\text{if } y^{k+1}\neq y^k,\\
            \text{equal to } \oxi_{i,j}^k,& \text{otherwise}. \end{cases}$
        \EndFor
   \EndFor
\end{algorithmic}
\end{algorithm}

The main algorithm in this section is Stochastically Shifted {\tt Local-SGD} ({\tt SS-Local-SVRG}, see Algorithm~\ref{alg:l_local_svrg}). We notice that the updates for $x_i^{k+1}$ and $y^{k+1}$ can be dependent, e.g., one can take $p = q$ and update $y^{k+1}$ as $x^k$ every time $x_i^{k+1}$ is updated by $x^{k+1}$. Moreover, with probability $q$ line $8$ implies a round of communication and computation of new stochastic gradient by each worker.

We emphasize that in expectation $y^k$ is updated only once per $\left\lceil\nicefrac{1}{q}\right\rceil$ iterations. Therefore, if $r = O\left(\nicefrac{1}{q}\right)$ and $q \le p$, then up to a constant numerical factor the overall expected number of oracle calls and communication rounds are the same as for {\tt Local-SGD} with either the same probability $p$ of communication or with constant local loop length $\tau = \left\lceil\nicefrac{1}{p}\right\rceil$.

Finally, we notice that due to independence of $\oxi_{i,1}^k,\oxi_{i,2}^k,\ldots,\oxi_{i,r}^k$ we have
\begin{equation}
	\EE\|\nabla f_{\oxi_i^k}(y^k)-\nabla f_i(y^k)\|^2 \overset{\eqref{eq:bounded_variance}}{\le} \frac{D_{1,i}}{r}. \label{eq:stoch_shifts_variance}
\end{equation}

\begin{lemma}\label{lem:loopless_local_svrg_second_moment}
	Let $f_i$ be convex and $L$-smooth for all $i\in[n]$. Then for all $k\ge 0$
	\begin{eqnarray}
		\frac{1}{n}\sum\limits_{i=1}^n \EE_k\left[g_i^k\right] &=& \frac{1}{n}\sum\limits_{i=1}^n\nabla f_i(x_i^k), \label{eq:unbiasedness_loopless_local_svrg}\\
		\frac{1}{n}\sum\limits_{i=1}^n \EE\left[\|\bar{g}_i^k\|^2\right] &\le& 8L\EE\left[f(x^k)-f(x^*)\right] + 2\EE[\sigma_k^2] + 4L^2\EE[V_k] + \frac{2\sigma^2}{r},\label{eq:second_moment_loopless_local_svrg}\\
		\frac{1}{n}\sum\limits_{i=1}^n \EE\left[\|g_i^k-\bar{g}_i^k\|^2\right] &\le& \sigma^2,\label{eq:variance_loopless_local_svrg}\\
		\EE\left[\left\|\frac{1}{n}\sum\limits_{i=1}^ng_i^k\right\|^2\right] &\le& 4L\EE\left[f(x^k) - f(x^*)\right] + 2L^2\EE\left[V_k\right] + \frac{\sigma^2}{n},\label{eq:second_moment_loopless_local_svrg_2}
	\end{eqnarray}
	where $\sigma_k^2 \eqdef \frac{1}{n}\sum\limits_{i=1}^n\left\|\nabla f_i(y^k) - \nabla f_i(x^*)\right\|^2$ and $\sigma^2 \eqdef \frac{1}{n}\sum_{i=1}^nD_{1,i}$.
\end{lemma}
\begin{proof}
	We start with unbiasedness:
	\begin{eqnarray*}
		\frac{1}{n}\sum\limits_{i=1}^n \EE_k\left[g_i^k\right] &=& \frac{1}{n}\sum\limits_{i=1}^n \EE_k\left[\nabla f_{\xi_i^k}(x_i^k) - \nabla f_{\oxi_i^k}(y^k) + \nabla f_{\oxi^k}(y^k)\right]\\
		&=& \frac{1}{n}\sum\limits_{i=1}^n\EE_k\left[\nabla f_{\xi_i^k}(x_i^k)\right] + \EE_k\left[\nabla f_{\oxi^k}(y^k) - \frac{1}{n}\sum\limits_{i=1}^n \nabla f_{\oxi_i^k}(y^k)\right] = \frac{1}{n}\sum\limits_{i=1}^n\nabla f_i(x_i^k).
	\end{eqnarray*}
	Using this we get
	\begin{eqnarray*}
		\frac{1}{n}\sum\limits_{i=1}^n\EE\left[\|\bar{g}_i^k\|^2\right] &\overset{\eqref{eq:a_b_norm_squared}}{\le}& \frac{2}{n}\sum\limits_{i=1}^n\EE\left[\|\nabla f_i(x_i^k) - \nabla f_i(x^*)\|^2\right]\\
		&&\quad + \frac{2}{n}\sum\limits_{i=1}^n\EE\left[\left\|\nabla f_{\oxi_i^k}(y^k) - \nabla f_i(x^*) - \left(\nabla f_{\oxi^k}(y^k) - \nabla f(x^*)\right)\right\|^2\right]\\
		&\overset{\eqref{eq:L_smoothness_cor},\eqref{eq:variance_decomposition}}{\le}& \frac{4L}{n}\sum\limits_{i=1}^n\EE\left[D_{f_i}(x_i^k,x^*)\right] + \frac{2}{n}\sum\limits_{i=1}^n\EE\left[\left\|\nabla f_{\oxi_i^k}(y^k) - \nabla f_i(x^*)\right\|^2\right]\\
		&\overset{\eqref{eq:poiouhnkj},\eqref{eq:variance_decomposition}}{\le}& 8L\EE\left[f(x^k)-f(x^*)\right] + 4L^2\EE[V_k] + \frac{2}{n}\sum\limits_{i=1}^n\EE\left[\left\|\nabla f_{i}(y^k) - \nabla f_i(x^*)\right\|^2\right]\\
		&&\quad + \frac{2}{n}\sum\limits_{i=1}^n\EE\left[\left\|\nabla f_{\oxi_i^k}(y^k) - \nabla f_i(y^k)\right\|^2\right]\\
		&\overset{\eqref{eq:stoch_shifts_variance}}{\le}& 8L\EE\left[f(x^k)-f(x^*)\right] + 2\EE[\sigma_k^2] + 4L^2\EE[V_k] + \frac{2\sigma^2}{r}
	\end{eqnarray*}
	and
	\begin{equation*}
		\frac{1}{n}\sum\limits_{i=1}^n \EE\left[\|g_i^k-\bar{g}_i^k\|^2\right] = \frac{1}{n}\sum\limits_{i=1}^n \EE\left[\|\nabla f_{\xi_i^k}(x_i^k)-\nabla f_i(x_i^k)\|^2\right] \overset{\eqref{eq:bounded_variance}}{\le} \sigma^2.
	\end{equation*}
	Finally, we use independence of $\nabla f_{\xi_1^k}(x_1^k),\ldots,\nabla f_{\xi_n^k}(x_n^k)$ and derive
	\begin{eqnarray*}
		\EE\left[\left\|\frac{1}{n}\sum\limits_{i=1}^n g_i^k\right\|^2\right] &=& \EE\left[\left\|\frac{1}{n}\sum\limits_{i=1}^n\nabla f_{\xi_i^k}(x_i^k)\right\|^2\right]\\
		&\overset{\eqref{eq:variance_decomposition}}{=}& \EE\left[\left\|\frac{1}{n}\sum\limits_{i=1}^n \nabla f_i(x_i^k)\right\|^2\right] + \EE\left[\left\|\frac{1}{n}\sum\limits_{i=1}^n\left(\nabla f_{\xi_i^k}(x_i^k) - \nabla f_i(x_i^k)\right)\right\|^2\right]\\
		&\overset{\eqref{eq:vdgasvgda}}{\le}& 4L\EE\left[f(x^k)-f(x^*)\right] + 2L^2\EE[V_k] + \frac{1}{n^2}\sum\limits_{i=1}^n\EE\left[\|\nabla f_{\xi_i^k}(x_i^k) - \nabla f_i(x_i^k)\|^2\right]\\
		&\overset{\eqref{eq:bounded_variance}}{\le}& 4L\EE\left[f(x^k) - f(x^*)\right] + 2L^2\EE\left[V_k\right] + \frac{\sigma^2}{n}
	\end{eqnarray*}
	which finishes the proof.
\end{proof}

\begin{lemma}\label{lem:loopless_local_svrg_sigma_k_bound}
	Let $f_i$ be convex and $L$-smooth for all $i\in[n]$. Then for all $k\ge 0$
	\begin{eqnarray}
		\EE\left[\sigma_{k+1}^2\right] &\le& (1-q)\EE\left[\sigma_k^2\right] + 2Lq\EE\left[f(x^k) - f(x^*)\right]\label{eq:loopless_local_svrg_sigma_k_bound}
	\end{eqnarray}
	where $\sigma_k^2 \eqdef \frac{1}{n}\sum\limits_{i=1}^n\left\|\nabla f_i(y^k) - \nabla f_i(x^*)\right\|^2$.
\end{lemma}
\begin{proof}
	By definition of $y^{k+1}$ we have
	\begin{eqnarray*}
		\EE\left[\sigma_{k+1}^2\mid x_1^k,\ldots, x_n^k\right] &=& \frac{1-q}{n}\sum\limits_{i=1}^n\|\nabla f_i(y^k) - \nabla f_i(x^*)\|^2 + \frac{q}{n}\sum\limits_{i=1}^n\|\nabla f_i(x^k) - \nabla f_i(x^*)\|^2\\
		&\overset{\eqref{eq:L_smoothness_cor}}{\le}& (1-q)\sigma_k^2 + 2Lq(f(x^k) - f(x^*)).
	\end{eqnarray*}
	Taking the full mathematical expectation on both sides of previous inequality and using the tower property \eqref{eq:tower_property} we get the result.
\end{proof}

Using Corollary~\ref{cor:rand_loop} we obtain the following theorem.
\begin{theorem}\label{thm:ss_local_sgd}
	Assume that $f_i(x)$ is $\mu$-strongly convex and $L$-smooth for every $i\in[n]$. Then {\tt SS-Local-SGD} satisfies Assumption~\ref{ass:hetero_second_moment} with
	\begin{gather*}
		\tA = 4L,\quad \hA = 0,\quad \tB = 2,\quad \hB = 0,\quad \tF = 4L^2,\quad \hF = 0, \quad \tD_1 = \frac{2\sigma^2}{r},\\
		\hD_1 = \sigma^2 = \frac{1}{n}\sum\limits_{i=1}^nD_{1,i},\quad A' = 2L,\quad B' = 0,\quad F' = 2L^2, \quad D_1' = \frac{\sigma^2}{n},\\
		\sigma_k^2 = \frac{1}{n}\sum\limits_{i=1}^n\left\|\nabla f_i(y^k) - \nabla f_i(x^*)\right\|^2,\quad \rho = q,\quad C = Lq,\quad G = 0,\quad D_2 = 0,\\
		H = \frac{128(1-p)(2+p)(2+q)\gamma^2}{3p^2q},\quad D_3 = \frac{8(1-p)}{p^2}\left(\frac{2(p+2)\sigma^2}{r} + p\sigma^2\right)
	\end{gather*}
	under assumption that
	\begin{eqnarray*}
		\gamma &\le& \min\left\{\frac{1}{4L}, \frac{p\sqrt{3}}{32L\sqrt{2(1-p)(2+p)\left(1+\nicefrac{1}{(1-q)}\right)}}\right\}.
	\end{eqnarray*}
	Moreover, for $\mu > 0$ we have
	\begin{eqnarray}
		\EE\left[f(\overline{x}^K) - f(x^*)\right] &\le& \left(1 - \min\left\{\gamma\mu,\frac{q}{4}\right\}\right)^K\frac{\Phi^0}{\gamma}\notag\\
		&&\quad +2\gamma\left(\frac{\sigma^2}{n} + \gamma \frac{16L(1-p)}{p^2}\left(\frac{2(p+2)\sigma^2}{r} + p\sigma^2\right)\right) \notag
	\end{eqnarray}
	and when $\mu = 0$ we have
	\begin{eqnarray}
		\EE\left[f(\overline{x}^K) - f(x^*)\right] &\le& \frac{\Phi^0}{\gamma K} +2\gamma\left(\frac{\sigma^2}{n} + \gamma \frac{16L(1-p)}{p^2}\left(\frac{2(p+2)\sigma^2}{r} + p\sigma^2\right)\right) \notag
	\end{eqnarray}
	where $\Phi^0 = 2\|x^0-x^*\|^2+ \frac{512L(1-p)(2+p)(2+q)\gamma^3\sigma_0^2}{3p^2q}$.
\end{theorem}

The theorem above together with Lemma~\ref{lem:lemma2_stich} implies the following result.
\begin{corollary}\label{cor:ss_local_sgd_str_cvx}
	Let assumptions of Theorem~\ref{thm:ss_local_sgd} hold with $\mu > 0$. Then for 
	\begin{eqnarray*}
		\gamma_0 &=& \min\left\{\frac{1}{4L}, \frac{p\sqrt{3}}{32L\sqrt{2(1-p)(2+p)\left(1+\nicefrac{1}{(1-q)}\right)}}\right\},\\
		\widetilde{\Phi}^0 &=& 2\|x^0-x^*\|^2+ \frac{512L(1-p)(2+p)(2+q)\gamma_0^3\sigma_0^2}{3p^2q},\quad q = p,\\
		\gamma &=& \min\left\{\gamma_0,\frac{\ln\left(\max\left\{2, \min\left\{\nicefrac{n\widetilde{\Phi}^0\mu^2K^2}{2\sigma^2 },\nicefrac{p\widetilde{\Phi}^0\mu^3K^3}{32L(1-p)(3p+4)\sigma^2}\right\}\right\}\right)}{\mu K}\right\},\quad r = \left\lceil\frac{1}{p}\right\rceil,
	\end{eqnarray*}
	for all $K$ such that 
	\begin{eqnarray*}
		\text{either} && \frac{\ln\left(\max\left\{2, \min\left\{\nicefrac{n\widetilde{\Phi}^0\mu^2K^2}{2\sigma^2 },\nicefrac{p\widetilde{\Phi}^0\mu^3K^3}{32L(1-p)(3p+4)\sigma^2}\right\}\right\}\right)}{ K} \le p\\
		\text{or} && \gamma_0 \le \frac{\ln\left(\max\left\{2, \min\left\{\nicefrac{n\widetilde{\Phi}^0\mu^2K^2}{2\sigma^2 },\nicefrac{p\widetilde{\Phi}^0\mu^3K^3}{32L(1-p)(3p+4)\sigma^2}\right\}\right\}\right)}{\mu K}
	\end{eqnarray*}
	we have that $\EE\left[f(\overline{x}^K)-f(x^*)\right]$ is of the order
	\begin{equation}
		 \widetilde\cO\left(\frac{\widetilde{\Phi}^0}{\gamma_0}\exp\left(- \min\left\{\frac{1}{p}, \gamma_0\mu\right\} K\right) + \frac{\sigma^2}{n\mu K}  + \frac{L(1-p)\sigma^2}{p\mu^2 K^2}\right).\notag
	\end{equation}
	That is, to achieve $\EE\left[f(\overline{x}^K)-f(x^*)\right] \le \varepsilon$ in this case {\tt SS-Local-SGD} requires
	\begin{equation*}
		\widetilde{\cO}\left(\frac{L}{p\mu} + \frac{\sigma^2}{n\mu\varepsilon} + \sqrt{\frac{L(1-p)\sigma^2}{p\mu^2\varepsilon}}\right)
	\end{equation*}
	iterations/oracle calls per node (in expectation) and $\nicefrac{1}{p}$ times less communication rounds.	
\end{corollary}

Combining Theorem~\ref{thm:ss_local_sgd} and Lemma~\ref{lem:lemma_technical_cvx} we derive the following result for the convergence of {\tt SS-Local-SGD} in the case when $\mu = 0$.
\begin{corollary}
	\label{cor:ss_local_sgd_cvx}
	Let assumptions of Theorem~\ref{thm:ss_local_sgd} hold with $\mu = 0$. Then for $q = p,$ $r = \lceil\nicefrac{1}{p}\rceil$ and
	\begin{eqnarray*}
		\gamma_0 &=& \min\left\{\frac{1}{4L}, \frac{p\sqrt{3}}{32L\sqrt{2(1-p)(2+p)\left(1+\nicefrac{1}{(1-q)}\right)}}\right\},\\	
		\gamma &=& \min\left\{\gamma_0, \sqrt[3]{\frac{3p^3R_0^2}{256L(1-p)(2+p)^2\sigma_0^2}}, \sqrt{\frac{nR_0^2}{\sigma^2 K}}, \sqrt[3]{\frac{pR_0^2}{16L(1-p)(3p+4)\sigma^2 K}}\right\},
	\end{eqnarray*}
	where $R_0 = \|x^0 - x^*\|$, we have that $\EE\left[f(\overline{x}^K)-f(x^*)\right]$ is of the order
	\begin{eqnarray*}
		\cO\left(\frac{LR_0^2 + \sqrt[3]{L(1-p)\sigma_0^2R_0^4}}{pK} + \sqrt{\frac{\sigma^2R_0^2}{n K}} + \frac{\sqrt[3]{LR_0^4(1-p)\sigma^2}}{p^{\nicefrac{1}{3}}K^{\nicefrac{2}{3}}} \right).
	\end{eqnarray*}
	That is, to achieve $\EE\left[f(\overline{x}^K)-f(x^*)\right] \le \varepsilon$ in this case {\tt SS-Local-SGD} requires
	\begin{eqnarray*}
		\cO\left(\frac{LR_0^2 + \sqrt[3]{L(1-p)\sigma_0^2R_0^4}}{p\varepsilon} + \frac{\sigma^2 R_0^2}{n\varepsilon^2} + \frac{R_0^2\sqrt{L(1-p)\sigma^2}}{p^{\nicefrac{1}{2}}\varepsilon^{\nicefrac{3}{2}}}\right)
	\end{eqnarray*}
	iterations/oracle calls per node (in expectation) and $\nicefrac{1}{p}$ times less communication rounds.
\end{corollary}

\begin{remark}
	To get the rate from Tbl.~\ref{tbl:special_cases_weakly_convex} it remains to apply the following inequality:
	\begin{eqnarray*}
		\sigma_0^2 = \frac{1}{n}\sum\limits_{i=1}^n\|\nabla f_{i}(x^0)-\nabla f_{i}(x^*)\|^2 \overset{\eqref{eq:L_smoothness}}{\le} L^2 \|x^0-x^*\|^2.
	\end{eqnarray*}
\end{remark}

\subsubsection{Expected Smoothness and Arbitrary Sampling}\label{sec:loopless_local_svrg_es}
In this section we consider the same method {\tt SS-Local-SGD}, but without assumption that the stochastic gradient has a uniformly bounded variance. Instead of this we consider the same setup as in Section~\ref{sec:sgd_es}, i.e.\ we assume that each worker $i\in [n]$ at any point $x\in\R^d$ has an access to the unbiased estimator $\nabla f_{\xi_i}(x)$ of $\nabla f_i(x)$ satisfying Assumption~\ref{ass:expected_smoothness}.

\begin{lemma}\label{lem:loopless_local_svrg_es_second_moment}
	Let $f_i$ be convex and $L$-smooth for all $i\in[n]$. Let Assumption~\ref{ass:expected_smoothness} holds. Then for all $k\ge 0$
	\begin{eqnarray}
		\frac{1}{n}\sum\limits_{i=1}^n \EE_k\left[g_i^k\right] &=& \frac{1}{n}\sum\limits_{i=1}^n\nabla f_i(x_i^k), \label{eq:unbiasedness_loopless_local_svrg_es}\\
		\frac{1}{n}\sum\limits_{i=1}^n \EE\left[\|\bar{g}_i^k\|^2\right] &\le& 8L\EE\left[f(x^k)-f(x^*)\right] + 2\EE[\sigma_k^2] + 4L^2\EE[V_k],\label{eq:second_moment_loopless_local_svrg_es}\\
		\frac{1}{n}\sum\limits_{i=1}^n \EE\left[\|g_i^k-\bar{g}_i^k\|^2\right] &\le& 8\cL\EE\left[f(x^k) - f(x^*)\right] + 4\cL L\EE[V_k] + 2\sigma_*^2,\label{eq:variance_loopless_local_svrg_es}\\
		\EE\left[\left\|\frac{1}{n}\sum\limits_{i=1}^ng_i^k\right\|^2\right] &\le& 4\left(\frac{2\cL}{n}+L\right)\EE\left[f(x^k)-f(x^*)\right] + 2L\left(\frac{2\cL}{n} + L\right)\EE[V_k]\notag\\
		&&\quad + \frac{2\sigma_*^2}{n},\label{eq:second_moment_loopless_local_svrg_es_2}
	\end{eqnarray}
	where $\sigma_k^2 \eqdef \frac{1}{n}\sum\limits_{i=1}^n\left\|\nabla f_{\oxi_i^k}(y^k) - \nabla f_i(x^*)\right\|^2$ and $\sigma_*^2 \eqdef \frac{1}{n}\sum_{i=1}^n\EE_{\xi_i}\|\nabla f_{\xi_i}(x^*) - \nabla f_i(x^*)\|^2$.
\end{lemma}
\begin{proof}
	First of all, \eqref{eq:unbiasedness_loopless_local_svrg_es} follows from \eqref{eq:unbiasedness_loopless_local_svrg}. Next, using $\bar{g}_i^k = \nabla f_i(x_i^k) - \nabla f_{\oxi_i^k}(y^k) + \nabla f_{\oxi^k}(y^k)$ we get
	\begin{eqnarray*}
		\frac{1}{n}\sum\limits_{i=1}^n\EE\left[\|\bar{g}_i^k\|^2\right] &\overset{\eqref{eq:a_b_norm_squared}}{\le}& \frac{2}{n}\sum\limits_{i=1}^n\EE\left[\|\nabla f_{i}(x_i^k) - \nabla f_{i}(x^*)\|^2\right]\\
		&&\quad + \frac{2}{n}\sum\limits_{i=1}^n\EE\left[\left\|\nabla f_{\oxi_i^k}(y^k) - \nabla f_i(x^*) - (\nabla f_{\oxi^k}(y^k) - \nabla f(x^*))\right\|^2\right]\\
		&\overset{\eqref{eq:L_smoothness_cor},\eqref{eq:variance_decomposition}}{\le}& \frac{4L}{n}\sum\limits_{i=1}^n\EE\left[D_{f_i}(x_i^k,x^*)\right] + \frac{2}{n}\sum\limits_{i=1}^n\EE\left[\|\nabla f_{\oxi_i^k}(y^k)- \nabla f_i(x^*)\|^2\right]\\
		&\overset{\eqref{eq:poiouhnkj}}{\le}& 8L\EE\left[f(x^k)-f(x^*)\right] + 2\EE[\sigma_k^2] + 4L^2\EE[V_k]
	\end{eqnarray*}
	and	
	\begin{eqnarray}
		\frac{1}{n}\sum\limits_{i=1}^n\EE\left[\|g_i^k-\bar{g}_i^k\|^2\right] &=& \frac{1}{n}\sum\limits_{i=1}^n\EE\left[\|\nabla f_{\xi_i^k}(x_i^k) - \nabla f_i(x_i^k)\|^2\right]\notag\\
		&\overset{\eqref{eq:variance_decomposition}}{\le}& \frac{1}{n}\sum\limits_{i=1}^n\EE\left[\|\nabla f_{\xi_i^k}(x_i^k) - \nabla f_i(x^*)\|^2\right]\notag\\
		&\overset{\eqref{eq:a_b_norm_squared}}{\le}& \frac{2}{n}\sum\limits_{i=1}^n\EE\left[\|\nabla f_{\xi_i^k}(x_i^k) - \nabla f_{\xi_i^k}(x^*)\|^2\right]\\
		&&\quad + \frac{2}{n}\sum\limits_{i=1}^n\EE\left[\|\nabla f_{\xi_i^k}(x^*) - \nabla f_{i}(x^*)\|^2\right]\notag \\
		&\overset{\eqref{eq:expected_smoothness_1}}{\le}& \frac{4\cL}{n}\sum\limits_{i=1}^n\EE\left[D_{f_i}(x_i^k,x^*)\right] + 2\sigma_*^2\notag\\
		&\overset{\eqref{eq:poiouhnkj}}{\le}& 8\cL\EE\left[f(x^k) - f(x^*)\right] + 4\cL L\EE[V_k] + 2\sigma_*^2. \label{eq:bshjbdhsbdhbucsb}
	\end{eqnarray}
	Finally, we use independence of $\xi_1^k,\ldots,\xi_n^k$ and derive
	\begin{eqnarray*}
		\EE\left[\left\|\frac{1}{n}\sum\limits_{i=1}^n g_i^k\right\|^2\right] &=& \EE\left[\left\|\frac{1}{n}\sum\limits_{i=1}^n \nabla f_{\xi_i^k}(x_i^k)\right\|^2\right]\\
		&\overset{\eqref{eq:tower_property},\eqref{eq:variance_decomposition}}{=}& \EE\left[\left\|\frac{1}{n}\sum\limits_{i=1}^n (\nabla f_{\xi_i^k}(x_i^k)-\nabla f_i(x_i^k))\right\|^2\right] + \EE\left[\left\|\frac{1}{n}\sum\limits_{i=1}^n \nabla f_{i}(x_i^k)\right\|^2\right]\\
		&=& \frac{1}{n^2}\sum\limits_{i=1}^n\EE\left[\|\nabla f_{\xi_i^k}(x_i^k) - \nabla f_i(x_i^k)\|^2\right] + \EE\left[\left\|\frac{1}{n}\sum\limits_{i=1}^n \nabla f_{i}(x_i^k)\right\|^2\right]\\
		&\overset{\eqref{eq:bshjbdhsbdhbucsb},\eqref{eq:vdgasvgda}}{\le}& 4\left(\frac{2\cL}{n}+L\right)\EE\left[f(x^k)-f(x^*)\right] + 2L\left(\frac{2\cL}{n} + L\right)\EE[V_k] + \frac{2\sigma_*^2}{n}
	\end{eqnarray*}
	which finishes the proof.
\end{proof}

\begin{lemma}\label{lem:loopless_local_svrg_es_sigma_k_bound}
	Let $f_i$ be convex and $L$-smooth for all $i\in[n]$ and Assumption~\ref{ass:expected_smoothness} holds. Then for all $k\ge 0$
	\begin{eqnarray}
		\EE\left[\sigma_{k+1}^2\right] &\le& (1-q)\EE\left[\sigma_k^2\right] + 2q\left(\frac{2\cL}{r} + L\right)\EE\left[f(x^k) - f(x^*)\right] + \frac{2q\sigma_*^2}{r}\label{eq:loopless_local_svrg_es_sigma_k_bound}
	\end{eqnarray}
	where $\sigma_k^2 \eqdef \frac{1}{n}\sum\limits_{i=1}^n\left\|\nabla f_{\oxi_i^k}(y^k) - \nabla f_i(x^*)\right\|^2$ and $\sigma_*^2 \eqdef \frac{1}{n}\sum_{i=1}^n\EE_{\xi_i}\|\nabla f_{\xi_i}(x^*) - \nabla f_i(x^*)\|^2$.
\end{lemma}
\begin{proof}
	By definition of $y^{k+1}$ we have
	\begin{eqnarray*}
		\EE\left[\sigma_{k+1}^2\mid x_1^k,\ldots, x_n^k\right] &=& \frac{1-q}{n}\sum\limits_{i=1}^n\|\nabla f_{\oxi_i^k}(y^k) - \nabla f_i(x^*)\|^2\\
		&&\quad + \frac{q}{n}\sum\limits_{i=1}^n\EE_{\oxi_i^{k+1}}\left[\|\nabla f_{\oxi_i^{k+1}}(x^k) - \nabla f_i(x^*)\|^2\right]\\
		&\overset{\eqref{eq:variance_decomposition}}{=}& (1-q)\sigma_k^2 + \frac{q}{n}\sum\limits_{i=1}^n\|\nabla f_{i}(x^k) - \nabla f_i(x^*)\|^2\\
		&&\quad + \frac{q}{n}\sum\limits_{i=1}^n\EE_{\oxi_i^{k+1}}\left[\|\nabla f_{\oxi_i^{k+1}}(x^k) - \nabla f_i(x^k)\|^2\right].
	\end{eqnarray*}
	Next, we use independence of $\oxi_{i,1}^{k+1}, \oxi_{i,2}^{k+1},\ldots, \oxi_{i,r}^{k+1}$ for all $i\in [n]$ and derive
	\begin{eqnarray*}
		\EE\left[\sigma_{k+1}^2\mid x_1^k,\ldots, x_n^k\right] &=& (1-q)\sigma_k^2 + \frac{q}{n}\sum\limits_{i=1}^n\|\nabla f_{i}(x^k) - \nabla f_i(x^*)\|^2\\
		&&\quad + \frac{q}{nr^2}\sum\limits_{i=1}^n\sum\limits_{j=1}^r\EE_{\oxi_{i,j}^{k+1}}\left[\|\nabla f_{\oxi_{i,j}^{k+1}}(x^k) - \nabla f_i(x^k)\|^2\right]\\
		&\overset{\eqref{eq:L_smoothness_cor},\eqref{eq:variance_decomposition}}{\le}& (1-q)\sigma_k^2 + 2Lq\left(f(x^k)-f(x^*)\right)\\
		&&\quad + \frac{q}{nr^2}\sum\limits_{i=1}^n\sum\limits_{j=1}^r\EE_{\oxi_{i,j}^{k+1}}\left[\|\nabla f_{\oxi_{i,j}^{k+1}}(x^k) - \nabla f_i(x^*)\|^2\right]\\
		&\overset{\eqref{eq:a_b_norm_squared}}{\le}& (1-q)\sigma_k^2 + 2Lq\left(f(x^k)-f(x^*)\right)\\
		&&\quad + \frac{2q}{nr^2}\sum\limits_{i=1}^n\sum\limits_{j=1}^r\EE_{\oxi_{i,j}^{k+1}}\left[\|\nabla f_{\oxi_{i,j}^{k+1}}(x^k) - \nabla f_{\oxi_{i,j}^{k+1}}(x^*)\|^2\right]\\
		&&\quad + \frac{2q}{nr^2}\sum\limits_{i=1}^n\sum\limits_{j=1}^r\EE_{\oxi_{i,j}^{k+1}}\left[\|\nabla f_{\oxi_{i,j}^{k+1}}(x^*) - \nabla f_{i}(x^*)\|^2\right]\\
		&\overset{\eqref{eq:expected_smoothness_1}}{\le}& (1-q)\sigma_k^2 + 2q\left(\frac{2\cL}{r} + L\right)\left(f(x^k) - f(x^*)\right) + \frac{2q\sigma_*^2}{r}.
	\end{eqnarray*}
	Taking the full mathematical expectation on both sides of previous inequality and using the tower property \eqref{eq:tower_property} we get the result.
\end{proof}

Using Corollary~\ref{cor:rand_loop} we obtain the following theorem.
\begin{theorem}\label{thm:ss_local_sgd_es}
	Assume that $f_i(x)$ is $\mu$-strongly convex and $L$-smooth for every $i\in[n]$. Let Assumption~\ref{ass:expected_smoothness} holds. Then {\tt SS-Local-SGD} satisfies Assumption~\ref{ass:hetero_second_moment} with
	\begin{gather*}
		\tA = 4L,\quad \hA = 4\cL,\quad \tB = 2,\quad \hB = 0,\quad \tF = 4L^2,\quad \hF = 4\cL L, \quad \tD_1 = 0,\quad B' = 0,\\
		\hD_1 = 2\sigma_*^2 = \frac{2}{n}\sum\limits_{i=1}^n\EE_{\xi_i}\|\nabla f_{\xi_i}(x^*) - \nabla f_i(x^*)\|^2,\quad A' = 2\left(\frac{2\cL}{n} + L\right),\quad F' = 2L\left(\frac{2\cL}{n} + L\right), \\
		D_1' = \frac{2\sigma_*^2}{n},\quad \sigma_k^2 = \frac{1}{n}\sum\limits_{i=1}^n\left\|\nabla f_{\oxi_i^k}(y^k) - \nabla f_i(x^*)\right\|^2,\quad \rho = q,\quad C = q\left(\frac{2\cL}{r} + L\right),\quad G = 0,\\
		D_2 = \frac{2q\sigma_*^2}{r},,\quad H = \frac{128(1-p)(2+p)(2+q)\gamma^2}{3p^2q},\quad D_3 = \frac{8(1-p)}{p^2}\left(2p\sigma_*^2+\frac{32(2+p)\sigma_*^2}{3r}\right)
	\end{gather*}
	under assumption that
	\begin{eqnarray*}
		\gamma &\le& \min\left\{\frac{1}{4\left(\frac{2\cL}{n}+L\right)}, \frac{p\sqrt{3}}{32\sqrt{2L(1-p)\left((2+p)L + p\cL + \frac{(2+p)\left(\nicefrac{2\cL}{r} + L\right)}{(1-q)}\right)}}\right\}.
	\end{eqnarray*}
	Moreover, for $\mu > 0$ we have
	\begin{eqnarray}
		\EE\left[f(\overline{x}^K) - f(x^*)\right] &\le& \left(1 - \min\left\{\gamma\mu,\frac{q}{4}\right\}\right)^K\frac{\Phi^0}{\gamma}\notag\\
		&&\quad +2\gamma\left(\frac{2\sigma_*^2}{n} + \gamma \frac{16L(1-p)}{p^2}\left(2p\sigma_*^2+\frac{32(2+p)\sigma_*^2}{3r}\right)\right) \notag
	\end{eqnarray}
	and when $\mu = 0$ we have
	\begin{eqnarray}
		\EE\left[f(\overline{x}^K) - f(x^*)\right] &\le& \frac{\Phi^0}{\gamma K} +2\gamma\left(\frac{2\sigma_*^2}{n} + \gamma \frac{16L(1-p)}{p^2}\left(2p\sigma_*^2+\frac{32(2+p)\sigma_*^2}{3r}\right)\right) \notag
	\end{eqnarray}
	where $\Phi^0 = 2\|x^0-x^*\|^2+ \frac{512L(1-p)(2+p)(2+q)\gamma^3\EE[\sigma_0^2]}{3p^2q}$.
\end{theorem}

The theorem above together with Lemma~\ref{lem:lemma2_stich} implies the following result.
\begin{corollary}\label{cor:ss_local_sgd_str_cvx_es}
	Let assumptions of Theorem~\ref{thm:ss_local_sgd_es} hold with $\mu > 0$. Then for 
	\begin{eqnarray*}
		\gamma_0 &=& \min\left\{\frac{1}{4\left(\frac{2\cL}{n}+L\right)}, \frac{p\sqrt{3}}{32\sqrt{2L(1-p)\left((2+p)L + p\cL + \frac{(2+p)\left(\nicefrac{2\cL}{r} + L\right)}{(1-q)}\right)}}\right\},\\
		\widetilde{\Phi}^0 &=& 2\|x^0-x^*\|^2+ \frac{512L(1-p)(2+p)(2+q)\gamma_0^3\EE[\sigma_0^2]}{p^2q},\quad q = p,\quad r = \left\lceil\frac{1}{p}\right\rceil,\\
		\gamma &=& \min\left\{\gamma_0,\frac{\ln\left(\max\left\{2, \min\left\{\nicefrac{n\widetilde{\Phi}^0\mu^2K^2}{4\sigma_*^2 },\nicefrac{p\widetilde{\Phi}^0\mu^3K^3}{64L(1-p)(1+\nicefrac{32(2+p)}{3})\sigma_*^2}\right\}\right\}\right)}{\mu K}\right\},
	\end{eqnarray*}
	for all $K$ such that 
	\begin{eqnarray*}
		\text{either} && \frac{\ln\left(\max\left\{2, \min\left\{\nicefrac{n\widetilde{\Phi}^0\mu^2K^2}{4\sigma_*^2 },\nicefrac{p\widetilde{\Phi}^0\mu^3K^3}{64L(1-p)(1+\nicefrac{32(2+p)}{3})\sigma_*^2}\right\}\right\}\right)}{ K} \le p\\
		\text{or} && \gamma_0 \le \frac{\ln\left(\max\left\{2, \min\left\{\nicefrac{n\widetilde{\Phi}^0\mu^2K^2}{4\sigma_*^2 },\nicefrac{p\widetilde{\Phi}^0\mu^3K^3}{64L(1-p)(1+\nicefrac{32(2+p)}{3})\sigma_*^2}\right\}\right\}\right)}{\mu K}
	\end{eqnarray*}
	we have that $\EE\left[f(\overline{x}^K)-f(x^*)\right]$ is of the order
	\begin{equation}
		 \widetilde\cO\left(\frac{\widetilde{\Phi}^0}{\gamma_0}\exp\left(- \min\left\{\frac{1}{p}, \gamma_0\mu\right\} K\right) + \frac{\sigma_*^2}{n\mu K}  + \frac{L(1-p)\sigma_*^2}{p\mu^2 K^2}\right).\notag
	\end{equation}
	That is, to achieve $\EE\left[f(\overline{x}^K)-f(x^*)\right] \le \varepsilon$ in this case {\tt SS-Local-SGD} requires
	\begin{equation*}
		\widetilde{\cO}\left(\frac{L}{p\mu} + \frac{\cL}{n\mu} + \frac{\sqrt{\cL L(1-p)}}{\sqrt{p}\mu} + \frac{\sigma_*^2}{n\mu\varepsilon} + \sqrt{\frac{L(1-p)\sigma_*^2}{p\mu^2\varepsilon}}\right)
	\end{equation*}
	iterations/oracle calls per node (in expectation) and $\nicefrac{1}{p}$ times less communication rounds.	
\end{corollary}

Combining Theorem~\ref{thm:ss_local_sgd_es} and Lemma~\ref{lem:lemma_technical_cvx} we derive the following result for the convergence of {\tt SS-Local-SGD} in the case when $\mu = 0$.
\begin{corollary}
	\label{cor:ss_local_sgd_cvx_es}
	Let assumptions of Theorem~\ref{thm:ss_local_sgd_es} hold with $\mu = 0$. Then for $q = p,$ $r = \lceil\nicefrac{1}{p}\rceil$ and
	\begin{eqnarray*}
		\gamma_0 &=& \min\left\{\frac{1}{4\left(\frac{2\cL}{n}+L\right)}, \frac{p\sqrt{3}}{32\sqrt{2L(1-p)\left((2+p)L + p\cL + \frac{(2+p)\left(\nicefrac{2\cL}{r} + L\right)}{(1-q)}\right)}}\right\},\\	
		\gamma &=& \min\left\{\gamma_0, \sqrt[3]{\frac{p^3R_0^2}{256L(1-p)(2+p)^2\EE[\sigma_0^2]}}, \sqrt{\frac{nR_0^2}{2\sigma_*^2 K}}, \sqrt[3]{\frac{pR_0^2}{32L(1-p)\left(1+\nicefrac{32(2+p)}{3}\right)\sigma_*^2 K}}\right\},
	\end{eqnarray*}
	where $R_0 = \|x^0 - x^*\|$, we have that $\EE\left[f(\overline{x}^K)-f(x^*)\right]$ is of the order
	\begin{eqnarray*}
		\cO\left(\frac{\left(L+\nicefrac{p\cL}{n} + \sqrt{p(1-p)\cL L}\right)R_0^2 + \sqrt[3]{L(1-p)\EE[\sigma_0^2]R_0^4}}{pK} + \sqrt{\frac{\sigma_*^2R_0^2}{n K}} + \frac{\sqrt[3]{LR_0^4(1-p)\sigma_*^2}}{p^{\nicefrac{1}{3}}K^{\nicefrac{2}{3}}} \right).
	\end{eqnarray*}
	That is, to achieve $\EE\left[f(\overline{x}^K)-f(x^*)\right] \le \varepsilon$ in this case {\tt SS-Local-SGD} requires
	\begin{eqnarray*}
		\cO\left(\frac{\left(L+\nicefrac{p\cL}{n} + \sqrt{p(1-p)\cL L}\right)R_0^2 + \sqrt[3]{L(1-p)\EE[\sigma_0^2]R_0^4}}{p\varepsilon} + \frac{\sigma_*^2 R_0^2}{n\varepsilon^2} + \frac{R_0^2\sqrt{L(1-p)\sigma_*^2}}{p^{\nicefrac{1}{2}}\varepsilon^{\nicefrac{3}{2}}}\right)
	\end{eqnarray*}
	iterations/oracle calls per node (in expectation) and $\nicefrac{1}{p}$ times less communication rounds.
\end{corollary}

\begin{remark}
	To get the rate from Tbl.~\ref{tbl:special_cases_weakly_convex} it remains to apply the following inequality:
	\begin{eqnarray*}
		\EE[\sigma_0^2] &=& \frac{1}{n}\sum\limits_{i=1}^n\EE_{\oxi_i^0}\left[\|\nabla f_{\oxi_i^0}(x^0)-\nabla f_{i}(x^*)\|^2\right]\\
		&\overset{\eqref{eq:variance_decomposition}}{=}& \frac{1}{n}\sum\limits_{i=1}^n\|\nabla f_{i}(x^0) - \nabla f_i(x^*)\|^2 + \frac{1}{n}\sum\limits_{i=1}^n\EE_{\oxi_i^0}\left[\|\nabla f_{\oxi_i^0}(x^0) - \nabla f_i(x^0)\|^2\right]\\
		&\overset{\eqref{eq:L_smoothness_cor}}{\le}& 2L(f(x^0)-f(x^*)) + \frac{1}{nr^2}\sum\limits_{i=1}^n\sum\limits_{j=1}^r\EE_{\oxi_{i,j}^0}\left[\|\nabla f_{\oxi_{i,j}^0}(x^0) - \nabla f_i(x^0)\|^2\right]\\
		&\overset{\eqref{eq:variance_decomposition}}{\le}& 2L(f(x^0)-f(x^*)) + \frac{1}{nr}\sum\limits_{i=1}^n\EE_{\xi_i}\left[\|\nabla f_{\xi_i}(x^0)-\nabla f_i(x^*)\|^2\right]\\
		&\overset{\eqref{eq:a_b_norm_squared}}{\le}& 2L(f(x^0)-f(x^*)) + \frac{2}{nr}\sum\limits_{i=1}^n\EE_{\xi_i}\left[\|\nabla f_{\xi_i}(x^0)-\nabla f_{\xi_i}(x^*)\|^2\right] \\
		&&\quad + \frac{2}{nr}\sum\limits_{i=1}^n\EE_{\xi_i}\left[\|\nabla f_{\xi_i}(x^*)-\nabla f_{i}(x^*)\|^2\right]\\
		&\overset{r = \lceil\nicefrac{1}{p}\rceil,\eqref{eq:expected_smoothness_1}}{\le}& 2\left(L + 2p\cL\right)(f(x^0) - f(x^*)) + 2p\sigma_*^2.
	\end{eqnarray*}
\end{remark}

\subsection{{\tt S*-Local-SGD*}} \label{sec:S*-Local-SGD*}
 In this section we present doubly idealized algorithm for solving problem \eqref{eq:main_problem}+\eqref{eq:f_i_sum}. Specifically, we choose $b_i^k$ to the optimal shift $\nabla f_i(x^*)$ as per Case II, while $a_i^k$ is selected as { \tt SGD-star} gradient estimator~\citep{gorbunov2019unified}, i.e., 
\[
a^{k}_i = \nabla f_{i,j_i}(x_i^k) - \nabla f_{i,j_i}(x^*) +  \nabla f_{i}(x^*), \qquad b^{k}_i = \nabla f_{i}(x^*).
\]

 Note that now $a_i^k$ serves as an ambitious target for the local variance reduced estimators, while $b_i^k$ serves as an ambitious goal for the local shift. The resulting instance of~\eqref{eq:local_sgd_def} is presented as Algorithm~\ref{alg:local_sgd_star_star} and called Star-Shifted {\tt Local-SGD-star} ({\tt S*-Local-SGD*}).

\begin{algorithm}[h]
   \caption{{\tt S*-Local-SGD*}}\label{alg:local_sgd_star_star}
\begin{algorithmic}[1]
   \Require learning rate $\gamma>0$, initial vector $x^0 \in \R^d$, communication period $\tau \ge 1$
	\For{$k=0,1,\dotsc$}
       	\For{$i=1,\dotsc,n$ in parallel}
            \State Set $ g^{k}_i = \nabla f_{i,j_i}(x_i^k) - \nabla f_{i,j_i}(x_*)$ where $1\leq j_i\leq m$ is sampled independently from all nodes
            \If{$k+1 \mod \tau = 0$}
            \State $x_i^{k+1} = x^{k+1} = \frac{1}{n}\sum\limits_{i=1}^n\left(x_i^k - \gamma g_i^k\right)$ \Comment{averaging}
            \Else
            \State $x_i^{k+1} = x_i^k - \gamma g_i^k$ \Comment{local update}
            \EndIf
        \EndFor
   \EndFor
\end{algorithmic}
\end{algorithm}

Let us next provide the details on the convergence rate. In order to do so, let us identify the parameters of Assumption~\ref{ass:sigma_k_original}.

\begin{lemma}\label{lem:s*_local_sgd*_lemma}
	Let $f_i$ be convex and $L$-smooth and $f_{i,j}$ be convex and $\max L_{ij}$-smooth for all $i\in[n]$, $j\in[m]$. Then for all $k\ge 0$
	\begin{eqnarray}
		\frac{1}{n}\sum\limits_{i=1}^n \EE_k\left[g_i^k\right] &=& \frac{1}{n}\sum\limits_{i=1}^n\nabla f_i(x_i^k), \label{eq:unbiasedness_s*_local_sgd*}\\
		\frac{1}{n}\sum\limits_{i=1}^n \EE\left[\|\bar{g}_i^k\|^2\right] &\le& 4L\EE\left[f(x^k)-f(x^*)\right] + 2L^2\EE[V_k],\label{eq:second_moment_s*_local_sgd*}\\
		\frac{1}{n}\sum\limits_{i=1}^n \EE\left[\|g_i^k-\bar{g}_i^k\|^2\right] &\le& 4\max L_{ij}\EE\left[f(x^k) - f(x^*)\right] + 2L\max L_{ij}\EE[V_k],\label{eq:variance_s*_local_sgd*}\\
		\EE\left[\left\|\frac{1}{n}\sum\limits_{i=1}^ng_i^k\right\|^2\right] &\le& 4\left(\frac{\max L_{ij}}{n}+L\right)\EE\left[f(x^k)-f(x^*)\right]\notag\\
		&&\quad + 2L\left(\frac{\max L_{ij}}{n} + L\right)\EE[V_k].\label{eq:second_moment_s*_local_sgd*_2}
	\end{eqnarray}
\end{lemma}
\begin{proof}
	First of all,
	\begin{eqnarray*}
		\frac{1}{n}\sum\limits_{i=1}^n \EE_k\left[g_i^k\right] &=& \frac{1}{nm}\sum\limits_{i=1}^n\sum\limits_{j=1}^m \left(\nabla f_{i,j}(x_i^k) - \nabla f_{i,j}(x^*)\right) = \frac{1}{n}\sum\limits_{i=1}^n\nabla f_i(x_i^k)
	\end{eqnarray*}
	and, in particular, $\bar{g}_i^k = \EE_k\left[g_i^k\right] = \nabla f_i(x_i^k) - \nabla f_i(x^*)$. Using this we derive
	\begin{eqnarray*}
		\frac{1}{n}\sum\limits_{i=1}^n \EE\left[\|\bar{g}_i^k\|^2\right] &=& \frac{1}{n}\sum\limits_{i=1}^n \EE\left[\|\nabla f_i(x_i^k) - \nabla f_i(x^*)\|^2\right]\\
		&\overset{\eqref{eq:L_smoothness_cor}}{\le}& \frac{2L}{n}\sum\limits_{i=1}^n\EE\left[D_{f_i}(x_i^k,x^*)\right] \overset{\eqref{eq:poiouhnkj}}{\le} 4L\EE\left[f(x^k)-f(x^*)\right] + 2L^2\EE[V_k]
	\end{eqnarray*}
	and
	\begin{eqnarray}
		\frac{1}{n}\sum\limits_{i=1}^n \EE\left[\|g_i^k-\bar{g}_i^k\|^2\right] &\overset{\eqref{eq:variance_decomposition}}{\le}& \frac{1}{n}\sum\limits_{i=1}^n \EE\left[\|g_i^k\|^2\right]\notag\\
		&=& \frac{1}{nm}\sum\limits_{i=1}^n\sum\limits_{j=1}^m \|\nabla f_{i,j}(x_i^k) - \nabla f_{i,j}(x^*)\|^2\notag \\
		&\overset{\eqref{eq:L_smoothness_cor}}{\le}& \frac{2\max L_{ij}}{n}\sum\limits_{i=1}^n\EE\left[D_{f_i}(x_i^k,x^*)\right]\notag\\
		&\overset{\eqref{eq:poiouhnkj}}{\le}& 4\max L_{ij}\EE\left[f(x^k)-f(x^*)\right] + 2L\max L_{ij}\EE[V_k].\label{eq:hbdsujbcvshvbcbu}
	\end{eqnarray}
	Finally, due to the independence of $j_1, j_2, \ldots, j_n$ we have
	\begin{eqnarray*}
		\EE\left[\left\|\frac{1}{n}\sum\limits_{i=1}^ng_i^k\right\|^2\right] &\overset{\eqref{eq:variance_decomposition},\eqref{eq:tower_property}}{=}& \EE\left[\left\|\frac{1}{n}\sum\limits_{i=1}^n\left(\nabla f_{i,j_i}(x_i^k) - \nabla f_{i,j_i}(x_*) - (\nabla f_i(x_i^k) - \nabla f_i(x^*))\right)\right\|^2\right]\\
		&&\quad + \EE\left[\left\|\frac{1}{n}\sum\limits_{i=1}^n\left(\nabla f_i(x_i^k) - \nabla f_i(x^*)\right)\right\|^2\right]\\
		&=& \frac{1}{n^2}\sum\limits_{i=1}^n\EE\left[\|\nabla f_{i,j_i}(x_i^k) - \nabla f_{i,j_i}(x_*) - (\nabla f_i(x_i^k) - \nabla f_i(x^*))\|^2\right]\\
		&&\quad + \EE\left[\left\|\frac{1}{n}\sum\limits_{i=1}^n\nabla f_i(x_i^k)\right\|^2\right]\\
		&\overset{\eqref{eq:variance_decomposition}}{\le}& \frac{1}{n^2m}\sum\limits_{i=1}^n\sum\limits_{j=1}^m \|\nabla f_{i,j}(x_i^k) - \nabla f_{i,j}(x^*)\|^2 + \EE\left[\left\|\frac{1}{n}\sum\limits_{i=1}^n\nabla f_i(x_i^k)\right\|^2\right]\\
		&\overset{\eqref{eq:hbdsujbcvshvbcbu},\eqref{eq:poiouhnkj}}{\le}& 4\left(\frac{\max L_{ij}}{n}+L\right)\EE\left[f(x^k)-f(x^*)\right] + 2L\left(\frac{\max L_{ij}}{n} + L\right)\EE[V_k].
	\end{eqnarray*}
\end{proof}

Using Corollary~\ref{cor:const_loop} we obtain the following theorem.
\begin{theorem}\label{thm:s*_local_sgd*}
	Assume that $f_i(x)$ is $\mu$-strongly convex and $L$-smooth and $f_{i,j}$ is convex and $\max L_{ij}$-smooth for every $i\in[n]$, $j\in[m]$. Then {\tt S*-Local-SGD*} satisfies Assumption~\ref{ass:hetero_second_moment} with
	\begin{gather*}
		\tA = 2L,\quad \hA = 2\max L_{ij},\quad \tB = \hB = 0,\quad \tF = 2L^2,\quad \hF = 2L\max L_{ij}, \quad \tD_1 = \hD_1 = 0,\\
		A' = 2\left(\frac{\max L_{ij}}{n} + L\right),\quad B' = 0,\quad F' = 2L\left(\frac{\max L_{ij}}{n} + L\right), \\
		D_1' = 0,\quad \sigma_k^2 \equiv 0,\quad \rho = 1,\quad C = 0,\quad G = 0,\quad D_2 = 0,\quad H = 0,\quad D_3 = 0
	\end{gather*}
	under assumption that
	\begin{eqnarray*}
		\gamma &\le& \min\left\{\frac{1}{4\left(\frac{\max L_{ij}}{n}+L\right)}, \frac{1}{8\sqrt{eL(\tau-1)\left(L(\tau-1)+\max L_{ij}\right)}}\right\}.
	\end{eqnarray*}
	Moreover, for $\mu > 0$ we have
	\begin{eqnarray}
		\EE\left[f(\overline{x}^K) - f(x^*)\right] &\le& \left(1 - \gamma\mu\right)^K\frac{2\|x^0-x^*\|^2}{\gamma} \notag
	\end{eqnarray}
	and when $\mu = 0$ we have
	\begin{eqnarray}
		\EE\left[f(\overline{x}^K) - f(x^*)\right] &\le& \frac{2\|x^0-x^*\|^2}{\gamma K}. \notag
	\end{eqnarray}
\end{theorem}

The theorem above together with Lemma~\ref{lem:lemma2_stich} implies the following result.
\begin{corollary}\label{cor:s*_local_sgd*_str_cvx}
	Let assumptions of Theorem~\ref{thm:s*_local_sgd*} hold with $\mu > 0$. Then for 
	\begin{eqnarray*}
		\gamma &=& \min\left\{\frac{1}{4\left(\frac{\max L_{ij}}{n}+L\right)}, \frac{1}{8\sqrt{eL(\tau-1)\left(L(\tau-1)+\max L_{ij}\right)}}\right\}
	\end{eqnarray*}
	and for all $K\ge 1$ we have $\EE\left[f(\overline{x}^K)-f(x^*)\right]$ of order
	\begin{equation}
		\cO\left(\left(L\tau + \frac{\max L_{ij}}{n} + \sqrt{(\tau-1)L\max L_{ij}}\right)R_0^2\exp\left(- \frac{\mu}{L\tau + \frac{\max L_{ij}}{n} + \sqrt{(\tau-1)L\max L_{ij}}}K\right)\right)\notag
	\end{equation}
	with $R_0 = \|x^0 - x^*\|$. That is, to achieve $\EE\left[f(\overline{x}^K)-f(x^*)\right] \le \varepsilon$ in this case {\tt S*-Local-SGD*} requires
	\begin{equation*}
		\cO\left(\left(\frac{L\tau}{\mu} + \frac{\max L_{ij}}{n\mu} + \frac{\sqrt{(\tau-1) L\max L_{ij}}}{\mu}\right)\log\frac{\left(L\tau + \frac{\max L_{ij}}{n} + \sqrt{(\tau-1)L\max L_{ij}}\right)R_0^2}{\varepsilon}\right)
	\end{equation*}
	iterations/oracle calls per node and $\tau$ times less communication rounds.	
\end{corollary}

Next, we derive the following result for the convergence of {\tt S*-Local-SGD*} in the case when $\mu = 0$.
\begin{corollary}
	\label{cor:s*_local_sgd*_cvx}
	Let assumptions of Theorem~\ref{thm:s*_local_sgd*} hold with $\mu = 0$. Then for
	\begin{eqnarray*}
		\gamma &=& \min\left\{\frac{1}{4\left(\frac{\max L_{ij}}{n}+L\right)}, \frac{1}{8\sqrt{eL(\tau-1)\left(L(\tau-1)+\max L_{ij}\right)}}\right\},
	\end{eqnarray*}
	we have that $\EE\left[f(\overline{x}^K)-f(x^*)\right]$ is of the order
	\begin{eqnarray*}
		\cO\left(\frac{\left(L\tau + \nicefrac{\max L_{ij}}{n} + \sqrt{(\tau-1)L\max L_{ij}}\right)R_0^2}{K}\right),
	\end{eqnarray*}
	where $R_0 = \|x^0 - x^*\|$. That is, to achieve $\EE\left[f(\overline{x}^K)-f(x^*)\right] \le \varepsilon$ in this case {\tt S*-Local-SGD*} requires
	\begin{eqnarray*}
		\cO\left(\frac{\left(L\tau + \nicefrac{\max L_{ij}}{n} + \sqrt{(\tau-1)L\max L_{ij}}\right)R_0^2}{\varepsilon}\right)
	\end{eqnarray*}
	iterations/oracle calls per node and $\tau$ times less communication rounds.
\end{corollary}

%
%
%
%
%

\subsection{{\tt S-Local-SVRG}}\label{sec:loopless_local_svrg_fs}

\begin{algorithm}[h]
   \caption{Shifted Local {\tt SVRG} ({\tt S-Local-SVRG}) for minimizing local finite sums}\label{alg:l_local_svrg_fs}
\begin{algorithmic}[1]
   \Require learning rate $\gamma>0$, initial vector $x^0 \in \R^d$, probability of communication $p\in(0,1]$, probability of local full gradient computation $q\in(0,1]$, initialization $y^0 = x^0$
	\For{$k=0,1,\dotsc$}
       	\For{$i=1,\dotsc,n$ in parallel}
       		\State Choose $j_i$ uniformly at random from $[m]$
            \State $g_i^k = \nabla f_{i,j_i}(x_i^k) - \nabla f_{i,j_i}(y^k) + \nabla f(y^k)$
            \State $x_i^{k+1} = \begin{cases}x^{k+1},&\text{w.p. } p,\\
            x_i^k - \gamma g_i^k,& \text{w.p. } 1 - p, \end{cases}$ where $x^{k+1} = \frac{1}{n}\sum\limits_{i=1}^n(x_i^k - \gamma g_i^k)$
            \State $y^{k+1} = \begin{cases}x^k,&\text{w.p. } q,\\
            y^k,& \text{w.p. } 1 - q \end{cases}$
        \EndFor
   \EndFor
\end{algorithmic}
\end{algorithm}

In this section we are interested in problem \eqref{eq:main_problem}+\eqref{eq:f_i_sum}. To solve this problem we propose a new method called Shifted {\tt Local-SVRG} ({\tt S-Local-SVRG}, see Algorithm~\ref{alg:l_local_svrg_fs}).

We note that our analysis works even when updates in lines $5$,$6$ are not independent. Moreover, in order for {\tt S-Local-SVRG} to be efficient, we shall require $q\leq p$.

\begin{remark}
Unlike all other special cases, the rate of {\tt S-Local-SVRG} can not be directly obtained from the theory of the local stochastic solver described in Section~\ref{sec:local_solver}. Specifically, we construct the sequence $l_i^k$ using $y^k$ in contrast to $x_i^k$ used in Section~\ref{sec:local_solver}. While we could construct $l_i^k$ from the local iterate sequences, setting it as the virtual iterates yields a tighter rate. We remark that such a choice is rather poor in general; we can implement it efficiently thanks to the specific structure of {\tt S-Local-SVRG}. 
\end{remark}

\begin{lemma}\label{lem:loopless_local_svrg_fs_second_moment}
	Let $f_{i}$ be convex and $L$-smooth and $f_{i,j}$ be convex and $\max L_{ij}$-smooth for all $i\in[n]$, $j\in[m]$. Then for all $k\ge 0$
	\begin{eqnarray}
		\frac{1}{n}\sum\limits_{i=1}^n \EE_k\left[g_i^k\right] &=& \frac{1}{n}\sum\limits_{i=1}^n\nabla f_i(x_i^k), \label{eq:unbiasedness_loopless_local_svrg_fs}\\
		\frac{1}{n}\sum\limits_{i=1}^n \EE\left[\|\bar{g}_i^k\|^2\right] &\le&8L\EE\left[f(x^k)-f(x^*)\right] + 2\EE[\sigma_k^2] + 4L^2\EE[V_k],\label{eq:second_moment_loopless_local_svrg_fs}\\
		\frac{1}{n}\sum\limits_{i=1}^n \EE\left[\|g_i^k-\bar{g}_i^k\|^2\right] &\le& 8\max L_{ij}\EE\left[f(x^k)-f(x^*)\right] + 2\EE[\sigma_k^2] + 4L\max L_{ij}\EE[V_k],\label{eq:variance_loopless_local_svrg_fs}\\
		\EE\left[\left\|\frac{1}{n}\sum\limits_{i=1}^ng_i^k\right\|^2\right] &\le& 4\left(\frac{2\max L_{ij}}{n} + L\right)\EE\left[f(x^k)-f(x^*)\right] + \frac{2}{n}\EE[\sigma_k^2]\notag \\
		&&\quad + 2L\left(\frac{2\max L_{ij}}{n} + L\right)\EE[V_k],\label{eq:second_moment_loopless_local_svrg_fs_2}
	\end{eqnarray}
	where $\sigma_k^2 \eqdef \frac{1}{nm}\sum\limits_{i=1}^n\sum\limits_{j=1}^m\left\|\nabla f_{i,j}(y^k) - \nabla f_{i,j}(x^*)\right\|^2 + \frac{1}{n}\sum\limits_{i=1}^n\left\|\nabla f_{i}(y^k) - \nabla f_{i}(x^*)\right\|^2$.
\end{lemma}
\begin{proof}
	First of all, we have
	\begin{eqnarray*}
		\frac{1}{n}\sum\limits_{i=1}^n \EE_k\left[g_i^k\right] &=& \frac{1}{n}\sum\limits_{i=1}^n \EE_k\left[\nabla f_{i,j^k}(x_i^k) - \nabla f_{i,j_i}(y^k)  + \nabla f(y^k)\right]\\
		&=& \frac{1}{nm}\sum\limits_{i=1}^n\sum\limits_{j=1}^m\left(\nabla f_{i,j}(x_i^k) - \nabla f_{i,j}(y^k) + \nabla f(y^k)\right)\\
		&=& \frac{1}{n}\sum\limits_{i=1}^n\nabla f_i(x_i^k)
	\end{eqnarray*}
	and, in particular, $\bar{g}_i^k = \EE_k[g_i^k] = \nabla f_i(x_i^k) - \nabla f_i(y^k) + \nabla f(y^k)$.
	Using this we get
	\begin{eqnarray*}
		\frac{1}{n}\sum\limits_{i=1}^n\EE\left[\|\bar{g}_i^k\|^2\right] &\overset{\eqref{eq:a_b_norm_squared}}{\le}& \frac{2}{n}\sum\limits_{i=1}^n\EE\left[\|\nabla f_i(x_i^k)-\nabla f_i(x^*)\|^2\right]\\
		&&\quad + \frac{2}{n}\sum\limits_{i=1}^n\EE\left[\|\nabla f_i(y^k)-\nabla f_i(x^*) - (\nabla f(y^k)-\nabla f(x^*))\|^2\right]\\
		&\overset{\eqref{eq:L_smoothness_cor},\eqref{eq:variance_decomposition}}{\le}& \frac{4L}{n}\sum\limits_{i=1}^n\EE\left[D_{f_i}(x_i^k,x^*)\right] + \frac{2}{n}\sum\limits_{i=1}^n\EE\left[\|\nabla f_i(y^k)-\nabla f_i(x^*)\|^2\right]\\
		&\overset{\eqref{eq:poiouhnkj}}{\le}& 8L\EE\left[f(x^k)-f(x^*)\right] + 2\EE[\sigma_k^2] + 4L^2\EE[V_k]
	\end{eqnarray*}
	and
	\begin{eqnarray}
		\frac{1}{n}\sum\limits_{i=1}^n \EE\left[\|g_i^k-\bar{g}_i^k\|^2\right] &=& \frac{1}{n}\sum\limits_{i=1}^n\EE\left[\|\nabla f_{i,j_i}(x_i^k) - \nabla f_{i,j_i}(y^k)-(\nabla f_i(x_i^k) - \nabla f_i(y^k))\|^2\right]\notag\\
		&\overset{\eqref{eq:variance_decomposition}}{\le}&\frac{1}{n}\sum\limits_{i=1}^n\EE\left[\|\nabla f_{i,j_i}(x_i^k) - \nabla f_{i,j_i}(y^k)\|^2\right]\notag\\
		&\overset{\eqref{eq:a_b_norm_squared}}{\le}& \frac{2}{nm}\sum\limits_{i=1}^n\sum\limits_{j=1}^m\EE\left[\|\nabla f_{i,j}(x_i^k) - \nabla f_{i,j}(x^*)\|^2\right]\notag\\
		&&\quad +\frac{2}{nm}\sum\limits_{i=1}^n\sum\limits_{j=1}^m\EE\left[\|\nabla f_{i,j}(y^k) - \nabla f_{i,j}(x^*)\|^2\right]\notag \\
		&\overset{\eqref{eq:L_smoothness_cor}}{\le}& \frac{4\max L_{ij}}{n}\sum\limits_{i=1}^n\EE\left[D_{f_i}(x_i^k,x^*)\right] + 2\EE[\sigma_k^2]\notag\\
		&\overset{\eqref{eq:poiouhnkj}}{\le}& 8\max L_{ij}\EE\left[f(x^k)-f(x^*)\right] + 2\EE[\sigma_k^2] + 4L\max L_{ij}\EE[V_k].\label{eq:hscdvgdvauaicna}
	\end{eqnarray}
	Finally, using independence of $j_1,j_2,\ldots,j_n$ we derive
	\begin{eqnarray*}
		\EE\left[\left\|\frac{1}{n}\sum\limits_{i=1}^n g_i^k\right\|^2\right] &\overset{\eqref{eq:variance_decomposition},\eqref{eq:unbiasedness_loopless_local_svrg_fs}}{=}&  \EE\left[\left\|\frac{1}{n}\sum\limits_{i=1}^n \nabla f_i(x_i^k)\right\|^2\right] \\
		&&\quad + \EE\left[\left\|\frac{1}{n}\sum\limits_{i=1}^n (\nabla f_{i,j_i}(x_i^k)-\nabla f_{i,j_i}(y^k) - (\nabla f_i(x_i^k) - \nabla f_i(y^k)))\right\|^2\right]\\
		&=& \EE\left[\left\|\frac{1}{n}\sum\limits_{i=1}^n \nabla f_i(x_i^k)\right\|^2\right] \\
		&&\quad + \frac{1}{n^2}\sum\limits_{i=1}^n \EE\left[\|(\nabla f_{i,j_i}(x_i^k)-\nabla f_{i,j_i}(y^k) - (\nabla f_i(x_i^k) - \nabla f_i(y^k)))\|^2\right]\\
		&\overset{\eqref{eq:vdgasvgda},\eqref{eq:hscdvgdvauaicna}}{\le}& 4\left(\frac{2\max L_{ij}}{n} + L\right)\EE\left[f(x^k)-f(x^*)\right] + \frac{2}{n}\EE[\sigma_k^2]\\
		&&\quad + 2L\left(\frac{2\max L_{ij}}{n} + L\right)\EE[V_k].
	\end{eqnarray*}
\end{proof}

\begin{lemma}\label{lem:loopless_local_svrg_fs_sigma_k_bound}
	Let $f_{i}$ be convex and $L$-smooth and $f_{i,j}$ be convex and $\max L_{ij}$-smooth for all $i\in[n]$, $j\in[m]$. Then for all $k\ge 0$
	\begin{eqnarray}
		\EE\left[\sigma_{k+1}^2\right] &\le& \left(1-q\right)\EE\left[\sigma_k^2\right] + 2(L+\max L_{ij})q\EE\left[f(x^k) - f(x^*)\right] \label{eq:loopless_local_svrg_fs_sigma_k_bound}
	\end{eqnarray}
	where  $\sigma_k^2 \eqdef \frac{1}{nm}\sum\limits_{i=1}^n\sum\limits_{j=1}^m\left\|\nabla f_{i,j}(y^k) - \nabla f_{i,j}(x^*)\right\|^2 + \frac{1}{n}\sum\limits_{i=1}^n\left\|\nabla f_{i}(y^k) - \nabla f_{i}(x^*)\right\|^2$.
\end{lemma}
\begin{proof}
	First of all, we introduce new notations:
	\begin{eqnarray*}
		\sigma_{k,1}^2 \eqdef \frac{1}{nm}\sum\limits_{i=1}^n\sum\limits_{j=1}^m\left\|\nabla f_{i,j}(y^k) - \nabla f_{i,j}(x^*)\right\|^2,\quad \sigma_{k,2}^2 = \frac{1}{n}\sum\limits_{i=1}^n\left\|\nabla f_{i}(y^k) - \nabla f_{i}(x^*)\right\|^2.
	\end{eqnarray*}
	Secondly, by definition of $y^{k+1}$ we have
	\begin{eqnarray*}
		\EE\left[\sigma_{k+1,1}^2\mid x_1^k,\ldots, x_n^k\right] &=& \frac{1-q}{nm}\sum\limits_{i=1}^n\sum\limits_{j=1}^m\left\|\nabla f_{i,j}(y^k) - \nabla f_{i,j}(x^*)\right\|^2\\
		&&\quad + \frac{q}{nm}\sum\limits_{i=1}^n\sum\limits_{j=1}^m\left\|\nabla f_{i,j}(x^k) - \nabla f_{i,j}(x^*)\right\|^2\\
		&\overset{\eqref{eq:L_smoothness_cor}}{\le}& (1-q)\sigma_{k,1}^2 + 2q\max L_{ij}(f(x^k) - f(x^*)),
	\end{eqnarray*}
	hence
	\begin{equation}
		\EE\left[\sigma_{k+1,1}^2\right] \le (1-q)\EE\left[\sigma_{k,1}^2\right] + 2q\max L_{ij}\EE\left[f(x^k) - f(x^*)\right]. \label{eq:loopless_local_svrg_fs_sigma_k_bound_tech1}
	\end{equation}
	Next, the definition of $y^{k+1}$ implies
	\begin{eqnarray*}
		\EE\left[\sigma_{k+1,2}^2\mid x_1^k,\ldots, x_n^k\right] &=& \frac{1-q}{n}\sum\limits_{i=1}^n\|\nabla f_i(y^k) - \nabla f_i(x^*)\|^2 + \frac{q}{n}\sum\limits_{i=1}^n\|\nabla f_i(x^k) - \nabla f_i(x^*)\|^2\\
		&\overset{\eqref{eq:L_smoothness_cor}}{\le}& (1-q)\sigma_k^2 + 2Lq(f(x^k) - f(x^*)),
	\end{eqnarray*}
	hence
	\begin{equation}
		\EE\left[\sigma_{k+1,2}^2\right] \le (1-q)\EE\left[\sigma_{k,2}^2\right] + 2Lq\EE\left[f(x^k) - f(x^*)\right]. \label{eq:loopless_local_svrg_fs_sigma_k_bound_tech2}
	\end{equation}
	Finally, we combine obtained inequalities and get
	\begin{eqnarray*}
		\EE\left[\sigma_{k+1}\right] &=& \EE\left[\sigma_{k+1,1}^2\right] + \EE\left[\sigma_{k+1,2}^2\right]\\
		&\overset{\eqref{eq:loopless_local_svrg_fs_sigma_k_bound_tech1},\eqref{eq:loopless_local_svrg_fs_sigma_k_bound_tech2}}{\le}& (1-q)\left(\EE\left[\sigma_{k,1}^2\right] + \EE\left[\sigma_{k,2}^2\right]\right) + 2(L+\max L_{ij})q\EE\left[f(x^k) - f(x^*)\right] \\
		&=& \left(1-q\right)\EE\left[\sigma_k^2\right] + 2(L+\max L_{ij})q\EE\left[f(x^k) - f(x^*)\right],
	\end{eqnarray*}
	which concludes the proof.
\end{proof}

Using Corollary~\ref{cor:rand_loop} we obtain the following theorem.
\begin{theorem}\label{thm:loopless_local_svrg_fs}
	Assume that $f_{i}$ is $\mu$-strongly convex and $L$-smooth and $f_{i,j}$ is convex and $\max L_{ij}$-smooth for all $i\in[n]$, $j\in[m]$. Then {\tt S-Local-SVRG} satisfies Assumption~\ref{ass:hetero_second_moment} with
	\begin{gather*}
		\tA = 4L,\quad \hA = 4\max L_{ij},\quad \tB = \hB = 2,\quad \tF = 4L^2,\quad \hF = 4L\max L_{ij} \quad \tD_1 = \hD_1 = 0,\\
		A' = \frac{4\max L_{ij}}{n} + 2L,\quad B' = \frac{2}{n},\quad F' = 2L\left(\frac{2\max L_{ij}}{n} + L\right), \quad D_1' = 0,\\
		\sigma_k^2 = \frac{1}{nm}\sum\limits_{i=1}^n\sum\limits_{j=1}^m\left\|\nabla f_{i,j}(y^k) - \nabla f_{i,j}(x^*)\right\|^2 + \frac{1}{n}\sum\limits_{i=1}^n\left\|\nabla f_{i}(y^k) - \nabla f_{i}(x^*)\right\|^2,\\
		 \rho = q,\quad C = (L+\max L_{ij})q,\quad G = 0,\quad D_2 = 0,\quad H = \frac{256(1-p^2)(2+q)\gamma^2}{3p^2q},\quad D_3 = 0
	\end{gather*}
	under assumption that
	\begin{eqnarray*}
		\gamma \le \min\left\{\frac{1}{\frac{56\max L_{ij}}{3n} + 4L + \frac{32L}{3n}}, \frac{p\sqrt{3}}{32\sqrt{2L(1-p)\left(L(2+p)+p\max L_{ij} + \frac{4(L+\max L_{ij})(1+p)}{(1-q)}\right)}}\right\}.
	\end{eqnarray*}
	Moreover, for $\mu > 0$ we have
	\begin{eqnarray}
		\EE\left[f(\overline{x}^K) - f(x^*)\right] &\le& \left(1 - \min\left\{\gamma\mu,\frac{q}{4}\right\}\right)^K\frac{2\|x^0-x^*\|^2 + \frac{16\gamma^2\sigma_0^2}{nq} + \frac{1024L(1-p^2)(2+q)\gamma^3\sigma_0^2}{3p^2q}}{\gamma}\notag
	\end{eqnarray}
	and when $\mu = 0$ we have
	\begin{eqnarray}
		\EE\left[f(\overline{x}^K) - f(x^*)\right] &\le& \frac{2\|x^0-x^*\|^2 + \frac{16\gamma^2\sigma_0^2}{nq} + \frac{1024L(1-p^2)(2+q)\gamma^3\sigma_0^2}{3p^2q}}{\gamma K}. \notag
	\end{eqnarray}
\end{theorem}

The theorem above together with Lemma~\ref{lem:lemma2_stich} implies the following result.
\begin{corollary}\label{cor:s_local_svrg_str_cvx}
	Let assumptions of Theorem~\ref{thm:loopless_local_svrg_fs} hold with $\mu > 0$. Then for $q = \nicefrac{1}{m}$, $m \ge \nicefrac{1}{p}$,
	\begin{eqnarray*}
		\gamma = \min\left\{\frac{1}{\frac{56\max L_{ij}}{3n} + 4L + \frac{32L}{3n}}, \frac{p\sqrt{3}}{32\sqrt{2L(1-p)\left(L(2+p)+p\max L_{ij} + \frac{4(L+\max L_{ij})(1+p)}{(1-q)}\right)}}\right\}
	\end{eqnarray*}
	and for all $K\ge 1$ we have $\EE\left[f(\overline{x}^K)-f(x^*)\right]$ of order
	\begin{equation}
		\cO\left(\left(\frac{L}{p} + \frac{\max L_{ij}}{n} + \frac{\sqrt{(1-p)L\max L_{ij}}}{p}\right)\Phi^0\exp\left(- \min\left\{\Lambda,\frac{1}{m}\right\}K\right)\right),\notag
	\end{equation}
	where $\Phi^0 = 2\|x^0-x^*\|^2 + \frac{16\gamma^2\sigma_0^2}{nq} + \frac{1024L(1-p^2)(2+q)\gamma^3\sigma_0^2}{3p^2q}$, $\Lambda = \frac{\mu}{\frac{L}{p} + \frac{\max L_{ij}}{n} + \frac{\sqrt{(1-p)L\max L_{ij}}}{p}}$. That is, to achieve $\EE\left[f(\overline{x}^K)-f(x^*)\right] \le \varepsilon$ in this case {\tt S-Local-SVRG} requires
	\begin{equation*}
		\cO\left(\left(m+ \frac{L}{p\mu} + \frac{\max L_{ij}}{n\mu} + \frac{\sqrt{(1-p) L\max L_{ij}}}{p\mu}\right)\log\frac{\left(\frac{L}{p} + \frac{\max L_{ij}}{n} + \frac{\sqrt{(1-p)L\max L_{ij}}}{p}\right)\Phi^0}{\varepsilon}\right)
	\end{equation*}
	iterations/oracle calls per node (in expectation) and $\nicefrac{1}{p}$ times less communication rounds.	
\end{corollary}

That is, {\tt S-Local-SVRG} is the first implementable linearly converging stochastic method with local updates with a convergence guarantee in terms of the number of communications that is not worse than that of {\tt GD} even in the arbitrary heterogeneous data regime.

Next, we derive the following result for the convergence of {\tt S-Local-SVRG} in the case when $\mu = 0$.
\begin{corollary}
	\label{cor:s_local_svrg_cvx}
	Let assumptions of Theorem~\ref{thm:loopless_local_svrg_fs} hold with $\mu = 0$. Then for $q = \nicefrac{1}{m}$, $m \ge \nicefrac{1}{p}$ and
	\begin{eqnarray*}
		\gamma_0 = \min\left\{\frac{1}{\frac{56\max L_{ij}}{3n} + 4L + \frac{32L}{3n}}, \frac{p\sqrt{3}}{32\sqrt{2L(1-p)\left(L(2+p)+p\max L_{ij} + \frac{4(L+\max L_{ij})(1+p)}{(1-q)}\right)}}\right\}
	\end{eqnarray*}
	\begin{eqnarray*}
		\gamma = \min\left\{\gamma_0,\sqrt{\frac{nR_0^2}{8m\sigma_0^2}},\sqrt[3]{\frac{3p^2R_0^2}{512L(1-p^2)(2m+1)\sigma_0^2}}\right\}
	\end{eqnarray*}
	we have that $\EE\left[f(\overline{x}^K)-f(x^*)\right]$ is of the order 
	\begin{eqnarray*}
        \cO\left(\frac{\left(L + \nicefrac{p\max L_{ij}}{n} + \sqrt{(1-p)L\max L_{ij}}\right)R_0^2}{pK} + \frac{\sqrt{m\sigma_0^2 R_0^2}}{\sqrt{n}K} + \frac{\sqrt[3]{Lm\sigma_0^2 R_0^4}}{p^{\nicefrac{2}{3}}K}\right),
	\end{eqnarray*}
	where $R_0 = \|x^0 - x^*\|$. That is, to achieve $\EE\left[f(\overline{x}^K)-f(x^*)\right] \le \varepsilon$ in this case {\tt S-Local-SVRG} requires
	\begin{eqnarray*}
	K=	\cO\left(\frac{\left(L + \nicefrac{p\max L_{ij}}{n} + \sqrt{(1-p)L\max L_{ij}}\right)R_0^2}{p\varepsilon} + \frac{\sqrt{m\sigma_0^2 R_0^2}}{\sqrt{n}\varepsilon} + \frac{\sqrt[3]{Lm\sigma_0^2 R_0^4}}{p^{\nicefrac{2}{3}}\varepsilon}\right)
	\end{eqnarray*}
	iterations/oracle calls per node (in expectation) and $\nicefrac{1}{p}$ times less communication rounds.
\end{corollary}

\begin{remark}
	To get the rate from Tbl.~\ref{tbl:special_cases_weakly_convex} it remains to apply the following inequality:
	\begin{eqnarray*}
		\sigma_0^2 &=& \frac{1}{nm}\sum\limits_{i=1}^n\sum\limits_{j=1}^m\left\|\nabla f_{i,j}(x^0) - \nabla f_{i,j}(x^*)\right\|^2 + \frac{1}{n}\sum\limits_{i=1}^n\left\|\nabla f_{i}(x^0) - \nabla f_{i}(x^*)\right\|^2\\
		&\overset{\eqref{eq:L_smoothness}}{\le}& 2\left(\max L_{ij}^2 + L^2\right)\|x^0 - x^*\|^2.
	\end{eqnarray*}
\end{remark}

\section{Experiments} \label{sec:exp}
We perform multiple experiments to verify the theoretical claims of this chapter. Due to space limitations, we only present a single experiment in the main body; the rest can be found in Section~\ref{sec:extra_exp} of the appendix.

We demonstrate the benefit of on-device variance reduction, which we introduce in this chapter. For that purpose, we compare standard {\tt Local-SGD} (Algorithm~\ref{alg:local_sgd}) with our {\tt Local-SVRG} (Algorithm~\ref{alg:local_svrg}) on a regularized logistic regression problem with LibSVM data~\citep{chang2011libsvm}. For each problem instance, we compare the two algorithms with the stepsize $\gamma\in \{1, 0.1, 0.01 \}$ (we have normalized the data so that $L=1$).  The remaining details for the setup are presented in Section~\ref{sec:extralog} of the appendix.

Our theory predicts that both {\tt Local-SGD} and {\tt Local-SVRG} have identical convergence rate early on. However, the neighborhood of the optimum to which {\tt Local-SVRG} converges is smaller comparing to {\tt Local-SGD}. For both methods, the neighborhood is controlled by the stepsize: the smaller the stepsize is, the smaller the optimum neighborhood is. The price to pay is a slower rate at the beginning.

\begin{figure}[!h]
\centering
\includegraphics[width =  0.49\textwidth ]{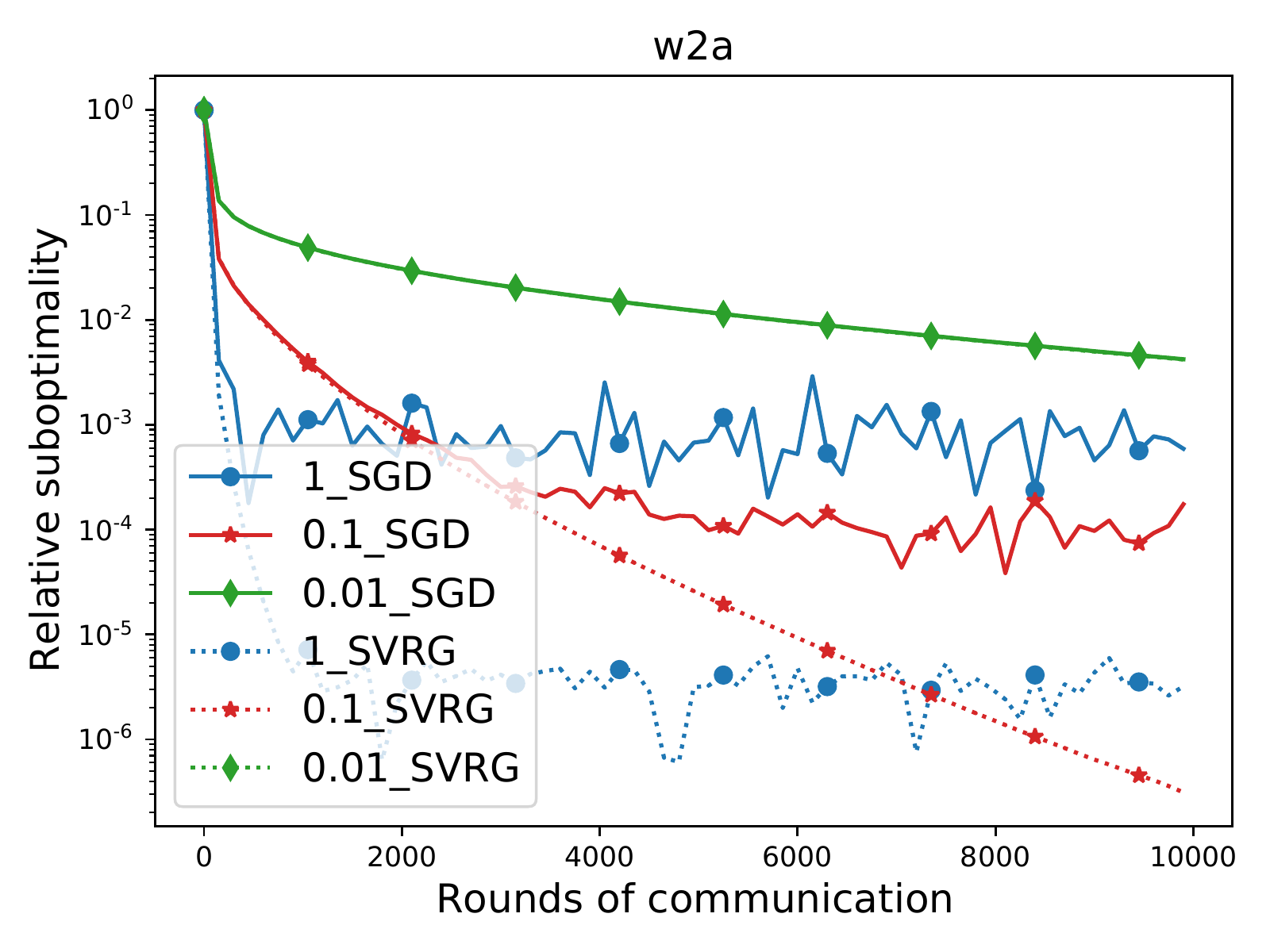}
\includegraphics[width =  0.49\textwidth ]{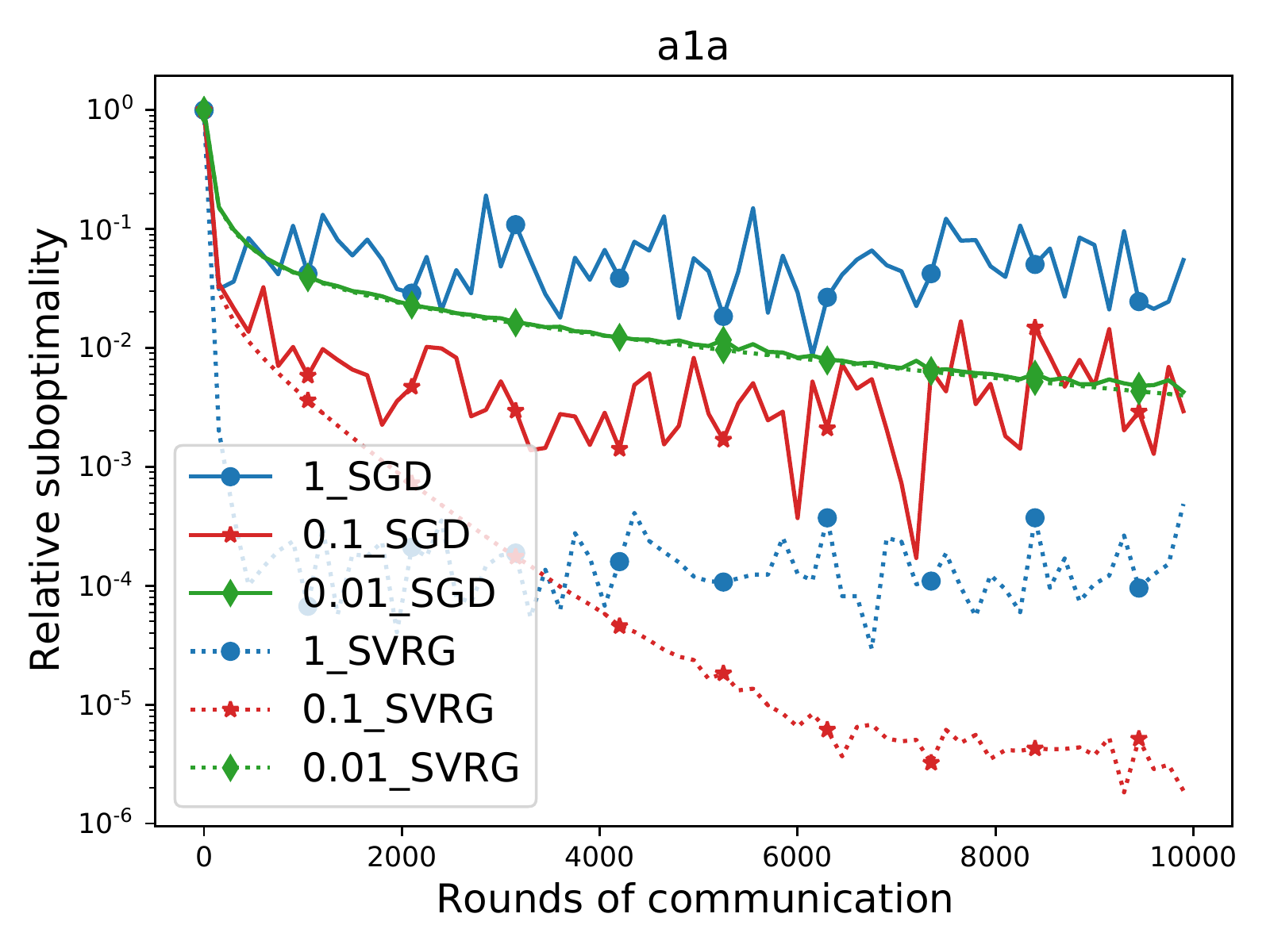}
\includegraphics[width =  0.49\textwidth ]{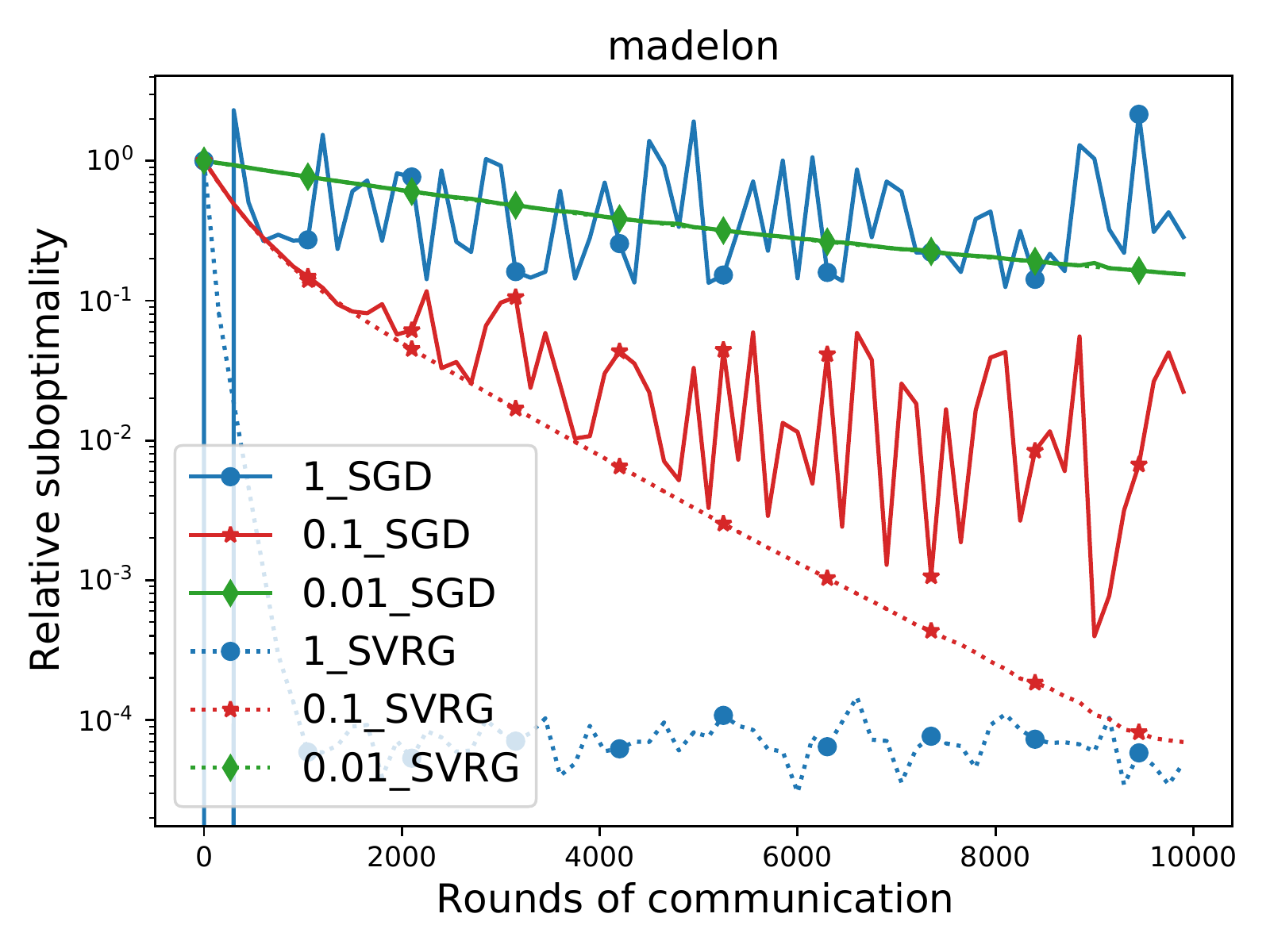}
\includegraphics[width =  0.49\textwidth ]{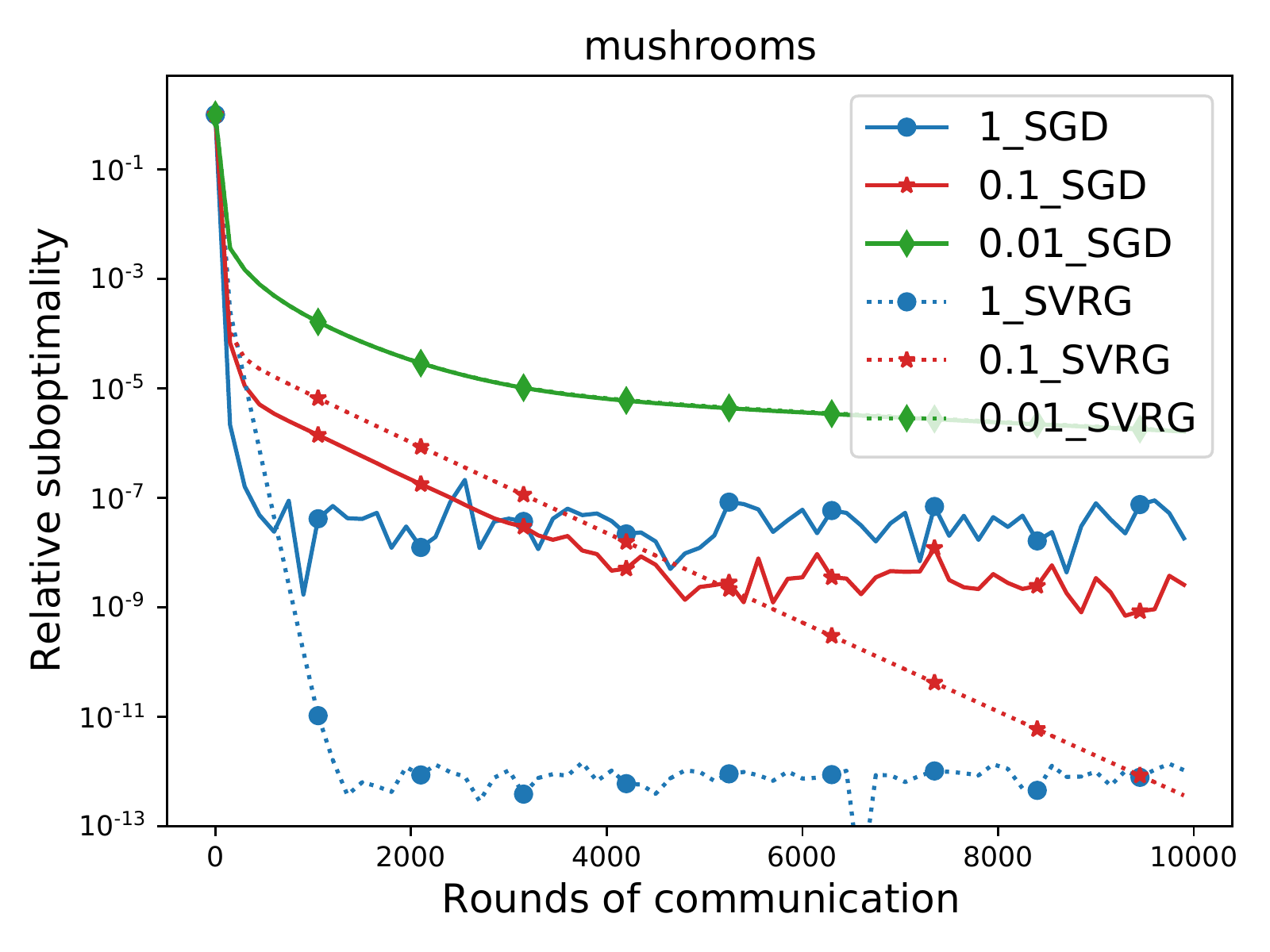}
\includegraphics[width =  0.49\textwidth ]{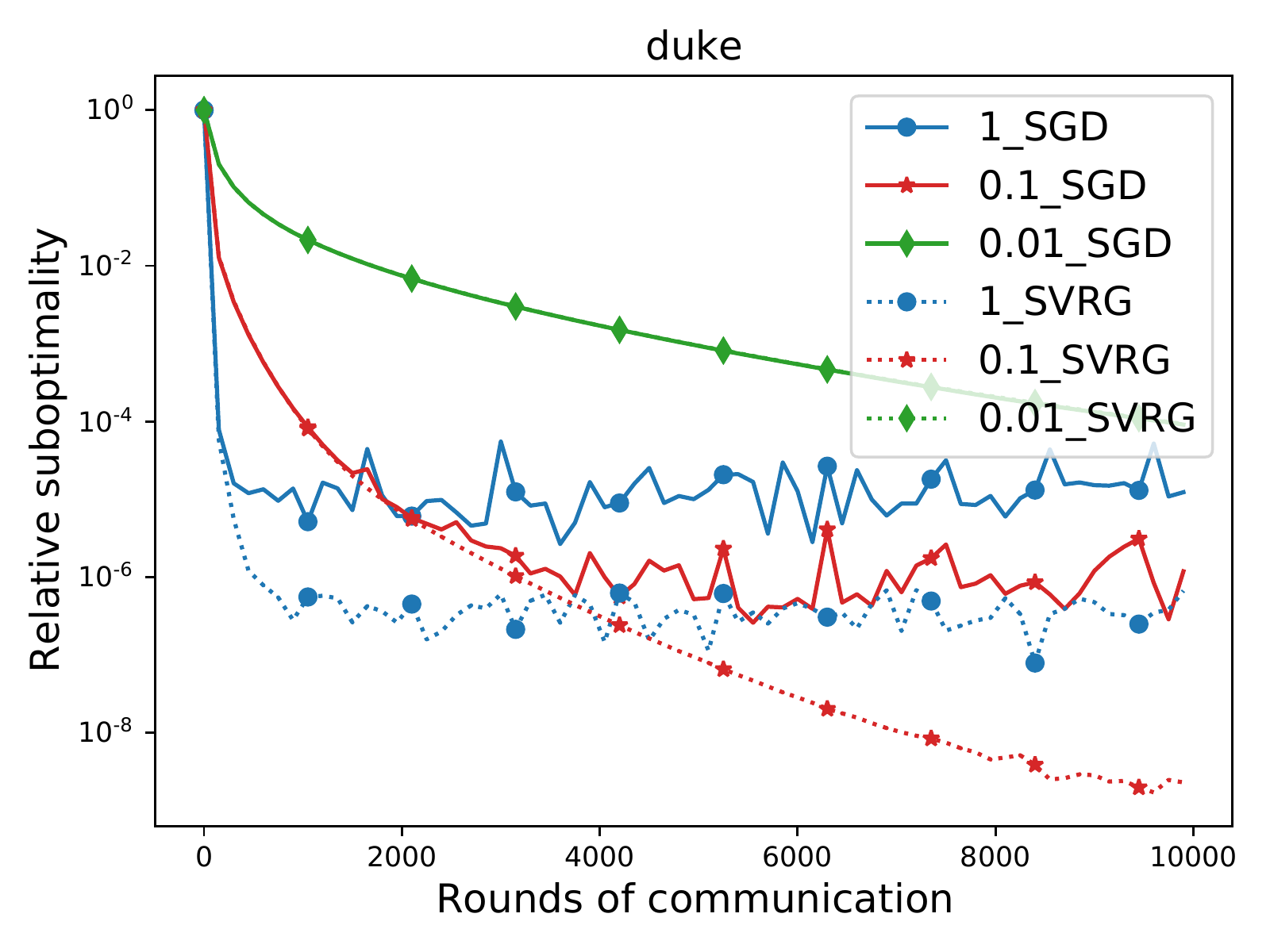}
\includegraphics[width =  0.49\textwidth ]{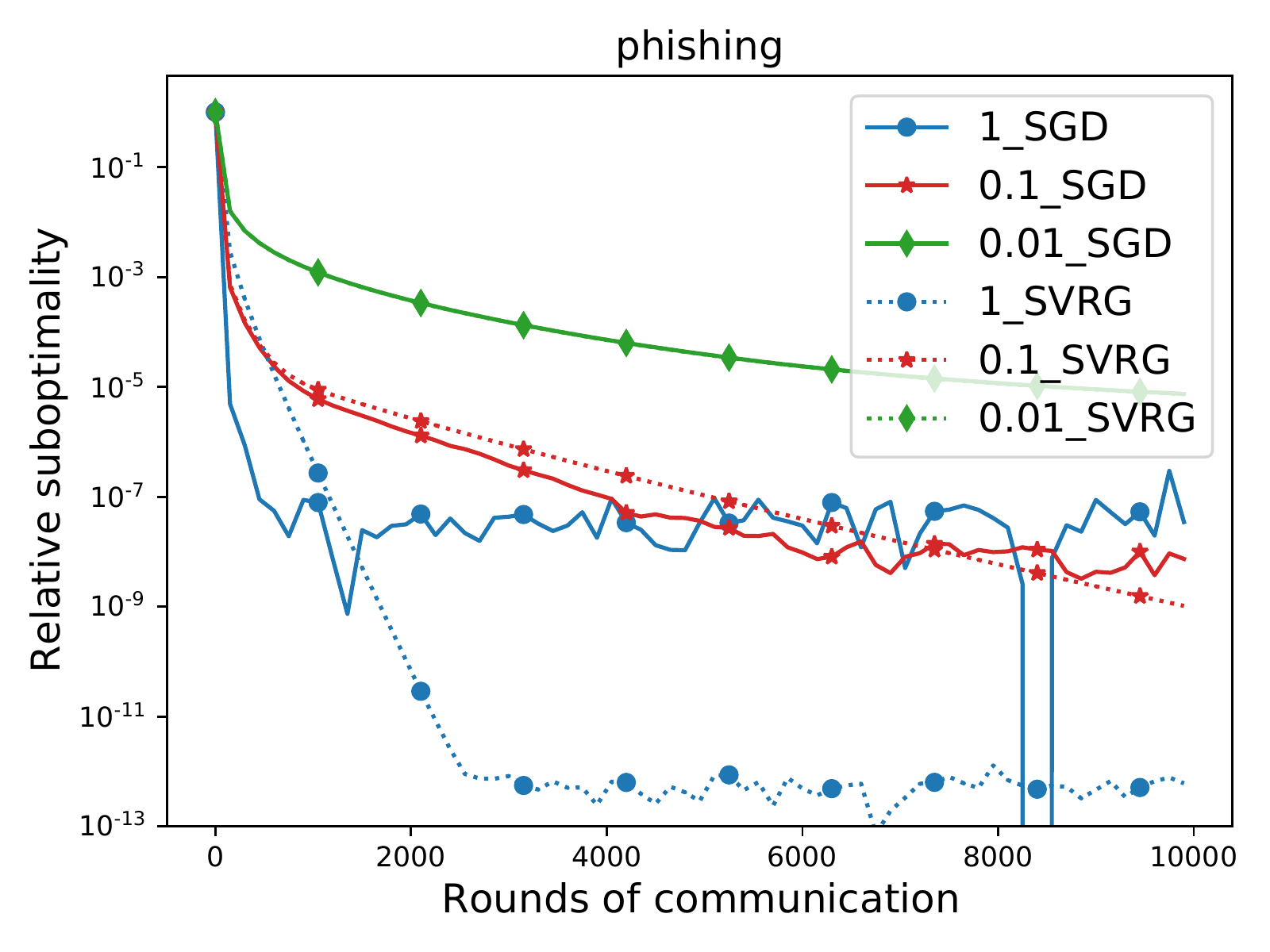}
\caption{Comparison of standard {\tt Local-SGD} (Alg.~\ref{alg:local_sgd}) and our {\tt Local-SVRG} (Alg.~\ref{alg:local_svrg}) for varying $\gamma$. Logistic regression applied on LibSVM~\citep{chang2011libsvm}.  Other parameters: $L=1, \mu=10^{-4}, \tau = 40$. Parameter $n$ chosen as per Tbl.~\ref{tbl:ns} in the appendix. }
\label{fig:sgd_svrg_hom}
\end{figure}

The results are presented in Figure~\ref{fig:sgd_svrg_hom}. As predicted, {\tt Local-SVRG} always outperforms {\tt Local-SGD} as it converges to a better neighborhood. Figure~\ref{fig:sgd_svrg_hom} also demonstrates that one can trade the smaller neighborhood for the slower convergence by modifying the stepsize.

\section{Conclusions and Future Work}
This chapter develops a unified approach to analyzing and designing a wide class of local stochastic first order algorithms. While our framework covers a broad range of methods, there are still some types of algorithms that we did not include but desire attention in future work. First, it would be interesting to study algorithms with  {\em biased} local stochastic gradients; these are popular for minimizing finite sums; see {\tt SAG}~\citep{schmidt2017minimizing} or {\tt SARAH}~\citep{nguyen2017sarah}. The second hitherto unexplored direction is including Nesterov's acceleration~\citep{nesterov1983method} in our framework. This idea is gaining traction in the area of  local methods already~\citep{pathak2020fedsplit, yuan2020federated}. However, it is not at all clear how this should be done and several attempts at achieving this unification goal failed. The third direction is allowing for a regularized local objective, which has been under-explored in the FL community so far. Other compelling directions that we do not cover are the local higher-order or proximal methods~\citep{li2018federated, pathak2020fedsplit} and methods supporting partial participation~\citep{mcmahan2016federated}.

%% file: ch5_marina.tex
\chapter{MARINA: Faster Non-Convex Distributed Learning with Compression}\label{ch:marina}

\section{Introduction}\label{sec:intro}
Non-convex\footnote{The results from this chapter were obtained while I was a research intern at KAUST. We thank Konstantin Mishchenko (KAUST) for a suggestion related to the experiments, Elena Bazanova (MIPT) for the suggestions about improving the text, and Slavom{\'{i}}r Hanzely (KAUST) and Egor Shulgin (KAUST) for spotting the typos.} optimization problems appear in various applications of machine learning, such as training deep neural networks \cite{Goodfellow-et-al-2016} and matrix completion and recovery \cite{ma2018implicit, bhojanapalli2016dropping}. Because of their practical importance, these problems gained much attention in recent years, which led to a rapid development of new efficient methods for non-convex optimization problems \cite{danilova2020recent}, and especially the training of deep learning models \cite{sun2019optimization}.

Training deep neural networks is notoriously computationally challenging and time-consuming. In the quest to improve the generalization performance of modern deep learning models, practitioners resort to using increasingly larger datasets in the training process, and to support such workloads, it is imperative to use advanced parallel and distributed hardware, systems, and algorithms. Distributed computing is often necessitated by the desire to train models from data naturally distributed across several edge devices, as is the case in federated learning \cite{FEDLEARN, FL2017-AISTATS}. However, even when this is not the case, distributed methods are often very efficient at reducing the training time \cite{goyal2017accurate,You2020Large}. Due to these and other reasons,   distributed optimization has gained immense popularity in recent years.

However, distributed methods almost invariably suffer from the so-called \textit{communication bottleneck}: the communication cost of information necessary for the workers  to jointly solve the problem  at hand is often very high, and depending on the particular compute architecture, workload, and algorithm used, it can be orders of magnitude higher than the computation cost. A popular technique for resolving this issue  is \textit{communication compression} \cite{seide20141, FEDLEARN, suresh2017distributed}, which is based on applying a lossy transformation/compression to the models, gradients, or tensors to be sent over the network to save on communication.  Since applying a lossy compression  generally decreases the utility of the exchanged messages, such an approach will typically lead to an increase in the number of communications, and the overall usefulness of this technique manifests itself in situations where the communication savings are larger compared to the increased need for the number of communication rounds \cite{Cnat}. 

The optimization and machine learning communities have exerted considerable effort in recent years to design  distributed methods  supporting compressed communication. From many methods proposed, we emphasize \algname{VR-DIANA} \cite{horvath2019stochastic},  \algname{FedCOMGATE} \cite{haddadpour2020federated}, and \algname{FedSTEPH} \cite{das2020improved} because these papers
contain the state-of-the-art results in the setup when the local loss functions can be arbitrary heterogeneous.


\begin{table}[H]
    \centering
    \scriptsize
	\caption{ Summary of the state-of-the-art results for finding an \textbf{$\varepsilon$-stationary point} for the problem \eqref{eq:main_problem_marina}, i.e., such a point $\hat x$ that $\EE\left[\|\nabla f(\hat x)\|^2\right] \le \varepsilon^2$. Dependences on the numerical constants, ``quality'' of the starting point, and smoothness constants are omitted in the complexity bounds.  Abbreviations: ``PP'' = partial participation; ``Communication complexity'' = the number of communications rounds needed to find an $\varepsilon$-stationary point; ``Oracle complexity'' = the number of (stochastic) first-order oracle calls needed to find an $\varepsilon$-stationary point. Notation: $\omega$ = the quantization parameter (see Def.~\ref{def:quantization}); $n$ = the number of nodes; $m$ = the size of the local dataset; $r$ = (expected) number of clients sampled at each iteration; $b'$ = the batchsize for \algnamex{VR-MARINA} at the iterations with compressed communication. To simplify the bounds, we assume that the expected density $\zeta_{\cQ}$ of the quantization operator $\cQ$ (see Def.~\ref{def:quantization}) satisfies $\omega+1 = \Theta(\nicefrac{d}{\zeta_{\cQ}})$ (e.g., this holds for RandK and $\ell_2$-quantization, see \cite{beznosikov2020biased}). We notice that \cite{haddadpour2020federated} and \cite{das2020improved} contain also better rates under different assumptions on clients' similarity.}
    \label{tab:comparison}    
    \begin{threeparttable}
    \begin{tabular}{|c|c|c c c|}
         \hline
         Setup & Method & Citation & Communication Complexity & Oracle Complexity\\ 
\hline\hline
    \multirow{5}{0.7cm}{\centering \eqref{eq:main_problem_marina}} &\algnamex{DIANA} &\makecell{\cite{mishchenko2019distributed}\\\cite{horvath2019stochastic}\\ \cite{li2020unified}} & $\frac{1+\left(1 + \omega\right)\sqrt{\nicefrac{\omega}{n}}}{\varepsilon^2}$ &  $\frac{1+\left(1 + \omega\right)\sqrt{\nicefrac{\omega}{n}}}{\varepsilon^2}$ \\
    & \algnamex{FedCOMGATE}\tnote{\color{red} (1)} & \cite{haddadpour2020federated} & $\frac{1 + \omega}{\varepsilon^2}$ & $\frac{1+\omega}{n\varepsilon^4}$\\
    & \algnamex{FedSTEPH}, $r=n$ & \cite{das2020improved} & $\frac{1 + \nicefrac{\omega}{n}}{\varepsilon^4}$ & $\frac{1 + \nicefrac{\omega}{n}}{\varepsilon^4}$ \\
    & \cellcolor{bgcolor2}\algnamex{MARINA} (Alg.~\ref{alg:marina}) &\cellcolor{bgcolor2} \begin{tabular}{c}
        Thm.~\ref{thm:main_result_non_cvx} \\ Cor.~\ref{cor:main_result_non_cvx}
    \end{tabular} &\cellcolor{bgcolor2} $\frac{1 + \nicefrac{\omega}{\sqrt{n}}}{\varepsilon^2}$ &\cellcolor{bgcolor2} $\frac{1 + \nicefrac{\omega}{\sqrt{n}}}{\varepsilon^2}$\\    
    \hline\hline
    \multirow{4.5}{0.7cm}{\centering\eqref{eq:main_problem_marina}, 
    \eqref{eq:f_i_finite_sum}}& \algnamex{DIANA} &\cite{li2020unified} & $\frac{1+\left(1 + \omega\right)\sqrt{\nicefrac{\omega}{n}}}{\varepsilon^2} + \frac{1+\omega}{n\varepsilon^4}$ &  $\frac{1+\left(1 + \omega\right)\sqrt{\nicefrac{\omega}{n}}}{\varepsilon^2} + \frac{1+\omega}{n\varepsilon^4}$\\
    & \algnamex{VR-DIANA} & \cite{horvath2019stochastic} & $\frac{\left(m^{\nicefrac{2}{3}} + \omega\right)\sqrt{1+\nicefrac{\omega}{n}}}{\varepsilon^2}$ &  $\frac{\left(m^{\nicefrac{2}{3}} + \omega\right)\sqrt{1+\nicefrac{\omega}{n}}}{\varepsilon^2}$\\
    &\cellcolor{bgcolor2} \begin{tabular}{c}
        \algnamex{VR-MARINA} (Alg.~\ref{alg:vr_marina})\\ $b'=1$\tnote{\color{red}(2)}
    \end{tabular} &\cellcolor{bgcolor2} \begin{tabular}{c}
        Thm.~\ref{thm:main_result_non_cvx_finite_sums} \\ Cor.~\ref{cor:main_result_non_cvx_finite_sums}
    \end{tabular} &\cellcolor{bgcolor2} $\frac{1 + \nicefrac{\max\left\{\omega,\sqrt{(1+\omega)m}\right\}}{\sqrt{n}}}{\varepsilon^2}$ &\cellcolor{bgcolor2} $\frac{1 + \nicefrac{\max\left\{\omega,\sqrt{(1+\omega)m}\right\}}{\sqrt{n}}}{\varepsilon^2}$\\
    \hline\hline
    \multirow{6.5}{0.7cm}{\centering\eqref{eq:main_problem_marina}, \eqref{eq:f_i_expectation}}& \algnamex{DIANA}\tnote{\color{red}(3)} &\makecell{\cite{mishchenko2019distributed}\\\cite{li2020unified}} & $\frac{1+\left(1 + \omega\right)\sqrt{\nicefrac{\omega}{n}}}{\varepsilon^2} + \frac{1+\omega}{n\varepsilon^4}$ &  $\frac{1+\left(1 + \omega\right)\sqrt{\nicefrac{\omega}{n}}}{\varepsilon^2} + \frac{1+\omega}{n\varepsilon^4}$\\
    & \algnamex{FedCOMGATE}\tnote{\color{red}(3)} & \cite{haddadpour2020federated} & $\frac{1 + \omega}{\varepsilon^2}$ & $\frac{1+\omega}{n\varepsilon^4}$\\
    &\cellcolor{bgcolor2} \begin{tabular}{c}
        \algnamex{VR-MARINA}  (Alg.~\ref{alg:vr_marina})\\ $b' = 1$
    \end{tabular} &\cellcolor{bgcolor2} \begin{tabular}{c}
        Thm.~\ref{thm:main_result_non_cvx_online} \\ Cor.~\ref{cor:main_result_non_cvx_online}
    \end{tabular} &\cellcolor{bgcolor2} $\frac{1 + \nicefrac{\omega}{\sqrt{n}}}{\varepsilon^2} + \frac{\sqrt{1+\omega}}{n\varepsilon^3}$ &\cellcolor{bgcolor2} $\frac{1 + \nicefrac{\omega}{\sqrt{n}}}{\varepsilon^2} + \frac{\sqrt{1+\omega}}{n\varepsilon^3}$\\
    &\cellcolor{bgcolor2} \begin{tabular}{c}
        \algnamex{VR-MARINA} (Alg.~\ref{alg:vr_marina})\\ $b' = \Theta\left(\frac{1}{n\varepsilon^2}\right)$
    \end{tabular} &\cellcolor{bgcolor2} \begin{tabular}{c}
        Thm.~\ref{thm:main_result_non_cvx_online} \\ Cor.~\ref{cor:main_result_non_cvx_online}
    \end{tabular} &\cellcolor{bgcolor2} $\frac{1 + \nicefrac{\omega}{\sqrt{n}}}{\varepsilon^2}$ &\cellcolor{bgcolor2} $\frac{1 + \nicefrac{\omega}{\sqrt{n}}}{n\varepsilon^4} + \frac{1+\omega}{n\varepsilon^3}$\\
    \hline\hline
    \multirow{3}{0.7cm}{\centering PP, \eqref{eq:main_problem_marina}}& \algnamex{FedSTEPH} & \cite{das2020improved} & $\frac{1+\nicefrac{\omega}{n}}{r\varepsilon^4} + \frac{(1+\omega)(n-r)}{r(n-1)\varepsilon^4}$ & $\frac{1+\nicefrac{\omega}{n}}{r\varepsilon^4} + \frac{(1+\omega)(n-r)}{r(n-1)\varepsilon^4}$\\
    &\cellcolor{bgcolor2} \algnamex{PP-MARINA} (Alg.~\ref{alg:pp_marina}) &\cellcolor{bgcolor2} \begin{tabular}{c}
        Thm.~\ref{thm:main_result_non_cvx_pp} \\ Cor.~\ref{cor:main_result_non_cvx_pp}
    \end{tabular} &\cellcolor{bgcolor2} $\frac{1+ \nicefrac{(1 + \omega)\sqrt{n}}{r}}{\varepsilon^2}$ &\cellcolor{bgcolor2} $\frac{1+ \nicefrac{(1 + \omega)\sqrt{n}}{r}}{\varepsilon^2}$\\
    \hline
    \end{tabular}
    \begin{tablenotes}
      {
        \item [{\color{red}(1)}] The results for \algnamex{FedCOMGATE} are derived under assumption that for all vectors $x_1,\ldots,x_n \in\R^d$ the quantization operator $\cQ$ satisfies $\EE\left[\left\|\frac{1}{n}\sum_{i=1}^n\cQ(x_j)\right\|^2 - \left\|\cQ\left(\frac{1}{n}\sum_{i=1}^n x_j\right)\right\|^2\right] \le G$ for some constant $G \ge 0$. In fact, this assumption does not hold for classical quantization operators like RandK and $\ell_2$-quantization on $\R^d$. The counterexample: $n=2$ and $x_1 = -x_2 = (t,t,\ldots,t)^\top$ with arbitrary large $t > 0$.
        \item [{\color{red}(2)}] One can even further improve the communication complexity by increasing $b'$.
        \item [{\color{red}(3)}] No assumptions on the smoothness of the stochastic realizations $f_{\xi}(x)$ are used.
      }
    \end{tablenotes}
    \end{threeparttable}
\end{table}

\subsection{Contributions}

We propose several new distributed optimization methods supporting compressed communication, specifically focusing on smooth but  nonconvex problems of the form
\begin{equation}
	\min\limits_{x\in\R^d}\left\{f(x) = \frac{1}{n}\sum\limits_{i=1}^n f_i(x)\right\}, \label{eq:main_problem_marina}
\end{equation}
where $n$ workers/devices/clients/peers are connected in a centralized way with a parameter-server, and client $i$ has an access to the local loss function $f_i$ only.   We establish strong complexity rates for them and show that they are better than previous state-of-the-art results.

$\bullet$ \textbf{MARINA.} The main contribution of our chapter is a new distributed method supporting communication compression called \algname{MARINA} (Alg~\ref{alg:marina}). In this algorithm, workers apply an unbiased compression operator to the {\em gradient differences} at each iteration with some probability and send them to the server that performs aggregation by averaging. Unlike all known methods operating with unbiased compression operators, this procedure leads to a {\em biased} gradient estimator.    We prove  convergence guarantees for \algname{MARINA}, which are strictly better than previous state-of-the-art methods (see Table~\ref{tab:comparison}). For example, \algname{MARINA}'s rate $\cO(\frac{1+\omega/\sqrt{n}}{\varepsilon^2})$ is $\cO(\sqrt{\omega})$ times better than that of the state-of-the-art method \algname{DIANA} \citep{mishchenko2019distributed}, where $\omega$ is the variance parameter associated with the deployed compressor. For example, in the case of the Rand1 sparsification compressor, we have $\omega=d-1$, and hence we get an improvement by the factor $\cO(\sqrt{d})$. Since the number $d$ of features can be truly very large when training modern models, this is a substantial improvement that can even amount to  {\em several orders of magnitude.} 

$\bullet$ \textbf{Variance Reduction on Nodes.} We generalize \algname{MARINA} to \algname{VR-MARINA}, which can handle the situation when the local functions $f_i$ have either a finite-sum (each $f_i$ is an average of $m$ functions) or an expectation form, and when it is more efficient to rely on local stochastic gradients rather than on local gradients. When compared with \algname{MARINA},  \algname{VR-MARINA} additionally performs {\em local variance reduction} on all nodes, progressively removing the variance coming from the stochastic approximation, leading to a better oracle complexity than previous state-of-the-art results (see Table~\ref{tab:comparison}). When no compression is used (i.e., $\omega=0$), the rate of \algname{VR-MARINA} is $\cO(\frac{\sqrt{m}}{\sqrt{n} \varepsilon^2})$, while the rate of the state-of-the-art method \algname{VR-DIANA} is $\cO(\frac{m^{2/3}}{\varepsilon^2})$. This is an improvement by the factor $\cO(\sqrt{n}m^{1/6})$. When much compression is applied, and $\omega$ is large, our method is faster by the factor  $\cO(\frac{m^{2/3} + \omega}{m^{1/2} + \omega^{1/2}})$. In the special case, when there is just a single node ($n=1$), and no compression is used, \algname{VR-MARINA} reduces to the \algname{PAGE} method of \cite{li2020page}; this is an optimal first-order algorithm for smooth non-convex finite-sum/online optimization problems.

$\bullet$ \textbf{Partial Participation.} We develop a modification of \algname{MARINA} allowing for {\em partial participation} of the clients, which is a feature critical in federated learning. The resulting method, \algname{PP-MARINA}, has  superior communication complexity to the existing methods developed for this settings (see Table~\ref{tab:comparison}).

$\bullet$ \textbf{Convergence Under the Polyak-{\L}ojasiewicz Condition.} We analyze all proposed methods for problems satisfying the Polyak-{\L}ojasiewicz condition \cite{polyak1963gradient,lojasiewicz1963topological}. Again, the obtained results are strictly better than previous ones (see Table~\ref{tab:comparison_pl}). Statements and proofs of all these results are in the Appendix.

$\bullet$ \textbf{Simple Analysis.} The simplicity and flexibility of our analysis offer several extensions. For example, one can easily generalize our analysis to the case of different quantization operators and different batch sizes used by clients. Moreover, one can combine the ideas of \algname{VR-MARINA} and \algname{PP-MARINA} and obtain a single distributed algorithm with compressed communications, variance reduction on nodes, and clients' sampling. We did not do this to keep the exposition simpler.

\begin{table}[H]
    \centering
    \scriptsize
	\caption{Summary of the state-of-the-art results for finding an $\varepsilon$-solution for the problem \eqref{eq:main_problem_marina} satifying \textbf{Polyak-{\L}ojasiewicz condition} (see As.~\ref{as:pl_condition}), i.e., such a point $\hat x$ that $\EE\left[f(\hat x) - f(x^*)\right] \le \varepsilon$. Dependences on the numerical constants and $\log(\nicefrac{1}{\varepsilon})$ factors are omitted and all smoothness constanst are denoted by $L$ in the complexity bounds.  Abbreviations: ``PP'' = partial participation; ``Communication complexity'' = the number of communications rounds needed to find an $\varepsilon$-stationary point; ``Oracle complexity'' = the number of (stochastic) first-order oracle calls needed to find an $\varepsilon$-stationary point. Notation: $\omega$ = the quantization parameter (see Def.~\ref{def:quantization}); $n$ = the number of nodes; $m$ = the size of the local dataset; $r$ = (expected) number of clients sampled at each iteration; $b'$ = the batchsize for \algnamex{VR-MARINA} at the iterations with compressed communication. To simplify the bounds, we assume that the expected density $\zeta_{\cQ}$ of the quantization operator $\cQ$ (see Def.~\ref{def:quantization}) satisfies $\omega+1 = \Theta(\nicefrac{d}{\zeta_{\cQ}})$ (e.g., this holds for RandK and $\ell_2$-quantization, see \cite{beznosikov2020biased}). We notice that \cite{haddadpour2020federated} and \cite{das2020improved} contain also better rates under different assumptions on clients' similarity.}
    \label{tab:comparison_pl}    
   \begin{threeparttable}
    \begin{tabular}{|c|c|c c c|}
         \hline
         Setup & Method & Citation & Communication Complexity & Oracle Complexity \\ 
\hline\hline
    \multirow{4}{0.7cm}{\centering \eqref{eq:main_problem_marina}} & \algnamex{DIANA} &\cite{li2020unified} & $\frac{L(1+\left(1 + \omega\right)\sqrt{\nicefrac{\omega}{n}})}{\mu}$ &  $\frac{L(1+\left(1 + \omega\right)\sqrt{\nicefrac{\omega}{n}})}{\mu}$ \\
    & \algnamex{FedCOMGATE}\tnote{\color{red} (1)} & \cite{haddadpour2020federated} & $\frac{L(1 + \omega)}{\mu}$ & $\frac{L(1+\omega)}{n\mu\varepsilon}$\\
    & \cellcolor{bgcolor2}\algnamex{MARINA} (Alg.~\ref{alg:marina}) &\cellcolor{bgcolor2} \begin{tabular}{c}
        Thm.~\ref{thm:main_result_pl}\\ Cor.~\ref{cor:main_result_pl_appendix}
    \end{tabular} &\cellcolor{bgcolor2} $\omega+\frac{L(1 + \nicefrac{\omega}{\sqrt{n}})}{\mu}$ &\cellcolor{bgcolor2} $\omega+\frac{L(1 + \nicefrac{\omega}{\sqrt{n}})}{\mu}$\\    
    \hline\hline
    \multirow{4.5}{0.7cm}{\centering\eqref{eq:main_problem_marina}, \eqref{eq:f_i_finite_sum}}& \algnamex{DIANA} &\cite{li2020unified} & \makecell{$\frac{L(1+\left(1 + \omega\right)\sqrt{\nicefrac{\omega}{n}})}{\mu}+$\quad\quad\quad\\ \quad\quad\quad$ + \frac{L(1+\omega)}{n\mu}\left(\frac{L}{\mu}+\frac{1}{\varepsilon}\right)$} &  \makecell{$\frac{L(1+\left(1 + \omega\right)\sqrt{\nicefrac{\omega}{n}})}{\mu}+$\quad\quad\quad\\ \quad\quad\quad$ + \frac{L(1+\omega)}{n\mu}\left(\frac{L}{\mu}+\frac{1}{\varepsilon}\right)$}\\
    & \algnamex{VR-DIANA}& \cite{li2020unified} & $\frac{L\left(m^{\nicefrac{2}{3}} + \omega\right)\sqrt{1+\nicefrac{\omega}{n}}}{\mu}$ &  $\frac{L\left(m^{\nicefrac{2}{3}} + \omega\right)\sqrt{1+\nicefrac{\omega}{n}}}{\mu}$\\
    &\cellcolor{bgcolor2} \begin{tabular}{c}
        \algnamex{VR-MARINA} (Alg.~\ref{alg:vr_marina})\\ $b'=1$\tnote{\color{red}(2)}
    \end{tabular} &\cellcolor{bgcolor2} \begin{tabular}{c}
        Thm.~\ref{thm:main_result_pl_finite_sums_appendix}\\ Cor.~\ref{cor:main_result_pl_finite_sums_appendix}
    \end{tabular} &\cellcolor{bgcolor2} \makecell{$\omega + m +$\quad\quad\quad\quad\quad\quad\quad\quad\quad\quad\\$+\frac{L(1 + \nicefrac{\max\left\{\omega,\sqrt{(1+\omega)m}\right\}}{\sqrt{n}})}{\mu}$} &\cellcolor{bgcolor2} \makecell{$\omega + m +$\quad\quad\quad\quad\quad\quad\quad\quad\quad\quad\\$+\frac{L(1 + \nicefrac{\max\left\{\omega,\sqrt{(1+\omega)m}\right\}}{\sqrt{n}})}{\mu}$}\\
    \hline\hline
    \multirow{6.5}{0.7cm}{\centering\eqref{eq:main_problem_marina}, \eqref{eq:f_i_expectation}}& \algnamex{DIANA}\tnote{\color{red}(3)} &\makecell{\cite{mishchenko2019distributed}\\\cite{li2020unified}} & $\frac{1+\left(1 + \omega\right)\sqrt{\nicefrac{\omega}{n}}}{\varepsilon^2} + \frac{1+\omega}{n\varepsilon^4}$ &  $\frac{1+\left(1 + \omega\right)\sqrt{\nicefrac{\omega}{n}}}{\varepsilon^2} + \frac{1+\omega}{n\varepsilon^4}$\\
    & \algnamex{FedCOMGATE}\tnote{\color{red}(3)} & \cite{haddadpour2020federated} & $\frac{L(1 + \omega)}{\mu}$ & $\frac{L(1+\omega)}{n\mu\varepsilon}$\\
    &\cellcolor{bgcolor2} \begin{tabular}{c}
        \algnamex{VR-MARINA} (Alg.~\ref{alg:vr_marina})\\ $b' = 1$
    \end{tabular} &\cellcolor{bgcolor2} \begin{tabular}{c}
        Thm.~\ref{thm:main_result_pl_online_appendix}\\ Cor.~\ref{cor:main_result_pl_online_appendix}
    \end{tabular} &\cellcolor{bgcolor2} $\omega+\frac{L(1 + \nicefrac{\omega}{\sqrt{n}})}{\mu} + \frac{L\sqrt{1+\omega}}{n\mu\varepsilon}$ &\cellcolor{bgcolor2} $\omega+\frac{L(1 + \nicefrac{\omega}{\sqrt{n}})}{\mu} + \frac{L\sqrt{1+\omega}}{n\mu\varepsilon}$\\
    &\cellcolor{bgcolor2} \begin{tabular}{c}
        \algnamex{VR-MARINA} (Alg.~\ref{alg:vr_marina})\\ $b' = \Theta\left(\frac{1}{n\mu\varepsilon}\right)$
    \end{tabular} &\cellcolor{bgcolor2} \begin{tabular}{c}
        Thm.~\ref{thm:main_result_pl_online_appendix}\\ Cor.~\ref{cor:main_result_pl_online_appendix}
    \end{tabular} &\cellcolor{bgcolor2} $\omega+\frac{L(1 + \nicefrac{\omega}{\sqrt{n}})}{\mu}$ &\cellcolor{bgcolor2} $\frac{1+\omega}{n\mu\varepsilon}+\frac{L(1 + \nicefrac{\omega}{\sqrt{n}})}{n\mu^2\varepsilon} + \frac{L(1+\omega)}{n\mu^2\sqrt{\varepsilon}}$\\
    \hline\hline
    \multirow{3}{0.7cm}{\centering PP, \eqref{eq:main_problem_marina}}& \algnamex{FedSTEPH}\tnote{\color{red}(4)} & \cite{das2020improved} & $\left(\frac{L}{\mu}\right)^{\nicefrac{3}{2}}$ & $\left(\frac{L}{\mu}\right)^{\nicefrac{3}{2}}$\\
    &\cellcolor{bgcolor2} \algnamex{PP-MARINA} (Alg.~\ref{alg:pp_marina}) &\cellcolor{bgcolor2} \begin{tabular}{c}
        Thm.~\ref{thm:main_result_pl_pp_appendix}\\ Cor.~\ref{cor:main_result_pl_pp_appendix}
    \end{tabular} &\cellcolor{bgcolor2} $\frac{(\omega+1)n}{r}+\frac{L(1+ \nicefrac{(1 + \omega)\sqrt{n}}{r})}{\mu}$ &\cellcolor{bgcolor2} $\frac{(\omega+1)n}{r}+\frac{L(1+ \nicefrac{(1 + \omega)\sqrt{n}}{r})}{\mu}$\\
    \hline
    \end{tabular}
    \begin{tablenotes}
      {
        \item [{\color{red}(1)}] The results for \algnamex{FedCOMGATE} are derived under assumption that for all vectors $x_1,\ldots,x_n \in\R^d$ the quantization operator $\cQ$ satisfies $\EE\left[\left\|\frac{1}{n}\sum_{i=1}^n\cQ(x_j)\right\|^2 - \left\|\cQ\left(\frac{1}{n}\sum_{i=1}^n x_j\right)\right\|^2\right] \le G$ for some constant $G \ge 0$. In fact, this assumption does not hold for classical quantization operators like RandK and $\ell_2$-quantization on $\R^d$. The counterexample: $n=2$ and $x_1 = -x_2 = (t,t,\ldots,t)^\top$ with arbitrary large $t > 0$.
        \item [{\color{red}(2)}] One can even further improve the communication complexity by increasing $b'$.
        \item [{\color{red}(3)}] No assumptions on the smoothness of the stochastic realizations $f_{\xi}(x)$ are used.
        \item [{\color{red}(4)}] The rate is derived under assumption that $r = \Omega((1+\omega)\sqrt{\nicefrac{L}{\mu}}\log(\nicefrac{1}{\varepsilon}))$.
      }
    \end{tablenotes}
    \end{threeparttable}
\end{table}

\subsection{Related Work}
\noindent\textbf{Non-Convex Optimization.} Since finding a global minimum of a non-convex function is, in general, an NP-hard problem \cite{murty1987some}, many researchers in non-convex optimization focus on relaxed goals such as finding an $\varepsilon$-stationary point. The theory of stochastic first-order methods for finding $\varepsilon$-stationary points is well-developed: it contains lower bounds for expectation minimization without smoothness of stochastic realizations \cite{arjevani2019lower} and for finite-sum/expectation minimization \cite{fang2018near,li2020page} as well as optimal methods matching the lower bounds (see \cite{danilova2020recent,li2020page} for the overview). Recently, distributed variants of such methods were proposed \cite{sun2020improving,sharma2019parallel,khanduri2020distributed}.

\noindent\textbf{Compressed Communications.} Works on  distributed methods supporting communication compression can be roughly split into two large groups: the first group focuses on methods using {\em unbiased} compression operators (which refer to as quantizations in this chapter), such as RandK, and the second one studies methods using {\em biased} compressors such as TopK. One can find a detailed summary of the most popular compression operators in \citep{UP2020, beznosikov2020biased}.

\noindent\textbf{Unbiased Compression.} In this line of work, the first convergence result in the non-convex case was obtained by \cite{alistarh2017qsgd} for  \algname{QSGD}, under assumptions that the local loss functions are the same for all workers, and the stochastic gradient has uniformly bounded second moment. After that, \cite{mishchenko2019distributed} proposed \algname{DIANA} (and its momentum version)  and proved its convergence rate for non-convex problems without any assumption on the boundedness of the second moment of the stochastic gradient, but under the assumption that the dissimilarity between local loss functions is bounded. This restriction was later eliminated by \cite{horvath2019stochastic} for the variance reduced version of \algname{DIANA} called  \algname{VR-DIANA}, and the analysis was extended to a large class of unbiased compressors. Finally, the results for  \algname{QSGD} and \algname{DIANA} were recently generalized and tightened by \cite{li2020unified} in a unifying framework that included many other methods as well.

\noindent\textbf{Biased Compression.} Biased compression operators are less ``optimization-friendly'' than unbiased ones. Indeed, one can construct a simple convex quadratic problem for which distributed \algname{SGD} with Top1 compression diverges exponentially fast  \cite{beznosikov2020biased}. However, this issue can be resolved using {\em error compensation} \cite{seide20141}. The first analysis of error-compensated \algname{SGD} (\algname{EC-SGD}) for non-convex problems was obtained by \cite{karimireddy2019error} for homogeneous problems under the assumption that the second moment of the stochastic gradient is uniformly  bounded. The last assumption was recently removed from the analysis of \algname{EC-SGD} by \cite{stich2020error, beznosikov2020biased}, while the first results without the homogeneity assumption were obtained by \cite{KoloskovaLSJ19decentralized} for \algname{Choco-SGD}, but still under the assumption that the second moment of the stochastic gradient is uniformly  bounded. This issue was resolved by \cite{beznosikov2020biased}. In general, the current understanding of optimization methods with biased compressors is far from complete: even in the strongly convex case, the first linearly converging \cite{gorbunov2020linearly} and accelerated \cite{qian2020error} error-compensated stochastic methods were proposed just recently.

\noindent\textbf{Other Approaches.} Besides communication compression, there are also different techniques aiming to reduce the overall communication cost of distributed methods. The most popular ones are based on decentralized communications and multiple local steps between communication rounds, where the second technique is very popular in federated learning \cite{FEDLEARN,kairouz2019advances}. One can find the state-of-the-art distributed optimization methods using these techniques and their combinations in \cite{lian2017can,karimireddy2020scaffold,li2019communication,koloskova2020unified}. Moreover, there exist results based on the combinations of communication compression with either decentralized communication, e.g., \algname{Choco-SGD} \cite{KoloskovaLSJ19decentralized}, or local updates, e.g., \algname{Qsparse-Local-SGD} \cite{basu2019qsparse}, \algname{FedCOMGATE} \cite{haddadpour2020federated}, \algname{FedSTEPH} \cite{das2020improved}, where in \cite{basu2019qsparse} the convergence rates were derived under an assumption that the stochastic gradient has uniformly bounded second moment and the results for \algname{Choco-SGD}, \algname{FedCOMGATE}, \algname{FedSTEPH} were described either earlier in the text, or in Table~\ref{tab:comparison}.


\subsection{Preliminaries}

We will rely on two key assumptions throughout the text.
\begin{assumption}[Uniform lower bound]\label{as:lower_bound}
	There exists  $f_*\in \R$ such that $f(x) \ge f_*$ for all $x\in\R^d$.
\end{assumption}
\begin{assumption}[$L$-smoothness]\label{as:L_smoothness}
	We assume that $f_i$ is $L_i$-smooth for all $i\in [n] = \{1,2,\dots,n\}$ meaning that the following inequality holds $\forall x,y\in \R^d$, $\forall i\in [n]$:
	\begin{equation}
		\left\|\nabla f_{i}(x) - \nabla f_{i}(y)\right\| \le L_i\|x-y\|.\label{eq:L_smoothness_local_marina}
	\end{equation}
\end{assumption}
This assumption implies that $f$ is $L_f$-smooth with $L_f^2 \le L^2 = \frac{1}{n}\sum_{i=1}^nL_i^2$.


\section{{\tt MARINA}: Compressing Gradient Differences}\label{sec:marina}
In this section, we describe the main algorithm of this work: \algname{MARINA}  (see Algorithm~\ref{alg:marina}). At each iteration of \algname{MARINA}, each  worker $i$ either sends to the server  the dense vector $\nabla f_i(x^{k+1})$ with probability $p$, or it sends the quantized gradient difference $\cQ\left(\nabla f_{i}(x^{k+1}) - \nabla f_{i}(x^k))\right)$ with probability $1-p$. In the first situation, the server just averages the vectors received from workers and gets $g^{k+1} = \nabla f(x^{k+1})$, whereas in the second case, the server averages the quantized differences from all workers and then adds the result to $g^k$ to get $g^{k+1}.$ Moreover, if $\cQ$ is identity quantization, i.e., $\cQ(x) = x$, then \algname{MARINA} reduces to Gradient Descent (\algname{GD}).

\begin{algorithm}[h]
   \caption{\algname{MARINA}}\label{alg:marina}
\begin{algorithmic}[1]
   \State {\bfseries Input:} starting point $x^0$, stepsize $\gamma$, probability $p\in(0,1]$, number of iterations $K$
   \State Initialize $g^0 = \nabla f(x^0)$
   \For{$k=0,1,\ldots,K-1$}
   \State Sample $c_k \sim \text{Be}(p)$
   \State Broadcast $g^k$ to all workers
   \For{$i = 1,\ldots,n$ in parallel} 
   \State $x^{k+1} = x^k - \gamma g^k$
   \State Set $g_i^{k+1} = \begin{cases}\nabla f_i(x^{k+1}),& \text{if } c_k = 1,\\ g^k + \cQ\left(\nabla f_{i}(x^{k+1}) - \nabla f_{i}(x^k))\right),& \text{if } c_k = 0\end{cases}$
   \EndFor
   \State $g^{k+1} = \frac{1}{n}\sum_{i=1}^ng_i^{k+1}$
   \EndFor
   \State {\bfseries Return:} $\hat x^K$ chosen uniformly at random from $\{x^k\}_{k=0}^{K-1}$
\end{algorithmic}
\end{algorithm}

However, for non-trivial quantizations, we have $\EE[g^{k+1}\mid x^{k+1}] \neq \nabla f(x^{k+1})$ unlike all other distributed methods using exclusively unbiased compressors we know of. That is, $g^{k+1}$ is a \textit{biased} stochastic estimator of $\nabla f(x^{k+1})$. However, \algname{MARINA} is an example of a rare phenomenon in stochastic optimization  when the {\em bias of the stochastic gradient helps to achieve better complexity.}

\subsection{Convergence Results for Generally Non-Convex Problems}
We start with the following result.
\begin{theorem}\label{thm:main_result_non_cvx}
	Let Assumptions~\ref{as:lower_bound}~and~\ref{as:L_smoothness} be satisfied. Then, after
	\begin{equation}
		K = \cO\left(\frac{\Delta_0 L}{\varepsilon^2}\left(1 + \sqrt{\frac{(1-p)\omega}{pn}}\right)\right) \notag
	\end{equation}
	iterations with $\Delta_0 = f(x^0)-f_*$, $L^2 = \frac{1}{n}\sum_{i=1}^nL_i^2$ and the stepsize 
	\begin{equation}
		\gamma \le \frac{1}{L\left(1 + \sqrt{\frac{(1-p)\omega}{pn}}\right)}\notag
	\end{equation}
	\algname{MARINA} produces  point $\hat x^K$ for which $\EE[\|\nabla f(\hat x^K)\|^2] \le \varepsilon^2$.
\end{theorem}
One can find the full statement of the theorem together with its proof in Section~\ref{sec:proof_of_thm_non_cvx} of the Appendix.

The following corollary provides the bounds on the number of iterations/communication rounds and estimates the total communication cost needed to achieve an $\varepsilon$-stationary point in expectation. Moreover, for simplicity, throughout the chapter we assume that the communication cost is proportional to the number of non-zero components of transmitted vectors from workers to the server.
\begin{corollary}\label{cor:main_result_non_cvx}
	Let the assumptions of Theorem~\ref{thm:main_result_non_cvx} hold and $p = \nicefrac{\zeta_{\cQ}}{d}$. If 
	\begin{equation*}
		\gamma \le \frac{1}{L\left(1 + \sqrt{\frac{\omega}{n}\left(\frac{d}{\zeta_{\cQ}}-1\right)}\right)},
	\end{equation*}
	then \algname{MARINA} requires 
	\begin{equation*}
		\cO\left(\frac{\Delta_0 L}{\varepsilon^2}\left(1 + \sqrt{\frac{\omega}{n}\left(\frac{d}{\zeta_{\cQ}}-1\right)}\right)\right)
	\end{equation*}
	iterations/communication rounds in order to achieve $\EE[\|\nabla f(\hat x^K)\|^2] \le \varepsilon^2$, and the expected total communication cost per worker is $\cO(d + \zeta_{\cQ}K)$.
\end{corollary}

Let us clarify the obtained result. First of all, if $\omega = 0$ (no quantization), then $\zeta_{\cQ} = 0$ and the rate coincides with the rate of Gradient Descent (\algname{GD}). Since \algname{GD} is optimal among first-order methods in terms of reducing the norm of the gradient \cite{carmon2019lower}, the dependence on $\varepsilon$ in our bound cannot be improved in general. Next, if $n$ is large enough, i.e., $n \ge \omega(\nicefrac{d}{\zeta_{\cQ}}-1)$, then\footnote{For $\ell_2$-quantization this requirement is satisfied when $n \ge d$.}  the iteration complexity of \algname{MARINA} (method with compressed communications) and \algname{GD} (method with dense communications) coincide. This means that in this regime,  \algname{MARINA} is able to reach a provably better communication complexity than \algname{GD}!

\begin{remark}
    When $p = \nicefrac{1}{(\omega + 1)}$ the complexity bound for {\tt MARINA} becomes
    \begin{equation*}
		\cO\left(\frac{\Delta_0 L}{\varepsilon^2}\left(1 + \frac{\omega}{\sqrt{n}}\right)\right).
	\end{equation*}
	Since the definition of quantization (Definition~\ref{def:quantization}) covers uniform coordinate-wise randomization and directional derivative oracle (directions are sampled from the uniform distribution on the unit Euclidean sphere) with $\omega = d-1$, the dependence on $\omega$ cannot be improved in general. One can prove this using the standard results for derivative-free methods from \cite{nemirovskij1983problem} that multiplicative dependence on $\cO(d)$ is unavoidable and approximating partial or directional derivative oracle using finite differences. Similar arguments hold for the methods from the next sections as well. 
\end{remark}

\subsection{Convergence Results Under Polyak-{\L}ojasiewicz Condition}
In this section, we provide a complexity bounds for \algname{MARINA} under the Polyak-{\L}ojasiewicz (P{\L}) condition.

\begin{assumption}[P{\L} condition]\label{as:pl_condition}
	Function $f$ satisfies Polyak-{\L}ojasiewicz (P{\L}) condition with parameter $\mu$, i.e., 	\begin{equation}
		\|\nabla f(x)\|^2 \ge 2\mu\left(f(x) - f(x^*)\right). \label{eq:pl_condition}
	\end{equation}
	holds for $x^*= \arg \min_{x\in\R^d} f(x)$ and for all $x\in\R^d$.

\end{assumption}

Under this and previously introduced assumptions, we derive the following result.
\begin{theorem}\label{thm:main_result_pl}
	Let Assumptions~\ref{as:lower_bound},~\ref{as:L_smoothness}~and~\ref{as:pl_condition} be satisfied. Then, after
	\begin{equation}
		K = \cO\left(\max\left\{\frac{1}{p},\frac{L}{\mu}\left(1 + \sqrt{\frac{(1-p)\omega}{pn}}\right)\right\}\log\frac{\Delta_0}{\varepsilon}\right) \notag
	\end{equation}
	iterations with $\Delta_0 = f(x^0)-f(x^*)$, $L^2 = \frac{1}{n}\sum_{i=1}^nL_i^2$ and the stepsize 
	\begin{equation}
		\gamma \le \min\left\{\frac{1}{L\left(1 + \sqrt{\frac{2(1-p)\omega}{pn}}\right)}, \frac{p}{2\mu}\right\}\notag
	\end{equation}
	\algname{MARINA} produces a point $x^K$ for which $\EE[f(x^K) - f(x^*)] \le \varepsilon$.
\end{theorem}
One can find the full statement of the theorem together with its proof in Section~\ref{sec:proof_of_thm_pl} of the Appendix.


\section{{\tt MARINA} and Variance Reduction}\label{sec:vr}
Throughout this section, we assume that the local loss on each node has either a finite-sum form (finite sum case), 
\begin{equation}
	f_i(x) = \frac{1}{m}\sum\limits_{j=1}^mf_{ij}(x), \label{eq:f_i_finite_sum}
\end{equation}
or an expectation form (online case),
\begin{equation}
	f_i(x) = \EE_{\xi_i\sim\cD_i}[f_{\xi_i}(x)]. \label{eq:f_i_expectation}
\end{equation}

\subsection{Finite Sum Case}
In this section, we generalize \algname{MARINA} to problems of the form \eqref{eq:main_problem_marina}+\eqref{eq:f_i_finite_sum}, obtaining  \algname{VR-MARINA} (see Algorithm~\ref{alg:vr_marina}).
\begin{algorithm}[h]
   \caption{\algname{VR-MARINA}: finite sum case}\label{alg:vr_marina}
\begin{algorithmic}[1]
   \State {\bfseries Input:} starting point $x^0$, stepsize $\gamma$, minibatch size $b'$, probability $p\in(0,1]$, number of iterations $K$
   \State Initialize $g^0 = \nabla f(x^0)$ 
   \For{$k=0,1,\ldots,K-1$}
   \State Sample $c_k \sim \text{Be}(p)$
   \State Broadcast $g^k$ to all workers
   \For{$i = 1,\ldots,n$ in parallel} 
   \State $x^{k+1} = x^k - \gamma g^k$
   \State Set $g_i^{k+1} = \begin{cases}\nabla f_i(x^{k+1}),& \text{if } c_k = 1,\\ g^k + \cQ\left(\frac{1}{b'}\sum_{j\in I_{i,k}'}(\nabla f_{ij}(x^{k+1}) - \nabla f_{ij}(x^k))\right),& \text{if } c_k = 0,\end{cases}$ where $I_{i,k}'$ is the set of the indices in the minibatch, $|I_{i,k}'| = b'$
   \EndFor
   \State $g^{k+1} = \frac{1}{n}\sum_{i=1}^ng_i^{k+1}$
   \EndFor
   \State {\bfseries Return:} $\hat x^K$ chosen uniformly at random from $\{x^k\}_{k=0}^{K-1}$
\end{algorithmic}
\end{algorithm}
At each iteration of \algname{VR-MARINA}, devices are to compute the full gradients $\nabla f_i(x^{k+1})$ and send them to the server with probability $p$. Typically, $p \le \nicefrac{1}{m}$ and $m$ is large, meaning that workers compute full gradients rarely (once per $\ge m$ iterations in expectation). At other iterations, workers compute minibatch stochastic gradients evaluated at the current and previous points, compress them using an unbiased compression operator, i.e., quantization/quantization operator, and send the resulting vectors $g_i^{k+1} - g^{k}$ to the server. Moreover, if $\cQ$ is the identity quantization, i.e., $\cQ(x) = x$, and $n=1$, then \algname{MARINA} reduces to the optimal method \algname{PAGE} \cite{li2020page}.

In this part, we will rely on the following average smoothness assumption.
\begin{assumption}[Average $\cL$-smoothness]\label{as:avg_smoothness}
	For all $k\ge 0$ and $i\in[n]$ the minibatch stochastic gradients difference $\widetilde{\Delta}_i^k = \frac{1}{b'}\sum_{j\in I_{i,k}'}(\nabla f_{ij}(x^{k+1}) - \nabla f_{ij}(x^k))$ computed on the $i$-th worker satisfies $\EE\left[\widetilde{\Delta}_i^k\mid x^k,x^{k+1}\right] = \Delta_i^k$ and
	\begin{equation}
		\EE\left[\left\|\widetilde{\Delta}_i^k - \Delta_i^k\right\|^2\mid x^k,x^{k+1}\right] \le \frac{\cL_i^2}{b'}\|x^{k+1}-x^k\|^2\label{eq:avg_L_smoothness}
	\end{equation}
	with some $\cL_i \ge 0$, where $\Delta_i^k = \nabla f_i(x^{k+1}) - \nabla f_i(x^k)$.
\end{assumption}
This assumption is satisfied in many standard minibatch regimes. In particular, if $I_{i,k}' = \{1,\ldots,m\}$, then $\cL_i = 0$, and if $I_{i,k}'$ consists of $b'$ i.i.d.\ samples from the uniform distributions on $\{1,\ldots,m\}$ and $f_{ij}$ are $L_{ij}$-smooth, then $\cL_i \le \max_{j\in [m]}L_{ij}$.

Under this and the previously introduced assumptions, we derive the following result.
\begin{theorem}\label{thm:main_result_non_cvx_finite_sums}
	Consider the finite sum case \eqref{eq:main_problem_marina}+\eqref{eq:f_i_finite_sum}. Let Assumptions~\ref{as:lower_bound},~\ref{as:L_smoothness}~and \ref{as:avg_smoothness} be satisfied. Then, after
	\begin{equation}
		K = \cO\left(\frac{\Delta_0}{\varepsilon^2}\left(L + \sqrt{\frac{1-p}{pn}\left(\omega L^2 + \frac{(1+\omega)\cL^2}{b'}\right)}\right)\right) \notag
	\end{equation}
	iterations with $\Delta_0 = f(x^0)-f_*$, $L^2 = \frac{1}{n}\sum_{i=1}^nL_i^2$, $\cL^2 = \frac{1}{n}\sum_{i=1}^n\cL_i^2$ and the stepsize 
	\begin{equation}
		\gamma \le \frac{1}{L + \sqrt{\frac{1-p}{pn}\left(\omega L^2 + \frac{(1+\omega)\cL^2}{b'}\right)}}\notag
	\end{equation}
	\algname{VR-MARINA} produces such a point $\hat x^K$ that $\EE[\|\nabla f(\hat x^K)\|^2] \le \varepsilon^2$.
\end{theorem}
One can find the full statement of the theorem together with its proof in Section~\ref{sec:proof_of_thm_non_cvx_fin_sums} of the Appendix.
\begin{corollary}\label{cor:main_result_non_cvx_finite_sums}
	Let the assumptions of Theorem~\ref{thm:main_result_non_cvx_finite_sums} hold and $p = \min\left\{\nicefrac{\zeta_{\cQ}}{d},\nicefrac{b'}{(m+b')}\right\}$, where $b'\le m$. If 
	\begin{equation*}
		\gamma \le \frac{1}{L + \sqrt{\frac{\max\left\{\nicefrac{d}{\zeta_{\cQ}} - 1,\nicefrac{m}{b'}\right\}}{n}\left(\omega L^2 + \frac{(1+\omega)\cL^2}{b'}\right)}},
	\end{equation*}
	then \algname{VR-MARINA} requires 
	\begin{eqnarray*}
	\cO\left(\frac{\Delta_0}{\varepsilon^2}\left(L\left(1 + \sqrt{\frac{\omega\max\left\{\nicefrac{d}{\zeta_{\cQ}} - 1,\nicefrac{m}{b'}\right\}}{n}}\right) + \cL\sqrt{\frac{(1+\omega)\max\left\{\nicefrac{d}{\zeta_{\cQ}} - 1,\nicefrac{m}{b'}\right\}}{nb'}}\right)\right)
	\end{eqnarray*}
	iterations/communication rounds and $\cO\left(m + b'K\right)$
	stochastic oracle calls per node in expectation in order to achieve $\EE[\|\nabla f(\hat x^K)\|^2] \le \varepsilon^2$, and the expected total communication cost per worker is $\cO(d + \zeta_{\cQ}K)$.
\end{corollary}

First of all, when workers quatize differences of the full gradients, then $I_{i,k}' = \{1,\ldots,m\}$ for all $i\in[n]$ and $k\ge 0$,  implying $\cL = 0$. In this case, the complexity bounds for \algname{VR-MARINA} recover the ones for \algname{MARINA}. Next, when $\omega = 0$ (no quantization) and $n = 1$, our bounds for iteration and oracle complexities for \algname{VR-MARINA} recover the bounds for \algname{PAGE} \cite{li2020unified}, which is optimal for finite-sum smooth non-convex optimization. This observation implies that the dependence on $\varepsilon$ and $m$ in the complexity bounds for \algname{VR-MARINA} cannot be improved in the class of first-order stochastic methods. Next, we notice that up to the differences in smoothness constants, the iteration and oracle complexities for \algname{VR-MARINA} benefit from the number of workers $n$. Finally, as Table~\ref{tab:comparison} shows, the rates for \algname{VR-MARINA} are strictly better than ones for the previous state-of-the-art method \algname{VR-DIANA}\cite{horvath2019stochastic}.

We provide the convergence results for \algname{VR-MARINA} in the finite-sum case under the Polyak-{\L}ojasiewicz condition,  together with complete proofs, in Section~\ref{sec:proof_of_thm_pl_fin_sums} of the Appendix.

\subsection{Online Case}
In this section, we focus on  problems of type \eqref{eq:main_problem_marina}+\eqref{eq:f_i_expectation}. For this type of problems, we consider a  slightly modified version of \algname{VR-MARINA}. That is, we replace line 8 in Algorithm~\ref{alg:vr_marina} with the following update rule: $g_i^{k+1} = \frac{1}{b}\sum_{j\in I_{i,k}}\nabla f_{\xi_{ij}^k}(x^{k+1})$ if $c_k = 1$, and $g_i^{k+1} = g^k + \cQ\left(\frac{1}{b'}\sum_{j\in I_{i,k}'}(\nabla f_{\xi_{ij}^k}(x^{k+1}) - \nabla f_{\xi_{ij}^k}(x^k))\right)$ otherwise, where $I_{i,k}, I_{i,k}'$ are the sets of the indices in the minibatches, $|I_{i,k}| = b$, $|I_{i,k}'| = b'$,  and $\xi_{ij}^k$ is independently sampled from $\cD_i$ for $i\in[n]$, $j\in[m]$ (see Algorithm~\ref{alg:vr_marina_online}).

\begin{algorithm}[h]
   \caption{\algname{VR-MARINA}: online case}\label{alg:vr_marina_online}
\begin{algorithmic}[1]
   \State {\bfseries Input:} starting point $x^0$, stepsize $\gamma$, minibatch sizes $b$, $b' < b$, probability $p\in(0,1]$, number of iterations $K$
   \State Initialize $g^0 = \frac{1}{nb}\sum_{i=1}^n\sum_{j\in I_{i,0}}\nabla f_{\xi_{ij}^0}(x^{k+1})$, where $I_{i,0}$ is the set of the indices in the minibatch, $|I_{i,0}| = b$, and $\xi_{ij}^0$ is independently sampled from $\cD_i$ for $i\in[n]$, $j\in[m]$
   \For{$k=0,1,\ldots,K-1$}
   \State Sample $c_k \sim \text{Be}(p)$
   \State Broadcast $g^k$ to all workers
   \For{$i = 1,\ldots,n$ in parallel} 
   \State $x^{k+1} = x^k - \gamma g^k$
   \State Set $g_{i}^{k+1} = \begin{cases}\frac{1}{b}\sum_{j\in I_{i,k}}\nabla f_{\xi_{ij}^k}(x^{k+1}),& \text{if } c_k = 1,\\ g^k + \cQ\left(\frac{1}{b'}\sum_{j\in I_{i,k}'}(\nabla f_{\xi_{ij}^k}(x^{k+1}) - \nabla f_{\xi_{ij}^k}(x^k))\right),& \text{if } c_k = 0, \end{cases}$  
    where $I_{i,k}, I_{i,k}'$ are the sets of the indices in the minibatches, $|I_{i,k}| = b$, $|I_{i,k}'| = b'$,  and $\xi_{ij}^k$ is independently sampled from $\cD_i$ for $i\in[n]$, $j\in[m]$
   \EndFor
   \State $g^{k+1} = \frac{1}{n}\sum_{i=1}^ng_i^{k+1}$
   \EndFor
   \State {\bfseries Return:} $\hat x^K$ chosen uniformly at random from $\{x^k\}_{k=0}^{K-1}$
\end{algorithmic}
\end{algorithm}

Before we provide our convergence results in this setup, we reformulate Assumption~\ref{as:avg_smoothness} for the online case.
\begin{assumption}[Average $\cL$-smoothness]\label{as:avg_smoothness_online}
	For all $k\ge 0$ and $i\in[n]$ the minibatch stochastic gradients difference $\widetilde{\Delta}_i^k = \frac{1}{b'}\sum_{j\in I_{i,k}'}(\nabla f_{\xi_{ij}^k}(x^{k+1}) - \nabla f_{\xi_{ij}^k}(x^k))$ computed on the $i$-th worker satisfies $\EE\left[\widetilde{\Delta}_i^k\mid x^k,x^{k+1}\right] = \Delta_i^k$ and
	\begin{equation}
		\EE\left[\left\|\widetilde{\Delta}_i^k - \Delta_i^k\right\|^2\mid x^k,x^{k+1}\right] \le \frac{\cL_i^2}{b'}\|x^{k+1}-x^k\|^2\label{eq:avg_L_smoothness_online}
	\end{equation}
	with some $\cL_i \ge 0$, where $\Delta_i^k = \nabla f_i(x^{k+1}) - \nabla f_i(x^k)$.
\end{assumption}
Moreover, we assume that the variance of the stochastic gradients on all nodes is uniformly upper bounded.
\begin{assumption}\label{as:bounded_var}
	We assume that for all $i \in [n]$ there exists such constant $\sigma_i \in [0,+\infty)$ that for all $x\in\R^d$
	\begin{eqnarray}
		\EE_{\xi_i \sim\cD_i}\left[\nabla f_{\xi_i}(x)\right]	&=& \nabla f_i(x), \label{eq:unbiasedness}\\
		\EE_{\xi_i \sim\cD_i}\left[\left\|\nabla f_{\xi_i}(x) - \nabla f_i(x)\right\|^2\right] &\le& \sigma_i^2. \label{eq:bounded_var}
	\end{eqnarray}
\end{assumption}

Under these and previously introduced assumptions, we derive the following result.
\begin{theorem}\label{thm:main_result_non_cvx_online}
	Consider the online case \eqref{eq:main_problem_marina}+\eqref{eq:f_i_expectation}. Let Assumptions~\ref{as:lower_bound},~\ref{as:L_smoothness},~\ref{as:avg_smoothness_online}~and \ref{as:bounded_var} be satisfied. Then, after
	\begin{equation}
		K = \cO\left(\frac{\Delta_0}{\varepsilon^2}\left(L + \sqrt{\frac{1-p}{pn}\left(\omega L^2 + \frac{(1+\omega)\cL^2}{b'}\right)}\right)\right) \notag
	\end{equation}
	iterations with $\Delta_0 = f(x^0)-f_*$, $L^2 = \frac{1}{n}\sum_{i=1}^nL_i^2$, $\cL^2 = \frac{1}{n}\sum_{i=1}^n\cL_i^2$, the stepsize 
	\begin{equation}
		\gamma \le \frac{1}{L + \sqrt{\frac{1-p}{pn}\left(\omega L^2 + \frac{(1+\omega)\cL^2}{b'}\right)}},\notag
	\end{equation}
	and $b = \Theta\left(\nicefrac{\sigma^2}{(n\varepsilon^2)}\right),$ $\sigma^2 = \frac{1}{n}\sum_{i=1}^n\sigma_i^2$,  
	\algname{VR-MARINA} produces a point $\hat x^K$ for which\newline  $\EE[\|\nabla f(\hat x^K)\|^2] \le \varepsilon^2$.
\end{theorem}
One can find the full statement of the theorem, together with its proof, in Section~\ref{sec:proof_of_thm_non_cvx_online} of the Appendix.
\begin{corollary}\label{cor:main_result_non_cvx_online}
	Let the assumptions of Theorem~\ref{thm:main_result_non_cvx_online} hold and choose $p = \min\left\{\nicefrac{\zeta_{\cQ}}{d},\nicefrac{b'}{(b+b')}\right\}$, where $b'\le b$, $b = \Theta\left(\nicefrac{\sigma^2}{(n\varepsilon^2)}\right)$. If 
	\begin{equation*}
		\gamma \le \frac{1}{L + \sqrt{\frac{\max\left\{\nicefrac{d}{\zeta_{\cQ}} - 1,\nicefrac{b}{b'}\right\}}{n}\left(\omega L^2 + \frac{(1+\omega)\cL^2}{b'}\right)}},
	\end{equation*}
	then \algname{VR-MARINA} requires 
	\begin{eqnarray*}
	\cO\left(\frac{\Delta_0}{\varepsilon^2}\left(L\left(1 + \sqrt{\frac{\omega}{n}\max\left\{\frac{d}{\zeta_{\cQ}} - 1,\frac{\sigma^2}{nb'\varepsilon^2}\right\}}\right) + \cL\sqrt{\frac{(1+\omega)}{nb'}\max\left\{\frac{d}{\zeta_{\cQ}} - 1,\frac{\sigma^2}{nb'\varepsilon^2}\right\}}\right)\right)
	\end{eqnarray*}
	iterations/communication rounds and $\cO(\zeta_{\cQ}K+ \nicefrac{\sigma^2}{(n\varepsilon^2)})$  
	stochastic oracle calls per node in expectation to achieve $\EE[\|\nabla f(\hat x^K)\|^2] \le \varepsilon^2$, and the expected total communication cost per worker is $\cO(d + \zeta_{\cQ}K)$.
\end{corollary}

Similarly to the finite-sum case, when $\omega = 0$ (no quantization) and $n = 1$, our bounds for iteration and oracle complexities for \algname{VR-MARINA} recover the bounds for \algname{PAGE} \cite{li2020unified}, which is optimal for online smooth non-convex optimization as well. That is, the dependence on $\varepsilon$ in the complexity bound for \algname{VR-MARINA} cannot be improved in the class of first-order stochastic methods. As previously, up to the differences in smoothness constants, the iteration and oracle complexities for \algname{VR-MARINA} benefit from an increase in the  number of workers $n$.

We provide the convergence results for \algname{VR-MARINA} in the online case under the Polyak-{\L}ojasiewicz condition, together with complete proofs, in Section~\ref{sec:proof_of_thm_pl_online} of the Appendix.

\section{{\tt MARINA} and Partial Participation}\label{sec:pp}
Finally, we propose another modification of \algname{MARINA}. In particular, we prove an option for {\em partial participation} of the clients - a feature important in federated learning. The resulting method is called \algname{PP-MARINA} (see Algorithm~\ref{alg:pp_marina}). At each iteration of \algname{PP-MARINA}, the server receives the quantized gradient differences from $r$ clients with probability $1-p$, and aggregates full gradients from all clients with probability $p$, i.e., \algname{PP-MARINA} coincides with \algname{MARINA} up to the following difference: $g_i^{k+1} = \nabla f_i(x^{k+1})$, $g^{k+1} = \frac{1}{n}\sum_{i=1}^ng_i^{k+1}$ if $c_k = 1$, and $g_i^{k+1} = g^k + \cQ\left(\nabla f_{i}(x^{k+1}) - \nabla f_{i}(x^k))\right)$, $g^{k+1} = \frac{1}{r}\sum_{i_k\in I_k'}g_{i_k}^{k+1}$ otherwise, where $I_k'$ is the set of $r$ i.i.d.\ samples from the uniform distribution over $\{1,\ldots,n\}$. That is, if the probability $p$ is chosen to be small enough, then with high probability the server receives only quantized vectors from a subset of clients at each iteration.

\begin{algorithm}[h]
   \caption{\algname{PP-MARINA}}\label{alg:pp_marina}
\begin{algorithmic}[1]
   \State {\bfseries Input:} starting point $x^0$, stepsize $\gamma$, probability $p\in(0,1]$, number of iterations $K$, clients-batchsize $r \le n$
   \State Initialize $g^0 = \nabla f(x^0)$
   \For{$k=0,1,\ldots,K-1$}
   \State Sample $c_k\sim \text{Be}(p)$
   \State Choose $I_k' = \{1,\ldots,n\}$ if $c_k = 1$, and choose $I_k'$ as the set of $r$ i.i.d.\ samples from the uniform distribution over $\{1,\ldots,n\}$ otherwise
   \State Broadcast $g^k$ to all workers
   \For{$i = 1,\ldots,n$ in parallel} 
   \State $x^{k+1} = x^k - \gamma g^k$
   \State Set $g_i^{k+1} = \begin{cases}\nabla f_i(x^{k+1})& \text{if } c_k = 1,\\ g^k + \cQ\left(\nabla f_{i}(x^{k+1}) - \nabla f_{i}(x^k)\right)& \text{if } c_k = 0. \end{cases}$ 
   \EndFor
   \State Set $g^{k+1} = \begin{cases}\nabla f(x^{k+1})& \text{if } c_k=1,\\ g^k + \frac{1}{r}\sum\limits_{i_k\in I'_{k}}\cQ\left(\nabla f_{i_k}(x^{k+1}) - \nabla f_{i_k}(x^k)\right)& \text{if } c_k=0. \end{cases}$ 
   \EndFor
   \State {\bfseries Return:} $\hat x^K$ chosen uniformly at random from $\{x^k\}_{k=0}^{K-1}$
\end{algorithmic}
\end{algorithm}

Below, we provide a convergence result for \algname{PP-MARINA} for smooth  non-convex problems.
\begin{theorem}\label{thm:main_result_non_cvx_pp}
	Let Assumptions~\ref{as:lower_bound}~and~\ref{as:L_smoothness} be satisfied. Then, after
	\begin{equation}
	K = \cO\left(\frac{\Delta_0 L}{\varepsilon^2}\left(1 + \sqrt{\frac{(1-p)(1+\omega)}{pr}}\right)\right) \notag
	\end{equation}
	iterations with $\Delta_0 = f(x^0)-f_*$, $L^2 = \frac{1}{n}\sum_{i=1}^nL_i^2$ and the stepsize 
	\begin{equation}
		\gamma \le \frac{1}{L\left(1 + \sqrt{\frac{(1-p)(1+\omega)}{pr}}\right)}\notag
	\end{equation}
	\algname{PP-MARINA} produces  a point $\hat x^K$ for which  $\EE[\|\nabla f(\hat x^K)\|^2] \le \varepsilon^2$.
\end{theorem}
One can find the full statement of the theorem together with its proof in Section~\ref{sec:proof_of_thm_non_cvx_pp} of the appendix.
\begin{corollary}\label{cor:main_result_non_cvx_pp}
	Let the assumptions of Theorem~\ref{thm:main_result_non_cvx_pp} hold and choose $p = \nicefrac{\zeta_{\cQ}r}{(dn)}$, where $r\le n$. If 
	\begin{equation*}
		\gamma \le \frac{1}{L\left(1 + \sqrt{\frac{1+\omega}{b'}\left(\frac{dn}{\zeta_{\cQ}r}-1\right)}\right)},
	\end{equation*}
	then \algname{PP-MARINA} requires 
	\begin{equation*}
	\cO\left(\frac{\Delta_0 L}{\varepsilon^2}\left(1 + \sqrt{\frac{1+\omega}{r}\left(\frac{dn}{\zeta_{\cQ}r}-1\right)}\right)\right)
	\end{equation*}
	iterations/communication rounds to achieve $\EE[\|\nabla f(\hat x^K)\|^2] \le \varepsilon^2$, and the expected total communication cost is $\cO\left(dn + \zeta_{\cQ}rK\right)$.
\end{corollary}

When $r=n$, i.e., all clients participate in communication with the server at each iteration, the rate for \algname{PP-MARINA} recovers the rate for \algname{MARINA} under the assumption that $(1+\omega)(\nicefrac{d}{\zeta_{\cQ}}-1) = \cO(\omega(\nicefrac{d}{\zeta_{\cQ}}-1))$, which holds for a wide class of quantization operators, e.g., for identical quantization, RandK, and $\ell_p$-quantization. In general, the derived complexity is strictly better than previous state-of-the-art one (see Table~\ref{tab:comparison}).

We provide the convergence results for \algname{PP-MARINA} under the Polyak-{\L}ojasiewicz condition, together with complete proofs, in Section~\ref{sec:proof_of_thm_pl_pp} of the Appendix.

\section{Numerical Experiments}
\subsection{Binary Classification with Non-Convex Loss}\label{sec:experiments}
We conduct several numerical experiments\footnote{Our code is available at \url{https://github.com/burlachenkok/marina}.} on binary classification problem involving non-convex loss \cite{zhao2010convex} (used for two-layer neural networks) with LibSVM data \cite{chang2011libsvm} to justify the theoretical claims of the chapter. That is, we consider the following optimization problem:
\begin{equation}
	\min\limits_{x\in\R^d}\left\{f(x) = \frac{1}{N}\sum\limits_{t=1}^N \ell(a_t^\top x, y_i)\right\},\label{eq:experiment_problem}
\end{equation}
where $\{a_t\} \in \R^d$, $y_i\in\{-1,1\}$ for all $t=1,\ldots,N$, and the function $\ell:\R^d \to \R$ is defined as
\begin{equation*}
	\ell(b,c) = \left(1 - \frac{1}{1+\exp(-bc)}\right)^2.
\end{equation*}
The datasets were taken from LibSVM \cite{chang2011libsvm} and split into five equal parts among five clients (we excluded $N - 5\cdot\lfloor\nicefrac{N}{5}\rfloor$ last datapoints from each dataset), see the summary in Table~\ref{tbl:ns}.
\begin{table}[H]
 \caption{Summary of the datasets and splitting of the data among clients (Figure~\ref{fig:full_batched_methods}).}
\label{tbl:ns}
\begin{center}
\begin{tabular}{|c|c|c|c|}
\hline
Dataset  & $n$ & $N$ (\# of datapoints) & $d$ (\# of features)   \\
 \hline
  \hline
\texttt{mushrooms} & 5 & 8 120 & 112   \\ \hline
\texttt{w8a} & 5  &49 745 & 300  \\ \hline
\texttt{phishing} & 5  &11 055 & 69  \\ \hline
\texttt{a9a} & 5  &32 560 & 124  \\ \hline
\end{tabular}
\end{center}
\end{table}

The code was written in Python 3.8 using \textsc{mpi4py} to emulate the distributed environment and then was executed on a machine with 48 cores, each is Intel(R) Xeon(R) Gold 6246 CPU 3.30GHz.

In our experiments, we compare \algname{MARINA} with the full-batch version of \algname{DIANA}, and then \algname{VR-MARINA} with \algname{VR-DIANA}. We exclude \algname{FedCOMGATE} and \algname{FedPATH} from this comparison since they have significantly worse oracle complexities (see Table~\ref{tab:comparison}). Since one of the main goals of our experiments is to justify the theoretical findings of the chapter, in the experiments, we used the stepsizes from the corresponding theoretical results for the methods (for \algname{DIANA} and \algname{VR-DIANA} the stepsizes were chosen according to \cite{horvath2019stochastic,li2020unified}). Next, to compute the stochastic gradients, we use batchsizes $= \max\{1, \nicefrac{m}{100}\}$ for \algname{VR-MARINA} and \algname{VR-DIANA}.

\begin{figure}[H]
\centering
\includegraphics[width=0.24\textwidth]{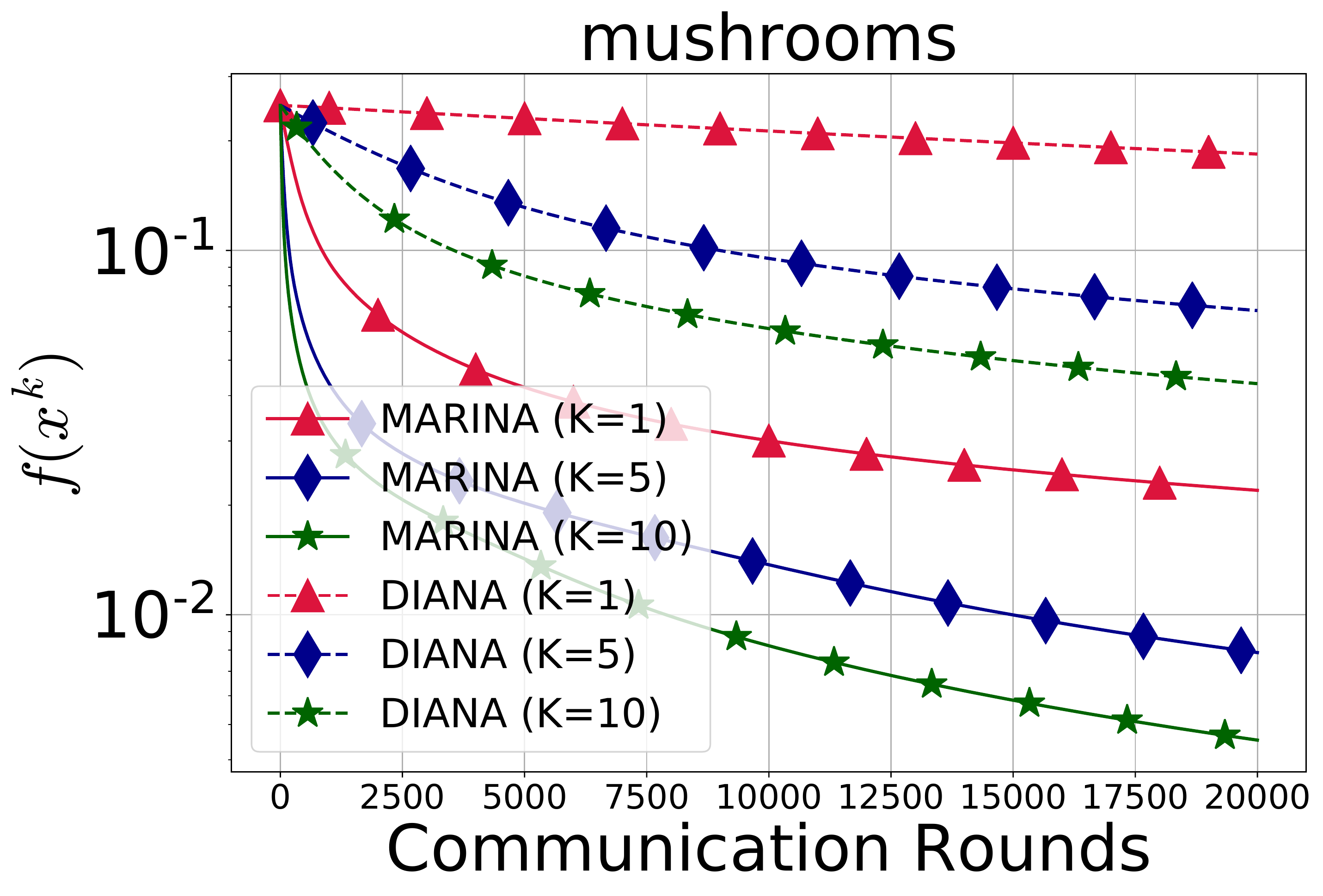}
\includegraphics[width=0.24\textwidth]{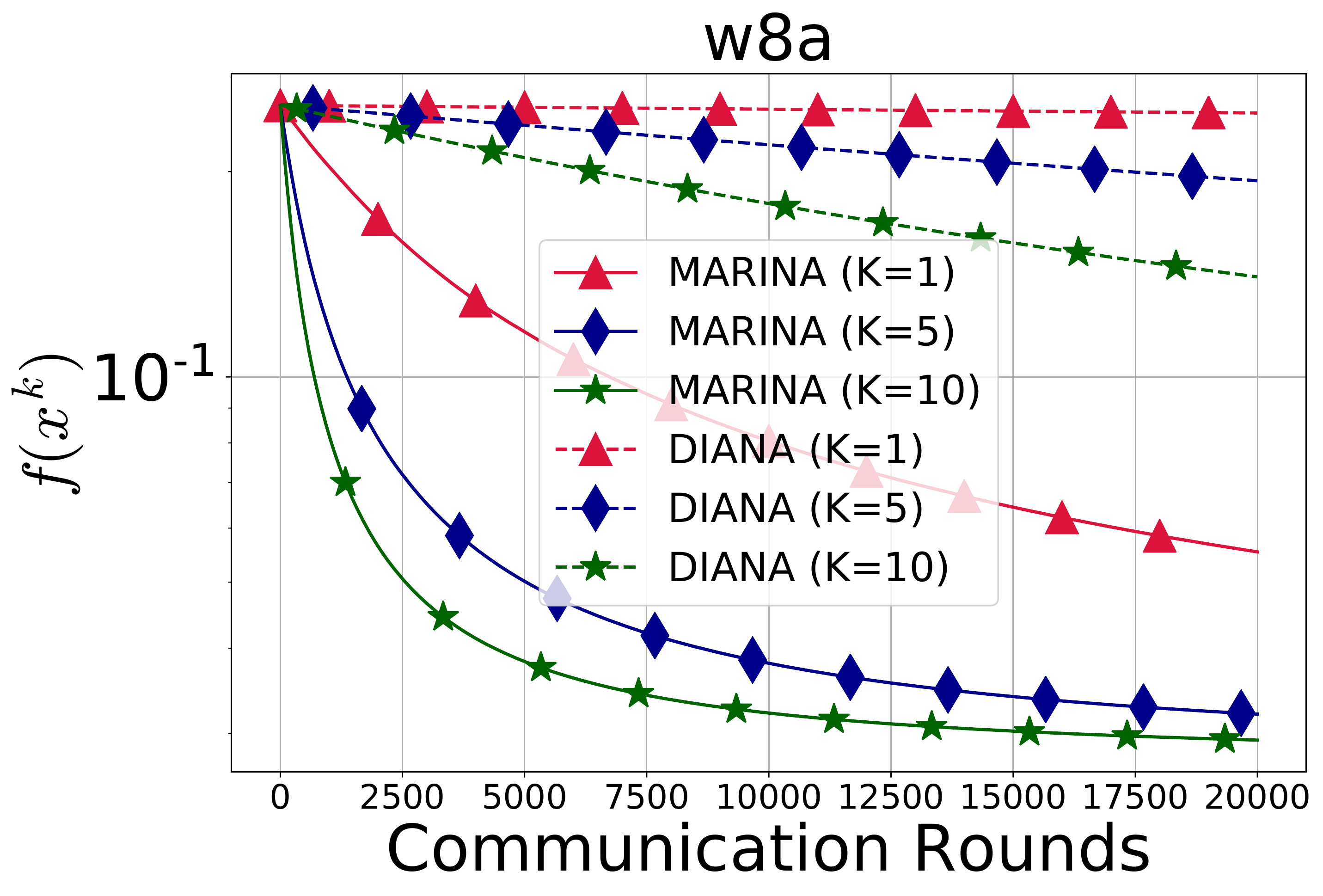}
\includegraphics[width=0.24\textwidth]{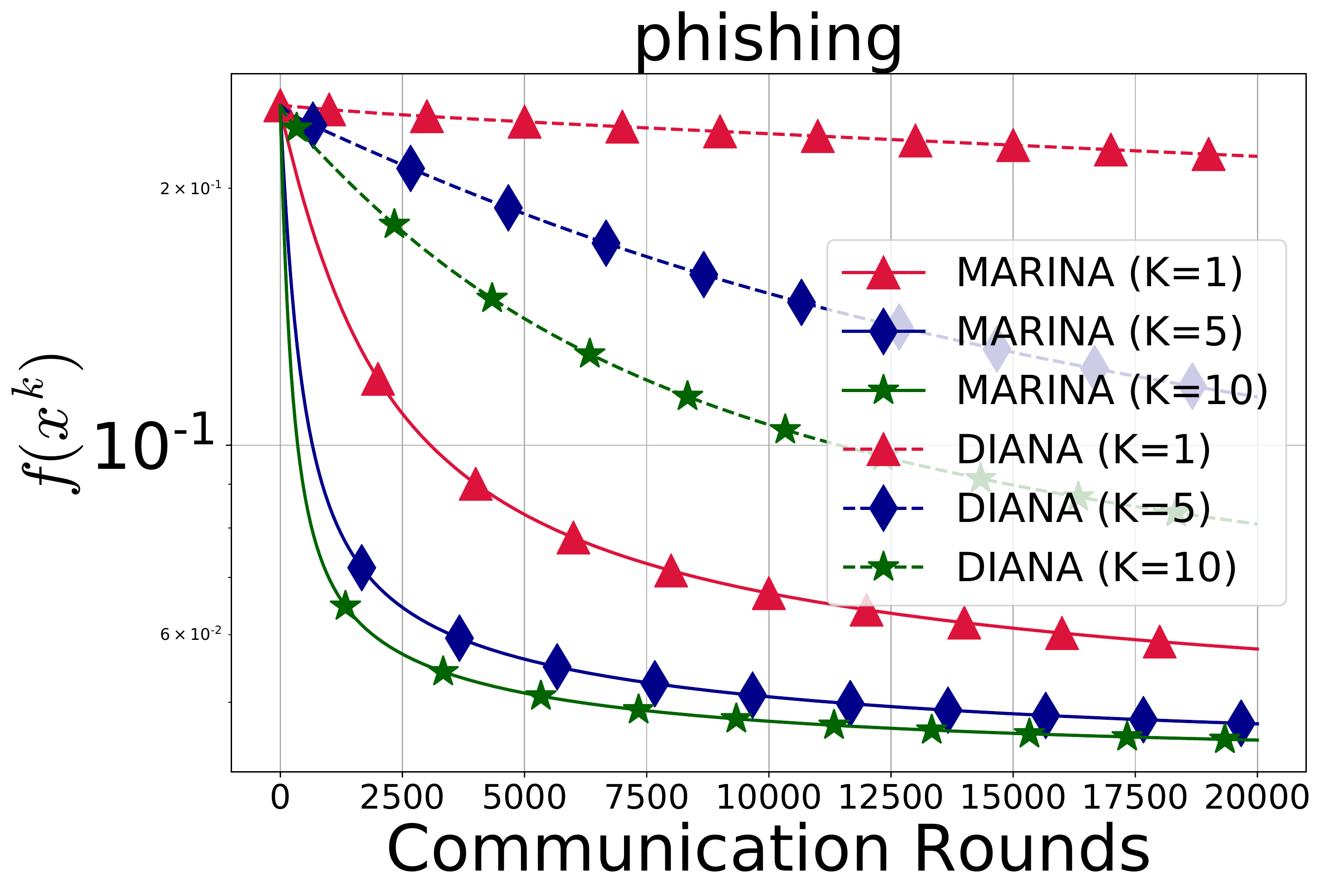}
\includegraphics[width=0.24\textwidth]{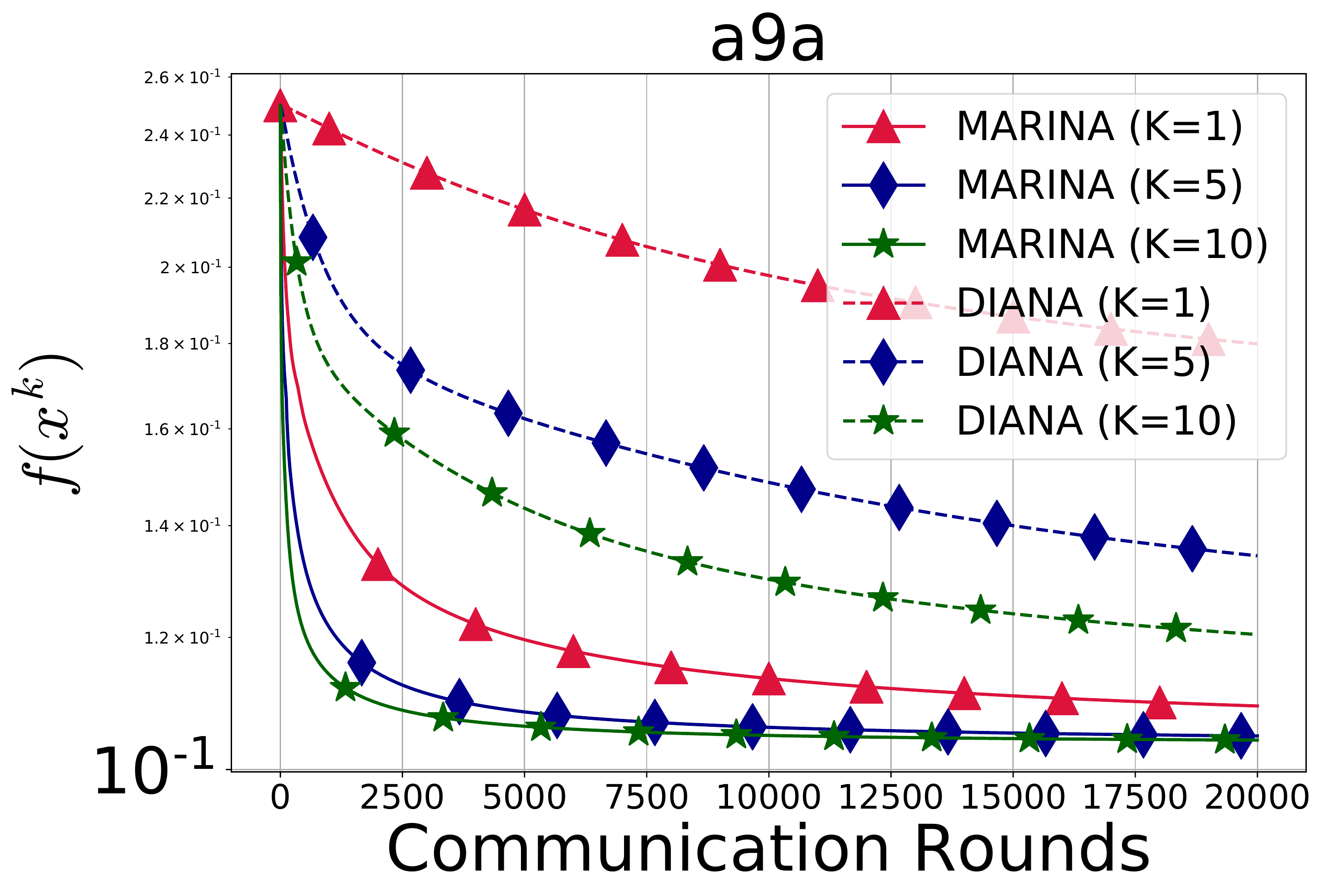}
\includegraphics[width=0.24\textwidth]{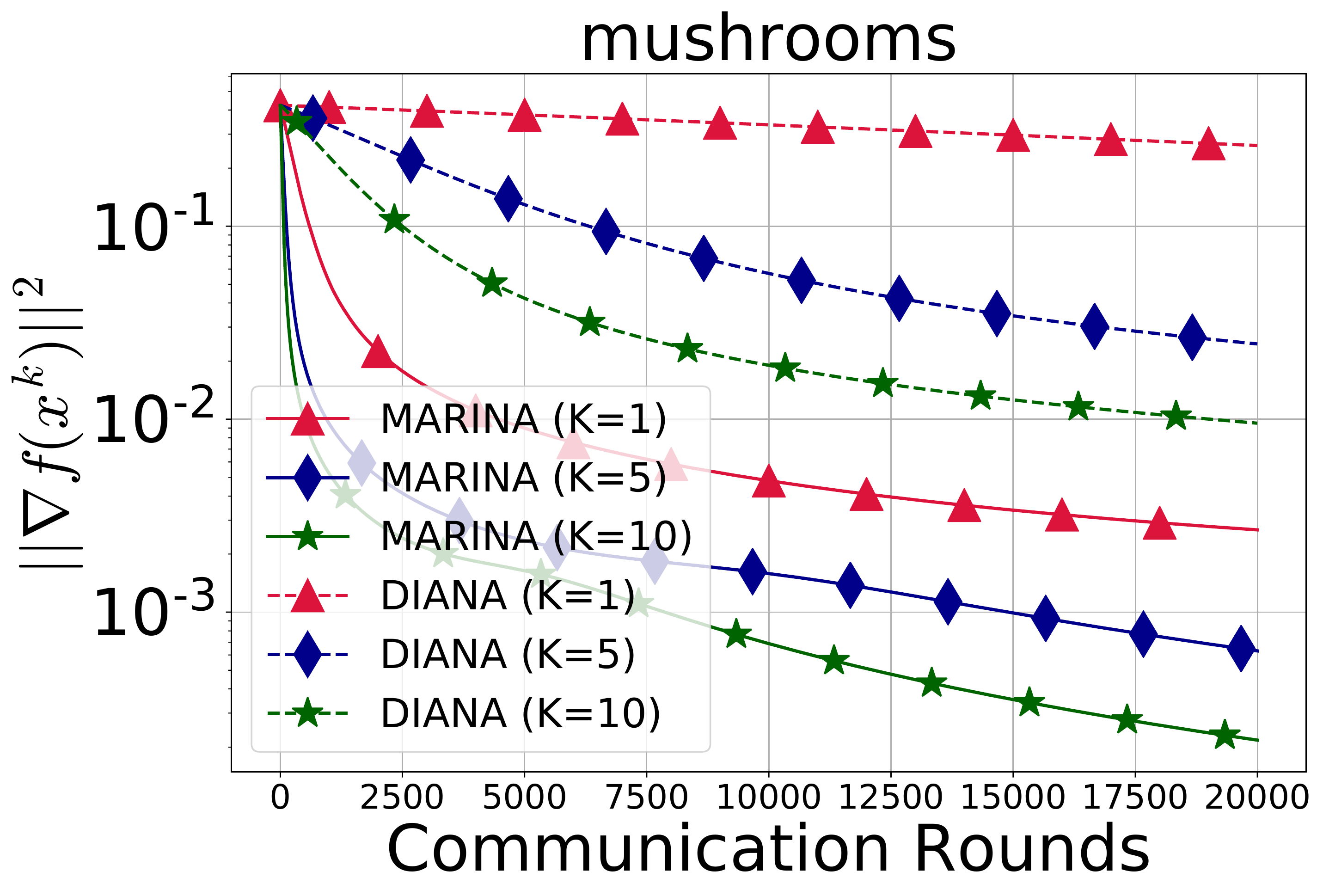}
\includegraphics[width=0.24\textwidth]{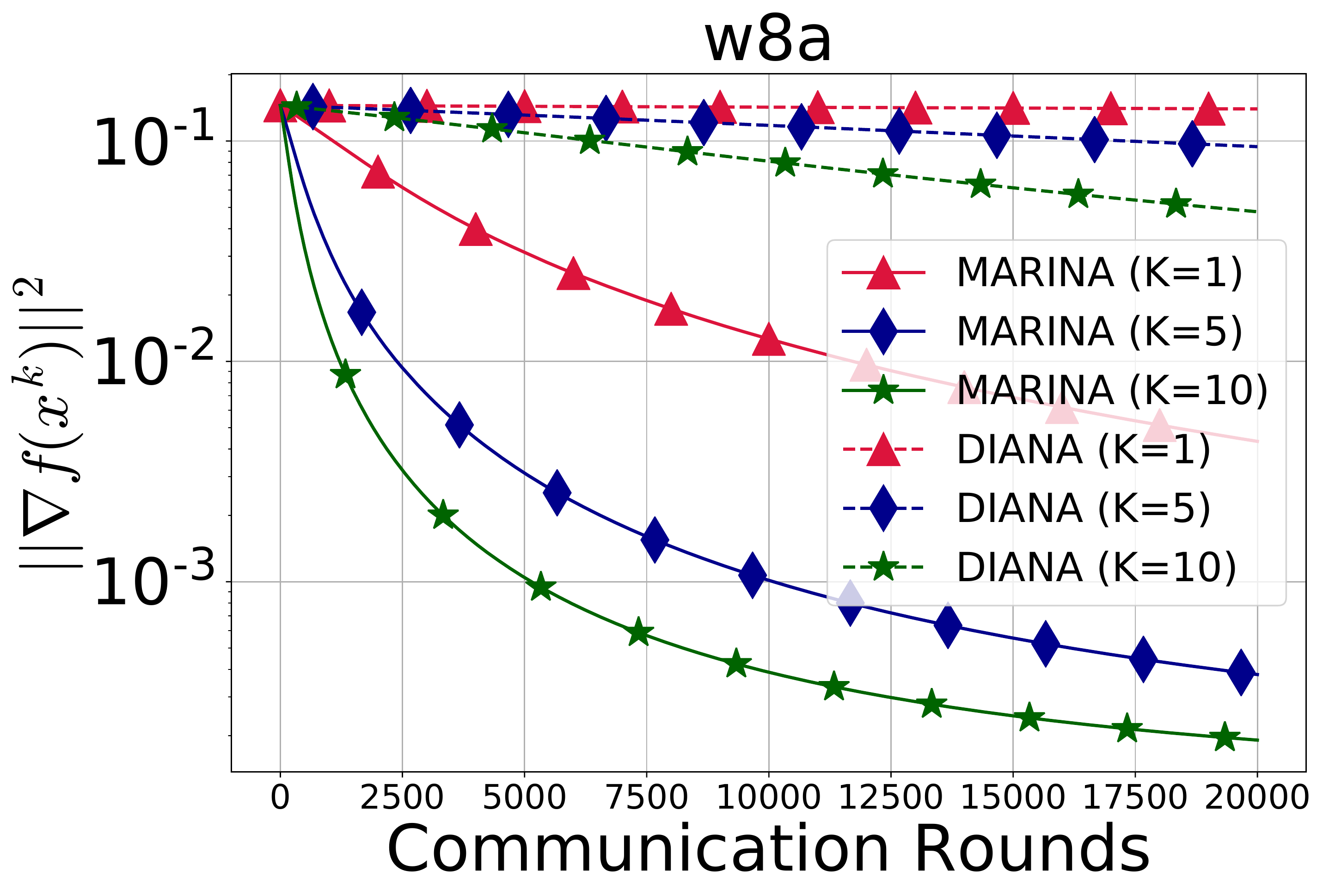}
\includegraphics[width=0.24\textwidth]{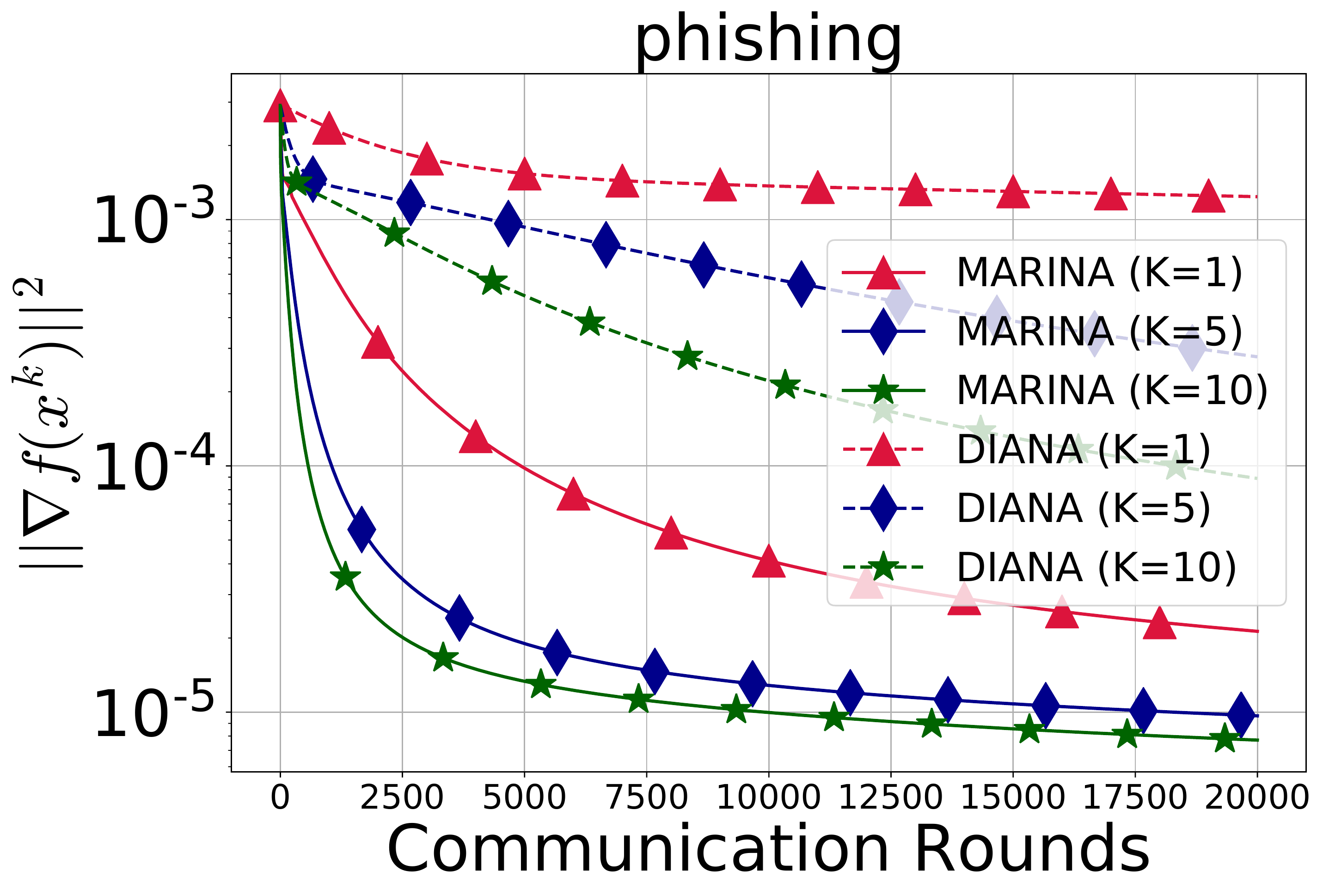}
\includegraphics[width=0.24\textwidth]{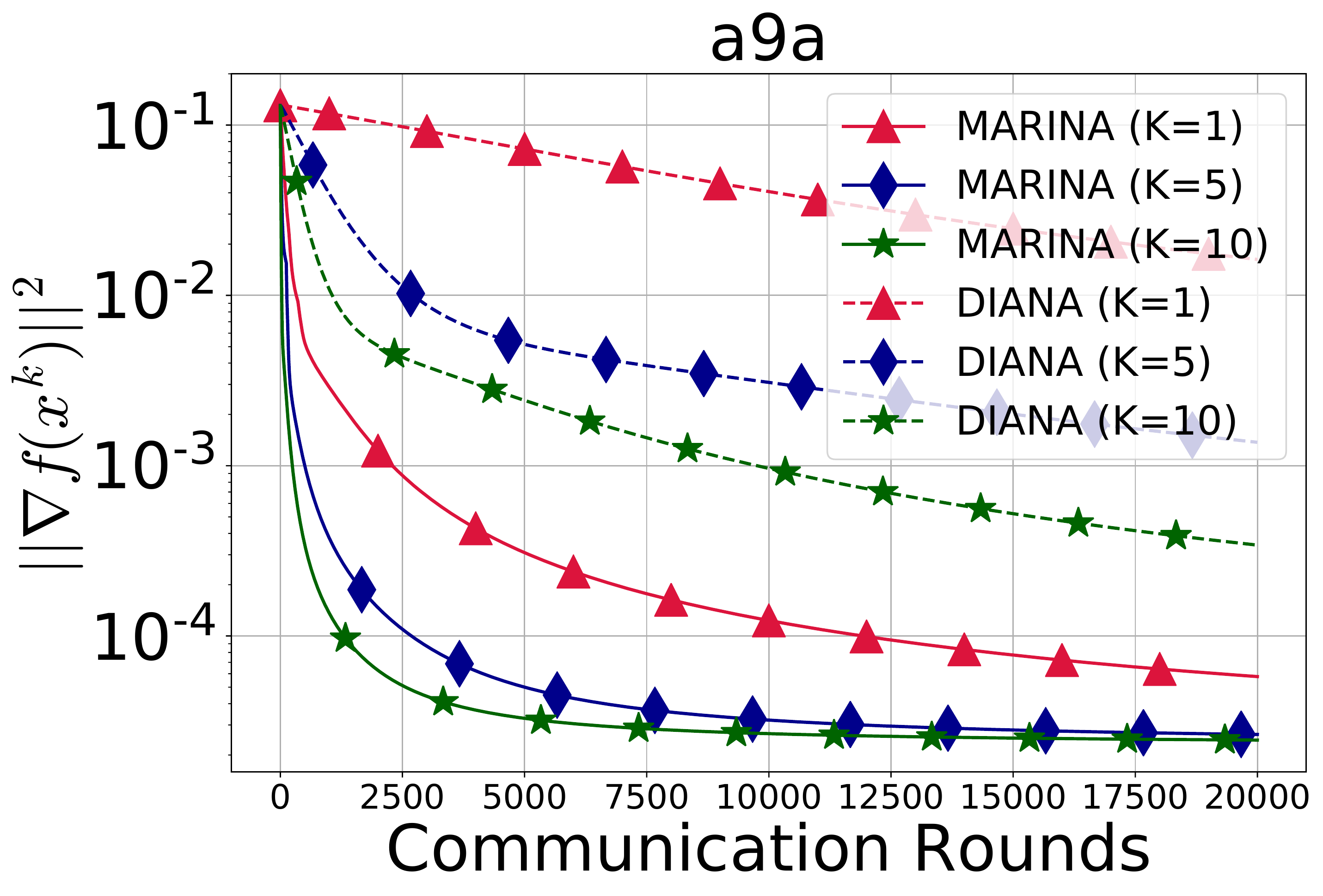}
\includegraphics[width=0.24\textwidth]{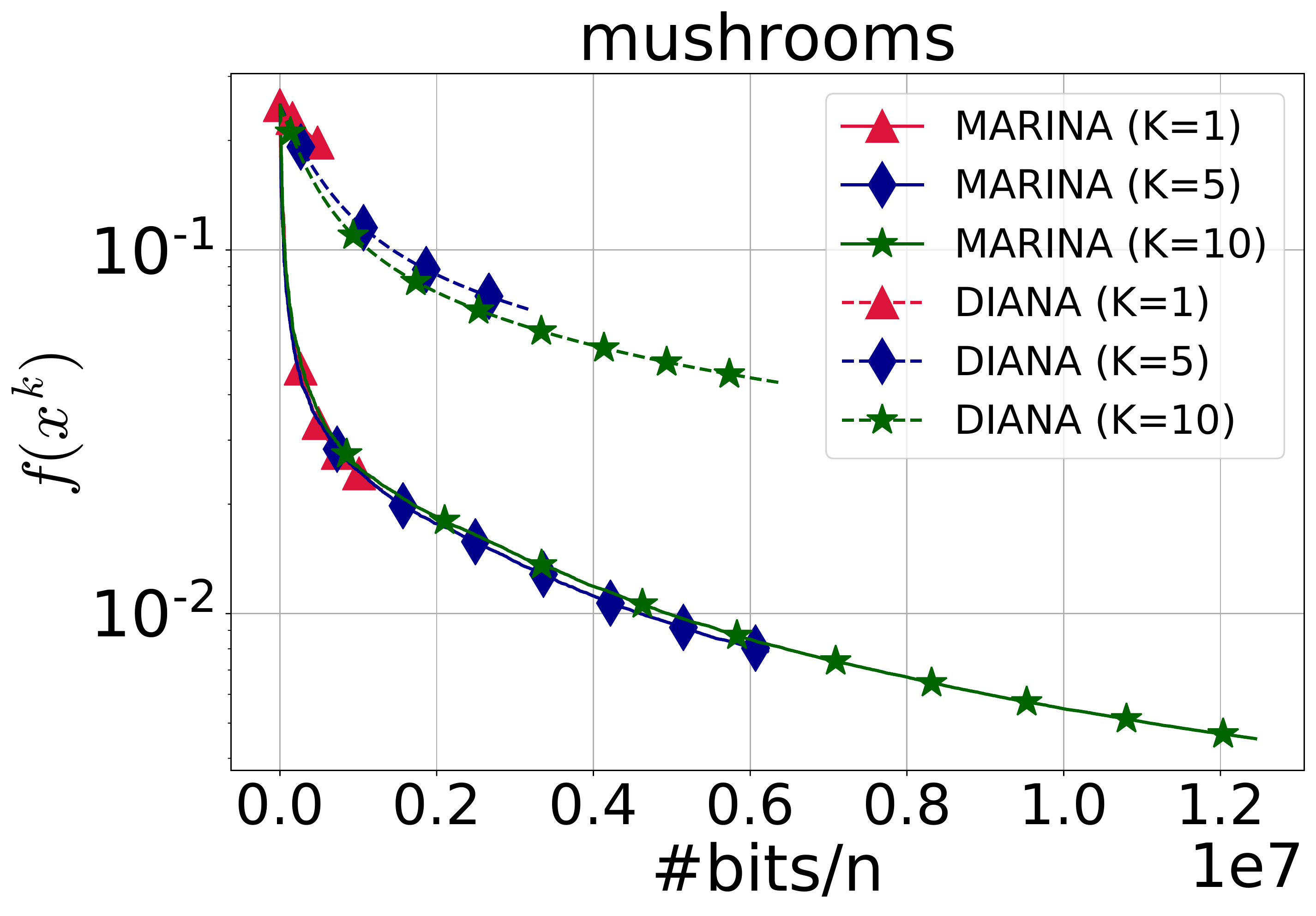}
\includegraphics[width=0.24\textwidth]{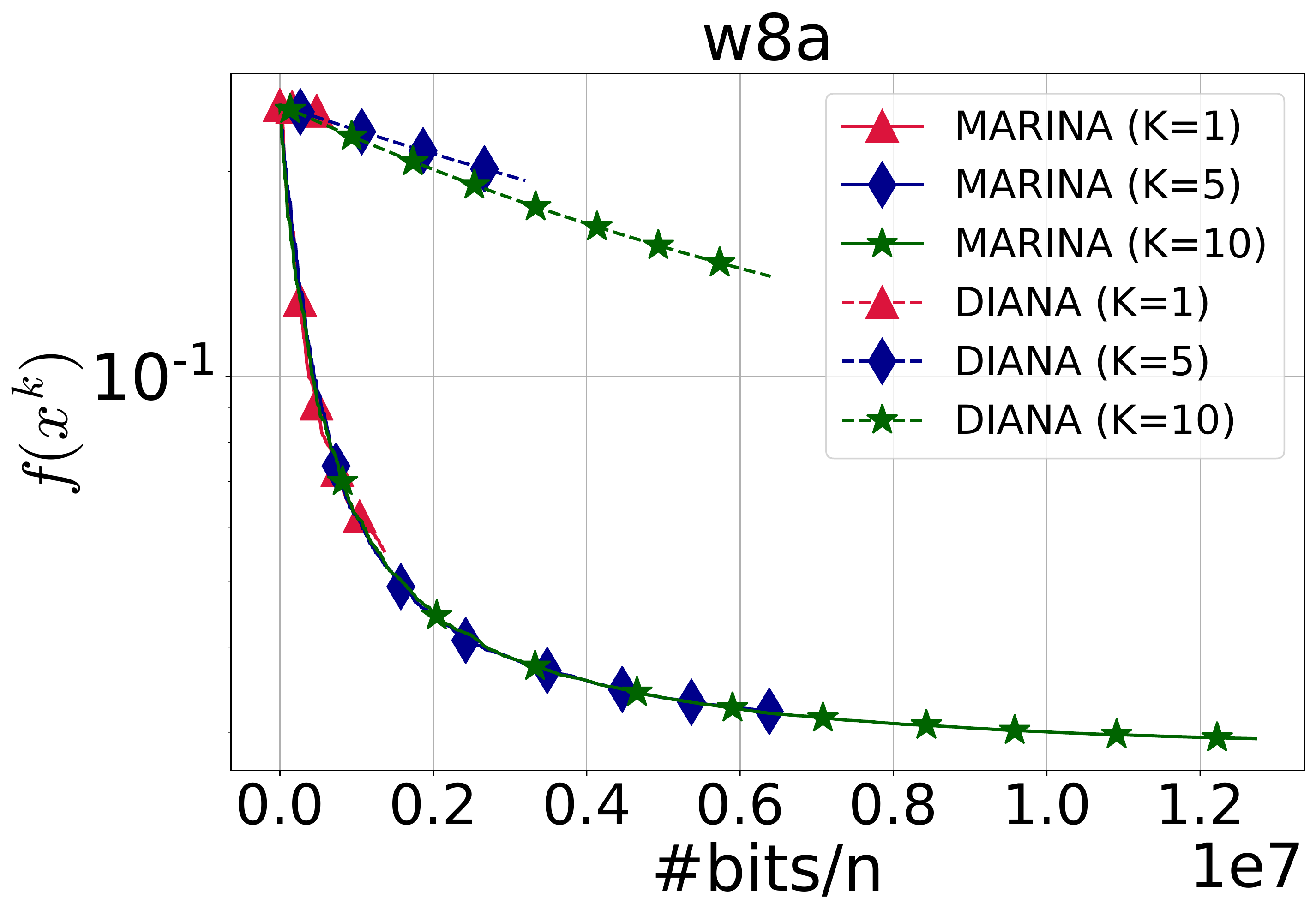}
\includegraphics[width=0.24\textwidth]{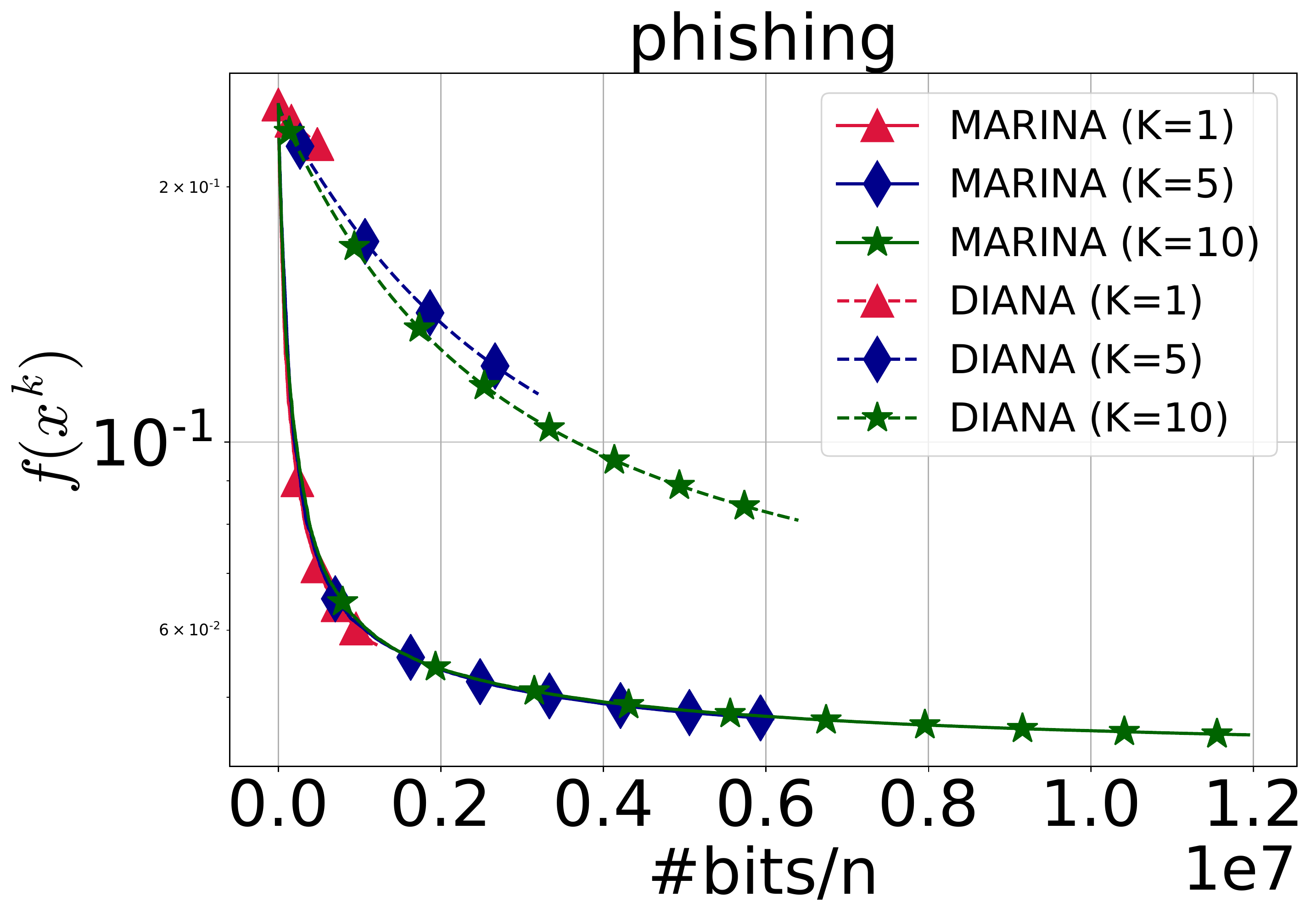}
\includegraphics[width=0.24\textwidth]{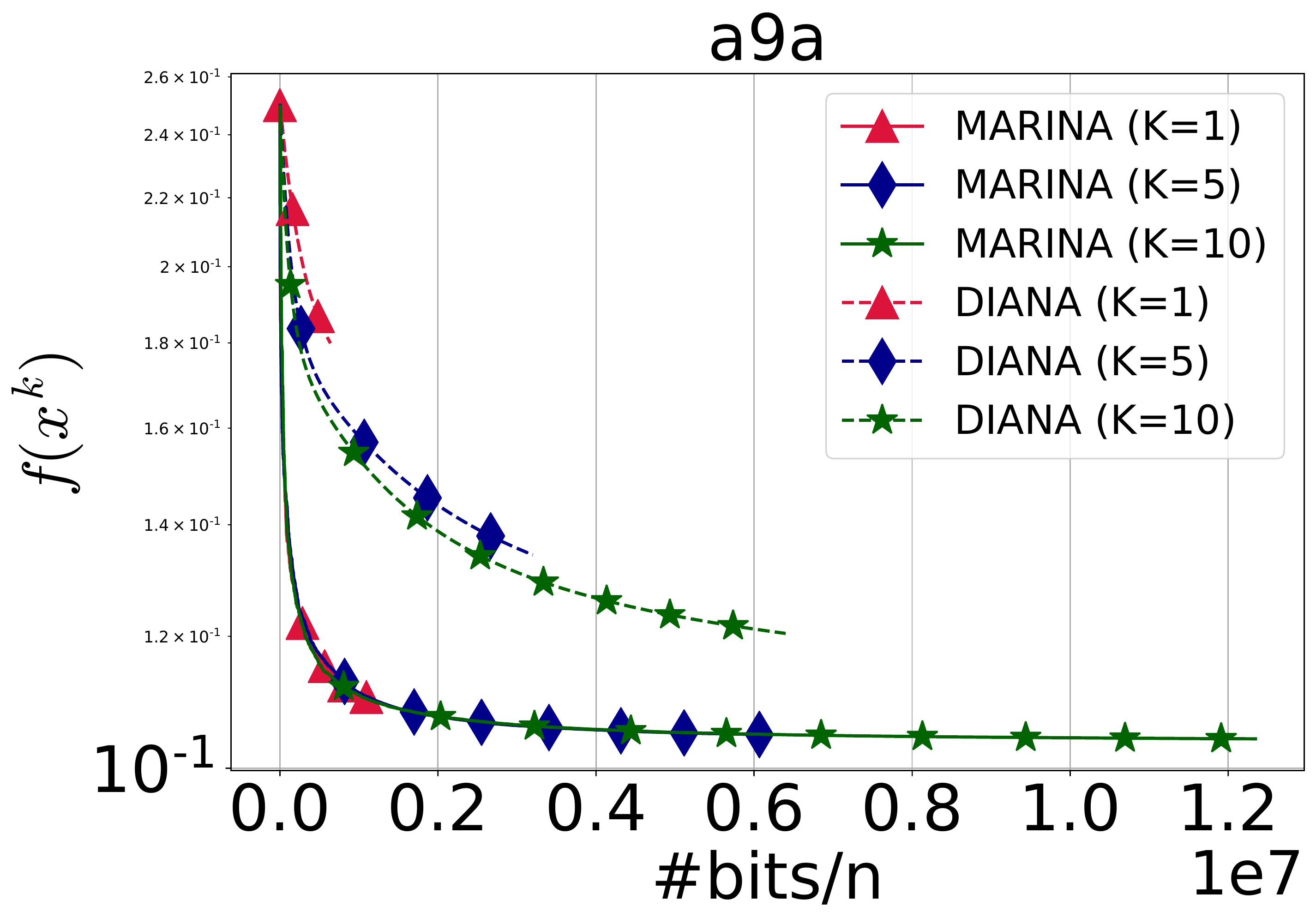}
\includegraphics[width=0.24\textwidth]{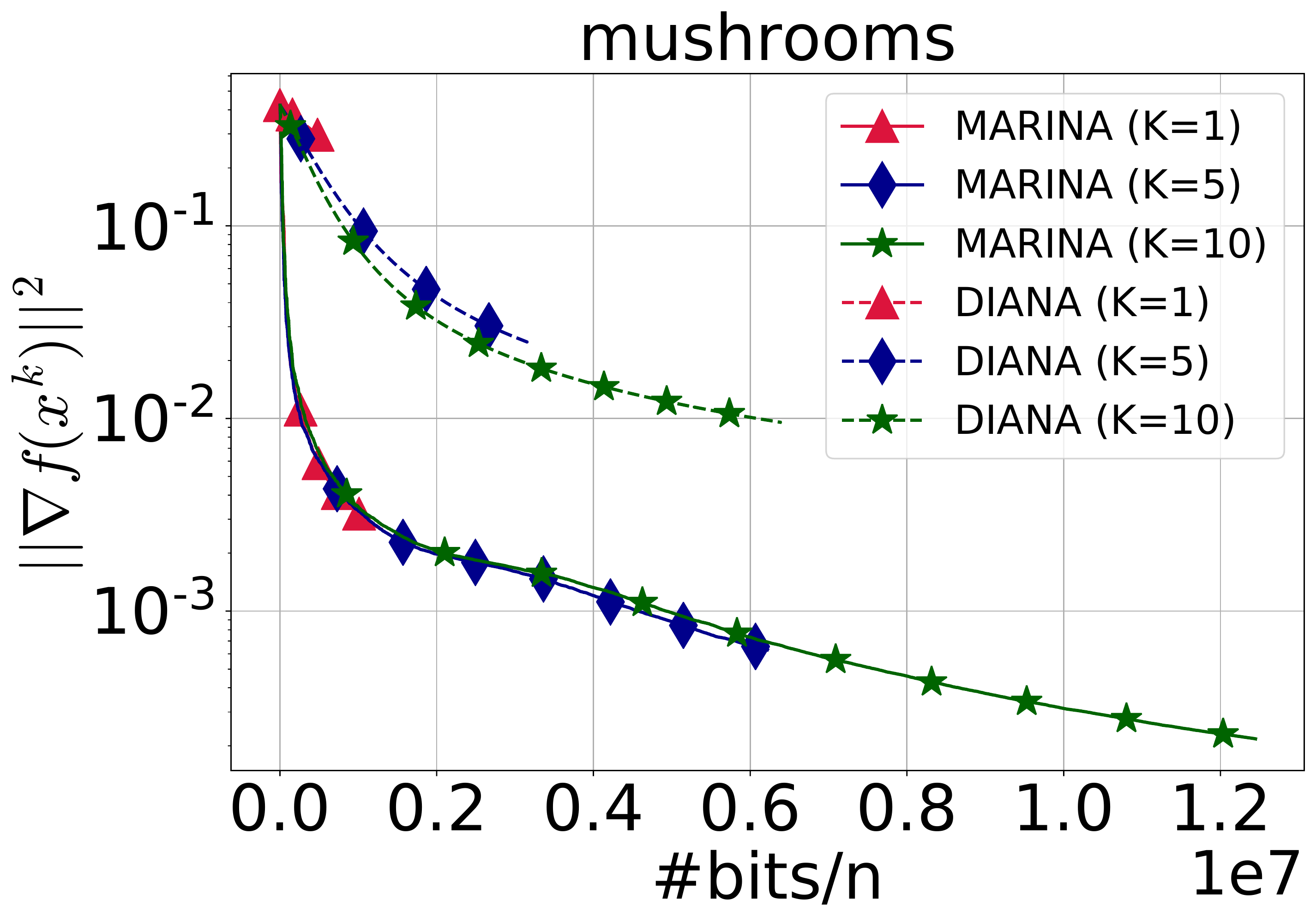}
\includegraphics[width=0.24\textwidth]{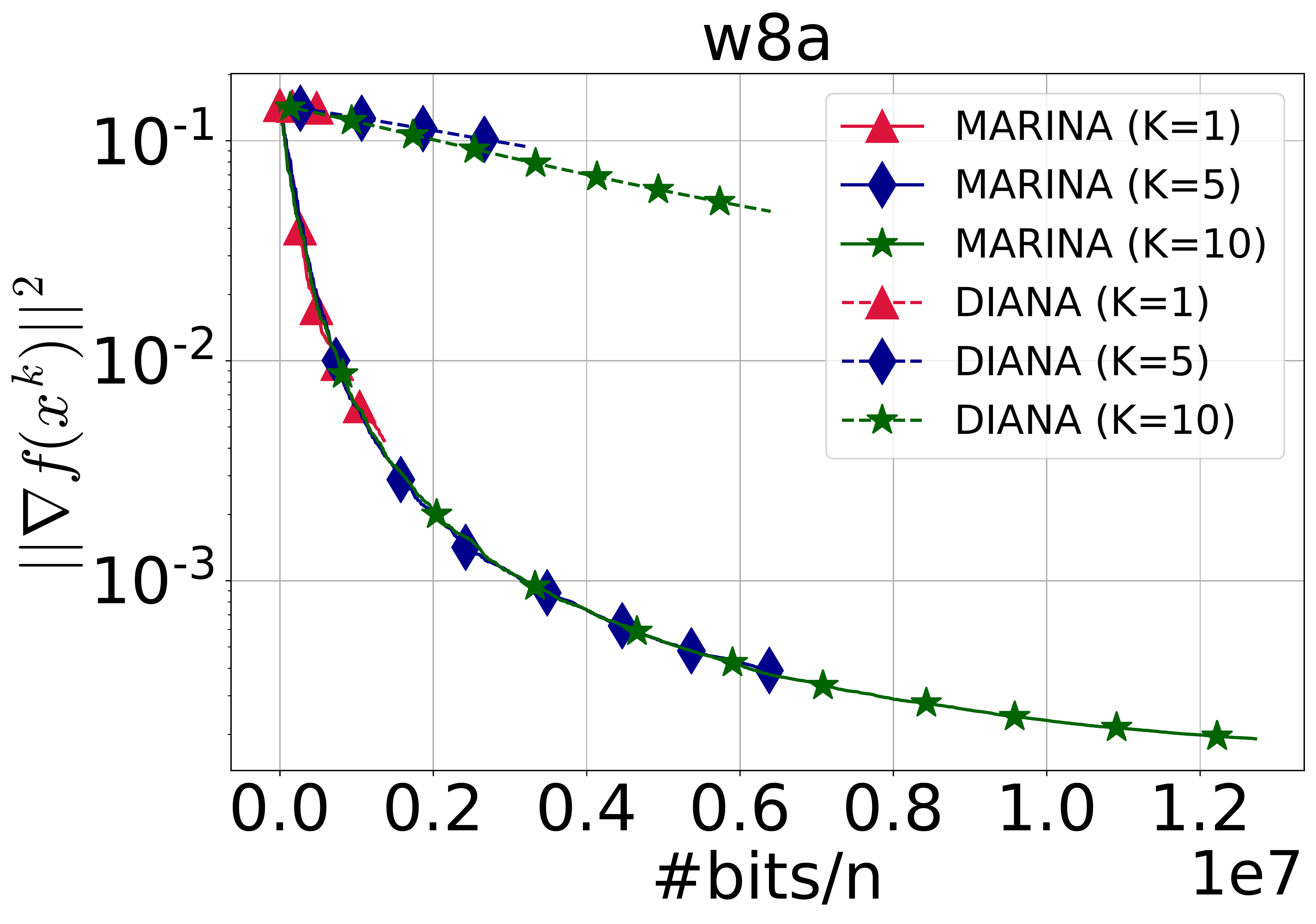}
\includegraphics[width=0.24\textwidth]{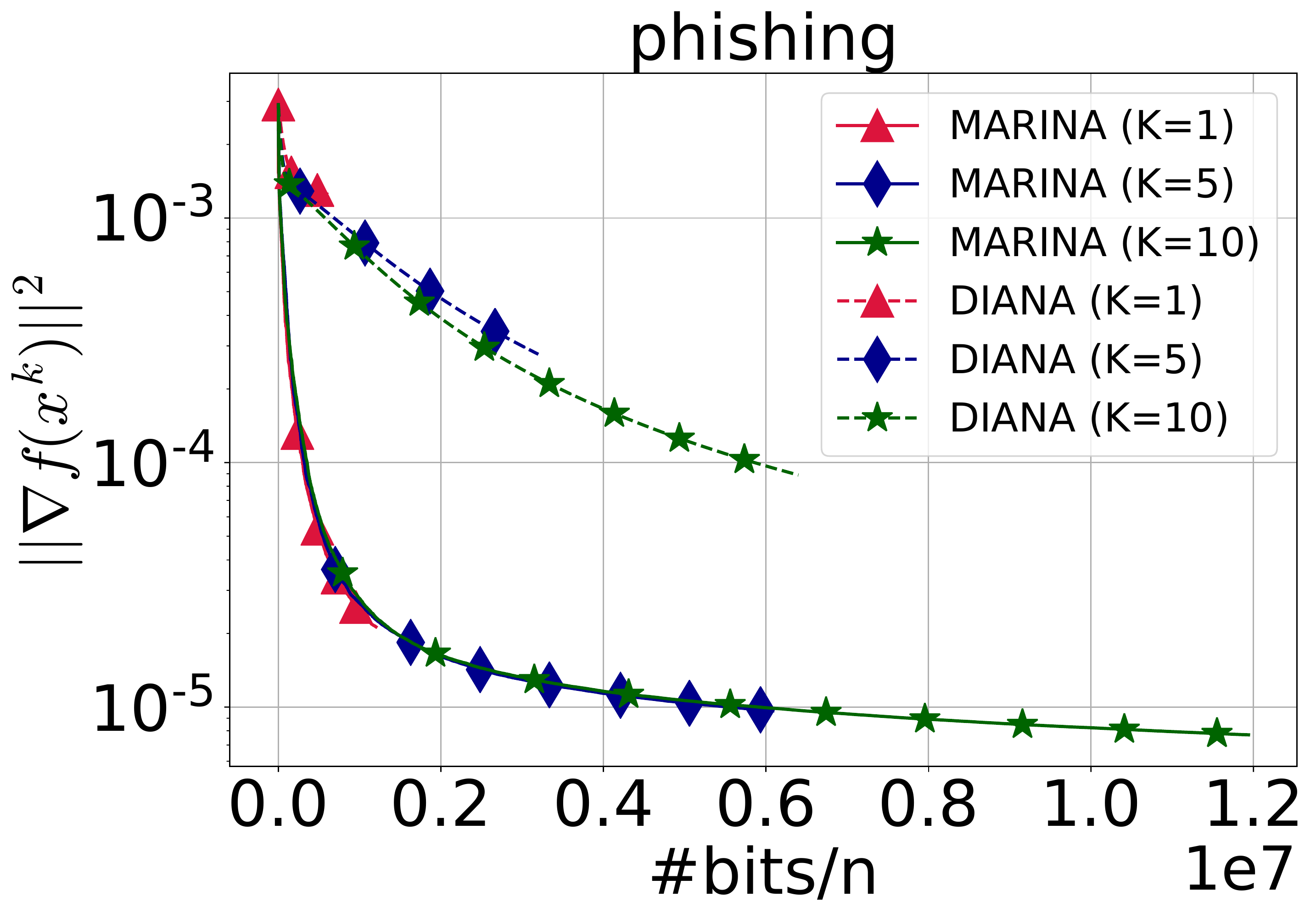}
\includegraphics[width=0.24\textwidth]{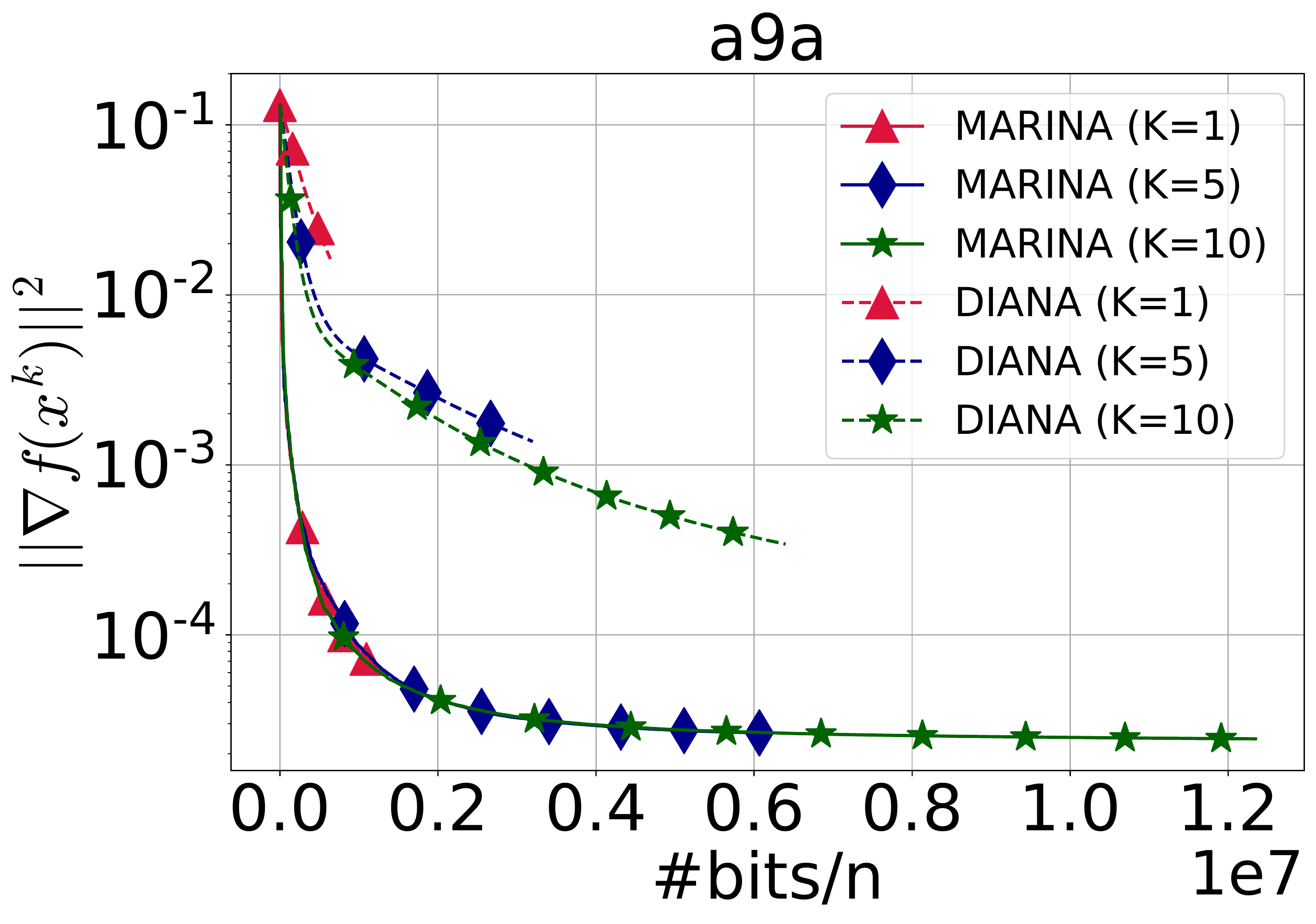}
\caption{Comparison of \algname{MARINA} with  \algname{DIANA} on binary classification problem involving non-convex loss \eqref{eq:experiment_problem} with LibSVM data \cite{chang2011libsvm}. Parameter $n$ is chosen as per Table~\ref{tbl:ns} ($n = 5$). Stepsizes for the methods are chosen according to the theory. In all cases, we used the RandK sparsification operator with K $\in \{1,5,10\}$.}
\label{fig:full_batched_methods}
\end{figure}
\begin{figure}[H]
\centering
\includegraphics[width=0.24\textwidth]{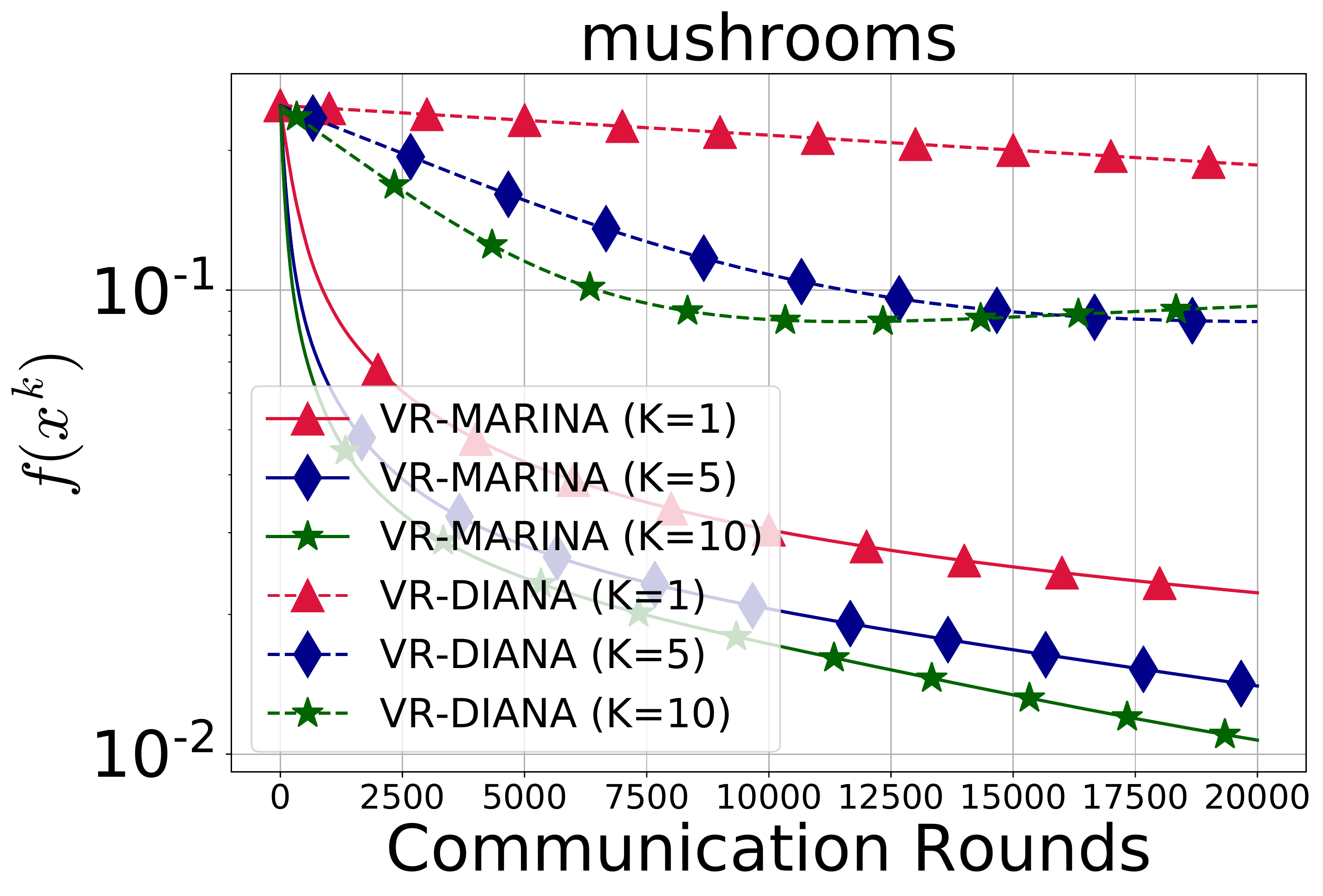}
\includegraphics[width=0.24\textwidth]{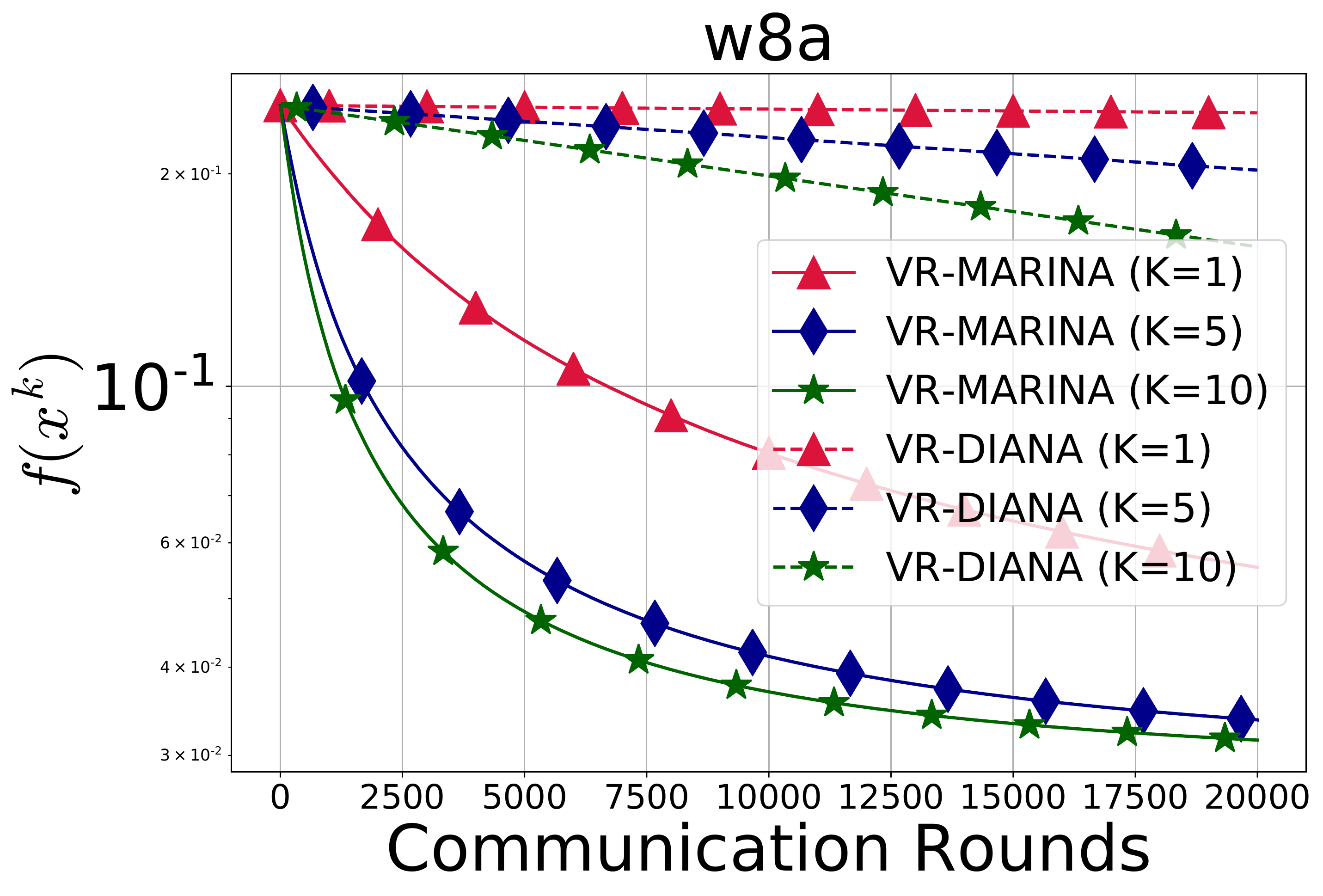}
\includegraphics[width=0.24\textwidth]{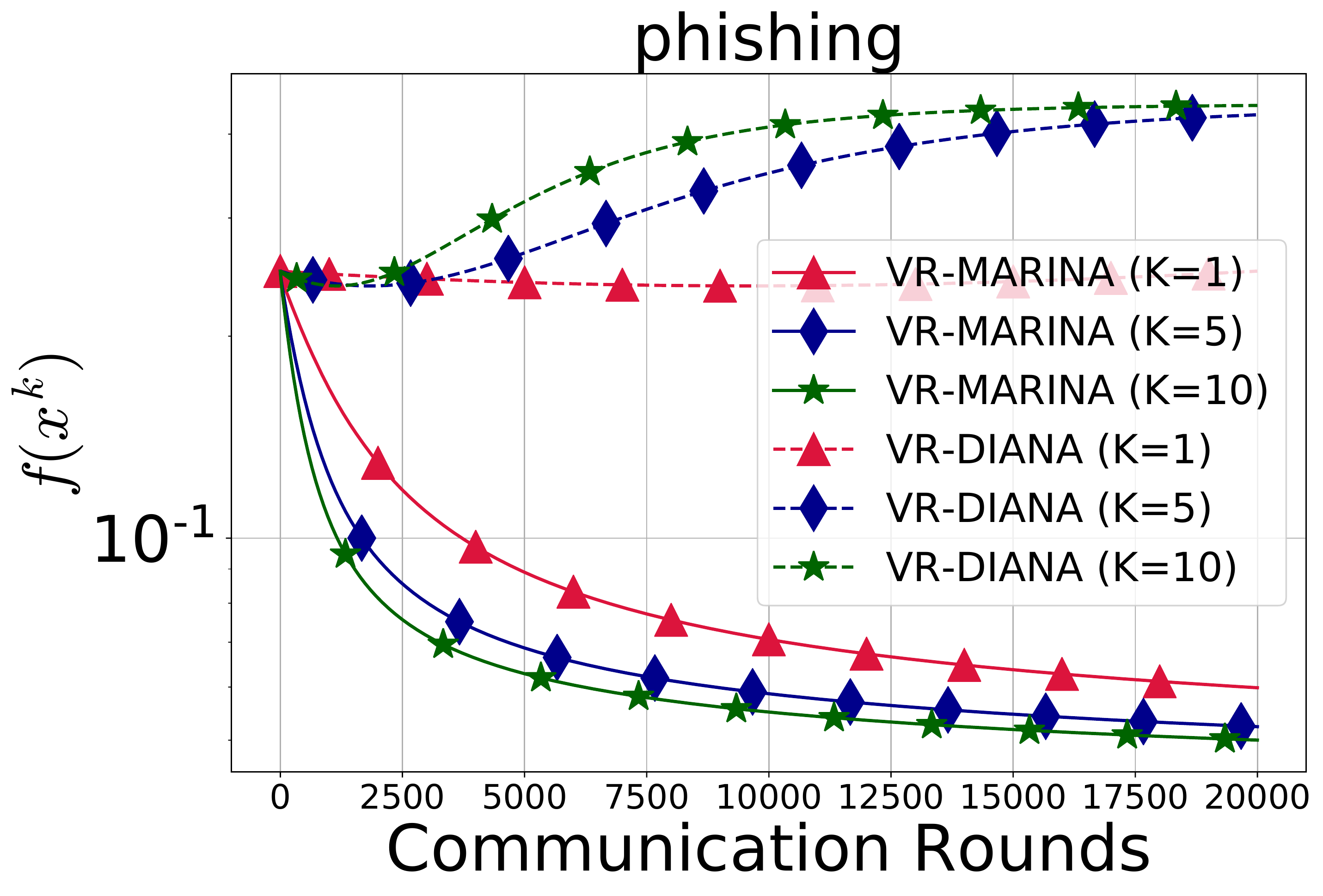}
\includegraphics[width=0.24\textwidth]{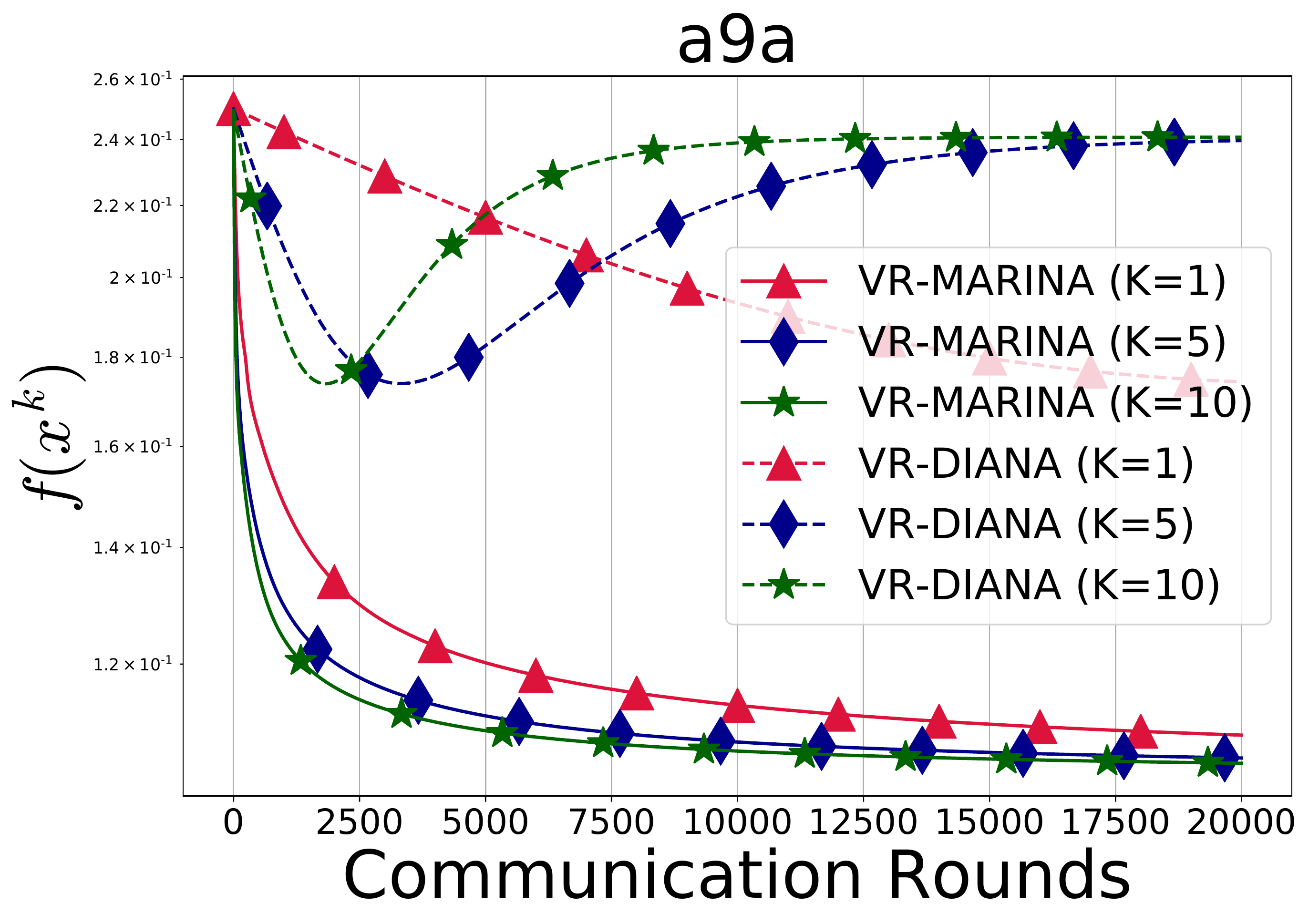}
\includegraphics[width=0.24\textwidth]{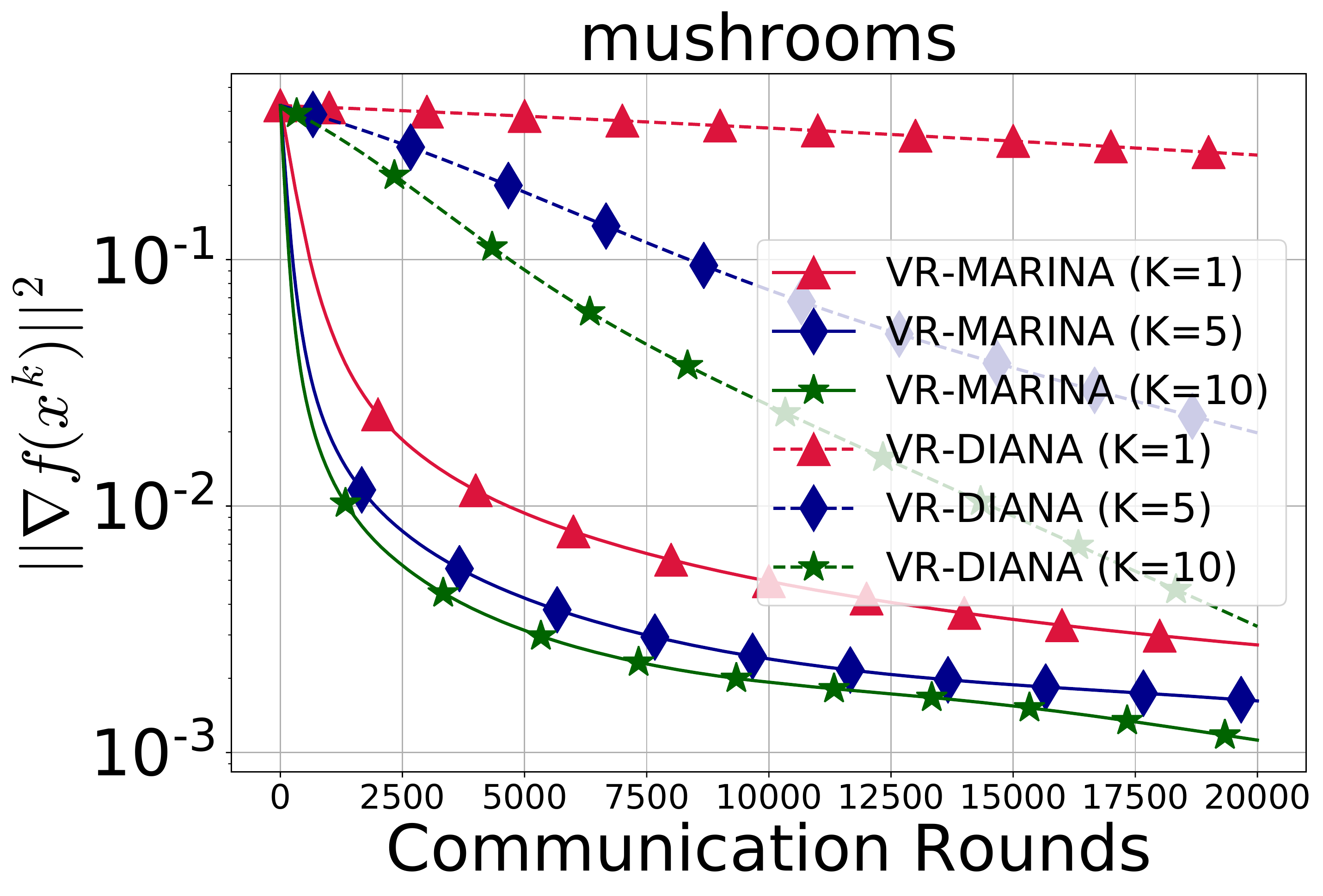}
\includegraphics[width=0.24\textwidth]{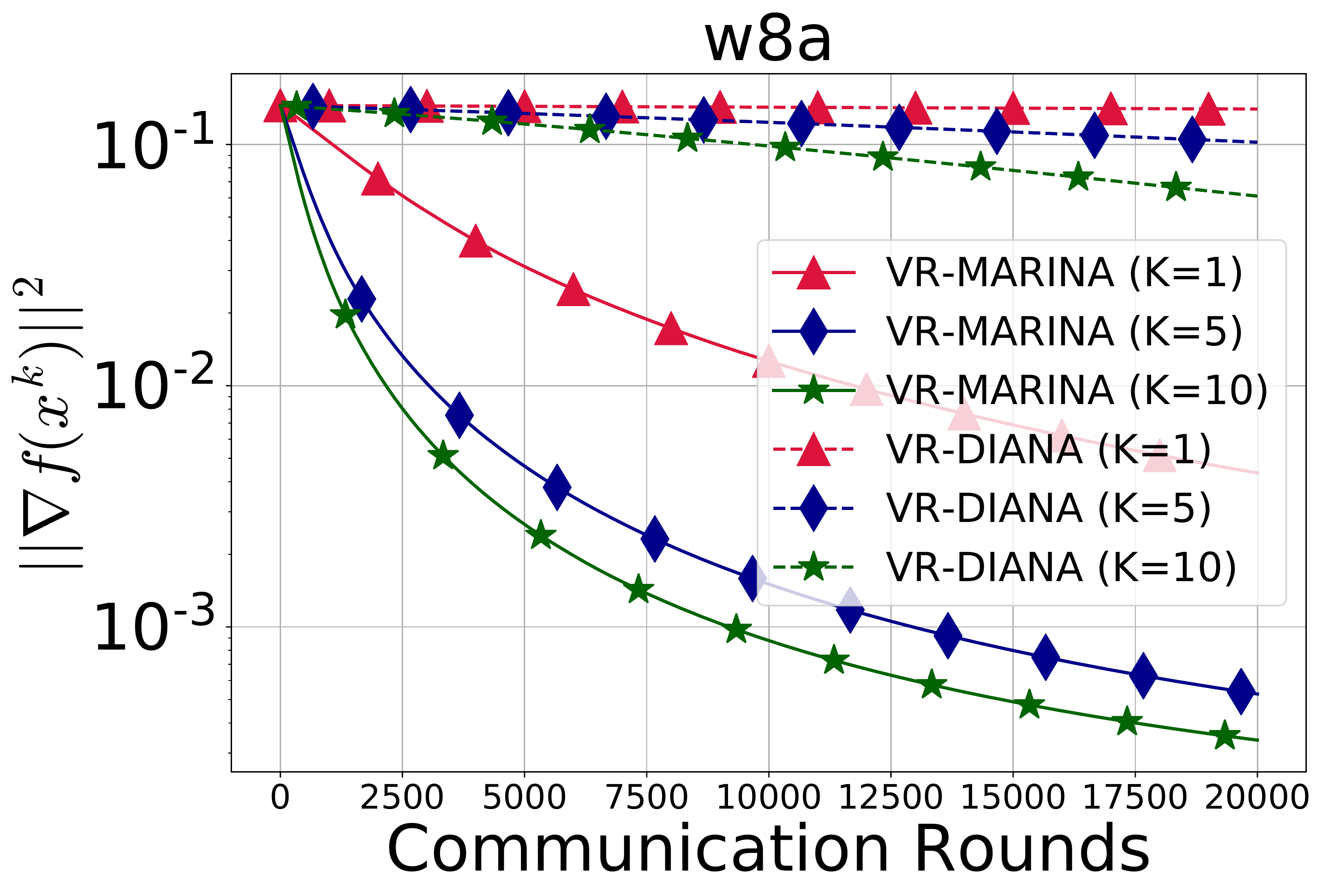}
\includegraphics[width=0.24\textwidth]{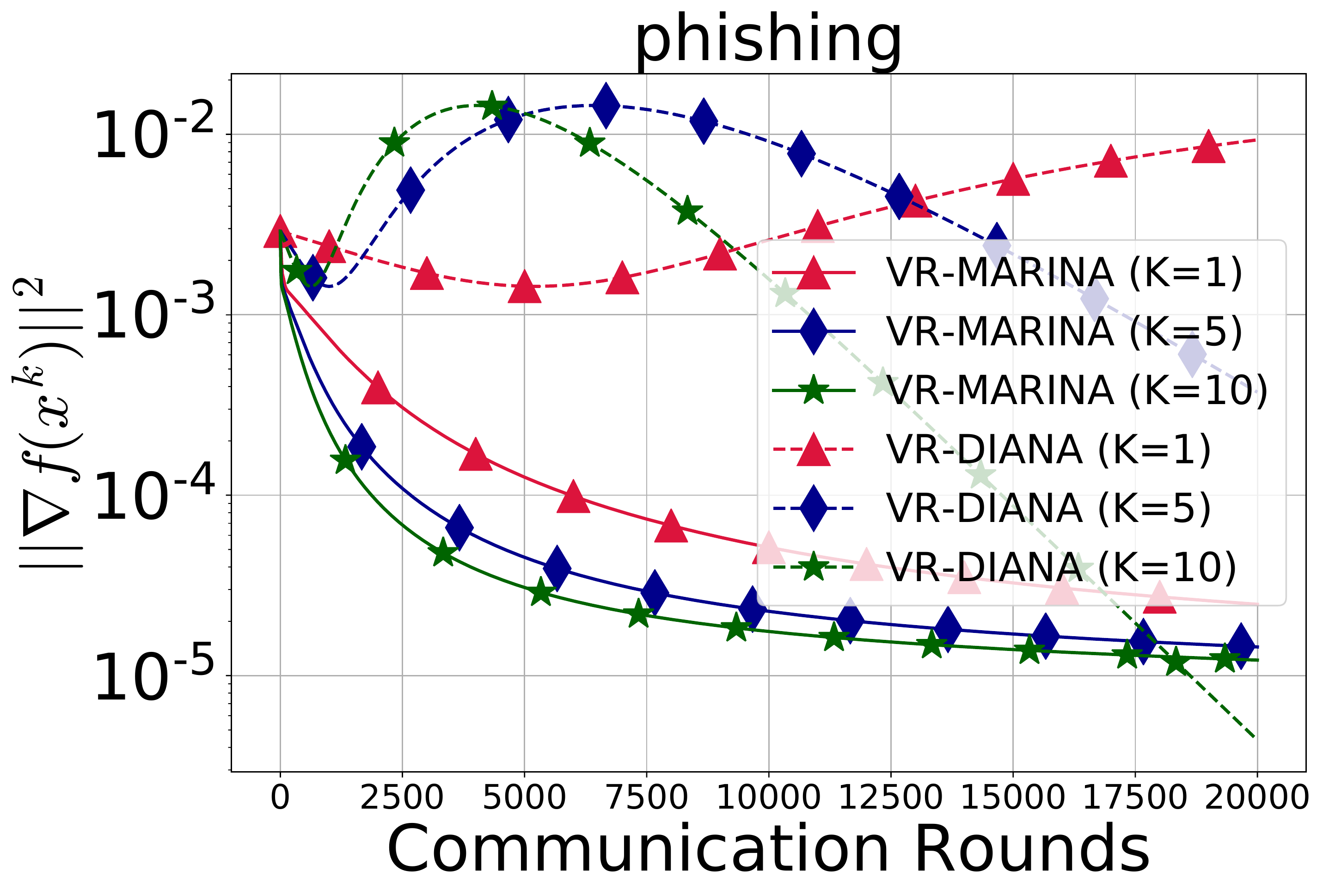}
\includegraphics[width=0.24\textwidth]{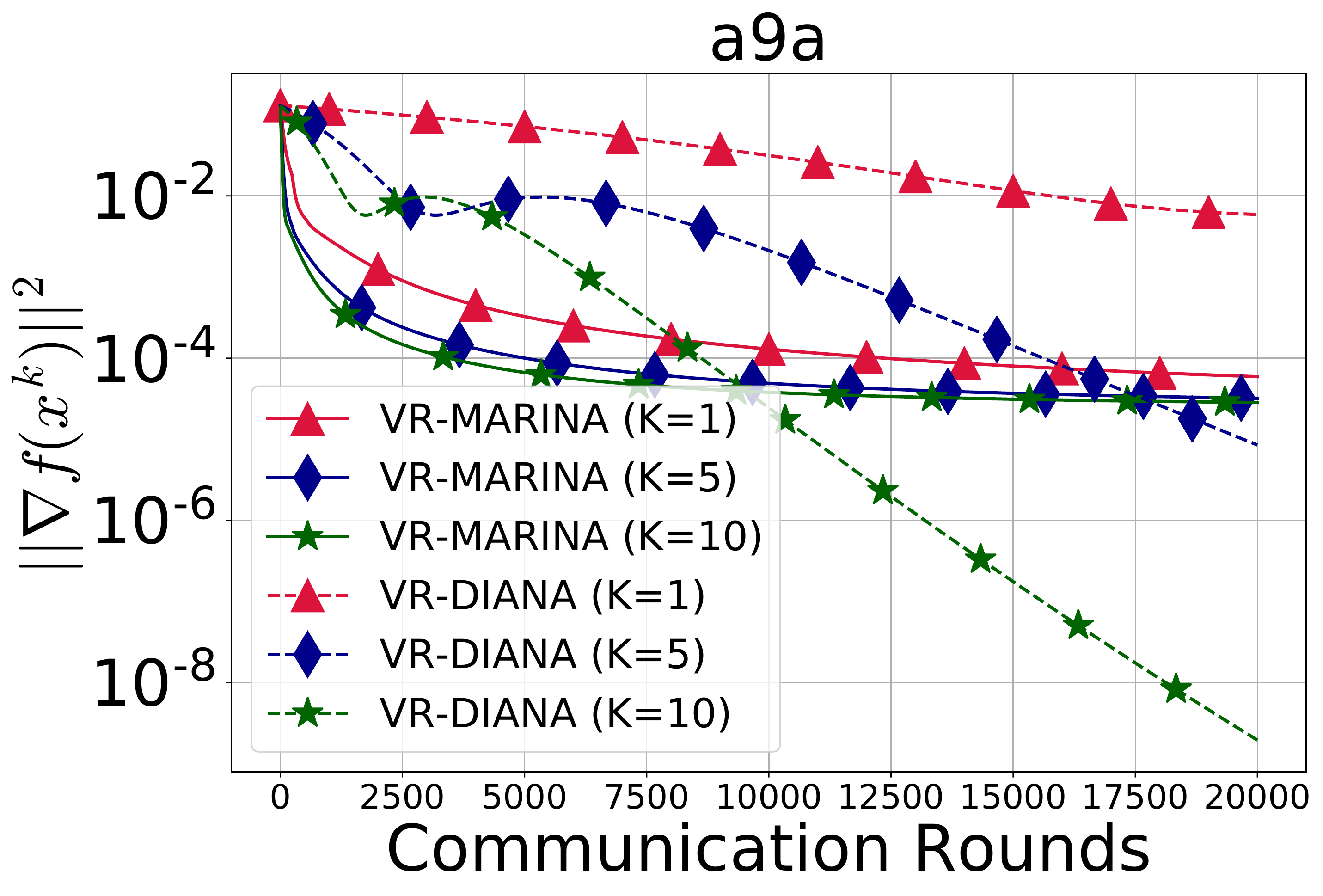}
\includegraphics[width=0.24\textwidth]{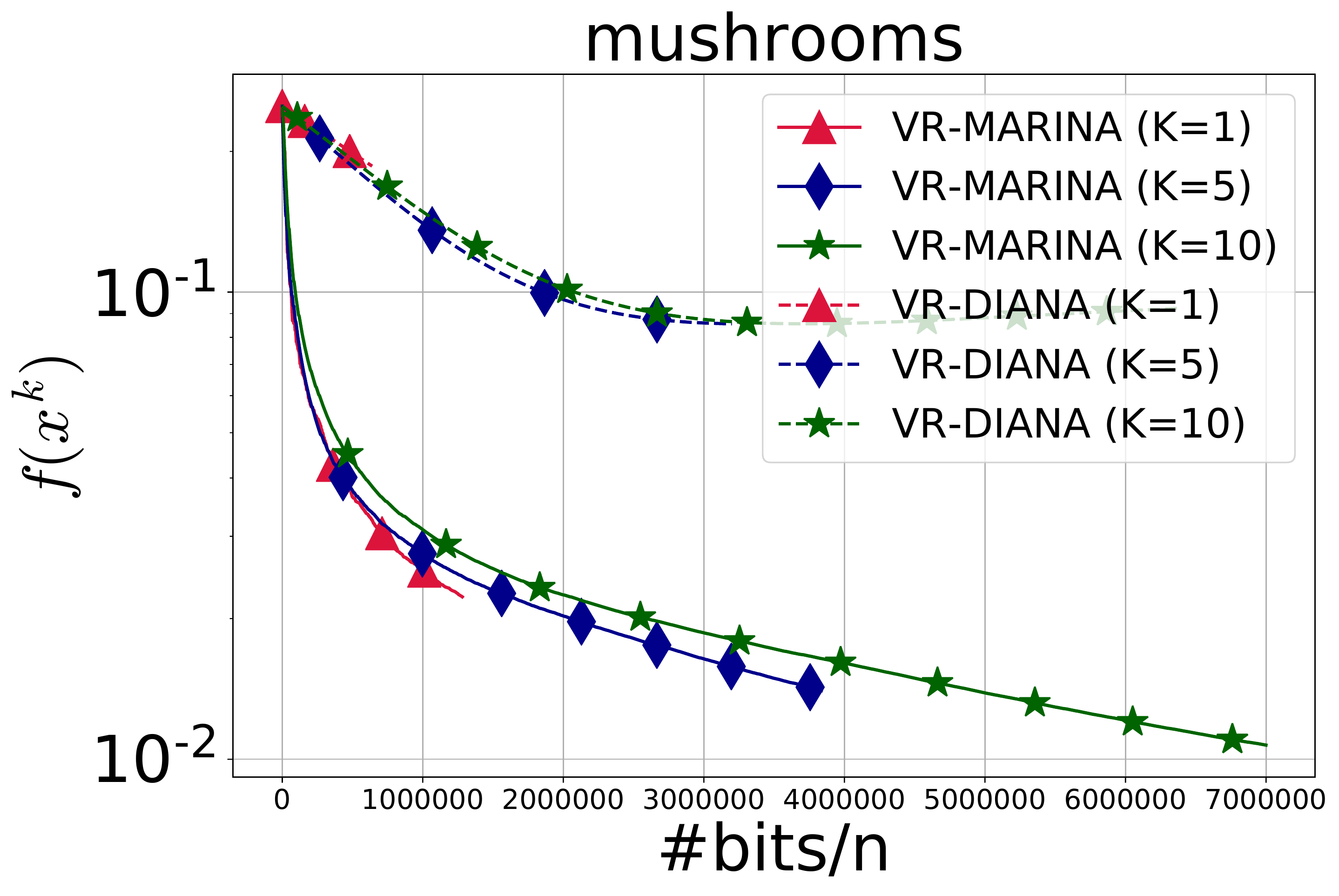}
\includegraphics[width=0.24\textwidth]{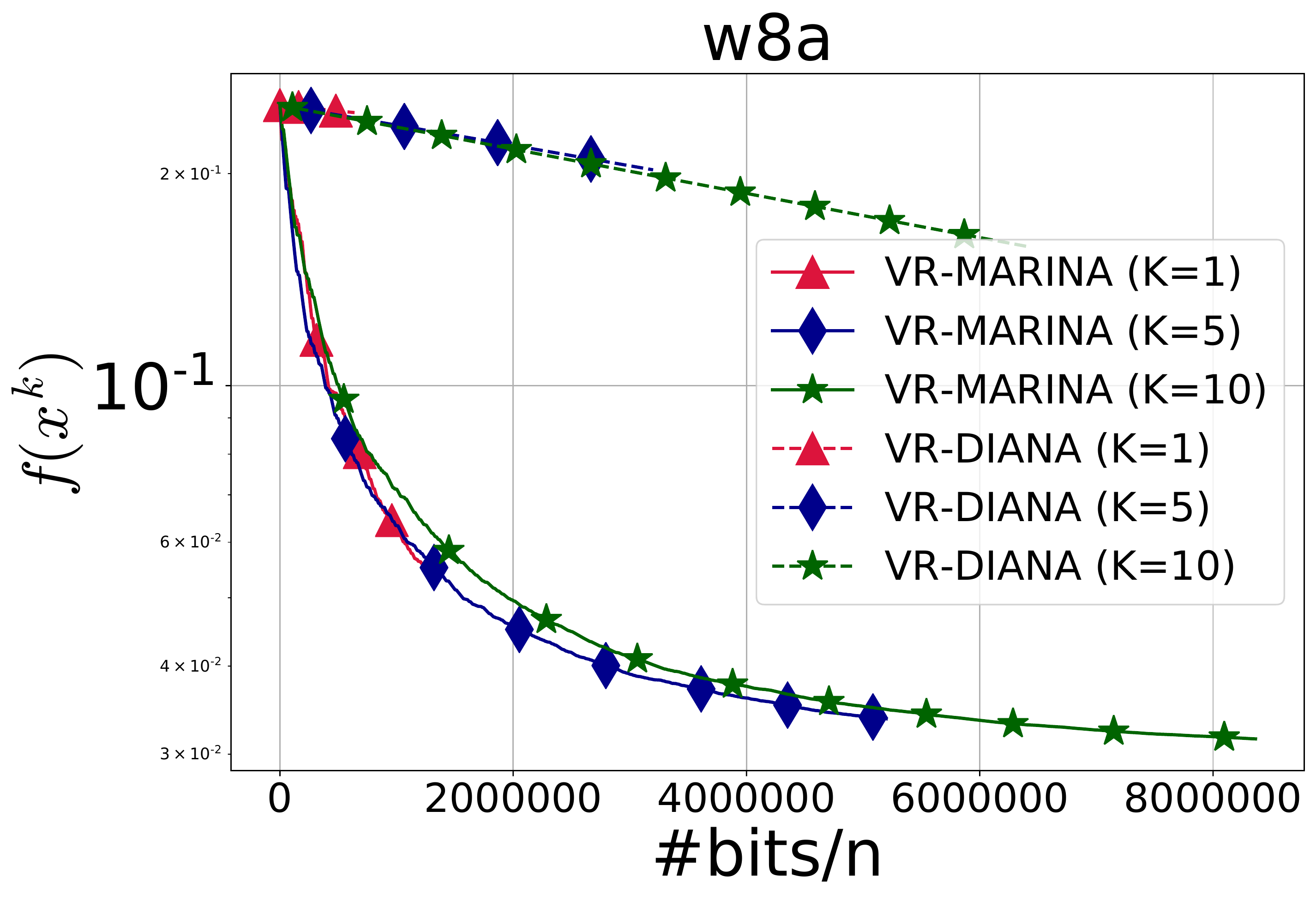}
\includegraphics[width=0.24\textwidth]{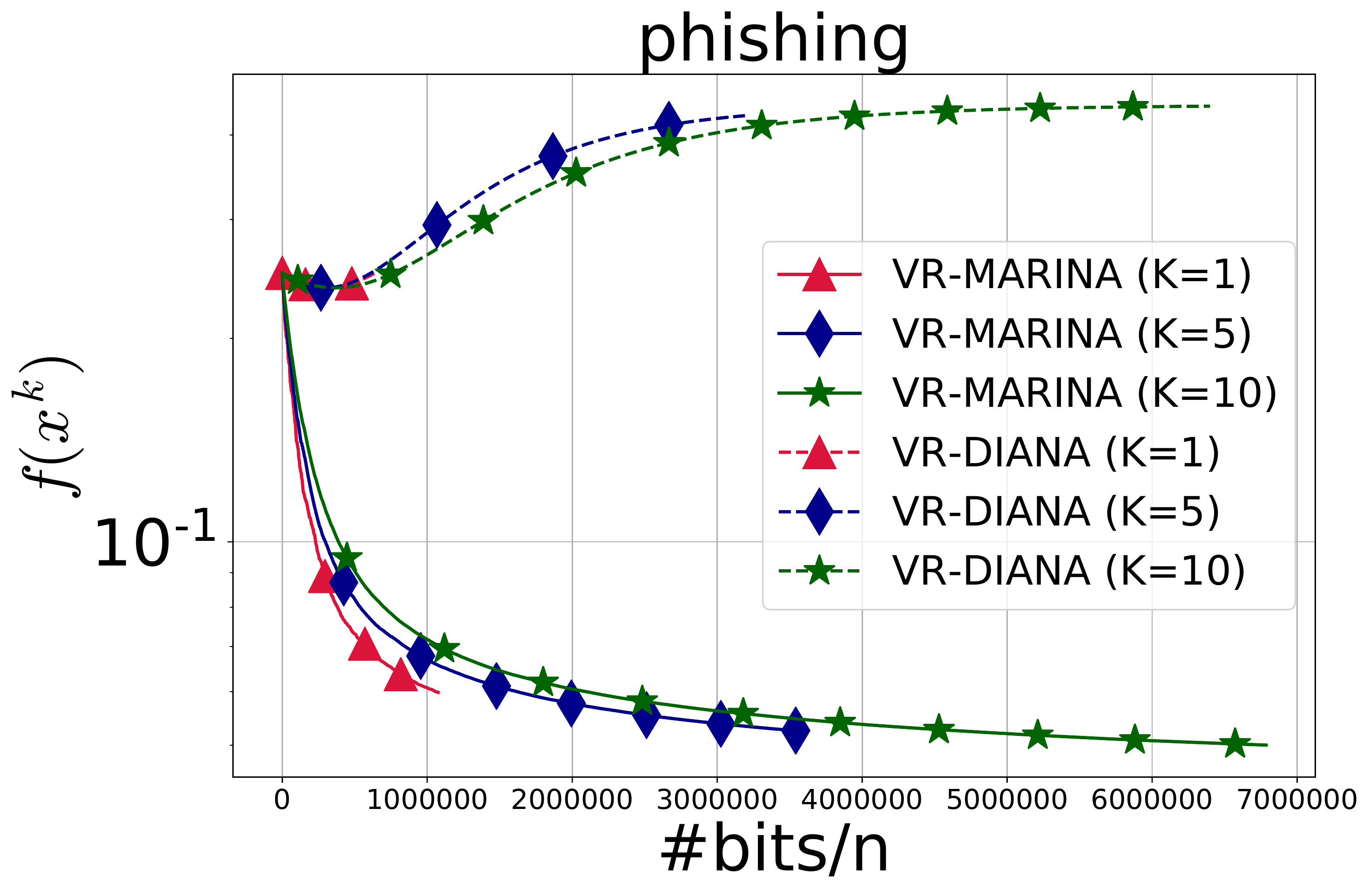}
\includegraphics[width=0.24\textwidth]{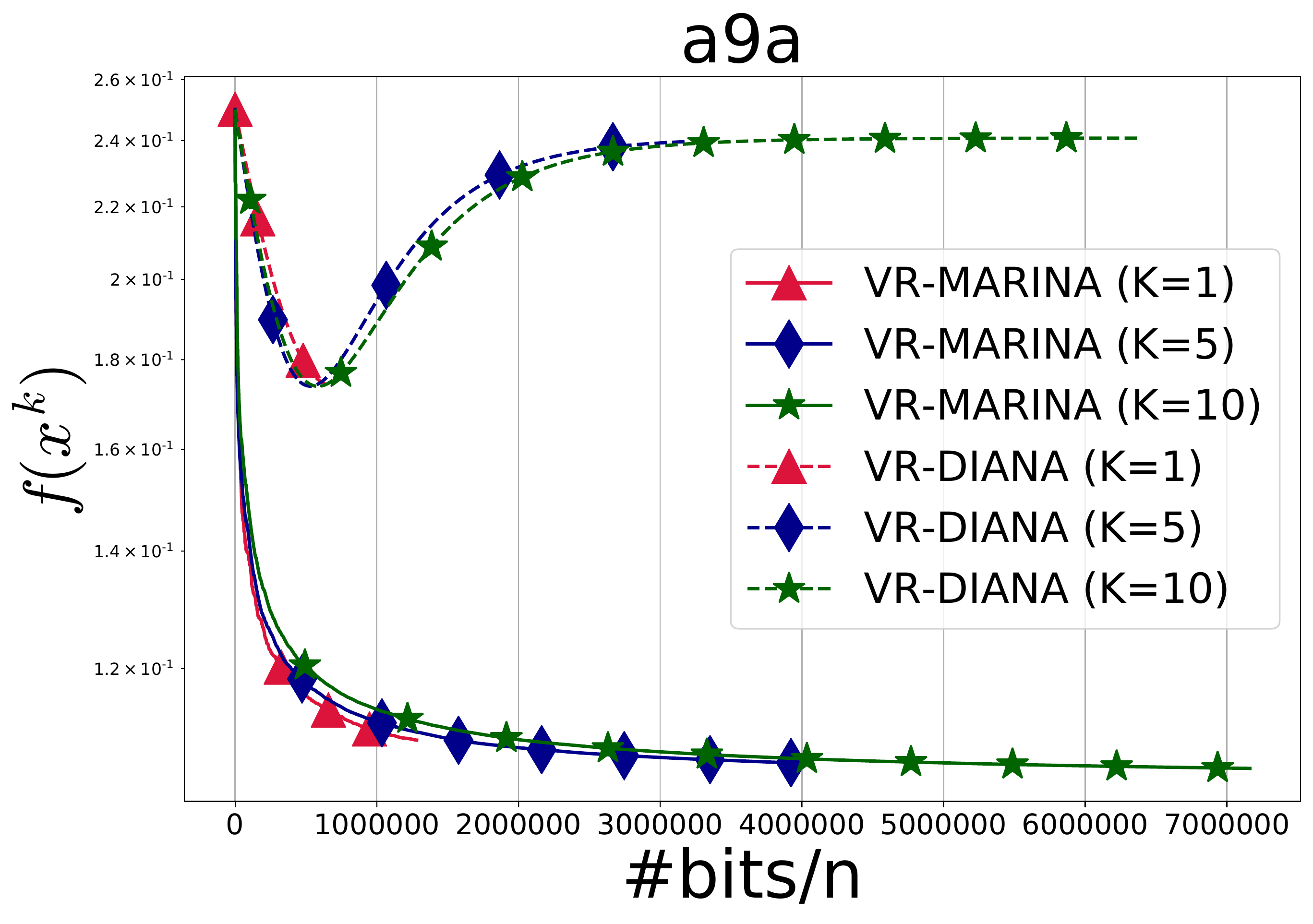}
\includegraphics[width=0.24\textwidth]{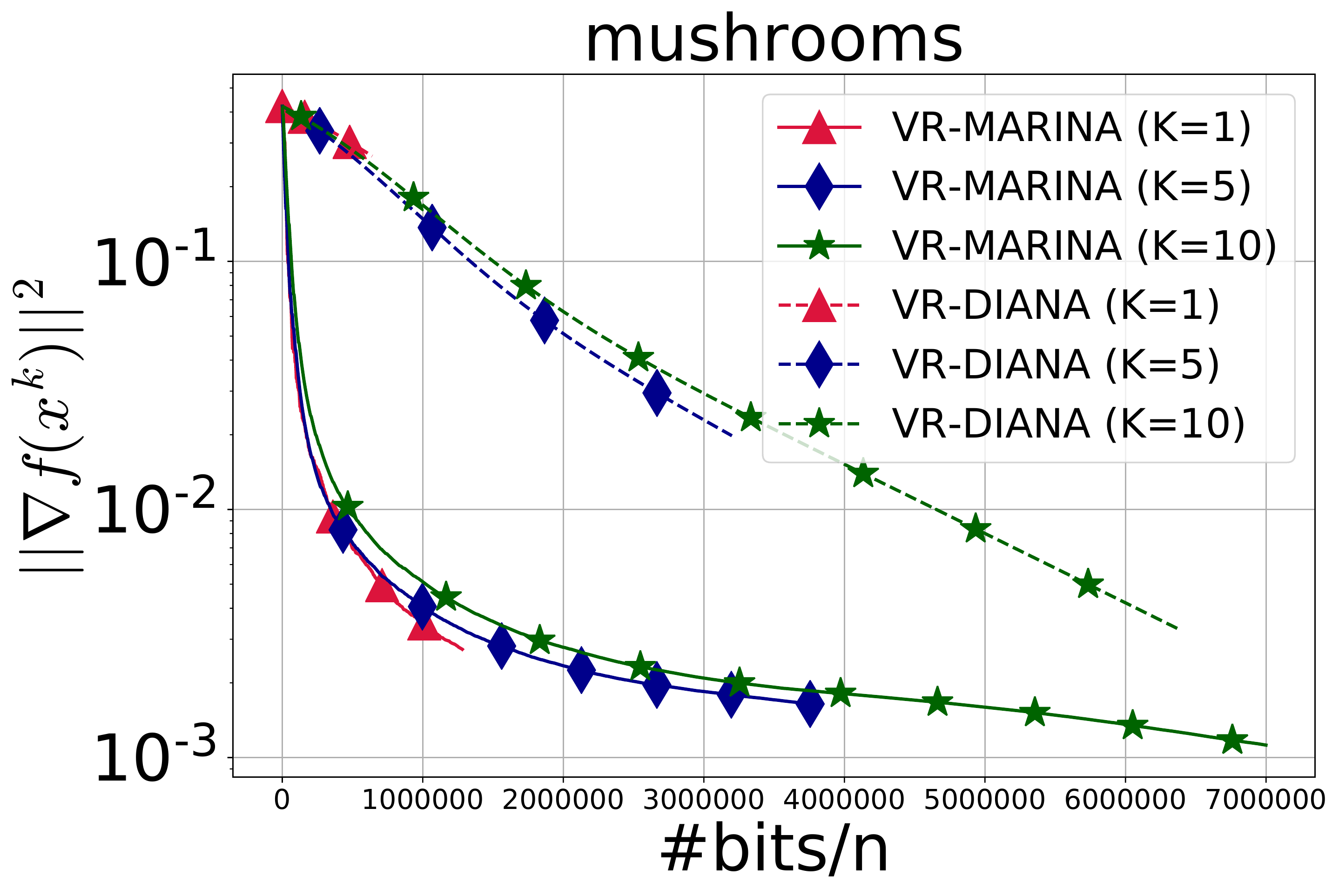}
\includegraphics[width=0.24\textwidth]{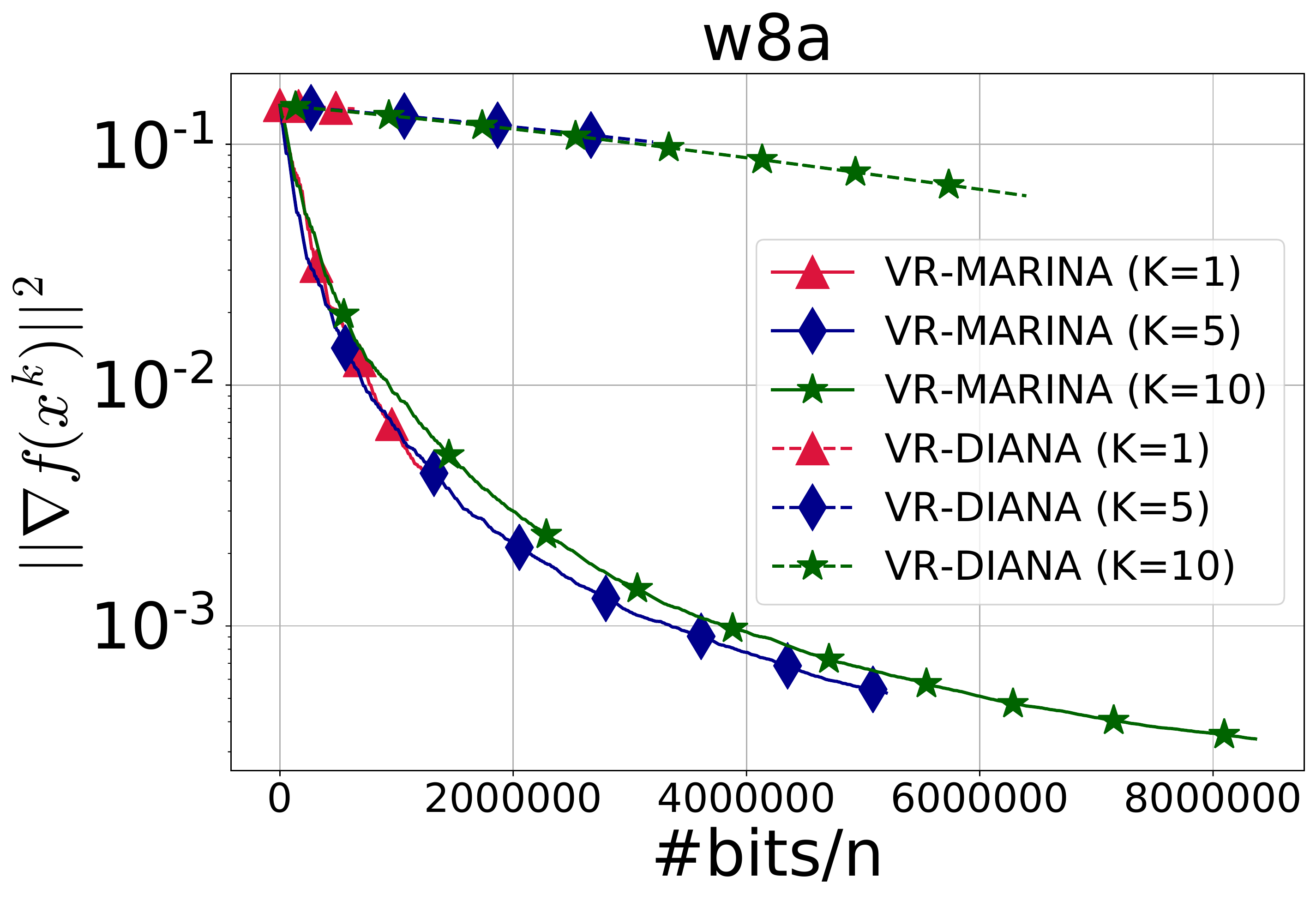}
\includegraphics[width=0.24\textwidth]{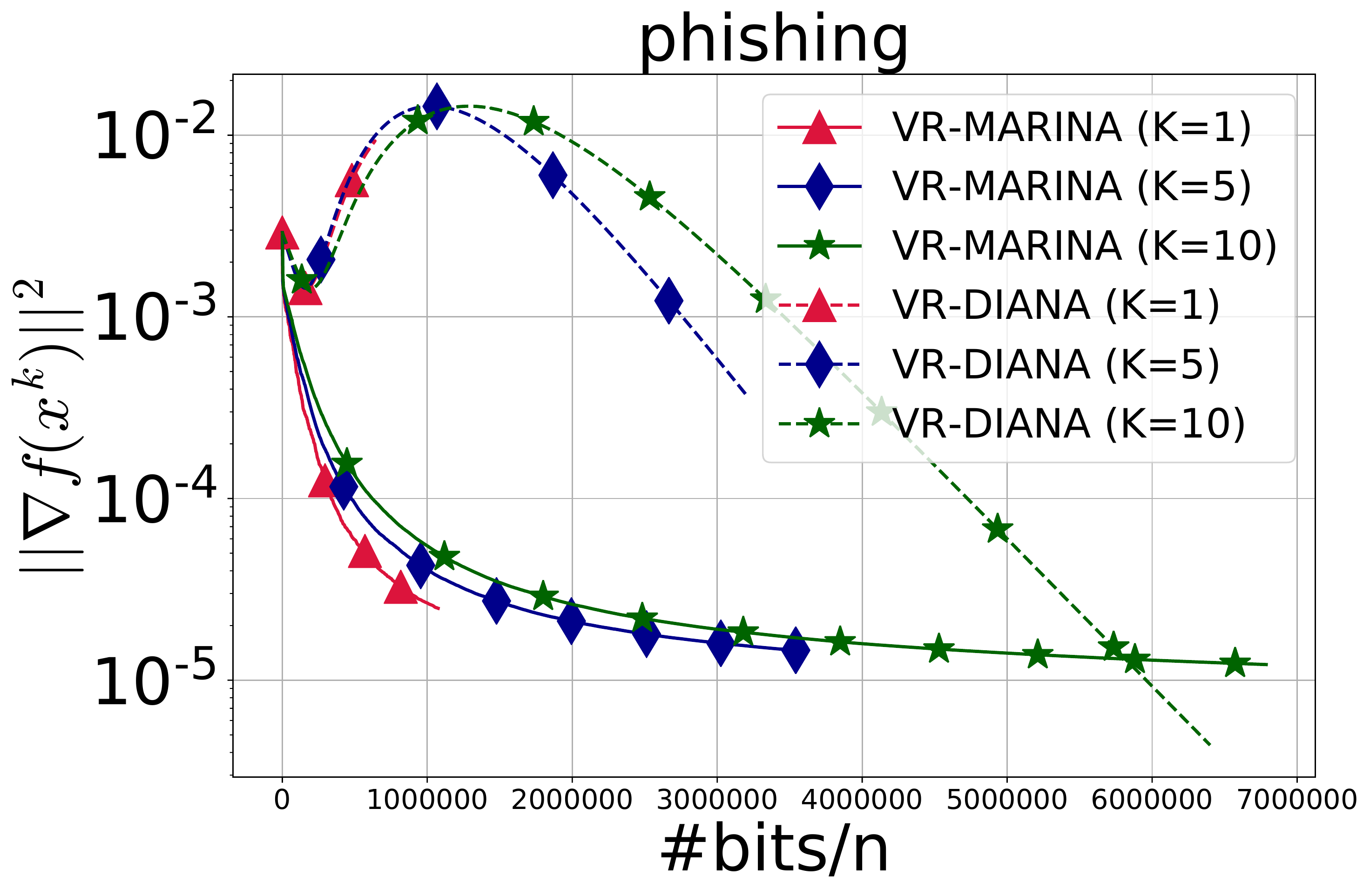}
\includegraphics[width=0.24\textwidth]{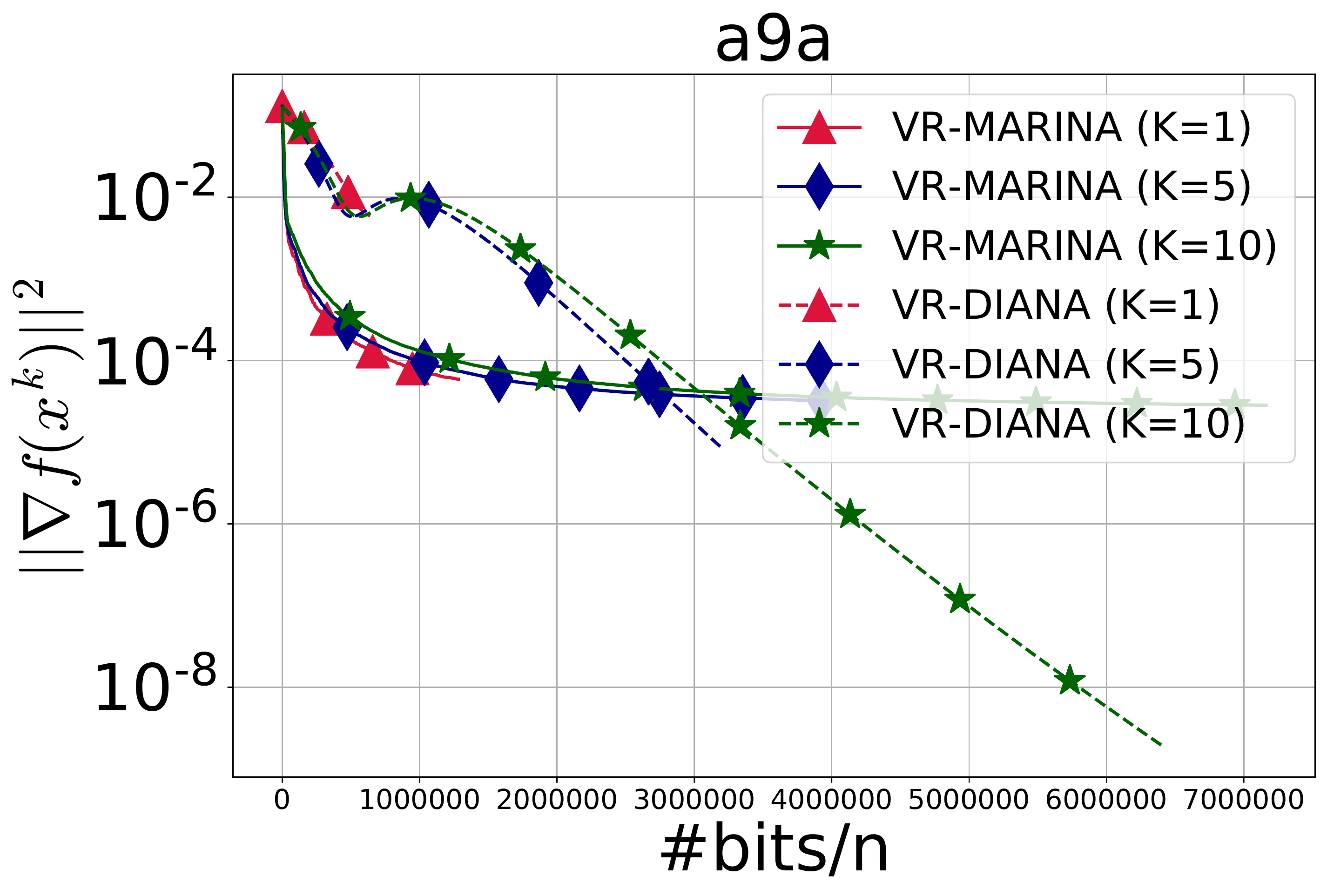}
\includegraphics[width=0.24\textwidth]{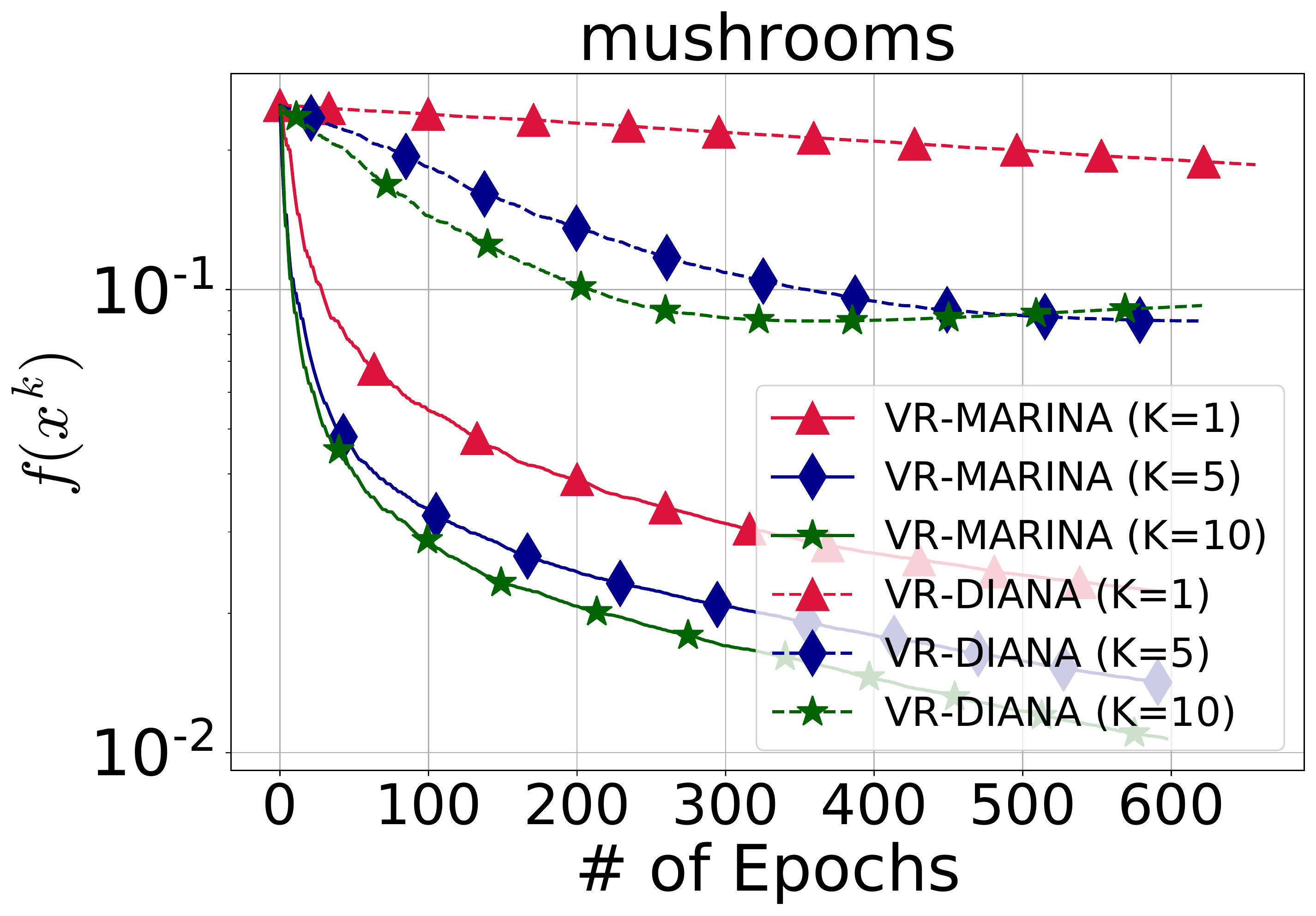}
\includegraphics[width=0.24\textwidth]{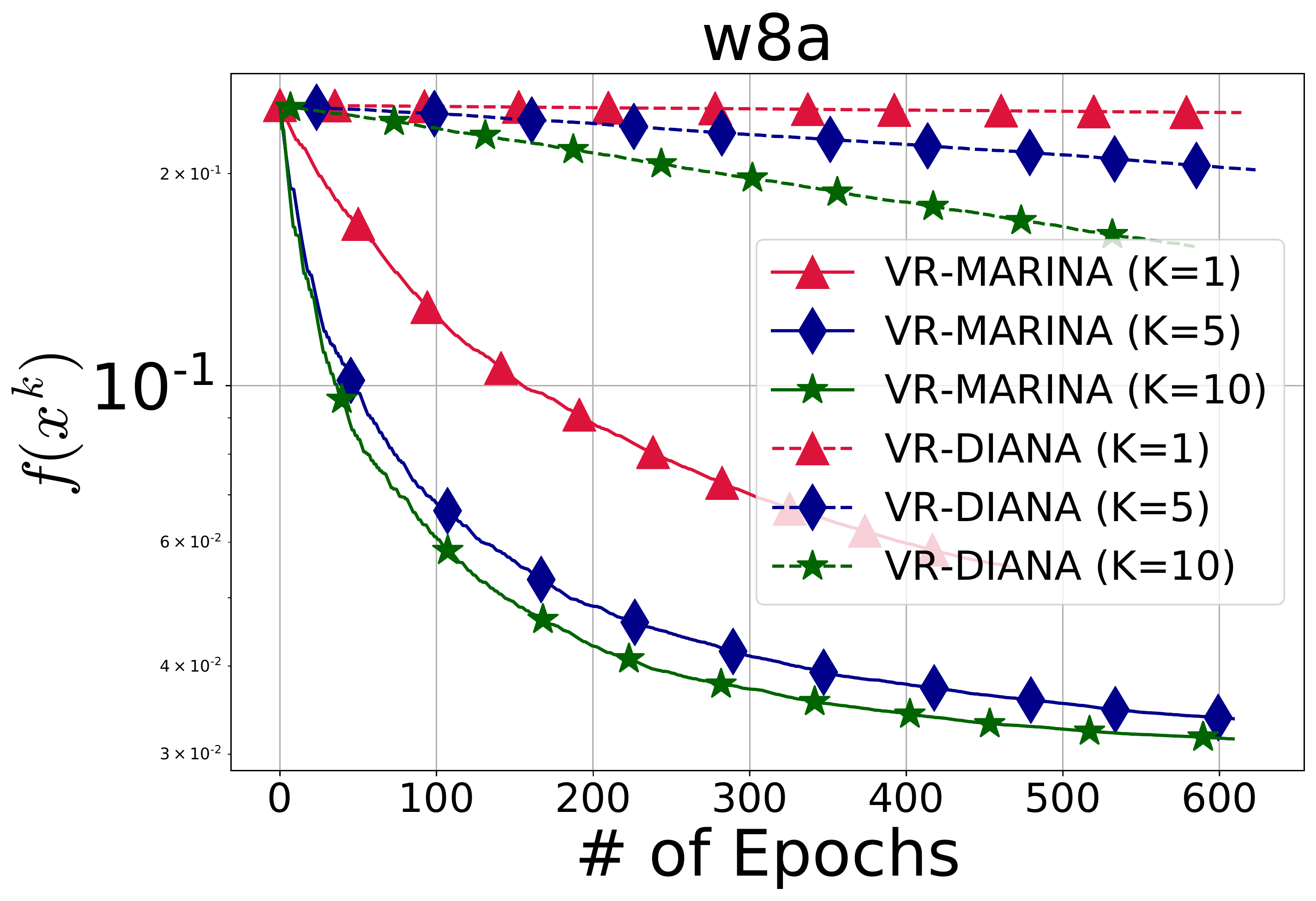}
\includegraphics[width=0.24\textwidth]{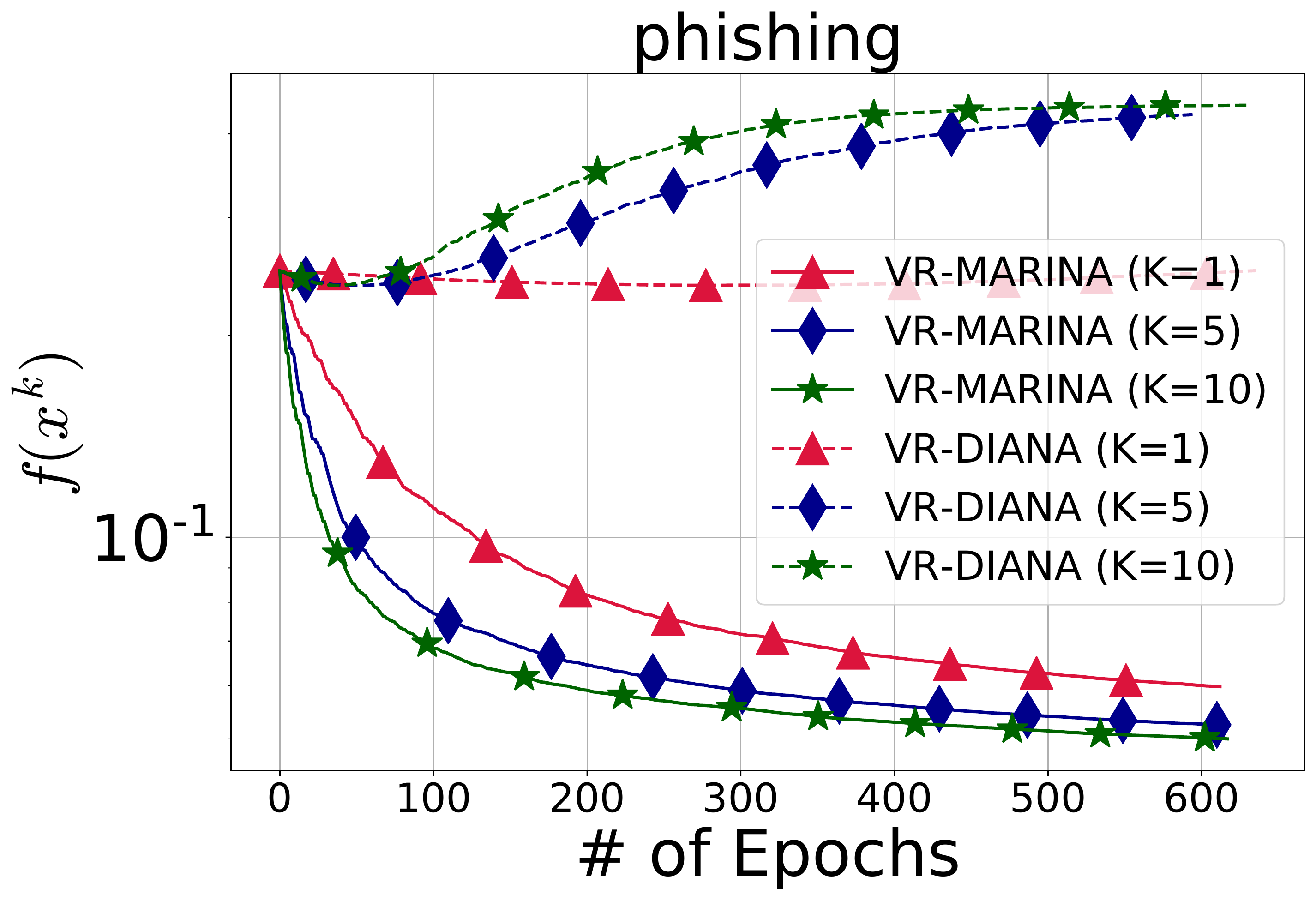}
\includegraphics[width=0.24\textwidth]{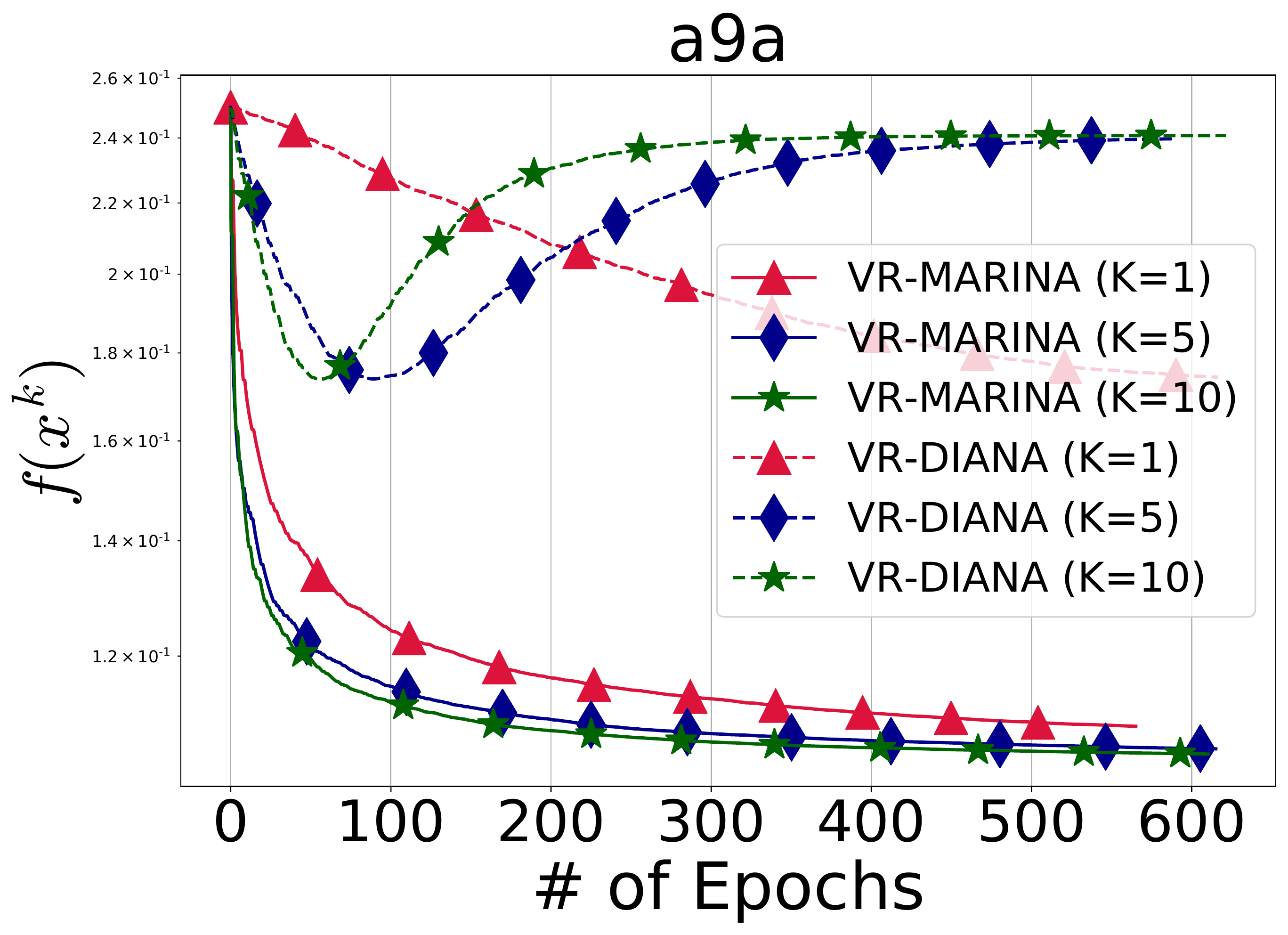}
\includegraphics[width=0.24\textwidth]{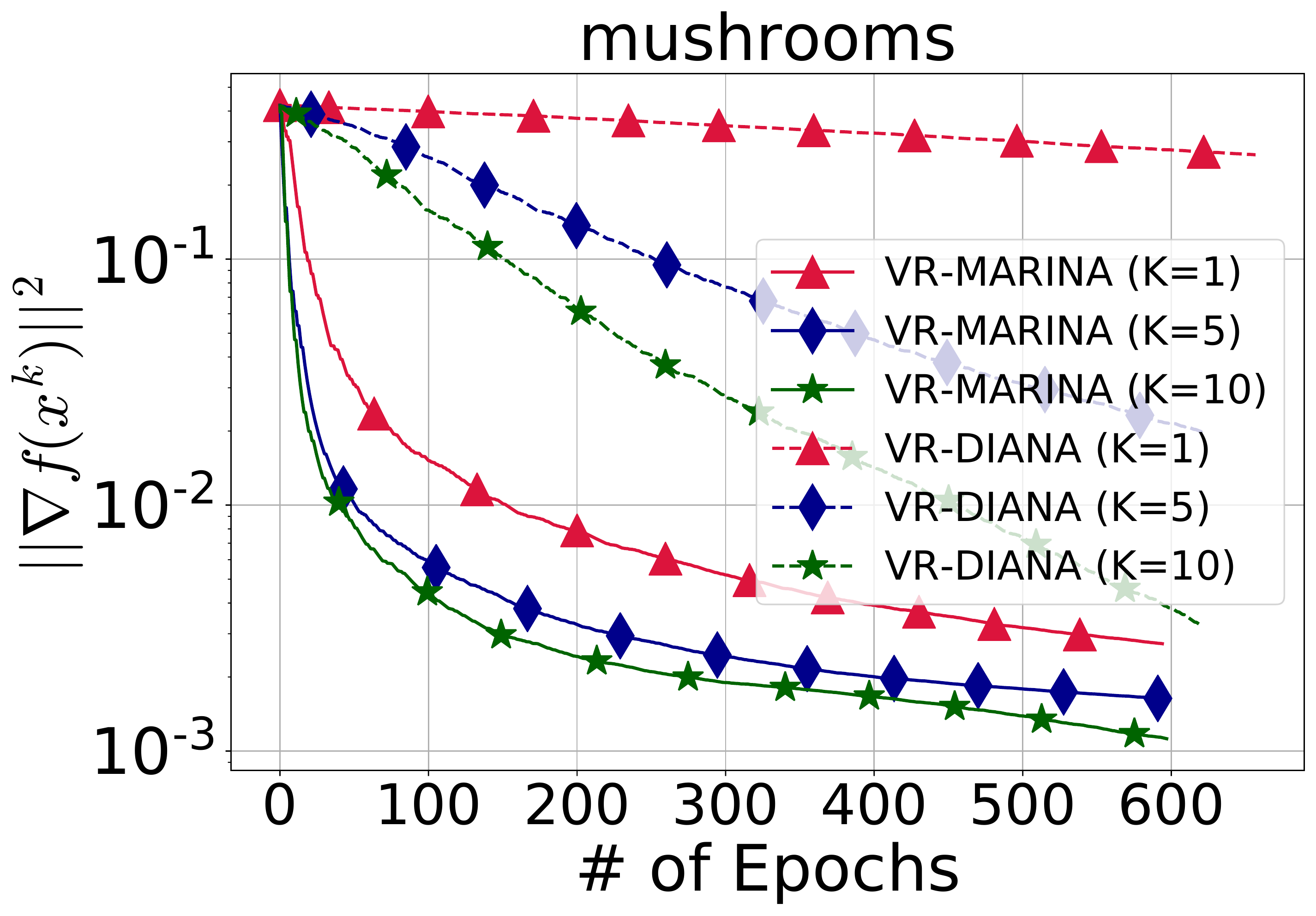}
\includegraphics[width=0.24\textwidth]{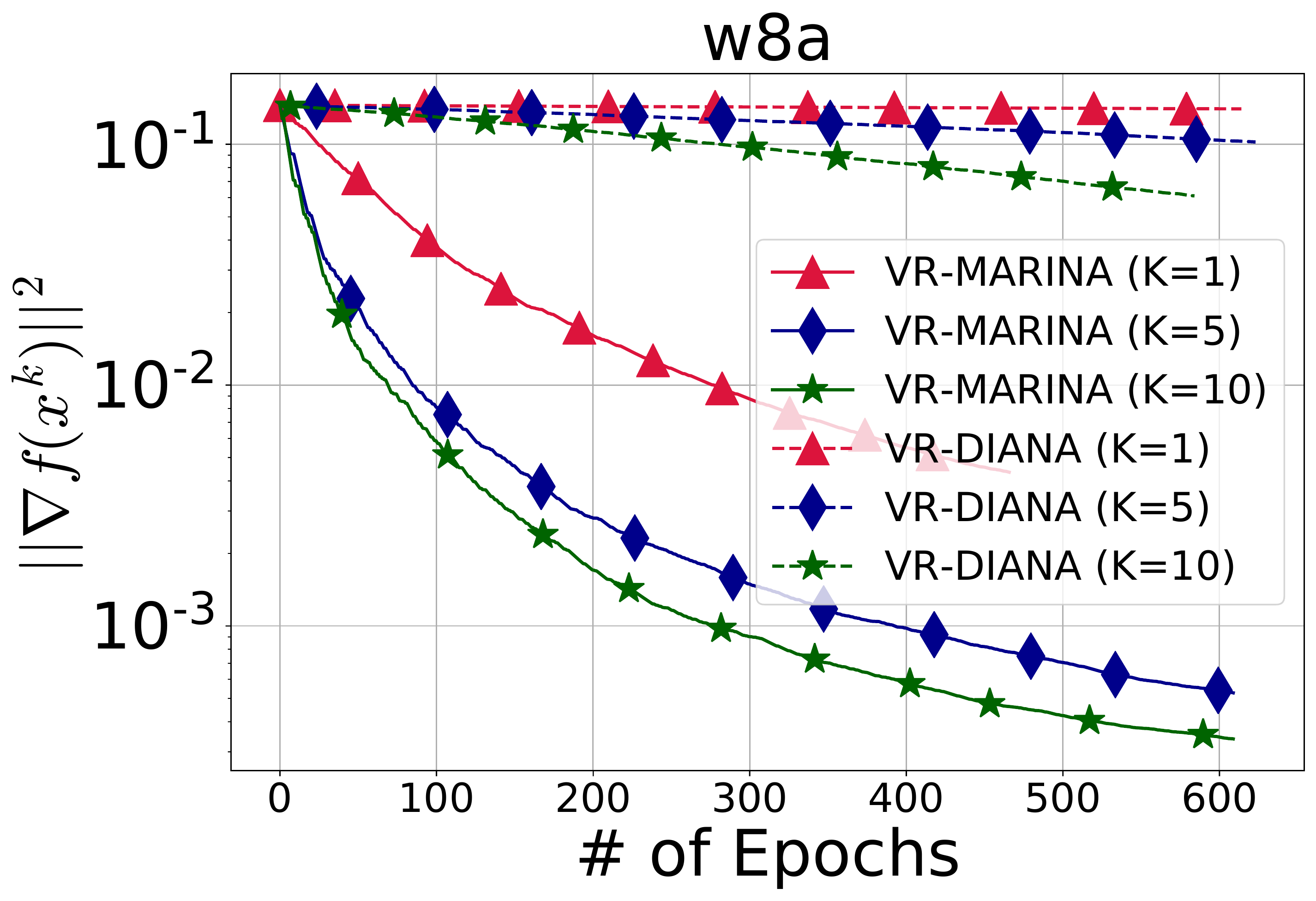}
\includegraphics[width=0.24\textwidth]{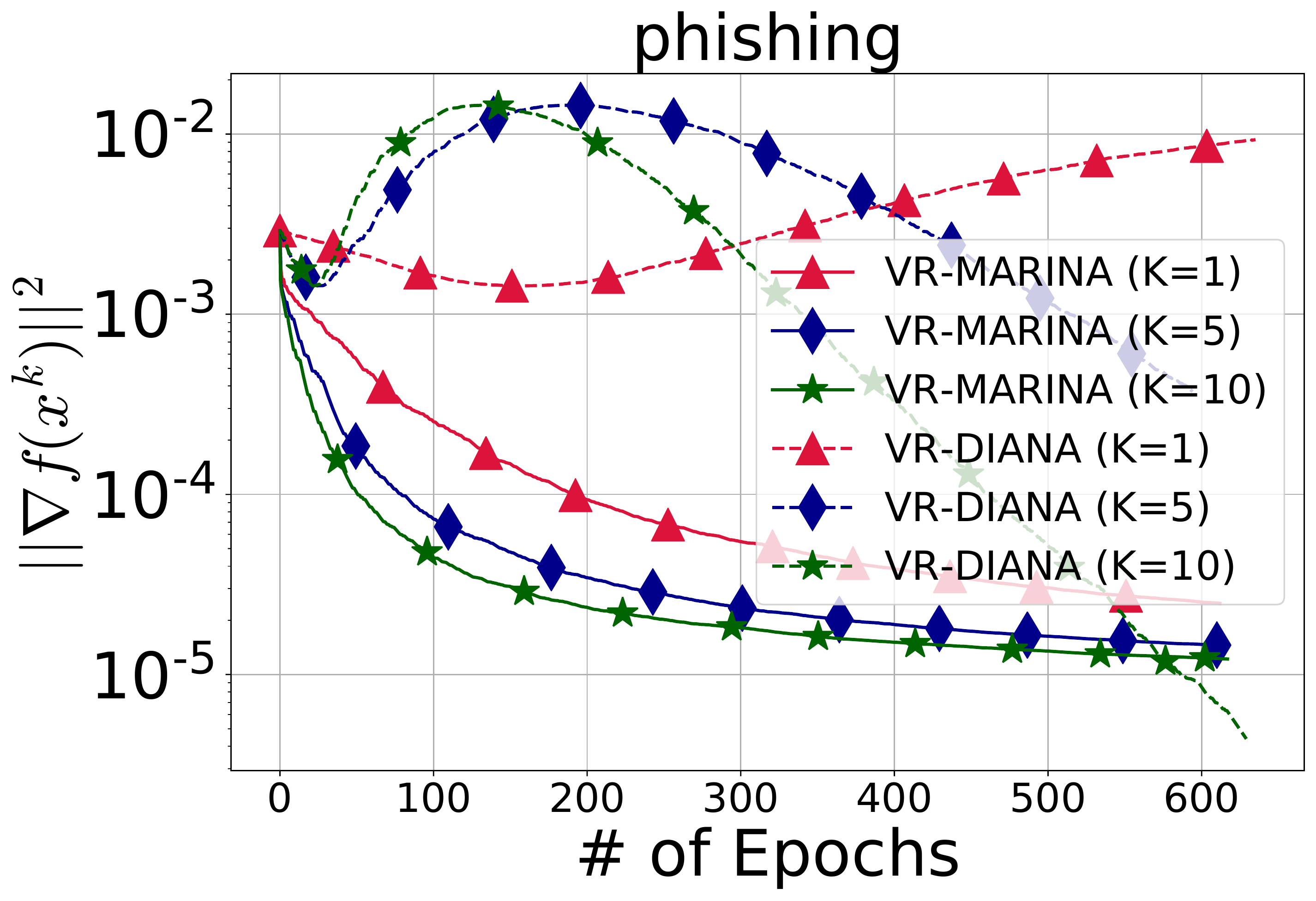}
\includegraphics[width=0.24\textwidth]{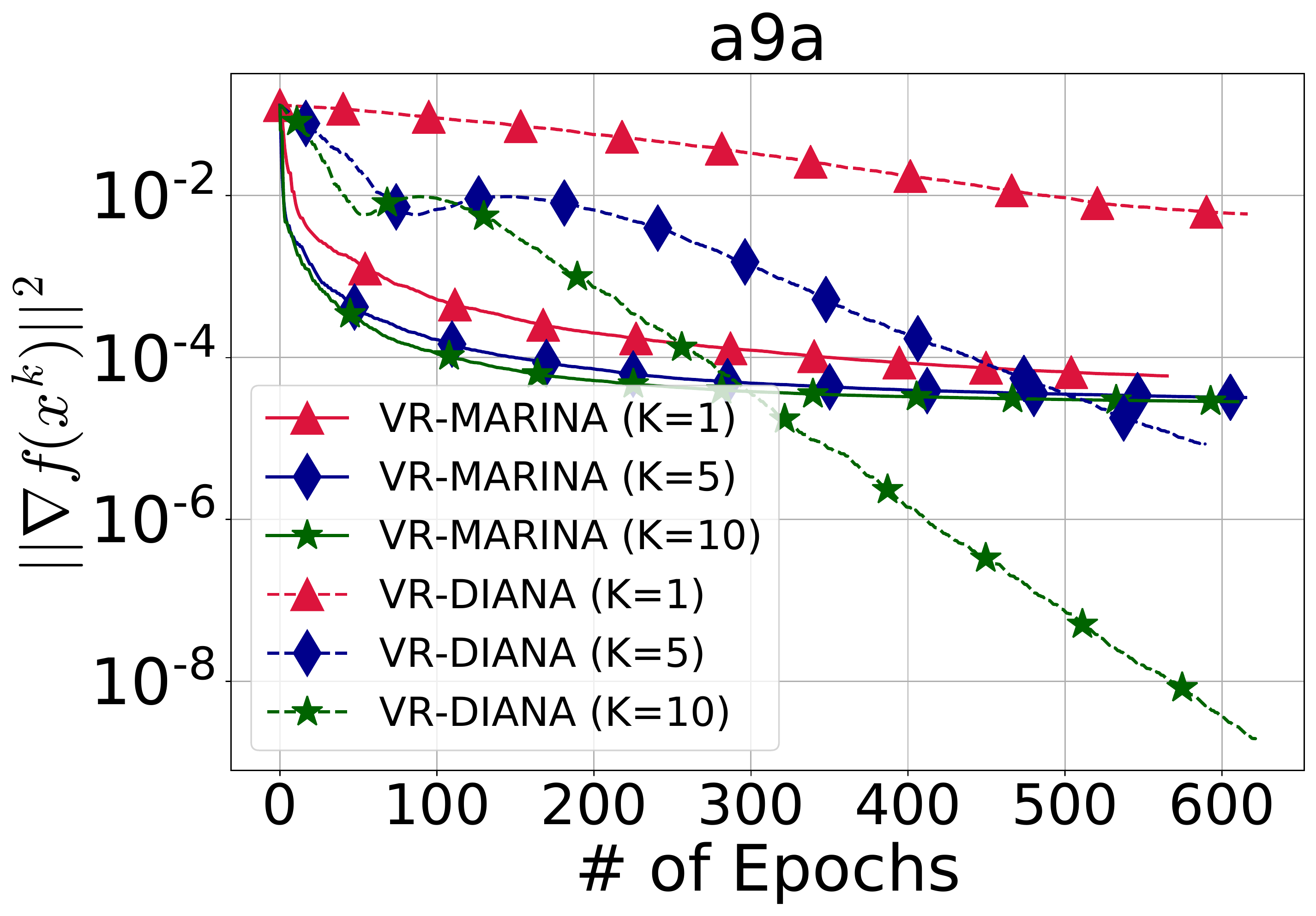}
\caption{Comparison of \algname{VR-MARINA} with  \algname{VR-DIANA} on binary classification problem involving non-convex loss \eqref{eq:experiment_problem} with LibSVM data \cite{chang2011libsvm}. Parameter $n$ is chosen as per Table~\ref{tbl:ns} ($n = 5$). Stepsizes for the methods are chosen according to the theory and the batchsizes are $\sim \nicefrac{m}{100}$. In all cases, we used the RandK sparsification operator with K $\in \{1,5,10\}$.}
\label{fig:vr_methods}
\end{figure}

The results for the full-batched methods are reported in Figure~\ref{fig:full_batched_methods}, and the comparison of \algname{VR-MARINA} and \algname{VR-DIANA} is given in Figure~\ref{fig:vr_methods}. Clearly, in both cases, \algname{MARINA} and \algname{VR-MARINA} show faster convergence than the previous state-of-the-art methods, \algname{DIANA} and \algname{VR-DIANA}, for distributed non-convex optimization with compression in terms of $\|\nabla f(x^k)\|^2$ and $f(x^k)$ decrease w.r.t.\ the number of communication rounds, oracle calls per node and the total number of transferred bits from workers to the master.

\begin{figure}[H]
\centering
\includegraphics[width=0.3\textwidth]{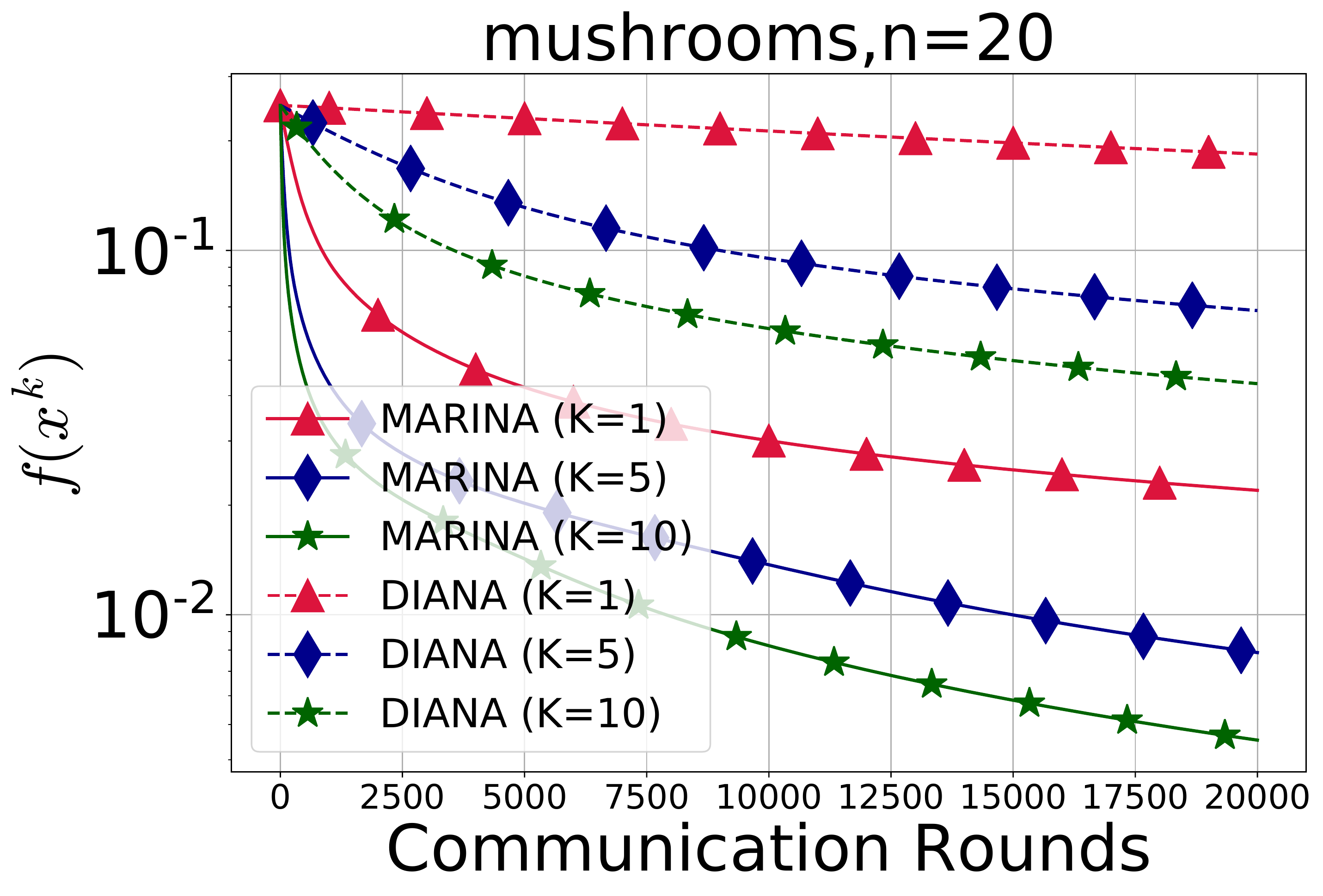}
\includegraphics[width=0.3\textwidth]{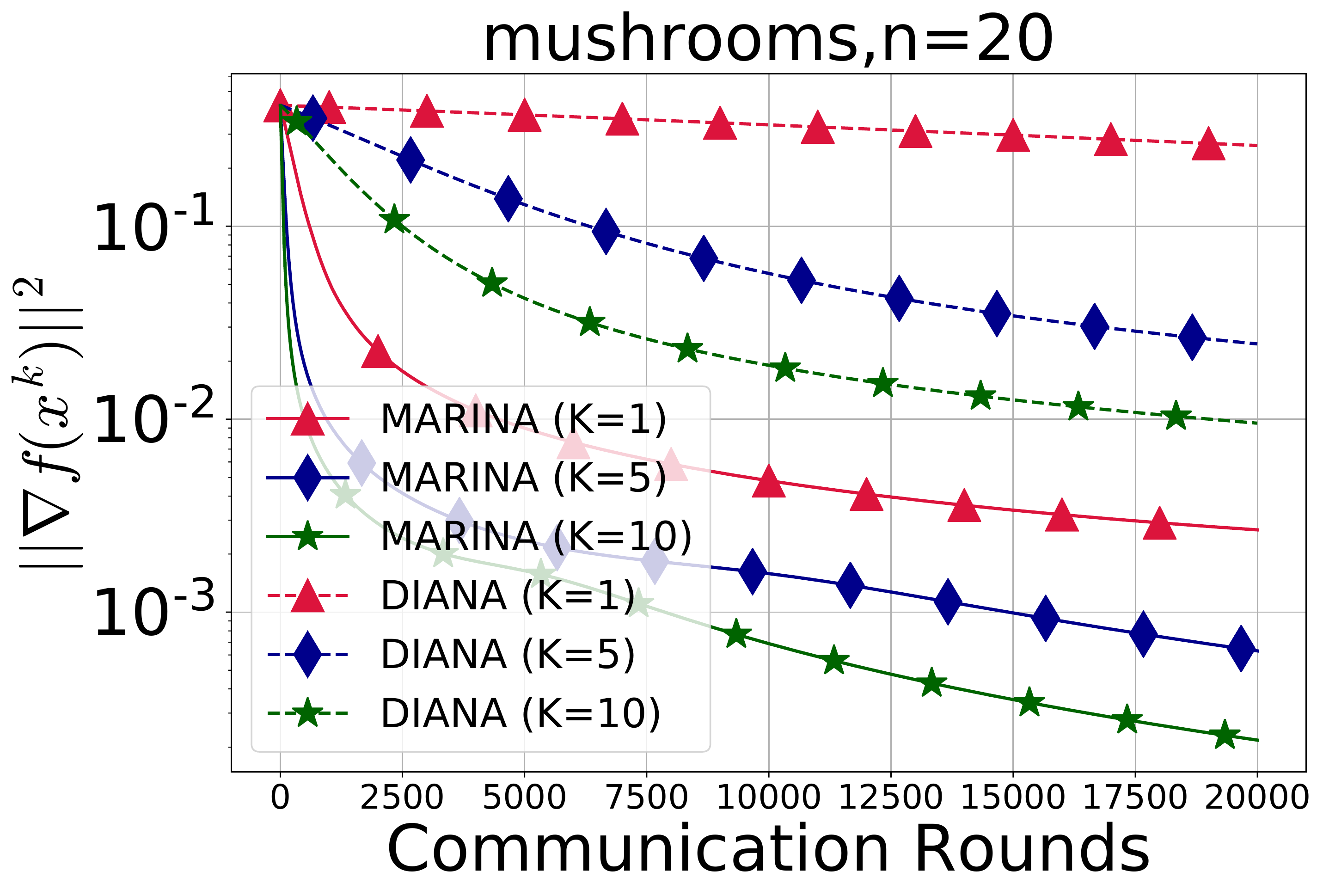}\\
\includegraphics[width=0.3\textwidth]{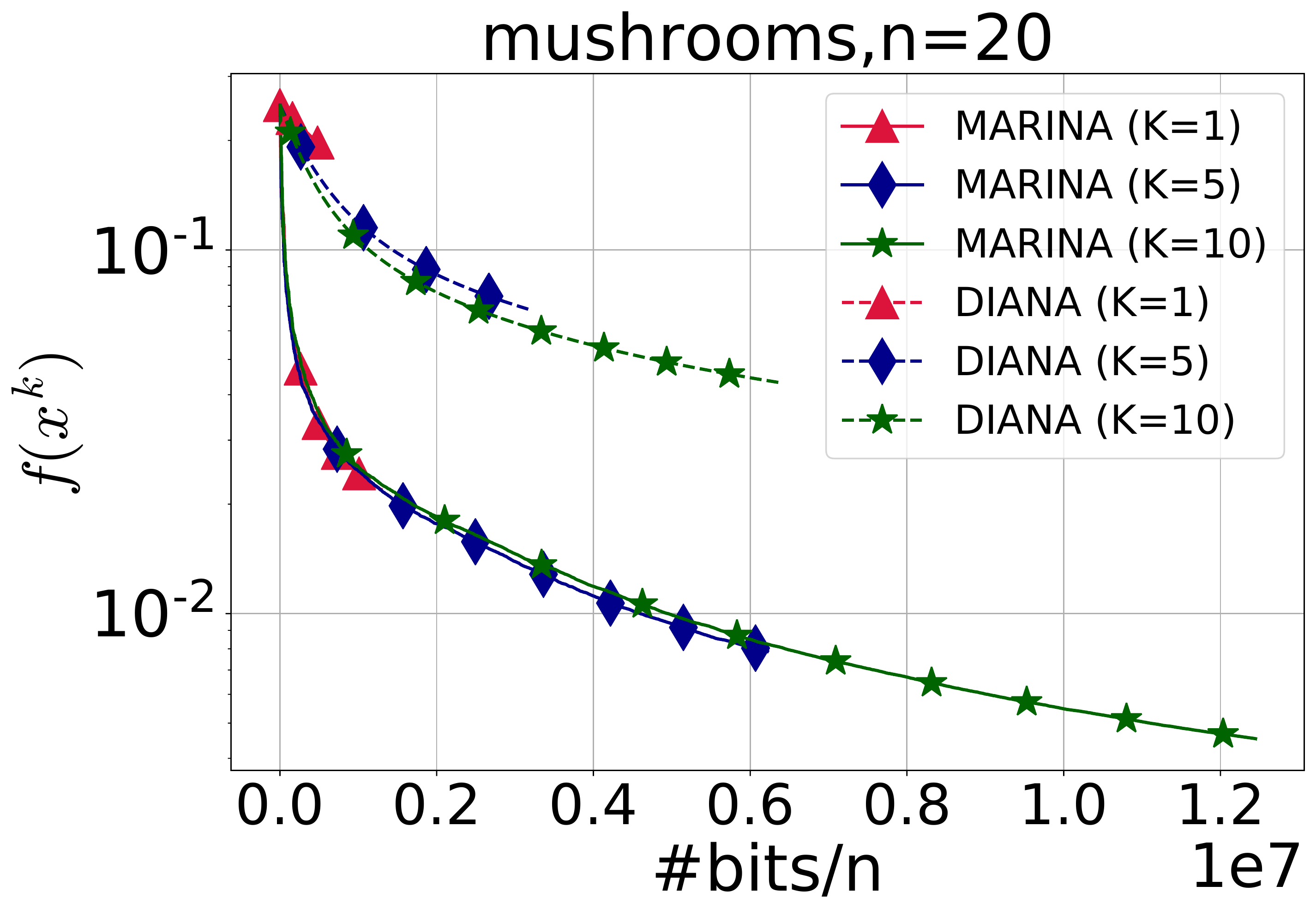}
\includegraphics[width=0.3\textwidth]{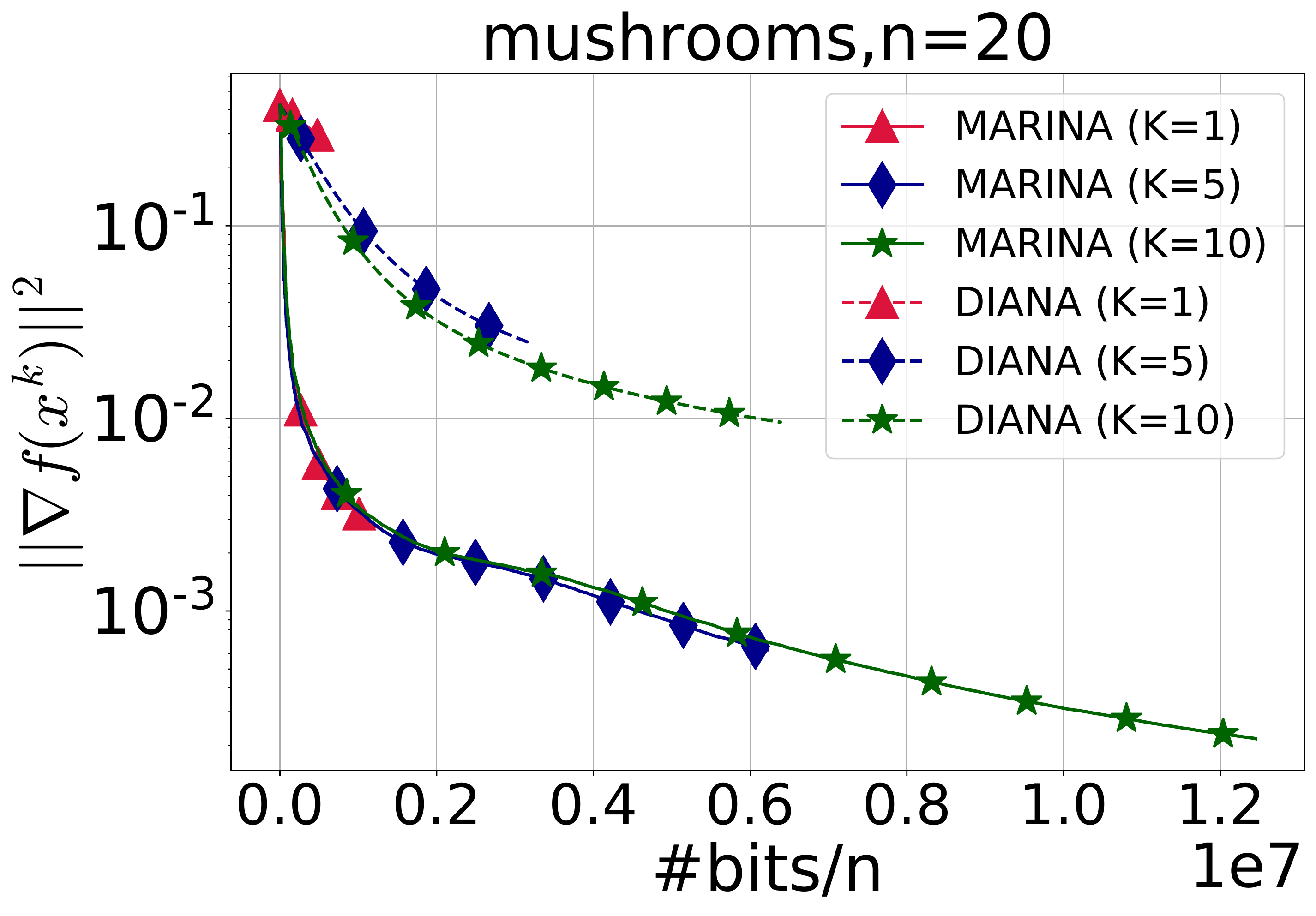}
\caption{Comparison of \algname{MARINA} with  \algname{DIANA} on binary classification problem involving non-convex loss \eqref{eq:experiment_problem} with \texttt{mushrooms} dataset and $n=20$ workers. Stepsizes for the methods are chosen according to the theory. In all cases, we used the RandK sparsification operator with K $\in \{1,5,10\}$.}
\label{fig:full_batched_methods_more_workers}
\end{figure}

We also tested \algname{MARINA} and \algname{DIANA} on \texttt{mushrooms} dataset with a bigger number of workers ($n=20$). The results are reported in Figure~\ref{fig:full_batched_methods_more_workers}. Similarly to the previous numerical tests, \algname{MARINA} shows its superiority to \algname{DIANA} with $n=20$ as well.

\subsection{Image Classification}\label{sec:NN_experiments}
We also compared the performance of \algname{VR-MARINA} and \algname{VR-DIANA} on the training {\tt ResNet-18} \cite{he2016deep} at {\tt CIFAR100} \cite{krizhevsky2009learning} dataset. Formally, the optimization problem is
\begin{equation}
	\min\limits_{x\in \R^d}\left\{f(x) = \frac{1}{N}\sum\limits_{i=1}^N\ell(p(f(a_i, x)), y_i)\right\}, \label{eq:NN_problem}
\end{equation} 
where $\{(a_i,y_i)\}_{i=1}^N$ encode images and labels from {\tt CIFAR100} dataset, $f(a_i,x)$ is the output of {\tt ResNet-18} on image $a_i$ with weights $x$, $p$ is softmax function, and $\ell(\cdot,\cdot)$ is cross-entropy loss. {\tt ResNet-18} has $d =$ 11~689~512 parameters to train and {\tt CIFAR100} contains $N =$ 50~000 colored images. The dataset is split into $5$ parts among $5$ workers in such a way that the first four workers get 10~112 samples and the fifth one get 9~552 samples. The code was written in Python 3.9 using \textsc{PyTorch 1.7} and then was executed on a machine with NVIDIA GPU Geforce RTX 2080 Ti with 11 GByte onboard global GPU memory.

In all experiments, we use batchsize $= 256$ on each worker and tune the stepsizes for each method separately. That is, for each method  and for each choice of $K$ for RandK operator we run the method with stepsize $\gamma \in \{10^{-6}, 0.1, 0.2, 0.5, 1.0, 5.0\}$ to find the interval containing the best stepsize. Next, the obtained interal is split into $\sim 10$ equal parts and the method is run with corresponding stepsizes. Other parameters of the methods are chosen according to the theory. The summary of used parameters is given in Table~\ref{tbl:params_resnet}.

\begin{table}[H]
 \caption{Summary of the parameters used in the experiments presented in Fig.~\ref{fig:resnet_at_cifar100} and Fig.~\ref{fig:resnet_at_cifar100_no_compr}. Stepsizes were tuned, batchsize $= 256$ on each worker, other parameters were picked according to the theory, except the last line, where $p$ for \algname{VR-MARINA} without compression was picked as for \algname{VR-MARINA} with RandK, $K = $ 100 000 compression operator.}
\label{tbl:params_resnet}
\begin{center}
\begin{tabular}{|c|c|c|c|}
\hline
Method  & RandK, $K = $ & $\gamma$ & $p$  \\
 \hline
  \hline
\algname{VR-MARINA} & 100 000 & $0.95$ & $0.008554$    \\ \hline
\algname{VR-MARINA} & 500 000 & $0.95$ & $0.024691$    \\ \hline
\algname{VR-MARINA} & 1 000 000 & $0.95$ & $0.024691$    \\ \hline
\algname{VR-DIANA} & 100 000 & $0.15$ & $0.025316$    \\ \hline
\algname{VR-DIANA} & 500 000 & $0.35$ & $0.025316$    \\ \hline
\algname{VR-DIANA} & 1 000 000 & $0.35$ & $0.025316$    \\ \hline
\algname{VR-MARINA} & 11 689 512 ($K = d$) & $3.5$ & $0.024691$    \\ \hline
\algname{VR-DIANA} & 11 689 512 ($K = d$) & $2.5$ & $0.025316$    \\ \hline
\algname{VR-MARINA} & 11 689 512 ($K = d$) & $3.5$ & $0.008554$    \\ \hline
\end{tabular}
\end{center}
\end{table}

\begin{figure}[H]
\centering
\includegraphics[width=0.23\textwidth]{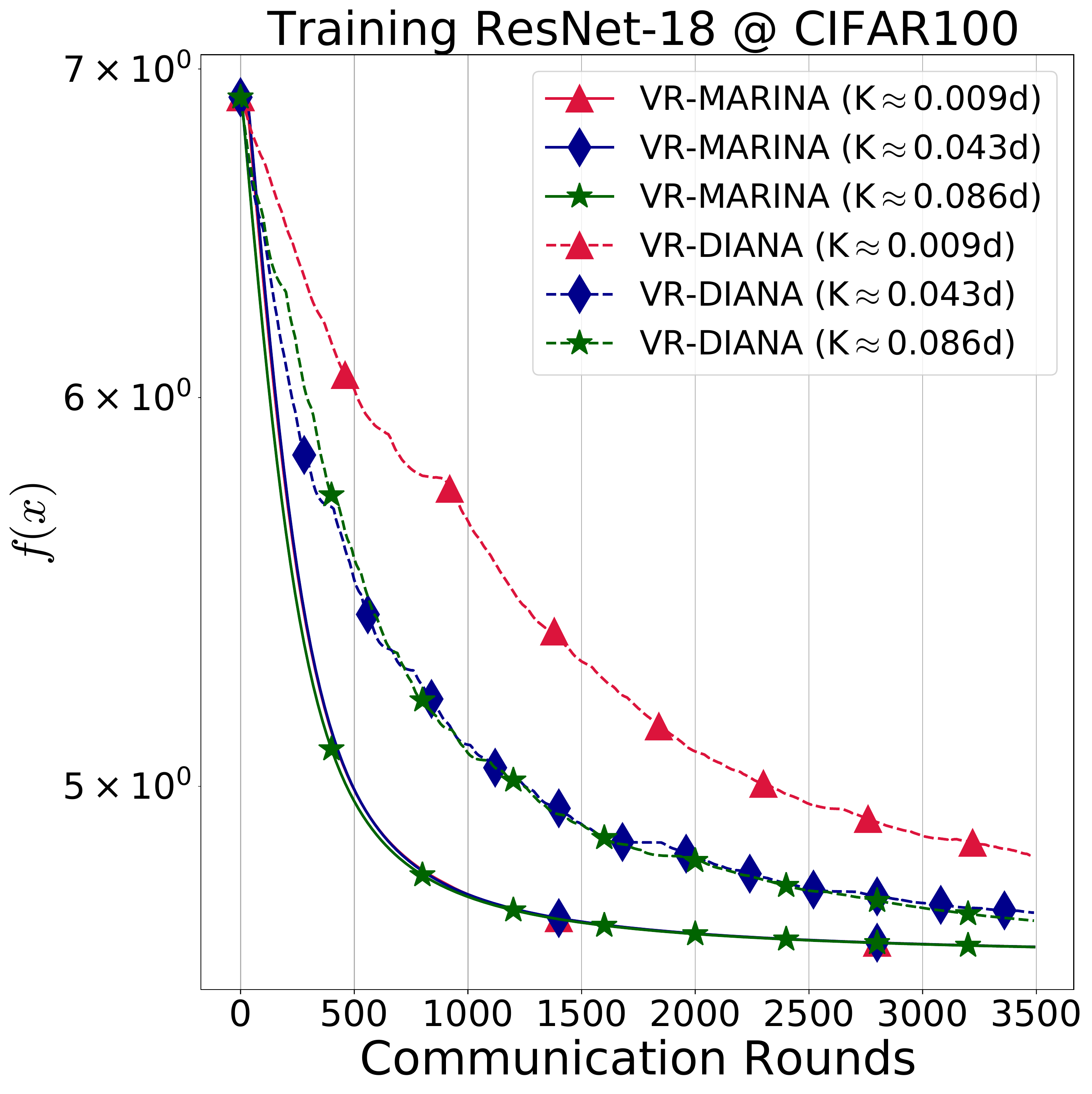}
\includegraphics[width=0.23\textwidth]{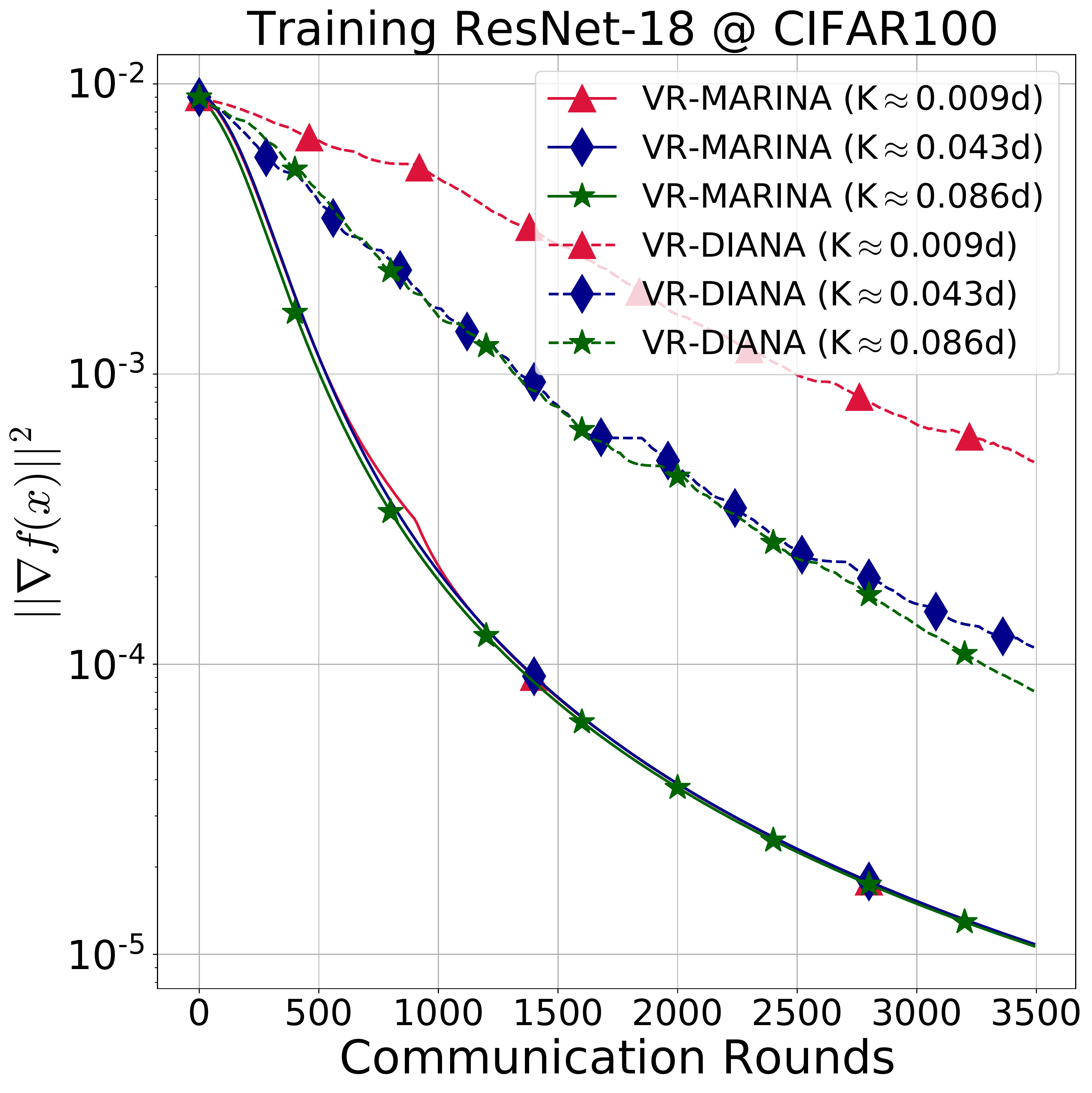}
\includegraphics[width=0.23\textwidth]{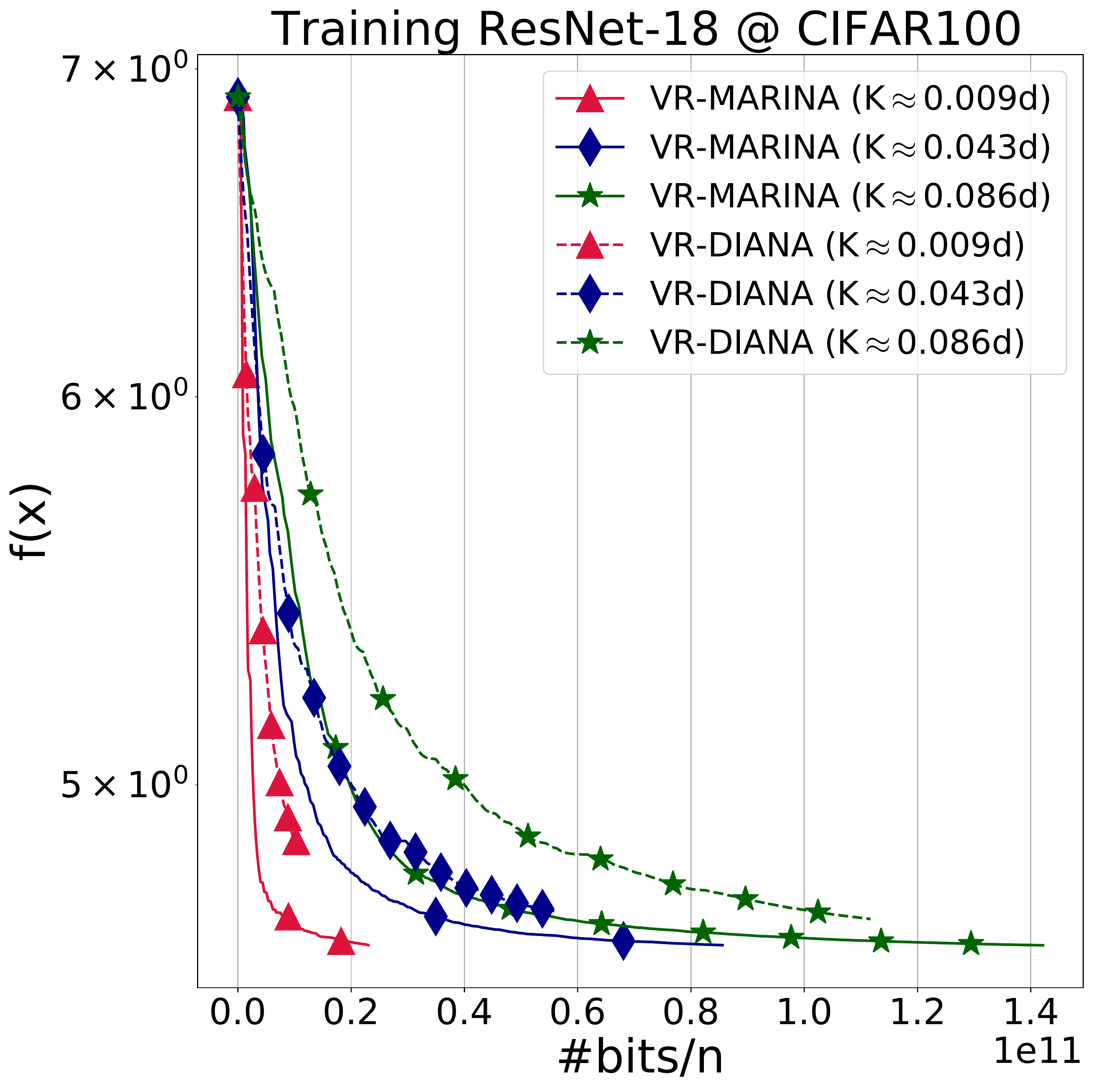}
\includegraphics[width=0.23\textwidth]{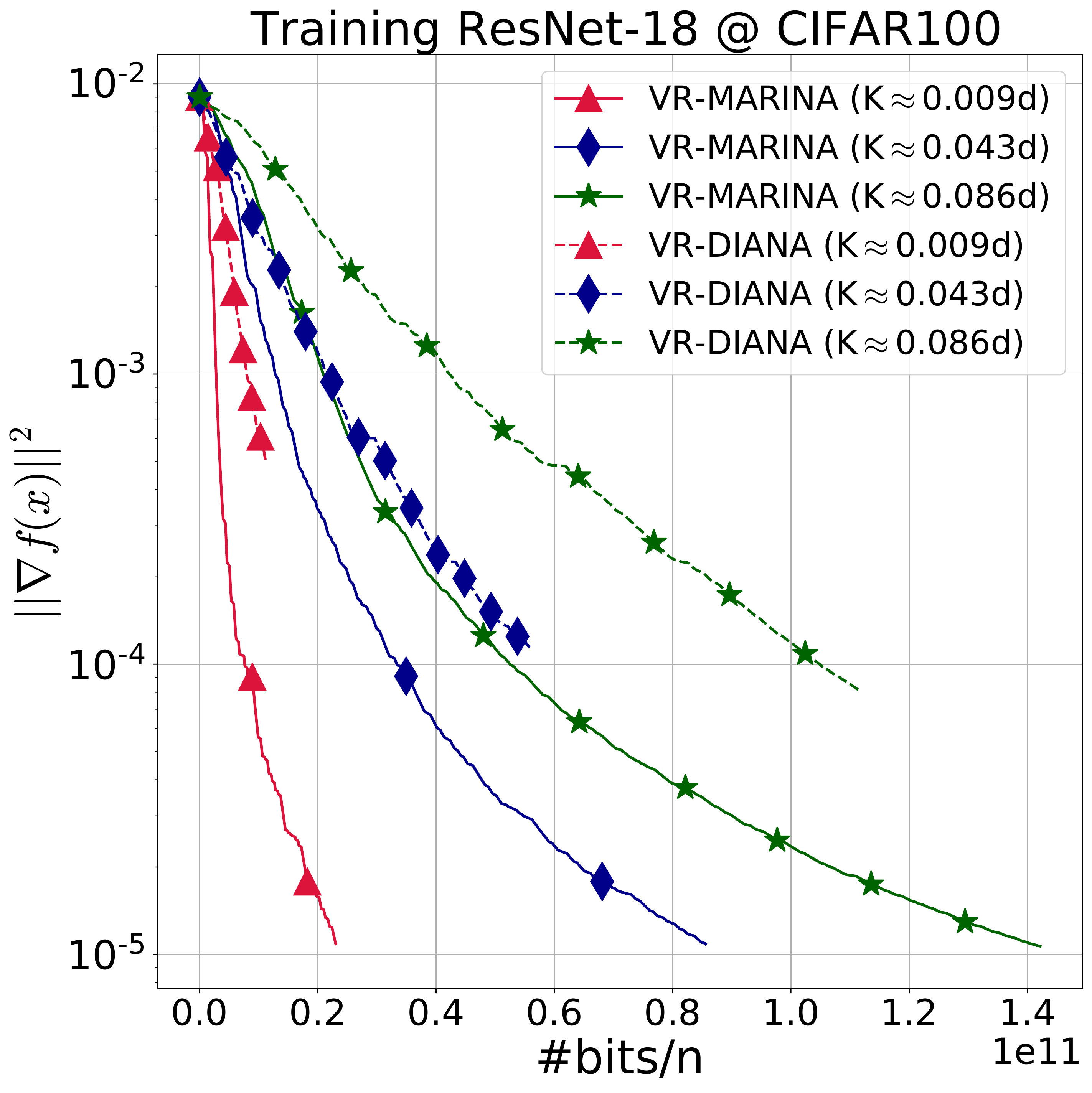}
\includegraphics[width=0.23\textwidth]{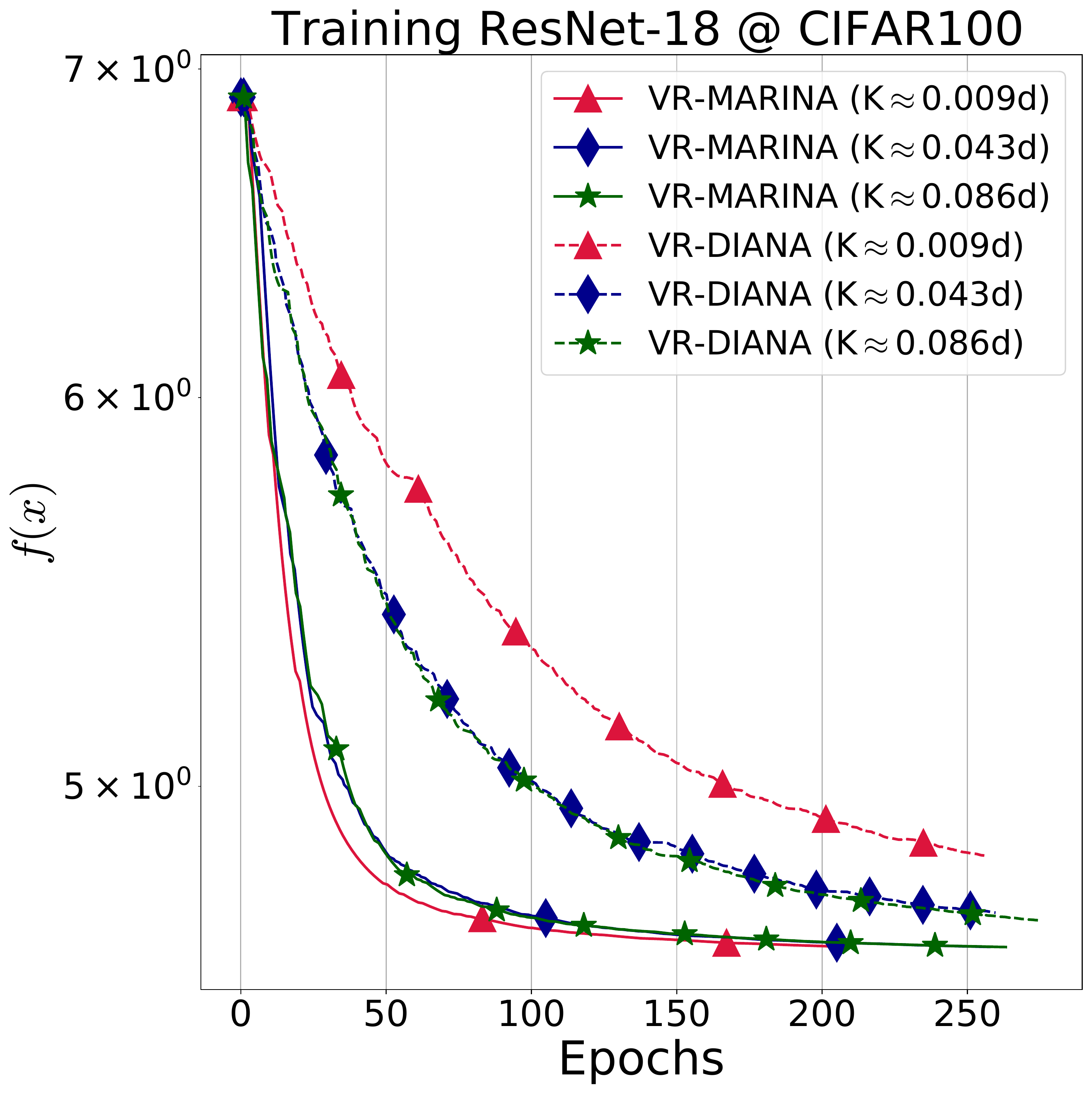}
\includegraphics[width=0.23\textwidth]{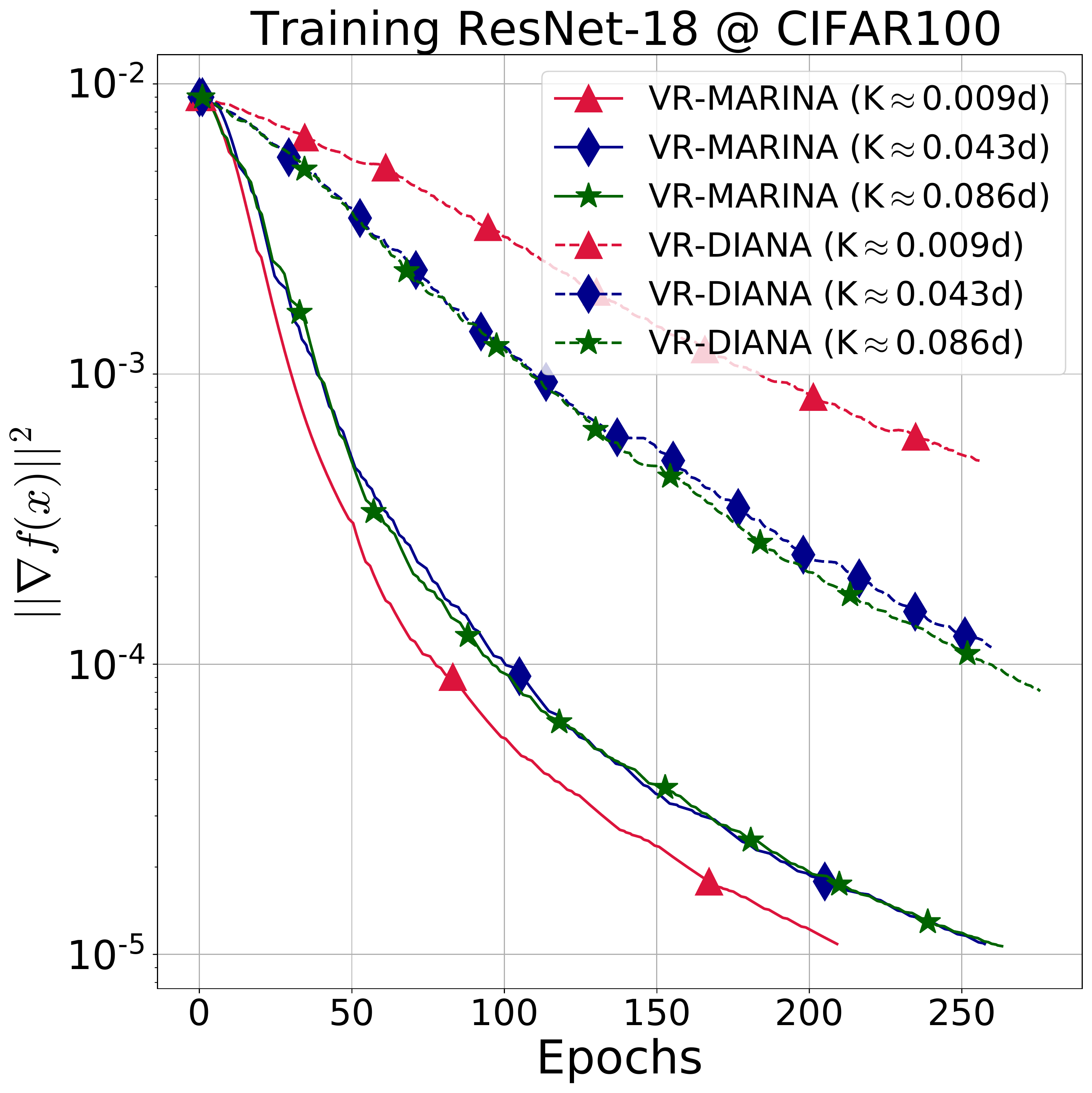}
\caption{Comparison of \algname{VR-MARINA} with \algname{VR-DIANA} on training {\tt ResNet-18} at {\tt CIFAR100} dataset. Number of workers equals $5$. Stepsizes for the methods were tuned and the batchsizes are $\sim \nicefrac{m}{50}$. In all cases, we used the RandK sparsification operator, the approximate values of $K$ are given in the legends ($d$ is dimension of the problem).}
\label{fig:resnet_at_cifar100}
\end{figure}

The results are presented in Fig.~\ref{fig:resnet_at_cifar100}. Again, \algname{VR-MARINA} converges significantly faster than \algname{VR-DIANA} both in terms of the oracle complexity and the total number of transmitted bits to achieve the given accuracy.

\begin{figure}[H]
\centering
\includegraphics[width=0.4\textwidth]{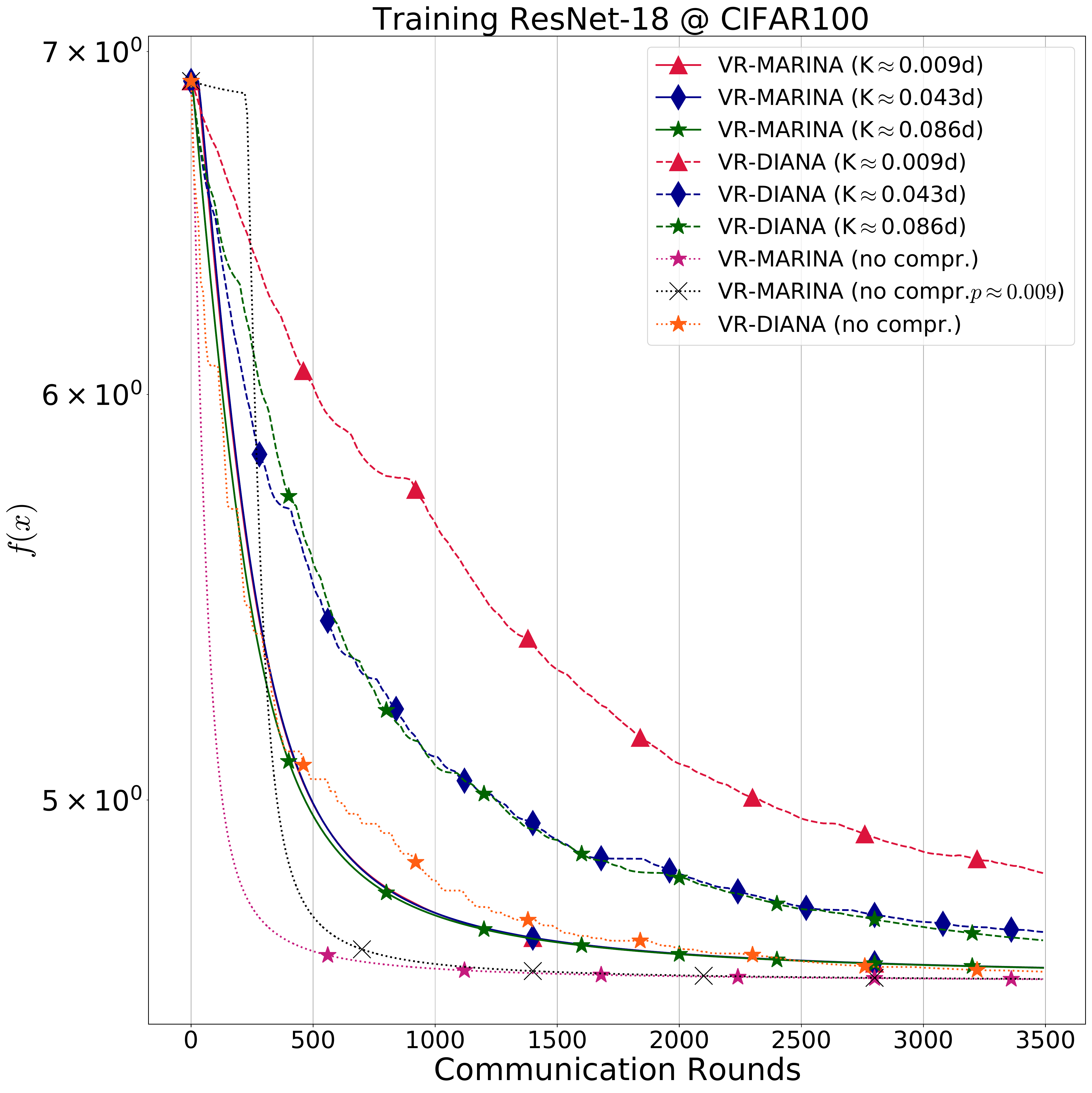}
\includegraphics[width=0.4\textwidth]{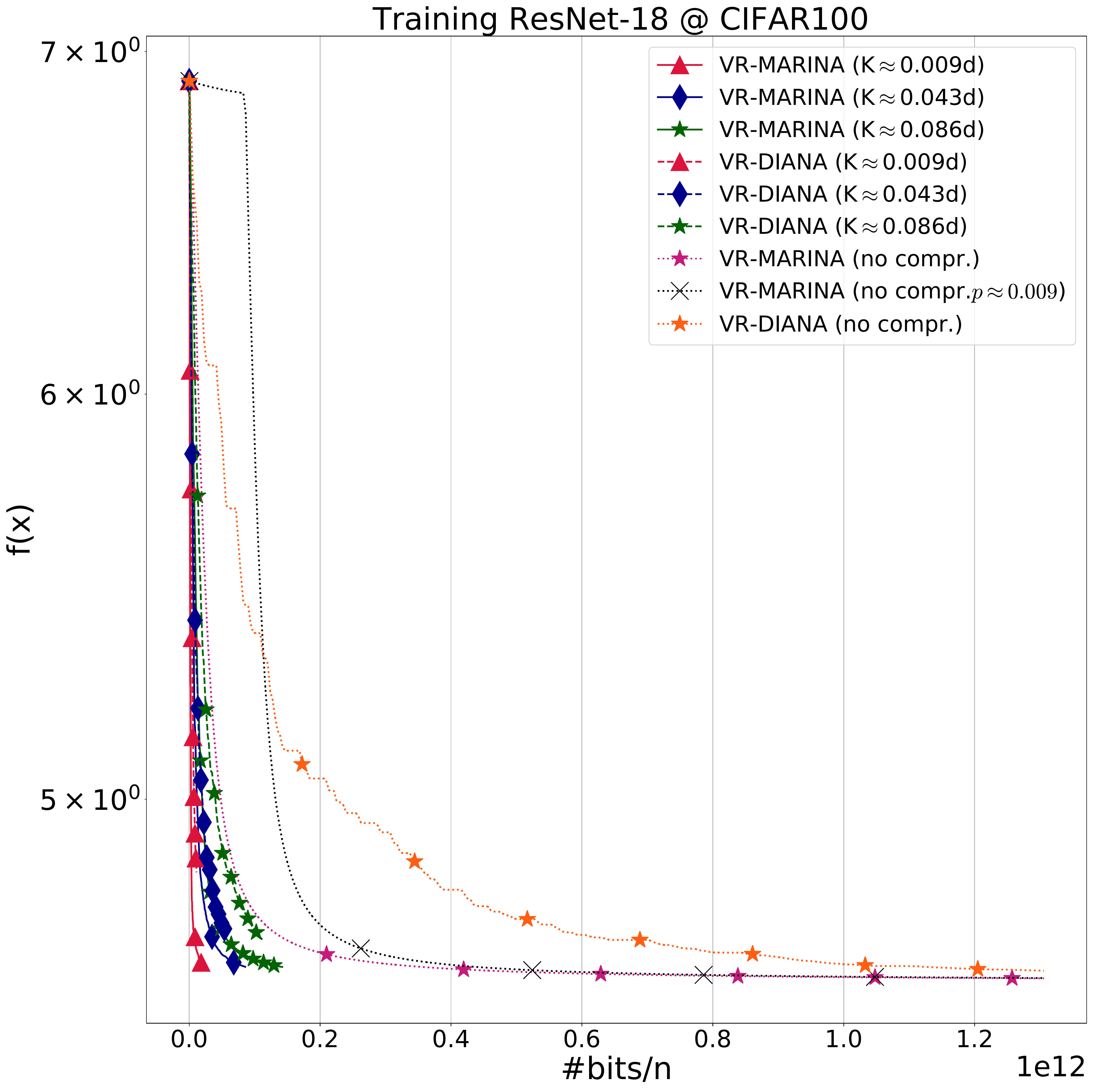}
\includegraphics[width=0.4\textwidth]{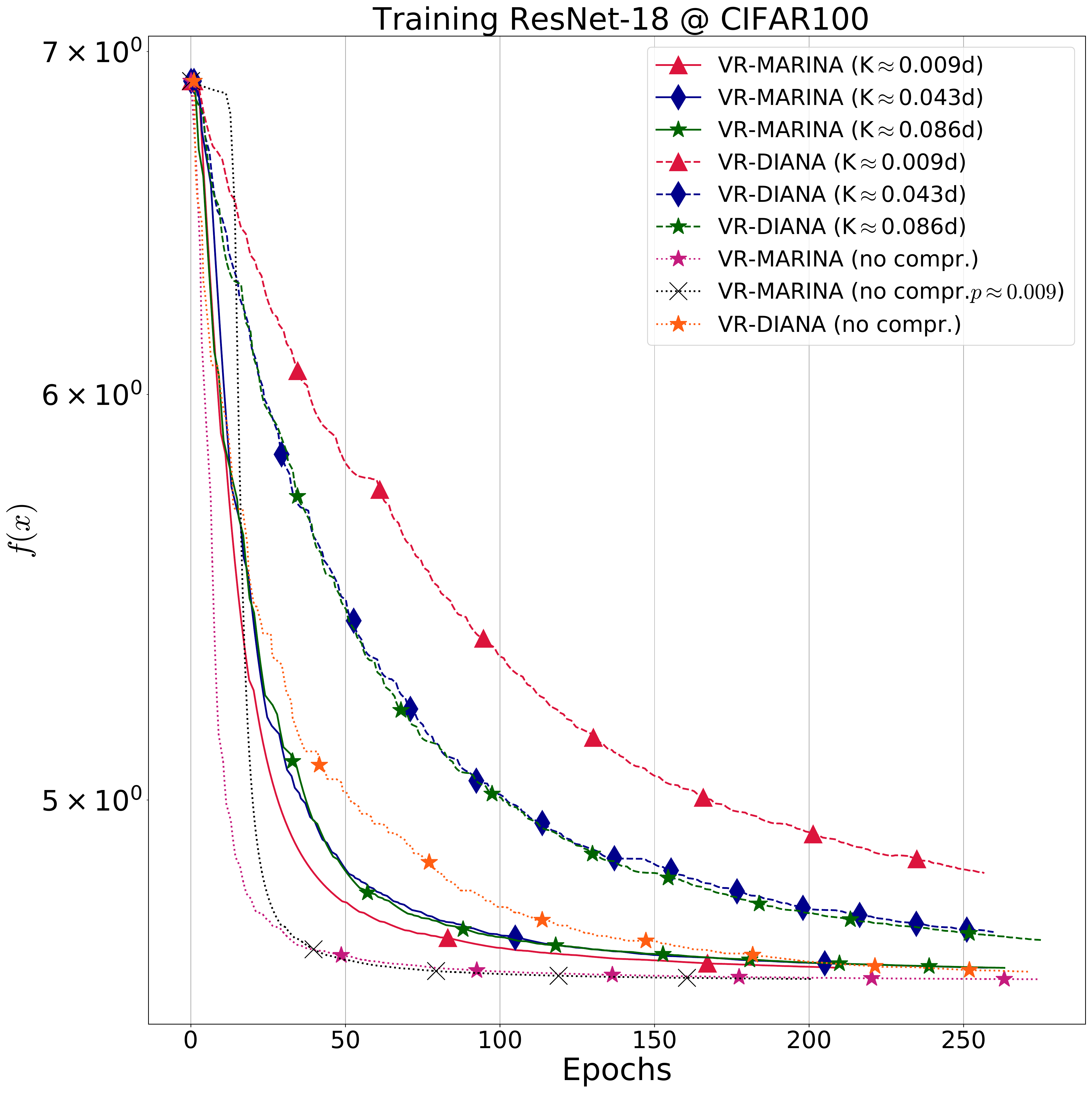}
\includegraphics[width=0.4\textwidth]{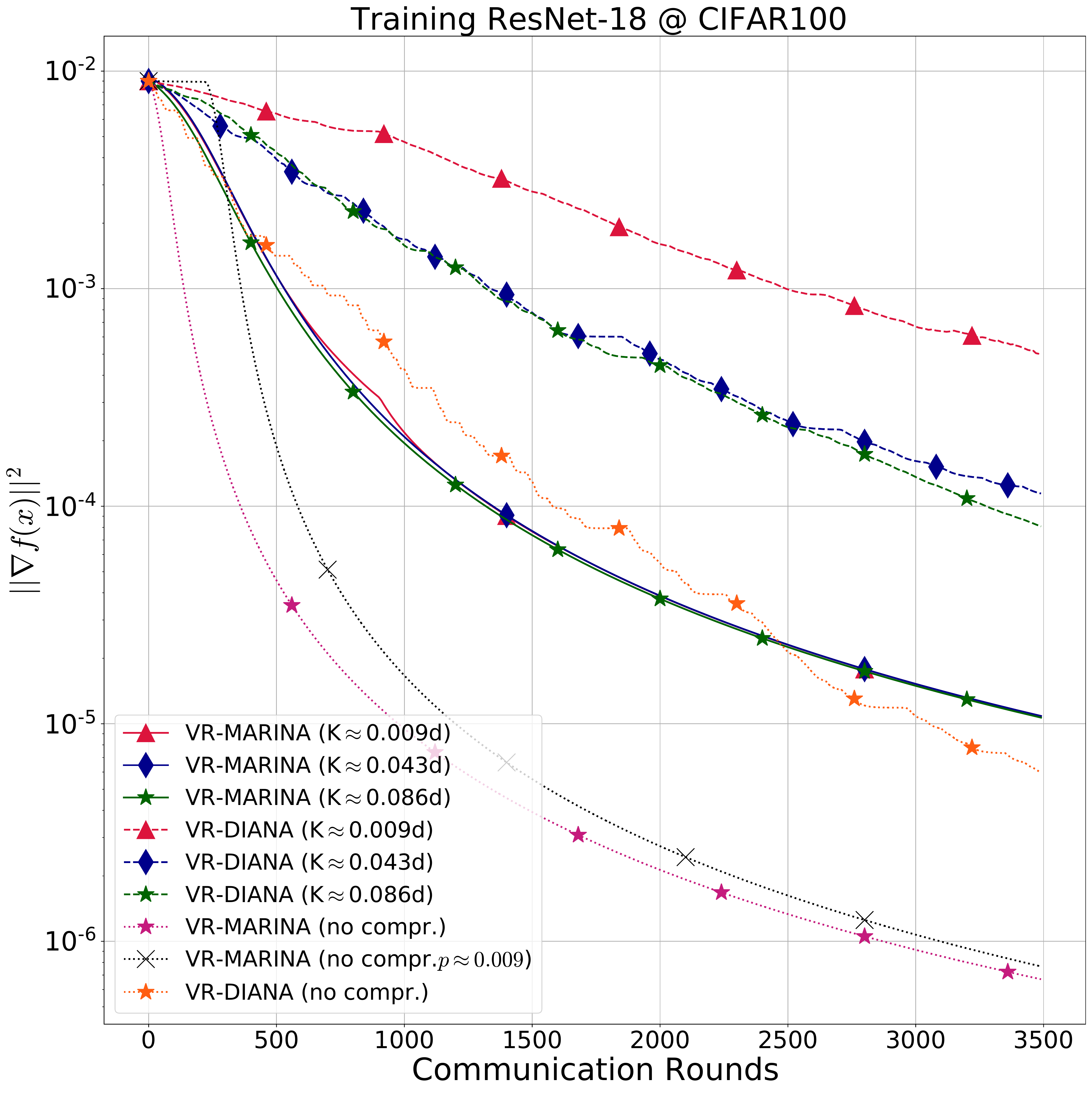}
\includegraphics[width=0.4\textwidth]{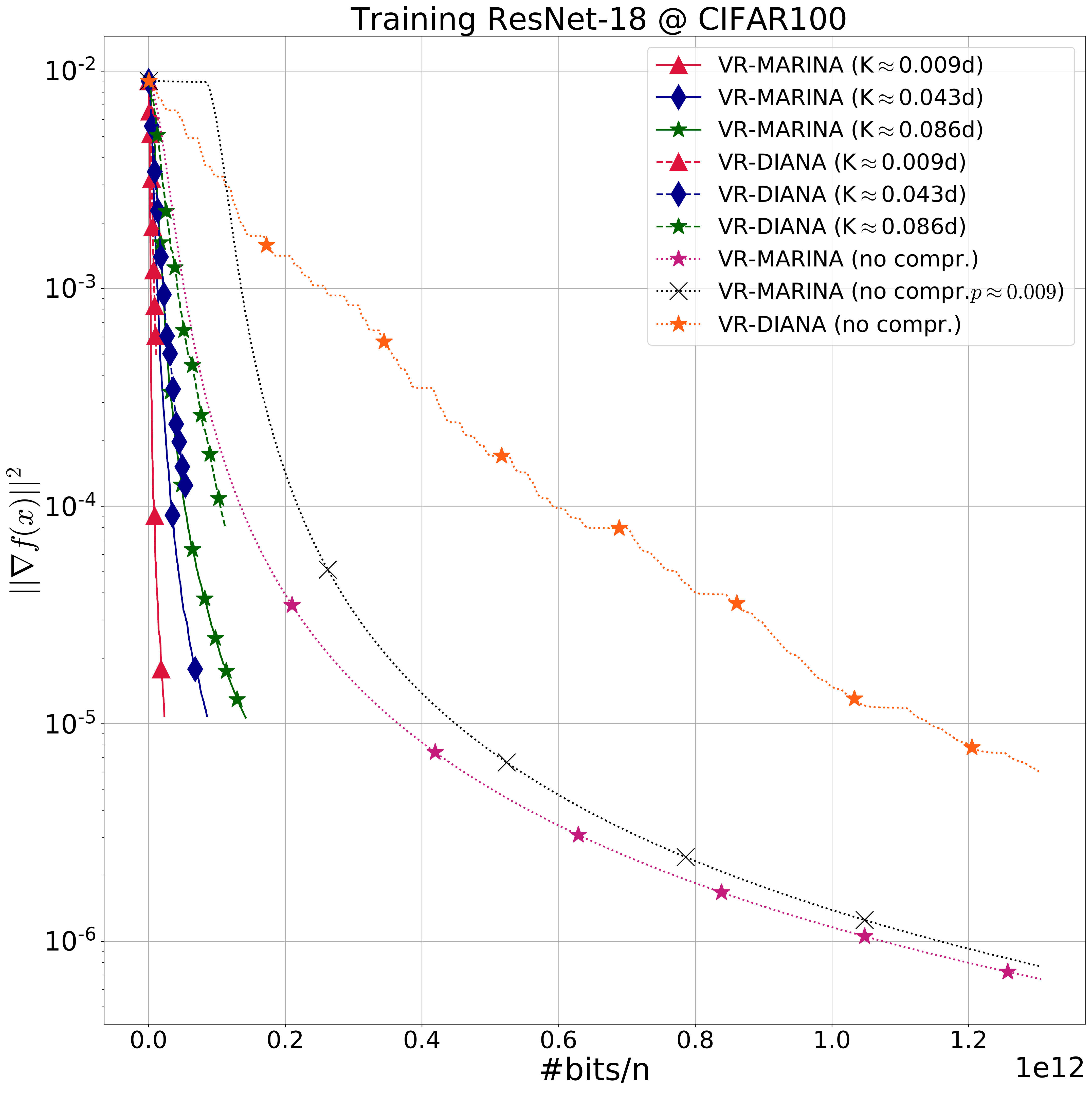}
\includegraphics[width=0.4\textwidth]{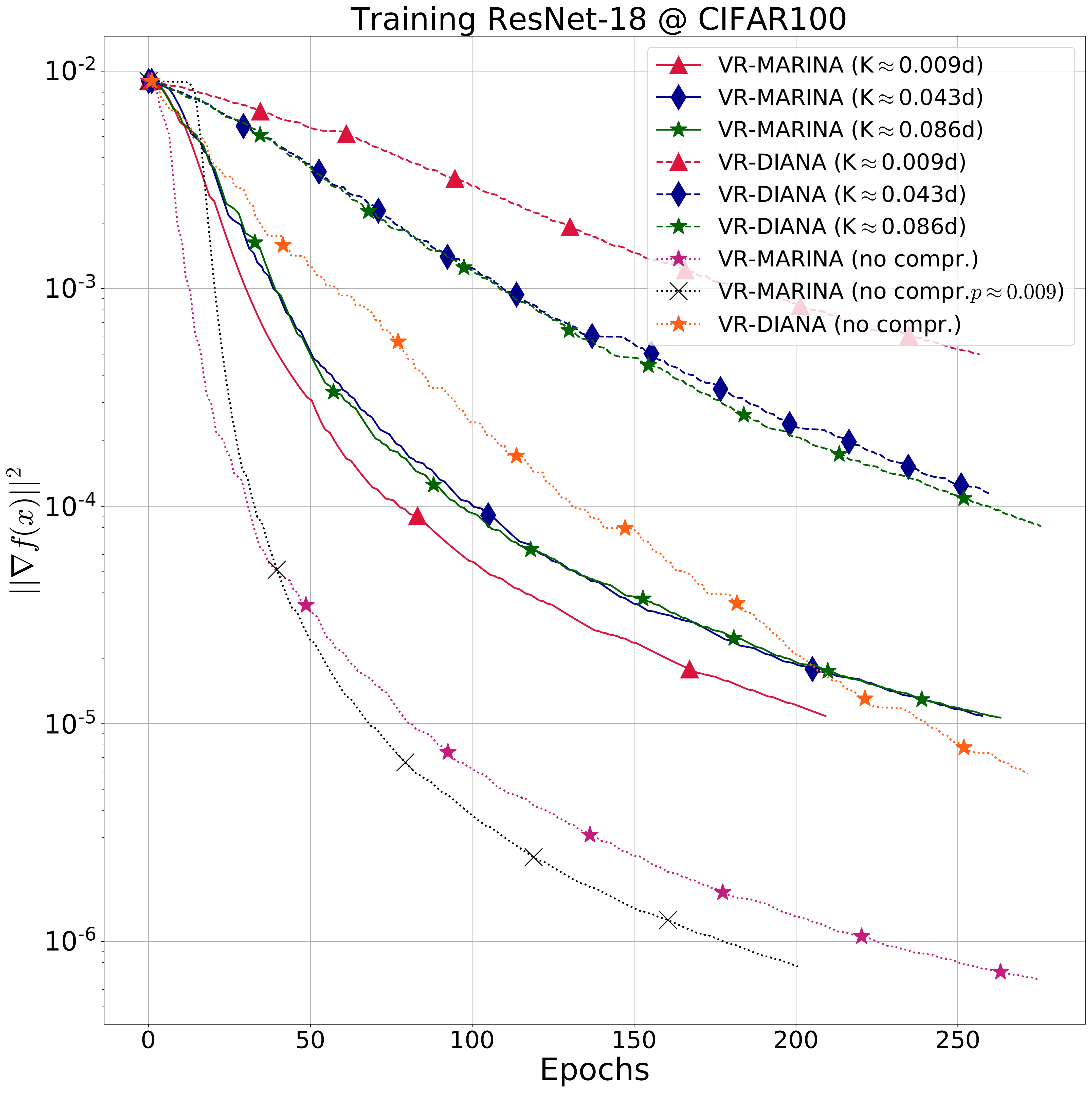}
\caption{Comparison of \algname{VR-MARINA} with \algname{VR-DIANA} on training {\tt ResNet-18} at {\tt CIFAR100} dataset. Number of workers equals $5$. Stepsizes for the methods were tuned and the batchsizes are $\sim \nicefrac{m}{50}$. We used the RandK sparsification operator, the approximate values of $K$ are given in the legends ($d$ is dimension of the problem). We also show the performance of \algname{VR-MARINA} and \algname{VR-DIANA} without compression. }
\label{fig:resnet_at_cifar100_no_compr}
\end{figure}

To emphasize the effect of compression we also run \algname{VR-MARINA} and \algname{VR-DIANA} without compression, see the results in Fig.~\ref{fig:resnet_at_cifar100_no_compr}. First of all, one con notice that the methods do benefit from compression: \algname{VR-MARINA} and \algname{VR-DIANA} with compression converge much faster than their non-comressed versions in terms of the total number of transmitted bits to achieve given accuracy.

Moreover, as Fig.~\ref{fig:resnet_at_cifar100} shows, \algname{VR-MARINA} with $K = $ 100 000 converges faster than \algname{VR-MARINA} with larger $K$ \textit{in terms of the epochs}. That is, the method with more aggresive compression requires less oracle calls to achieve the same accuracy. The reason of such an unusual behavior is the choice of $p$: when $K = $ 100 000 the theoretical choice of $p$ is much smaller than for $K = $ 500 000 and $K = $ 1 000 000. Therefore, in \algname{VR-MARINA} with $K = $ 100 000, the workers compute the full gradients more rarely than in the case of larger $K$. As the result, it turns out, that the total number of oracle calls needed to achieve given accuracy also smaller for $K =$ 100 000 than for larger $K$. Moreover, we see this phenomenon even without applying compression: \algname{VR-MARINA} without compression and with $p$ as in the experiment with \algname{VR-MARINA} with $K = $ 100 000 converges faster than \algname{VR-MARINA} without compression and with theoretical choice of $p$, which is the same as in the case when $K =$ 500 000, 1 000 000, see Table~\ref{tbl:params_resnet}.

%% file: ch6_moshpit.tex
\chapter{Moshpit SGD: Communication-Efficient Decentralized Training on Heterogeneous Unreliable Devices}\label{ch:moshpit}

\section{Introduction}\label{sect:intro}

Many\footnote{We would like to thank Anastasia Koloskova, Liudmila Prokhorenkova and Anton Osokin for helpful feedback and discussions. Finally, we would like to thank Dmitry Afanasiev, Vladimir Aliev, Anand Jayarajan and Michael Solotky for their suggestions on the technical aspects of our study. The computational resources for the experiments were provided by the Amazon Research Awards program and Yandex.} recent influential discoveries in deep learning were enabled by the trend of scaling model and dataset size.
Over the last decade, computer vision has grown from training models with 60 million parameters~\cite{alexnet} on 1.3 million images~\cite{imagenet_cvpr09} to 15 times more parameters~\cite{Kolesnikov2020BigT} and 200 times more training data~\cite{jft-300m}. In natural language processing, the state-of-the-art language models~\cite{gpt3} with 175 billion parameters are trained on over 570GB of texts, and even this does not saturate the model quality~\cite{kaplan2020scaling}.
Training these large models can take years even with a top-of-the-line GPU server~\cite{gpt3costlambda}. As a result, researchers and practitioners often have to run distributed training with multiple machines~\cite{mlperf}.

The dominant approach to distributed deep learning is data-parallel training~\cite{valiant1990bridging}, where each worker processes a fraction of the training batch and then exchanges its gradients with peers. If done naïvely, the gradient exchange can overload the network as the number of workers increases. To combat this issue, modern distributed training algorithms take advantage of communication-efficient protocols, such as all-reduce~\cite{bandwidth_optimal_allreduce}. These protocols 
allow workers to collectively compute the global average gradient with a constant communication overhead, regardless of the total number of peers.
However, this efficiency makes the protocols more fragile: if any single participant fails or takes too long to process its batch, all other nodes will be stalled.

Therefore, scaling all-reduce protocols beyond a couple of servers requires specialized infrastructure with dedicated ultra-high bandwidth networking~\cite{mlperf}.
This kind of infrastructure is notoriously expensive compared to regular
GPU servers or preemptible cloud VMs (see Appendix~\ref{sect:cloud_costs}). Hence, it is tempting to consider distributed training with cheap unreliable instances as a cost-efficient alternative. A similar scenario arises in federated learning~\cite{FL2017-AISTATS}, where one must run distributed training with heterogeneous devices due to privacy concerns.

In both scenarios, participants use a shared network, where both latency and bandwidth can vary drastically due to interference from other users~\cite{variability_azure}\nocite{variability_aws}. Furthermore, compute nodes are also subject to failure (or preemption) caused by factors beyond the protocol's control.

Running large-scale distributed training in these circumstances requires fault- and latency-tolerant algorithms~\cite{lian2017can,sgpush}. Most of these algorithms replace all-reduce averaging with \textbf{gossip}: each participant periodically downloads the latest parameters from his neighbors in a sparsely connected communication graph and averages the results. The updates gradually propagate through the graph over multiple rounds of averaging.
However, the communication required to perform gossip grows linearly with the number of neighbors. Hence, when scaling to hundreds of peers, decentralized SGD has to keep the communication graph sparse, slowing down the convergence.

In this work, we propose an alternative approach. Instead of relying on a predefined communication graph, participants dynamically organize themselves into groups using a fully decentralized matchmaking algorithm which we call {\tt Moshpit All-Reduce}. This strategy allows us to use communication-efficient all-reduce protocols that significantly reduce the network load compared to gossip-based averaging, while still being able to operate in unreliable hardware and network conditions.

Our contributions can be summarized as follows:
\begin{itemize}\vspace{-6px}
    \item We propose {\tt Moshpit All-Reduce} --- a novel decentralized averaging protocol for large-scale training with unreliable communication-constrained devices. According to our analysis, this method has exponential convergence independent of network topology.
    \item Armed with this averaging protocol, we develop {\tt Moshpit SGD} for distributed optimization. We derive convergence rates for this algorithm and establish its equivalence to Centralized (Local) SGD in terms of iteration complexity under realistic assumptions.
    \item Our experiments demonstrate that {\tt Moshpit All-Reduce} is significantly more efficient under network latency. In particular, we train ResNet-50 on ImageNet to 75\% accuracy 1.3 times faster than existing decentralized training algorithms and train ALBERT-large from scratch 1.5 times faster on preemptible cloud VMs.
    \item We release the reference implementation of {\tt Moshpit SGD} and the code for all experiments.\footnote{\href{https://github.com/yandex-research/moshpit-sgd}{\texttt{github.com/yandex-research/moshpit-sgd}}}
\end{itemize}

\vspace{-10px}%
\section{Related Work}\label{sect:related}
\vspace{-4px}
\subsection{Data Parallel Training}\label{sect:related_data_parallel}
\vspace{-4px}

The most popular way to accelerate neural network training with multiple devices is data-parallel training~\cite{valiant1990bridging,goyal2017accurate,You2020Large}. On each optimization step, this strategy splits the training batch among participants. Each participant then runs forward and backward passes to obtain gradients of the objective function on their part of the training batch. After that, we can aggregate the gradients from workers and perform an optimization step. There are two main strategies for this aggregation.

Historically, the first solution to gradient aggregation was to use Parameter Server (PS)~\cite{parameter_server_first}: a separate process or a dedicated server that keeps track of model parameters and optimizer statistics. After each round, the PS accumulates the gradients from each worker and updates the model parameters using SGD or any other optimizer, such as Adam~\cite{adam}. Finally, the server distributes the updated model parameters to workers.

This strategy is robust and easy to implement, but it requires the server to regularly download full model gradients from every single worker. As a result, the parameter server can quickly become a bottleneck for large-scale training~\cite{survey_distributed2}\nocite{survey_distributed}. Since the original PS, researchers have proposed several modifications that reduce the communication load: accumulating multiple batches~\cite{localsgd_first}, compression~\cite{deep,pmlr-v97-koloskova19a}, server sharding~\cite{sharded_ps_first,byteps}. A more detailed overview is given in Appendix~\ref{sect:post_related}.

In turn, many practical distributed training systems have instead switched to averaging with All-Reduce ~\cite{goyal2017accurate,mikami2019massively,shoeybi2019megatron,You2020Large}. This name refers to a collection of protocols originally developed for HPC applications. Workers can follow these protocols to collectively compute the average\footnote{All-Reduce works with any commutative associative operation, such as min, max, or product.} gradient more efficiently than with a central server.

\subsection{Communication-Efficient All-Reduce}\label{sect:related_allreduce}
There are several all-reduce protocols optimized for different network topologies. The simplest one is known as Butterfly All-Reduce~\cite{bandwidth_optimal_allreduce}. Each of $n$ participants splits its local vector into $n$ chunks. Then, $i$-th worker aggregates $i$-th chunk of data from all peers and sends back the averaged chunk.

\vspace{-4pt}
\begin{figure}[h!]
    \centering
    \includegraphics[width=0.6\linewidth]{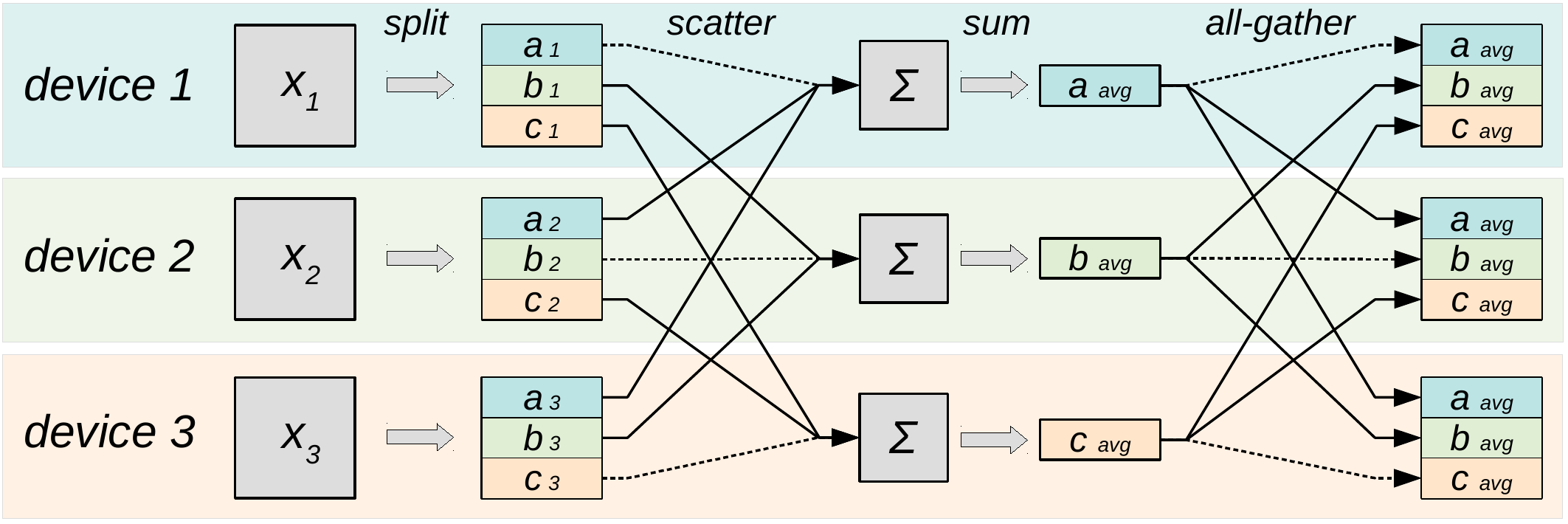}\vspace{-6pt}
    \caption{A schematic illustration of Butterfly All-Reduce.}
    \label{fig:butterfly_allreduce}
\end{figure}
\vspace{-6pt}

As long as the vector size $s$ is greater than $n$, this protocol uses $\cO\left(s \times \frac{n - 1}{n}\right)$ total bandwidth on each worker. However, it requires all-to-all communication, which is not always practical for the HPC infrastructure. Real-world systems typically use Ring or Tree All-Reduce, where each worker only communicates with a small subset of its peers. 

These protocols enable highly efficient and scalable averaging with $\cO(1)$ or $\cO(\log n)$ total communication per worker, but they also share a common drawback: they cannot tolerate node failures or network instability. If any single participant fails to execute its part or takes long to respond, this paralyzes all other workers.

\subsection{Distributed Training in Unstable Conditions}\label{sect:related_unreliable}
Some distributed training applications must deal with unstable network bandwidth and/or unreliable workers. This issue is most prevalent in federated learning~\cite{FL2017-AISTATS,secure_aggregation,federatedlearningatscale}. When dealing with privacy-sensitive data distributed across multiple actors, such as hospital servers~\cite{fed_intel,fed_nvidia} or mobile phones~\cite{fed_google1,fed_google2}, one must train the model using whichever hardware and network available to those actors.

Another important motivational factor is cost: HPC-grade infrastructure can be prohibitively expensive, pushing researchers and practitioners towards commodity servers or preemptible cloud VMs that are significantly cheaper (see Appendix~\ref{sect:cloud_costs}). Another solution is to use volunteer computing~\cite{volunteer_dl_async, mryab} with abundant, but even less reliable, compute resources.

Training under these conditions requires specialized strategies. At a small scale, one can deploy one or a few reliable parameter servers to aggregate the updates from workers. This strategy can tolerate individual node failures~\cite{proteus}, but scales poorly due to the reasons discussed in Section~\ref{sect:related_data_parallel}.

\subsection{Decentralized Training}\label{sect:related_decentralized_training}
If there are too many participants for PS, it can be advantageous to use decentralized SGD via \textbf{gossip-based} averaging \cite{boyd2006randomized,tsitsiklis1984problems,lian2017can}. In this scenario, participants form a sparse graph: each worker periodically downloads parameters from its neighbors and mixes them with local parameters.

In essence, gossip-based averaging removes the communication bottlenecks of PS at the cost of using different local parameters on each peer. That said, gossip-based optimization algorithms can match, and sometimes even outperform, their centralized counterparts in terms of training speed~\cite{scaman2017optimal,scaman2018optimal,scaman2019optimal,lian2017can,sgpush}. However, the convergence properties of gossip averaging and gossip-based optimization methods significantly depend on the communication graph through the spectral properties of the mixing matrix~\cite{xiao2004fast,scaman2019optimal} or the Laplacian matrix of the network~\cite{merris1994laplacian,uribe2020dual}. 

Consequently, as the number of peers increases, gossip-based averaging has to either increase the number of neighbors (hence more communication) or accept slower convergence speed. Because of this, gossip is less communication-efficient than all-reduce algorithms reviewed in Section~\ref{sect:related_allreduce}. However, gossip-based algorithms are more robust to the changes, which makes them applicable to time-varying networks~\cite{nedic2014distributed,nedic2016stochastic,nedic2018network,rogozin2019projected} and federated learning~\cite{ram2009asynchronous,yan2012distributed,yuan2016convergence}.

\section{Method Description}\label{sect:method}


Large-scale training with unreliable participants requires a protocol that is both communication-efficient and fault-tolerant. Unfortunately, existing methods have only  provide one of these properties. To better address our conditions, we propose {\tt Moshpit All-Reduce} --- a fully decentralized averaging protocol that combines the efficiency of all-reduce and the fault tolerance of gossip-based averaging. 

The rest of this section is organized as follows:
\vspace{-8px}
\begin{itemize}
    \item Section~\ref{sect:method_algorithm} describes the protocol and proves its correctness and communication efficiency;
    \item Section~\ref{sect:method_convergence} provides the analysis of the proposed protocol and proves exponential convergence rate for averaging and linear convergence rate for optimization;
    \item Section~\ref{sect:method_implementation_details} contains implementation details for training with heterogeneous compute nodes.
\end{itemize}\vspace{-8pt}

\subsection{Moshpit Averaging}
\label{sect:method_algorithm}

The core idea of {\tt Moshpit All-Reduce} is that workers perform averaging in small independent groups. That way, a single failed participant would only affect his current group. In turn, the composition of each group should be chosen dynamically to converge in the least number of steps.
Ideally, if there are 16 peers with local parameters $\x$, we can average them in 2 rounds, as demonstrated in Figure~\ref{fig:square_allreduce}.

\vspace{-4pt}
\noindent
\begin{minipage}{.5\textwidth}
\centering
\includegraphics[width=\textwidth]{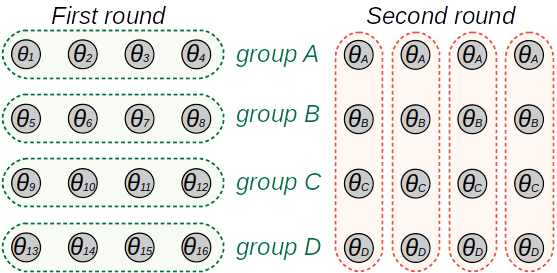}
\captionof{figure}{Example averaging order for 16 peers in 2 rounds. On each round, peers are split into 4 groups that run All-Reduce in parallel.}
\label{fig:square_allreduce}
\end{minipage}
\begin{minipage}{0.5\textwidth}
\begin{algorithm}[H]
\caption{{\tt Moshpit All-Reduce} (for $i$-th peer)}
   \label{alg:ungossip}
\begin{algorithmic}[H]
   \State {\bfseries Input:} parameters $\{\x_j\}_{j=1}^n$, number of peers $n$, $N$, $M$, number of iterations $T$, peer index $i$

   $\x_{i}^0 := \x_i$
   
   $C^0_i :=\texttt{get\_initial\_index(i)}$
   
   \For{$t \in 1 \dots T$}
     \State $\texttt{DHT}[C^{t-1}_i, t].\texttt{add}(\texttt{address}_i)$
     
     \State \texttt{/* wait for peers to assemble */}
     
     \State $\texttt{peers}_t := \texttt{DHT}.\texttt{get}([C^{t-1}_i, t])$ 
     
     \State $\x_{i}^t, c^t_i := \texttt{AllReduce}(\x_{i}^{t - 1}, \texttt{peers}_t)$
     
     \State $C^t_i := (C^{t-1}_i\texttt{[1:]}, c^t_i)$  // same as eq. (1)
   \EndFor
   \State {\bfseries Return} $\x^T_i$
\end{algorithmic}
\end{algorithm}
\end{minipage}

To achieve this in a decentralized system, we use Distributed Hash Tables (DHT) --- a decentralized key-value storage; \autoref{sect:post_related} contains a more detailed description. On each averaging round:
\begin{itemize}\vspace{-8pt}
    \item Each worker computes his group key $C_i$;
    \item Workers add their network addresses to the DHT key corresponding to $C_i$;
    \item Each worker can now fetch a full list of peers that have the same $C_i$ and run All-Reduce with those peers.
\end{itemize}\vspace{-8pt}

Unfortunately, the averaging structure from Figure~\ref{fig:square_allreduce} is impossible to maintain when participants are constantly joining, leaving, and failing. However, we can achieve equivalent results without global structure using a simple rule: \textit{if two peers were in the same group in round $t$, they must choose different groups in round $t {+} 1$.}

A natural way to enforce this rule is to take advantage of the chunk indices from Butterfly All-Reduce (see Figure~\ref{fig:butterfly_allreduce}). Recall that each worker accumulates a \textit{unique} chunk of parameters defined by an index $c_i$. By setting $C_i := c_i$, we can guarantee that any workers that were in the same group at a round $t$ will have different group indices in round $t {+} 1$.

This averaging scheme can be generalized to more than two dimensions in order to fit a larger number of peers or reduce the group size. For a $N$-dimensional hypercube, nodes should find groups of peers that they have not communicated with during $N {-} 1$ previous rounds. To that end, we define $C_i$ as tuples containing chunk indices from $N{-}1$ previous rounds ($t$ denotes the communication round):
\vspace{-2pt}
\begin{equation}
    C^t_i := (c^{t-N+1}_i, c^{t-N+2}_i, \ldots, c^{t }_i).
    \label{eq:group}
\end{equation}
The above intuition can be formalized with Algorithm \ref{alg:ungossip}.
Here, $n$ peers form a virtual $N$-dimensional grid with $M$ peers per row and average their parameters $\x_i$ over $T$ rounds. $\texttt{DHT}[\cdot]$ is a shortcut for using the DHT to add or retrieve values for a given key. In turn, \texttt{AllReduce} denotes running all-reduce to compute the average $\x$ in a given group of peers. The \texttt{get\_initial\_index} function takes the peer index $i$ and returns $N{-}1$ integers
in range $[0, M)$ such as the size of initial groups does not exceed $M$.
That way, the groups formed on all subsequent rounds will also have at most $M$ participants. One possible strategy is:

\vspace{-14pt}
\begin{equation}
    \texttt{get\_initial\_index}(i) = 
    \begin{pmatrix}
           \lfloor i / M^{N{-}1} \rfloor \Mod M \\
         \end{pmatrix}_{j\in \{1,\ \ldots,\ N\}}
    \label{eq:get_initial_index}
\end{equation}

If $n {=} M^N$ and there are no node/network failures, Algorithm~\ref{alg:ungossip} is equivalent to Torus All-Reduce~\cite{torus_allreduce}, achieving the exact average after $N$ rounds of communication (see Appendix~\ref{sect:equiv_to_torus}).
However, our typical use case is far from this perfect scenario; for example, some groups can have less than $M$ members. Furthermore, a peer might fail during all-reduce, causing its groupmates to skip a round of averaging. 
Still, {\tt Moshpit All-Reduce} is applicable even in these conditions:
\begin{theorem}[Correctness]\label{thm:quality_of_avg_deterministic_vectors_0}
If all workers have a non-zero probability of successfully running a communication round and the order of $\texttt{peers}_t$ is random, then all local vectors $\x^t_i$ converge to the global average with probability 1:

\begin{equation}
    \forall i, \Big|\Big|\x^t_i - \frac1n \sum_i \x^0_i\Big|\Big|^2_2 \xrightarrow[t\to\infty]{} 0.
\end{equation}
\end{theorem}
\begin{proof}[Proof (sketch, complete in Appendix~\ref{sect:correctness_proof})]
Running all-reduce with a subset of peers preserves the invariant $\frac1n \sum_i \x^t_i=\frac1n \sum_i \x^{t-1}_i$ and reduces the deviation of $\x^t_i$ from the overall average.
\end{proof}\vspace{-6pt}

\textbf{Complexity.} The matchmaking protocol is implemented over Kademlia DHT~\cite{kademlia}, meaning that each read and write operation needs at most $\cO(\log n)$ requests and $\cO(M)$ bandwidth to load $\texttt{peers}_t$.

After the matchmaking is over, each group runs a single all-reduce round to compute the average. In principle, Moshpit Averaging can use any general-purpose all-reduce protocol. We opted for a butterfly-like version (Figure~\ref{fig:butterfly_allreduce}), as it is simpler than Ring All-Reduce while still being communication-efficient. The communication complexity of this algorithm is $\cO\left(\max(s, M) \times \frac{M - 1}{M}\right)$, where $s$ is the size of vector $\x$. Thus, the total time complexity of Algorithm \ref{alg:ungossip} becomes:
\begin{equation}
    \cO\left(T \times \left[\log_2{n} + M + \max(s, M) \times {\frac{M - 1}{M}}\right]\right).
\end{equation}
This compares favorably to gossip, where network load grows linearly with the number of neighbors.

\vspace{-2pt}
\subsection{Convergence Analysis}\label{sect:method_convergence}
\subsubsection{Mixing Properties of Moshpit Averaging}\label{sect:theory_about_avg}
As stated in the previous section, {\tt Moshpit All-Reduce} computes the exact average when $n = M^N$, which cannot be guaranteed in practice. Therefore, additional analysis is needed to establish how quickly Moshpit Averaging approximates the actual average of $n$ vectors stored on peers.

In the following theorem, we provide such analysis for a simplified version of Moshpit Averaging. One can find the full proof in Appendix~\ref{sec:proof_quality_of_avg_deterministic_vectors}.
\begin{theorem}\label{thm:quality_of_avg_deterministic_vectors}
    Consider a modification of {\tt Moshpit All-Reduce} that works as follows: at each iteration $k\ge 1$, 1) peers are randomly split in $r$ disjoint groups of sizes $M_1^k,\ldots, M_r^k$ in such a way that $\sum_{i=1}^r M_i^k = n$ and $M_i^k \ge 1$ for all $i = 1,\ldots,r$ and 2) peers from each group compute their group average via All-Reduce. Let $\x_1,\ldots,\x_n$ be the input vectors of this procedure and $\x_1^T,\ldots,\x_n^T$ be the outputs after $T$ iterations. Also, let $\overline{\x} = \frac{1}{n}\sum_{i=1}^n\x_i$ Then,
    \begin{equation}
         \hspace{-0.1cm}\EE\left[\frac{1}{n}\sum\limits_{i=1}^n\|\x_i^T - \overline{\x}\|^2\right]= \left(\frac{r-1}{n} + \frac{r}{n^2}\right)^T\frac{1}{n}\sum\limits_{i=1}^n\|\x_i - \overline{\x}\|^2. \label{eq:determ_quality_of_avg}
    \end{equation}
\end{theorem}

\begin{algorithm}[h]
   \caption{{\tt Moshpit SGD}}
   \label{alg:ungossip_Local_SGD}
\begin{algorithmic}[1]
   \State {\bfseries Input:} starting point $\x^0$, learning rate $\gamma > 0$, communication period $\tau \ge 1$
   \For{$k = 0, 1, \ldots$}
   \For{each peer $i\in P_{k+1}$ in parallel}
   \State Compute the stochastic gradient $g_i^k$ at the current point $\x_i^k$
   \If{$k+1 \mod \tau = 0$}
   \State $\x_i^{k+1} = \texttt{Moshpit All-Reduce}_{j\in P_{k+1}}(\x_j^k - \gamma g_j^k)$ for $i$-th peer (Algorithm~\ref{alg:ungossip})
   \Else
   \State $\x_i^{k+1} = \x_i^k - \gamma g_i^k$
   \EndIf
   \EndFor
   \EndFor
\end{algorithmic}
\end{algorithm}\setlength{\textfloatsep}{8pt}

In particular, this result implies that even if workers are randomly split into pairs at each iteration, the simplified version of Moshpit Averaging makes the average distortion (the left-hand side of Equation~\ref{eq:determ_quality_of_avg}) less than $\varepsilon$ in expectation after $\cO\left(\log(\nicefrac{1}{\varepsilon})\right)$ iterations. That is, this algorithm finds $\varepsilon$-accurate average on each node with the rate that \textit{does not} depend on the spectral properties of the communication graph. Since Moshpit Averaging prevents two peers from participating in the same groups during successive iterations, the actual algorithm should find $\varepsilon$-accurate averages on participating peers even faster than Equation~\ref{eq:determ_quality_of_avg} predicts. Moreover, in Appendix~\ref{sec:proof_quality_of_avg_deterministic_vectors} we explain how this result can be generalized to the case when $\{M_i^k\}_{i=1}^n$ and $r$ depends on $k$ or even is random. In Appendix~\ref{sec:mix_rand_proof}, we also provide the guarantees measuring how fast Algorithm~\ref{alg:ungossip} reduces the variance when averaging random vectors.

\subsubsection{Moshpit SGD}\label{sect:optim_theory}
We consider a classical distributed optimization problem
\begin{equation}
    \min\limits_{\x\in\R^d}\left\{f(\x) = \frac{1}{n}\sum\limits_{i=1}^n f_i(\x)\right\}, \label{eq:main_problem}
\end{equation}
where $n$ is the number of workers and worker $i$ has access only to the function $f_i$.

We propose a new algorithm called {\tt Moshpit SGD} to solve this problem (see Algorithm~\ref{alg:ungossip_Local_SGD}). In this algorithm, workers perform independent local SGD steps and periodically synchronize their parameters $\x_i^k$ with other peers using {\tt Moshpit All-Reduce}. Moreover, we define the indices of participating nodes at iteration $k$ as $P_{k+1}$ ($P_0 = \{1,\ldots,n\}$) allowing peers to vanish.

First of all, we list the key assumptions that we use in the convergence analysis of {\tt Moshpit SGD}.
\begin{assumption}[Bounded variance]\label{as:bounded_var}
    We assume that for all $k\ge 0$ and $i=1,\ldots, n$ stochastic gradients $g_i^k$ satisfy $\EE\left[g_i^k\mid \x_i^k\right] = \nabla f_i(\x_i^k)$ and
    \begin{eqnarray}
        \EE\left[\|g_i^k - \nabla f_i(\x_i^k)\|^2\mid \x_i^k\right] &\le& \sigma^2.\label{eq:bounded_variance}
    \end{eqnarray}
\end{assumption}\vspace{-6px}
This assumption is classical in the stochastic optimization literature \cite{Nemirovski-Juditsky-Lan-Shapiro-2009,ghadimi2013stochastic}. We notice that our analysis can be generalized to the settings when the stochastic gradients satisfy less restrictive assumptions such as expected smoothness \cite{gower2019sgd} or have more sophisticated structure similar to \cite{karimireddy2020scaffold} using the theoretical framework from \cite{gorbunov2020local}.

The following assumption controls the averaging properties and the effect of the peers' vanishing.
\begin{assumption}[Averaging quality \& peers' vanishing]\label{as:averaging_quality}
    We assume that the vanishing of peers does not change the global average of the iterates of {\tt Moshpit SGD} too much, i.e., $P_{k+1}\subseteq P_{k}$ and $|P_k| \ge n_{\min}$ for all $k\ge 0$, $|P_{a\tau}| \le 2|P_{a(\tau+1)}|$ for all non-negative integers $a\ge 0$, and there exist such $\widetilde{\x}\in \R^d$ and a sequence of non-negative numbers $\{\Delta_{pv}^k\}_{k\ge 0}$ that $\forall k \ge 0$
    \begin{align}
        \EE\left[\langle\x^{k+1} - \widehat{\x}^{k+1}, \x^{k+1}+\widehat{\x}^{k+1} - 2\widetilde\x\rangle\right] \!\le\! \Delta_{pv}^k\label{eq:stationary_avg_almost}&,f\text{ convex;}\\
        \EE\!\left[\langle\nabla f(\x^k), \x^{k+1}-\widehat{\x}^{k+1}\rangle + L\|\widehat{\x}^{k+1} - \x^{k+1}\|^2\right] \!\le\! \Delta_{pv}^k\label{eq:stationary_avg_almost_2}&,f\text{ non-convex, $L$-smooth,}
    \end{align}
    where $n_k = |P_k|$, $\x^{k+1} = \frac{1}{n_{k+1}}\sum_{i\in P_{k+1}}\x_i^{k+1}$, and $\widehat \x^{k+1} = \frac{1}{n_{k}}\sum_{i\in P_{k}}(\x_i^{k}-\gamma g_i^k)$ for $k\ge 0$. Moreover, we assume that for some $\delta_{aq} \ge 0$ and for all non-negative integers $a\ge 0$
    \begin{eqnarray}
        \EE\left[\frac{1}{n_{a\tau}}\sum\limits_{i\in P_{a\tau}}\|\x_i^{a\tau} - \x^{a\tau}\|^2\right] &\le& \gamma^2\delta_{aq}^2.\label{eq:quality_of_avg}
    \end{eqnarray}
\end{assumption}
If $P_k = P_{k+1} = \{1,\ldots,n\}$ for all $k\ge 0$, i.e., peers do not vanish, then $\x^{k} = \widehat{\x}^{k}$ and properties (\ref{eq:stationary_avg_almost}, \ref{eq:stationary_avg_almost_2}) hold with $\Delta_{pv}^k \equiv 0$ for all $k\ge 0$. Moreover, according to the mixing properties of Moshpit Averaging established in Theorem~\ref{thm:quality_of_avg_deterministic_vectors}, inequality \ref{eq:quality_of_avg} holds after $\cO\left(\log\left(\nicefrac{1}{\gamma^2\delta_{aq}^2}\right)\right)$ iterations of Algorithm~\ref{alg:ungossip}. Therefore, the assumption above is natural and well-motivated.

Under these assumptions, we derive the convergence rates both for convex and non-convex problems. The full statements and complete proofs are deferred to Appendix~\ref{sect:missing_proofs_local_sgd}.
\begin{theorem}[Convex case]\label{thm:cvx_convergence}
    Let $f_1 = \ldots = f_n = f$, function $f$ be $\mu$-strongly convex (Def.~\ref{def:str_cvx}) and $L$-smooth (see Def.~\ref{def:L_smoothness}), and Assumptions~\ref{as:bounded_var}~and~\ref{as:averaging_quality} hold with $\Delta_{pv}^k = \delta_{pv,1}\gamma\mu\EE[\|\x^k-\x^*\|^2] + \gamma^2\delta_{pv,2}^2$ and $\widetilde{\x} = \x^*$, where $\x^* \in \argmin_{\x\in\R^d} f(\x)$ and $\delta_{pv,1}\in [0,1)$, $\delta_{pv,2}\ge 0$. Then there exists a choice of $\gamma$ such that $\EE\left[f(\overline{\x}^K) - f(\x^*)\right]\le \varepsilon$ after $K$ iterations of {\tt Moshpit SGD}, where $K$ equals
    \vspace{-2pt}
    \begin{align}
        \widetilde{\cO}\!\left(\!\frac{L}{(1\!-\!\delta_{pv,1})\mu}\! +\! \frac{\delta_{pv,2}^2\!+\!\nicefrac{\sigma^2}{n_{\min}}}{(1-\delta_{pv,1})\mu\varepsilon}\! +\! \sqrt{\frac{L((\tau\!-\!1)\sigma^2\!+\!\delta_{aq}^2)}{(1\!-\!\delta_{pv,1})^2\mu^2\varepsilon}}\!\right)&,\ \mu>0;\\
        \cO\!\left(\!\frac{LR_0^2}{\varepsilon}\!+\! \frac{R_0^2(\delta_{pv,2}^2\!+\!\nicefrac{\sigma^2}{n_{\min}})}{\varepsilon^2}\!+\! \frac{R_0^2\!\sqrt{L\!(\!(\tau\!-\!1)\!\sigma^2\!+\!\delta_{aq}^2)}}{\varepsilon^{\nicefrac{3}{2}}}\!\right)&,\ \mu=0,
    \end{align}
    \vspace{-2pt}
    where $\overline{\x}^K = \frac{1}{W_K}\sum\limits_{k=0}^K\frac{1}{n_k}\sum\limits_{i\in P_k} w_k \x_i^k$, $w_k = (1-\gamma\mu)^{-(k+1)}$, $W_K = \sum_{k=0}^Kw_k$, $R_0 = \|\x^0 - \x^*\|$ and $\widetilde{\cO}(\cdot)$ hides constant and $\log(\nicefrac{1}{\varepsilon})$ factors.
\end{theorem}
That is, if $\delta_{pv,1} \le \nicefrac{1}{2}$, $n_{\min} = \Omega(n)$, $\delta_{pv,2}^2 = \cO(\nicefrac{\sigma^2}{n_{\min}})$, and $\delta_{aq}^2 = \cO((\tau-1)\sigma)$, then {\tt Moshpit SGD} has the same iteration complexity as Local-SGD in the homogeneous case \cite{khaled2020tighter,woodworth2020local}. However, the averaging steps of {\tt Moshpit SGD} are much faster than those of the parameter-server architecture when the number of peers is large. Also, unlike the state-of-the-art convergence guarantees for Decentralized Local-SGD \cite{koloskova2020unified}, our bounds do not depend on the spectral properties of the communication graph.

\begin{theorem}[Non-convex case]\label{thm:non_cvx_convergence}
    Let $f_1 = \ldots = f_n = f$, function $f$ be $L$-smooth and bounded from below by $f_*$, and Assumptions~\ref{as:bounded_var}~and~\ref{as:averaging_quality} hold with $\Delta_{pv}^k = \delta_{pv,1}\gamma\EE[\|\nabla f(\x^k)\|^2] + L\gamma^2\delta_{pv,2}^2$, $\delta_{pv,1}\in [0,\nicefrac{1}{2})$, $\delta_{pv,2}\ge 0$. Then there exists such choice of $\gamma$ that $\EE\left[\|\nabla f(\x_{\text{rand}}^K)\|^2\right]\le \varepsilon^2$ after $K$ iterations of {\tt Moshpit SGD}, where $K$ equals
    {\begin{eqnarray*}
        \cO\Bigg(\tfrac{L\Delta_0}{(\!1\!-\!2\delta_{pv,1}\!)^2\varepsilon^2}\!\Bigg[\!1\! +\!\tau\sqrt{1\!-\!2\delta_{pv,1}}\! +\! \tfrac{\delta_{pv,2}^2 + \nicefrac{\sigma^2}{n_{\min}}}{\varepsilon^2}\!+\! \tfrac{\sqrt{(1-2\delta_{pv,1})(\delta_{aq}^2+(\tau-1)\sigma^2)}}{\varepsilon}\!\Bigg]\!\Bigg),
    \end{eqnarray*}}
    $\Delta_0 = f(\x^0) - f(\x^*)$ and $\x_{\text{rand}}^K$ is chosen uniformly from $\{\x^0,\x^1,\ldots,\x^{K-1}\}$ defined in As.~\ref{as:averaging_quality}.
\end{theorem}
Again, if $\delta_{pv,1} \le \nicefrac{1}{3}$, $n_{\min} = \Omega(n)$, $\delta_{pv,2}^2 = \cO(\nicefrac{\sigma^2}{n_{\min}})$, and $\delta_{aq}^2 = \cO((\tau-1)\sigma)$, then the above theorem recovers the state-of-the-art results in the non-convex case for Local-SGD \cite{li2019communication,koloskova2020unified}. 

\subsection{Implementation Details}
\label{sect:method_implementation_details}

Training on heterogeneous unreliable hardware also poses a number of engineering challenges. The most obvious one is that the system must be able to recover from node failures. To address these challenges, we use a fully decentralized infrastructure where all information is replicated. When a new worker joins midway through training, it can download the latest model parameters and metadata from any other peer (see \autoref{sect:load_state_from_peers}). Another challenge arises when devices in a group have uneven network bandwidth. In that case, we dynamically adjust the communication load of each peer to avoid being bottlenecked. More information on this procedure can be found in \autoref{sect:load_balancing}.

\vspace{-10pt}
\section{Experiments}\label{sect:experiments}
\vspace{-2pt}
In this section, we
first check the theoretical properties of {\tt Moshpit All-Reduce} in a controlled setup (Section~\ref{sect:experiments_averaging}). Then, we compare {\tt Moshpit SGD} with other distributed methods on practical tasks of image classification and masked language model pretraining (Sections~\ref{sect:experiments_vision} and~\ref{sect:experiments_nlp}).
\vspace{-6pt}
\subsection{Decentralized Averaging}
\label{sect:experiments_averaging}
We aim to verify the convergence and fault tolerance properties proven in Section~\ref{sect:method_convergence}.
To achieve this, we initialize vectors of $512{-}1024$ peers with standard Gaussian noise and run Moshpit Averaging for up to $18$ steps.
We report the average squared difference between the worker parameters and the true average parameters for a $32{\times}32$ grid with varying density and failure rate. We simulate failures by randomly shutting down peers with probability $p$. Failed peers return in the next round of averaging.

\begin{figure}[h]
    \vspace{-4pt}
    \centering
    \includegraphics[width=0.6\linewidth]{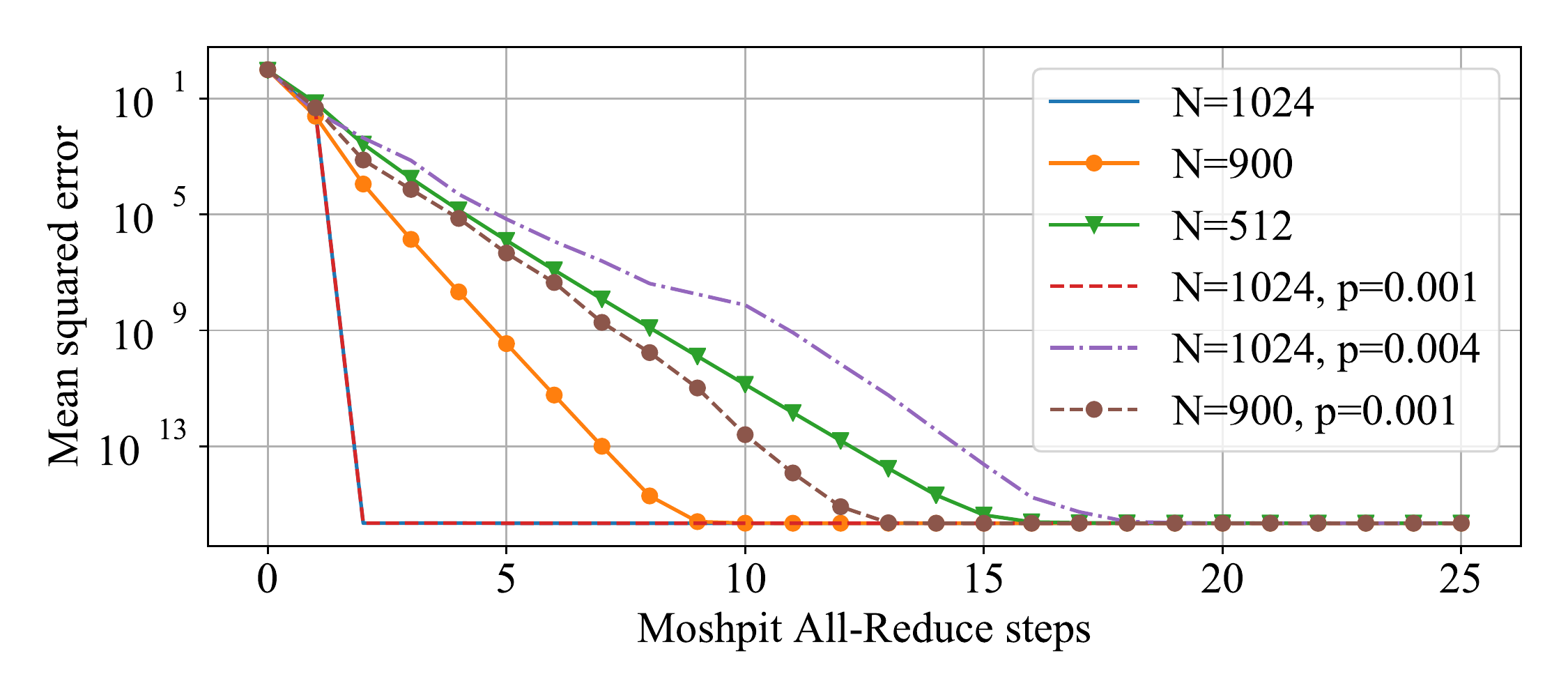}
    \vspace{-12pt}
    \caption{Averaging error for {\tt Moshpit All-Reduce}.}
    \label{fig:log10_std_graph}\vspace{-8pt}
\end{figure}

The results in Figure~\ref{fig:log10_std_graph} outperform the theoretical estimate (Theorem~\ref{thm:quality_of_avg_deterministic_vectors}) in all but one scenario: when $n{=}1024$, the algorithm finds the exact average (within 32-bit precision) in $2$ steps. We also verified that despite worker failures, the global average vector among all participants remains constant throughout each run.
We report additional grid configurations in Appendix~\ref{sect:extra_averaging}.


\subsection{ImageNet Training}\label{sect:experiments_vision}
Here, we evaluate the performance of {\tt Moshpit SGD} in distributed training. More specifically, we train ResNet-50~\cite{he2016deep} on the ILSVRC~\cite{imagenet_cvpr09} dataset, following the training protocol of~\cite{goyal2017accurate}. Trainers use SGD with Nesterov momentum with a batch size of 256 and 32-bit precision regardless of the GPU type\footnote{For GPUs that cannot fit this into memory, we accumulate gradients over 2 batches of 128 examples.}. We evaluate the following training strategies:
\begin{itemize}
    \item \textbf{All-Reduce SGD (AR-SGD)} --- traditional distributed training with all-reduce gradient averaging;\vspace{1px}
    \item \textbf{Asynchronous Decentralized Parallel SGD (AD-PSGD)} --- parallel SGD that runs gossip communication in a cycle: each worker averages parameters with 2 neighbors. Communication rounds are performed in background while the algorithm trains;\vspace{1px}
    \item \textbf{Stochastic Gradient Push (SGP)} --- a more advanced algorithm with an exponential communication graph and push-based communication~\cite{sgpush}. \vspace{1px}
    \item \textbf{Moshpit SGD} --- similar to \textbf{SGP}, but with 1 round of Moshpit Averaging instead of PushSum.
\end{itemize}

We report top-1 validation accuracy as a function of training time in two experimental setups:
\begin{itemize}
    \item \textbf{Homogeneous}: 16 servers with a single Tesla V100-PCIe GPU, 6 CPU cores, and 64GB RAM.
    \item \textbf{Heterogeneous}: a total of 81 GPUs (V100, 1080Ti, and P40) across 64 servers and workstations.\footnote{We provide a detailed configuration in Appendix~\ref{sect:detailed_setup}.}
\end{itemize}

All servers and workstations communicate over the network with 1Gb/s Ethernet (non-dedicated symmetric bandwidth). The machines are located in two data centers and one office within 300 km of one another. The communication latency is 1--6ms depending on the location. To simulate shared usage, at the beginning of each communication round we inject additional latency sampled from the exponential distribution~\cite{sukhov2016generating} with the mean of 100ms.

For {\tt Moshpit SGD}, we use a two-dimensional ``grid'' with 4 and 8 groups for homogeneous and heterogeneous setups respectively. For AD-PSGD, we attempt to compensate for slow convergence by training for 60 more epochs without changing the learning rate schedule. Finally, we only report AR-SGD in the first setup, as it is unsuitable for heterogeneous hardware.

The results in Figure~\ref{fig:all} (Left) demonstrate that the two most efficient strategies for our setting are {\tt Moshpit SGD} and SGP. In the \textbf{homogeneous} setup, Moshpit is only slightly more efficient than SGP, likely due to higher efficiency of all-reduce. This advantage increases to over 30\% for the \textbf{heterogeneous} setup with 64 servers. In turn, AR-SGD demonstrates the best performance per iteration, but its training time is by far the longest due to network latency ($1.5{\times}$ of {\tt Moshpit SGD}). Finally, AD-PSGD predictably shows the best throughput (time per epoch), but achieves lower accuracy even after training for 150 epochs. We report results for smaller setups in Appendix~\ref{sect:extra_classification}. 

\subsection{Masked Language Model Training}
\label{sect:experiments_nlp}
Finally, we evaluate {\tt Moshpit All-Reduce} training performance in the wild with preemptible cloud instances. For this experiment, we perform one of the most resource-demanding tasks in modern deep learning --- unsupervised pretraining of Transformers~\cite{bert,roberta,radford2019language,gpt3}.
We opt for the ALBERT model~\cite{albert} to make better use of communication-constrained devices. This model has fewer trainable parameters due to layer-wise weight sharing.

\begin{figure*}[t]
    \noindent
    \centering
    \vspace{-10pt}
    \includegraphics[width=\textwidth]{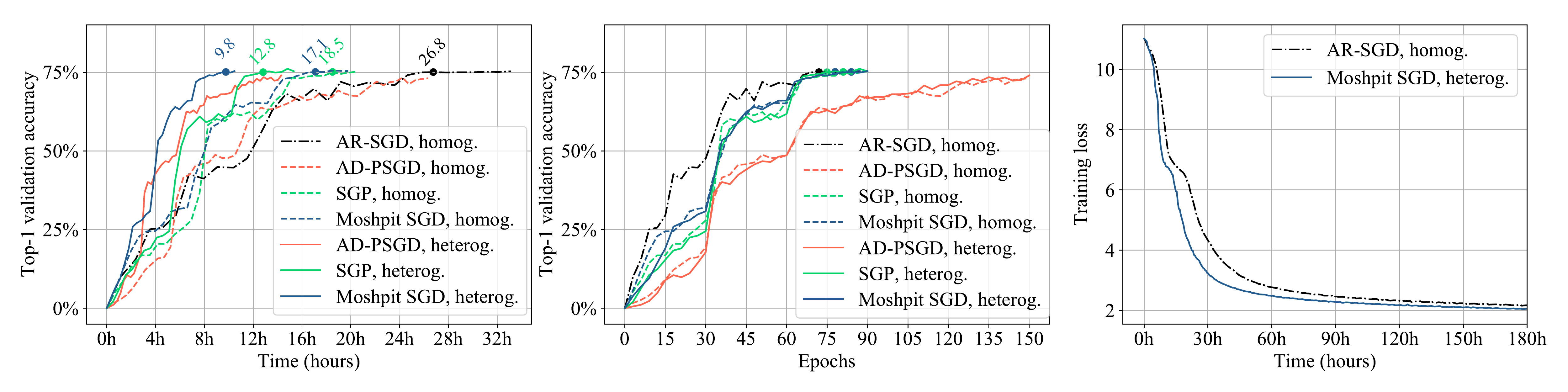}
    \vspace{-20pt}
    \caption{\textbf{(Left, Middle)} ResNet-50 top-1 validation accuracy for ImageNet as a function of training time (left) and epochs (middle). \textbf{(Right)} Full training objective (MLM + SOP) of ALBERT-large on BookCorpus as a function of training time.}
    \label{fig:all}\
\end{figure*}

Specifically, we train ALBERT-large (18M parameters) on the BookCorpus~\cite{bookcorpus} dataset, following the training setup from the original paper. We minimize the masked language modeling loss (MLM) along with the sentence order prediction loss (SOP) using the LAMB optimizer~\cite{You2020Large} with a global batch size of 4096 and sequence length 512. We measure convergence in terms of full training loss~\cite{lin2020multinode,fedus2021switch}. Similarly to Section~\ref{sect:experiments_vision}, we use two training setups:
\begin{itemize}
    \item \textbf{Homogeneous:} a single cloud instance with $8$ Tesla V100-PCIe GPUs and 56 vCPUs;
    \item \textbf{Heterogeneous:} a total of 66 preemptible GPUs, 32 of which are cloud T4, and the remaining 34 are various devices rented on a public marketplace.
\end{itemize}

Despite the fact that the latter setup has almost $3{\times}$ more raw compute\footnote{Based on official performance benchmarks~\cite{nvidia_perf}.}, its hourly rent costs less than the homogeneous setup due to relying on preemptible instances\footnote{Please refer to Appendix~\ref{sect:detailed_setup} for full experimental setups.}. This instance type is much cheaper than regular cloud instances, but it can be interrupted at any time. As a side-effect, the participants in \textbf{heterogeneous} setup are also spread across 3 continents with uneven network bandwidth, ranging from 100Mb/s to 1500Mb/s per worker. These limitations make it impractical to deploy conventional all-reduce protocols. By contrast, the fully decentralized nature of {\tt Moshpit SGD} allows it to operate on unreliable nodes.

In this setup, the participants accumulate gradients over multiple local batches and use DHT to track the global batch size. Once the swarm collectively accumulates gradients over 4096 training samples, it runs 2 rounds of {\tt Moshpit All-Reduce} with $M{=}8$ and $N{=}2$. Unfortunately, training with simple parameter averaging does not converge, likely due to diverging LAMB statistics. To mitigate this issue, workers recover ``pseudo-gradients''~\cite{reddi2021adaptive,chen2020toward} after averaging to update the optimizer statistics.



Figure~\ref{fig:all} (right) demonstrates that {\tt Moshpit SGD} with a fully preemptible fleet of machines trains 1.5 times faster than the traditional data-parallel setup.
The final loss achieved by two training strategies is the same within the margin of error.
A closer investigation reveals that this speedup is entirely explained by the reduced iteration time.
An interesting observation is that the iteration time of {\tt Moshpit SGD} varies between {10--22} seconds, while AR-SGD consistently spends {25}s per step. This can be explained by natural variation in the preemptible fleet size: there were 30--66 active participants depending on resource availability.
\vspace{-6pt}

\section{Conclusion}
In this work, we propose {\tt Moshpit All-Reduce} --- a decentralized averaging protocol intended for distributed optimization. It has favorable theoretical properties when compared to gossip-based approaches and achieves considerable distributed training speedups for image classification and masked language modeling.

Our approach was primarily designed for cloud-based training and federated learning, as well as for distributed training on unreliable instances; future work might explore additional settings, such as collaborative training of neural networks.
Another perspective research direction is to study the combination of the proposed protocol with other techniques that aim for communication efficiency in distributed optimization, such as gradient compression.

%% file: Appendix_Basic_ineqs.tex
\chapter{Basic Facts, Technical Lemmas, and Auxiliary Results}

\section{Standard Definitions from Optimization Theory}
In this section, we provide the most frequently used definitions and simple facts from optimization theory. The proofs of the facts mentioned below are given in \cite{nesterov2018lectures}.

\noindent{\bf Notation.} We use the following notation. $\langle x, y\rangle \eqdef \sum_i x_i y_i$ is the standard Euclidean inner product, and $\norm{x}\eqdef \langle x, x\rangle ^{1/2} $ is the induced $\ell_2$ norm. For simplicity we assume that \eqref{eq:problem_gen} has a unique minimizer, which we denote $x^*$. Let $D_f(x,y)$ denote the \textit{Bregman divergence} associated with $f$: $D_f(x,y) \eqdef f(x) - f(y) - \<\nabla f(y), x-y>$. We  often write $[n]\eqdef \{1,2,\dots,n\}$.

\begin{definition}[$L$-smoothness]\label{def:L_smoothness}
A function $f:\R^n \to \R$ is called $L$-smooth if for all $x,y\in \R^n$, the following inequality holds:
\begin{equation}
    \|\nabla f(x) - \nabla f(y)\| \le L\|x-y\|.\label{eq:L_smoothness_def}
\end{equation}
\end{definition}
If the function $f$ is $L$-smooth, then for all $x,y\in\R^n$
\begin{equation}
    f(y) \le f(x) + \langle\nabla f(x), y-x \rangle + \frac{L}{2}\|y-x\|^2. \label{eq:L_smoothness_cor}
\end{equation}
Next, if $f$ is additionally lower bounded by $f_*$, then for all $x\in\R^d$
\begin{equation}
    \|\nabla f(x)\|^2 \le 2L\left(f(x) - f_*\right). \label{eq:L_smoothness_cor_2}
\end{equation}
Finally, if $f$ is additionally convex, then for all $x,y\in\R^d$
\begin{equation}
    \|\nabla f(x) - \nabla f(y)\|^2 \le 2LD_f(x,y). \label{eq:L_smoothness_cor_3}
\end{equation}

\begin{definition}[$\mu$-strong convexity]\label{def:str_cvx}
    A differentiable function $f:\R^n \to\R$ is called $\mu$-strongly convex if there exists a constant $\mu \ge 0$ such that for all $x,y\in \R^n$
    \begin{equation}
        f(y) \ge f(x) + \langle\nabla f(x), y-x \rangle + \frac{\mu}{2}\|y-x\|^2. \label{eq:str_cvx_def}
    \end{equation}
\end{definition}

\section{Compression and Quantization Operators}

\begin{definition}[Quantization]\label{def:quantization}
	We say that a stochastic mapping $\cQ:\R^d \to \R^d$ is a quantization operator/quantization if there exists  $\omega > 0$ such that for any $x\in\R^d$ , we have
	\begin{equation}
		\EE\left[\cQ(x)\right] = x,\quad \EE\left[\|\cQ(x) - x\|^2\right] \le \omega\|x\|^2. \label{eq:quantization_def}
	\end{equation}
	For the given quantization operator $\cQ(x)$, we define the the expected density as $\zeta_{\cQ} = \sup_{x\in\R^d}\EE\left[\left\|\cQ(x)\right\|_0\right],$ 
	where $\|y\|_0$ is the number of non-zero components of $y\in\R^d$.
\end{definition}
Notice that the expected density is well-defined for any quantization operator since $\left\|\cQ(x)\right\|_0 \le d$.

Below we enumerate some classical compression and quantization operators (see more in \cite{beznosikov2020biased}).
\begin{enumerate}
\item \textbf{TopK sparsification.} This compression operator is defined as follows:
\begin{equation*}
	\cC(x) = \sum\limits_{i=1}^K x_{(i)}e_{(i)}
\end{equation*}
where $|x_{(1)}| \ge |x_{(2)}| \ge \ldots \ge |x_{(d)}|$ are components of $x$ sorted in the decreasing order of their absolute values, $e_1,\ldots,e_d$ is the standard basis in $\R^d$ and $K$ is some number from $[d]$. Clearly, TopK is a biased compression operator. One can show that TopK satisfies \eqref{eq:compression_def} with $\delta = \frac{K}{d}$ \cite{beznosikov2020biased}.
\item \textbf{RandK sparsification} operator is defined as 
\begin{equation*}
	\cQ(x) = \frac{d}{K}\sum\limits_{i\in S} x_{i}e_{i}
\end{equation*}
where $S$ is a random subset of $[d]$ sampled from the uniform distribution on the all subset of $[d]$ with cardinality $K$. RandK is an unbiased compression operator satisfying \eqref{eq:quantization_def} with $\omega = \frac{d}{K}$.
\item \textbf{$\ell_p$-quantization.} By $\ell_2$-quantization we mean the following random operator:
\begin{equation*}
	\cQ(x) = \|x\|_p\text{sign}(x)\circ\xi
\end{equation*}
where $\|x\|_p = \left(\sum_{i=1}^d|x_i|^p\right)^{\nicefrac{1}{p}}$ is an $\ell_p$-norm of vector $x$, $\text{sign}(x)$ is a component-wise sign of vector $x$, $a\circ b$ defines a component-wise product of vectors $a$ and $b$ and $\xi = (\xi_1,\ldots,\xi_d)^\top$ is a random vector such that
\begin{equation*}
\xi_i = \begin{cases}1,&\text{with probability } \frac{|x_i|}{\|x\|_p},\\ 0,&\text{with probability } 1-\frac{|x_i|}{\|x\|_p}. \end{cases}
\end{equation*}
One can show that this operator satisfies \eqref{eq:quantization_def}. In particular, if $p = 2$ it satisfies \eqref{eq:quantization_def} with $\omega = \sqrt{d}-1$ and if $p = \infty$, then $\omega = \frac{1+\sqrt{d}}{2}-1$ (see \cite{mishchenko2019distributed}).
\end{enumerate}

We assume that $\cC$ is any operator which enjoys the following contractive property:  there exists a constant $0< \delta \leq 1$ such that \[\EE\left[\|x-\cC(x)\|^2\right] \leq (1-\delta) \|x\|^2, \qquad \forall x\in \R^d .\]

\section{Basic Inequalities}

For all $a,b,x_1,\ldots,x_n\in\R^d$, $\beta > 0$ and $p\in(0,1]$ the following inequalities hold
\begin{equation}\label{eq:fenchel_young}
	\langle a,b\rangle \le \frac{\|a\|^2}{2\beta} + \frac{\beta\|b\|^2}{2},
\end{equation}
\begin{equation}\label{eq:a-b_a+b}
	\langle a-b,a+b\rangle = \|a\|^2 - \|b\|^2,
\end{equation}
\begin{equation}\label{eq:1/2a_minus_b}
    \frac{1}{2}\|a\|^2 - \|b\|^2 \le \|a+b\|^2,
\end{equation}
\begin{equation}\label{eq:a+b_norm_beta}
    \|a+b\|^2 \le (1+\beta)\|a\|^2 + (1+\nicefrac{1}{\beta})\|b\|^2,
\end{equation}
\begin{equation}\label{eq:a_b_norm_squared}
	\left\|\sum\limits_{i=1}^n x_n\right\|^2 \le n\sum\limits_{i=1}^n\|x_i\|^2,
\end{equation}
\begin{equation}
	\left(1 - \frac{p}{2}\right)^{-1} \le 1 + p, \label{eq:1-p/2_inequality}
\end{equation}
\begin{equation}
	\left(1 + \frac{p}{2}\right)(1 - p) \le 1 - \frac{p}{2}. \label{eq:1+p/2_inequality}
\end{equation}

\section{Identities and Inequalities Involving Random Variables}

\textbf{Variance decomposition.} For a random vector $\xi \in \R^d$ and any deterministic vector $x \in \R^d$ the variance can be decomposed as
\begin{equation}\label{eq:variance_decomposition}
	\EE\left[\left\|\xi - \EE\xi\right\|^2\right] = \EE\left[\|\xi-x\|^2\right] - \left\|\EE\xi - x\right\|^2
\end{equation}

\textbf{Tower property of mathematical expectation.} For random variables $\xi,\eta\in \R^d$ we have
\begin{equation}
	\EE\left[\xi\right] = \EE\left[\EE\left[\xi\mid \eta\right]\right]\label{eq:tower_property}
\end{equation}
under assumption that all expectations in the expression above are well-defined.

\section{Auxiliary Results and Technical Lemmas}\label{sec:tech_lemmas}

The next lemma is used in the analysis of methods with delayed gradients (see Section~\ref{sec:d_sgd}).
\begin{lemma}[Lemma 14 from \cite{stich2020error}]\label{lem:lemma14_stich}
	For any $\tau$ vectors $a_1,\ldots,a_\tau\in\R^d$ and $\xi_1, \ldots, \xi_\tau$ zero-mean random vectors in $\R^d$, each $\xi_t$ conditionally independent of $\{\xi_i\}_{i=1}^{t-1}$ for all $1\le t \le \tau$ the following inequality holds
	\begin{equation}
		\EE\left[\left\|\sum\limits_{t=1}^\tau (a_t + \xi_t)\right\|^2\right] \le \tau\sum\limits_{t=1}^\tau\|a_t\|^2 + \sum\limits_{t=1}^\tau\EE\|\xi_t\|^2. \label{eq:lemma14_stich}
	\end{equation}
\end{lemma}

However, the above lemma is not applicable in the analysis of methods with local steps. To overcome this issue, we propose a generalized version of this result.
\begin{lemma}\label{lem:lemma14_stich_general}
	For any $\tau$ random vectors $\xi_1,\ldots,\xi_\tau\in\R^d$ such that for all $t=2,\ldots,\tau$ random vector $\xi_t$ depends on $\xi_{1},\ldots,\xi_{t-1}$ and does not depend on $\xi_{t+1},\ldots,\xi_{\tau}$ the following inequality holds
	\begin{equation}
		\EE\left[\left\|\sum\limits_{t=1}^\tau\xi_t\right\|^2\right] \le e\tau\sum\limits_{t=1}^\tau\EE\left[\left\|\EE_t[\xi_{t}]\right\|^2\right] + e\sum\limits_{t=1}^\tau\EE\left[\left\|\xi_t-\EE_t[\xi_{t}]\right\|^2\right], \label{eq:lemma14_stich_general}
	\end{equation}
	where $\EE_t[\cdot]$ denotes the conditional expectation $\EE[\cdot\mid \xi_{t-1},\ldots,\xi_1]$.
\end{lemma}
\begin{proof}
	First of all, if $\tau = 1$ then \eqref{eq:lemma14_stich} immediately follows from variance decompostion \eqref{eq:variance_decomposition}. Otherwise ($\tau > 1$) for all $l=1,\ldots,\tau$ we have
	\begin{eqnarray*}
		\EE_l\left[\left\|\sum\limits_{t=1}^l\xi_t\right\|^2\right] &\overset{\eqref{eq:variance_decomposition}}{=}& \left\|\EE_l[\xi_l] + \sum\limits_{t=1}^{l-1}\xi_t\right\|^2 + \EE_l\left[\|\xi_l - \EE_l[\xi_l]\|^2\right]\\
		&\overset{\eqref{eq:a+b_norm_beta}}{\le}& \left(1 + \frac{1}{\tau-1}\right)\left\|\sum\limits_{t=1}^{l-1}\xi_t\right\|^2 + \tau\left\|\EE_l[\xi_l]\right\|^2 + \EE_l\left[\|\xi_l - \EE_l[\xi_l]\|^2\right].
	\end{eqnarray*}
	Taking full mathematical expectation and using tower property \eqref{eq:tower_property} we derive
	\begin{equation*}
		\EE\left[\left\|\sum\limits_{t=1}^l\xi_t\right\|^2\right] \le \left(1 + \frac{1}{\tau-1}\right)\EE\left[\left\|\sum\limits_{t=1}^{l-1}\xi_t\right\|^2\right] + \tau\EE\left[\left\|\EE_l[\xi_l]\right\|^2\right] + \EE\left[\|\xi_l - \EE_l[\xi_l]\|^2\right]
	\end{equation*}
	for all $l=1,\ldots,\tau$. Unrolling the recurrence for $\EE\left[\left\|\sum\limits_{t=1}^l\xi_t\right\|^2\right]$ we obtain
	\begin{eqnarray*}
		\EE\left[\left\|\sum\limits_{t=1}^\tau\xi_t\right\|^2\right] &\le& \tau\sum\limits_{t=1}^\tau \left(1 + \frac{1}{\tau-1}\right)^{\tau-t}\EE\left[\left\|\EE_t[\xi_{t}]\right\|^2\right] + \sum\limits_{t=1}^\tau \left(1 + \frac{1}{\tau-1}\right)^{\tau-t}\EE\left[\left\|\xi_t-\EE_t[\xi_{t}]\right\|^2\right].
	\end{eqnarray*}
	Since $\left(1 + \frac{1}{\tau-1}\right)^{\tau-t} \le \left(1 + \frac{1}{\tau-1}\right)^{\tau-1} \le e$ for all $t=1,\ldots,\tau$ we get \eqref{eq:lemma14_stich}.
\end{proof}

We use the following lemma to derive the final complexity results from Chapter~\ref{ch:ef_sigma_k} in the strongly convex case.
\begin{lemma}[see also Lemma 2 from \cite{stich2019unified}]\label{lem:lemma2_stich}
	Let $\{r_k\}_{k\ge 0}$ satisfy
	\begin{equation}
		r_K \le \frac{a}{\gamma W_K} + c_1\gamma + c_2\gamma^2 \label{eq:lemma2_stich_tech_1}
	\end{equation}
	for all $K\ge 0$ with some constants $a, c_2\ge 0$, $c_1 \ge 0$ where $\{w_k\}_{k\ge 0}$ and $\{W_K\}_{K\ge 0}$ are defined in \eqref{eq:w_k_definition_new}, $\gamma \le \frac{1}{d}$. Then for all $K$ such that
		\begin{eqnarray*}
		\text{either} && \frac{\ln\left(\max\{2,\min\{\nicefrac{a\mu^2K^2}{c_1},\nicefrac{a\mu^3K^3}{c_2}\}\}\right)}{K}\le \min\{\rho_1,\rho_2\}\\
		\text{or} && \frac{1}{h}\le \frac{\ln\left(\max\{2,\min\{\nicefrac{a\mu^2K^2}{c_1},\nicefrac{a\mu^3K^3}{c_2}\}\}\right)}{\mu K}
	\end{eqnarray*}	
	and
	\begin{equation}
		\gamma = \min\left\{\frac{1}{d}, \frac{\ln\left(\max\{2,\min\{\nicefrac{a\mu^2K^2}{c_1},\nicefrac{a\mu^3K^3}{c_2}\}\}\right)}{\mu K}\right\} \label{eq:lemma2_stich_gamma}
	\end{equation}
	we have that
	\begin{equation}
		r_K = \widetilde\cO\left(da\exp\left(-\min\left\{\frac{\mu}{d}, \rho_1, \rho_2\right\}K\right) + \frac{c_1}{\mu K} + \frac{c_2}{\mu^2 K^2}\right). \label{eq:lemma2_stich}
	\end{equation}
\end{lemma}
\begin{proof}
	Since $W_K \ge w_K = (1-\eta)^{-(K+1)}$ we have
	\begin{eqnarray}
		r_K &\le& (1-\eta)^{K+1}\frac{a}{\gamma} + c_1\gamma + c_2\gamma^2 \le \frac{a}{\gamma}\exp\left(-\eta(K+1)\right) + c_1\gamma + c_2\gamma^2.\label{eq:lemma2_stich_tech_2}
	\end{eqnarray}
	Next we consider two possible situations.
	\begin{enumerate}
		\item If $\frac{1}{d} \ge \frac{\ln\left(\max\{2,\min\{\nicefrac{a\mu^2K^2}{c_1},\nicefrac{a\mu^3K^3}{c_2}\}\}\right)}{\mu K}$ then we choose $\gamma = \frac{\ln\left(\max\{2,\min\{\nicefrac{a\mu^2K^2}{c_1},\nicefrac{a\mu^3K^3}{c_2}\}\}\right)}{\mu K}$ and get that
		\begin{eqnarray*}
			r_K &\overset{\eqref{eq:lemma2_stich_tech_2}}{\le}& \frac{a}{\gamma}\exp\left(-\eta(K+1)\right) + c_1\gamma + c_2\gamma^2 \\
			&=& \widetilde\cO\left(a\mu K\exp\left(-\min\left\{\rho_1,\rho_2, \frac{\ln\left(\max\{2,\min\{\nicefrac{a\mu^2K^2}{c_1},\nicefrac{a\mu^3K^3}{c_2}\}\}\right)}{K}\right\}K\right)\right) \\
			&&\quad + \widetilde\cO\left(\frac{c_1}{\mu K} + \frac{c_2}{\mu^2 K^2}\right).
		\end{eqnarray*}
		Since $\frac{\ln\left(\max\{2,\min\{\nicefrac{a\mu^2K^2}{c_1},\nicefrac{a\mu^3K^3}{c_2}\}\}\right)}{K}\le \min\{\rho_1,\rho_2\}$ we have
		\begin{eqnarray*}
			r_K &=& \widetilde\cO\left(a\mu K\exp\left(-\ln\left(\max\left\{2,\min\left\{\frac{a\mu^2K^2}{c_1},\frac{a\mu^3K^3}{c_2}\right\}\right\}\right)\right)\right)\\
			&&\quad + \widetilde\cO\left(\frac{c_1}{\mu K} + \frac{c_2}{\mu^2 K^2}\right)\\
			&=& \widetilde\cO\left(\frac{c_1}{\mu K} + \frac{c_2}{\mu^2 K^2}\right).
		\end{eqnarray*}
		\item If $\frac{1}{d} \le \frac{\ln\left(\max\{2,\min\{\nicefrac{a\mu^2K^2}{c_1},\nicefrac{a\mu^3K^3}{c_2}\}\}\right)}{\mu K}$ then we choose $\gamma = \frac{1}{d}$ which implies that
		\begin{eqnarray*}
			r_K &\overset{\eqref{eq:lemma2_stich_tech_2}}{\le}& da\exp\left(-\min\left\{\frac{\mu}{d},\frac{\rho_1}{4},\frac{\rho_2}{4}\right\}(K+1)\right) + \frac{c_1}{d} + \frac{c_2}{d^2} \\
			&=& \widetilde\cO\left(da\exp\left(-\min\left\{\frac{\mu}{d}, \rho_1, \rho_2\right\}K\right) + \frac{c_1}{\mu K} + \frac{c_2}{\mu^2K^2}\right). 
		\end{eqnarray*}
	\end{enumerate}
	Combining the obtained bounds we get the result. 
\end{proof}

In Chapter~\ref{ch:local_sigma_k}, we apply slightly different result in the strongly convex case.
\begin{lemma}[see also Lemma 2 from \citep{stich2019unified}]\label{lem:lemma2_stich_local}
	Let $\{r_k\}_{k\ge 0}$ satisfy
	\begin{equation}
		r_K \le \frac{a}{\gamma W_K} + c_1\gamma + c_2\gamma^2 \label{eq:lemma2_stich_tech_1}
	\end{equation}
	for all $K\ge 0$ with some constants $a, c_2\ge 0$, $c_1 \ge 0$ where $\{w_k\}_{k\ge 0}$ and $\{W_K\}_{K\ge 0}$ are defined in \eqref{eq:w_k_definition}, $\gamma \le \frac{1}{h}$. Then for all $K$ such that
	\begin{eqnarray*}
		\text{either} && \frac{\ln\left(\max\{2,\min\{\nicefrac{a\mu^2K^2}{c_1},\nicefrac{a\mu^3K^3}{c_2}\}\}\right)}{K}\le \rho\\
		\text{or} && \frac{1}{h}\le \frac{\ln\left(\max\{2,\min\{\nicefrac{a\mu^2K^2}{c_1},\nicefrac{a\mu^3K^3}{c_2}\}\}\right)}{\mu K}
	\end{eqnarray*}		
	 and
	\begin{equation}
		\gamma = \min\left\{\frac{1}{h}, \frac{\ln\left(\max\{2,\min\{\nicefrac{a\mu^2K^2}{c_1},\nicefrac{a\mu^3K^3}{c_2}\}\}\right)}{\mu K}\right\} \label{eq:lemma2_stich_gamma}
	\end{equation}
	we have that
	\begin{equation}
		r_K = \widetilde\cO\left(ha\exp\left(-\min\left\{\frac{\mu}{h}, \rho\right\}K\right) + \frac{c_1}{\mu K} + \frac{c_2}{\mu^2 K^2}\right). \label{eq:lemma2_stich}
	\end{equation}
\end{lemma}
\begin{proof}
	Since $W_K \ge w_K = (1-\eta)^{-(K+1)}$ we have
	\begin{eqnarray}
		r_K &\le& (1-\eta)^{K+1}\frac{a}{\gamma} + c_1\gamma + c_2\gamma^2 \le \frac{a}{\gamma}\exp\left(-\eta(K+1)\right) + c_1\gamma + c_2\gamma^2.\label{eq:lemma2_stich_tech_2}
	\end{eqnarray}
	Next we consider two possible situations.
	\begin{enumerate}
		\item If $\frac{1}{h} \ge \frac{\ln\left(\max\{2,\min\{\nicefrac{a\mu^2K^2}{c_1},\nicefrac{a\mu^3K^3}{c_2}\}\}\right)}{\mu K}$ then we choose $\gamma = \frac{\ln\left(\max\{2,\min\{\nicefrac{a\mu^2K^2}{c_1},\nicefrac{a\mu^3K^3}{c_2}\}\}\right)}{\mu K}$ and get that
		\begin{eqnarray*}
			r_K &\overset{\eqref{eq:lemma2_stich_tech_2}}{\le}& \frac{a}{\gamma}\exp\left(-\eta(K+1)\right) + c_1\gamma + c_2\gamma^2 \\
			&=& \widetilde\cO\left(a\mu K\exp\left(-\min\left\{\rho, \frac{\ln\left(\max\{2,\min\{\nicefrac{a\mu^2K^2}{c_1},\nicefrac{a\mu^3K^3}{c_2}\}\}\right)}{K}\right\}K\right)\right) \\
			&&\quad + \widetilde\cO\left(\frac{c_1}{\mu K} + \frac{c_2}{\mu^2 K^2}\right).
		\end{eqnarray*}
		Since $\frac{\ln\left(\max\{2,\min\{\nicefrac{a\mu^2K^2}{c_1},\nicefrac{a\mu^3K^3}{c_2}\}\}\right)}{K}\le \rho$ we have
		\begin{eqnarray*}
			r_K &=& \widetilde\cO\left(a\mu K\exp\left(-\ln\left(\max\left\{2,\min\left\{\frac{a\mu^2K^2}{c_1},\frac{a\mu^3K^3}{c_2}\right\}\right\}\right)\right)\right)\\
			&&\quad + \widetilde\cO\left(\frac{c_1}{\mu K} + \frac{c_2}{\mu^2 K^2}\right)\\
			&=& \widetilde\cO\left(\frac{c_1}{\mu K} + \frac{c_2}{\mu^2 K^2}\right).
		\end{eqnarray*}
		\item If $\frac{1}{h} \le \frac{\ln\left(\max\{2,\min\{\nicefrac{a\mu^2K^2}{c_1},\nicefrac{a\mu^3K^3}{c_2}\}\}\right)}{\mu K}$ then we choose $\gamma = \frac{1}{h}$ which implies that
		\begin{eqnarray*}
			r_K &\overset{\eqref{eq:lemma2_stich_tech_2}}{\le}& ha\exp\left(-\min\left\{\frac{\mu}{h},\frac{\rho}{4}\right\}(K+1)\right) + \frac{c_1}{h} + \frac{c_2}{h^2} \\
			&=& \widetilde\cO\left(ha\exp\left(-\min\left\{\frac{\mu}{h}, \rho\right\}K\right) + \frac{c_1}{\mu K} + \frac{c_2}{\mu^2K^2}\right). 
		\end{eqnarray*}
	\end{enumerate}
	Combining the obtained bounds we get the result. 
\end{proof}

In the analysis of {\tt Moshpit-SGD}, we also use the following lemma that follows from the previous one.
\begin{lemma}\label{lem:lemma_i_2_gorbunov}
    Let $\{r_k\}_{k\ge 0}$ satisfy
    \begin{equation*}
        r_K \le \frac{a}{\gamma W_K} + c_1\gamma + c_2\gamma^2
    \end{equation*}
    for all $K \ge 0$ with some constants $a,c_2 \ge 0$, $c_1 \ge 0$, where $w_k = (1-\gamma\mu(1-\delta_{pv,1}))^{-(k+1)}$, $W_K = \sum_{k=0}^Kw_k$, $\mu > 0$, $\delta_{pv,1}\in [0,1)$ and $\gamma \le \gamma_0$ for some $\gamma_0 > 0$, $\gamma_0 \le \nicefrac{1}{\mu(1-\delta_{pv,1})}$. Then, for all $K$ such that
    \begin{align*}
        \text{either  } & \frac{\ln\left(\max\left\{2, \min\left\{\nicefrac{a\mu^2(1-\delta_{pv,1})^2K^2}{c_1},\nicefrac{a\mu^3(1-\delta_{pv,1})^3K^3}{c_2}\right\}\right\}\right)}{K} \le 1\\
        \text{or  } & \gamma_0 \le \frac{\ln\left(\max\left\{2, \min\left\{\nicefrac{a\mu^2(1-\delta_{pv,1})^2K^2}{c_1},\nicefrac{a\mu^3(1-\delta_{pv,1})^3K^3}{c_2}\right\}\right\}\right)}{(1-\delta_{pv,1})\mu K}
    \end{align*}
    and
    \begin{equation*}
        \gamma = \min\left\{\gamma_0, \frac{\ln\left(\max\left\{2, \min\left\{\nicefrac{a\mu^2(1-\delta_{pv,1})^2K^2}{c_1},\nicefrac{a\mu^3(1-\delta_{pv,1})^3K^3}{c_2}\right\}\right\}\right)}{(1-\delta_{pv,1})\mu K}\right\}
    \end{equation*}
    we have that
    \begin{equation*}
        r_K = \widetilde{\cO}\left(\frac{a}{\gamma_0}\exp\left(-\gamma_0\mu(1-\delta_{pv,1})K\right) + \frac{c_1}{(1-\delta_{pv,1})\mu K} + \frac{c_2}{(1-\delta_{pv,1})^2\mu^2 K^2}\right).
    \end{equation*}
\end{lemma}

To establish the complexity bounds in the convex case, we apply the lemma below.
\begin{lemma}\label{lem:lemma_technical_cvx}
	Let $\{r_k\}_{k\ge 0}$ satisfy
	\begin{equation}
		r_K \le \frac{a}{\gamma K} + \frac{b_1\gamma}{K} + \frac{b_2\gamma^2}{K} + c_1\gamma + c_2\gamma^2 \label{eq:lemma_technical_cvx_1}
	\end{equation}
	for all $K\ge 0$ with some constants $a> 0$, $b_1, b_2, c_1, c_2 \ge 0$ where $\gamma \le \gamma_0$. Then for all $K$ and
	\begin{equation*}
		\gamma = \min\left\{\gamma_0, \sqrt{\frac{a}{b_1}}, \sqrt[3]{\frac{a}{b_2}}, \sqrt{\frac{a}{c_1 K}}, \sqrt[3]{\frac{a}{c_2 K}}\right\}
	\end{equation*}
	we have that
	\begin{equation}
		r_K = \cO\left(\frac{a}{\gamma_0 K} + \frac{\sqrt{ab_1}}{K} + \frac{\sqrt[3]{a^2b_2}}{K} + \sqrt{\frac{ac_1}{K}} + \frac{\sqrt[3]{a^2c_2}}{K^{\nicefrac{2}{3}}} \right). \label{eq:lemma_technical_cvx_2}
	\end{equation}
\end{lemma}
\begin{proof}
	We have
	\begin{eqnarray*}
		r_K &\le& \frac{a}{\gamma K} + \frac{b_1\gamma}{K} + \frac{b_2\gamma^2}{K} + c_1\gamma + c_2\gamma^2\\
		&\le& \frac{a}{\min\left\{\gamma_0, \sqrt{\frac{a}{b_1}}, \sqrt[3]{\frac{a}{b_2}}, \sqrt{\frac{a}{c_1 K}}, \sqrt[3]{\frac{a}{c_2 K}}\right\}K} + \frac{b_1}{K}\cdot\sqrt{\frac{a}{b_1}} + \frac{b_2}{K}\cdot\sqrt[3]{\frac{a}{b_2}}\\
		&&\quad + c_1\cdot\sqrt{\frac{a}{c_1 K}} + c_2 \left(\sqrt[3]{\frac{a}{c_2 K}}\right)^2\\
		 &=& \cO\left(\frac{a}{\gamma_0 K} + \frac{\sqrt{ab_1}}{K} + \frac{\sqrt[3]{a^2b_2}}{K} + \sqrt{\frac{ac_1}{K}} + \frac{\sqrt[3]{a^2c_2}}{K^{\nicefrac{2}{3}}} \right).
	\end{eqnarray*}
\end{proof}

Next, we use the following result in the analysis of methods presented in Chapter~\ref{ch:marina}.
\begin{lemma}[Lemma~2 from \cite{li2020page}]\label{lem:lemma_2_page}
	Assume that function $f$ is $L$-smooth and $x^{k+1} = x^k - \gamma g^k$. Then 
	\begin{equation}
		f(x^{k+1}) \le f(x^k) - \frac{\gamma}{2}\|\nabla f(x^k)\|^2 - \left(\frac{1}{2\gamma} - \frac{L}{2}\right)\|x^{k+1}-x^k\|^2 + \frac{\gamma}{2}\|g^k - \nabla f(x^k)\|^2. \label{eq:key_inequality}
	\end{equation}
\end{lemma}

Finally, in the analysis of {\tt Moshpit-SGD}, we use the following classical result establishing contractiveness of the gradient descent step.
\begin{lemma}[Lemma 6 from \cite{karimireddy2020scaffold}]\label{lem:gd_contraction}
    For any $L$-smooth and $\mu$-strongly convex function $f:\R^n\to\R$, points $x,y\in \R^n$, and stepsize $\gamma \in (0,\nicefrac{1}{L}]$, the following inequality holds:
    \begin{equation}
        \|x - \gamma\nabla f(x) - y + \gamma\nabla f(y)\|^2 \le (1-\gamma\mu)\|x-y\|^2. \label{eq:gd_contraction}
    \end{equation}
\end{lemma}


%% file: Appendix_ef_sigma_k.tex
\chapter{Appendix for Chapter \ref{ch:ef_sigma_k}}\label{app:ef_sigma_k}
\section{Missing Plots} \label{sec:moreexperiments}
\subsection{Compressing Stochastic Gradients}
\begin{figure}[H]
    \centering
    \includegraphics[width=0.32\textwidth]{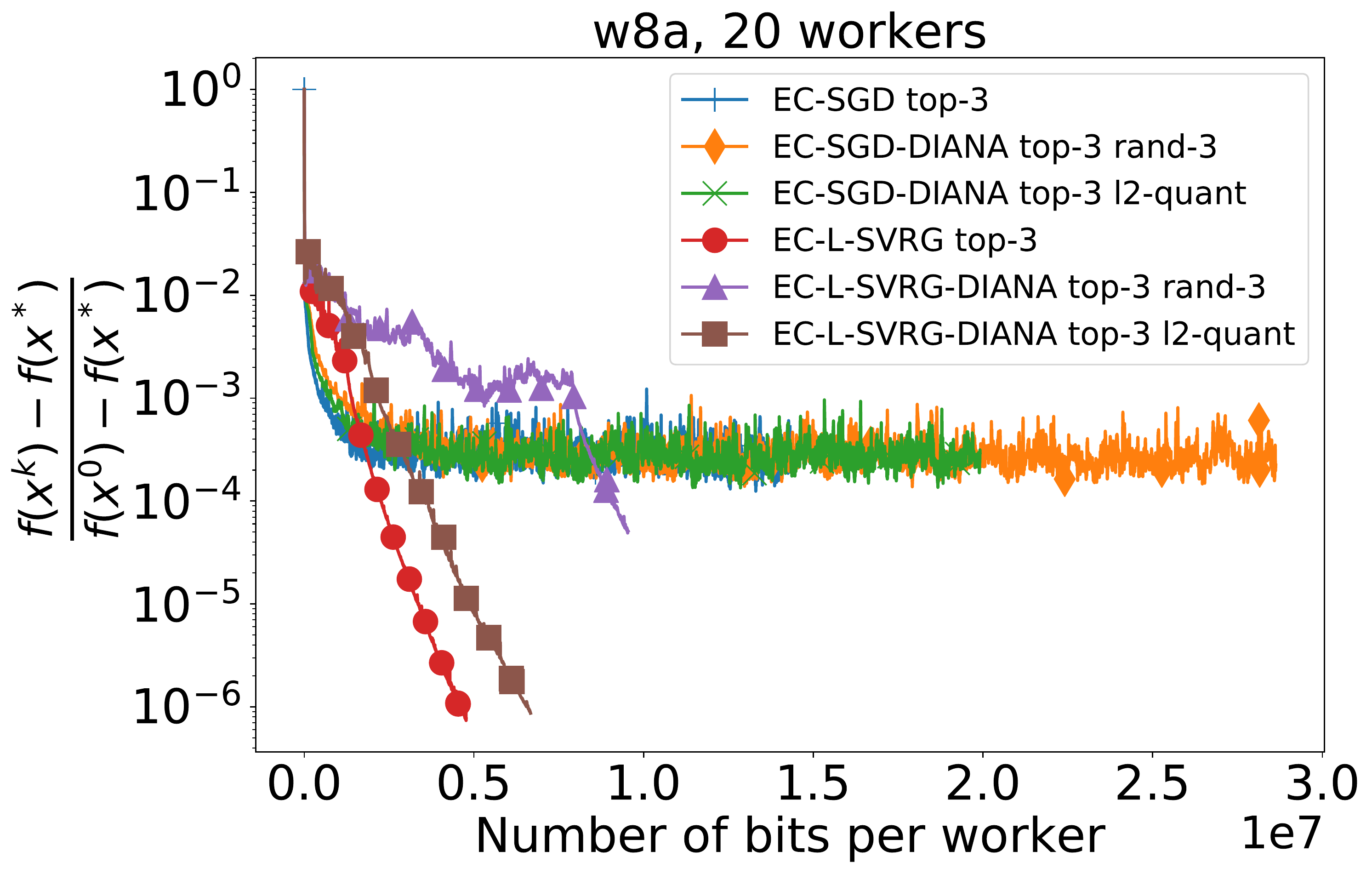}
	\includegraphics[width=0.32\textwidth]{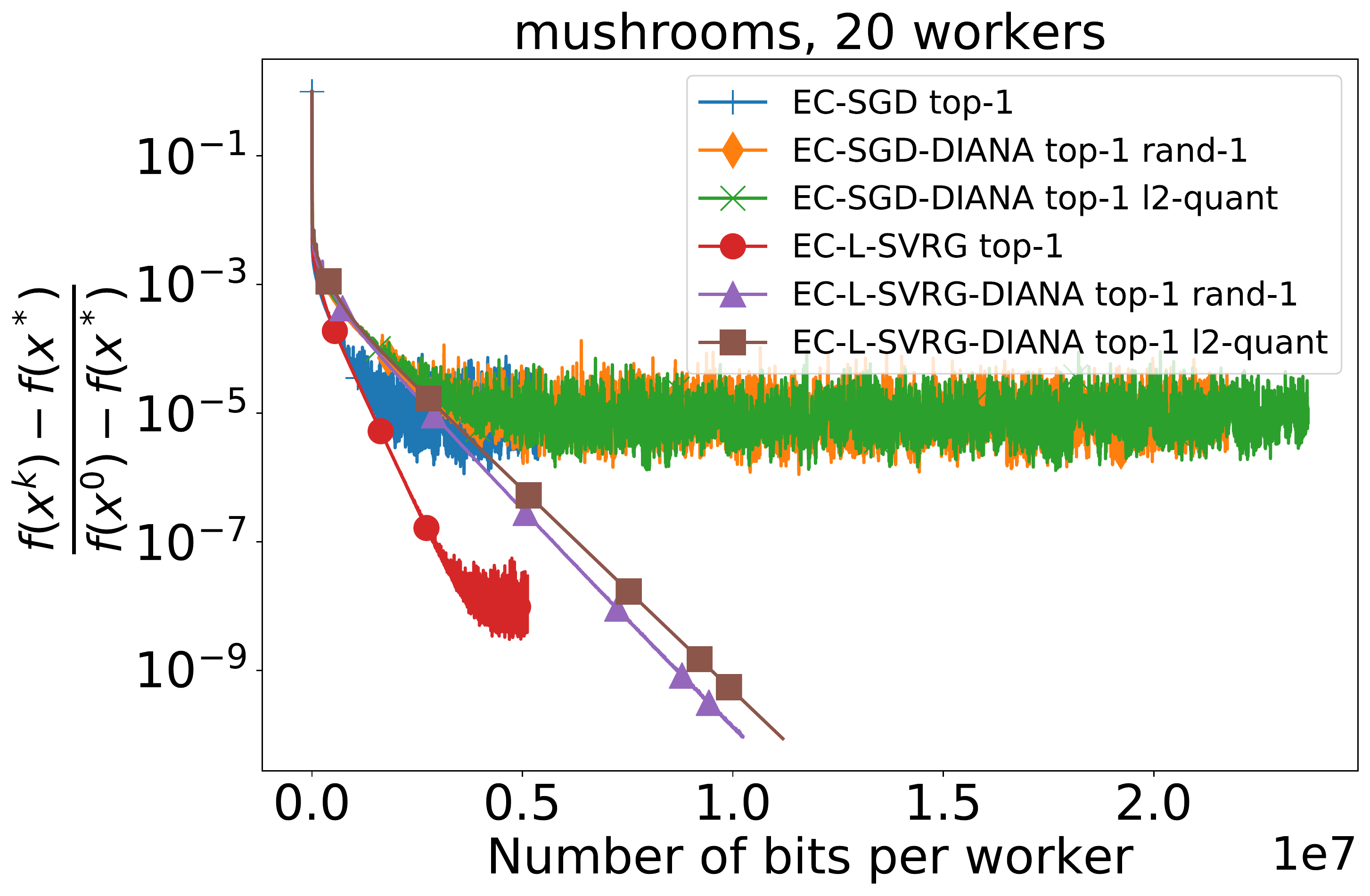}    
	\includegraphics[width=0.32\textwidth]{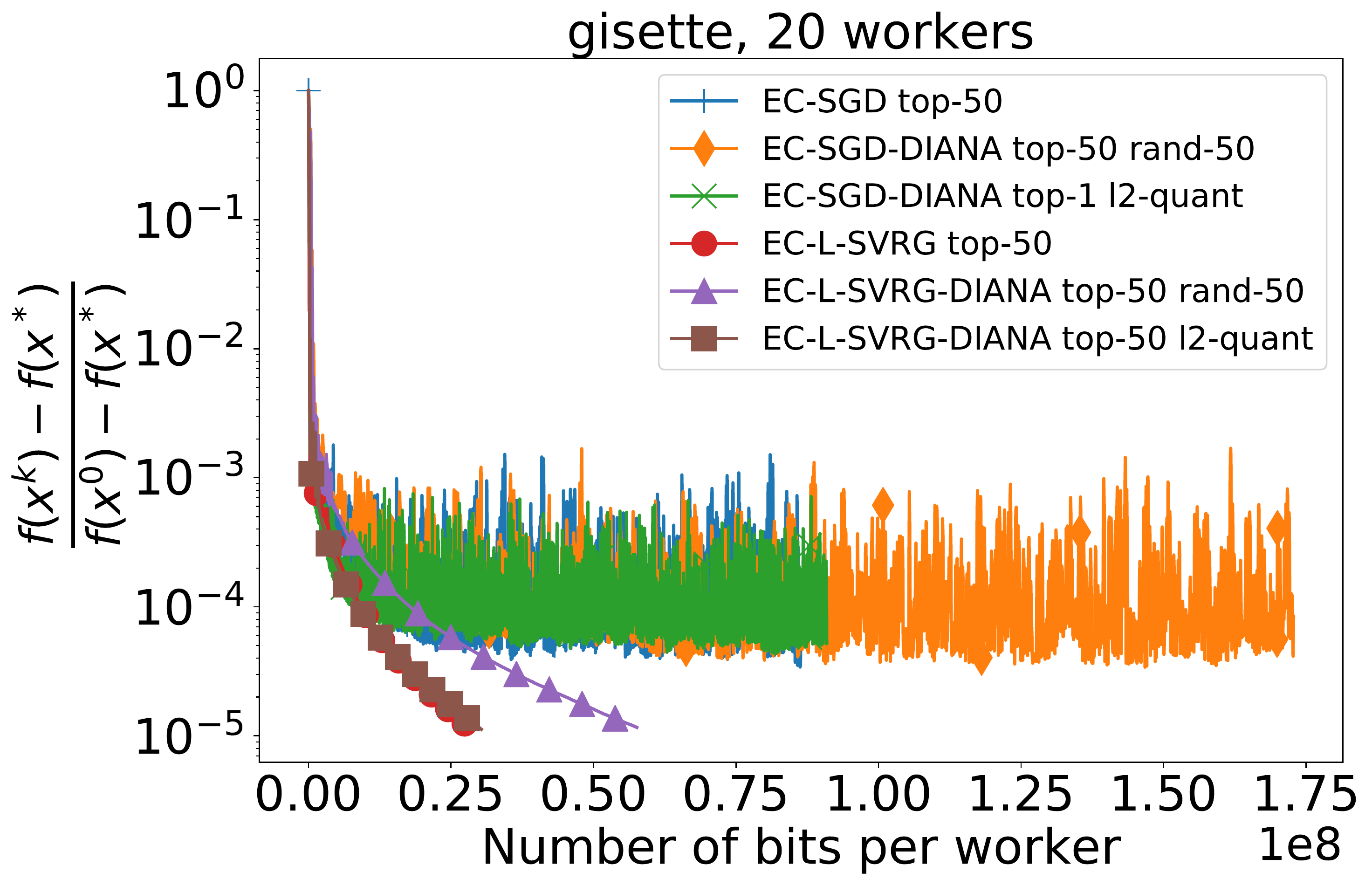}    
    \\
    \includegraphics[width=0.32\textwidth]{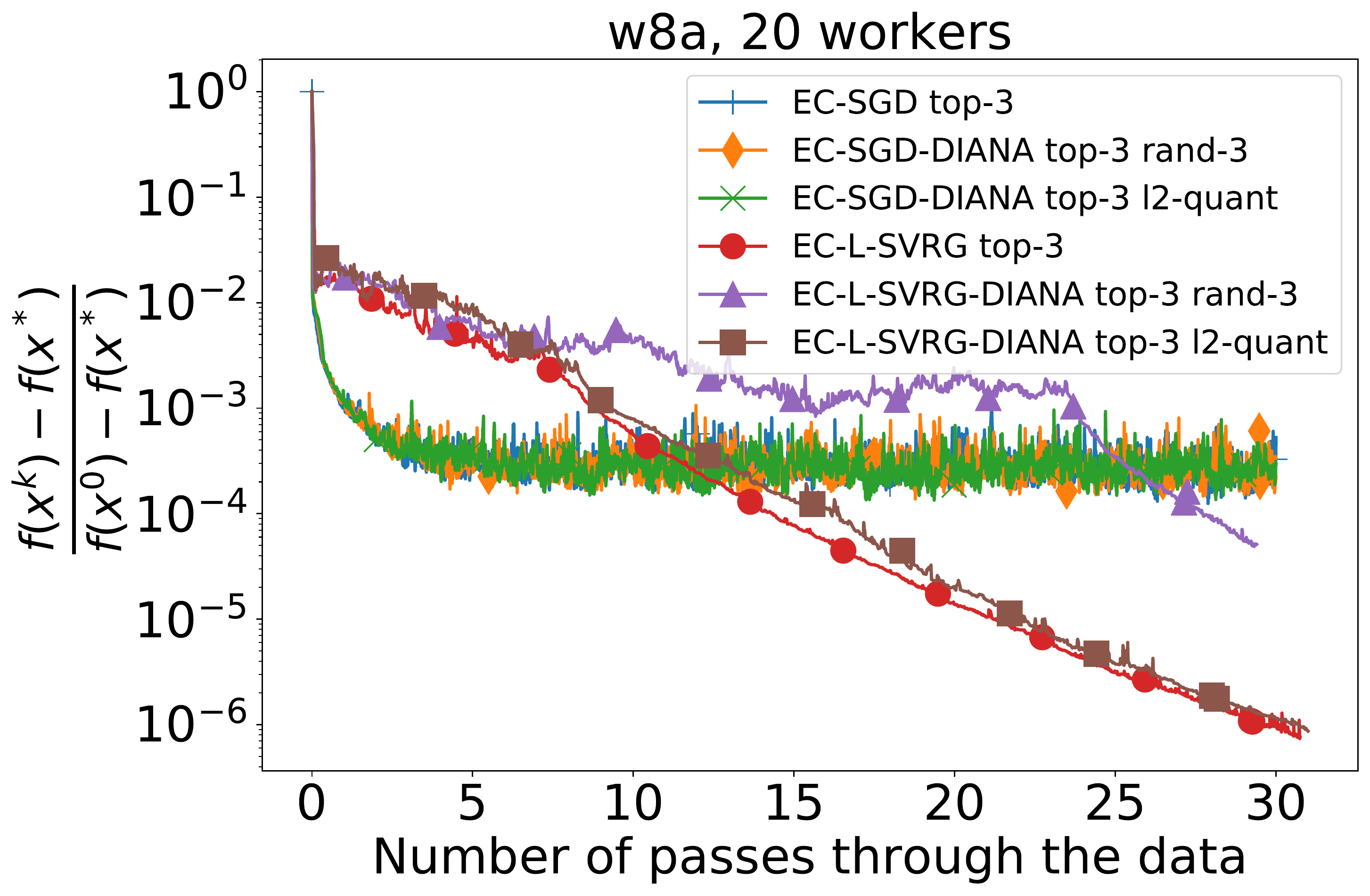}
	\includegraphics[width=0.32\textwidth]{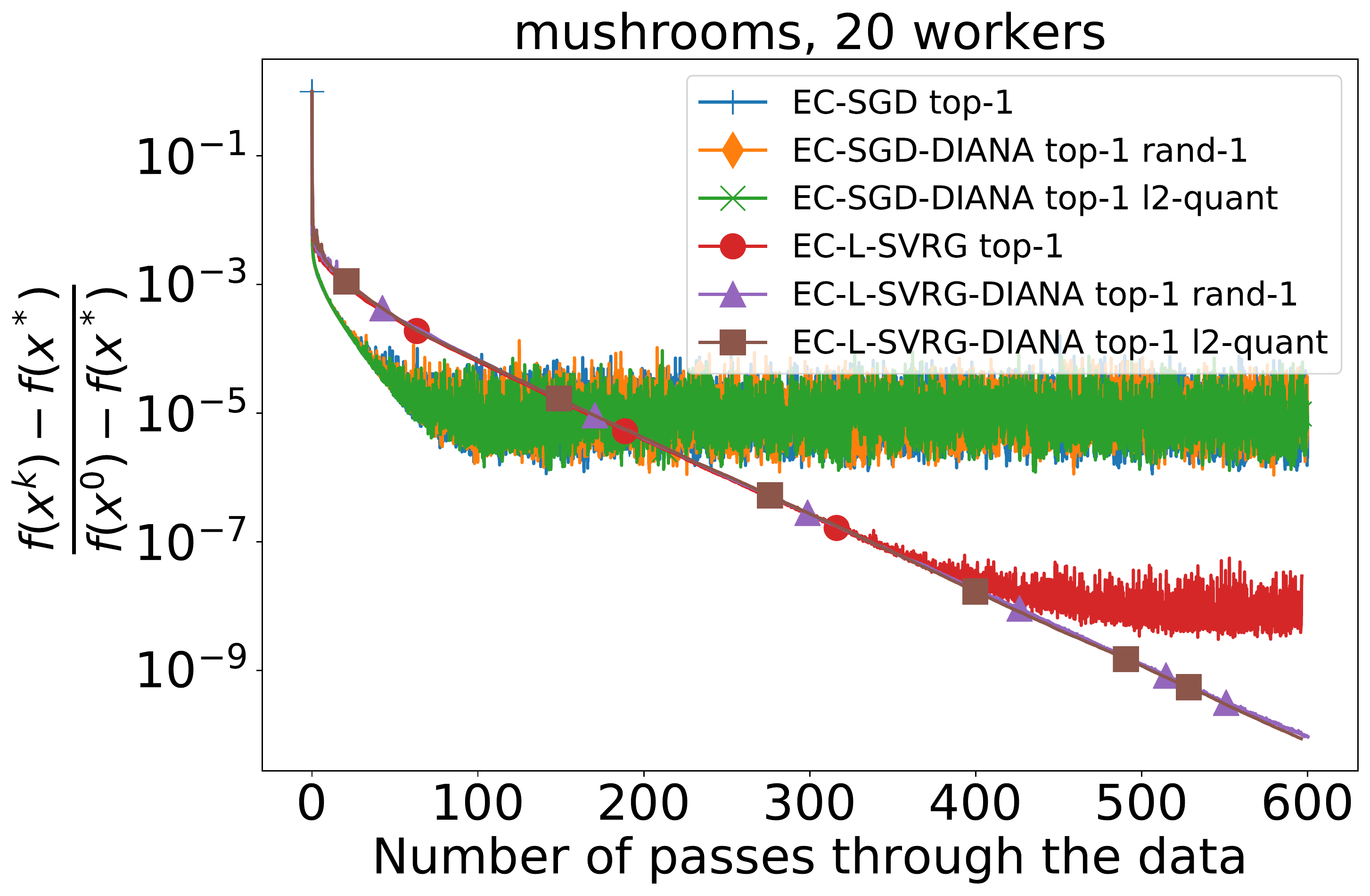}        
	\includegraphics[width=0.32\textwidth]{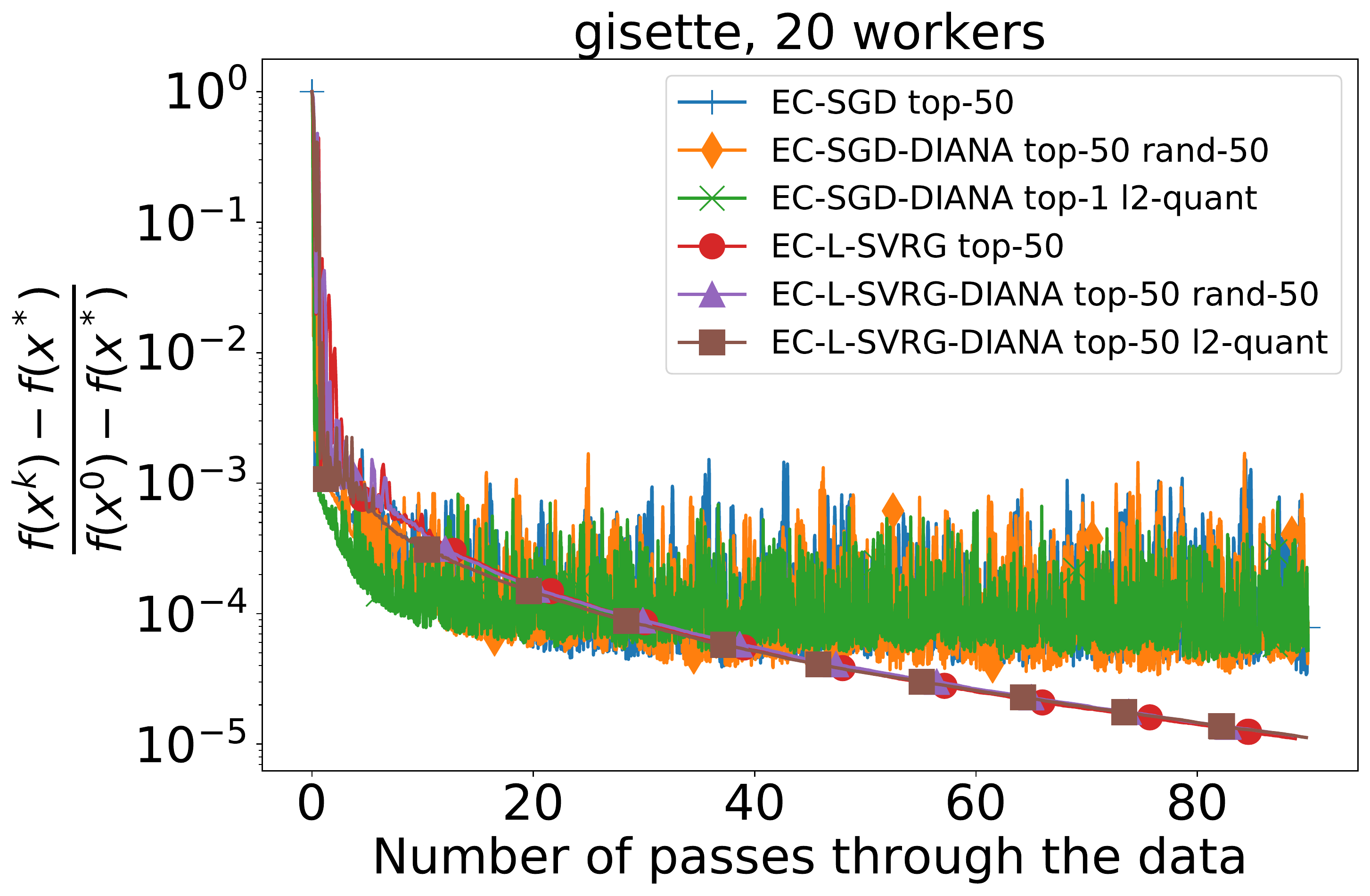}          
    \\
	\includegraphics[width=0.32\textwidth]{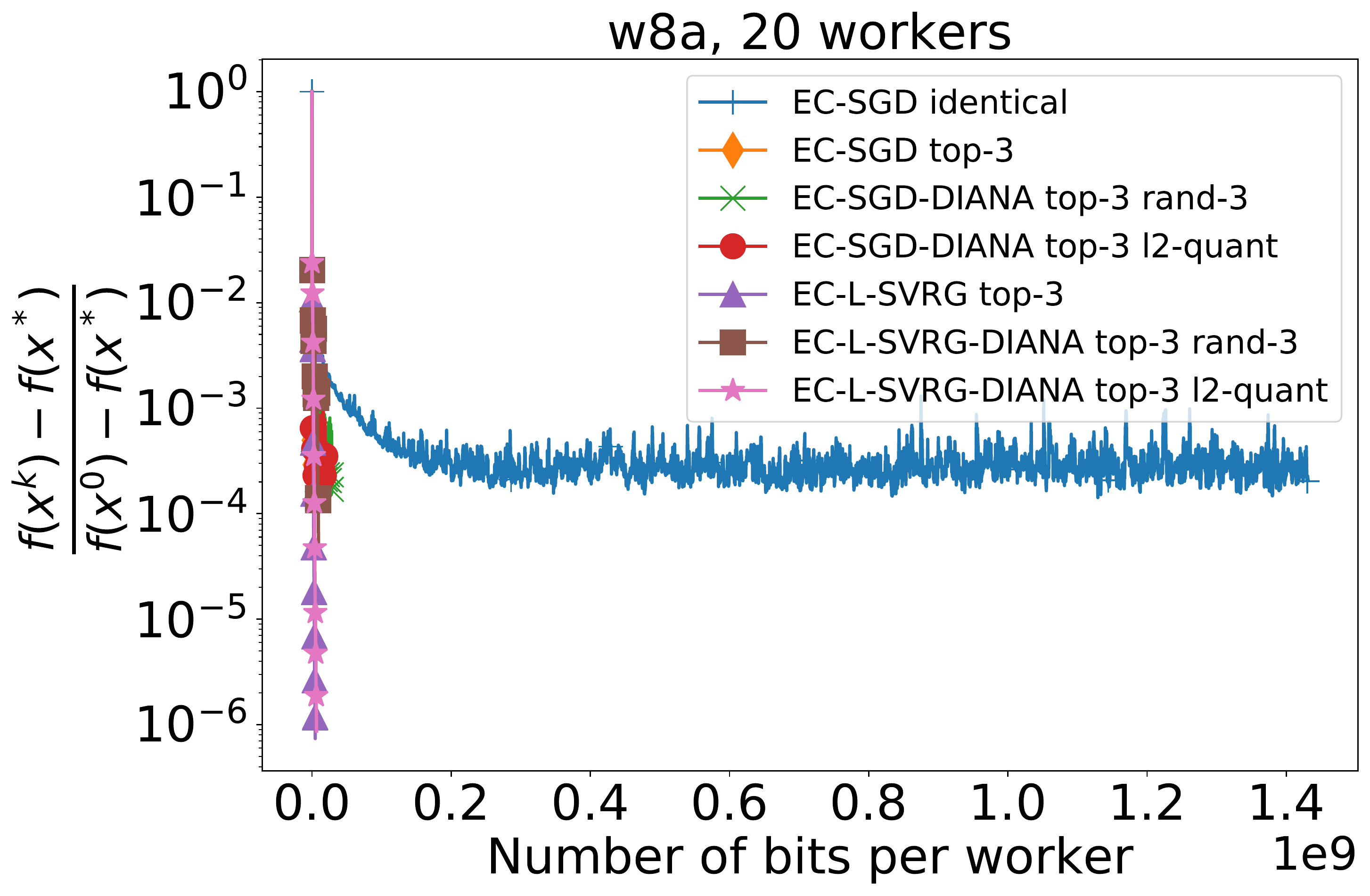}
	\includegraphics[width=0.32\textwidth]{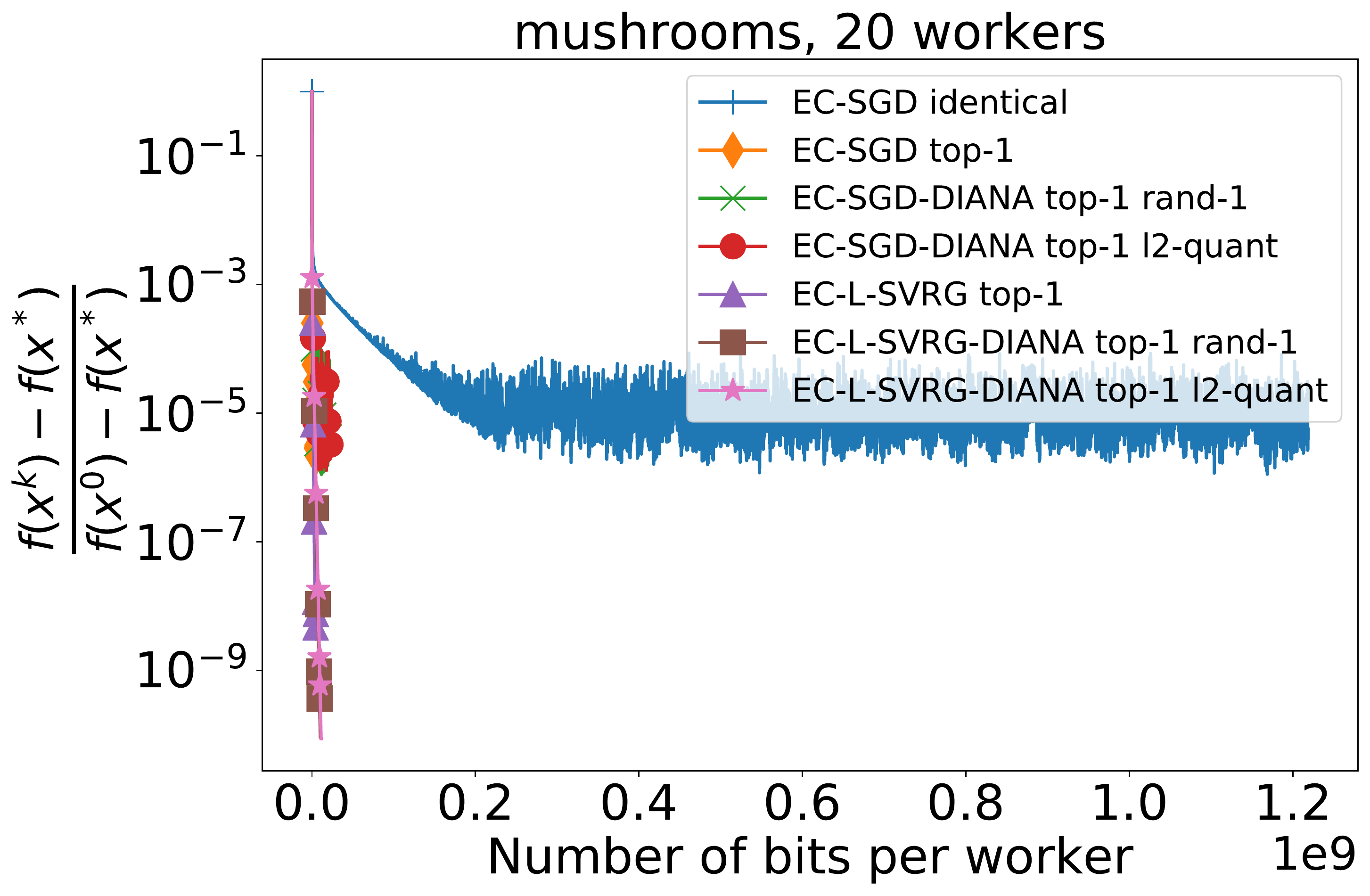}    
	\includegraphics[width=0.32\textwidth]{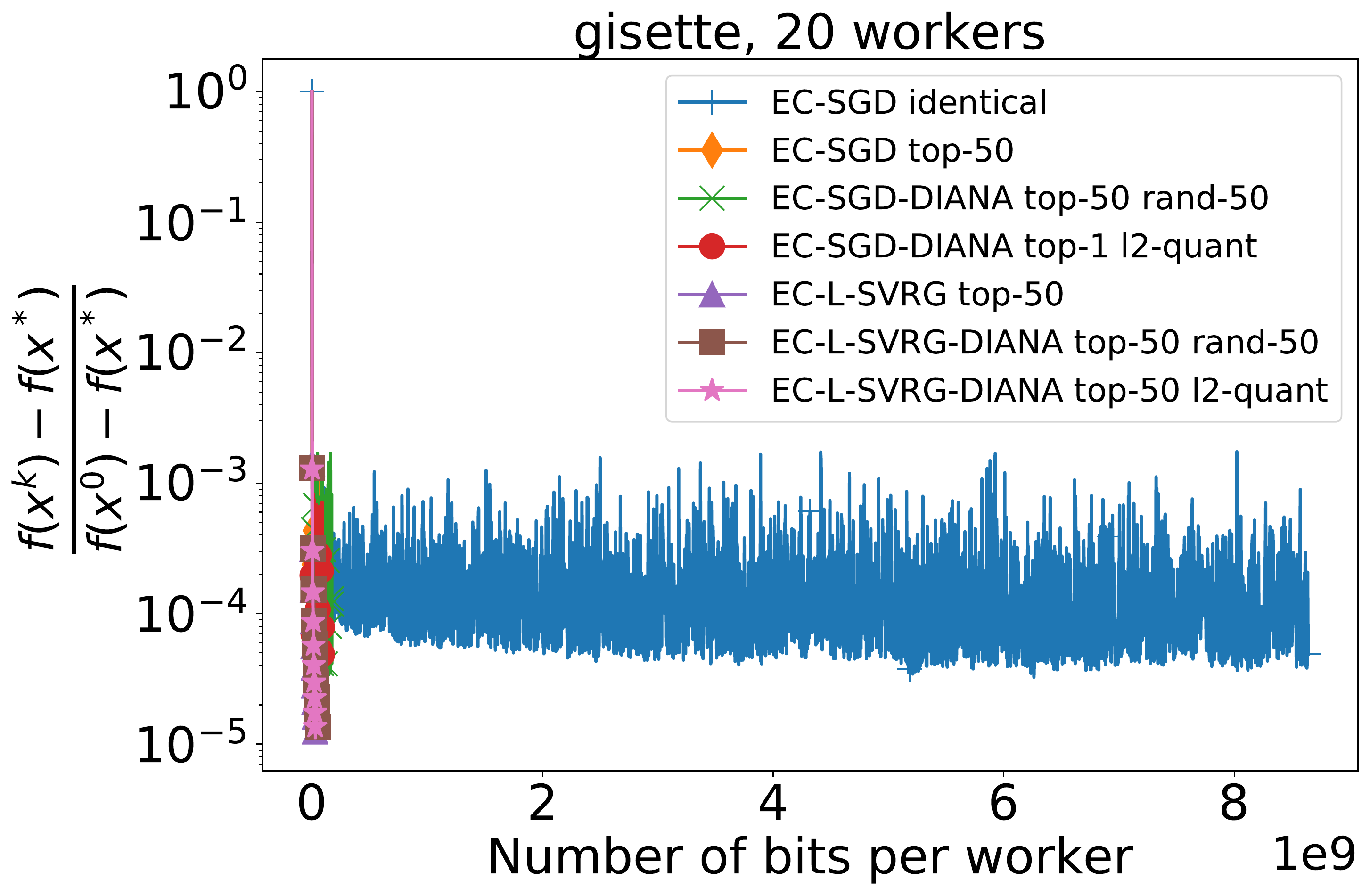}           
    \\
    \includegraphics[width=0.32\textwidth]{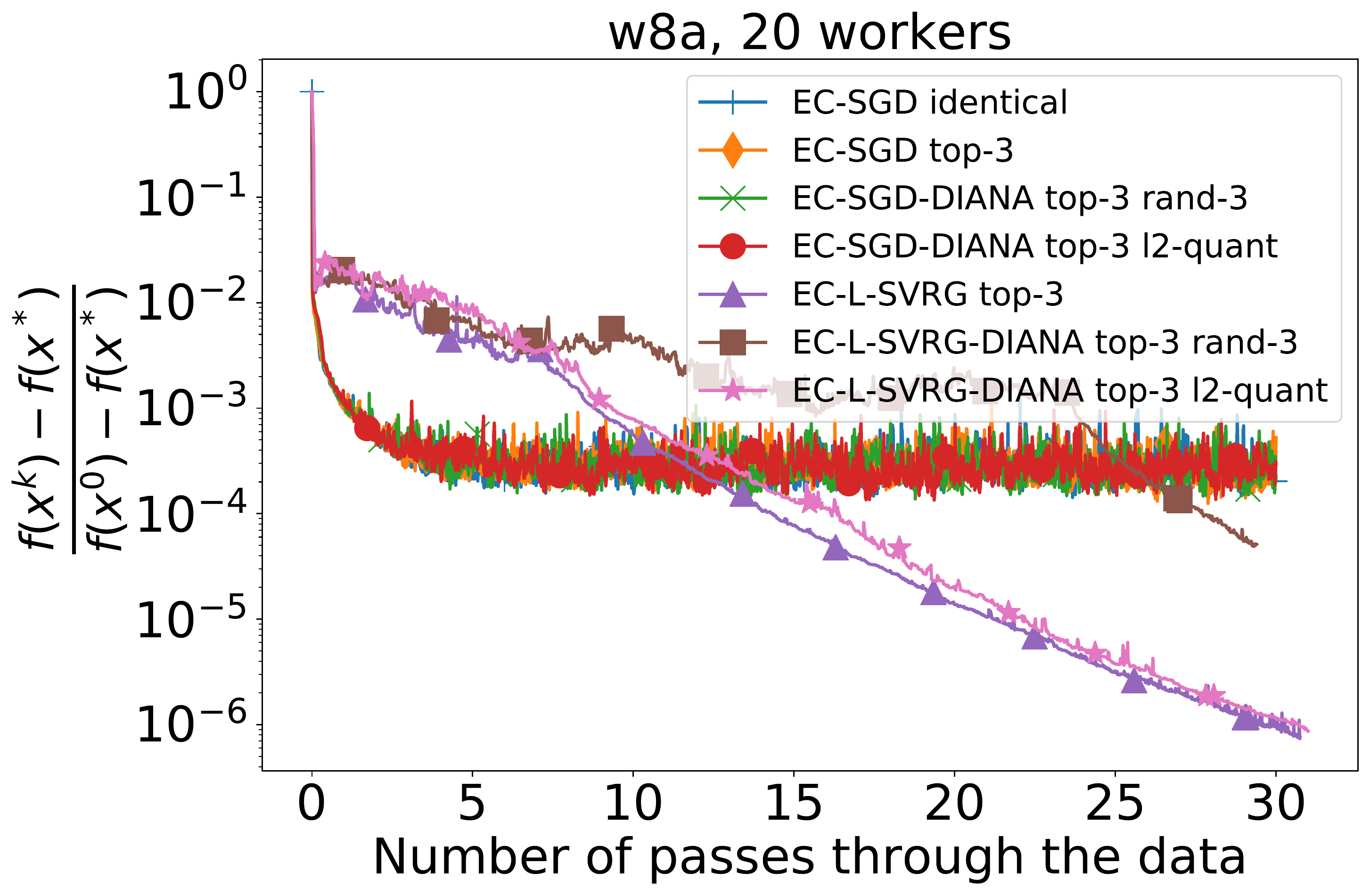}
    \includegraphics[width=0.32\textwidth]{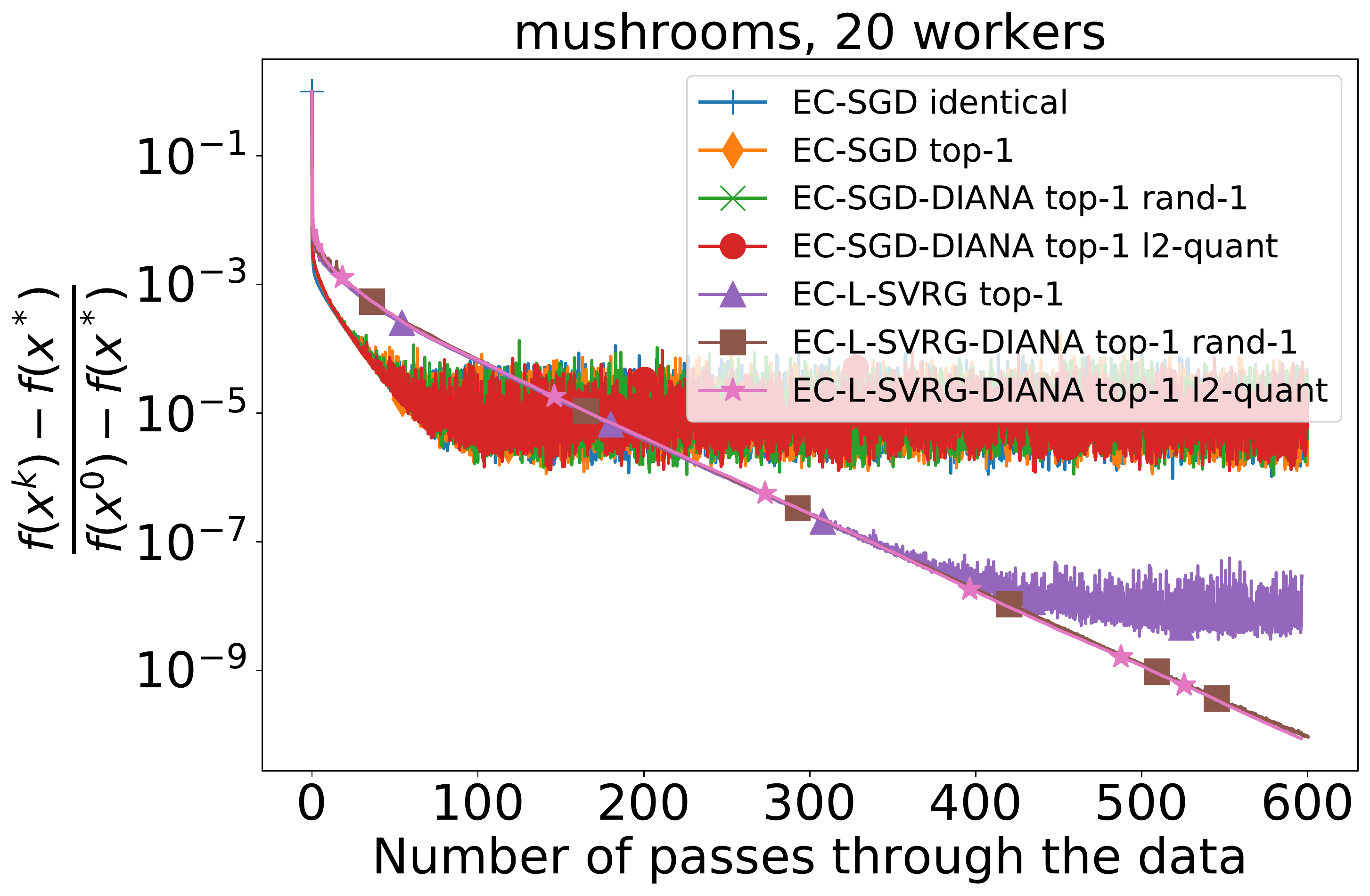}        
	\includegraphics[width=0.32\textwidth]{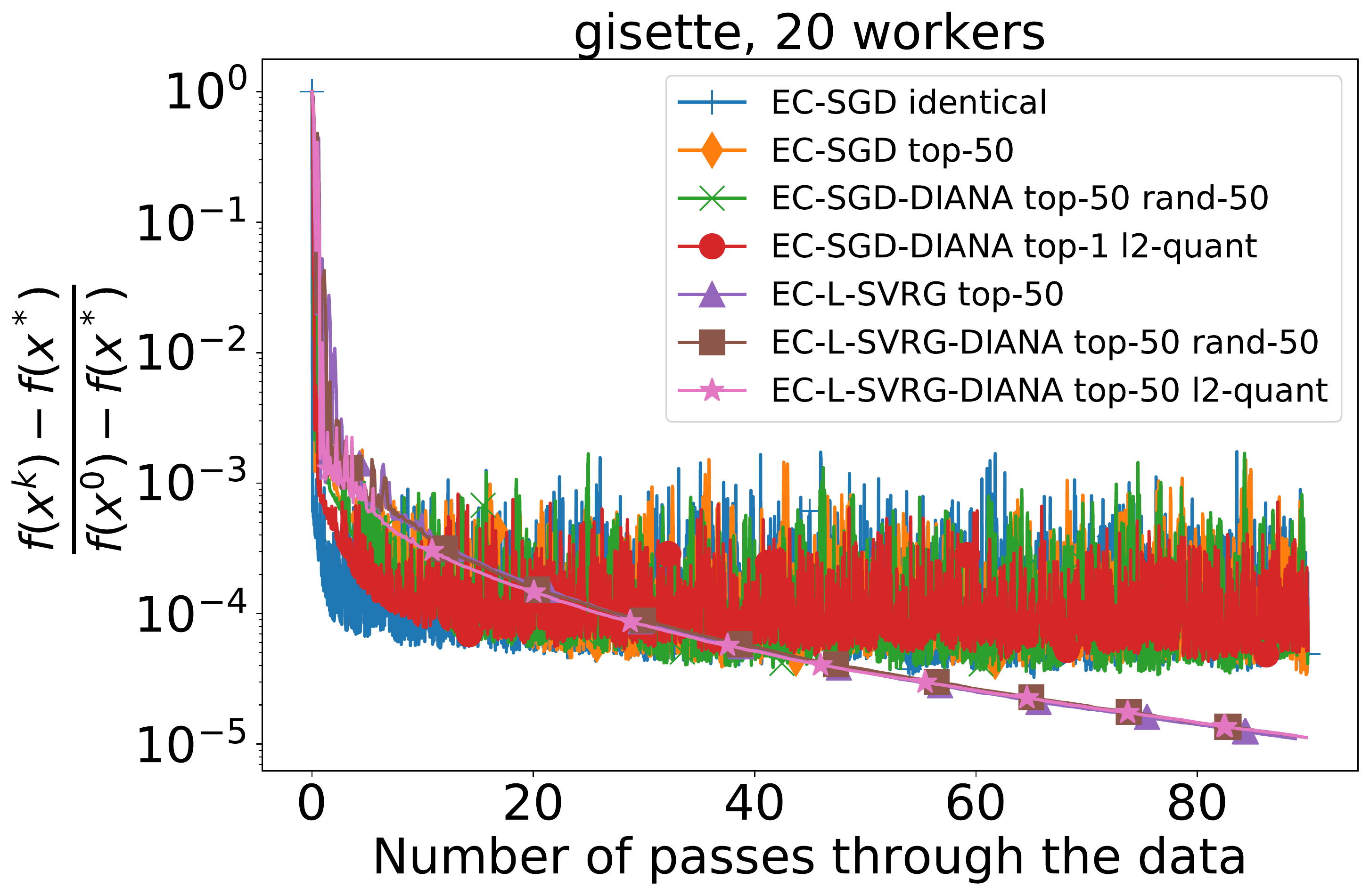}        
    \caption{Trajectories of {\tt EC-SGD}, {\tt EC-SGD-DIANA}, {\tt EC-LSVRG} and {\tt EC-LSVRG-DIANA} applied to solving logistic regression problem with $20$ workers. {\tt EC-SGD identical} corresponds to {\tt SGD} with error compensation with the identity compression operator $\cC(x) = x$, i.e., it is just parallel {\tt SGD}.}
    \label{fig:sgd_logreg_extra}
\end{figure}

\subsection{Compressing Full Gradients}

\begin{figure}[H]
    \centering
	\includegraphics[width=0.32\textwidth]{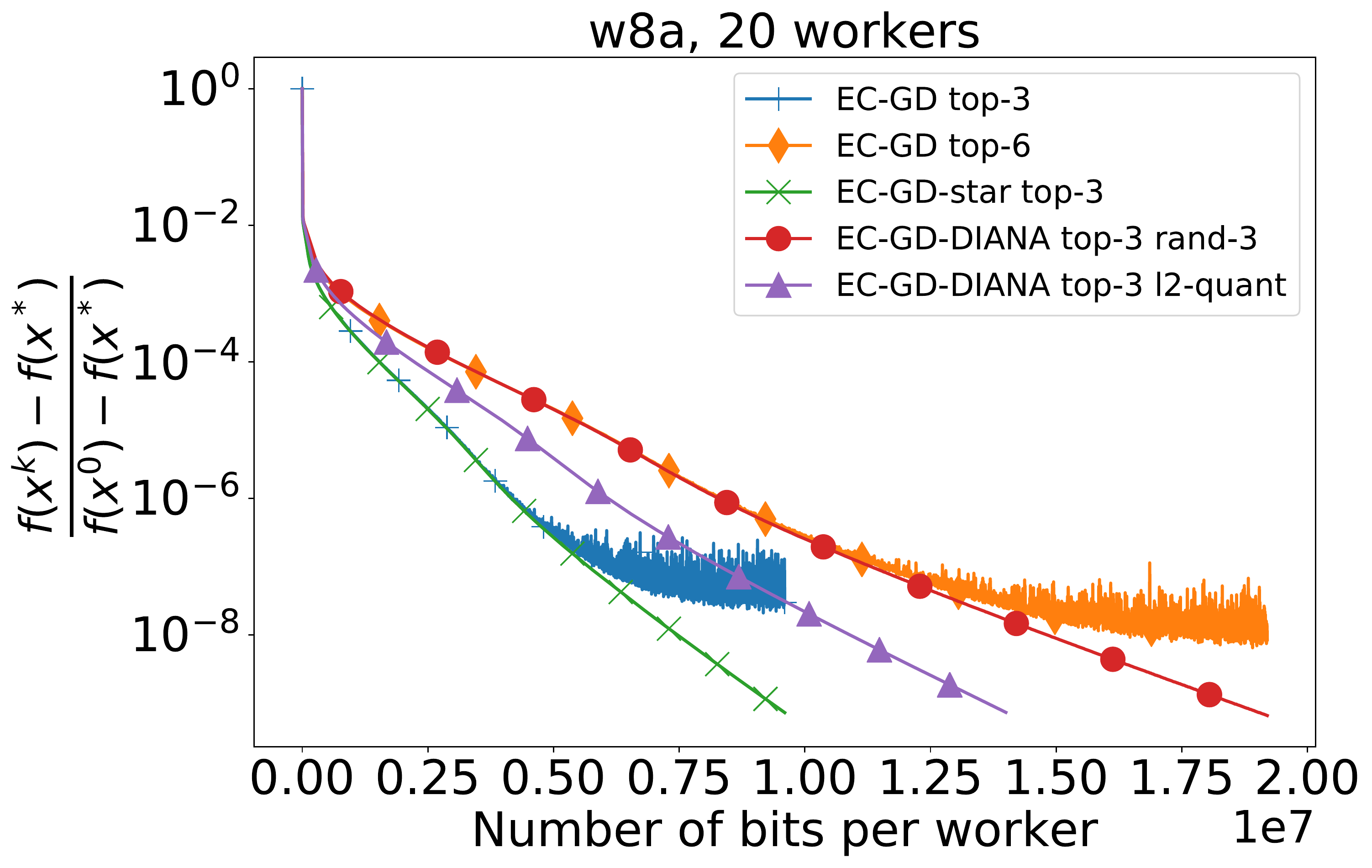}	\includegraphics[width=0.32\textwidth]{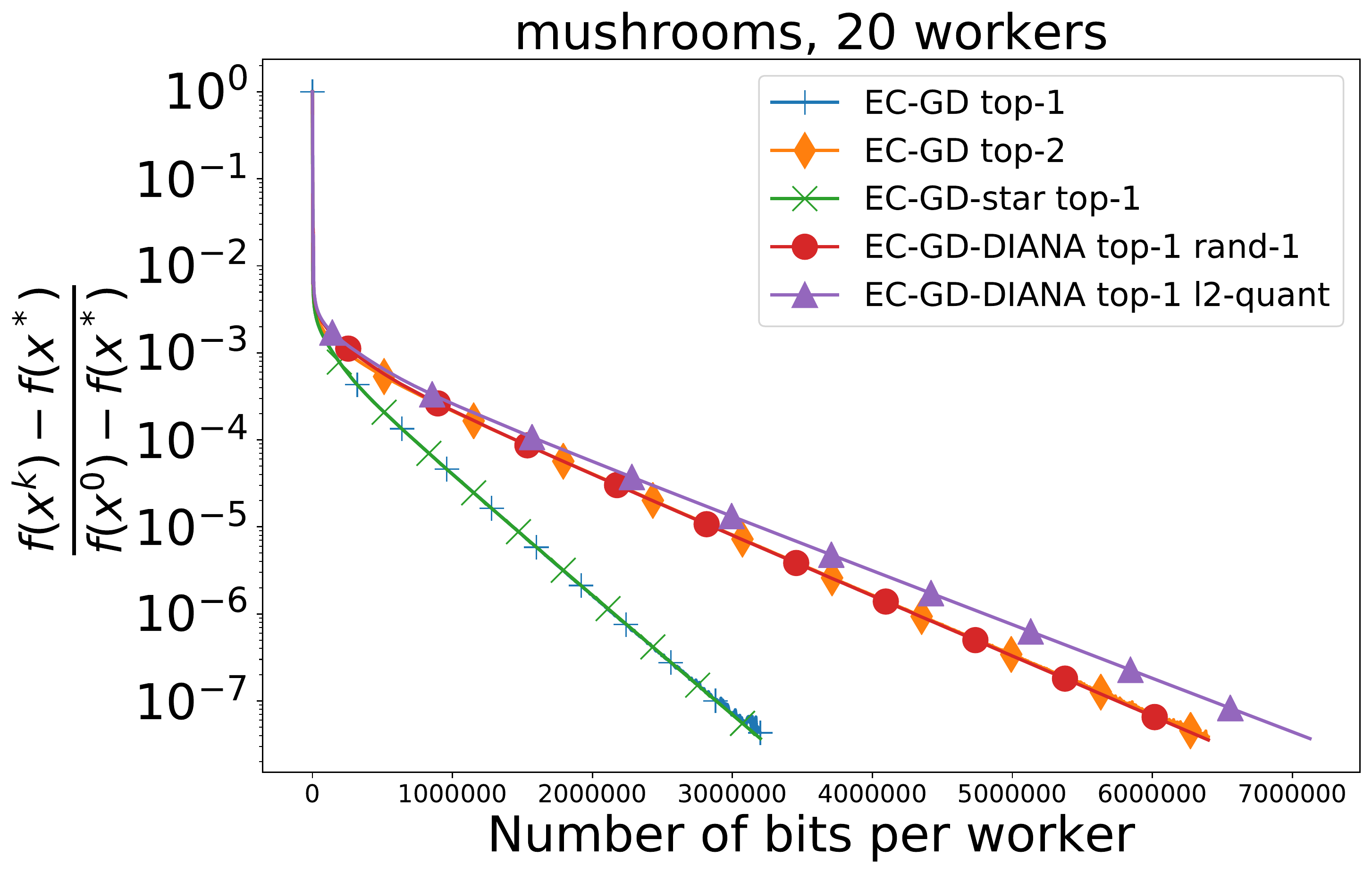}
	\includegraphics[width=0.32\textwidth]{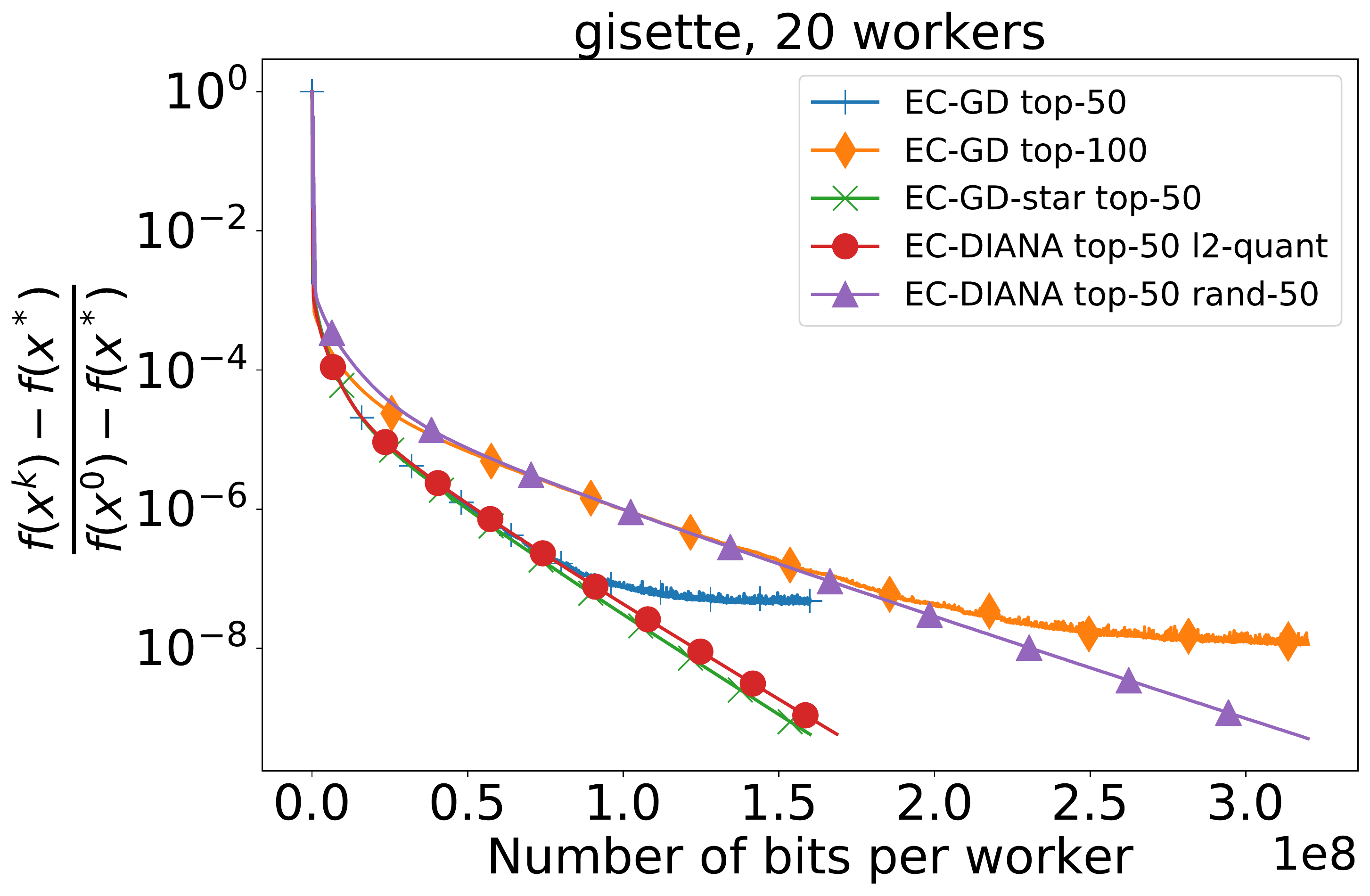}
	\\
	\includegraphics[width=0.32\textwidth]{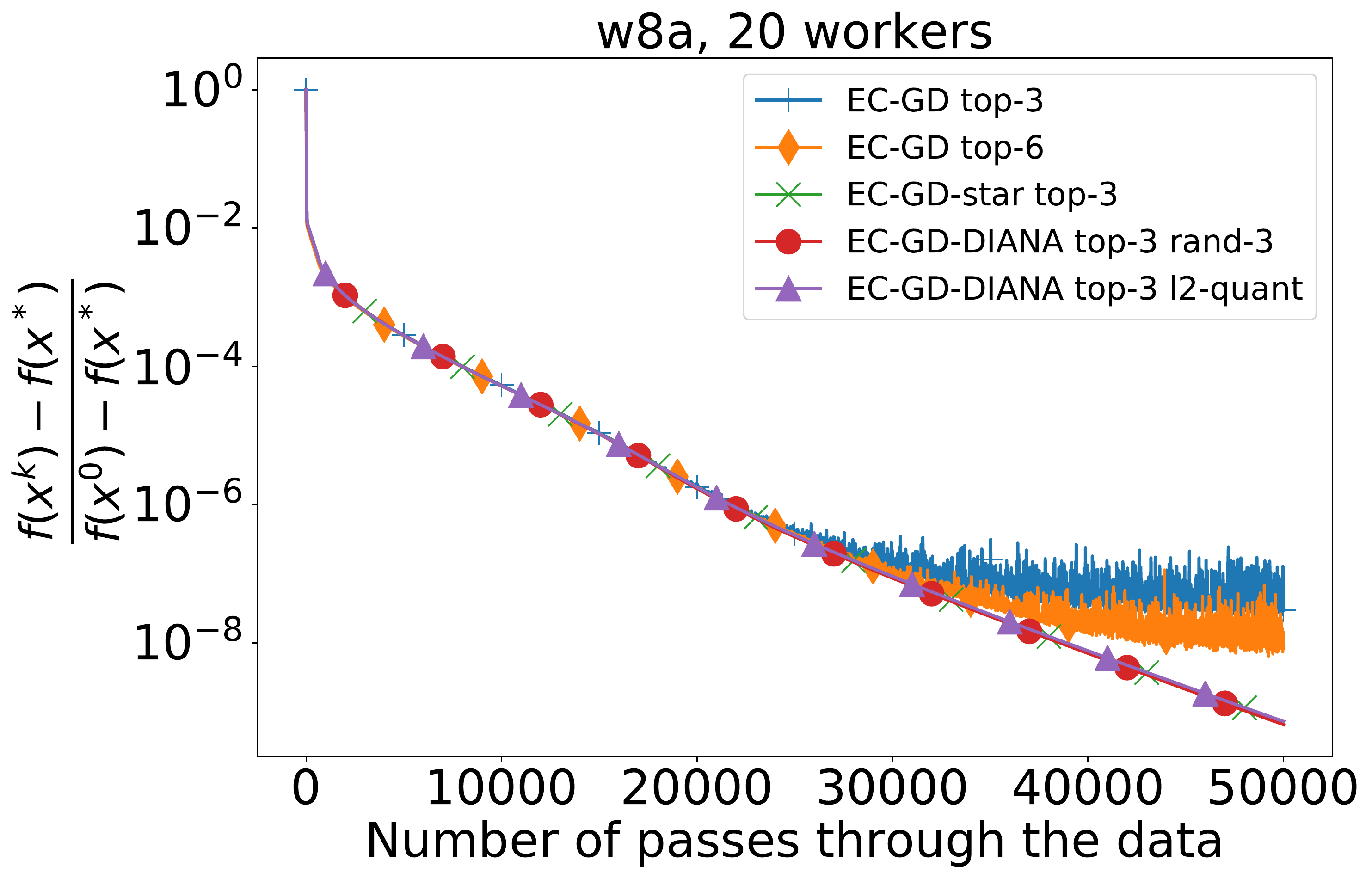}
	\includegraphics[width=0.32\textwidth]{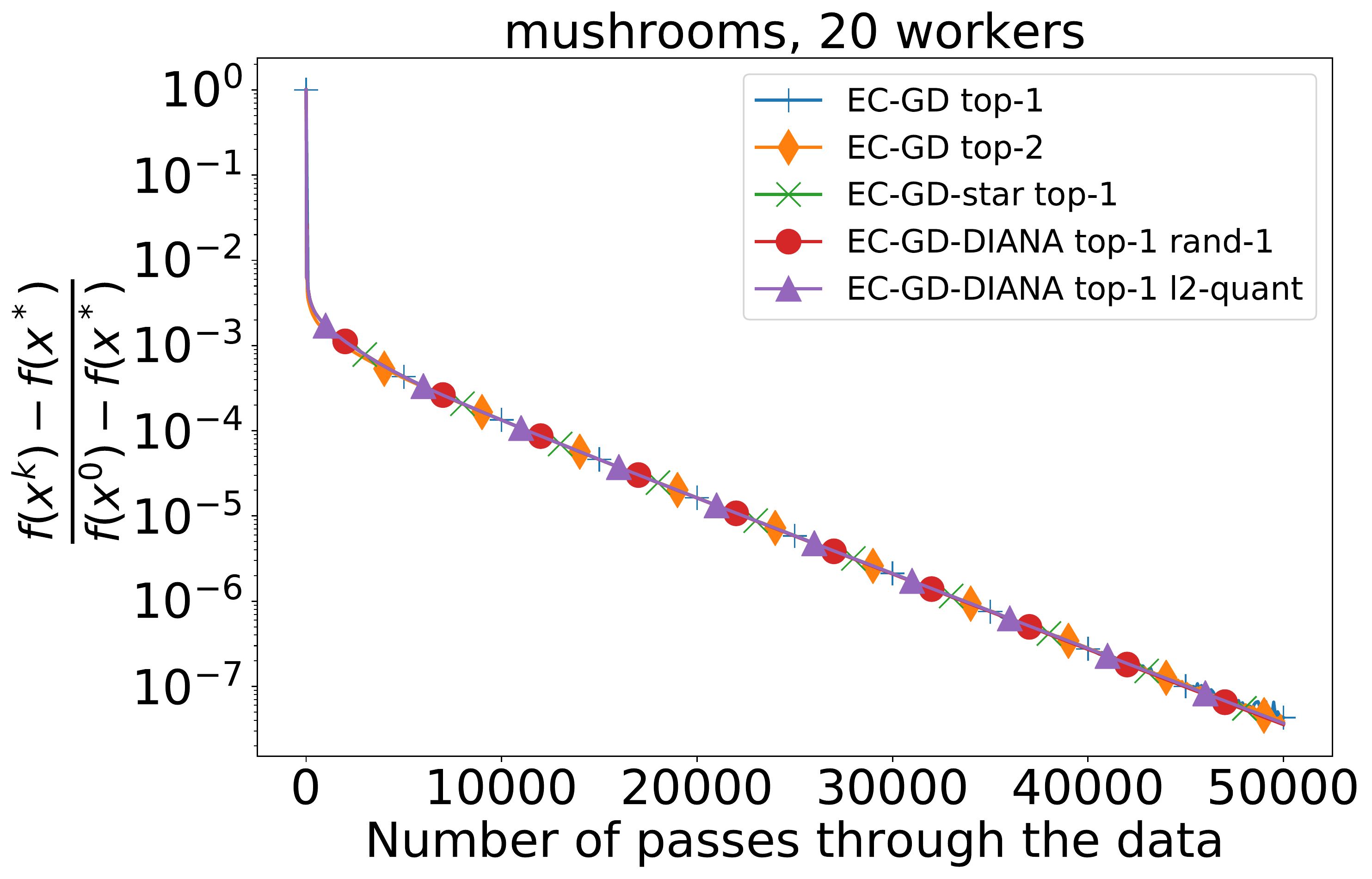}	\includegraphics[width=0.32\textwidth]{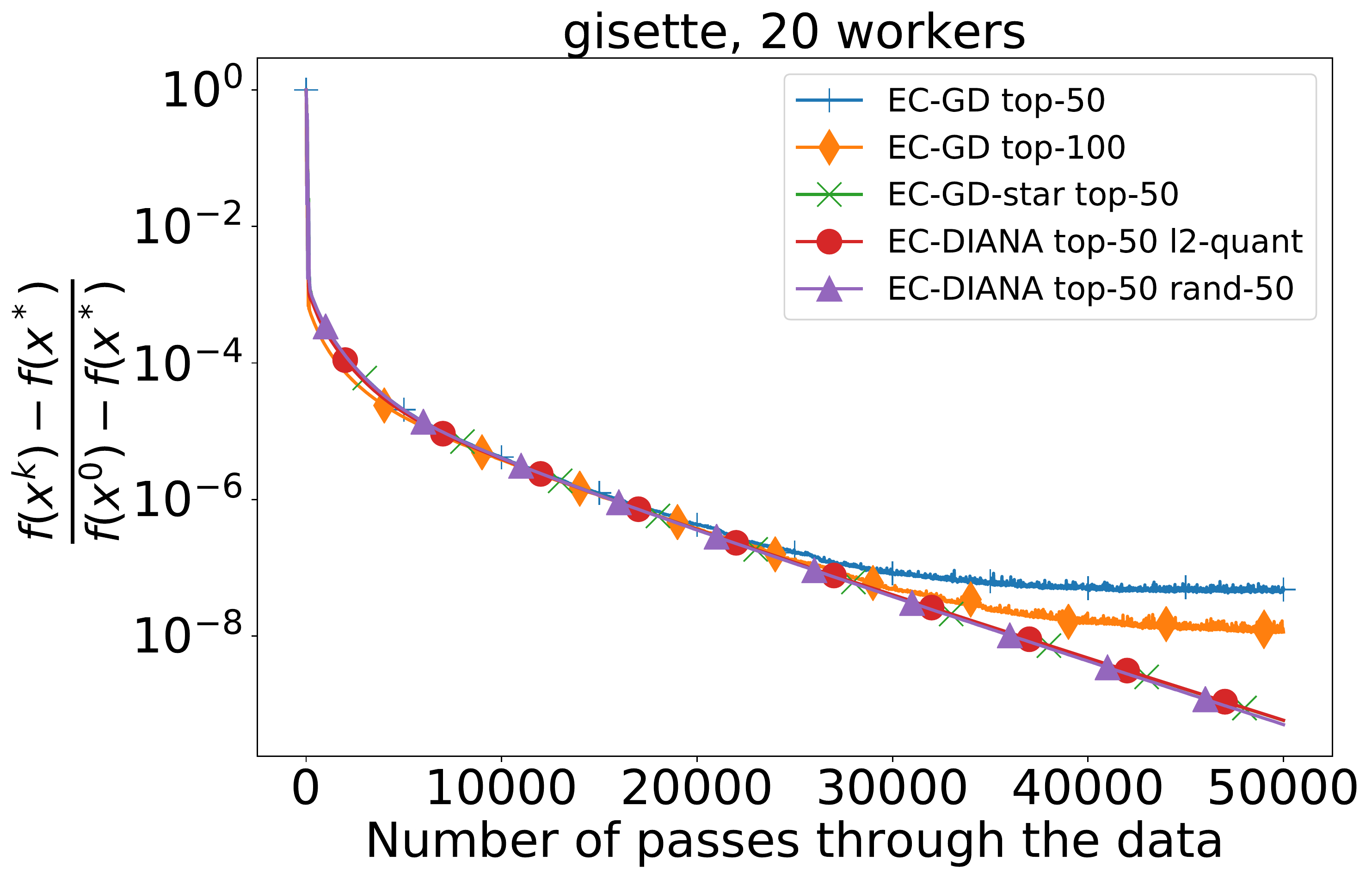}
	\caption{Trajectories of {\tt EC-GD}, {\tt EC-GD-star} and {\tt EC-DIANA} applied to solving logistic regression problem with $20$ workers.}
    \label{fig:gd_logreg_20_workers_appendix}
\end{figure}
\begin{figure}[H]
    \centering
    \includegraphics[width=0.32\textwidth]{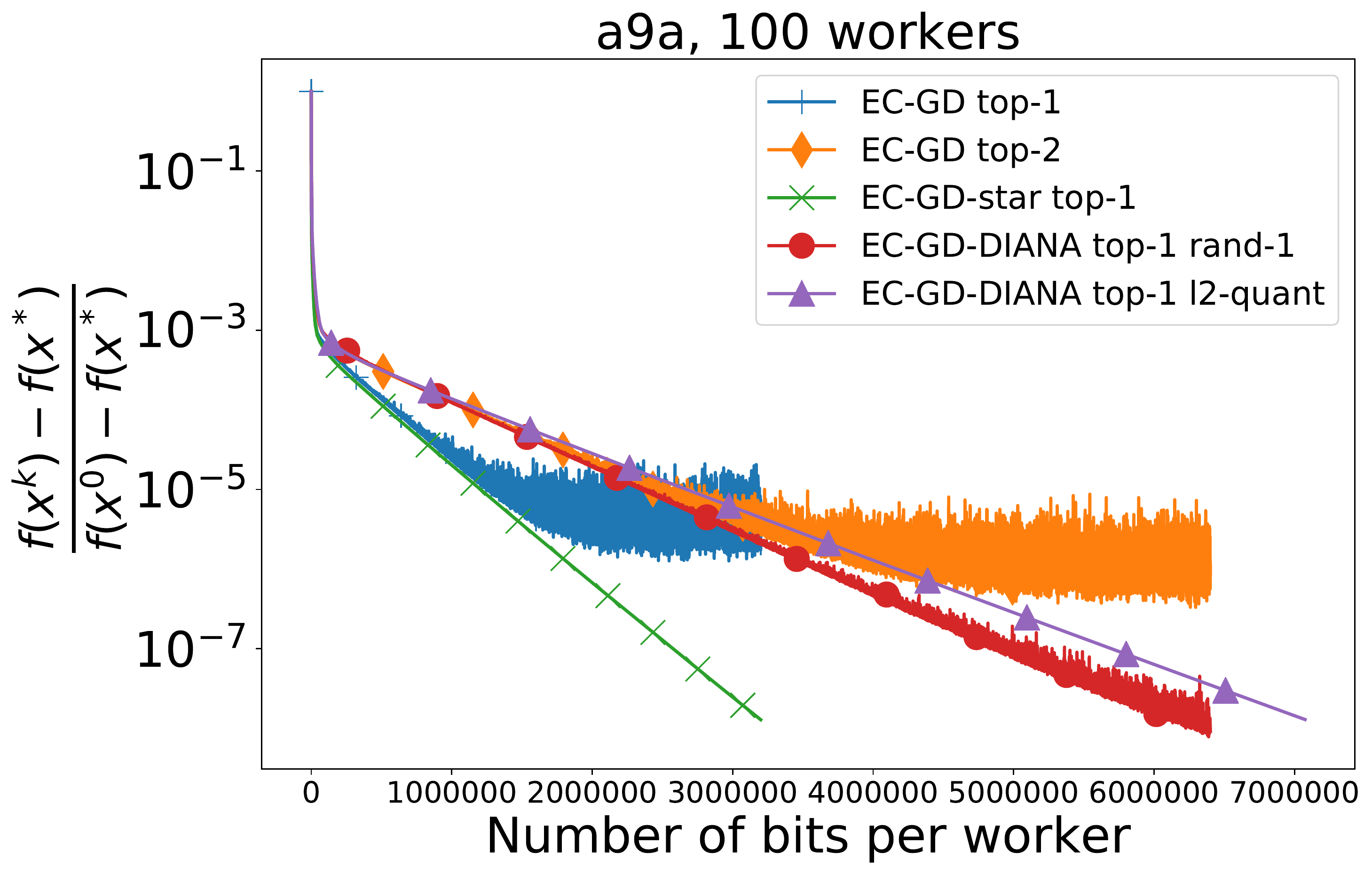}
	\includegraphics[width=0.32\textwidth]{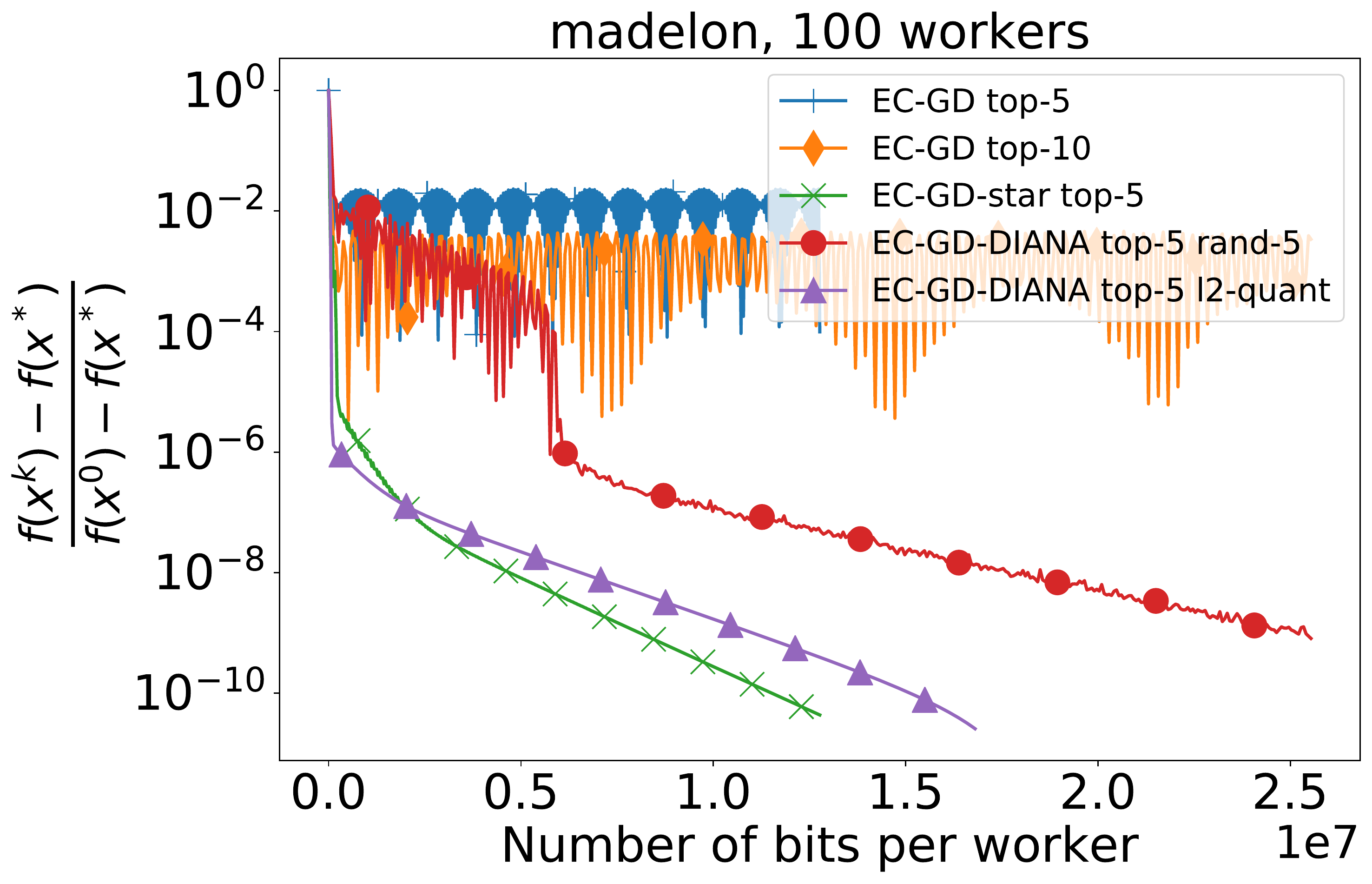}    
	\includegraphics[width=0.32\textwidth]{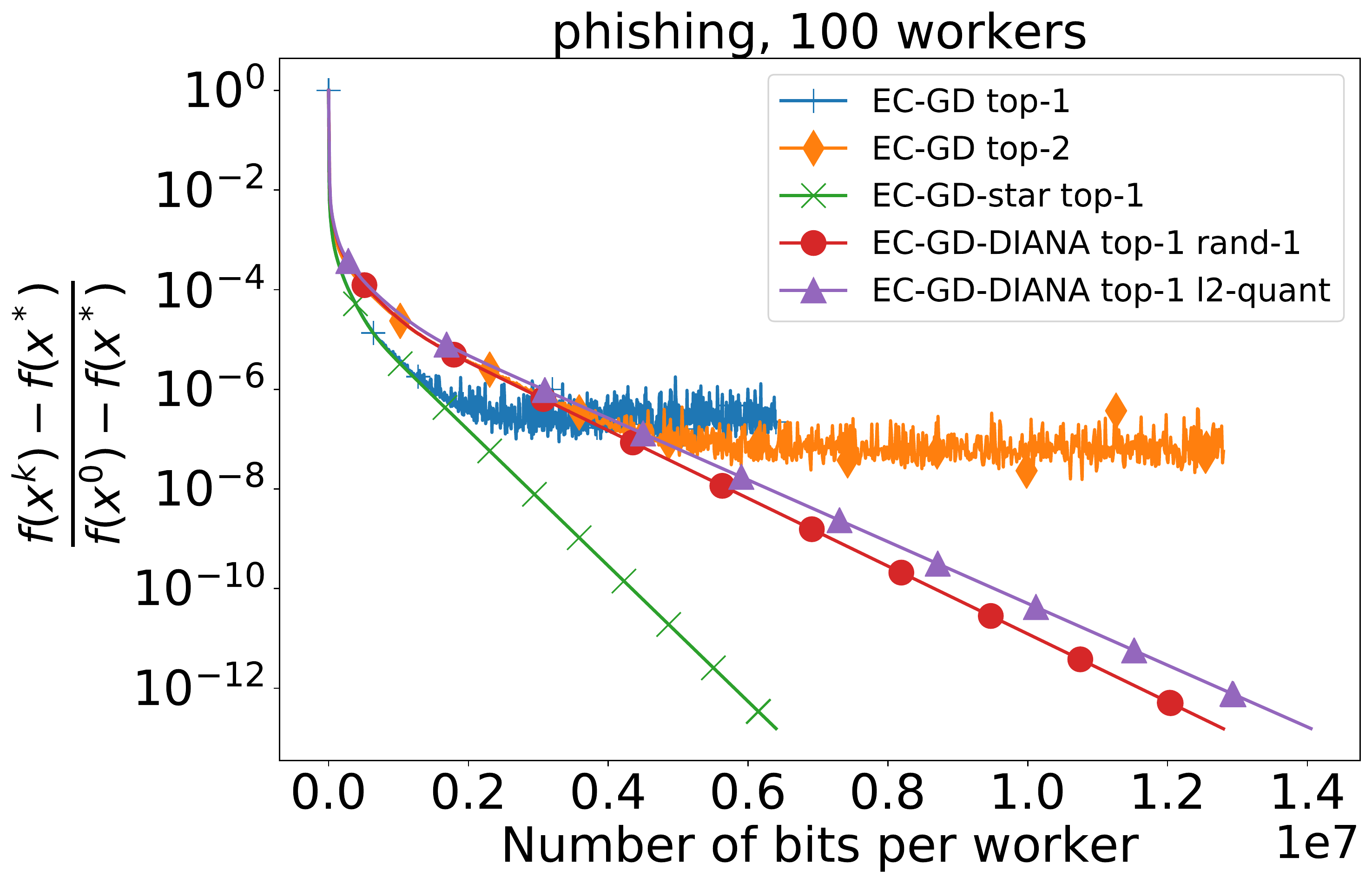}    
    \\
    \includegraphics[width=0.32\textwidth]{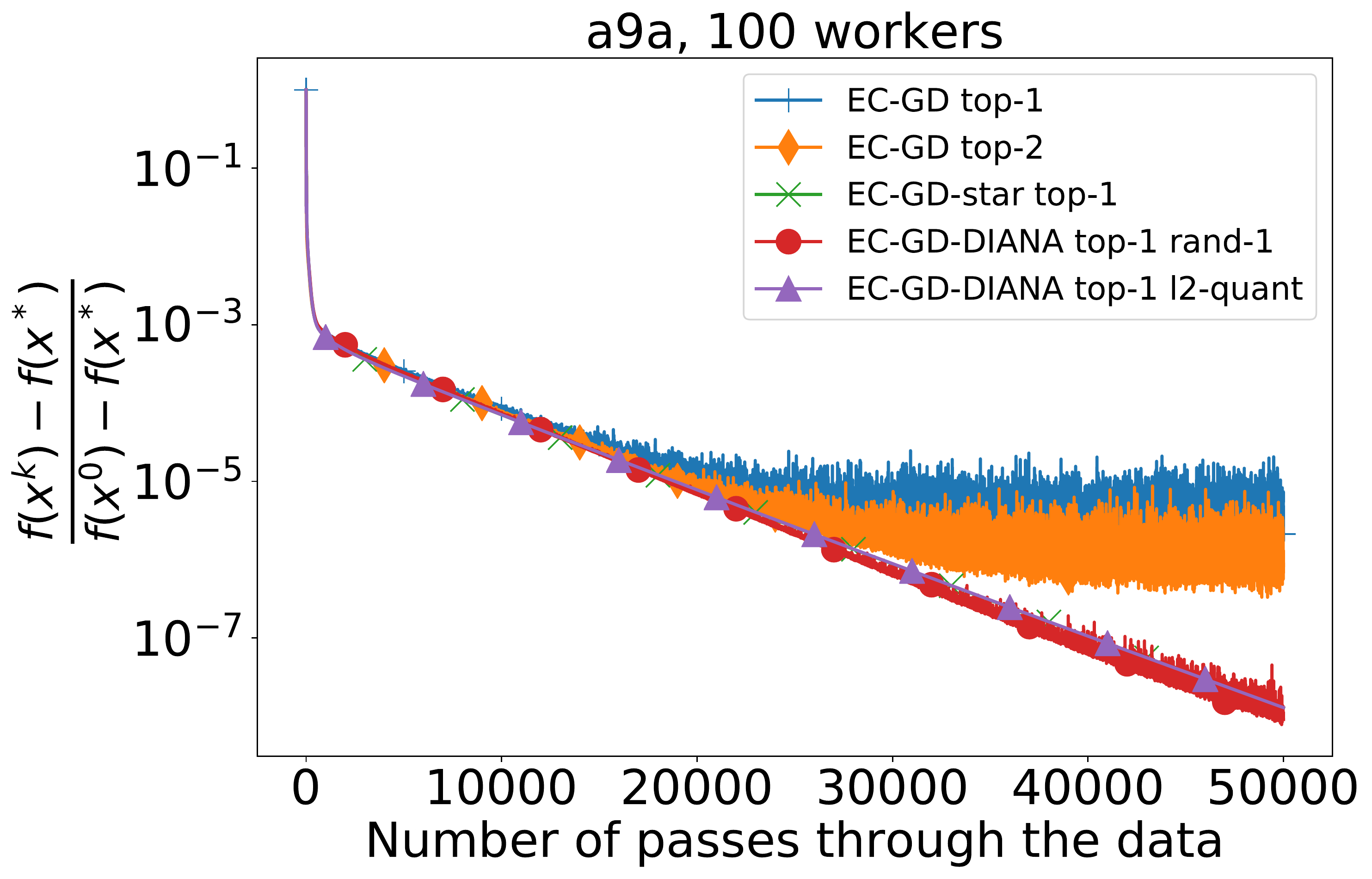}    
	\includegraphics[width=0.32\textwidth]{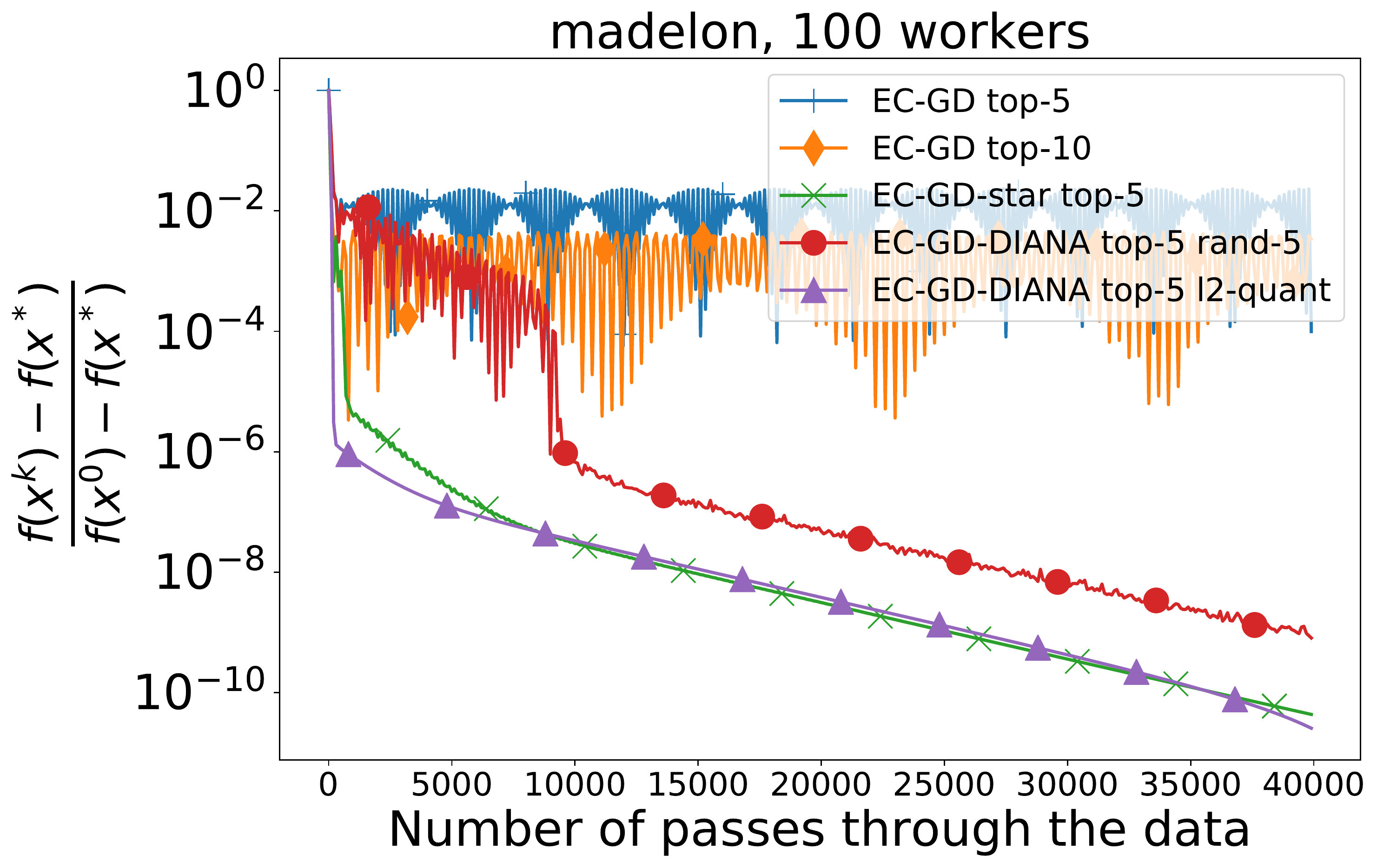}    
	\includegraphics[width=0.32\textwidth]{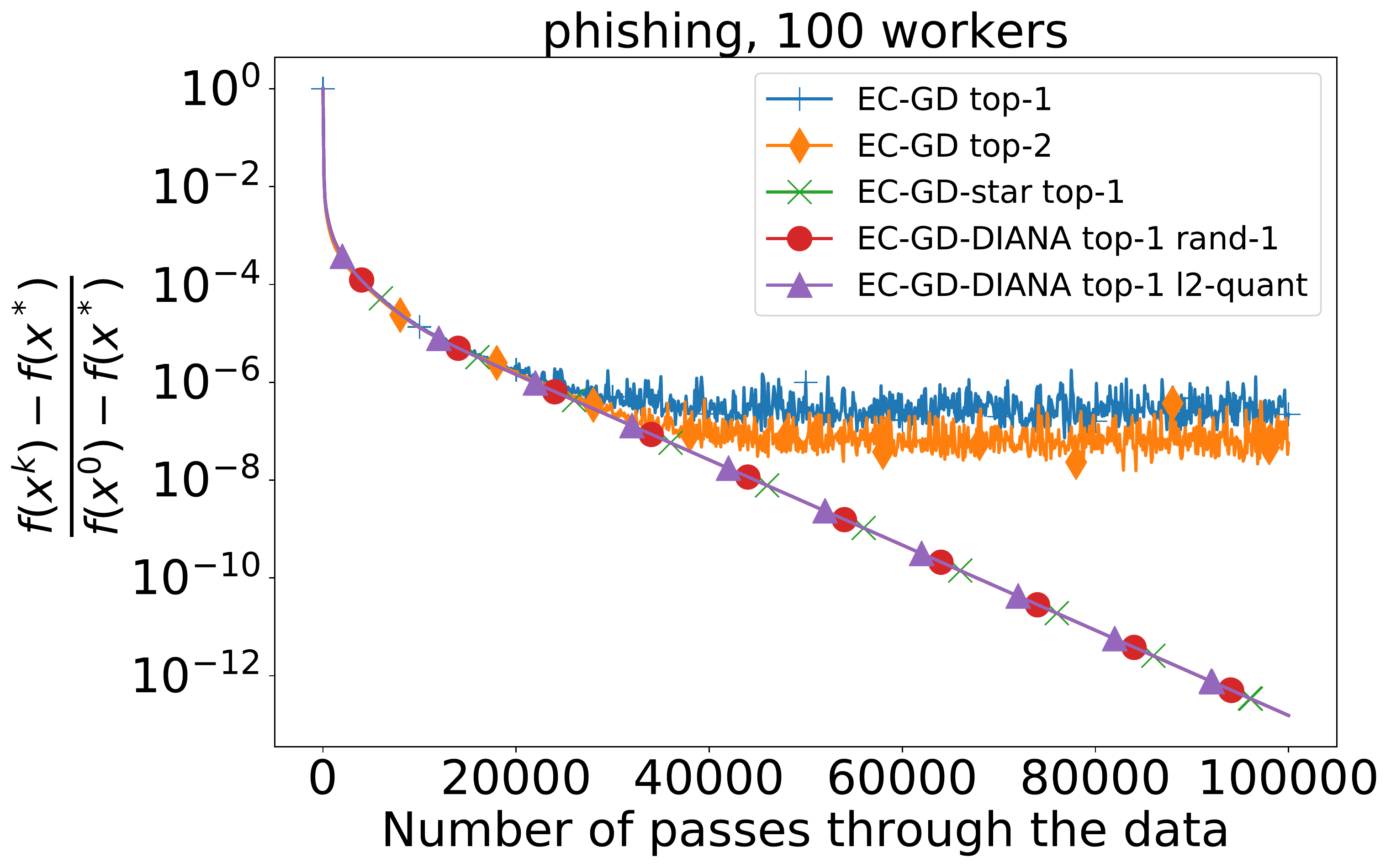}    	
	\\
	\includegraphics[width=0.32\textwidth]{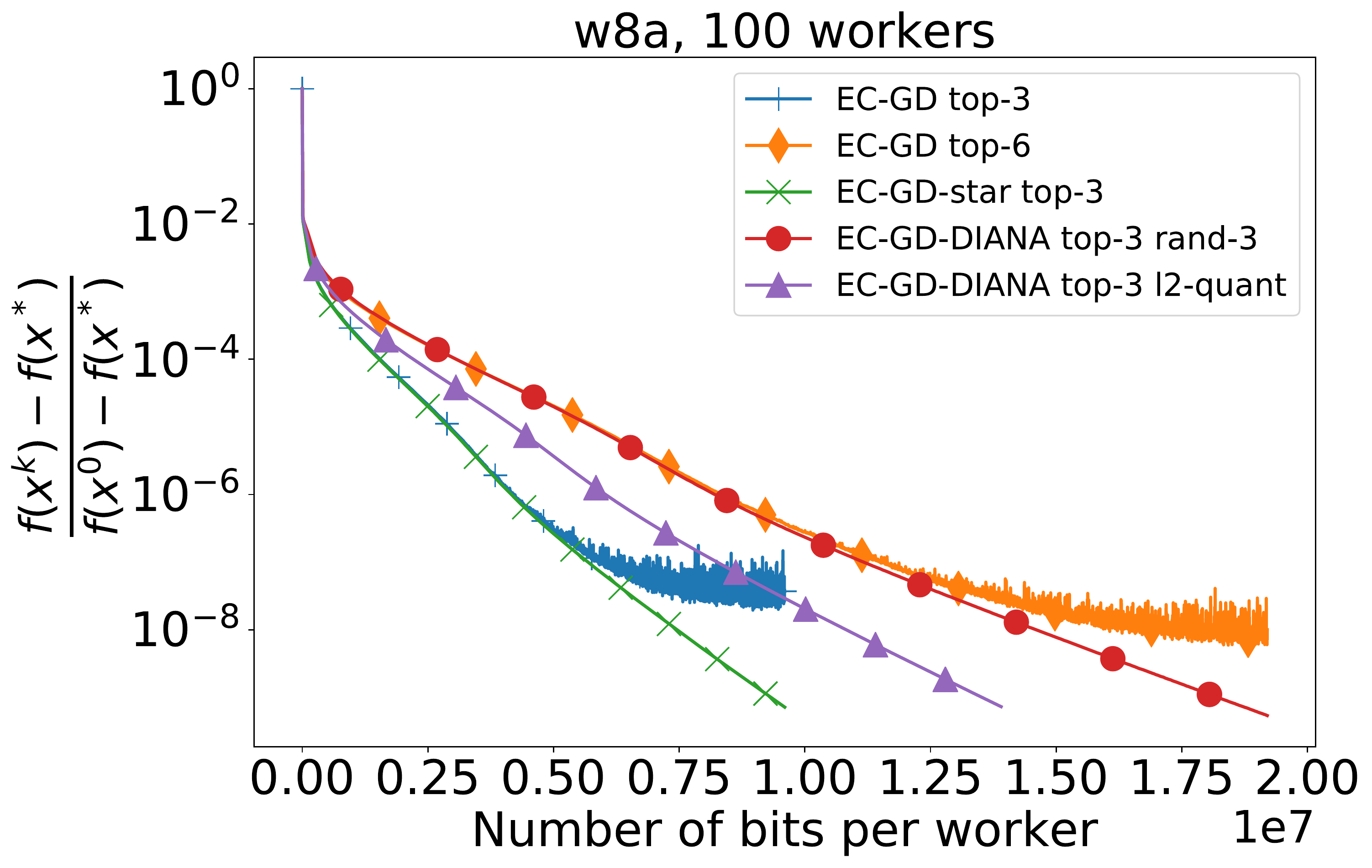}	\includegraphics[width=0.32\textwidth]{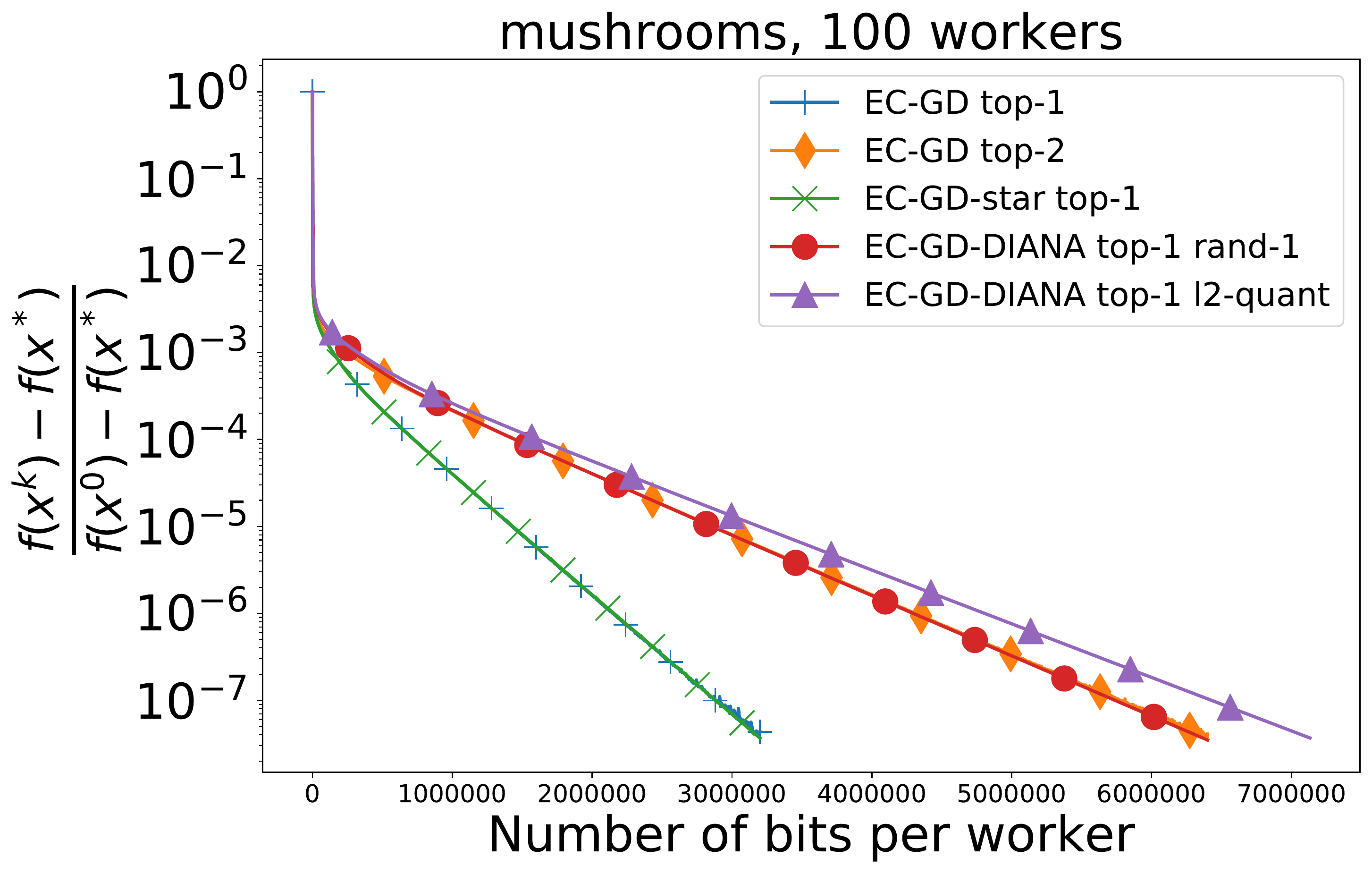}
	\includegraphics[width=0.32\textwidth]{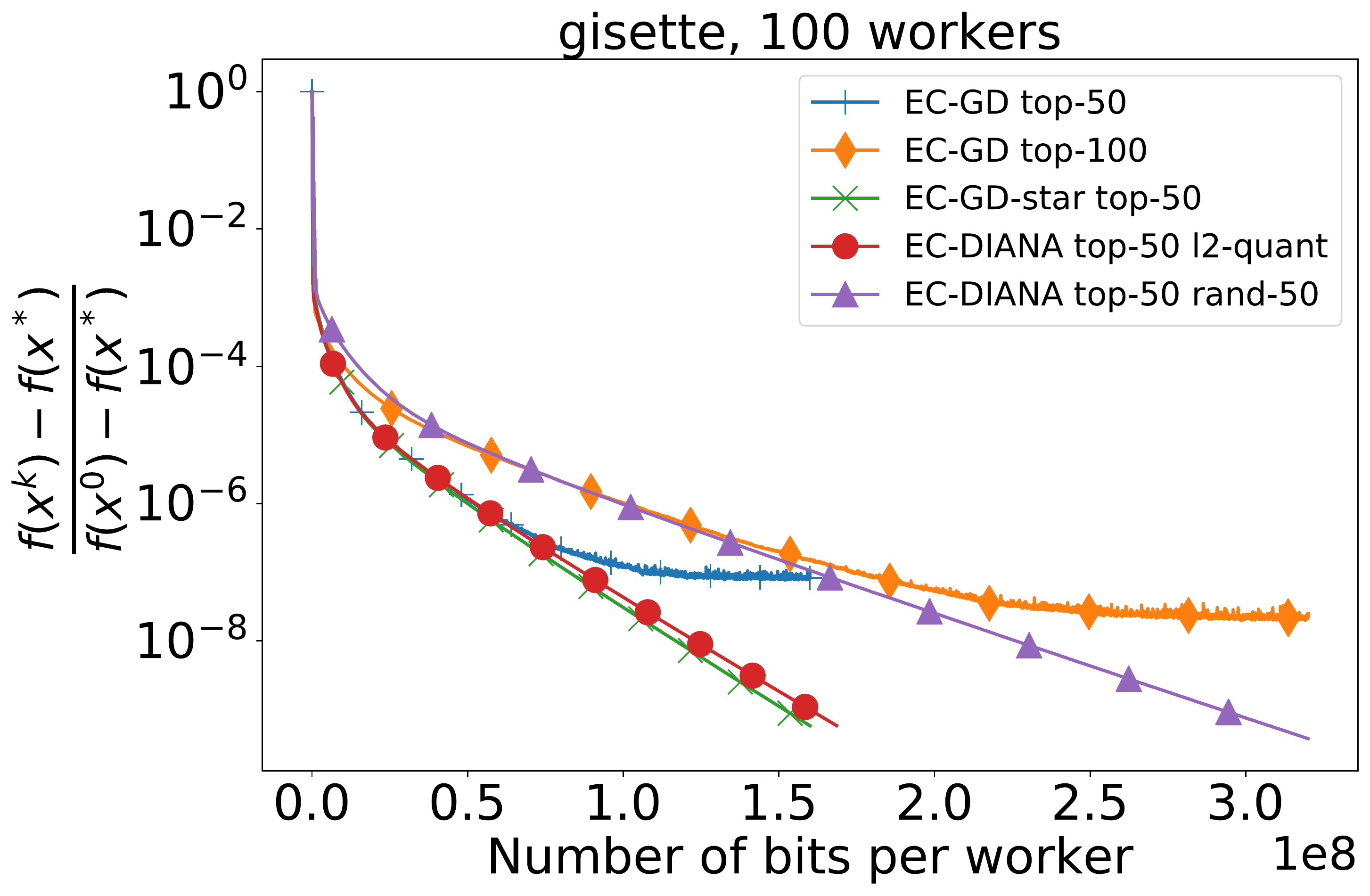}
	\\
	\includegraphics[width=0.32\textwidth]{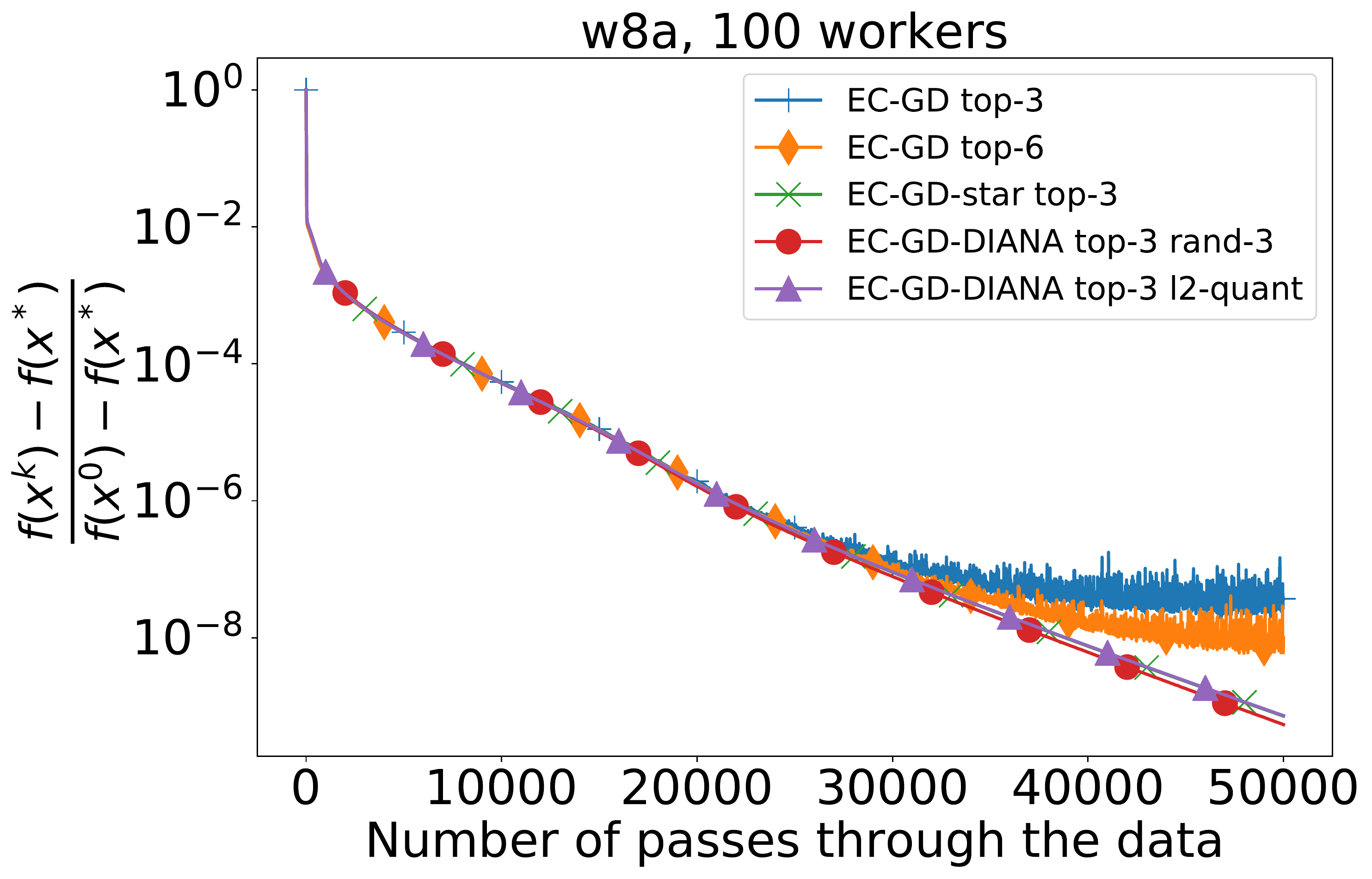}
	\includegraphics[width=0.32\textwidth]{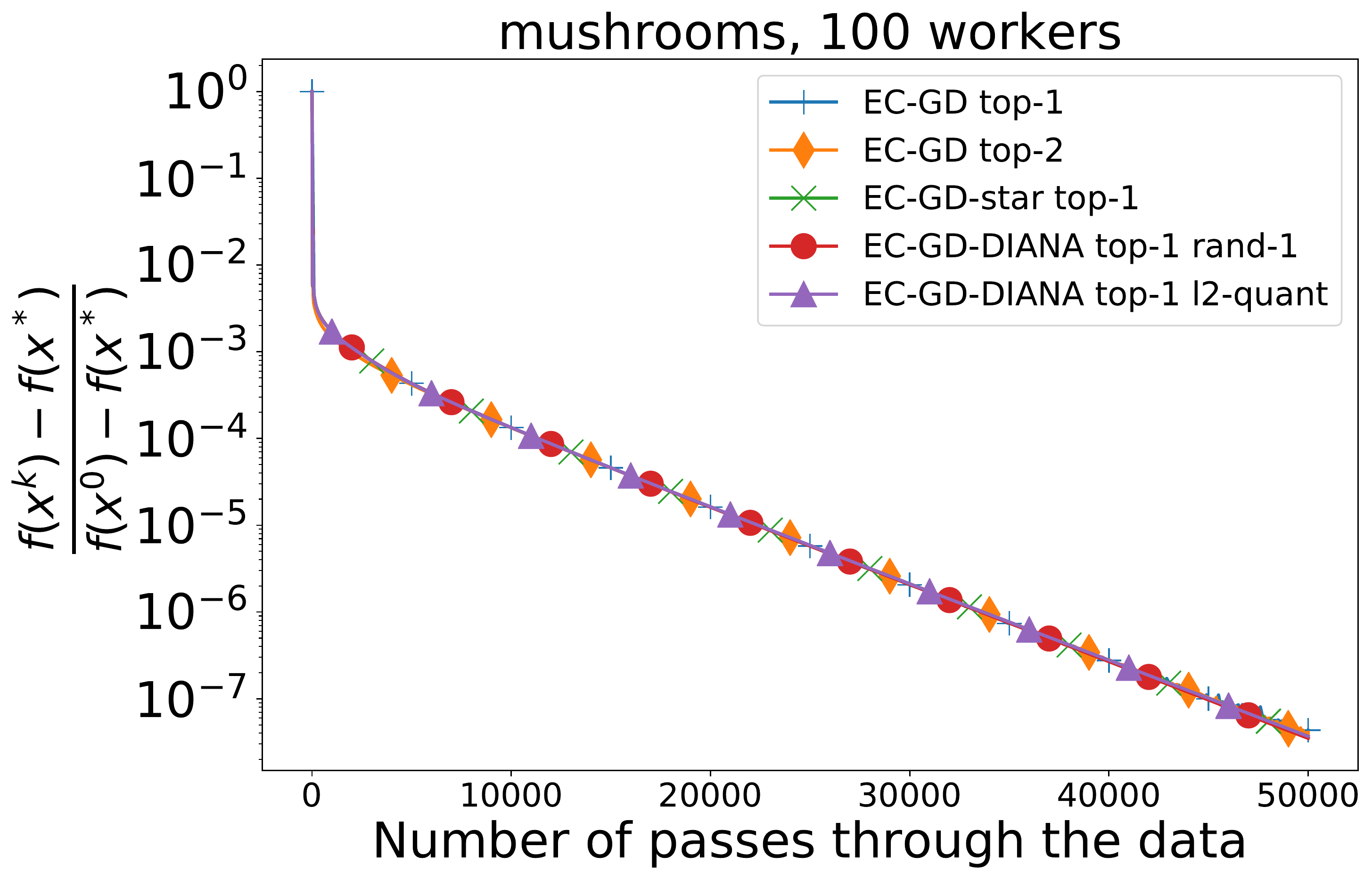}
	\includegraphics[width=0.32\textwidth]{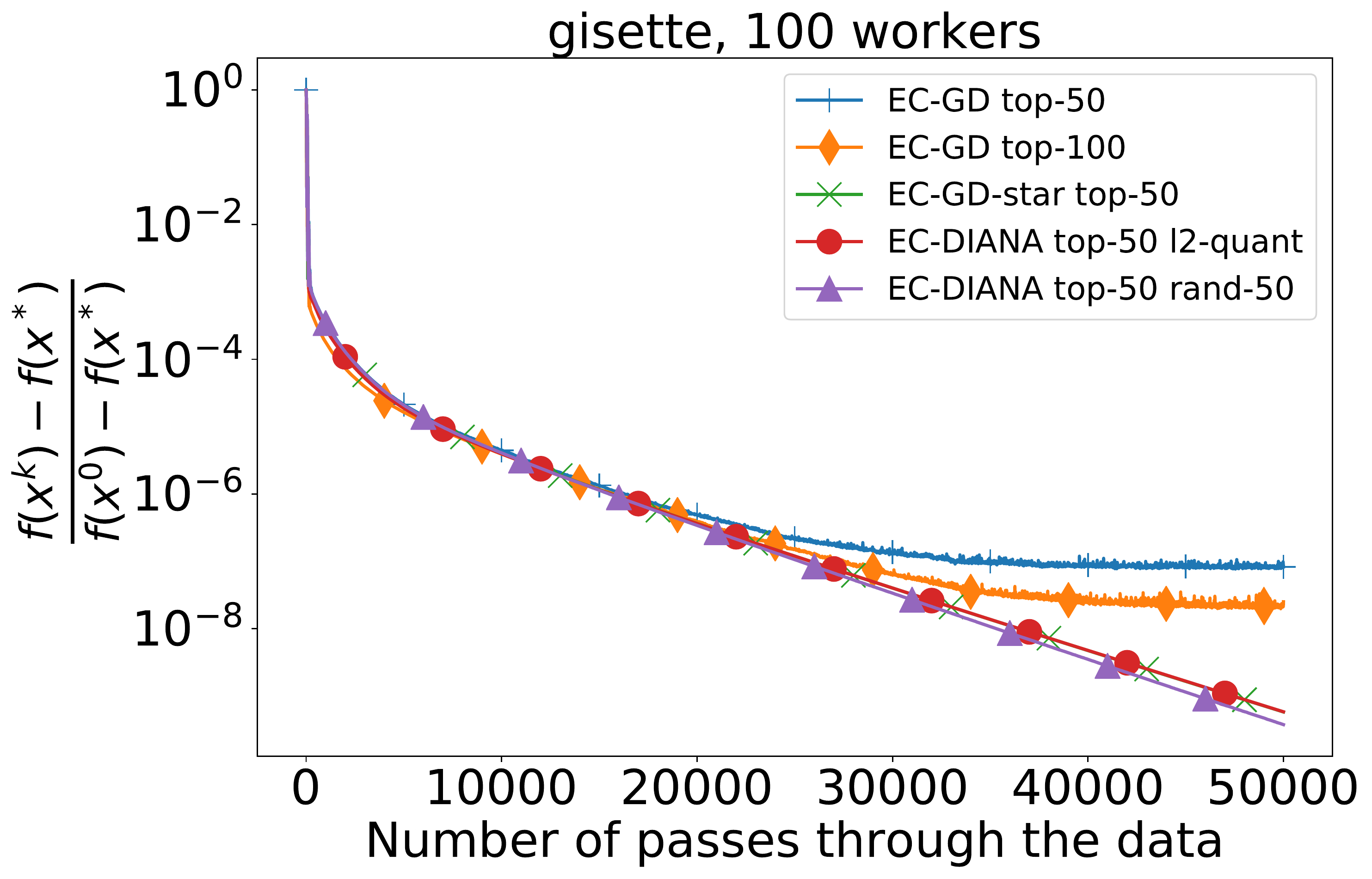}	    
	    \caption{Trajectories of {\tt EC-GD}, {\tt EC-GD-star} and {\tt EC-DIANA} applied to solving logistic regression problem with $100$ workers.}
    \label{fig:gd_logreg_100_workers}
\end{figure}

\begin{figure}[H]
    \centering
    \includegraphics[width=0.32\textwidth]{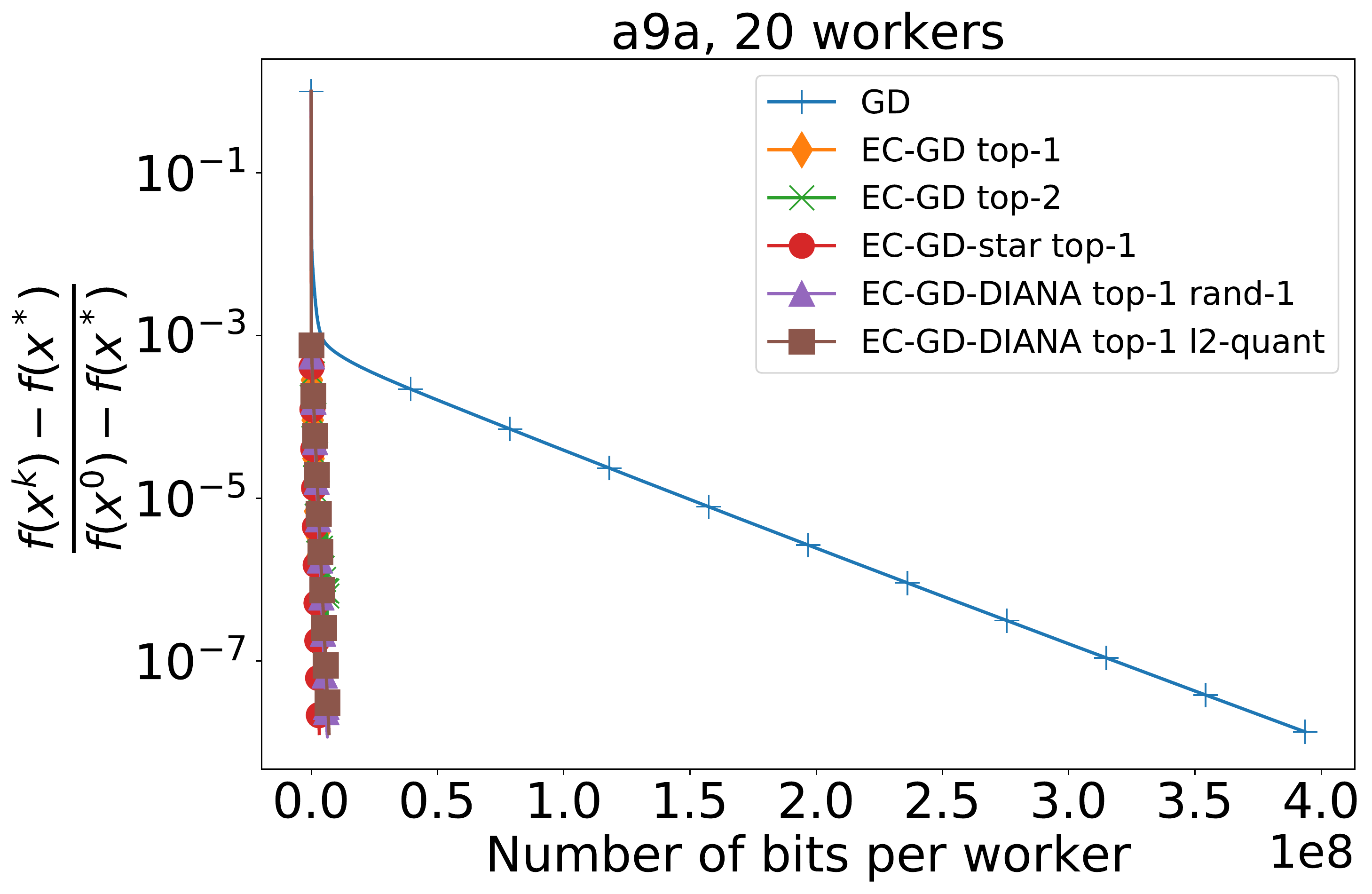}
	\includegraphics[width=0.32\textwidth]{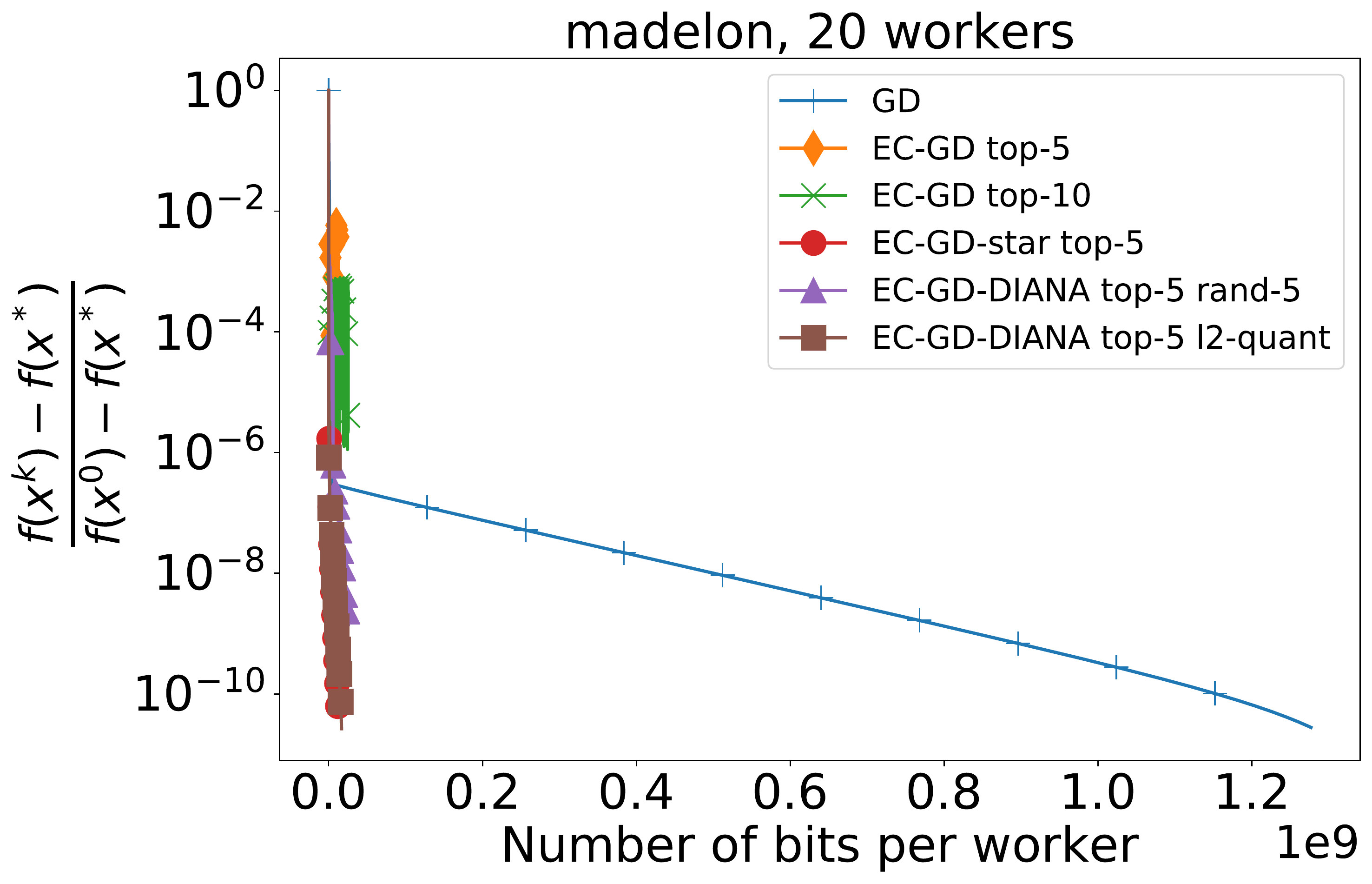}    
	\includegraphics[width=0.32\textwidth]{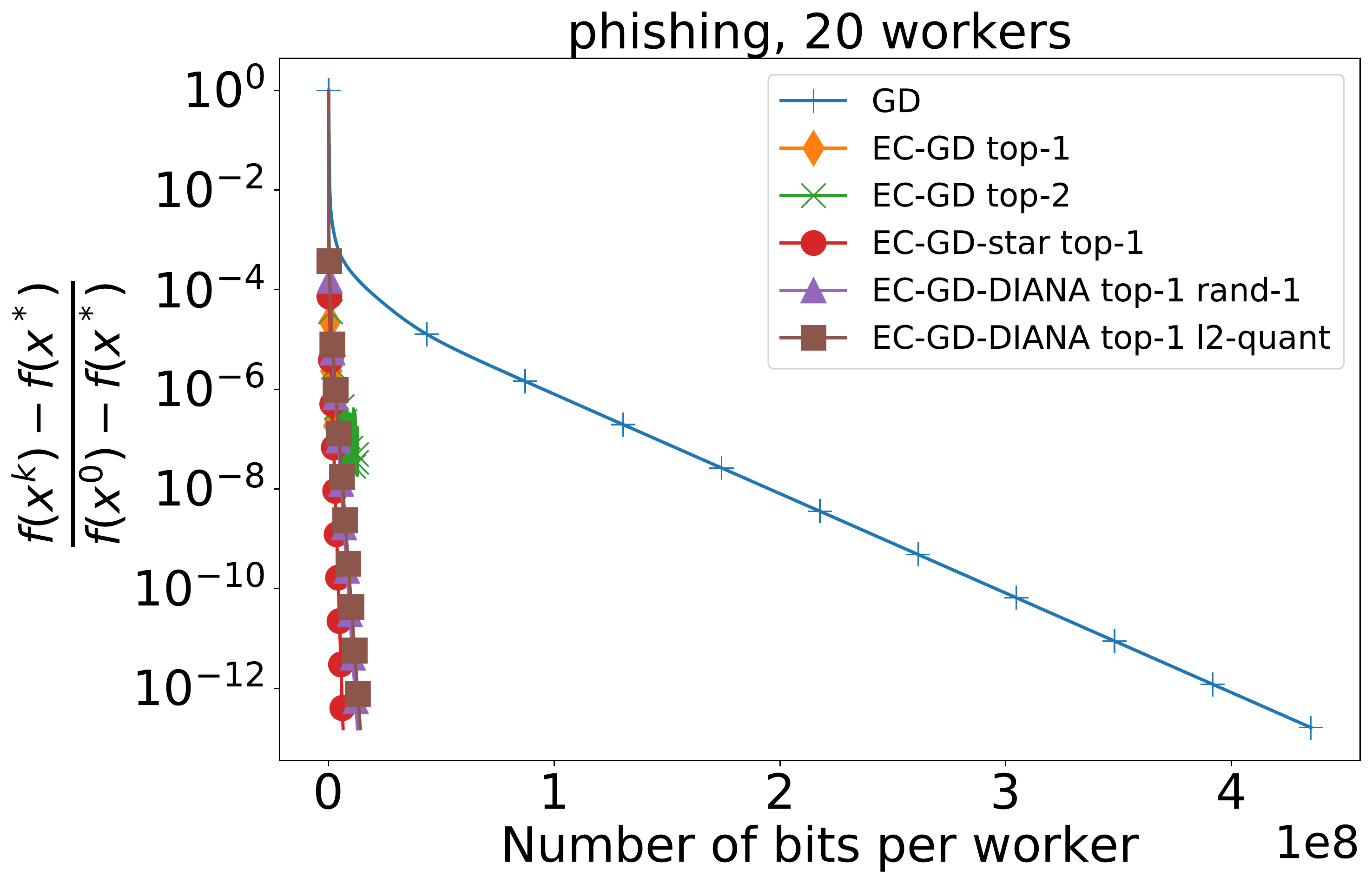}    
    \\
    \includegraphics[width=0.32\textwidth]{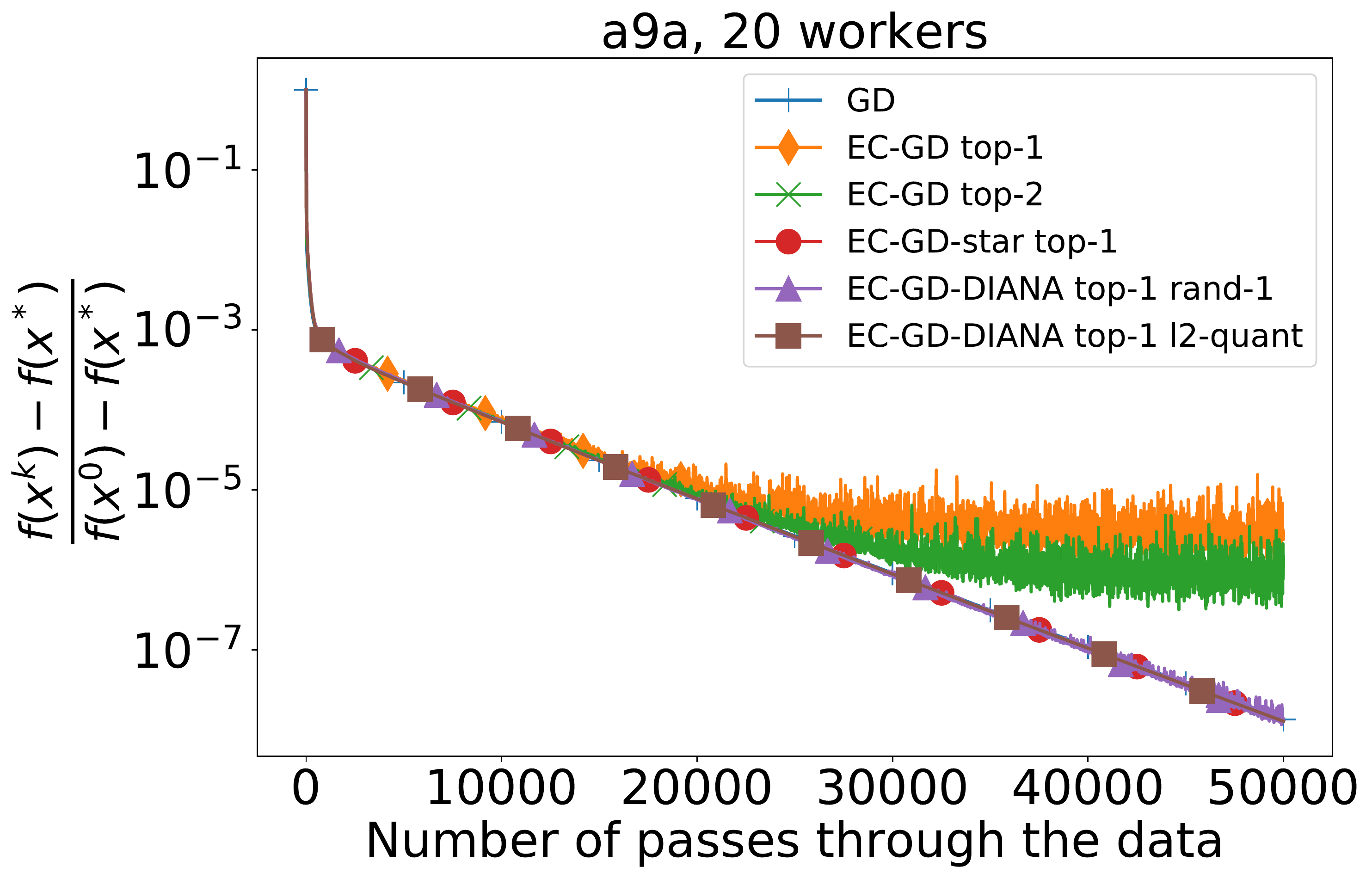}    
	\includegraphics[width=0.32\textwidth]{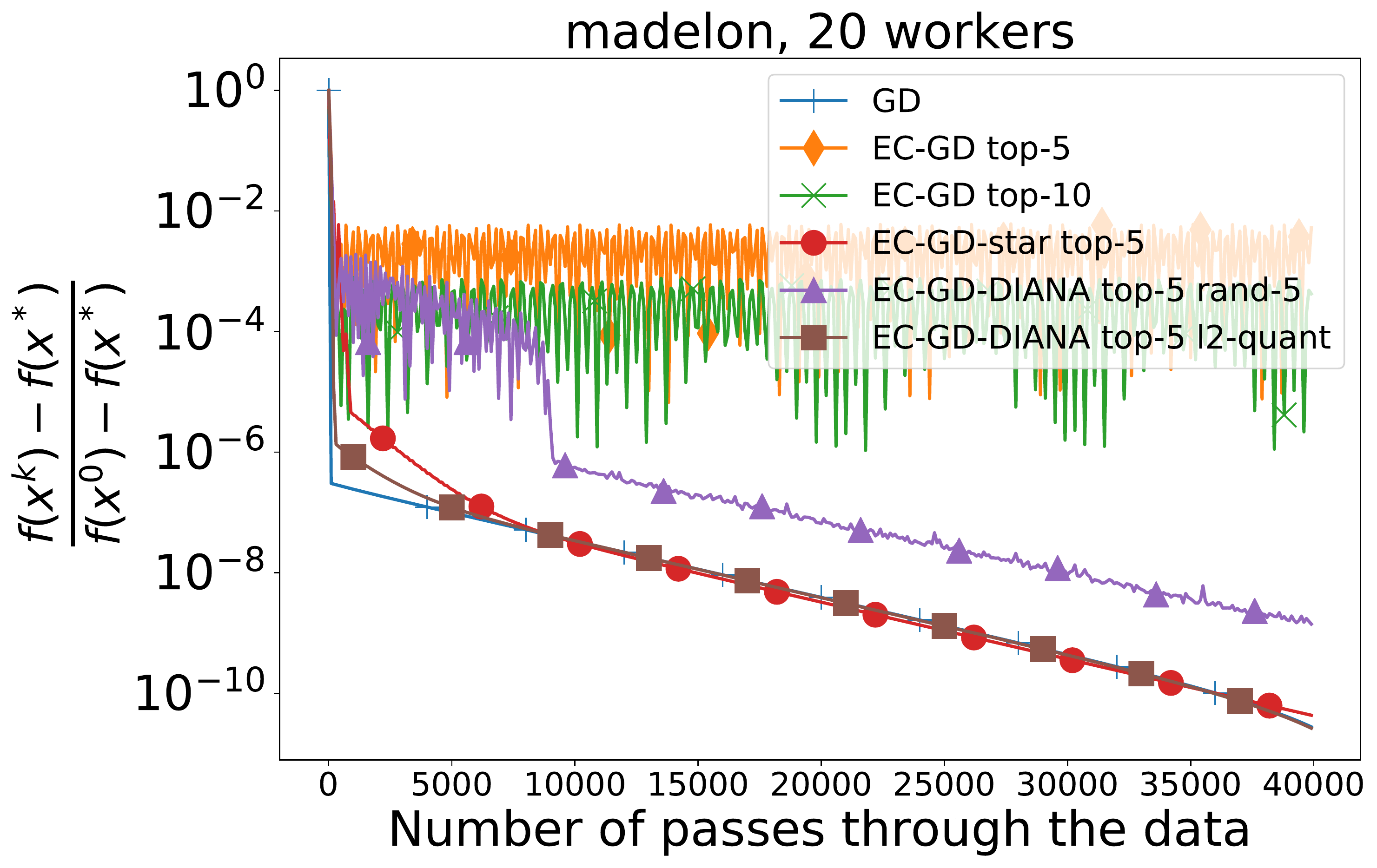}    
	\includegraphics[width=0.32\textwidth]{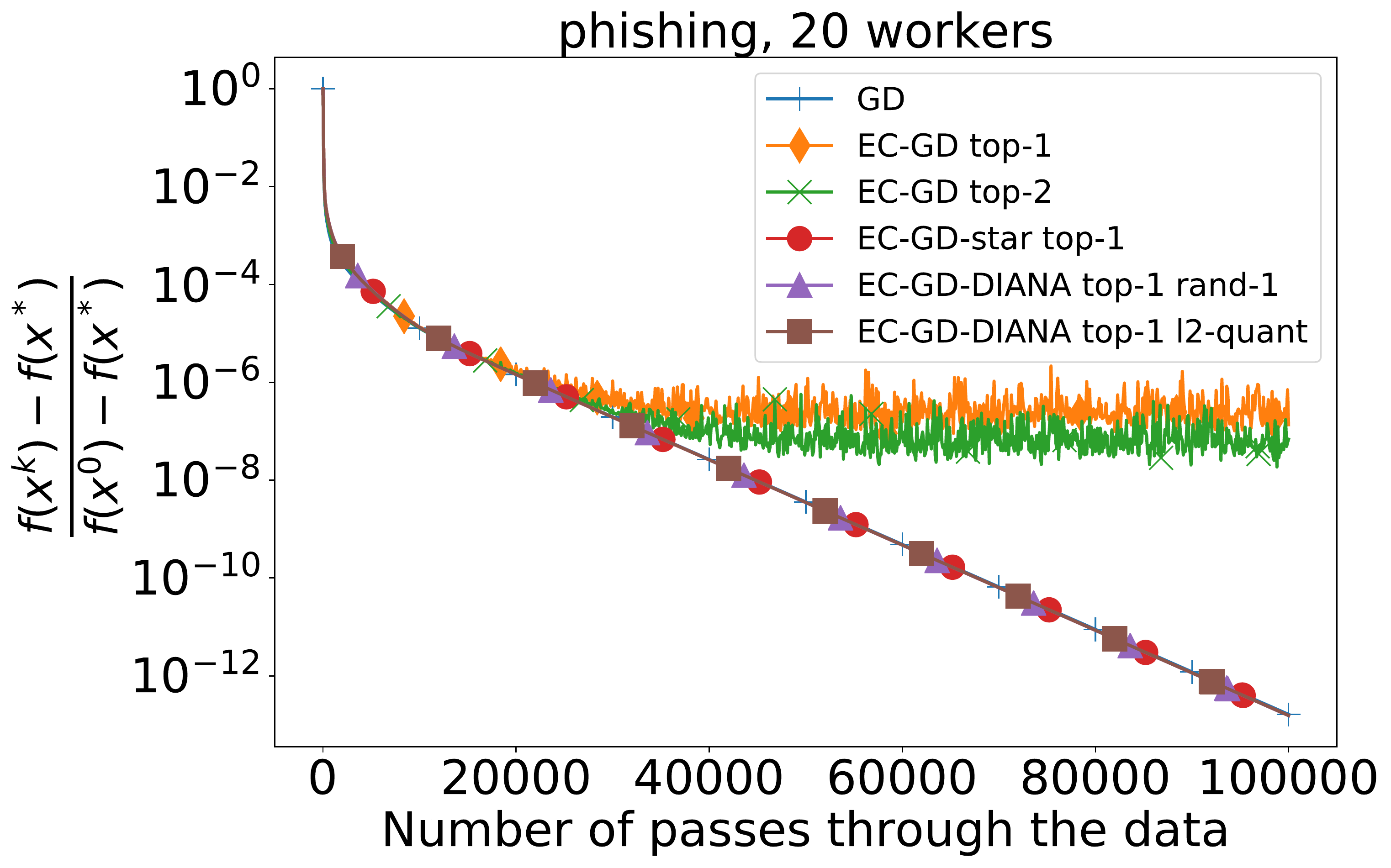}    	
	\\
	\includegraphics[width=0.32\textwidth]{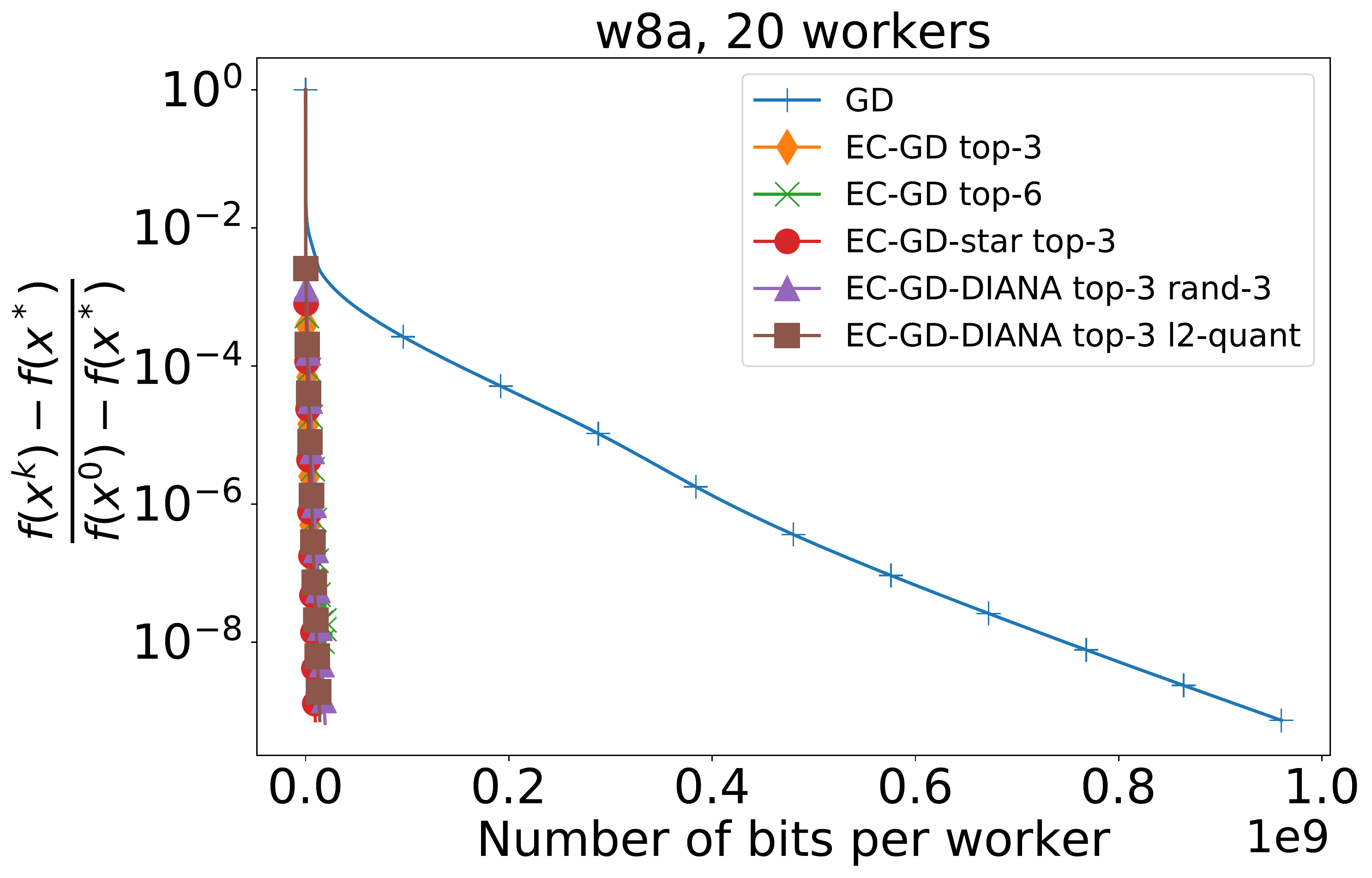}	\includegraphics[width=0.32\textwidth]{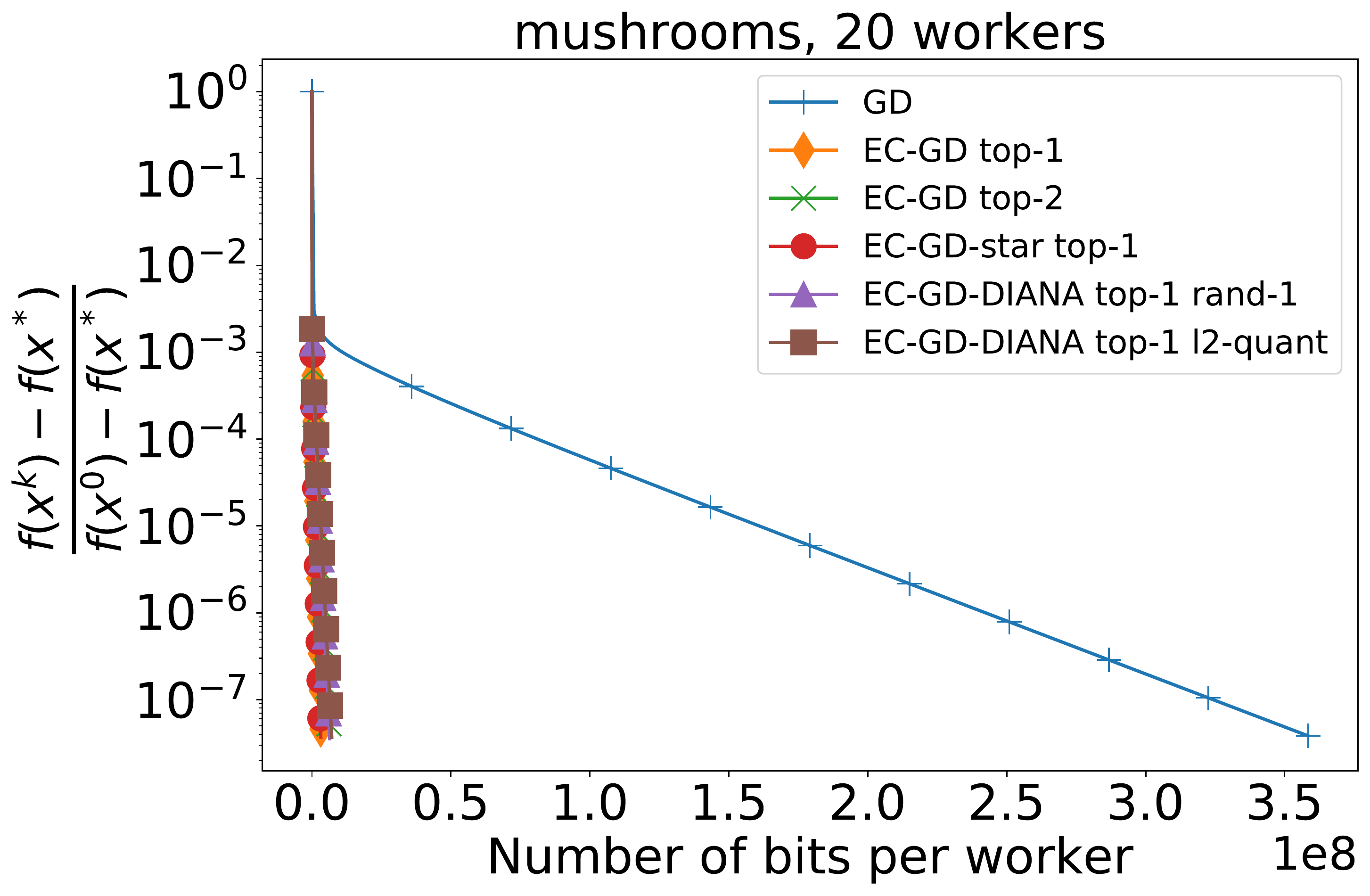}
	\includegraphics[width=0.32\textwidth]{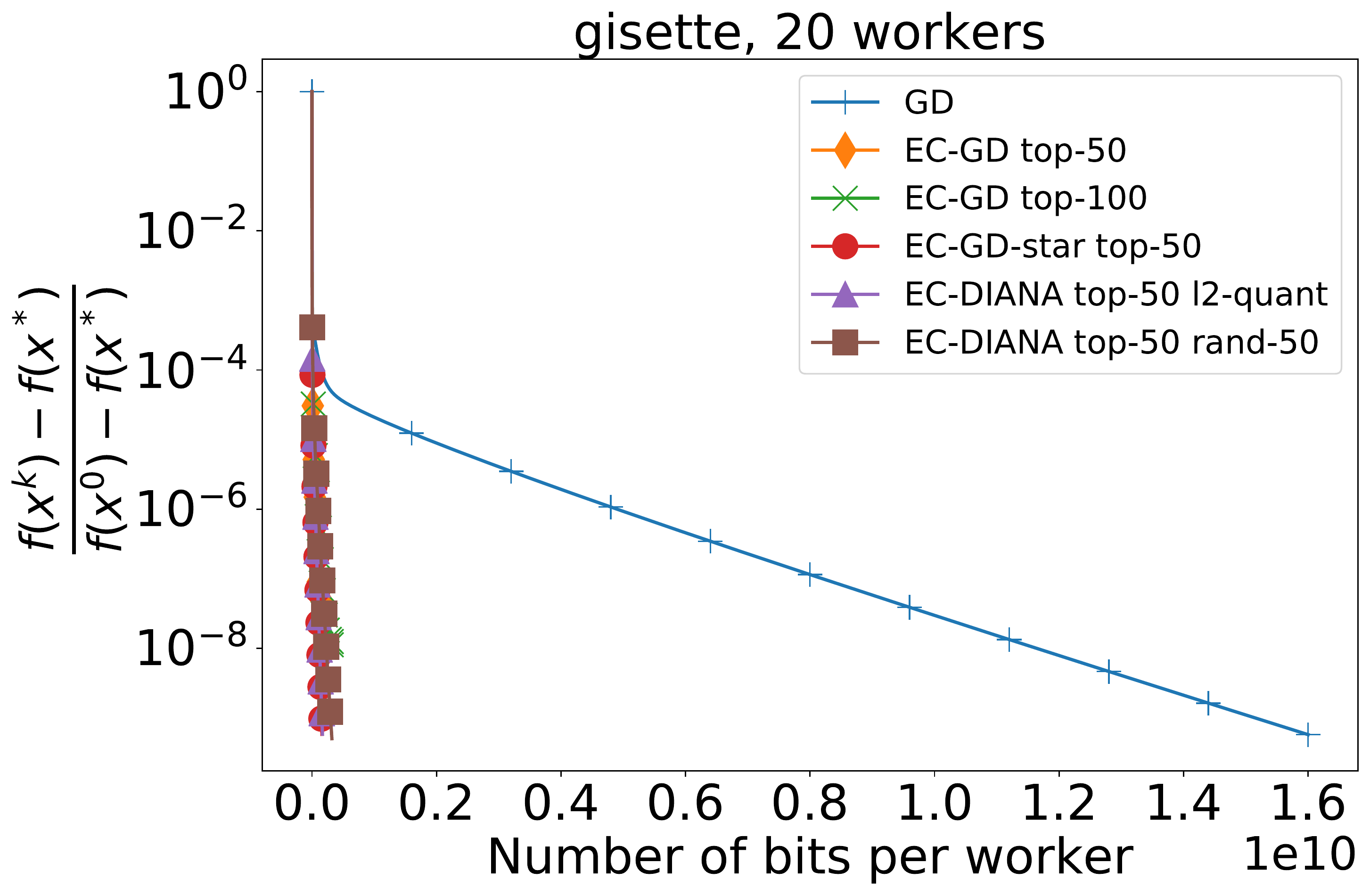}
	\\
	\includegraphics[width=0.32\textwidth]{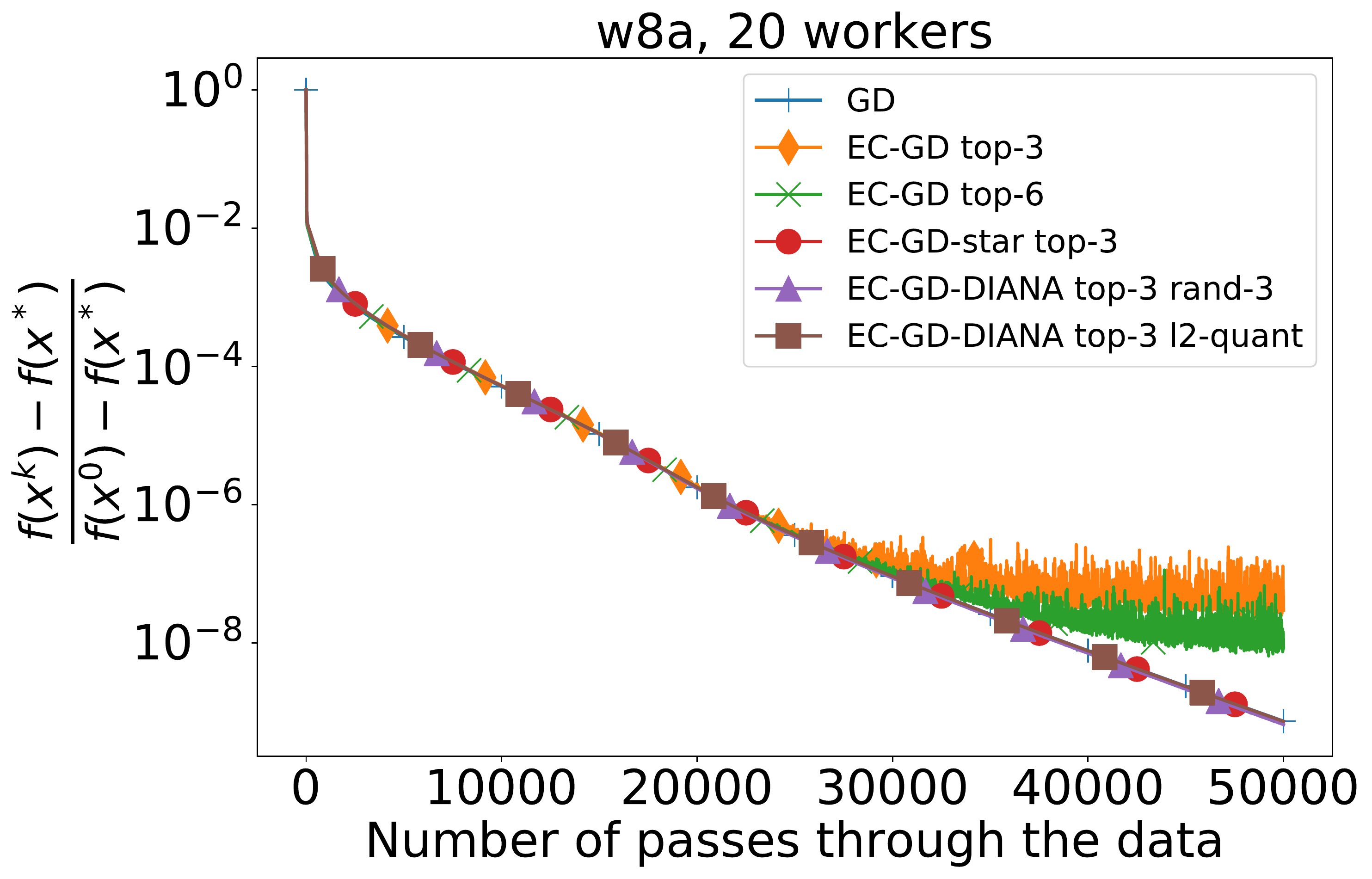}
	\includegraphics[width=0.32\textwidth]{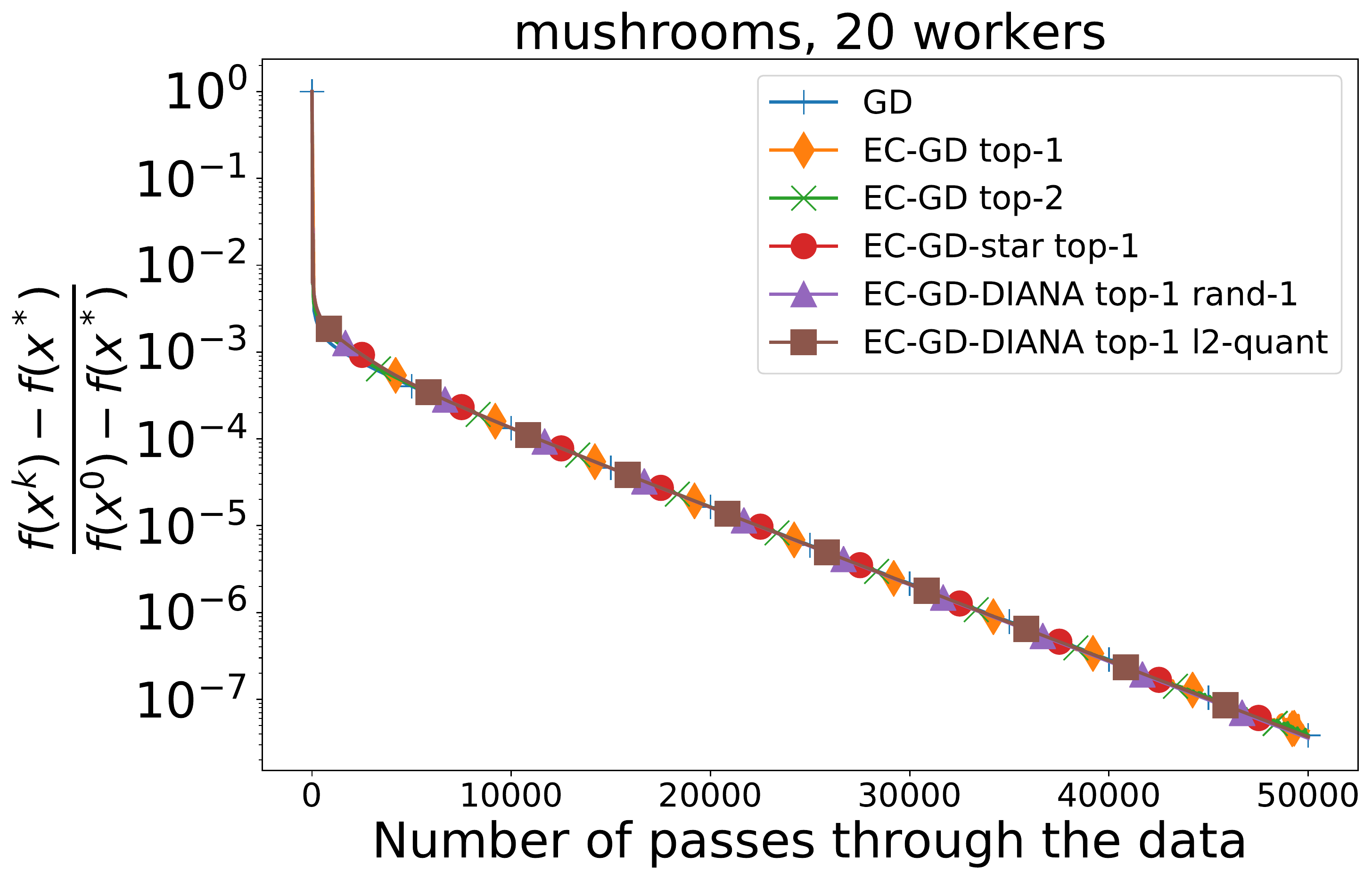}	\includegraphics[width=0.32\textwidth]{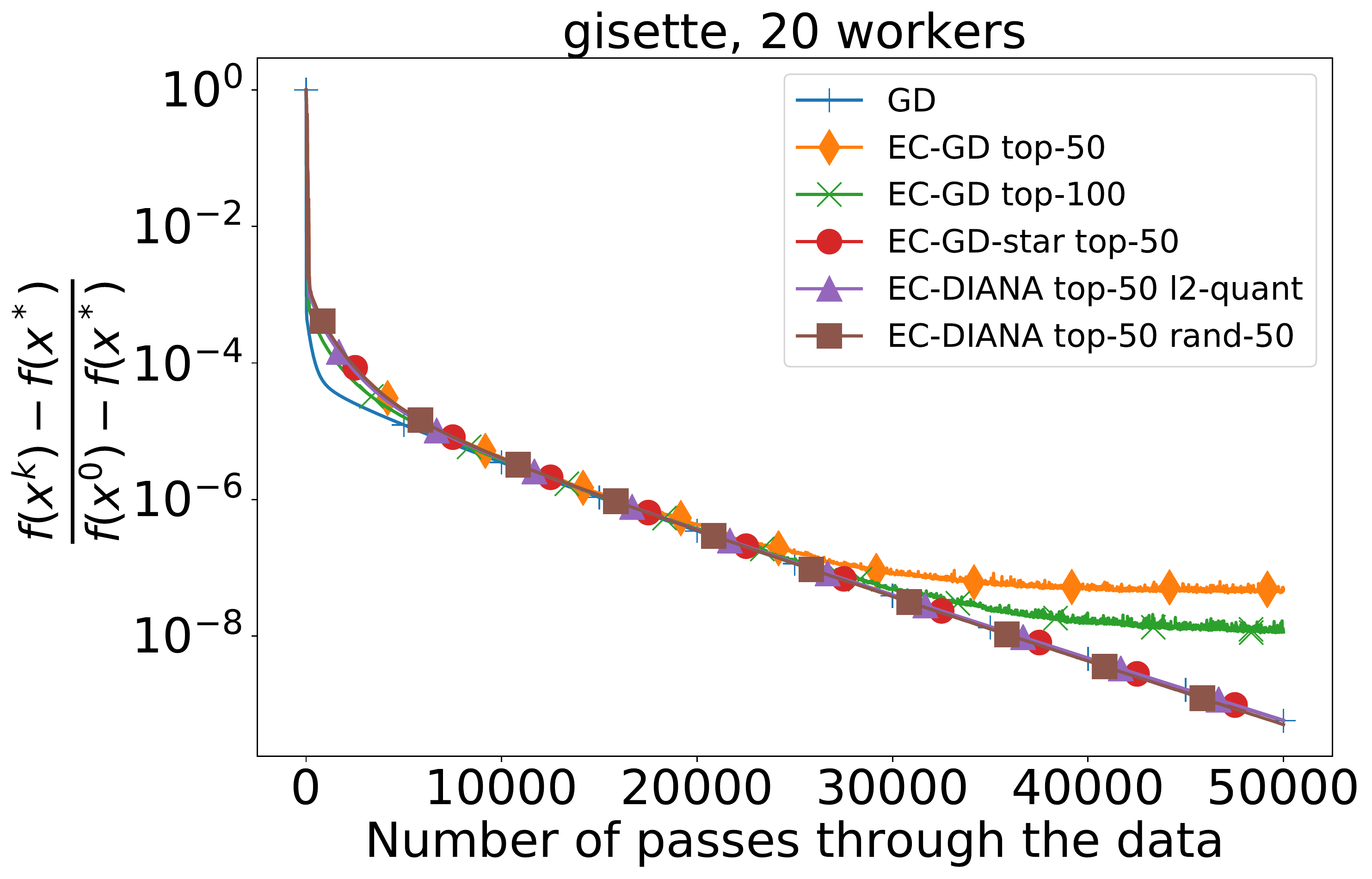}
	\caption{Trajectories of {\tt EC-GD}, {\tt EC-GD-star}, {\tt EC-DIANA} and {\tt GD} applied to solving logistic regression problem with $20$ workers.}
    \label{fig:gd_logreg_20_workers_id}
\end{figure}

\begin{figure}[H]
    \centering
    \includegraphics[width=0.32\textwidth]{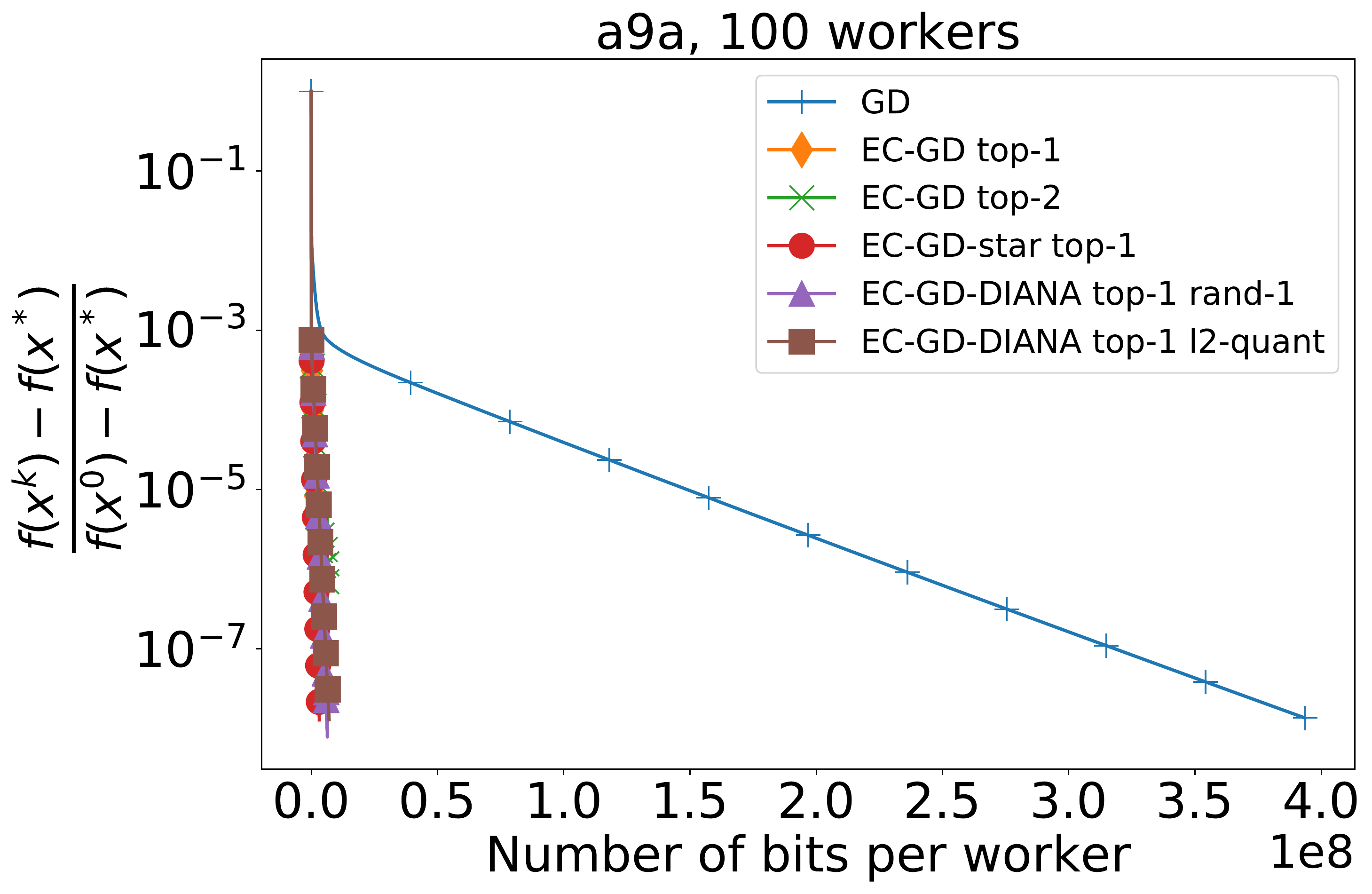}
	\includegraphics[width=0.32\textwidth]{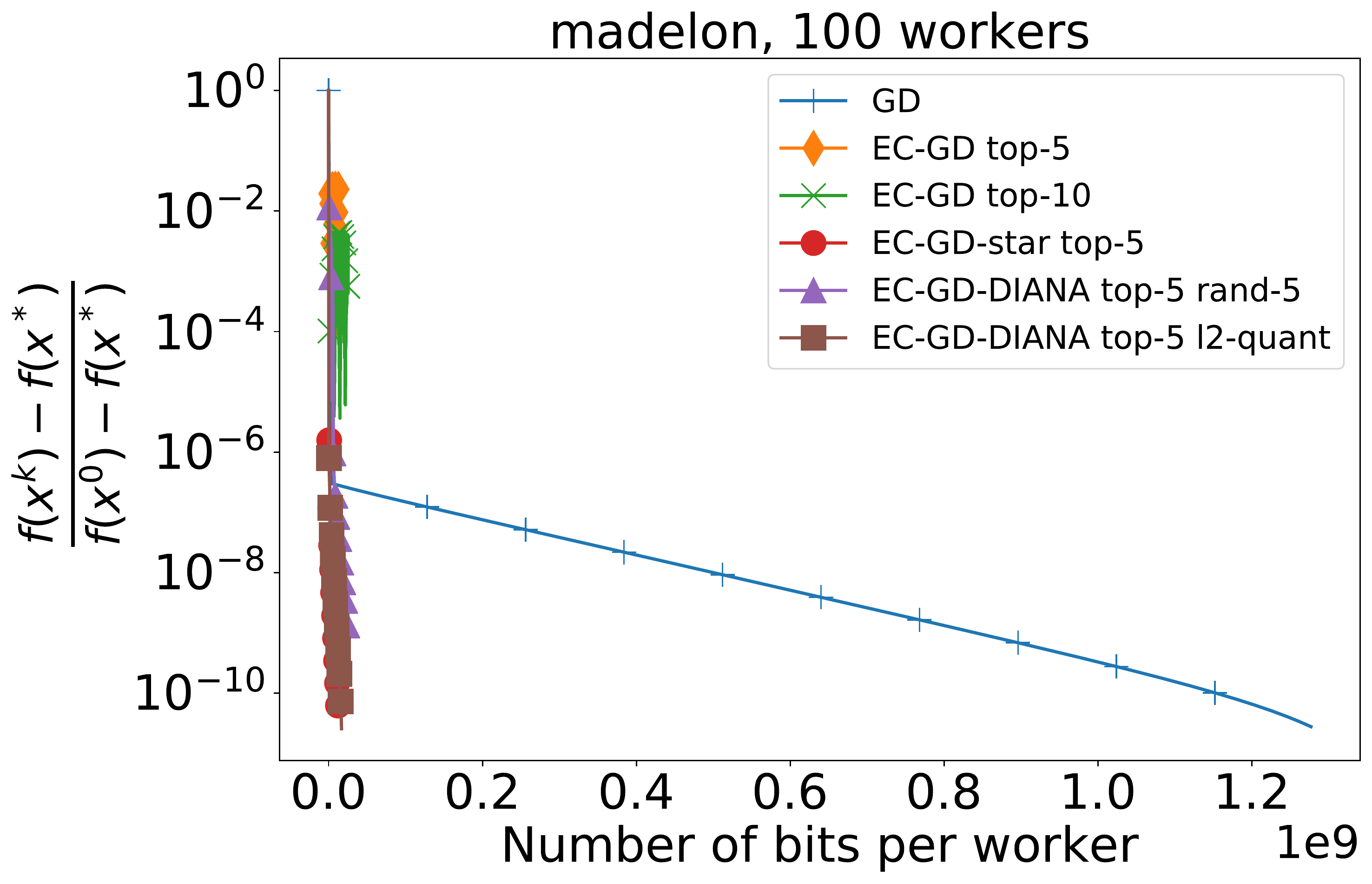}    
	\includegraphics[width=0.32\textwidth]{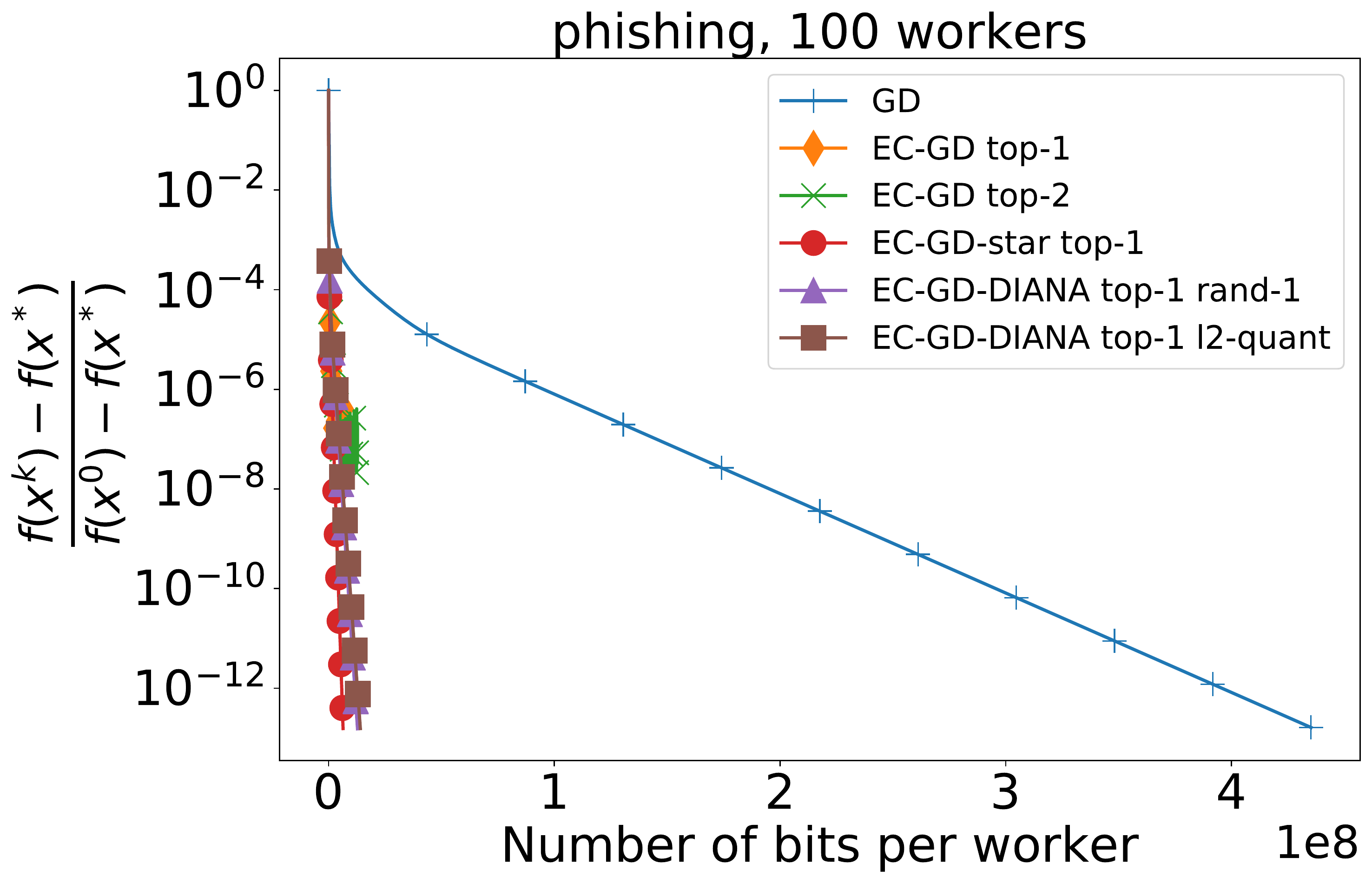}    
    \\
    \includegraphics[width=0.32\textwidth]{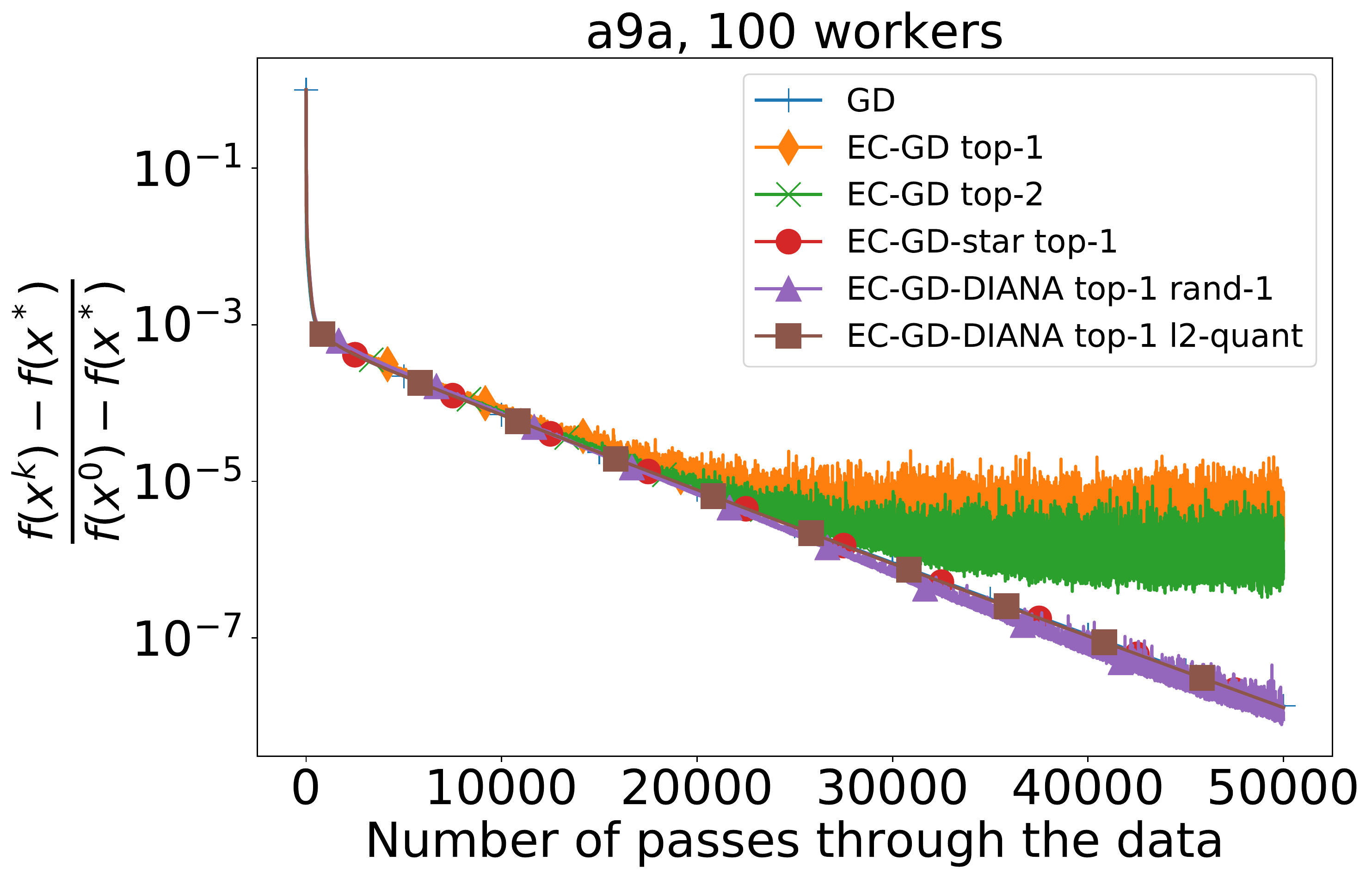}    
	\includegraphics[width=0.32\textwidth]{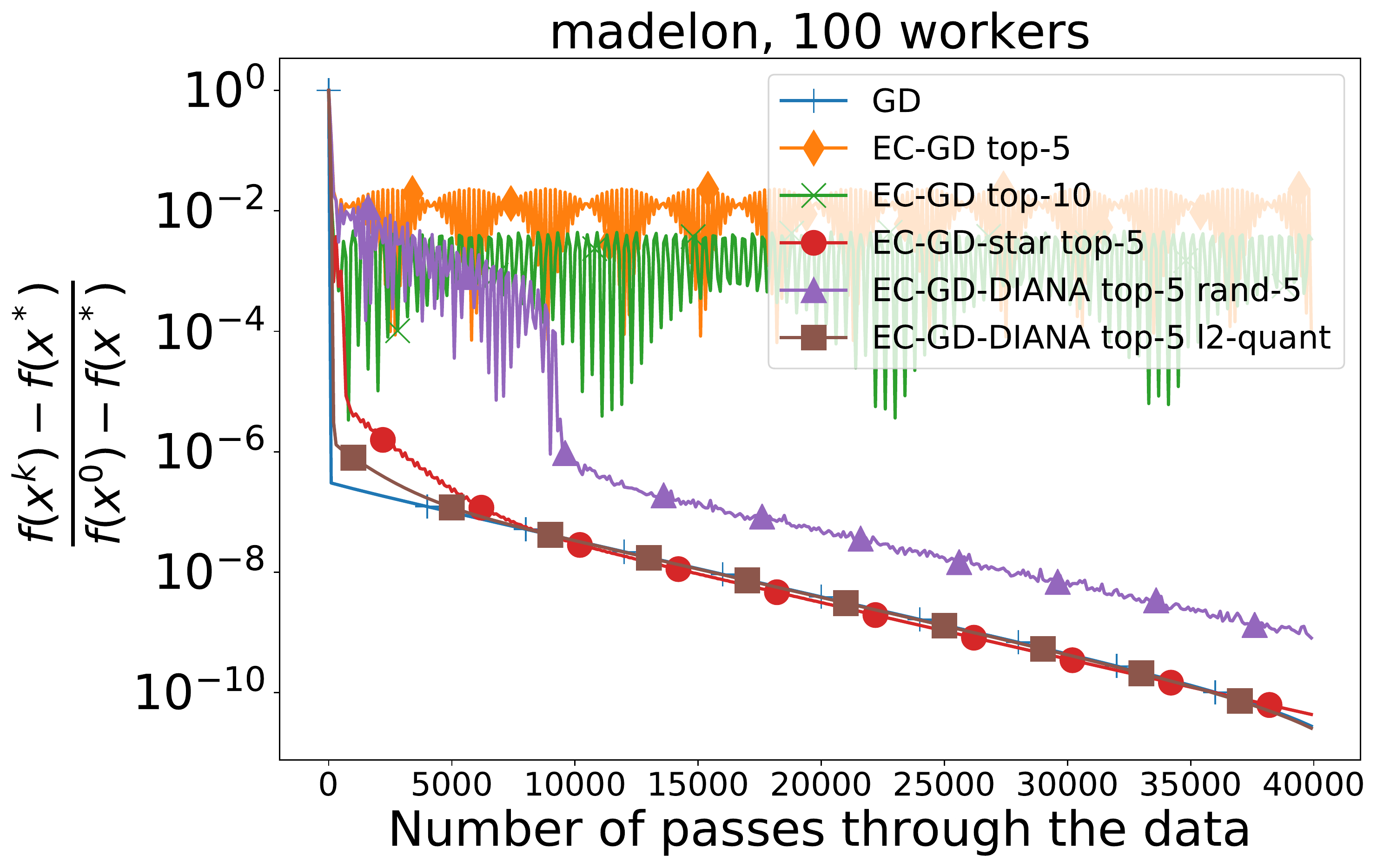}    
	\includegraphics[width=0.32\textwidth]{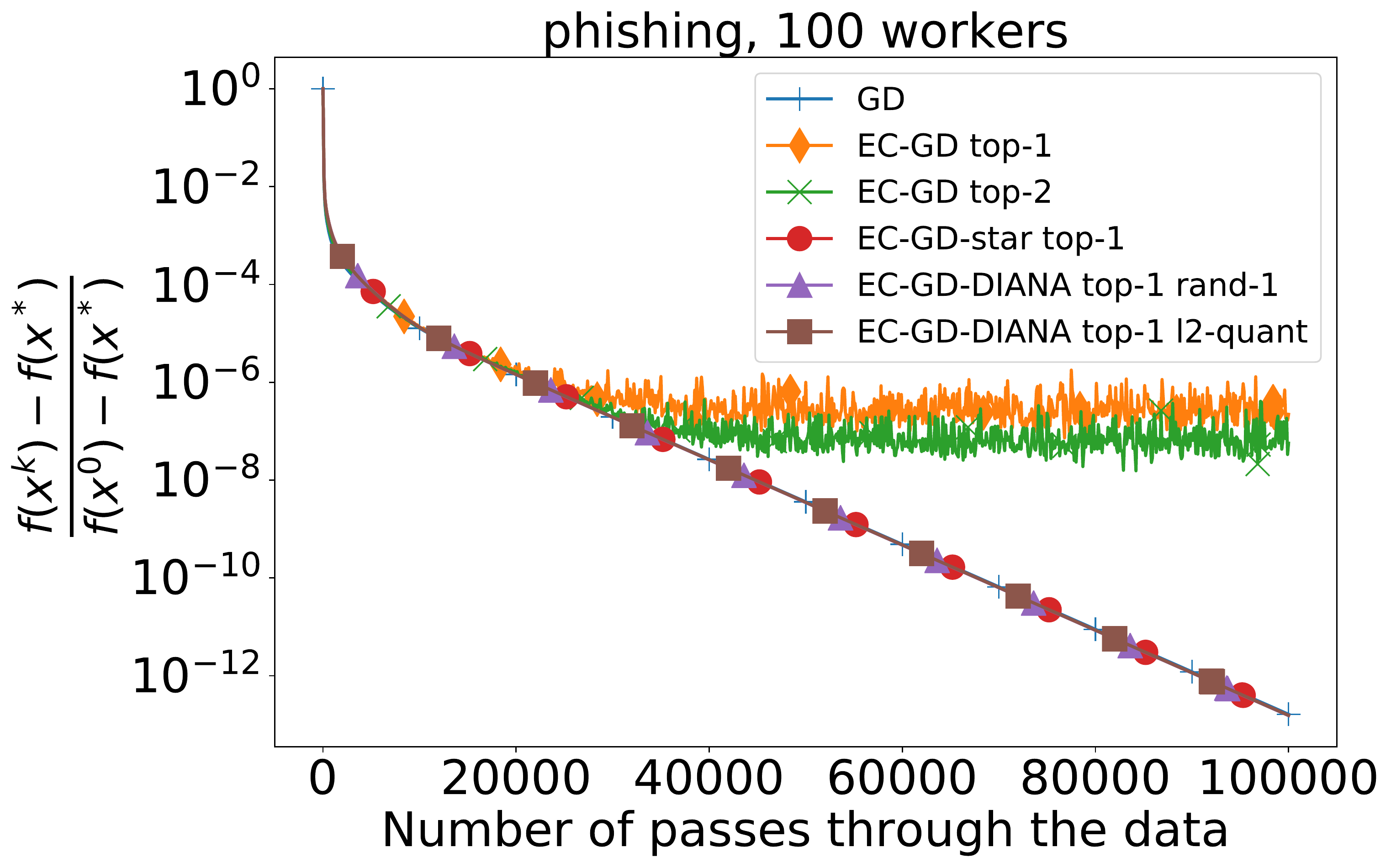}    	
	\\
	\includegraphics[width=0.32\textwidth]{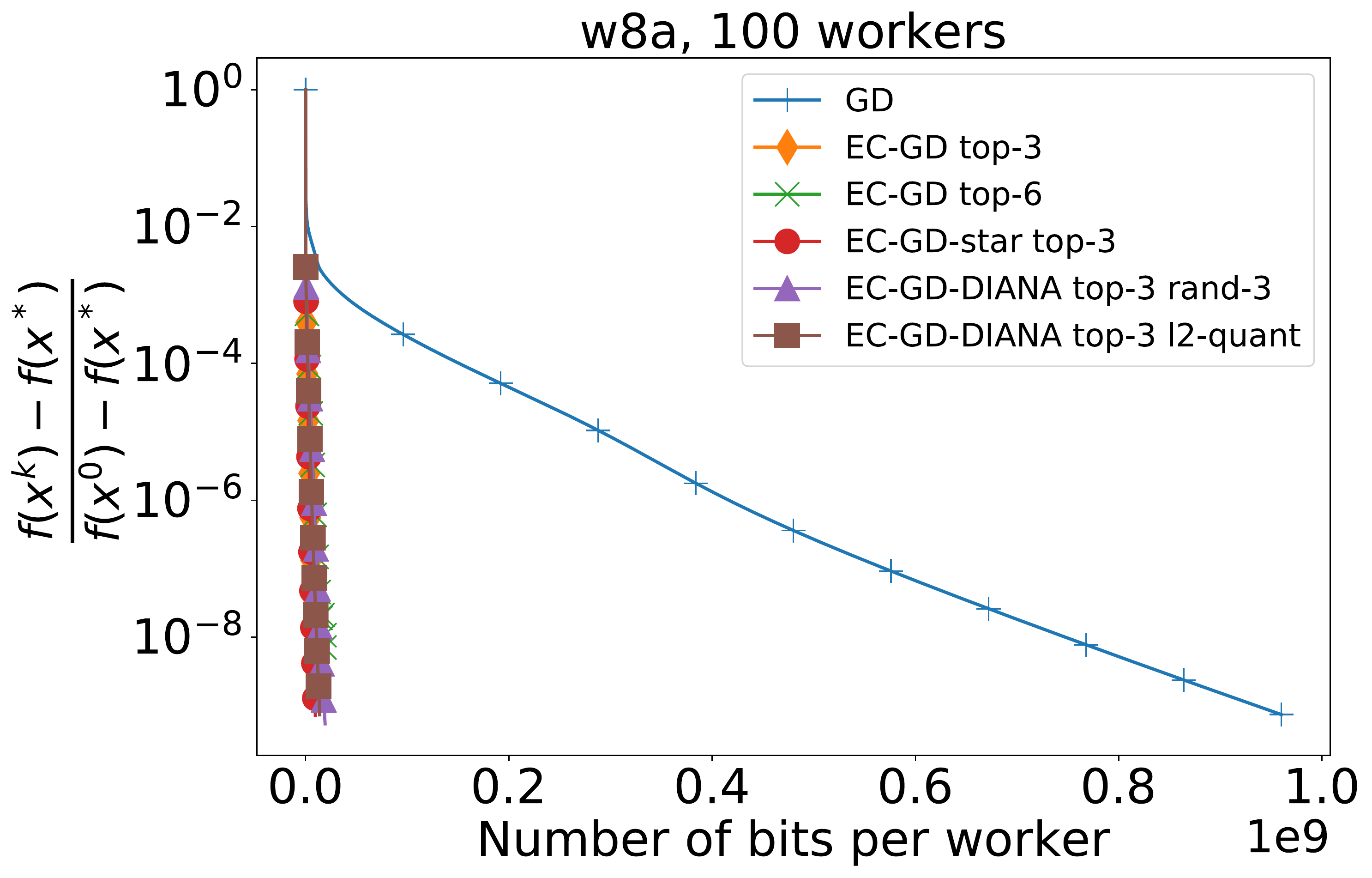}	\includegraphics[width=0.32\textwidth]{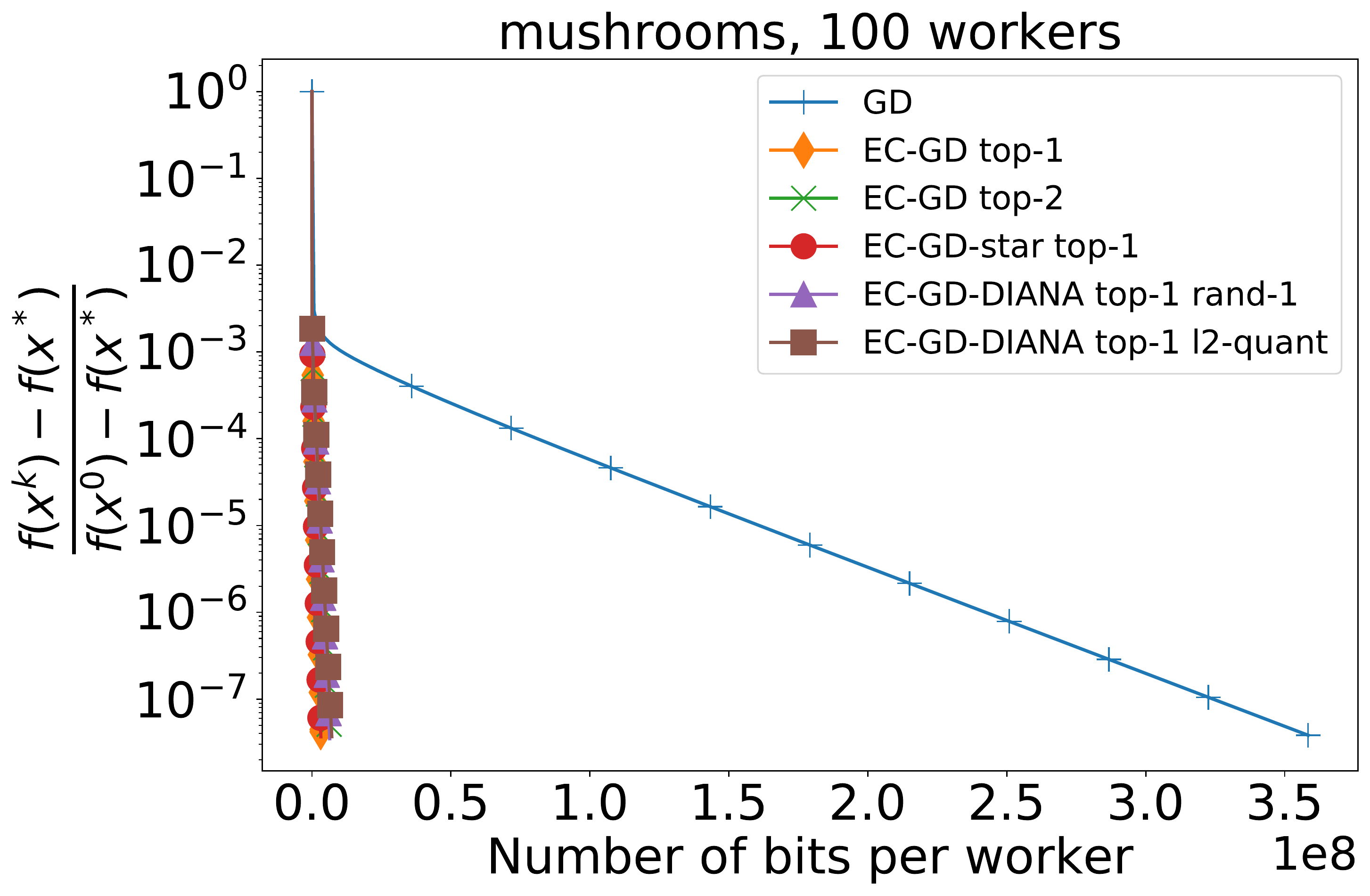}
	\includegraphics[width=0.32\textwidth]{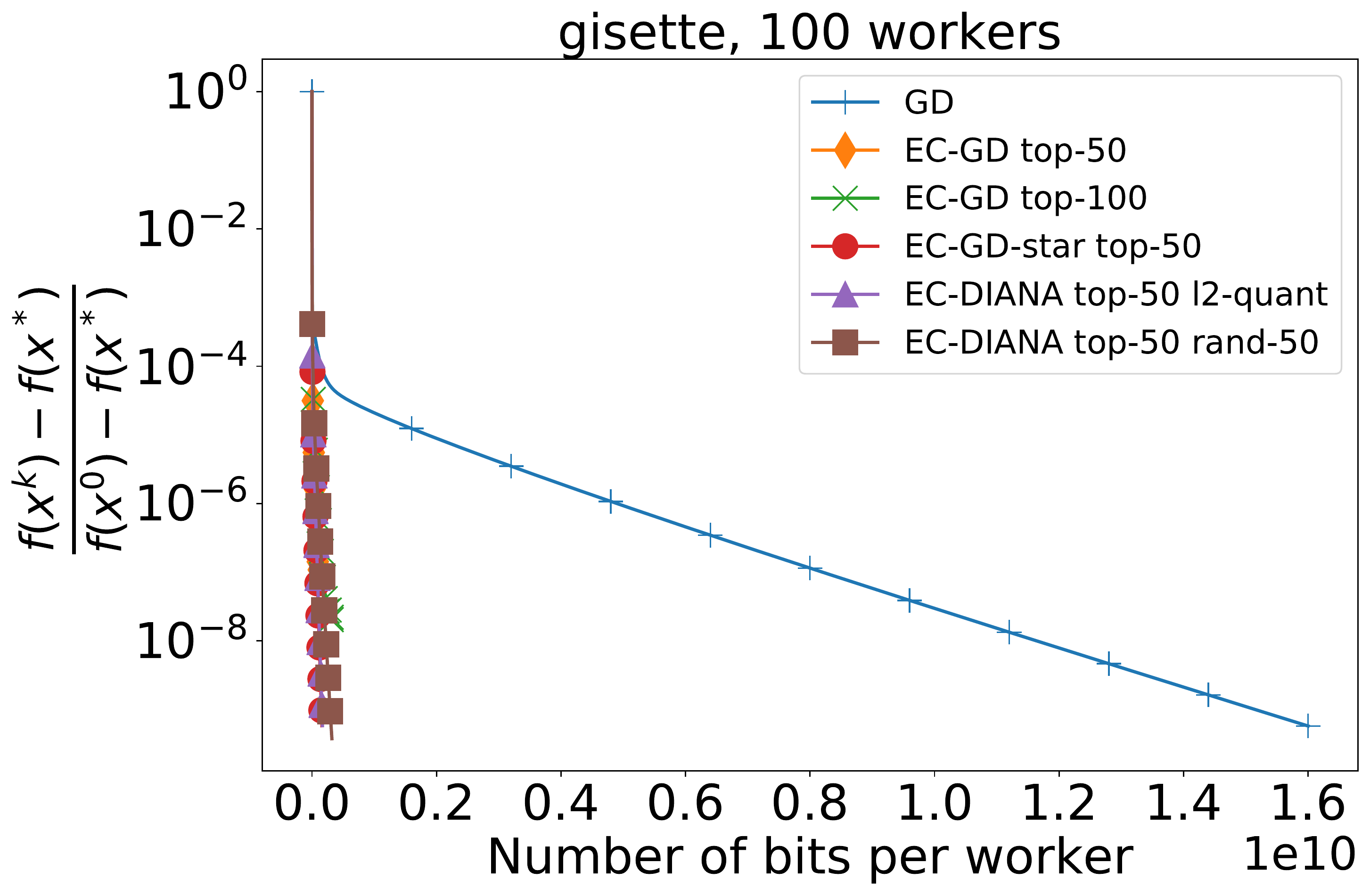}
	\\
	\includegraphics[width=0.32\textwidth]{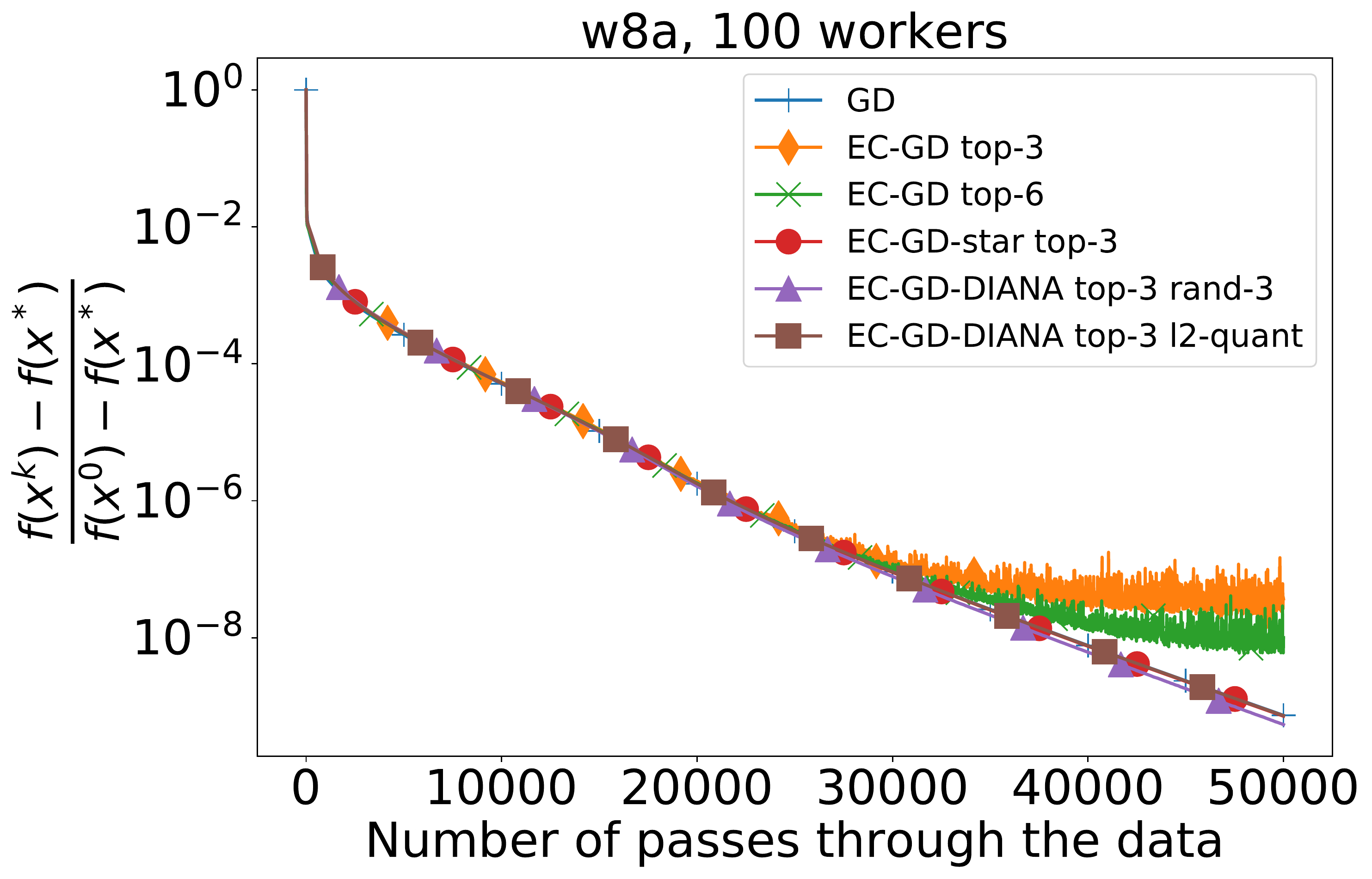}
	\includegraphics[width=0.32\textwidth]{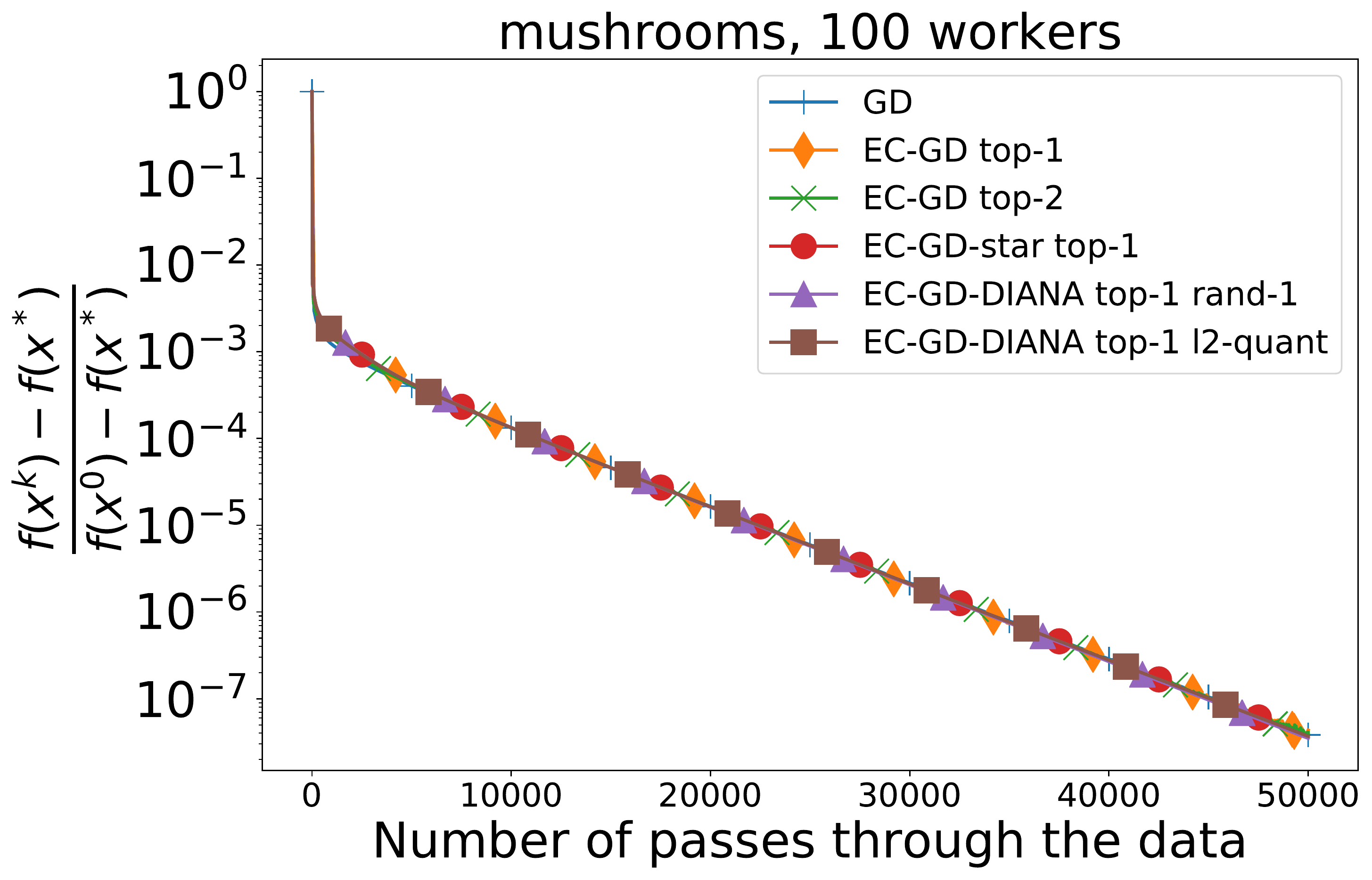}	\includegraphics[width=0.32\textwidth]{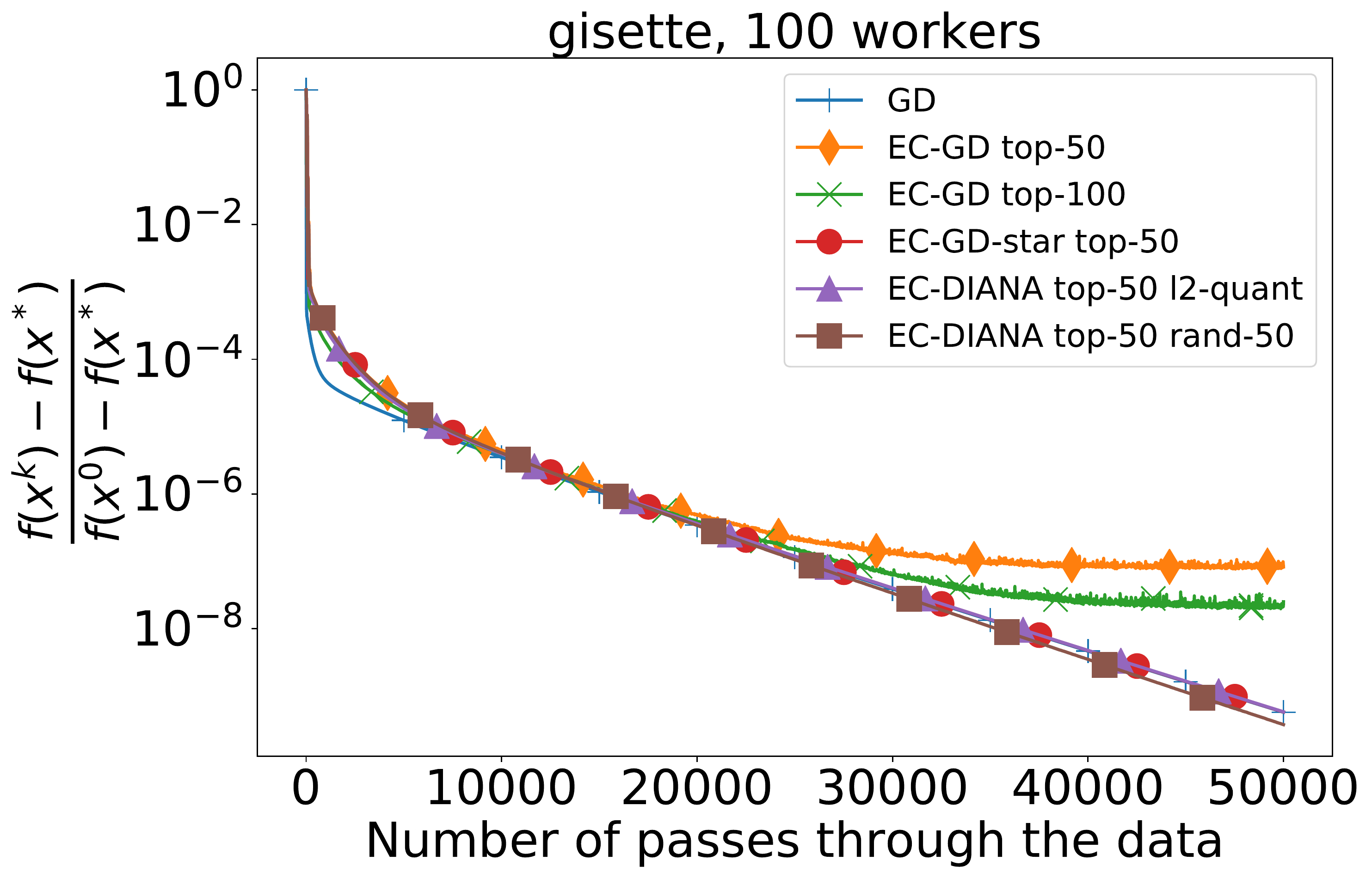}
	\caption{Trajectories of {\tt EC-GD}, {\tt EC-GD-star}, {\tt EC-DIANA} and {\tt GD} applied to solve logistic regression problem with $100$ workers.}
    \label{fig:gd_logreg_100_workers_id}
\end{figure}


\section{Compression Operators: Extra Commentary} \label{sec:compressions}

Communication efficient distributed {\tt SGD} methods based on the idea of communication compression exists in two distinct varieties: i) methods based on unbiased compression operators, and ii) methods based on biased compression operators. The first class of methods is much mire developed than the latter since it is easier to theoretically analyze unbiased operators. The subject of this chapter is the study of the latter and dramatically less developed and understood class.

\subsection{Unbiased Compressors} 

By unbiased compression operators we mean randomized mappings $\cQ:\R^d\to \R$ satisfying the relations
\[\Exp{\cQ(x)} = x \qquad \text{and} \qquad \Exp \|\cQ(x)-x\|^2 \leq \omega \|x\|^2, \qquad \forall x\in \R^d\]
for some $\omega \geq 0$.
While operators satisfying the above relations are often in the literature called \textit{quantization operators}, this class includes compressors which perform sparsification as well.

 Among the first methods using unbiased compressors developed in this field are {\tt QSGD} \cite{alistarh2017qsgd}, {\tt TernGrad} \cite{wen2017terngrad} and {\tt DQGD} \cite{khirirat2018distributed}. The first analysis of {\tt QSGD} and {\tt TernGrad} without bounded gradients assumptions was proposed in \cite{mishchenko2019distributed}, which contains the best known results for {\tt QSGD} and {\tt TernGrad}. However, existing guarantees in the strongly convex case for {\tt QGSD}, {\tt TernGrad}, and {\tt DQGD} establish linear convergence  to some neighborhood of the solution only, even if the workers quantize the full gradients of their functions. This problem was resolved by \cite{mishchenko2019distributed}, who proposed the first method, called {\tt DIANA}, which uses quantization for communication and enjoys the linear rate of convergence to the exact optimum asymptotically in the strongly convex case when workers compute the full gradients of their functions in each iteration. Unlike all previous approaches, {\tt DIANA} is based on the quantization of gradient differences rather than iterates or gradients. In essence, {\tt DIANA} is a technique for reducing the variance  introduced by quantization.  \cite{horvath2019stochastic} generalized the {\tt DIANA} method to the case of more general quantization operators. Moreover, the same authors developed a new method called {\tt VR-DIANA} specially designed to solve problems  \eqref{eq:main_problem_ef} with the individual functions having the finite sum structure \eqref{eq:f_i_sum_ef}.

\subsection{Biased Compressors}

By biased compressors we mean (possibly) randomized mappings $\cC:\R^d\to \R$ satisfying the average contraction relation
\[		\EE\left[\|\cC(x) - x\|^2\right] \le (1 - \delta)\|x\|^2, \qquad \forall x\in \R^d \]
for some $\delta>0$.
	
Perhaps the most popular biased  compression operator is  TopK, which takes vector $x$ as input and substitutes all coordinates of $x$ by zero except the $k$ components with the largest absolute values. However, such a greedy approach applied to simple distributed {\tt SGD} and even distributed {\tt GD} can break the convergence of the method even when applied to simple functions in small dimensions, and may even lead to exponential divergence \cite{beznosikov2020biased}. The \textit{error-feedback} framework described in \cite{karimireddy2019error,stich2020error,stich2018sparsified} and studies in this chapter can fix this problem, and it remains the only known mechanism that does so for all compressors described in \eqref{eq:compression_def}.
This is one of the main motivations for the study of the error-feedback mechanism. For instance, error feedback can fix convergence issues with   methods like {\tt sign-SGD} \cite{Bernstein2019signSGDWM}. The analysis of error feedback by \cite{karimireddy2019error,stich2020error,stich2018sparsified} works either under the assumption that the second moment of the stochastic gradient is uniformly bounded or only for the single-worker case. Recently Beznosikov et al. \cite{beznosikov2020biased} proposed the first analysis of {\tt SGD} with error feedback for the general case of multiple workers without bounded second moment assumption. There is another line of works \cite{pmlr-v97-koloskova19a, KoloskovaLSJ19decentralized} where authors apply arbitrary compressions in the decentralized setup. This approach has better potential than a centralized one in terms of reducing the communication cost. However, in this chapter, we study only centralized architecture.



\section{Proofs for Section~\ref{sec:main_res}}\label{sec:proofs_sec_3}

\subsection{A Lemma}
\begin{lemma}[See also Lemma~8 from \cite{stich2020error}]\label{lem:main_lemma_new}
	Let Assumptions~\ref{ass:quasi_strong_convexity},~\ref{ass:key_assumption_new}~and~\ref{ass:L_smoothness} be satisfied and $\gamma \le \nicefrac{1}{4(A'+C_1M_1+C_2M_2)}$. Then for all $k\ge 0$ we have
	\begin{equation}
		\frac{\gamma}{2}\EE\left[f(x^k) - f(x^*)\right] \le (1-\eta)\EE T^{k} - \EE T^{k+1} + \gamma^2(D_1' + M_1D_2) + 3L\gamma \EE\|e^k\|^2, \label{eq:main_lemma_new}
	\end{equation}	
	where $T^k \eqdef \|\tx^k - x^*\|^2 + M_1\gamma^2 \sigma_{1,k}^2 + M_2\gamma^2 \sigma_{2,k}^2$ and $M_1 = \frac{4B_1'}{3\rho_1}$, $M_2 = \frac{4\left(B_2' + \frac{4}{3}G\right)}{3\rho_2}$. 
\end{lemma}
\begin{proof}
	We start with the upper bound for $\EE\|\tx^{k+1} - x^*\|^2$. First of all, by definition of $\tx^k$ we have
	\begin{eqnarray}
		\|\tx^{k+1} - x^*\|^2 &\overset{\eqref{eq:perturbed_iterates_key_relation}}{=}& \|\tx^k - x^* - \gamma g^k\|^2\notag\\
		&=& \|\tx^k - x^*\|^2 -2\gamma\langle \tx^k - x^*, g^k\rangle + \gamma^2\|g^k\|^2\notag\\
		&=& \|\tx^k - x^*\|^2 -2\gamma\langle x^k - x^*, g^k\rangle + \gamma^2\|g^k\|^2 + 2\gamma\langle x^k - \tx^k, g^k\rangle.\notag
	\end{eqnarray}
	Taking conditional expectation $\EE\left[\cdot\mid x^k\right]$ from the both sides of the previous inequality we get
	\begin{eqnarray}
		\EE\left[\|\tx^{k+1} - x^*\|^2\mid x^k\right] &\overset{\eqref{eq:unbiasedness_new},\eqref{eq:second_moment_bound_new}}{\le}& \|\tx^k - x^*\|^2 -2\gamma\langle x^k - x^*, \nabla f(x^k)\rangle \notag\\
		&&\quad+ \gamma^2\left(2A'(f(x^k) - f(x^*)) + B_1'\sigma_{1,k}^2 + B_2'\sigma_{2,k}^2 + D_1'\right)\notag\\
		&&\quad + 2\gamma\langle x^k - \tx^k, \nabla f(x^k)\rangle\notag\\
		&\overset{\eqref{eq:str_quasi_cvx}}{\le}& \|\tx^k - x^*\|^2 - \gamma\mu\|x^k - x^*\|^2 - \gamma (2 - 2A'\gamma)(f(x^k) - f(x^*)) \notag\\
		&&\quad + \gamma^2 B_1'\sigma_{1,k}^2 + \gamma^2 B_2'\sigma_{2,k}^2 + \gamma^2 D_1' \notag\\
		&&\quad + 2\gamma\langle x^k - \tx^k, \nabla f(x^k)\rangle. \label{eq:main_lemma_technical_1_new}
	\end{eqnarray}
	Next,
	\begin{equation}
		-\|x^k - x^*\|^2 = -\| \tx^k - x^* + x^k - \tx^k\|^2 \overset{\eqref{eq:1/2a_minus_b}}{\le} -\frac{1}{2}\|\tx^k - x^*\|^2 + \|x^k - \tx^k\|^2. \label{eq:main_lemma_technical_2_new}
	\end{equation}
	Using Fenchel-Young inequality we derive an upper bound for the inner product from \eqref{eq:main_lemma_technical_1_new}:
	\begin{equation}
		\langle x^k - \tx^k, \nabla f(x^k)\rangle \overset{\eqref{eq:fenchel_young}}{\le} L\|x^k - \tx^k\|^2 + \frac{1}{4L}\|\nabla f(x^k)\|^2 \overset{\eqref{eq:L_smoothness_cor_3}}{\le} L\|x^k - \tx^k\|^2 + \frac{1}{2}(f(x^k) - f(x^*)).\label{eq:main_lemma_technical_3_new}
	\end{equation}
	Combining previous three inequalities we get
	\begin{eqnarray}
		\EE\left[\|\tx^{k+1} - x^*\|^2\mid x^k\right] &\overset{\eqref{eq:main_lemma_technical_1_new}-\eqref{eq:main_lemma_technical_3_new}}{\le}& \left(1 - \frac{\gamma\mu}{2}\right)\|\tx^k - x^*\|^2 - \gamma\left(1 - 2A'\gamma\right)(f(x^k) - f(x^*))\notag\\
		&&\quad + \gamma^2 B_1'\sigma_{1,k}^2 + \gamma^2 B_2'\sigma_{2,k}^2 + \gamma^2 D_1'\notag\\
		&&\quad + \gamma (2L + \mu)\|x^k - \tx^k\|^2. \label{eq:main_lemma_norm_squared_ub_new}
	\end{eqnarray}
	Taking into account that $T^k = \|\tx^k - x^*\|^2 + M_1\gamma^2 \sigma_{1,k}^2 + M_2\gamma^2 \sigma_{2,k}^2$ with $M_1 = \frac{4B_1'}{3\rho_1}$ and $M_2 = \frac{4\left(B_2' + \frac{4}{3}G\right)}{3\rho_2}$, using the tower property \eqref{eq:tower_property} of mathematical expectation together with $\gamma \le \frac{1}{4(A'+C_1M_1 + C_2M_2)}$, we conclude
	\begin{eqnarray*}
		\EE\left[T^{k+1}\right] &\overset{\eqref{eq:main_lemma_norm_squared_ub_new}}{\le}& \left(1 - \frac{\gamma\mu}{2}\right)\EE\|\tx^k - x^*\|^2 - \gamma\left(1 - 2A'\gamma\right)\EE\left[f(x^k) - f(x^*)\right] + M_1\gamma^2\EE\left[\sigma_{1,k+1}^2\right]\\
		&&\quad M_2\gamma^2\EE\left[\sigma_{2,k+1}^2\right] + \gamma^2 B_1'\sigma_{1,k}^2 + \gamma^2 B_2'\sigma_{2,k}^2 + \gamma^2 D_1' + \gamma (2L + \mu)\EE\|x^k - \tx^k\|^2\\
		&\overset{\eqref{eq:sigma_k+1_bound_1},\eqref{eq:sigma_k+1_bound_2}}{\le}& \left(1 - \frac{\gamma\mu}{2}\right)\EE\|\tx^k - x^*\|^2 + \left(1 + \frac{B_1'}{M_1} - \rho_1\right)M_1\gamma^2\EE\left[\sigma_{1,k}^2\right] \\
		&&\quad + \left(1 + \frac{B_2'+M_1G\rho_1}{M_2} - \rho_2\right)M_2\gamma^2\EE\left[\sigma_{2,k}^2\right]+ \gamma^2(D_1' + M_1D_2)\notag\\
		&&\quad - \gamma\left(1 - 2(A'+C_1M_1+C_2M_2)\gamma\right)\EE\left[f(x^k) - f(x^*)\right] \\
		&&\quad + \gamma(2L+\mu)\EE\|x^k - \tx^k\|^2\\
		&\le& \left(1 - \frac{\gamma\mu}{2}\right)\EE\|\tx^k - x^*\|^2 + \left(1 - \frac{\rho_1}{4}\right)M_1\gamma^2\EE\left[\sigma_{1,k}^2\right] + \left(1 - \frac{\rho_2}{4}\right)M_2\gamma^2\EE\left[\sigma_{2,k}^2\right] \notag\\
		&&\quad - \frac{\gamma}{2}\EE\left[f(x^k) - f(x^*)\right] + \gamma(2L+\mu)\EE\|x^k - \tx^k\|^2 + \gamma^2(D_1' + M_1D_2).
	\end{eqnarray*}
	Since $L\ge \mu$, $\tx^k = x^k - e^k$ and $\eta \eqdef \min\{\frac{\gamma\mu}{2},\frac{\rho_1}{4},\frac{\rho_2}{4}\}$ the last inequality implies
	\begin{equation*}
		\frac{\gamma}{2}\EE\left[f(x^k) - f(x^*)\right] \le (1-\eta)\EE T^{k} - \EE T^{k+1} + \gamma^2(D_1' + M_1D_2) + 3L\gamma \EE\|e^k\|^2,
	\end{equation*}
	which concludes the proof.
\end{proof}

\subsection{Proof of Theorem~\ref{thm:main_result_new}}
\begin{proof}
	Form Lemma~\ref{lem:main_lemma_new} we have
	\begin{equation*}
		\frac{\gamma}{2}\EE\left[f(x^k) - f(x^*)\right] \le (1-\eta)\EE T^{k} - \EE T^{k+1} + \gamma^2(D_1' + M_1D_2) + 3L\gamma \EE\|e^k\|^2.
	\end{equation*}
	Summing up these inequalities for $k=0,\ldots,K$ with weights $w_k = (1-\eta)^{-(k+1)}$ we get
	\begin{eqnarray*}
		\frac{1}{2}\sum\limits_{k=0}^K w_k\EE\left[f(x^k) - f(x^*)\right] &\le& \sum\limits_{k=0}^K\left(\frac{w_k(1-\eta)}{\gamma}\EE T^k - \frac{w_k}{\gamma}\EE T^{k+1}\right) + \gamma(D_1' + M_1D_2)\sum\limits_{k=0}^K w_k\\
		&&\quad +3L\sum\limits_{k=0}^Kw_k\EE\|e^k\|^2 \\
		&\overset{\eqref{eq:sum_of_errors_bound_new},\eqref{eq:w_k_definition_new}}{\le}& \sum\limits_{k=0}^K\left(\frac{w_{k-1}}{\gamma}\EE T^k - \frac{w_k}{\gamma}\EE T^{k+1}\right) + F_1\sigma_{1,0}^2 + F_2\sigma_{2,0}^2  \\
		&&\quad+ \gamma^2(D_1' + M_1D_2 + D_3)W_K\\
		&&\quad + \frac{1}{4}\sum\limits_{k=0}^K w_k\EE\left[f(x^k) - f(x^*)\right].
	\end{eqnarray*}
	Rearranging the terms and using $\bar{x}^K = \frac{1}{W_K}\sum_{k=0}^K w_k x^k$ together with Jensen's inequality we obtain
	\begin{eqnarray*}
		\EE\left[f(\bar x^K) - f(x^*)\right] &\le& \frac{4(T^0 + \gamma F_1 \sigma_{1,0}^2+ \gamma F_2 \sigma_{2,0}^2)}{\gamma W_K} + 4\gamma\left(D_1' + M_1D_2 + D_3\right).
	\end{eqnarray*}
	Finally, using the definition of the sequences $\{W_K\}_{K\ge 0}$ and $\{w_k\}_{k\ge 0}$ we derive that if $\mu > $, then $W_K \ge w_K \ge (1-\eta)^{-K}$ and we get \eqref{eq:main_result_new}. In the case when $\mu = 0$ we have $w_k = 1$ and $W_K = K$ which implies \eqref{eq:main_result_new_cvx}.
\end{proof}


\section{{Distributed \tt SGD} with Compression and Error Compensation: Missing Proofs}\label{sec:ec_sgd_appendix}

\begin{lemma}[Lemma~\ref{lem:ec_sgd_key_lemma_new}]\label{lem:ec_sgd_key_lemma_new_appendix}
	Let Assumptions~\ref{ass:quasi_strong_convexity}~and~\ref{ass:L_smoothness} be satisfied,  Assumption~\ref{ass:key_assumption_finite_sums_new} holds and\footnote{When $\rho_1 = 1$ and $\rho_2=1$ one can always set the parameters in such a way that $B_1 = \widetilde{B}_1 = B_2 = \widetilde{B}_2 = C_1 = C_2 = 0$, $D_2 = 0$. In this case we assume that $\frac{2}{1-\rho_1}\left(\frac{C_1}{\rho_1}+\frac{2GC_2}{\rho_2(1-\rho_2)}\right)\left(\frac{2B_1}{\delta}+\widetilde{B}_1\right) + \frac{2C_2\left(\frac{2B_2}{\delta}+\widetilde{B}_2\right)}{\rho_2(1-\rho_2)} = 0$.}
	\begin{equation}
		\gamma \le \min\left\{\frac{\delta}{4\mu}, \sqrt{\frac{\delta}{96L\left(\frac{2A}{\delta} + \widetilde{A} + \frac{2}{1-\rho_1}\left(\frac{C_1}{\rho_1}+\frac{2GC_2}{\rho_2(1-\rho_2)}\right)\left(\frac{2B_1}{\delta}+\widetilde{B}_1\right) + \frac{2C_2\left(\frac{2B_2}{\delta}+\widetilde{B}_2\right)}{\rho_2(1-\rho_2)}\right)}}\right\},\label{eq:gamma_condition_ec_sgd_new_appendix}
	\end{equation}
	where $ M_1 = \frac{4B_1'}{3\rho_1}$ and $M_2 = \frac{4\left(B_2' + \frac{4}{3}G\right)}{3\rho_2}$. Then {\tt EC-SGD} satisfies Assumption~\ref{ass:key_assumption_new}, i.e., inequality \eqref{eq:sum_of_errors_bound_new} holds with the following parameters:
	\begin{equation}
		F_1 = \frac{24L\gamma^2}{\delta\rho_1(1-\eta)}\left(\frac{2B_1}{\delta}+\widetilde{B}_1\right),\quad F_2 = \frac{24L\gamma^2}{\delta\rho_2(1-\eta)}\left(\frac{2G}{1-\rho_1}\left(\frac{2B_1}{\delta}+\tilde{B}_1\right)+\frac{2B_2}{\delta} + \widetilde{B}_2\right), \label{eq:ec_sgd_parameters_new_appendix}
	\end{equation}
	\begin{equation}
		D_3 = \frac{6L\gamma}{\delta}\left(\frac{D_2}{\rho_1}\left(\frac{2B_1}{\delta}+\widetilde{B}_1\right) + \frac{2D_1}{\delta} + \widetilde{D}_1\right).\label{eq:ec_sgd_parameters_new_2_appendix}
	\end{equation}
\end{lemma}
\begin{proof}
	First of all, we derive an upper bound for the second moment of $e_i^{k+1}$:
	\begin{eqnarray}
		\EE\|e_i^{k+1}\|^2 &\overset{\eqref{eq:e^k_def_ec_sgd},\eqref{eq:tower_property}}{=}& \EE\left[\EE\left[\|e_i^k + \gamma g_i^k - C(e_i^k + \gamma g_i^k)\|^2\mid e_i^k,g_i^k\right]\right] \notag\\
		&\overset{\eqref{eq:compression_def}}{\le}& (1-\delta)\EE\|e_i^k + \gamma g_i^k\|^2 \notag\\
		&\overset{\eqref{eq:tower_property},\eqref{eq:variance_decomposition}}{=}& (1-\delta)\EE\|e_i^k + \gamma \bar{g}_i^k\|^2 + (1-\delta)\gamma^2\EE\|g_i^k-\bar{g}_i^k\|^2 \notag\\
		&\overset{\eqref{eq:a+b_norm_beta}}{\le}& (1-\delta)(1+\beta)\EE\|e_i^k\|^2 + (1-\delta)\left(1+\frac{1}{\beta}\right)\gamma^2\EE\|\bar{g}_i^k\|^2\notag\\
		&&\quad + (1-\delta)\gamma^2\EE\|g_i^k-\bar{g}_i^k\|^2. \notag
	\end{eqnarray}
	Summing up these inequalities for $i=1,\ldots, n$ we get
	\begin{eqnarray}
		\frac{1}{n}\sum\limits_{i=1}^n\EE\|e_i^{k+1}\|^2 &\le& (1-\delta)(1+\beta)\frac{1}{n}\sum\limits_{i=1}^n\EE\|e_i^{k}\|^2\notag\\
		&&\quad + (1-\delta)\left(1+\frac{1}{\beta}\right)\gamma^2\frac{1}{n}\sum\limits_{i=1}^n\EE\|\bar{g}_i^k\|^2\notag\\
		&&\quad + (1-\delta)\gamma^2\frac{1}{n}\sum\limits_{i=1}^n\EE\|g_i^k-\bar{g}_i^k\|^2.\label{eq:ec_sgd_technical_1_new}
	\end{eqnarray}
	Consider $\beta = \frac{\delta}{2(1-\delta)}$. For this choice of $\beta$ we have
	\begin{eqnarray*}
		(1-\delta)(1+\beta) &=& (1-\delta)\left(1 + \frac{\delta}{2(1-\delta)}\right) = 1 - \frac{\delta}{2}\\
		(1-\delta)\left(1+\frac{1}{\beta}\right) &=& (1-\delta)\left(1 + \frac{2(1-\delta)}{\delta}\right) = \frac{(1-\delta)(2-\delta)}{\delta} \le \frac{2(1-\delta)}{\delta}.
	\end{eqnarray*}
	Using this we continue our derivations:
	\begin{eqnarray}
		\frac{1}{n}\sum\limits_{i=1}^n\EE\|e_i^{k+1}\|^2 &\le& \left(1 - \frac{\delta}{2}\right)\frac{1}{n}\sum\limits_{i=1}^n\EE\|e_i^{k}\|^2 + \frac{2\gamma^2(1-\delta)}{\delta}\frac{1}{n}\sum\limits_{i=1}^n\EE\|\bar{g}_i^k\|^2\notag\\
		&&\quad + (1-\delta)\gamma^2\frac{1}{n}\sum\limits_{i=1}^n\EE\|g_i^k-\bar{g}_i^k\|^2\notag\\
		&\overset{\eqref{eq:second_moment_bound_g_i^k_new},\eqref{eq:variance_bound_g_i^k_new}}{\le}& \left(1 - \frac{\delta}{2}\right)\frac{1}{n}\sum\limits_{i=1}^n\EE\|e_i^{k}\|^2 + 2\gamma^2(1-\delta)\left(\frac{2A}{\delta}+\widetilde{A}\right)\EE\left[f(x^k) - f(x^*)\right]\notag\\
		&&\quad + \gamma^2(1-\delta)\left(\frac{2B_1}{\delta}+\widetilde{B}_1\right)\EE\sigma_{1,k}^2 + \gamma^2(1-\delta)\left(\frac{2B_2}{\delta}+\widetilde{B}_2\right)\EE\sigma_{2,k}^2\notag\\
		&&\quad + \gamma^2(1-\delta)\left(\frac{2D_1}{\delta}+\widetilde{D}_1\right). \label{eq:ec_sgd_technical_2_new} 
	\end{eqnarray}
	Unrolling the recurrence above we get
	\begin{eqnarray}
		\frac{1}{n}\sum\limits_{i=1}^n\EE\|e_i^{k+1}\|^2 &\overset{\eqref{eq:ec_sgd_technical_2_new}}{\le}& 2\gamma^2(1-\delta)\left(\frac{2A}{\delta}+\widetilde{A}\right)\sum\limits_{l=0}^k\left(1 - \frac{\delta}{2}\right)^{k-l}\EE\left[f(x^l) - f(x^*)\right]\notag\\
		&&\quad + \gamma^2(1-\delta)\left(\frac{2B_1}{\delta}+\widetilde{B}_1\right)\sum\limits_{l=0}^k\left(1 - \frac{\delta}{2}\right)^{k-l}\EE\sigma_{1,l}^2 \notag\\
		&&\quad + \gamma^2(1-\delta)\left(\frac{2B_2}{\delta}+\widetilde{B}_2\right)\sum\limits_{l=0}^k\left(1 - \frac{\delta}{2}\right)^{k-l}\EE\sigma_{2,l}^2 \notag\\
		&&\quad + \gamma^2(1-\delta)\left(\frac{2D_1}{\delta}+\widetilde{D}_1\right)\sum\limits_{l=0}^k\left(1 - \frac{\delta}{2}\right)^{k-l}\label{eq:ec_sgd_technical_3_new}
	\end{eqnarray}
	which implies
	\begin{eqnarray}
		3L\sum\limits_{k=0}^Kw_k\EE\|e^k\|^2 &\overset{\eqref{eq:e^k_def_ec_sgd}}{=}& 3L\sum\limits_{k=0}^Kw_k\EE\left\|\frac{1}{n}\sum\limits_{i=1}^n e_i^k\right\|^2 \overset{\eqref{eq:a_b_norm_squared}}{\le} 3L\sum\limits_{k=0}^K w_k\frac{1}{n}\sum\limits_{i=1}^n\EE\left\|e_i^k\right\|^2\notag\\
		&\overset{\eqref{eq:ec_sgd_technical_3_new}}{\le}& \frac{6L\gamma^2(1-\delta)}{1-\frac{\delta}{2}}\left(\frac{2A}{\delta}+\widetilde{A}\right)\sum\limits_{k=0}^K\sum\limits_{l=0}^k w_k\left(1 - \frac{\delta}{2}\right)^{k-l}\EE\left[f(x^l) - f(x^*)\right]\notag\\
		&&\quad + \frac{3L\gamma^2(1-\delta)}{1-\frac{\delta}{2}}\left(\frac{2B_1}{\delta}+\widetilde{B}_1\right)\sum\limits_{k=0}^K\sum\limits_{l=0}^k w_k\left(1 - \frac{\delta}{2}\right)^{k-l}\EE\sigma_{1,l}^2 \notag\\
		&&\quad + \frac{3L\gamma^2(1-\delta)}{1-\frac{\delta}{2}}\left(\frac{2B_2}{\delta}+\widetilde{B}_2\right)\sum\limits_{k=0}^K\sum\limits_{l=0}^k w_k\left(1 - \frac{\delta}{2}\right)^{k-l}\EE\sigma_{2,l}^2 \notag\\
		&&\quad +\frac{3L\gamma^2(1-\delta)}{1-\frac{\delta}{2}}\left(\frac{2D_1}{\delta}+\widetilde{D}_1\right)\sum\limits_{k=0}^K\sum\limits_{l=0}^k w_k\left(1 - \frac{\delta}{2}\right)^{k-l}.\label{eq:ec_sgd_technical_4_new}
	\end{eqnarray}
	In the remaining part of the proof we derive upper bounds for three terms in the right-hand side of the previous inequality. First of all, recall that $w_k = (1 - \eta)^{-(k+1)}$ and $\eta = \min\left\{\frac{\gamma\mu}{2}, \frac{\rho_1}{4}, \frac{\rho_2}{4}\right\}$. It implies that for all $0 \le i < k$ we have
	\begin{eqnarray}
		w_k &=& (1 - \eta)^{-(k-j+1)}\left(1 - \eta\right)^{-j} \overset{\eqref{eq:1-p/2_inequality}}{\le} w_{k-j}\left(1 + 2\eta\right)^{j} \notag\\
		&\le& w_{k-j}\left(1 + \gamma\mu\right)^{j} \overset{\eqref{eq:gamma_condition_ec_sgd_new}}{\le} w_{k-j}\left(1 + \frac{\delta}{4}\right)^j, \label{eq:ec_sgd_technical_5_new}\\
		w_k &=& \left(1 - \eta\right)^{-(k-j+1)}\left(1 - \eta\right)^{-j} \overset{\eqref{eq:1-p/2_inequality}}{\le} w_{k-j}\left(1 + 2\eta\right)^j \notag\\
		&\le& w_{k-j}\left(1 + \frac{\min\{\rho_1,\rho_2\}}{2}\right)^j. \label{eq:ec_sgd_technical_6_new}
	\end{eqnarray}
	For simplicity, we introduce new notation: $r_k \eqdef \EE\left[f(x^k) - f(x^*)\right]$. Using this we get
	\begin{eqnarray}
		\sum\limits_{k=0}^K\sum\limits_{l=0}^k w_k\left(1 - \frac{\delta}{2}\right)^{k-l}r_l &\overset{\eqref{eq:ec_sgd_technical_5_new}}{\le}& \sum\limits_{k=0}^K\sum\limits_{l=0}^k w_l r_l\left(1 + \frac{\delta}{4}\right)^{k-l}\left(1 - \frac{\delta}{2}\right)^{k-l}\notag\\
		&\overset{\eqref{eq:1+p/2_inequality}}{\le}& \sum\limits_{k=0}^K\sum\limits_{l=0}^k w_l r_l\left(1 - \frac{\delta}{4}\right)^{k-l}\notag\\
		&\le& \left(\sum\limits_{k=0}^K w_k r_k\right)\left(\sum\limits_{k=0}^\infty \left(1 - \frac{\delta}{4}\right)^{k}\right) = \frac{4}{\delta}\sum\limits_{k=0}^K w_k r_k. \label{eq:ec_sgd_technical_7_new}
	\end{eqnarray}
	Next, we apply our assumption on $\sigma_{2,k}^2$ and derive that
	\begin{eqnarray}
		\EE\sigma_{2,k+1}^2 &\overset{\eqref{eq:sigma_k+1_bound_2}}{\le}& (1 - \rho_2)\EE\sigma_{2,k}^2 + 2C_2 \underbrace{\EE\left[f(x^k) - f(x^*)\right]}_{r_k}\notag\\
		&\le& (1-\rho_2)^{k+1}\sigma_{2,0}^2 + 2C_2\sum\limits_{l=0}^{k}(1-\rho_2)^{k-l}r_l,\label{eq:sigma_2_k_useful_recurrence_new}
	\end{eqnarray}
	hence
	\begin{eqnarray}
		\sum\limits_{k=0}^K\sum\limits_{l=0}^k w_k\left(1 - \frac{\delta}{2}\right)^{k-l}\EE\sigma_{2,l}^2 &\le& \sum\limits_{k=0}^K\sum\limits_{l=0}^k w_k\left(1 - \frac{\delta}{2}\right)^{k-l}(1-\rho_2)^l\sigma_{2,0}^2\notag\\
		&&\quad + \frac{2C_2}{1-\rho_2}\sum\limits_{k=0}^K\sum\limits_{l=0}^k\sum\limits_{t=0}^lw_k\left(1 - \frac{\delta}{2}\right)^{k-l}(1-\rho_2)^{l-t}r_t\notag.
	\end{eqnarray}
	Using this and 
	\begin{eqnarray*}
		w_k\left(1 - \frac{\delta}{2}\right)^{k-l}(1-\rho_2)^{l-t} &\overset{\eqref{eq:ec_sgd_technical_5_new}}{\le}& w_{l}\left(1+\frac{\delta}{4}\right)^{k-l}\left(1 - \frac{\delta}{2}\right)^{k-l}(1-\rho_2)^{l-t}\\
		&\overset{\eqref{eq:1+p/2_inequality},\eqref{eq:ec_sgd_technical_6_new}}{\le}& \left(1-\frac{\delta}{4}\right)^{k-l}\left(1+\frac{\rho_2}{2}\right)^{l-t}(1-\rho_2)^{l-t}w_t\\
		&\overset{\eqref{eq:1+p/2_inequality}}{\le}& \left(1-\frac{\delta}{4}\right)^{k-l}\left(1-\frac{\rho_2}{2}\right)^{l-t}w_t
	\end{eqnarray*}
	we derive
	\begin{eqnarray}
		\sum\limits_{k=0}^K\sum\limits_{l=0}^k w_k\left(1 - \frac{\delta}{2}\right)^{k-l}\EE\sigma_{2,l}^2 &\le& \sum\limits_{k=0}^K\sum\limits_{l=0}^k w_k\left(1 - \frac{\delta}{4}\right)^{k-l}\left(1-\frac{\rho_2}{2}\right)^lw_0\sigma_{2,0}^2\notag\\
		&&\quad + \frac{2C_2}{1-\rho_2}\sum\limits_{k=0}^K\sum\limits_{l=0}^k\sum\limits_{t=0}^l\left(1 - \frac{\delta}{4}\right)^{k-l}\left(1-\frac{\rho_2}{2}\right)^{l-t}w_tr_t\notag\\
		&\le& w_0\sigma_{2,0}^2 \left(\sum\limits_{k=0}^{\infty}\left(1-\frac{\delta}{4}\right)^k\right)\left(\sum\limits_{k=0}^{\infty}\left(1-\frac{\rho_2}{2}\right)^k\right)\notag\\
		&&\quad \frac{2C_2}{1-\rho_2}\left(\sum\limits_{k=0}^K w_kr_k\right)\left(\sum\limits_{k=0}^{\infty}\left(1-\frac{\delta}{4}\right)^k\right)\left(\sum\limits_{k=0}^{\infty}\left(1-\frac{\rho_2}{2}\right)^k\right)\notag\\
		&=& \frac{8\sigma_{2,0}^2}{\delta\rho_2(1-\eta)} + \frac{16C_2}{\delta\rho_2(1-\rho_2)}\sum\limits_{k=0}^K w_kr_k.\label{eq:ec_sgd_technical_8_new}
	\end{eqnarray}
	Similarly, we estimate $\sigma_{1,k}^2$:
	\begin{eqnarray}
		\EE\sigma_{1,k+1}^2 &\overset{\eqref{eq:sigma_k+1_bound_1}}{\le}& (1 - \rho_1)\EE\sigma_{1,k}^2 + 2C_1 \underbrace{\EE\left[f(x^k) - f(x^*)\right]}_{r_k} + G\rho_1\EE\sigma_{2,k}^2 + D_2\notag\\
		&\le& (1-\rho_1)^{k+1}\sigma_{1,0}^2 + 2C_1\sum\limits_{l=0}^{k}(1-\rho_1)^{k-l}r_l + G\rho_1\sum\limits_{l=0}^k(1-\rho_1)^{k-l}\EE\sigma_{2,k}^2 \notag\\
		&&+ D_2\sum\limits_{l=0}^{k}(1-\rho_1)^l\notag\\
		&\le& (1-\rho_1)^{k+1}\sigma_{1,0}^2 + 2C_1\sum\limits_{l=0}^{k}(1-\rho_1)^{k-l}r_l + G\rho_1\sum\limits_{l=0}^k(1-\rho_1)^{k-l}\EE\sigma_{2,k}^2 \notag\\
		&&\quad + D_2\sum\limits_{l=0}^{\infty}(1-\rho_1)^l\notag\\
		&=& (1-\rho_1)^{k+1}\sigma_{1,0}^2 + 2C_1\sum\limits_{l=0}^{k}(1-\rho_1)^{k-l}r_l + G\rho_1\sum\limits_{l=0}^k(1-\rho_1)^{k-l}\EE\sigma_{2,k}^2 \notag\\
		&&\quad + \frac{D_2}{\rho_1}.\label{eq:sigma_1_k_useful_recurrence_new}
	\end{eqnarray}
	Using this we get
	\begin{eqnarray}
		\sum\limits_{k=0}^K\sum\limits_{l=0}^k w_k\left(1 - \frac{\delta}{2}\right)^{k-l}\EE\sigma_{1,l}^2 &\le& \sigma_{1,0}^2\sum\limits_{k=0}^K\sum\limits_{l=0}^kw_k\left(1-\frac{\delta}{2}\right)^{k-l}(1-\rho_1)^{l}\notag\\
		&&\quad + \frac{2C_1}{1-\rho_1}\sum\limits_{k=0}^K\sum\limits_{l=0}^k\sum\limits_{t=0}^{l} w_k\left(1 - \frac{\delta}{2}\right)^{k-l}(1-\rho_1)^{l-t}r_t\notag\\
		&&\quad + \frac{G\rho_1}{1-\rho_1}\sum\limits_{k=0}^K\sum\limits_{l=0}^k\sum\limits_{t=0}^{l} w_k\left(1 - \frac{\delta}{2}\right)^{k-l}(1-\rho_1)^{l-t}\EE\sigma_{2,t}^2\notag\\
		&&\quad + \frac{D_2}{\rho_1}\sum\limits_{k=0}^K\sum\limits_{l=0}^k\sum\limits_{t=0}^{l} w_k\left(1 - \frac{\delta}{2}\right)^{k-l}(1-\rho_1)^{l-t}. \label{eq:ec_sgd_technical_9_new}
	\end{eqnarray}
	Moreover,
	\begin{eqnarray*}
		w_k\left(1 - \frac{\delta}{2}\right)^{k-l}(1-\rho_1)^{l-t} &\overset{\eqref{eq:ec_sgd_technical_5_new}}{\le}& w_{l}\left(1+\frac{\delta}{4}\right)^{k-l}\left(1 - \frac{\delta}{2}\right)^{k-l}(1-\rho_1)^{l-t}\\
		&\overset{\eqref{eq:1+p/2_inequality},\eqref{eq:ec_sgd_technical_6_new}}{\le}& \left(1-\frac{\delta}{4}\right)^{k-l}\left(1+\frac{\rho_1}{2}\right)^{l-t}(1-\rho_1)^{l-t}w_t\\
		&\overset{\eqref{eq:1+p/2_inequality}}{\le}& \left(1-\frac{\delta}{4}\right)^{k-l}\left(1-\frac{\rho_1}{2}\right)^{l-t}w_t,
	\end{eqnarray*}
	hence
	\begin{eqnarray}
		\sum\limits_{k=0}^K\sum\limits_{l=0}^k w_k\left(1 - \frac{\delta}{2}\right)^{k-l}\EE\sigma_{1,l}^2 &\overset{\eqref{eq:ec_sgd_technical_9_new}}{\le}& w_0\sigma_{1,0}^2\sum\limits_{k=0}^K\sum\limits_{l=0}^k\left(1-\frac{\delta}{4}\right)^{k-l}\left(1-\frac{\rho_1}{2}\right)^{l}\notag\\
		&&\quad + \frac{2C_1}{1-\rho_1}\sum\limits_{k=0}^K\sum\limits_{l=0}^k\sum\limits_{t=0}^{l} \left(1 - \frac{\delta}{4}\right)^{k-l}\left(1-\frac{\rho_1}{2}\right)^{l-t}w_tr_t\notag\\
		&&\quad + \frac{G\rho_1}{1-\rho_1}\sum\limits_{k=0}^K\sum\limits_{l=0}^k\sum\limits_{t=0}^{l} \left(1 - \frac{\delta}{4}\right)^{k-l}\left(1-\frac{\rho_1}{2}\right)^{l-t}w_t\EE\sigma_{2,t}^2\notag\\
		&&\quad + \frac{D_2}{\rho_1}\left(\sum\limits_{k=0}^Kw_k\right)\left(\sum\limits_{k=0}^\infty\left(1 - \frac{\delta}{2}\right)^{k}\right)\left(\sum\limits_{k=0}^\infty(1-\rho_1)^{k}\right)\notag\\
		&\le& w_0\sigma_{1,0}^2\left(\sum\limits_{k=0}^\infty\left(1 - \frac{\delta}{4}\right)^{k}\right)\left(\sum\limits_{k=0}^\infty\left(1-\frac{\rho_1}{2}\right)^{k}\right)\notag\\
		&&\quad + \frac{2C_1}{1-\rho_1}\left(\sum\limits_{k=0}^Kw_kr_k\right)\left(\sum\limits_{k=0}^\infty\left(1 - \frac{\delta}{4}\right)^{k}\right)\left(\sum\limits_{k=0}^\infty\left(1-\frac{\rho_1}{2}\right)^{k}\right)\notag\\
		&&\quad+ \frac{G\rho_1}{1-\rho_1}\left(\sum\limits_{k=0}^Kw_k\EE\sigma_{2,k}^2\right)\left(\sum\limits_{k=0}^\infty\left(1 - \frac{\delta}{4}\right)^{k}\right)\left(\sum\limits_{k=0}^\infty\left(1-\frac{\rho_1}{2}\right)^{k}\right)\notag\\
		&&\quad + \frac{2D_2}{\delta\rho_1}W_K\notag\\
		&=& \frac{8\sigma_{1,0}^2}{\delta\rho_1(1-\eta)} + \frac{16C_1}{\delta\rho_1(1-\rho_1)}\sum\limits_{k=0}^K w_kr_k + \frac{8G}{\delta(1-\rho_1)}\sum\limits_{k=0}^Kw_k\EE\sigma_{2,k}^2\notag\\
		&&\quad + \frac{2D_2}{\delta\rho_1}W_K.\label{eq:ec_sgd_technical_10_new}
	\end{eqnarray}
	For the third term in the right-hand side of previous inequality we have
	\begin{eqnarray}
		\frac{8G}{\delta(1-\rho_1)}\sum\limits_{k=0}^Kw_k\EE\sigma_{2,k}^2 &\overset{\eqref{eq:sigma_2_k_useful_recurrence_new}}{\le}& \frac{8G\sigma_{2,0}^2}{\delta(1-\rho_1)}\sum\limits_{k=0}^Kw_k(1-\rho_2)^k\notag\\
		&&\quad + \frac{16GC_2}{\delta(1-\rho_1)(1-\rho_2)}\sum\limits_{k=0}^K\sum\limits_{l=0}^kw_k(1-\rho_2)^{k-l}r_l\notag\\
		&\overset{\eqref{eq:ec_sgd_technical_6_new}}{\le}&\frac{8G\sigma_{2,0}^2w_0}{\delta(1-\rho_1)}\sum\limits_{k=0}^K\left(1+\frac{\rho_2}{2}\right)^k(1-\rho_2)^k\notag\\
		&&\quad + \frac{16GC_2}{\delta(1-\rho_1)(1-\rho_2)}\sum\limits_{k=0}^K\sum\limits_{l=0}^k\left(1+\frac{\rho_2}{2}\right)^{k-l}(1-\rho_2)^{k-l}w_lr_l\notag\\
		&\overset{\eqref{eq:1+p/2_inequality}}{\le}&\frac{8G\sigma_{2,0}^2w_0}{\delta(1-\rho_1)}\sum\limits_{k=0}^\infty\left(1-\frac{\rho_2}{2}\right)^k\notag\\
		&&\quad + \frac{16GC_2}{\delta(1-\rho_1)(1-\rho_2)}\sum\limits_{k=0}^K\sum\limits_{l=0}^k\left(1-\frac{\rho_2}{2}\right)^{k-l}w_lr_l\notag\\
		&\le& \frac{16G\sigma_{2,0}^2w_0}{\delta\rho_2(1-\rho_1)}+ \frac{16GC_2}{\delta(1-\rho_1)(1-\rho_2)}\left(\sum\limits_{k=0}^Kw_kr_k\right)\left(\sum\limits_{k=0}^\infty\left(1-\frac{\rho_2}{2}\right)^k\right)\notag\\
		&=&\frac{16G\sigma_{2,0}^2}{\delta\rho_2(1-\rho_1)(1-\eta)}+ \frac{32GC_2}{\delta\rho_2(1-\rho_1)(1-\rho_2)}\sum\limits_{k=0}^Kw_kr_k\label{eq:ec_sgd_technical_11_new}
	\end{eqnarray}
	Combining inequalities \eqref{eq:ec_sgd_technical_10_new} and \eqref{eq:ec_sgd_technical_11_new} we get
	\begin{eqnarray}
		\sum\limits_{k=0}^K\sum\limits_{l=0}^k w_k\left(1 - \frac{\delta}{2}\right)^{k-l}\EE\sigma_{1,l}^2 &\le& \frac{8\sigma_{1,0}^2}{\delta\rho_1(1-\eta)} + \frac{16}{\delta(1-\rho_1)}\left(\frac{C_1}{\rho_1} + \frac{2GC_2}{\rho_2(1-\rho_2)}\right) \sum\limits_{k=0}^K w_kr_k \notag\\
		&&\quad + \frac{16G\sigma_{2,0}^2}{\delta\rho_2(1-\rho_1)(1-\eta)}+ \frac{2D_2}{\delta\rho_1}W_K\label{eq:ec_sgd_technical_12_new}
	\end{eqnarray}
	Finally, we estimate the last term in the right-hand side of \eqref{eq:ec_sgd_technical_4_new}:
	\begin{eqnarray}
		\sum\limits_{k=0}^K\sum\limits_{l=0}^k w_k\left(1 - \frac{\delta}{2}\right)^{k-l} &\le& \left(\sum\limits_{k=0}^K w_k\right)\left(\sum\limits_{k=0}^\infty \left(1 - \frac{\delta}{2}\right)^k\right) = \frac{2}{\delta}W_K. \label{eq:ec_sgd_technical_13_new}  
	\end{eqnarray}
	Plugging inequalities \eqref{eq:ec_sgd_technical_7_new}, \eqref{eq:ec_sgd_technical_8_new}, \eqref{eq:ec_sgd_technical_12_new}, \eqref{eq:ec_sgd_technical_13_new} and $\frac{1-\delta}{1-\frac{\delta}{2}} \le 1$ in \eqref{eq:ec_sgd_technical_4_new} we obtain
	\begin{eqnarray}
		3L\sum\limits_{k=0}^Kw_k\EE\|e^k\|^2 &\le& \frac{24L\left(\frac{2A}{\delta} + \widetilde{A} + \frac{2}{1-\rho_1}\left(\frac{C_1}{\rho_1}+\frac{2GC_2}{\rho_2(1-\rho_2)}\right)\left(\frac{2B_1}{\delta}+\widetilde{B}_1\right) + \frac{2C_2\left(\frac{2B_2}{\delta}+\widetilde{B}_2\right)}{\rho_2(1-\rho_2)}\right)\gamma^2}{\delta}\sum\limits_{k=0}^K w_k r_k\notag\\
		&&\quad + \frac{24L\gamma^2}{\delta\rho_1(1-\eta)}\left(\frac{2B_1}{\delta}+\widetilde{B}_1\right)\sigma_{1,0}^2\notag\\
		&&\quad + \frac{24L\gamma^2}{\delta\rho_2(1-\eta)}\left(\frac{2G}{1-\rho_1}\left(\frac{2B_1}{\delta}+\tilde{B}_1\right)+\frac{2B_2}{\delta} + \widetilde{B}_2\right)\sigma_{2,0}^2\notag\\
		&&\quad + \frac{6L\gamma^2}{\delta}\left(\frac{D_2}{\rho_1}\left(\frac{2B_1}{\delta}+\widetilde{B}_1\right) + \frac{2D_1}{\delta} + \widetilde{D}_1\right)W_K. \notag
	\end{eqnarray}
	Taking into account that $\gamma \le \sqrt{\frac{\delta}{96L\left(\frac{2A}{\delta} + \widetilde{A} + \frac{2}{1-\rho_1}\left(\frac{C_1}{\rho_1}+\frac{2GC_2}{\rho_2(1-\rho_2)}\right)\left(\frac{2B_1}{\delta}+\widetilde{B}_1\right) + \frac{2C_2\left(\frac{2B_2}{\delta}+\widetilde{B}_2\right)}{\rho_2(1-\rho_2)}\right)}}$, $F_1 = \frac{24L\gamma^2}{\delta\rho_1(1-\eta)}\left(\frac{2B_1}{\delta}+\widetilde{B}_1\right)$, $F_2 = \frac{24L\gamma^2}{\delta\rho_2(1-\eta)}\left(\frac{2G}{1-\rho_1}\left(\frac{2B_1}{\delta}+\tilde{B}_1\right)+\frac{2B_2}{\delta} + \widetilde{B}_2\right)$ and\newline $D_3 = \frac{6L\gamma}{\delta}\left(\frac{D_2}{\rho_1}\left(\frac{2B_1}{\delta}+\widetilde{B}_1\right) + \frac{2D_1}{\delta} + \widetilde{D}_1\right)$ we get
	\begin{eqnarray*}
		3L\sum\limits_{k=0}^Kw_k\EE\|e^k\|^2 &\le& \frac{1}{4}\sum\limits_{k=0}^K w_k r_k + F_1\sigma_{1,0}^2 + F_2\sigma_{2,0}^2 + \gamma D_3.
	\end{eqnarray*}
\end{proof}

\section{{\tt SGD} with Delayed Updates}\label{sec:d_sgd}
In this section we consider the {\tt SGD} with delayed updates ({\tt D-SGD}) \cite{agarwal2011distributed,lian2015asynchronous,feyzmahdavian2016asynchronous,arjevani2018tight,stich2020error}. This method has updates of the form \eqref{eq:x^k+1_update}-\eqref{eq:error_update} with
\begin{eqnarray}
	g^k &=& \frac{1}{n}\sum\limits_{i=1}^ng_i^k\notag\\
	v^k &=& \frac{1}{n}\sum\limits_{i=1}^n v_i^k,\quad v_i^k = \begin{cases}\gamma g_i^{k-\tau},&\text{if } t\ge\tau,\\ 0,&\text{if } t < \tau \end{cases}\label{eq:v^k_def_d_sgd}\\
	e^k &=& \frac{1}{n}\sum\limits_{i=1}^n e_i^k,\quad e_i^{k+1} = e_i^k + \gamma g_i^k - v_i^k = \gamma\sum\limits_{t=1}^\tau g_i^{k+1-t},\label{eq:e^k_def_d_sgd}
\end{eqnarray}
where the summation is performed only for non-negative indices. Moreover, we assume that $e_i^0 = 0$ for $i = 1,\ldots,n$.

For convenience we also introduce new constant:
\begin{equation}
	\hat A = A' + L\tau. \label{eq:d_sgd_hat_A}
\end{equation}
\begin{lemma}\label{lem:d_sgd_key_lemma_new}
	Let Assumptions~\ref{ass:quasi_strong_convexity}~and~\ref{ass:L_smoothness} be satisfied, inequalities \eqref{eq:second_moment_bound_new}, \eqref{eq:sigma_k+1_bound_1} and \eqref{eq:sigma_k+1_bound_2} hold and\footnote{When $\rho_1 = 1$ and $\rho_2=1$ one can always set the parameters in such a way that $B_1 = B_1' = B_2 = B_2' = C_1 = C_2 = 0$, $D_2 = 0$. In this case we assume that $\frac{2B_1'C_1}{\rho_1(1-\rho_1)} = \frac{2B_2'C_2}{\rho_2(1-\rho_2)} = 0$.}
	\begin{equation}
		\gamma \le \min\left\{\frac{1}{2\tau\mu}, \frac{1}{8\sqrt{L\tau\left(\hat A + \frac{2 B_1'C_1}{\rho_1(1-\rho_1)} + \frac{2 B_2'C_2}{\rho_2(1-\rho_2)} + \frac{4B_1'GC_2}{\rho_2(1-\rho_1)(1-\rho_2)}\right)}}\right\}, \label{eq:gamma_condition_d_sgd_new}
	\end{equation}
	where $M_1 = \frac{4B_1'}{3\rho_1}$ and  $M_2 = \frac{4\left(B_2' + \frac{4}{3}G\right)}{3\rho_2}$. Then {\tt D-SGD} satisfies Assumption~\ref{ass:key_assumption_new}, i.e., inequality \eqref{eq:sum_of_errors_bound_new} holds with the following parameters:
	\begin{equation}
		F_1 = \frac{6\gamma^2L B_1' \tau(2+\rho_1)}{\rho_1},\quad F_2=\frac{6\gamma^2\tau L(2+\rho_2)}{\rho_2}\left(\frac{2B_1' G}{1-\rho_1} + B_2'\right), \label{eq:d_sgd_parameters_new_1}
	\end{equation}
	\begin{equation}
		D_3 = 3\gamma\tau L\left( D_1' + \frac{2 B_1' D_2}{\rho_1}\right). \label{eq:d_sgd_parameters_new_2}
	\end{equation}
\end{lemma}
\begin{proof}
	First of all, we derive an upper bound for the second moment of $e_i^{k}$:
	\begin{eqnarray}
		\EE\|e^{k}\|^2 &\overset{\eqref{eq:e^k_def_d_sgd}}{=}& \gamma^2\EE\left[\left\|\sum\limits_{t=1}^\tau g^{k-t}\right\|^2\right] \notag\\
		&\overset{\eqref{eq:lemma14_stich}}{\le}& \gamma^2\tau\sum\limits_{t=1}^\tau\EE\left[\left\|\nabla f(x^{k-t})\right\|^2\right] + \gamma^2\sum\limits_{t=1}^\tau\EE\left[\left\|g^{k-t} - \nabla f(x^{k-t})\right\|^2\right]\notag\\
		&\overset{\eqref{eq:variance_decomposition}}{\le}& \gamma^2\tau\sum\limits_{t=1}^\tau\EE\left[\left\|\nabla f(x^{k-t})\right\|^2\right] + \gamma^2\sum\limits_{t=1}^\tau\EE\left[\left\|g^{k-t}\right\|^2\right]\notag\\
		&\overset{\eqref{eq:second_moment_bound_new},\eqref{eq:L_smoothness_cor_3}}{\le}& 2\gamma^2\underbrace{(A'+L\tau)}_{\hat A}\sum\limits_{t=1}^\tau\EE\left[f(x^{k-t}) - f(x^*)\right] + \gamma^2 B_1'\sum\limits_{t=1}^\tau \EE\sigma_{1,k-t}^2\notag\\
		&&\quad + \gamma^2 B_2'\sum\limits_{t=1}^\tau \EE\sigma_{2,k-t}^2 + \gamma^2\tau D_1' \label{eq:d_sgd_technical_1_new}
	\end{eqnarray}
	which implies
	\begin{eqnarray}
		3L\sum\limits_{k=0}^Kw_k\EE\|e^k\|^2 &\overset{\eqref{eq:d_sgd_technical_1_new}}{\le}& 6\gamma^2L\hat A\sum\limits_{k=0}^K\sum\limits_{t=1}^\tau w_k\EE\left[f(x^{k-t}) - f(x^*)\right]\notag\\
		&&\quad + 3\gamma^2L B_1'\sum\limits_{k=0}^K\sum\limits_{t=1}^\tau w_k\EE\sigma_{1,k-t}^2\notag\\
		&&\quad + 3\gamma^2L B_2'\sum\limits_{k=0}^K\sum\limits_{t=1}^\tau w_k\EE\sigma_{2,k-t}^2 + 3\gamma^2\tau L D_1' W_K \label{eq:d_sgd_technical_2_new}
	\end{eqnarray}
	In the remaining part of the proof we derive upper bounds for four terms in the right-hand side of the previous inequality. First of all, recall that $w_k = (1 - \eta)^{-(k+1)}$ and $\eta = \min\left\{\frac{\gamma\mu}{2}, \frac{\rho_1}{4}, \frac{\rho_2}{4}\right\}$. It implies that for all $0 \le i < k$ and $0\le t \le \tau$ we have
	\begin{eqnarray}
		w_k &=& (1 - \eta)^{-(k-t+1)}\left(1 - \eta\right)^{-t} \overset{\eqref{eq:1-p/2_inequality}}{\le} w_{k-t}\left(1 + 2\eta\right)^{t} \notag\\
		&\le& w_{k-t}\left(1 + \gamma\mu\right)^{t} \overset{\eqref{eq:gamma_condition_d_sgd_new}}{\le} w_{k-t}\left(1 + \frac{1}{2\tau}\right)^t \le w_{k-t}\exp\left(\frac{t}{2\tau}\right) \le 2w_{k-t}, \label{eq:d_sgd_technical_3_new}\\
		w_k &=& \left(1 - \eta\right)^{-(k-j+1)}\left(1 - \eta\right)^{-j} \overset{\eqref{eq:1-p/2_inequality}}{\le} w_{k-j}\left(1 + 2\eta\right)^j \le w_{k-j}\left(1 + \frac{\min\{\rho_1,\rho_2\}}{2}\right)^j. \label{eq:d_sgd_technical_4_new}
	\end{eqnarray}
	For simplicity, we introduce new notation: $r_k \eqdef \EE\left[f(x^k) - f(x^*)\right]$. Using this we get
	\begin{eqnarray}
		\sum\limits_{k=0}^K\sum\limits_{t=1}^\tau w_kr_{k-t} &\overset{\eqref{eq:d_sgd_technical_3_new}}{\le}& \sum\limits_{k=0}^K\sum\limits_{t=1}^\tau 2w_{k-t} r_{k-t} \le 2\tau\sum\limits_{k=0}^K w_kr_k \label{eq:d_sgd_technical_6_new}
	\end{eqnarray}
	Similarly, we estimate the second term in the right-hand side of \eqref{eq:d_sgd_technical_4_new}:
	\begin{eqnarray}
		\sum\limits_{k=0}^K\sum\limits_{t=1}^\tau w_k\EE\sigma_{1,k-t}^2 &\le& \sum\limits_{k=0}^K\sum\limits_{t=1}^\tau 2w_{k-t}\EE\sigma_{1,k-t}^2 \le 2\tau\sum\limits_{k=0}^K w_k \EE\sigma_{1,k}^2\notag\\
		&\overset{\eqref{eq:sigma_1_k_useful_recurrence_new}}{\le}& 2\tau\sigma_{1,0}^2\sum\limits_{k=0}^K w_k(1-\rho_1)^k + \frac{4C_1\tau}{1-\rho_1}\sum\limits_{k=0}^K\sum\limits_{l=0}^k w_k (1-\rho_1)^{k-l}r_l \notag\\
		&&\quad + \frac{2G\rho_1\tau}{1-\rho_1}\sum\limits_{k=0}^K\sum\limits_{l=0}^k w_k (1-\rho_1)^{k-l}\EE\sigma_{2,l}^2 + \frac{2\tau D_2}{\rho}W_K. \label{eq:d_sgd_technical_7_new}
	\end{eqnarray}
	For the first term in the right-hand side of previous inequality we have
	\begin{eqnarray}
		2\tau\sigma_{1,0}^2\sum\limits_{k=0}^K w_k(1-\rho_1)^{k} &\overset{\eqref{eq:d_sgd_technical_4_new}}{\le}& 2\tau\sigma_{1,0}^2\sum\limits_{k=0}^K \left(1 + \frac{\rho_1}{2}\right)^{k+1}(1-\rho_1)^k \notag\\
		&\overset{\eqref{eq:1+p/2_inequality}}{\le}& 2\tau\left(1+\frac{\rho_1}{2}\right)\sigma_{1,0}^2\sum\limits_{k=0}^K\left(1 - \frac{\rho_1}{2}\right)^k \notag\\
		&\le& \tau\left(2+\rho_1\right)\sigma_{1,0}^2\sum\limits_{k=0}^\infty\left(1 - \frac{\rho_1}{2}\right)^k \le \frac{2\tau\left(2 + \rho_1\right)\sigma_{1,0}^2}{\rho_1}.\label{eq:d_sgd_technical_8_new}
	\end{eqnarray}
	The second term in the right-hand side of \eqref{eq:d_sgd_technical_7_new} can be upper bounded in the following way:
	\begin{eqnarray}
		\frac{4C_1\tau}{1-\rho_1}\sum\limits_{k=0}^K\sum\limits_{l=0}^k w_k (1-\rho_1)^{k-l}r_l &\overset{\eqref{eq:d_sgd_technical_4_new}}{\le}& \frac{4C_1\tau}{1-\rho_1}\sum\limits_{k=0}^K\sum\limits_{l=0}^k w_l r_l \left(1 + \frac{\rho_1}{2}\right)^{k-l}(1-\rho_1)^{k-l}\notag\\
		&\overset{\eqref{eq:1+p/2_inequality}}{\le}& \frac{4C_1\tau}{1-\rho_1}\sum\limits_{k=0}^K\sum\limits_{l=0}^k w_l r_l \left(1 - \frac{\rho_1}{2}\right)^{k-l}\notag\\
		&\le& \frac{4C_1\tau}{1-\rho_1}\left(\sum\limits_{k=0}^K w_k r_k\right)\left(\sum\limits_{k=0}^\infty\left(1 - \frac{\rho_1}{2}\right)^k\right)\notag\\
		&\le& \frac{8C_1\tau}{\rho_1(1-\rho_1)}\sum\limits_{k=0}^K w_k r_k.\label{eq:d_sgd_technical_9_new}
	\end{eqnarray}
	Repeating similar steps we estimate the third term in the right-hand side of \eqref{eq:d_sgd_technical_7_new}:
	\begin{eqnarray}
		\frac{2G\rho_1\tau}{1-\rho_1}\sum\limits_{k=0}^K\sum\limits_{l=0}^k w_k (1-\rho_1)^{k-l}\EE\sigma_{2,l}^2 &\le& \frac{4G\tau}{1-\rho_1}\sum\limits_{k=0}^Kw_k\EE\sigma_{2,k}^2\notag\\
		&\overset{\eqref{eq:sigma_2_k_useful_recurrence_new}}{\le}& \frac{4G\tau\sigma_{2,0}^2}{1-\rho_1}\sum\limits_{k=0}^Kw_k(1-\rho_2)^k\notag\\
		&& + \frac{8GC_2}{(1-\rho_1)(1-\rho_2)}\sum\limits_{k=0}^K\sum\limits_{l=0}^kw_k(1-\rho_2)^{k-l}r_l\notag\\
		&\overset{\eqref{eq:d_sgd_technical_4_new}}{\le}& \frac{4G\tau\sigma_{2,0}^2}{1-\rho_1}\sum\limits_{k=0}^K\left(1+\frac{\rho_2}{2}\right)^{k+1}(1-\rho_2)^k\notag\\
		&&\hspace{-2cm}+ \frac{8GC_2\tau}{(1-\rho_1)(1-\rho_2)}\sum\limits_{k=0}^K\sum\limits_{l=0}^k\left(1+\frac{\rho_2}{2}\right)^{k-l}(1-\rho_2)^{k-l}w_lr_l\notag\\
		&\overset{\eqref{eq:1+p/2_inequality}}{\le}& \frac{2G\tau(2+\rho_2)\sigma_{2,0}^2}{1-\rho_1}\sum\limits_{k=0}^\infty\left(1-\frac{\rho_2}{2}\right)^{k}\notag\\
		&&+ \frac{8GC_2\tau}{(1-\rho_1)(1-\rho_2)}\sum\limits_{k=0}^K\sum\limits_{l=0}^k\left(1-\frac{\rho_2}{2}\right)^{k-l}w_lr_l\notag\\
		&\le& \frac{4G\tau(2+\rho_2)\sigma_{2,0}^2}{\rho_2(1-\rho_1)}\notag\\
		&&+\frac{8GC_2\tau}{(1-\rho_1)(1-\rho_2)}\left(\sum\limits_{k=0}^Kw_kr_k\right)\left(\sum\limits_{k=0}^\infty\left(1-\frac{\rho_2}{2}\right)^k\right)\notag\\
		&=& \frac{4G\tau(2+\rho_2)\sigma_{2,0}^2}{\rho_2(1-\rho_1)} \notag\\
		&&\quad+ \frac{16GC_2\tau}{\rho_2(1-\rho_1)(1-\rho_2)}\sum\limits_{k=0}^Kw_kr_k\label{eq:d_sgd_technical_9_1_new}
	\end{eqnarray}
	Combining inequalities \eqref{eq:d_sgd_technical_7_new}, \eqref{eq:d_sgd_technical_8_new}, \eqref{eq:d_sgd_technical_9_new} and \eqref{eq:d_sgd_technical_9_1_new} we get
	\begin{eqnarray}
		\sum\limits_{k=0}^K\sum\limits_{t=1}^\tau w_k\EE\sigma_{1,k-t}^2 &\le& \frac{2\tau\left(2 + \rho_1\right)\sigma_{1,0}^2}{\rho_1} + \frac{8\tau}{1-\rho_1}\left(\frac{C_1}{\rho_1}+\frac{2GC_2}{\rho_2(1-\rho_2)}\right)\sum\limits_{k=0}^K w_k r_k \notag\\
		&&\quad + \frac{4G\tau(2+\rho_2)\sigma_{2,0}^2}{\rho_2(1-\rho_1)} + \frac{2\tau D_2}{\rho} W_K.\label{eq:d_sgd_technical_10_new}
	\end{eqnarray}
	Next, we derive
	\begin{eqnarray}
		\sum\limits_{k=0}^K\sum\limits_{t=1}^\tau w_k\EE\sigma_{2,k-t}^2 &\le& \sum\limits_{k=0}^K\sum\limits_{t=1}^\tau 2w_{k-t}\EE\sigma_{2,k-t}^2 \le 2\tau\sum\limits_{k=0}^K w_k \EE\sigma_{2,k}^2\notag\\
		&\overset{\eqref{eq:sigma_2_k_useful_recurrence_new}}{\le}& 2\tau\sigma_{2,0}^2\sum\limits_{k=0}^K w_k(1-\rho_1)^k \notag\\
		&&\quad+ \frac{4C_2\tau}{1-\rho_2}\sum\limits_{k=0}^K\sum\limits_{l=0}^k w_k (1-\rho_2)^{k-l}r_l.\label{eq:d_sgd_technical_11_new}
	\end{eqnarray}
	For the first term in the right-hand side of previous inequality we have
	\begin{eqnarray}
		2\tau\sigma_{2,0}^2\sum\limits_{k=0}^K w_k(1-\rho_2)^{k} &\overset{\eqref{eq:d_sgd_technical_4_new}}{\le}& 2\tau\sigma_{2,0}^2\sum\limits_{k=0}^K \left(1 + \frac{\rho_2}{2}\right)^{k+1}(1-\rho_2)^k \notag\\
		&\overset{\eqref{eq:1+p/2_inequality}}{\le}& 2\tau\left(1+\frac{\rho_2}{2}\right)\sigma_{2,0}^2\sum\limits_{k=0}^K\left(1 - \frac{\rho_2}{2}\right)^k \notag\\
		&\le& \tau\left(2+\rho_2\right)\sigma_{2,0}^2\sum\limits_{k=0}^\infty\left(1 - \frac{\rho_2}{2}\right)^k \le \frac{2\tau\left(2 + \rho_2\right)\sigma_{2,0}^2}{\rho_2}.\notag
	\end{eqnarray}
	The second term in the right-hand side of \eqref{eq:d_sgd_technical_11_new} can be upper bounded in the following way:
	\begin{eqnarray}
		\frac{4C_2\tau}{1-\rho_2}\sum\limits_{k=0}^K\sum\limits_{l=0}^k w_k (1-\rho_2)^{k-l}r_l &\overset{\eqref{eq:d_sgd_technical_4_new}}{\le}& \frac{4C_2\tau}{1-\rho_2}\sum\limits_{k=0}^K\sum\limits_{l=0}^k w_l r_l \left(1 + \frac{\rho_2}{2}\right)^{k-l}(1-\rho_2)^{k-l}\notag\\
		&\overset{\eqref{eq:1+p/2_inequality}}{\le}& \frac{4C_2\tau}{1-\rho_2}\sum\limits_{k=0}^K\sum\limits_{l=0}^k w_l r_l \left(1 - \frac{\rho_2}{2}\right)^{k-l}\notag\\
		&\le& \frac{4C_2\tau}{1-\rho_2}\left(\sum\limits_{k=0}^K w_k r_k\right)\left(\sum\limits_{k=0}^\infty\left(1 - \frac{\rho_2}{2}\right)^k\right)\notag\\
		&\le& \frac{8C_2\tau}{\rho_2(1-\rho_2)}\sum\limits_{k=0}^K w_k r_k,\notag
	\end{eqnarray}	
	hence
	\begin{eqnarray}
		\sum\limits_{k=0}^K\sum\limits_{t=1}^\tau w_k\EE\sigma_{2,k-t}^2 &\overset{\eqref{eq:d_sgd_technical_11_new}}{\le}& \frac{2\tau\left(2 + \rho_2\right)\sigma_{2,0}^2}{\rho_2} + \frac{8C_2\tau}{\rho_2(1-\rho_2)}\sum\limits_{k=0}^K w_k r_k. \label{eq:d_sgd_technical_12_new}
	\end{eqnarray}
	Plugging inequalities \eqref{eq:d_sgd_technical_6_new}, \eqref{eq:d_sgd_technical_10_new} and \eqref{eq:d_sgd_technical_12_new} in \eqref{eq:d_sgd_technical_2_new} we obtain
	\begin{eqnarray}
		3L\sum\limits_{k=0}^Kw_k\EE\|e^k\|^2 &\le& 12\gamma^2L\tau\left(\hat A + \frac{2 B_1'C_1}{\rho_1(1-\rho_1)} + \frac{2 B_2'C_2}{\rho_2(1-\rho_2)} + \frac{4B_1'GC_2}{\rho_2(1-\rho_1)(1-\rho_2)}\right)\sum\limits_{k=0}^K w_k r_k\notag\\
		&&\quad + \frac{6\gamma^2L B_1' \tau(2+\rho_1)}{\rho_1}\sigma_0^2 + \frac{6\gamma^2\tau L(2+\rho_2)}{\rho_2}\left(\frac{2B_1' G}{1-\rho_1} + B_2'\right)\sigma_{2,0}^2\notag\\
		&&\quad + 3\gamma^2\tau L\left(D_1' + \frac{2 B_1' D_2}{\rho}\right)W_K. \notag
	\end{eqnarray}
	Taking into account that $\gamma \le \frac{1}{4\sqrt{4L\tau\left(\hat A + \frac{2 B_1'C_1}{\rho_1(1-\rho_1)} + \frac{2 B_2'C_2}{\rho_2(1-\rho_2)} + \frac{4B_1'GC_2}{\rho_2(1-\rho_1)(1-\rho_2)}\right)}}$, $F_1 = \frac{6\gamma^2L B_1' \tau(2+\rho_1)}{\rho_1}$, $F_2=\frac{6\gamma^2\tau L}{\rho_2}\left(\frac{2B_1' G(2+\rho_2)}{1-\rho_1} + B_2'\right)$ and $D_3 = 3\gamma\tau L\left( D_1' + \frac{2 B_1' D_2}{\rho}\right)$ we get
	\begin{eqnarray*}
		3L\sum\limits_{k=0}^Kw_k\EE\|e^k\|^2 &\le& \frac{1}{4}\sum\limits_{k=0}^K w_k r_k + F_1\sigma_{1,0}^2 + F_2\sigma_{2,0}^2 + \gamma D_3.
	\end{eqnarray*}
\end{proof}

As a direct application of Lemma~\ref{lem:d_sgd_key_lemma_new} and Theorem~\ref{thm:main_result_new} we get the following result.
\begin{theorem}\label{thm:d_sgd_main_result_new}
	Let Assumptions~\ref{ass:quasi_strong_convexity}~and~\ref{ass:L_smoothness} be satisfied, inequalities \eqref{eq:second_moment_bound_new}, \eqref{eq:sigma_k+1_bound_1} and \eqref{eq:sigma_k+1_bound_2} hold and
	\begin{equation*}
		\gamma \le \min\left\{\frac{1}{4(A'+C_1M_1 + C_2M_2)},\frac{1}{2\tau\mu}, \frac{1}{8\sqrt{L\tau\left(\hat A + \frac{2 B_1'C_1}{\rho_1(1-\rho_1)} + \frac{2 B_2'C_2}{\rho_2(1-\rho_2)} + \frac{4B_1'GC_2}{\rho_2(1-\rho_1)(1-\rho_2)}\right)}}\right\},
	\end{equation*}
	where $M_1 = \frac{4B_1'}{3\rho_1}$ and  $M_2 = \frac{4\left(B_2' + \frac{4}{3}G\right)}{3\rho_2}$. Then for all $K\ge 0$ we have
	\begin{equation*}
		\EE\left[f(\bar x^K) - f(x^*)\right] \le \left(1 - \eta\right)^K\frac{4(T^0 + \gamma F_1 \sigma_{1,0}^2+ \gamma F_2 \sigma_{2,0}^2)}{\gamma} + 4\gamma\left(D_1' + MD_2 + D_3\right) 
	\end{equation*}
	when $\mu > 0$ and
	\begin{equation*}
		\EE\left[f(\bar x^K) - f(x^*)\right] \le \frac{4(T^0 + \gamma F_1 \sigma_{1,0}^2+ \gamma F_2 \sigma_{2,0}^2)}{\gamma K} + 4\gamma\left(D_1' + MD_2 + D_3\right) 
	\end{equation*}
	when $\mu = 0$, where $\eta = \min\left\{\nicefrac{\gamma\mu}{2},\nicefrac{\rho_1}{4},\nicefrac{\rho_2}{4}\right\}$, $T^k \eqdef \|\tx^k - x^*\|^2 + M_1\gamma^2 \sigma_{1,k}^2 + M_2\gamma^2 \sigma_{2,k}^2$ and 
	\begin{equation*}
		F_1 = \frac{6\gamma^2L B_1' \tau(2+\rho_1)}{\rho_1},\quad F_2=\frac{6\gamma^2\tau L(2+\rho_2)}{\rho_2}\left(\frac{2B_1' G}{1-\rho_1} + B_2'\right),
	\end{equation*}
	\begin{equation*}
		D_3 = 3\gamma\tau L\left( D_1' + \frac{2 B_1' D_2}{\rho_1}\right).
	\end{equation*}
\end{theorem}


\section{Special Cases: Delayed Updates Methods}\label{sec:special_cases2}

\begin{table*}[!t]
\caption{Complexity of {\tt SGD} methods with delayed updates established in this chapter. Symbols: $\varepsilon = $ error tolerance; $\delta = $ contraction factor of compressor $\cC$; $\omega = $ variance parameter of compressor $\cQ$; $\kappa = \nicefrac{L}{\mu}$; $\cL =$ expected smoothness constant; $\sigma_*^2 = $ variance of the stochastic gradients in the solution; $\zeta_*^2 =$ average of $\|\nabla f_i(x^*)\|^2$; $\sigma^2 =$ average of the uniform bounds for the variances of stochastic gradients of workers; $\cM_{2,q} = (\omega+1)\sigma^2 + \omega\zeta_*^2$; $\sigma^2_q = (1+\omega)\left(1+\frac{\omega}{n}\right)\sigma^2$. $^\dagger${\tt D-QGDstar} is a special case of {\tt D-QSGDstar} where each worker $i$ computes the full gradient  $\nabla f_i(x^k)$; $^\ddagger${\tt D-GD-DIANA} is a special case of {\tt D-SGD-DIANA} where each worker $i$ computes the full gradient  $\nabla f_i(x^k)$.
}
\label{tbl:special_cases_delayed_methods}
\begin{center}
\footnotesize
\begin{tabular}{|c|l|c|c|c|c|}
\hline
\bf Problem & \bf Method &   \bf Alg \# &  \bf Citation &  \bf  Sec \#  
& \bf Rate (constants ignored)\\
\hline
\eqref{eq:main_problem_ef}+\eqref{eq:f_i_sum_ef} & {\tt D-SGDsr}  & Alg \ref{alg:d-SGDsr} & {\color{red}\bf new} & \ref{sec:d_SGDsr} 
& {\color{red}$\widetilde{\cO}\left(\frac{\cL + \sqrt{L^2\tau^2 + L\cL\tau}}{\mu} + \frac{\sigma_*^2}{n\mu\varepsilon} + \frac{\sqrt{L\tau \sigma_*^2}}{\mu\sqrt{n\varepsilon}}\right)$}\\
\hline
\eqref{eq:main_problem_ef}+\eqref{eq:f_i_expectation_ef} & {\tt D-SGD}  & Alg \ref{alg:d-sgd} & {\cite{stich2020error}}  & \ref{sec:d_sgd_pure} 
& $\widetilde{\cO}\left(\tau\kappa + \frac{\sigma_*^2}{n\mu\varepsilon} + \frac{\sqrt{L\tau\sigma_*^2}}{\mu\sqrt{n\varepsilon}}\right)$\\
\hline
\eqref{eq:main_problem_ef}+\eqref{eq:f_i_expectation_ef} & {\tt D-QSGD}  & Alg \ref{alg:d-qsgd} & {\color{red}\bf new}
& \ref{sec:d_qsgd} 
& {\color{red} $\widetilde{\cO}\left( \kappa\left(\tau + \frac{\omega}{n}\right) + \frac{\cM_{2,q}}{n\mu\varepsilon} + \frac{\sqrt{L\tau\cM_{2,q}}}{\mu\sqrt{n\varepsilon}} \right)$ } \\
\hline
\eqref{eq:main_problem_ef}+\eqref{eq:f_i_expectation_ef} & {\tt D-QSGDstar}  & Alg \ref{alg:d-qSGDstar} & {\color{red}\bf new}
& \ref{sec:d_qsgd_star} 
& {\color{red} $\widetilde{\cO}\left( \kappa\left(\tau + \frac{\omega}{n}\right) + \frac{\sigma^2}{n\mu\varepsilon} + \frac{\sqrt{L\tau\sigma^2}}{\mu\sqrt{n\varepsilon}} \right)$ } \\
\hline
\eqref{eq:main_problem_ef}+\eqref{eq:f_i_expectation_ef} & {\tt D-QGDstar}$^\dagger$  & Alg \ref{alg:d-qSGDstar} & {\color{red}\bf new}
& \ref{sec:d_qsgd_star} 
& {\color{red} $\cO\left( \kappa\left(\tau + \frac{\omega}{n}\right)\log\frac{1}{\varepsilon}\right)$ } \\
\hline
\eqref{eq:main_problem_ef}+\eqref{eq:f_i_expectation_ef} & {\tt D-SGD-DIANA}  & Alg \ref{alg:d-diana} & {\color{red}\bf new}
& \ref{sec:d_diana} 
& {\color{red}$\widetilde{\cO}\left(\omega +\kappa\left(\tau + \frac{\omega}{n}\right) + \frac{\sigma^2}{n\mu\varepsilon} + \frac{\sqrt{L\tau\sigma_q^2}}{\mu\sqrt{n\varepsilon}}\right)$}\\
\hline
\eqref{eq:main_problem_ef}+\eqref{eq:f_i_expectation_ef} & {\tt D-GD-DIANA}$^\ddagger$  & Alg \ref{alg:d-diana} & {\color{red}\bf new}
& \ref{sec:d_diana} 
& {\color{red}$\cO\left(\left(\omega + \kappa\left(\tau + \frac{\omega}{n}\right)\right) \log \frac{1}{\varepsilon}\right)$} \\
\hline
\eqref{eq:main_problem_ef}+\eqref{eq:f_i_sum_ef} & {\tt D-LSVRG}  & Alg \ref{alg:d-LSVRG} & {\color{red}\bf new}
& \ref{sec:d_LSVRG} 
& {\color{red}$\cO\left(\left(m + \kappa\tau\right)\log\frac{1}{\varepsilon}\right)$}\\
\hline
\eqref{eq:main_problem_ef}+\eqref{eq:f_i_sum_ef} & {\tt D-QLSVRG}  & Alg \ref{alg:d-qLSVRG} & {\color{red}\bf new}
& \ref{sec:d_qLSVRG} 
& {\color{red}$\widetilde{\cO}\left(m + \kappa\left(\tau+\frac{\omega}{n}\right) + \frac{\zeta_*^2}{n\mu\varepsilon} + \frac{\sqrt{L\tau\zeta_*^2}}{\mu\sqrt{n\varepsilon}} \right)$}\\
\hline
\eqref{eq:main_problem_ef}+\eqref{eq:f_i_sum_ef} & {\tt D-QLSVRGstar}  & Alg \ref{alg:d-qLSVRGstar} & {\color{red}\bf new}
 & \ref{sec:d_qLSVRGstar} 
 & {\color{red}$\cO\left(\left(m + \kappa\left(\tau+\frac{\omega}{n}\right)\right)\log\frac{1}{\varepsilon}\right)$} \\
\hline
\eqref{eq:main_problem_ef}+\eqref{eq:f_i_sum_ef} & {\tt D-LSVRG-DIANA}  & Alg \ref{alg:d-LSVRG-diana} & {\color{red}\bf new}
& \ref{sec:d_LSVRG-diana} 
& {\color{red}$\cO \left(\left( \omega + m + \kappa\left(\tau+\frac{\omega}{n}\right) \right) \log \frac{1}{\varepsilon}\right)$} \\
\hline
\end{tabular}
\end{center}
\end{table*}

\subsection{{\tt D-SGD}}\label{sec:d_sgd_pure}
In this section we consider the same setup as in Section~\ref{sec:ec_sgd_pure}.
\begin{algorithm}[t]
   \caption{{\tt D-SGD}}\label{alg:d-sgd}
\begin{algorithmic}[1]
   \Require learning rate $\gamma>0$, initial vector $x^0 \in \R^d$
	\State Set $e_i^0 = 0$ for all $i=1,\ldots, n$   
   \For{$k=0,1,\dotsc$}
       \State Broadcast $x^{k}$ to all workers
        \For{$i=1,\dotsc,n$ in parallel}
            \State Sample $g^{k}_i = \nabla f_{\xi_i}(x^k) - \nabla f_i(x^*)$
            \State $v_i^k = \begin{cases}\gamma g_i^{k-\tau},& \text{if } k \ge \tau,\\ 0,& \text{if } k < \tau \end{cases}$
            \State $e_i^{k+1} = e_i^k + \gamma g_i^k - v_i^k$
        \EndFor
        \State $e^k = \frac{1}{n}\sum_{i=1}^n e_i^k$, $g^k = \frac{1}{n}\sum_{i=1}^ng_i^k$, $v^k = \frac{1}{n}\sum_{i=1}^nv_i^k = \frac{1}{n}\sum_{i=1}^n \nabla f_{\xi_i}(x^{k-\tau})$
       \State $x^{k+1} = x^k - v^k$
   \EndFor
\end{algorithmic}
\end{algorithm}
We notice that vectors $e_i^k$ appear only in the analysis and there is no need to compute them. Moreover, we use $\nabla f_i(x^*)$ in the definition of $g_i^k$ which is problematic at the firt glance. Indeed, workers do not know $\nabla f_i(x^*)$. However, since $0 = \nabla f(x^*) = \frac{1}{n}\nabla f_i(x^*)$ and master node uses averages of $g_i^k$ for the updates one can ignore $\nabla f_i(x^*)$ in $g_i^k$ in the implementation of {\tt D-SGD} and get exactly the same method. We define $g_i^k$ in such a way only for the theoretical analysis.
\begin{lemma}[see also Lemmas 1,2 from~\cite{nguyen2018sgd}]\label{lem:lemma_d_sgd}
    Assume that $f_{\xi_i}(x)$ are convex in $x$ for every $\xi_i$, $i=1,\ldots,n$. Then for every $x\in\R^d$ and $i=1,\ldots, n$
    \begin{equation}\label{eq:lemma_d_sgd_1}
        \EE\left[\|g^k\|^2\mid x^k\right] \le 4L(f(x^k) - f(x^*)) + \frac{2}{n^2}\sum\limits_{i=1}^n \Var\left[\nabla f_{\xi_i}(x^*)\right].
    \end{equation}
    If further $f(x)$ is $\mu$-quasi strongly convex with possibly non-convex $f_i,f_{\xi_i}$ and $\mu > 0$, then for every $x\in\R^d$ and $i = 1,\ldots, n$
    \begin{equation}\label{eq:lemma_d_sgd_2}
        \EE\left[\|g^k\|^2\mid x^k\right] \le 4L\kappa(f(x^k) - f(x^*)) + \frac{2}{n^2}\sum\limits_{i=1}^n \Var\left[\nabla f_{\xi_i}(x^*)\right],
    \end{equation}
    where $\kappa = \frac{L}{\mu}$.
\end{lemma}
\begin{proof}
	By definition of $g^k$ we have
	\begin{eqnarray}
		\EE\left[\|g^k\|^2\mid x^k\right] &=& \EE\left[\left\|\frac{1}{n}\sum\limits_{i=1}^n\left(\nabla f_{\xi_i}(x^k) - \nabla f_{\xi_i}(x^*) + \nabla f_{\xi_i}(x^*) - \nabla f_i(x^*)\right)\right\|^2\mid x^k\right]\notag \\
		&\overset{\eqref{eq:a_b_norm_squared}}{\le}& 2\EE\left[\left\|\frac{1}{n}\sum\limits_{i=1}^n\left(\nabla f_{\xi_i}(x^k) - \nabla f_{\xi_i}(x^*)\right)\right\|^2\mid x^k\right]\notag\\
		&&\quad + 2\underbrace{\EE\left[\left\|\frac{1}{n}\sum\limits_{i=1}^n\left(\nabla f_{\xi_i}(x^*) - \nabla f_i(x^*)\right)\right\|^2\right]}_{\Var\left[\frac{1}{n}\sum\limits_{i=1}^n\nabla f_{\xi_i}(x^*)\right]}\notag\\
		&\overset{\eqref{eq:a_b_norm_squared}}{\le}& \frac{2}{n}\sum\limits_{i=1}^n\EE\left[\|\nabla f_{\xi_i}(x^k)-\nabla f_{\xi_i}(x^*)\|^2\mid x^k\right]\notag\\
		&&\quad + \frac{2}{n^2}\sum\limits_{i=1}^n \underbrace{\EE\left[\|\nabla f_{\xi_i}(x^*) - \nabla f_i(x^*)\|^2\right]}_{\Var\left[\nabla f_{\xi_i}(x^*)\right]}, \label{eq:d_sgd_pure_technical_1}
	\end{eqnarray}
	where in the last inequality we use independence of $\nabla f_{\xi_i}(x^*)$, $i=1,\ldots,n$. Using this we derive inequality \eqref{eq:lemma_d_sgd_1}:
	\begin{eqnarray*}
		\EE\left[\|g^k\|^2\mid x^k\right] &\overset{\eqref{eq:d_sgd_pure_technical_1},\eqref{eq:L_smoothness_cor_3}}{\le}& \frac{4L}{n}\sum\limits_{i=1}^n \EE\left[D_{f_{\xi_i}}(x^k,x^*)\mid x^k\right] + \frac{2}{n^2}\sum\limits_{i=1}^n\Var\left[\nabla f_{\xi_i}(x^*)\right]\\
		&=& \frac{4L}{n}\sum\limits_{i=1}^n D_{f_i}(x^k,x^*) + \frac{2}{n^2}\sum\limits_{i=1}^n\Var\left[\nabla f_{\xi_i}(x^*)\right]\\
		&=& 4L\left(f(x^k) - f(x^*)\right) + \frac{2}{n^2}\sum\limits_{i=1}^n\Var\left[\nabla f_{\xi_i}(x^*)\right].
	\end{eqnarray*}
	Next, if $f(x)$ is $\mu$-quasi strongly convex, but $f_i,f_{\xi_i}$ are not necessary convex, we obtain
	\begin{eqnarray*}
		\EE\left[\|g^k\|^2\mid x^k\right] &\overset{\eqref{eq:d_sgd_pure_technical_1},\eqref{eq:L_smoothness_def}}{\le}& \frac{2L^2}{n}\sum\limits_{i=1}^n \|x^k - x^*\|^2 + \frac{2}{n^2}\sum\limits_{i=1}^n\Var\left[\nabla f_{\xi_i}(x^*)\right]\\
		&\overset{\eqref{eq:str_quasi_cvx}}{\le}& \frac{4L^2}{\mu}\left(f(x^k) - f(x^*)\right) + \frac{2}{n^2}\sum\limits_{i=1}^n\Var\left[\nabla f_{\xi_i}(x^*)\right].
	\end{eqnarray*}
\end{proof}

\begin{theorem}\label{thm:d_sgd_pure}
	Assume that $f_\xi(x)$ is convex in $x$ for every $\xi$. Then {\tt D-SGD} satisfies Assumption~\ref{ass:key_assumption_new} with
	\begin{gather*}
		A' = 2L,\quad B_1' = B_2' = 0,\quad D_1' = \frac{2}{n^2}\sum\limits_{i=1}^n\Var\left[\nabla f_{\xi_i}(x^*)\right],\quad \sigma_{1,k}^2 \equiv \sigma_{2,k}^2 \equiv 0\\
		\rho_1 = \rho_2 = 1,\quad C_1 = C_2 = 0,\quad D_2 = 0\\
		F_1 = F_2 = 0,\quad D_3 = \frac{6\gamma\tau L}{n^2}\sum\limits_{i=1}^n\Var\left[\nabla f_{\xi_i}(x^*)\right]	
	\end{gather*}
	with $\gamma$ satisfying
	\begin{equation*}
		\gamma \le \frac{1}{8L\sqrt{2\tau\left(\tau + 2\right)}}
	\end{equation*}
	and for all $K \ge 0$
	\begin{equation*}
		\EE\left[f(\bar{x}^K) - f(x^*)\right] \le \left(1 - \frac{\gamma\mu}{2}\right)^K\frac{4\|x^0 - x^*\|^2}{\gamma} + \frac{8\gamma}{n^2}\left(1 + 3L\gamma\tau\right)\sum\limits_{i=1}^n\Var\left[\nabla f_{\xi_i}(x^*)\right]
	\end{equation*}
	when $\mu > 0$ and
	\begin{equation*}
		\EE\left[f(\bar{x}^K) - f(x^*)\right] \le \frac{4\|x^0 - x^*\|^2}{\gamma K} + \frac{8\gamma}{n^2}\left(1 + 3L\gamma\tau\right)\sum\limits_{i=1}^n\Var\left[\nabla f_{\xi_i}(x^*)\right]
	\end{equation*}
	when $\mu = 0$. If further $f_i(x)$ are $\mu$-strongly convex with possibly non-convex $f_{\xi_i}$ and $\mu > 0$, then
	{\tt D-SGD} satisfies Assumption~\ref{ass:key_assumption_new} with
	\begin{gather*}
		A' = 2\kappa L,\quad B_1' = B_2' = 0,\quad D_1' = \frac{2}{n^2}\sum\limits_{i=1}^n\Var\left[\nabla f_{\xi_i}(x^*)\right],\quad \sigma_{1,k}^2 \equiv \sigma_{2,k}^2 \equiv 0,\\
		\rho_1 = \rho_2 = 1,\quad C_1 = C_2 = 0,\quad D_2 = 0,\quad G = 0,\\
		F_1 = F_2 = 0,\quad D_3 = \frac{6\gamma\tau L}{n^2}\sum\limits_{i=1}^n\Var\left[\nabla f_{\xi_i}(x^*)\right]	
	\end{gather*}
	with $\gamma$ satisfying
	\begin{equation*}
		\gamma \le \min\left\{\frac{1}{8\kappa L}, \frac{1}{8L\sqrt{2\tau\left(\tau + 2\kappa\right)}}\right\}
	\end{equation*}
	and for all $K \ge 0$
	\begin{equation*}
		\EE\left[f(\bar{x}^K) - f(x^*)\right] \le \left(1 - \frac{\gamma\mu}{2}\right)^K\frac{4\|x^0 - x^*\|^2}{\gamma} + \frac{8\gamma}{n^2}\left(1 + 3L\gamma\tau\right)\sum\limits_{i=1}^n\Var\left[\nabla f_{\xi_i}(x^*)\right].
	\end{equation*}
\end{theorem}

In other words, {\tt D-SGD} converges with linear rate $\cO\left(\tau\kappa\ln\frac{1}{\varepsilon}\right)$ to the neighbourhood of the solution when $\mu > 0$. Applying Lemma~\ref{lem:lemma2_stich} we establish the rate of convergence to $\varepsilon$-solution.
\begin{corollary}\label{cor:d_SGD_pure_str_cvx_cor}
	Let the assumptions of Theorem~\ref{thm:d_sgd_pure} hold, $f_{\xi}(x)$ are convex for each $\xi$ and $\mu > 0$. Then after $K$ iterations of {\tt D-SGD} with the stepsize
	\begin{equation*}
		\gamma = \min\left\{\frac{1}{8L\sqrt{2\tau\left(\tau + 2\right)}}, \frac{\ln\left(\max\left\{2,\min\left\{\frac{\|x^0-x^*\|^2\mu^2K^2}{D_1'}, \frac{\|x^0-x^*\|^2\mu^3K^3}{3\tau LD_1}\right\}\right\}\right)}{\mu K}\right\}
	\end{equation*}	 
	we have
	\begin{equation*}
		\EE\left[f(\bar{x}^K) - f(x^*)\right] = \widetilde\cO\left(L\tau\|x^0 - x^*\|^2\exp\left(-\frac{\mu}{\tau L}K\right) + \frac{D_1'}{\mu K} + \frac{L\tau D_1'}{\mu^2 K^2}\right).
	\end{equation*}
	That is, to achive $\EE\left[f(\bar{x}^K) - f(x^*)\right] \le \varepsilon$ {\tt D-SGD} requires
	\begin{equation*}
		\widetilde{\cO}\left(\frac{\tau L}{\mu} + \frac{D_1'}{\mu\varepsilon} + \frac{\sqrt{L\tau D_1'}}{\mu\sqrt{\varepsilon}}\right) \text{ iterations.}
	\end{equation*}
\end{corollary}
\begin{corollary}\label{cor:d_SGD_pure_str_cvx_cor_2}
	Let the assumptions of Theorem~\ref{thm:d_sgd_pure} hold and $f(x)$ is $\mu$-strongly convex with $\mu > 0$ and possibly non-convex $f_i,f_{\xi_i}$. Then after $K$ iterations of {\tt D-SGD} with the stepsize
	\begin{equation*}
		\gamma = \min\left\{\frac{1}{8\kappa L}, \frac{1}{8L\sqrt{2\tau\left(\tau + 2\kappa\right)}}, \frac{\ln\left(\max\left\{2,\min\left\{\frac{\|x^0-x^*\|^2\mu^2K^2}{D_1'}, \frac{\|x^0-x^*\|^2\mu^3K^3}{L\tau D_1'}\right\}\right\}\right)}{\mu K}\right\}
	\end{equation*}	 
	we have $\EE\left[f(\bar{x}^K) - f(x^*)\right]$ of order
	\begin{equation*}
		\widetilde\cO\left(L\left(\kappa+\tau\sqrt{\kappa}\right)\|x^0 - x^*\|^2\exp\left(-\min\left\{\frac{\mu}{\tau L\sqrt{\kappa}}, \frac{1}{\kappa^2}\right\}K\right) + \frac{D_1'}{\mu K} + \frac{L\tau D_1'}{\mu^2 K^2}\right).
	\end{equation*}
	That is, to achive $\EE\left[f(\bar{x}^K) - f(x^*)\right] \le \varepsilon$ {\tt D-SGD} requires
	\begin{equation*}
		\widetilde{\cO}\left(\kappa^2 + \tau\kappa^{\nicefrac{3}{2}} + \frac{D_1'}{\mu\varepsilon} + \frac{\sqrt{L\tau D_1'}}{\mu\sqrt{\varepsilon}}\right) \text{ iterations.}
	\end{equation*}
\end{corollary}

Applying Lemma~\ref{lem:lemma_technical_cvx} we get the complexity result in the case when $\mu = 0$.
\begin{corollary}\label{cor:d_sgd_cvx_cor}
	Let the assumptions of Theorem~\ref{thm:d_sgd_pure} hold, $f_{\xi}(x)$ are convex for each $\xi$ and $\mu = 0$. Then after $K$ iterations of {\tt D-SGD} with the stepsize
	\begin{eqnarray*}
		\gamma &=& \min\left\{\frac{1}{8L\sqrt{2\tau\left(\tau + 2\right)}}, \sqrt{\frac{\|x^0 - x^*\|^2}{D_1' K}}, \sqrt[3]{\frac{\|x^0 - x^*\|^2}{3L\tau D_1' K}}\right\}	\end{eqnarray*}		
	we have $\EE\left[f(\bar{x}^K) - f(x^*)\right]$ of order
	\begin{equation*}
		\cO\left(\frac{\tau LR_0^2}{K} + \sqrt{\frac{R_0^2 \tau D_1'}{K}} + \frac{\sqrt[3]{LR_0^4\tau D_1'}}{K^{\nicefrac{2}{3}}}\right)
	\end{equation*}
	where $R_0 = \|x^0 - x^*\|$. That is, to achive $\EE\left[f(\bar{x}^K) - f(x^*)\right] \le \varepsilon$ {\tt D-SGD} requires
	\begin{equation*}
		\cO\left(\frac{\tau LR_0^2}{\varepsilon} + \frac{R_0^2 D_1'}{\varepsilon^2} + \frac{R_0^2\sqrt{L\tau D_1'}}{\varepsilon^{\nicefrac{3}{2}}}\right)
	\end{equation*}
	iterations.
\end{corollary}

\subsection{{\tt D-QSGD}}\label{sec:d_qsgd}
In this section we show how one can combine delayed updates with quantization using our scheme.
\begin{algorithm}[h!]
   \caption{{\tt D-QSGD}}\label{alg:d-qsgd}
\begin{algorithmic}[1]
   \Require learning rate $\gamma>0$, initial vector $x^0 \in \R^d$
	\State Set $e_i^0 = 0$ for all $i=1,\ldots, n$   
   \For{$k=0,1,\dotsc$}
       \State Broadcast $x^{k-\tau}$ to all workers
        \For{$i=1,\dotsc,n$}
			\State Sample $\hat g_i^{k-\tau}$ independently from other nodes such that $\EE[\hat g_i^{k-\tau}\mid x^{k-\tau}] = \nabla f_i(x^{k-\tau})$ and $\EE\left[\|\hat g_i^{k-\tau} - \nabla f_i(x^{k-\tau})\|^2\mid x^{k-\tau}\right] \le D_{i}$            
            \State $g^{k-\tau}_i = Q(\hat g_i^{k-\tau}) - \nabla f_i(x^*)$ (quantization is performed independently from other nodes)
            \State $v_i^k = \gamma g_i^{k-\tau}$
            \State $e_i^{k+1} = e_i^k + \gamma g_i^k - v_i^k$
        \EndFor
        \State $e^k = \frac{1}{n}\sum_{i=1}^ne_i^k$, $g^k = \frac{1}{n}\sum_{i=1}^ng_i^k$, $v^k = \frac{1}{n}\sum_{i=1}^nv_i^k = \frac{\gamma}{n}\sum_{i=1}^ng_i^{k-\tau} = \frac{\gamma}{n}\sum_{i=1}^nQ(\hat g_i^{k-\tau})$
       \State $x^{k+1} = x^k - v^k$
   \EndFor
\end{algorithmic}
\end{algorithm}
\begin{lemma}\label{lem:d_qsgd_second_moment_bound}
	Assume that $f_i(x)$ is convex and $L$-smooth for all $i=1,\ldots,n$. Then, for all $k\ge 0$ we have
	\begin{eqnarray*}
		\EE\left[g^k\mid x^k\right] &=& \nabla f(x^k), 
		\\
		\EE\left[\|g^k\|^2\mid x^k\right] &\le& 2L\left(1 + \frac{2\omega}{n}\right)\left(f(x^k) - f(x^*)\right) + \frac{(\omega+1)D}{n} + \frac{2\omega}{n^2}\sum\limits_{i=1}^n\|\nabla f_i(x^*)\|^2 
	\end{eqnarray*}
	where $D = \frac{1}{n}\sum_{i=1}^n D_{i}$.
\end{lemma}
\begin{proof}
	First of all, we show unbiasedness of $g^k$:
	\begin{eqnarray*}
		\EE\left[g^k\mid x^k\right] &=& \frac{1}{n}\sum\limits_{i=1}^n\EE\left[g_i^k\mid x^k\right] = \frac{1}{n}\sum\limits_{i=1}^n\EE\left[\EE_Q\left[Q(\hat g_i^k) - \nabla f_i(x^*)\right]\mid x^k\right]\\
		&\overset{\eqref{eq:quantization_def}}{=}& \frac{1}{n}\sum\limits_{i=1}^n\left(\nabla f_i(x^k) - \nabla f_i(x^*)\right) = \nabla f(x^k),
	\end{eqnarray*}
	where $\EE_{Q}\left[\cdot\right]$ denotes mathematical expectation w.r.t.\ the randomness coming only from the quantization. Next, we derive the upper bound for the second moment of $g^k$:
	\begin{eqnarray}
		\EE_Q\left[\|g^k\|^2\right] &=& \EE_Q\left[\left\|\frac{1}{n}\sum\limits_{i=1}^n\left(Q(\hat g_i^k) - \nabla f_i(x^*)\right)\right\|^2\right]\notag\\
		&\overset{\eqref{eq:variance_decomposition}}{=}& \EE_Q\left[\left\|\frac{1}{n}\sum\limits_{i=1}^n\left(Q(\hat g_i^k) - \hat g_i^k\right)\right\|^2\right] + \left\|\frac{1}{n}\sum\limits_{i=1}^n\left(\hat g_i^k - \nabla f_i(x^*)\right)\right\|^2.\label{eq:d_qsgd_technical_1}
	\end{eqnarray}
	Since $Q(\hat g_1^k),\ldots, Q(\hat g_n^k)$ are independent quantizations, we get
	\begin{eqnarray*}
		\EE_Q\left[\|g^k\|^2\right] &\overset{\eqref{eq:d_qsgd_technical_1}}{\le}& \frac{1}{n^2}\sum\limits_{i=1}^n\EE_Q\left[\left\|Q(\hat g_i^k) - \hat g_i^k\right\|^2\right] + \left\|\frac{1}{n}\sum\limits_{i=1}^n\left(\hat g_i^k - \nabla f_i(x^*)\right)\right\|^2\\
		&\overset{\eqref{eq:quantization_def}}{\le}& \frac{\omega}{n^2}\sum\limits_{i=1}^n\|\hat g_i^k\|^2 + \left\|\frac{1}{n}\sum\limits_{i=1}^n\left(\hat g_i^k - \nabla f_i(x^*)\right)\right\|^2.
	\end{eqnarray*}
	Taking conditional expectation $\EE\left[\cdot\mid x^k\right]$ from the both sides of the previous inequality we obtain
	\begin{eqnarray}
		\EE\left[\|g^k\|^2\mid x^k\right] &\le& \frac{\omega}{n^2}\sum\limits_{i=1}^n\EE\left[\|\hat g_i^k\|^2\mid x^k\right] + \EE\left[\left\|\frac{1}{n}\sum\limits_{i=1}^n\left(\hat g_i^k - \nabla f_i(x^*)\right)\right\|^2\mid x^k\right]\notag\\
		&\overset{\eqref{eq:variance_decomposition}}{\le}& \frac{\omega}{n^2}\sum\limits_{i=1}^n\|\nabla f_i(x^k)\|^2 + \frac{\omega}{n^2}\sum\limits_{i=1}^n\EE\left[\|\hat g_i^k - \nabla f_i(x^k)\|^2\mid x^k\right]\notag\\
		&&\quad + \underbrace{\left\|\frac{1}{n}\sum\limits_{i=1}^n\left(\nabla f_i(x^k) - \nabla f_i(x^*)\right)\right\|^2}_{\|\nabla f(x^k)-\nabla f(x^*)\|^2} + \EE\left[\left\|\frac{1}{n}\sum\limits_{i=1}^n\left(\hat g_i^k - \nabla f_i(x^k)\right)\right\|^2\mid x^k\right].\notag
	\end{eqnarray}
	It remains to estimate terms in the second and the third lines of the previous inequality:
	\begin{eqnarray*}
		\frac{\omega}{n^2}\sum\limits_{i=1}^n\|\nabla f_i(x^k)\|^2 &\overset{\eqref{eq:a_b_norm_squared}}{\le}& \frac{2\omega}{n^2}\sum\limits_{i=1}^n \|\nabla f_i(x^k) - \nabla f_i(x^*)\|^2 + \frac{2\omega}{n^2}\sum\limits_{i=1}^n \|\nabla f_i(x^*)\|^2\\
		&\overset{\eqref{eq:L_smoothness_cor_3}}{\le}& \frac{4\omega L}{n}\left(f(x^k) - f(x^*)\right) + \frac{2\omega}{n^2}\sum\limits_{i=1}^n \|\nabla f_i(x^*)\|^2,\\
		\frac{\omega}{n}\sum\limits_{i=1}^n\EE\left[\|\hat g_i^k - \nabla f_i(x^k)\|^2\mid x^k\right] &\le& \frac{\omega}{n^2}\sum\limits_{i=1}^nD_{i} = \frac{\omega D}{n},\\
		\|\nabla f(x^k) - \nabla f(x^*)\|^2 &\overset{\eqref{eq:L_smoothness_cor_3}}{\le}& 2L\left(f(x^k) - f(x^*)\right),\\
		\EE\left[\left\|\frac{1}{n}\sum\limits_{i=1}^n\left(\hat g_i^k - \nabla f_i(x^k)\right)\right\|^2\mid x^k\right] &=& \frac{1}{n^2}\sum\limits_{i=1}^n\EE\left[\|\hat g_i^k - \nabla f_i(x^k)\|^2\mid x^k\right]\\
		&\le& \frac{1}{n^2}\sum\limits_{i=1}^n D_i = \frac{D}{n}.
	\end{eqnarray*}
	Putting all together we get
	\begin{eqnarray*}
		\EE\left[\|g^k\|^2\mid x^k\right] &\le& 2L\left(1 + \frac{2\omega}{n}\right)\left(f(x^k) - f(x^*)\right) + \frac{(\omega+1)D}{n} + \frac{2\omega}{n^2}\sum\limits_{i=1}^n \|\nabla f_i(x^*)\|^2.
	\end{eqnarray*}
\end{proof}

\begin{theorem}\label{thm:d_qsgd}
	Assume that $f_i(x)$ is convex and $L$-smooth for all $i=1,\ldots, n$ and $f(x)$ is $\mu$-quasi strongly convex. Then {\tt D-QSGD} satisfies Assumption~\ref{ass:key_assumption_new} with
	\begin{gather*}
		A' = L\left(1 + \frac{2\omega}{n}\right),\quad B_1' = B_2' = 0,\quad D_1' = \frac{(\omega+1)D}{n} + \frac{2\omega}{n^2}\sum\limits_{i=1}^n \|\nabla f_i(x^*)\|^2,\\
		\sigma_{1,k}^2 \equiv \sigma_{2,k}^2 \equiv 0,\quad \rho_1 = \rho_2 = 1,\quad C_1 = C_2 = 0,\quad D_2 = 0\\
		F_1 = F_2 = 0,\quad G = 0,\quad D_3 = \frac{3\gamma\tau L}{n}\left((\omega+1)D + \frac{2\omega}{n}\sum\limits_{i=1}^n\|\nabla f_i(x^*)\|^2\right)
	\end{gather*}
	with $\gamma$ satisfying
	\begin{equation*}
		\gamma \le \min\left\{\frac{1}{4L(1+\nicefrac{2\omega}{n})}, \frac{1}{8L\sqrt{2\tau\left(\tau + 1 + \nicefrac{2\omega}{n}\right)}}\right\}
	\end{equation*}
	and for all $K \ge 0$
	\begin{equation*}
		\EE\left[f(\bar x^K) - f(x^*)\right] \le \left(1 - \frac{\gamma\mu}{2}\right)^K\frac{4\|x^0 - x^*\|^2}{\gamma} + \gamma\left(D_1' + D_3\right)
	\end{equation*}
	when $\mu > 0$ and
	\begin{equation*}
		\EE\left[f(\bar x^K) - f(x^*)\right] \le \frac{4\|x^0 - x^*\|^2}{\gamma K} + \gamma\left(D_1' + D_3\right)
	\end{equation*}
	when $\mu = 0$.
\end{theorem}
In other words, {\tt D-QSGD} converges with the linear rate
\begin{equation*}
	\cO\left(\left(\kappa\left(1+\frac{\omega}{n}\right) + \kappa\sqrt{\tau\left(\tau + \frac{\omega}{n}\right)}\right)\ln\frac{1}{\varepsilon}\right)
\end{equation*}
to the neighbourhood of the solution when $\mu > 0$. Applying Lemma~\ref{lem:lemma2_stich} we establish the rate of convergence to $\varepsilon$-solution.
\begin{corollary}\label{cor:d_QSGD_str_cvx_cor}
	Let the assumptions of Theorem~\ref{thm:d_qsgd} hold, $f_{\xi}(x)$ are convex for each $\xi$ and $\mu > 0$. Then after $K$ iterations of {\tt D-QSGD} with the stepsize
	\begin{eqnarray*}
		\gamma_0 &=& \min\left\{\frac{1}{4L(1+\nicefrac{2\omega}{n})}, \frac{1}{8L\sqrt{2\tau\left(\tau + 1 + \nicefrac{2\omega}{n}\right)}}\right\},\quad R_0 = \|x^0 - x^*\|,\\
		\gamma &=& \min\left\{\gamma_0, \frac{\ln\left(\max\left\{2,\min\left\{\frac{R_0^2\mu^2K^2}{D_1'}, \frac{R_0^2\mu^3K^3}{3\tau LD_1'}\right\}\right\}\right)}{\mu K}\right\}
	\end{eqnarray*}
	we have $\EE\left[f(\bar{x}^K) - f(x^*)\right]$ of order
	\begin{equation*}
		 \widetilde\cO\left(LR_0^2\left(1+\frac{\omega}{n}+\sqrt{\tau\left(\tau + \frac{\omega}{n}\right)}\right)\exp\left(-\frac{\mu}{L\left(1+\frac{\omega}{n}+\sqrt{\tau\left(\tau + \frac{\omega}{n}\right)}\right)}K\right) + \frac{D_1'}{\mu K} + \frac{L\tau D_1'}{\mu^2 K^2}\right).
	\end{equation*}
	That is, to achive $\EE\left[f(\bar{x}^K) - f(x^*)\right] \le \varepsilon$ {\tt D-QSGD} requires
	\begin{equation*}
		\widetilde{\cO}\left(\frac{L}{\mu}\left(1+\frac{\omega}{n}\right) + \frac{L}{\mu}\sqrt{\tau\left(\tau + \frac{\omega}{n}\right)} + \frac{D_1'}{\mu\varepsilon} + \frac{\sqrt{L\tau D_1'}}{\mu\sqrt{\varepsilon}}\right) \text{ iterations.}
	\end{equation*}
\end{corollary}

Applying Lemma~\ref{lem:lemma_technical_cvx} we get the complexity result in the case when $\mu = 0$.
\begin{corollary}\label{cor:d_QSGD_cvx_cor}
	Let the assumptions of Theorem~\ref{thm:d_qsgd} hold and $\mu = 0$. Then after $K$ iterations of {\tt D-QSGD} with the stepsize
	\begin{eqnarray*}
		\gamma_0 &=& \min\left\{\frac{1}{4L(1+\nicefrac{2\omega}{n})}, \frac{1}{8L\sqrt{2\tau\left(\tau + 1 + \nicefrac{2\omega}{n}\right)}}\right\},\\	
		\gamma &=& \min\left\{\gamma_0, \sqrt{\frac{\|x^0 - x^*\|^2}{D_1' K}}, \sqrt[3]{\frac{\|x^0 - x^*\|^2}{3L\tau D_1' K}}\right\}	
	\end{eqnarray*}		
	we have $\EE\left[f(\bar{x}^K) - f(x^*)\right]$ of order
	\begin{equation*}
		\cO\left(\frac{LR_0^2\left(1+\frac{\omega}{n}\right)}{K} + \frac{LR_0^2\sqrt{\tau\left(\tau + \frac{\omega}{n}\right)}}{K} + \sqrt{\frac{R_0^2 D_1'}{K}} + \frac{\sqrt[3]{LR_0^4\tau D_1'}}{K^{\nicefrac{2}{3}}}\right)
	\end{equation*}
	where $R_0 = \|x^0 - x^*\|$. That is, to achive $\EE\left[f(\bar{x}^K) - f(x^*)\right] \le \varepsilon$ {\tt D-QSGD} requires
	\begin{equation*}
		\cO\left(\frac{LR_0^2\left(1+\frac{\omega}{n}\right)}{\varepsilon} + \frac{LR_0^2\sqrt{\tau\left(\tau + \frac{\omega}{n}\right)}}{\varepsilon} + \frac{R_0^2 D_1'}{\varepsilon^2} + \frac{R_0^2\sqrt{L\tau D_1'}}{\varepsilon^{\nicefrac{3}{2}}}\right)
	\end{equation*}
	iterations.
\end{corollary}

\subsection{{\tt D-QSGDstar}}\label{sec:d_qsgd_star}
As we saw in Section~\ref{sec:d_qsgd} {\tt D-QSGD} fails to converge to the exact optimum asymptotically even if $\hat g_i^k = \nabla f_i(x^k)$ for all $i=1,\ldots,n$ almost surely, i.e., all $D_i = 0$ for all $i=1,\ldots,n$. As for {\tt EC-GDstar} we assume now that $i$-th worker has an access to $\nabla f_i(x^*)$. Using this one can construct the method with delayed updates that converges asymptotically to the exact solution when the full gradients are available.
\begin{algorithm}[h!]
   \caption{{\tt D-QSGDstar}}\label{alg:d-qSGDstar}
\begin{algorithmic}[1]
   \Require learning rate $\gamma>0$, initial vector $x^0 \in \R^d$
	\State Set $e_i^0 = 0$ for all $i=1,\ldots, n$   
   \For{$k=0,1,\dotsc$}
       \State Broadcast $x^{k-\tau}$ to all workers
        \For{$i=1,\dotsc,n$}
			\State Sample $\hat g_i^{k-\tau}$ independently from other nodes such that $\EE[\hat g_i^{k-\tau}\mid x^{k-\tau}] = \nabla f_i(x^{k-\tau})$ and $\EE\left[\|\hat g_i^{k-\tau} - \nabla f_i(x^{k-\tau})\|^2\mid x^{k-\tau}\right] \le D_{i}$            
            \State $g^{k-\tau}_i = Q(\hat g_i^{k-\tau}- \nabla f_i(x^*))$ (quantization is performed independently from other nodes)
            \State $v_i^k = \gamma g_i^{k-\tau}$
            \State $e_i^{k+1} = e_i^k + \gamma g_i^k - v_i^k$
        \EndFor
        \State $e^k = \frac{1}{n}\sum_{i=1}^ne_i^k$, $g^k = \frac{1}{n}\sum_{i=1}^ng_i^k$, $v^k = \frac{1}{n}\sum_{i=1}^nv_i^k = \frac{\gamma}{n}\sum_{i=1}^ng_i^{k-\tau} = \frac{\gamma}{n}\sum_{i=1}^nQ(\hat g_i^{k-\tau}-\nabla f_i(x^*))$
       \State $x^{k+1} = x^k - v^k$
   \EndFor
\end{algorithmic}
\end{algorithm}
\begin{lemma}\label{lem:d_qsgd_star_second_moment_bound}
	Assume that $f_i(x)$ is convex and $L$-smooth for all $i=1,\ldots,n$. Then, for all $k\ge 0$ we have
	\begin{eqnarray}
		\EE\left[g^k\mid x^k\right] &=& \nabla f(x^k), \label{eq:d_qsgd_star_unbiasedness}\\
		\EE\left[\|g^k\|^2\mid x^k\right] &\le& 2L\left(1 + \frac{\omega}{n}\right)\left(f(x^k) - f(x^*)\right) + \frac{(\omega+1)D}{n} \label{eq:d_qsgd_star_second_moment_bound}
	\end{eqnarray}
	where $D = \frac{1}{n}\sum_{i=1}^n D_{i}$.
\end{lemma}
\begin{proof}
	First of all, we show unbiasedness of $g^k$:
	\begin{eqnarray*}
		\EE\left[g^k\mid x^k\right] &=& \frac{1}{n}\sum\limits_{i=1}^n\EE\left[g_i^k\mid x^k\right] = \frac{1}{n}\sum\limits_{i=1}^n\EE\left[\EE_Q\left[Q(\hat g_i^k - \nabla f_i(x^*))\right]\mid x^k\right]\\
		&\overset{\eqref{eq:quantization_def}}{=}& \frac{1}{n}\sum\limits_{i=1}^n\left(\nabla f_i(x^k) - \nabla f_i(x^*)\right) = \nabla f(x^k),
	\end{eqnarray*}
	where $\EE_{Q}\left[\cdot\right]$ denotes mathematical expectation w.r.t.\ the randomness coming only from the quantization. Next, we derive the upper bound for the second moment of $g^k$:
	\begin{eqnarray}
		\EE_Q\left[\|g^k\|^2\right] &=& \EE_Q\left[\left\|\frac{1}{n}\sum\limits_{i=1}^n\left(Q\left(\hat g_i^k - \nabla f_i(x^*)\right)\right)\right\|^2\right]\notag\\
		&\overset{\eqref{eq:variance_decomposition}}{=}& \EE_Q\left[\left\|\frac{1}{n}\sum\limits_{i=1}^n\left(Q\left(\hat g_i^k - \nabla f_i(x^*)\right) - \left(\hat g_i^k - \nabla f_i(x^*)\right)\right)\right\|^2\right]\notag\\
		&&\quad + \left\|\frac{1}{n}\sum\limits_{i=1}^n\left(\hat g_i^k - \nabla f_i(x^*)\right)\right\|^2.\label{eq:d_qsgd_star_technical_1}
	\end{eqnarray}
	Since $Q\left(\hat g_1^k - \nabla f_1(x^*)\right),\ldots, Q\left(\hat g_n^k - \nabla f_n(x^*)\right)$ are independent quantizations, we get
	\begin{eqnarray*}
		\EE_Q\left[\|g^k\|^2\right] &\overset{\eqref{eq:d_qsgd_star_technical_1}}{\le}& \frac{1}{n^2}\sum\limits_{i=1}^n\EE_Q\left[\left\|Q\left(\hat g_i^k - \nabla f_i(x^*)\right) - \left(\hat g_i^k - \nabla f_i(x^*)\right)\right\|^2\right]\notag\\
		&&\quad + \left\|\frac{1}{n}\sum\limits_{i=1}^n\left(\hat g_i^k - \nabla f_i(x^*)\right)\right\|^2\\
		&\overset{\eqref{eq:quantization_def}}{\le}& \frac{\omega}{n^2}\sum\limits_{i=1}^n\|\hat g_i^k - \nabla f_i(x^*)\|^2 + \left\|\frac{1}{n}\sum\limits_{i=1}^n\left(\hat g_i^k - \nabla f_i(x^*)\right)\right\|^2.
	\end{eqnarray*}
	Taking conditional expectation $\EE\left[\cdot\mid x^k\right]$ from the both sides of the previous inequality we obtain
	\begin{eqnarray}
		\EE\left[\|g^k\|^2\mid x^k\right] &\le& \frac{\omega}{n^2}\sum\limits_{i=1}^n\EE\left[\|\hat g_i^k - \nabla f_i(x^*)\|^2\mid x^k\right] + \EE\left[\left\|\frac{1}{n}\sum\limits_{i=1}^n\left(\hat g_i^k - \nabla f_i(x^*)\right)\right\|^2\mid x^k\right]\notag\\
		&\overset{\eqref{eq:variance_decomposition}}{\le}& \frac{\omega}{n^2}\sum\limits_{i=1}^n\|\nabla f_i(x^k)- \nabla f_i(x^*)\|^2 + \frac{\omega}{n^2}\sum\limits_{i=1}^n\EE\left[\|\hat g_i^k - \nabla f_i(x^k)\|^2\mid x^k\right]\notag\\
		&&\quad + \underbrace{\left\|\frac{1}{n}\sum\limits_{i=1}^n\left(\nabla f_i(x^k) - \nabla f_i(x^*)\right)\right\|^2}_{\|\nabla f(x^k)-\nabla f(x^*)\|^2} + \EE\left[\left\|\frac{1}{n}\sum\limits_{i=1}^n\left(\hat g_i^k - \nabla f_i(x^k)\right)\right\|^2\mid x^k\right].\notag
	\end{eqnarray}
	It remains to estimate terms in the second and the third lines of the previous inequality:
	\begin{eqnarray*}
		\frac{\omega}{n^2}\sum\limits_{i=1}^n\|\nabla f_i(x^k)-\nabla f_i(x^*)\|^2 &\overset{\eqref{eq:L_smoothness_cor_3}}{\le}& \frac{2\omega L}{n}\left(f(x^k) - f(x^*)\right),\\
		\frac{\omega}{n}\sum\limits_{i=1}^n\EE\left[\|\hat g_i^k - \nabla f_i(x^k)\|^2\mid x^k\right] &\le& \frac{\omega}{n^2}\sum\limits_{i=1}^nD_{i} = \frac{\omega D}{n},\\
		\|\nabla f(x^k) - \nabla f(x^*)\|^2 &\overset{\eqref{eq:L_smoothness_cor_3}}{\le}& 2L\left(f(x^k) - f(x^*)\right),\\
		\EE\left[\left\|\frac{1}{n}\sum\limits_{i=1}^n\left(\hat g_i^k - \nabla f_i(x^k)\right)\right\|^2\mid x^k\right] &=& \frac{1}{n^2}\sum\limits_{i=1}^n\EE\left[\|\hat g_i^k - \nabla f_i(x^k)\|^2\mid x^k\right]\\
		&\le& \frac{1}{n^2}\sum\limits_{i=1}^n D_i = \frac{D}{n}.
	\end{eqnarray*}
	Putting all together we get
	\begin{eqnarray*}
		\EE\left[\|g^k\|^2\mid x^k\right] &\le& 2L\left(1 + \frac{\omega}{n}\right)\left(f(x^k) - f(x^*)\right) + \frac{(\omega+1)D}{n}.
	\end{eqnarray*}
\end{proof}

\begin{theorem}\label{thm:d_qsgd_star}
	Assume that $f_i(x)$ is convex and $L$-smooth for all $i=1,\ldots, n$ and $f(x)$ is $\mu$-quasi strongly convex. Then {\tt D-QSGDstar} satisfies Assumption~\ref{ass:key_assumption_new} with
	\begin{gather*}
		A' = L\left(1 + \frac{\omega}{n}\right),\quad B_1' = B_2' = 0,\quad D_1' = \frac{(\omega+1)D}{n},\quad \sigma_{1,k}^2 \equiv \sigma_{2,k}^2 \equiv 0,\\
		\rho_1 = \rho_2 = 1,\quad C_1 = C_2 = 0,\quad D_2 = 0,\quad G = 0,\\
		F_1 = F_2 = 0,\quad D_3 = \frac{3\gamma\tau L(\omega+1)D}{n}
	\end{gather*}
	with $\gamma$ satisfying
	\begin{equation*}
		\gamma \le \min\left\{\frac{1}{4L(1+\nicefrac{\omega}{n})}, \frac{1}{8L\sqrt{\tau\left(\tau + 1 + \nicefrac{\omega}{n}\right)}}\right\}.
	\end{equation*}
	and for all $K \ge 0$
	\begin{equation*}
		\EE\left[f(\bar x^K) - f(x^*)\right] \le \left(1 - \frac{\gamma\mu}{2}\right)^K\frac{4\|x^0 - x^*\|^2}{\gamma} + 4\gamma\left(D_1' + D_3\right)
	\end{equation*}	
	when $\mu > 0$ and
	\begin{equation*}
		\EE\left[f(\bar x^K) - f(x^*)\right] \le \frac{4\|x^0 - x^*\|^2}{\gamma K} + 4\gamma\left(D_1' + D_3\right)
	\end{equation*}
	when $\mu = 0$.
\end{theorem}
In other words, {\tt D-QSGDstar} converges with the linear rate
\begin{equation*}
	\cO\left(\left(\tau + \kappa\left(1+\frac{\omega}{n}\right) + \kappa\sqrt{\tau\left(\tau + \frac{\omega}{n}\right)}\right)\ln\frac{1}{\varepsilon}\right)
\end{equation*}
to the exact solution when $\mu > 0$ and $D = 0$, i.e., $\hat g_i^k = \nabla f_i(x^k)$ for all $i=1,\ldots,n$ almost surely. Applying Lemma~\ref{lem:lemma2_stich} we establish the rate of convergence to $\varepsilon$-solution.
\begin{corollary}\label{cor:d_QSGDstar_str_cvx_cor}
	Let the assumptions of Theorem~\ref{thm:d_qsgd_star} hold and $\mu > 0$. Then after $K$ iterations of {\tt D-QSGDstar} with the stepsize
	\begin{eqnarray*}
		\gamma_0 &=& \min\left\{\frac{1}{4L(1+\nicefrac{\omega}{n})}, \frac{1}{8L\sqrt{\tau\left(\tau + 1 + \nicefrac{\omega}{n}\right)}}\right\},\quad R_0 = \|x^0 - x^*\|,\\
		\gamma &=& \min\left\{\gamma_0, \frac{\ln\left(\max\left\{2,\min\left\{\frac{nR_0^2\mu^2K^2}{D}, \frac{nR_0^2\mu^3K^3}{3\tau LD}\right\}\right\}\right)}{\mu K}\right\}
	\end{eqnarray*}	 
	we have $\EE\left[f(\bar{x}^K) - f(x^*)\right]$ of order
	\begin{equation*}
		 \widetilde\cO\left(LR_0^2\left(1+\frac{\omega}{n}+\sqrt{\tau\left(\tau + \frac{\omega}{n}\right)}\right)\exp\left(-\frac{\mu}{L\left(1+\frac{\omega}{n}+\sqrt{\tau\left(\tau + \frac{\omega}{n}\right)}\right)}K\right) + \frac{D}{n\mu K} + \frac{L\tau D}{n\mu^2 K^2}\right).
	\end{equation*}
	That is, to achive $\EE\left[f(\bar{x}^K) - f(x^*)\right] \le \varepsilon$ {\tt D-QSGDstar} requires
	\begin{equation*}
		\widetilde{\cO}\left(\frac{L}{\mu}\left(1+\frac{\omega}{n}\right) + \frac{L}{\mu}\sqrt{\tau\left(\tau + \frac{\omega}{n}\right)} + \frac{D}{n\mu\varepsilon} + \frac{\sqrt{L\tau D}}{\mu\sqrt{n\varepsilon}}\right) \text{ iterations.}
	\end{equation*}
\end{corollary}

Applying Lemma~\ref{lem:lemma_technical_cvx} we get the complexity result in the case when $\mu = 0$.
\begin{corollary}\label{cor:d_QSGDstar_cvx_cor}
	Let the assumptions of Theorem~\ref{thm:d_qsgd_star} hold and $\mu = 0$. Then after $K$ iterations of {\tt D-QSGDstar} with the stepsize
	\begin{eqnarray*}
		\gamma_0 &=& \min\left\{\frac{1}{4L(1+\nicefrac{2\omega}{n})}, \frac{1}{8L\sqrt{\tau\left(\tau + 1 + \nicefrac{\omega}{n}\right)}}\right\},\\	
		\gamma &=& \min\left\{\gamma_0, \sqrt{\frac{n\|x^0 - x^*\|^2}{D K}}, \sqrt[3]{\frac{n\|x^0 - x^*\|^2}{3L\tau D K}}\right\}	
	\end{eqnarray*}		
	we have $\EE\left[f(\bar{x}^K) - f(x^*)\right]$ of order
	\begin{equation*}
		\cO\left(\frac{LR_0^2\left(1+\frac{\omega}{n}\right)}{K} + \frac{LR_0^2\sqrt{\tau\left(\tau + \frac{\omega}{n}\right)}}{K} + \sqrt{\frac{R_0^2 D}{nK}} + \frac{\sqrt[3]{LR_0^4\tau D}}{n^{\nicefrac{1}{3}}K^{\nicefrac{2}{3}}}\right)
	\end{equation*}
	where $R_0 = \|x^0 - x^*\|$. That is, to achive $\EE\left[f(\bar{x}^K) - f(x^*)\right] \le \varepsilon$ {\tt D-QSGDstar} requires
	\begin{equation*}
		\cO\left(\frac{LR_0^2\left(1+\frac{\omega}{n}\right)}{\varepsilon} + \frac{LR_0^2\sqrt{\tau\left(\tau + \frac{\omega}{n}\right)}}{\varepsilon} + \frac{R_0^2 D}{n\varepsilon^2} + \frac{R_0^2\sqrt{L\tau D}}{\sqrt{n}\varepsilon^{\nicefrac{3}{2}}}\right)
	\end{equation*}
	iterations.
\end{corollary}

\subsection{{\tt D-SGD-DIANA}}\label{sec:d_diana}
In this section we present a practical version of {\tt D-QSGDstar}: {\tt D-SGD-DIANA}.
\begin{algorithm}[h!]
   \caption{{\tt D-SGD-DIANA}}\label{alg:d-diana}
\begin{algorithmic}[1]
   \Require learning rates $\gamma>0, \alpha\in(0,1]$, initial vectors $x^0, h_1^0,\ldots,h_n^0 \in \R^d$
	\State Set $e_i^0 = 0$ for all $i=1,\ldots, n$   
	\State Set $h^0 = \frac{1}{n}\sum_{i=1}^n h_i^0$   
   \For{$k=0,1,\dotsc$}
       \State Broadcast $x^{k-\tau}$ to all workers
        \For{$i=1,\dotsc,n$}
			\State Sample $\hat g_i^{k-\tau}$ independently from other nodes such that $\EE[\hat g_i^{k-\tau}\mid x^{k-\tau}] = \nabla f_i(x^{k-\tau})$ and $\EE\left[\|\hat g_i^{k-\tau} - \nabla f_i(x^{k-\tau})\|^2\mid x^{k-\tau}\right] \le D_{i}$            
			\State $\hat \Delta_i^{k-\tau} = Q(\hat g_i^{k-\tau}- h_i^{k-\tau})$ (quantization is performed independently from other nodes)            
            \State $g^{k-\tau}_i = h_i^{k-\tau} + \hat \Delta_i^{k-\tau}$
            \State $v_i^k = \gamma g_i^{k-\tau}$
            \State $e_i^{k+1} = e_i^k + \gamma g_i^k - v_i^k$
            \State $h_i^{k-\tau+1} = h_i^{k-\tau} + \alpha\hat \Delta_i^{k-\tau}$
        \EndFor
        \State $h^{k-\tau} = \frac{1}{n}\sum_{i=1}^nh_i^{k-\tau}$, $e^k = \frac{1}{n}\sum_{i=1}^ne_i^k$, $g^k = \frac{1}{n}\sum_{i=1}^ng_i^k$, $v^k = \frac{1}{n}\sum_{i=1}^nv_i^k = \frac{\gamma}{n}\sum_{i=1}^ng_i^{k-\tau} = \gamma h^{k-\tau} + \frac{\gamma}{n}\sum_{i=1}^n\hat{\Delta}_i^{k-\tau}$
       \State $x^{k+1} = x^k - v^k$
       \State $h^{k-\tau+1} = h^{k-\tau}+\frac{\alpha}{n}\sum_{i=1}^n\hat{\Delta}_i^{k-\tau}$
   \EndFor
\end{algorithmic}
\end{algorithm}

\begin{lemma}[Lemmas 1 and 2 from \cite{horvath2019stochastic}]\label{lem:d_diana_second_moment_bound}
	Assume that $f_i(x)$ is convex and $L$-smooth for all $i=1,\ldots,n$ and $\alpha \le \nicefrac{1}{(\omega+1)}$. Then, for all $k\ge 0$ we have
	\begin{eqnarray}
		\EE\left[g^k\mid x^k\right] &=& \nabla f(x^k), \label{eq:d_diana_unbiasedness}\\
		\EE\left[\|g^k\|^2\mid x^k\right] &\le& 2L\left(1 + \frac{2\omega}{n}\right)\left(f(x^k) - f(x^*)\right) + \frac{2\omega\sigma_k^2}{n} + \frac{(\omega+1)D}{n} \label{eq:d_diana_second_moment_bound}\\
		\EE\left[\sigma_{k+1}^2\mid x^k\right] &\le& (1-\alpha)\sigma_k^2 + 2L\alpha\left(f(x^k) - f(x^*)\right) + \alpha D \label{eq:d_diana_sigma_k+1_bound}
	\end{eqnarray}
	where $\sigma_k^2 = \frac{1}{n}\sum_{i=1}^n\|h_i^k - \nabla f_i(x^*)\|^2$ and $D = \frac{1}{n}\sum_{i=1}^n D_{i}$.
\end{lemma}

\begin{theorem}\label{thm:d_diana}
	Assume that $f_i(x)$ is convex and $L$-smooth for all $i=1,\ldots, n$ and $f(x)$ is $\mu$-quasi strongly convex. Then {\tt D-SGD-DIANA} satisfies Assumption~\ref{ass:key_assumption_new} with
	\begin{gather*}
		A' = L\left(1 + \frac{2\omega}{n}\right),\quad B_1' = \frac{2\omega}{n},\quad D_1' = \frac{(\omega+1)D}{n},\quad \sigma_{1,k}^2 = \sigma_k^2 = \frac{1}{n}\sum_{i=1}^n\|h_i^k - \nabla f_i(x^*)\|^2,\\
		B_2' = 0,\quad \rho_1 = \alpha,\quad \rho_2 = 1,\quad C_1 = L\alpha,\quad C_2 = 0,\quad D_2 = \frac{\alpha(\omega+1)D}{n},\quad G = 0,\\
		F_1 = \frac{12\gamma^2L\omega \tau(2+\alpha)}{n\alpha},\quad F_2 = 0,\quad D_3 = 3\gamma\tau L\left(1 + \frac{4\omega}{n}\right)\frac{(\omega+1)D}{n}
	\end{gather*}
	with $\gamma$ and $\alpha$ satisfying
	\begin{equation*}
		\gamma \le \min\left\{\frac{1}{4L(1+\nicefrac{14\omega}{3n})}, \frac{1}{8L\sqrt{2\tau\left(1+\tau + \nicefrac{2\omega}{n} + \nicefrac{4\omega}{n(1-\alpha)} \right)}}\right\},\quad \alpha \le \frac{1}{\omega+1}, \quad M_1 = \frac{8\omega}{3n\alpha}
	\end{equation*}
	and for all $K \ge 0$
	\begin{equation*}
		\EE\left[f(\bar x^K) - f(x^*)\right] \le \left(1 - \min\left\{\frac{\gamma\mu}{2},\frac{\alpha}{4}\right\}\right)^K\frac{4(T^0 + \gamma F_1 \sigma_0^2)}{\gamma} + 4\gamma\left(D_1' + M_1D_2 + D_3\right)
	\end{equation*}	
	when $\mu > 0$ and 
	\begin{equation*}
		\EE\left[f(\bar x^K) - f(x^*)\right] \le \frac{4(T^0 + \gamma F_1 \sigma_0^2)}{\gamma K} + 4\gamma\left(D_1' + M_1D_2 + D_3\right)
	\end{equation*}
	when $\mu = 0$, where $T^k \eqdef \|\tx^k - x^*\|^2 + M_1\gamma^2 \sigma_k^2$.	
\end{theorem}
In other words, if
\begin{equation*}
	\gamma \le \min\left\{\frac{1}{4L(1+\nicefrac{14\omega}{3n})}, \frac{1}{8L\sqrt{2\tau\left(1+\tau + \nicefrac{10\omega}{n} \right)}}\right\},\quad \alpha \le \min\left\{\frac{1}{\omega+1},\frac{1}{2}\right\}
\end{equation*}
then {\tt D-SGD-DIANA} converges with the linear rate
\begin{equation*}
	\cO\left(\left(\omega + \kappa\left(1+\frac{\omega}{n}\right) + \kappa\sqrt{\tau\left(\tau + \frac{\omega}{n}\right)}\right)\ln\frac{1}{\varepsilon}\right)
\end{equation*}
to the exact solution when $\mu > 0$. Applying Lemma~\ref{lem:lemma2_stich} we establish the rate of convergence to $\varepsilon$-solution.
\begin{corollary}\label{cor:d_diana_str_cvx_cor}
	Let the assumptions of Theorem~\ref{thm:d_diana} hold and $\mu > 0$. Then after $K$ iterations of {\tt D-SGD-DIANA} with the stepsize
	\begin{eqnarray*}
		\gamma_0 &=& \min\left\{\frac{1}{4L(1+\nicefrac{14\omega}{3n})}, \frac{1}{8L\sqrt{2\tau\left(1+\tau + \nicefrac{10\omega}{n} \right)}}\right\},\quad R_0 = \|x^0 - x^*\|,\\
		 \tilde{F}_1 &=& \frac{12L\omega\tau(2+\alpha)\gamma_0^2}{n\alpha},\quad \tilde{T}^0 = R_0^2 + M_1\gamma_0^2\sigma_0^2,\\
		\gamma &=& \min\left\{\gamma_0, \frac{\ln\left(\max\left\{2,\min\left\{\frac{\left(\tilde{T}^0+\gamma_0\tilde{F}_1\sigma_0^2\right)\mu^2K^2}{D_1'+M_1D_2}, \frac{\left(\tilde{T}^0 +\gamma_0\tilde{F}_1\sigma_0^2\right)\mu^3K^3}{3\tau L\left(D_1' + \frac{2B_1'D_2}{\alpha}\right)}\right\}\right\}\right)}{\mu K}\right\}
	\end{eqnarray*}
	and $\alpha \le \min\left\{\frac{1}{\omega+1},\frac{1}{2}\right\}$ we have $\EE\left[f(\bar{x}^K) - f(x^*)\right]$ of order
	\begin{gather*}
		\widetilde\cO\left(LR_0^2\left(1+\frac{\omega}{n}+\sqrt{\tau\left(\tau + \frac{\omega}{n}\right)}\right)\exp\left(-\min\left\{\frac{\mu}{L\left(1+\frac{\omega}{n}+\sqrt{\tau\left(\tau + \frac{\omega}{n}\right)}\right)},\frac{1}{1+\omega}\right\}K\right)\right)\\
		 + \widetilde{O}\left(\frac{D_1'+M_1D_2}{\mu K} + \frac{\tau L\left(D_1' + \frac{B_1'D_2}{\alpha}\right)}{\mu^2 K^2}\right).
	\end{gather*}
	That is, to achive $\EE\left[f(\bar{x}^K) - f(x^*)\right] \le \varepsilon$ {\tt D-SGD-DIANA} requires
	\begin{equation*}
		\widetilde{\cO}\left(\omega+\frac{L}{\mu}\left(1+\frac{\omega}{n}\right) + \frac{L}{\mu}\sqrt{\tau\left(\tau + \frac{\omega}{n}\right)} + \frac{(\omega+1)\left(1+\frac{\omega}{n}\right)D}{n\mu\varepsilon} + \frac{\sqrt{L\tau(\omega+1)\left(1+\frac{\omega}{n}\right) D}}{\mu\sqrt{n\varepsilon}}\right)
	\end{equation*}
	iterations.
\end{corollary}

Applying Lemma~\ref{lem:lemma_technical_cvx} we get the complexity result in the case when $\mu = 0$.
\begin{corollary}\label{cor:d_diana_cvx_cor}
	Let the assumptions of Theorem~\ref{thm:d_diana} hold and $\mu = 0$. Then after $K$ iterations of {\tt D-SGD-DIANA} with the stepsize
	\begin{eqnarray*}
		\gamma_0 &=& \min\left\{\frac{1}{4L(1+\nicefrac{14\omega}{3n})}, \frac{1}{8L\sqrt{2\tau\left(1+\tau + \nicefrac{10\omega}{n} \right)}}\right\},\quad R_0 = \|x^0 - x^*\|,\\
		\gamma &=& \min\left\{\gamma_0, \sqrt{\frac{R_0^2}{M_1\sigma_0^2}}, \sqrt[3]{\frac{R_0^2n\alpha}{12L\omega\tau(2+\alpha)\sigma_0^2}}, \sqrt{\frac{R_0^2}{(D_1' + M_1D_2)K}}, \sqrt[3]{\frac{R_0^2}{3\tau L\left(D_1' + \frac{2B_1'D_2}{\alpha}\right)K}}\right\}	
	\end{eqnarray*}		
	we have $\EE\left[f(\bar{x}^K) - f(x^*)\right]$ of order
	\begin{gather*}
		\cO\left(\frac{L\left(1+\frac{\omega}{n}\right)R_0^2}{K} + \frac{L\sqrt{\tau\left(\tau + \frac{\omega}{n}\right)}R_0^2}{K} + \frac{\sqrt{R_0^2\omega(1+\omega)\sigma_0^2}}{\sqrt{n}K} + \frac{\sqrt[3]{R_0^4 L\tau\omega(1+\omega)\sigma_0^2}}{\sqrt[3]{n}K}\right)\\
		+\cO\left(\sqrt{\frac{(1+\omega)\left(1+\frac{\omega}{n}\right)R_0^2D}{nK}} + \frac{\sqrt[3]{R_0^4\tau L(1+\omega)\left(1+\frac{\omega}{n}\right)D}}{n^{\nicefrac{1}{3}}K^{\nicefrac{2}{3}}}\right).
	\end{gather*}
	That is, to achive $\EE\left[f(\bar{x}^K) - f(x^*)\right] \le \varepsilon$ {\tt D-SGD-DIANA} requires
	\begin{gather*}
		\cO\left(\frac{L\left(1+\frac{\omega}{n}\right)R_0^2}{\varepsilon} + \frac{L\sqrt{\tau\left(\tau + \frac{\omega}{n}\right)}R_0^2}{\varepsilon} + \frac{\sqrt{R_0^2\omega(1+\omega)\sigma_0^2}}{\sqrt{n}\varepsilon} + \frac{\sqrt[3]{R_0^4 L\tau\omega(1+\omega)\sigma_0^2}}{\sqrt[3]{n}\varepsilon}\right)\\
		+\cO\left(\frac{(1+\omega)\left(1+\frac{\omega}{n}\right)R_0^2D}{n\varepsilon^2} + \frac{R_0^2\sqrt{\tau L(1+\omega)\left(1+\frac{\omega}{n}\right)D}}{n^{\nicefrac{1}{2}}\varepsilon^{\nicefrac{3}{2}}}\right)\quad \text{iterations.}
	\end{gather*}
\end{corollary}

\subsection{{\tt D-SGDsr}}\label{sec:d_SGDsr}
In this section we consider the same settings as in Section~\ref{sec:ec_SGDsr}, but this time we consider delayed updates. Moreover, in this section we need slightly weaker assumption.
\begin{assumption}[Expected smoothness]\label{ass:exp_smoothness_f}
	We assume that function $f$ is $\cL$-smooth in expectation w.r.t.\ distribution $\cD$, i.e., there exists constant $\cL = \cL(f,\cD)$ such that
	\begin{equation}
		\EE_{\cD}\left[\|\nabla f_{\xi}(x) - \nabla f_{\xi}(x^*)\|^2\right] \le 2\cL\left(f(x) - f(x^*)\right) \label{eq:exp_smoothness_f}
	\end{equation}
	for all $i\in [n]$ and $x\in\R^d$.
\end{assumption}
\begin{algorithm}[t]
   \caption{{\tt D-SGDsr}}\label{alg:d-SGDsr}
\begin{algorithmic}[1]
   \Require learning rate $\gamma>0$, initial vector $x^0 \in \R^d$
	\State Set $e_i^0 = 0$ for all $i=1,\ldots, n$   
   \For{$k=0,1,\dotsc$}
       \State Broadcast $x^{k-\tau}$ to all workers
        \For{$i=1,\dotsc,n$ in parallel}
            \State Sample $g^{k-\tau}_i = \nabla f_{\xi_i}(x^{k-\tau}) - \nabla f_i(x^*)$
            \State $v_i^k = \gamma g_i^{k-\tau}$
            \State $e_i^{k+1} = e_i^k + \gamma g_i^k - v_i^k$
        \EndFor
        \State $e^k = \frac{1}{n}\sum_{i=1}^ne_i^k$, $g^k = \frac{1}{n}\sum_{i=1}^ng_i^k$, $v^k = \frac{1}{n}\sum_{i=1}^nv_i^k = \frac{1}{n}\sum_{i=1}^n\nabla f_{\xi_i}(x^{k-\tau})$
       \State $x^{k+1} = x^k - v^k$
   \EndFor
\end{algorithmic}
\end{algorithm}

\begin{lemma}\label{lem:key_lemma_d-SGDsr}
	For all $k\ge 0$ we have
	\begin{equation}
		\EE\left[\|g^k\|^2\mid x^k\right] \le 4\cL\left(f(x^k) - f(x^*)\right) + 2\EE_{\cD}\left[\|\nabla f_{\xi}(x^*)\|^2\right]. \label{eq:key_lemma_d-SGDsr} 
	\end{equation}
\end{lemma}
\begin{proof}
	Applying straightforward inequality $\|a+b\|^2 \le 2\|a\|^2 + 2\|b\|^2$ for $a,b\in\R^d$ we get
	\begin{eqnarray*}
		\EE\left[\|g^k\|^2\mid x^k\right] &=& \EE\left[\left\|\frac{1}{n}\sum\limits_{i=1}^n\left(\nabla f_{\xi_i}(x^k) - \nabla f_i(x^*)\right)\right\|^2\mid x^k\right] \\
		&\overset{\eqref{eq:a_b_norm_squared}}{\le}& 2\EE_{\cD}\left[\|\nabla f_{\xi}(x^k) - \nabla f_{\xi}(x^*)\|^2\right] + 2\EE_{\cD}\left[\|\nabla f_{\xi}(x^*) - \nabla f(x^*)\|^2\right]\\
		&\overset{\eqref{eq:exp_smoothness_f}}{\le}& 4\cL\left(f(x^k) - f(x^*)\right) + 2\EE_{\cD}\left[\|\nabla f_{\xi}(x^*)\|^2\right].
	\end{eqnarray*}
\end{proof}

\begin{theorem}\label{thm:d_SGDsr}
	Assume that $f(x)$ is $\mu$-quasi strongly convex, $L$-smooth and Assumption~\ref{ass:exp_smoothness_f} holds. Then {\tt D-SGDsr} satisfies Assumption~\ref{ass:key_assumption_new} with
	\begin{gather*}
		A' = 2\cL,\quad B_1' = B_2' = 0,\quad D_1' = 2\EE_{\cD}\|\nabla f_{\xi}(x^*)\|^2,\quad \sigma_{1,k}^2 \equiv \sigma_{2,k}^2 \equiv 0\\
		\rho_1 = \rho_2 = 1,\quad C_1 = C_2 = 0,\quad D_2 = 0,\quad G = 0,\\
		F_1 = F_2 = 0,\quad D_3 = 6\gamma\tau L\EE_{\cD}\|\nabla f_{\xi}(x^*)\|^2
	\end{gather*}
	with $\gamma$ satisfying
	\begin{equation*}
		\gamma \le \min\left\{\frac{1}{8\cL}, \frac{1}{8\sqrt{L\tau\left(L\tau + 2\cL\right)}}\right\}
	\end{equation*}
	and for all $K \ge 0$
	\begin{equation*}
		\EE\left[f(\bar{x}^K) - f(x^*)\right] \le \left(1 - \frac{\gamma\mu}{2}\right)^K\frac{4\|x^0 - x^*\|^2}{\gamma} + 8\gamma(1 + 3\gamma\tau L)\EE_{\cD}\|\nabla f_\xi(x^*)\|^2
	\end{equation*}
	when $\mu > 0$ and
	\begin{equation*}
		\EE\left[f(\bar{x}^K) - f(x^*)\right] \le \frac{4\|x^0 - x^*\|^2}{\gamma K} + 8\gamma(1 + 3\gamma\tau L)\EE_{\cD}\|\nabla f_\xi(x^*)\|^2
	\end{equation*}
	when $\mu = 0$.
\end{theorem}
In other words, {\tt D-SGDsr} converges with linear rate $\cO\left(\left(\frac{\cL}{\mu} + \frac{\sqrt{L\cL\tau + L^2\tau^2}}{\mu}\right)\ln\frac{1}{\varepsilon}\right)$ to the neighbourhood of the solution when $\mu > 0$. Applying Lemma~\ref{lem:lemma2_stich} we establish the rate of convergence to $\varepsilon$-solution.
\begin{corollary}\label{cor:d_SGDsr_str_cvx_cor}
	Let the assumptions of Theorem~\ref{thm:d_SGDsr} hold and $\mu > 0$. Then after $K$ iterations of {\tt D-SGDsr} with the stepsize
	\begin{eqnarray*}
		\gamma_0 &=& \min\left\{\frac{1}{8\cL}, \frac{1}{8\sqrt{L\tau\left(L\tau + 2\cL\right)}}\right\},\quad R_0 = \|x^0 - x^*\|,\\
		\gamma &=& \min\left\{\gamma_0, \frac{\ln\left(\max\left\{2,\min\left\{\frac{R_0^2\mu^2K^2}{D_1'}, \frac{R_0^2\mu^3K^3}{3\tau LD_1'}\right\}\right\}\right)}{\mu K}\right\}
	\end{eqnarray*}	 
	we have $\EE\left[f(\bar{x}^K) - f(x^*)\right]$ of order
	\begin{equation*}
		\widetilde\cO\left(R_0^2\left(\cL + \sqrt{L^2\tau^2 + L\cL\tau}\right)\exp\left(-\frac{\mu}{\tau L}K\right) + \frac{\EE_{\cD}\|\nabla f_{\xi}(x^*)\|^2}{\mu K} + \frac{L\tau \EE_{\cD}\|\nabla f_{\xi}(x^*)\|^2}{\mu^2 K^2}\right).
	\end{equation*}
	That is, to achive $\EE\left[f(\bar{x}^K) - f(x^*)\right] \le \varepsilon$ {\tt D-SGDsr} requires
	\begin{equation*}
		\widetilde{\cO}\left(\frac{\cL + \sqrt{L^2\tau^2 + L\cL\tau}}{\mu} + \frac{\EE_{\cD}\|\nabla f_{\xi}(x^*)\|^2}{\mu\varepsilon} + \frac{\sqrt{L\tau \EE_{\cD}\|\nabla f_{\xi}(x^*)\|^2}}{\mu\sqrt{\varepsilon}}\right) \text{ iterations.}
	\end{equation*}
\end{corollary}

Applying Lemma~\ref{lem:lemma_technical_cvx} we get the complexity result in the case when $\mu = 0$.
\begin{corollary}\label{cor:d_sgdsr_cvx_cor}
	Let the assumptions of Theorem~\ref{thm:d_SGDsr} hold and $\mu = 0$. Then after $K$ iterations of {\tt D-SGDsr} with the stepsize
	\begin{eqnarray*}
		\gamma &=& \min\left\{\frac{1}{8\cL}, \frac{1}{8\sqrt{L\tau\left(L\tau + 2\cL\right)}}, \sqrt{\frac{\|x^0 - x^*\|^2}{D_1' K}}, \sqrt[3]{\frac{\|x^0 - x^*\|^2}{3L\tau D_1' K}}\right\}
		\end{eqnarray*}		
	we have $\EE\left[f(\bar{x}^K) - f(x^*)\right]$ of order
	\begin{equation*}
		\cO\left(\frac{\cL R_0^2}{K} + \frac{\sqrt{L^2\tau^2 + L\cL\tau} R_0^2}{K} + \sqrt{\frac{R_0^2 \tau \EE_\cD\|\nabla f_{\xi}(x^*)\|^2}{K}} + \frac{\sqrt[3]{LR_0^4\tau \EE_\cD\|\nabla f_{\xi}(x^*)\|^2}}{K^{\nicefrac{2}{3}}}\right)
	\end{equation*}
	where $R_0 = \|x^0 - x^*\|$. That is, to achive $\EE\left[f(\bar{x}^K) - f(x^*)\right] \le \varepsilon$ {\tt D-SGDsr} requires
	\begin{equation*}
		\cO\left(\frac{\cL R_0^2}{\varepsilon} + \frac{\sqrt{L^2\tau^2 + L\cL\tau} R_0^2}{\varepsilon} + \frac{R_0^2 \EE_\cD\|\nabla f_{\xi}(x^*)\|^2}{\varepsilon^2} + \frac{R_0^2\sqrt{L\tau \EE_\cD\|\nabla f_{\xi}(x^*)\|^2}}{\varepsilon^{\nicefrac{3}{2}}}\right)
	\end{equation*}
	iterations.
\end{corollary}

\subsection{{\tt D-LSVRG}}\label{sec:d_LSVRG}
In the same settings as in Section~\ref{sec:ec_LSVRG} we now consider a new method called {\tt D-LSVRG} which is another modification of {\tt LSVRG} that works with delayed updates. 
\begin{algorithm}[t]
   \caption{{\tt D-LSVRG}}\label{alg:d-LSVRG}
\begin{algorithmic}[1]
   \Require learning rate $\gamma>0$, initial vector $x^0 \in \R^d$
   \State Set $e_i^0 = 0$ for all $i=1,\ldots, n$   
   \For{$k=0,1,\dotsc$}
       \State Broadcast $x^{k-\tau}$ to all workers
        \For{$i=1,\dotsc,n$ in parallel}
        	\State Pick $l$ uniformly at random from $[m]$
            \State Set $g^{k-\tau}_i = \nabla f_{il}(x^{k-\tau}) - \nabla f_{il}(w_i^{k-\tau}) + \nabla f_i(w_i^{k-\tau})$
            \State $v_i^k = \gamma g_{i}^{k-\tau}$
            \State $e_i^{k+1} = e_i^k + \gamma g_i^k - v_i^k$
            \State $w_i^{k-\tau+1} = \begin{cases}x^{k-\tau},& \text{with probability } p,\\ w_i^{k-\tau},& \text{with probability } 1-p\end{cases}$
        \EndFor
        \State $e^k = \frac{1}{n}\sum_{i=1}^ne_i^k$, $g^k = \frac{1}{n}\sum_{i=1}^ng_i^k$, $v^k = \frac{1}{n}\sum_{i=1}^nv_i^k$
       \State $x^{k+1} = x^k - v^k$
   \EndFor
\end{algorithmic}
\end{algorithm}

\begin{lemma}\label{lem:second_moment_bound_d-LSVRG}
	For all $k\ge 0$, $i\in [n]$ we have
	\begin{equation}
		\EE\left[g_i^k\mid x^k\right] = \nabla f_i(x^k) \label{eq:unbiasedness_g_i^k_d-LSVRG}
	\end{equation}		
	and
	\begin{equation}
		\EE\left[\|g^k\|^2\mid x^k\right] \le 4L\left(f(x^k) - f(x^*)\right) + 2\sigma_k^2, \label{eq:second_moment_bound_d-LSVRG} 
	\end{equation}
	where $\sigma_k^2 = \frac{1}{nm}\sum_{i=1}^n\sum_{j=1}^n\|\nabla f_{ij}(w_i^k) - \nabla f_{ij}(x^*)\|^2$.
\end{lemma}
\begin{proof}
	First of all, we derive unbiasedness of $g_i^k$:
	\begin{equation*}
		\EE\left[g_i^k\mid x^k\right] = \frac{1}{m}\sum\limits_{j=1}^m\left(\nabla f_{ij}(x^k) - \nabla f_{ij}(w_i^k) + \nabla f_i(w_i^k)\right) = \nabla f_i(x^k).
	\end{equation*}
	Next, we estimate the second moment of $g^k$:
	\begin{eqnarray*}
		\EE\left[\|g^k\|^2\mid x^k\right] &=& \EE\left[\left\|\frac{1}{n}\sum\limits_{i=1}^n\left(\nabla f_{il}(x^k) - \nabla f_{il}(w_i^k) + \nabla f_{i}(w_i^k)\right)\right\|^2\right]\\
		&=& \EE\left[\left\|\frac{1}{n}\sum\limits_{i=1}^n\left(\nabla f_{il}(x^k)  - \nabla f_{il}(w_i^k) + \nabla f_{i}(w_i^k) - \nabla f_i(x^*)\right)\right\|^2\right]\\
		&\overset{\eqref{eq:a_b_norm_squared}}{\le}& \frac{2}{n}\sum\limits_{i=1}^n\EE\left[\|\nabla f_{il}(x^k) - \nabla f_{il}(x^*)\|^2\mid x^k\right]\\
		&&\quad + \frac{2}{n}\sum\limits_{i=1}^n\EE\left[\left\|\nabla f_{il}(w_i^k)- \nabla f_{il}(x^*) - \left(\nabla f_{i}(w_i^k) - \nabla f_i(x^*)\right)\right\|^2\mid x^k\right]\\
		&\overset{\eqref{eq:variance_decomposition}}{\le}& \frac{2}{nm}\sum\limits_{i=1}^n\sum\limits_{j=1}^m\|\nabla f_{ij}(x^k) - \nabla f_{ij}(x^*)\|^2 + \frac{2}{n}\EE\left[\left\|\nabla f_{il}(w_i^k)- \nabla f_{il}(x^*)\right\|^2\mid x^k\right]\\
		&\overset{\eqref{eq:L_smoothness_cor_3}}{\le}& \frac{4L}{nm}\sum\limits_{i=1}^n\sum\limits_{j=1}^mD_{f_{ij}}(x^k,x^*) + \frac{2}{nm}\sum\limits_{i=1}^n\sum\limits_{j=1}^m\|\nabla f_{ij}(w_i^k) - \nabla f_{ij}(x^*)\|^2\\
		&=& 4L\left(f(x^k) - f(x^*)\right) + 2\sigma_k^2.
	\end{eqnarray*}
\end{proof}

\begin{lemma}\label{lem:sigma_k+1_bound_d-LSVRG}
	For all $k\ge 0$, $i\in [n]$ we have
	\begin{equation}
		\EE\left[\sigma_{k+1}^2\mid x^k\right] \le (1-p)\sigma_k^2 + 2Lp\left(f(x^k) - f(x^*)\right), \label{eq:sigma_k+1_d-LSVRG} 
	\end{equation}
	where $\sigma_k^2 = \frac{1}{nm}\sum_{i=1}^n\sum_{j=1}^n\|\nabla f_{ij}(w_i^k) - \nabla f_{ij}(x^*)\|^2$.
\end{lemma}
\begin{proof}
	The proof is identical to the proof of Lemma~\ref{lem:sigma_k+1_bound_ec-LSVRG}.
\end{proof}

\begin{theorem}\label{thm:d_LSVRG}
	Assume that $f(x)$ is $\mu$-quasi strongly convex and functions $f_{ij}$ are convex and $L$-smooth for all $i\in[n],j\in[m]$. Then {\tt D-LSVRG} satisfies Assumption~\ref{ass:key_assumption_new} with
	\begin{gather*}
		A' = 2L,\quad B_1' = 0,\quad B_2' = 2,\quad D_1' = 0,\quad \sigma_{2,k}^2 = \sigma_k^2 = \frac{1}{nm}\sum\limits_{i=1}^n\sum\limits_{j=1}^m\|\nabla f_{ij}(w_i^{k}) - \nabla f_{ij}(x^*)\|^2,\\
		\sigma_{1,k}^2 \equiv 0,\quad \rho_1 = 1,\quad\rho_2 = p,\quad C_1 = 0,\quad C_2 = Lp,\quad D_2 = 0,\\
		G = 0,\quad F_1 = 0,\quad F_2 = \frac{12\gamma^2L \tau(2+p)}{p},\quad D_3 = 0
	\end{gather*}
	with $\gamma$ satisfying
	\begin{equation*}
		\gamma \le \min\left\{\frac{3}{56L}, \frac{1}{8L\sqrt{\tau\left(2+\tau + \nicefrac{4}{(1-p)}\right)}}\right\}, \quad M_2 = \frac{8}{3p}
	\end{equation*}
	and for all $K \ge 0$
	\begin{equation*}
		\EE\left[f(\bar x^K) - f(x^*)\right] \le \left(1 - \min\left\{\frac{\gamma\mu}{2},\frac{p}{4}\right\}\right)^K\frac{4(T^0 + \gamma F_2 \sigma_0^2)}{\gamma}
	\end{equation*}	
	when $\mu > 0$ and
	\begin{equation*}
		\EE\left[f(\bar x^K) - f(x^*)\right] \le \frac{4(T^0 + \gamma F_2 \sigma_0^2)}{\gamma K}
	\end{equation*}		
	when $\mu = 0$, where $T^k \eqdef \|\tx^k - x^*\|^2 + M_2\gamma^2 \sigma_k^2$.
\end{theorem}

In other words, {\tt D-LSVRG} converges with linear rate $\cO\left(\left(\frac{1}{p} + \kappa\sqrt{\tau\left(\tau + \frac{1}{(1-p)}\right)}\right)\ln\frac{1}{\varepsilon}\right)$ to the exact solution when $\mu > 0$. If $m\ge 2$ then taking $p = \frac{1}{m}$ we get that in expectation the sample complexity of one iteration of {\tt D-LSVRG} is $\cO(1)$ gradients calculations per node as for {\tt D-SGDsr} with standard sampling and the rate of convergence to the exact solution becomes $\cO\left(\left(m + \kappa\tau\right)\ln\frac{1}{\varepsilon}\right)$.

Applying Lemma~\ref{lem:lemma_technical_cvx} we get the complexity result in the case when $\mu = 0$.
\begin{corollary}\label{cor:d_lsvrg_cvx_cor}
	Let the assumptions of Theorem~\ref{thm:d_LSVRG} hold and $\mu = 0$. Then after $K$ iterations of {\tt D-LSVRG} with the stepsize
	\begin{eqnarray*}
		\gamma &=& \min\left\{\frac{3}{56L}, \frac{1}{8L\sqrt{\tau\left(2+\tau + \nicefrac{4}{(1-p)}\right)}}, \sqrt{\frac{\|x^0 - x^*\|^2}{M_2\sigma_0^2}}, \sqrt[3]{\frac{\|x^0 - x^*\|^2p}{12L\tau(2+p)\sigma_0^2}}\right\}
	\end{eqnarray*}		
	and $p = \frac{1}{m}$, $m\ge 2$ we have $\EE\left[f(\bar{x}^K) - f(x^*)\right]$ of order
	\begin{equation*}
		\cO\left(\frac{L\tau R_0^2}{K} + \frac{\sqrt{R_0^2m\sigma_0^2}}{K} + \frac{\sqrt[3]{R_0^4 L\tau \sigma_0^2}}{K}\right)
	\end{equation*}
	where $R_0 = \|x^0 - x^*\|$. That is, to achive $\EE\left[f(\bar{x}^K) - f(x^*)\right] \le \varepsilon$ {\tt D-LSVRG} requires
	\begin{equation*}
		\cO\left(\frac{L\tau R_0^2}{\varepsilon} + \frac{\sqrt{R_0^2m\sigma_0^2}}{\varepsilon} + \frac{\sqrt[3]{R_0^4 L\tau \sigma_0^2}}{\varepsilon}\right)
	\end{equation*}
	iterations.
\end{corollary}

\subsection{{\tt D-QLSVRG}}\label{sec:d_qLSVRG}
In this section we add a quantization to {\tt D-LSVRG}.  
\begin{algorithm}[t]
   \caption{{\tt D-QLSVRG}}\label{alg:d-qLSVRG}
\begin{algorithmic}[1]
   \Require learning rate $\gamma>0$, initial vector $x^0 \in \R^d$
   \State Set $e_i^0 = 0$ for all $i=1,\ldots, n$   
   \For{$k=0,1,\dotsc$}
       \State Broadcast $x^{k-\tau}$ to all workers
        \For{$i=1,\dotsc,n$ in parallel}
        	\State Pick $l$ uniformly at random from $[m]$
            \State Set $\hat g^{k-\tau}_i = \nabla f_{il}(x^{k-\tau}) - \nabla f_{il}(w_i^{k-\tau}) + \nabla f_i(w_i^{k-\tau})$
			\State Set $g_i^{k-\tau} = Q(\hat g_i^{k-\tau})$ (quantization is performed independently from other nodes)        
            \State $v_i^k = \gamma g_{i}^{k-\tau}$
            \State $e_i^{k+1} = e_i^k + \gamma g_i^k - v_i^k$
            \State $w_i^{k-\tau+1} = \begin{cases}x^{k-\tau},& \text{with probability } p,\\ w_i^{k-\tau},& \text{with probability } 1-p\end{cases}$
        \EndFor
        \State $e^k = \frac{1}{n}\sum_{i=1}^ne_i^k$, $g^k = \frac{1}{n}\sum_{i=1}^ng_i^k$, $v^k = \frac{1}{n}\sum_{i=1}^nv_i^k$
       \State $x^{k+1} = x^k - v^k$
   \EndFor
\end{algorithmic}
\end{algorithm}

\begin{lemma}\label{lem:second_moment_bound_d-qLSVRG}
	For all $k\ge 0$, $i\in [n]$ we have
	\begin{equation*}
		\EE\left[g_i^k\mid x^k\right] = \nabla f_i(x^k)
	\end{equation*}		
	and
	\begin{equation*}
		\EE\left[\|g^k\|^2\mid x^k\right] \le 4L\left(1 + \frac{2\omega}{n}\right)\left(f(x^k) - f(x^*)\right) + 2\left(1 + \frac{2\omega}{n}\right)\sigma_k^2 + \frac{2\omega}{n^2}\sum\limits_{i=1}^n\|\nabla f_i(x^*)\|^2,  
	\end{equation*}
	where $\sigma_k^2 = \frac{1}{nm}\sum_{i=1}^n\sum_{j=1}^n\|\nabla f_{ij}(w_i^k) - \nabla f_{ij}(x^*)\|^2$.
\end{lemma}
\begin{proof}
	First of all, we derive unbiasedness of $g_i^k$:
	\begin{eqnarray*}
		\EE\left[g_i^k\mid x^k\right] &\overset{\eqref{eq:tower_property}}{=}& \EE\left[\EE_{Q}\left[Q(\hat g_i^k)\right]\mid x^k\right] \overset{\eqref{eq:quantization_def}}{=} \EE\left[\hat g_i^k\mid x^k\right]\\
		&=& \frac{1}{m}\sum\limits_{j=1}^m\left(\nabla f_{ij}(x^k) - \nabla f_{ij}(w_i^k) + \nabla f_i(w_i^k)\right) = \nabla f_i(x^k).
	\end{eqnarray*}
	Next, we estimate the second moment of $g^k$:
	\begin{eqnarray*}
		\EE_{Q}\left[\|g^k\|^2\right] &=& \EE_{Q}\left[\left\|\frac{1}{n}\sum\limits_{i=1}^nQ(\hat g_i^k)\right\|^2\right]\\
		&\overset{\eqref{eq:variance_decomposition}}{=}& \EE_{Q}\left[\left\|\frac{1}{n}\sum\limits_{i=1}^n\left(Q(\hat g_i^k)-\hat g_i^k\right)\right\|^2\right] + \left\|\frac{1}{n}\sum\limits_{i=1}^n\hat g_i^k\right\|^2.
	\end{eqnarray*}
	Since quantization on nodes is performed independently we can decompose the first term from the last row of the previous inequality into the sum of variances:
	\begin{eqnarray*}
		\EE_{Q}\left[\|g^k\|^2\right] &=& \frac{1}{n^2}\sum\limits_{i=1}^n\EE_{Q}\left\|Q(\hat g_i^k)-\hat g_i^k\right\|^2 + \left\|\frac{1}{n}\sum\limits_{i=1}^n\hat g_i^k\right\|^2\\
		&\overset{\eqref{eq:quantization_def}}{\le}& \frac{\omega}{n^2}\sum\limits_{i=1}^n\|\hat g_i^k\|^2 + \left\|\frac{1}{n}\sum\limits_{i=1}^n\left(\hat g_i^k-\nabla f_i(x^*)\right)\right\|^2\\
		&\overset{\eqref{eq:a_b_norm_squared}}{\le}& \left(1 + \frac{2\omega}{n}\right)\frac{1}{n}\sum\limits_{i=1}^n\|\hat g_i^k - \nabla f_i(x^*)\|^2 + \frac{2\omega}{n^2}\sum\limits_{i=1}^n\|\nabla f_i(x^*)\|^2.
	\end{eqnarray*}
	Taking conditional mathematical expectation $\EE\left[\cdot\mid x^k\right]$ from the both sides of previous inequality we get
	\begin{eqnarray*}
		\EE\left[\|g^k\|^2\mid x^k\right] &\le& \left(1 + \frac{2\omega}{n}\right)\frac{2}{n}\sum\limits_{i=1}^n\EE\left[\|\nabla f_{il}(x^k) - \nabla f_{il}(x^*)\|^2\mid x^k\right]\\
		&&\quad + \left(1 + \frac{2\omega}{n}\right)\frac{2}{n}\sum\limits_{i=1}^n\EE\left[\left\|\nabla f_{il}(w_i^k) - \nabla f_{il}(x^*) - \left(\nabla f_i(w_i^k) - \nabla f_i(x^*)\right) \right\|^2\mid x^k\right]\\
		&&\quad  + \frac{2\omega}{n^2}\sum\limits_{i=1}^n\|\nabla f_i(x^*)\|^2\\
		&\le& \left(1 + \frac{2\omega}{n}\right)\frac{2}{nm}\sum\limits_{i=1}^n\sum\limits_{j=1}^m\|\nabla f_{ij}(x^k) - \nabla f_{ij}(x^*)\|^2\\
		&&\quad + \left(1 + \frac{2\omega}{n}\right)\frac{2}{n}\sum\limits_{i=1}^n\EE\left[\left\|\nabla f_{il}(w_i^k) - \nabla f_{il}(x^*)\right\|^2\mid x^k\right] + \frac{2\omega}{n^2}\sum\limits_{i=1}^n\|\nabla f_i(x^*)\|^2\\
		&\overset{\eqref{eq:L_smoothness_cor_3}}{\le}& \left(1 + \frac{2\omega}{n}\right)\frac{4L}{nm}\sum\limits_{i=1}^n\sum\limits_{j=1}^mD_{f_{ij}}(x^k,x^*)\\
		&&\quad + \left(1 + \frac{2\omega}{n}\right)\frac{2}{nm}\sum\limits_{i=1}^n\sum\limits_{j=1}^m\left\|\nabla f_{ij}(w_i^k) - \nabla f_{ij}(x^*)\right\|^2 + \frac{2\omega}{n^2}\sum\limits_{i=1}^n\|\nabla f_i(x^*)\|^2\\
		&=& 4L\left(1 + \frac{2\omega}{n}\right)\left(f(x^k) - f(x^*)\right) + 2\left(1 + \frac{2\omega}{n}\right)\sigma_k^2 + \frac{2\omega}{n^2}\sum\limits_{i=1}^n\|\nabla f_i(x^*)\|^2.
	\end{eqnarray*}
\end{proof}

\begin{lemma}\label{lem:sigma_k+1_bound_d-qLSVRG}
	For all $k\ge 0$, $i\in [n]$ we have
	\begin{equation}
		\EE\left[\sigma_{k+1}^2\mid x^k\right] \le (1-p)\sigma_k^2 + 2Lp\left(f(x^k) - f(x^*)\right), \label{eq:sigma_k+1_d-qLSVRG} 
	\end{equation}
	where $\sigma_k^2 = \frac{1}{nm}\sum_{i=1}^n\sum_{j=1}^n\|\nabla f_{ij}(w_i^k) - \nabla f_{ij}(x^*)\|^2$.
\end{lemma}
\begin{proof}
	The proof is identical to the proof of Lemma~\ref{lem:sigma_k+1_bound_ec-LSVRG}.
\end{proof}

\begin{theorem}\label{thm:d_qLSVRG}
	Assume that $f(x)$ is $\mu$-quasi strongly convex and functions $f_{ij}$ are convex and $L$-smooth for all $i\in[n],j\in[m]$. Then {\tt D-QLSVRG} satisfies Assumption~\ref{ass:key_assumption_new} with
	\begin{gather*}
		A' = 2L\left(1+\frac{2\omega}{n}\right),\quad B_1' = 0,\quad B_2' = 2\left(1+\frac{2\omega}{n}\right),\quad D_1' = \frac{2\omega}{n^2}\sum\limits_{i=1}^n\|\nabla f_i(x^*)\|^2,\quad \sigma_{1,0}^2 \equiv 0,\\
		\sigma_{2,k}^2 = \sigma_k^2 = \frac{1}{nm}\sum\limits_{i=1}^n\sum\limits_{j=1}^m\|\nabla f_{ij}(w_i^{k}) - \nabla f_{ij}(x^*)\|^2,\quad \rho_1 = 1,\quad \rho_2 = p,\quad C_2 = Lp,\quad D_2 = 0,\\
		C_1 = 0,\quad G = 0,\quad F_1 = 0,\quad F_2 = \frac{12\gamma^2L \tau\left(1+\frac{2\omega}{n}\right)(2+p)}{p},\quad D_3 = \frac{6\gamma\tau L\omega}{n^2}\sum\limits_{i=1}^n\|\nabla f_i(x^*)\|^2
	\end{gather*}
	with $\gamma$ satisfying
	\begin{equation*}
		\gamma \le \min\left\{\frac{3}{56L(1+\nicefrac{2\omega}{n})}, \frac{1}{8L\sqrt{\tau\left(\tau+2\left(1+\nicefrac{2\omega}{n}\right)\left(1+\nicefrac{2}{(1-p)}\right)\right)}}\right\}, \quad M_2 = \frac{8\left(1+\frac{2\omega}{n}\right)}{3p}
	\end{equation*}
	and for all $K \ge 0$
	\begin{equation*}
		\EE\left[f(\bar x^K) - f(x^*)\right] \le \left(1 - \min\left\{\frac{\gamma\mu}{2},\frac{p}{4}\right\}\right)^K\frac{4(T^0 + \gamma F_2 \sigma_0^2)}{\gamma} + 4\gamma\left(D_1' + D_3\right)
	\end{equation*}	
	when $\mu > 0$ and
	\begin{equation*}
		\EE\left[f(\bar x^K) - f(x^*)\right] \le \frac{4(T^0 + \gamma F_2 \sigma_0^2)}{\gamma K} + 4\gamma\left(D_1' + D_3\right)
	\end{equation*}
	when $\mu = 0$, where $T^k \eqdef \|\tx^k - x^*\|^2 + M_2\gamma^2 \sigma_k^2$.
\end{theorem}

In other words, {\tt D-QLSVRG} converges with linear rate 
$$
\cO\left(\left(\frac{1}{p} + \kappa\left(1+\frac{\omega}{n}\right) + \kappa\sqrt{\tau\left(\tau + \left(1+\frac{\omega}{n}\right)\left(1+\frac{1}{(1-p)}\right)\right)}\right)\ln\frac{1}{\varepsilon}\right)
$$
to neighbourhood the solution when $\mu > 0$. If $m\ge 2$ then taking $p = \frac{1}{m}$ we get that in expectation the sample complexity of one iteration of {\tt D-QLSVRG} is $\cO(1)$ gradients calculations per node as for {\tt D-QSGDsr} with standard sampling and the rate of convergence to the neighbourhood of the solution becomes
$$
\cO\left(\left(m + \kappa\left(1+\frac{\omega}{n}\right) + \kappa\sqrt{\tau\left(\tau + \frac{\omega}{n}\right)}\right)\ln\frac{1}{\varepsilon}\right).
$$
Applying Lemma~\ref{lem:lemma2_stich} we establish the rate of convergence to $\varepsilon$-solution.
\begin{corollary}\label{cor:d_QLSVRG_str_cvx_cor}
	Let the assumptions of Theorem~\ref{thm:d_qLSVRG} hold, $f_{\xi}(x)$ are convex for each $\xi$ and $\mu > 0$. Then after $K$ iterations of {\tt D-QLSVRG} with the stepsize
	\begin{eqnarray*}
		\gamma_0 &=& \min\left\{\frac{3}{56L(1+\nicefrac{2\omega}{n})}, \frac{1}{8L\sqrt{\tau\left(\tau+2\left(1+\nicefrac{2\omega}{n}\right)\left(1+\nicefrac{2}{(1-p)}\right)\right)}}\right\},\quad R_0 = \|x^0 - x^*\|,\\
		\gamma &=& \min\left\{\gamma_0, \frac{\ln\left(\max\left\{2,\min\left\{\frac{R_0^2\mu^2K^2}{D_1'}, \frac{R_0^2\mu^3K^3}{3\tau LD_1'}\right\}\right\}\right)}{\mu K}\right\}
	\end{eqnarray*}
	and $p = \frac{1}{m}$, $m \ge 2$ we have $\EE\left[f(\bar{x}^K) - f(x^*)\right]$ of order
	\begin{equation*}
		 \widetilde\cO\left(LR_0^2\left(1+\frac{\omega}{n}+\sqrt{\tau\left(\tau + \frac{\omega}{n}\right)}\right)\exp\left(-\frac{\mu}{L\left(1+\frac{\omega}{n}+\sqrt{\tau\left(\tau + \frac{\omega}{n}\right)}\right)}K\right) + \frac{D_1'}{\mu K} + \frac{L\tau D_1'}{\mu^2 K^2}\right).
	\end{equation*}
	That is, to achive $\EE\left[f(\bar{x}^K) - f(x^*)\right] \le \varepsilon$ {\tt D-QLSVRG} requires
	\begin{equation*}
		\widetilde{\cO}\left(\frac{L}{\mu}\left(1+\frac{\omega}{n}\right) + \frac{L}{\mu}\sqrt{\tau\left(\tau + \frac{\omega}{n}\right)} + \frac{D_1'}{\mu\varepsilon} + \frac{\sqrt{L\tau D_1'}}{\mu\sqrt{\varepsilon}}\right) \text{ iterations.}
	\end{equation*}
\end{corollary}

Applying Lemma~\ref{lem:lemma_technical_cvx} we get the complexity result in the case when $\mu = 0$.
\begin{corollary}\label{cor:d_QLSVRG_cvx_cor}
	Let the assumptions of Theorem~\ref{thm:d_qLSVRG} hold and $\mu = 0$. Then after $K$ iterations of {\tt D-QLSVRG} with the stepsize
	\begin{eqnarray*}
		\gamma_0 &=& \min\left\{\frac{3}{56L(1+\nicefrac{2\omega}{n})}, \frac{1}{8L\sqrt{\tau\left(\tau+2\left(1+\nicefrac{2\omega}{n}\right)\left(1+\nicefrac{2}{(1-p)}\right)\right)}}\right\},\quad R_0 = \|x^0 - x^*\|,\\	
		\gamma &=& \min\left\{\gamma_0, \sqrt{\frac{R_0^2}{M_2\sigma_{0}^2}}, \sqrt[3]{\frac{R_0^2p}{12L\tau\left(1+\frac{2\omega}{n}\right)(2+p)}}, \sqrt{\frac{R_0^2}{D_1' K}}, \sqrt[3]{\frac{R_0^2}{3L\tau D_1' K}}\right\}	
	\end{eqnarray*}		
	and $p = \frac{1}{m}$, $m\ge 2$ we have $\EE\left[f(\bar{x}^K) - f(x^*)\right]$ of order
	\begin{gather*}
		\cO\left(\frac{LR_0^2\left(1+\frac{\omega}{n}+\sqrt{\tau\left(\tau + \frac{\omega}{n}\right)}\right)}{K} + \frac{\sqrt{R_0^2 m\left(1+\frac{\omega}{n}\right)\sigma_0^2}}{K} + \frac{\sqrt[3]{R_0^4 L\tau m\left(1+\frac{\omega}{n}\right)}}{K}\right)\\
		 + \cO\left(\sqrt{\frac{R_0^2 D_1'}{K}} + \frac{\sqrt[3]{LR_0^4\tau D_1'}}{K^{\nicefrac{2}{3}}}\right).
	\end{gather*}
	That is, to achive $\EE\left[f(\bar{x}^K) - f(x^*)\right] \le \varepsilon$ {\tt D-QLSVRG} requires
	\begin{gather*}
		\cO\left(\frac{LR_0^2\left(1+\frac{\omega}{n}+\sqrt{\tau\left(\tau + \frac{\omega}{n}\right)}\right)}{\varepsilon} + \frac{\sqrt{R_0^2 m\left(1+\frac{\omega}{n}\right)\sigma_0^2}}{\varepsilon} + \frac{\sqrt[3]{R_0^4 L\tau m\left(1+\frac{\omega}{n}\right)}}{\varepsilon}\right)\\
		 + \cO\left(\frac{R_0^2 D_1'}{\varepsilon^2} + \frac{R_0^2\sqrt{L\tau D_1'}}{\varepsilon^{\nicefrac{3}{2}}}\right)
	\end{gather*}
	iterations.
\end{corollary}

\subsection{{\tt D-QLSVRGstar}}\label{sec:d_qLSVRGstar}
Now we assume that $i$-th node has an access to $\nabla f_i(x^*)$ and modify {\tt D-QLSVRG} in order to get convergence asymptotically to the exact optimum.   
\begin{algorithm}[t]
   \caption{{\tt D-QLSVRGstar}}\label{alg:d-qLSVRGstar}
\begin{algorithmic}[1]
   \Require learning rate $\gamma>0$, initial vector $x^0 \in \R^d$
   \State Set $e_i^0 = 0$ for all $i=1,\ldots, n$   
   \For{$k=0,1,\dotsc$}
       \State Broadcast $x^{k-\tau}$ to all workers
        \For{$i=1,\dotsc,n$ in parallel}
        	\State Pick $l$ uniformly at random from $[m]$
            \State Set $\hat g^{k-\tau}_i = \nabla f_{il}(x^{k-\tau}) - \nabla f_{il}(w_i^{k-\tau}) + \nabla f_i(w_i^{k-\tau})$
			\State Set $g_i^{k-\tau} = Q(\hat g_i^{k-\tau} - \nabla f_i(x^*))$ (quantization is performed independently from other nodes)        
            \State $v_i^k = \gamma g_{i}^{k-\tau}$
            \State $e_i^{k+1} = e_i^k + \gamma g_i^k - v_i^k$
            \State $w_i^{k-\tau+1} = \begin{cases}x^{k-\tau},& \text{with probability } p,\\ w_i^{k-\tau},& \text{with probability } 1-p\end{cases}$
        \EndFor
        \State $e^k = \frac{1}{n}\sum_{i=1}^ne_i^k$, $g^k = \frac{1}{n}\sum_{i=1}^ng_i^k$, $v^k = \frac{1}{n}\sum_{i=1}^nv_i^k$
       \State $x^{k+1} = x^k - v^k$
   \EndFor
\end{algorithmic}
\end{algorithm}

\begin{lemma}\label{lem:second_moment_bound_d-qLSVRGstar}
	For all $k\ge 0$, $i\in [n]$ we have
	\begin{equation}
		\EE\left[g^k\mid x^k\right] = \nabla f(x^k) \label{eq:unbiasedness_g^k_d-qLSVRGstar}
	\end{equation}		
	and
	\begin{equation}
		\EE\left[\|g^k\|^2\mid x^k\right] \le 2L\left(1 + \frac{\omega}{n}\right)\left(f(x^k) - f(x^*)\right) + 2\left(1 + \frac{\omega}{n}\right)\sigma_k^2, \label{eq:second_moment_bound_d-qLSVRGstar} 
	\end{equation}
	where $\sigma_k^2 = \frac{1}{nm}\sum_{i=1}^n\sum_{j=1}^n\|\nabla f_{ij}(w_i^k) - \nabla f_{ij}(x^*)\|^2$.
\end{lemma}
\begin{proof}
	First of all, we derive unbiasedness of $g_i^k$:
	\begin{eqnarray*}
		\EE\left[g^k\mid x^k\right] &\overset{\eqref{eq:tower_property}}{=}& \EE\left[\EE_{Q}\left[\frac{1}{n}\sum\limits_{i=1}^nQ(\hat g_i^k - \nabla f_i(x^*))\right]\mid x^k\right] \overset{\eqref{eq:quantization_def}}{=} \EE\left[\frac{1}{n}\sum\limits_{i=1}^n\left(\hat g_i^k-\nabla f_i(x^*)\right)\mid x^k\right]\\
		&=& \frac{1}{nm}\sum\limits_{i=1}^n\sum\limits_{j=1}^m\left(\nabla f_{ij}(x^k) - \nabla f_{ij}(w_i^k) + \nabla f_i(w_i^k)\right) = \nabla f(x^k).
	\end{eqnarray*}
	Next, we estimate the second moment of $g^k$:
	\begin{eqnarray*}
		\EE_{Q}\left[\|g^k\|^2\right] &=& \EE_{Q}\left[\left\|\frac{1}{n}\sum\limits_{i=1}^nQ(\hat g_i^k - \nabla f_i(x^*))\right\|^2\right]\\
		&\overset{\eqref{eq:variance_decomposition}}{=}& \EE_{Q}\left[\left\|\frac{1}{n}\sum\limits_{i=1}^n\left(Q(\hat g_i^k- \nabla f_i(x^*))-\left(\hat g_i^k- \nabla f_i(x^*)\right)\right)\right\|^2\right] + \left\|\frac{1}{n}\sum\limits_{i=1}^n\hat g_i^k- \nabla f_i(x^*)\right\|^2.
	\end{eqnarray*}
	Since quantization on nodes is performed independently we can decompose the first term from the last row of the previous inequality into the sum of variances:
	\begin{eqnarray*}
		\EE_{Q}\left[\|g^k\|^2\right] &=& \frac{1}{n^2}\sum\limits_{i=1}^n\EE_{Q}\left\|Q(\hat g_i^k- \nabla f_i(x^*))-\left(\hat g_i^k- \nabla f_i(x^*)\right)\right\|^2 + \left\|\frac{1}{n}\sum\limits_{i=1}^n\hat g_i^k- \nabla f_i(x^*)\right\|^2\\
		&\overset{\eqref{eq:quantization_def}}{\le}& \frac{\omega}{n^2}\sum\limits_{i=1}^n\|\hat g_i^k- \nabla f_i(x^*)\|^2 + \left\|\frac{1}{n}\sum\limits_{i=1}^n\left(\hat g_i^k-\nabla f_i(x^*)\right)\right\|^2\\
		&\overset{\eqref{eq:a_b_norm_squared}}{\le}& \left(1 + \frac{\omega}{n}\right)\frac{1}{n}\sum\limits_{i=1}^n\|\hat g_i^k - \nabla f_i(x^*)\|^2.
	\end{eqnarray*}
	Taking conditional mathematical expectation $\EE\left[\cdot\mid x^k\right]$ from the both sides of previous inequality and using the bound
	\begin{equation*}
		\frac{1}{n}\sum\limits_{i=1}^n\EE\left[\|\hat g_i^k - \nabla f_i(x^*)\|^2\mid x^k\right] \le 4L\left(f(x^k)- f(x^*)\right) + 2\sigma_k^2
	\end{equation*}
	implicitly obtained in the proof of Lemma~\ref{lem:second_moment_bound_d-qLSVRG} we get \eqref{eq:second_moment_bound_d-qLSVRGstar}.
\end{proof}

\begin{lemma}\label{lem:sigma_k+1_bound_d-qLSVRGstar}
	For all $k\ge 0$, $i\in [n]$ we have
	\begin{equation}
		\EE\left[\sigma_{k+1}^2\mid x^k\right] \le (1-p)\sigma_k^2 + 2Lp\left(f(x^k) - f(x^*)\right), \label{eq:sigma_k+1_d-qLSVRGstar} 
	\end{equation}
	where $\sigma_k^2 = \frac{1}{nm}\sum_{i=1}^n\sum_{j=1}^n\|\nabla f_{ij}(w_i^k) - \nabla f_{ij}(x^*)\|^2$.
\end{lemma}
\begin{proof}
	The proof is identical to the proof of Lemma~\ref{lem:sigma_k+1_bound_ec-LSVRG}.
\end{proof}

\begin{theorem}\label{thm:d_qLSVRGstar}
	Assume that $f(x)$ is $\mu$-quasi strongly convex and functions $f_{ij}$ are convex and $L$-smooth for all $i\in[n],j\in[m]$. Then {\tt D-QLSVRGstar} satisfies Assumption~\ref{ass:key_assumption_new} with
	\begin{gather*}
		A' = 2L\left(1+\frac{2\omega}{n}\right),\quad B_1' = 0,\quad B_2' = 2\left(1+\frac{2\omega}{n}\right),\quad D_1' = 0,\quad \sigma_{1,0}^2 \equiv 0,\\
		\sigma_{2,k}^2 = \sigma_k^2 = \frac{1}{nm}\sum\limits_{i=1}^n\sum\limits_{j=1}^m\|\nabla f_{ij}(w_i^{k}) - \nabla f_{ij}(x^*)\|^2,\quad \rho_1 = 1,\quad \rho_2 = p,\quad C_2 = Lp,\quad D_2 = 0,\\
		C_1 = 0,\quad G = 0,\quad F_1 = 0,\quad F_2 = \frac{12\gamma^2L \tau\left(1+\frac{2\omega}{n}\right)(2+p)}{p},\quad D_3 = 0
	\end{gather*}
	with $\gamma$ satisfying
	\begin{equation*}
		\gamma \le \min\left\{\frac{3}{56L(1+\nicefrac{2\omega}{n})}, \frac{1}{8L\sqrt{\tau\left(\tau+2\left(1+\nicefrac{2\omega}{n}\right)\left(1+\nicefrac{2}{(1-p)}\right)\right)}}\right\}, \quad M_2 = \frac{8\left(1+\frac{2\omega}{n}\right)}{3p}
	\end{equation*}
	and for all $K \ge 0$
	\begin{equation*}
		\EE\left[f(\bar x^K) - f(x^*)\right] \le \left(1 - \min\left\{\frac{\gamma\mu}{2},\frac{p}{4}\right\}\right)^K\frac{4(T^0 + \gamma F_2 \sigma_0^2)}{\gamma}
	\end{equation*}	
	when $\mu > 0$ and
	\begin{equation*}
		\EE\left[f(\bar x^K) - f(x^*)\right] \le \frac{4(T^0 + \gamma F_2 \sigma_0^2)}{\gamma K}
	\end{equation*}
	when $\mu = 0$, where $T^k \eqdef \|\tx^k - x^*\|^2 + M_2\gamma^2 \sigma_k^2$.
\end{theorem}

In other words, {\tt D-QLSVRGstar} converges with linear rate 
$$
\cO\left(\left(\frac{1}{p} + \kappa\left(1+\frac{\omega}{n}\right) + \kappa\sqrt{\tau\left(\tau + \left(1+\frac{\omega}{n}\right)\left(1+\frac{1}{(1-p)}\right)\right)}\right)\ln\frac{1}{\varepsilon}\right)
$$
to the exact solution when $\mu > 0$. If $m\ge 2$ then taking $p = \frac{1}{m}$ we get that in expectation the sample complexity of one iteration of {\tt D-QLSVRGstar} is $\cO(1)$ gradients calculations per node as for {\tt D-QSGDsr} with standard sampling and the rate of convergence to the exact solution becomes
$$
\cO\left(\left(m + \kappa\left(1+\frac{\omega}{n}\right) + \kappa\sqrt{\tau\left(\tau + \frac{\omega}{n}\right)}\right)\ln\frac{1}{\varepsilon}\right).
$$

Applying Lemma~\ref{lem:lemma_technical_cvx} we get the complexity result in the case when $\mu = 0$.
\begin{corollary}\label{cor:d_QLSVRGstar_cvx_cor}
	Let the assumptions of Theorem~\ref{thm:d_qLSVRGstar} hold and $\mu = 0$. Then after $K$ iterations of {\tt D-QLSVRGstar} with the stepsize
	\begin{eqnarray*}
		\gamma_0 &=& \min\left\{\frac{3}{56L(1+\nicefrac{2\omega}{n})}, \frac{1}{8L\sqrt{\tau\left(\tau+2\left(1+\nicefrac{2\omega}{n}\right)\left(1+\nicefrac{2}{(1-p)}\right)\right)}}\right\},\quad R_0 = \|x^0 - x^*\|,\\	
		\gamma &=& \min\left\{\gamma_0, \sqrt{\frac{R_0^2}{M_2\sigma_{0}^2}}, \sqrt[3]{\frac{R_0^2p}{12L\tau\left(1+\frac{2\omega}{n}\right)(2+p)}}\right\}	
	\end{eqnarray*}		
	and $p = \frac{1}{m}$, $m\ge 2$ we have $\EE\left[f(\bar{x}^K) - f(x^*)\right]$ of order
	\begin{gather*}
		\cO\left(\frac{LR_0^2\left(1+\frac{\omega}{n}+\sqrt{\tau\left(\tau + \frac{\omega}{n}\right)}\right)}{K} + \frac{\sqrt{R_0^2 m\left(1+\frac{\omega}{n}\right)\sigma_0^2}}{K} + \frac{\sqrt[3]{R_0^4 L\tau m\left(1+\frac{\omega}{n}\right)}}{K}\right).
	\end{gather*}
	That is, to achive $\EE\left[f(\bar{x}^K) - f(x^*)\right] \le \varepsilon$ {\tt D-QLSVRGstar} requires
	\begin{gather*}
		\cO\left(\frac{LR_0^2\left(1+\frac{\omega}{n}+\sqrt{\tau\left(\tau + \frac{\omega}{n}\right)}\right)}{\varepsilon} + \frac{\sqrt{R_0^2 m\left(1+\frac{\omega}{n}\right)\sigma_0^2}}{\varepsilon} + \frac{\sqrt[3]{R_0^4 L\tau m\left(1+\frac{\omega}{n}\right)}}{\varepsilon}\right)
	\end{gather*}
	iterations.
\end{corollary}

However, such convergence guarantees are obtained under very restrictive assumption: the method requires to know vectors $\nabla f_i(x^*)$.

\subsection{{\tt D-LSVRG-DIANA}}\label{sec:d_LSVRG-diana}
In the setup of Section~\ref{sec:d_LSVRG} we construct a new method with delayed updates and quantization called {\tt D-LSVRG-DIANA} which does not require to know $\nabla f_i(x^*)$ and has linear convergence to the exact solution.

\begin{algorithm}[t]
   \caption{{\tt D-LSVRG-DIANA}}\label{alg:d-LSVRG-diana}
\begin{algorithmic}[1]
   \Require learning rates $\gamma>0$, $\alpha \in (0,1]$, initial vectors $x^0, h_1^0,\ldots, h_n^0 \in \R^d$
	\State Set $e_i^0 = 0$ for all $i=1,\ldots, n$   
	\State Set $h^0 = \frac{1}{n}\sum_{i=1}^n h_i^0$   
   \For{$k=0,1,\dotsc$}
       \State Broadcast $x^{k-\tau}$ to all workers
        \For{$i=1,\dotsc,n$ in parallel}
			\State Pick $l$ uniformly at random from $[m]$
            \State Set $\hat g^{k-\tau}_i = \nabla f_{il}(x^{k-\tau}) - \nabla f_{il}(w_i^{k-\tau}) + \nabla f_i(w_i^{k-\tau})$           
			\State $\hat \Delta_i^{k-\tau} = Q(\hat g^{k-\tau}_i - h_i^{k-\tau})$ (quantization is performed independently from other nodes)           
			\State $g_i^{k-\tau} = h_i^{k-\tau} + \hat \Delta_i^{k-\tau}$            
            \State $v_i^k = \gamma g_i^{k-\tau}$
            \State $e_i^{k+1} = e_i^k + \gamma g_i^k - v_i^k$
            \State $h_i^{k-\tau+1} = h_i^{k-\tau} + \alpha \hat \Delta_i^{k-\tau}$
        \EndFor
        \State $e^k = \frac{1}{n}\sum_{i=1}^ne_i^k$, $g^k = \frac{1}{n}\sum_{i=1}^ng_i^k = h^k + \frac{1}{n}\sum\limits_{i=1}^n\hat\Delta_i^k$, $v^k = \frac{1}{n}\sum_{i=1}^nv_i^k = \gamma h^{k-\tau} + \frac{\gamma}{n}\sum_{i=1}^n \hat\Delta_i^{k-\tau}$
        \State $h^{k-\tau+1} = \frac{1}{n}\sum\limits_{i=1}^n h_i^{k-\tau+1} = h^{k-\tau} + \alpha\frac{1}{n}\sum\limits_{i=1}^n \hat \Delta_i^{k-\tau}$
       \State $x^{k+1} = x^k - v^k$
   \EndFor
\end{algorithmic}
\end{algorithm}

\begin{lemma}\label{lem:d_LSVRG-diana_second_moment_bound}
	Assume that $f_{ij}(x)$ is convex and $L$-smooth for all $i=1,\ldots,n$, $j=1,\ldots,m$. Then, for all $k\ge 0$ we have
	\begin{eqnarray}
		\EE\left[g^k\mid x^k\right] &=& \nabla f(x^k), \label{eq:d_LSVRG-diana_unbiasedness}\\
		\EE\left[\|g^k\|^2\mid x^k\right] &\le& 4L\left(1+\frac{2\omega}{n}\right)\left(f(x^k) - f(x^*)\right) + \frac{2\omega}{n}\sigma_{1,k}^2 + 2\left(1+\frac{2\omega}{n}\right)\sigma_{2,k}^2 \label{eq:d_LSVRG-diana_second_moment_bound}
	\end{eqnarray}
	where $\sigma_{1,k}^2 = \frac{1}{n}\sum_{i=1}^n\|h_i^k - \nabla f(x^*)\|^2$ and $\sigma_{2,k}^2 = \frac{1}{nm}\sum_{i=1}^n\sum_{j=1}^m\|\nabla f_{ij}(w_i^k) - \nabla f_{ij}(x^*)\|^2$.
\end{lemma}
\begin{proof}
	First of all, we show unbiasedness of $g^k$:
	\begin{eqnarray*}
		\EE\left[g^k\mid x^k\right] &\overset{\eqref{eq:tower_property}}{=}& h^k + \frac{1}{n}\sum\limits_{i=1}^n\EE\left[\EE_Q\left[\hat \Delta_i^k\right]\mid x^k\right] \overset{\eqref{eq:quantization_def}}{=} h^k + \frac{1}{n}\sum\limits_{i=1}^n\EE\left[\hat g_i^k - h_i^k\mid x^k\right]	\\
		&=& \frac{1}{nm}\sum\limits_{i=1}^n\sum\limits_{j=1}^m\left(\nabla f_{ij}(x^k) - \nabla f_{ij}(w_i^k) + \nabla f_i(w_i^k)\right) = \nabla f(x^k).
	\end{eqnarray*}
	Next, we derive the upper bound for the second moment of $g^k$:
	\begin{eqnarray*}
		\EE_{Q}\left[\|g^k\|^2\right] &=& \EE_{Q}\left[\left\|h^k + \frac{1}{n}\sum\limits_{i=1}^n\hat{\Delta}_i^k\right\|^2\right]\\
		&\overset{\eqref{eq:variance_decomposition}}{=}& \EE_{Q}\left[\left\|\frac{1}{n}\sum\limits_{i=1}^n\left(\hat{\Delta}_i^k - \hat g_i^k + h_i^k\right)\right\|^2\right] + \left\|\frac{1}{n}\sum\limits_{i=1}^n\hat g_i^k\right\|^2.
	\end{eqnarray*}
	Since quantization on nodes is performed independently we can decompose the first term from the last row of the previous inequality into the sum of variances:
	\begin{eqnarray*}
		\EE_{Q}\left[\|g^k\|^2\right] &\le& \frac{1}{n^2}\sum\limits_{i=1}^n\EE_{Q}\left[\|\hat{\Delta}_i^k - \hat g_i^k + h_i^k\|^2\right] + \left\|\frac{1}{n}\sum\limits_{i=1}^n\left(\hat g_i^k - \nabla f_i(x^*)\right)\right\|^2\\
		&\overset{\eqref{eq:quantization_def},\eqref{eq:a_b_norm_squared}}{\le}& \frac{\omega}{n^2}\sum\limits_{i=1}^n\|\hat g_i^k - h_i^k\|^2 + \frac{1}{n}\sum\limits_{i=1}^n\|\hat g_i^k - \nabla f_i(x^*)\|^2\\
		&\overset{\eqref{eq:a_b_norm_squared}}{\le}& \left(1+\frac{2\omega}{n}\right)\frac{1}{n}\sum\limits_{i=1}^n\|\hat g_i^k - \nabla f_i(x^*)\|^2 + \frac{2\omega}{n^2}\sum\limits_{i=1}^n\|h_i^k - f_i(x^*)\|^2.
	\end{eqnarray*}
	Taking mathematical expectation $\EE\left[\cdot\mid x^k\right]$ from the both sides of the previous inequality and using the bound
	\begin{equation*}
		\frac{1}{n}\sum\limits_{i=1}^n\EE\left[\|\hat g_i^k - \nabla f_i(x^*)\|^2\mid x^k\right] \le 4L\left(f(x^k)- f(x^*)\right) + \frac{2}{nm}\sum\limits_{i=1}^n\sum\limits_{j=1}^m\|\nabla f_{ij}(w_i^k) - \nabla f_{ij}(x^*)\|^2
	\end{equation*}
	implicitly obtained in the proof of Lemma~\ref{lem:second_moment_bound_d-qLSVRG} we get \eqref{eq:d_LSVRG-diana_second_moment_bound}.
\end{proof}

\begin{lemma}\label{lem:d_LSVRG-diana_sigma_k+1_bound}
	Assume that $\alpha \le \nicefrac{1}{(\omega+1)}$. Then, for all $k\ge 0$ we have
	\begin{equation*}
		\EE\left[\sigma_{1,k+1}^2\mid x^k\right] \le (1 - \alpha)\sigma_{1,k}^2 + 6L\alpha(f(x^k) - f(x^*)) + 2\alpha\sigma_{2,k}^2,
	\end{equation*}
	\begin{equation*}
		\EE\left[\sigma_{2,{k+1}}^2\mid x^k\right] \le (1 - p)\sigma_{k,2}^2 + 2Lp\left(f(x^k)-f(x^*)\right)
	\end{equation*}
	where $\sigma_{1,k}^2 = \frac{1}{n}\sum_{i=1}^n\|h_i^k - \nabla f_i(x^*)\|^2$ and $\sigma_{2,k}^2= \frac{1}{nm}\sum_{i=1}^n\sum_{j=1}^m\|\nabla f_{ij}(w_i^k) - \nabla f_{ij}(x^*)\|^2$.
\end{lemma}
\begin{proof}
	The proof is identical to the proof of Lemma~\ref{lem:ec_LSVRG-diana_sigma_k+1_bound}.
\end{proof}

\begin{theorem}\label{thm:d_LSVRG-diana}
	Assume that $f_{ij}(x)$ is convex and $L$-smooth for all $i=1,\ldots, n$, $j=1,\ldots,m$ and $f(x)$ is $\mu$-quasi strongly convex. Then {\tt D-LSVRG-DIANA} satisfies Assumption~\ref{ass:key_assumption_new} with
	\begin{gather*}
		A' = 2L\left(1+\frac{2\omega}{n}\right),\quad B_1' = \frac{2\omega}{n},\quad B_2' = 2\left(1+\frac{2\omega}{n}\right),\quad D_1' = 0,\\
		\sigma_{1,k}^2 = \frac{1}{n}\sum\limits_{i=1}^n\|h_i^k - \nabla f_i(x^*)\|^2,\quad \sigma_{2,k}^2 = \frac{1}{nm}\sum\limits_{i=1}^n\sum_{j=1}^m\|\nabla f_{ij}(w_i^k) - \nabla f_{ij}(x^*)\|^2,\\
		\rho_1 = \alpha,\quad \rho_2 = p,\quad C_1 = 3L\alpha,\quad C_2 = Lp,\quad D_2 = 0,\quad G = 2,\\
		F_1 = \frac{12\gamma^2 L\omega\tau(2+\alpha)}{n\alpha},\quad F_2 = \frac{12\gamma^2\tau L(2+p)}{p}\left(\frac{4\omega}{n(1-\alpha)} + 1 + \frac{2\omega}{n}\right),\quad D_3 = 0
	\end{gather*}
	with $\gamma$ and $\alpha$ satisfying
	\begin{equation*}
		\gamma \le \min\left\{\frac{1}{8L\left(\frac{37}{9} + \frac{24\omega}{3n}\right)}, \frac{1}{8L\sqrt{\tau\left(2+\tau+ \frac{4}{1-p}+\frac{4\omega}{n}\left(1+\frac{3}{1-\alpha}+\frac{2}{1-p}+\frac{4}{(1-\alpha)(1-p)}\right)\right)}}\right\},
	\end{equation*}
	\begin{equation*}
		\alpha \le \frac{1}{\omega+1},\quad M_1 = \frac{8\omega}{3n\alpha},\quad M_2 = \frac{8\left(7 + \frac{6\omega}{n}\right)}{9p}.
	\end{equation*}
	and for all $K \ge 0$
	\begin{equation*}
		\EE\left[f(\bar x^K) - f(x^*)\right] \le \left(1 - \min\left\{\frac{\gamma\mu}{2},\frac{\alpha}{4},\frac{p}{4}\right\}\right)^K\frac{4(T^0 + \gamma F_1 \sigma_{1,0}^2 + \gamma F_2 \sigma_{2,0}^2)}{\gamma}
	\end{equation*}	
	when $\mu > 0$ and
	\begin{equation*}
		\EE\left[f(\bar x^K) - f(x^*)\right] \le \frac{4(T^0 + \gamma F_1 \sigma_{1,0}^2 + \gamma F_2 \sigma_{2,0}^2)}{\gamma K}
	\end{equation*}	
	when $\mu = 0$, where $T^k \eqdef \|\tx^k - x^*\|^2 + M_1\gamma^2 \sigma_{1,k}^2 + M_2\gamma^2 \sigma_{2,k}^2$.
\end{theorem}
In other words, if $m\ge 2$, $p = \nicefrac{1}{m}$, $\alpha = \min\left\{\frac{1}{\omega+1},\frac{1}{2}\right\}$ and 
\begin{equation*}
		\gamma \le \min\left\{\frac{1}{8L\left(\frac{37}{9} + \frac{24\omega}{3n}\right)}, \frac{1}{8L\sqrt{\tau\left(2+\tau+ \frac{4}{1-p}+\frac{4\omega}{n}\left(1+\frac{3}{1-\alpha}+\frac{2}{1-p}+\frac{4}{(1-\alpha)(1-p)}\right)\right)}}\right\},
\end{equation*}
{\tt D-LSVRG-DIANA} converges with the linear rate
\begin{equation*}
	\cO\left(\left(\omega + m + \kappa\left(1+\frac{\omega}{n}\right) + \kappa\sqrt{\tau\left(\tau+\frac{\omega}{n}\right)}\right)\ln\frac{1}{\varepsilon}\right)
\end{equation*}
to the exact solution when $\mu > 0$.

Applying Lemma~\ref{lem:lemma_technical_cvx} we get the complexity result in the case when $\mu = 0$.
\begin{corollary}\label{cor:d_LSVRG_DIANA_cvx_cor}
	Let the assumptions of Theorem~\ref{thm:d_LSVRG-diana} hold and $\mu = 0$. Then after $K$ iterations of {\tt D-LSVRG-DIANA} with the stepsize
	\begin{eqnarray*}
		\gamma_0 &=& \min\left\{\frac{1}{8L\left(\frac{37}{9} + \frac{24\omega}{3n}\right)}, \frac{1}{8L\sqrt{\tau\left(2+\tau+ \frac{4}{1-p}+\frac{4\omega}{n}\left(1+\frac{3}{1-\alpha}+\frac{2}{1-p}+\frac{4}{(1-\alpha)(1-p)}\right)\right)}}\right\},\\	
		\gamma &=& \min\left\{\gamma_0, \sqrt{\frac{R_0^2}{M_1\sigma_{1,0}^2 + M_2\sigma_{2,0}^2}}, \sqrt[3]{\frac{R_0^2}{12\tau L\left(\frac{\omega(2+\alpha)}{n\alpha} + \frac{2+p}{p}\left(1+\frac{2\omega}{n}+\frac{4\omega}{n(1-\alpha)}\right)\right)}}\right\},
	\end{eqnarray*}		
	where $R_0 = \|x^0 - x^*\|$, $\alpha = \min\left\{\frac{1}{\omega+1},\frac{1}{2}\right\}$ and $p = \frac{1}{m}$, $m\ge 2$ we have $\EE\left[f(\bar{x}^K) - f(x^*)\right]$ of order
	\begin{gather*}
		\cO\left(\frac{LR_0^2\left(1+\frac{\omega}{n}+\sqrt{\tau\left(\tau + \frac{\omega}{n}\right)}\right)}{K} + \frac{\sqrt{R_0^2\omega(\omega+1)\sigma_{1,0}^2}}{\sqrt{n}K} + \frac{\sqrt{R_0^2m\left(1+\frac{\omega}{n}\right)\sigma_{2,0}^2}}{K}\right)\\
		+\cO\left(\frac{\sqrt[3]{R_0^4\tau L \omega(\omega+1)\sigma_{1,0}^2}}{\sqrt[3]{n}K} + \frac{\sqrt[3]{R_0^4\tau L m\left(1+\frac{\omega}{n}\right)\sigma_{2,0}^2}}{K}\right)
	\end{gather*}
	That is, to achive $\EE\left[f(\bar{x}^K) - f(x^*)\right] \le \varepsilon$ {\tt D-LSVRG-DIANA} requires
	\begin{gather*}
		\cO\left(\frac{LR_0^2\left(1+\frac{\omega}{n}+\sqrt{\tau\left(\tau + \frac{\omega}{n}\right)}\right)}{\varepsilon} + \frac{\sqrt{R_0^2\omega(\omega+1)\sigma_{1,0}^2}}{\sqrt{n}\varepsilon} + \frac{\sqrt{R_0^2m\left(1+\frac{\omega}{n}\right)\sigma_{2,0}^2}}{\varepsilon}\right)\\
		+\cO\left(\frac{\sqrt[3]{R_0^4\tau L \omega(\omega+1)\sigma_{1,0}^2}}{\sqrt[3]{n}\varepsilon} + \frac{\sqrt[3]{R_0^4\tau L m\left(1+\frac{\omega}{n}\right)\sigma_{2,0}^2}}{\varepsilon}\right)
	\end{gather*}
	iterations.
\end{corollary}

\begin{table*}[!t]
\caption{The parameters for which the methods from Tables~\ref{tbl:special_cases2_ef} and \ref{tbl:special_cases_delayed_methods} satisfy Assumption~\ref{ass:key_assumption_new}. The meaning of the expressions appearing in the table, as well as their justification is defined in details in the Sections~\ref{sec:special_cases_ef} and \ref{sec:special_cases2}. Symbols: $\varepsilon = $ error tolerance; $\delta = $ contraction factor of compressor $\cC$; $\omega = $ variance parameter of compressor $\cQ$; $\kappa = \nicefrac{L}{\mu}$; $\cL =$ expected smoothness constant; $\sigma_*^2 = $ variance of the stochastic gradients in the solution; $\zeta_*^2 =$ average of $\|\nabla f_i(x^*)\|^2$; $\sigma^2 =$ average of the uniform bounds for the variances of stochastic gradients of workers.}
\label{tbl:special_cases-parameters}
\begin{center}
 \tiny
\begin{adjustbox}{angle=90}
\begin{tabular}{|c|c|c|c|c|c|c|c|c|c|c|}
\hline
 Method &   $A'$ & $B_1'$ & $B_2'$ & $\rho_1$ & $\rho_2$ & $C_1$ & $C_2$ & $F_1,\quad F_2$ & $G$ & $D_1'$, $D_2$, $D_3$\\
 \hline
\hline
 {\tt EC-SGDsr}   &  $2\cL $ & $0$ & $0$ & $1$ & $1$ & $0$ & $0$ & $0,\quad 0$ & $0$ & $\frac{2\sigma_*^2}{n},\quad 0,\quad\frac{6L\gamma}{\delta}\left(\frac{4\zeta_*^2}{\delta}+3\sigma_*^2\right)$\\
 \hline
 {\tt EC-SGD}   &  $2L $ & $0$ & $0$ & $1$ & $1$ & $0$ & $0$ & $0,\quad 0$ & $0$ & $\frac{2\sigma_*^2}{n},\quad 0,\quad\frac{12L\gamma}{\delta}\left(\frac{2\zeta_*^2}{\delta}+\sigma_*^2\right)$\\
 \hline
 {\tt EC-GDstar}   &  $L $ & $0$ & $0$ & $1$ & $1$ & $0$ & $0$ & $0,\quad 0$ & $0$ & $0,\quad 0,\quad 0$\\
 \hline
 {\tt EC-SGD-DIANA}   &  $L $ & $0$ & $0$ & $\alpha$ & $1$ & $L\alpha$ & $0$ & $\frac{96L\gamma^2}{\delta^2\alpha(1-\eta)},\quad 0$ & $0$ & \makecell{$\frac{\sigma^2}{n},\quad \alpha^2(\omega+1)\sigma^2,$\\ $\frac{6L\gamma}{\delta}\left(\frac{4\alpha(\omega+1)}{\delta}+1\right)\sigma^2$}\\
 \hline
 {\tt EC-SGDsr-DIANA}   &  $2\cL $ & $0$ & $0$ & $\alpha$ & $1$ & $2\alpha(3\cL + 4L)$ & $0$ & $\frac{96L\gamma^2}{\delta^2\alpha(1-\eta)},\quad 0$ & $0$ & \makecell{$\frac{2\sigma_*^2}{n},\quad \alpha^2(\omega+1)\sigma_*^2,$\\ $\frac{18L\gamma}{\delta}\left(\frac{4\alpha(\omega+1)}{\delta}+1\right)\sigma_*^2$}\\
 \hline
 {\tt EC-LSVRG}   &  $2L $ & $0$ & $2$ & $1$ & $p$ & $0$ & $Lp$ & $0,\quad \frac{72L\gamma^2}{\delta p(1-\eta)}$ & $0$ & $0,\quad 0,\quad \frac{24L\gamma}{\delta^2}\zeta_*^2$\\
 \hline
  {\tt EC-LSVRGstar}   &  $2L $ & $0$ & $2$ & $1$ & $p$ & $0$ & $Lp$ & $0,\quad \frac{48L\gamma^2}{\delta p}$ & $0$ & $0,\quad 0,\quad 0$\\
 \hline
  {\tt EC-LSVRG-DIANA}   &  $2L $ & $0$ & $2$ & $\alpha$ & $p$ & $3L\alpha$ & $Lp$ & \makecell{$\frac{24L\gamma^2\left(\frac{4}{\delta}+3\right)}{\delta\alpha(1-\eta)},$\\ $\frac{24L\gamma^2\left(\frac{4}{1-\alpha}\left(\frac{4}{\delta}+3\right)+3\right)}{\delta p(1-\eta)}$} & $2$ & $0,\quad 0,\quad 0$\\
 \hline
 {\tt D-SGDsr}   &  $2\cL $ & $0$ & $0$ & $1$ & $1$ & $0$ & $0$ & $0,\quad 0$ & $0$ & $\frac{2\sigma_*^2}{n},\quad 0,\quad \frac{6L\tau\gamma\sigma_*^2}{n}$\\
 \hline
 {\tt D-SGD}   &  $2L $ & $0$ & $0$ & $1$ & $1$ & $0$ & $0$ & $0,\quad 0$ & $0$ & $\frac{2\sigma_*^2}{n},\quad 0,\quad \frac{6L\tau\gamma\sigma_*^2}{n}$\\
 \hline
 {\tt D-QSGD}   &  $L\left(1+\frac{2\omega}{n}\right) $ & $0$ & $0$ & $1$ & $1$ & $0$ & $0$ & $0,\quad 0$ & $0$ & \makecell{$\frac{(\omega+1)\sigma^2}{n}+\frac{2\omega\zeta_*^2}{n},\quad 0,$\\ $\frac{3\gamma\tau L}{n}\left((\omega+1)\sigma^2 + 2\omega\zeta_*^2\right)$}\\
 \hline
 {\tt D-QSGDstar}   &  $L\left(1+\frac{\omega}{n}\right) $ & $0$ & $0$ & $1$ & $1$ & $0$ & $0$ & $0,\quad 0$ & $0$ & $\frac{(\omega+1)\sigma^2}{n},\quad 0,\quad \frac{3\gamma\tau L(\omega+1)\sigma^2}{n}$\\
 \hline
 {\tt D-QGDstar}   &  $L\left(1+\frac{\omega}{n}\right) $ & $0$ & $0$ & $1$ & $1$ & $0$ & $0$ & $0,\quad 0$ & $0$ & $0,\quad 0,\quad 0$\\
 \hline
 {\tt D-SGD-DIANA}   &  $L\left(1+\frac{2\omega}{n}\right) $ & $\frac{2\omega}{n}$ & $0$ & $\alpha$ & $1$ & $L\alpha$ & $0$ & $\frac{12\gamma^2 L\omega\tau(2+\alpha)}{n\alpha},\quad 0$ & $0$ & \makecell{$\frac{(\omega+1)\sigma^2}{n},\quad \frac{\alpha(\omega+1)\sigma^2}{n},$\\ $3\gamma\tau L\left(1+\frac{4\omega}{n}\right)\frac{(\omega+1)\sigma^2}{n}$}\\
 \hline
 {\tt D-LSVRG}   &  $2L$ & $0$ & $2$ & $1$ & $p$ & $0$ & $Lp$ & $0,\quad \frac{12\gamma^2 L\tau(2+p)}{np}$ & $0$ & $0,\quad 0,\quad 0$\\
 \hline
 {\tt D-QLSVRG}   &  $2L\left(1+\frac{2\omega}{n}\right)$ & $0$ & $2\left(1+\frac{2\omega}{n}\right)$ & $1$ & $p$ & $0$ & $Lp$ & $0,\quad \frac{12\gamma^2 L\tau\left(1+\frac{2\omega}{n}\right)\tau(2+p)}{p}$ & $0$ & \makecell{$\frac{2\omega\zeta_*^2}{n},\quad 0,$\\ $\frac{6\gamma\tau L\omega\zeta_*^2}{n}$}\\
 \hline
 {\tt D-QLSVRGstar}   &  $2L\left(1+\frac{2\omega}{n}\right)$ & $0$ & $2\left(1+\frac{2\omega}{n}\right)$ & $1$ & $p$ & $0$ & $Lp$ & $0,\quad \frac{12\gamma^2 L\left(1+\frac{2\omega}{n}\right)\tau(2+p)}{p}$ & $0$ & $0,\quad 0,\quad 0$\\
 \hline
 {\tt D-LSVRG-DIANA}   &  $2L\left(1+\frac{2\omega}{n}\right)$ & $\frac{2\omega}{n}$ & $2\left(1+\frac{2\omega}{n}\right)$ & $\alpha$ & $p$ & $3L\alpha$ & $Lp$ & \begin{tabular}{c}
 $\frac{12\gamma^2 L\omega\tau(2+\alpha)}{n\alpha}$,\\$\frac{12\gamma^2\tau L(2+p)}{p}\left(1 + \frac{2\omega(3-\alpha)}{n(1-\alpha)}\right)$
 \end{tabular} & $0$ & $0,\quad 0,\quad 0$\\
 \hline
\end{tabular}
\end{adjustbox}
\end{center}
\end{table*}

%% file: Appendix_local_sigma_k.tex
\chapter{Appendix for Chapter~\ref{ch:local_sigma_k}}\label{app:local_sigma_k}
\section{Table of Frequently Used Notation} \label{sec:notation_table}
\begin{table}[H]
\scriptsize
\caption{Summary of frequently used notation.}
\label{tbl:notation}
\begin{center}
\begin{tabular}{|c|l|c|}
\hline
\multicolumn{3}{|c|}{{\bf Main notation} }\\
\hline
\hline
  $f: \R^d \rightarrow \R$ & Objective to be minimized&  \eqref{eq:main_problem} \\ \hline
  $f_i: \R^d \rightarrow \R$ & Local objective  owned by device/worker $i$ & \eqref{eq:f_i_expectation} or \eqref{eq:f_i_sum} \\ \hline
 $x^*$ & Global optimum of \eqref{eq:main_problem};  $x^* \in \R^d$  & \\ \hline
  $d$ & Dimensionality of the problem space  & \eqref{eq:main_problem} \\ \hline
    $n$ & Number of clients/devices/nodes/workers & \eqref{eq:main_problem} \\ \hline
            $x_i^k$ & Local iterate;  $x_i^k \in \R^d$  & \eqref{eq:local_sgd_def} \\ \hline
        $g_i^k$ & Local stochastic direction;  $g_i^k \in \R^d$  & \eqref{eq:local_sgd_def} \\ \hline
           $\gamma$ & Stepsize/learning rate; $\gamma \geq 0$ & \eqref{eq:local_sgd_def} \\ \hline
         $c_k$ & Indicator of the communication; $c_k \in \{0,1\}$& \eqref{eq:local_sgd_def} \\ \hline
          $\mu$ &Strong quasi-convexity of the local objective; $\mu \geq 0$ & ~\eqref{eq:str_quasi_cvx} \\ \hline
    $L$ &  Smoothness of the local objective; $L\geq \mu$ & ~\eqref{eq:L_smoothness} \\  \hline
          $x^k$ & Virtual iterate;  $x^k \in \R^d$  & Sec~\ref{sec:main_res} \\ \hline
      $V^k$ & Discrepancy between local and virtual iterates;  $V^k \geq 0$  & Sec~\ref{sec:main_res}\\ \hline
         $    \overline{x}^K$ & Weighted average of historical iterates;  $\overline{x}^K \in \R^d$  & Thm~\ref{thm:main_result} \\ \hline
         $\zeta$ &Heterogeneity parameter; $\zeta \geq 0$ & ~\eqref{eq:bounded_data_dissimilarity} \\ \hline
       $\tau$ &  Size of the fixed local loop $\tau \geq 0$ & Sec~\ref{sec:data_and_loop}  \\ \hline
        $p$ &  Probability of aggregation fixed for the random local loop $p \in [0,1] $ & Sec~\ref{sec:data_and_loop}  \\ \hline
  $a_i^k$ & Unbiased local gradient; $a_i^k \in \R^d $ & Sec~\ref{sec:local_solver}  \\ \hline
    $b_i^k$ & Local shift; $b_i^k \in \R^d $ & Sec~\ref{sec:local_solver}  \\ \hline
        $h_i^k$ & Delayed local gradient estimator used to construct $b_i^k$; $h_i^k \in \R^d $ & Sec~\ref{sec:local_solver}  \\ \hline
       $l_i^k$ & Unbiased local gradient estimator used to construct $b_i^k$; $l_i^k \in \R^d $ & Sec~\ref{sec:local_solver}  \\ \hline
       $\cL$ & Expected smoothness of local objectives; $\cL \geq 0 $ & \eqref{eq:expected_smoothness_1}  \\ \hline
      $\max L_{ij}$ & Smoothness constant of local summands; $\max L_{ij}\geq 0$ & Sec~\eqref{sec:llsvrg}  \\ \hline
     $\sigma^2$ & Averaged upper bound for the variance of local stochastic gradient  & Tab~\eqref{tbl:instances}  \\ \hline
          $\sigma^2_*$ & Averaged variance of local stochastic gradients at the solution  & Tab~\eqref{tbl:instances}  \\ \hline
           $\zeta_*^2$ &  $\eqdef \frac1n\sum_{i=1}^n \| \nabla f_i(x^*)\|^2$  & Tab~\eqref{tbl:instances}  \\         
\hline
\hline
\multicolumn{3}{|c|}{{\bf Parametric Assumptions} }\\
\hline
\hline
\makecell{ $A, A', B, B', C, C', F, F',$ \\ $G, H, D_1, D_1',D_2, D_3, \rho$} & \multicolumn{2}{|c|}{ Parameters of Assumption~\ref{ass:key_assumption} } \\
\hline
 $A_i, B_i\, D_{1,i} , \rho_i , C_i ,  D_{2,i}$ & \multicolumn{2}{|c|}{ Parameters of Assumption~\ref{ass:sigma_k_original} } \\
 \hline
  $A'_i, D_{3,i}$ & \multicolumn{2}{|c|}{ Parameters of Assumption~\ref{ass:bik} } \\
  \hline
   $\sigma^2_k, \sigma^2_{i,k}$ & \multicolumn{2}{|c|}{ Possibly random non-negative sequences from Assumptions~\ref{ass:key_assumption},~\ref{ass:sigma_k_original},~\ref{ass:hetero_second_moment} } \\
 \hline
\hline
\multicolumn{3}{|c|}{{\bf Standard} }\\
\hline
\hline
$\EE[\cdot]$ & Expectation & \\
\hline
$\EE\left[\cdot\mid x^k\right]$  & $\eqdef \EE\left[\cdot\mid x_1^k,\ldots,x_n^k\right]$; expectation conditioned on $k$-th local iterates  & \\
\hline
$D_h(x,y)$ & $\eqdef h(x)-h(y)- \langle \nabla h(y), x-y \rangle$; Bregman distance of $x,y$ w.r.t. $h$  & As~\ref{ass:sigma_k_original}\\
 \hline
\end{tabular}
 \end{center}
\end{table}


\section{Extra Experiments \label{sec:extra_exp}}

\subsection{Missing Details from Section~\ref{sec:exp} and an Extra Figure\label{sec:extralog}}

In Section~\ref{sec:exp} we study the effect of  local variance reduction on the communication complexity of local methods. We consider the regularized logistic regression objective, i.e., we choose
\[
f_i(x) \eqdef \frac1m \sum_{j=1}^m \log \left(1+\exp\left( \langle a_{(i-1)m+j}, x \rangle\cdot  b_{(i-1)m+j}\right) \right) + \frac{\mu}{2}\| x\|^2,
\]
where $a_{j}\in \R^d, b_j\in \{-1, 1\}$ for $j\leq nm$ are the training data and labels.

\paragraph{Number of the clients.} We select a different number of clients for each dataset in order to capture a variety of scenarios. See Table~\ref{tbl:ns} for details.

 \begin{table}[H]
 \caption{Number of clients per dataset (Figures~\ref{fig:sgd_svrg_hom} and~\ref{fig:sgd_svrg_het}). }
\label{tbl:ns}
\begin{center}
\small
\begin{tabular}{|c|c|c|c|}
\hline
{\bf Dataset}  & $n$ & {\bf \# datapoints} ($=mn$) & $d$   \\
 \hline
  \hline
\texttt{a1a} & 5 & 1 605	& 123  \\ \hline
\texttt{mushrooms} & 12 & 8 124 & 112   \\ \hline
\texttt{phishing} & 11   & 11 055	&	68 \\ \hline
\texttt{madelon} & 50 & 2 000& 500 \\ \hline
\texttt{duke} & 4  &44 & 7 129  \\ \hline
\texttt{w2a} & 10  &3 470 & 300  \\ \hline
\end{tabular}
\end{center}
\end{table}

\begin{figure}[H]
\centering
\begin{minipage}{0.3\textwidth}
  \centering
\includegraphics[width =  \textwidth ]{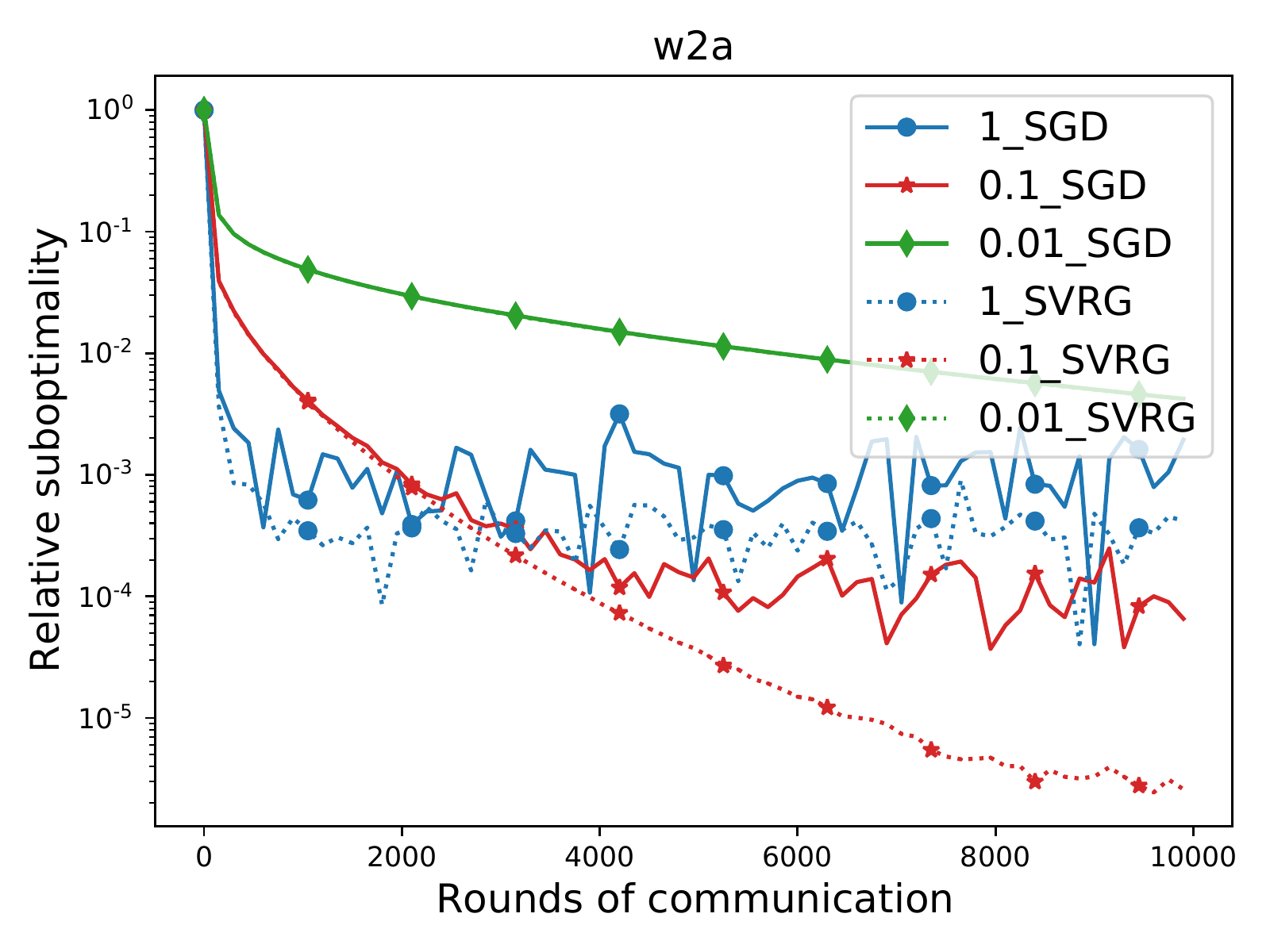}
\end{minipage}
\begin{minipage}{0.3\textwidth}
  \centering
\includegraphics[width =  \textwidth ]{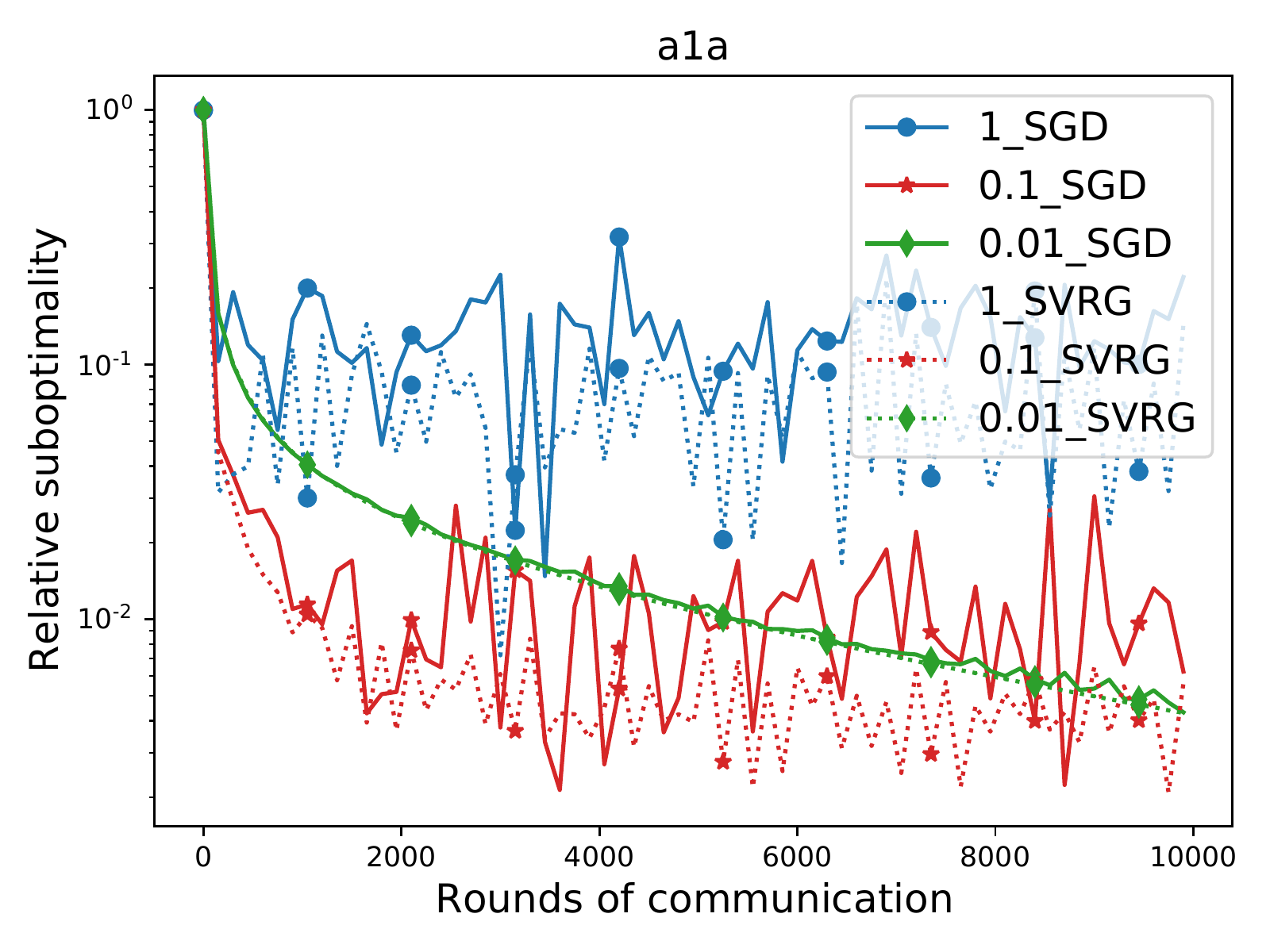}
\end{minipage}
\begin{minipage}{0.3\textwidth}
  \centering
\includegraphics[width =  \textwidth ]{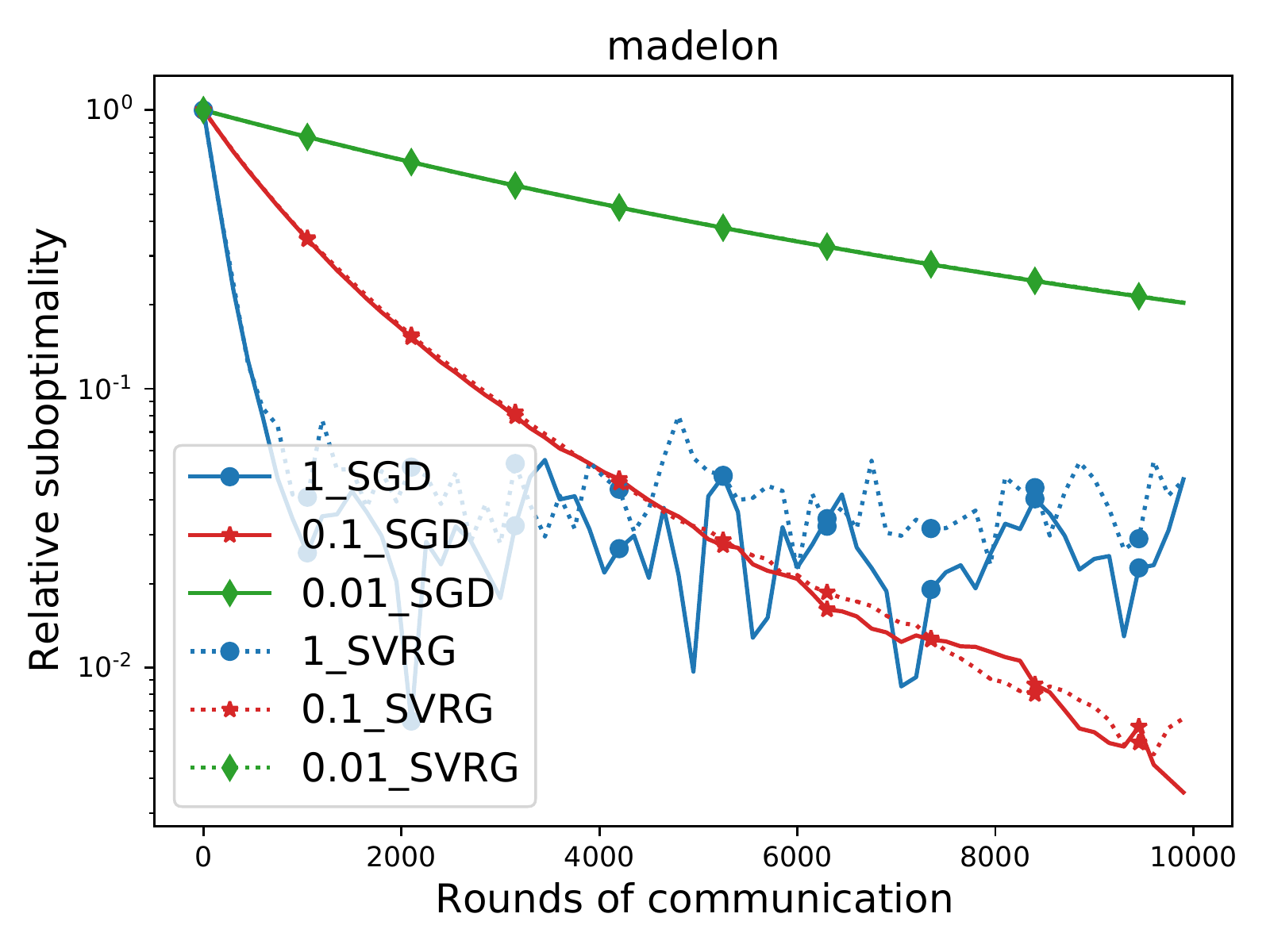}
\end{minipage}
\\
\begin{minipage}{0.3\textwidth}
  \centering
\includegraphics[width =  \textwidth ]{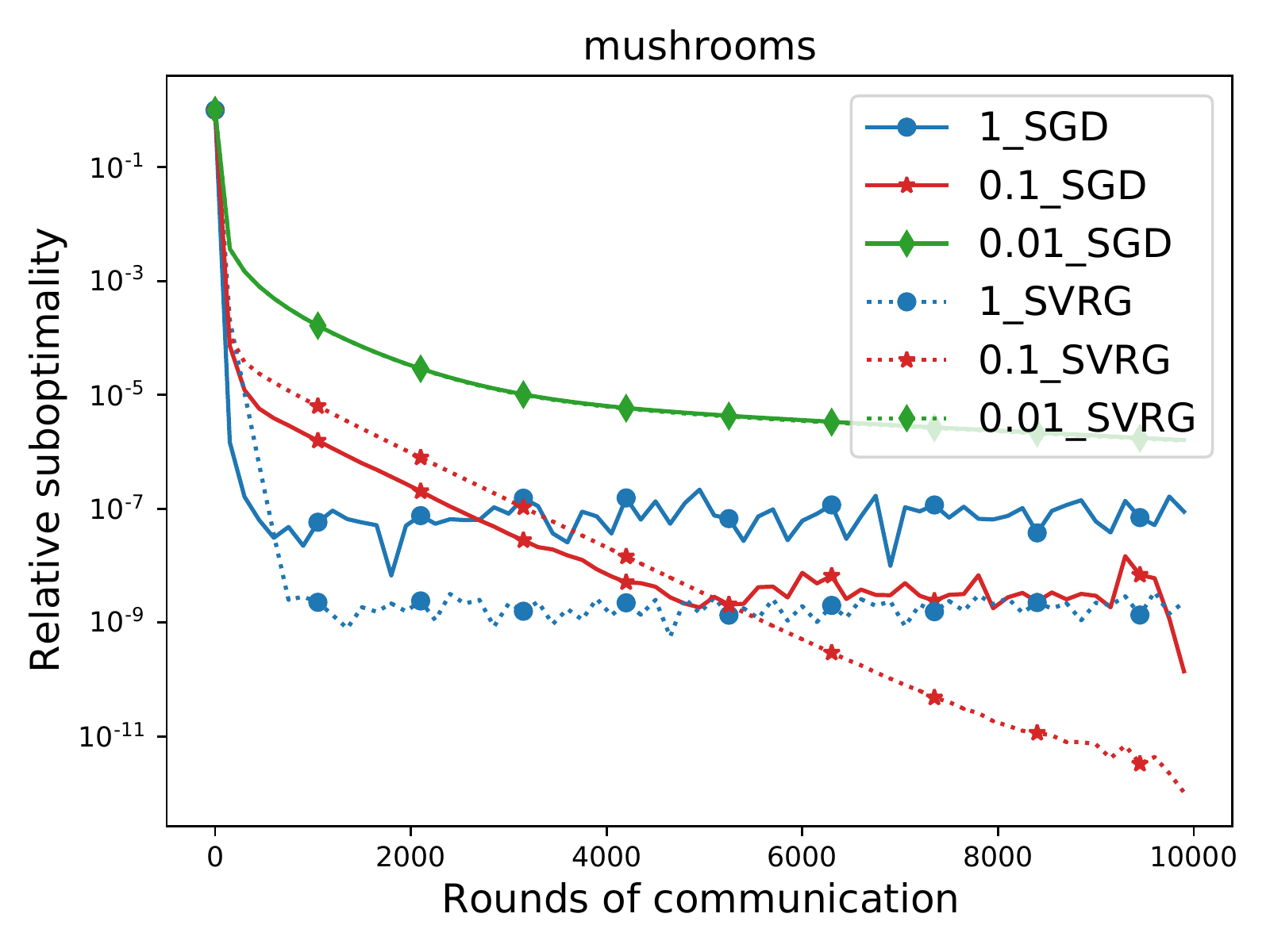}
\end{minipage}
\begin{minipage}{0.3\textwidth}
  \centering
\includegraphics[width =  \textwidth ]{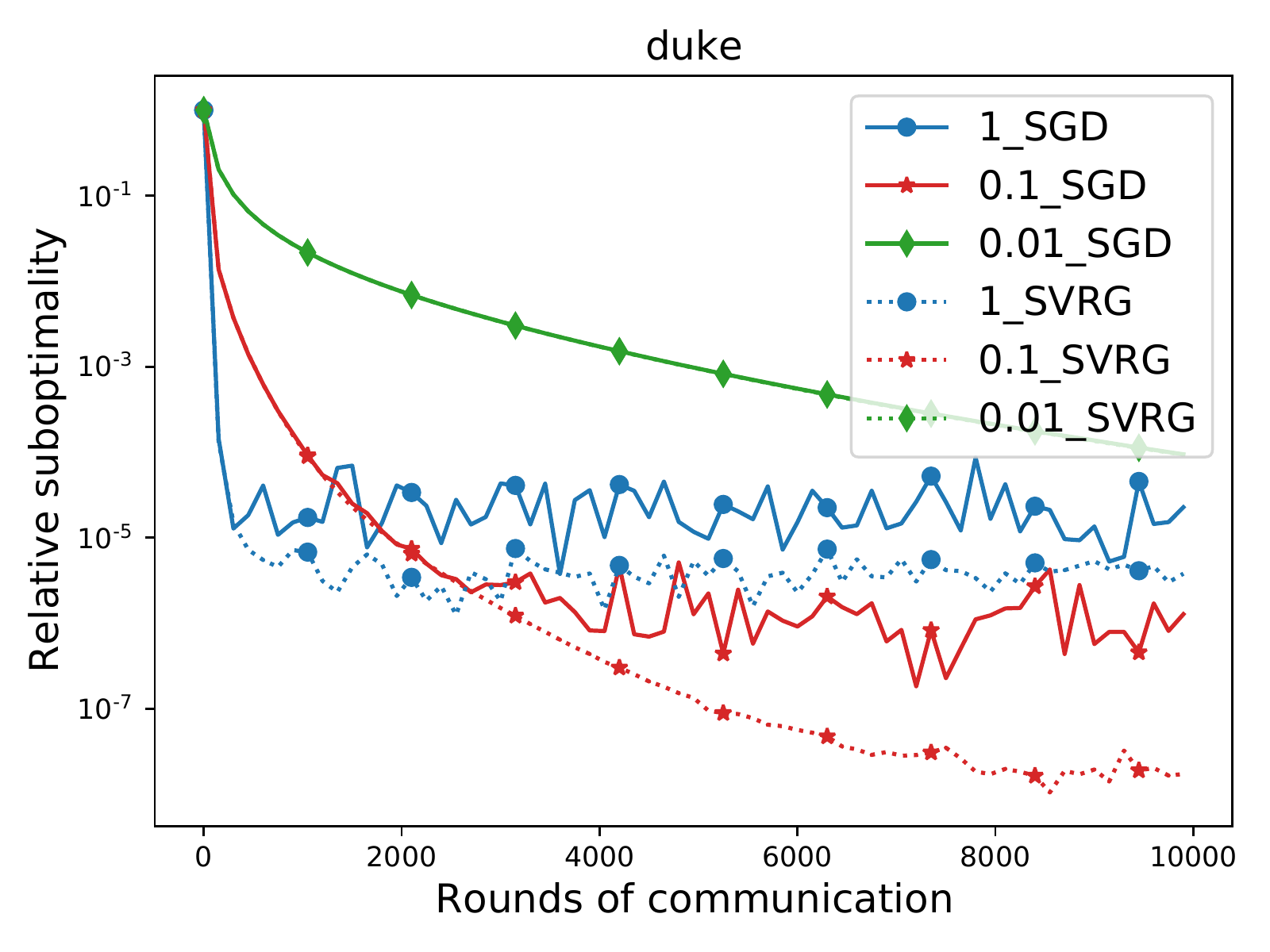}
\end{minipage}
\begin{minipage}{0.3\textwidth}
  \centering
\includegraphics[width =  \textwidth ]{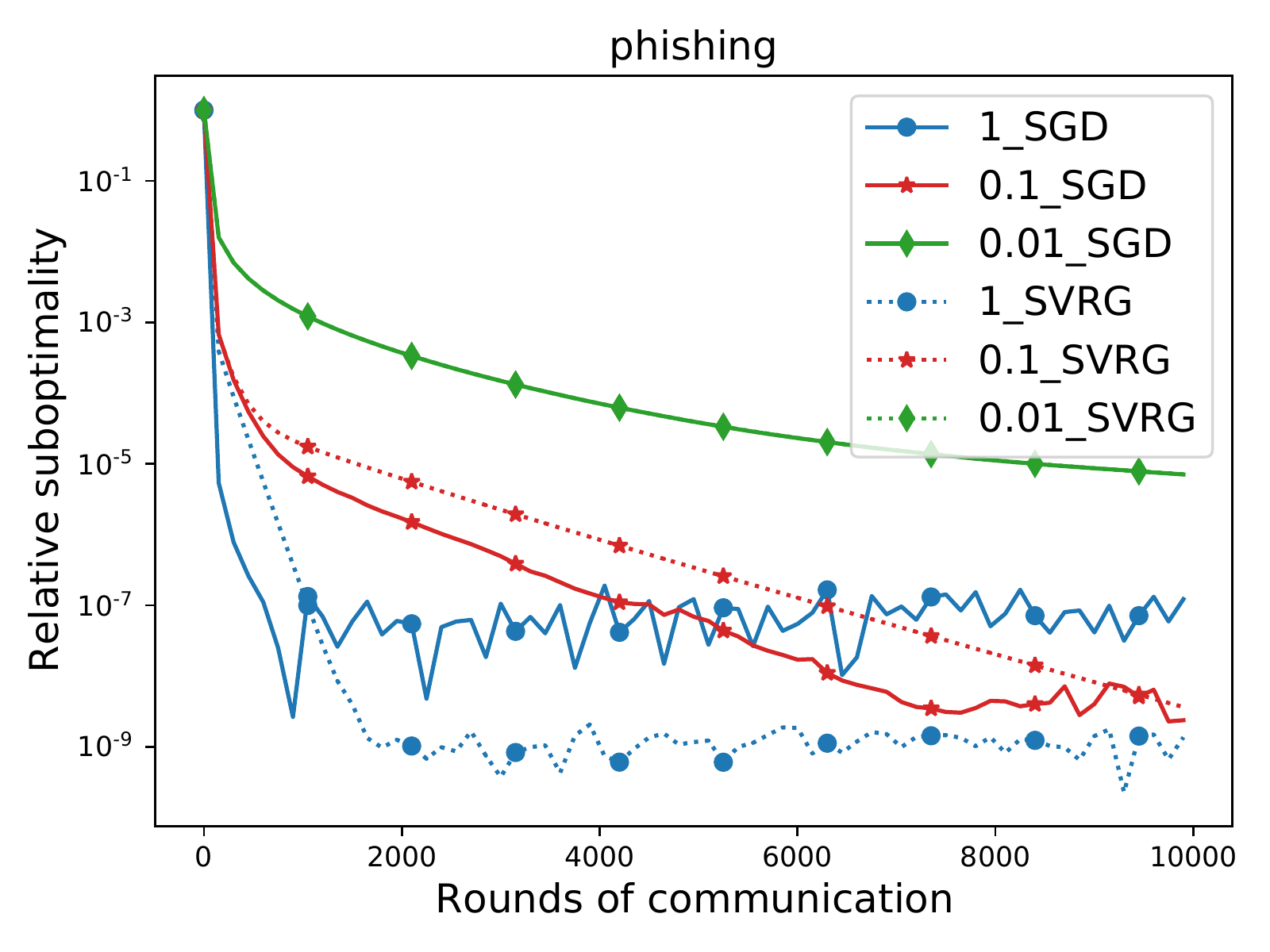}
\end{minipage}
\caption{Comparison of standard {\tt Local-SGD} (Algorithm~\ref{alg:local_sgd}), and {\tt Local-SVRG} (Algorithm~\ref{alg:local_svrg})  with various stepsizes $\gamma$. Logistic regression applied on LibSVM data~\citep{chang2011libsvm} with heterogenously splitted data.  Other parameters: $L=1, \mu=10^{-4}, \tau = 40$. Parameter $n$ chosen as per Table~\ref{tbl:ns}. (Same as Figure~\ref{fig:sgd_svrg_hom}, but with the heterogenous data split)}
\label{fig:sgd_svrg_het}
\end{figure}

\paragraph{Data split.} The experiment from Figure~\ref{fig:sgd_svrg_hom} in the main body of the paper splits the data among the clients  uniformly at random (i.e., split according to the the order given by a random permutation). However, in a typical FL scenario, the local data might significantly differ from the population average. For this reason, we also test on a different split of the data: we first sort the data according to the labels, and then split them among the clients. Figure~\ref{fig:sgd_svrg_het} shows the results. We draw a conclusions identical to Figure~\ref{fig:sgd_svrg_hom}. We see that {\tt Local-SVRG} was at least as good as {\tt Local-SGD} for every stepsize choice and every dataset. Further, the prediction that the smaller stepsize yields the smaller of the optimum neighborhood for the price of slower convergence was confirmed.

\paragraph{Environment.} All experiments were performed in a simulated environment on a single machine.

\subsection{The Effect of  Local Shift/Drifts}

The experiment presented in Section~\ref{sec:exp} examined the effect of the noise on the performance of local methods and demonstrated that control variates can be efficiently employed to reduce that noise. In this section, we study the second factor that influences the neighborhood to which {\tt Local-SGD} converges: non-stationarity of {\tt Local-GD}.

We have already shown that the mentioned non-stationarity of {\tt Local-GD} can be fixed using a carefully designed idealized/optimal shift that depends on the solution $x^*$ (see Algorithm~\ref{alg:local_sgd_star}). Furthermore, we have shown that this idealized shift can be learned on-the-fly at the small price of slightly slower convergence rate (see Algorithm~\ref{alg:l_local_svrg} -- {\tt SS-Local-SGD}/{\tt SCAFFOLD}).\footnote{In fact, {\tt SCAFFOLD} can be coupled together with {\tt Local-SVRG} given that the local objectives are of a finite-sum structure, resulting in Algorithm~\ref{alg:l_local_svrg_fs}.}

\begin{figure}[!h]
\centering
\begin{minipage}{0.3\textwidth}
  \centering
\includegraphics[width =  \textwidth ]{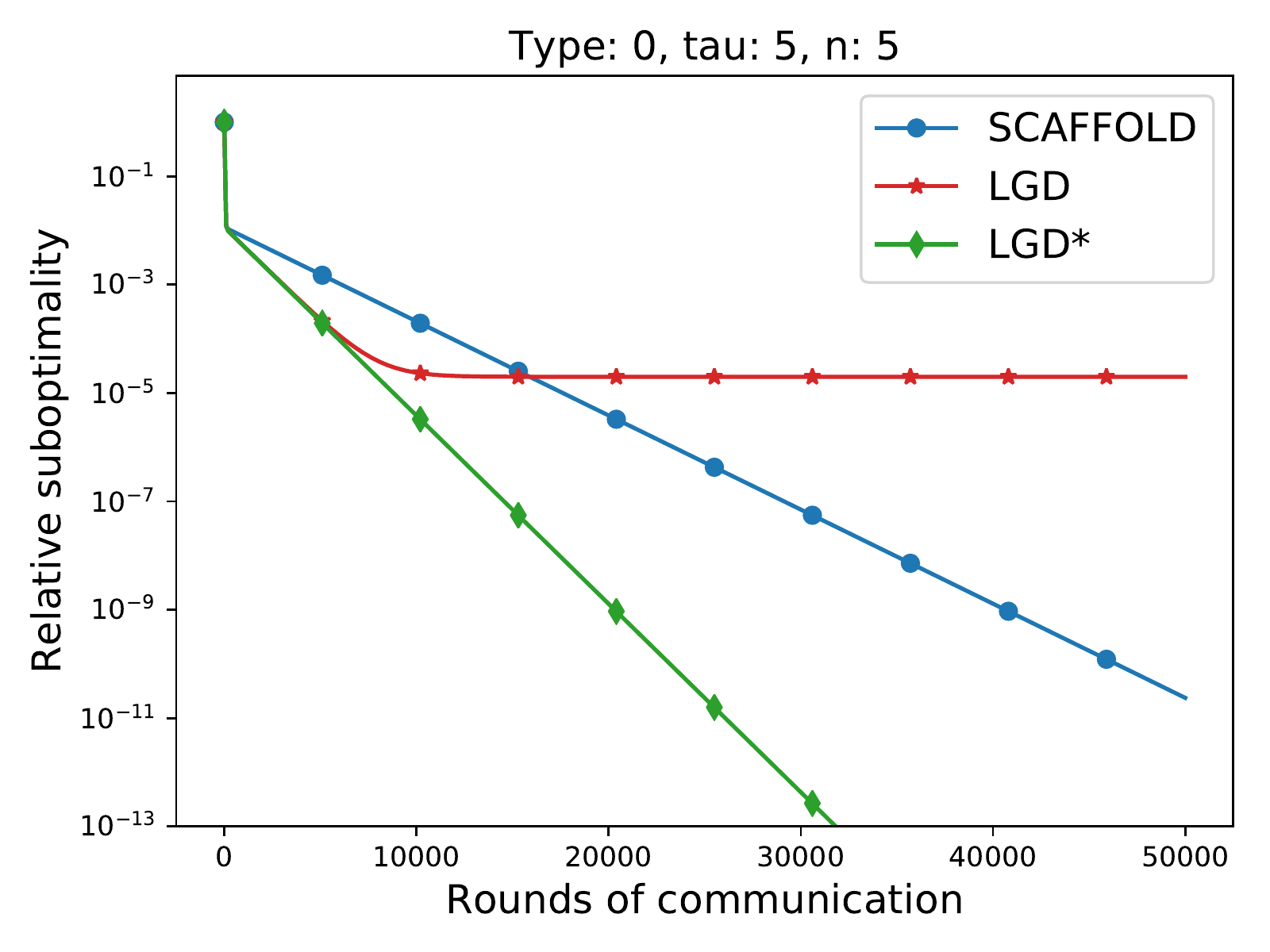}
\end{minipage}
\begin{minipage}{0.3\textwidth}
  \centering
\includegraphics[width =  \textwidth ]{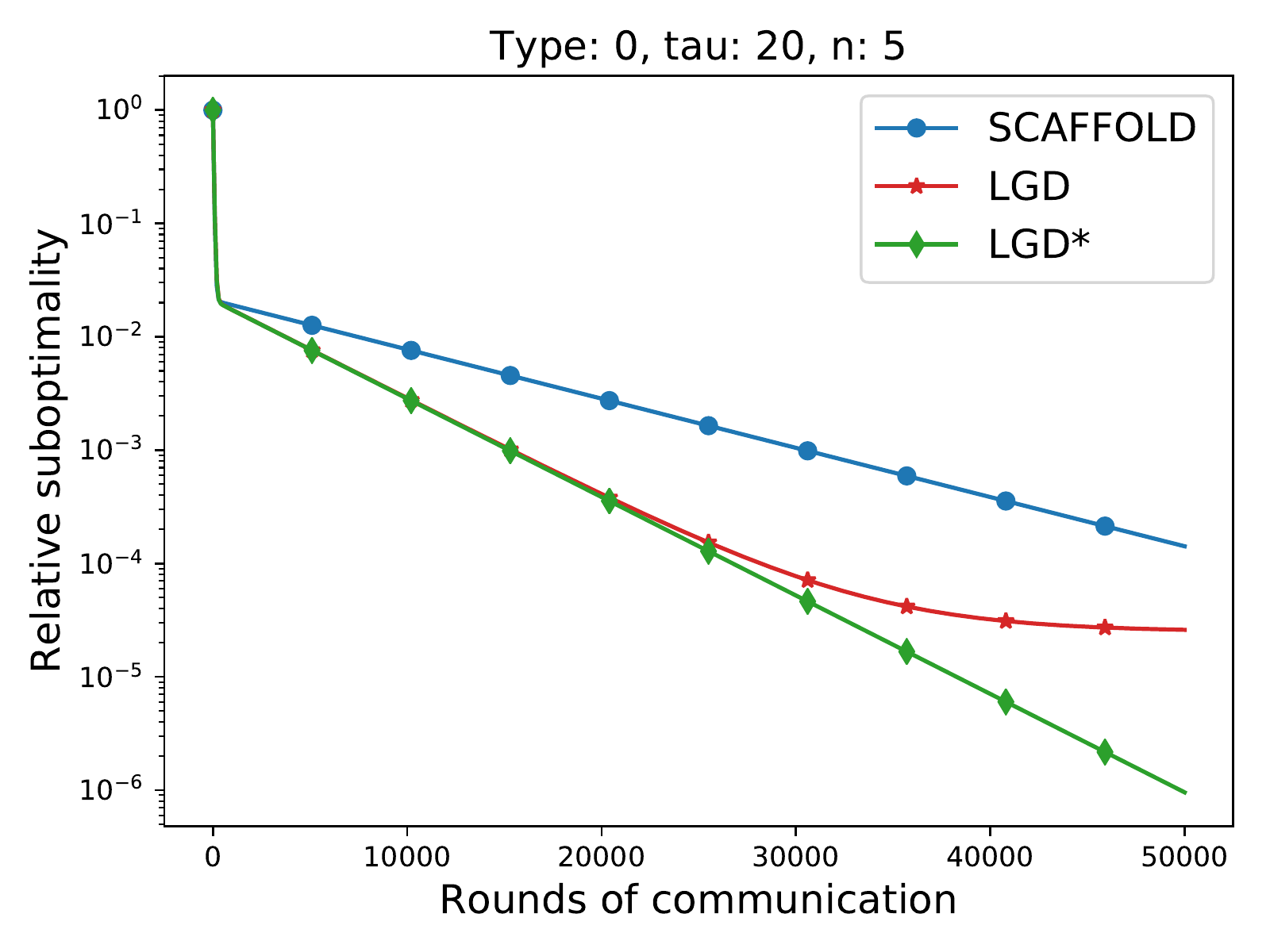}
\end{minipage}
\begin{minipage}{0.3\textwidth}
  \centering
\includegraphics[width =  \textwidth ]{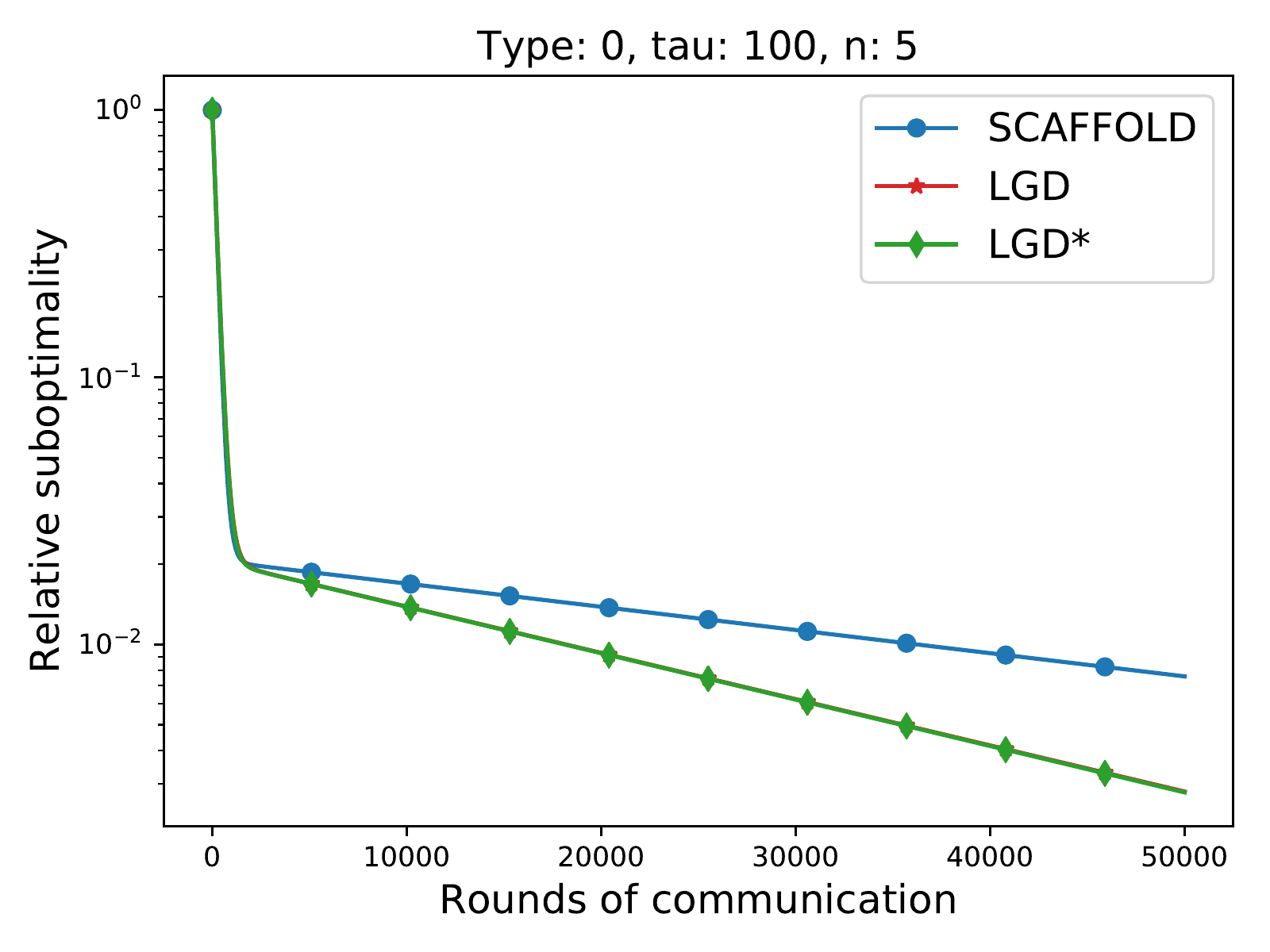}
\end{minipage}
\\
\begin{minipage}{0.3\textwidth}
  \centering
\includegraphics[width =  \textwidth ]{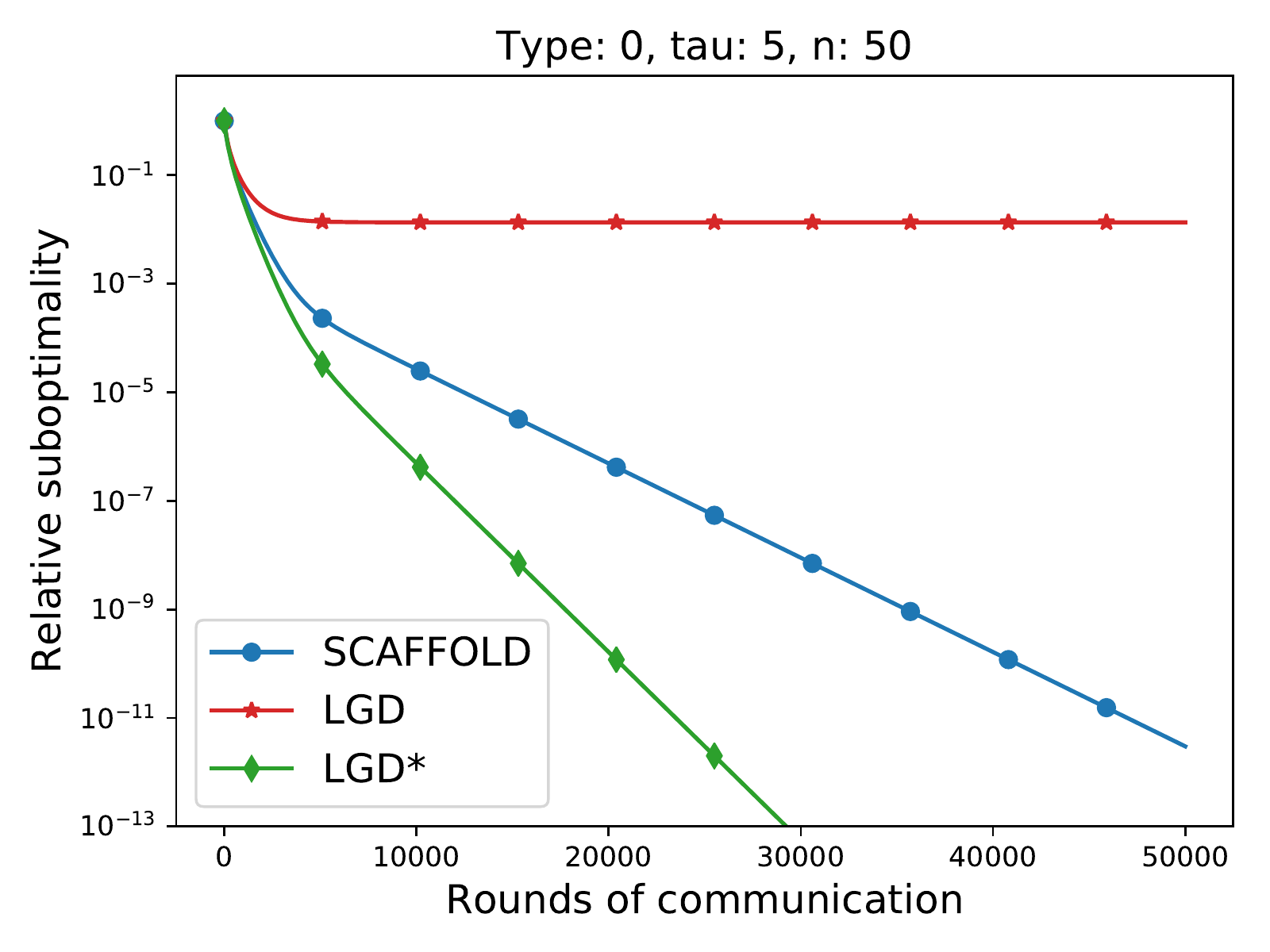}
\end{minipage}
\begin{minipage}{0.3\textwidth}
  \centering
\includegraphics[width =  \textwidth ]{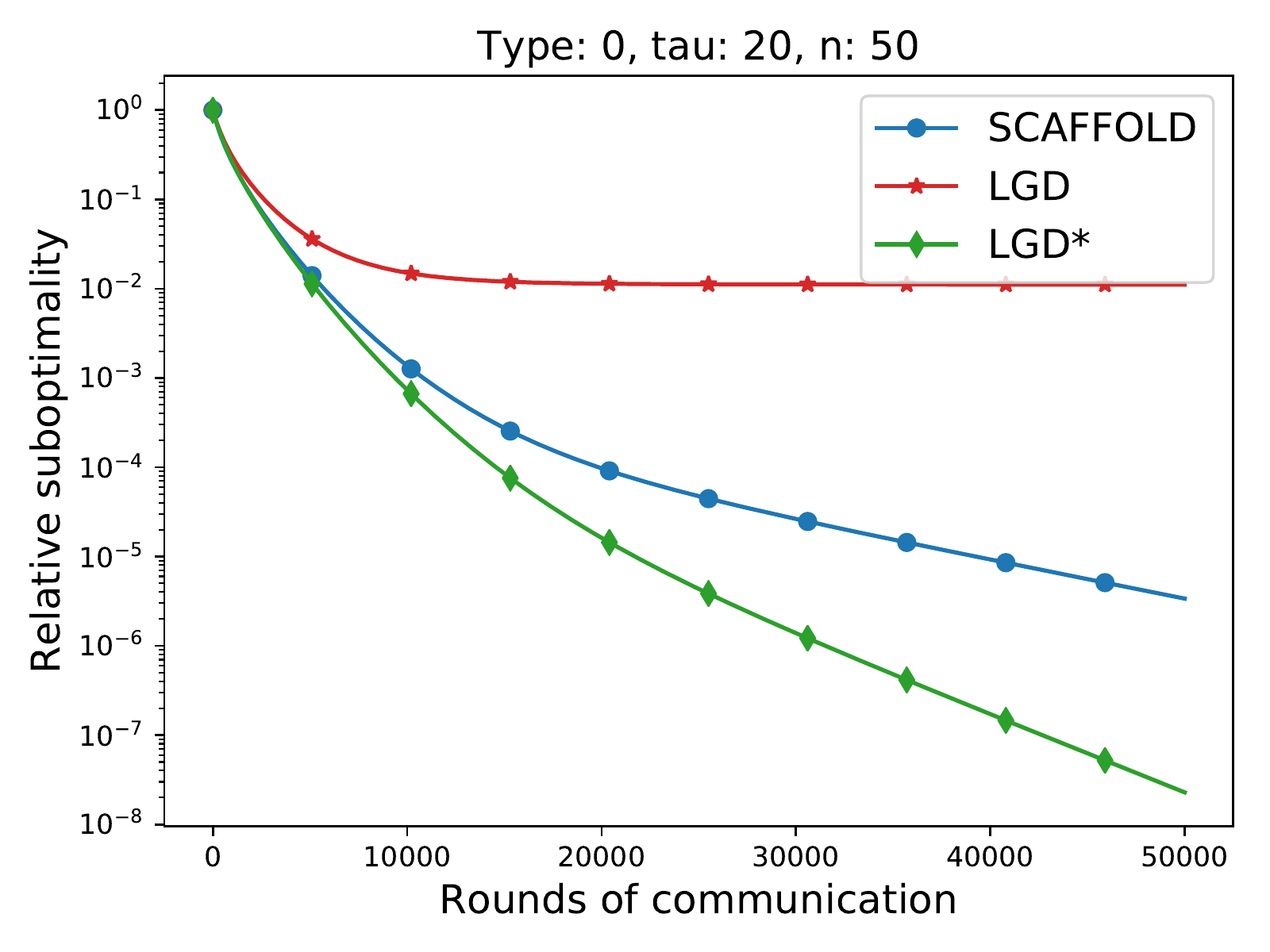}
\end{minipage}
\begin{minipage}{0.3\textwidth}
  \centering
\includegraphics[width =  \textwidth ]{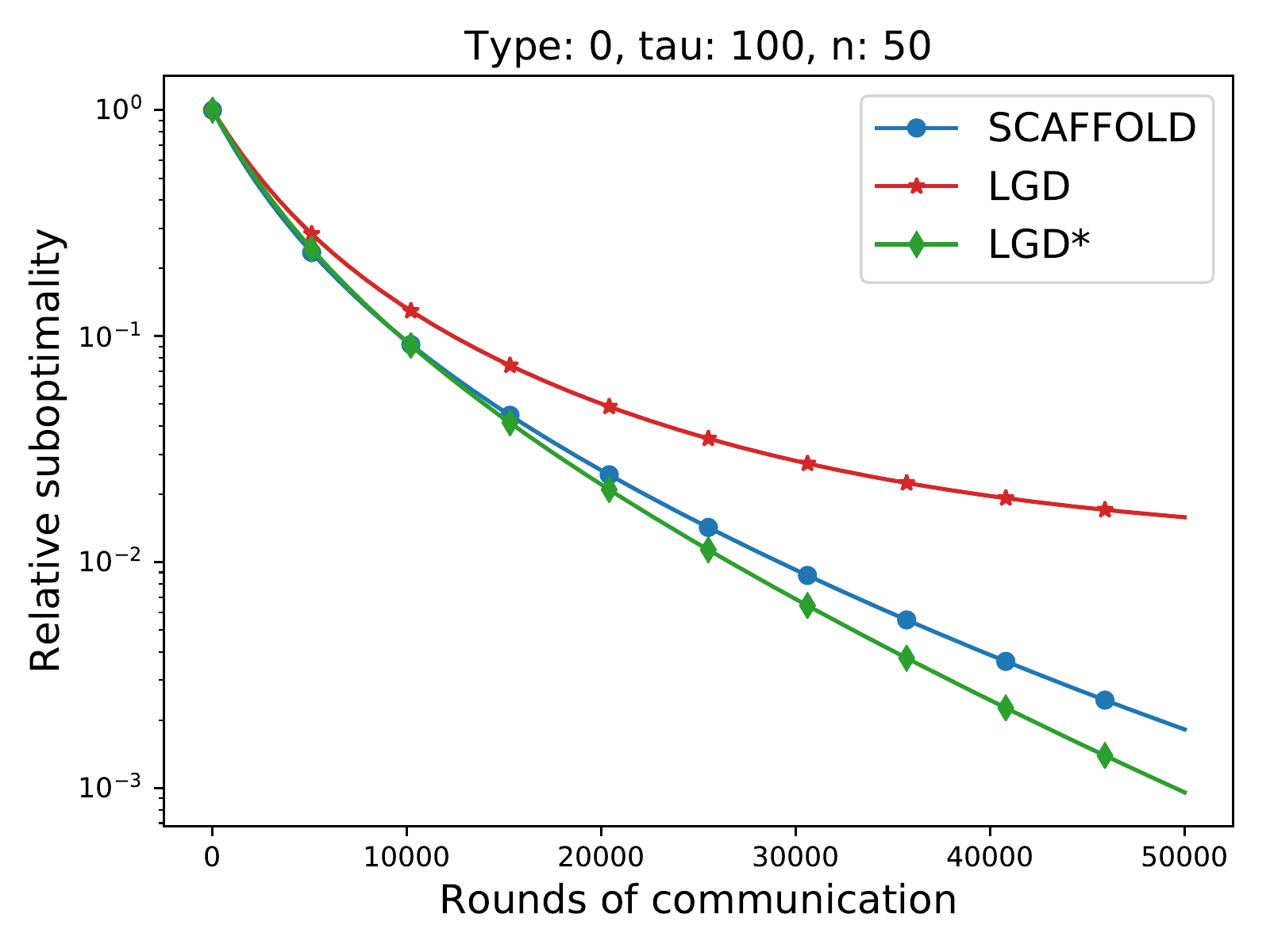}
\end{minipage}
\\
\begin{minipage}{0.3\textwidth}
  \centering
\includegraphics[width =  \textwidth ]{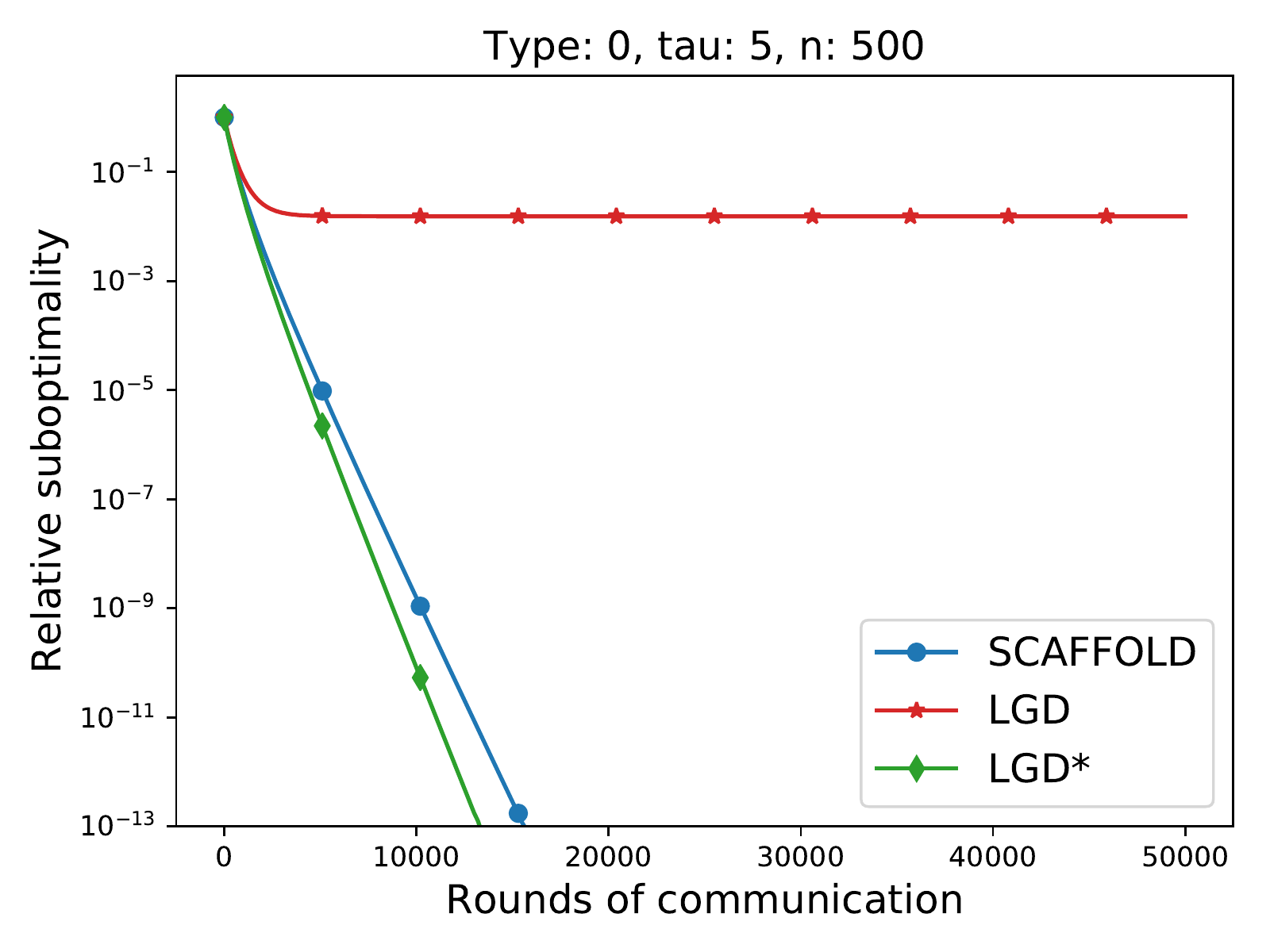}
\end{minipage}
\begin{minipage}{0.3\textwidth}
  \centering
\includegraphics[width =  \textwidth ]{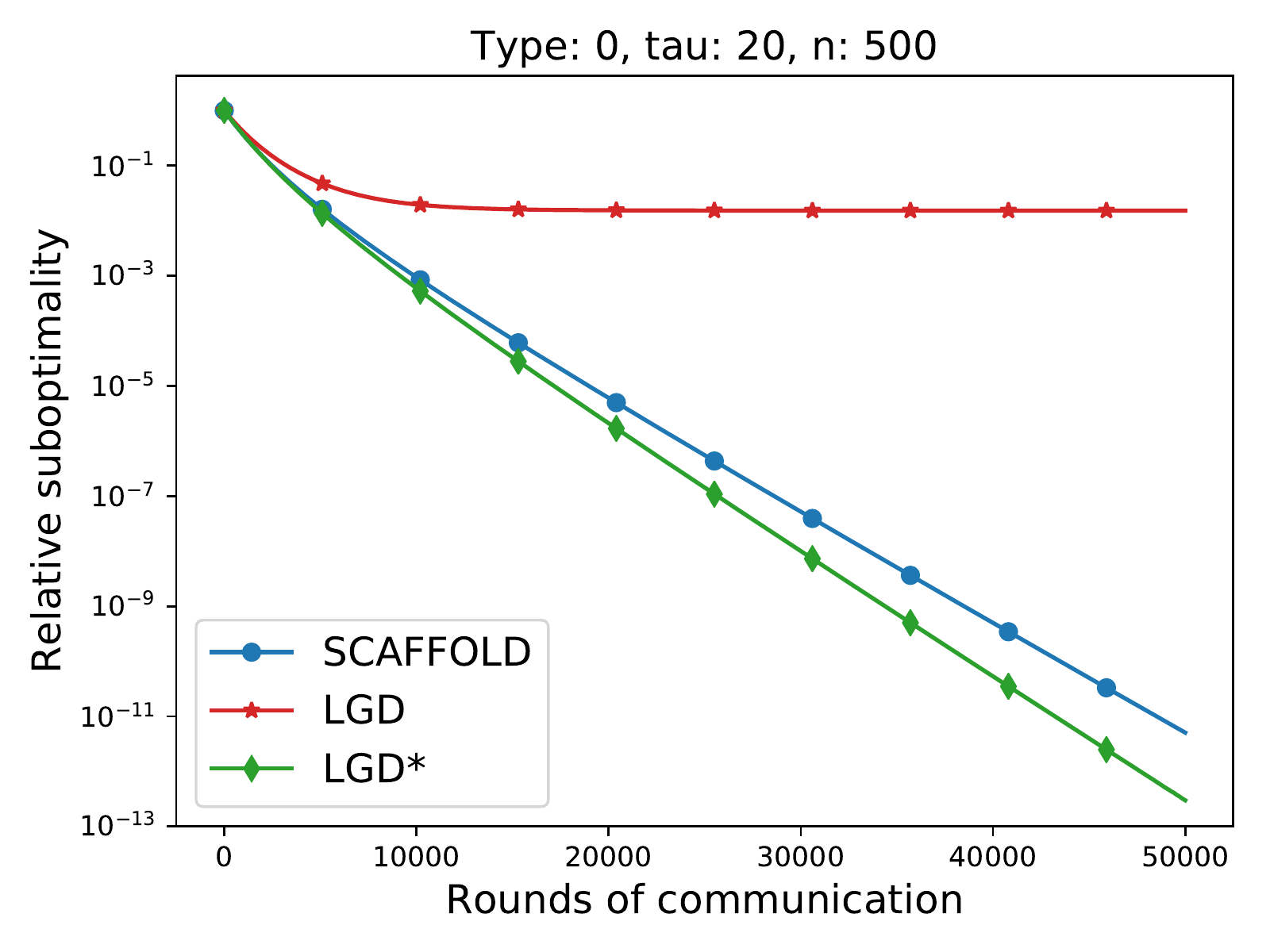}
\end{minipage}
\begin{minipage}{0.3\textwidth}
  \centering
\includegraphics[width =  \textwidth ]{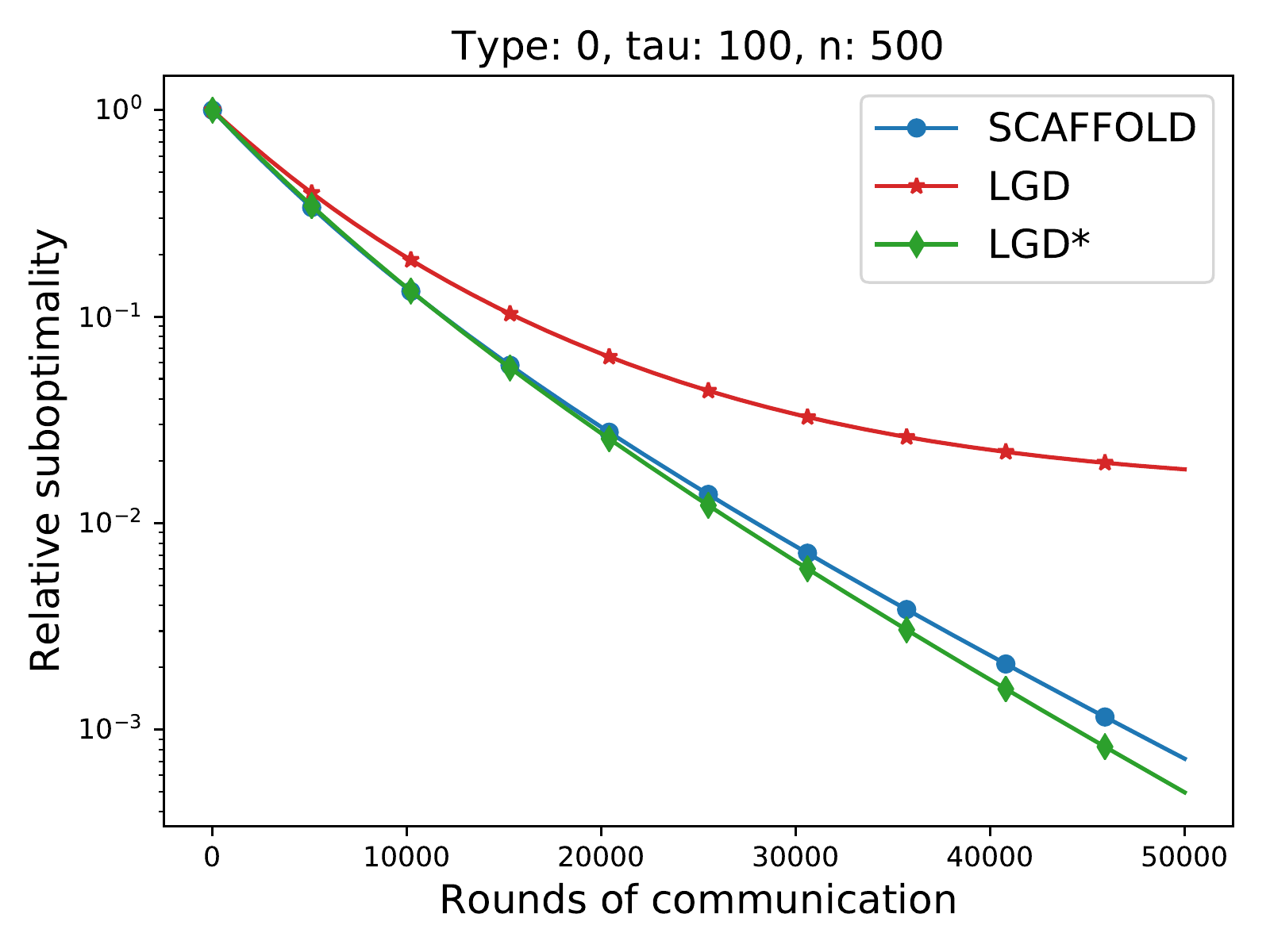}
\end{minipage}
\caption{Comparison of the following noiseless algorithms  {\tt Local-SGD} ({\tt LGD}, Algorithm~\ref{alg:local_sgd} with no local noise) and {\tt SCAFFOLD}~\citep{karimireddy2020scaffold} (Algorithm~\ref{alg:l_local_svrg} without ``Loopless'') and {\tt S*-Local-SGD} ({\tt LGD*}, Algorithm~\ref{alg:local_sgd_star}). Quadratic minimization, problem type 0 (see Table~\ref{tbl:instances}). }
\label{fig:artif1}
\end{figure}

\begin{figure}[!h]
\centering
\begin{minipage}{0.3\textwidth}
  \centering
\includegraphics[width =  \textwidth ]{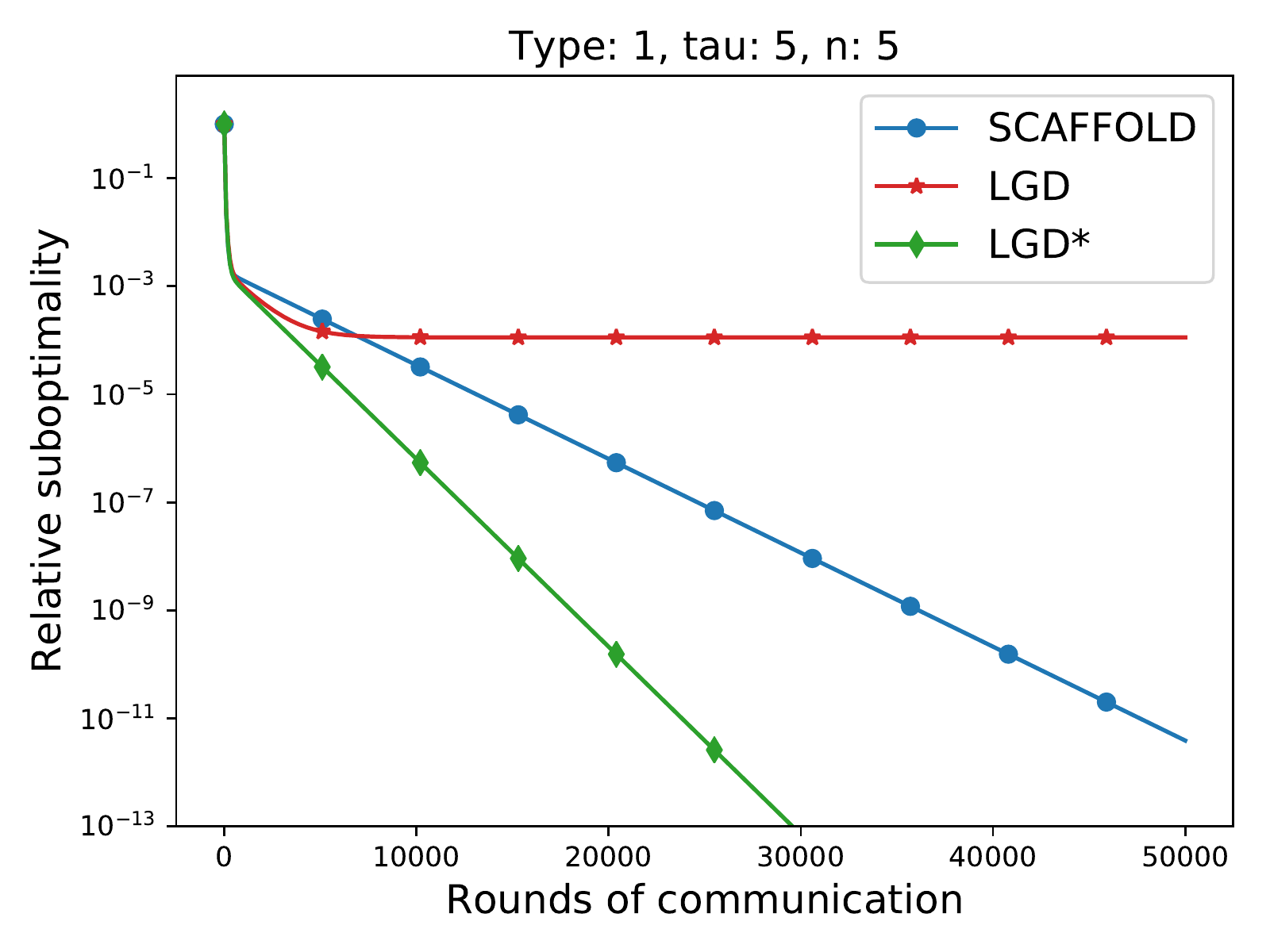}
\end{minipage}
\begin{minipage}{0.3\textwidth}
  \centering
\includegraphics[width =  \textwidth ]{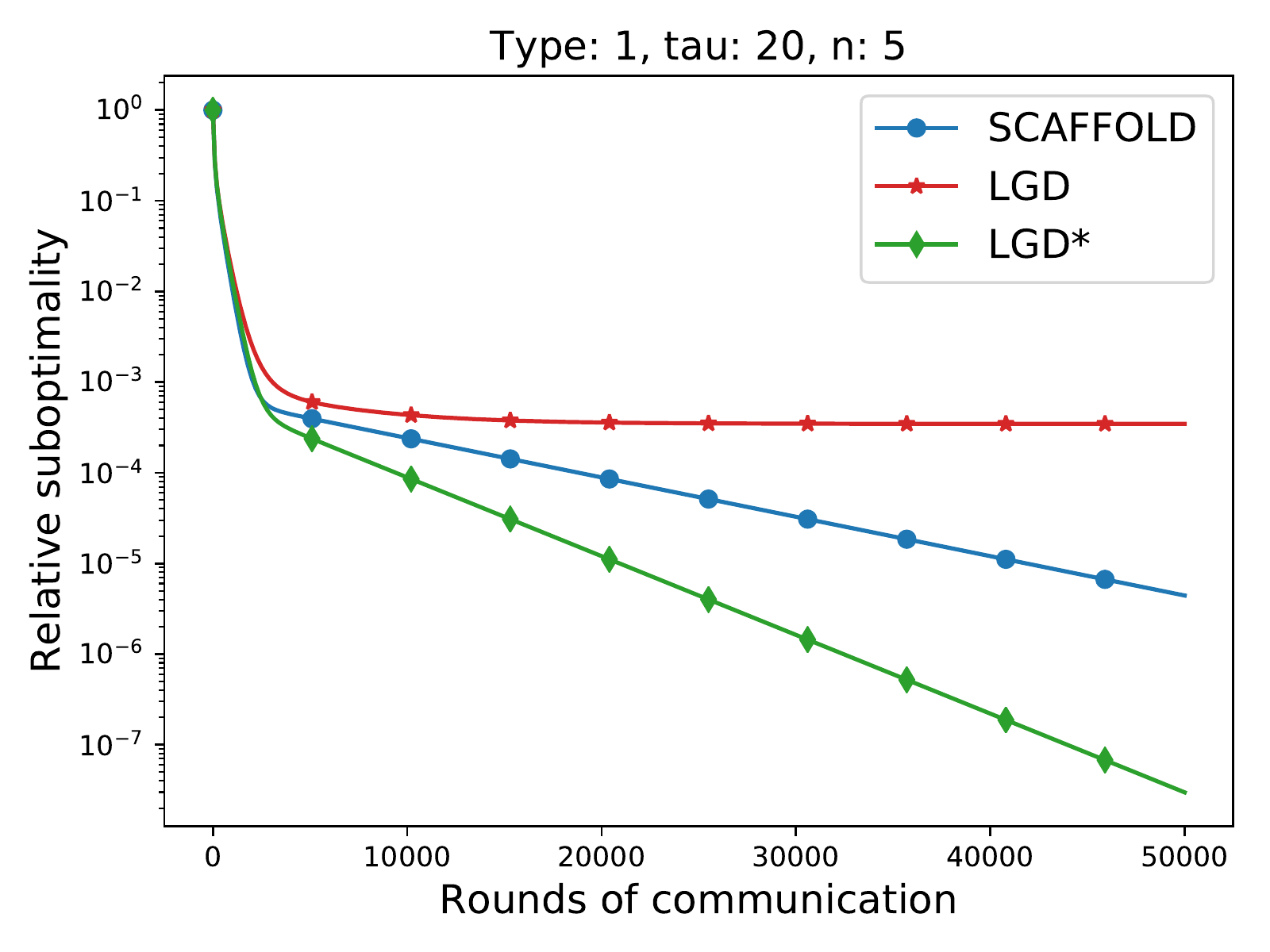}
\end{minipage}
\begin{minipage}{0.3\textwidth}
  \centering
\includegraphics[width =  \textwidth ]{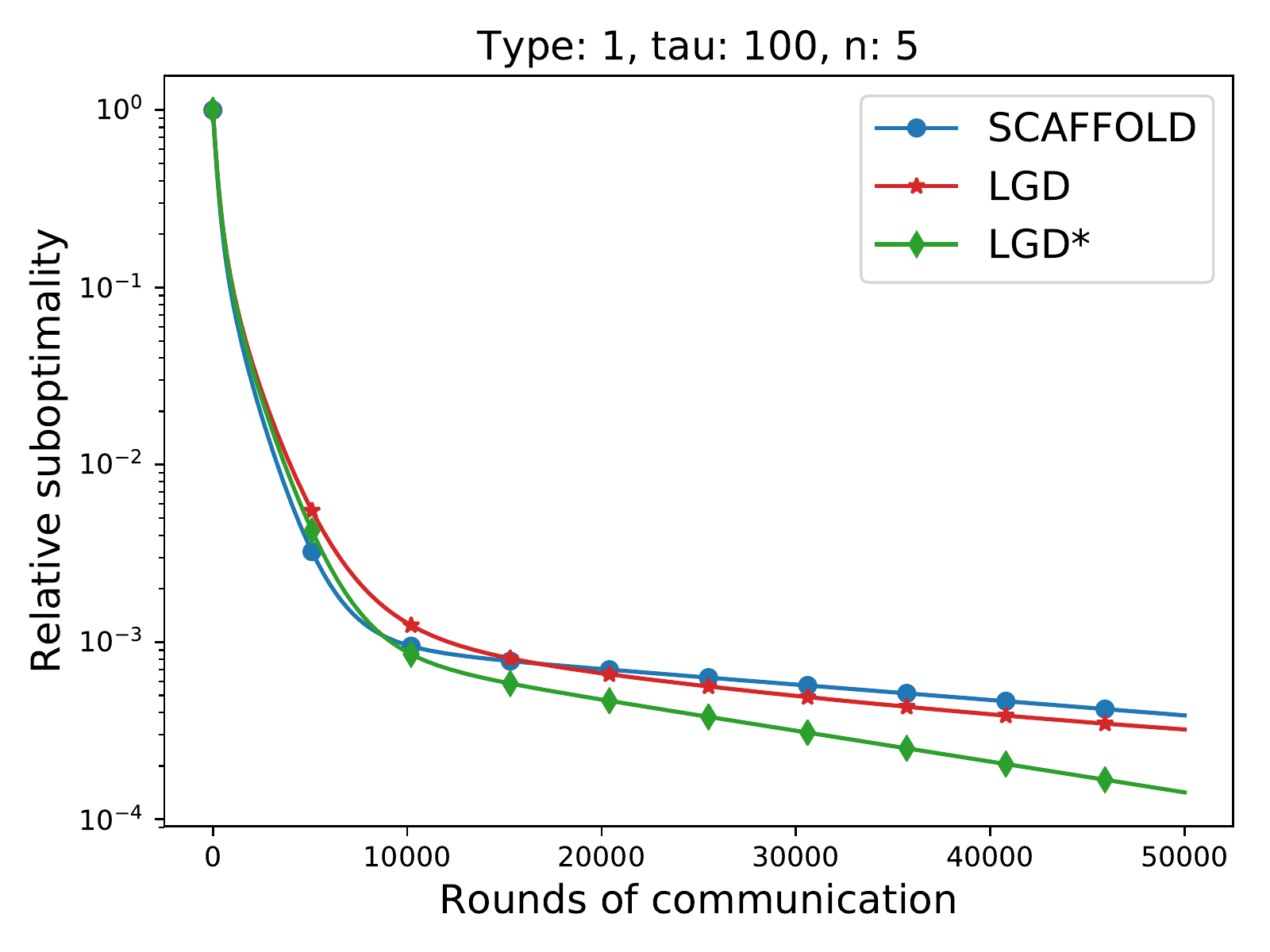}
\end{minipage}
\\
\begin{minipage}{0.3\textwidth}
  \centering
\includegraphics[width =  \textwidth ]{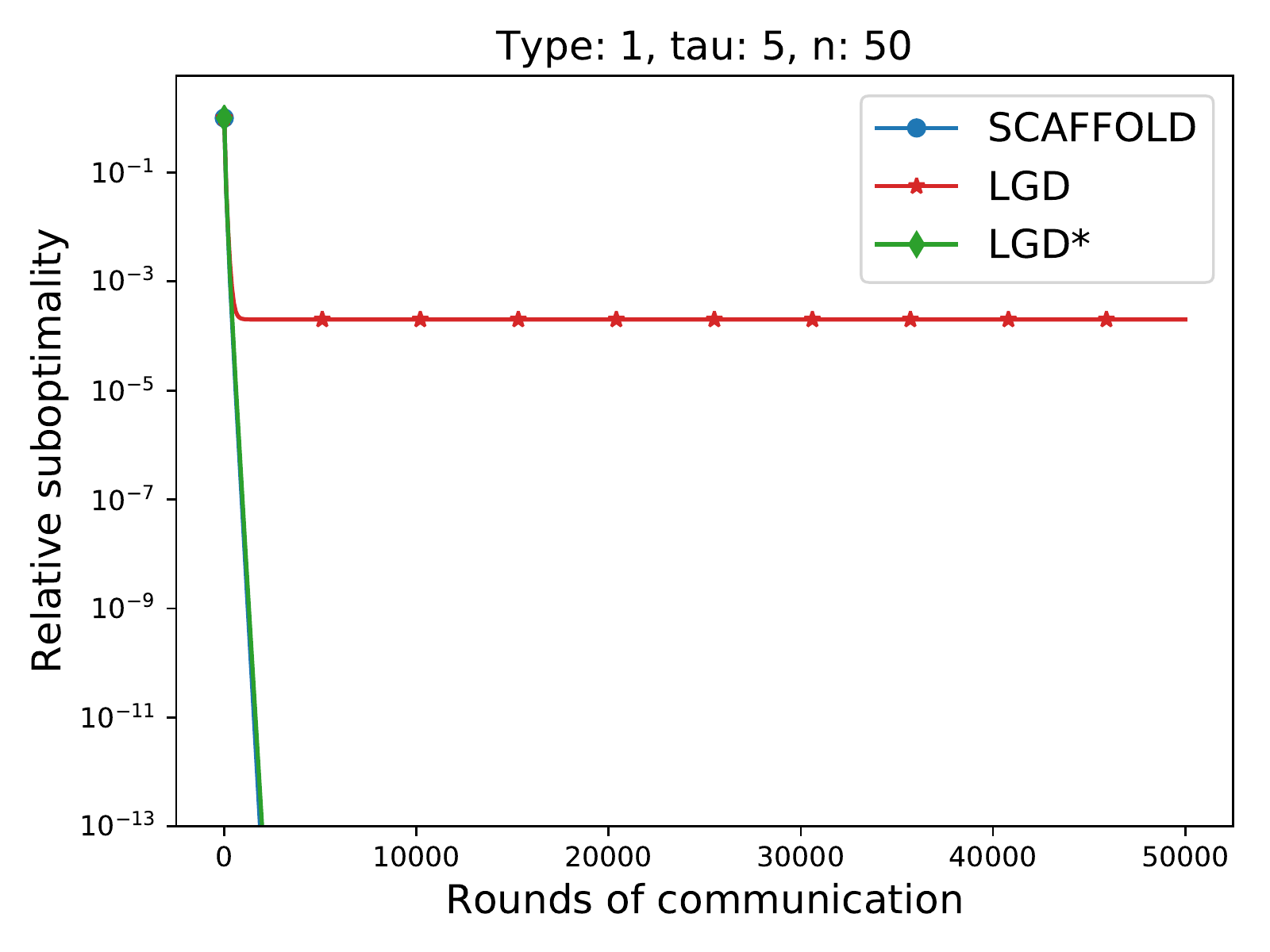}
\end{minipage}
\begin{minipage}{0.3\textwidth}
  \centering
\includegraphics[width =  \textwidth ]{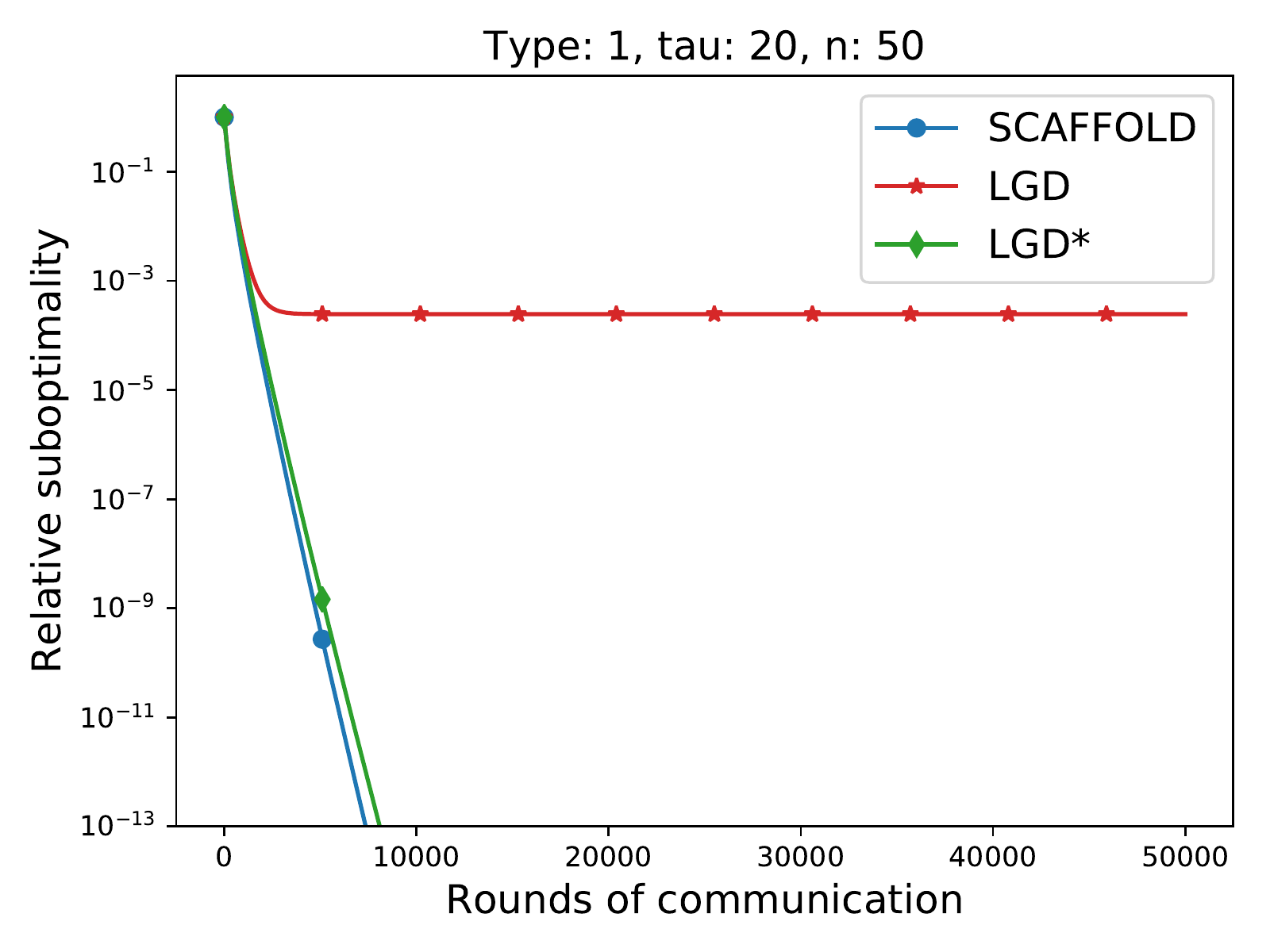}
\end{minipage}
\begin{minipage}{0.3\textwidth}
  \centering
\includegraphics[width =  \textwidth ]{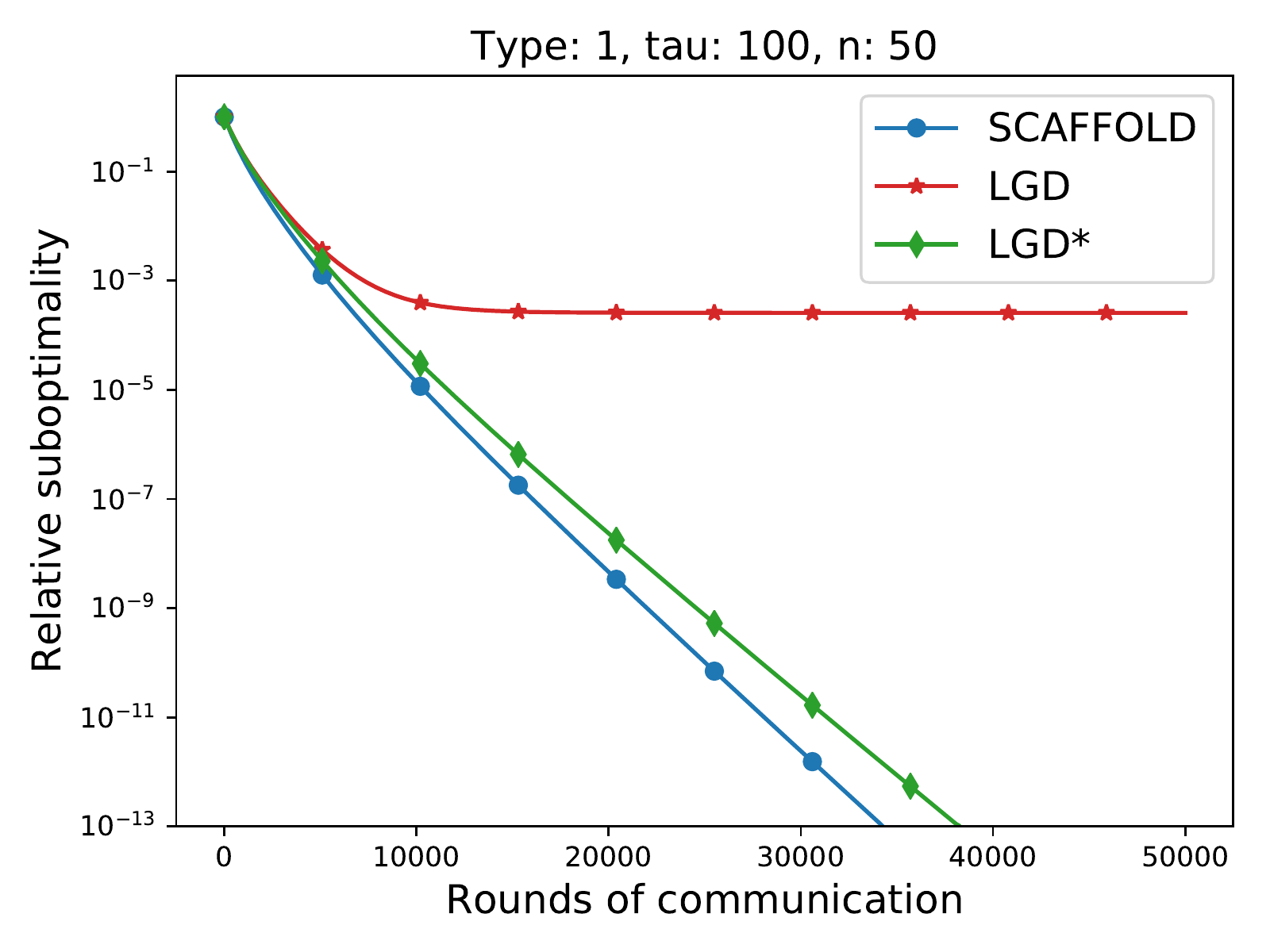}
\end{minipage}
\\
\begin{minipage}{0.3\textwidth}
  \centering
\includegraphics[width =  \textwidth ]{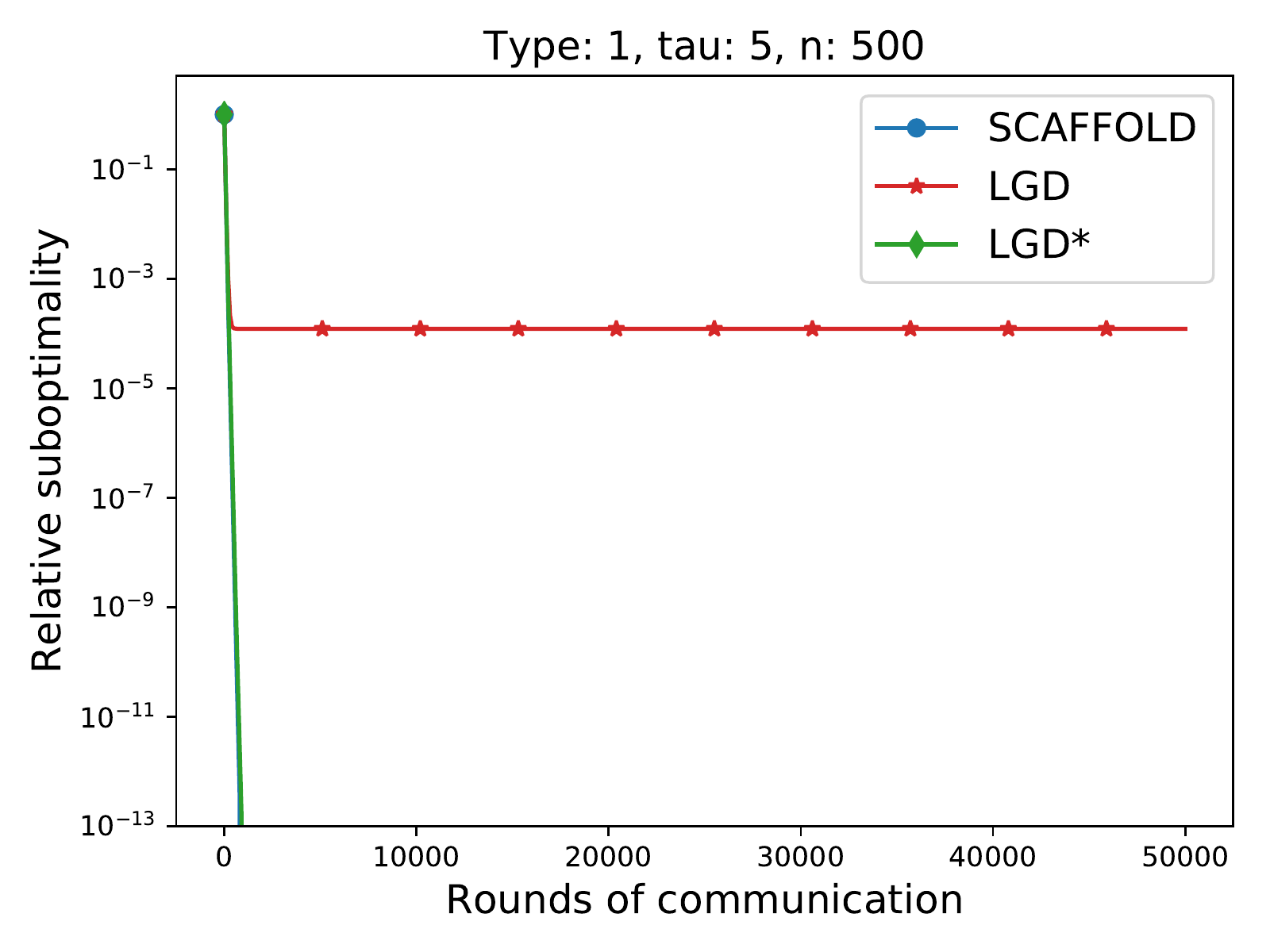}
\end{minipage}
\begin{minipage}{0.3\textwidth}
  \centering
\includegraphics[width =  \textwidth ]{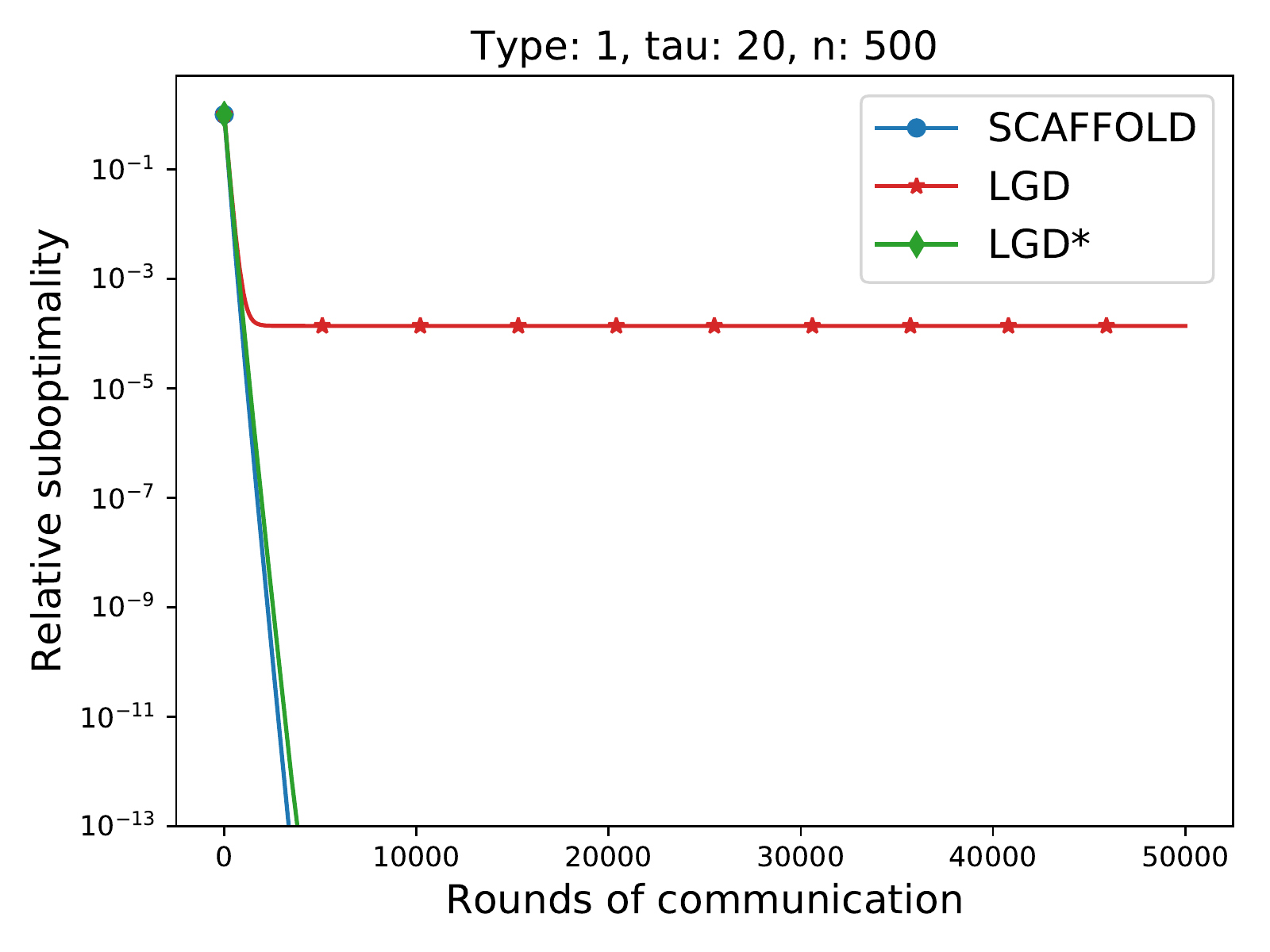}
\end{minipage}
\begin{minipage}{0.3\textwidth}
  \centering
\includegraphics[width =  \textwidth ]{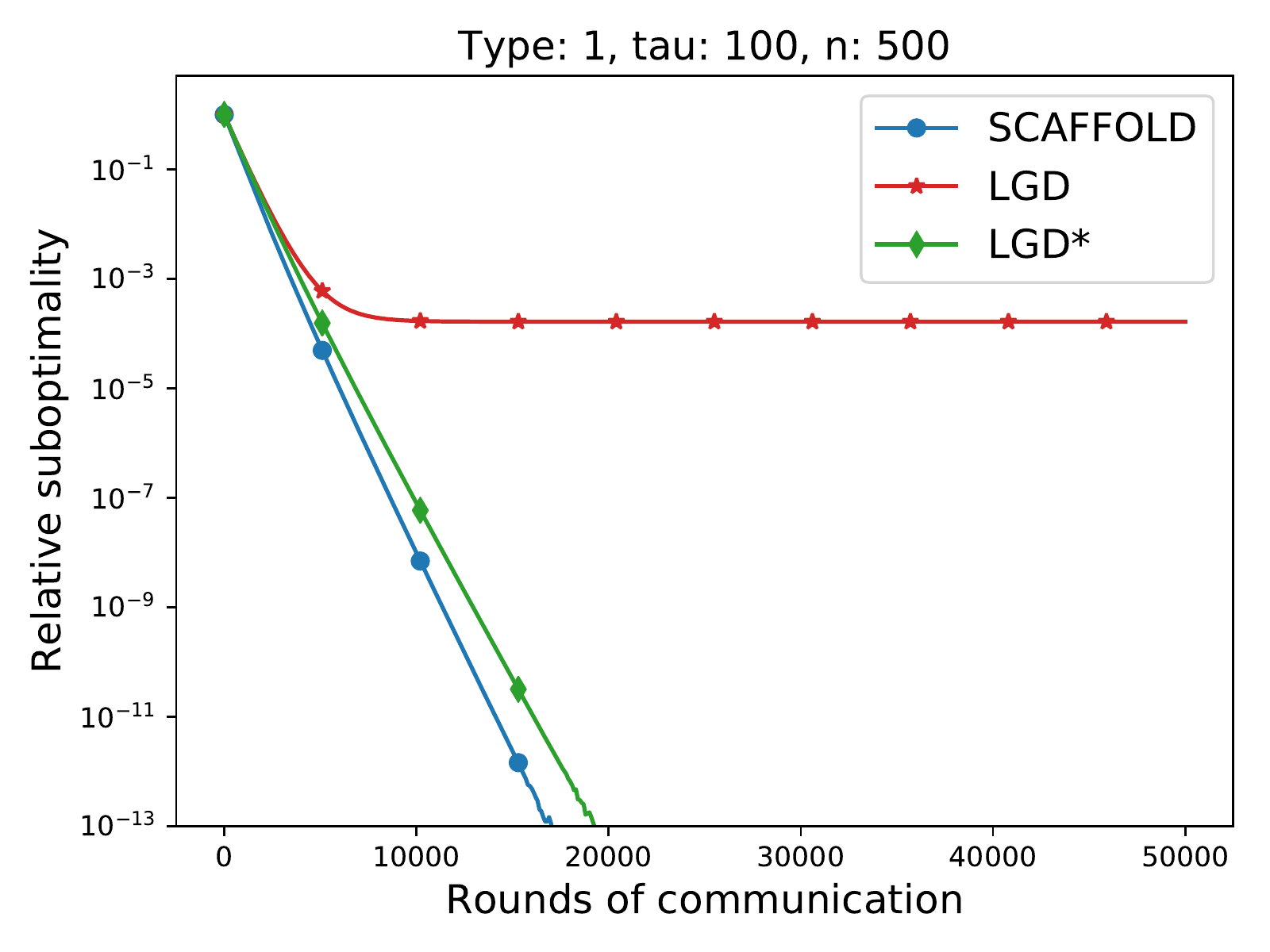}
\end{minipage}
\caption{Comparison of the following noiseless algorithms  {\tt Local-SGD} ({\tt LGD}, Algorithm~\ref{alg:local_sgd} with no local noise) and {\tt SCAFFOLD}~\citep{karimireddy2020scaffold} (Algorithm~\ref{alg:l_local_svrg} without ``Loopless'') and {\tt S*-Local-SGD} ({\tt LGD*}, Algorithm~\ref{alg:local_sgd_star}). Quadratic minimization, problem type 1 (see Table~\ref{tbl:instances}). }
\label{fig:artif2}
\end{figure}

\begin{figure}[!h]
\centering
\begin{minipage}{0.3\textwidth}
  \centering
\includegraphics[width =  \textwidth ]{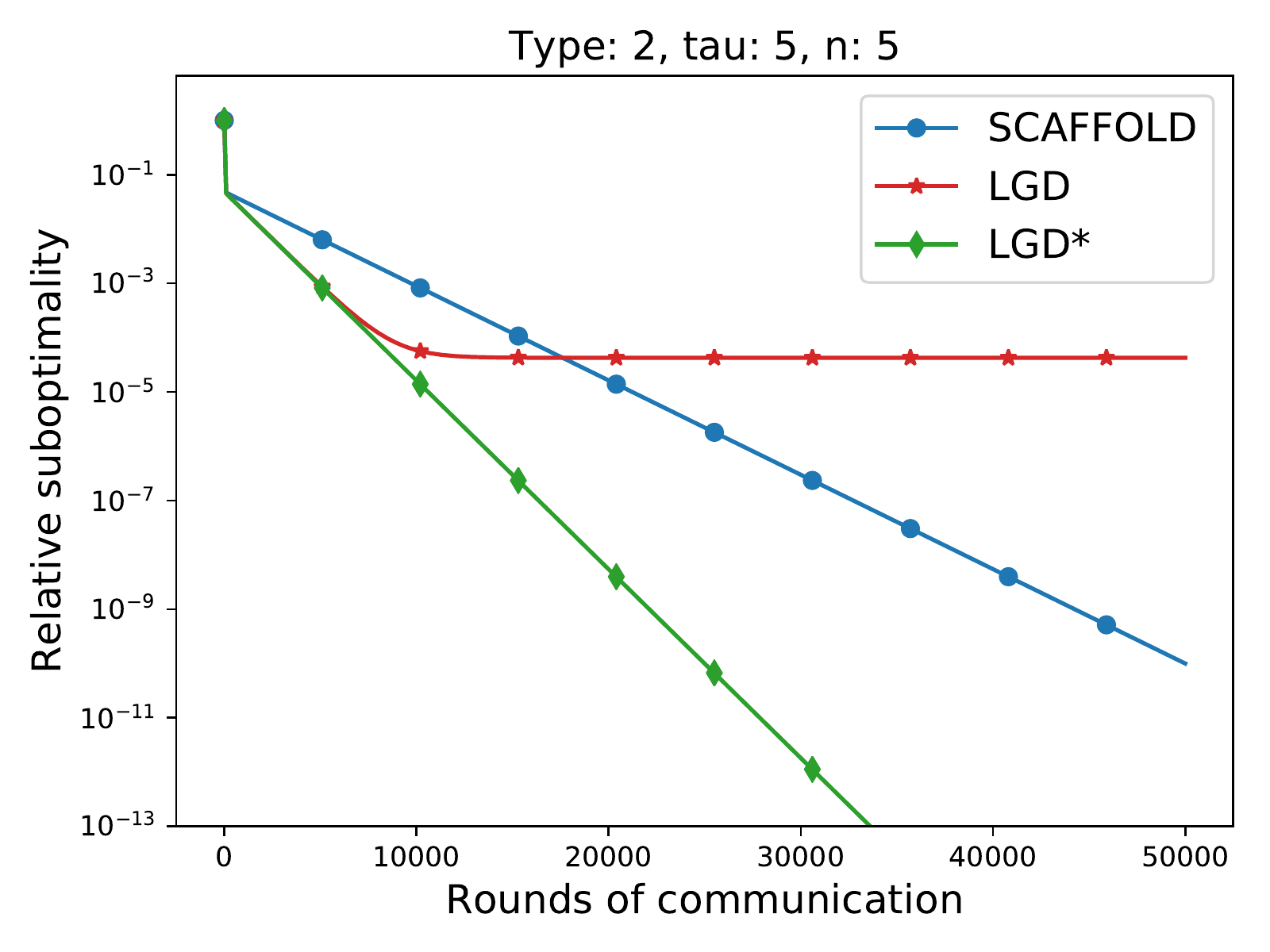}
\end{minipage}
\begin{minipage}{0.3\textwidth}
  \centering
\includegraphics[width =  \textwidth ]{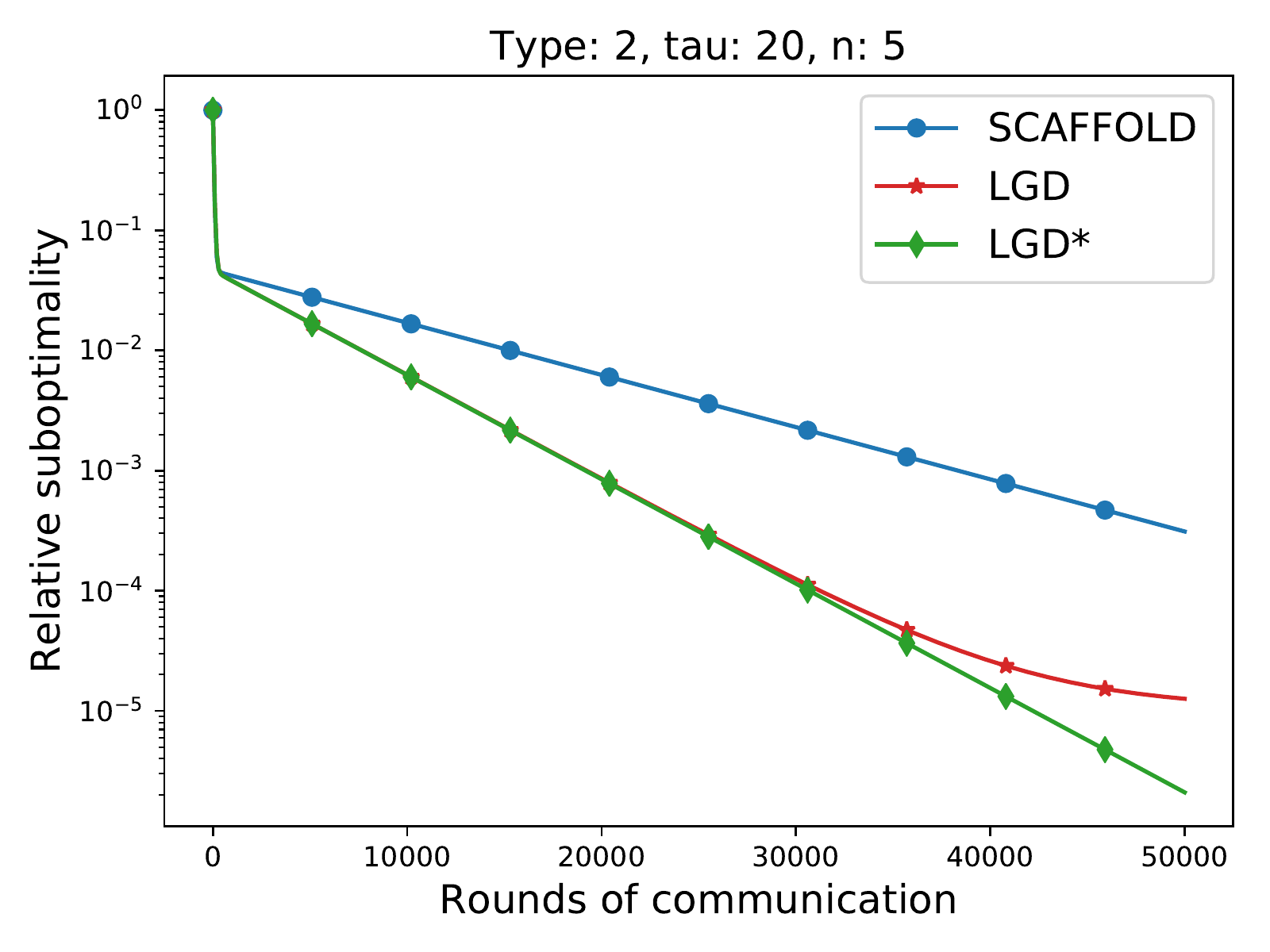}
\end{minipage}
\begin{minipage}{0.3\textwidth}
  \centering
\includegraphics[width =  \textwidth ]{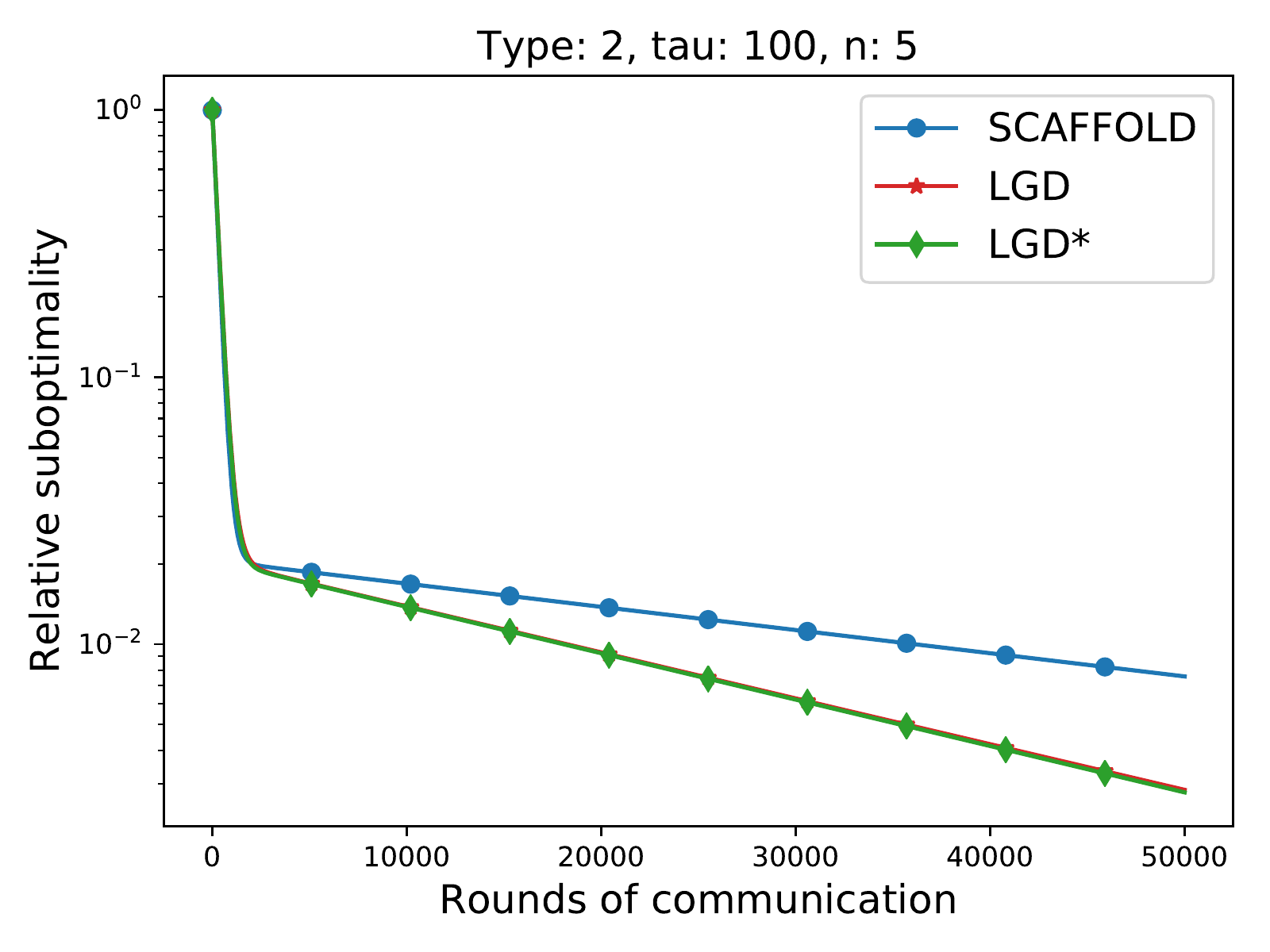}
\end{minipage}
\\
\begin{minipage}{0.3\textwidth}
  \centering
\includegraphics[width =  \textwidth ]{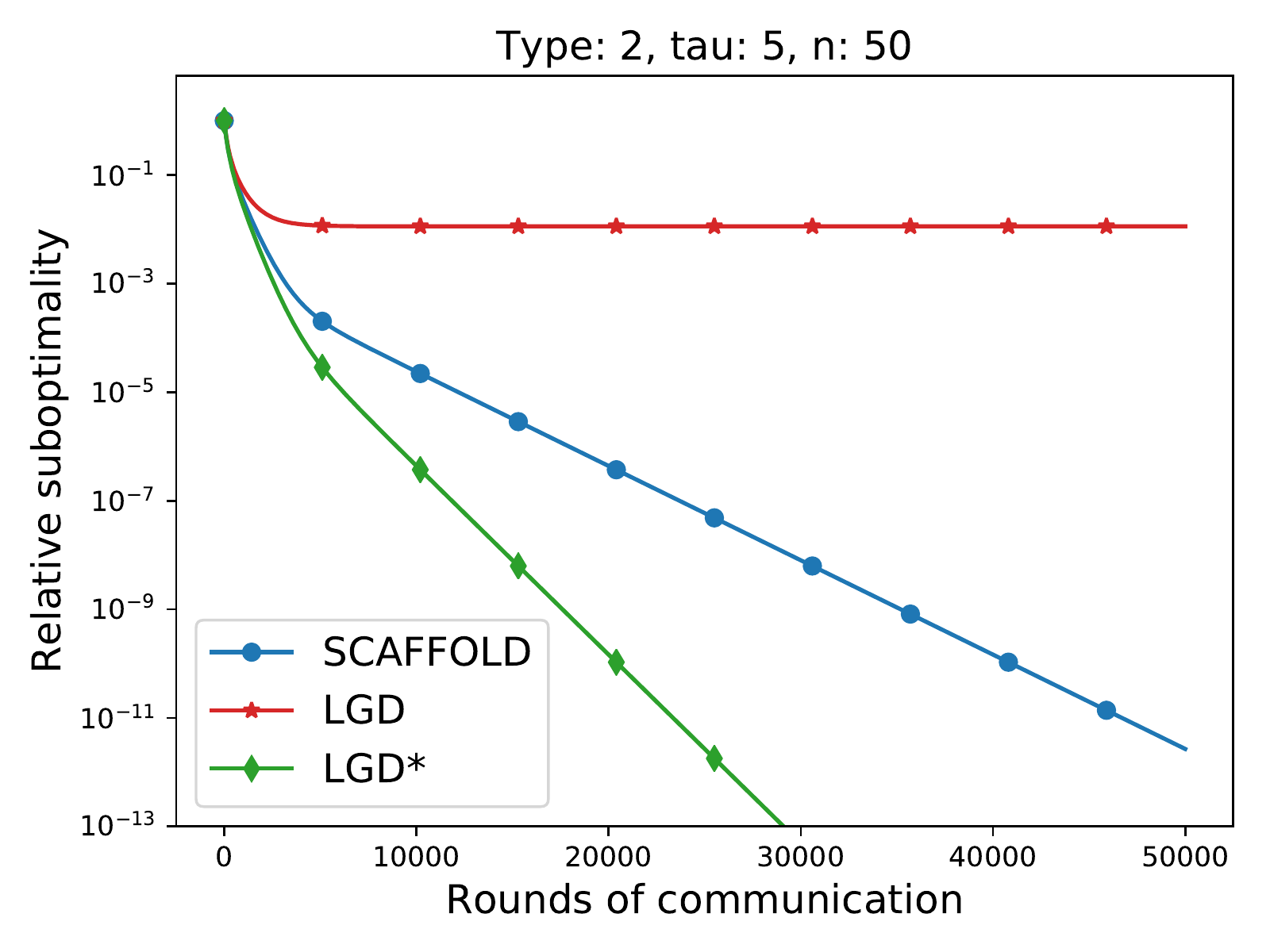}
\end{minipage}
\begin{minipage}{0.3\textwidth}
  \centering
\includegraphics[width =  \textwidth ]{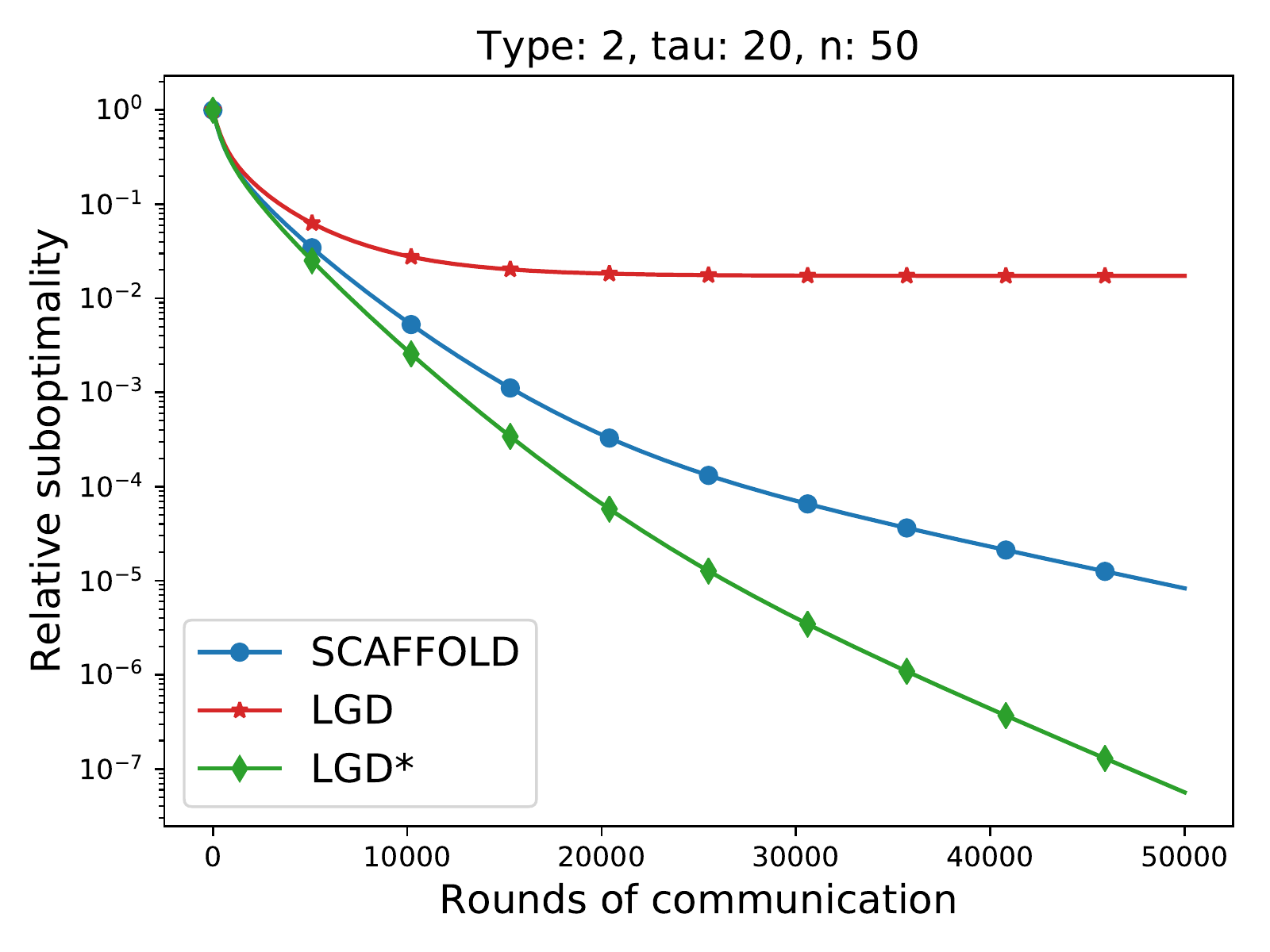}
\end{minipage}
\begin{minipage}{0.3\textwidth}
  \centering
\includegraphics[width =  \textwidth ]{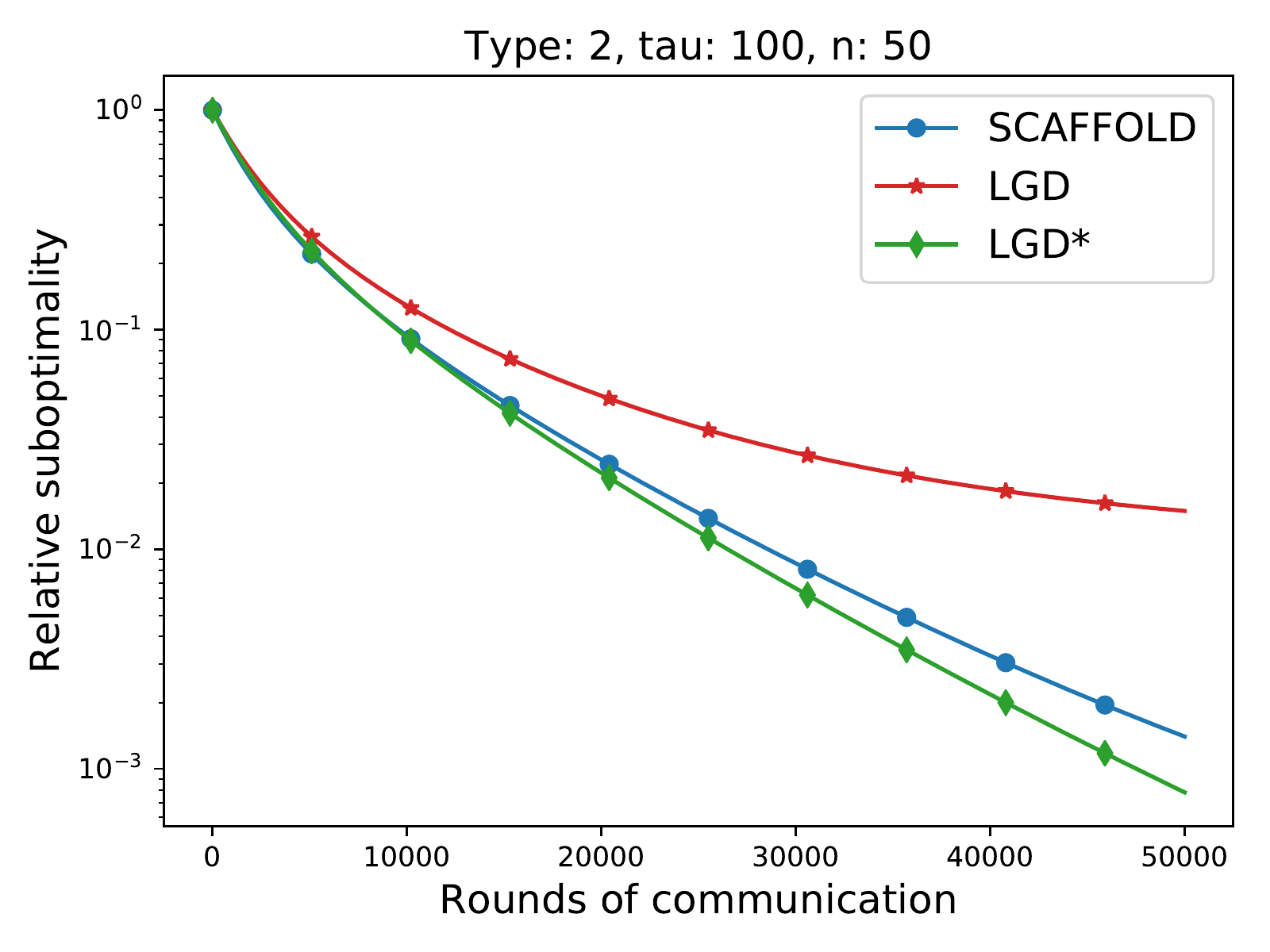}
\end{minipage}
\\
\begin{minipage}{0.3\textwidth}
  \centering
\includegraphics[width =  \textwidth ]{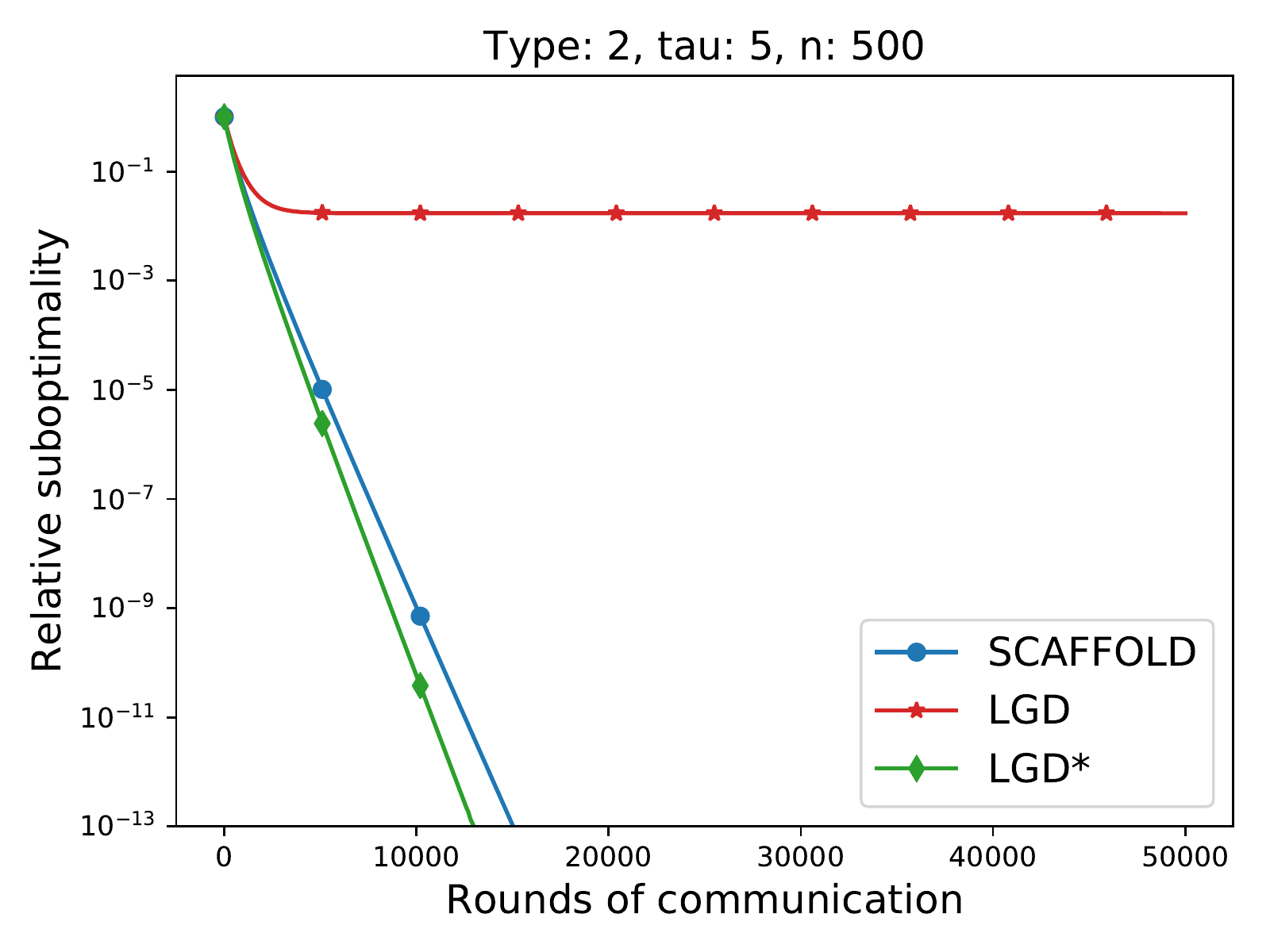}
\end{minipage}
\begin{minipage}{0.3\textwidth}
  \centering
\includegraphics[width =  \textwidth ]{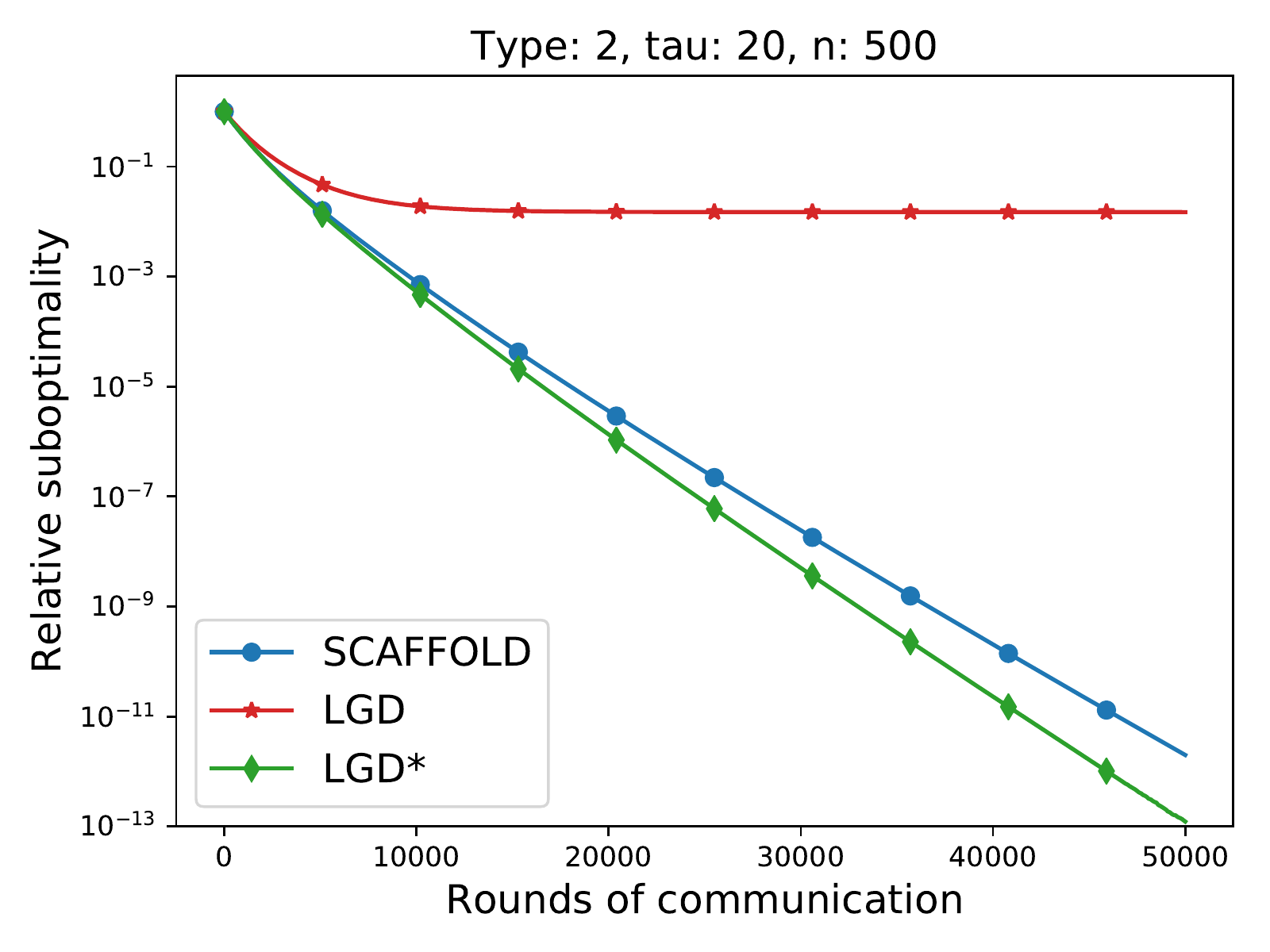}
\end{minipage}
\begin{minipage}{0.3\textwidth}
  \centering
\includegraphics[width =  \textwidth ]{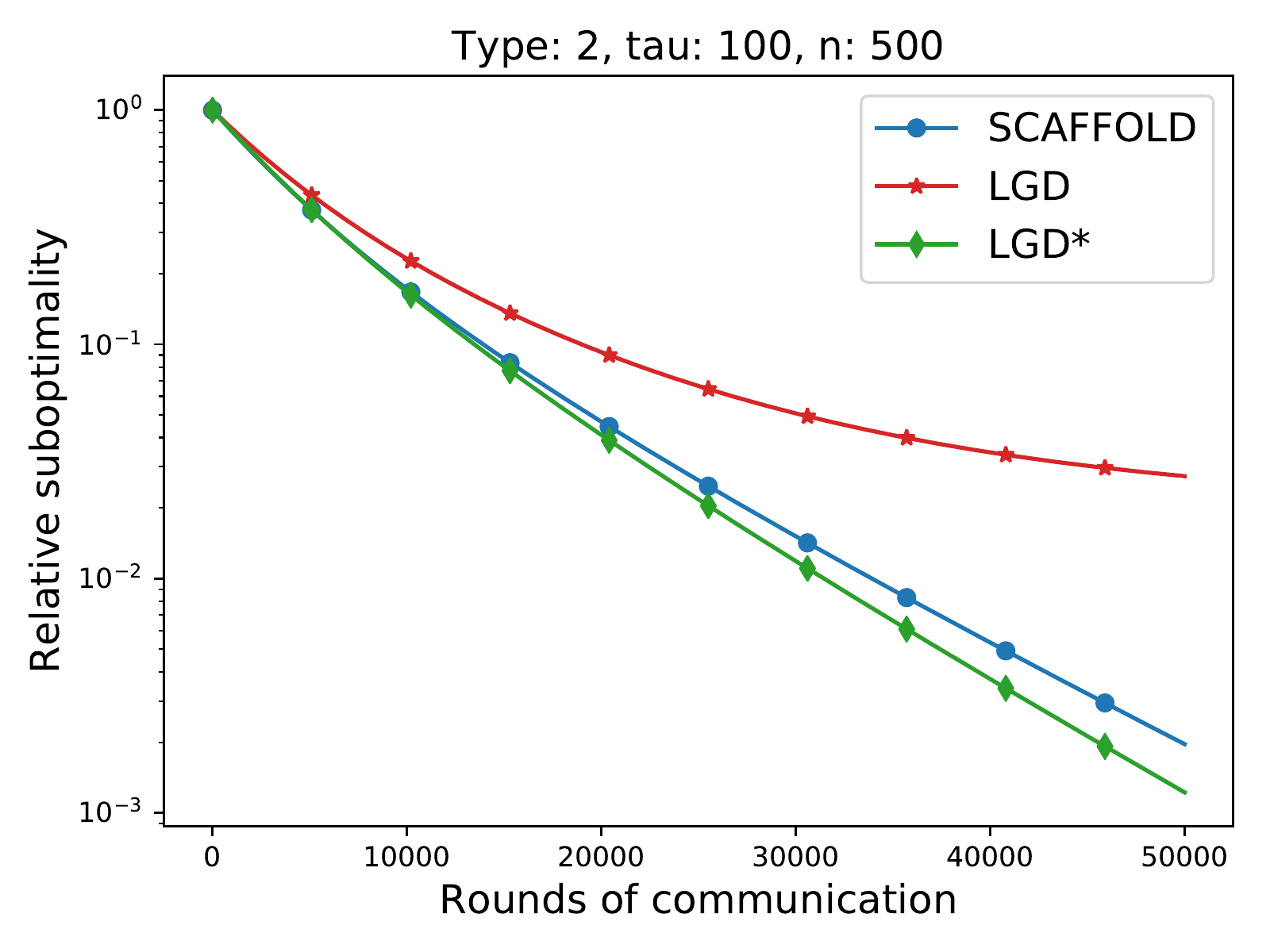}
\end{minipage}
\caption{Comparison of the following noiseless algorithms  {\tt Local-SGD} ({\tt LGD}, Algorithm~\ref{alg:local_sgd} with no local noise) and {\tt SCAFFOLD}~\citep{karimireddy2020scaffold} (Algorithm~\ref{alg:l_local_svrg} without ``Loopless'') and {\tt S*-Local-SGD} ({\tt LGD*}, Algorithm~\ref{alg:local_sgd_star}). Quadratic minimization, problem type 2 (see Table~\ref{tbl:instances}). }
\label{fig:artif3}
\end{figure}

\begin{figure}[!h]
\centering
\begin{minipage}{0.3\textwidth}
  \centering
\includegraphics[width =  \textwidth ]{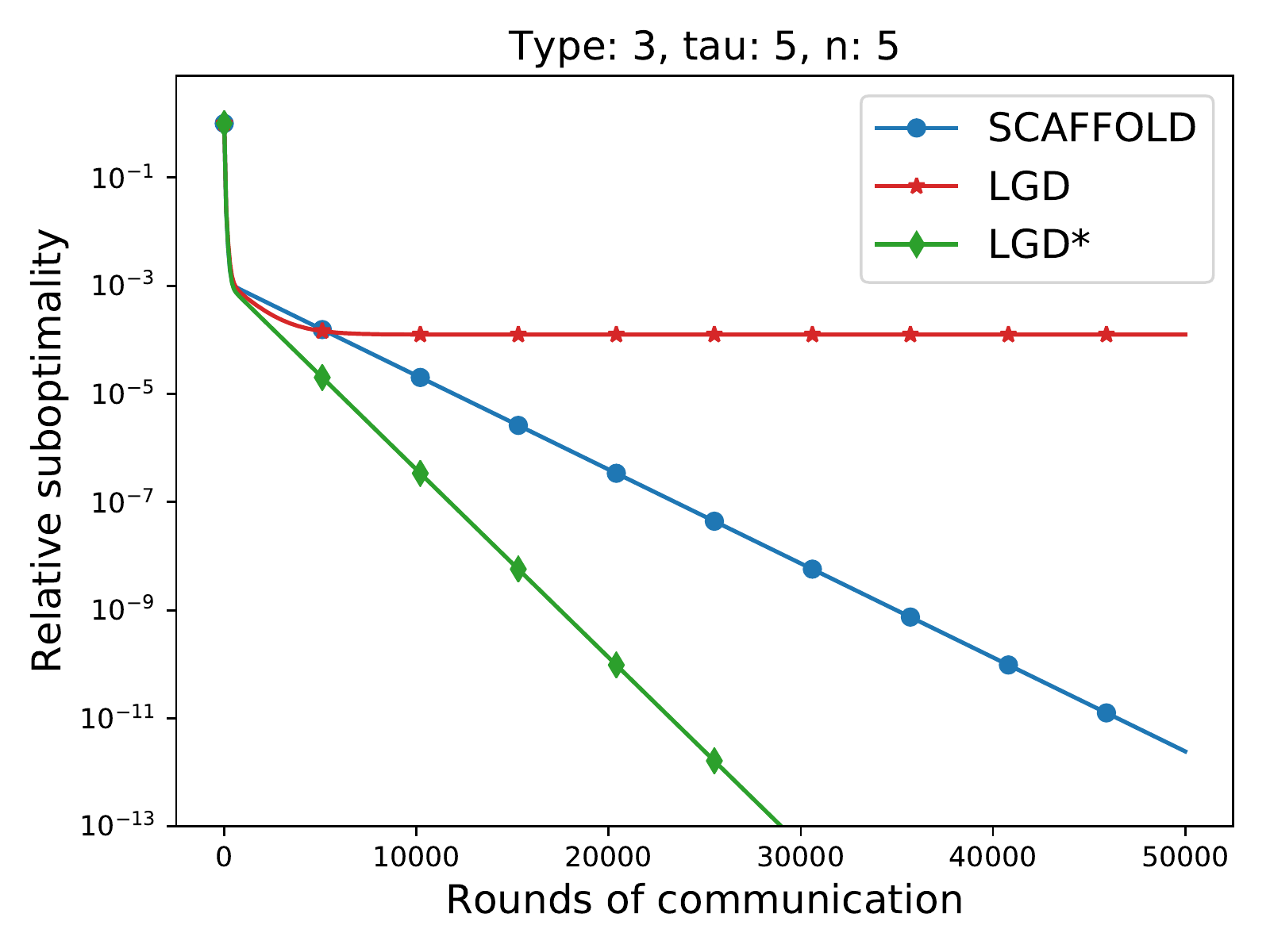}
\end{minipage}
\begin{minipage}{0.3\textwidth}
  \centering
\includegraphics[width =  \textwidth ]{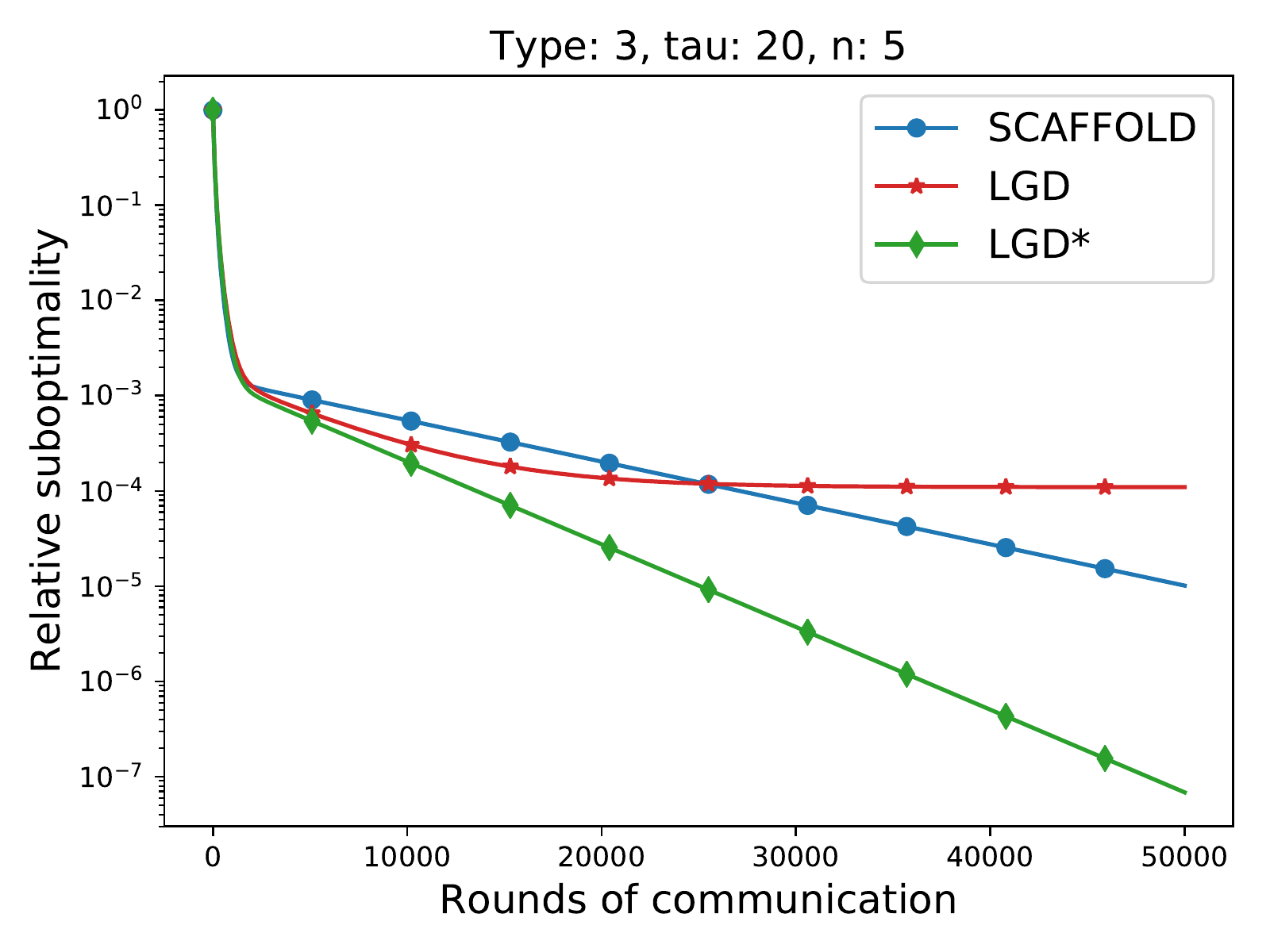}
\end{minipage}
\begin{minipage}{0.3\textwidth}
  \centering
\includegraphics[width =  \textwidth ]{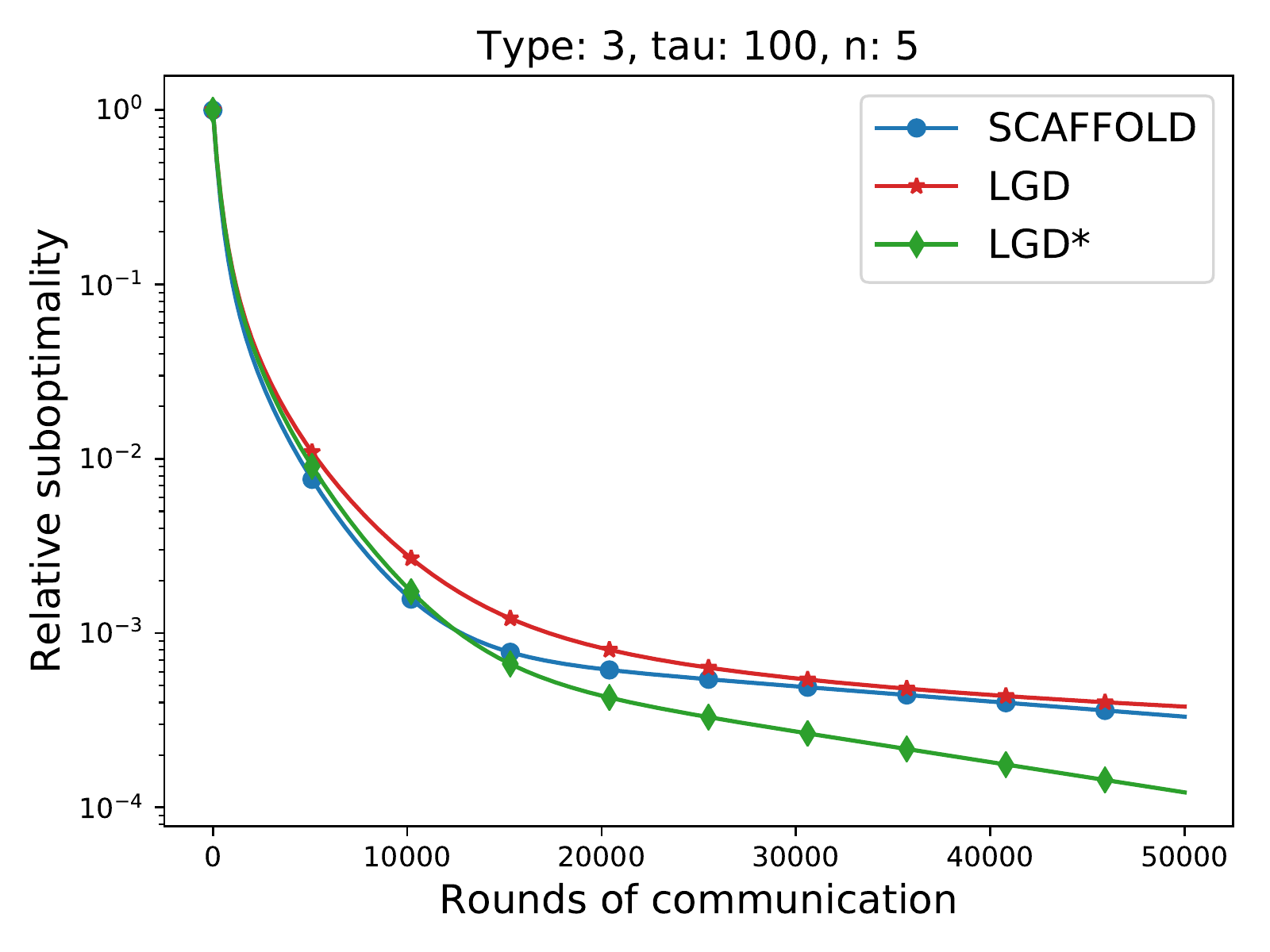}
\end{minipage}
\\
\begin{minipage}{0.3\textwidth}
  \centering
\includegraphics[width =  \textwidth ]{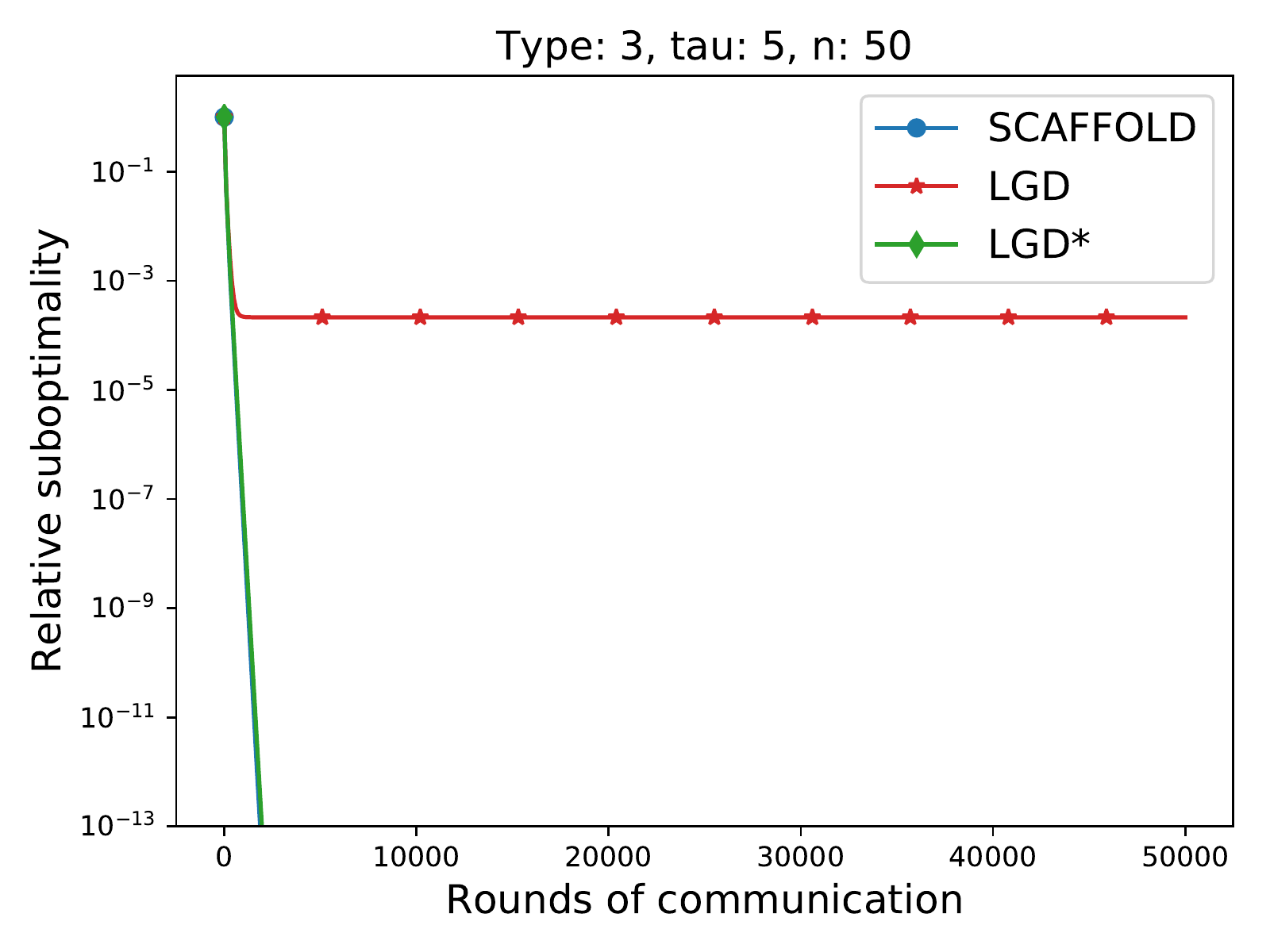}
\end{minipage}
\begin{minipage}{0.3\textwidth}
  \centering
\includegraphics[width =  \textwidth ]{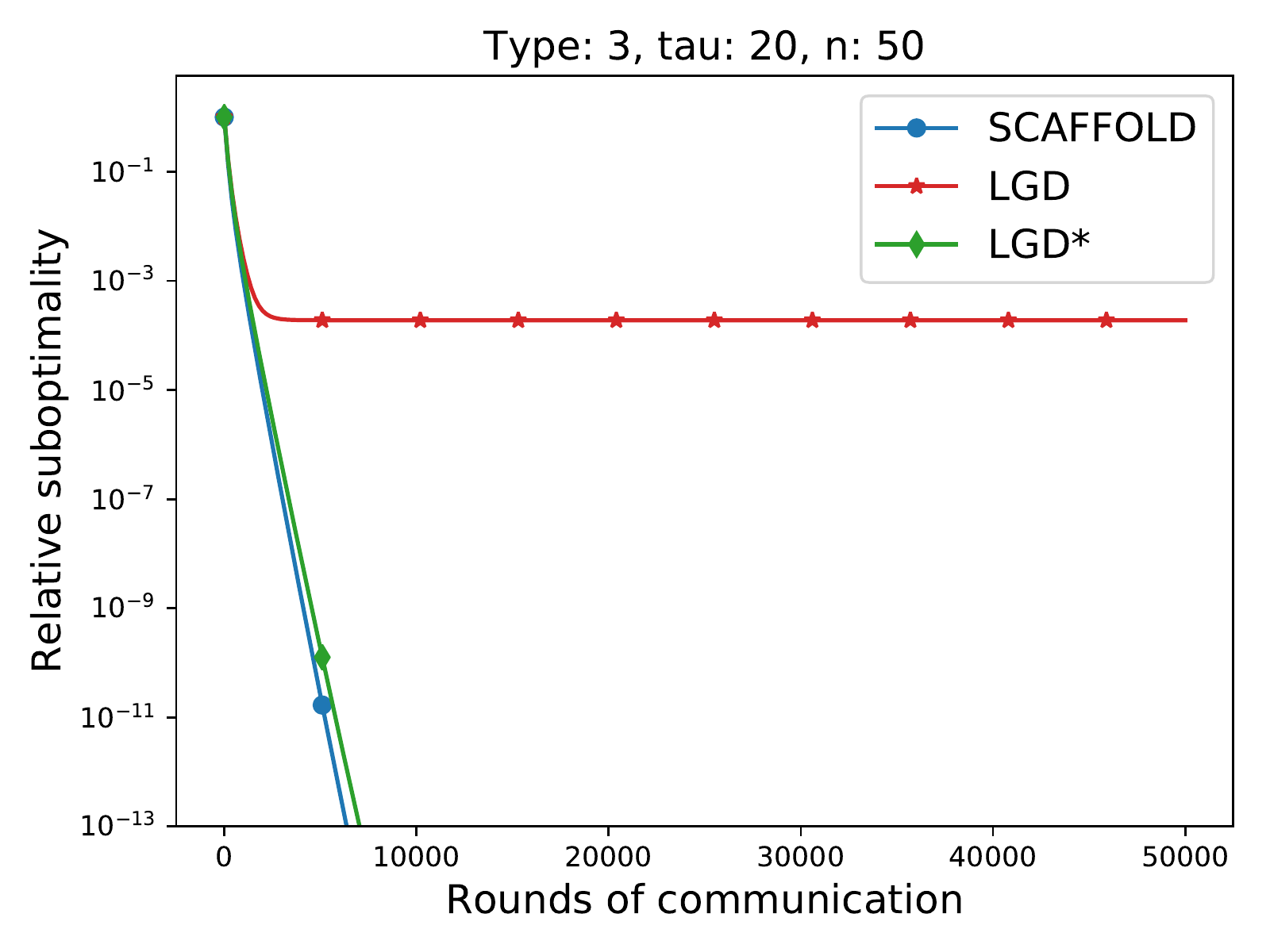}
\end{minipage}
\begin{minipage}{0.3\textwidth}
  \centering
\includegraphics[width =  \textwidth ]{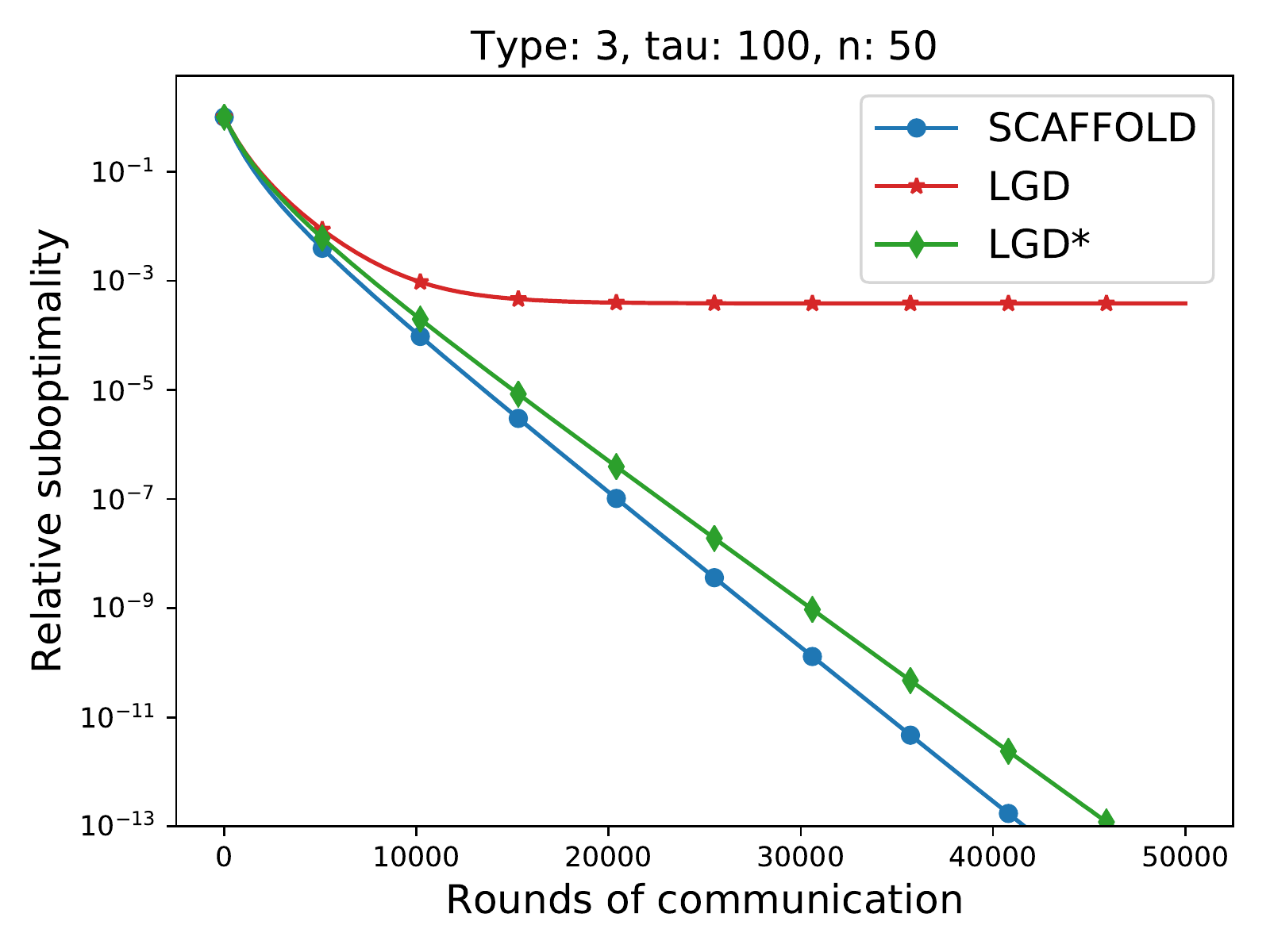}
\end{minipage}
\\
\begin{minipage}{0.3\textwidth}
  \centering
\includegraphics[width =  \textwidth ]{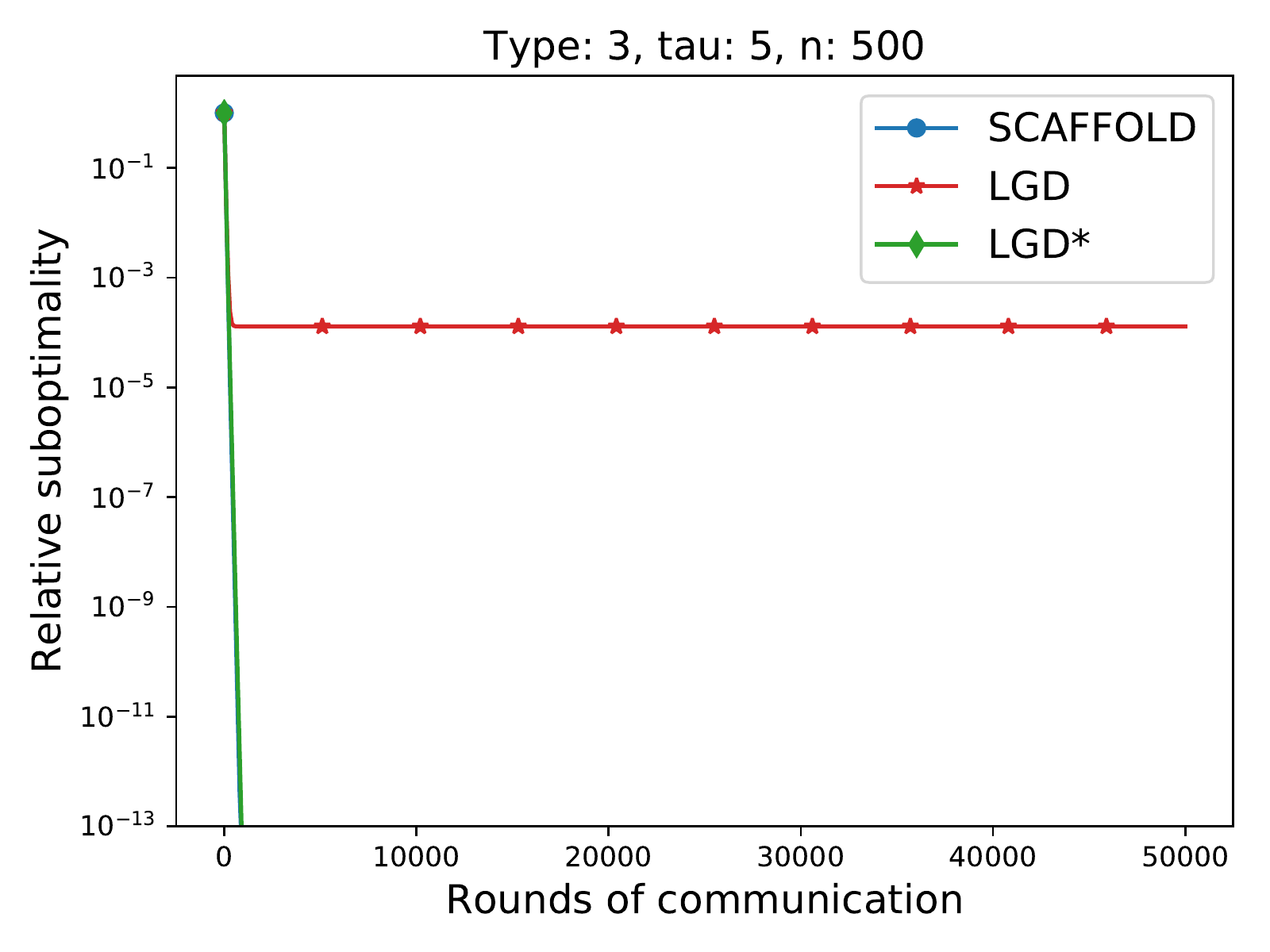}
\end{minipage}
\begin{minipage}{0.3\textwidth}
  \centering
\includegraphics[width =  \textwidth ]{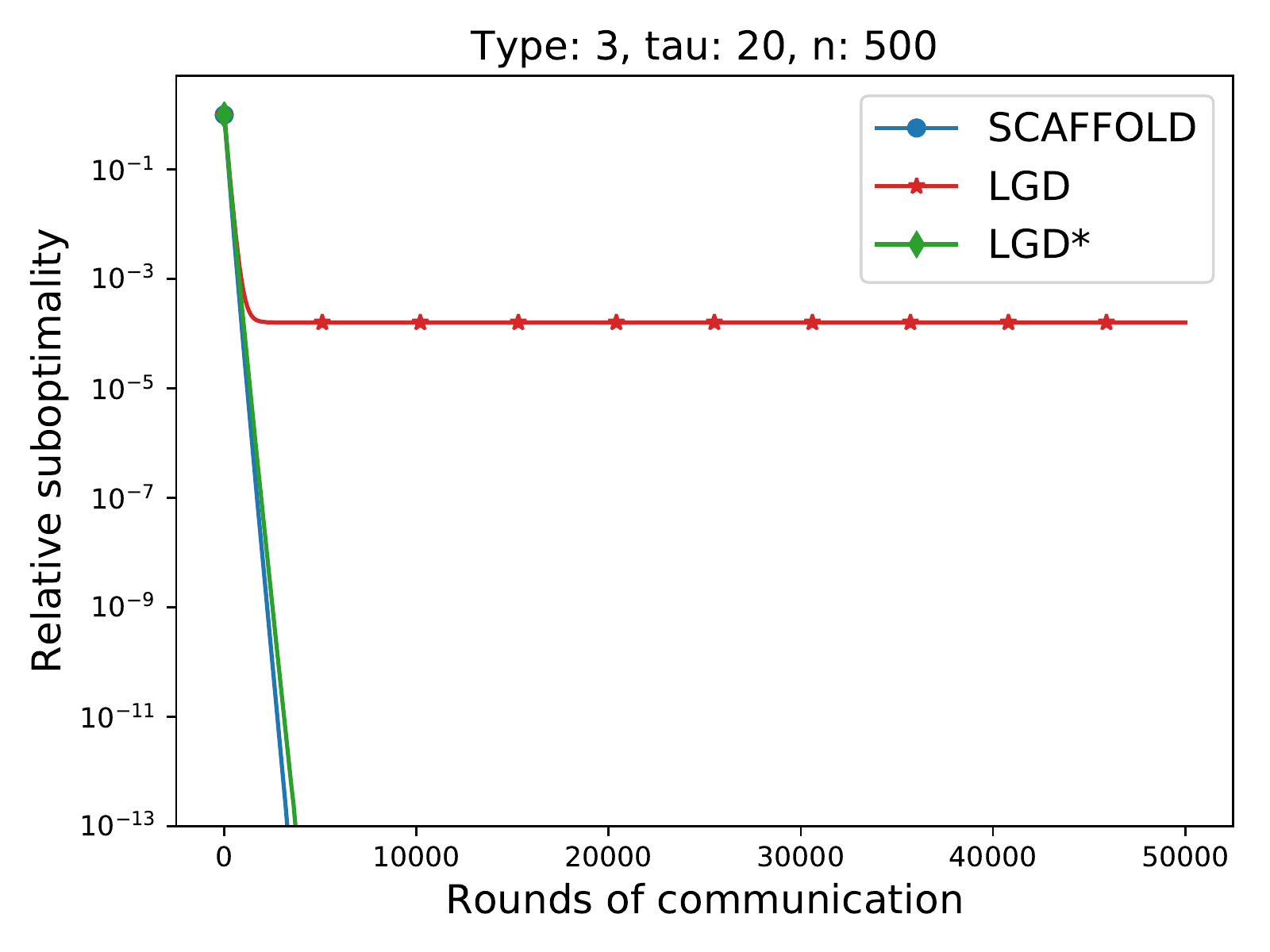}
\end{minipage}
\begin{minipage}{0.3\textwidth}
  \centering
\includegraphics[width =  \textwidth ]{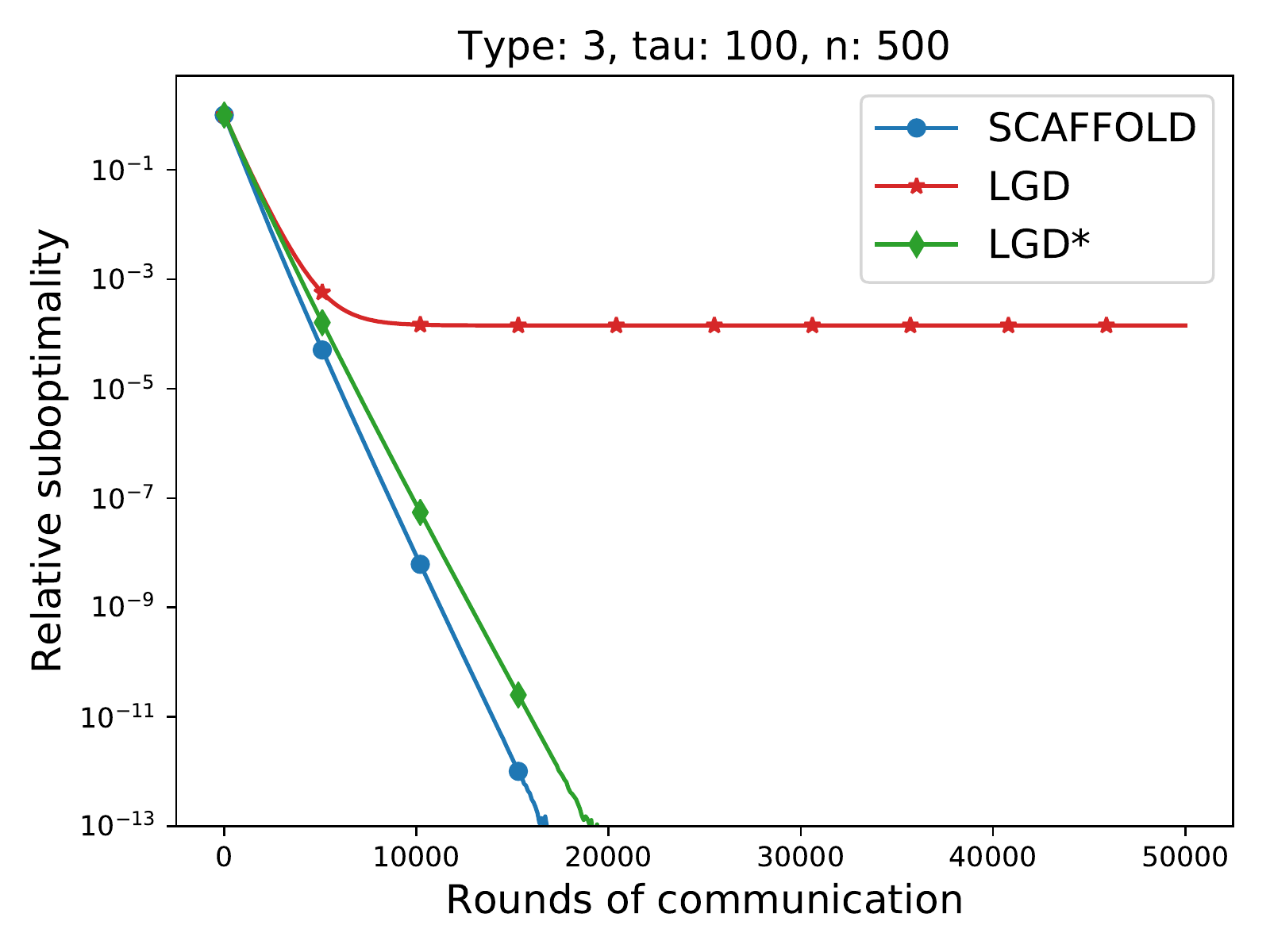}
\end{minipage}
\caption{Comparison of the following noiseless algorithms:  {\tt Local-SGD} ({\tt LGD}, Algorithm~\ref{alg:local_sgd} with no local noise) and {\tt SCAFFOLD}~\citep{karimireddy2020scaffold} (Algorithm~\ref{alg:l_local_svrg} without ``Loopless'') and {\tt S*-Local-SGD} ({\tt LGD*}, Algorithm~\ref{alg:local_sgd_star}). Quadratic minimization, problem type 3 (see Table~\ref{tbl:instances}). }
\label{fig:artif4}
\end{figure}

In this experiment, we therefore compare {\tt Local-SGD},  {\tt S*-Local-SGD} and {\tt SCAFFOLD}. In order to decouple the local variance with the non-stationarity of the local methods, we let each algorithm access the full local gradients. Next, in order to have a full control of the setting, we let the local objectives to be artificially generated quadratic problems. Specifically, we set
\begin{equation}\label{eq:quadproblem}
f_i (x)= \frac{\mu}{2} \| x\|^2 + \frac{1-\mu}{2} (x-z_i^*)^\top \left( \sum_{j=1}^m a_i a_i^\top \right)(x-z_i^*),
\end{equation}
where $a_i$ are mutually orthogonal vectors of norm 1 with $m<d$ (generated by orthogonalizing Gaussian vectors), $z_i^*$ are Gaussian vectors and $\mu = 10^{-3}$. We consider four different instances of~\eqref{eq:quadproblem} given by Table~\ref{eq:quadproblem}. Figures~\ref{fig:artif1},~\ref{fig:artif2},~\ref{fig:artif3},~\ref{fig:artif4} show the result.

Through most of the plots across all combinations of type, $\tau$, $n$, we can see that {\tt Local-SGD}  suffers greatly from the fact that it is attracted to an incorrect fixed point and as a result, it never converges to the exact optimum. On the other hand, both {\tt S*-Local-SGD}  and {\tt SCAFFOLD} converge to the exact optimum and therefore outperform {\tt Local-SGD}  in most examples. We shall note that the rate of {\tt SCAFFOLD} involves  slightly worse constants than those in {\tt Local-SGD}  and {\tt S*-Local-SGD}, and therefore it sometimes performs worse in the early stages of the optimization process when compared to the other methods. Furthermore, notice that our method {\tt S*-Local-SGD}  always performed best.

To summarize, our results demonstrate that \begin{itemize}
\item [(i)] the incorrect fixed point of used by standard local methods is an issue not only theory but also in practice, and should be addressed if better performance is required, 
\item [(ii)] the theoretically optimal shift employed by {\tt S*-Local-SGD} is ideal from a performance perspective if it was available (however, this strategy is impractical to implement as the optimal shift presumes the knowledge of the optimal solution), and 
\item [(iii)] {\tt SCAFFOLD}/{\tt SS-Local-SGD} is a practical solution to fixing the incorrect fixed point problem -- it converges to the exact optimum at a price of a slightly worse initial convergence speed. 
\end{itemize}

 \begin{table}[!t]
 \caption{Instances of~\eqref{eq:quadproblem}. }
\label{tbl:instances}
\begin{center}
\small
\begin{tabular}{|c|c|c|}
\hline
 {\bf Type}  & $m$ &$z_i^*$   \\
 \hline
  \hline
0 & 1 & $\sim \cN(0,\mI)$ \\
\hline
1 & 10 & $\sim \cN(0,\mI)$ \\
\hline
2 & 1 & $\sim \cN(0,\mI)$ \\
\hline
3 & 10 & $\sim \cN(0,\mI)$ \\
\hline
\end{tabular}
\end{center}
\end{table}

\section{Missing Proofs for Section~\ref{sec:main_res}}
Let us first state some well-known consequences of $L$-smoothness. Specifically, if $f_i$ is $L$-smooth, we must have
\begin{equation}
	f_i(y) \le f_i(x) + \langle\nabla f_i(x), y-x \rangle + \frac{L}{2}\|x-y\|^2, \qquad \forall x,y\in \R^d.\label{eq:L_smoothness_cor_1}
\end{equation}
If in addition to this we assume that $f_i$ is  convex, the following bound holds:
\begin{equation}
	\|\nabla f_i(x) - \nabla f_i(y)\|^2 \le 2L(f_i(x) - f_i(y) - \langle\nabla f_i(y), x-y\rangle) \eqdef 2LD_{f_i}(x,y), \qquad \forall x,y\in \R^d \label{eq:L_smoothness_cor}
\end{equation}

We next proceed with the proof of Theorem~\ref{thm:main_result}. Following the technique of virtual iterates from~\citep{stich2020error,khaled2020tighter}, notice that the sequence $\{x^k\}_{k\ge 0}$ satisfies the recursion
\begin{equation}
	x^{k+1} = x^k - \frac{\gamma}{n}\sum\limits_{i=1}^ng_i^k .\label{eq:x^k_recurrsion}
\end{equation}
 This observation forms the backbone of the key lemma of our paper, which we present next. 

\begin{lemma}\label{lem:main_lemma}
	Let Assumption~\ref{ass:quasi_strong_convexity_local},~\ref{ass:L_smoothness}~and~\ref{ass:key_assumption} be satisfied and $\gamma \le \min\left\{\nicefrac{1}{2(A'+MC)},\nicefrac{L}{(F'+MG)}\right\}$, where $M = \frac{4B'}{3\rho}$. Let  $\eta \eqdef\min\left\{\gamma\mu, \frac{\rho}{4}\right\} $. Then for all $k\ge 0$ we have
	\begin{equation}
		\gamma\EE\left[f(x^k) - f(x^*)\right] \le (1-\eta)\EE T^{k} - \EE T^{k+1} + \gamma^2(D_1' + MD_2) + 2L\gamma \EE V_k, \label{eq:main_lemma}
	\end{equation}
	where $\eta \eqdef \min\left\{\gamma\mu,\frac{\rho}{4}\right\}$, $T^k \eqdef \|x^k - x^*\|^2 + M\gamma^2 \sigma_k^2$.
\end{lemma}
\begin{proof}
	First of all, to simplify the proofs we introduce new notation: $g^k \eqdef \frac{1}{n}\sum_{i=1}^n g_i^k$. Using this and \eqref{eq:x^k_recurrsion} we get
	\begin{eqnarray*}
		\|x^{k+1}-x^*\|^2 &\overset{\eqref{eq:x^k_recurrsion}}{=}& \left\|x^k - x^* - \gamma g^k\right\|^2\\
		&=& \|x^k - x^*\|^2 -2\gamma\langle x^k - x^*, g^k\rangle + \gamma^2\|g^k\|^2.
	\end{eqnarray*}
	Taking conditional mathematical expectation $\EE_k[\cdot] = \EE[\cdot\mid x^k] \eqdef \EE[\cdot\mid x_1^k,\ldots, x_n^k]$ on both sides of the previous inequality we get
	\begin{eqnarray*}
		\EE\left[\|x^{k+1} - x^*\|^2\mid x^k\right] &\overset{\eqref{eq:unbiasedness}}{=}& \|x^k-x^*\|^2 -\frac{2\gamma}{n}\sum\limits_{i=1}^n\left\langle x^k-x^*,\nabla f_i(x_i^k) \right\rangle + \gamma^2\EE\left[\|g^k\|^2\mid x^k\right],
	\end{eqnarray*}
	hence
	\begin{eqnarray}
		\EE\left[\|x^{k+1}-x^*\|^2\right] &\overset{\eqref{eq:tower_property}}{\le}& \EE\left[\|x^k-x^*\|^2\right] -\frac{2\gamma}{n}\sum\limits_{i=1}^n\EE\left[\left\langle x^k-x^*,\nabla f_i(x_i^k) \right\rangle\right] + \gamma^2\EE\left[\|g^k\|^2\right]\notag\\
		&\overset{\eqref{eq:second_moment_bound}}{\le}& \EE\left[\|x^k-x^*\|^2\right] -\frac{2\gamma}{n}\sum\limits_{i=1}^n\EE\left[\left\langle x^k-x^*,\nabla f_i(x_i^k) \right\rangle\right] + B'\gamma^2\EE\left[\sigma_k^2\right] \notag\\
		&&\quad +2A'\gamma^2\EE\left[f(x^k) - f(x^*)\right] + F'\gamma^2\EE\left[V_k\right] + \gamma^2D_1'. \label{eq:main_lemma_technical_1}
	\end{eqnarray}
	Next, we derive an upper bound for the second term on the right-hand side of the previous inequality:
	\begin{eqnarray}
		-\frac{2\gamma}{n}\sum\limits_{i=1}^n\left\langle x^k-x^*,\nabla f_i(x_i^k) \right\rangle &=& \frac{2\gamma}{n}\sum\limits_{i=1}^n\left(\left\langle x^* - x_i^k, \nabla f_i(x_i^k)\right\rangle + \left\langle x_i^k - x^k, \nabla f_i(x_i^k)\right\rangle\right)\notag\\
		&\overset{\eqref{eq:str_quasi_cvx},\eqref{eq:L_smoothness_cor_1}}{\le}& \frac{2\gamma}{n}\sum\limits_{i=1}^n \left(f_i(x^*) - f_i(x_i^k) - \frac{\mu}{2}\|x_i^k - x^*\|^2\right)\notag\\
		&&\quad + \frac{2\gamma}{n}\sum\limits_{i=1}^n\left(f_i(x_i^k) - f_i(x^k) + \frac{L}{2}\|x^k - x_i^k\|^2\right)\notag\\
		&\overset{\eqref{eq:a_b_norm_squared}}{\le}& -2\gamma\left(f(x^k) - f(x^*)\right) -\mu\gamma\|x^k - x^*\|^2 + L\gamma V_k.\label{eq:main_lemma_technical_2}
	\end{eqnarray}
	Plugging \eqref{eq:main_lemma_technical_2} in \eqref{eq:main_lemma_technical_1}, we obtain
	\begin{eqnarray}
		\EE\left[\|x^{k+1} - x^*\|^2\right] &\overset{\eqref{eq:main_lemma_technical_1},\eqref{eq:main_lemma_technical_2}}{\le}& (1-\gamma\mu)\EE\left[\|x^k - x^*\|^2\right] -2\gamma\left(1 - A'\gamma\right)\EE\left[f(x^k) - f(x^*)\right]\notag\\
		&&\quad + B'\gamma^2\EE\left[\sigma_k^2\right] + \gamma\left(L+F'\gamma\right)\EE\left[V_k\right] + \gamma^2D_1'.\label{eq:main_lemma_technical_3}
	\end{eqnarray}
	It implies that
	\begin{eqnarray*}
		\EE T^{k+1} &=& \EE\left[\|x^{k+1}-x^*\|^2\right] + M\gamma^2\EE\left[\sigma_{k+1}^2\right]\\
		&\overset{\eqref{eq:main_lemma_technical_3},\eqref{eq:sigma_k+1_bound}}{\le}& (1-\gamma\mu)\EE\|x^k - x^*\|^2 + \left(1+\frac{B'}{M}-\rho\right)M\gamma^2\EE\sigma_k^2 \\
		&&\quad -2\gamma\left(1 - \left(A'+MC\right)\gamma\right)\EE\left[f(x^k) - f(x^*)\right]\\
		&&\quad + \gamma\left(L + (F'+MG)\gamma\right)\EE V_k + \gamma^2\left(D_1' + MD_2\right).
	\end{eqnarray*}
	Since $M = \frac{4B'}{3\rho}$, $\eta = \min\left\{\gamma\mu, \frac{\rho}{4}\right\}$ and $\gamma \le \min\left\{\nicefrac{1}{2(A'+MC)},\nicefrac{L}{(F'+MG)}\right\}$, we get
	\begin{eqnarray*}
		\EE T^{k+1} &\le& (1-\gamma\mu)\EE\|x^k - x^*\|^2 + \left(1-\frac{\rho}{4}\right)M\gamma^2\EE\sigma_k^2 -\gamma\EE\left[f(x^k) - f(x^*)\right]\\
		&&\quad + 2L\gamma \EE V_k + \gamma^2\left(D_1' + MD_2\right)\\
		&\le& (1-\eta)\EE T^k  -\gamma\EE\left[f(x^k) - f(x^*)\right] +  2L\gamma \EE V_k + \gamma^2\left(D_1' + MD_2\right).
	\end{eqnarray*}
	Rearranging the terms we get \eqref{eq:main_lemma}.
\end{proof}

Using the above lemma we derive the main complexity result.

\subsection{Proof of Theorem~\ref{thm:main_result}}
	From Lemma~\ref{lem:main_lemma} we have that
	\begin{eqnarray*}
		\gamma\EE\left[f(x^k) - f(x^*)\right] \le (1-\eta)\EE T^{k} - \EE T^{k+1} + \gamma^2(D_1' + MD_2) + 2L\gamma \EE V_k.
	\end{eqnarray*}
	Summing up previous inequalities for $k=0,\ldots,K$ with weights $w_k$ defined in \eqref{eq:w_k_definition} we derive
	\begin{eqnarray*}
		\gamma\sum\limits_{k=0}^K w_k\EE\left[f(x^k) - f(x^*)\right] &\le& \sum\limits_{k=0}^K \left(w_k(1-\eta)\EE T^{k} - w_k\EE T^{k+1}\right) + \gamma^2(D_1' + MD_2)W_K\\
		&&\quad + 2L\gamma\sum\limits_{k=0}^K w_k\EE V_k\\
		&\overset{\eqref{eq:w_k_definition},\eqref{eq:sum_V_k_bounds}}{\le}& \sum\limits_{k=0}^K \left(w_{k-1}\EE T^{k} - w_k\EE T^{k+1}\right) + \gamma^2\left(D_1' + MD_2\right)W_K\\
		&&\quad + \frac{\gamma}{2}\sum\limits_{k=0}^Kw_k\EE\left[f(x^k)-f(x^*)\right] + 2LH\gamma\EE\sigma_0^2 +2L\gamma^3 D_3 W_K.
	\end{eqnarray*}
	Relations $T^k \ge 0$ and $w_{-1} = 1$ imply that
	\begin{eqnarray*}
		\frac{\gamma}{2}\sum\limits_{k=0}^Kw_k\EE\left[f(x^k)-f(x^*)\right] &\le& T^0 + 2LH\gamma\EE\sigma_0^2 + \gamma^2\left(D_1' + MD_2 + 2L\gamma D_3\right)W_K.
	\end{eqnarray*}
	Using the definition of $\overline{x}^K$ and convexity of $f$, we get 
		\begin{eqnarray}
		\EE\left[f(\overline{x}^K) - f(x^*)\right] &\le& \frac{2T^0 + 4LH\gamma\EE\sigma_0^2}{\gamma W_K} + 2\gamma\left(D_1' + MD_2 + 2L\gamma D_3\right). \label{eq:main_result}
	\end{eqnarray}
	 It remains to consider two cases: $\mu > 0$ and $\mu = 0$. If $\mu > 0$ we have $W_K \ge w_K \ge (1-\eta)^{-K}$, where $\eta \eqdef \min\left\{\gamma\mu, \frac{\rho}{4}\right\}$ which implies \eqref{eq:main_result_1}. Finally, when $\mu = 0$, we have $w_k = 1$ for all $k\ge 0$, which implies $W_K = K+1 \ge K$ and \eqref{eq:main_result_2}.

\subsection{Corollaries}\label{sec:corollaries}
We state the full complexity results that can be obtained from Theorem~\ref{thm:main_result}. These results can be obtained as a direct consequence of Lemmas~\ref{lem:lemma2_stich_local}~and~\ref{lem:lemma_technical_cvx}.

\begin{corollary}\label{cor:app_complexity_cor_str_cvx}
	Consider the setup from Theorem~\ref{thm:main_result} and denote $\frac1h$ to be the resulting upper bound on $\gamma$\footnote{In order to obtain tight estimate of parameters $D_3$ and $H$, we shall impose further bounds on $\gamma$ (see Section~\ref{sec:data_and_loop} and Table~\ref{tbl:data_loop} therein). }
	  and $\mu > 0$.
	\begin{enumerate}
		\item If $D_3$ does not depend on $\gamma$, then for all $K$ such that
		\begin{eqnarray*}
		\text{either} && \frac{\ln\left(\max\{2,\min\{\nicefrac{a\mu^2K^2}{c_1},\nicefrac{a\mu^3K^3}{c_2}\}\}\right)}{K}\le \rho\\
		\text{or} && \frac{1}{h}\le \frac{\ln\left(\max\{2,\min\{\nicefrac{a\mu^2K^2}{c_1},\nicefrac{a\mu^3K^3}{c_2}\}\}\right)}{\mu K},
	\end{eqnarray*}			
	$a = 2\|x^0 - x^*\|^2 +   \frac{8B'\EE\sigma_0^2}{3h^2\rho} + \frac{4LH\EE\sigma_0^2}{h}$, $c_1 = 2D_1' + \frac{4B'D_2}{3\rho}$, $c_2 = 4LD_3$ and
	\begin{eqnarray*}
		\gamma &=& \min\left\{\frac{1}{h}, \gamma_K\right\},\\
		\gamma_K &=& \frac{\ln\left(\max\left\{2,\min\left\{\frac{a\mu^2K^2}{c_1},\frac{a\mu^3K^3}{c_2}\right\}\right\}\right)}{\mu K},
	\end{eqnarray*}
	we have\footnote{$\widetilde{\cO}$ hides numerical constants and logarithmical factors depending on $K$ and parameters of the problem.}
	\begin{equation*}
	\EE\left[f(\overline{x}^K) \right] - f(x^*) =	\widetilde\cO\left(ha\exp\left(-\min\left\{\frac{\mu}{h}, \rho\right\}K\right) + \frac{c_1}{\mu K} + \frac{c_2}{\mu^2 K^2}\right).
	\end{equation*}
	That is, to achieve $\EE\left[f(\overline{x}^K) \right] - f(x^*)\le \varepsilon$, the method requires\footnote{If $c_1 = c_2 = 0$, then one can replace $\widetilde{\cO}$ by $\cO$.}:
	\begin{equation*}
		K = \widetilde\cO\left(\left(\frac{1}{\rho} + \frac{h}{\mu}\right)\log\left(\frac{ha}{\varepsilon}\right) + \frac{c_1}{\mu \varepsilon} + \sqrt{\frac{c_2}{\mu^2 \varepsilon}}\right).
	\end{equation*}
	\item If $D_3 = D_{3,1} + \frac{D_{3,2}}{\gamma}$, then the same bounds hold with $c_1 = 2D_1' + \frac{4B'D_2}{3\rho} + 2LD_{3,2}$ and $c_2 = 4LD_{3,1}$. 
	\end{enumerate}
\end{corollary}

\begin{corollary}\label{cor:app_complexity_cor_cvx}
	Let assumptions of Theorem~\ref{thm:main_result} be satisfied with any $\gamma \le \frac{1}{h}$ and $\mu = 0$.
	\begin{enumerate}
		\item If $D_3$ does not depend on $\gamma$, then for all $K$ and
	\begin{eqnarray*}
		\gamma &=& \min\left\{\frac{1}{h}, \sqrt{\frac{a}{b_1}}, \sqrt[3]{\frac{a}{b_2}}, \sqrt{\frac{a}{c_1 K}}, \sqrt[3]{\frac{a}{c_2 K}}\right\},
	\end{eqnarray*}
	where $a = 2\|x^0 - x^*\|^2$, $b_1 = 4LH\EE\sigma_0^2$, $b_2 = \frac{8B'\EE\sigma_0^2}{3\rho}$, $c_1 = 2D_1' + \frac{4B'D_2}{3\rho}$, $c_2 = 4LD_3$, we have 	\begin{equation*}
\EE\left[f(\overline{x}^K) \right] - f(x^*) =		\cO\left(\frac{ha}{K} + \frac{\sqrt{ab_1}}{K} + \frac{\sqrt[3]{a^2b_2}}{K} + \sqrt{\frac{ac_1}{K}} + \frac{\sqrt[3]{a^2c_2}}{K^{\nicefrac{2}{3}}} \right).
	\end{equation*}
	That is, to achieve $\EE\left[f(\overline{x}^K) \right] - f(x^*)\le \varepsilon$, the method requires 
	\begin{equation*}
		K = \cO\left(\frac{ha}{\varepsilon} + \frac{\sqrt{ab_1}}{\varepsilon} + \frac{\sqrt[3]{a^2b_2}}{\varepsilon} + \frac{ac_1}{\varepsilon^2} + \frac{a\sqrt{c_2}}{\varepsilon^{\nicefrac{3}{2}}} \right).
	\end{equation*}
	\item If $D_3 = D_{3,1} + \frac{D_{3,2}}{\gamma}$, then the same bounds hold with $c_1 = 2D_1' + \frac{4B'D_2}{3\rho} + 2LD_{3,2}$ and $c_2 = 4LD_{3,1}$. 
	\end{enumerate}
\end{corollary}

\section{Missing Proofs and Details for Section~\ref{sec:data_and_loop} \label{sec:a_data_and_loop}}

\subsection{Constant Local Loop}\label{sec:constant_loop}
In this section we show how our results can be applied to analyze \eqref{eq:local_sgd_def} in the case when
\begin{equation*}
	c_{k} = \begin{cases}1,& \text{if } k \mod \tau = 0,\\ 0,& \text{if } k\mod \tau \neq 0, \end{cases}
\end{equation*}
where $\tau$ is number of local steps between two neighboring rounds of communications. This corresponds to the setting in which the local loop size on each device has a fixed length.

\subsubsection{Heterogenous Data}\label{sec:const_loop_hetero}
First of all, we need to assume more about $g_i^k$.
\begin{assumption}\label{ass:hetero_second_moment}
	We assume that inequalities \eqref{eq:second_moment_bound}-\eqref{eq:sigma_k+1_bound} hold and additionally there exist such non-negative constants $\tA, \hA, \tB, \hB, \tF, \hF, \tD_1, \hD_{1}$ that for all $k \ge 0$
	\begin{eqnarray}
		\frac{1}{n}\sum\limits_{i=1}^n\EE\left[\|\bar{g}_i^k\|^2\right] &\le & 2\tA\EE\left[f(x^k) - f(x^*)\right] + \tB\EE\left[\sigma_k^2\right] + \tF\EE\left[V_k\right] + \tD_{1}, \label{eq:hetero_second_moment_bound}\\
		\frac{1}{n}\sum\limits_{i=1}^n\EE\left[\|g_i^k-\bar{g}_i^k\|^2\right] &\le & 2\hA\EE\left[f(x^k) - f(x^*)\right] + \hB\EE\left[\sigma_k^2\right] + \hF\EE\left[V_k\right] + \hD_{1}, \label{eq:hetero_var_bound}
	\end{eqnarray}
	where $\bar{g}_i^k = \EE\left[g_i^k\mid x_1^k,\ldots,x_n^k\right]$.
\end{assumption}
We notice that inequalities \eqref{eq:hetero_second_moment_bound}-\eqref{eq:hetero_var_bound} imply \eqref{eq:second_moment_bound} and vice versa. Indeed, if \eqref{eq:hetero_second_moment_bound}-\eqref{eq:hetero_var_bound} hold then inequality \eqref{eq:second_moment_bound} holds with $A = \tA + \hA$, $B = \tB+\hB$, $F = \tF + \hF$, $D_1 = \tD_{1}+\hD_{1}$ due to variance decomposition formula \eqref{eq:variance_decomposition}, and if \eqref{eq:second_moment_bound} is true then \eqref{eq:hetero_second_moment_bound}-\eqref{eq:hetero_var_bound} also hold with $\tA = \hA = A$, $\tB = \hB = B$, $\tF = \hF = F$, $\tD_{1} = \hD_{1} = D_1$.

We start our analysis without making any assumption on homogeneity of data that workers have an access to. Next lemma provides an upper bound for the weighted sum of $\EE V_k$.
\begin{lemma}\label{lem:V_k_lemma}
	Let Assumption~\ref{ass:quasi_strong_convexity_local},~\ref{ass:L_smoothness}~and~\ref{ass:hetero_second_moment} hold and\footnote{When $\rho = 1$ one can always set the parameters in such a way that $\tB = \hB = C = G = 0$, $D_2 = 0$. In this case we assume that $\frac{2\tB C}{\rho(1-\rho)} = \frac{2\hB C}{\rho(1-\rho)} = \frac{2\tB G}{\rho(1-\rho)} = \frac{2\hB G}{\rho(1-\rho)} = 0$.}
	\begin{eqnarray*}
		\gamma &\le& \min\left\{\frac{1}{4(\tau-1)\mu}, \frac{1}{2\sqrt{e(\tau-1)\left(\tF(\tau-1) + \hF + \frac{2G(\tB(\tau-1)+\hB)}{\rho(1-\rho)}\right)}}\right\},\\
		\gamma &\le& \frac{1}{4\sqrt{2eL(\tau-1)\left(\tA(\tau-1)+\hA +\frac{2C(\tB(\tau-1)+\hB)}{\rho(1-\rho)}\right)}}
	\end{eqnarray*}
	Then \eqref{eq:sum_V_k_bounds} holds with
	\begin{eqnarray}
        H &=& \frac{4e(\tau-1)(\tB(\tau-1)+\hB)(2+\rho)\gamma^2}{\rho},\notag\\
        D_3 &=& 2e(\tau-1)\left(\tD_{1}(\tau-1)+\hD_{1} + \frac{2D_2(\tB(\tau-1)+\hB)}{\rho}\right).\label{eq:V_k_bound}	
	\end{eqnarray}
\end{lemma}
\begin{proof}
	Consider some integer $k\ge 0$. There exists such integer $t\ge 0$ that $\tau t \le k \le \tau(t+1)-1$. Using this and Lemma~\ref{lem:lemma14_stich_general} we get
	\begin{eqnarray*}
		\EE[V_k] &\overset{\eqref{eq:local_sgd_def},\eqref{eq:x^k_recurrsion}}{=}& \frac{1}{n}\sum\limits_{i=1}^n \EE\left[\left\|x_i^{\tau t} - \gamma\sum\limits_{l=\tau t}^{k-1} g_i^l - x^{\tau t} + \gamma\sum\limits_{l=\tau t}^{k-1} g^l\right\|^2\right]\\
		&=& \frac{\gamma^2}{n}\sum\limits_{i=1}^n\EE\left[\left\|\sum\limits_{l=\tau t}^{k-1} \left(g_i^l - g^l\right)\right\|^2\right]\\
		&\overset{\eqref{eq:lemma14_stich}}{\le}& \frac{e\gamma^2(k-\tau t)}{n}\sum\limits_{i=1}^n\sum\limits_{l=\tau t}^{k-1}\EE\left[\left\|\bar{g}_i^l - \bar{g}^l\right\|^2\right] + \frac{e\gamma^2}{n}\sum\limits_{i=1}^n\sum\limits_{l=\tau t}^{k-1}\EE\left[\left\|g_i^l - \bar{g}_i^l - \left(g^l - \bar{g}^l\right)\right\|^2\right]\\
		&\overset{\eqref{eq:variance_decomposition}}{\le}& \frac{e\gamma^2(\tau - 1)}{n}\sum\limits_{i=1}^n\sum\limits_{l=\tau t}^{k-1}\EE\left[\left\|\bar{g}_i^l\right\|^2\right] + \frac{e\gamma^2}{n}\sum\limits_{i=1}^n\sum\limits_{l=\tau t}^{k-1}\EE\left[\left\|g_i^l - \bar{g}_i^l\right\|^2\right],
	\end{eqnarray*}
	where $\bar{g}^k = \frac{1}{n}\sum\limits_{i=1}^n\bar{g}_i^k$.
	Applying Assumption~\ref{ass:hetero_second_moment}, we obtain
	\begin{eqnarray*}
		\EE V_k &\overset{\eqref{eq:hetero_second_moment_bound},\eqref{eq:hetero_var_bound}}{\le}& 2e\left(\tA(\tau - 1)+\hA\right)\gamma^2\sum\limits_{l=\tau t}^{k-1}\EE\left[f(x^l) - f(x^*)\right] + e\left(\tB(\tau-1)+\hB\right)\gamma^2\sum\limits_{l=\tau t}^{k-1}\EE\sigma_l^2 \\
		&&\quad + e\left(\tF(\tau-1)+\hF\right)\gamma^2\sum\limits_{l=\tau t}^{k-1} \EE V_l + e(\tau- 1)\left(\tD_{1}(\tau-1)+\hD_{1}\right)\gamma^2,
	\end{eqnarray*}
	hence
	\begin{eqnarray}
		\sum\limits_{j=\tau t}^k w_j \EE V_j &\le& 2e\left(\tA(\tau - 1)+\hA\right)\gamma^2\sum\limits_{j=\tau t}^k\sum\limits_{l=\tau t}^{j-1}w_j\EE\left[f(x^l) - f(x^*)\right]\notag\\
		&&\quad + e\left(\tB(\tau-1)+\hB\right)\gamma^2\sum\limits_{j=\tau t}^k\sum\limits_{l=\tau t}^{j-1}w_j\EE\sigma_l^2 \notag \\
		&&\quad + e\left(\tF(\tau-1)+\hF\right)\gamma^2\sum\limits_{j=\tau t}^k\sum\limits_{l=\tau t}^{j-1} w_j\EE V_l\notag\\
		&&\quad+ e(\tau- 1)\left(\tD_{1}(\tau-1)+\hD_{1}\right)\gamma^2\sum\limits_{j=\tau t}^kw_j. \label{eq:V_k_lemma_technical_1}
	\end{eqnarray}
	Recall that $w_k = (1 - \eta)^{-(k+1)}$ and $\eta = \min\left\{\gamma\mu, \frac{\rho}{4}\right\}$. Together with our assumption on $\gamma$ it implies that for all $0 \le i < k$, $0\le j \le \tau-1$ we have
	\begin{eqnarray}
		w_k &=& (1 - \eta)^{-(k-j+1)}\left(1 - \eta\right)^{-j} \overset{\eqref{eq:1-p/2_inequality}}{\le} w_{k-j}\left(1 + 2\eta\right)^{j} \notag\\
		&\le& w_{k-j}\left(1 + 2\gamma\mu\right)^{j} \le w_{k-j}\left(1 + \frac{1}{2(\tau-1)}\right)^j \le w_{k-j}\exp\left(\frac{j}{2(\tau-1)}\right)\notag\\
		&\le& w_{k-j}\exp\left(\frac{1}{2}\right) \le 2w_{k-j}, \label{eq:V_k_lemma_technical_2}\\
		w_k &=& \left(1 - \eta\right)^{-(k-i+1)}\left(1 - \eta\right)^{-i} \overset{\eqref{eq:1-p/2_inequality}}{\le} w_{k-i}\left(1 + 2\eta\right)^i \le w_{k-i}\left(1 + \frac{\rho}{2}\right)^i, \label{eq:V_k_lemma_technical_3}\\
		w_k &\overset{\eqref{eq:1-p/2_inequality}}{\le}& \left(1 + 2\eta\right)^{k+1} \le \left(1 + \frac{\rho}{2}\right)^{k+1}. \label{eq:V_k_lemma_technical_4}
	\end{eqnarray}
	For simplicity, we introduce new notation: $r_k \eqdef \EE\left[f(x^k) - f(x^*)\right]$. Using this we get
	\begin{eqnarray*}
		\sum\limits_{j=\tau t}^k\sum\limits_{l=\tau t}^{j-1}w_jr_l &\overset{\eqref{eq:V_k_lemma_technical_2}}{\le}& \sum\limits_{j=\tau t}^k\sum\limits_{l=\tau t}^{j-1}2w_l r_l \le 2(k-\tau t)\sum\limits_{j=\tau t}^kw_jr_j \le 2(\tau - 1)\sum\limits_{j=\tau t}^k w_jr_j,\\
		\sum\limits_{j=\tau t}^k\sum\limits_{l=\tau t}^{j-1}w_j\EE\sigma_l^2 &\overset{\eqref{eq:V_k_lemma_technical_2}}{\le}& \sum\limits_{j=\tau t}^k\sum\limits_{l=\tau t}^{j-1}2w_l \EE\sigma_l^2 \le 2(k-\tau t)\sum\limits_{j=\tau t}^kw_j\EE\sigma_j^2 \le 2(\tau - 1)\sum\limits_{j=\tau t}^kw_j\EE\sigma_j^2,\\
		\sum\limits_{j=\tau t}^k\sum\limits_{l=\tau t}^{j-1}w_j\EE V_l &\overset{\eqref{eq:V_k_lemma_technical_2}}{\le}& \sum\limits_{j=\tau t}^k\sum\limits_{l=\tau t}^{j-1}2w_l \EE V_l \le 2(k-\tau t)\sum\limits_{j=\tau t}^kw_j\EE V_j \le 2(\tau - 1)\sum\limits_{j=\tau t}^kw_j \EE V_j.
	\end{eqnarray*}
	Plugging these inequalities in \eqref{eq:V_k_lemma_technical_1} we derive
	\begin{eqnarray*}
		\sum\limits_{j=\tau t}^k w_j \EE V_j &\le& 4e(\tau - 1)(\tA(\tau-1)+\hA)\gamma^2\sum\limits_{j=\tau t}^kw_j r_j + 2e(\tau - 1)(\tB(\tau-1)+\hB)\gamma^2\sum\limits_{j=\tau t}^kw_j\EE\sigma_j^2 \notag \\
		&&\quad + 2e(\tau - 1)(\tF(\tau-1)+\hF)\gamma^2\sum\limits_{j=\tau t}^k w_j\EE V_j + e\left(\tD_{1}(\tau - 1)+\hD_{1}\right)\gamma^2\sum\limits_{j=\tau t}^kw_j.
	\end{eqnarray*}
	Since $V_{\tau t} = 0$ for all integer $t \ge 0$ we obtain
	\begin{eqnarray}
		\sum\limits_{k=0}^K w_k \EE V_k &\le& 4e(\tau - 1)(\tA(\tau-1)+\hA)\gamma^2\sum\limits_{k=0}^K w_k r_k + 2e(\tau - 1)(\tB(\tau-1)+\hB)\gamma^2\sum\limits_{k=0}^Kw_k\EE\sigma_k^2 \notag \\
		&&\quad + 2e(\tau - 1)(\tF(\tau-1)+\hF)\gamma^2\sum\limits_{k=0}^K w_k\EE V_k\notag\\
		&&\quad + e\left(\tD_{1}(\tau - 1)+\hD_{1}\right)\gamma^2 \sum\limits_{k=0}^K w_k\label{eq:V_k_lemma_technical_5}
	\end{eqnarray}
	It remains to estimate the second term in the right-hand side of the previous inequality. First of all,
	\begin{eqnarray}
		\EE\sigma_{k+1}^2 &\overset{\eqref{eq:sigma_k+1_bound}}{\le}& (1 - \rho)\EE\sigma_{k}^2 + 2C \underbrace{\EE\left[f(x^k) - f(x^*)\right]}_{r_k} + G\EE V_k + D_2\notag\\
		&\le& (1-\rho)^{k+1}\EE\sigma_0^2 + 2C\sum\limits_{l=0}^{k}(1-\rho)^{k-l}r_l + G\sum\limits_{l=0}^k(1-\rho)^{k-l}\EE V_l + D_2\sum\limits_{l=0}^{k}(1-\rho)^l\notag\\
		&\le& (1-\rho)^{k+1}\EE\sigma_0^2 + 2C\sum\limits_{l=0}^{k}(1-\rho)^{k-l}r_l + G\sum\limits_{l=0}^k(1-\rho)^{k-l}\EE V_l + D_2\sum\limits_{l=0}^{\infty}(1-\rho)^l\notag\\
		&=& (1-\rho)^{k+1}\EE\sigma_0^2 + 2C\sum\limits_{l=0}^{k}(1-\rho)^{k-l}r_l + G\sum\limits_{l=0}^k(1-\rho)^{k-l}\EE V_l + \frac{D_2}{\rho}.\label{eq:sigma_k_useful_recurrence}
	\end{eqnarray}
	It implies that
	\begin{eqnarray}
		\sum\limits_{k=0}^K w_k\EE\sigma_k^2 &\overset{\eqref{eq:sigma_k_useful_recurrence}}{\le}& \EE\sigma_0^2\sum\limits_{k=0}^K w_k(1-\rho)^{k} + \frac{2C}{1-\rho}\sum\limits_{k=0}^K\sum\limits_{l=0}^k w_k(1-\rho)^{k-l}r_l\notag\\
		&&\quad + \frac{G}{1-\rho}\sum\limits_{k=0}^K\sum\limits_{l=0}^k w_k(1-\rho)^{k-l}\EE V_l +  \frac{D_2 W_K}{\rho}\notag\\
		&\overset{\eqref{eq:V_k_lemma_technical_3},\eqref{eq:V_k_lemma_technical_4}}{\le}& \EE\sigma_0^2\left(1+\frac{\rho}{2}\right)\sum\limits_{k=0}^K\left(1+\frac{\rho}{2}\right)^k(1-\rho)^k\notag\\
		&&\quad + \frac{2C}{1-\rho}\sum\limits_{k=0}^K\sum\limits_{l=0}^k w_l\left(1+\frac{\rho}{2}\right)^{k-l}(1-\rho)^{k-l}r_l\notag\\
		&&\quad + \frac{G}{1-\rho}\sum\limits_{k=0}^K\sum\limits_{l=0}^k w_l\left(1+\frac{\rho}{2}\right)^{k-l}(1-\rho)^{k-l}\EE V_l +  \frac{D_2 W_K}{\rho}\notag\\
		&\overset{\eqref{eq:1+p/2_inequality}}{\le}& \EE\sigma_0^2\left(1+\frac{\rho}{2}\right)\sum\limits_{k=0}^K\left(1-\frac{\rho}{2}\right)^k + \frac{2C}{1-\rho}\sum\limits_{k=0}^K\sum\limits_{l=0}^k w_lr_l\left(1-\frac{\rho}{2}\right)^{k-l}\notag\\
		&&\quad + \frac{G}{1-\rho}\sum\limits_{k=0}^K\sum\limits_{l=0}^k w_l\EE V_l\left(1-\frac{\rho}{2}\right)^{k-l} +  \frac{D_2 W_K}{\rho}\notag\\
		&\le& \EE\sigma_0^2\left(1+\frac{\rho}{2}\right)\sum\limits_{k=0}^\infty\left(1-\frac{\rho}{2}\right)^k + \frac{2C}{1-\rho}\left(\sum\limits_{k=0}^Kw_kr_k\right)\left(\sum\limits_{l=0}^\infty \left(1-\frac{\rho}{2}\right)^{l}\right)\notag\\
		&&\quad + \frac{G}{1-\rho}\left(\sum\limits_{k=0}^Kw_k\EE V_k\right)\left(\sum\limits_{l=0}^\infty \left(1-\frac{\rho}{2}\right)^{l}\right) +  \frac{D_2 W_K}{\rho}\notag\\
		&=& \frac{\EE\sigma_0^2(2+\rho)}{\rho} + \frac{4C}{\rho(1-\rho)}\sum\limits_{k=0}^Kw_kr_k + \frac{2G}{\rho(1-\rho)}\sum\limits_{k=0}^Kw_k\EE V_k\notag \\
		&&\quad +  \frac{D_2 W_K}{\rho}.\label{eq:sigma_k_technical_bound}
	\end{eqnarray}
	Plugging this inequality in \eqref{eq:V_k_lemma_technical_5} we get
	\begin{eqnarray*}
		\sum\limits_{k=0}^K w_k \EE V_k &\le& 4e(\tau-1)\gamma^2\left(\tA(\tau-1)+\hA +\frac{2C(\tB(\tau-1)+\hB)}{\rho(1-\rho)}\right)\sum\limits_{k=0}^K w_k r_k\\
		&&\quad + \frac{2e(\tau-1)(\tB(\tau-1)+\hB)\EE\sigma_0^2(2+\rho)\gamma^2}{\rho}\\
		&&\quad + 2e(\tau-1)\gamma^2\left(\tF(\tau-1) + \hF + \frac{2G(\tB(\tau-1)+\hB)}{\rho(1-\rho)}\right)\sum\limits_{k=0}^K w_k\EE V_k \\
		&&\quad + e(\tau-1)\gamma^2\left(\tD_{1}(\tau-1)+\hD_{1} + \frac{2D_2(\tB(\tau-1)+\hB)}{\rho}\right)W_K.
	\end{eqnarray*}
	Our choice of $\gamma$ implies
	\begin{equation*}
		4e(\tau-1)\gamma^2\left(\tA(\tau-1)+\hA +\frac{2C(\tB(\tau-1)+\hB)}{\rho(1-\rho)}\right) \le \frac{1}{8L}
	\end{equation*}
	and
	\begin{equation*}
		2e(\tau-1)\gamma^2\left(\tF(\tau-1) + \hF + \frac{2G(\tB(\tau-1)+\hB)}{\rho(1-\rho)}\right) \le \frac{1}{2}.
	\end{equation*}		
	Using these inequalities we continue our derivations
	\begin{eqnarray*}
		\frac{1}{2}\sum\limits_{k=0}^K w_k\EE V_k &\le& \frac{1}{8L}\sum\limits_{k=0}^K w_k r_k + \frac{2e(\tau-1)(\tB(\tau-1)+\hB)\EE\sigma_0^2(2+\rho)\gamma^2}{\rho} \\
		&&\quad+ e(\tau-1)\gamma^2\left(\tD_{1}(\tau-1)+\hD_{1} + \frac{2D_2(\tB(\tau-1)+\hB)}{\rho}\right)W_K.
	\end{eqnarray*}
	Multiplying both sides by $4L$ we get the result.
\end{proof}

Clearly, this lemma and Theorem~\ref{thm:main_result} imply the following result.
\begin{corollary}\label{cor:const_loop}
	Let the assumptions of Lemma~\ref{lem:V_k_lemma} are satisfied. Then Assumption~\ref{ass:key_assumption} holds and, in particular, if
	\begin{eqnarray*}
		\gamma &\le& \min\left\{\frac{1}{2\left(A'+\frac{4B'C}{3\rho}\right)}, \frac{L}{F'+\frac{4B'G}{3\rho}}\right\},\\
		\gamma &\le& \min\left\{\frac{1}{4(\tau-1)\mu}, \frac{1}{2\sqrt{e(\tau-1)\left(\tF(\tau-1) + \hF + \frac{2G(\tB(\tau-1)+\hB)}{\rho(1-\rho)}\right)}}\right\},\\
		\gamma &\le& \frac{1}{4\sqrt{2eL(\tau-1)\left(\tA(\tau-1)+\hA +\frac{2C(\tB(\tau-1)+\hB)}{\rho(1-\rho)}\right)}},
	\end{eqnarray*}
	then for all $K\ge 0$ we have
	\begin{eqnarray}
		\EE\left[f(\overline{x}^K) - f(x^*)\right] &\le& \frac{2\|x^0 - x^*\|^2 +   \frac{8B'}{3\rho}\gamma^2 \EE\sigma_0^2 + 4LH\gamma\EE\sigma_0^2}{\gamma W_K}\notag\\
		&&\quad + 2\gamma\left(D_1' + \frac{4B'D_2}{3\rho} + 2L\gamma D_3\right), \label{eq:main_result_const_loop}
	\end{eqnarray}
	where $\overline{x}^K \eqdef \frac{1}{W_K}\sum_{k=0}^K w_k x^k$ and
	\begin{eqnarray*}
	    H &=& \frac{4e(\tau-1)(\tB(\tau-1)+\hB)(2+\rho)\gamma^2}{\rho},\\ D_3 &=& 2e(\tau-1)\left(\tD_{1}(\tau-1)+\hD_{1} + \frac{2D_2(\tB(\tau-1)+\hB)}{\rho}\right).
	\end{eqnarray*}
	 Moreover, if $\mu > 0$, then
	\begin{eqnarray}
		\EE\left[f(\overline{x}^K) - f(x^*)\right] &\le& \left(1 - \min\left\{\gamma\mu,\frac{\rho}{4}\right\}\right)^K\frac{2\|x^0 - x^*\|^2 +   \frac{8B'}{3\rho}\gamma^2 \EE\sigma_0^2 + 4LH\gamma\EE\sigma_0^2}{\gamma}\notag\\
		&&\quad + 2\gamma\left(D_1' + \frac{4B'D_2}{3\rho} + 2L\gamma D_3\right), \label{eq:main_result_1_const_loop}
	\end{eqnarray}
	and in the case when $\mu = 0$, we have
	\begin{eqnarray}
		\EE\left[f(\overline{x}^K) - f(x^*)\right] &\le& \frac{2\|x^0 - x^*\|^2 +   \frac{8B'}{3\rho}\gamma^2 \EE\sigma_0^2 + 4LH\gamma\EE\sigma_0^2}{\gamma K}\notag\\
		&&\quad + 2\gamma\left(D_1' + \frac{4B'D_2}{3\rho} + 2L\gamma D_3\right). \label{eq:main_result_2_const_loop}
	\end{eqnarray}
\end{corollary}


\subsubsection{$\zeta$-Heterogeneous Data}\label{sec:cont_loop_homo}
In this section we assume that $f_1, f_2, \ldots, f_n$ are $\zeta$-heterogeneous (see Definition~\ref{def:zeta_hetero}). Moreover, we additionally assume that $\EE\left[g_i^k\mid x_i^k\right] = \nabla f_i(x_i^k)$ and that the functions $f_i$ for $i\in[n]$ are  $\mu$-strongly convex,
\begin{equation}
	f_i(x) \ge f_i(y) + \langle\nabla f_i(y), x-y\rangle + \frac{\mu}{2}\|x-y\|^2\qquad \forall x,y\in\R^d \label{eq:strong_convexity}
\end{equation}
which implies (e.g., see \citep{nesterov2018lectures})
\begin{equation}
	\langle \nabla f_i(x) - \nabla f_i(y), x-y\rangle \ge \mu\|x-y\|^2\qquad \forall x,y\in\R^d. \label{eq:coercivity}
\end{equation}
\begin{lemma}\label{lem:V_k_lemma_homo}
	Let Assumption~\ref{ass:L_smoothness} be satisfied, inequalities \eqref{eq:unbiasedness}-\eqref{eq:sigma_k+1_bound} hold and\footnote{When $\rho = 1$ one can always set the parameters in such a way that $B = C = G = 0$, $D_2 = 0$. In this case we assume that $\frac{2BC}{\rho(1-\rho)} = \frac{2BG}{\rho(1-\rho)} = 0$.}
	\begin{equation*}
		\gamma \le \min\left\{\frac{1}{4(\tau-1)\mu}, \frac{1}{2\sqrt{(\tau-1)\left(F+\frac{2BG}{\rho(1-\rho)}\right)}}, \frac{1}{4\sqrt{2L(\tau-1)\left(A + \frac{2BC}{\rho(1-\rho)}\right)}}\right\}.
	\end{equation*}
	Moreover, assume that $f_1, f_2, \ldots, f_n$ are $\zeta$-heterogeneous and $\mu$-strongly convex, and\newline  $\EE\left[g_i^k\mid x_i^k\right] = \nabla f_i	(x_i^k)$ for all $i\in [n]$. Then \eqref{eq:sum_V_k_bounds} holds with
	\begin{equation}
		H = \frac{4B(\tau-1)\gamma^2(2+\rho)}{\rho},\quad D_3 = 2(\tau-1)\left(D_1 + \frac{\zeta^2}{\gamma\mu} + \frac{2BD_2}{\rho}\right).\label{eq:V_k_bound_homo}
	\end{equation}
\end{lemma}
\begin{proof}
	First of all, if $k \mod \tau = 0$, then $V_k = 0$ by definition. Otherwise, we have
	\begin{eqnarray*}
		V_k &\overset{\eqref{eq:local_sgd_def},\eqref{eq:x^k_recurrsion}}{=}& \frac{1}{n}\sum\limits_{i=1}^n \left\|x_i^{k-1} - x^{k-1} -\gamma g_i^{k-1} + \gamma g^{k-1}\right\|^2\\
		&=& \frac{1}{n}\sum\limits_{i=1}^n\|x_i^{k-1} - x^{k-1}\|^2 + \frac{2\gamma}{n}\sum\limits_{i=1}^n\left\langle x_i^{k-1} - x^{k-1}, g^{k-1} - g_i^{k-1} \right\rangle  + \frac{\gamma^2}{n}\sum\limits_{i=1}^n\|g_i^{k-1} - g^{k-1}\|^2\\
		&=& V_{k-1} + 2\gamma\left\langle\frac{1}{n}\sum\limits_{i=1}^nx_i^{k-1} - x^{k-1}, g^{k-1}\right\rangle + \frac{2\gamma}{n}\sum\limits_{i=1}^n\left\langle x^{k-1}-x_i^{k-1}, g_i^{k-1} \right\rangle\\
		&&\quad  + \frac{\gamma^2}{n}\sum\limits_{i=1}^n\|g_i^{k-1} - g^{k-1}\|^2\\
		&=& V_{k-1} + \frac{2\gamma}{n}\sum\limits_{i=1}^n\left\langle x^{k-1}-x_i^{k-1}, g_i^{k-1} \right\rangle  + \frac{\gamma^2}{n}\sum\limits_{i=1}^n\|g_i^{k-1} - g^{k-1}\|^2.
	\end{eqnarray*}
	Next, we take the conditional expectation $\EE\left[\cdot\mid x^{k-1}\right] \eqdef \EE\left[\cdot\mid x_1^{k-1},\ldots, x_n^{k-1}\right]$ on both sides of the obtained inequality and get
	\begin{eqnarray*}
		\EE\left[V_k\mid x^{k-1}\right] &=& V_{k-1} + \frac{2\gamma}{n}\sum\limits_{i=1}^n\left\langle x^{k-1} - x_i^{k-1}, \nabla f_i(x_i^{k-1}) \right\rangle + \frac{\gamma^2}{n}\sum\limits_{i=1}^n\EE\left[\|g_i^{k-1} - g^{k-1}\|^2\mid x^{k-1}\right]\\
		&\overset{\eqref{eq:variance_decomposition}}{\le}& V_{k-1} + \frac{2\gamma}{n}\sum\limits_{i=1}^n \left\langle x^{k-1} - x_i^{k-1}, \nabla f_i(x_i^{k-1}) - \nabla f_i(x^{k-1}) \right\rangle\\
		&&\quad + \frac{2\gamma}{n}\sum\limits_{i=1}^n\left\langle x^{k-1} - x_i^{k-1}, \nabla f_i(x^{k-1}) \right\rangle + \frac{\gamma^2}{n}\sum\limits_{i=1}^n\EE\left[\|g_i^{k-1}\|^2\mid x^{k-1}\right].
	\end{eqnarray*}
	Since $\frac{1}{n}\sum_{i=1}^n\langle x^{k-1}-x_i^{k-1},\nabla f(x^{k-1}) \rangle = 0$, we can continue as follows:
	\begin{eqnarray*}
		\EE\left[V_k\mid x^{k-1}\right] &\overset{\eqref{eq:coercivity}}{\le}& V_{k-1} - \frac{2\gamma\mu}{n}\sum\limits_{i=1}^n \| x^{k-1} - x_i^{k-1}\|^2 + \frac{\gamma^2}{n}\sum\limits_{i=1}^n\EE\left[\|g_i^{k-1}\|^2\mid x^{k-1}\right]\\
		&&\quad + \frac{2\gamma}{n}\sum\limits_{i=1}^n\left\langle x^{k-1} - x_i^{k-1}, \nabla f_i(x^{k-1}) - \nabla f(x^{k-1}) \right\rangle\\
		&\overset{\eqref{eq:fenchel_young}}{\le}& (1-2\gamma\mu)V_{k-1} + \frac{\gamma^2}{n}\sum\limits_{i=1}^n\EE\left[\|g_i^{k-1}\|^2\mid x^{k-1}\right]\\
		&&\quad +\frac{2\gamma}{n}\sum\limits_{i=1}^n\left(\frac{\mu}{2}\|x^{k-1}-x_i^{k-1}\|^2 + \frac{1}{2\mu}\|\nabla f_i(x^{k-1})-\nabla f(x^{k-1})\|^2\right)\\
		&\overset{\eqref{eq:bounded_data_dissimilarity}}{\le}& (1-\gamma\mu)V_{k-1} + \frac{\gamma^2}{n}\sum\limits_{i=1}^n\EE\left[\|g_i^{k-1}\|^2\mid x^{k-1}\right] + \frac{\gamma\zeta^2}{\mu}.
	\end{eqnarray*}
	Taking full expectation on both sides of previous inequality, we obtain
	\begin{eqnarray*}
		\EE V_k &\overset{\eqref{eq:tower_property}}{\le}& \EE\left[V_{k-1}\right] + \frac{\gamma^2}{n}\sum\limits_{i=1}^n\EE\left[\|g_i^{k-1}\|^2\right] + \frac{\gamma\zeta^2}{\mu}.
	\end{eqnarray*}
	Let $t$ be a non-negative integer for which $\tau t \le k < \tau(t+1)$. Using this and $V_{\tau t} = 0$, we unroll the recurrence and derive
	\begin{eqnarray*}
		\EE[V_k] &\le& \frac{\gamma^2}{n}\sum\limits_{l = \tau t}^{k-1}\sum\limits_{i=1}^n\EE\left[\|g_i^{l}\|^2\right] + \frac{\gamma\zeta^2(k-\tau t)}{\mu}\\
		&\overset{\eqref{eq:second_moment_bound}}{\le}& \gamma^2\sum\limits_{l=\tau t}^{k-1}\left(2A\EE\left[f(x^l) - f(x^*)\right] + B\EE[\sigma_l^2] + F\EE[V_l] + D_1\right) + \frac{\gamma\zeta^2(k-\tau t)}{\mu},
	\end{eqnarray*}
	whence
	\begin{eqnarray}
		\sum\limits_{j=\tau t}^k w_j \EE V_j &\le& 2A\gamma^2\sum\limits_{j=\tau t}^k\sum\limits_{l=\tau t}^{j-1}w_j\EE\left[f(x^l) - f(x^*)\right] + B\gamma^2\sum\limits_{j=\tau t}^k\sum\limits_{l=\tau t}^{j-1}w_j\EE\sigma_l^2 \notag \\
		&&\quad + F\gamma^2\sum\limits_{j=\tau t}^k\sum\limits_{l=\tau t}^{j-1} w_j\EE V_l + (\tau-1)\left(\gamma^2 D_1+\frac{\gamma\zeta^2}{\mu}\right)\sum\limits_{j=\tau t}^kw_j. \notag
	\end{eqnarray}
	If we substitute $A$ with $e(\tA(\tau-1)+\hA)$, $B$ with $e(\tB(\tau-1)+\hB)$, $F$ with $e(\tF(\tau-1)+\hF)$, and $\left(\gamma^2 D_1+\frac{\gamma\zeta^2}{\mu}\right)$ with $e\gamma^2(\tD_1(\tau-1) + \hD_1)$ in the inequality above, we will get inequality \eqref{eq:V_k_lemma_technical_1}.
	Following the same steps as in the proof of Lemma~\ref{lem:V_k_lemma}, we get
	\begin{eqnarray*}
		\sum\limits_{k=0}^K w_k \EE V_k &\le& 4(\tau-1)\gamma^2\left(A +\frac{2BC}{\rho(1-\rho)}\right)\sum\limits_{k=0}^K w_k r_k + \frac{2B\EE\sigma_0^2(2+\rho)(\tau-1)\gamma^2}{\rho}\\
		&&\quad + 2(\tau-1)\gamma^2\left(F + \frac{2BG}{\rho(1-\rho)}\right)\sum\limits_{k=0}^K w_k\EE V_k + (\tau-1)\gamma^2\left(D_1 + \frac{\zeta^2}{\gamma\mu} + \frac{2BD_2}{\rho}\right)W_K.
	\end{eqnarray*}
	Our choice of $\gamma$ implies that
	\begin{equation*}
		4(\tau-1)\gamma^2\left(A + \frac{2BC}{\rho(1-\rho)}\right) \le \frac{1}{8L}\quad\text{and}\quad 2(\tau-1)\gamma^2\left(F + \frac{2BG}{\rho(1-\rho)}\right) \le \frac{1}{2}.
	\end{equation*}
	Using these inequalities we continue our derivations
	\begin{eqnarray*}
		\frac{1}{2}\sum\limits_{k=0}^K w_k\EE V_k &\le& \frac{1}{8L}\sum\limits_{k=0}^K w_k r_k + \frac{2B\EE\sigma_0^2(2+\rho)(\tau-1)\gamma^2}{\rho}\\
		&&\quad + (\tau-1)\gamma^2\left(D_1 + \frac{\zeta^2}{\gamma\mu} + \frac{2BD_2}{\rho}\right)W_K.
	\end{eqnarray*}
	Multiplying both sides by $4L$ we get the result.
\end{proof}

Clearly, this lemma and Theorem~\ref{thm:main_result} imply the following result.

\begin{corollary}\label{cor:const_loop_homo}
	Let the assumptions of Lemma~\ref{lem:V_k_lemma_homo} be satisfied. Then Assumption~\ref{ass:key_assumption} holds and, in particular, if
	\begin{eqnarray*}
		\gamma &\le& \min\left\{\frac{1}{2(A'+CM)}, \frac{L}{F'+GM}\right\},\quad M = \frac{4B'}{3\rho},\\
		\gamma &\le& \min\left\{\frac{1}{4(\tau-1)\mu}, \frac{1}{2\sqrt{(\tau-1)\left(F+\frac{2BG}{\rho(1-\rho)}\right)}}, \frac{1}{4\sqrt{2L(\tau-1)\left(A + \frac{2BC}{\rho(1-\rho)}\right)}}\right\},
	\end{eqnarray*}
	then for all $K\ge 0$ we have
	\begin{eqnarray}
		\EE\left[f(\overline{x}^K) - f(x^*)\right] &\le& \frac{2T^0 + 4LH\gamma\EE\sigma_0^2}{\gamma W_K} + 2\gamma\left(D_1' + MD_2 + 2L\gamma D_3\right), \label{eq:main_result_const_loop_homo}
	\end{eqnarray}
	where $\overline{x}^K \eqdef \frac{1}{W_K}\sum_{k=0}^K w_k x^k$ and
	\begin{equation*}
		H = \frac{4B(\tau-1)\gamma^2(2+\rho)}{\rho},\quad D_3 = 2(\tau-1)\left(D_1 + \frac{\zeta^2}{\gamma\mu} + \frac{2BD_2}{\rho}\right).
	\end{equation*}
	 Moreover, if $\mu > 0$, then
	\begin{eqnarray}
		\EE\left[f(\overline{x}^K) - f(x^*)\right] &\le& \left(1 - \min\left\{\gamma\mu,\frac{\rho}{4}\right\}\right)^K\frac{2T^0 + 4LH\gamma\EE\sigma_0^2}{\gamma}\notag\\
		&&\quad+ 2\gamma\left(D_1' + MD_2 + 2L\gamma D_3\right), \label{eq:main_result_1_const_loop_homo}
	\end{eqnarray}
	and in the case when $\mu = 0$, we have
	\begin{eqnarray}
		\EE\left[f(\overline{x}^K) - f(x^*)\right] &\le& \frac{2T^0 + 4LH\gamma\EE\sigma_0^2}{\gamma K} + 2\gamma\left(D_1' + MD_2 + 2L\gamma D_3\right). \label{eq:main_result_2_const_loop_homo}
	\end{eqnarray}
\end{corollary}

\subsection{Random Local Loop}\label{sec:random_local_loop}
In this section we show how our results can be applied to analyze \eqref{eq:local_sgd_def} in the case when
\begin{equation*}
	c_{k} = \begin{cases}1,& \text{with probability } p,\\ 0,& \text{with probability } 1-p, \end{cases}
\end{equation*}
where $p$ encodes the probability of initiating communication. This choice in effect leads to a method using a random-length local loop on all devices.

\subsubsection{Heterogeneous Data}\label{sec:random_local_loop_hetero}
As in Section~\ref{sec:const_loop_hetero}, our analysis of \eqref{eq:local_sgd_def} with random length of the local loop relies on Assumption~\ref{ass:hetero_second_moment}. Next lemma provides an upper bound for the weighted sum of $\EE\left[ V_k \right]$ in this case.
\begin{lemma}\label{lem:V_k_lemma_random}
	Let Assumptions~\ref{ass:quasi_strong_convexity_local},~\ref{ass:L_smoothness}~and~\ref{ass:hetero_second_moment} be satisfied and\footnote{When $\rho = 1$ one can always set the parameters in such a way that $\tB = \hB = C = G = 0$, $D_2 = 0$. In this case we assume that $\frac{2\tB C}{\rho(1-\rho)} = \frac{2\hB C}{\rho(1-\rho)} = \frac{2\tB G}{\rho(1-\rho)} = \frac{2\hB G}{\rho(1-\rho)} = 0$.}
	\begin{eqnarray*}
		\gamma &\le& \min\left\{\frac{p}{16\mu}, \frac{p}{2\sqrt{(1-p)((2+p)\tF+p\hF)}}, \frac{p\sqrt{3\rho(1-\rho)}}{8\sqrt{2G(1-p)\left((p+2)\tB + p\hB\right)}}\right\},\\
		\gamma &\le& \frac{p\sqrt{3}}{16\sqrt{2L(1-p)\left((2+p)\tA + p\hA + \frac{2C\left((p+2)\tB + p\hB\right)}{\rho(1-\rho)}\right)}}.
	\end{eqnarray*}
	Then \eqref{eq:sum_V_k_bounds} holds with
	\begin{eqnarray}
	    H &=& \frac{64(1-p)\left((p+2)\tB + p\hB\right)(2+\rho)\gamma^2}{3p^2\rho},\notag\\
	    D_3 &=& \frac{8(1-p)}{p^2}\left((p+2)\tD_{1} + p \hD_{1} + \frac{8D_2\left((p+2)\tB+p\hB\right)}{3\rho}\right).\label{eq:V_k_bound_random}
	\end{eqnarray}
\end{lemma}
\begin{proof}
	First of all, we introduce new notation: $\EE[\cdot\mid x^{k},g^{k}]\eqdef \EE[\cdot\mid x_1^{k},\ldots,x_n^{k},g_1^{k},\ldots,g_n^{k}]$, $\EE[\cdot\mid x^{k}]\eqdef \EE[\cdot\mid x_1^{k},\ldots,x_n^{k}]$. By definition of $V_k$, we have
	\begin{eqnarray*}
		\EE\left[V_{k+1}\mid x^k\right] &\overset{\eqref{eq:tower_property}}{=}& \frac{1}{n}\sum\limits_{i=1}^n\EE\left[\EE\left[\|x_i^{k+1} - x^{k+1}\|^2\mid x^k, g^k\right]\mid x^k\right]\\
		&=& \frac{1-p}{n}\sum\limits_{i=1}^n\EE\left[\|x_i^k - x^k - \gamma g_i^k + \gamma g^k\|^2\mid x^k\right]\\
		&\overset{\eqref{eq:variance_decomposition}}{=}& \frac{1-p}{n}\sum\limits_{i=1}^n\|x_i^k - x^k - \gamma \bar{g}_i^k + \gamma \bar{g}^k\|^2 \\
		&&\quad + \frac{(1-p)\gamma^2}{n}\sum\limits_{i=1}^n\EE\left[\|g_i^k - \bar{g}_i^k - (g^k - \bar{g}^k)\|^2\mid x^k\right]\\
		&\overset{\eqref{eq:a+b_norm_beta},\eqref{eq:variance_decomposition}}{\le}& \frac{(1-p)\left(1+\frac{p}{2}\right)}{n}\sum\limits_{i=1}^n\|x_i^k - x^k\|^2 + \frac{(1-p)\left(1+\frac{2}{p}\right)\gamma^2}{n}\sum\limits_{i=1}^n\|\bar{g}_i^k - \bar{g}^k\|^2\\
		&&\quad + \frac{(1-p)\gamma^2}{n}\sum\limits_{i=1}^n\EE\left[\|g_i^k - \bar{g}_i^k\|^2\mid x^k\right]\\
		&\overset{\eqref{eq:1+p/2_inequality},\eqref{eq:variance_decomposition}}{\le}& \left(1 - \frac{p}{2}\right)V_k + \frac{(1-p)(2+p)\gamma^2}{pn}\sum\limits_{i=1}^n\left\|\bar{g}_i^k\right\|^2\\
		&&\quad+ \frac{(1-p)\gamma^2}{n}\sum\limits_{i=1}^n\EE\left[\|g_i^k - \bar{g}_i^k\|^2\mid x^k\right],
	\end{eqnarray*}
	where $\bar{g}^k = \EE[g^k\mid x^k]$. Taking the full expectation we derive
	\begin{eqnarray*}
		\EE\left[V_{k+1}\right] &\le& \left(1 - \frac{p}{2}\right)\EE\left[V_k\right] + \frac{(1-p)(2+p)\gamma^2}{pn}\sum\limits_{i=1}^n\EE\left[\left\|\bar{g}_i^k\right\|^2\right]\\
		&&\quad + \frac{(1-p)\gamma^2}{n}\sum\limits_{i=1}^n\EE\left[\|g_i^k - \bar{g}_i^k\|^2\right]\\
		&\overset{\eqref{eq:hetero_second_moment_bound},\eqref{eq:hetero_var_bound}}{\le}& \left(1 - \frac{p}{2}\right)\EE\left[V_k\right] + 2(1-p)\gamma^2\left(\frac{2+p}{p}\tA + \hA\right)\EE\left[f(x^k)-f(x^*)\right]\\
		&&\quad + (1-p)\gamma^2\left(\left(\frac{2+p}{p}\tB + \hB\right)\EE\sigma_k^2+\left(\frac{2+p}{p}\tF + \hF\right)\EE V_k\right)\\
		&&\quad + (1-p)\gamma^2\left(\frac{2+p}{p}\tD_{1} + \hD_{1}\right).
	\end{eqnarray*}
	This inequality together with $\gamma \le \frac{p}{2\sqrt{(1-p)((2+p)\tF+p\hF)}}$  imply
	\begin{eqnarray*}
		\EE\left[V_{k+1}\right] &\le& \left(1 - \frac{p}{4}\right)\EE\left[V_k\right] + 2(1-p)\gamma^2\left(\frac{2+p}{p}\tA + \hA\right)\EE\left[f(x^k)-f(x^*)\right]\\
		&&\quad + (1-p)\gamma^2\left(\frac{2+p}{p}\tB + \hB\right)\EE\sigma_k^2 + (1-p)\gamma^2\left(\frac{2+p}{p}\tD_{1} + \hD_{1}\right).
	\end{eqnarray*}
	Unrolling the recurrence, we obtain
	\begin{eqnarray*}
		\EE\left[V_{k+1}\right] &\le& 2(1-p)\gamma^2\left(\frac{2+p}{p}\tA + \hA\right)\sum\limits_{l=0}^k \left(1 - \frac{p}{4}\right)^{k-l}\EE\left[f(x^l) - f(x^*)\right]\\
		&&\quad + (1-p)\gamma^2\left(\frac{2+p}{p}\tB + \hB\right)\sum\limits_{l=0}^k \left(1 - \frac{p}{4}\right)^{k-l}\EE\sigma_l^2\\
		&&\quad + (1-p)\gamma^2\left(\frac{2+p}{p}\tD_{1} + \hD_{1}\right)\sum\limits_{l=0}^k\left(1-\frac{p}{4}\right)^{k-l}.
	\end{eqnarray*}
	As a consequence, we derive
	\begin{eqnarray}
		\sum\limits_{k=0}^K w_k\EE\left[V_k\right] &\le& \frac{2(1-p)\left((2+p)\tA + p\hA\right)\gamma^2}{p\left(1-\frac{p}{4}\right)}\sum\limits_{k=0}^K\sum\limits_{l=0}^k\left(1 - \frac{p}{4}\right)^{k-l}w_kr_l\notag\\
		&&\quad + \frac{(1-p)\left((2+p)\tB + p\hB\right)\gamma^2}{p\left(1-\frac{p}{4}\right)}\sum\limits_{k=0}^K\sum\limits_{l=0}^k\left(1 - \frac{p}{4}\right)^{k-l}w_k\EE\left[\sigma_l^2\right]\notag\\
		&&\quad + \frac{(1-p)\left((2+p)\tD_{1} + p\hD_{1}\right)\gamma^2}{p}\sum\limits_{k=0}^K\sum\limits_{l=0}^{k-1}\left(1 - \frac{p}{4}\right)^{k-1-l}w_k,\label{eq:V_k_bound_rand_tech_1}
	\end{eqnarray}
	where we use new notation: $r_l = \EE\left[f(x^l) - f(x^*)\right]$. Recall that $w_k = (1 - \eta)^{-(k+1)}$ and $\eta = \min\left\{\gamma\mu, \frac{\rho}{4}\right\}$. Together with our assumption on $\gamma$ it implies that for all $0 \le i < k$ we have
	\begin{eqnarray}
		w_k &=& (1 - \eta)^{-(k-i+1)}\left(1 - \eta\right)^{-i} \overset{\eqref{eq:1-p/2_inequality}}{\le} w_{k-i}\left(1 + 2\eta\right)^{i} \notag\\
		&\le& w_{k-i}\left(1 + 2\gamma\mu\right)^{i} \le w_{k-i}\left(1+\frac{p}{8}\right)^i, \label{eq:V_k_bound_rand_tech_2}\\
		w_k &=& \left(1 - \eta\right)^{-(k-i+1)}\left(1 - \eta\right)^{-i} \overset{\eqref{eq:1-p/2_inequality}}{\le} w_{k-i}\left(1 + 2\eta\right)^i \le w_{k-i}\left(1 + \frac{\rho}{2}\right)^i , \label{eq:V_k_bound_rand_tech_3}\\
		w_k &\overset{\eqref{eq:1-p/2_inequality}}{\le}& \left(1 + 2\eta\right)^{k+1} \le \left(1 + \frac{\rho}{2}\right)^{k+1}. \label{eq:V_k_bound_rand_tech_4}
	\end{eqnarray}
	Having these inequalities in hand we obtain
	\begin{eqnarray*}
		\sum\limits_{k=0}^K\sum\limits_{l=0}^k\left(1 - \frac{p}{4}\right)^{k-l}w_kr_l &\overset{\eqref{eq:V_k_bound_rand_tech_2}}{\le}& \sum\limits_{k=0}^K\sum\limits_{l=0}^k\left(1 - \frac{p}{4}\right)^{k-l}\left(1 + \frac{p}{8}\right)^{k-l}w_lr_l\\
		&\overset{\eqref{eq:1+p/2_inequality}}{\le}& \sum\limits_{k=0}^K\sum\limits_{l=0}^k\left(1 - \frac{p}{8}\right)^{k-l}w_lr_l \le \left(\sum\limits_{k=0}^K w_k r_k\right)\left(\sum\limits_{k=0}^{\infty}\left(1 - \frac{p}{8}\right)^{k}\right)\\
		&=&\frac{8}{p}\sum\limits_{k=0}^K w_k r_k,
	\end{eqnarray*}
	\begin{eqnarray*}
		\sum\limits_{k=0}^K\sum\limits_{l=0}^k\left(1 - \frac{p}{4}\right)^{k-l}w_k\EE\left[\sigma_l^2\right] &\overset{\eqref{eq:V_k_bound_rand_tech_2}}{\le}& \sum\limits_{k=0}^K\sum\limits_{l=0}^k\left(1 - \frac{p}{4}\right)^{k-l}\left(1 + \frac{p}{8}\right)^{k-l}w_l\EE\left[\sigma_l^2\right]\\
		&\overset{\eqref{eq:1+p/2_inequality}}{\le}& \sum\limits_{k=0}^K\sum\limits_{l=0}^k\left(1 - \frac{p}{8}\right)^{k-l}w_l\EE\left[\sigma_l^2\right]\\
		&\le& \left(\sum\limits_{k=0}^K w_k \EE\left[\sigma_k^2\right]\right)\left(\sum\limits_{k=0}^{\infty}\left(1 - \frac{p}{8}\right)^{k}\right) = \frac{8}{p}\sum\limits_{k=0}^K w_k \EE\left[\sigma_k^2\right],
	\end{eqnarray*}
	and
	\begin{eqnarray*}
		\sum\limits_{k=0}^K\sum\limits_{l=0}^{k-1}\left(1 - \frac{p}{4}\right)^{k-1-l}w_k &\le& \left(\sum\limits_{k=0}^K w_k \right)\left(\sum\limits_{k=0}^{\infty}\left(1 - \frac{p}{4}\right)^{k}\right) = \frac{4W_K}{p}.
	\end{eqnarray*}
	Plugging these inequalities together with $1-\frac{p}{4}\ge \frac{3}{4}$ in \eqref{eq:V_k_bound_rand_tech_1}, we derive
	\begin{eqnarray}
		\sum\limits_{k=0}^K w_k\EE\left[V_k\right] &\le& \frac{64(1-p)\left((2+p)\tA + p\hA\right)\gamma^2}{3p^2}\sum\limits_{k=0}^Kw_kr_k\notag\\
		&&\quad + \frac{32(1-p)\left((2+p)\tB + p\hB\right)\gamma^2}{3p^2}\sum\limits_{k=0}^Kw_k\EE\left[\sigma_k^2\right]\notag\\
		&&\quad + \frac{4(1-p)\left((2+p)\tD_{1} + p\hD_{1}\right)\gamma^2}{p^2}W_K.\label{eq:V_k_bound_rand_tech_5}
	\end{eqnarray}
	It remains to estimate the second term on the right-hand side of this inequality. We notice that an analogous term appears in the proof of Lemma~\ref{lem:V_k_lemma}. In particular, in that proof inequality \eqref{eq:sigma_k_technical_bound} was shown via inequalities \eqref{eq:sigma_k+1_bound}, \eqref{eq:V_k_bound_rand_tech_3}, \eqref{eq:V_k_bound_rand_tech_4} and \eqref{eq:1+p/2_inequality} which hold in this case too. Therefore, we get that
	\begin{eqnarray}
	\sum\limits_{k=0}^Kw_k\EE\left[\sigma_k^2\right] &\overset{\eqref{eq:sigma_k_technical_bound}}{\le}& \frac{\EE\sigma_0^2(2+\rho)}{\rho} + \frac{4C}{\rho(1-\rho)}\sum\limits_{k=0}^Kw_kr_k + \frac{2G}{\rho(1-\rho)}\sum\limits_{k=0}^Kw_k\EE V_k +  \frac{D_2 W_K}{\rho},\notag
	\end{eqnarray}
	whence
	\begin{eqnarray}
		\sum\limits_{k=0}^K w_k\EE\left[V_k\right] &\overset{\eqref{eq:V_k_bound_rand_tech_5}}{\le}& \frac{64(1-p)\gamma^2\left((2+p)\tA + p\hA + \frac{2C\left((p+2)\tB + p\hB\right)}{\rho(1-\rho)}\right)}{3p^2}\sum\limits_{k=0}^Kw_kr_k \notag\\
		&&\quad + \frac{32(1-p)\left((p+2)\tB + p\hB\right)(2+\rho)\gamma^2\EE\sigma_0^2}{3p^2\rho}\notag\\
		&&\quad + \frac{64G(1-p)\left((p+2)\tB + p\hB\right)\gamma^2}{3p^2\rho(1-\rho)}\sum\limits_{k=0}^K w_k\EE\left[V_k\right]\notag\\
		&&\quad + \frac{4(1-p)\gamma^2}{p^2}\left((p+2)\tD_{1} + p \hD_{1} + \frac{8D_2\left((p+2)\tB+p\hB\right)}{3\rho}\right)W_K.\notag
	\end{eqnarray}
	Our assumptions on $\gamma$ imply
	\begin{eqnarray*}
		\frac{64(1-p)\gamma^2\left((2+p)\tA + p\hA + \frac{2C\left((p+2)\tB + p\hB\right)}{\rho(1-\rho)}\right)}{3p^2} \le \frac{1}{8L},\quad \frac{64G(1-p)\left((p+2)\tB + p\hB\right)\gamma^2}{3p^2\rho(1-\rho)} \le \frac{1}{2}.
	\end{eqnarray*}
	Next, we introduce new notation as follows:
	\begin{eqnarray*}
	    H &=& \frac{64(1-p)\left((p+2)\tB + p\hB\right)(2+\rho)\gamma^2}{3p^2\rho},\\
	    D_3 &=& \frac{8(1-p)}{p^2}\left((p+2)\tD_{1} + p \hD_{1} + \frac{8D_2\left((p+2)\tB+p\hB\right)}{3\rho}\right).
	\end{eqnarray*}
	Putting all together, we get
	\begin{equation*}
		\frac{1}{2}\sum\limits_{k=0}^Kw_k\EE\left[V_k\right] \le \frac{1}{8L}\sum\limits_{k=0}^K w_k r_k + \frac{H}{2}\EE\sigma_0^2 + \frac{D_3}{2}\gamma^2 W_K,
	\end{equation*}
	which concludes the proof.
\end{proof}

This lemma and Theorem~\ref{thm:main_result} imply the following result.
\begin{corollary}\label{cor:rand_loop}
	Let the assumptions of Lemma~\ref{lem:V_k_lemma_random} be satisfied. Then Assumption~\ref{ass:key_assumption} holds and, in particular, if
	\begin{eqnarray*}
		\gamma &\le& \min\left\{\frac{1}{2\left(A'+\frac{4B'C}{3\rho}\right)}, \frac{L}{F'+\frac{4B'G}{3\rho}}, \frac{p}{16\mu}, \frac{p\sqrt{3\rho(1-\rho)}}{8\sqrt{2G(1-p)\left((p+2)\tB + p\hB\right)}}\right\},\\
		\gamma &\le& \min\left\{\frac{p}{2\sqrt{(1-p)((2+p)\tF+p\hF)}}, \frac{p\sqrt{3}}{16\sqrt{2L(1-p)\left((2+p)\tA + p\hA + \frac{2C\left((p+2)\tB + p\hB\right)}{\rho(1-\rho)}\right)}}\right\},
	\end{eqnarray*}
	then for all $K\ge 0$ we have
	\begin{eqnarray}
		\EE\left[f(\overline{x}^K) - f(x^*)\right] &\le& \frac{2\|x^0 - x^*\|^2 +   \frac{8B'}{3\rho}\gamma^2 \EE\sigma_0^2 + 4LH\gamma\EE\sigma_0^2}{\gamma W_K}\notag\\
		&&\quad + 2\gamma\left(D_1' + \frac{4B'D_2}{3\rho} + 2L\gamma D_3\right), \label{eq:main_result_rand_loop}
	\end{eqnarray}
	where $\overline{x}^K \eqdef \frac{1}{W_K}\sum_{k=0}^K w_k x^k$ and
	\begin{eqnarray*}
	    H &=& \frac{64(1-p)\left((p+2)\tB + p\hB\right)(2+\rho)\gamma^2}{3p^2\rho},\\
	    D_3 &=& \frac{8(1-p)}{p^2}\left((p+2)\tD_{1} + p \hD_{1} + \frac{8D_2\left((p+2)\tB+p\hB\right)}{3\rho}\right).
	\end{eqnarray*}
	 Moreover, if $\mu > 0$, then
	\begin{eqnarray}
		\EE\left[f(\overline{x}^K) - f(x^*)\right] &\le& \left(1 - \min\left\{\gamma\mu,\frac{\rho}{4}\right\}\right)^K\frac{2\|x^0 - x^*\|^2 +   \frac{8B'}{3\rho}\gamma^2 \EE\sigma_0^2 + 4LH\gamma\EE\sigma_0^2}{\gamma}\notag\\
		&&\quad + 2\gamma\left(D_1' + \frac{4B'D_2}{3\rho} + 2L\gamma D_3\right), \label{eq:main_result_1_rand_loop}
	\end{eqnarray}
	and in the case when $\mu = 0$, we have
	\begin{eqnarray}
		\EE\left[f(\overline{x}^K) - f(x^*)\right] &\le& \frac{2\|x^0 - x^*\|^2 +   \frac{8B'}{3\rho}\gamma^2 \EE\sigma_0^2 + 4LH\gamma\EE\sigma_0^2}{\gamma K}\notag\\
		&&\quad + 2\gamma\left(D_1' + \frac{4B'D_2}{3\rho} + 2L\gamma D_3\right). \label{eq:main_result_2_rand_loop}
	\end{eqnarray}
\end{corollary}

\subsubsection{$\zeta$-Heterogeneous Data}\label{sec:random_loop_homo}
In this section we assume that $f_1, f_2, \ldots, f_n$ are $\zeta$-heterogeneous (see Definition~\ref{def:zeta_hetero}). Moreover, we additionally assume that $\EE\left[g_i^k\mid x_i^k\right] = \nabla f_i(x_i^k)$ and we also assume $\mu$-strong convexity of the functions $f_i$ for $i\in[n]$.

\begin{lemma}\label{lem:V_k_lemma_random_homo}
	Let Assumption \ref{ass:L_smoothness} be satisfied, inequalities \eqref{eq:unbiasedness}-\eqref{eq:sigma_k+1_bound} hold and\footnote{When $\rho = 1$ one can always set the parameters in such a way that $B = C = G = 0$, $D_2 = 0$. In this case we assume that $\frac{2BC}{\rho(1-\rho)} = \frac{2BG}{\rho(1-\rho)} = 0$.}
	\begin{eqnarray*}
		\gamma &\le& \min\left\{\frac{p}{8\mu},\sqrt{\frac{p}{2F(1-p)}}, \sqrt{\frac{p\rho(1-\rho)}{32BG(1-p)}}, \sqrt{\frac{p}{128L(1-p)\left(A + \frac{2BC}{\rho(1-\rho)}\right)}}\right\}.
	\end{eqnarray*}
	Moreover, assume that $f_1, f_2, \ldots, f_n$ are $\zeta$-heterogeneous and $\mu$-strongly convex, and\newline $\EE\left[g_i^k\mid x_i^k\right] = \nabla f_i	(x_i^k)$ for all $i\in [n]$. Then \eqref{eq:sum_V_k_bounds} holds with
	\begin{equation}
		H = \frac{16B(1-p)(2+\rho)\gamma^2}{p\rho},\quad D_3 = \frac{4(1-p)}{p}\left(D_1 + \frac{\zeta^2}{\gamma\mu} + \frac{4BD_2}{\rho}\right).\label{eq:V_k_bound_random_homo}
	\end{equation}
\end{lemma}
\begin{proof}
	First of all, we introduce new notation: $\EE[\cdot\mid x^{k},g^{k}]\eqdef \EE[\cdot\mid x_1^{k},\ldots,x_n^{k},g_1^{k},\ldots,g_n^{k}]$. By definition of $V_k$ for all $k\ge 1$ we have
	\begin{eqnarray*}
		\EE[V_k\mid x^{k-1},g^{k-1}] &\overset{\eqref{eq:local_sgd_def},\eqref{eq:x^k_recurrsion}}{=}& \frac{1-p}{n}\sum\limits_{i=1}^n \left\|x_i^{k-1} - x^{k-1} -\gamma g_i^{k-1} + \gamma g^{k-1}\right\|^2\\
		&=& \frac{1-p}{n}\sum\limits_{i=1}^n\|x_i^{k-1} - x^{k-1}\|^2 + \frac{2\gamma(1-p)}{n}\sum\limits_{i=1}^n\left\langle x_i^{k-1} - x^{k-1}, g^{k-1} - g_i^{k-1} \right\rangle\\
		&&\quad  + \frac{\gamma^2(1-p)}{n}\sum\limits_{i=1}^n\|g_i^{k-1} - g^{k-1}\|^2\\
		&=& (1-p)V_{k-1} + 2\gamma(1-p)\left\langle\frac{1}{n}\sum\limits_{i=1}^nx_i^{k-1} - x^{k-1}, g^{k-1}\right\rangle\\
		&&\quad + \frac{2\gamma(1-p)}{n}\sum\limits_{i=1}^n\left\langle x^{k-1}-x_i^{k-1}, g_i^{k-1} \right\rangle + \frac{\gamma^2(1-p)}{n}\sum\limits_{i=1}^n\|g_i^{k-1} - g^{k-1}\|^2\\
		&=& (1-p)V_{k-1} + \frac{2\gamma(1-p)}{n}\sum\limits_{i=1}^n\left\langle x^{k-1}-x_i^{k-1}, g_i^{k-1} \right\rangle \\
		&&\quad + \frac{\gamma^2(1-p)}{n}\sum\limits_{i=1}^n\|g_i^{k-1} - g^{k-1}\|^2.
	\end{eqnarray*}
	Next, we take the conditional expectation $\EE\left[\cdot\mid x^{k-1}\right] \eqdef \EE\left[\cdot\mid x_1^{k-1},\ldots, x_n^{k-1}\right]$ on both sides of the obtained inequality and get
	\begin{eqnarray*}
		\EE\left[V_k\mid x^{k-1}\right] &=& (1-p)V_{k-1} + \frac{2\gamma(1-p)}{n}\sum\limits_{i=1}^n\left\langle x^{k-1} - x_i^{k-1}, \nabla f_i(x_i^{k-1}) \right\rangle\\
		&&\quad + \frac{\gamma^2(1-p)}{n}\sum\limits_{i=1}^n\EE\left[\|g_i^{k-1} - g^{k-1}\|^2\mid x^{k-1}\right]\\
		&\overset{\eqref{eq:variance_decomposition}}{\le}& (1-p)V_{k-1} + \frac{2\gamma(1-p)}{n}\sum\limits_{i=1}^n\left\langle x^{k-1} - x_i^{k-1}, \nabla f_i(x_i^{k-1}) - \nabla f_i(x^{k-1}) \right\rangle\\
		&&\quad + \frac{2\gamma(1-p)}{n}\sum\limits_{i=1}^n\left\langle x^{k-1} - x_i^{k-1}, \nabla f_i(x^{k-1}) \right\rangle\\
		&&\quad + \frac{\gamma^2(1-p)}{n}\sum\limits_{i=1}^n\EE\left[\|g_i^{k-1}\|^2\mid x^{k-1}\right].
	\end{eqnarray*}
	Since $\frac{1}{n}\sum_{i=1}^n\langle x^{k-1}-x_i^{k-1},\nabla f(x^{k-1}) \rangle = 0$, we can continue as follows:
	\begin{eqnarray*}
		\EE\left[V_k\mid x^{k-1}\right] &\overset{\eqref{eq:coercivity}}{\le}& (1-p)V_{k-1} - \frac{2\gamma\mu(1-p)}{n}\sum\limits_{i=1}^n\|x^{k-1} - x_i^{k-1}\|^2\\
		&&\quad + \frac{2\gamma(1-p)}{n}\sum\limits_{i=1}^n\left\langle x^{k-1} - x_i^{k-1}, \nabla f_i(x^{k-1}) - \nabla f(x^{k-1}) \right\rangle\\
		&&\quad + \frac{\gamma^2(1-p)}{n}\sum\limits_{i=1}^n\EE\left[\|g_i^{k-1}\|^2\mid x^{k-1}\right]\\
		&\overset{\eqref{eq:fenchel_young}}{\le}& (1-p)(1-2\gamma\mu)V_{k-1} + \frac{\gamma^2(1-p)}{n}\sum\limits_{i=1}^n\EE\left[\|g_i^{k-1}\|^2\mid x^{k-1}\right]\\
		&&\quad + \frac{2\gamma(1-p)}{n}\sum\limits_{i=1}^n\left(\frac{\mu}{2}\|x^{k-1}-x_i^{k-1}\|^2 + \frac{1}{2\mu}\|\nabla f_i(x^{k-1}) - \nabla f(x^{k-1})\|^2\right)\\
		&\overset{\eqref{eq:bounded_data_dissimilarity}}{\le}& (1-p)(1-\gamma\mu)V_{k-1} + \frac{\gamma^2(1-p)}{n}\sum\limits_{i=1}^n\EE\left[\|g_i^{k-1}\|^2\mid x^{k-1}\right] + \frac{(1-p)\gamma\zeta^2}{\mu}.
	\end{eqnarray*}
	Taking full mathematical expectation on both sides of previous inequality and using $1-\gamma\mu \le 1$ we obtain
	\begin{eqnarray*}
		\EE V_k &\overset{\eqref{eq:tower_property}}{\le}& (1-p)\EE\left[V_{k-1}\right] + \frac{\gamma^2(1-p)}{n}\sum\limits_{i=1}^n\EE\left[\|g_i^{k-1}\|^2\right] + \frac{(1-p)\gamma\zeta^2}{\mu}\\
		&\overset{\eqref{eq:second_moment_bound}}{\le}&(1-p)\EE[V_{k-1}] + (1-p)\gamma^2\left(2A\EE[f(x^{k-1})-f(x^*)] + B\EE[\sigma_k^2] + F\EE[V_{k-1}] + D_1\right)\\
		&&\quad + \frac{(1-p)\gamma\zeta^2}{\mu}.
	\end{eqnarray*}
	Since $\gamma \le \sqrt{\frac{p}{2F(1-p)}}$ we have $(1-p)\gamma^2 F \le \frac{p}{2}$ and
	\begin{eqnarray*}
		\EE V_k &\le& \left(1-\frac{p}{2}\right)\EE[V_{k-1}] + (1-p)\gamma^2\left(2A\EE[f(x^{k-1})-f(x^*)] + B\EE[\sigma_k^2] + D_1 + \frac{\zeta^2}{\gamma\mu}\right).
	\end{eqnarray*}
	Unrolling the recurrence we obtain
	\begin{eqnarray*}
		\EE\left[V_{k}\right] &\le& (1-p)\gamma^2\sum\limits_{l=0}^{k-1} \left(1 - \frac{p}{2}\right)^{k-1-l}\left(2A\EE\left[f(x^l) - f(x^*)\right] + B\EE\left[\sigma_l^2\right] + D_1 + \frac{\zeta^2}{\gamma\mu}\right).
	\end{eqnarray*}
	As a consequence, we derive
	\begin{eqnarray}
		\sum\limits_{k=0}^K w_k\EE\left[V_k\right] &\le& \frac{2A(1-p)\gamma^2}{1-\frac{p}{2}}\sum\limits_{k=0}^K\sum\limits_{l=0}^k\left(1 - \frac{p}{2}\right)^{k-l}w_kr_l\notag\\
		&&\quad + \frac{B(1-p)\gamma^2}{1-\frac{p}{2}}\sum\limits_{k=0}^K\sum\limits_{l=0}^k\left(1 - \frac{p}{2}\right)^{k-l}w_k\EE\left[\sigma_l^2\right]\notag\\
		&&\quad + \left(D_1 + \frac{\zeta^2}{\gamma\mu}\right)(1-p)\gamma^2\sum\limits_{k=0}^{K}\sum\limits_{l=0}^{k-1}\left(1 - \frac{p}{2}\right)^{k-1-l}w_k,\label{eq:V_k_bound_rand_tech_1_homo}
	\end{eqnarray}
	where we use new notation: $r_l = \EE\left[f(x^l) - f(x^*)\right]$. Recall that $w_k = (1 - \eta)^{-(k+1)}$ and $\eta = \min\left\{\gamma\mu, \frac{\rho}{4}\right\}$. Together with our assumption on $\gamma$ it implies that for all $0 \le i < k$ we have
	\begin{eqnarray}
		w_k &=& (1 - \eta)^{-(k-i+1)}\left(1 - \eta\right)^{-i} \overset{\eqref{eq:1-p/2_inequality}}{\le} w_{k-i}\left(1 + 2\eta\right)^{i} \notag\\
		&\le& w_{k-i}\left(1 + 2\gamma\mu\right)^{i} \le w_{k-i}\left(1+\frac{p}{4}\right)^i, \label{eq:V_k_bound_rand_tech_2_homo}\\
		w_k &=& \left(1 - \eta\right)^{-(k-i+1)}\left(1 - \eta\right)^{-i} \overset{\eqref{eq:1-p/2_inequality}}{\le} w_{k-i}\left(1 + 2\eta\right)^i \le w_{k-i}\left(1 + \frac{\rho}{2}\right)^i , \label{eq:V_k_bound_rand_tech_3_homo}\\
		w_k &\overset{\eqref{eq:1-p/2_inequality}}{\le}& \left(1 + 2\eta\right)^{k+1} \le \left(1 + \frac{\rho}{2}\right)^{k+1}. \label{eq:V_k_bound_rand_tech_4_homo}
	\end{eqnarray}
	Having these inequalities in hand we obtain
	\begin{eqnarray*}
		\sum\limits_{k=0}^K\sum\limits_{l=0}^k\left(1 - \frac{p}{2}\right)^{k-l}w_kr_l &\overset{\eqref{eq:V_k_bound_rand_tech_2_homo}}{\le}& \sum\limits_{k=0}^K\sum\limits_{l=0}^k\left(1 - \frac{p}{2}\right)^{k-l}\left(1 + \frac{p}{4}\right)^{k-l}w_lr_l\\
		&\overset{\eqref{eq:1+p/2_inequality}}{\le}& \sum\limits_{k=0}^K\sum\limits_{l=0}^k\left(1 - \frac{p}{4}\right)^{k-l}w_lr_l \le \left(\sum\limits_{k=0}^K w_k r_k\right)\left(\sum\limits_{k=0}^{\infty}\left(1 - \frac{p}{4}\right)^{k}\right)\\
		&=&\frac{4}{p}\sum\limits_{k=0}^K w_k r_k,
	\end{eqnarray*}
	\begin{eqnarray*}
		\sum\limits_{k=0}^K\sum\limits_{l=0}^k\left(1 - \frac{p}{2}\right)^{k-l}w_k\EE\left[\sigma_l^2\right] &\overset{\eqref{eq:V_k_bound_rand_tech_2_homo}}{\le}& \sum\limits_{k=0}^K\sum\limits_{l=0}^k\left(1 - \frac{p}{2}\right)^{k-l}\left(1 + \frac{p}{4}\right)^{k-l}w_l\EE\left[\sigma_l^2\right]\\
		&\overset{\eqref{eq:1+p/2_inequality}}{\le}& \sum\limits_{k=0}^K\sum\limits_{l=0}^k\left(1 - \frac{p}{4}\right)^{k-l}w_l\EE\left[\sigma_l^2\right]\\
		&\le& \left(\sum\limits_{k=0}^K w_k \EE\left[\sigma_k^2\right]\right)\left(\sum\limits_{k=0}^{\infty}\left(1 - \frac{p}{4}\right)^{k}\right) = \frac{4}{p}\sum\limits_{k=0}^K w_k \EE\left[\sigma_k^2\right],
	\end{eqnarray*}
	and
	\begin{eqnarray*}
		\sum\limits_{k=0}^K\sum\limits_{l=0}^{k-1}\left(1 - \frac{p}{2}\right)^{k-1-l}w_k &\le& \left(\sum\limits_{k=0}^K w_k \right)\left(\sum\limits_{k=0}^{\infty}\left(1 - \frac{p}{2}\right)^{k}\right) = \frac{2W_K}{p}.
	\end{eqnarray*}
	Plugging these inequalities together with $1-\frac{p}{2}\ge \frac{1}{2}$ in \eqref{eq:V_k_bound_rand_tech_1_homo} we derive
	\begin{eqnarray}
		\sum\limits_{k=0}^K w_k\EE\left[V_k\right] &\le& \frac{16A(1-p)\gamma^2}{p}\sum\limits_{k=0}^Kw_kr_k + \frac{8B(1-p)\gamma^2}{p}\sum\limits_{k=0}^Kw_k\EE\left[\sigma_k^2\right]\notag\\
		&&\quad + \frac{2\left(D_1 + \frac{\zeta^2}{\gamma\mu}\right)(1-p)\gamma^2}{p}W_K.\label{eq:V_k_bound_rand_tech_5_homo}
	\end{eqnarray}
	It remains to estimate the second term in the right-hand side of this inequality. We notice that an analogous term appear in the proof of Lemma~\ref{lem:V_k_lemma}. In particular, in that proof inequality \eqref{eq:sigma_k_technical_bound} was shown via inequalities \eqref{eq:sigma_k+1_bound}, \eqref{eq:V_k_bound_rand_tech_3}, \eqref{eq:V_k_bound_rand_tech_4} and \eqref{eq:1+p/2_inequality} which hold in this case too. Therefore, we get that
	\begin{eqnarray}
	\sum\limits_{k=0}^Kw_k\EE\left[\sigma_k^2\right] &\overset{\eqref{eq:sigma_k_technical_bound}}{\le}& \frac{\EE\sigma_0^2(2+\rho)}{\rho} + \frac{4C}{\rho(1-\rho)}\sum\limits_{k=0}^Kw_kr_k + \frac{2G}{\rho(1-\rho)}\sum\limits_{k=0}^Kw_k\EE V_k +  \frac{D_2 W_K}{\rho},\notag
	\end{eqnarray}
	hence
	\begin{eqnarray}
		\sum\limits_{k=0}^K w_k\EE\left[V_k\right] &\overset{\eqref{eq:V_k_bound_rand_tech_5}}{\le}& \frac{16(1-p)\gamma^2\left(A + \frac{2BC}{\rho(1-\rho)}\right)}{p}\sum\limits_{k=0}^Kw_kr_k \notag\\
		&&\quad + \frac{8B(1-p)(2+\rho)\gamma^2\EE\sigma_0^2}{p\rho}+ \frac{16BG(1-p)\gamma^2}{p\rho(1-\rho)}\sum\limits_{k=0}^K w_k\EE\left[V_k\right]\notag\\
		&&\quad + \frac{2(1-p)\gamma^2}{p}\left(D_1 + \frac{\zeta^2}{\gamma\mu} + \frac{4BD_2}{\rho}\right)W_K.\notag
	\end{eqnarray}
	Our assumption on $\gamma$ imply
	\begin{eqnarray*}
		\frac{16(1-p)\gamma^2\left(A + \frac{2BC}{\rho(1-\rho)}\right)}{p} \le \frac{1}{8L},\quad \frac{16BG(1-p)\gamma^2}{p\rho(1-\rho)} \le \frac{1}{2}.
	\end{eqnarray*}
	Next, we introduce new notation as follows:
	\begin{equation*}
		H = \frac{16B(1-p)(2+\rho)\gamma^2}{p\rho},\quad D_3 = \frac{4(1-p)}{p}\left(D_1 + \frac{\zeta^2}{\gamma\mu} + \frac{4BD_2}{\rho}\right).
	\end{equation*}
	Putting all together we get
	\begin{equation*}
		\frac{1}{2}\sum\limits_{k=0}^Kw_k\EE\left[V_k\right] \le \frac{1}{8L}\sum\limits_{k=0}^K w_k r_k + \frac{H}{2}\EE\sigma_0^2 + \frac{D_3}{2}\gamma^2 W_K
	\end{equation*}
	which concludes the proof.
\end{proof}

This lemma and Theorem~\ref{thm:main_result} imply the following result.
\begin{corollary}\label{cor:rand_loop_homo}
	Let the assumptions of Lemma~\ref{lem:V_k_lemma_random_homo} are satisfied. Then Assumption~\ref{ass:key_assumption} holds and, in particular, if
	\begin{eqnarray*}
		\gamma &\le& \min\left\{\frac{1}{2(A'+CM)}, \frac{L}{F'+GM}, \frac{p}{8\mu}\right\},\quad M = \frac{4B'}{3\rho},\\
		\gamma &\le& \min\left\{\sqrt{\frac{p}{2F(1-p)}}, \sqrt{\frac{p\rho(1-\rho)}{32BG(1-p)}}, \sqrt{\frac{p}{128L(1-p)\left(A + \frac{2BC}{\rho(1-\rho)}\right)}}\right\},
	\end{eqnarray*}
	then for all $K\ge 0$ we have
	\begin{eqnarray}
		\EE\left[f(\overline{x}^K) - f(x^*)\right] &\le& \frac{2T^0 + 4LH\gamma\EE\sigma_0^2}{\gamma W_K} + 2\gamma\left(D_1' + MD_2 + 2L\gamma D_3\right), \label{eq:main_result_rand_loop_homo}
	\end{eqnarray}
	where $\overline{x}^K \eqdef \frac{1}{W_K}\sum_{k=0}^K w_k x^k$ and
	\begin{equation*}
		H = \frac{16B(1-p)(2+\rho)\gamma^2}{p\rho},\quad D_3 = \frac{4(1-p)}{p}\left(D_1 + \frac{\zeta^2}{\gamma\mu} + \frac{4BD_2}{\rho}\right).
	\end{equation*}
	 Moreover, if $\mu > 0$, then
	\begin{eqnarray}
		\EE\left[f(\overline{x}^K) - f(x^*)\right] &\le& \left(1 - \min\left\{\gamma\mu,\frac{\rho}{4}\right\}\right)^K\frac{2T^0 + 4LH\gamma\EE\sigma_0^2}{\gamma}\notag\\
		&&\quad + 2\gamma\left(D_1' + MD_2 + 2L\gamma D_3\right), \label{eq:main_result_1_rand_loop_homo}
	\end{eqnarray}
	and in the case when $\mu = 0$ we have
	\begin{eqnarray}
		\EE\left[f(\overline{x}^K) - f(x^*)\right] &\le& \frac{2T^0 + 4LH\gamma\EE\sigma_0^2}{\gamma K} + 2\gamma\left(D_1' + MD_2 + 2L\gamma D_3\right). \label{eq:main_result_2_rand_loop_homo}
	\end{eqnarray}
\end{corollary}

\section{Missing Parts from Section~\ref{sec:local_solver}}
Let us start with an useful Lemma that bounds the Bregman distance between the local iterate $x_i^k$  and the optimum $x^*$ by the Bregman distance between the virtual iterate $x^k$ and the optimum.

\begin{lemma}
Assume $f_i$ is $L$-smooth for all $i\in [n]$. Then
\begin{equation} \label{eq:poiouhnkj}
D_{f_i} (x^k_i, x^*) \leq
2D_{f_i} (x^k, x^*) +
L \| x_i^k-x^k\|^2\quad \forall i\in[n].
\end{equation}
\end{lemma}
\begin{proof}
Using corollaries of $L$-smoothness and Young's inequality, we derive
\begin{eqnarray*}
 D_{f_i} (x^k_i, x^*)
&\overset{\eqref{eq:L_smoothness_cor_1}}{\leq} &
D_{f_i} (x^k, x^*) + \langle \nabla f_i(x^k)- \nabla f_i(x^*) , x_i^k-x^k\rangle + \frac{L}{2}\| x_i^k-x^k \|^2
 \\
&\overset{\eqref{eq:fenchel_young}}{\leq}&
D_{f_i} (x^k, x^*) +
\frac{1}{2L} \|  \nabla f_i(x^k)- \nabla f_i(x^*) \|^2  + L\| x_i^k-x^k\|^2
\\
&\overset{\eqref{eq:L_smoothness}}{\leq}&
2D_{f_i} (x^k, x^*) +
L \| x_i^k-x^k\|^2.
\end{eqnarray*}
\end{proof}

\subsection{Proof of Lemma~\ref{lem:local_solver}}

 Let us bound $\frac{1}{n}\sum\limits_{i=1}^n \EE_k\left[\| g_i^k \|^2\right] $ first:
\begin{eqnarray*}
	\frac{1}{n}\sum\limits_{i=1}^n \EE_k\left[\| g_i^k \|^2\right]
	&=& 	\frac{1}{n}\sum\limits_{i=1}^n \EE_k\left[\| a_i^k - b_i^k \|^2\right] \\
	&=& 	\frac{1}{n}\sum\limits_{i=1}^n \EE_k\left[\| a_i^k - \nabla f_i(x^*)- (b_i^k - \nabla f_i(x^*))\|^2\right] \\
		&\leq& 	\frac{2}{n}\sum\limits_{i=1}^n \EE_k\left[\| a_i^k - \nabla f_i(x^*)\|^2 + \|b_i^k - \nabla f_i(x^*)\|^2\right] \\
						&\leq& 	\frac{2}{n}\sum\limits_{i=1}^n \left(2A_iD_{f_i} (x^k_i, x^*) + B_i\sigma_{i,k}^2 + D_{1,i} +\EE_k\left[\|b_i^k - \nabla f_i(x^*)\|^2\right]\right)
						\\
		&\stackrel{\eqref{eq:poiouhnkj}}{\leq}& \frac{2}{n}\sum\limits_{i=1}^n \left(4A_iD_{f_i} (x^k, x^*)
				+ 2A_iL \| x_i^k-x^k\|^2
		+ B_i\sigma_{i,k}^2\right)
		 \\
		&&\quad + \frac{2}{n}\sum\limits_{i=1}^n \left(D_{1,i} + \EE_k\left[\|b_i^k - \nabla f_i(x^*)\|^2\right]\right)
		 \\
		&\leq &
		8\max_{i} \{A_i\} (f(x^k)-f(x^*))
				+ 4\max_{i} \{A_i\}  L V_k\\
				&&\quad +
		\frac{2}{n}\sum\limits_{i=1}^n  \left( B_i\sigma_{i,k}^2 + D_{1,i} + \EE_k\left[\|b_i^k - \nabla f_i(x^*)\|^2\right] \right).
\end{eqnarray*}

Taking the full expectation, we arrive at
\begin{eqnarray}
    \frac{1}{n}\sum\limits_{i=1}^n \EE\left[\| g_i^k \|^2\right]
&\leq&
		8\max_{i} \{A_i\} \EE(f(x^k)-f(x^*))
				+ 4\max_{i} \{A_i\}  L \EE V_k\notag\\
				&&\quad +
		\frac{2}{n}\sum\limits_{i=1}^n  \left( B_i\EE\sigma_{i,k}^2 + D_{1,i} +\EE \|b_i^k - \nabla f_i(x^*)\|^2 \right).\label{eq:ndjkabhdavgjda}
\end{eqnarray}

Next, we have

\begin{eqnarray*}
	\EE_k\left[ \left \|  \frac{1}{n}\sum\limits_{i=1}^n  g_i^k \right\|^2\right]
	&=&
	\EE_k\left[\left\| \frac{1}{n}\sum\limits_{i=1}^n  a_i^k - b_i^k \right\|^2\right]
	\\
	&=&
	\EE_k\left[\left\| \frac{1}{n}\sum\limits_{i=1}^n  a_i^k - \nabla f_i(x^*)\right\|^2\right]
		\\
	&=&
	\Var\left[\frac{1}{n}\sum\limits_{i=1}^n  a_i^k - \nabla f_i(x^*)\right] + 
	\left\|\frac{1}{n}\sum\limits_{i=1}^n \nabla f_i(x_i^k)- \nabla f_i(x^*)\right\|^2
	\\
	&\leq &
		\Var\left[\frac{1}{n}\sum\limits_{i=1}^n  a_i^k - \nabla f_i(x^*)\right] + 
 \frac{1}{n}\sum\limits_{i=1}^n	\left\| \nabla f_i(x_i^k)- \nabla f_i(x^*)\right\|^2
 	\\
	&\leq &
		\Var\left[\frac{1}{n}\sum\limits_{i=1}^n  a_i^k - \nabla f_i(x^*)\right] + 
\frac{2L}{n}\sum\limits_{i=1}^nD_{f_i}(x_i^k,x^*)
\\
	&= &
	\frac{1}{n^2}	\sum\limits_{i=1}^n \Var\left[  a_i^k - \nabla f_i(x^*)\right] + \frac{2L}{n}\sum\limits_{i=1}^nD_{f_i}(x_i^k,x^*)
		\\
	&\leq &
	\frac{1}{n^2} \sum\limits_{i=1}^n   \EE_k\left[\left\| a_i^k - \nabla f_i(x^*)\right\|^2\right]
	+ 
\frac{2L}{n}\sum\limits_{i=1}^nD_{f_i}(x_i^k,x^*)
	\\
		&\le &
	\frac{1}{n^2} \sum\limits_{i=1}^n  \left(
	 2A_iD_{f_i} (x^k_i, x^*) + B_i\sigma_{i,k}^2 + D_{1,i}\right)
	 + 
\frac{2L}{n}\sum\limits_{i=1}^nD_{f_i}(x_i^k,x^*)
	 \\
	 	&\le &
\frac{1}{n^2} \sum\limits_{i=1}^n  \left(  	2\left(\max_{i} \{A_i\}+nL\right)D_{f_i} (x^k_i, x^*) +   B_i\sigma_{i,k}^2 + D_{1,i}\right)
		 \\
	 	&\stackrel{\eqref{eq:poiouhnkj}}{\leq} &
  	\left( \frac{4\max_{i} \{A_i\}}{n} + 2L\right) D_{f} (x^k, x^*)\\
  	&&\quad+ \frac{1}{n^2} \sum\limits_{i=1}^n  \left( 2(\max_{i} \{A_i\} L+nL^2) \| x_i^k - x^*\|^2+   B_i\sigma_{i,k}^2 + D_{1,i}\right)
	 \\
	 	&=&
  	\left( \frac{4\max_{i} \{A_i\}}{n} + 2L\right)\left(f(x^k)- f(x^*)\right) +2\left(\frac{\max_{i}\{A_i\} L  }{n} + L^2\right) V_k\\
  	&&\quad+ \frac{1}{n^2}\sum\limits_{i=1}^n  \left(   B_i\sigma_{i,k}^2 + D_{1,i}\right).
\end{eqnarray*}

 Further, we define
\begin{equation}\label{eq:omegak_eq_dnakjdjsgajvgdhasvgjlds}
\omega_k^2 \eqdef \frac{2}{n}\sum\limits_{i=1}^n   B_i\sigma_{i,k}^2
\end{equation}
and consequently, we get
\begin{eqnarray}
\EE\left[\omega_{k+1}^2  \right]
&=&
\frac{2}{n}\sum\limits_{i=1}^n B_i \EE\left[   \sigma_{i,k+1}^2  \right]
\nonumber
\\
&\leq&
 (1-\rho) \omega_{k}^2+ \frac{2}{n} \sum\limits_{i=1}^n B_iC_i D_{f_i}(x^k_i, x^*) +  \frac{2}{n}\sum\limits_{i=1}^n B_iD_{2,i}
\nonumber
\\
&\stackrel{\eqref{eq:poiouhnkj}}{\leq}&
 (1-\rho) \omega_{k}^2+ \frac{4}{n} \sum\limits_{i=1}^n B_iC_iD_{f_i}(x^k, x^*) + \frac{2}{n} \sum\limits_{i=1}^n B_iC_i L\|x^k_i - x^k\|^2 +\frac{2}{n} \sum\limits_{i=1}^n B_iD_{2,i}
 \nonumber
\\
&\leq &
 (1-\rho) \omega_{k}^2+ 4 \max_{i}\{B_iC_i \} D_{f}(x^k, x^*) + 2 \max_{i}\{B_iC_i \} LV_k + \frac2n\sum\limits_{i=1}^n B_iD_{2,i}.
\notag
\end{eqnarray}

We will provide a bound on $\EE\|b_i^k - \nabla f_i(x^*)\|^2 $ based on the choices of $b_i^k$:

\begin{itemize}
\item[Case I.] The choice $b_i^k = 0$ yields
$
 \EE\|b_i^k - \nabla f_i(x^*)\|^2 =  \|\nabla f_i(x^*)\|^2.
$
\item[Case II.] The choice $b_i^k = \nabla f_i(x^*)$ yields $ \EE\|b_i^k - \nabla f_i(x^*)\|^2 =  0$.
 Overall, for both Case I and II we have
\[
\EE\sigma^2_{k+1} \leq
(1-\rho)\EE\sigma^2_k
+ 4 \max_{i}\{B_iC_i \}  D_{f}(x^k, x^*) + 2 \max_{i}\{B_iC_i \}L V_k + \frac2n\sum\limits_{i=1}^n B_iD_{2,i}
\]
as desired, where $\sigma_k = \omega_k$.
\item[Case III.] The choice $b_i^k = h_i^k - \frac1n  \sum_{i=1}^n h_i^k $ yields
\[
\frac1n \sum_{i=1}^n  \|b_i^k - \nabla f_i(x^*)\|^2 = \frac1n \sum_{i=1}^n  \left\| h_i^k  - \frac1n \sum_{i=1}^n  h_i^k  - \nabla f_i(x^*)\right \|^2 \leq  \frac1n \sum_{i=1}^n  \| h_i^k   - \nabla f_i(x^*)\|^2
\]
where
\begin{eqnarray*}
\EE_k \left[ \|  h_i^{k+1}   - \nabla f_i(x^*)\|^2\right]
 &=&
 (1-\rho_i')\|  h_i^k   - \nabla f_i(x^*)\|^2  + \rho_i' \EE_k\|  l_i^k   - \nabla f_i(x^*)\|^2
 \\
  &\stackrel{\eqref{eq:bdef}}{\leq } &
   (1-\rho_i')\|h_i^k   - \nabla f_i(x^*)\|^2  +   2\rho_i' A'_iD_{f_i} (x^k_i, x^*)+ \rho_i' D_{3,i}.
\end{eqnarray*}
Next, set $\sigma^2_k \eqdef \omega^2_k +  \| h_i^k  - \nabla f_i(x^*)\|^2$ for this case.  Consequently, we have
\begin{eqnarray*}
	\EE_k\sigma^2_{k+1} &\leq& 
(1-\rho)\sigma^2_k
+ 4( \max_{i}\{B_iC_i \} + \max_i\{\rho_i' A'_i\})  D_{f}(x^k, x^*) +  2(\max_{i}\{B_iC_i \} \\
&&\quad+\max_{i}\{\rho_i'A_i' \}) LV_k + \frac1n\sum\limits_{i=1}^n \left( 2 B_iD_{2,i} + \rho_i' D_{3,i} \right),
\end{eqnarray*}
where $\rho = \min_i\min\{\rho_i,\rho_i'\}$.
\end{itemize}

 It remains to plug everything back to~\eqref{eq:second_moment_bound}, \eqref{eq:second_moment_bound_2} and \eqref{eq:sigma_k+1_bound}.

\newpage
\begin{table*}[!t]
\caption{The parameters for which the methods from Table~\ref{tbl:special_cases} satisfy Assumption~\ref{ass:key_assumption}/\ref{ass:hetero_second_moment}. Absolute  constants were omitted. The meaning of the expressions appearing in the table, as well as their justification, is detailed in Section~\ref{sec:special_cases}. UBV stands for the ``Uniform Bound on the Variance'' of local stochastic gradient, which is often assumed when $f_i$ is of the form~\eqref{eq:f_i_expectation}. ES stands for the ``Expected Smoothness'' inequality~\citep{gower2019sgd}, which does not impose any extra assumption on the objective/noise, but rather can be derived given the sampling strategy and the smoothness structure of $f_i$. Consequently, such a setup allows us to obtain local methods with importance sampling. Next, the simple setting is a special case of ES when we uniformly sample a single index on each node each iteration.}
\label{tbl:special_cases-parameters}
\begin{center}
\footnotesize
\begin{adjustbox}{angle=90}
\begin{tabular}{|c|c|c|c|c|c|c|c|}
\hline
 Method, Setting &   $A$, $\tA$, $\hA$, $A'$ & $B$, $\tB$, $\hB$, $B'$ & $\rho$ & $C$ & $F$, $\tF$, $\hF$, $F'$ & $G$ & $D_1'$, $D_1$, $\tD_1$, $\hD_1$, $D_2$, $D_3$ \\
\hline
 \begin{tabular}{c}
 	{\tt Local-SGD}\\ 
 	UBV, $\zeta$-Het.
 \end{tabular}   &  $L$, $-$, $-$, $L$ & $0$, $-$, $-$, $0$ & $1$ & $0$ & $L^2$, $-$, $-$, $L^2$ & $0$ & \makecell{$\frac{\sigma^2}{n}$, $\sigma^2+\zeta_*^2$, $-$, $-$, $0$,\\ $\tau\sigma^2+\tau^2\zeta^2$}\\
 \hline
 \begin{tabular}{c}
 	{\tt Local-SGD}\\ 
 	UBV, Het.
 \end{tabular}   &  $-$, $L$, $0$, $L$ & $-$, $0$, $0$, $0$ & $1$ & $0$ & $-$, $L^2$, $0$, $L^2$ & $0$ & \makecell{$\frac{\sigma^2}{n}$, $-$, $\zeta_*^2$, $\sigma^2$, $0$,\\ $(\tau-1)\sigma^2+(\tau-1)^2\zeta_*^2$}\\
 \hline
 \begin{tabular}{c}
 	{\tt Local-SGD}\\ 
 	ES, $\zeta$-Het.
 \end{tabular}   &  $\cL$, $-$, $-$, $\frac{\cL}{n}+L$ & $0$, $-$, $-$, $0$ & $1$ & $0$ & $\cL L$, $-$, $-$, $\frac{\cL L}{n}+L^2$ & $0$ & \makecell{$\frac{\sigma_*^2}{n}$, $\sigma_*^2 + \zeta_*^2$, $-$, $-$, $0$,\\ $(\tau-1)\left(\sigma_*^2+\zeta_*^2+\frac{\zeta^2}{\gamma\mu}\right)$}\\
 \hline
 \begin{tabular}{c}
 	{\tt Local-SGD}\\ 
 	ES, Het.
 \end{tabular}   &  $-$, $L$, $\cL$, $\frac{\cL}{n}+L$ & $-$, $0$, $0$, $0$ & $1$ & $0$ & $-$, $L^2$, $\cL L$, $\frac{\cL L}{n}+L^2$ & $0$ & \makecell{$\frac{\sigma_*^2}{n}$, $-$, $\zeta_*^2$, $\sigma_*^2$, $0$,\\ $(\tau-1)\sigma_*^2+(\tau-1)^2\zeta_*^2$}\\
 \hline
 \begin{tabular}{c}
	{\tt Local-SVRG}\\
	simple, $\zeta$-Het.  
\end{tabular}   &  \makecell{$\max L_{ij}$, $-$, $-$,\\ $\frac{\max L_{ij}}{n}+L$} & $1$, $-$, $-$, $\frac1n$ & $\psvrg$ & $\max L_{ij}\psvrg$ & \makecell{$\max L_{ij} L$, $-$, $-$,\\ $\frac{\max L_{ij} L}{n} + L^2$} & $\max L_{ij} L q$ & \makecell{$0$, $\zeta_*^2$, $-$, $-$, $0$,\\ $(\tau-1)\left(\zeta_*^2 + \frac{\zeta^2}{\gamma\mu}\right)$}\\
 \hline
 \begin{tabular}{c}
	{\tt Local-SVRG}\\
	simple, Het.  
\end{tabular}   &  \makecell{$-$, $L$, $\max L_{ij}$,\\ $\frac{\max L_{ij}}{n}+L$} & $-$, $0$, $1$, $\frac1n$ & $\psvrg$ & $\max L_{ij}\psvrg$ & \makecell{$-$, $L^2$, $\max L_{ij} L$,\\ $\frac{\max L_{ij} L}{n} + L^2$} & $\max L_{ij} L q$ & \makecell{$0$, $-$, $\zeta_*^2$, $0$, $0$, $(\tau-1)^2\zeta_*^2$}\\
 \hline
  \begin{tabular}{c}
  	{\tt S*-Local-SGD}\\
  	UBV, Het.
	\end{tabular}    &  $-$, $L$, $0$, $L$ & $-$, $0$, $0$, $0$ & $1$ & $0$ & $-$, $L^2$, $0$, $l^2$ & $0$ & $\frac{\sigma^2}{n}$, $-$, $0$, $\sigma^2$, $(\tau-1)\sigma^2$\\
 \hline
 \begin{tabular}{c}
 	{\tt SS-Local-SGD}\\
 	UBV, Het.,\\
 	$p=q$, $r=\lceil\nicefrac{1}{p} \rceil$
\end{tabular}    & $-$, $L$, $0$, $L$ & $-$, $1$, $0$, $0$ & $p$ & $Lp$ & $-$, $L^2$, $0$, $L^2$ & $0$ & \makecell{$\frac{\sigma^2}{n}$, $-$, $p\sigma^2$, $\sigma^2$, $0$, $\frac{(1-p)\sigma^2}{p}$ }  \\
 \hline
 \begin{tabular}{c}
 	{\tt SS-Local-SGD}\\
 	ES, Het.,\\
 	$p=q$, $r=\lceil\nicefrac{1}{p} \rceil$
\end{tabular}    & $-$, $L$, $\cL$, $\frac{\cL}{n}+L$ & $-$, $1$, $0$, $0$ & $p$ & $Lp+\cL p^2$ & $-$, $L^2$, $\cL L$, $\frac{\cL L}{n}+L^2$ & $0$ & \makecell{$\frac{\sigma_*^2}{n}$, $-$, $0$, $\sigma_*^2$, $p^2\sigma_*^2$, $\frac{(1-p)\sigma_*^2}{p}$ }  \\
 \hline
 \begin{tabular}{c}
	{\tt S*-Local-SGD*}\\
	simple, Het. 
\end{tabular}    &  \makecell{$-$, $L$, $\max L_{ij}$,\\ $\frac{\max L_{ij}}{n}+L$ }& $-$, $0$, $0$, $0$ & $p$ & $0$ & \makecell{$-$, $L^2$, $\max L_{ij} L$,\\ $\frac{L\max L_{ij}}{n}+L^2$} & $0$ & $0$, $-$, $0$, $0$, $0$, $0$ \\
 \hline
 \begin{tabular}{c}
	{\tt S-Local-SVRG}\\
	simple, Het.,\\
	$q = \frac{1}{m}$, $m \ge \frac{1}{p}$
\end{tabular}    &  \makecell{$-$, $L$, $\max L_{ij}$,\\ $\frac{\max L_{ij}}{n}+L$ }& $-$, $1$, $1$, $\frac{1}{n}$ & $\frac{1}{m}$ & $\frac{L+\max L_{ij}}{m}$ & \makecell{$-$, $L^2$, $\max L_{ij} L$,\\ $\frac{L\max L_{ij}}{n}+L^2$} & $0$ & $0$, $-$, $0$, $0$, $0$, $0$ \\
 \hline
\end{tabular}
\end{adjustbox}
\end{center}
\end{table*}

%% file: Appendix_marina.tex
\chapter{Appendix for Chapter~\ref{ch:marina}}\label{app:marina}

\section{Missing Proofs for \algname{MARINA}}\label{sec:marina_proofs}
\subsection{Generally Non-Convex Problems}\label{sec:proof_of_thm_non_cvx}
In this section, we provide the full statement of Theorem~\ref{thm:main_result_non_cvx} together with the proof of this result.
\begin{theorem}[Theorem~\ref{thm:main_result_non_cvx}]\label{thm:main_result_non_cvx_appendix}
	Let Assumptions~\ref{as:lower_bound}~and~\ref{as:L_smoothness} be satisfied and 
	\begin{equation}
		\gamma \le \frac{1}{L\left(1 + \sqrt{\frac{(1-p)\omega}{pn}}\right)},\label{eq:gamma_bound_non_cvx_appendix}
	\end{equation}
	where $L^2 = \frac{1}{n}\sum_{i=1}^nL_i^2$. Then after $K$ iterations of \algname{MARINA} we have
	\begin{equation}
		\EE\left[\left\|\nabla f(\hat x^K)\right\|^2\right] \le \frac{2\Delta_0}{\gamma K}, \label{eq:main_res_non_cvx_appendix}
	\end{equation}
	where $\hat{x}^K$ is chosen uniformly at random from $x^0,\ldots,x^{K-1}$ and $\Delta_0 = f(x^0)-f_*$. That is, after
	\begin{equation}
		K = \cO\left(\frac{\Delta_0 L}{\varepsilon^2}\left(1 + \sqrt{\frac{(1-p)\omega}{pn}}\right)\right) \label{eq:main_res_2_non_cvx_appendix}
	\end{equation}
	iterations \algname{MARINA} produces such a point $\hat x^K$ that $\EE[\|\nabla f(\hat x^K)\|^2] \le \varepsilon^2$.
	Moreover, under an assumption that the communication cost is proportional to the number of non-zero components of transmitted vectors from workers to the server, we have that the expected total communication cost per worker equals
	\begin{equation}
		d + K(pd + (1-p)\zeta_{\cQ}) =  \cO\left(d+\frac{\Delta_0 L}{\varepsilon^2}\left(1 + \sqrt{\frac{(1-p)\omega}{pn}}\right)(pd + (1-p)\zeta_{\cQ})\right),\label{eq:main_res_4_non_cvx_appendix}
	\end{equation}
	where $\zeta_{\cQ}$ is the expected density of the quantization (see Def.~\ref{def:quantization}).
\end{theorem}
\begin{proof}[Proof of Theorem~\ref{thm:main_result_non_cvx}]
	The scheme of the proof is similar to the proof of Theorem~1 from \cite{li2020page}. From Lemma~\ref{lem:lemma_2_page}, we have
	\begin{equation}
		\EE[f(x^{k+1})] \le \EE[f(x^k)] - \frac{\gamma}{2}\EE\left[\|\nabla f(x^k)\|^2\right] - \left(\frac{1}{2\gamma} - \frac{L}{2}\right)\EE\left[\|x^{k+1}-x^k\|^2\right] + \frac{\gamma}{2}\EE\left[\|g^k - \nabla f(x^k)\|^2\right]. \label{eq:non_cvx_technical_1}
	\end{equation}
	Next, we need to derive an upper bound for $\EE\left[\|g^{k+1}-\nabla f(x^{k+1})\|^2\right]$. By definition of $g^{k+1}$, we have
	\begin{equation}
		g^{k+1} = \begin{cases}\nabla f(x^{k+1})& \text{with probability } p,\\ g^k + \frac{1}{n}\sum\limits_{i = 1}^n\cQ\left(\nabla f_{i}(x^{k+1}) - \nabla f_{i}(x^k)\right)& \text{with probability } 1-p. \end{cases}\notag
	\end{equation}
	Using this, variance decomposition \eqref{eq:variance_decomposition} and tower property \eqref{eq:tower_property}, we derive:
	\begin{eqnarray}
		\EE\left[\|g^{k+1}-\nabla f(x^{k+1})\|^2\right] &\overset{\eqref{eq:tower_property}}{=}& (1-p)\EE\left[\left\|g^k + \frac{1}{n}\sum\limits_{i=1}^n \cQ\left(\nabla f_{i}(x^{k+1}) - \nabla f_{i}(x^k)\right) - \nabla f(x^{k+1})\right\|^2\right]\notag\\
		&\overset{\eqref{eq:tower_property},\eqref{eq:variance_decomposition}}{=}& (1-p)\EE\left[\left\|\frac{1}{n}\sum\limits_{i=1}^n \cQ\left(\nabla f_{i}(x^{k+1}) - \nabla f_{i}(x^k)\right) - \nabla f(x^{k+1}) + \nabla f(x^k)\right\|^2\right]\notag\\
		&&\quad + (1-p)\EE\left[\left\|g^k - \nabla f(x^k)\right\|^2\right].\notag
	\end{eqnarray}
	Since $\cQ\left(\nabla f_{1}(x^{k+1}) - \nabla f_{1}(x^k)\right),\ldots,\cQ\left(\nabla f_{n}(x^{k+1}) - \nabla f_{n}(x^k)\right)$ are independent random vectors for fixed $x^k$ and $x^{k+1}$ we have
	\begin{eqnarray*}
		\EE\left[\|g^{k+1}-\nabla f(x^{k+1})\|^2\right] &=& (1-p)\EE\left[\left\|\frac{1}{n}\sum\limits_{i=1}^n \left(\cQ\left(\nabla f_{i}(x^{k+1}) - \nabla f_{i}(x^k)\right) - \nabla f_i(x^{k+1}) + \nabla f_i(x^k)\right)\right\|^2\right]\\
		&&\quad + (1-p)\EE\left[\left\|g^k - \nabla f(x^k)\right\|^2\right]\\
		&=& \frac{1-p}{n^2}\sum\limits_{i=1}^n\EE\left[\left\|\cQ\left(\nabla f_{i}(x^{k+1}) - \nabla f_{i}(x^k)\right) - \nabla f_i(x^{k+1}) + \nabla f_i(x^k)\right\|^2\right]\\
		&&\quad + (1-p)\EE\left[\left\|g^k - \nabla f(x^k)\right\|^2\right]\\
		&\overset{\eqref{eq:quantization_def}}{\le}& \frac{(1-p)\omega}{n^2}\sum\limits_{i=1}^n\EE\left[\left\|\nabla f_i(x^{k+1}) - \nabla f_i(x^k)\right\|^2\right] + (1-p)\EE\left[\left\|g^k - \nabla f(x^k)\right\|^2\right].
	\end{eqnarray*}
	Using $L$-smoothness \eqref{eq:L_smoothness_local_marina} of $f_i$ together with the tower property \eqref{eq:tower_property}, we obtain
	\begin{eqnarray}
		\EE\left[\|g^{k+1}-\nabla f(x^{k+1})\|^2\right] &\le& \frac{(1-p)\omega}{n^2}\sum\limits_{i=1}^nL_i^2\EE\left[\|x^{k+1} - x^k\|^2\right] + (1-p)\EE\left[\left\|g^k - \nabla f(x^k)\right\|^2\right]\notag\\
		&=&\frac{(1-p)\omega L^2}{n}\EE\left[\|x^{k+1}-x^k\|^2\right] + (1-p)\EE\left[\left\|g^k - \nabla f(x^k)\right\|^2\right].\label{eq:non_cvx_technical_2}
	\end{eqnarray}
	Next, we introduce a new notation: $\Phi_k = f(x^k) - f_* + \frac{\gamma}{2p}\|g^k - \nabla f(x^k)\|^2$. Using this and inequalities \eqref{eq:non_cvx_technical_1} and \eqref{eq:non_cvx_technical_2}, we establish the following inequality:
	\begin{eqnarray}
		\EE\left[\Phi_{k+1}\right] &\le& \EE\left[f(x^k) - f_* - \frac{\gamma}{2}\|\nabla f(x^k)\|^2 - \left(\frac{1}{2\gamma} - \frac{L}{2}\right)\|x^{k+1}-x^k\|^2 + \frac{\gamma}{2}\|g^k - \nabla f(x^k)\|^2\right]\notag\\
		&&\quad + \frac{\gamma}{2p}\EE\left[\frac{(1-p)\omega L^2}{n}\|x^{k+1}-x^k\|^2 + (1-p)\left\|g^k - \nabla f(x^k)\right\|^2\right] \notag\\
		&=& \EE\left[\Phi_k\right] - \frac{\gamma}{2}\EE\left[\|\nabla f(x^k)\|^2\right] + \left(\frac{\gamma(1-p)\omega L^2}{2pn} - \frac{1}{2\gamma} + \frac{L}{2}\right)\EE\left[\|x^{k+1}-x^k\|^2\right]\notag\\
		&\overset{\eqref{eq:gamma_bound_non_cvx_appendix}}{\le}& \EE\left[\Phi_k\right] - \frac{\gamma}{2}\EE\left[\|\nabla f(x^k)\|^2\right],\label{eq:non_cvx_technical_3}
	\end{eqnarray}
	where in the last inequality, we use $\frac{\gamma(1-p)\omega L^2}{2pn} - \frac{1}{2\gamma} + \frac{L}{2} \le 0$ following from \eqref{eq:gamma_bound_non_cvx_appendix}. Summing up inequalities \eqref{eq:non_cvx_technical_3} for $k=0,1,\ldots,K-1$ and rearranging the terms, we derive
	\begin{eqnarray}
		\frac{1}{K}\sum\limits_{k=0}^{K-1}\EE\left[\|\nabla f(x^k)\|^2\right] &\le& \frac{2}{\gamma K}\sum\limits_{k=0}^{K-1}\left(\EE[\Phi_k]-\EE[\Phi_{k+1}]\right) = \frac{2\left(\EE[\Phi_0]-\EE[\Phi_{K}]\right)}{\gamma K} = \frac{2\Delta_0}{\gamma K},\notag
	\end{eqnarray}
	since $g^0 = \nabla f(x^0)$ and $\Phi_{k+1} \ge 0$. Finally, using the tower property \eqref{eq:tower_property} and the definition of $\hat x^K$, we obtain \eqref{eq:main_res_non_cvx_appendix} that implies \eqref{eq:main_res_2_non_cvx_appendix} and \eqref{eq:main_res_4_non_cvx_appendix}.
\end{proof}

\begin{corollary}[Corollary~\ref{cor:main_result_non_cvx}]\label{cor:main_result_non_cvx_appendix}
	Let the assumptions of Theorem~\ref{thm:main_result_non_cvx} hold and $p = \frac{\zeta_{\cQ}}{d}$, where $\zeta_{\cQ}$ is the expected density of the quantization (see Def.~\ref{def:quantization}). If 
	\begin{equation*}
		\gamma \le \frac{1}{L\left(1 + \sqrt{\frac{\omega}{n}\left(\frac{d}{\zeta_{\cQ}}-1\right)}\right)},
	\end{equation*}
	then \algname{MARINA} requires 
	\begin{equation*}
		K = \cO\left(\frac{\Delta_0 L}{\varepsilon^2}\left(1 + \sqrt{\frac{\omega}{n}\left(\frac{d}{\zeta_{\cQ}}-1\right)}\right)\right)
	\end{equation*}
	iterations/communication rounds to achieve $\EE[\|\nabla f(\hat x^K)\|^2] \le \varepsilon^2$, and the expected total communication cost per worker is
	\begin{equation*}
		\cO\left(d+\frac{\Delta_0 L}{\varepsilon^2}\left(\zeta_{\cQ} + \sqrt{\frac{\omega\zeta_{\cQ}}{n}\left(d-\zeta_{\cQ}\right)}\right)\right)
	\end{equation*}
	 under an assumption that the communication cost is proportional to the number of non-zero components of transmitted vectors from workers to the server.
\end{corollary}
\begin{proof}[Proof of Corollary~\ref{cor:main_result_non_cvx}]
	The choice of $p = \frac{\zeta_{\cQ}}{d}$ implies
	\begin{eqnarray*}
		\frac{1-p}{p} &=& \frac{d}{\zeta_{\cQ}}-1,\\
		pd + (1-p)\zeta_{\cQ} &\le& \zeta_{\cQ} + \left(1 - \frac{\zeta_{\cQ}}{d}\right)\cdot\zeta_{\cQ} \le 2\zeta_{\cQ}.
	\end{eqnarray*}	
	Plugging these relations in \eqref{eq:gamma_bound_non_cvx_appendix}, \eqref{eq:main_res_2_non_cvx_appendix}, and \eqref{eq:main_res_4_non_cvx_appendix}, we get that if
	\begin{equation*}
		\gamma \le \frac{1}{L\left(1 + \sqrt{\frac{\omega}{n}\left(\frac{d}{\zeta_{\cQ}}-1\right)}\right)},
	\end{equation*}
	then \algname{MARINA} requires 
	\begin{eqnarray*}
		K &=& \cO\left(\frac{\Delta_0 L}{\varepsilon^2}\left(1 + \sqrt{\frac{(1-p)\omega}{pn}}\right)\right)\\
		&=& \cO\left(\frac{\Delta_0 L}{\varepsilon^2}\left(1 + \sqrt{\frac{\omega}{n}\left(\frac{d}{\zeta_{\cQ}}-1\right)}\right)\right)
	\end{eqnarray*}
	iterations/communication rounds in order to achieve $\EE[\|\nabla f(\hat x^K)\|^2] \le \varepsilon^2$, and the expected total communication cost per worker is
	\begin{eqnarray*}
		d + K(pd + (1-p)\zeta_{\cQ}) &=&  \cO\left(d+\frac{\Delta_0 L}{\varepsilon^2}\left(1 + \sqrt{\frac{(1-p)\omega}{pn}}\right)(pd + (1-p)\zeta_{\cQ})\right)\\
		&=&\cO\left(d+\frac{\Delta_0 L}{\varepsilon^2}\left(\zeta_{\cQ} + \sqrt{\frac{\omega\zeta_{\cQ}}{n}\left(d-\zeta_{\cQ}\right)}\right)\right)
	\end{eqnarray*}
	 under an assumption that the communication cost is proportional to the number of non-zero components of transmitted vectors from workers to the server.
\end{proof}

\subsection{Convergence Results Under Polyak-{\L}ojasiewicz Condition}\label{sec:proof_of_thm_pl}
In this section, we provide the full statement of Theorem~\ref{thm:main_result_pl} together with the proof of this result.
\begin{theorem}[Theorem~\ref{thm:main_result_pl}]\label{thm:main_result_pl_appendix}
	Let Assumptions~\ref{as:lower_bound},~\ref{as:L_smoothness}~and~\ref{as:pl_condition} be satisfied and 
	\begin{equation}
		\gamma \le \min\left\{\frac{1}{L\left(1 + \sqrt{\frac{2(1-p)\omega}{pn}}\right)}, \frac{p}{2\mu}\right\},\label{eq:gamma_bound_pl_appendix}
	\end{equation}
	where $L^2 = \frac{1}{n}\sum_{i=1}^nL_i^2$. Then after $K$ iterations of \algname{MARINA} we have
	\begin{equation}
		\EE\left[f(x^K) - f(x^*)\right] \le (1-\gamma\mu)^K\Delta_0, \label{eq:main_res_pl_appendix}
	\end{equation}
	where $\Delta_0 = f(x^0)-f(x^*)$. That is, after
	\begin{equation}
		K = \cO\left(\max\left\{\frac{1}{p},\frac{L}{\mu}\left(1 + \sqrt{\frac{(1-p)\omega}{pn}}\right)\right\}\log\frac{\Delta_0}{\varepsilon}\right) \label{eq:main_res_2_pl_appendix}
	\end{equation}
	iterations \algname{MARINA} produces such a point $x^K$ that $\EE[f(x^K) - f(x^*)] \le \varepsilon$.
	Moreover, under an assumption that the communication cost is proportional to the number of non-zero components of transmitted vectors from workers to the server, we have that the expected total communication cost per worker equals
	\begin{equation}
		d + K(pd + (1-p)\zeta_{\cQ}) =  \cO\left(d+\max\left\{\frac{1}{p},\frac{L}{\mu}\left(1 + \sqrt{\frac{(1-p)\omega}{pn}}\right)\right\}(pd + (1-p)\zeta_{\cQ})\log\frac{\Delta_0}{\varepsilon}\right),\label{eq:main_res_4_pl_appendix}
	\end{equation}
	where $\zeta_{\cQ}$ is the expected density of the quantization (see Def.~\ref{def:quantization}).
\end{theorem}
\begin{proof}[Proof of Theorem~\ref{thm:main_result_pl}]
	The proof is very similar to the proof of Theorem~\ref{thm:main_result_non_cvx}. From Lemma~\ref{lem:lemma_2_page} and P{\L} condition, we have
	\begin{eqnarray}
		\EE[f(x^{k+1}) - f(x^*)] &\le& \EE[f(x^k) - f(x^*)] - \frac{\gamma}{2}\EE\left[\|\nabla f(x^k)\|^2\right] - \left(\frac{1}{2\gamma} - \frac{L}{2}\right)\EE\left[\|x^{k+1}-x^k\|^2\right]\notag\\
		&&\quad + \frac{\gamma}{2}\EE\left[\|g^k - \nabla f(x^k)\|^2\right]\notag\\
		&\overset{\eqref{eq:pl_condition}}{\le}& (1-\gamma\mu)\EE\left[f(x^k) - f(x^*)\right] - \left(\frac{1}{2\gamma} - \frac{L}{2}\right)\EE\left[\|x^{k+1}-x^k\|^2\right]\\
		&&\quad + \frac{\gamma}{2}\EE\left[\|g^k - \nabla f(x^k)\|^2\right]. \notag
	\end{eqnarray}
	Using the same arguments as in the proof of \eqref{eq:non_cvx_technical_2}, we obtain
	\begin{eqnarray}
		\EE\left[\|g^{k+1}-\nabla f(x^{k+1})\|^2\right] &\le& \frac{(1-p)\omega L^2}{n}\EE\left[\|x^{k+1}-x^k\|^2\right] + (1-p)\EE\left[\left\|g^k - \nabla f(x^k)\right\|^2\right].\notag
	\end{eqnarray}
	Putting all together, we derive that the sequence $\Phi_k = f(x^k) - f(x^*) + \frac{\gamma}{p}\|g^k - \nabla f(x^k)\|^2$ satisfies
	\begin{eqnarray}
		\EE\left[\Phi_{k+1}\right] &\le& \EE\left[(1-\gamma\mu)(f(x^k) - f(x^*)) - \left(\frac{1}{2\gamma} - \frac{L}{2}\right)\|x^{k+1}-x^k\|^2 + \frac{\gamma}{2}\|g^k - \nabla f(x^k)\|^2\right]\notag\\
		&&\quad + \frac{\gamma}{p}\EE\left[\frac{(1-p)\omega L^2}{n}\|x^{k+1}-x^k\|^2 + (1-p)\left\|g^k - \nabla f(x^k)\right\|^2 \right] \notag\\
		&=& \EE\left[(1-\gamma\mu)(f(x^k) - f(x^*)) + \left(\frac{\gamma}{2} + \frac{\gamma}{p}(1-p)\right)\left\|g^k - \nabla f(x^k)\right\|^2\right]\notag\\
		&&\quad + \left(\frac{\gamma(1-p)\omega L^2}{pn} - \frac{1}{2\gamma} + \frac{L}{2}\right)\EE\left[\|x^{k+1}-x^k\|^2\right]\notag\\
		&\overset{\eqref{eq:gamma_bound_pl_appendix}}{\le}& (1-\gamma\mu)\EE[\Phi_k],\notag
	\end{eqnarray}
	where in the last inequality, we use $\frac{\gamma(1-p)\omega L^2}{pn} - \frac{1}{2\gamma} + \frac{L}{2} \le 0$ and $\frac{\gamma}{2} + \frac{\gamma}{p}(1-p) \le (1-\gamma\mu)\frac{\gamma}{p}$ following from \eqref{eq:gamma_bound_pl_appendix}. Unrolling the recurrence and using $g^0 = \nabla f(x^0)$, we obtain 
	\begin{eqnarray*}
		\EE\left[f(x^{K}) - f(x^*)\right] \le \EE[\Phi_{K}] &\le& (1-\gamma\mu)^{K}\Phi_0 = (1-\gamma\mu)^{K}(f(x^0) - f(x^*))
	\end{eqnarray*}
	that implies \eqref{eq:main_res_2_pl_appendix} and \eqref{eq:main_res_4_pl_appendix}.
\end{proof}

\begin{corollary}\label{cor:main_result_pl_appendix}
	Let the assumptions of Theorem~\ref{thm:main_result_pl} hold and $p = \frac{\zeta_{\cQ}}{d}$, where $\zeta_{\cQ}$ is the expected density of the quantization (see Def.~\ref{def:quantization}). If 
	\begin{equation*}
		\gamma \le \min\left\{\frac{1}{L\left(1 + \sqrt{\frac{2\omega}{n}\left(\frac{d}{\zeta_{\cQ}}-1\right)}\right)}, \frac{p}{2\mu}\right\},
	\end{equation*}
	then \algname{MARINA} requires 
	\begin{equation*}
		K = \cO\left(\max\left\{\frac{d}{\zeta_{\cQ}},\frac{L}{\mu}\left(1 + \sqrt{\frac{\omega}{n}\left(\frac{d}{\zeta_{\cQ}}-1\right)}\right)\right\}\log\frac{\Delta_0}{\varepsilon}\right)
	\end{equation*}
	iterations/communication rounds to achieve $\EE[f(x^K) - f(x^*)] \le \varepsilon$, and the expected total communication cost per worker is
	\begin{equation*}
		\cO\left(d+\max\left\{d,\frac{L}{\mu}\left(\zeta_{\cQ} + \sqrt{\frac{\omega\zeta_{\cQ}}{n}\left(d-\zeta_{\cQ}\right)}\right)\right\}\log\frac{\Delta_0}{\varepsilon}\right)
	\end{equation*}
	 under an assumption that the communication cost is proportional to the number of non-zero components of transmitted vectors from workers to the server.
\end{corollary}
\begin{proof}
	The choice of $p = \frac{\zeta_{\cQ}}{d}$ implies
	\begin{eqnarray*}
		\frac{1-p}{p} &=& \frac{d}{\zeta_{\cQ}}-1,\\
		pd + (1-p)\zeta_{\cQ} &\le& \zeta_{\cQ} + \left(1 - \frac{\zeta_{\cQ}}{d}\right)\cdot\zeta_{\cQ} \le 2\zeta_{\cQ}.
	\end{eqnarray*}	
	Plugging these relations in \eqref{eq:gamma_bound_pl_appendix}, \eqref{eq:main_res_2_pl_appendix}, and \eqref{eq:main_res_4_pl_appendix}, we get that if
	\begin{equation*}
		\gamma \le \min\left\{\frac{1}{L\left(1 + \sqrt{\frac{2\omega}{n}\left(\frac{d}{\zeta_{\cQ}}-1\right)}\right)}, \frac{p}{2\mu}\right\},
	\end{equation*}
	then \algname{MARINA} requires 
	\begin{eqnarray*}
		K &=& \cO\left(\max\left\{\frac{1}{p},\frac{L}{\mu}\left(1 + \sqrt{\frac{(1-p)\omega}{pn}}\right)\right\}\log\frac{\Delta_0}{\varepsilon}\right)\\
		&=& \cO\left(\max\left\{\frac{d}{\zeta_{\cQ}},\frac{L}{\mu}\left(1 + \sqrt{\frac{\omega}{n}\left(\frac{d}{\zeta_{\cQ}}-1\right)}\right)\right\}\log\frac{\Delta_0}{\varepsilon}\right)
	\end{eqnarray*}
	iterations/communication rounds in order to achieve $\EE[f(x^K)-f(x^*)] \le \varepsilon$, and the expected total communication cost per worker is
	\begin{eqnarray*}
		d + K(pd + (1-p)\zeta_{\cQ}) &=&  \cO\left(d+\max\left\{\frac{1}{p},\frac{L}{\mu}\left(1 + \sqrt{\frac{(1-p)\omega}{pn}}\right)\right\}(pd + (1-p)\zeta_{\cQ})\log\frac{\Delta_0}{\varepsilon}\right)\\
		&=&\cO\left(d+\max\left\{d,\frac{L}{\mu}\left(\zeta_{\cQ} + \sqrt{\frac{\omega\zeta_{\cQ}}{n}\left(d-\zeta_{\cQ}\right)}\right)\right\}\log\frac{\Delta_0}{\varepsilon}\right)
	\end{eqnarray*}
	 under an assumption that the communication cost is proportional to the number of non-zero components of transmitted vectors from workers to the server.
\end{proof}

\section{Missing Proofs for \algname{VR-MARINA}}\label{sec:missing_proofs}

\subsection{Finite Sum Case}

\subsubsection{Generally Non-Convex Problems}\label{sec:proof_of_thm_non_cvx_fin_sums}
In this section, we provide the full statement of Theorem~\ref{thm:main_result_non_cvx_finite_sums} together with the proof of this result.
\begin{theorem}[Theorem~\ref{thm:main_result_non_cvx_finite_sums}]\label{thm:main_result_non_cvx_finite_sums_appendix}
	Consider the finite sum case \eqref{eq:main_problem_marina}+\eqref{eq:f_i_finite_sum}. Let Assumptions~\ref{as:lower_bound},~\ref{as:L_smoothness}~and~\ref{as:avg_smoothness} be satisfied and 
	\begin{equation}
		\gamma \le \frac{1}{L + \sqrt{\frac{1-p}{pn}\left(\omega L^2 + \frac{(1+\omega)\cL^2}{b'}\right)}},\label{eq:gamma_bound_non_cvx_finite_sums_appendix}
	\end{equation}
	where $L^2 = \frac{1}{n}\sum_{i=1}^nL_i^2$ and $\cL^2 = \frac{1}{n}\sum_{i=1}^n\cL_i^2$. Then after $K$ iterations of \algname{VR-MARINA} we have
	\begin{equation}
		\EE\left[\left\|\nabla f(\hat x^K)\right\|^2\right] \le \frac{2\Delta_0}{\gamma K}, \label{eq:main_res_non_cvx_finite_sums_appendix}
	\end{equation}
	where $\hat{x}^K$ is chosen uniformly at random from $x^0,\ldots,x^{K-1}$ and $\Delta_0 = f(x^0)-f_*$. That is, after
	\begin{equation}
		K = \cO\left(\frac{\Delta_0}{\varepsilon^2}\left(L + \sqrt{\frac{1-p}{pn}\left(\omega L^2 + \frac{(1+\omega)\cL^2}{b'}\right)}\right)\right) \label{eq:main_res_2_non_cvx_finite_sums_appendix}
	\end{equation}
	iterations \algname{VR-MARINA} produces such a point $\hat x^K$ that $\EE[\|\nabla f(\hat x^K)\|^2] \le \varepsilon^2$, and the expected total number of stochastic oracle calls per node equals
	\begin{equation}
		m + K(pm + 2(1-p)b') = \cO\left(m + \frac{\Delta_0}{\varepsilon^2}\left(L + \sqrt{\frac{1-p}{pn}\left(\omega L^2 + \frac{(1+\omega)\cL^2}{b'}\right)}\right)(pm + (1-p)b')\right). \label{eq:main_res_3_non_cvx_finite_sums_appendix}
	\end{equation}
	Moreover, under an assumption that the communication cost is proportional to the number of non-zero components of transmitted vectors from workers to the server, we have that the expected total communication cost per worker equals
	\begin{equation}
		d + K(pd + (1-p)\zeta_{\cQ}) =  \cO\left(d+\frac{\Delta_0}{\varepsilon^2}\left(L + \sqrt{\frac{1-p}{pn}\left(\omega L^2 + \frac{(1+\omega)\cL^2}{b'}\right)}\right)(pd + (1-p)\zeta_{\cQ})\right),\label{eq:main_res_4_non_cvx_finite_sums_appendix}
	\end{equation}
	where $\zeta_{\cQ}$ is the expected density of the quantization (see Def.~\ref{def:quantization}).
\end{theorem}
\begin{proof}[Proof of Theorem~\ref{thm:main_result_non_cvx_finite_sums}]
	The proof of this theorem is a generalization of the proof of Theorem~\ref{thm:main_result_non_cvx}. From Lemma~\ref{lem:lemma_2_page}, we have
	\begin{equation}
		\EE[f(x^{k+1})] \le \EE[f(x^k)] - \frac{\gamma}{2}\EE\left[\|\nabla f(x^k)\|^2\right] - \left(\frac{1}{2\gamma} - \frac{L}{2}\right)\EE\left[\|x^{k+1}-x^k\|^2\right] + \frac{\gamma}{2}\EE\left[\|g^k - \nabla f(x^k)\|^2\right]. \label{eq:non_cvx_finite_sums_technical_1}
	\end{equation}
	Next, we need to derive an upper bound for $\EE\left[\|g^{k+1}-\nabla f(x^{k+1})\|^2\right]$. Since $g^{k+1} = \frac{1}{n}\sum\limits_{i=1}^ng_i^{k+1}$, we get the following representation of $g^{k+1}$:
	\begin{equation}
		g^{k+1} = \begin{cases}\nabla f(x^{k+1})& \text{with probability } p,\\ g^k + \frac{1}{n}\sum\limits_{i=1}^n \cQ\left(\frac{1}{b'}\sum\limits_{j\in I'_{i,k}}(\nabla f_{ij}(x^{k+1}) - \nabla f_{ij}(x^k))\right)& \text{with probability } 1-p. \end{cases}\notag
	\end{equation}
	Using this, variance decomposition \eqref{eq:variance_decomposition} and tower property \eqref{eq:tower_property}, we derive:
	\begin{eqnarray}
		\EE\left[\|g^{k+1}-\nabla f(x^{k+1})\|^2\right]&&\notag\\ 
		&&\hspace{-4.5cm}\overset{\eqref{eq:tower_property}}{=} (1-p)\EE\left[\left\|g^k + \frac{1}{n}\sum\limits_{i=1}^n \cQ\left(\frac{1}{b'}\sum\limits_{j\in I'_{i,k}}(\nabla f_{ij}(x^{k+1}) - \nabla f_{ij}(x^k))\right) - \nabla f(x^{k+1})\right\|^2\right]\notag\\
		&&\hspace{-5cm}\overset{\eqref{eq:tower_property},\eqref{eq:variance_decomposition}}{=} (1-p)\EE\left[\left\|\frac{1}{n}\sum\limits_{i=1}^n \cQ\left(\frac{1}{b'}\sum\limits_{j\in I'_{i,k}}(\nabla f_{ij}(x^{k+1}) - \nabla f_{ij}(x^k))\right) - \nabla f(x^{k+1}) + \nabla f(x^k)\right\|^2\right]\notag\\
		&&\hspace{-3.5cm} + (1-p)\EE\left[\left\|g^k - \nabla f(x^k)\right\|^2\right].\notag
	\end{eqnarray}
	Next, we use the notation: $\widetilde{\Delta}_i^k = \frac{1}{b'}\sum\limits_{j\in I'_{i,k}}(\nabla f_{ij}(x^{k+1}) - \nabla f_{ij}(x^k))$ and $\Delta_i^k = \nabla f_i(x^{k+1}) - \nabla f_i(x^k)$. These vectors satisfy $\EE\left[\widetilde{\Delta}_i^k \mid x^k,x^{k+1}\right] = \Delta_i^k$ for all $i\in [n]$. Moreover, $\cQ(\widetilde{\Delta}_1^k),\ldots,\cQ(\widetilde{\Delta}_n^k)$ are independent random vectors for fixed $x^k$ and $x^{k+1}$. These observations imply
	\begin{eqnarray}
		\EE\left[\|g^{k+1}-\nabla f(x^{k+1})\|^2\right] &=& (1-p)\EE\left[\left\|\frac{1}{n}\sum\limits_{i=1}^n \left(\cQ(\widetilde{\Delta}_i^k) - \Delta_i^k\right)\right\|^2\right]\notag\\
		&&\quad +(1-p)\EE\left[\left\|g^k - \nabla f(x^k)\right\|^2\right]\notag\\
		&=& \frac{1-p}{n^2}\sum\limits_{i=1}^n\EE\left[\left\|\cQ(\widetilde{\Delta}_i^k) - \widetilde{\Delta}_i^k + \widetilde{\Delta}_i^k - \Delta_i^k\right\|^2\right]\notag\\
		&&\quad + (1-p)\EE\left[\left\|g^k - \nabla f(x^k)\right\|^2\right]\notag\\
		&\overset{\eqref{eq:tower_property},\eqref{eq:variance_decomposition}}{=}& \frac{1-p}{n^2}\sum\limits_{i=1}^n\left(\EE\left[\left\|\cQ(\widetilde{\Delta}_i^k) - \widetilde{\Delta}_i^k\right\|^2\right] + \EE\left[\left\|\widetilde{\Delta}_i^k - \Delta_i^k\right\|^2\right]\right)\notag\\
		&&\quad + (1-p)\EE\left[\left\|g^k - \nabla f(x^k)\right\|^2\right]\notag\\
		&\overset{\eqref{eq:tower_property},\eqref{eq:quantization_def}}{=}& \frac{1-p}{n^2}\sum\limits_{i=1}^n\left(\omega\EE\left[\left\|\widetilde{\Delta}_i^k\right\|^2\right] + \EE\left[\left\|\widetilde{\Delta}_i^k - \Delta_i^k\right\|^2\right]\right)\notag\\
		&&\quad + (1-p)\EE\left[\left\|g^k - \nabla f(x^k)\right\|^2\right]\notag\\
		&\overset{\eqref{eq:tower_property},\eqref{eq:variance_decomposition}}{=}& \frac{1-p}{n^2}\sum\limits_{i=1}^n\left(\omega\EE\left[\left\|\Delta_i^k\right\|^2\right] + (1+\omega)\EE\left[\left\|\widetilde{\Delta}_i^k - \Delta_i^k\right\|^2\right]\right)\notag\\
		&&\quad + (1-p)\EE\left[\left\|g^k - \nabla f(x^k)\right\|^2\right].\notag
	\end{eqnarray}
	Using $L$-smoothness \eqref{eq:L_smoothness_local_marina} and average $\cL$-smoothness \eqref{eq:avg_L_smoothness} of $f_i$ together with the tower property \eqref{eq:tower_property}, we get
	\begin{eqnarray}
		\EE\left[\|g^{k+1}-\nabla f(x^{k+1})\|^2\right] &\le& \frac{1-p}{n^2}\sum\limits_{i=1}^n\left(\omega L_i^2 + \frac{(1+\omega)\cL_i^2}{b'}\right)\EE\left[\|x^{k+1} - x^k\|^2\right]\notag\\
		&&\quad + (1-p)\EE\left[\left\|g^k - \nabla f(x^k)\right\|^2\right]\notag\\
		&=&\frac{1-p}{n}\left(\omega L^2 + \frac{(1+\omega)\cL^2}{b'}\right)\EE\left[\|x^{k+1}-x^k\|^2\right]\notag\\
		&&\quad + (1-p)\EE\left[\left\|g^k - \nabla f(x^k)\right\|^2\right].\label{eq:non_cvx_finite_sums_technical_2}
	\end{eqnarray}
	Next, we introduce new notation: $\Phi_k = f(x^k) - f_* + \frac{\gamma}{2p}\|g^k - \nabla f(x^k)\|^2$. Using this and inequalities \eqref{eq:non_cvx_finite_sums_technical_1} and \eqref{eq:non_cvx_finite_sums_technical_2}, we establish the following inequality:
	\begin{eqnarray}
		\EE\left[\Phi_{k+1}\right] &\le& \EE\left[f(x^k) - f_* - \frac{\gamma}{2}\|\nabla f(x^k)\|^2 - \left(\frac{1}{2\gamma} - \frac{L}{2}\right)\|x^{k+1}-x^k\|^2 + \frac{\gamma}{2}\|g^k - \nabla f(x^k)\|^2\right]\notag\\
		&&\quad + \frac{\gamma}{2p}\EE\left[\frac{1-p}{n}\left(\omega L^2 + \frac{(1+\omega)\cL^2}{b'}\right)\|x^{k+1}-x^k\|^2 + (1-p)\left\|g^k - \nabla f(x^k)\right\|^2\right] \notag\\
		&=& \EE\left[\Phi_k\right] - \frac{\gamma}{2}\EE\left[\|\nabla f(x^k)\|^2\right]\notag\\
		&&\quad + \left(\frac{\gamma(1-p)}{2pn}\left(\omega L^2 + \frac{(1+\omega)\cL^2}{b'}\right) - \frac{1}{2\gamma} + \frac{L}{2}\right)\EE\left[\|x^{k+1}-x^k\|^2\right]\notag\\
		&\overset{\eqref{eq:gamma_bound_non_cvx_finite_sums_appendix}}{\le}& \EE\left[\Phi_k\right] - \frac{\gamma}{2}\EE\left[\|\nabla f(x^k)\|^2\right],\label{eq:non_cvx_finite_sums_technical_3}
	\end{eqnarray}
	where in the last inequality, we use $\frac{\gamma(1-p)}{2pn}\left(\omega L^2 + \frac{(1+\omega)\cL^2}{b'}\right) - \frac{1}{2\gamma} + \frac{L}{2} \le 0$ following from \eqref{eq:gamma_bound_non_cvx_finite_sums_appendix}. Summing up inequalities \eqref{eq:non_cvx_finite_sums_technical_3} for $k=0,1,\ldots,K-1$ and rearranging the terms, we derive
	\begin{eqnarray}
		\frac{1}{K}\sum\limits_{k=0}^{K-1}\EE\left[\|\nabla f(x^k)\|^2\right] &\le& \frac{2}{\gamma K}\sum\limits_{k=0}^{K-1}\left(\EE[\Phi_k]-\EE[\Phi_{k+1}]\right) = \frac{2\left(\EE[\Phi_0]-\EE[\Phi_{K}]\right)}{\gamma K} = \frac{2\Delta_0}{\gamma K},\notag
	\end{eqnarray}
	since $g^0 = \nabla f(x^0)$ and $\Phi_{k+1} \ge 0$. Finally, using the tower property \eqref{eq:tower_property} and the definition of $\hat x^K$, we obtain \eqref{eq:main_res_non_cvx_finite_sums_appendix} that implies \eqref{eq:main_res_2_non_cvx_finite_sums_appendix}, \eqref{eq:main_res_3_non_cvx_finite_sums_appendix}, and \eqref{eq:main_res_4_non_cvx_finite_sums_appendix}.
\end{proof}

\begin{remark}[About batchsizes dissimilarity]
	We notice that our analysis can be easily extended to handle the version of \algname{VR-MARINA} with different batchsizes $b'_1,\ldots,b'_n$ on different workers, i.e., when $|I_{i,k}'| = b_i'$ and $\widetilde{\Delta}_i^k = \frac{1}{b_i'}\sum_{j\in I_{i,k}'}(\nabla f_{ij}(x^{k+1}) - \nabla f_{ij}(x^k))$. In this case, the statement of Theorem~\ref{thm:main_result_non_cvx_finite_sums} remains the same with the small modificiation: instead of $\frac{\cL^2}{b'}$ the complexity bounds will have $\frac{1}{n}\sum_{i=1}^n\frac{\cL_i^2}{b_i'}$.
\end{remark}

\begin{corollary}[Corollary~\ref{cor:main_result_non_cvx_finite_sums}]\label{cor:main_result_non_cvx_finite_sums_appendix}
	Let the assumptions of Theorem~\ref{thm:main_result_non_cvx_finite_sums} hold and $p = \min\left\{\frac{\zeta_{\cQ}}{d},\frac{b'}{m+b'}\right\}$, where $b' \le m$ and $\zeta_{\cQ}$ is the expected density of the quantization (see Def.~\ref{def:quantization}). If 
	\begin{equation*}
		\gamma \le \frac{1}{L + \sqrt{\frac{\max\left\{\nicefrac{d}{\zeta_{\cQ}} - 1,\nicefrac{m}{b'}\right\}}{n}\left(\omega L^2 + \frac{(1+\omega)\cL^2}{b'}\right)}},
	\end{equation*}
	then \algname{VR-MARINA} requires 
	\begin{equation*}
		\cO\left(\frac{\Delta_0}{\varepsilon^2}\left(L\left(1 + \sqrt{\frac{\omega\max\left\{\nicefrac{d}{\zeta_{\cQ}} - 1,\nicefrac{m}{b'}\right\}}{n}}\right) + \cL\sqrt{\frac{(1+\omega)\max\left\{\nicefrac{d}{\zeta_{\cQ}} - 1,\nicefrac{m}{b'}\right\}}{nb'}}\right)\right)
	\end{equation*}
	iterations/communication rounds, 
	\begin{equation*}
		\cO\left(m+\frac{\Delta_0}{\varepsilon^2}\left(L\left(b' + \sqrt{\frac{\omega\max\left\{(\nicefrac{d}{\zeta_{\cQ}} - 1)(b')^2,mb'\right\}}{n}}\right) + \cL\sqrt{\frac{(1+\omega)\max\left\{(\nicefrac{d}{\zeta_{\cQ}} - 1) b',m\right\}}{n}}\right)\right)
	\end{equation*}
	stochastic oracle calls per node in expectation in order to achieve $\EE[\|\nabla f(\hat x^K)\|^2] \le \varepsilon^2$, and the expected total communication cost per worker is
	\begin{equation*}
		\cO\left(d+\frac{\Delta_0\zeta_{\cQ}}{\varepsilon^2}\left(L\left(1 + \sqrt{\frac{\omega\max\left\{\nicefrac{d}{\zeta_{\cQ}} - 1,\nicefrac{m}{b'}\right\}}{n}}\right) + \cL\sqrt{\frac{(1+\omega)\max\left\{\nicefrac{d}{\zeta_{\cQ}} - 1,\nicefrac{m}{b'}\right\}}{nb'}}\right)\right)
	\end{equation*}
	 under an assumption that the communication cost is proportional to the number of non-zero components of transmitted vectors from workers to the server.
\end{corollary}
\begin{proof}[Proof of Corollary~\ref{cor:main_result_non_cvx_finite_sums}]
	The choice of $p = \min\left\{\frac{\zeta_{\cQ}}{d},\frac{b'}{m+b'}\right\}$ implies
	\begin{eqnarray*}
		\frac{1-p}{p} &=& \max\left\{\frac{d}{\zeta_{\cQ}}-1,\frac{m}{b'}\right\},\\
		pm + (1-p)b' &\le& \frac{2mb'}{m+b'} \le 2b',\\
		pd + (1-p)\zeta_{\cQ} &\le& \frac{\zeta_{\cQ}}{d}\cdot d + \left(1 - \frac{\zeta_{\cQ}}{d}\right)\cdot\zeta_{\cQ} \le 2\zeta_{\cQ}.
	\end{eqnarray*}	
	Plugging these relations in \eqref{eq:gamma_bound_non_cvx_finite_sums_appendix}, \eqref{eq:main_res_2_non_cvx_finite_sums_appendix}, \eqref{eq:main_res_3_non_cvx_finite_sums_appendix} and \eqref{eq:main_res_4_non_cvx_finite_sums_appendix} and using $\sqrt{a+b} \le \sqrt{a} + \sqrt{b}$, we get that if
	\begin{equation*}
		\gamma \le \frac{1}{L + \sqrt{\frac{\max\left\{\nicefrac{d}{\zeta_{\cQ}} - 1,\nicefrac{m}{b'}\right\}}{n}\left(\omega L^2 + \frac{(1+\omega)\cL^2}{b'}\right)}},
	\end{equation*}
	then \algname{VR-MARINA} requires 
	\begin{eqnarray*}
		K &=& \cO\left(\frac{\Delta_0}{\varepsilon^2}\left(L + \sqrt{\frac{1-p}{pn}\left(\omega L^2 + \frac{(1+\omega)\cL^2}{b'}\right)}\right)\right)\\
		&=& \cO\left(\frac{\Delta_0}{\varepsilon^2}\left(L + \sqrt{L^2\frac{\omega\max\left\{\nicefrac{d}{\zeta_{\cQ}} - 1,\nicefrac{m}{b'}\right\}}{n} + \cL^2\frac{(1+\omega)\max\left\{\nicefrac{d}{\zeta_{\cQ}} - 1,\nicefrac{m}{b'}\right\}}{nb'}}\right)\right)\\
		&=& \cO\left(\frac{\Delta_0}{\varepsilon^2}\left(L\left(1 + \sqrt{\frac{\omega\max\left\{\nicefrac{d}{\zeta_{\cQ}} - 1,\nicefrac{m}{b'}\right\}}{n}}\right) + \cL\sqrt{\frac{(1+\omega)\max\left\{\nicefrac{d}{\zeta_{\cQ}} - 1,\nicefrac{m}{b'}\right\}}{nb'}}\right)\right)
	\end{eqnarray*}
	iterations/communication rounds and 
	\begin{eqnarray*}
		m + K(pm + 2(1-p)b') &=& \cO\left(m + \frac{\Delta_0}{\varepsilon^2}\left(L + \sqrt{\frac{1-p}{pn}\left(\omega L^2 + \frac{(1+\omega)\cL^2}{b'}\right)}\right)(pm + (1-p)b')\right)\\	
		&=& \cO\Bigg(m+\frac{\Delta_0}{\varepsilon^2}\Bigg(L\Bigg(1 + \sqrt{\frac{\omega\max\left\{\nicefrac{d}{\zeta_{\cQ}} - 1,\nicefrac{m}{b'}\right\}}{n}}\Bigg)\\
		&&\hspace{4.5cm} + \cL\sqrt{\frac{(1+\omega)\max\left\{\nicefrac{d}{\zeta_{\cQ}} - 1,\nicefrac{m}{b'}\right\}}{nb'}}\Bigg)b'\Bigg)\\
		&=& \cO\Bigg(m+\frac{\Delta_0}{\varepsilon^2}\Bigg(L\Bigg(b' + \sqrt{\frac{\omega\max\left\{(\nicefrac{d}{\zeta_{\cQ}} - 1)(b')^2, mb'\right\}}{n}}\Bigg)\\
		&&\hspace{4.5cm} + \cL\sqrt{\frac{(1+\omega)\max\left\{(\nicefrac{d}{\zeta_{\cQ}} - 1)b',m\right\}}{n}}\Bigg)\Bigg)
	\end{eqnarray*}
	stochastic oracle calls per node in expectation in order to achieve $\EE[\|\nabla f(\hat x^K)\|^2] \le \varepsilon^2$, and the expected total communication cost per worker is
	\begin{eqnarray*}
		d + K(pd + (1-p)\zeta_{\cQ}) &=&  \cO\left(d+\frac{\Delta_0}{\varepsilon^2}\left(L + \sqrt{\frac{1-p}{pn}\left(\omega L^2 + \frac{(1+\omega)\cL^2}{b'}\right)}\right)(pd + (1-p)\zeta_{\cQ})\right)\\
		&=&\cO\Bigg(d+\frac{\Delta_0\zeta_{\cQ}}{\varepsilon^2}\Bigg(L\Bigg(1 + \sqrt{\frac{\omega\max\left\{\nicefrac{d}{\zeta_{\cQ}} - 1,\nicefrac{m}{b'}\right\}}{n}}\Bigg) \\
		&&\hspace{4.5cm} + \cL\sqrt{\frac{(1+\omega)\max\left\{\nicefrac{d}{\zeta_{\cQ}} - 1,\nicefrac{m}{b'}\right\}}{nb'}}\Bigg)\Bigg)
	\end{eqnarray*}
	 under an assumption that the communication cost is proportional to the number of non-zero components of transmitted vectors from workers to the server.
\end{proof}

\subsubsection{Convergence Results Under Polyak-{\L}ojasiewicz condition}\label{sec:proof_of_thm_pl_fin_sums}
In this section, we provide an analysis of \algname{VR-MARINA} under the Polyak-{\L}ojasiewicz condition in the finite sum case.
\begin{theorem}\label{thm:main_result_pl_finite_sums_appendix}
	Consider the finite sum case \eqref{eq:main_problem_marina}+\eqref{eq:f_i_finite_sum}. Let Assumptions~\ref{as:lower_bound},~\ref{as:L_smoothness},~\ref{as:avg_smoothness} and \ref{as:pl_condition} be satisfied and 
	\begin{equation}
		\gamma \le \min\left\{\frac{1}{L + \sqrt{\frac{2(1-p)}{pn}\left(\omega L^2 + \frac{(1+\omega)\cL^2}{b'}\right)}},\frac{p}{2\mu}\right\},\label{eq:gamma_bound_pl_finite_sums_appendix}
	\end{equation}
	where $L^2 = \frac{1}{n}\sum_{i=1}^nL_i^2$ and $\cL^2 = \frac{1}{n}\sum_{i=1}^n\cL_i^2$. Then after $K$ iterations of \algname{VR-MARINA}, we have
	\begin{equation}
		\EE\left[f(x^K) - f(x^*)\right] \le (1-\gamma\mu)^K\Delta_0, \label{eq:main_res_pl_finite_sums_appendix}
	\end{equation}
	where $\Delta_0 = f(x^0)-f(x^*)$. That is, after
	\begin{equation}
		K = \cO\left(\max\left\{\frac{1}{p}, \frac{L + \sqrt{\frac{1-p}{pn}\left(\omega L^2 + \frac{(1+\omega)\cL^2}{b'}\right)}}{\mu}\right\}\log\frac{\Delta_0}{\varepsilon}\right) \label{eq:main_res_2_pl_finite_sums_appendix}
	\end{equation}
	iterations \algname{VR-MARINA} produces such a point $x^K$ that $\EE\left[f(x^K) - f(x^*)\right] \le \varepsilon$, and the expected total number of stochastic oracle calls per node $m + K(pm + 2(1-p)b')$ equals
	\begin{equation}
		 \cO\left(m + \max\left\{\frac{1}{p}, \frac{L + \sqrt{\frac{1-p}{pn}\left(\omega L^2 + \frac{(1+\omega)\cL^2}{b'}\right)}}{\mu}\right\}(pm + (1-p)b')\log\frac{\Delta_0}{\varepsilon}\right). \label{eq:main_res_3_pl_finite_sums_appendix}
	\end{equation}
	Moreover, under an assumption that the communication cost is proportional to the number of non-zero components of transmitted vectors from workers to the server we have that the expected total communication cost per worker $d + K(pd + (1-p)\zeta_{\cQ})$ equals
	\begin{equation}
		 \cO\left(d+\max\left\{\frac{1}{p}, \frac{L + \sqrt{\frac{1-p}{pn}\left(\omega L^2 + \frac{(1+\omega)\cL^2}{b'}\right)}}{\mu}\right\}(pd + (1-p)\zeta_{\cQ})\log\frac{\Delta_0}{\varepsilon}\right),\label{eq:main_res_4_pl_finite_sums_appendix}
	\end{equation}
	where $\zeta_{\cQ}$ is the expected density of the quantization (see Def.~\ref{def:quantization}).
\end{theorem}
\begin{proof}
	The proof is very similar to the proof of Theorem~\ref{thm:main_result_non_cvx_finite_sums}. From Lemma~\ref{lem:lemma_2_page} and P{\L} condition, we have
	\begin{eqnarray}
		\EE[f(x^{k+1}) - f(x^*)] &\le& \EE[f(x^k) - f(x^*)] - \frac{\gamma}{2}\EE\left[\|\nabla f(x^k)\|^2\right] - \left(\frac{1}{2\gamma} - \frac{L}{2}\right)\EE\left[\|x^{k+1}-x^k\|^2\right]\notag\\
		&&\quad + \frac{\gamma}{2}\EE\left[\|g^k - \nabla f(x^k)\|^2\right]\notag\\
		&\overset{\eqref{eq:pl_condition}}{\le}& (1-\gamma\mu)\EE\left[f(x^k) - f(x^*)\right] - \left(\frac{1}{2\gamma} - \frac{L}{2}\right)\EE\left[\|x^{k+1}-x^k\|^2\right]\notag\\
		&&\quad + \frac{\gamma}{2}\EE\left[\|g^k - \nabla f(x^k)\|^2\right]. \notag
	\end{eqnarray}
	Using the same arguments as in the proof of \eqref{eq:non_cvx_finite_sums_technical_2}, we obtain
	\begin{eqnarray}
		\EE\left[\|g^{k+1}-\nabla f(x^{k+1})\|^2\right] &\le&\frac{1-p}{n}\left(\omega L^2 + \frac{(1+\omega)\cL^2}{b'}\right)\EE\left[\|x^{k+1}-x^k\|^2\right] \notag\\
		&&\quad + (1-p)\EE\left[\left\|g^k - \nabla f(x^k)\right\|^2\right].\notag
	\end{eqnarray}
	Putting all together we derive that the sequence $\Phi_k = f(x^k) - f(x^*) + \frac{\gamma}{p}\|g^k - \nabla f(x^k)\|^2$ satisfies
	\begin{eqnarray}
		\EE\left[\Phi_{k+1}\right] &\le& \EE\left[(1-\gamma\mu)(f(x^k) - f(x^*)) - \left(\frac{1}{2\gamma} - \frac{L}{2}\right)\|x^{k+1}-x^k\|^2 + \frac{\gamma}{2}\|g^k - \nabla f(x^k)\|^2\right]\notag\\
		&&\quad + \frac{\gamma}{p}\EE\left[\frac{1-p}{n}\left(\omega L^2 + \frac{(1+\omega)\cL^2}{b'}\right)\|x^{k+1}-x^k\|^2 + (1-p)\left\|g^k - \nabla f(x^k)\right\|^2\right] \notag\\
		&=& \EE\left[(1-\gamma\mu)(f(x^k) - f(x^*)) + \left(\frac{\gamma}{2} + \frac{\gamma}{p}(1-p)\right)\left\|g^k - \nabla f(x^k)\right\|^2\right]\notag\\
		&&\quad + \left(\frac{\gamma(1-p)}{pn}\left(\omega L^2 + \frac{(1+\omega)\cL^2}{b'}\right) - \frac{1}{2\gamma} + \frac{L}{2}\right)\EE\left[\|x^{k+1}-x^k\|^2\right]\notag\\
		&\overset{\eqref{eq:gamma_bound_pl_finite_sums_appendix}}{\le}& (1-\gamma\mu)\EE[\Phi_k], \notag
	\end{eqnarray}
	where in the last inequality we use $\frac{\gamma(1-p)}{pn}\left(\omega L^2 + \frac{(1+\omega)\cL^2}{b'}\right) - \frac{1}{2\gamma} + \frac{L}{2} \le 0$ and $\frac{\gamma}{2} + \frac{\gamma}{p}(1-p) \le (1-\gamma\mu)\frac{\gamma}{p}$ following from \eqref{eq:gamma_bound_pl_finite_sums_appendix}. Unrolling the recurrence and using $g^0 = \nabla f(x^0)$, we obtain 
	\begin{equation*}
		\EE\left[f(x^{k+1}) - f(x^*)\right] \le \EE[\Phi_{k+1}] \le (1-\gamma\mu)^{k+1}\Phi_0 = (1-\gamma\mu)^{k+1}(f(x^0) - f(x^*))
	\end{equation*}
	that implies \eqref{eq:main_res_2_pl_finite_sums_appendix}, \eqref{eq:main_res_3_pl_finite_sums_appendix}, and \eqref{eq:main_res_4_pl_finite_sums_appendix}.
\end{proof}

\begin{corollary}\label{cor:main_result_pl_finite_sums_appendix}
	Let the assumptions of Theorem~\ref{thm:main_result_pl_finite_sums_appendix} hold and $p = \min\left\{\frac{\zeta_{\cQ}}{d},\frac{b'}{m+b'}\right\}$, where $b' \le m$ and $\zeta_{\cQ}$ is the expected density of the quantization (see Def.~\ref{def:quantization}). If 
	\begin{equation*}
		\gamma \le \min\left\{\frac{1}{L + \sqrt{\frac{2\max\left\{\nicefrac{d}{\zeta_{\cQ}} - 1,\nicefrac{m}{b'}\right\}}{n}\left(\omega L^2 + \frac{(1+\omega)\cL^2}{b'}\right)}},\frac{p}{2\mu}\right\},
	\end{equation*}
	then \algname{VR-MARINA} requires 
	\begin{equation*}
		\cO\left(\max\left\{\frac{1}{p}, \frac{L}{\mu}\left(1 + \sqrt{\frac{\omega\max\left\{\nicefrac{d}{\zeta_{\cQ}} - 1,\nicefrac{m}{b'}\right\}}{n}}\right) + \frac{\cL}{\mu}\sqrt{\frac{(1+\omega)\max\left\{\nicefrac{d}{\zeta_{\cQ}} - 1,\nicefrac{m}{b'}\right\}}{nb'}}\right\}\log\frac{\Delta_0}{\varepsilon}\right)
	\end{equation*}
	iterations/communication rounds, 
	\begin{eqnarray*}
	    \cO\Bigg(m+\max\Bigg\{\frac{b'}{p}, \frac{L}{\mu}\left(b' + \sqrt{\frac{\omega\max\left\{(\nicefrac{d}{\zeta_{\cQ}} - 1)(b')^2,mb'\right\}}{n}}\right)&\\
	    &\hspace{-2cm}+ \frac{\cL}{\mu}\sqrt{\frac{(1+\omega)\max\left\{(\nicefrac{d}{\zeta_{\cQ}} - 1) b',m\right\}}{n}}\Bigg\}\log\frac{\Delta_0}{\varepsilon}\Bigg)
	\end{eqnarray*}
	stochastic oracle calls per node in expectation to achieve $\EE[f(x^K)-f(x^*)] \le \varepsilon$, and the expected total communication cost per worker is
	\begin{eqnarray*}
	    \cO\Bigg(d+\zeta_{\cQ}\max\Bigg\{\frac{1}{p}, \frac{L}{\mu}\left(1 + \sqrt{\frac{\omega\max\left\{\nicefrac{d}{\zeta_{\cQ}} - 1,\nicefrac{m}{b'}\right\}}{n}}\right)&\\
	    &\hspace{-2cm}+ \frac{\cL}{\mu}\sqrt{\frac{(1+\omega)\max\left\{\nicefrac{d}{\zeta_{\cQ}} - 1,\nicefrac{m}{b'}\right\}}{nb'}}\Bigg\}\log\frac{\Delta_0}{\varepsilon}\Bigg)
	\end{eqnarray*}
	 under an assumption that the communication cost is proportional to the number of non-zero components of transmitted vectors from workers to the server.
\end{corollary}
\begin{proof}
	The choice of $p = \min\left\{\frac{\zeta_{\cQ}}{d},\frac{b'}{m+b'}\right\}$ implies
	\begin{eqnarray*}
		\frac{1-p}{p} &=& \max\left\{\frac{d}{\zeta_{\cQ}}-1,\frac{m}{b'}\right\},\\
		pm + (1-p)b' &\le& \frac{2mb'}{m+b'} \le 2b',\\
		pd + (1-p)\zeta_{\cQ} &\le& \frac{\zeta_{\cQ}}{d}\cdot d + \left(1 - \frac{\zeta_{\cQ}}{d}\right)\cdot\zeta_{\cQ} \le 2\zeta_{\cQ}.
	\end{eqnarray*}	
	Plugging these relations in \eqref{eq:gamma_bound_pl_finite_sums_appendix}, \eqref{eq:main_res_2_pl_finite_sums_appendix}, \eqref{eq:main_res_3_pl_finite_sums_appendix} and \eqref{eq:main_res_4_pl_finite_sums_appendix} and using $\sqrt{a+b} \le \sqrt{a} + \sqrt{b}$, we get that if
	\begin{equation*}
		\gamma \le \min\left\{\frac{1}{L + \sqrt{\frac{2\max\left\{\nicefrac{d}{\zeta_{\cQ}} - 1,\nicefrac{m}{b'}\right\}}{n}\left(\omega L^2 + \frac{(1+\omega)\cL^2}{b'}\right)}},\frac{p}{2\mu}\right\},
	\end{equation*}
	then \algname{VR-MARINA} requires 
	\begin{eqnarray*}
		K &=& \cO\left(\max\left\{\frac{1}{p}, \frac{L + \sqrt{\frac{1-p}{pn}\left(\omega L^2 + \frac{(1+\omega)\cL^2}{b'}\right)}}{\mu}\right\}\log\frac{\Delta_0}{\varepsilon}\right)\\
		&=& \cO\left(\max\left\{\frac{1}{p}, \frac{L + \sqrt{L^2\frac{\omega\max\left\{\nicefrac{d}{\zeta_{\cQ}} - 1,\nicefrac{m}{b'}\right\}}{n} + \cL^2\frac{(1+\omega)\max\left\{\nicefrac{d}{\zeta_{\cQ}} - 1,\nicefrac{m}{b'}\right\}}{nb'}}}{\mu}\right\}\log\frac{\Delta_0}{\varepsilon}\right)\\
		&=& \cO\Bigg(\max\Bigg\{\frac{1}{p}, \frac{L}{\mu}\left(1 + \sqrt{\frac{\omega\max\left\{\nicefrac{d}{\zeta_{\cQ}} - 1,\nicefrac{m}{b'}\right\}}{n}}\right)\\
		&&\hspace{5cm}+ \frac{\cL}{\mu}\sqrt{\frac{(1+\omega)\max\left\{\nicefrac{d}{\zeta_{\cQ}} - 1,\nicefrac{m}{b'}\right\}}{nb'}}\Bigg\}\log\frac{\Delta_0}{\varepsilon}\Bigg)
	\end{eqnarray*}
	iterations/communication rounds and 
	\begin{eqnarray*}
		m + K(pm + 2(1-p)b') &\\
		&\hspace{-3cm}= \cO\left(m + \max\left\{\frac{1}{p}, \frac{L + \sqrt{\frac{1-p}{pn}\left(\omega L^2 + \frac{(1+\omega)\cL^2}{b'}\right)}}{\mu}\right\}(pm + (1-p)b')\log\frac{\Delta_0}{\varepsilon}\right)\\	
		&\hspace{-6.8cm}= \cO\Bigg(m+\max\Bigg\{\frac{1}{p}, \frac{L}{\mu}\Bigg(1 + \sqrt{\frac{\omega\max\left\{\nicefrac{d}{\zeta_{\cQ}} - 1,\nicefrac{m}{b'}\right\}}{n}}\Bigg)\\
		&\hspace{2cm} + \frac{\cL}{\mu}\sqrt{\frac{(1+\omega)\max\left\{\nicefrac{d}{\zeta_{\cQ}} - 1,\nicefrac{m}{b'}\right\}}{nb'}}\Bigg\}b'\log\frac{\Delta_0}{\varepsilon}\Bigg)\\
		&\hspace{-5.9cm}= \cO\Bigg(m+\max\Bigg\{\frac{b'}{p}, \frac{L}{\mu}\Bigg(b' + \sqrt{\frac{\omega\max\left\{(\nicefrac{d}{\zeta_{\cQ}} - 1)(b')^2, mb'\right\}}{n}}\Bigg)\\
		&\hspace{2cm} + \frac{\cL}{\mu}\sqrt{\frac{(1+\omega)\max\left\{(\nicefrac{d}{\zeta_{\cQ}} - 1)b',m\right\}}{n}}\Bigg\}\log\frac{\Delta_0}{\varepsilon}\Bigg)
	\end{eqnarray*}
	stochastic oracle calls per node in expectation in order to achieve $\EE[f(x^K) - f(x^*)] \le \varepsilon$, and the expected total communication cost per worker is
	\begin{eqnarray*}
		d + K(pd + (1-p)\zeta_{\cQ}) &\\
		&\hspace{-2cm}=  \cO\left(d+\max\left\{\frac{1}{p}, \frac{L + \sqrt{\frac{1-p}{pn}\left(\omega L^2 + \frac{(1+\omega)\cL^2}{b'}\right)}}{\mu}\right\}(pd + (1-p)\zeta_{\cQ})\log\frac{\Delta_0}{\varepsilon}\right)\\
		&\hspace{-5.2cm}=\cO\Bigg(d+\zeta_{\cQ}\max\Bigg\{\frac{1}{p}, \frac{L}{\mu}\Bigg(1 + \sqrt{\frac{\omega\max\left\{\nicefrac{d}{\zeta_{\cQ}} - 1,\nicefrac{m}{b'}\right\}}{n}}\Bigg) \\
		&\hspace{3cm} + \frac{\cL}{\mu}\sqrt{\frac{(1+\omega)\max\left\{\nicefrac{d}{\zeta_{\cQ}} - 1,\nicefrac{m}{b'}\right\}}{nb'}}\Bigg\}\log\frac{\Delta_0}{\varepsilon}\Bigg)
	\end{eqnarray*}
	 under an assumption that the communication cost is proportional to the number of non-zero components of transmitted vectors from workers to the server.
\end{proof}

\subsection{Online Case}

\subsubsection{Generally Non-Convex Problems}\label{sec:proof_of_thm_non_cvx_online}

In this section, we provide the full statement of Theorem~\ref{thm:main_result_non_cvx_online} together with the proof of this result.
\begin{theorem}[Theorem~\ref{thm:main_result_non_cvx_online}]\label{thm:main_result_non_cvx_online_appendix}
	Consider the finite sum case \eqref{eq:main_problem_marina}+\eqref{eq:f_i_expectation}. Let Assumptions~\ref{as:lower_bound},~\ref{as:L_smoothness}~and~\ref{as:avg_smoothness_online} be satisfied and 
	\begin{equation}
		\gamma \le \frac{1}{L + \sqrt{\frac{1-p}{pn}\left(\omega L^2 + \frac{(1+\omega)\cL^2}{b'}\right)}},\label{eq:gamma_bound_non_cvx_online_appendix}
	\end{equation}
	where $L^2 = \frac{1}{n}\sum_{i=1}^nL_i^2$ and $\cL^2 = \frac{1}{n}\sum_{i=1}^n\cL_i^2$. Then after $K$ iterations of \algname{VR-MARINA}, we have
	\begin{equation}
		\EE\left[\left\|\nabla f(\hat x^K)\right\|^2\right] \le \frac{2\Delta_0}{\gamma K} + \frac{\sigma^2}{nb}, \label{eq:main_res_non_cvx_online_appendix}
	\end{equation}
	where $\hat{x}^K$ is chosen uniformly at random from $x^0,\ldots,x^{K-1}$ and $\Delta_0 = f(x^0)-f_*$. That is, after
	\begin{equation}
		K = \cO\left(\frac{\Delta_0}{\varepsilon^2}\left(L + \sqrt{\frac{1-p}{pn}\left(\omega L^2 + \frac{(1+\omega)\cL^2}{b'}\right)}\right)\right) \label{eq:main_res_2_non_cvx_online_appendix}
	\end{equation}
	iterations with $b = \Theta(\frac{\sigma^2}{n\varepsilon^2})$ \algname{VR-MARINA} produces such a point $\hat x^K$ that $\EE[\|\nabla f(\hat x^K)\|^2] \le \varepsilon^2$, and the expected total number of stochastic oracle calls per node $b + K(pb + 2(1-p)b')$ equals
	\begin{equation}
	    \cO\left(\frac{\sigma^2}{n\varepsilon^2} + \frac{\Delta_0}{\varepsilon^2}\left(L + \sqrt{\frac{1-p}{pn}\left(\omega L^2 + \frac{(1+\omega)\cL^2}{b'}\right)}\right)\left(p\frac{\sigma^2}{n\varepsilon^2} + (1-p)b'\right)\right). \label{eq:main_res_3_non_cvx_online_appendix}
	\end{equation}
	Moreover, under an assumption that the communication cost is proportional to the number of non-zero components of transmitted vectors from workers to the server we have that the expected total communication cost per worker $d + K(pd + (1-p)\zeta_{\cQ})$ equals
	\begin{equation}
		\cO\left(d+\frac{\Delta_0}{\varepsilon^2}\left(L + \sqrt{\frac{1-p}{pn}\left(\omega L^2 + \frac{(1+\omega)\cL^2}{b'}\right)}\right)(pd + (1-p)\zeta_{\cQ})\right),\label{eq:main_res_4_non_cvx_online_appendix}
	\end{equation}
	where $\zeta_{\cQ}$ is the expected density of the quantization (see Def.~\ref{def:quantization}).
\end{theorem}
\begin{proof}[Proof of Theorem~\ref{thm:main_result_non_cvx_online}]
	The proof follows the same steps as the proof of Theorem~\ref{thm:main_result_non_cvx_finite_sums}. From Lemma~\ref{lem:lemma_2_page}, we have
	\begin{equation}
		\EE[f(x^{k+1})] \le \EE[f(x^k)] - \frac{\gamma}{2}\EE\left[\|\nabla f(x^k)\|^2\right] - \left(\frac{1}{2\gamma} - \frac{L}{2}\right)\EE\left[\|x^{k+1}-x^k\|^2\right] + \frac{\gamma}{2}\EE\left[\|g^k - \nabla f(x^k)\|^2\right]. \label{eq:non_cvx_online_technical_1}
	\end{equation}
	Next, we need to derive an upper bound for $\EE\left[\|g^{k+1}-\nabla f(x^{k+1})\|^2\right]$. Since $g^{k+1} = \frac{1}{n}\sum\limits_{i=1}^ng_i^{k+1}$, we get the following representation of $g^{k+1}$:
	\begin{equation}
		g^{k+1} = \begin{cases}\frac{1}{nb}\sum\limits_{i=1}^n\sum\limits_{j\in I_{i,k}}\nabla f_{\xi_{ij}^k}(x^{k+1})& \text{with probability } p,\\ g^k + \frac{1}{n}\sum\limits_{i=1}^n \cQ\left(\frac{1}{b'}\sum\limits_{j\in I'_{i,k}}(\nabla f_{\xi_{ij}^k}(x^{k+1}) - \nabla f_{\xi_{ij}^k}(x^k))\right)& \text{with probability } 1-p. \end{cases}\notag
	\end{equation}
	Using this, variance decomposition \eqref{eq:variance_decomposition}, tower property \eqref{eq:tower_property}, and independence of $\xi_{ij}^k$ for $i\in[n]$, $j\in I_{i,k}$, we derive:
	\begin{eqnarray}
		\EE\left[\|g^{k+1}-\nabla f(x^{k+1})\|^2\right] & \notag\\ &\hspace{-4.9cm}\overset{\eqref{eq:tower_property}}{=} (1-p)\EE\left[\left\|g^k + \frac{1}{n}\sum\limits_{i=1}^n \cQ\left(\frac{1}{b'}\sum\limits_{j\in I'_{i,k}}(\nabla f_{ij}(x^{k+1}) - \nabla f_{ij}(x^k))\right) - \nabla f(x^{k+1})\right\|^2\right]\notag\\
		&\hspace{-6cm} + \frac{p}{n^2b^2}\EE\left[\left\|\sum\limits_{i=1}^n\sum\limits_{j\in I_{i,k}}\left(\nabla f_{\xi_{ij}^k}(x^{k+1}) - \nabla f(x^{k+1})\right)\right\|^2\right] \notag\\
		&\hspace{-4cm}\overset{\eqref{eq:tower_property},\eqref{eq:variance_decomposition}}{=} (1-p)\EE\left[\left\|\frac{1}{n}\sum\limits_{i=1}^n \cQ\left(\frac{1}{b'}\sum\limits_{j\in I'_{i,k}}(\nabla f_{ij}(x^{k+1}) - \nabla f_{ij}(x^k))\right) - \nabla f(x^{k+1}) + \nabla f(x^k)\right\|^2\right]\notag\\
		& \hspace{-1.5cm}+ (1-p)\EE\left[\left\|g^k - \nabla f(x^k)\right\|^2\right] + \frac{p}{n^2b^2}\sum\limits_{i=1}^n\sum\limits_{j\in I_{i,k}}\EE\left[\left\|\nabla f_{\xi_{ij}^k}(x^{k+1}) - \nabla f(x^{k+1})\right\|^2\right]\notag\\
		&\hspace{-4cm}\overset{\eqref{eq:tower_property},\eqref{eq:bounded_var}}{=} (1-p)\EE\left[\left\|\frac{1}{n}\sum\limits_{i=1}^n \cQ\left(\frac{1}{b'}\sum\limits_{j\in I'_{i,k}}(\nabla f_{ij}(x^{k+1}) - \nabla f_{ij}(x^k))\right) - \nabla f(x^{k+1}) + \nabla f(x^k)\right\|^2\right]\notag\\
		&\hspace{-8.5cm} + (1-p)\EE\left[\left\|g^k - \nabla f(x^k)\right\|^2\right] + \frac{p\sigma^2}{nb},\notag
	\end{eqnarray}
	where $\sigma^2 = \frac{1}{n}\sum_{i=1}^n\sigma_i^2$. Applying the same arguments as in the proof of inequality \eqref{eq:non_cvx_finite_sums_technical_2}, we obtain
	\begin{eqnarray}
		\EE\left[\|g^{k+1}-\nabla f(x^{k+1})\|^2\right]	&\le&\frac{1-p}{n}\left(\omega L^2 + \frac{(1+\omega)\cL^2}{b'}\right)\EE\left[\|x^{k+1}-x^k\|^2\right]\notag\\
		&&\quad + (1-p)\EE\left[\left\|g^k - \nabla f(x^k)\right\|^2\right] + \frac{p\sigma^2}{nb}.\label{eq:non_cvx_online_technical_2}
	\end{eqnarray}
	Next, we introduce new notation: $\Phi_k = f(x^k) - f_* + \frac{\gamma}{2p}\|g^k - \nabla f(x^k)\|^2$. Using this and inequalities \eqref{eq:non_cvx_online_technical_1} and \eqref{eq:non_cvx_online_technical_2}, we establish the following inequality:
	\begin{eqnarray}
		\EE\left[\Phi_{k+1}\right] &\le& \EE\left[f(x^k) - f_* - \frac{\gamma}{2}\|\nabla f(x^k)\|^2 - \left(\frac{1}{2\gamma} - \frac{L}{2}\right)\|x^{k+1}-x^k\|^2 + \frac{\gamma}{2}\|g^k - \nabla f(x^k)\|^2\right]\notag\\
		&&\quad + \frac{\gamma}{2p}\EE\left[\frac{1-p}{n}\left(\omega L^2 + \frac{(1+\omega)\cL^2}{b'}\right)\|x^{k+1}-x^k\|^2\right] \notag\\
		&&\quad + \frac{\gamma}{2p}\EE\left[(1-p)\left\|g^k - \nabla f(x^k)\right\|^2 + \frac{p\sigma^2}{nb}\right] \notag\\
		&=& \EE\left[\Phi_k\right] - \frac{\gamma}{2}\EE\left[\|\nabla f(x^k)\|^2\right]\notag\\
		&&\quad+ \left(\frac{\gamma(1-p)}{2pn}\left(\omega L^2 + \frac{(1+\omega)\cL^2}{b'}\right) - \frac{1}{2\gamma} + \frac{L}{2}\right)\EE\left[\|x^{k+1}-x^k\|^2\right] + \frac{\gamma\sigma^2}{2nb}\notag\\
		&\overset{\eqref{eq:gamma_bound_non_cvx_online_appendix}}{\le}& \EE\left[\Phi_k\right] - \frac{\gamma}{2}\EE\left[\|\nabla f(x^k)\|^2\right] + \frac{\gamma\sigma^2}{2nb},\label{eq:non_cvx_online_technical_3}
	\end{eqnarray}
	where in the last inequality, we use $\frac{\gamma(1-p)}{2pn}\left(\omega L^2 + \frac{(1+\omega)\cL^2}{b'}\right) - \frac{1}{2\gamma} + \frac{L}{2} \le 0$ following from \eqref{eq:gamma_bound_non_cvx_online_appendix}. Summing up inequalities \eqref{eq:non_cvx_online_technical_3} for $k=0,1,\ldots,K-1$ and rearranging the terms, we derive
	\begin{eqnarray}
		\frac{1}{K}\sum\limits_{k=0}^{K-1}\EE\left[\|\nabla f(x^k)\|^2\right] &\le& \frac{2}{\gamma K}\sum\limits_{k=0}^{K-1}\left(\EE[\Phi_k]-\EE[\Phi_{k+1}]\right) + \frac{\sigma^2}{nb}\notag\\
		&=& \frac{2\left(\EE[\Phi_0]-\EE[\Phi_{K}]\right)}{\gamma K} + \frac{\sigma^2}{nb} = \frac{2\Delta_0}{\gamma K} + \frac{\sigma^2}{nb},\notag
	\end{eqnarray}
	since $g^0 = \nabla f(x^0)$ and $\Phi_{k+1} \ge 0$. Finally, using the tower property \eqref{eq:tower_property} and the definition of $\hat x^K$, we obtain \eqref{eq:main_res_non_cvx_online_appendix} that implies \eqref{eq:main_res_2_non_cvx_online_appendix}, \eqref{eq:main_res_3_non_cvx_online_appendix}, and \eqref{eq:main_res_4_non_cvx_online_appendix}.
\end{proof}

\begin{remark}[About batchsizes dissimilarity]
	Similarly to the finite sum case, our analysis can be easily extended to handle the version of \algname{VR-MARINA} with different batchsizes $b_1,\ldots,b_n$ and $b'_1,\ldots,b'_n$ on different workers, i.e., when $|I_{i,k}| = b_i$, $|I_{i,k}'| = b_i'$ for $i\in[n]$. In this case, the statement of Theorem~\ref{thm:main_result_non_cvx_online} remains the same with the small modificiation: instead of $\frac{\cL^2}{b'}$ the complexity bounds will have $\frac{1}{n}\sum_{i=1}^n\frac{\cL_i^2}{b_i'}$, and instead of the requirement $b = \Theta\left(\frac{\sigma^2}{n\varepsilon}\right)$ it will have $\frac{1}{n^2}\sum_{i=1}^n\frac{\sigma_i^2}{b_i} = \Theta(\varepsilon^2)$.
\end{remark}

\begin{corollary}[Corollary~\ref{cor:main_result_non_cvx_online}]\label{cor:main_result_non_cvx_online_appendix}
	Let the assumptions of Theorem~\ref{thm:main_result_non_cvx_online} hold and $p = \min\left\{\frac{\zeta_{\cQ}}{d},\frac{b'}{b+b'}\right\}$, where $b'\le b$, $b = \Theta\left(\nicefrac{\sigma^2}{(n\varepsilon^2)}\right)$ and $\zeta_{\cQ}$ is the expected density of the quantization (see Def.~\ref{def:quantization}). If 
	\begin{equation*}
		\gamma \le \frac{1}{L + \sqrt{\frac{\max\left\{\nicefrac{d}{\zeta_{\cQ}} - 1,\nicefrac{b}{b'}\right\}}{n}\left(\omega L^2 + \frac{(1+\omega)\cL^2}{b'}\right)}},
	\end{equation*}
	then \algname{VR-MARINA} requires 
	\begin{eqnarray*}
		\cO\left(\frac{\Delta_0}{\varepsilon^2}\left(L\left(1 + \sqrt{\frac{\omega}{n}\max\left\{\frac{d}{\zeta_{\cQ}} - 1,\frac{\sigma^2}{nb'\varepsilon^2}\right\}}\right) + \cL\sqrt{\frac{(1+\omega)}{nb'}\max\left\{\frac{d}{\zeta_{\cQ}} - 1,\frac{\sigma^2}{nb'\varepsilon^2}\right\}}\right)\right)
	\end{eqnarray*}
	iterations/communication rounds and 
	\begin{eqnarray*}
		\cO\Bigg(\frac{\sigma^2}{n\varepsilon^2}+\frac{\Delta_0 L b'}{\varepsilon^2}+ \frac{\Delta_0 L}{\varepsilon^2}\sqrt{\frac{\omega b'}{n}\max\left\{\left(\frac{d}{\zeta_{\cQ}} - 1\right)b',\frac{\sigma^2}{n\varepsilon^2}\right\}}&\\
		&\hspace{-2cm}+ \frac{\Delta_0\cL}{\varepsilon^2}\sqrt{\frac{1+\omega}{n}\max\left\{\left(\frac{d}{\zeta_{\cQ}} - 1\right) b',\frac{\sigma^2}{n\varepsilon^2}\right\}}\Bigg)
	\end{eqnarray*}
	stochastic oracle calls per node in expectation to achieve $\EE[\|\nabla f(\hat x^K)\|^2] \le \varepsilon^2$, and the expected total communication cost per worker is
	\begin{eqnarray*}
		\cO\left(d+\frac{\Delta_0\zeta_{\cQ}}{\varepsilon^2}\left(L\left(1 + \sqrt{\frac{\omega}{n}\max\left\{\frac{d}{\zeta_{\cQ}} - 1,\frac{\sigma^2}{nb'\varepsilon^2}\right\}}\right)+ \cL\sqrt{\frac{1+\omega}{nb'}\max\left\{\frac{d}{\zeta_{\cQ}} - 1,\frac{\sigma^2}{nb'\varepsilon^2}\right\}}\right)\right)
	\end{eqnarray*}
	 under an assumption that the communication cost is proportional to the number of non-zero components of transmitted vectors from workers to the server.
\end{corollary}
\begin{proof}[Proof of Corollary~\ref{cor:main_result_non_cvx_finite_sums}]
	The choice of $p = \min\left\{\frac{\zeta_{\cQ}}{d},\frac{b'}{b+b'}\right\}$ implies
	\begin{eqnarray*}
		\frac{1-p}{p} &=& \max\left\{\frac{d}{\zeta_{\cQ}}-1,\frac{b}{b'}\right\},\\
		pm + (1-p)b' &\le& \frac{2mb'}{m+b'} \le 2b',\\
		pd + (1-p)\zeta_{\cQ} &\le& \frac{\zeta_{\cQ}}{d}\cdot d + \left(1 - \frac{\zeta_{\cQ}}{d}\right)\cdot\zeta_{\cQ} \le 2\zeta_{\cQ}.
	\end{eqnarray*}	
	Plugging these relations in \eqref{eq:gamma_bound_non_cvx_online_appendix}, \eqref{eq:main_res_2_non_cvx_online_appendix}, \eqref{eq:main_res_3_non_cvx_online_appendix} and \eqref{eq:main_res_4_non_cvx_online_appendix} and using $\sqrt{a+b} \le \sqrt{a} + \sqrt{b}$, we get that if
	\begin{equation*}
		\gamma \le \frac{1}{L + \sqrt{\frac{\max\left\{\nicefrac{d}{\zeta_{\cQ}} - 1,\nicefrac{b}{b'}\right\}}{n}\left(\omega L^2 + \frac{(1+\omega)\cL^2}{b'}\right)}},
	\end{equation*}
	then \algname{VR-MARINA} requires 
	\begin{eqnarray*}
		K &=& \cO\left(\frac{\Delta_0}{\varepsilon^2}\left(L + \sqrt{\frac{1-p}{pn}\left(\omega L^2 + \frac{(1+\omega)\cL^2}{b'}\right)}\right)\right)\\
		&=& \cO\left(\frac{\Delta_0}{\varepsilon^2}\left(L + \sqrt{L^2\frac{\omega\max\left\{\nicefrac{d}{\zeta_{\cQ}} - 1,\nicefrac{b}{b'}\right\}}{n} + \cL^2\frac{(1+\omega)\max\left\{\nicefrac{d}{\zeta_{\cQ}} - 1,\nicefrac{b}{b'}\right\}}{nb'}}\right)\right)\\
		&=& \cO\left(\frac{\Delta_0}{\varepsilon^2}\left(L\left(1 + \sqrt{\frac{\omega}{n}\max\left\{\frac{d}{\zeta_{\cQ}} - 1,\frac{\sigma^2}{nb'\varepsilon^2}\right\}}\right) + \cL\sqrt{\frac{(1+\omega)}{nb'}\max\left\{\frac{d}{\zeta_{\cQ}} - 1,\frac{\sigma^2}{nb'\varepsilon^2}\right\}}\right)\right)
	\end{eqnarray*}
	iterations/communication rounds and 
	\begin{eqnarray*}
		b + K(pb + 2(1-p)b') &=& \cO\left(b + \frac{\Delta_0}{\varepsilon^2}\left(L + \sqrt{\frac{1-p}{pn}\left(\omega L^2 + \frac{(1+\omega)\cL^2}{b'}\right)}\right)(pb + (1-p)b')\right)\\
		&=& \cO\Bigg(b+\frac{\Delta_0}{\varepsilon^2}\Bigg(L\Bigg(1 + \sqrt{\frac{\omega\max\left\{\nicefrac{d}{\zeta_{\cQ}} - 1,\nicefrac{b}{b'}\right\}}{n}}\Bigg)\\
		&&\hspace{4cm} + \cL\sqrt{\frac{(1+\omega)\max\left\{\nicefrac{d}{\zeta_{\cQ}} - 1,\nicefrac{b}{b'}\right\}}{nb'}}\Bigg)b'\Bigg)\\
		&=& \cO\Bigg(\frac{\sigma^2}{n\varepsilon^2}+\frac{\Delta_0}{\varepsilon^2}\Bigg(L\Bigg(b' + \sqrt{\frac{\omega b'}{n}\max\left\{\left(\frac{d}{\zeta_{\cQ}} - 1\right)b', \frac{\sigma^2}{n\varepsilon^2}\right\}}\Bigg) \\
		&&\hspace{4cm}+ \cL\sqrt{\frac{1+\omega}{n}\max\left\{\left(\frac{d}{\zeta_{\cQ}} - 1\right)b',\frac{\sigma^2}{n\varepsilon^2}\right\}}\Bigg)\Bigg)
	\end{eqnarray*}
	stochastic oracle calls per node in expectation to achieve $\EE[\|\nabla f(\hat x^K)\|^2] \le \varepsilon^2$, and the expected total communication cost per worker is
	\begin{eqnarray*}
		d + K(pd + (1-p)\zeta_{\cQ}) &=&  \cO\left(d+\frac{\Delta_0}{\varepsilon^2}\left(L + \sqrt{\frac{1-p}{pn}\left(\omega L^2 + \frac{(1+\omega)\cL^2}{b'}\right)}\right)(pd + (1-p)\zeta_{\cQ})\right)\\
		&=&\cO\Bigg(d+\frac{\Delta_0\zeta_{\cQ}}{\varepsilon^2}\Bigg(L\Bigg(1 + \sqrt{\frac{\omega}{n}\max\left\{\frac{d}{\zeta_{\cQ}} - 1,\frac{\sigma^2}{nb'\varepsilon^2}\right\}}\Bigg)\\
		&&\hspace{4cm}+ \cL\sqrt{\frac{1+\omega}{nb'}\max\left\{\frac{d}{\zeta_{\cQ}} - 1,\frac{\sigma^2}{nb'\varepsilon^2}\right\}}\Bigg)\Bigg)
	\end{eqnarray*}
	 under an assumption that the communication cost is proportional to the number of non-zero components of transmitted vectors from workers to the server.
\end{proof}

\subsubsection{Convergence Results Under Polyak-{\L}ojasiewicz condition}\label{sec:proof_of_thm_pl_online}
In this section, we provide an analysis of \algname{VR-MARINA} under Polyak-{\L}ojasiewicz condition in the online case.
\begin{theorem}\label{thm:main_result_pl_online_appendix}
	Consider the finite sum case \eqref{eq:main_problem_marina}+\eqref{eq:f_i_expectation}. Let Assumptions~\ref{as:lower_bound},~\ref{as:L_smoothness},~\ref{as:avg_smoothness_online}, \ref{as:pl_condition}~and~\ref{as:bounded_var} be satisfied and 
	\begin{equation}
		\gamma \le \min\left\{\frac{1}{L + \sqrt{\frac{2(1-p)}{pn}\left(\omega L^2 + \frac{(1+\omega)\cL^2}{b'}\right)}},\frac{p}{2\mu}\right\},\label{eq:gamma_bound_pl_online_appendix}
	\end{equation}
	where $L^2 = \frac{1}{n}\sum_{i=1}^nL_i^2$ and $\cL^2 = \frac{1}{n}\sum_{i=1}^n\cL_i^2$. Then after $K$ iterations of \algname{VR-MARINA}, we have
	\begin{equation}
		\EE\left[f(x^K) - f(x^*)\right] \le (1-\gamma\mu)^K\Delta_0 + \frac{\sigma^2}{nb\mu}, \label{eq:main_res_pl_online_appendix}
	\end{equation}
	where $\Delta_0 = f(x^0)-f(x^*)$. That is, after
	\begin{equation}
		K = \cO\left(\max\left\{\frac{1}{p}, \frac{L + \sqrt{\frac{1-p}{pn}\left(\omega L^2 + \frac{(1+\omega)\cL^2}{b'}\right)}}{\mu}\right\}\log\frac{\Delta_0}{\varepsilon}\right) \label{eq:main_res_2_pl_online_appendix}
	\end{equation}
	iterations with $b = \Theta\left(\frac{\sigma^2}{n\mu\varepsilon}\right)$ \algname{VR-MARINA} produces such a point $x^K$ that $\EE\left[f(x^K) - f(x^*)\right] \le \varepsilon$, and the expected total number of stochastic oracle calls per node $b + K(pb + 2(1-p)b')$ equals
	\begin{equation}
		\cO\left(m + \max\left\{\frac{1}{p}, \frac{L + \sqrt{\frac{1-p}{pn}\left(\omega L^2 + \frac{(1+\omega)\cL^2}{b'}\right)}}{\mu}\right\}(pb + (1-p)b')\log\frac{\Delta_0}{\varepsilon}\right). \label{eq:main_res_3_pl_online_appendix}
	\end{equation}
	Moreover, under an assumption that the communication cost is proportional to the number of non-zero components of transmitted vectors from workers to the server, we have that the expected total communication cost per worker $d + K(pd + (1-p)\zeta_{\cQ})$ equals
	\begin{equation}
		\cO\left(d+\max\left\{\frac{1}{p}, \frac{L + \sqrt{\frac{1-p}{pn}\left(\omega L^2 + \frac{(1+\omega)\cL^2}{b'}\right)}}{\mu}\right\}(pd + (1-p)\zeta_{\cQ})\log\frac{\Delta_0}{\varepsilon}\right),\label{eq:main_res_4_pl_online_appendix}
	\end{equation}
	where $\zeta_{\cQ}$ is the expected density of the quantization (see Def.~\ref{def:quantization}).
\end{theorem}
\begin{proof}
	The proof is very similar to the proof of Theorem~\ref{thm:main_result_non_cvx_online}. From Lemma~\ref{lem:lemma_2_page} and P{\L} condition, we have
	\begin{eqnarray}
		\EE[f(x^{k+1}) - f(x^*)] &\le& \EE[f(x^k) - f(x^*)] - \frac{\gamma}{2}\EE\left[\|\nabla f(x^k)\|^2\right] - \left(\frac{1}{2\gamma} - \frac{L}{2}\right)\EE\left[\|x^{k+1}-x^k\|^2\right]\notag\\
		&&\quad + \frac{\gamma}{2}\EE\left[\|g^k - \nabla f(x^k)\|^2\right]\notag\\
		&\overset{\eqref{eq:pl_condition}}{\le}& (1-\gamma\mu)\EE\left[f(x^k) - f(x^*)\right] - \left(\frac{1}{2\gamma} - \frac{L}{2}\right)\EE\left[\|x^{k+1}-x^k\|^2\right]\notag\\
		&&\quad+ \frac{\gamma}{2}\EE\left[\|g^k - \nabla f(x^k)\|^2\right]. \notag
	\end{eqnarray}
	Using the same arguments as in the proof of \eqref{eq:non_cvx_online_technical_2}, we obtain
	\begin{eqnarray}
		\EE\left[\|g^{k+1}-\nabla f(x^{k+1})\|^2\right] &\le&\frac{1-p}{n}\left(\omega L^2 + \frac{(1+\omega)\cL^2}{b'}\right)\EE\left[\|x^{k+1}-x^k\|^2\right]\notag\\
		&&+ (1-p)\EE\left[\left\|g^k - \nabla f(x^k)\right\|^2\right] + \frac{p\sigma^2}{nb}.
	\end{eqnarray}
	Putting all together, we derive that the sequence $\Phi_k = f(x^k) - f(x^*) + \frac{\gamma}{p}\|g^k - \nabla f(x^k)\|^2$ satisfies
	\begin{eqnarray}
		\EE\left[\Phi_{k+1}\right] &\le& \EE\left[(1-\gamma\mu)(f(x^k) - f(x^*)) - \left(\frac{1}{2\gamma} - \frac{L}{2}\right)\|x^{k+1}-x^k\|^2 + \frac{\gamma}{2}\|g^k - \nabla f(x^k)\|^2\right]\notag\\
		&&\quad + \frac{\gamma}{p}\EE\left[\frac{1-p}{n}\left(\omega L^2 + \frac{(1+\omega)\cL^2}{b'}\right)\|x^{k+1}-x^k\|^2\right] \notag\\
		&&\quad + \frac{\gamma}{p}\EE\left[(1-p)\left\|g^k - \nabla f(x^k)\right\|^2 + \frac{p\sigma^2}{nb}\right] \notag\\
		&=& \EE\left[(1-\gamma\mu)(f(x^k) - f(x^*)) + \left(\frac{\gamma}{2} + \frac{\gamma}{p}(1-p)\right)\left\|g^k - \nabla f(x^k)\right\|^2\right] + \frac{\gamma\sigma^2}{nb}\notag\\
		&&\quad + \left(\frac{\gamma(1-p)}{pn}\left(\omega L^2 + \frac{(1+\omega)\cL^2}{b'}\right) - \frac{1}{2\gamma} + \frac{L}{2}\right)\EE\left[\|x^{k+1}-x^k\|^2\right]\notag\\
		&\overset{\eqref{eq:gamma_bound_pl_finite_sums_appendix}}{\le}& (1-\gamma\mu)\EE[\Phi_k] + \frac{\gamma\sigma^2}{nb},\notag
	\end{eqnarray}
	where in the last inequality we use $\frac{\gamma(1-p)}{pn}\left(\omega L^2 + \frac{(1+\omega)\cL^2}{b'}\right) - \frac{1}{2\gamma} + \frac{L}{2} \le 0$ and $\frac{\gamma}{2} + \frac{\gamma}{p}(1-p) \le (1-\gamma\mu)\frac{\gamma}{p}$ following from \eqref{eq:gamma_bound_pl_online_appendix}. Unrolling the recurrence and using $g^0 = \nabla f(x^0)$, we obtain 
	\begin{eqnarray*}
		\EE\left[f(x^{K}) - f(x^*)\right] \le \EE[\Phi_{K}] &\le& (1-\gamma\mu)^{K}\Phi_0 + \frac{\gamma\sigma^2}{nb}\sum\limits_{k=0}^{K-1}(1-\gamma\mu)^k \\
		&\le& (1-\gamma\mu)^{K}(f(x^0) - f(x^*)) + \frac{\gamma\sigma^2}{nb}\sum\limits_{k=0}^{\infty}(1-\gamma\mu)^k\\
		&\le& (1-\gamma\mu)^{K}(f(x^0) - f(x^*)) + \frac{\sigma^2}{nb\mu}.
	\end{eqnarray*}
	Together with $b = \Theta\left(\frac{\sigma^2}{n\mu\varepsilon}\right)$ it implies \eqref{eq:main_res_2_pl_online_appendix}, \eqref{eq:main_res_3_pl_online_appendix}, and \eqref{eq:main_res_4_pl_online_appendix}.
\end{proof}

\begin{corollary}\label{cor:main_result_pl_online_appendix}
	Let the assumptions of Theorem~\ref{thm:main_result_pl_online_appendix} hold and $p = \min\left\{\frac{\zeta_{\cQ}}{d},\frac{b'}{b+b'}\right\}$, where $b' \le b$ and $\zeta_{\cQ}$ is the expected density of the quantization (see Def.~\ref{def:quantization}). If 
	\begin{equation*}
		\gamma \le \min\left\{\frac{1}{L + \sqrt{\frac{2\max\left\{\nicefrac{d}{\zeta_{\cQ}} - 1,\nicefrac{b}{b'}\right\}}{n}\left(\omega L^2 + \frac{(1+\omega)\cL^2}{b'}\right)}},\frac{p}{2\mu}\right\}
	\end{equation*}
	and
	\begin{equation*}
		b = \Theta\left(\frac{\sigma^2}{n\mu\varepsilon}\right),\quad \sigma^2 = \frac{1}{n}\sum\limits_{i=1}^n\sigma_i^2,
	\end{equation*}
	then \algname{VR-MARINA} requires 
	\begin{equation*}
		\cO\left(\max\left\{\frac{1}{p}, \frac{L}{\mu}\left(1 + \sqrt{\frac{\omega}{n}\max\left\{\frac{d}{\zeta_{\cQ}} - 1,\frac{\sigma^2}{nb'\mu}\right\}}\right) + \frac{\cL}{\mu}\sqrt{\frac{1+\omega}{nb'}\max\left\{\frac{d}{\zeta_{\cQ}} - 1,\frac{\sigma^2}{nb'\mu}\right\}}\right\}\log\frac{\Delta_0}{\varepsilon}\right)
	\end{equation*}
	iterations/communication rounds, 
	\begin{eqnarray*}
		\cO\Bigg(\frac{\sigma^2}{n\mu\varepsilon}+\max\Bigg\{\frac{b'}{p}, \frac{L}{\mu}\Bigg(b' + \sqrt{\frac{\omega b'}{n}\max\left\{\left(\frac{d}{\zeta_{\cQ}} - 1\right)b',\frac{\sigma^2}{n\mu\varepsilon}\right\}}\Bigg)&\\
		&\hspace{-3cm} + \frac{\cL}{\mu}\sqrt{\frac{1+\omega}{n}\max\left\{\left(\frac{d}{\zeta_{\cQ}} - 1\right)b',\frac{\sigma^2}{n\mu\varepsilon}\right\}}\Bigg\}\log\frac{\Delta_0}{\varepsilon}\Bigg)
	\end{eqnarray*}
	stochastic oracle calls per node in expectation to achieve $\EE[f(x^K)-f(x^*)] \le \varepsilon$, and the expected total communication cost per worker is
	\begin{eqnarray*}
	    \cO\Bigg(d+\zeta_{\cQ}\max\Bigg\{\frac{1}{p}, \frac{L}{\mu}\left(1 + \sqrt{\frac{\omega}{n}\max\left\{\frac{d}{\zeta_{\cQ}} - 1,\frac{\sigma^2}{nb'\mu}\right\}}\right)&\\
	    &\hspace{-2cm}+ \frac{\cL}{\mu}\sqrt{\frac{1+\omega}{nb'}\max\left\{\frac{d}{\zeta_{\cQ}} - 1,\frac{\sigma^2}{nb'\mu}\right\}}\Bigg\}\log\frac{\Delta_0}{\varepsilon}\Bigg)
	\end{eqnarray*}
	 under an assumption that the communication cost is proportional to the number of non-zero components of transmitted vectors from workers to the server.
\end{corollary}
\begin{proof}
	The choice of $p = \min\left\{\frac{\zeta_{\cQ}}{d},\frac{b'}{b+b'}\right\}$ implies
	\begin{eqnarray*}
		\frac{1-p}{p} &=& \max\left\{\frac{d}{\zeta_{\cQ}}-1,\frac{b}{b'}\right\},\\
		pm + (1-p)b' &\le& \frac{2bb'}{b+b'} \le 2b',\\
		pd + (1-p)\zeta_{\cQ} &\le& \frac{\zeta_{\cQ}}{d}\cdot d + \left(1 - \frac{\zeta_{\cQ}}{d}\right)\cdot\zeta_{\cQ} \le 2\zeta_{\cQ}.
	\end{eqnarray*}	
	Plugging these relations in \eqref{eq:gamma_bound_pl_online_appendix}, \eqref{eq:main_res_2_pl_online_appendix}, \eqref{eq:main_res_3_pl_online_appendix} and \eqref{eq:main_res_4_pl_online_appendix} and using $\sqrt{a+b} \le \sqrt{a} + \sqrt{b}$, we get that if
	\begin{equation*}
		\gamma \le \min\left\{\frac{1}{L + \sqrt{\frac{2\max\left\{\nicefrac{d}{\zeta_{\cQ}} - 1,\nicefrac{b}{b'}\right\}}{n}\left(\omega L^2 + \frac{(1+\omega)\cL^2}{b'}\right)}},\frac{p}{2\mu}\right\},
	\end{equation*}
	then \algname{VR-MARINA} requires 
	\begin{eqnarray*}
		K &=& \cO\left(\max\left\{\frac{1}{p}, \frac{L + \sqrt{\frac{1-p}{pn}\left(\omega L^2 + \frac{(1+\omega)\cL^2}{b'}\right)}}{\mu}\right\}\log\frac{\Delta_0}{\varepsilon}\right)\\
		&=& \cO\left(\max\left\{\frac{1}{p}, \frac{L + \sqrt{L^2\frac{\omega\max\left\{\nicefrac{d}{\zeta_{\cQ}} - 1,\nicefrac{b}{b'}\right\}}{n} + \cL^2\frac{(1+\omega)\max\left\{\nicefrac{d}{\zeta_{\cQ}} - 1,\nicefrac{b}{b'}\right\}}{nb'}}}{\mu}\right\}\log\frac{\Delta_0}{\varepsilon}\right)\\
		&=& \cO\left(\max\left\{\frac{1}{p}, \frac{L}{\mu}\left(1 + \sqrt{\frac{\omega}{n}\max\left\{\frac{d}{\zeta_{\cQ}} - 1,\frac{\sigma^2}{nb'\mu}\right\}}\right) + \frac{\cL}{\mu}\sqrt{\frac{1+\omega}{nb'}\max\left\{\frac{d}{\zeta_{\cQ}} - 1,\frac{\sigma^2}{nb'\mu}\right\}}\right\}\log\frac{\Delta_0}{\varepsilon}\right)
	\end{eqnarray*}
	iterations/communication rounds and 
	\begin{eqnarray*}
		b + K(pb + 2(1-p)b') &=& \cO\left(b + \max\left\{\frac{1}{p}, \frac{L + \sqrt{\frac{1-p}{pn}\left(\omega L^2 + \frac{(1+\omega)\cL^2}{b'}\right)}}{\mu}\right\}(pb + (1-p)b')\log\frac{\Delta_0}{\varepsilon}\right)\\	
		&=& \cO\Bigg(b+\max\Bigg\{\frac{1}{p}, \frac{L}{\mu}\Bigg(1 + \sqrt{\frac{\omega\max\left\{\nicefrac{d}{\zeta_{\cQ}} - 1,\nicefrac{b}{b'}\right\}}{n}}\Bigg)\\
		&&\hspace{3cm} + \frac{\cL}{\mu}\sqrt{\frac{(1+\omega)\max\left\{\nicefrac{d}{\zeta_{\cQ}} - 1,\nicefrac{b}{b'}\right\}}{nb'}}\Bigg\}b'\log\frac{\Delta_0}{\varepsilon}\Bigg)\\
		&=& \cO\Bigg(\frac{\sigma^2}{n\mu\varepsilon}+\max\Bigg\{\frac{b'}{p}, \frac{L}{\mu}\Bigg(b' + \sqrt{\frac{\omega b'}{n}\max\left\{\left(\frac{d}{\zeta_{\cQ}} - 1\right)b',\frac{\sigma^2}{n\mu\varepsilon}\right\}}\Bigg)\\
		&&\hspace{3cm} + \frac{\cL}{\mu}\sqrt{\frac{1+\omega}{n}\max\left\{\left(\frac{d}{\zeta_{\cQ}} - 1\right)b',\frac{\sigma^2}{n\mu\varepsilon}\right\}}\Bigg\}\log\frac{\Delta_0}{\varepsilon}\Bigg)
	\end{eqnarray*}
	stochastic oracle calls per node in expectation to achieve $\EE[f(x^K) - f(x^*)] \le \varepsilon$, and the expected total communication cost per worker is
	\begin{eqnarray*}
		d + K(pd + (1-p)\zeta_{\cQ}) &=&  \cO\left(d+\max\left\{\frac{1}{p}, \frac{L + \sqrt{\frac{1-p}{pn}\left(\omega L^2 + \frac{(1+\omega)\cL^2}{b'}\right)}}{\mu}\right\}(pd + (1-p)\zeta_{\cQ})\log\frac{\Delta_0}{\varepsilon}\right)\\
		&=&\cO\Bigg(d+\zeta_{\cQ}\max\Bigg\{\frac{1}{p}, \frac{L}{\mu}\left(1 + \sqrt{\frac{\omega}{n}\max\left\{\frac{d}{\zeta_{\cQ}} - 1,\frac{\sigma^2}{nb'\mu}\right\}}\right)\\
		&&\hspace{3cm} + \frac{\cL}{\mu}\sqrt{\frac{1+\omega}{nb'}\max\left\{\frac{d}{\zeta_{\cQ}} - 1,\frac{\sigma^2}{nb'\mu}\right\}}\Bigg\}\log\frac{\Delta_0}{\varepsilon}\Bigg)
	\end{eqnarray*}
	 under an assumption that the communication cost is proportional to the number of non-zero components of transmitted vectors from workers to the server.
\end{proof}


\section{Missing Proofs for \algname{PP-MARINA}}\label{sec:pp_marina_proofs}

\subsection{Generally Non-Convex Problems}\label{sec:proof_of_thm_non_cvx_pp}
In this section, we provide the full statement of Theorem~\ref{thm:main_result_non_cvx_pp} together with the proof of this result.
\begin{theorem}[Theorem~\ref{thm:main_result_non_cvx_pp}]\label{thm:main_result_non_cvx_pp_appendix}
	Let Assumptions~\ref{as:lower_bound}~and~\ref{as:L_smoothness} be satisfied and 
	\begin{equation}
		\gamma \le \frac{1}{L\left(1 + \sqrt{\frac{(1-p)(1+\omega)}{pr}}\right)},\label{eq:gamma_bound_non_cvx_pp_appendix}
	\end{equation}
	where $L^2 = \frac{1}{n}\sum_{i=1}^nL_i^2$. Then after $K$ iterations of \algname{PP-MARINA}, we have
	\begin{equation}
		\EE\left[\left\|\nabla f(\hat x^K)\right\|^2\right] \le \frac{2\Delta_0}{\gamma K}, \label{eq:main_res_non_cvx_pp_appendix}
	\end{equation}
	where $\hat{x}^K$ is chosen uniformly at random from $x^0,\ldots,x^{K-1}$ and $\Delta_0 = f(x^0)-f_*$. That is, after
	\begin{equation}
		K = \cO\left(\frac{\Delta_0 L}{\varepsilon^2}\left(1 + \sqrt{\frac{(1-p)(1+\omega)}{pr}}\right)\right) \label{eq:main_res_2_non_cvx_pp_appendix}
	\end{equation}
	iterations \algname{PP-MARINA} produces such a point $\hat x^K$ that $\EE[\|\nabla f(\hat x^K)\|^2] \le \varepsilon^2$.
	Moreover, under an assumption that the communication cost is proportional to the number of non-zero components of transmitted vectors from workers to the server, we have that the expected total communication cost (for all workers) equals
	\begin{equation}
		dn + K(pdn + (1-p)\zeta_{\cQ}r) =  \cO\left(dn+\frac{\Delta_0 L}{\varepsilon^2}\left(1 + \sqrt{\frac{(1-p)(1+\omega)}{pr}}\right)(pdn + (1-p)\zeta_{\cQ}r)\right),\label{eq:main_res_4_non_cvx_pp_appendix}
	\end{equation}
	where $\zeta_{\cQ}$ is the expected density of the quantization (see Def.~\ref{def:quantization}).
\end{theorem}
\begin{proof}[Proof of Theorem~\ref{thm:main_result_non_cvx_pp}]
	The proof is very similar to the proof of Theorem~\ref{thm:main_result_non_cvx_finite_sums}. From Lemma~\ref{lem:lemma_2_page}, we have
	\begin{equation}
		\EE[f(x^{k+1})] \le \EE[f(x^k)] - \frac{\gamma}{2}\EE\left[\|\nabla f(x^k)\|^2\right] - \left(\frac{1}{2\gamma} - \frac{L}{2}\right)\EE\left[\|x^{k+1}-x^k\|^2\right] + \frac{\gamma}{2}\EE\left[\|g^k - \nabla f(x^k)\|^2\right]. \label{eq:non_cvx_pp_technical_1}
	\end{equation}
	Next, we need to derive an upper bound for $\EE\left[\|g^{k+1}-\nabla f(x^{k+1})\|^2\right]$. By definition of $g^{k+1}$, we have
	\begin{equation}
		g^{k+1} = \begin{cases}\nabla f(x^{k+1})& \text{with probability } p,\\ g^k + \frac{1}{r}\sum\limits_{i_k\in I'_{k}}\cQ\left(\nabla f_{i_k}(x^{k+1}) - \nabla f_{i_k}(x^k)\right)& \text{with probability } 1-p. \end{cases}\notag
	\end{equation}
	Using this, variance decomposition \eqref{eq:variance_decomposition} and tower property \eqref{eq:tower_property}, we derive:
	\begin{eqnarray}
		\EE\left[\|g^{k+1}-\nabla f(x^{k+1})\|^2\right]&\notag\\ &\hspace{-4.3cm}\overset{\eqref{eq:tower_property}}{=} (1-p)\EE\left[\left\|g^k + \frac{1}{r}\sum\limits_{i_k\in I_k'} \cQ\left(\nabla f_{i_k}(x^{k+1}) - \nabla f_{i_k}(x^k)\right) - \nabla f(x^{k+1})\right\|^2\right]\notag\\
		&\hspace{-3.5cm}\overset{\eqref{eq:tower_property},\eqref{eq:variance_decomposition}}{=} (1-p)\EE\left[\left\|\frac{1}{r}\sum\limits_{i_k\in I_k'} \cQ\left(\nabla f_{i_k}(x^{k+1}) - \nabla f_{i_k}(x^k)\right) - \nabla f(x^{k+1}) + \nabla f(x^k)\right\|^2\right]\notag\\
		&\hspace{-5cm}+ (1-p)\EE\left[\left\|g^k - \nabla f(x^k)\right\|^2\right].\notag
	\end{eqnarray}
	Next, we use the notation: $\Delta_i^k = \nabla f_{i}(x^{k+1}) - \nabla f_{i}(x^k)$ for $i\in [n]$ and $\Delta^k = \nabla f(x^{k+1}) - \nabla f(x^k)$. These vectors satisfy $\EE\left[\Delta_{i_k}^k \mid x^k,x^{k+1}\right] = \Delta^k$ for all $i_k\in I_{k}'$. Moreover, $\cQ(\Delta_{i_k}^k)$ for $i_k\in I_k'$ are independent random vectors for fixed $x^k$ and $x^{k+1}$. These observations imply
	\begin{eqnarray}
		\EE\left[\|g^{k+1}-\nabla f(x^{k+1})\|^2\right] &=& (1-p)\EE\left[\left\|\frac{1}{r}\sum\limits_{i_k\in I_k'} \left(\cQ(\Delta_{i_k}^k) - \Delta^k\right)\right\|^2\right]\notag\\
		&&\quad +(1-p)\EE\left[\left\|g^k - \nabla f(x^k)\right\|^2\right]\notag\\
		&=& \frac{1-p}{r}\EE\left[\left\|\cQ(\Delta_{i_k}^k) - \Delta_{i_k}^k + \Delta_{i_k}^k - \Delta^k\right\|^2\right]\notag\\
		&&\quad + (1-p)\EE\left[\left\|g^k - \nabla f(x^k)\right\|^2\right]\notag\\
		&\overset{\eqref{eq:tower_property},\eqref{eq:variance_decomposition}}{=}& \frac{1-p}{r}\left(\EE\left[\left\|\cQ(\Delta_{i_k}^k) - \Delta_{i_k}^k\right\|^2\right] + \EE\left[\left\|\Delta_{i_k}^k - \Delta^k\right\|^2\right]\right)\notag\\
		&&\quad + (1-p)\EE\left[\left\|g^k - \nabla f(x^k)\right\|^2\right]\notag\\
		&\overset{\eqref{eq:tower_property},\eqref{eq:quantization_def}}{=}& \frac{1-p}{r}\left(\omega\EE\left[\left\|\Delta_{i_k}^k\right\|^2\right] + \EE\left[\left\|\Delta_{i_k}^k - \Delta^k\right\|^2\right]\right)\notag\\
		&&\quad + (1-p)\EE\left[\left\|g^k - \nabla f(x^k)\right\|^2\right]\notag\\
		&\overset{\eqref{eq:tower_property},\eqref{eq:variance_decomposition}}{=}& \frac{(1-p)(1+\omega)}{r}\EE\left[\left\|\Delta_{i_k}^k\right\|^2\right]+ (1-p)\EE\left[\left\|g^k - \nabla f(x^k)\right\|^2\right].\notag
	\end{eqnarray}
	Using $L$-smoothness \eqref{eq:L_smoothness_local_marina} of $f_i$ together with the tower property \eqref{eq:tower_property}, we get
	\begin{eqnarray}
		\EE\left[\|g^{k+1}-\nabla f(x^{k+1})\|^2\right] &\le& \frac{(1-p)(1+\omega)}{nr}\sum\limits_{i=1}^nL_i^2\EE\left[\|x^{k+1} - x^k\|^2\right] \notag\\
		&&\quad + (1-p)\EE\left[\left\|g^k - \nabla f(x^k)\right\|^2\right]\notag\\
		&=&\frac{(1-p)(1+\omega)L^2}{r}\EE\left[\|x^{k+1}-x^k\|^2\right]\notag\\
		&&\quad + (1-p)\EE\left[\left\|g^k - \nabla f(x^k)\right\|^2\right].\label{eq:non_cvx_pp_technical_2}
	\end{eqnarray}
	Next, we introduce new notation: $\Phi_k = f(x^k) - f_* + \frac{\gamma}{2p}\|g^k - \nabla f(x^k)\|^2$. Using this and inequalities \eqref{eq:non_cvx_pp_technical_1} and \eqref{eq:non_cvx_pp_technical_2}, we establish the following inequality:
	\begin{eqnarray}
		\EE\left[\Phi_{k+1}\right] &\le& \EE\left[f(x^k) - f_* - \frac{\gamma}{2}\|\nabla f(x^k)\|^2 - \left(\frac{1}{2\gamma} - \frac{L}{2}\right)\|x^{k+1}-x^k\|^2 + \frac{\gamma}{2}\|g^k - \nabla f(x^k)\|^2\right]\notag\\
		&&\quad + \frac{\gamma}{2p}\EE\left[\frac{(1-p)(1+\omega)L^2}{r}\|x^{k+1}-x^k\|^2 + (1-p)\left\|g^k - \nabla f(x^k)\right\|^2\right] \notag\\
		&=& \EE\left[\Phi_k\right] - \frac{\gamma}{2}\EE\left[\|\nabla f(x^k)\|^2\right] + \left(\frac{\gamma(1-p)(1+\omega)L^2}{2pn} - \frac{1}{2\gamma} + \frac{L}{2}\right)\EE\left[\|x^{k+1}-x^k\|^2\right]\notag\\
		&\overset{\eqref{eq:gamma_bound_non_cvx_pp_appendix}}{\le}& \EE\left[\Phi_k\right] - \frac{\gamma}{2}\EE\left[\|\nabla f(x^k)\|^2\right],\label{eq:non_cvx_pp_technical_3}
	\end{eqnarray}
	where in the last inequality we use $\frac{\gamma(1-p)(1+\omega)L^2}{2pn} - \frac{1}{2\gamma} + \frac{L}{2} \le 0$ following from \eqref{eq:gamma_bound_non_cvx_pp_appendix}. Summing up inequalities \eqref{eq:non_cvx_finite_sums_technical_3} for $k=0,1,\ldots,K-1$ and rearranging the terms, we derive
	\begin{eqnarray}
		\frac{1}{K}\sum\limits_{k=0}^{K-1}\EE\left[\|\nabla f(x^k)\|^2\right] &\le& \frac{2}{\gamma K}\sum\limits_{k=0}^{K-1}\left(\EE[\Phi_k]-\EE[\Phi_{k+1}]\right) = \frac{2\left(\EE[\Phi_0]-\EE[\Phi_{K}]\right)}{\gamma K} = \frac{2\Delta_0}{\gamma K},\notag
	\end{eqnarray}
	since $g^0 = \nabla f(x^0)$ and $\Phi_{k+1} \ge 0$. Finally, using the tower property \eqref{eq:tower_property} and the definition of $\hat x^K$, we obtain \eqref{eq:main_res_non_cvx_pp_appendix} that implies \eqref{eq:main_res_2_non_cvx_pp_appendix} and \eqref{eq:main_res_4_non_cvx_pp_appendix}.
\end{proof}

\begin{corollary}[Corollary~\ref{cor:main_result_non_cvx_pp}]\label{cor:main_result_non_cvx_pp_appendix}
	Let the assumptions of Theorem~\ref{thm:main_result_non_cvx_pp} hold and $p = \frac{\zeta_{\cQ}r}{dn}$, where $r \le n$ and $\zeta_{\cQ}$ is the expected density of the quantization (see Def.~\ref{def:quantization}). If 
	\begin{equation*}
		\gamma \le \frac{1}{L\left(1 + \sqrt{\frac{1+\omega}{r}\left(\frac{dn}{\zeta_{\cQ}r}-1\right)}\right)},
	\end{equation*}
	then \algname{PP-MARINA} requires 
	\begin{equation*}
		K = \cO\left(\frac{\Delta_0 L}{\varepsilon^2}\left(1 + \sqrt{\frac{1+\omega}{r}\left(\frac{dn}{\zeta_{\cQ}r}-1\right)}\right)\right)
	\end{equation*}
	iterations/communication rounds to achieve $\EE[\|\nabla f(\hat x^K)\|^2] \le \varepsilon^2$, and the expected total communication cost is
	\begin{equation*}
		\cO\left(dn+\frac{\Delta_0 L}{\varepsilon^2}\left(\zeta_{\cQ}r + \sqrt{(1+\omega)\zeta_{\cQ}\left(dn-\zeta_{\cQ}r\right)}\right)\right)
	\end{equation*}
	 under an assumption that the communication cost is proportional to the number of non-zero components of transmitted vectors from workers to the server.
\end{corollary}
\begin{proof}[Proof of Corollary~\ref{cor:main_result_non_cvx_pp}]
	The choice of $p = \frac{\zeta_{\cQ}r}{dn}$ implies
	\begin{eqnarray*}
		\frac{1-p}{p} &=& \frac{dn}{\zeta_{\cQ}r}-1,\\
		pdn + (1-p)\zeta_{\cQ}r &\le& \zeta_{\cQ}r + \left(1 - \frac{\zeta_{\cQ}r}{dn}\right)\cdot\zeta_{\cQ}r \le 2\zeta_{\cQ}r.
	\end{eqnarray*}	
	Plugging these relations in \eqref{eq:gamma_bound_non_cvx_pp_appendix}, \eqref{eq:main_res_2_non_cvx_pp_appendix}, and \eqref{eq:main_res_4_non_cvx_pp_appendix}, we get that if
	\begin{equation*}
		\gamma \le \frac{1}{L\left(1 + \sqrt{\frac{1+\omega}{r}\left(\frac{dn}{\zeta_{\cQ}r}-1\right)}\right)},
	\end{equation*}
	then \algname{PP-MARINA} requires 
	\begin{eqnarray*}
		K &=& \cO\left(\frac{\Delta_0 L}{\varepsilon^2}\left(1 + \sqrt{\frac{(1-p)(1+\omega)}{pr}}\right)\right)\\
		&=& \cO\left(\frac{\Delta_0 L}{\varepsilon^2}\left(1 + \sqrt{\frac{1+\omega}{r}\left(\frac{dn}{\zeta_{\cQ}r}-1\right)}\right)\right)
	\end{eqnarray*}
	iterations/communication rounds in order to achieve $\EE[\|\nabla f(\hat x^K)\|^2] \le \varepsilon^2$, and the expected total communication cost is
	\begin{eqnarray*}
		dn + K(pdn + (1-p)\zeta_{\cQ}r) &=&  \cO\left(dn+\frac{\Delta_0 L}{\varepsilon^2}\left(1 + \sqrt{\frac{(1-p)(1+\omega)}{pr}}\right)(pdn + (1-p)\zeta_{\cQ}r)\right)\\
		&=&\cO\left(dn+\frac{\Delta_0 L}{\varepsilon^2}\left(\zeta_{\cQ}r + \sqrt{(1+\omega)\zeta_{\cQ}\left(dn-\zeta_{\cQ}r\right)}\right)\right)
	\end{eqnarray*}
	 under an assumption that the communication cost is proportional to the number of non-zero components of transmitted vectors from workers to the server.
\end{proof}

\subsection{Convergence Results Under Polyak-{\L}ojasiewicz Condition}\label{sec:proof_of_thm_pl_pp}
In this section, we provide an analysis of \algname{PP-MARINA} under Polyak-{\L}ojasiewicz condition.
\begin{theorem}\label{thm:main_result_pl_pp_appendix}
	Let Assumptions~\ref{as:lower_bound},~\ref{as:L_smoothness}~and~\ref{as:pl_condition} be satisfied and 
	\begin{equation}
		\gamma \le \min\left\{\frac{1}{L\left(1 + \sqrt{\frac{2(1-p)(1+\omega)}{pr}}\right)}, \frac{p}{2\mu}\right\},\label{eq:gamma_bound_pl_pp_appendix}
	\end{equation}
	where $L^2 = \frac{1}{n}\sum_{i=1}^nL_i^2$. Then after $K$ iterations of \algname{PP-MARINA}, we have
	\begin{equation}
		\EE\left[f(x^K) - f(x^*)\right] \le (1-\gamma\mu)^K\Delta_0, \label{eq:main_res_pl_pp_appendix}
	\end{equation}
	where $\Delta_0 = f(x^0)-f(x^*)$. That is, after
	\begin{equation}
		K = \cO\left(\max\left\{\frac{1}{p},\frac{L}{\mu}\left(1 + \sqrt{\frac{(1-p)(1+\omega)}{pr}}\right)\right\}\log\frac{\Delta_0}{\varepsilon}\right) \label{eq:main_res_2_pl_pp_appendix}
	\end{equation}
	iterations \algname{PP-MARINA} produces such a point $x^K$ that $\EE[f(x^K) - f(x^*)] \le \varepsilon$.
	Moreover, under an assumption that the communication cost is proportional to the number of non-zero components of transmitted vectors from workers to the server, we have that the expected total communication cost (for all workers) $dn + K(pdn + (1-p)\zeta_{\cQ}r)$ equals
	\begin{equation}
		 \cO\left(dn+\max\left\{\frac{1}{p},\frac{L}{\mu}\left(1 + \sqrt{\frac{(1-p)(1+\omega)}{pr}}\right)\right\}(pdn + (1-p)\zeta_{\cQ}r)\log\frac{\Delta_0}{\varepsilon}\right),\label{eq:main_res_4_pl_pp_appendix}
	\end{equation}
	where $\zeta_{\cQ}$ is the expected density of the quantization (see Def.~\ref{def:quantization}).
\end{theorem}
\begin{proof}
	The proof is very similar to the proof of Theorem~\ref{thm:main_result_non_cvx_pp}. From Lemma~\ref{lem:lemma_2_page} and P{\L} condition we have
	\begin{eqnarray}
		\EE[f(x^{k+1}) - f(x^*)] &\le& \EE[f(x^k) - f(x^*)] - \frac{\gamma}{2}\EE\left[\|\nabla f(x^k)\|^2\right] - \left(\frac{1}{2\gamma} - \frac{L}{2}\right)\EE\left[\|x^{k+1}-x^k\|^2\right]\notag\\
		&&\quad + \frac{\gamma}{2}\EE\left[\|g^k - \nabla f(x^k)\|^2\right]\notag\\
		&\overset{\eqref{eq:pl_condition}}{\le}& (1-\gamma\mu)\EE\left[f(x^k) - f(x^*)\right] - \left(\frac{1}{2\gamma} - \frac{L}{2}\right)\EE\left[\|x^{k+1}-x^k\|^2\right]\notag\\
		&&\quad + \frac{\gamma}{2}\EE\left[\|g^k - \nabla f(x^k)\|^2\right]. \notag
	\end{eqnarray}
	Using the same arguments as in the proof of \eqref{eq:non_cvx_pp_technical_2}, we obtain
	\begin{eqnarray}
		\EE\left[\|g^{k+1}-\nabla f(x^{k+1})\|^2\right] &\le& \frac{(1-p)(1+\omega)L^2}{r}\EE\left[\|x^{k+1}-x^k\|^2\right] + (1-p)\EE\left[\left\|g^k - \nabla f(x^k)\right\|^2\right].\notag
	\end{eqnarray}
	Putting all together, we derive that the sequence $\Phi_k = f(x^k) - f(x^*) + \frac{\gamma}{p}\|g^k - \nabla f(x^k)\|^2$ satisfies
	\begin{eqnarray}
		\EE\left[\Phi_{k+1}\right] &\le& \EE\left[(1-\gamma\mu)(f(x^k) - f(x^*)) - \left(\frac{1}{2\gamma} - \frac{L}{2}\right)\|x^{k+1}-x^k\|^2 + \frac{\gamma}{2}\|g^k - \nabla f(x^k)\|^2\right]\notag\\
		&&\quad + \frac{\gamma}{p}\EE\left[\frac{(1-p)(1+\omega)L^2}{r}\|x^{k+1}-x^k\|^2 + (1-p)\left\|g^k - \nabla f(x^k)\right\|^2 \right] \notag\\
		&=& \EE\left[(1-\gamma\mu)(f(x^k) - f(x^*)) + \left(\frac{\gamma}{2} + \frac{\gamma}{p}(1-p)\right)\left\|g^k - \nabla f(x^k)\right\|^2\right]\notag\\
		&&\quad + \left(\frac{\gamma(1-p)(1+\omega)L^2}{pr} - \frac{1}{2\gamma} + \frac{L}{2}\right)\EE\left[\|x^{k+1}-x^k\|^2\right]\notag\\
		&\overset{\eqref{eq:gamma_bound_pl_pp_appendix}}{\le}& (1-\gamma\mu)\EE[\Phi_k],\notag
	\end{eqnarray}
	where in the last inequality we use $\frac{\gamma(1-p)(1+\omega)L^2}{pr} - \frac{1}{2\gamma} + \frac{L}{2} \le 0$ and $\frac{\gamma}{2} + \frac{\gamma}{p}(1-p) \le (1-\gamma\mu)\frac{\gamma}{p}$ following from \eqref{eq:gamma_bound_pl_pp_appendix}. Unrolling the recurrence and using $g^0 = \nabla f(x^0)$, we obtain 
	\begin{eqnarray*}
		\EE\left[f(x^{K}) - f(x^*)\right] \le \EE[\Phi_{K}] &\le& (1-\gamma\mu)^{K}\Phi_0 = (1-\gamma\mu)^{K}(f(x^0) - f(x^*))
	\end{eqnarray*}
	that implies \eqref{eq:main_res_2_pl_pp_appendix} and \eqref{eq:main_res_4_pl_pp_appendix}.
\end{proof}

\begin{corollary}\label{cor:main_result_pl_pp_appendix}
	Let the assumptions of Theorem~\ref{thm:main_result_pl_pp_appendix} hold and $p = \frac{\zeta_{\cQ}r}{dn}$, where $r \le n$ and $\zeta_{\cQ}$ is the expected density of the quantization (see Def.~\ref{def:quantization}). If 
	\begin{equation*}
		\gamma \le \min\left\{\frac{1}{L\left(1 + \sqrt{\frac{2(1+\omega)}{r}\left(\frac{dn}{\zeta_{\cQ}r}-1\right)}\right)}, \frac{p}{2\mu}\right\},
	\end{equation*}
	then \algname{PP-MARINA} requires 
	\begin{equation*}
		K = \cO\left(\max\left\{\frac{dn}{\zeta_{\cQ}r}\frac{L}{\mu}\left(1 + \sqrt{\frac{1+\omega}{r}\left(\frac{dn}{\zeta_{\cQ}r}-1\right)}\right)\right\}\log\frac{\Delta_0}{\varepsilon}\right)
	\end{equation*}
	iterations/communication rounds to achieve $\EE[f(x^K) - f(x^*)] \le \varepsilon$, and the expected total communication cost is
	\begin{equation*}
		\cO\left(dn+\max\left\{dn,\frac{L}{\mu}\left(\zeta_{\cQ}r + \sqrt{(1+\omega)\zeta_{\cQ}\left(dn-\zeta_{\cQ}r\right)}\right)\right\}\log\frac{\Delta_0}{\varepsilon}\right)
	\end{equation*}
	 under an assumption that the communication cost is proportional to the number of non-zero components of transmitted vectors from workers to the server.
\end{corollary}
\begin{proof}
	The choice of $p = \frac{\zeta_{\cQ}r}{dn}$ implies
	\begin{eqnarray*}
		\frac{1-p}{p} &=& \frac{dn}{\zeta_{\cQ}r}-1,\\
		pdn + (1-p)\zeta_{\cQ}r &\le& \zeta_{\cQ}r + \left(1 - \frac{\zeta_{\cQ}r}{dn}\right)\cdot\zeta_{\cQ}r \le 2\zeta_{\cQ}r.
	\end{eqnarray*}	
	Plugging these relations in \eqref{eq:gamma_bound_pl_pp_appendix}, \eqref{eq:main_res_2_pl_pp_appendix}, and \eqref{eq:main_res_4_pl_pp_appendix}, we get that if
	\begin{equation*}
		\gamma \le \min\left\{\frac{1}{L\left(1 + \sqrt{\frac{2(1+\omega)}{r}\left(\frac{dn}{\zeta_{\cQ}r}-1\right)}\right)}, \frac{p}{2\mu}\right\},
	\end{equation*}
	then \algname{PP-MARINA} requires 
	\begin{eqnarray*}
		K &=& \cO\left(\max\left\{\frac{1}{p},\frac{L}{\mu}\left(1 + \sqrt{\frac{(1-p)(1+\omega)}{pr}}\right)\right\}\log\frac{\Delta_0}{\varepsilon}\right)\\
		&=& \cO\left(\max\left\{\frac{dn}{\zeta_{\cQ}r}\frac{L}{\mu}\left(1 + \sqrt{\frac{1+\omega}{r}\left(\frac{dn}{\zeta_{\cQ}r}-1\right)}\right)\right\}\log\frac{\Delta_0}{\varepsilon}\right)
	\end{eqnarray*}
	iterations/communication rounds to achieve $\EE[f(x^K)-f(x^*)] \le \varepsilon$, and the expected total communication cost is
	\begin{eqnarray*}
		dn + K(pdn + (1-p)\zeta_{\cQ}r)& \\
		&\hspace{-2cm}=  \cO\left(dn+\max\left\{\frac{1}{p},\frac{L}{\mu}\left(1 + \sqrt{\frac{(1-p)(1+\omega)}{pr}}\right)\right\}(pdn + (1-p)\zeta_{\cQ}r)\log\frac{\Delta_0}{\varepsilon}\right)\\
		&\hspace{-3.1cm}=\cO\left(dn+\max\left\{dn,\frac{L}{\mu}\left(\zeta_{\cQ}r + \sqrt{(1+\omega)\zeta_{\cQ}\left(dn-\zeta_{\cQ}r\right)}\right)\right\}\log\frac{\Delta_0}{\varepsilon}\right)
	\end{eqnarray*}
	 under an assumption that the communication cost is proportional to the number of non-zero components of transmitted vectors from workers to the server.
\end{proof}

%% file: Appendix_moshpit.tex
\chapter{Appendix for Chapter~\ref{ch:moshpit}}\label{app:moshpit}

\section{GPU Instance Costs}
\label{sect:cloud_costs}

This section provides a brief cost analysis of typical deep learning compute resources both in the cloud and on-premises.
For brevity, we limit this analysis to the popular GPUs available at the time of submission. Note that the exact costs will depend on a variety of factors such as the cloud provider, the region, electricity costs, and market fluctuations. Therefore, we warn the reader to consider this analysis only as a rough estimate. 

Specifically, we estimate the compute costs for the occasional usage scenario: running a single set of experiments over several weeks or conducting infrequent experiments. This scenario covers most research scientists and small organizations. The most straightforward way to provision a GPU server in such a scenario is to rent it from a cloud provider (e.g., GCP or AWS) or a public marketplace (e.g., Vast.ai or Golem).

While the exact server specifications vary from one provider to another, there are two broad categories of GPU machines: regular and preemptible. Regular instance types typically offer 1--8 GPUs per node with tight uptime guarantees (typically $99.99\%$) and a high-bandwidth network (tens of Gb/s). In turn, preemptible instances provide the same resource type at a significant discount with the condition that the machine can be terminated at any time after short notice.

To account for individual variations, we report the average rent price over three popular cloud providers.
We consider three popular instance types: two high-end instances with 8 Tesla V100 or A100 GPUs and a low-end instance with a single Tesla T4 GPU.
We also describe several low-end servers and workstations available on a public marketplace. Unlike cloud VMs, these instances are hosted on non-curated hardware with less uptime guarantees (typically 95\% -- 99.9\%), slower network and significant variation in performance. However, marketplace instances are the cheapest in terms of cost per TFLOPS. To quantify this, we report the average over three most affordable instances that fit the chosen minimum requirements.

As a point of comparison, we also measure each system's training performance for BERT-Large~\cite{bert} fine-tuning on SQuAD v1.1~\cite{squad} in PyTorch with mixed precision. We follow the official benchmarking protocol by~\cite{nvidia_perf} and reuse the official performance results for V100, A100, and T4 instances. The only exception is GTX 1080Ti, where we use full 32-bit precision because that device does not support efficient half-precision operations.

\begin{table}[h]
\scriptsize
 \centering
\caption{Cloud and marketplace GPU instance pricing for short-term usage.}
\label{fig:cloud_costs}
\begin{tabular}{@{}ccccccc@{}}
\toprule
\multicolumn{4}{c}{Minimum system specifications} & \multicolumn{2}{c}{Average cost, \$/hour} & \multirow{2}[2]{*}{\shortstack{BERT-Large\\ training samples/s}} \\
\cmidrule(lr){1-4}\cmidrule(lr){5-6}
GPU & CPU cores & CPU type & RAM, GB & Regular & Preemptible &  \\ \midrule
\multicolumn{7}{c}{Cloud instances} \\ \midrule
8$\times$ V100 & 64 & Intel Xeon Broadwell & 480 & 23.47 & 7.13 & 354 \\
8$\times$  A100 & 96 & AMD Epyc ROME & 960 & 30.65 & 10.18 & 755 \\
1$\times$  T4 & 4 & Intel Xeon Cascade Lake & 16 & 0.46 & 0.18 & 18 \\ \midrule
\multicolumn{7}{c}{Marketplace instances} \\ \midrule
6$\times$ 3090 & 32 & AMD Epyc Rome & 480 & 5.04 & 4.17 & 154 \\
4$\times$  2080Ti & 16 & Intel Xeon Haswell & 240 & 0.96 & 0.84 & 83.4 \\
1$\times$  RTX 1080Ti & 8 & Intel Xeon Haswell & 16 & 0.22 & 0.16 & 12 \\ \bottomrule
\end{tabular}
\end{table}

Table~\ref{fig:cloud_costs} shows two main tendencies. First, preemptible \textit{cloud} instances are, on average, three times cheaper than their non-preemptible counterparts\footnote{The cost can be up to $11{\times}$ cheaper for some instance types, e.g. Azure V100 instances in the central US region at the time of writing.}. Second, the high-end HPC-grade servers that offer the highest raw performance are less cost-effective than lower-tier servers and marketplace instances. In theory, one could match the raw floating-point performance of a $8{\times}$V100 instance at a fraction of its cost using multiple lower-tier workstations, such as $4{\times}$ RTX 2080Ti, with a smaller total cost.
However, in practice, running distributed training with these workstations is challenging due to their unreliability and slow network connection.

Note that this analysis does not represent the cloud costs for sustained GPU usage. If an organization plans to constantly use GPU resources over a period of multiple years, they can reduce the costs by deploying their own compute infrastructure or relying on the sustained usage discounts reaching up to 60--70\%. Thus, the long-term compute costs are much harder to analyze and depend on a number of additional factors, such as local electricity prices for on-premise infrastructure. However, this scenario offers similar trade-offs: HPC-grade infrastructure offers greater interconnectivity, but requires expensive network interface cards, high-end switches and a more complex setup process.

\section{Additional Related Work}
\label{sect:post_related}

In this section, we review some of the papers relevant to our work, but omitted from the main part due to space constraints. 

\subsection{Decentralized Training}
In this subsection, we give additional details about the dependence of gossip-based optimization methods on the spectral properties on the communication graph through the spectral properties of the mixing matrix~\cite{xiao2004fast,scaman2019optimal} or the Laplacian matrix \cite{merris1994laplacian,uribe2020dual} of the network. 
That is, gossip finds approximate average on nodes with accuracy $\varepsilon$ after $\cO\left((1-\lambda_2(\mM))^{-1}\log(\varepsilon^{-1})\right)$ iterations, where $\mM$ is the mixing matrix and $\lambda_2(\mM)$ is the second largest eigenvalue of $\mM$ when sorted by absolute value. 
The quantity $\eta = 1-\lambda_2(\mM)$ is called the spectral gap of the mixing matrix $\mM$, and $\eta^{-1}$ is typically a polynomial of the total number of nodes $n$ when the maximal degree of the node is $\cO(1)$. For example, for uniformly averaging $\mM$ one can show that $\eta^{-1} = \cO(n^2)$ for the ring topology (node degree $2$), $\eta^{-1} = \cO(n)$ for the two-dimensional torus topology (node degree  $2$), and $\eta^{-1} = \cO(1)$ for the fully connected graph (node degree $n-1$); one can find more examples in~\cite{aldous2002reversible}. Similarly, the communication complexity of decentralized optimization methods often has multiplicative dependence on either $\cO(\eta^{-1})$ (see \cite{xu2020distributed} and references therein) or $\cO(\eta^{-\nicefrac{1}{2}})$ \cite{scaman2019optimal,uribe2020dual,fallah2019robust,kovalev2020optimal}, which is not improvable for gossip-based methods~\cite{arjevani2015communication,scaman2017optimal}.

Contrary to this, {\tt Moshpit All-Reduce} does not depend on a fixed communication graph and the properties of its mixing matrix.
However, it depends on the number of averaging groups and the total number of peers (see Theorem~\ref{thm:quality_of_avg_deterministic_vectors}), which can be viewed as properties of a time-varying random communication graph. Fortunately, this dependence is often much better than in gossip: as we mentioned in the main part of the paper, even if workers are randomly split into pairs at each iteration, the simplified version of {\tt Moshpit All-Reduce} makes the average distortion (the left-hand side of Equation~\ref{eq:determ_quality_of_avg}) at least $2$ times smaller after each round on average.

\subsection{Compressed Communication}
Another popular approach to addressing the communication bottleneck is communication compression~\cite{seide20141,alistarh2017qsgd,suresh2017distributed}: before sending any information (e.g., iterates, gradients, Hessians or more sophisticated data) over the network, peers compress this information by applying some (possibly random) transformation. As the result, peers send fewer bits for each communication round, but the total number of communication rounds needed to achieve the predefined accuracy of the solution increases. However, communication compression is very useful in the situations when the reduction in communication costs of one round is more important than the increase in the number of these rounds~\cite{Cnat}.

There are two distinct groups of works on distributed training with compressed communication: ones that focus on unbiased compression operators (e.g., Rand-K, $\ell_p$-quantization) and ones studying algorithms with biased compressors (e.g., Top-K); see a detailed summary of  popular compression operators in~\cite{beznosikov2020biased}). 
Quantized SGD (QSGD)~\cite{alistarh2017qsgd} and TernGrad~\cite{wen2017terngrad} were among the first compression methods with convergence guarantees. Next, the convergence analysis of these methods was generalized and tightened in the (strongly) convex case in~\cite{mishchenko2019distributed}. Moreover, the authors of \cite{mishchenko2019distributed} proposed a modification of QSGD called DIANA: this algorithm is based on the quantization of gradients' differences, which helps it achieve linear convergence in the strongly convex case when peers compute full gradients. Next, DIANA was generalized to arbitrary unbiased compression in~\cite{horvath2019stochastic}, where authors also developed and analyzed the variance-reduced version of DIANA. After that, several further modifications, such as Accelerated DIANA~\cite{li2020acceleration} and DIANA with bidirectional compression~\cite{gorbunov2020linearly,philippenko2020artemis}, were proposed. Finally, we refer the reader to~\cite{li2020unified,haddadpour2020federated,das2020improved} for state-of-the-art results for distributed methods with unbiased compression in the non-convex case.

However, naïve application of biased compression operators can lead to significantly worse performance in practice. For instance, as it was shown recently in~\cite{beznosikov2020biased}, parallel SGD with Top-1 compression can diverge exponentially fast. Therefore, biased compressors are used jointly with so-called error-compensation~\cite{seide20141}. The first analysis of Error-Compensated SGD (EC-SGD) was proposed in~\cite{stich2018sparsified,karimireddy2019error} which then was generalized and tightened in~\cite{beznosikov2020biased}. Next, several further improvements, such as an accelerated version of EC-SGD~\cite{qian2020error} and linearly converging EC-SGD~\cite{gorbunov2020linearly}, were recently proposed. However, current theory does not show any superiority of distributed methods with biased compressors to the ones with unbiased compression operators.
In addition, one can combine decentralized communication with compression. Such combinations with unbiased compression operators were studied in~\cite{reisizadeh2019exact,kovalev2020linearly} and with biased operators in~\cite{pmlr-v97-koloskova19a,KoloskovaLSJ19decentralized}.
In this paper, we do not study the interaction of different compression methods and Moshpit Averaging, leaving this promising direction to future work.

\subsection{Multiple Local Steps}
Alternatively, to reduce the impact of the communication bottleneck, it is possible to perform several local optimization steps on each peer between the communication rounds.
This approach is based on the idea that the increased computational load of peers will decrease the number of communication rounds required to obtain the optimal parameters; it is frequently used in federated learning~\cite{FEDLEARN,kairouz2019advances}. In particular, one of the most popular methods with multiple local steps is called Local-SGD or Federated Averaging~\cite{FEDLEARN,Stich18local}. The first results on its convergence were given in \cite{Stich18local,LinSPJ2018local}, and later they were tightened and generalized both for homogeneous~\cite{khaled2020tighter,woodworth2020local} and heterogeneous  cases~\cite{khaled2020tighter,woodworth2020minibatch}. Recently, further modifications of Local-SGD were proposed and analyzed: these modifications include acceleration \cite{yuan2020federated}, variance reduction \cite{gorbunov2020local}, communication compression \cite{basu2019qsparse,haddadpour2020federated,das2020improved}, decentralization \cite{li2019communication,koloskova2020unified}, adaptive and proximal methods \cite{reddi2021adaptive,yuan2020federated_comp}, and resistance to client drift \cite{karimireddy2020scaffold}.
{\tt Moshpit SGD} can perform multiple local gradient steps before synchronization by design, as shown in Algorithm~\ref{alg:ungossip_Local_SGD}.

\subsection{Asynchronous Methods}
In the previous subsections, we mostly discussed synchronous distributed methods, since they are more widespread and better studied than asynchronous ones. Mainly, this is because asynchronous methods are more difficult to implement, debug and analyze under general assumptions. However, such methods can be more efficient in terms of using computational resources, which leads to faster wall-clock convergence \cite{assran2020advances}. In recent years, several asynchronous stochastic methods~\cite{recht2011hogwild,zhao2016fast,leblond2017asaga}, methods with no shared memory \cite{peng2016arock,mishchenko2018delay}, and methods with delayed updates~\cite{agarwal2011distributed,feyzmahdavian2016asynchronous,arjevani2020tight,gorbunov2020linearly} were proposed and analyzed. One can find more details in a recent survey of asynchronous distributed methods~\cite{assran2020advances}.
{\tt Moshpit SGD} belongs to this family of asynchronous approaches as well, because the averaging steps happen in smaller groups and can be interleaved with local parameter updates.

\subsection{Distributed Hash Tables}
\label{sect:related_dht}

In this work, we set out to improve distributed averaging with a dynamic matchmaking protocol. Without a central server, this protocol relies on decentralized data structures to organize peers. The main data structure we use is the Distributed Hash Table, or DHT. On a high level, DHT is a distributed fault-tolerant ``dictionary'' that can be accessed by every participant. Each key-value pair is stored on a subset of peers determined by the $\mathrm{hash}$ function of the key.

Each participant has a unique identifier (ID) sampled uniformly from the $\mathrm{hash}$ function output range. When storing a $(key,\ value)$ pair, one must find $k$ peers whose IDs are nearest to $\mathrm{hash}(key)$ according to a chosen metric. After that, the participant requests each of those peers to store $(key,\ value)$. When retrieving a value for a key, one should compute $\mathrm{hash}(key)$, search for peers with IDs nearest to that $\mathrm{hash}$ value and request the value from those peers.

Specific DHT versions, such as Chord~\cite{chord} or Kademlia~\cite{kademlia}, employ different hash types and algorithms for finding nearest peers. For instance, Kademlia DHT sorts peers based on the XOR distance function: $d(x, y) = \mathrm{int}(x \oplus y)$.

In DHT, each participant is directly aware of only a small subset of peers. When storing or retrieving a key, the participant requests additional peers from its neighbors in a semi-greedy search, minimizing the XOR distance until it finds $k$ nearest peers. In Kademlia, nodes form a special navigable graph structure that lets them find nearest peers in at most $\cO(k + \log n)$ requests to other peers, where $n$ is the total number of participants. Due to their scalability and fault-tolerance, DHTs found numerous applications including BitTorrent, Ethereum, I2P and even deep learning with Mixtures-of-Experts~\cite{mryab}.

\section{Proofs of Mixing Properties of {\tt Moshpit All-Reduce}}\label{sect:missing_proofs}
Here we formally state the theorems about mixing properties of Moshpit Averaging along with their proofs.

\textbf{Notation.} Throughout the following sections, we use the standard notation from the literature on stochastic optimization. That is, for any $n$-dimensional vectors $x = (x_1,\ldots,x_n)^\top,y = (y_1,\ldots,y_n)^\top\in\R^d$ we use $\langle x,y\rangle$ to denote the standard inner product: $\langle x, y\rangle = x_1y_1 + \ldots + x_ny_n$. Next, we use $\|x\|$ to denote the $\ell_2$=norm of $x$ ($\|x\| = \sqrt{\langle x, x\rangle}$), $\EE[\xi]$ to denote an expectation of a random variable $\xi$, $\EE[\xi\mid \eta]$ is used for the conditional expectation of $\xi$ given $\eta$, and $\PP\{E\}$ denotes the probability of an event $E$.

\subsection{Computing Exact Average in a Full Grid}\label{sect:equiv_to_torus}
As discussed in Section~\ref{sect:method_algorithm}, {\tt Moshpit All-Reduce} obtains the exact average of parameter vectors from $n$ peers arranged in a grid with $N$ coordinates and $M$ positions per coordinate when $n\equiv M^N$. That is, when the grid is full and each step averages $M$ parameter values along a single grid coordinate without repetitions, the algorithm needs only $N$ steps to compute the actual average across all nodes. In this section, we give a proof of this fact.

First, let us formally define the setting and the averaging steps of {\tt Moshpit All-Reduce} in this specific case. Let $\x_{i_1 i_2\ldots i_N}$ be the parameter vector of the worker with coordinates $i_1, i_2,\ldots, i_N$; each coordinate $i_k$ takes values from $1$ to $M$, because the hypercube of peers is completely full (thus, due to the pigeonhole principle, there are no unoccupied coordinates). Next, arrange the coordinates of these vector according to the order of averaging iterations: namely, at iteration 1
\begin{equation}
    \overline{\x}_{i_1 i_2\ldots  i_N}^1=\frac{1}{M}\sum_{j_1=1}^M \x_{j_1 i_2\ldots i_N},\quad  i_1\in\{1,\ldots,M\},
\end{equation}
which means that for the first iteration, we take the average across the first axis $\overline{\x}^1$ and replicate it across all $M$ resulting vectors regardless of their index $i_1$. The next averaging steps can be expressed similarly with a simple recurrence relation:
\begin{equation}
\label{eqn:avg_recurrence}
    \overline{\x}_{i_1 i_2 \ldots i_N}^t=\frac{1}{M}\sum_{j_t=1}^M \overline{\x}_{i_1\ldots i_{t-1} j_t i_{t+1}\ldots i_N}^{t-1}.
\end{equation}
Given this formal definition, we can now state and prove the exact averaging result:
\begin{theorem}[Exact average in a full $N$-dimensional hypercube after $N$ steps]
Assume that $M^N$ peers are arranged in a $N$-dimensional hypercube with $M$ positions in each dimension. Also, assume that each peer fully participates in every averaging step and $M$-sized groups for each averaging iteration are determined based on the hypercube coordinates. Then, if {\tt Moshpit All-Reduce} is ran in the above setup for $N$ iterations without repeating groups (i.e. averaging across each dimension exactly once), its result for each participant is the average value of $\x$ across all $M^N$ peers.
\end{theorem}
\begin{proof}
We can directly obtain the expression for the average by expanding the recurrence and rearranging the sums:
\begin{gather*}
    \overline{\x}_{i_1 i_2\ldots i_N}^N=\frac{1}{M}\sum_{j_N=1}^M\overline{\x}_{i_1\ldots i_{N-1} j_N}^{N-1}=\frac{1}{M}\sum_{j_N=1}^M\left(\frac{1}{M}\sum_{j_{N-1}=1}^M \overline{\x}_{i_1 i_2\ldots j_{N-1}j_N}\right)=\ldots\\
    =\frac{1}{M}\Bigg(\underbrace{\sum_{j_N=1}^M\Bigg(\frac{1}{M}\sum_{j_{N-1}=1}^M\ldots\sum_{j_2=1}^M\Bigg(\frac{1}{M}\sum_{j_1=1}^M}_{N\textrm{ summations}} \x_{j_1 \ldots j_N}\Bigg)\Bigg)\Bigg)=\frac{1}{M^N}\sum_{j_N=1}^M\sum_{j_{N-1}=1}^M\ldots\sum_{j_2=1}^M\sum_{j_1=1}^M \x_{j_1 \ldots j_N}=\\
    =\frac{1}{M^N}\sum_{j_1, \ldots, j_N=1}^M  \x_{j_1 \ldots j_N}.
\end{gather*}
But this is exactly the global average of all $\x$, since there are $M^N$ participants and each vector is represented in the sum because of summation over all possible indices.
\end{proof}

Notice that for a given grid of peers, if some of its indices do not have corresponding parameter vectors, Equation~\eqref{eqn:avg_recurrence} may result in different average vectors on different workers due to different numbers of peers along a coordinate for different indices. For example, running two iterations of Moshpit Averaging with $N=2,\ M=2$ and three parameter vectors $\x_{11},\ \x_{21},\ \x_{22}$ results in $\frac{\x_{11}+\x_{21}}{2}$ on the first worker and $\frac{\x_{11}+\x_{21}}{4}+\x_{22}$ on other workers, so neither of the values is equal to the global average. However, the variance of the averaged vectors does decrease, which is formally proven in Section~\ref{sec:proof_quality_of_avg_deterministic_vectors}.

\subsection{Proof of Theorem~\ref{thm:quality_of_avg_deterministic_vectors_0}}\label{sect:correctness_proof}
Below we provide the complete proof of Theorem~\ref{thm:quality_of_avg_deterministic_vectors_0}. For the readers' convenience, we restate the theorem.
\begin{theorem}[Theorem~\ref{thm:quality_of_avg_deterministic_vectors_0}]\label{thm:quality_of_avg_deterministic_vectors_0_supp}
If all workers have non-zero probability of successfully running a communication round in Moshpit Averaging and the order of $\texttt{peers}_t$ is random, then all local vectors $\x^t_i$ converge to the global average with probability $1$:
\begin{equation}
    \forall i = 1,\ldots, n\quad \left\|\x^t_i - \frac1n \sum_{i=1}^n \x^0_i\right\|^2 \xrightarrow[t\to\infty]{} 0.\label{eq:quality_of_avg_deterministic_vectors_0_supp}
\end{equation}
\end{theorem}
\begin{proof}[Proof of Theorem~\ref{thm:quality_of_avg_deterministic_vectors_0}]
    First of all, we notice that \eqref{eq:quality_of_avg_deterministic_vectors_0_supp} is equivalent to
    \begin{equation}
    \forall i = 1,\ldots, n,\;\forall j=1,\ldots,n\quad \left(\x^t_i(j) - \frac1n \sum_{i=1}^n \x^0_i(j)\right)^2 \xrightarrow[t\to\infty]{} 0,\label{eq:quality_of_avg_deterministic_vectors_0_supp_tech_1}
    \end{equation}
    where $\x_i^t(j)$ denotes $j$-th component of $\x_i^t$. Consider an arbitrary component $j \in \{1,\ldots,n\}$ and the sequence of intervals $\{I_{j,t}\}_{t\ge 0}$ where $I_{j,t} = \text{conv}\{\x_1^t(j),\x_2^t(j),\ldots, \x_n^t(j)\}$. Then, $\{I_{j,t}\}_{t\ge 0}$ is a sequence of nested intervals ($I_{j,t+1} \subseteq I_{j,t} \forall t\ge 0$), since averaging in groups does not expand the convex hull of $\{\x_1^t,\x_2^t,\ldots, \x_n^t\}$. For convenience, we specify the bounds of the intervals: $I_{j,t} = [a_{j,t}, b_{j,t}]$. Using the Cantor's intersection theorem, we conclude that
    \begin{equation*}
        \bigcap\limits_{t=0}^\infty I_{j,t} = I_j = [a_j, b_j],
    \end{equation*}
    where $\overline{\x}(j) = \frac{1}{n}\sum_{i=1}^n\x_i^0(j) \in [a_j, b_j]$. If $[a_j, b_j] = \{\overline{\x}(j)\}$ with probability $1$, then \eqref{eq:quality_of_avg_deterministic_vectors_0_supp_tech_1} holds with probability $1$ as well. Suppose the opposite: there exist such $j \in \{1,\ldots,n\}$, $[a,b]$ and $\delta,\Delta > 0$ that $\overline{\x}(j) \in [a,b]$, $b-a = \Delta$ and
    \begin{equation*}
        \PP\Bigg\{\underbrace{[a,b] \subseteq \bigcap\limits_{t=0}^\infty I_{j,t}}_{E}\Bigg\} = \delta > 0\quad \text{ and }\quad \forall \varepsilon > 0\; \PP\Bigg\{\underbrace{[a-\varepsilon,b+\varepsilon] \subseteq \bigcap\limits_{t=0}^\infty I_{j,t}}_{E_{\varepsilon}}\Bigg\} < \delta.
    \end{equation*}
    This implies that for all $\varepsilon > 0$ there exists such $T_{\varepsilon} > 0$ that
    \begin{equation*}
        \PP\Big\{\underbrace{\forall t \ge T_{\varepsilon}\;\; a_{j,t}\in [a-\varepsilon,a], b_{j,t}\in[b,b+\varepsilon]}_{E_{\varepsilon}'}\Big\} = \delta_{\varepsilon} > 0.
    \end{equation*}
    Consider $\varepsilon = \frac{\Delta}{(2n+100)^{2n}}$ and assume that the event $E_{\varepsilon}'$ holds. Next, we introduce new notation: $J_{\text{left}}^t = \{i \in \{1,\ldots, n\}\mid \x_{i}^t(j) \in [a-\varepsilon,a]\}$ and $J_{\text{right}}^t = \{i \in \{1,\ldots, n\}\mid \x_{i}^t(j) \in [b,b+\varepsilon]\}$. Since $E_{\varepsilon}'$ holds the sets $J_{\text{left}}^t$ and $J_{\text{right}}^t$ are non-empty for all $t\ge T_{\varepsilon}$ with probability $\delta_{\varepsilon} > 0$:
    \begin{equation}
        \PP\left\{\forall t \ge T_{\varepsilon}\;\; J_{\text{left}}^t \neq \varnothing\text{ and }  J_{\text{right}}^t \neq \varnothing\right\} = \delta_{\varepsilon} > 0. \label{eq:quality_of_avg_deterministic_vectors_0_supp_tech_2}
    \end{equation}
    We notice that every pair of workers $i_1,i_2$ has a non-zero probability of taking part in the averaging inside the common group at each iteration since all workers have a non-zero probability of successfully running a communication round and the order of $\texttt{peers}_t$ is random. This implies that every pair of workers $i_1,i_2$ with probability $1$ take part in the averaging inside the common group infinitely many times when $t$ goes to the infinity.
    
    Next, we choose some $t_0 \ge T_{\varepsilon}$. Let $J_{\text{left}}^{t_0} = \{i_{l,1},\ldots, i_{l,q_l}\}$ and $J_{\text{right}}^{t_0} = \{i_{r,1},\ldots, i_{r,q_r}\}$. Consider the event $E_{\varepsilon,0}' \subseteq E_{\varepsilon}'$ such that in $E_{\varepsilon,0}'$ peer $i_{l,1}$ computes an average in the group containing any peer from $J_{\text{right}}^{t_0}$ at some iteration $t_1 > t_0$. Our observations above imply that $\PP\{E_{\varepsilon,0}'\} = \PP\{E_{\varepsilon}'\} = \delta_{\varepsilon} > 0$. Then, $\x_{i_{l,1}}^{t_1}(j) \ge \frac{n-1}{n}(a-\varepsilon) + \frac{1}{n}b = a-\varepsilon + \frac{1}{n}(\Delta + \varepsilon) = a - \frac{\Delta}{(2n+100)^{2n}} + \frac{1}{n}\left(\Delta + \frac{\Delta}{(2n+100)^{2n}}\right) > a + \frac{\Delta}{2n}$, i.e., $\x_{i_{l,1}}^{t_1}(j) \in (a,b]$ meaning that $i_{l,1} \not\in J_{\text{left}}^{t_1}$. The last part of the proof shows that for any $t\ge t_1$, the peer $i_{l,1}$ will never be the part of $J_{\text{left}}^t$ and after a finite number of iterations $J_{\text{left}}^t = \varnothing$ with probability $\delta_{\varepsilon} > 0$ when $E_{\varepsilon,0}'$ holds, implying the contradiction with \eqref{eq:quality_of_avg_deterministic_vectors_0_supp_tech_2}.
    
    To show that, we consider the following set of peers: $\widehat{J}_{\text{left}}^{t_1} = \{i\in\{1,\ldots,n\}\mid \exists t \ge t_1:\; \x_i^{t}(j)\in [a-\varepsilon, a+\frac{\Delta}{2n})\}$. Next, we consider the event $E_{\varepsilon,1}'\subseteq E_{\varepsilon,0}'$ such that in $E_{\varepsilon,1}'$ peer $i_{l,1}$ computes an average in the group containing some peer $i_{l,avg,1}$ from $\widehat{J}_{\text{left}}^{t_1}$ at some iteration $t_2 > t_1$ (and $t_2$ is the first such moment after $t_1$). Again, our observations imply $\PP\{E_{\varepsilon,1}'\} = \PP\{E_{\varepsilon,0}'\} = \delta_{\varepsilon}>0$. Then, $\x_{i_{l,1}}^{t_2}(j) = \x_{i_{l,avg,1}}^{t_2}(j) > \frac{n-1}{n}(a-\varepsilon) + \frac{1}{n}\left(a+\frac{\Delta}{2n}\right) = a + \frac{\Delta}{2n^2} - \frac{(n-1)\Delta}{n(2n+100)^{2n}} > a + \frac{\Delta}{4n^2}$. After that, we consider the event $E_{\varepsilon,2}'\subseteq E_{\varepsilon,1}'$ such that in $E_{\varepsilon,2}'$ peer $i_{l,1}$ or $i_{l,avg,1}$ computes an average in the group containing a peer $i_{l,avg,2}\neq i_{l,avg,1}$ from $\widehat{J}_{\text{left}}^{t_1}$ at an iteration $t_3 > t_2$ (and $t_3$ is the first such moment after $t_2$). Then, $\x_{i_{l,1}}^{t_3}(j), \x_{i_{l,avg,1}}^{t_3}(j)$ and $\x_{i_{l,avg,2}}^{t_3}(j)$ are greater than $\frac{n-1}{n}(a-\varepsilon) + \frac{1}{n}\left(a + \frac{\Delta}{4n^2}\right) = a + \frac{\Delta}{4n^3} - \frac{(n-1)\Delta}{n(2n+100)^{2n}} > a + \frac{\Delta}{8n^3}$.
    
    Therefore, after at least $n-1$ of such averaging iterations, with probability $\delta_\varepsilon$ all $\x_i^t(j)$ will be greater than $a + \frac{\Delta}{(2n)^n} > a$ while $E_{\varepsilon}'$ holds. This contradicts \eqref{eq:quality_of_avg_deterministic_vectors_0_supp_tech_2}. Therefore, 
    \begin{equation*}
        \bigcap\limits_{t=0}^\infty I_{j,t} = \{\overline{\x}(j)\}
    \end{equation*}
    with probability $1$, which concludes the proof.
\end{proof}

\subsection{Proof of Theorem~\ref{thm:quality_of_avg_deterministic_vectors}}\label{sec:proof_quality_of_avg_deterministic_vectors}
In this section, we provide the complete proof of Theorem~\ref{thm:quality_of_avg_deterministic_vectors}. For convenience, we restate the theorem below.
\begin{theorem}[Theorem~\ref{thm:quality_of_avg_deterministic_vectors}, averaging convergence rate]\label{thm:quality_of_avg_deterministic_vectors_supp}
    Consider the modification of {\tt Moshpit All-Reduce} that works as follows: at each iteration $k\geq 1$ 1) peers are randomly split into $r$ disjoint groups of sizes $M_1^k,\ldots, M_r^k$ in such a way that $\sum_{i=1}^r M_i^k = n$ and $M_i^k \ge 1\  \forall i = 1,\ldots,r$ and 2) peers from each group compute their group average via All-Reduce. Let $\x_1,\ldots,\x_n$ be the input vectors of this procedure and $\x_1^T,\ldots,\x_n^T$ be the outputs after $T$ iterations. Then,
    \begin{eqnarray}
         \EE\left[\frac{1}{n}\sum\limits_{i=1}^n\|\x_i^T - \overline{\x}\|^2\right] = \left(\frac{r-1}{n} + \frac{r}{n^2}\right)^T\cdot\frac{1}{n}\sum\limits_{i=1}^n\|\x_i - \overline{\x}\|^2, \label{eq:determ_quality_of_avg_supp}
    \end{eqnarray}
    where $\overline{\x} = \frac{1}{n}\sum_{i=1}^n\x_i$.
\end{theorem}
\begin{proof}
First of all, let us clarify the procedure of random splitting of peers in $r$ groups. We assume that at iteration $k$ of the modified algorithm we generate a random permutation $\pi^k = (\pi_1^k,\ldots,\pi_n^k)$ of $1,\ldots, n$. Next, $J_1^k = \{\pi_1^k,\ldots,\pi_{M_1^k}^k\}$ form the indices of the first group of workers, $J_2^k = \{\pi_{M_1^k+1}^k,\ldots,\pi_{M_2^k}^k\}$ are the indices of the second group, and $J_r^k = \{\pi_{M_1^k+M_2^k+\ldots+M_{r-1}^k+1}^k,\ldots,\pi_{n}^k\}$ are the indices of group $r$. In other words, we generate a random permutation and take contiguous subgroups of indices corresponding to predefined group sizes $M_i^k$, starting from the first group.

By definition, we have $\bigsqcup_{i=1}^r J_i^k = \{1,2,\ldots,n\}$, where $\sqcup$ defines the disjoint union operator. Moreover, notice that group sizes $M_1^k,\ldots,M_r^k$ can depend on $k$ and even be random: for our analysis, it is sufficient that the randomness defining the permutation is independent from $M_1^k,\ldots,M_r^k$. Next, vectors $\x_1^k,\ldots,\x_n^k$ are obtained by the following formula:
\begin{equation*}
    \forall j=1,\ldots,n,\quad \x_j^k = \frac{1}{M_i^k}\sum\limits_{t\in J_i^k}\x_t^{k-1},\quad \text{where } J_i^k \text{ is the group for which } j\in J_i^k.
\end{equation*}
Using this, we show that the average of vectors $\{\x_i^k\}_{i=1}^n$ remains the same throughout the iterations of {\tt Moshpit All-Reduce}:
\begin{equation*}
    \frac{1}{n}\sum\limits_{j=1}^n\x_j^k = \frac{1}{n}\sum\limits_{i=1}^rM_i^k\cdot\frac{1}{M_i^k}\sum\limits_{t\in J_i^k}\x_t^{k-1} = \frac{1}{n}\sum\limits_{i=1}^r\sum\limits_{t\in J_i^k}\x_t^{k-1} = \frac{1}{n}\sum\limits_{j=1}^n\x_j^{k-1}.
\end{equation*}
Therefore, the quantity $\frac{1}{n}\sum_{j=1}^n\|\x_j^k - \overline{\x}\|^2$ (average distortion) measures the quality of averaging. For this quantity, we can derive the following expression:
\begin{eqnarray}
    \frac{1}{n}\sum\limits_{j=1}^n\|\x_j^k - \overline{\x}\|^2 &=& \frac{1}{n}\sum\limits_{i=1}^r M_i^k\left\|\frac{1}{M_i^k}\sum\limits_{t\in J_i^k}\x_t^{k-1} - \overline{\x}\right\|^2\notag\\
    &=& \frac{1}{n}\sum\limits_{i=1}^r\frac{1}{M_i^k}\left(\sum\limits_{t\in J_i^k}\|\x_t^{k-1} - \overline{\x}\|^2 + 2\sum\limits_{t,l\in J_i^k, t < l}\langle \x_t^{k-1} - \overline{\x}, \x_l^{k-1} - \overline{\x} \rangle\right).\notag
\end{eqnarray}
Taking the expectation $\EE_{\pi^k}[\cdot]$ with respect to the randomness coming from the choice of $\pi^k$ we get
\begin{eqnarray}
    \EE_{\pi^k}\left[\frac{1}{n}\sum\limits_{j=1}^n\|\x_j^k - \overline{\x}\|^2\right] \!=\! \frac{1}{n}\sum\limits_{i=1}^r\frac{1}{M_i^k}\left(\EE_{\pi^k}\left[\sum\limits_{t\in J_i^k}\|\x_t^{k-1} - \overline{\x}\|^2\!\right] \!+\! 2\EE_{\pi^k}\!\left[\sum\limits_{t,l\in J_i^k, t < l}\langle \x_t^{k-1} - \overline{\x}, \x_l^{k-1} - \overline{\x} \rangle\right]\right).\notag
\end{eqnarray}
Since $\forall j,j_1,j_2 \in\{1,\ldots,n\},j_1\neq j_2$ and for all $i=1,\ldots,r$
\begin{equation*}
    \PP\left\{j\in J_i^k\right\} = \frac{M_i^k}{n},\quad \PP\left\{j_1,j_2 \in J_i^k\right\} = \frac{M_{i}^k(M_i^k - 1)}{n^2},
\end{equation*}
we have
\begin{eqnarray*}
    \EE_{\pi^k}\left[\frac{1}{n}\sum\limits_{j=1}^n\|\x_j^k - \overline{\x}\|^2\right] &=& \frac{1}{n}\sum\limits_{i=1}^r\frac{1}{M_i^k}\Bigg(\frac{M_i^k}{n}\sum\limits_{j=1}^n\|\x_j^{k-1} - \overline{\x}\|^2\\
    &&\quad+ 2\frac{M_{i}^k(M_i^k - 1)}{n^2}\sum\limits_{1 \le j_1 < j_2 \le n}\langle \x_{j_1}^{k-1} - \overline{\x}, \x_{j_2}^{k-1} - \overline{\x}\rangle\Bigg)\\
    &=& \frac{r}{n^2}\sum\limits_{j=1}^n\|\x_j^{k-1} - \overline{\x}\|^2 + 2\frac{n-r}{n^3}\sum\limits_{1 \le j_1 < j_2 \le n}\langle \x_{j_1}^{k-1} - \overline{\x}, \x_{j_2}^{k-1} - \overline{\x}\rangle\\
    &=& \left(\frac{r}{n^2} - \frac{n-r}{n^3}\right)\sum\limits_{j=1}^n\|\x_j^{k-1} - \overline{\x}\|^2\\
    &&\quad+\frac{n-r}{n^3}\left(\sum\limits_{j=1}^n\|\x_j^{k-1} - \overline{\x}\|^2 + 2\sum\limits_{1 \le j_1 < j_2 \le n}\langle \x_{j_1}^{k-1} - \overline{\x}, \x_{j_2}^{k-1} - \overline{\x}\rangle\right)\\
    &=& \frac{n(r-1)+r}{n^3}\sum\limits_{j=1}^n\|\x_j^{k-1} - \overline{\x}\|^2 + \frac{n-r}{n^3}\underbrace{\left\|\sum\limits_{j=1}^n(\x_j^{k-1} - \overline{\x})\right\|^2}_{\|n\overline{\x} - n\overline{\x}\|^2 = 0}\\
    &=& \left(\frac{r-1}{n} + \frac{r}{n^2}\right)\cdot\frac{1}{n}\sum\limits_{j=1}^n\|\x_j^{k-1} - \overline{\x}\|^2.
\end{eqnarray*}
Finally, we take the full expectation from the both sides of the above equation and apply the tower property $\EE\left[\EE_{\pi^k}\left[\cdot\right]\right] = \EE\left[\cdot\right]$:
\begin{equation*}
    \EE\left[\frac{1}{n}\sum\limits_{j=1}^n\|\x_j^k - \overline{\x}\|^2\right] = \left(\frac{r-1}{n} + \frac{r}{n^2}\right)\EE\left[\frac{1}{n}\sum\limits_{j=1}^n\|\x_j^{k-1} - \overline{\x}\|^2\right].
\end{equation*}
Unrolling the recurrence for $k=T$, we establish \eqref{eq:determ_quality_of_avg_supp}.
\end{proof}

\begin{remark}
    The result implies that increasing the group size $\alpha > 1$ times implies almost $\alpha$ times faster convergence to the average.
\end{remark}

\begin{remark}
    Our analysis can be easily generalized to the case when number of groups $r$ can depend on $k$ and be a random variable independent from the choice of permutations and the number of groups at previous steps. In this case, \eqref{eq:determ_quality_of_avg_supp} transforms into
    \begin{equation}
        \EE\left[\frac{1}{n}\sum\limits_{i=1}^n\|\x_i^T - \overline{\x}\|^2\right] = \frac{1}{n}\sum\limits_{i=1}^n\|\x_i - \overline{\x}\|^2\cdot\prod_{k=1}^T\left(\frac{\EE[r_k]-1}{n} + \frac{\EE[r_k]}{n^2}\right), \label{eq:determ_quality_of_avg_generalized_supp}
    \end{equation}
    where $r_k$ is the number of groups at iteration $k$.
\end{remark}

\subsection{Additional Guarantees For Moshpit Averaging}\label{sec:mix_rand_proof}
In this section,  we derive the result measuring the rate of variance reduction when averaging random vectors with Algorithm~\ref{alg:ungossip}. We start with the following technical lemma:
\begin{lemma}\label{lem:ode_lemma}
    Let $\xi \sim \text{Binom}(M,p)$ have a binomial distribution with parameters $M$ (number of trials) and $p$ (probability of success for each trial). Then
    \begin{eqnarray}
        m_1(M,p) := \EE\left[\min\left\{\frac{1}{\xi},1\right\}\right] &=& (1-p)^M + \sum\limits_{i=1}^M\frac{1}{i}\left((1-p)^{M-i} - (1-p)^M\right), \label{eq:binom_first_inverse_moment}\\
        m_2(M,p) := \EE\left[\min\left\{\frac{1}{\xi^2},1\right\}\right] &=& (1-p)^M + \sum\limits_{i=1}^M\frac{1}{i}\left((1-p)^{M-i} - (1-p)^M\right)\sum\limits_{j=i}^M\frac{1}{j}. \label{eq:binom_second_inverse_moment}
    \end{eqnarray}
\end{lemma}
\begin{proof}
    We start with the proof of \eqref{eq:binom_first_inverse_moment}. By definition of the expectation, we have
    \begin{eqnarray*}
        \EE\left[\min\left\{\frac{1}{\xi},1\right\}\right] &=& (1-p)^M + \sum\limits_{i=1}^M \frac{1}{i}p^i(1-p)^{M-i}\binom{M}{i}.
    \end{eqnarray*}
    For simplicity of further derivations, we introduce the following notation: $m_1(M,p) = \EE\left[\min\left\{\frac{1}{\xi},1\right\}\right]$ and $m_2(M,p) = \EE\left[\min\left\{\frac{1}{\xi^2},1\right\}\right]$. Taking the derivative of $m_1(M,p)$ by $p$, we obtain
    \begin{eqnarray*}
        m_1'(M,p) &=& -M(1-p)^{M-1} + \sum\limits_{i=1}^Mp^{i-1}(1-p)^{M-i}\binom{M}{i} - \sum\limits_{i=1}^M\frac{M-i}{i}p^i(1-p)^{M-i-1}\binom{M}{i}\\
        &=& -M(1-p)^{M-1} + \frac{1}{p}\left(-(1-p)^M + \sum\limits_{i=0}^Mp^{i}(1-p)^{M-i}\binom{M}{i}\right)\\
        && - \frac{M}{1-p}\sum\limits_{i=1}^M\frac{1}{i}p^i(1-p)^{M-i}\binom{M}{i} + \frac{1}{1-p}\left(-(1-p)^M + \sum\limits_{i=0}^Mp^i(1-p)^{M-i}\binom{M}{i}\right)\\
        &=& -M(1-p)^{M-1} + \frac{1}{p}\left(1 - (1-p)^M\right) - \frac{M}{1-p}\left(m_1(M,p) - (1-p)^M\right)\\
        &&\quad + \frac{1}{1-p}\left(1- (1-p)^M\right)\\
        &=& \frac{1}{p(1-p)} - \frac{(1-p)^{M-1}}{p} - \frac{M}{1-p}m_1(M,p).
    \end{eqnarray*}
    Rearranging the terms, we get the following linear first-order ODE
    \begin{equation}
        m_1'(M,p) + \frac{M}{1-p}m_1(M,p) = \frac{1}{p(1-p)} - \frac{(1-p)^{M-1}}{p}. \label{eq:first_moment_ODE}
    \end{equation}
    To solve it, we consider the following homogeneous ODE:
    \begin{equation*}
        m_1'(M,p) + \frac{M}{1-p}m_1(M,p) = 0.
    \end{equation*}
    The solution of this ODE is $m_1(M,p) = C(1-p)^M$, where $C\in\R$ is an arbitrary real constant. Next, we go back to the initial ODE \eqref{eq:first_moment_ODE} and try to find a solution of the form $m_1(M,p) = C(p)(1-p)^M$, where $C(p):\R \to \R$ is a differentiable function:
    \begin{eqnarray*}
        \left(C(p)(1-p)^M\right)' + \frac{M}{1-p}C(p)(1-p)^M &=& \frac{1}{p(1-p)} - \frac{(1-p)^{M-1}}{p}\\
        &\Downarrow&\\
        C'(p)(1-p)^M &=& \frac{1}{p(1-p)} - \frac{(1-p)^{M-1}}{p}\\
        &\Downarrow&\\
        C'(p) &=& \frac{1}{p(1-p)^{M+1}} - \frac{1}{p(1-p)}.
    \end{eqnarray*}
    Since 
    \begin{equation}
        \frac{1}{x(1-x)^{k+1}} = \frac{1}{x(1-x)^{k}} + \frac{1}{(1-x)^{k+1}}\label{eq:technical_expansion}
    \end{equation}
    for all $x\not\in \{0,1\}$ and all non-negative integers $k$, we have
    \begin{eqnarray*}
        C'(p) &=& \frac{1}{p} + \frac{1}{1-p} + \frac{1}{(1-p)^2} + \ldots + \frac{1}{(1-p)^{M+1}} - \frac{1}{p} - \frac{1}{1-p}\\
        &\Downarrow&\\
        C'(p) &=& \sum\limits_{i=1}^M(1-p)^{-i-1},
    \end{eqnarray*}
    hence
    \begin{eqnarray*}
        C(p) = \hat{C} + \sum\limits_{i=1}^M\frac{1}{i}(1-p)^{-i},
    \end{eqnarray*}
    where $\hat{C}$ is a real constant. Putting all together, we obtain
    \begin{eqnarray*}
        m_1(M,p) &=& C(p)(1-p)^M = \hat{C}(1-p)^M + \sum\limits_{i=1}^M\frac{1}{i}(1-p)^{M-i}.
    \end{eqnarray*}
    Taking $m_1(M,0) = 1$ into account, we conclude that $\hat{C} = 1 - \sum_{i=1}^M\frac{1}{i}$ and obtain \eqref{eq:binom_first_inverse_moment}.
    
    Using a similar technique, we derive \eqref{eq:binom_second_inverse_moment}. By definition of the expectation, we have
    \begin{eqnarray*}
        m_2(M,p) &=& (1-p)^M + \sum\limits_{i=1}^M \frac{1}{i^2}p^i(1-p)^{M-i}\binom{M}{i}.
    \end{eqnarray*}
    Taking the derivative of $m_2(M,p)$ by $p$, we obtain
    \begin{eqnarray*}
        m_2'(M,p) &=& -M(1-p)^{M-1} + \sum\limits_{i=1}^M\frac{1}{i}p^{i-1}(1-p)^{M-i}\binom{M}{i} - \sum\limits_{i=1}^M\frac{M-i}{i^2}p^i(1-p)^{M-i-1}\binom{M}{i}\\
        &=& -M(1-p)^{M-1} + \frac{1}{p} \sum\limits_{i=1}^M\frac{1}{i}p^{i}(1-p)^{M-i}\binom{M}{i}\\
        && - \frac{M}{1-p}\sum\limits_{i=1}^M\frac{1}{i^2}p^i(1-p)^{M-i}\binom{M}{i} + \frac{1}{1-p}\sum\limits_{i=1}^M\frac{1}{i}p^i(1-p)^{M-i}\binom{M}{i}\\
        &=& -M(1-p)^{M-1} + \frac{1}{p}\left(m_1(M,p) - (1-p)^M\right) \\
        &&\quad + \frac{1}{1-p}\left(-M m_2(M,p) + M(1-p)^M + m_1(M,p) - (1-p)^M\right)\\
        &=& \frac{m_1(M,p)}{p(1-p)} - \frac{(1-p)^{M-1}}{p} - \frac{M}{1-p}m_2(M,p).
    \end{eqnarray*}
    Rearranging the terms, we get the following linear first-order ODE
    \begin{equation}
        m_2'(M,p) + \frac{M}{1-p}m_2(M,p) = \frac{m_1(M,p)}{p(1-p)} - \frac{(1-p)^{M-1}}{p}. \label{eq:second_moment_ODE}
    \end{equation}
    To solve this ODE, we consider the homogeneous ODE:
    \begin{equation*}
        m_2'(M,p) + \frac{M}{1-p}m_2(M,p) = 0.
    \end{equation*}
    The solution of this ODE is $m_2(M,p) = C(1-p)^M$, where $C\in\R$ is an arbitrary real constant. Next, we go back to the initial ODE \eqref{eq:second_moment_ODE} and try to find a solution of the form $m_2(M,p) = C(p)(1-p)^M$, where $C(p):\R \to \R$ is a differentiable function:
    \begin{eqnarray*}
        \left(C(p)(1-p)^M\right)' + \frac{M}{1-p}C(p)(1-p)^M &=& \frac{m_1(M,p)}{p(1-p)} - \frac{(1-p)^{M-1}}{p}\\
        &\Downarrow&\\
        C'(p)(1-p)^M &=& \frac{m_1(M,p)}{p(1-p)} - \frac{(1-p)^{M-1}}{p}\\
        &\Downarrow&\\
        C'(p) &=& \frac{m_1(M,p)}{p(1-p)^{M+1}} - \frac{1}{p(1-p)}.
    \end{eqnarray*}
    Using \eqref{eq:technical_expansion} and \eqref{eq:binom_first_inverse_moment}, we derive
    \begin{eqnarray*}
        C'(p) &\overset{\eqref{eq:binom_first_inverse_moment}}{=}& -\frac{\sum\limits_{i=1}^M\frac{1}{i}}{p(1-p)} + \frac{\sum\limits_{i=1}^M\frac{1}{i}(1-p)^{M-i}}{p(1-p)^{M+1}}\\
        &=& -\sum\limits_{i=1}^M \frac{1}{ip(1-p)} + \sum\limits_{i=1}^M\frac{1}{ip(1-p)^{i+1}}\\
        &\overset{\eqref{eq:technical_expansion}}{=}& -\sum\limits_{i=1}^M\frac{1}{i}\left(\frac{1}{p} + \frac{1}{1-p}\right) + \sum\limits_{i=1}^M\frac{1}{i}\left(\frac{1}{p} + \frac{1}{1-p} + \frac{1}{(1-p)^2} + \ldots + \frac{1}{(1-p)^{i+1}}\right)\\
        &=& \sum\limits_{i=1}^M\frac{1}{i}\left(\frac{1}{(1-p)^2} + \ldots + \frac{1}{(1-p)^{i+1}}\right) = \sum\limits_{i=1}^M \frac{1}{(1-p)^{i+1}}\sum\limits_{j=i}^M\frac{1}{j},
    \end{eqnarray*}
    hence 
    \begin{eqnarray*}
        C(p) = \hat{C} + \sum\limits_{i=1}^M\frac{1}{i}(1-p)^{-i}\sum\limits_{j=i}^M\frac{1}{j},
    \end{eqnarray*}
    where $\hat{C}$ is a real constant. Putting all together, we obtain
    \begin{eqnarray*}
        m_2(M,p) &=& C(p)(1-p)^M = \hat{C}(1-p)^M + \sum\limits_{i=1}^M\frac{1}{i}(1-p)^{M-i}\sum\limits_{j=i}^M\frac{1}{j}.
    \end{eqnarray*}
    Taking $m_2(M,0) = 1$ into account, we conclude that $\hat{C} = 1 - \sum_{i=1}^M\frac{1}{i}\sum_{j=i}^M\frac{1}{j}$ and obtain \eqref{eq:binom_second_inverse_moment}.
\end{proof}

Using this lemma, we derive the following result:
\begin{theorem}\label{thm:quality_of_avg_supp}
    Assume that peers participating in Moshpit Averaging have independent random vectors $\x_1,\ldots,\x_n$ with means $\overline{\x}_1,\ldots,\overline{\x}_n$ and variances bounded by $\sigma^2$ before the averaging. Let $\x_1^T,\ldots,\x_n^T$ be the outputs of Moshpit Averaging after $T$ iterations. Finally, we assume that each peer from the grid can be dropped out for the whole averaging process before averaging independently from other peers, i.e., $n \sim \text{Binom}(M^N,p)$. Then, for all $i = 1,\ldots,n$ we have
    \begin{equation}
        \EE\left[\left\|\x_i^T - \EE_{\x}\left[\x_i^T\right]\right\|^2\right] \leq M^{T-1}\sigma^2 m_1(M-1,p)\left(m_2(M-1,p)\right)^{T-1},\label{eq:variance_bound_supp}
    \end{equation}
    where functions $m_1(M,p)$ and $m_2(M,p)$ are defined in \eqref{eq:binom_first_inverse_moment} and \eqref{eq:binom_second_inverse_moment} respectively, and $\EE_\x\left[\cdot\right]$ denotes the expectation w.r.t.\ the randomness from $\x_1,\ldots,\x_n$. Moreover, if $p \ge \frac{2}{3}$ and $M \ge 11$, then $m_1(M-1,p) \le \frac{2}{M}$, $m_2(M-1,p) \le \frac{3}{M^2}$ and 
    \begin{equation}
        \EE\left[\left\|\x_i^T - \EE_{\x}\left[\x_i^T\right]\right\|^2\right] \leq \frac{2\sigma^2}{M(\nicefrac{M}{3})^{T-1}}.\label{eq:variance_bound_2_supp}
    \end{equation}
\end{theorem}
\begin{proof}
First of all, we recall an equivalent formulation of Moshpit Averaging. Consider a hypercube $\{1,\ldots,M\}^N$. One can consider the elements of this hypercube as hyperindices and assign a unique hyperindex to each peer so that peers can be viewed as vertices in the hypercube. Then, during the $k$-th iteration of {\tt Moshpit All-Reduce}, each worker computes the average among those peers that have hyperindices with the same values except the $k$-th index; in other words, peers compute averages along the $k$-th dimension of the hypercube. Next, if $n = 0$, we assume that $\x_i^T = \EE_{\x}\left[\x_i^T\right]$ and \eqref{eq:variance_bound_supp} holds for free. Therefore, to derive \eqref{eq:variance_bound_supp}, we assume that $n > 0$.

More formally, we use the following notation: $\x_{C_i} = \x_i$ for all $i= 1,\ldots,n$, where $C_{i} = (c_{1}^i, c_2^i,\ldots, c_N^i)$, $c_{j}^i \in \{1,\ldots,M\}$ for all $j = 1,\ldots,M$, and $C_{i} \neq C_k$ for $i\neq k$. Let $\cC$ be the set of hyperindices corresponding to all peers. Next, we use $\x_{C_i}^t$ to define the vector stored on $i$-th peer after $t$ iterations of Moshpit Averaging. Then, for all $i = 1,\ldots,n$ we have $\x_{C_i}^0 = \x_{C_i}$ and for all $t = 1,\ldots,N$
\begin{equation*}
    \x_{C_i}^{t} = \frac{1}{b_{i,t}}\sum\limits_{k\in J_{i,t}}\x_{C_k}^{t-1},
\end{equation*}
where $J_{i,t} = \{k \in n\mid C_k = (c_1^k,\ldots,c_N^k) \in \cC \text{ and } c_j^k = c_j^i\; \forall j \neq t\}$ and $b_{i,t} = |J_{i,t}|$. Using this, we derive the following formula for $\x_{C_i}^t$:
\begin{equation*}
    \x_i^T \equiv \x_{C_i}^T = \frac{1}{b_{i,T}}\sum\limits_{i_1\in J_{i,T}}\frac{1}{b_{i_1,T-1}}\sum\limits_{i_2\in J_{i_1,T-1}}\frac{1}{b_{i_2,T-2}}\sum\limits_{i_3\in J_{i_2,T-1}}\ldots\frac{1}{b_{i_{T-1},1}}\sum\limits_{i_T\in J_{i_{T-1},1}}\x_{i_{T}}.
\end{equation*}
Taking the expectation w.r.t. $\x_1,\ldots,\x_n$, we get
\begin{equation*}
    \EE_{\x}\left[\x_i^T\right] = \frac{1}{b_{i,T}}\sum\limits_{i_1\in J_{i,T}}\frac{1}{b_{i_1,T-1}}\sum\limits_{i_2\in J_{i_1,T-1}}\frac{1}{b_{i_2,T-2}}\sum\limits_{i_3\in J_{i_2,T-1}}\ldots\frac{1}{b_{i_{T-1},1}}\sum\limits_{i_T\in J_{i_{T-1},1}}\overline{\x}_{i_{T}}.
\end{equation*}
Using the independence of $\x_1,\ldots,\x_n$, we derive
\begin{eqnarray*}
    \EE_\x\left[\left\|\x_i^T - \EE_{\x}\left[\x_i^T\right]\right\|^2\right] &=& \EE_\x\left[\left\|\sum\limits_{i_1\in J_{i,T}}\sum\limits_{i_2\in J_{i_1,T-1}}\ldots \sum\limits_{i_{T}\in J_{i_{T-1},1}}\frac{\x_{i_T} - \overline{\x}_{i_T}}{b_{i,T} b_{i_1,T-1}\ldots b_{i_{T-1},1}}\right\|^2\right]\\
    &=& \sum\limits_{i_1\in J_{i,T}}\sum\limits_{i_2\in J_{i_1,T-1}}\ldots \sum\limits_{i_{T}\in J_{i_{T-1},1}}\frac{\EE_\x\left[\|\x_{i_T} - \overline{\x}_{i_T}\|^2\right]}{b_{i,T}^2 b_{i_1,T-1}^2\ldots b_{i_{T-1},1}^2}\\
    &\le& \sum\limits_{i_1\in J_{i,T}}\sum\limits_{i_2\in J_{i_1,T-1}}\ldots \sum\limits_{i_{T}\in J_{i_{T-1},1}}\frac{\sigma^2}{b_{i,T}^2 b_{i_1,T-1}^2\ldots b_{i_{T-1},1}^2}\\
    &=& \sum\limits_{i_1\in J_{i,T}}\sum\limits_{i_2\in J_{i_1,T-1}}\ldots \sum\limits_{i_{T-1}\in J_{i_{T-2},2}}\frac{\sigma^2}{b_{i,T}^2 b_{i_1,T-1}^2\ldots b_{i_{T-2},2}^2b_{i_{T-1},1}}.
\end{eqnarray*}
Next, taking the full expectation from the both sides of the previous inequality and using the tower property, we obtain
\begin{equation}
     \EE\left[\left\|\x_i^T - \EE_{\x}\left[\x_i^T\right]\right\|^2\right] \le \EE\left[\sum\limits_{i_1\in J_{i,T}}\sum\limits_{i_2\in J_{i_1,T-1}}\ldots \sum\limits_{i_{T-1}\in J_{i_{T-2},2}}\frac{\sigma^2}{b_{i,T}^2 b_{i_1,T-1}^2\ldots b_{i_{T-2},2}^2b_{i_{T-1},1}}\right]. \label{eq:rand_mix_thm_technical_1}
\end{equation}
Notice that $J_{i_k,T-k} \cap J_{i_{k+1},T-k-1} = \{i_{k+1}\}$ for all $k=0,\ldots,T-1$, where $i_0 = i$. Moreover, for $k_1, k_2 \in\{0,1,\ldots,T\}$, $k_1 < k_2$ either $J_{i_{k_1},T-k_1} \cap J_{i_{k_2},T-k_2} = \{k_2\}$ or $J_{i_{k_1},T-k_1} \cap J_{i_{k_2},T-k_2} = \varnothing$. The first situation is possible iff $i_{k_1} = i_{k_1+1} = \ldots i_{k_2-1}$.

Taking these observations about sets $J_{i_{k}, T-k}$ into account, we consider the sets $J_{i_k,T-k}' = J_{i_k,T-k}\setminus\{i_{k}\}$ for $k = 0, 1, \ldots, T-1$. These sets are pairwise disjoint and their cardinalities $b_{i_k,T-k}' = |J_{i_k,T-k}'|$ satisfy the following relations: $b_{i_k,T-k} = 1 + b_{i_k,T-k}' \ge \max\{1, b_{i_k,T-k}'\} =: \hat{b}_{i_k,T-k}$ for $k = 1, 2, \ldots, T-1$. Moreover, $b_{i,T}', b_{i_1,T-1}',\ldots, b_{i_{T-1},1}'$ are independent random variables from the binomial distribution $\text{Binom}(M-1, p)$. Finally, we notice that the number of terms in \eqref{eq:rand_mix_thm_technical_1} is upper-bounded by $M^{T-1}$, since $|J_{i,t}| \le M$ for all $i = 1,\ldots,n$ and $t=0,\ldots,T$.

Putting all together, we obtain
\begin{eqnarray*}
    \EE\left[\left\|\x_i^T - \EE_{\x}\left[\x_i^T\right]\right\|^2\right] &\le& \EE\left[\sum\limits_{i_1\in J_{i,T}}\sum\limits_{i_2\in J_{i_1,T-1}}\ldots \sum\limits_{i_{T-1}\in J_{i_{T-2},2}}\frac{\sigma^2}{\hat b_{i,T}^2 \hat b_{i_1,T-1}^2\ldots \hat b_{i_{T-2},2}^2\hat b_{i_{T-1},1}}\right]\\
    &\le& M^{T-1}\sigma^2\EE\left[\frac{1}{\hat\xi_{1}^2 \hat\xi_{2}^2\ldots \hat\xi_{T-1}^2\hat\xi_{T}}\right]\\
    &=& M^{T-1}\sigma^2\EE\left[\frac{1}{\hat\xi_{1}^2}\right]\EE\left[\frac{1}{\hat\xi_{2}^2}\right]\ldots \EE\left[\frac{1}{\hat\xi_{T-1}^2}\right]\EE\left[\frac{1}{\hat\xi_{T}}\right],
\end{eqnarray*}
where $\hat \xi_k^2 = \max\{1,\xi_1^2\}$ for $k=1,\ldots,T$ and $\xi_1,\ldots,\xi_T$ are i.i.d.\ random variables having the binomial distribution $\text{Binom}(M-1, p)$. Then one can simplify the inequality above using Lemma~\ref{lem:ode_lemma} and get
\begin{eqnarray*}
    \EE\left[\left\|\x_i^T - \EE_{\x}\left[\x_i^T\right]\right\|^2\right] &\le& M^{T-1}\sigma^2 m_1(M-1,p)\left(m_2(M-1,p)\right)^{T-1},
\end{eqnarray*}
where functions $m_1(M,p)$ and $m_2(M,p)$ are defined in \eqref{eq:binom_first_inverse_moment} and \eqref{eq:binom_second_inverse_moment} respectively.

Next, we simplify the obtained upper bound under the assumption that $M$ and $p$ are not too small; specifically, $M\ge 11$ and $p\ge \nicefrac{2}{3}$. From \eqref{eq:binom_first_inverse_moment}, we have
\begin{eqnarray*}
    m_1(M-1,p) &=& (1-p)^{M-1} + \sum\limits_{i=1}^{M-1}\frac{1}{i}\left((1-p)^{M-1-i} - (1-p)^{M-1}\right)\\
    &\le& (1-p)^{M-1}\sum\limits_{i=1}^{M-1}\frac{1}{i(1-p)^{i}}.
\end{eqnarray*}
Since
\begin{equation*}
    \frac{1}{(k+1)(1-p)^{k+1}}\cdot\frac{k(1-p)^k}{1} = \frac{k}{(k+1)(1-p)} \xrightarrow[k\to\infty]{}\frac{1}{1-p} \ge 3,
\end{equation*}
we have
\begin{equation*}
    (1-p)^{M-1}\sum\limits_{i=1}^{M-1}\frac{1}{i(1-p)^{i}} = \x\left((1-p)^M\cdot\frac{1}{M(1-p)^M}\right) = \x\left(\frac{1}{M}\right).
\end{equation*}
Using simple algebra, one can prove that for $M\ge 11$ and $p \ge\nicefrac{2}{3}$ the following inequality holds:
\begin{equation*}
    m_1(M-1,p)\le (1-p)^{M-1}\sum\limits_{i=1}^{M-1}\frac{1}{i(1-p)^{i}} \le \frac{2}{M}.
\end{equation*}
Similarly, we analyze $m_2(M-1, p)$:
\begin{eqnarray*}
    m_2(M-1,p) &=& (1-p)^{M-1} + \sum\limits_{i=1}^{M-1}\frac{1}{i}\left((1-p)^{M-1-i} - (1-p)^{M-1}\right)\sum\limits_{j=i}^{M-1}\frac{1}{j}\\
    &\le& (1-p)^{M-1}\sum\limits_{i=1}^{M-1}\frac{1}{i(1-p)^i}\sum\limits_{j=i}^{M-1}\frac{1}{j}.
\end{eqnarray*}
Since
\begin{eqnarray*}
    \frac{\frac{1}{k(1-p)^k}\sum\limits_{j=k}^{M-1}\frac{1}{j}}{\frac{1}{(k-1)(1-p)^{k-1}}\sum\limits_{j=k-1}^{M-1}\frac{1}{j}} = \frac{(k-1)\sum\limits_{j=k}^{M-1}\frac{1}{j}}{k(1-p)\left(\frac{1}{k-1} + \sum\limits_{j=k}^{M-1}\frac{1}{j}\right)} \ge \frac{3(k-1)\cdot\frac{1}{k}}{k\left(\frac{1}{k-1}+\frac{1}{k}\right)} = \frac{3(k-1)^2}{k(2k-1)}\xrightarrow[k\to\infty]{}  \frac{3}{2},
\end{eqnarray*}
we have
\begin{equation*}
    (1-p)^{M-1}\sum\limits_{i=1}^{M-1}\frac{1}{i(1-p)^i}\sum\limits_{j=i}^{M-1}\frac{1}{j} = \x\left((1-p)^M\cdot\frac{1}{M^2(1-p)^M}\right) = \x\left(\frac{1}{M^2}\right).
\end{equation*}
Next, one can prove with simple algebra that for $M\ge 11$ and $p \ge\nicefrac{2}{3}$ the following inequality holds:
\begin{equation*}
    m_2(M-1,p) \le (1-p)^{M-1}\sum\limits_{i=1}^{M-1}\frac{1}{i(1-p)^i}\sum\limits_{j=i}^{M-1}\frac{1}{j} \le \frac{3}{M^2}.
\end{equation*}
Plugging the obtained upper bounds for $m_1(M-1,p)$ and $m_2(M-1,p)$ in \eqref{eq:variance_bound_supp}, we obtain \eqref{eq:variance_bound_2_supp}.
\end{proof}

\section{Convergence Proofs of {\tt Moshpit SGD}}\label{sect:missing_proofs_local_sgd}
In this section, we provide the complete statements of the theorems establishing the convergence of {\tt Moshpit SGD} together with the full proofs. First, we introduce all necessary definitions, basic inequalities and auxiliary lemmas; then we prove the convergence in strongly convex and convex cases; lastly, we provide the proofs for the non-convex case.

\subsection{Convex Case}
In this section, we give the full proof of Theorem~\ref{thm:cvx_convergence} about the convergence of {\tt Moshpit SGD} for convex and strongly convex problems. The scheme of the proof follows the similar steps as in the state-of-the-art analysis of Local-SGD \cite{khaled2020tighter,woodworth2020local,gorbunov2020local}. We start with the following lemma:
\begin{lemma}\label{lem:key_lemma_cvx}
    Let $f_1 = \ldots = f_n = f$, function $f$ be $\mu$-strongly convex (Def.~\ref{def:str_cvx}) and $L$-smooth (see Def.~\ref{def:L_smoothness}), and Assumptions~\ref{as:bounded_var}~and~\ref{as:averaging_quality} hold with $\Delta_{pv}^k = \delta_{pv,1}\gamma\mu\EE[\|\x^k-\x^*\|^2] + \gamma^2\delta_{pv,2}^2$ and $\widetilde{\x} = \x^*$, where $\x^* \in \argmin_{\x\in\R^d} f(\x)$ and $\delta_{pv,1}\in [0,1)$, $\delta_{pv,2}\ge 0$. Then, for any $k \ge 0$ the iterates produced by {\tt Moshpit SGD} with $\gamma \le \nicefrac{1}{4L}$ satisfy
    \begin{eqnarray}
        \gamma\EE\left[f(\x^k) - f(\x^*)\right] &\le& (1-\gamma\mu(1-\delta_{pv,1}))\EE\left[\|\x^k - \x^*\|^2\right] - \EE\left[\|\x^{k+1} - \x^*\|^2\right]\notag\\
        &&\quad+ \frac{3L\gamma}{2}\EE[V_k] + \gamma^2\left(\frac{\sigma^2}{n_{\min}} + \delta_{pv,2}^2\right),\label{eq:key_lemma_cvx}
    \end{eqnarray}
    where $V_k = \frac{1}{n_k}\sum_{i\in P_k}\|\x_i^k - \x^k\|^2$ and $\x^k = \frac{1}{n_k}\sum_{i\in P_k}\x_i^k$.
\end{lemma}
\begin{proof}
Recall that Assumption~\ref{as:averaging_quality} with $\Delta_{pv}^k = \delta_{pv,1}\gamma\mu\EE[\|\x^k-\x^*\|^2] + \gamma^2\delta_{pv,2}^2$ and $\widetilde{\x} = \x^*$ states
\begin{equation}
    \EE\left[\langle\x^{k+1} - \widehat{\x}^{k+1}, \x^{k+1}+\widehat{\x}^{k+1} - 2\x^*\rangle\right] \le \delta_{pv,1}\gamma\mu\EE[\|\x^k-\x^*\|^2] + \gamma^2\delta_{pv,2}^2, \label{eq:key_lemma_cvx_tech_1}
\end{equation}
where $\widehat \x^{k+1} = \frac{1}{n_{k}}\sum_{i\in P_{k}}(\x_i^{k}-\gamma g_i^k)$. Next, the definition of $\widehat \x^{k+1}$ implies
\begin{equation}
    \widehat \x^{k+1} = \frac{1}{n_k}\sum\limits_{i\in P_{k}}\x_i^{k} - \frac{\gamma}{n_k}\sum\limits_{i\in P_{k}} g_i^k = \x^k - \gamma g^k,\notag
\end{equation}
where $g^k = \frac{1}{n_k}\sum_{i\in P_k}g_i^k$. Using this, we derive
\begin{eqnarray}
    \|\x^{k+1} - \x^*\|^2 &=& \|\widehat{\x}^{k+1} - \x^*\|^2 + 2\langle \x^{k+1} - \widehat{\x}^{k+1}, \widehat{\x}^{k+1} - \x^* \rangle + \|\x^{k+1} - \widehat{\x}^{k+1}\|^2\notag\\
    &=& \|\x^k - \x^* - \gamma g^k\|^2 +  \langle\x^{k+1} - \widehat{\x}^{k+1}, \x^{k+1}+\widehat{\x}^{k+1} - 2\x^*\rangle \notag\\
    &=& \|\x^k - \x^*\|^2 -2\gamma\langle\x^k - \x^*, g^k\rangle + \gamma^2\|g^k\|^2 +  \langle\x^{k+1} - \widehat{\x}^{k+1}, \x^{k+1}+\widehat{\x}^{k+1} - 2\x^*\rangle. \notag
\end{eqnarray}
Taking the conditional expectation $\EE\left[\ \cdot \mid \x^k\right] := \EE\left[\ \cdot \mid P_k, \x_i^k, i\in P_k\right]$ from the both sides of the previous equation and using Assumption~\ref{as:bounded_var}, we obtain
\begin{eqnarray}
    \EE\left[\|\x^{k+1} - \x^*\|^2\mid \x^k\right] &=& \|\x^k - \x^*\|^2 -2\gamma\left\langle\x^k - \x^*, \frac{1}{n_k}\sum\limits_{i\in P_k}\nabla f(\x_i^k)\right\rangle + \gamma^2\EE\left[\left\|\frac{1}{n_k}\sum\limits_{i\in P_k}g_i^k\right\|^2\mid \x^k\right] \notag\\
    &&\quad +  \EE\left[\langle\x^{k+1} - \widehat{\x}^{k+1}, \x^{k+1}+\widehat{\x}^{k+1} - 2\x^*\rangle\mid \x^k\right]. \label{eq:key_lemma_cvx_tech_2}
\end{eqnarray}
Next, we estimate the second and the third terms in the right-hand side of \eqref{eq:key_lemma_cvx_tech_2}. First,
\begin{eqnarray}
    -2\gamma\left\langle\x^k - \x^*, \frac{1}{n_k}\sum\limits_{i\in P_k}\nabla f(\x_i^k)\right\rangle &=& \frac{2\gamma}{n_k}\sum\limits_{i\in P_k}\left(\langle\x^* - \x_i^k, \nabla f(\x_i^k) \rangle + \langle\x_i^k - \x^k, \nabla f(\x_i^k) \rangle \right)\notag\\
    &\overset{\eqref{eq:str_cvx_def},\eqref{eq:L_smoothness_cor}}{\le}& \frac{2\gamma}{n_k}\sum\limits_{i\in P_k}\left( f(\x^*) - f(\x_i^k) - \frac{\mu}{2}\|\x_i^k - \x^*\|^2\right)\notag\\
    &&\quad + \frac{2\gamma}{n_k}\sum\limits_{i\in P_k}\left(f(\x_i^k) - f(\x^k) + \frac{L}{2}\|\x_i^k - \x^k\|^2\right)\notag\\
    &\overset{\eqref{eq:a_b_norm_squared}}{\le}& 2\gamma\left(f(\x^*) - f(\x^k)\right) -\gamma\mu\|\x^k - \x^*\|^2 \notag\\
    &&\quad + L\gamma V_k, \label{eq:key_lemma_cvx_tech_3}
\end{eqnarray}
where $V_k = \frac{1}{n_k}\sum_{i\in P_k}\|\x_i^k - \x^k\|^2$. Secondly, since stochastic gradients $\{g_i^k\}_{i\in P_k}$ are computed independently, we get
\begin{eqnarray}
    \gamma^2\EE\left[\left\|\frac{1}{n_k}\sum\limits_{i\in P_k}g_i^k\right\|^2\mid \x^k\right] &\overset{\eqref{eq:variance_decomposition}}{=}& \gamma^2\left\|\frac{1}{n_k}\sum\limits_{i\in P_k}\nabla f(\x_i^k)\right\|^2\notag\\
    &&\quad + \gamma^2\EE\left[\left\|\frac{1}{n_k}\sum\limits_{i\in P_k}(g_i^k-\nabla f(\x_i^k))\right\|^2\mid \x^k\right]\notag\\
    &\overset{\eqref{eq:a_b_norm_squared}}{\le}& 2\gamma^2 \left\|\frac{1}{n_k}\sum\limits_{i\in P_k}(\nabla f(\x_i^k)-\nabla f(\x^k))\right\|^2 + 2\gamma^2\|\nabla f(\x^k)\|^2 \notag\\
    &&\quad + \frac{\gamma^2}{n_k^2}\sum\limits_{i\in P_k}\EE\left[\|g_i^k - \nabla f(\x_i^k)\|^2\mid \x^k\right]\notag\\
    &\overset{\eqref{eq:a_b_norm_squared},\eqref{eq:L_smoothness_cor_2},\eqref{eq:bounded_variance}}{\le}& \frac{2\gamma^2}{n_k}\sum\limits_{i\in P_k}\|\nabla f(\x_i^k)-\nabla f(\x^k)\|^2\notag\\
    &&\quad+ 4L\gamma^2\left(f(\x^k) - f(\x^*)\right) + \frac{\gamma^2\sigma^2}{n_k}\notag\\
    &\overset{\eqref{eq:L_smoothness_def}}{\le}& \underbrace{\frac{2L^2\gamma^2}{n_k}\sum\limits_{i\in P_k}\|\x_i^k - \x^k\|^2}_{2L^2\gamma^2 V_k}\notag\\
    &&\quad+ 4L\gamma^2\left(f(\x^k) - f(\x^*)\right) + \frac{\gamma^2\sigma^2}{n_{\min}}. \label{eq:key_lemma_cvx_tech_4}
\end{eqnarray}
Plugging \eqref{eq:key_lemma_cvx_tech_3} and \eqref{eq:key_lemma_cvx_tech_4} in \eqref{eq:key_lemma_cvx_tech_2}, we obtain
\begin{eqnarray}
    \EE\left[\|\x^{k+1} - \x^*\|^2\mid \x^k\right] &\le& (1-\gamma\mu)\|\x^k - \x^*\|^2 - 2\gamma\left(1 - 2L\gamma\right)\left(f(\x^k) - f(\x^*)\right)\notag\\
    &&\quad  + L\gamma\left(1+2L\gamma\right)V_k + \frac{\gamma^2\sigma^2}{n_{\min}} \notag\\
    &&\quad +  \EE\left[\langle\x^{k+1} - \widehat{\x}^{k+1}, \x^{k+1}+\widehat{\x}^{k+1} - 2\x^*\rangle\mid \x^k\right], \notag
\end{eqnarray}
and
\begin{eqnarray}
    \EE\left[\|\x^{k+1} - \x^*\|^2\right] &\overset{\eqref{eq:key_lemma_cvx_tech_1}}{\le}& (1-\gamma\mu(1-\delta_{pv,1}))\EE\left[\|\x^k - \x^*\|^2\right] - 2\gamma\left(1 - 2L\gamma\right)\EE\left[f(\x^k) - f(\x^*)\right]\notag\\
    &&\quad+ L\gamma\left(1+2L\gamma\right)\EE[V_k] + \gamma^2\left(\frac{\sigma^2}{n_{\min}} + \delta_{pv,2}^2\right)\notag\\
    &\le& (1-\gamma\mu(1-\delta_{pv,1}))\EE\left[\|\x^k - \x^*\|^2\right] - \gamma\EE\left[f(\x^k) - f(\x^*)\right]\notag\\
    &&\quad+ \frac{3L\gamma}{2}\EE[V_k] + \gamma^2\left(\frac{\sigma^2}{n_{\min}} + \delta_{pv,2}^2\right),\notag
\end{eqnarray}
where in the last inequality we use $\gamma \le \nicefrac{1}{4L}$.
\end{proof}

Next, we estimate the term $\EE[V_k]$ measuring the expected dissimilarity between local iterates and their global average at iteration $k$.

\begin{lemma}\label{lem:V_k_lemma_cvx}
    Let $f_1 = \ldots = f_n = f$, function $f$ be $\mu$-strongly convex (Def.~\ref{def:str_cvx}) and $L$-smooth (see Def.~\ref{def:L_smoothness}), and Assumptions~\ref{as:bounded_var}~and~\ref{as:averaging_quality} hold with $\Delta_{pv}^k = \delta_{pv,1}\gamma\mu\EE[\|\x^k-\x^*\|^2] + \gamma^2\delta_{pv,2}^2$ and $\widetilde{\x} = \x^*$, where $\x^* \in \argmin_{\x\in\R^d} f(\x)$ and $\delta_{pv,1}\in [0,1)$, $\delta_{pv,2}\ge 0$. Then, for any $k \ge 0$ the iterates produced by {\tt Moshpit SGD} with $\gamma \le \nicefrac{1}{4L}$ satisfy
    \begin{equation}
        \EE[V_k] \le 2\gamma^2\left(4\delta_{aq}^2 + (\tau-1)\sigma^2\right), \label{eq:V_k_bound_cvx}
    \end{equation}
    where $V_k = \frac{1}{n_k}\sum_{i\in P_k}\|\x_i^k - \x^k\|^2$ and $\x^k = \frac{1}{n_k}\sum_{i\in P_k}\x_i^k$.
\end{lemma}
\begin{proof}
    First of all, if $k = a\tau$ for some integer $a\ge 0$, then \eqref{eq:V_k_bound_cvx} follows from Assumption~\ref{as:averaging_quality} (eq.~\eqref{eq:quality_of_avg}). Therefore, we consider such $k$ that $k = a\tau + t'$ for some $t'\in (0,\tau)$. Then, for any $i,j \in P_{k}$, $i\neq j$
    \begin{eqnarray*}
        \EE\left[\|\x_i^k - \x_j^k\|^2\mid \x^{k-1}\right] &=& \EE\left[\|\x_i^{k-1} - \gamma g_i^{k-1} - \x_j^{k-1} + \gamma g_{j}^{k-1}\|^2\mid \x^{k-1}\right]\\
        &\overset{\eqref{eq:variance_decomposition}}{=}& \|\x_i^{k-1} - \gamma \nabla f(\x_i^{k-1}) - \x_j^{k-1} + \gamma \nabla f(\x_j^{k-1})\|^2\\
        &&\quad +\gamma^2\EE\left[\|g_i^{k-1} - \nabla f(\x_i^{k-1}) + g_{j}^{k-1} - \nabla f(\x_j^{k-1})\|^2\mid \x^{k-1}\right].
    \end{eqnarray*}
    Using Lemma~\ref{lem:gd_contraction} and independence of $g_i^{k-1}$ and $g_j^{k-1}$ for given $\x_i^{k-1}, \x_j^{k-1}$, $i\neq j$ we derive
    \begin{eqnarray*}
        \EE\left[\|\x_i^k - \x_j^k\|^2\mid \x^{k-1}\right] &\overset{\eqref{eq:gd_contraction}}{\le}& (1-\gamma\mu)\|\x_i^{k-1} - \x_j^{k-1}\|^2 +\gamma^2\EE\left[\|g_i^{k-1} - \nabla f(\x_i^{k-1})\|^2\mid \x^{k-1}\right]\\
        &&\quad +\gamma^2\EE\left[\|g_j^{k-1} - \nabla f(\x_j^{k-1})\|^2\mid \x^{k-1}\right]\\
        &\overset{\eqref{eq:bounded_variance}}{\le}& (1-\gamma\mu)\|\x_i^{k-1} - \x_j^{k-1}\|^2 + 2\gamma^2\sigma^2,
    \end{eqnarray*}
    from which we get the following: 
    \begin{equation}
        \EE_g\left[\|\x_i^k - \x_j^k\|^2\right] \le (1-\gamma\mu)\EE_g\left[\|\x_i^{k-1} - \x_j^{k-1}\|^2\right] + 2\gamma^2\sigma^2 \le \EE_g\left[\|\x_i^{k-1} - \x_j^{k-1}\|^2\right] + 2\gamma^2\sigma^2.\notag 
    \end{equation}
    Here, $\EE_g[\cdot]$ denotes the expectation conditioned on $\{P_k\}_{k = a\tau}^{(a+1)\tau-1}$. Unrolling the recurrence, we get
    \begin{eqnarray}
        \EE_g\left[\|\x_i^k - \x_j^k\|^2\right] &\le& \EE_g\left[\|\x_i^{a\tau} - \x_j^{a\tau}\|^2\right] + 2(k-a\tau)\gamma^2\sigma^2\notag\\
        &\le& \EE_g\left[\|\x_i^{a\tau} - \x_j^{a\tau}\|^2\right] + 2(\tau-1)\gamma^2\sigma^2.\label{eq:V_k_lemma_technical_1}
    \end{eqnarray}
    Using this, we estimate $\EE_{g}[V_k]$:
    \begin{eqnarray*}
        \EE_g[V_k] &=& \frac{1}{n_k}\sum\limits_{i\in P_k}\EE_g\left[\left\|\x_i^k - \frac{1}{n_k}\sum\limits_{j\in P_k}\x_j^k\right\|^2\right] \overset{\eqref{eq:a_b_norm_squared}}{\le} \frac{1}{n_k^2}\sum\limits_{i,j \in P_k}\EE_g\left[\|\x_i^k - \x_j^k\|^2\right]\\
        &\overset{\eqref{eq:V_k_lemma_technical_1}}{\le}& \frac{1}{n_k^2}\sum\limits_{i,j \in P_k}\EE_g\left[\|\x_i^{a\tau} - \x_j^{a\tau}\|^2\right] + 2(\tau-1)\gamma^2\sigma^2 \\
        &\overset{\eqref{eq:a+b_norm_beta}}{\le}& \frac{2}{n_k^2}\sum\limits_{i,j \in P_k}\left(\EE_g\left[\|\x_i^{a\tau} - \x^{a\tau}\|^2\right] + \EE_g\left[\|\x_j^{a\tau} - \x^{a\tau}\|^2\right]\right) + 2(\tau-1)\gamma^2\sigma^2\\
        &=& \frac{4}{n_k}\sum\limits_{i\in P_k}\EE_g\left[\|\x_i^{a\tau} - \x^{a\tau}\|^2\right]+ 2(\tau-1)\gamma^2\sigma^2\\
        &\le& \frac{4}{n_{a\tau}}\cdot\frac{n_{a\tau}}{n_k}\sum\limits_{i\in P_{a\tau}}\EE_g\left[\|\x_i^{a\tau} - \x^{a\tau}\|^2\right]+ 2(\tau-1)\gamma^2\sigma^2\\
        &\le& \EE_g\left[\frac{8}{n_{a\tau}}\sum\limits_{i\in P_{a\tau}}\|\x_i^{a\tau} - \x^{a\tau}\|^2\right]+ 2(\tau-1)\gamma^2\sigma^2,
    \end{eqnarray*}
    where in the last inequality we use $2n_{(a+1)\tau} = 2|P_{(a+1)\tau}| \ge |P_{a\tau}| = n_{a\tau}$ and $|n_k|\le |n_{k-1}|$ following from Assumption~\ref{as:averaging_quality}. Finally, we take the full expectation from the previous inequality and derive
    \begin{eqnarray*}
        \EE[V_k] &\overset{\eqref{eq:tower_property}}{\le}& 8\EE\left[\frac{1}{n_{a\tau}}\sum\limits_{i\in P_{a\tau}}\|\x_i^{a\tau} - \x^{a\tau}\|^2\right]+ 2(\tau-1)\gamma^2\sigma^2 \overset{\eqref{eq:quality_of_avg}}{\le} 2\gamma^2\left(4\delta_{aq}^2 + (\tau-1)\sigma^2\right),
    \end{eqnarray*}
    which finishes the proof.
\end{proof}

Combining Lemmas~\ref{lem:key_lemma_cvx}~and~\ref{lem:V_k_lemma_cvx}, we get the following result:
\begin{theorem}[Theorem~\ref{thm:cvx_convergence}, convergence in the convex case]\label{thm:cvx_convergence_supp}
    Let $f_1 = \ldots = f_n = f$ be $\mu$-strongly convex (Def.~\ref{def:str_cvx}) and $L$-smooth (see Def.~\ref{def:L_smoothness}), and Assumptions~\ref{as:bounded_var}~and~\ref{as:averaging_quality} hold with $\Delta_{pv}^k = \delta_{pv,1}\gamma\mu\EE[\|\x^k-\x^*\|^2] + \gamma^2\delta_{pv,2}^2$ and $\widetilde{\x} = \x^*$, where $\x^* \in \argmin_{\x\in\R^d} f(\x)$ and $\delta_{pv,1}\in [0,1)$, $\delta_{pv,2}\ge 0$. Then, for any $K \ge 0$, the iterates produced by {\tt Moshpit SGD} with $\gamma \le \nicefrac{1}{4L}$ satisfy
    \begin{equation}
        \EE\left[f(\overline{\x}^K) - f(\x^*)\right] \le (1-\gamma\mu(1-\delta_{pv,1}))^K\frac{R_0^2}{\gamma} + \gamma\left(\frac{\sigma^2}{n_{\min}} + \delta_{pv,2}^2 + 3L\gamma\left(4\delta_{aq}^2 + (\tau-1)\sigma^2\right)\right), \label{eq:str_cvx_bound_supp}
    \end{equation}
    when $\mu > 0$, and
    \begin{equation}
        \EE\left[f(\overline{\x}^K) - f(\x^*)\right] \le \frac{R_0^2}{\gamma K} + \gamma\left(\frac{\sigma^2}{n_{\min}} + \delta_{pv,2}^2 + 3L\gamma\left(4\delta_{aq}^2 + (\tau-1)\sigma^2\right)\right), \label{eq:cvx_bound_supp}
    \end{equation}
    when $\mu = 0$, where $R_0 = \|\x^0 - \x^*\|$, $\overline{\x}^K = \frac{1}{W_K}\sum_{k=0}^Kw_k\x^k = \frac{1}{W_K}\sum_{k=0}^K\frac{w_k}{n_k}\sum_{i\in P_k}\x_i^k$, $w_k = (1-\gamma\mu(1-\delta_{pv,1}))^{-(k+1)}$, and $W_K = \sum_{k=0}^Kw_k$. That is, {\tt Moshpit SGD} achieves $\EE[f(\overline{\x}^K) - f(\x^*)] \le \varepsilon$ after 
    \begin{equation}
        K = \widetilde{\cO}\left(\frac{L}{(1-\delta_{pv,1})\mu} +  \frac{\sigma^2}{n_{\min}(1-\delta_{pv,1})\mu\varepsilon} + \frac{\delta_{pv,2}^2}{(1-\delta_{pv,1})\mu\varepsilon} + \sqrt{\frac{L((\tau-1)\sigma^2+\delta_{aq}^2)}{(1-\delta_{pv,1})^2\mu^2\varepsilon}}\right)\label{eq:str_cvx_bound_2_supp}
    \end{equation}
    iterations with
    \begin{equation*}
        \gamma = \min\left\{\frac{1}{4L}, \frac{\ln\left(\max\left\{2, \min\left\{\frac{R_0^2\mu^2(1-\delta_{pv,1})^2K^2}{(\delta_{pv,2}^2 + \nicefrac{\sigma^2}{n_{\min}}) },\frac{R_0^2\mu^3(1-\delta_{pv,1})^3K^3}{3L\left(4\delta_{aq}^2 + (\tau-1)\sigma^2\right)}\right\}\right\}\right)}{(1-\delta_{pv,1})\mu K}\right\}
    \end{equation*}
    when $\mu > 0$, and after
    \begin{equation}
        K = \cO\left(\frac{LR_0^2}{\varepsilon} +  \frac{R_0^2\sigma^2}{n_{\min}\varepsilon^2} + \frac{R_0^2\delta_{pv,2}^2}{\varepsilon^2} + \frac{R_0^2\sqrt{L((\tau-1)\sigma^2+\delta_{aq}^2)}}{\varepsilon^{\nicefrac{3}{2}}}\right)\label{eq:cvx_bound_2_supp}
    \end{equation}
    iterations with
    \begin{equation*}
       \gamma = \min\left\{\frac{1}{4L} \sqrt{\frac{R_0}{(\delta_{pv,2}^2 + \nicefrac{\sigma^2}{n_{\min}})K}}, \sqrt[3]{\frac{R_0^2}{3L\left(4\delta_{aq}^2 + (\tau-1)\sigma^2\right) K}}\right\}
    \end{equation*}
    when $\mu = 0$.
\end{theorem}
\begin{proof}
    Plugging the result of Lemma~\ref{lem:V_k_lemma_cvx} in inequality \eqref{eq:key_lemma_cvx} from Lemma~\ref{lem:key_lemma_cvx}, we obtain
    \begin{eqnarray}
        \gamma\EE\left[f(\x^k) - f(\x^*)\right] &\le& (1-\gamma\mu(1-\delta_{pv,1}))\EE\left[\|\x^k - \x^*\|^2\right] - \EE\left[\|\x^{k+1} - \x^*\|^2\right]\notag\\
        &&\quad+ 3L\gamma^3\left(4\delta_{aq}^2 + (\tau-1)\sigma^2\right) + \gamma^2\left(\frac{\sigma^2}{n_{\min}} + \delta_{pv,2}^2\right).\notag
    \end{eqnarray}
    Next, we sum up these inequalities for $k=0,\ldots, K$ with weights $w_k = (1-\gamma\mu(1-\delta_{pv,1}))^{-(k+1)}$ and divide both sides by $\gamma W_K$, where $W_K = \sum_{k=0}^Kw_k$:
    \begin{eqnarray*}
        \frac{1}{W_K}\sum\limits_{k=0}^K w_k\EE\left[f(\x^k) - f(\x^*)\right] &\le& \frac{1}{\gamma W_K}\sum\limits_{k=0}^K(1-\gamma\mu(1-\delta_{pv,1}))w_k\EE\left[\|\x^k - \x^*\|^2\right]\notag\\
        &&\quad - \frac{1}{\gamma W_K}\sum\limits_{k=0}^Kw_k\EE\left[\|\x^{k+1} - \x^*\|^2\right]\notag\\
        &&\quad+ \gamma\left(\frac{\sigma^2}{n_{\min}} + \delta_{pv,2}^2 + 3L\gamma\left(4\delta_{aq}^2 + (\tau-1)\sigma^2\right)\right)\frac{1}{W_K}\sum\limits_{k=0}^Kw_k\\
        &=& \frac{1}{\gamma W_K}\sum\limits_{k=0}^K\left(w_{k-1}\EE\left[\|\x^k - \x^*\|^2\right] - w_k\EE\left[\|\x^{k+1} - \x^*\|^2\right]\right)\notag\\
        &&\quad+ \gamma\left(\frac{\sigma^2}{n_{\min}} + \delta_{pv,2}^2 + 3L\gamma\left(4\delta_{aq}^2 + (\tau-1)\sigma^2\right)\right)\\
        &=& \frac{w_{-1}\|\x^0 - \x^*\|^2 - w_K\EE\left[\|\x^{K+1}-\x^*\|^2\right]}{\gamma W_K}\\
        &&\quad+ \gamma\left(\frac{\sigma^2}{n_{\min}} + \delta_{pv,2}^2 + 3L\gamma\left(4\delta_{aq}^2 + (\tau-1)\sigma^2\right)\right)\\
        &\le& \frac{\|\x^0 - \x^*\|^2}{\gamma W_K} + \gamma\left(\frac{\sigma^2}{n_{\min}} + \delta_{pv,2}^2 + 3L\gamma\left(4\delta_{aq}^2 + (\tau-1)\sigma^2\right)\right).
    \end{eqnarray*}
    Since $f$ is convex, we apply the Jensen's inquality
    \begin{eqnarray*}
        f\left(\frac{1}{W_K}\sum\limits_{k=0}^K w_k\x^k\right) &\le& \frac{1}{W_K}\sum\limits_{k=0}^K w_k f(\x^k)
    \end{eqnarray*}
    to the previous result and get
    \begin{eqnarray*}
        \EE\left[f(\overline{\x}^K) - f(\x^*)\right] &\le& \frac{R_0^2}{\gamma W_K} + \gamma\left(\frac{\sigma^2}{n_{\min}} + \delta_{pv,2}^2 + 3L\gamma\left(4\delta_{aq}^2 + (\tau-1)\sigma^2\right)\right),
    \end{eqnarray*}
    where $R_0 = \|\x^0 - \x^*\|$ and $\overline{\x}^K = \frac{1}{W_K}\sum_{k=0}^Kw_k\x^k = \frac{1}{W_K}\sum_{k=0}^K\frac{w_k}{n_k}\sum_{i\in P_k}\x_i^k$. If $\mu > 0$, then $W_K \ge w_K \ge (1-\gamma\mu(1-\delta_{pv,1}))^{-K}$, implying \eqref{eq:str_cvx_bound_supp}. Next, $w_k = 1$ and $W_K = K$ when $\mu = 0$ gives \eqref{eq:cvx_bound_supp}. It remains to estimate the total number of iterations $K$ required by {\tt Moshpit SGD} to find an $\varepsilon$-solution, i.e., to achieve $\EE[f(\overline{\x}^K) - f(\x^*)] \le \varepsilon$. Applying Lemma~\ref{lem:lemma_i_2_gorbunov} to \eqref{eq:str_cvx_bound_supp}, we get the following result: if $\mu > 0$ and 
    \begin{equation*}
        \gamma = \min\left\{\frac{1}{4L}, \frac{\ln\left(\max\left\{2, \min\left\{\frac{R_0^2\mu^2(1-\delta_{pv,1})^2K^2}{\delta_{pv,2}^2 + \nicefrac{\sigma^2}{n_{\min}} },\frac{R_0^2\mu^3(1-\delta_{pv,1})^3K^3}{3L\left(4\delta_{aq}^2 + (\tau-1)\sigma^2\right)}\right\}\right\}\right)}{(1-\delta_{pv,1})\mu K}\right\},
    \end{equation*}
    then
    \begin{equation*}
        \EE\left[f(\overline{\x}^K) - f(\x^*)\right] = \widetilde{\cO}\left(LR_0^2\exp\left(-\frac{\mu}{L}(1-\delta_{pv,1})K\right) + \frac{\delta_{pv,2}^2 + \nicefrac{\sigma^2}{n_{\min}}}{(1-\delta_{pv,1})\mu K} + \frac{L\left(\delta_{aq}^2 + (\tau-1)\sigma^2\right)}{(1-\delta_{pv,1})^2\mu^2 K^2}\right),
    \end{equation*}
    implying \eqref{eq:str_cvx_bound_2_supp}. Similarly, we apply Lemma~\ref{lem:lemma_technical_cvx} to \eqref{eq:cvx_bound_supp} and get that for $\mu = 0$ and 
    \begin{equation*}
        \gamma = \min\left\{\frac{1}{4L} \sqrt{\frac{R_0}{(\delta_{pv,2}^2 + \nicefrac{\sigma^2}{n_{\min}})K}}, \sqrt[3]{\frac{R_0^2}{3L\left(4\delta_{aq}^2 + (\tau-1)\sigma^2\right) K}}\right\},
    \end{equation*}
    \begin{equation*}
        \EE\left[f(\overline{\x}^K) - f(\x^*)\right] = \cO\left(\frac{LR_0^2}{K} + \sqrt{\frac{R_0^2(\delta_{pv,2}^2 + \nicefrac{\sigma^2}{n_{\min}})}{K}} + \frac{\sqrt[3]{R_0^4L\left(\delta_{aq}^2 + (\tau-1)\sigma^2\right)}}{K^{\nicefrac{2}{3}}}\right),
    \end{equation*}
    implying \eqref{eq:cvx_bound_2_supp}.
\end{proof}

\subsection{Non-Convex Case}
In this section, we give the full proof of Theorem~\ref{thm:non_cvx_convergence} about convergence of {\tt Moshpit SGD} for general non-convex problems. The proof follows the similar steps as in the state-of-the-art analysis of Local-SGD in non-convex case~\cite{li2019communication,koloskova2020unified}. We start with the following lemma:
\begin{lemma}\label{lem:key_lemma_non_cvx}
    Let $f_1 = \ldots = f_n = f$, function $f$ be $L$-smooth and bounded from below by $f_*$, and Assumptions~\ref{as:bounded_var}~and~\ref{as:averaging_quality} hold with $\Delta_{pv}^k = \delta_{pv,1}\gamma\EE[\|\nabla f(\x^k)\|^2] + L\gamma^2\delta_{pv,2}^2$, $\delta_{pv,1}\in [0,\nicefrac{1}{2})$, $\delta_{pv,2}\ge 0$. Then, for any $K \ge 0$ the iterates produced by {\tt Moshpit SGD} with $\gamma \le \nicefrac{(1-2\delta_{pv,1})}{8L}$ satisfy
    \begin{eqnarray}
         \frac{(1-2\delta_{pv,1})\gamma}{4}\sum\limits_{k=0}^{K-1}\EE\left[\|\nabla f(\x^k)\|^2\right] &\le& f(\x^0) - f_* + \gamma L^2\sum\limits_{k=0}^{K-1} \EE[V_k] \notag\\
         &&\quad + KL\gamma^2\left(\frac{\sigma^2}{n_{\min}} + \delta_{pv,2}^2\right),\label{eq:key_lemma_non_cvx}
    \end{eqnarray}
    where $V_k = \frac{1}{n_k}\sum_{i\in P_k}\|\x_i^k - \x^k\|^2$ and $\x^k = \frac{1}{n_k}\sum_{i\in P_k}\x_i^k$.
\end{lemma}
\begin{proof}
    Recall that Assumption~\ref{as:averaging_quality} with $\Delta_{pv}^k = \delta_{pv,1}\gamma\EE[\|\nabla f(\x^k)\|^2] + L\gamma^2\delta_{pv,2}^2$ states
\begin{equation}
    \EE\left[\langle\nabla f(\x^k), \x^{k+1}-\widehat{\x}^{k+1}\rangle + L\|\widehat{\x}^{k+1} - \x^{k+1}\|^2\right] \le \delta_{pv,1}\gamma\EE[\|\nabla f(\x^k)\|^2] + L\gamma^2\delta_{pv,2}^2, \label{eq:key_lemma_non_cvx_tech_1}
\end{equation}
where $\widehat \x^{k+1} = \frac{1}{n_{k}}\sum_{i\in P_{k}}(\x_i^{k}-\gamma g_i^k)$. As for the convex case, we notice that the definition of $\widehat \x^{k+1}$ implies
\begin{equation}
    \widehat \x^{k+1} = \frac{1}{n_k}\sum\limits_{i\in P_{k}}\x_i^{k} - \frac{\gamma}{n_k}\sum\limits_{i\in P_{k}} g_i^k = \x^k - \gamma g^k,\notag
\end{equation}
where $g^k = \frac{1}{n_k}\sum_{i\in P_k}g_i^k$. Using this and L-smoothness of $f$, we derive
    \begin{eqnarray*}
        f(\x^{k+1}) - f(\x^k) &\overset{\eqref{eq:L_smoothness_cor}}{\le}& \langle\nabla f(\x^k), \x^{k+1} - \x^k \rangle + \frac{L}{2}\|\x^{k+1} - \x^k\|^2\\
        &\overset{\eqref{eq:a+b_norm_beta}}{\le}& \langle\nabla f(\x^k), \widehat{\x}^{k+1} - \x^k \rangle + \langle\nabla f(\x^k), \x^{k+1} - \widehat{\x}^{k+1} \rangle + L\|\widehat{\x}^{k+1} - \x^k\|^2\\
        &&\quad + L\|\x^{k+1} - \widehat{\x}^{k+1}\|^2\\
        &=& - \gamma\langle\nabla f(\x^k), g^k\rangle + L\gamma^2\|g^k\|^2 + \langle\nabla f(\x^k), \x^{k+1} - \widehat{\x}^{k+1} \rangle\\
        &&\quad + L\|\x^{k+1} - \widehat{\x}^{k+1}\|^2,
    \end{eqnarray*}
    from which it follows that
    \begin{eqnarray}
        \EE\left[f(\x^{k+1}) - f(\x^k)\mid \x^k\right] &\le& -\gamma\left\langle\nabla f(\x^k), \frac{1}{n_k}\sum\limits_{i\in P_k}\nabla f(\x_i^k) \right\rangle + L\gamma^2\EE\left[\left\|\frac{1}{n_k}\sum\limits_{i\in P_k}g_i^k\right\|^2\mid \x^k\right]\notag\\
        &&\quad + \EE\left[\langle\nabla f(\x^k), \x^{k+1} - \widehat{\x}^{k+1} \rangle + L\|\x^{k+1} - \widehat{\x}^{k+1}\|^2\mid \x^k\right],\label{eq:key_lemma_non_cvx_tech_2}
    \end{eqnarray}
    where $\EE\left[\ \cdot \mid \x^k\right] := \EE\left[\ \cdot \mid P_k, \x_i^k, i\in P_k\right]$. Next, we estimate the second and third terms in the right-hand side of \eqref{eq:key_lemma_non_cvx_tech_2}. First of all,
\begin{eqnarray}
    -\gamma\left\langle\nabla f(\x^k), \frac{1}{n_k}\sum\limits_{i\in P_k}\nabla f(\x_i^k)\right\rangle &=& -\gamma\|\nabla f(\x^k)\|^2 - \gamma\left\langle\nabla f(\x^k), \frac{1}{n_k}\sum\limits_{i\in P_k}\nabla f(\x_i^k) - \nabla f(\x^k)\right\rangle \notag\\
    &\overset{\eqref{eq:fenchel_young}}{\le}& -\gamma\|\nabla f(\x^k)\|^2 + \frac{\gamma}{2}\|\nabla f(\x^k)\|^2\notag\\
    &&\quad + \frac{\gamma}{2}\left\|\frac{1}{n_k}\sum\limits_{i\in P_k}(\nabla f(\x_i^k) - \nabla f(\x^k))\right\|^2\notag\\
    &\overset{\eqref{eq:a_b_norm_squared}}{\le}& - \frac{\gamma}{2}\|\nabla f(\x^k)\|^2 + \frac{\gamma}{2n_k}\sum\limits_{i\in P_k}\|\nabla f(\x_i^k) - \nabla f(\x^k)\|^2\notag\\
    &\overset{\eqref{eq:L_smoothness_def}}{\le}& - \frac{\gamma}{2}\|\nabla f(\x^k)\|^2 + \frac{\gamma L^2}{2}V_k, \label{eq:key_lemma_non_cvx_tech_3}
\end{eqnarray}
where $V_k = \frac{1}{n_k}\sum_{i\in P_k}\|\x_i^k - \x^k\|^2$. Secondly, since the stochastic gradients $\{g_i^k\}_{i\in P_k}$ are computed independently, we derive
\begin{eqnarray}
    L\gamma^2\EE\left[\left\|\frac{1}{n_k}\sum\limits_{i\in P_k}g_i^k\right\|^2\mid \x^k\right] &\overset{\eqref{eq:variance_decomposition}}{=}& L\gamma^2\left\|\frac{1}{n_k}\sum\limits_{i\in P_k}\nabla f(\x_i^k)\right\|^2\notag\\
    &&\quad + L\gamma^2\EE\left[\left\|\frac{1}{n_k}\sum\limits_{i\in P_k}(g_i^k-\nabla f(\x_i^k))\right\|^2\mid \x^k\right]\notag\\
    &\overset{\eqref{eq:a_b_norm_squared}}{\le}& 2L\gamma^2 \left\|\frac{1}{n_k}\sum\limits_{i\in P_k}(\nabla f(\x_i^k)-\nabla f(\x^k))\right\|^2 \notag\\
    &&\quad + 2L\gamma^2\|\nabla f(\x^k)\|^2 + \frac{\gamma^2L}{n_k^2}\sum\limits_{i\in P_k}\EE\left[\|g_i^k - \nabla f(\x_i^k)\|^2\mid \x^k\right]\notag\\
    &\overset{\eqref{eq:a_b_norm_squared},\eqref{eq:bounded_variance}}{\le}& \frac{2\gamma^2L}{n_k}\sum\limits_{i\in P_k}\|\nabla f(\x_i^k)-\nabla f(\x^k)\|^2 + 2L\gamma^2\|\nabla f(\x^k)\|^2\notag\\
    &&\quad + \frac{\gamma^2L\sigma^2}{n_k}\notag\\
    &\overset{\eqref{eq:L_smoothness_def}}{\le}& \underbrace{\frac{2L^3\gamma^2}{n_k}\sum\limits_{i\in P_k}\|\x_i^k - \x^k\|^2}_{2L^3\gamma^2 V_k} + 2L\gamma^2\|\nabla f(\x^k)\|^2\notag\\
    &&\quad + \frac{\gamma^2L\sigma^2}{n_{\min}}. \label{eq:key_lemma_non_cvx_tech_4}
\end{eqnarray}
Plugging \eqref{eq:key_lemma_non_cvx_tech_3} and \eqref{eq:key_lemma_non_cvx_tech_4} in \eqref{eq:key_lemma_non_cvx_tech_2}, we obtain
\begin{eqnarray}
    \EE\left[f(\x^{k+1}) - f(\x^k)\mid \x^k\right] &\le& -\frac{\gamma}{2}\left(1 - 4L\gamma\right)\|\nabla f(\x^k)\|^2 + \frac{\gamma L^2}{2}\left(1 + 4L\gamma\right)V_k + \frac{L\gamma^2\sigma^2}{n_{\min}}\notag\\
    &&\quad + \EE\left[\langle\nabla f(\x^k), \x^{k+1} - \widehat{\x}^{k+1} \rangle + L\|\x^{k+1} - \widehat{\x}^{k+1}\|^2\mid \x^k\right].\notag
\end{eqnarray}
Next, we take the full expectation from the both sides of the above inequality, apply the tower property \eqref{eq:tower_property} and take into account that $\gamma \le \nicefrac{(1-2\delta_{pv,1})}{8L}$:
\begin{eqnarray*}
    \EE\left[f(\x^{k+1}) - f(\x^k)\right] &\le& -\frac{\gamma}{2}\left(1 - 4L\gamma\right)\EE\left[\|\nabla f(\x^k)\|^2\right] + \frac{\gamma L^2}{2}\left(1 + 4L\gamma\right)\EE[V_k] + \frac{L\gamma^2\sigma^2}{n_{\min}}\\
    &&\quad + \EE\left[\langle\nabla f(\x^k), \x^{k+1} - \widehat{\x}^{k+1} \rangle + L\|\x^{k+1} - \widehat{\x}^{k+1}\|^2\right]\\
    &\overset{\eqref{eq:key_lemma_non_cvx_tech_1}}{\le}& -\frac{\gamma}{2}\left(1 - 2\delta_{pv,1} - 4L\gamma\right)\EE\left[\|\nabla f(\x^k)\|^2\right] + \frac{\gamma L^2}{2}\left(1 + 4L\gamma\right)\EE[V_k] \\
    &&\quad + L\gamma^2\left(\frac{\sigma^2}{n_{\min}} + \delta_{pv,2}^2\right)\\
    &\le& -\frac{(1-2\delta_{pv,1})\gamma}{4}\EE\left[\|\nabla f(\x^k)\|^2\right] + \gamma L^2 \EE[V_k] + L\gamma^2\left(\frac{\sigma^2}{n_{\min}} + \delta_{pv,2}^2\right).
\end{eqnarray*}
Summing up the obtained inequalities for $k = 0,\ldots, K-1$ and rearranging the terms, we derive
\begin{eqnarray*}
    \frac{(1-2\delta_{pv,1})\gamma}{4}\sum\limits_{k=0}^{K-1}\EE\left[\|\nabla f(\x^k)\|^2\right] &\le& \sum\limits_{k=0}^{K-1} \EE\left[f(\x^k) - f(\x^{k+1})\right] + \gamma L^2\sum\limits_{k=0}^{K-1} \EE[V_k] \\
    &&\quad + KL\gamma^2\left(\frac{\sigma^2}{n_{\min}} + \delta_{pv,2}^2\right)\\
    &=& f(\x^0) - \EE[f(\x^{K})] + \gamma L^2\sum\limits_{k=0}^{K-1} \EE[V_k] + KL\gamma^2\left(\frac{\sigma^2}{n_{\min}} + \delta_{pv,2}^2\right)\\
    &\le& f(\x^0) - f_* + \gamma L^2\sum\limits_{k=0}^{K-1} \EE[V_k] + KL\gamma^2\left(\frac{\sigma^2}{n_{\min}} + \delta_{pv,2}^2\right),
\end{eqnarray*}
where $f_*$ is a uniform lower bound for $f$.
\end{proof}
The next step towards completing the proof of Theorem~\ref{thm:non_cvx_convergence} gives the upper bound for $\sum_{k=0}^{K-1} \EE[V_k]$ that appeared in \eqref{eq:key_lemma_non_cvx}.

\begin{lemma}\label{lem:V_k_lemma_non_cvx}
    Let $f_1 = \ldots = f_n = f$ be $L$-smooth and bounded from below by $f_*$, and Assumptions~\ref{as:bounded_var}~and~\ref{as:averaging_quality} hold with $\Delta_{pv}^k = \delta_{pv,1}\gamma\EE[\|\nabla f(\x^k)\|^2] + L\gamma^2\delta_{pv,2}^2$, $\delta_{pv,1}\in [0,\nicefrac{1}{2})$, $\delta_{pv,2}\ge 0$. Then, for any $K \ge 0$ the iterates produced by {\tt Moshpit SGD} with $\gamma \le \nicefrac{1}{\left(4\sqrt{e}L(\tau-1)\right)}$ satisfy
    \begin{eqnarray}
        \sum\limits_{k=0}^{K-1}\EE[V_k] &\le& 8e\gamma^2(\tau-1)^2\sum\limits_{k=0}^{K-1}\EE[\|\nabla f(\x^k)\|^2] + 4\gamma^2K\left(2\delta_{aq}^2 + e(\tau-1)\sigma^2\right) ,\label{eq:V_k_lemma_non_cvx}
    \end{eqnarray}
    where $V_k = \frac{1}{n_k}\sum_{i\in P_k}\|\x_i^k - \x^k\|^2$ and $\x^k = \frac{1}{n_k}\sum_{i\in P_k}\x_i^k$.
\end{lemma}
\begin{proof}
    First of all, consider $k$ such that $k = a\tau + t'$ for some $t'\in [0,\tau)$. Let $\EE_g[\cdot]$ denote the expectation conditioned on $\{P_t\}_{t=a\tau}^{(a+1)\tau-1}$. Then
     \begin{eqnarray}
         \EE_g[V_k] &=& \frac{1}{n_k}\sum\limits_{i\in P_k}\EE_g\left[\|\x_i^k - \x^k\|^2\right] \overset{\eqref{eq:variance_decomposition}}{\le} \frac{1}{n_k}\sum\limits_{i\in P_k}\EE_g\left[\|\x_i^k - \x^{a\tau}\|^2\right] \notag\\
         &=& \frac{1}{n_k}\sum\limits_{i\in P_k}\EE_g\left[\left\|\x_i^{a\tau} - \x^{a\tau} - \gamma\sum\limits_{t=a\tau}^{k-1} g_i^t\right\|^2\right]\notag\\
         &\overset{\eqref{eq:a+b_norm_beta}}{\le}& \frac{2}{n_k} \sum\limits_{i\in P_k}\EE_g\left[\|\x_i^{a\tau} - \x^{a\tau}\|^2\right] + \frac{2\gamma^2}{n_k}\sum\limits_{i\in P_k}\EE_g\left[\left\|\sum\limits_{t=a\tau}^{k-1} g_i^t\right\|^2\right]. \label{eq:V_k_lemma_non_cvx_tech_1}
     \end{eqnarray}
     Next, we estimate the second term in the right-hand side of \eqref{eq:V_k_lemma_non_cvx_tech_1} using Lemma~\ref{lem:lemma14_stich_general}:
     \begin{eqnarray}
         \frac{2\gamma^2}{n_k}\sum\limits_{i\in P_k}\EE_g\left[\left\|\sum\limits_{t=a\tau}^{k-1} g_i^t\right\|^2\right] &\overset{\eqref{eq:lemma14_stich_general}}{\le}& \frac{2e\gamma^2(k - a\tau)}{n_k} \sum\limits_{i\in P_k} \sum\limits_{t=a\tau}^{k-1}\EE_g[\|\nabla f(\x_i^t)\|^2]\notag\\
         &&\quad + \frac{2e\gamma^2}{n_k}\sum\limits_{i\in P_k} \sum\limits_{t=a\tau}^{k-1}\EE_g[\|g_i^t - \nabla f(\x_i^t)\|^2]\notag\\
         &\overset{\eqref{eq:a+b_norm_beta},\eqref{eq:bounded_variance}}{\le}& 4e\gamma^2(\tau-1) \sum\limits_{t=a\tau}^{k-1}\EE_g[\|\nabla f(\x^t)\|^2] \notag\\
         &&\quad+ 4e\gamma^2(\tau-1) \sum\limits_{t=a\tau}^{k-1}\frac{1}{n_k}\sum\limits_{i\in P_k}\EE_g[\|\nabla f(\x_i^t) - \nabla f(\x^t)\|^2]\notag\\
         &&\quad + 2e\gamma^2 (k - a\tau)\sigma^2\notag\\
         &\overset{\eqref{eq:L_smoothness_def}}{\le}& 4e\gamma^2(\tau-1) \sum\limits_{t=a\tau}^{k-1}\EE_g[\|\nabla f(\x^t)\|^2]\notag\\
         &&\quad + 4e\gamma^2L^2(\tau-1) \sum\limits_{t=a\tau}^{k-1}\frac{n_t}{n_k}\cdot\frac{1}{n_t}\sum\limits_{i\in P_t}\EE_g[\|\x_i^t - \x^t\|^2]\notag\\
         &&\quad + 2e\gamma^2(\tau-1)\sigma^2\notag\\
         &\le& 4e\gamma^2(\tau-1) \sum\limits_{t=a\tau}^{k-1}\EE_g[\|\nabla f(\x^t)\|^2] + 8e\gamma^2L^2(\tau-1) \sum\limits_{t=a\tau}^{k-1}\EE_g[V_t]\notag\\
         &&\quad+ 2e\gamma^2(\tau-1)\sigma^2,\notag
     \end{eqnarray}
     where in the last two inequalities we use $n_k = |P_k| \le |P_{k-1}| = n_{k-1}$ for all $k\ge 1$ and $n_{a\tau} \le 2 n_{(a+1)\tau}$ for all integer $a \ge 0$. Plugging this inequality in \eqref{eq:V_k_lemma_non_cvx_tech_1} and taking the full expectation from the result, we get
     \begin{eqnarray}
         \EE[V_k] &\le& 2\EE\left[\frac{1}{n_k}\sum\limits_{i\in P_k}\|\x_i^{a\tau} - \x^{a\tau}\|^2\right] + 4e\gamma^2(\tau-1) \sum\limits_{t=a\tau}^{k-1}\EE[\|\nabla f(\x^t)\|^2] \notag\\
         &&\quad + 8e\gamma^2L^2(\tau-1) \sum\limits_{t=a\tau}^{k-1}\EE[V_t] + 2e\gamma^2(\tau-1)\sigma^2\notag\\
         &\le& 4\EE\left[\frac{1}{n_{a\tau}}\sum\limits_{i\in P_{a\tau}}\|\x_i^{a\tau} - \x^{a\tau}\|^2\right] + 4e\gamma^2(\tau-1) \sum\limits_{t=a\tau}^{k-1}\EE[\|\nabla f(\x^t)\|^2]\notag\\
         &&\quad + 8e\gamma^2L^2(\tau-1) \sum\limits_{t=a\tau}^{k-1}\EE[V_t] + 2e\gamma^2(\tau-1)\sigma^2\notag\\
         &\overset{\eqref{eq:quality_of_avg}}{\le}& 4e\gamma^2(\tau-1) \sum\limits_{t=a\tau}^{k-1}\EE[\|\nabla f(\x^t)\|^2] + 8e\gamma^2L^2(\tau-1) \sum\limits_{t=a\tau}^{k-1}\EE[V_t] + 2\gamma^2\left(2\delta_{aq}^2 + e(\tau-1)\sigma^2\right),\notag
     \end{eqnarray}
     where in the second inequality we also use $n_k = |P_k| \le |P_{k-1}| = n_{k-1}$ for all $k\ge 1$ and $n_{a\tau} \le 2 n_{(a+1)\tau}$ for all integer $a \ge 0$. Summing up the obtained inequalities for $k = a\tau, a\tau+1,\ldots, K'$ for some $K' \in[a\tau, (a+1)\tau-1]$ we derive
     \begin{eqnarray*}
         \sum\limits_{k=a\tau}^{K'}\EE[V_k] &\le& 4e\gamma^2(\tau-1)\sum\limits_{k=a\tau}^{K'} \sum\limits_{t=a\tau}^{k-1}\EE[\|\nabla f(\x^t)\|^2] + 8e\gamma^2L^2(\tau-1) \sum\limits_{k=a\tau}^{K'}\sum\limits_{t=a\tau}^{k-1}\EE[V_t]\\
         &&\quad + 2\gamma^2(K'-a\tau+1)\left(2\delta_{aq}^2 + e(\tau-1)\sigma^2\right)\\
         &\le& 4e\gamma^2(\tau-1)^2\sum\limits_{k=a\tau}^{K'} \EE[\|\nabla f(\x^k)\|^2] + 8e\gamma^2L^2(\tau-1)^2 \sum\limits_{k=a\tau}^{K'}\EE[V_k]\\
         &&\quad + 2\gamma^2(K'-a\tau+1)\left(2\delta_{aq}^2 + e(\tau-1)\sigma^2\right)\\
         &\le& 4e\gamma^2(\tau-1)^2\sum\limits_{k=a\tau}^{K'} \EE[\|\nabla f(\x^k)\|^2] + \frac{1}{2} \sum\limits_{k=a\tau}^{K'}\EE[V_k]\\
         &&\quad + 2\gamma^2(K'-a\tau+1)\left(2\delta_{aq}^2 + e(\tau-1)\sigma^2\right),
     \end{eqnarray*}
     where in the last inequality we use $\gamma \le \nicefrac{1}{\left(4\sqrt{e}L(\tau-1)\right)}$. Rearranging the terms, we get that for $K' \ge 0$
     \begin{eqnarray*}
         \sum\limits_{k=a\tau}^{K'} \EE[V_k] &\le& 8e\gamma^2(\tau-1)^2\sum\limits_{k=a\tau}^{K'}\EE[\|\nabla f(\x^k)\|^2] + 4\gamma^2(K'-a\tau+1)\left(2\delta_{aq}^2 + e(\tau-1)\sigma^2\right),
     \end{eqnarray*}
     where $a\ge 0$ is an integer such that $a\tau \le K' \le (a+1)\tau - 1$. Summing up the obtained inequalities for $K' = \tau-1, 2\tau-1,\ldots, \tau\lfloor\nicefrac{(K-1)}{\tau}\rfloor - 1, K-1$, we derive \eqref{eq:V_k_lemma_non_cvx}.
\end{proof}

Combining Lemmas~\ref{lem:key_lemma_non_cvx}~and~\ref{lem:V_k_lemma_non_cvx}, we get the following result:
\begin{theorem}[Theorem~\ref{thm:non_cvx_convergence}]
    Let $f_1 = \ldots = f_n = f$, function $f$ be $L$-smooth and bounded from below by $f_*$, and Assumptions~\ref{as:bounded_var}~and~\ref{as:averaging_quality} hold with $\Delta_{pv}^k = \delta_{pv,1}\gamma\EE[\|\nabla f(\x^k)\|^2] + L\gamma^2\delta_{pv,2}^2$, $\delta_{pv,1}\in [0,\nicefrac{1}{2})$, $\delta_{pv,2}\ge 0$. Then, for any $K \ge 0$ the iterates produced by {\tt Moshpit SGD} with
    \begin{equation*}
        \gamma \le \min\left\{\frac{1-2\delta_{pv,1}}{8L},\frac{\sqrt{1-2\delta_{pv,1}}}{8\sqrt{e}L(\tau-1)}\right\}
    \end{equation*}
    satisfy
    \begin{eqnarray}
        \EE\left[\|\nabla f(\x_{\text{rand}}^K)\|^2\right] &\le& \frac{8\Delta_0}{(1-2\delta_{pv,1})K\gamma}\notag\\
        &&\quad + \frac{8L\gamma}{1-2\delta_{pv,1}}\left(\frac{\sigma^2}{n_{\min}} + \delta_{pv,2}^2 + 4\gamma L\left(2\delta_{aq}^2 + e(\tau-1)\sigma^2\right)\right), \label{eq:non_cvx_bound_supp}
    \end{eqnarray}
    where $\Delta_0 = f(\x^0) - f_*$ and $\x_{\text{rand}}^K$ is chosen uniformly at random from $\{\x^0,\x^1,\ldots,\x^{K-1}\}$. That is, {\tt Moshpit SGD} achieves $\EE\left[\|\nabla f(\x_{\text{rand}}^K)\|^2\right] \le \varepsilon^2$ after 
    \begin{eqnarray}
        K = \cO\Bigg(\frac{L\Delta_0}{(1-2\delta_{pv,1})^2\varepsilon^2}\Bigg[1 +(\tau-1)\sqrt{1-2\delta_{pv,1}} + \frac{\delta_{pv,2}^2 + \nicefrac{\sigma^2}{n_{\min}}}{\varepsilon^2}&\notag\\
        &\hspace{-2.5cm}+ \frac{\sqrt{(1-2\delta_{pv,1})(\delta_{aq}^2+(\tau-1)\sigma^2)}}{\varepsilon}\Bigg]\Bigg)\label{eq:non_cvx_bound_2_supp}
    \end{eqnarray}
    iterations with
    \begin{equation*}
        \gamma = \min\left\{\frac{1-2\delta_{pv,1}}{8L},\frac{\sqrt{1-2\delta_{pv,1}}}{8\sqrt{e}L(\tau-1)}, \sqrt{\frac{\Delta_0}{LK\left(\delta_{pv,2}^2 + \nicefrac{\sigma^2}{n_{\min}}\right)}}, \sqrt[3]{\frac{\Delta_0}{4L^2\left(2\delta_{aq}^2 + e(\tau-1)\sigma^2\right)}}\right\}.
    \end{equation*}
\end{theorem}
\begin{proof}[Proof of Theorem~\ref{thm:non_cvx_convergence}]
    Plugging the result of Lemma~\ref{lem:V_k_lemma_non_cvx} in the inequality \eqref{eq:key_lemma_non_cvx} from Lemma~\ref{lem:key_lemma_non_cvx}, we obtain
    \begin{eqnarray*}
        \frac{(1-2\delta_{pv,1})\gamma}{4}\sum\limits_{k=0}^{K-1}\EE\left[\|\nabla f(\x^k)\|^2\right] &\le& f(\x^0) - f_* + 8e\gamma^3L^2\tau(\tau-1)\sum\limits_{k=0}^{K-1}\EE[\|\nabla f(\x^k)\|^2] \\
        &&\quad + KL\gamma^2\left(\frac{\sigma^2}{n_{\min}} + \delta_{pv,2}^2 + 4\gamma L\left(2\delta_{aq}^2 + e(\tau-1)\sigma^2\right)\right)\\
        &\le& f(\x^0) - f_* + \frac{(1-2\delta_{pv,1})\gamma}{8}\sum\limits_{k=0}^{K-1}\EE\left[\|\nabla f(\x^k)\|^2\right] \\
        &&\quad + KL\gamma^2\left(\frac{\sigma^2}{n_{\min}} + \delta_{pv,2}^2 + 4\gamma L\left(2\delta_{aq}^2 + e(\tau-1)\sigma^2\right)\right).
    \end{eqnarray*}
    Next,
    \begin{eqnarray*}
        \frac{1}{K}\sum\limits_{k=0}^K\EE\left[\|\nabla f(\x^k)\|^2\right] &\le& \frac{8\Delta_0}{(1-2\delta_{pv,1})K\gamma}\\
        &&\quad + \frac{8L\gamma}{1-2\delta_{pv,1}}\left(\frac{\sigma^2}{n_{\min}} + \delta_{pv,2}^2 + 4\gamma L\left(2\delta_{aq}^2 + e(\tau-1)\sigma^2\right)\right),
    \end{eqnarray*}
    where $\Delta_0 = f(\x^0) - f_*$. Since $\x_{\text{rand}}^K$ is chosen uniformly at random from $\{\x^0,\x^1,\ldots,\x^{K-1}\}$, we have
    \begin{equation*}
        \EE\left[\|\nabla f(\x_{\text{rand}}^K)\|^2\right] \overset{\eqref{eq:tower_property}}{=} \frac{1}{K}\sum\limits_{k=0}^K\EE\left[\|\nabla f(\x^k)\|^2\right]
    \end{equation*}
    and \eqref{eq:non_cvx_bound_supp} holds. Applying Lemma~\ref{lem:lemma_technical_cvx} to \eqref{eq:non_cvx_bound_supp}, we get the following result: if
    \begin{equation*}
        \gamma = \min\left\{\frac{1-2\delta_{pv,1}}{8L},\frac{\sqrt{1-2\delta_{pv,1}}}{8\sqrt{e}L(\tau-1)}, \sqrt{\frac{\Delta_0}{LK\left(\delta_{pv,2}^2 + \nicefrac{\sigma^2}{n_{\min}}\right)}}, \sqrt[3]{\frac{\Delta_0}{4L^2\left(2\delta_{aq}^2 + e(\tau-1)\sigma^2\right)}}\right\},
    \end{equation*}
    then $\EE\left[\|\nabla f(\x_{\text{rand}}^K)\|^2\right]$ is of the order
    \begin{equation*}
        \cO\left(\frac{L\Delta_0\left(1+ (\tau-1)\sqrt{1-2\delta_{pv,1}}\right)}{(1-2\delta_{pv,1})^2K} + \sqrt{\frac{L\Delta_0\left(\delta_{pv,2}^2 + \nicefrac{\sigma^2}{n_{\min}}\right)}{(1-2\delta_{pv,1})^2K}} + \frac{\sqrt[3]{L^2\Delta_0^2(\delta_{aq}^2 + (\tau-1)\sigma^2)}}{(1-2\delta_{pv,1})K^{\nicefrac{2}{3}}}\right),
    \end{equation*}
    which implies the desired convergence result from \eqref{eq:non_cvx_bound_2_supp}.
\end{proof}

\section{Training with a Dynamic Number of Peers}
\label{sect:load_state_from_peers}

Many practical setups with unreliable devices allow peers to join or leave at any time, which can produce undesirable side-effects. For instance, consider a participant that joins the ``swarm'' midway through the training process. If this participant starts with the initial model parameters, it can undo some of the progress made by other peers.

To circumvent this issue, we require each new participant to download the latest parameters from a random up-to-date peer discovered through DHT. The same technique is used to synchronize the optimizer statistics and the learning rate schedule. This protocol is also triggered if a peer becomes desynchronized with others, e.g., after a network freeze.

\section{Load Balancing via Linear Programming}
\label{sect:load_balancing}

When running Moshpit Averaging on heterogeneous devices, one must regularly perform Butterfly All-Reduce among peers with uneven network bandwidth.
In order to speed up the protocol, we can make low-throughput peers receive, average, and send smaller partitions of the averaged vector; conversely, the high-throughput peers can process greater fractions of the input vector.
To compute the optimal partitioning, peers must solve an optimization problem that minimizes the total time spent on communication during all-reduce.

Consider a group of $M$ peers with network bandwidths $b_1, ..., b_M$, defined for simplicity as the minimum of the upload and download speed for each peer. Our objective is to find $w_i$ --- a fraction of all input vectors to be processed by the $i$-th peer.

In Butterfly All-Reduce, each peer $i$ splits its vector into parts and sends these parts to corresponding peers. Since there is no need to send $w_i$ to itself, $i$-th peer will upload a total of $1 - w_i$ of the vector to its peers.
On the receiving side, peer $i$ will average $w_i$ of the vector from all peers in its group. To do so, it must download $M-1$ vector parts of size $w_i$ from all other peers.
After that, peers distribute the averaged parts by running the same procedure in reverse (see Figure~\ref{fig:butterfly_allreduce}).

Thus, the communication time for each peer is proportional to $t_i = (1-w_i+(M-1) w_i) \cdot \frac{1}{b_i}$ and the total runtime of Butterfly All-Reduce is the maximum communication time over all peers: $T = \max_i t_i=\max_i (1-w_i+(M-1) w_i) \cdot \frac{1}{b_i}$. Formally, we minimize $T$ with respect to $w_i$ with two constraints on the fraction weights:
\begin{alignat*}{3}
\min_w&\quad &\max_i &(1-w_i +&(M-1)w_i)\cdot\frac{1}{b_i}&\\
\text{subject to}&\quad& \sum_{i=1}^M w_i = 1&&&\\
&&w_i \geq 0 &&&\forall i=1,\ldots,M
\end{alignat*}

Because the functions being maximized and the constraints are linear in $w_i$, this problem can be reduced to linear programming~\cite{kaplan1974application}. Namely, we can minimize a surrogate variable $\xi$ such that $\forall i, \ \xi \geq (1-w_i+(M-1)\cdot w_i) \cdot \frac{1}{b_i}$. The resulting linear program is formulated as follows:

\begin{alignat*}{3}
\min_{w,\xi}&\quad& \xi && &\\
\text{subject to}&\quad& \sum_{i=1}^M w_i& = 1 &&\\
&\quad& w_i& \geq 0 &&\quad \forall i=1,\ldots,M\\
&\quad&\xi&\geq (1-&w_i+(M-1)w_i)\cdot\frac{1}{b_i}&\quad\forall i=1,\ldots,M
\end{alignat*}

We solve this problem using the interior point method~\cite{andersen} implemented as part of the SciPy package (\texttt{scipy.optimize.linprog}).
Note that depending on the conditions given by participant bandwidth, optimal weights of specific peers might be equal to 0 in some cases. In essence, this allows our method to smoothly interpolate between data parallelism~\cite{valiant1990bridging}, parameter server~\cite{parameter_server_first} and sharded parameter server~\cite{sharded_ps_first} in manner similar to BytePS~\cite{byteps}.

\section{Detailed Experimental Setup}
\label{sect:detailed_setup}

In this section, we provide the detailed hardware configuration of servers used for each of our distributed training experiments.

\subsection{ImageNet Training}\label{sect:detailed_setup_resnet}

Both homogeneous and heterogeneous training setups for ImageNet are provisioned in our on-premise infrastructure across multiple data centers and an office space (for the heterogeneous setup only).

\paragraph{Homogeneous.}For the homogeneous setup, we use 16 identical instances with the following specifications:
\begin{itemize}
    \item \textbf{GPU:} V100-PCIe,
    \item \textbf{CPU:} 6 vCPUs (Xeon E5-2650v4),
    \item \textbf{RAM:} 64GB.
\end{itemize}

\paragraph{Heterogeneous.}In turn, the heterogeneous setup contains multiple instance types listed in Table~\ref{fig:tab_setup_resnet}:
\begin{table}[h]
\centering
\caption{\textbf{Heterogeneous} setup for ImageNet training.}
\label{fig:tab_setup_resnet}
\renewcommand{\arraystretch}{1}
\begin{tabular}{@{}cccccc@{}}
\toprule
Instances & GPUs & GPU type & Cores & RAM, GB & CPU type \\ 
\midrule
4            & 1      & V100-PCIe  & 6        & 64     & E5-2650v4 \\
17           & 2      & GTX 1080Ti & 8        & 64     & E5-2650v4 \\
7            & 1      & GTX 1080Ti & 4        & 32     & E5-2650v4 \\
16           & 1      & P40  & 4        & 32     & E5-2667v2 \\
20           & 1      & M40-24GB  & 4        & 32     & E5-2667v2 \\

\bottomrule
\end{tabular}
\end{table}

\subsection{ALBERT Training}\label{sect:detailed_setup_albert}

\paragraph{Homogeneous.}For the homogeneous setup, we use a single virtual machine with the following specifications:
\begin{itemize}
    \item \textbf{GPU:} $8{\times}$ V100-PCIe,
    \item \textbf{CPU:} 48 vCPUs (Xeon E5-2650v4),
    \item \textbf{RAM:} 488GB.
\end{itemize}

At the time of writing, the cloud rent cost for this instance is \textbf{\$24.48} per hour.

\paragraph{Heterogeneous.}Our heterogeneous setup is composed of two parts: AWS EC2 Spot instances and crowdsourced machines from the \texttt{Vast.ai} marketplace. For spot instances, we picked the smallest suitable instance size available from the cloud provider and further limited their bandwidth to 1Gb/s\footnote{We use \texttt{tc qdisc} Linux utility to artificially limit the network throughput, similarly to~\cite{MLSYS2019_d09bf415}}. As for marketplace instances, we report the hardware specifications for each worker gathered 1 hour after the start of ALBERT training.

Since both cloud and marketplace instances are preemptible, the actual cost of the server fleet will vary based on the current price. For simplicity, we report the maximum hourly price we ended up paying for this instance (enforced via maximum bid). Finally, some marketplace instances have missing specifications, such as unknown CPU type. This is likely caused by non-standard virtualization configured by the device owner. The resulting fleet configuration, shown in Table~\ref{fig:tab_setup}, costs up to \$15.43/hour, depending on the number of active instances.

\begin{table*}[h!]
\centering
\caption{\textbf{Heterogeneous} setup for ALBERT training.}
\label{fig:tab_setup}
\scriptsize
\begin{tabular}{@{}ccccccc@{}}
\toprule
GPU           & Cores & RAM, GB & CPU type                       & Download, Mb/s & Upload, Mb/s &
Cost, \$/hour \\ 
\midrule
\multicolumn{7}{c}{Preemptible \texttt{g4dn.xlarge} instances ($32{\times}$)} \\
\midrule
T4            & 4         & 16     & Xeon Platinum 8259CL           & 1000          & 1000        & 0.1578         \\

\midrule
\multicolumn{7}{c}{Marketplace instances} \\    
\midrule
GTX 1070Ti    & 6         & 16     & E5-2640                        & 425           & 255         & 0.036         \\
GTX 1070Ti    & 6         & 16     & i3-6100T                       & 121           & 36          & 0.06          \\
GTX 1080Ti    & 4         & 20     & i3-6096P                       & 817           & 308         & 0.101         \\
GTX 1080Ti    & 20        & 129    & E5-2630v4                      & 660           & 475         & 0.182         \\
GTX 1080Ti    & 1         & 16     & i7-7700K                       & 245           & 210         & 0.302         \\
GTX 1080Ti    & 48        & 97     & Xeon Platinum 8124             & 583           & 539         & 0.217         \\
GTX 1080Ti    & 10        & 16     & Unknown                        & n/a           & n/a           & 0.15          \\
GTX 1080Ti    & 4         & 16     & Xeon Gold 6149                 & 98            & 100         & 0.2           \\ 
GTX 1080Ti    & 4         & 16     & Xeon Gold 6149                 & 99            & 98          & 0.2           \\ 
GTX 1080Ti    & 4         & 16     & Xeon Gold 6149                 & 99            & 99          & 0.2           \\ 
GTX 1080Ti    & 4         & 16     & Xeon Gold 6149                 & 99            & 99          & 0.2           \\ 
RTX 2070S     & 24        & 32     & E5-2620v2                      & 199           & 25          & 0.199         \\
RTX 2070S     & 32        & 97     & E5-2650                        & 162           & 64          & 0.285         \\
RTX 2080      & 6         & 16     & E5-2620v3                      & 271           & 287         & 0.25          \\
RTX 2080      & 24        & 32     & E5-2630v3                      & 199           & 25          & 0.302         \\
RTX 2080S     & 4         & 32     & E5-2697v4                      & 101           & 99          & 0.292         \\ 
RTX 2080S     & 4         & 32     & E5-2697v4                      & 93            & 99          & 0.292         \\ 
RTX 2080S     & 4         & 32     & E5-2697v4                      & 94            & 98          & 0.292         \\ 
RTX 2080S     & 4         & 32     & E5-2697v4                      & 94            & 98          & 0.292         \\ 
RTX 2080S     & 4         & 32     & E5-2697v4                      & 100           & 99          & 0.292         \\ 
RTX 2080Ti   & 4         & 16     & Ryzen Threadripper 3960x       & 279           & 271          & 0.35          \\
RTX 2080Ti   & 8         & 129    & E5-2670v3                      & 616           & 672          & 0.201         \\
RTX 2080Ti   & 6         & 32     & E5-2620v3                      & 217           & 61           & 0.22          \\
RTX 2080Ti   & 8         & 16     & E5-2697v2                      & 100           & 58           & 0.3           \\
RTX 2080Ti   & 8         & 21     & E5-2697v2                      & 145           & 49           & 0.243         \\
RTX 2080Ti    & 12        & 32     & Unknown                        & 111          & 92          & 0.326         \\
RTX 2080Ti    & 12        & 64     & E5-2690v3                      & 205          & 61          & 0.549         \\
RTX 3080      & 16        & 16     & i7-10700K                      & 69           & 49          & 0.462         \\
RTX 3090      & 14        & 32     & E5-2695v3                      & 93           & 37          & 0.498         \\
RTX 3090      & 16        & 32     & Ryzen 9 3950X                  & 338          & 38          & 0.511         \\
Titan RTX     & 4         & 32     & Xeon W-3223                   & 321           & 115          & 1             \\
Titan RTX     & 4         & 32     & Xeon Gold 6149                 & 99           & 100         & 0.702         \\ 
Titan V       & 8         & 32     & i7-7700K                       & 97           & 50          & 0.282         \\
V100-FHHL     & 8         & 60     & Xeon Gold 6148                 & 544          & 584         & 0.39          \\
\midrule
\multicolumn{6}{c}{Total hourly cost (as listed):} &\bf 15.43 \\    
\bottomrule
\end{tabular}
\end{table*}

\section{Additional Averaging Experiments}
\label{sect:extra_averaging}

\begin{figure}[h]
    \centering
    \includegraphics[width=\linewidth]{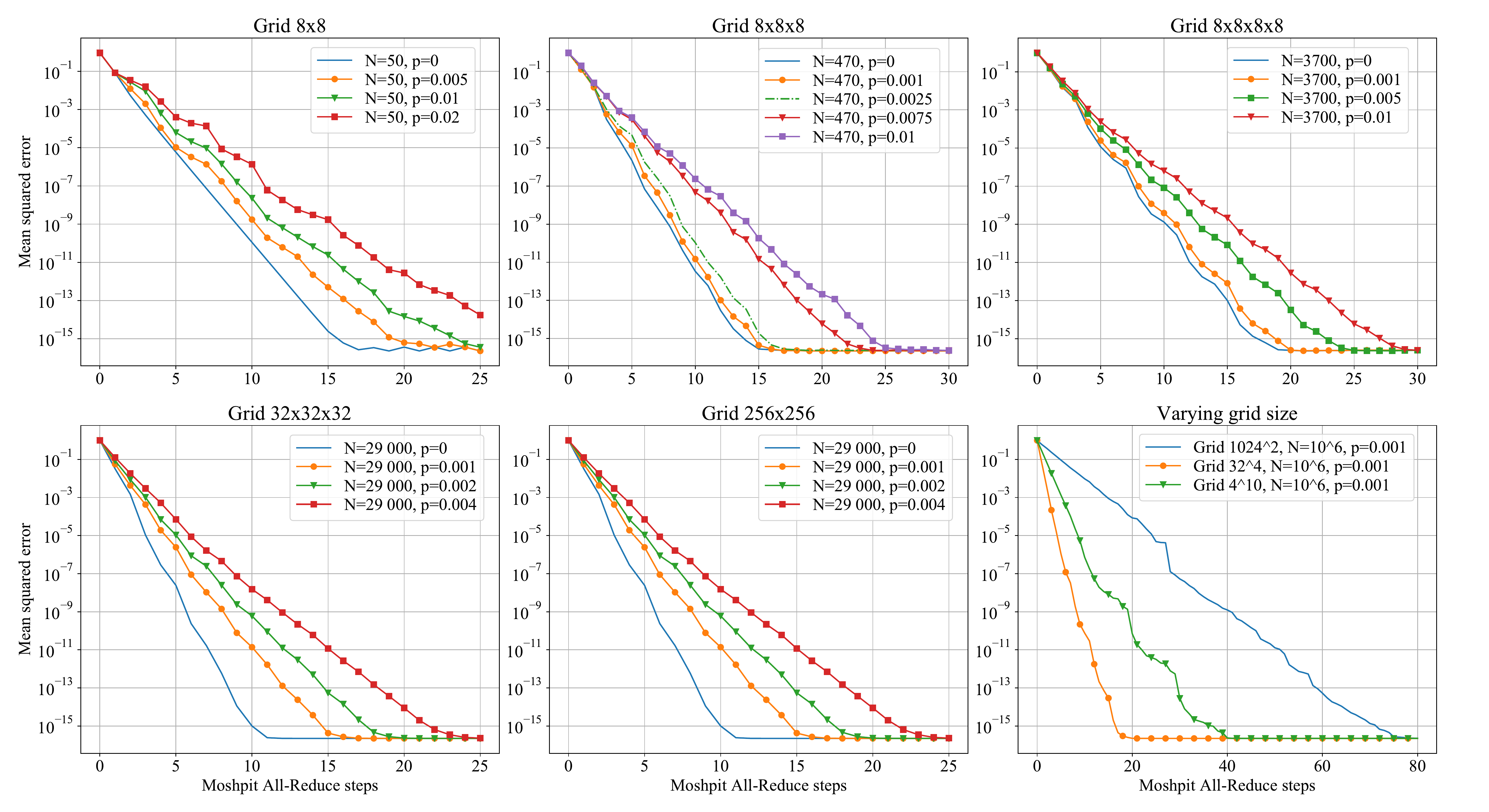}
    \vspace{-20pt}
    \caption{Averaging error of {\tt Moshpit All-Reduce} as a function of the iteration number for different configurations and failure rates.}
    \label{fig:many_averagings}
\end{figure}

In this section, we evaluate the averaging precision with the same methodology as in~\ref{sect:experiments_averaging}, but for different worker configurations. In Figure~\ref{fig:many_averagings}, plots 1--5 explore several combinations of grid sizes and failure rates, whereas plot 6 (bottom right) demonstrates a setup with the same number of peers ($10^6$) arranged into several different grid sizes and its relation to convergence. Note that $M{=}32$ outperforms the alternatives only for the specific failure rate of $0.001$.

\section{Additional Image Classification Experiments}
\label{sect:extra_classification}

Aside from the two evaluation scenarios provided in~\ref{sect:experiments_vision}, we also measure the performance of Moshpit-SGD in a non-distributed setup, i.e. on a single server with multiple GPUs. We conduct this experiment on the same $8{\times}$ V100 machine that was used in the \textbf{homogeneous} setup for training ALBERT (see Appendix~\ref{sect:detailed_setup_albert}).

\begin{figure}[h]
    \centering
    \begin{tabular}{cc}
    \hspace{-10pt}
        \includegraphics[width=0.5\textwidth]{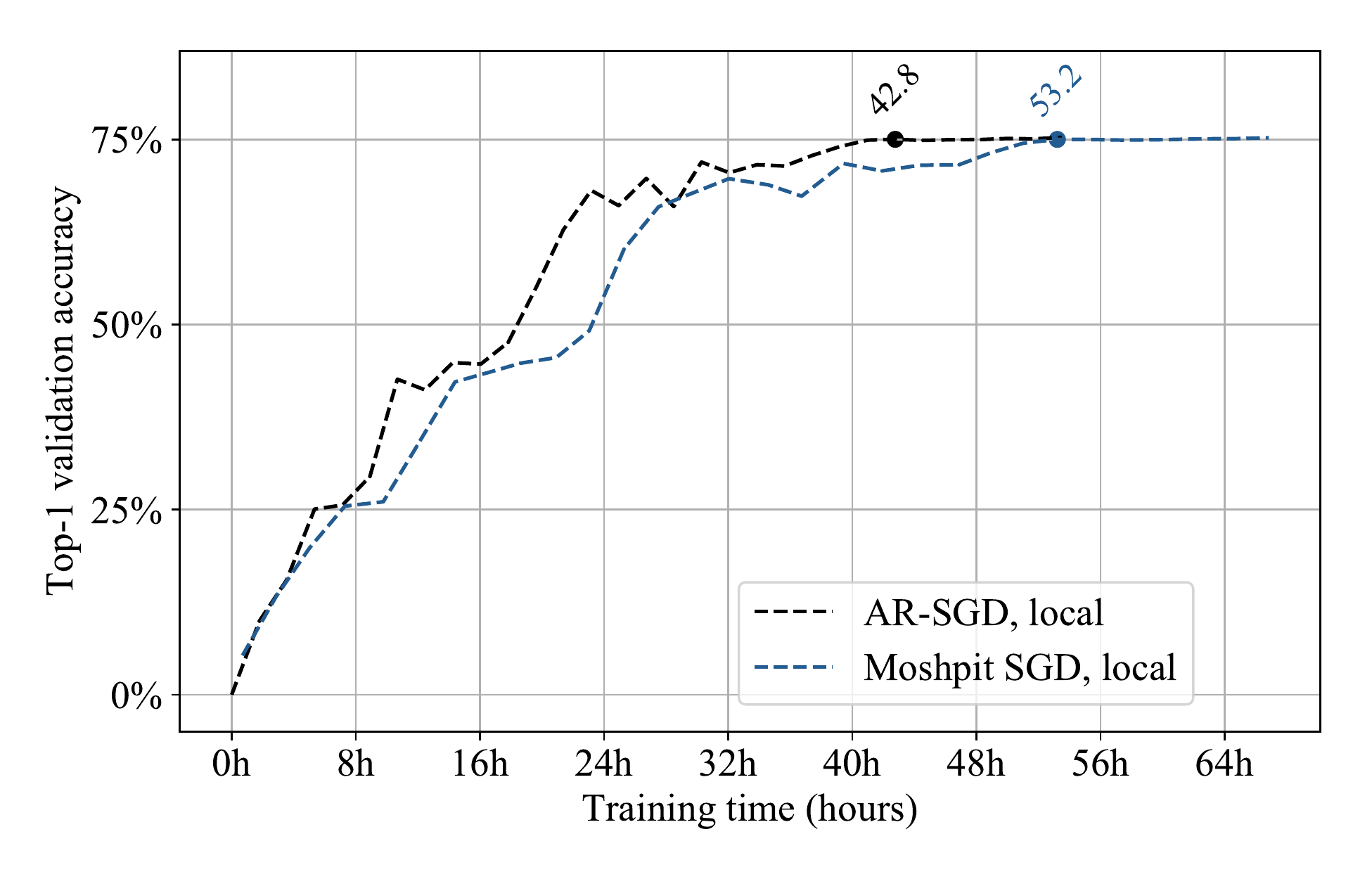} &
        \includegraphics[width=0.5\textwidth]{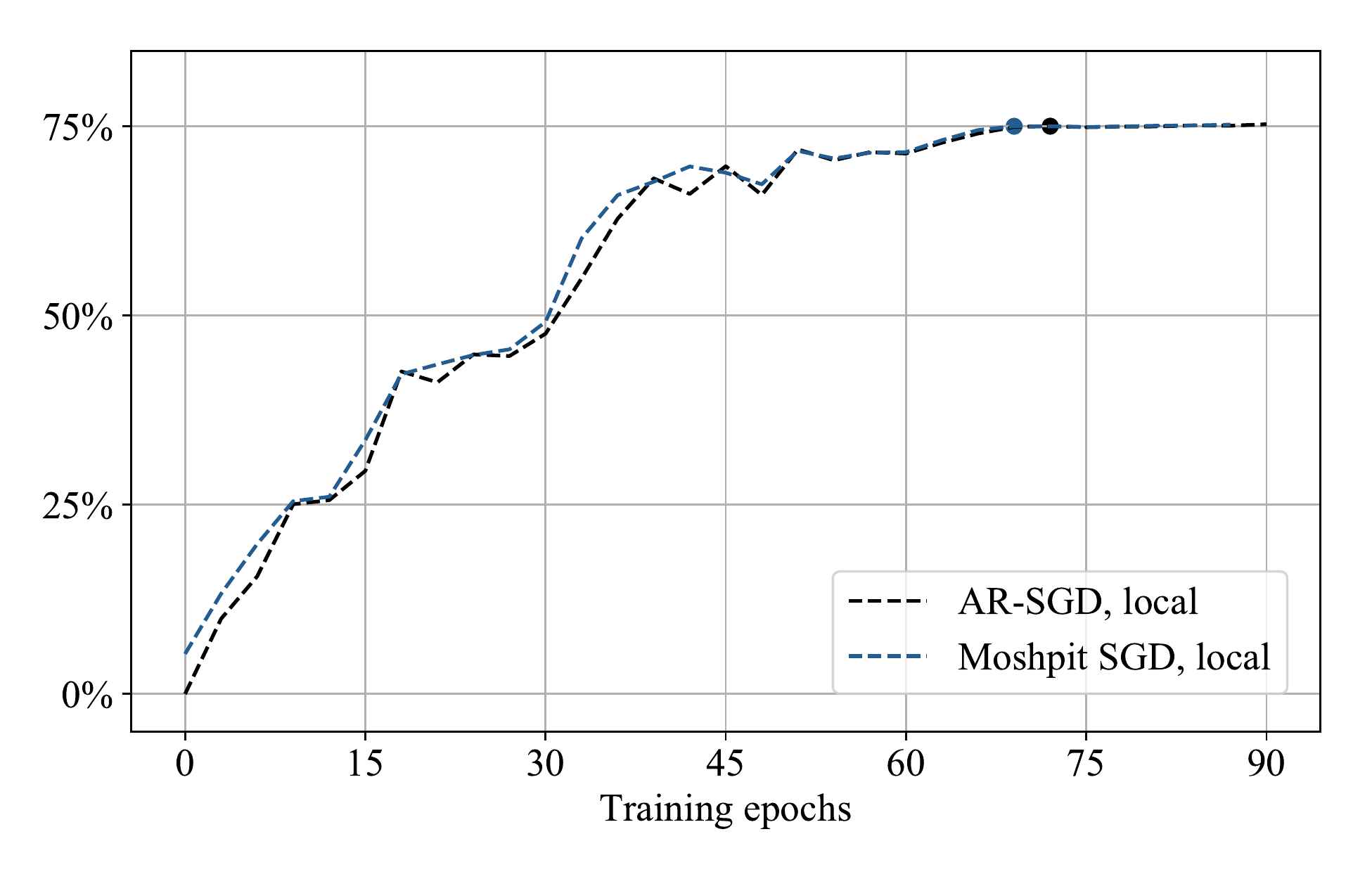}
    \end{tabular}
    \caption{
    ResNet-50 top-1 validation accuracy on ImageNet when training on a single node with $8{\times}$ V100-PCIe GPUs.
    \textbf{(Left)} Convergence in terms of training time, \textbf{(Right)} Convergence in terms of training epochs}
    \label{fig:resnet_local}\vspace{-8pt}
\end{figure}

As Figure~\ref{fig:resnet_local} demonstrates, {\tt Moshpit SGD} is slower than AR-SGD by approximately $25\%$. This result is expected, since our implementation of {\tt Moshpit All-Reduce} is more general and communicates over a TCP connection, whereas AR-SGD uses direct peer-to-peer GPU communication over PCIe. On average, this incurs a slowdown of $27\%$ in terms of training time.